\documentclass[10pt]{book}
 \usepackage[margin=1in]{geometry} 
\usepackage{amsmath,amsthm,amssymb,amsfonts,graphicx, fullpage, amsmath, url, verbatim, amssymb}
\usepackage{tikz}
\usepackage{xcolor}
\usepackage{mathtools}
\usepackage{enumitem}

\usepackage{enumitem}
\usepackage{hyperref}
\usepackage{bbm}
\usepackage{color}
\usepackage{afterpage}

\newcommand{\vphi}{\varphi}
\newcommand{\R}{\mathbb{R}} 

\newcommand{\Z}{\mathbb{Z}}

\newcommand{\ii}{\text{i}}
\renewcommand{\epsilon}{\varepsilon}

\newcommand{\lepe}{\le_{\mathtt{pe}}}
\newcommand{\lip}{\text{lip}}
\newcommand{\Lip}{\text{Lip}}
\newcommand{\ihat}{\hat{\textbf{\i}}}
\newcommand{\D}{\mathcal{D}}

\usepackage{csquotes}

\MakeOuterQuote{"}
\numberwithin{equation}{section}
\newtheorem{theorem}{Theorem}[section]
\newtheorem{lemma}[theorem]{Lemma}
\newtheorem{corollary}[theorem]{Corollary}
\newtheorem{proposition}[theorem]{Proposition}
\newtheorem{definition}[theorem]{Definition}
\newtheorem{remark}[theorem]{Remark}

\def\XXint#1#2#3{{\setbox0=\hbox{$#1{#2#3}{\int}$ }
\vcenter{\hbox{$#2#3$ }}\kern-.6\wd0}}

\definecolor{purple}{rgb}{0.65, 0, 1}
\definecolor{orange}{rgb}{1,.5,0}
\definecolor{purple}{rgb}{0.5, 0, 0.9}
\definecolor{green}{rgb}{0, 1, 0}
\definecolor{orange}{rgb}{1,.5,0}
\definecolor{gray}{rgb}{.6,.6,.6}

\newcommand{\colvec}[1]{\begin{pmatrix} #1 \end{pmatrix}}

\newenvironment{proofthm}[1]{\begin{trivlist} \item[] {\em Proof of Theorem \ref{#1}:}}{\hfill $\Box$
                       \end{trivlist}}

\newenvironment{proofprop}[1]{\begin{trivlist} \item[] {\em Proof of Proposition \ref{#1}:}}{\hfill $\Box$
                       \end{trivlist}}

\begin{document}

\author{Javier G\'omez-Serrano, Alexandru D. Ionescu, Jaemin Park}

\title{Quasiperiodic solutions of the generalized SQG equation}

\maketitle

\chapter*{Preface}
  
 The motivation for this book stems from a desire to understand the long-time behavior and structural stability of patch-type solutions in active scalar equations, particularly the generalized surface quasi-geostrophic (gSQG) equation. The gSQG model, a natural cousin of the 2D Euler equations, presents significant mathematical challenges due to its nonlocality and the supercritical nature of its dynamics when the parameter \( \alpha \in (1, 2) \).

When studying long-time or global-in-time behavior of the gSQG equation, there are several natural directions one might pursue—for example, classical questions of global well-posedness or finite-time singularity formation. Unfortunately, these questions remain largely open and notoriously difficult. A slightly more flexible line of inquiry could be constructing explicit global-in-time solutions and analyzing the typical behaviors arising from generic initial data. In this view, steady and time-periodic solutions offer natural candidates, as they are by definition global in time. Quasiperiodic solutions fall into a similar category but exhibit a far richer dynamical structure, revealing deeper aspects of the equation’s behavior.

This monograph focuses on the existence of linearly stable, quasiperiodic patch solutions of the gSQG equation for all \( \alpha \in (1, 2) \), and for a large class of initial data near the rotating disk solution. The key novelty lies in the construction of global, time-quasiperiodic solutions without the use of external parameters. Instead, we exploit the genericity of the initial data to generate rich families of such solutions.

While the techniques originate in Hamiltonian dynamical systems and perturbative analysis, their adaptation to the nonlocal, nonlinear PDE setting of gSQG patch dynamics requires substantial developments. We build upon and extend the foundational work of Berti and collaborators, combining tools from KAM theory, normal form transformations, pseudo-differential calculus, and the Nash--Moser iteration scheme.

This is a timely subject. Recent breakthroughs have clarified many aspects of linear stability and persistence in vortex dynamics, yet a comprehensive framework for quasiperiodic solutions in the gSQG setting remains largely unexplored. This book aims to fill that gap and serve as a reference for researchers interested in both the qualitative dynamics and the rigorous analysis of active scalar models.

This book is intended primarily for graduate students and researchers in partial differential equations, mathematical fluid dynamics, and dynamical systems. It may also appeal to mathematical physicists working on Hamiltonian PDEs and quasiperiodic phenomena in infinite dimensions.

We assume that the reader has a solid foundation in real and functional analysis, some familiarity with Sobolev spaces, and basic PDE theory. Prior exposure to Hamiltonian mechanics and KAM theory is helpful but not essential, as we include background material where appropriate. Throughout the book, we aim to balance rigor with intuition, and technical results with motivating ideas.

There are several excellent books and research articles on KAM theory and its applications to PDEs, including the following works:
\begin{itemize}
    \item M.~Berti and P.~Bolle, \emph{A Nash--Moser approach to KAM theory}, in \emph{Hamiltonian Partial Differential Equations and Applications}, vol.~75 of \emph{Fields Inst. Commun.}, pp.~255--284, Fields Inst. Res. Math. Sci., Toronto, ON, 2015.

    \item P.~Baldi, M.~Berti, and R.~Montalto, \emph{KAM for autonomous quasi-linear perturbations of KdV}, \emph{Ann.~Inst.~H.~Poincar\'e Anal.~Non~Lin\'eaire}, \textbf{33}(6):1589--1638, 2016.

    \item M.~Berti, Z.~Hassainia, and N.~Masmoudi, \emph{Time quasi-periodic vortex patches of Euler equation in the plane}, \emph{Invent.~Math.}, \textbf{233}(3):1279--1391, 2023.

    \item Z.~Hassainia, T.~Hmidi, and N.~Masmoudi, \emph{KAM theory for active scalar equations}, arXiv preprint \texttt{arXiv:2110.08615}, 2021.
\end{itemize}

Our approach to the gSQG equation builds on the general framework for quasi-linear Hamiltonian systems described in \emph{A Nash--Moser approach to KAM theory} by Berti and Bolle, as do many recent works in this direction. However, this book differs in that the problem posed by the gSQG equation possesses a more intricate nonlinear structure. Moreover, it targets a nonlocal, genuinely nonlinear PDE model without the use of external tuning parameters. The emphasis is on exploiting the structure of the initial data space rather than relying on auxiliary bifurcation parameters. In this sense, our approach complements existing treatments while extending the applicability of KAM theory into broader territory—namely, the theory of patch dynamics in geophysical fluid models.

\newpage

\chapter*{Acknowledgements}
JGS was partially supported by NSF through Grants DMS-1763356, DMS-2245017, DMS-2247537 and DMS-2434314; and by the AGAUR
project 2021-SGR-0087 (Catalunya). JGS and JP were partially supported by MICINN
(Spain) research grant number PID2021-125021NA-I00 and by the European Research Council through ERC-StG-852741-CAPA. AI was partially supported by NSF grants DMS-2007008 and DMS-2407694. JP was partially supported by ERC Starting Grant 676675-FLIRT, the SNF grant 212573-FLUTURA, Ambizione
fellowship project PZ00P2-216083, the Yonsei University Research Fund of 2024-22-0500, and the
POSCO Science Fellowship of POSCO TJ Park Foundation. 
Part of this work was carried out while JGS and JP were visiting the Department of Mathematics at Princeton University. We wholeheartedly thank them for their warm hospitality.
This material is based upon work supported by the National Science
Foundation under Grant No. DMS-1929284 while the authors were in residence at the Institute for Computational and Experimental Research in Mathematics in Providence, RI, during the program ``Hamiltonian Methods in Dispersive and Wave Evolution Equations''. This work is supported by the Spanish State Research Agency, through the Severo Ochoa and Mar\'ia de Maeztu Program for Centers and Units of Excellence in R\&D (CEX2020-001084-M).

\tableofcontents

\chapter{Introduction}

The {generalized surface quasi-geostrophic}\index{generalized surface quasi-geostrophic} (gSQG\index{gSQG}) equations describing the evolution of the potential temperature $\omega$ read as

\begin{align}\label{gSQG_1}
\begin{cases}
\omega_t + u\cdot \nabla \omega = 0, &\text{ for }(t,x)\in [0,\infty)\times \mathbb{R}^2,  \\
u = -\nabla^{\perp}(-\Delta)^{-1+\alpha/2}\omega, &\alpha \in [0,2]\\
\omega(0,\cdot) =\omega_0(\cdot).
\end{cases}
\end{align}
Formally, these equations interpolate between the case of the Euler equation\index{Euler equations} $(\alpha = 0)$ and the case of stationary solutions $(\alpha=2)$. The case $(\alpha = 1)$ is known as the SQG\index{SQG} equation.

The SQG equation models the evolution of the temperature from a general quasi-geostrophic system for atmospheric and oceanic flows (see  \cite{Constantin-Majda-Tabak:formation-fronts-qg,Held-Pierrehumbert-Garner-Swanson:sqg-dynamics,Pedlosky:geophysical,Majda-Bertozzi:vorticity-incompressible-flow} for more details).
The first rigorous mathematical study of the SQG equation was done by Constantin--Majda--Tabak \cite{Constantin-Majda-Tabak:formation-fronts-qg} where its mathematical importance due to its analogy with the incompressible 3D Euler equations\index{3D Euler equations} was highlighted and the first numerical and analytical study of the equation was carried out. 
C\'ordoba--Fontelos--Mancho--Rodrigo in \cite{Cordoba-Fontelos-Mancho-Rodrigo:evidence-singularities-contour-dynamics} proposed the gSQG or (SQG)$_\alpha$  model \eqref{gSQG_1} as an interpolation between the Euler and surface quasi--geostrophic equations. Nevertheless, very little is known for this family of equations, and specifically the question of global existence\index{global existence} versus finite-time singularities\index{finite-time singularities} is still open, for all $\alpha > 0$. In this monograph we aim to prove the existence of a large class of initial data for which there is time quasiperiodic behaviour and thus global existence in the more singular case $\alpha \in (1,2)$.

\section{Motivation of the problem: From  Euler to gSQG}
One of the fundamental equations in fluid mechanics is the three-dimensional Euler equations\index{Euler equations}, which are given by:

\[
\partial_t \mathbf{v} + (\mathbf{v} \cdot \nabla) \mathbf{v} = - \nabla \mathbf{p}, \quad \nabla \cdot \mathbf{v} = 0, \quad \mathbf{v}(0,x) = \mathbf{v}_0(x), \quad \text{for } (t,x) \in \mathbb{R}_+ \times \mathbb{R}^3,
\]
where \(\mathbf{v} = (v_1, v_2, v_3)\) is the velocity vector of the fluid, and \(\nabla \cdot \mathbf{v}\) denotes the divergence\index{divergence} of the velocity field. The Euler equations describe the motion of an incompressible fluid\index{incompressible} with no viscosity\index{viscosity} (inviscid flow) and constant density. Specifically, in the case of incompressible flow (constant density), these equations govern the evolution of the fluid's velocity field in three-dimensional space. The system consists of three components:

\begin{itemize}
    \item \textbf{Momentum equation\index{Momentum equation}}: \(\partial_t \mathbf{v} + (\mathbf{v} \cdot \nabla) \mathbf{v} = - \nabla \mathbf{p}\). This equation represents the conservation of momentum and describes the acceleration of the fluid due to pressure gradients.
    
    \item \textbf{Continuity equation\index{Continuity equation}}: \(\nabla \cdot \mathbf{v} = 0\). This equation ensures the incompressibility of the fluid, meaning that the volume of any fluid element remains constant over time. It guarantees that there is no net volumetric change in the fluid as it moves.
    
    \item \textbf{Initial condition}: \(\mathbf{v}(0,x) = \mathbf{v}_0(x)\). This specifies the initial velocity field of the fluid. The evolution of the fluid's velocity can be completely determined if the initial condition, together with the momentum and continuity equations, can be uniquely solved.
\end{itemize}


The Euler equations describe idealized, frictionless flows. In real fluids, viscosity introduces additional terms, leading to the more general \textit{Navier-Stokes equations\index{Navier-Stokes equations}}, which will not be covered in this monograph.
  
  \subsection{Singularity Formation vs Global Well-Posedness}
One of the key challenges in understanding the behavior of solutions to the Euler equations is the possibility of \textit{singularity formation}, where the solution or its derivatives become infinite in finite time, and the question of \textit{global well-posedness}, which concerns whether a smooth solution exists for all time. These two concepts are central to the study of fluid mechanics as they determine whether mathematical models accurately reflect physical phenomena and whether they can reliably predict the behavior of real-world fluids over extended periods.

{Singularity formation} refers to the possibility that the velocity field \(\mathbf{v}\) (or its derivatives) becomes unbounded or develops infinite values after a finite amount of time. Specifically, a singularity may manifest as a \textit{blow-up} of the solution, where the velocity becomes infinite, leading to a breakdown of the mathematical model. Despite several potential scenarios in which singularities might arise, a rigorous proof of finite-time singularity formation for the 3D Euler equations remains an open problem. This is an area of ongoing research, with the famous \textit{Navier-Stokes existence and smoothness problem} being one of the seven Millennium Prize Problems \cite{Fefferman:clay-statement}. Although the Euler equations are known to have smooth solutions in certain cases, 
a general proof of their behavior for arbitrary initial conditions remains elusive.






In contrast to singularity formation, 
in the context of the 3D Euler equations, \textit{global well-posedness} refers to the existence of a unique, smooth solution \(\mathbf{v}(t,x)\) for all times \(t \geq 0\), given an initial velocity field \(\mathbf{v}(0,x)\) that is smooth and incompressible. The question of global well-posedness remains an open problem.
Although {local well-posedness} has been established under certain conditions, which means that for sufficiently smooth initial data, solutions exist and remain smooth for a short time, the potential for singularity formation suggests that global well-posedness may not hold for arbitrary initial conditions.

 \subsection{Generalized Surface Quasi-Geostrophic Equation}
The primary mathematical model that motivates the main theorem in this monograph is the {Surface Quasi-Geostrophic (SQG)} equation. This equation is widely used in fluid dynamics, particularly in atmospheric and oceanic sciences, to describe the evolution of temperature or potential vorticity on the surface of a fluid under the influence of the Coriolis force\index{Coriolis force}. The SQG equation is given by:
\[
\partial_t \omega + \mathbf{u} \cdot \nabla \omega = 0, \quad \mathbf{u} := -\nabla^\perp (-\Delta)^{-1/2} \omega, \quad \omega(0, x) = \omega_0(x), \quad \text{for } (t,x) \in \mathbb{R}_+ \times \mathbb{R}^2,
\]
where \(\omega\) represents the temperature or potential vorticity\index{vorticity} in the quasi-geostrophic system for atmospheric flows, and \(\mathbf{u}\) is the velocity field.

The {3D Euler equations} and the (SQG) equation are both fundamental models in fluid dynamics. While the 3D Euler equations govern the motion of an incompressible, inviscid fluid in three-dimensional space, the SQG equation is a simplified model that describes the evolution of a scalar field on a surface, which is two-dimensional. Despite their differences, these two equations share several important similarities in their mathematical structures.

By formally taking the curl operator in the momentum equation of the 3D Euler equations, we can derive the vorticity equation for \(\boldsymbol{\Theta} := \nabla \times \mathbf{v}\):
\[
\partial_t \boldsymbol{\Theta} + (\mathbf{v} \cdot \nabla) \boldsymbol{\Theta} = (\boldsymbol{\Theta} \cdot \nabla) \mathbf{v}.
\]
We can list several features of the vorticity equation:
\begin{enumerate}[label=\textbf{E}-(\alph*)]
\item \(\nabla \cdot \boldsymbol{\Theta} = \nabla \cdot \mathbf{v} = 0\): This follows from the fact that \(\boldsymbol{\Theta}\) is the curl of a vector field.
\item \(\nabla \mathbf{v} = T(\boldsymbol{\Theta})\) for some singular integral operator \(T\) of order zero.
\item Conservation of kinetic energy\index{kinetic energy}: \(\| \mathbf{v}(t) \|_{L^2} = \| \mathbf{v}_0 \|_{L^2}\).
\end{enumerate}
On the other hand, the two-dimensional vector field \(\nabla^\perp \omega\) for a solution to the SQG equation satisfies:
\[
\partial_t (\nabla^\perp \omega) + (\mathbf{u} \cdot \nabla) (\nabla^\perp \omega) = (\nabla^\perp \omega \cdot \nabla) \mathbf{u},
\]
and exhibits similar features:
\begin{enumerate}[label=\textbf{S}-(\alph*)]
\item \(\nabla \cdot (\nabla^\perp \omega) = \nabla \cdot \mathbf{u} = 0\).
\item \(\nabla \mathbf{u} = S(\nabla^\perp \omega)\) for some singular integral operator \(S\) of order zero.
\item Conservation of the \(L^2\)-norm of \(\omega\): \(\| \omega \|_{L^2} = \| \omega_0 \|_{L^2}\).
\end{enumerate}
We also note that vortex lines\index{vortex line} in the Euler equation move along the flow, while the level curves of \(\omega\) in the SQG equation also move along the flow. This observation suggests that the Euler and SQG equations share many structural similarities, and their behaviors are expected to resemble each other. As in the Euler equation, the question of finite-time singularity formation versus global well-posedness for the SQG equation remains open. Consequently, the well-posedness question for the SQG equation has attracted significant attention and is an active area of research.

Despite the close relationship between the 3D Euler and SQG equations, there are also important differences. For instance, while the 3D Euler equations do not conserve \(\| \mathbf{v}(t) \|_{L^p}\) for \(p \neq 2\), the transport nature of the SQG equation immediately guarantees the conservation of \(\| \omega \|_{L^p}\) for all \(p \in [1, \infty]\). Moreover, the SQG equation is spatially two-dimensional, which simplifies certain technical computations compared to the 3D Euler equations. These features are reminiscent of the vorticity form of the {two-dimensional Euler equations}:
\[
\partial_t \theta + \mathbf{v} \cdot \nabla \theta = 0, \quad \mathbf{v} = -\nabla^\perp (-\Delta)^{-1} \theta, \quad \text{for } (t,x) \in \mathbb{R}_+ \times \mathbb{R}^2.
\]
The difference between the 2D Euler equations\index{2D Euler equations} and the SQG equation is that the velocity \(\mathbf{v}\) is related to the transported scalar \(\theta\) through an integral operator of order \(-1\) in the 2D Euler case, while it is an integral operator of order zero in the SQG equation. This distinction suggests that the SQG equation encodes some similarities and differences with the 2D Euler equation, motivating the generalized SQG equations, as described in equation~\ref{gSQG_1}, which mathematically interpolate the relationship between the scalar \(\omega\) and the velocity field \(\mathbf{u}\).
However, in contrast to the 2D Euler equation, for which global well-posedness is known, it is not yet known whether smooth solutions to the generalized SQG equations are globally well-posed or whether they may develop a singularity in finite time for any range of \(\alpha > 0\).

\subsection{Patch Problems\index{Patch Problems}}
As mentioned previously, the motivation for studying the generalized Surface Quasi-Geostrophic (gSQG) equations arises from their analogies with the Euler equations. For the 2D Euler equations, an important class of solutions is known as {vortex patch solutions}. A patch solution takes the form:
\[
\omega(t,x) := 1_{D(t)}(x), \quad \text{for a bounded domain } D(t) \subset \mathbb{R}^2.
\]
Here, \(1_{D(t)}\) denotes the characteristic function\index{characteristic function} of the domain \(D(t)\), and by a solution, we mean that \(\omega\) satisfies the 2D Euler equations in a distributional sense\index{distributional sense}:
\begin{align*}
& \int_{\mathbb{R}^2} \omega(T,x) \eta(T,x) \, dx - \int_{\mathbb{R}^2} \omega_0(x) \eta(0,x) \, dx \\
= & \int_0^T \int_{\mathbb{R}^2} \omega(t,x) \left( \partial_t \eta(t,x) + \mathbf{v}(t,x) \cdot \nabla \eta(t,x) \right) \, dx \, dt,
\end{align*}
for all smooth, compactly supported test functions \(\eta\). The fact that the solution remains a characteristic function relies  on the fact that the 2D Euler equation is a scalar transport equation\index{scalar transport equation}. This property suggests that such patch solutions can also be naturally formulated for the gSQG equations.

The question of finite-time singularity formation in the gSQG equations can thus be rephrased in the context of patch solutions as: Does the boundary \(\partial D(t)\) maintain smoothness throughout the evolution? Unfortunately, a complete answer to this question remains unresolved for any \(\alpha > 0\), except in certain cases \cite{Kiselev-Ryzhik-Yao-Zlatos:singularity-alpha-patch-boundary, Gancedo-Patel:local-existence-blowup-gsqg,Zlatos:local-regularity-singularity-gsqg-halfplane}, where the authors constructed finite-time singularities in domains with boundaries (rather than \(\mathbb{R}^2\)).

In both the smooth and patch cases, determining whether a finite-time singularity can occur in the gSQG equation remains an open and challenging problem. However, a potentially more accessible question is whether a global solution (or more generally, a large family of global solutions) can be constructed. This is the primary focus of investigation in this monograph.

To approach the construction of a global solution, it is useful to first consider a steady solution\index{steady solution}. A steady solution refers to a time-independent solution in a certain reference frame, typically under a Galilean transformation\index{Galilean transformation}. Radial functions are known to be steady solutions to the gSQG equations. For example, a patch solution with the domain \(D\) being a disk does not alter the shape of the patch during its evolution. However, one might argue that such a stationary solution is too trivial to capture the more complex features of a global solution. Therefore, the next step is to attempt the construction of a global solution that does not remain unchanged but instead exhibits mild evolution over time. To this end, we will explore the Hamiltonian structure of the gSQG equations.

 \subsection{Hamiltonian Systems and the gSQG Equations}

Let us briefly digress from our discussion of the gSQG equations and review the basic notions of a Hamiltonian system.

For a manifold\index{manifold} \( X \), a symplectic form\index{symplectic form} \( \Omega \) is a non-degenerate\index{non-degenerate}, skew-symmetric 2-form\index{skew-symmetric} on the tangent bundle\index{tangent bundle} \( TX \). This means that for every point \( p \in X \), the following properties hold:
\[
\Omega_p(V, W) = -\Omega_p(W, V) \quad \text{for all } V, W \in T_pX,
\]
and
\[
\Omega_p(V, W) = 0 \quad \text{for all } W \in T_pX \implies V = 0.
\]
Let us consider a functional \( H: X \to \mathbb{R} \cup \{ \infty \} \), which may take the value \( \infty \) at some point \( p \in X \). By the non-degeneracy of the symplectic form, there exists a unique vector field \( X_H \), called the Hamiltonian vector field\index{Hamiltonian vector field}, such that
\[
\Omega_p(X_H(p), W) = d_pH(W), \quad \text{for all } W \in T_pX,
\]
where \( d_pH \in T^*X \) denotes the gradient\index{gradient} of \( H \) at \( p \). When the Hamiltonian vector field is well-defined, the differential equation
\[
\partial_t r = X_H(r(t))
\]
is called a Hamiltonian system\index{Hamiltonian system}.

\textbf{Example 1: A Simple Harmonic Oscillator\index{Harmonic Oscillator}.} 
Consider \( X := \mathbb{R}^2 = \mathbb{R} \times \mathbb{R} \ni (y, p) \), with constants \( m, k \in \mathbb{R} \). The Hamiltonian is given by
\[
H(y, p) := \frac{p^2}{2m} + \frac{1}{2} k y^2,
\]
which represents the sum of the kinetic and potential energies of a mass \( m \) attached to a spring with spring constant\index{spring constant} \( k \). Here, \( (y, p) \) denotes the displacement and momentum of the mass. With the natural symplectic form
\[
\Omega(V, W) := W \cdot \begin{pmatrix} 0 & -1 \\ 1 & 0 \end{pmatrix} V, \quad \text{for } V, W \in \mathbb{R}^2,
\]
the associated Hamiltonian vector field is given by 
\[
X_H(y, p) := \begin{pmatrix} 0 & \frac{1}{m} \\ -k & 0 \end{pmatrix} \begin{pmatrix} y \\ p \end{pmatrix} = \begin{pmatrix} \frac{p}{m} \\ -ky \end{pmatrix}.
\]
Therefore, the Hamiltonian system is described by the equation
\[
\frac{d}{dt} \begin{pmatrix} y(t) \\ p(t) \end{pmatrix} = \begin{pmatrix} \frac{p(t)}{m} \\ -k y(t) \end{pmatrix} = \begin{pmatrix} \text{velocity} \\ \text{spring force} \end{pmatrix},
\]
which describes a simple harmonic oscillator.

\textbf{Example 2: The Airy Equation\index{Airy equation}.} 
Let us now consider \( X := L^2_0(\mathbb{T}) \), the set of square-integrable functions\index{square-integrable functions} with zero average on the torus \( \mathbb{T} \). The Airy equation is a simple linear PDE\index{linear PDE} given by
\[
\partial_t f(t, x) + \partial_{xxx} f(t, x)=0, \quad f(0, x) = f_0(x) \quad \text{for } (t, x) \in \mathbb{R}_+ \times \mathbb{T}.
\]
To formulate this as a Hamiltonian system, we consider a symplectic form \( \Omega \) and a Hamiltonian functional \( H \) given by
\[
\Omega(f, g) := \int (\partial_x^{-1} f)(x) g(x) \, dx, \quad H(f) := \frac{1}{2} \int_{\mathbb{T}} |\partial_x f(x)|^2 \, dx.
\]
The associated Hamiltonian vector field is \( X_H(f) := -\partial_{xxx} f \), which corresponds to the Airy equation. This equation can be solved explicitly using the Fourier transform, yielding the solution
\[
\widehat{f}(t, j) = \sum_{j \ne 0} \widehat{f}_0(j) e^{i(j^3 t + jx)},
\]
where \( \widehat{f}(j) \) denotes the \( j \)-th Fourier mode.

One key observation from these examples is that the linear operators associated with the Hamiltonian vector fields, \( \begin{pmatrix} 0 & \frac{1}{m} \\ -k & 0 \end{pmatrix} \) and \( \partial_{xxx} \), have purely imaginary eigenvalues\index{imaginary eigenvalues}, which lead to oscillations\index{oscillations} in motion. This is characteristic of Hamiltonian systems, where near equilibrium, solutions tend to exhibit oscillatory behavior. Such behavior is typical in many Hamiltonian systems, where the linearized model approximates the full system near an equilibrium\index{equilibrium}. In addition to periodic motion\index{periodic motion}, as shown in the simple harmonic oscillator, systems with many particles or waves can exhibit quasiperiodic motion, where different components oscillate with different frequencies.

As we will see in Chapter~\ref{11jshsdsss2wxowd2}, the gSQG equations can be formulated as a Hamiltonian system in an infinite-dimensional phase space\index{infinite-dimensional phase space}. In light of our earlier discussion on constructing global solutions, this observation leads us to investigate the existence of periodic and quasiperiodic solutions near a steady state\index{steady state}. To this end, we will first review several key research works related to global well-posedness, periodic and quasiperiodic solutions not only to the gSQG equations but also to various mathematical models.  Subsequently, we will present our main theorem and proof strategy, which leverages KAM theory, a systematic method for constructing quasiperiodic solutions in general Hamiltonian systems.

\section{Overview of related works}\color{black}

\subsection{Patch problems}

In this monograph we will work in the \emph{patch} setting, where $\omega(\cdot,t)=1_{D(t)}$ is an indicator function of a simply-connected\index{simply-connected}, bounded set that moves with the fluid. In such a situation, we parametrize $\partial D(t)$ as $z(\theta,t), \, \theta \in [0,2\pi]$ and the evolution equations read:

\begin{align}\label{sqg-patch}
\partial_t z(\theta,t) = \int_{0}^{2\pi} \frac{\partial_{\theta} z(\theta,t) - \partial_{\theta} z(\theta-\eta,t) }{|z(\theta,t) - z(\theta-\eta,t)|^\alpha} d\eta + c(\theta,t) \partial_{\theta} z(\theta,t),
\end{align}
where $c(\theta,t)$ accounts for the reparametrization freedom of the curve.

 Concerning well-posedness results for patch solutions, Rodrigo (in a $C^\infty$ space) \cite{Rodrigo:evolution-sharp-fronts-qg} and Gancedo \cite{Gancedo:existence-alpha-patch-sobolev} and Chae--Constantin--C\'ordoba--Gancedo--Wu \cite{Chae-Constantin-Cordoba-Gancedo-Wu:gsqg-singular-velocities} (in a Sobolev space) proved local existence for the case $0 < \alpha \leq 1$ and $1 < \alpha$ respectively. See also \cite{Kiselev-Yao-Zlatos:local-regularity-sqg-patch-boundary,Gancedo-Patel:local-existence-blowup-gsqg,Ai-Avadanei:well-posedness-sqg-front}.

\subsection{Steady solutions and global existence of the gSQG equation}
The construction of nontrivial global solutions for the generalized SQG equations is a very challenging open problem for all parameters $\alpha\in(0,2)$, both in the smooth case and in the patch case. For $\alpha = 0$ (the 2D Euler equations), global regularity of solutions was well-understood a long time ago, both in the smooth case and in the patch case. See for example the classical papers of Wolibner \cite{Wolibner:theoreme-existence-fluide-parfait-temps-infiniment-long}, Yudovich \cite{Yudovich:Nonstationary-ideal-incompressible}, Burbea \cite{Burbea:motions-vortex-patches}, Chemin \cite{Chemin:persistance-structures-fluides-incompressibles}, and Bertozzi-Constantin \cite{Bertozzi-Constantin:global-regularity-vortex-patches}. However, the construction of global solutions in the case of $\alpha\in(0,2)$ is much more challenging than when $\alpha=0$, since the velocity is more singular, and only partial results have been obtained in recent years. We review some of these results below.

Most of the results around global existence of the gSQG equation have revolved around solutions that exhibited some rigid character (steady, uniformly rotating\index{Uniformly rotating} --V-states--, traveling). In the case where $0 < \alpha < 1$, Hassainia--Hmidi \cite{Hassainia-Hmidi:v-states-generalized-sqg} proved the existence of V-states\index{V-states} with $C^k$ boundary regularity. Castro--C\'ordoba--G\'omez-Serrano then expanded upon this result in \cite{Castro-Cordoba-GomezSerrano:existence-regularity-vstates-gsqg} by showing that V-states also exist with $C^{\infty}$ boundary regularity in the remaining open cases of $\alpha \in [1,2)$ for existence and $\alpha \in (0,2)$ for regularity. This boundary regularity\index{boundary regularity} was later refined to be analytic in \cite{Castro-Cordoba-GomezSerrano:analytic-vstates-ellipses}. Other notable works on rotating solutions include \cite{Cuba-Ferreira:existence-asymmetric-vstates-sqg,Garcia:vortex-patch-choreography,Hmidi-Mateu:existence-corotating-counter-rotating,Hassainia-Wheeler:multipole-vstates-active-scalars}, which discuss other families of rotating solutions or even more steady states, \cite{delaHoz-Hassainia-Hmidi:doubly-connected-vstates-gsqg,Renault:relative-equlibria-holes-sqg} which address the doubly connected case, and \cite{Castro-Cordoba-GomezSerrano:global-smooth-solutions-sqg} which presents a construction in the smooth setting. 

In \cite{delaHoz-Hassainia-Hmidi:doubly-connected-vstates-gsqg}, de la Hoz--Hassainia--Hmidi showed that there exist non-radial patches bifurcating from annuli at negative angular velocities and G\'omez-Serrano \cite{GomezSerrano:stationary-patches} constructed non-radial, doubly connected stationary patches. Garc\'ia \cite{Garcia:Karman-vortex-street} proved the existence of a K\'arm\'an vortex street structure by desingularizing an infinite array of point vortices in the case $\alpha \in [0,1)$. In \cite{Castro-Cordoba-GomezSerrano-MartinZamora:remarks-geometric-properties-sqg} it was ruled out by Castro--C\'ordoba--G\'omez-Serrano--Mart\'in Zamora that ellipses could be rotating solutions\index{rotating solutions} for $\alpha > 0$, as opposed to the case $\alpha = 0$. Gravejat--Smets \cite{Gravejat-Smets:travelling-waves-smooth-sqg}, in the case $\alpha = 1$, constructed smooth translating solutions. Ao--D\'avila--del Pino--Musso--Wei \cite{Ao-Davila-delPino-Musso-Wei:travelling-rotating-solutions-gsqg}, expanded the range to $\alpha \in (0,2)$ as well as to rotating solutions. See also \cite{GodardCadillac:smooth-traveling-waves-sqg,GodardCadillac-Gravejat-Smets:corotating-nfold-sqg} and \cite{Cao-Qin-Zhan-Zou:corotating-traveling-gsqg,Cao-Qin-Zhan-Zou:vstates-gsqg,Cao-Qin-Zhan-Zou:corotating-traveling-gsqg-doubly-connected} for alternative constructions. In  \cite{GomezSerrano-Park-Shi-Yao:radial-symmetry-stationary-solutions}, G\'omez-Serrano--Park--Shi--Yao proved that any smooth, non-negative rotating solution with simply-connected superlevel sets can only rotate with positive angular velocity, and in the case of a patch of fixed area derived moreover a sharp upper bound on the angular velocity.

The drawback of the aforementioned solutions is that they are \textit{special} in the sense that general solutions will not have such behavior. Concerning results for general solutions, C\'ordoba--G\'omez-Serrano--Ionescu \cite{Cordoba-GomezSerrano-Ionescu:global-generalized-sqg-patch} proved global existence for small patch data close to a halfplane in the case $\alpha \in (1,2)$, using a different mechanism based on dispersion and decay. This was extended in \cite{Hunter-Shu-Zhang:global-model-front-sqg,Hunter-Shu-Zhang:global-gsqg}. The main idea was to show that general initial data that are small perturbations of the halfplane stationary patch solution lead to global solutions that decay in time (at an optimal rate of $t^{-1/2}$), thus converging back to the halfplane stationary patch. Unfortunately, the mechanism of dispersion and decay seems to require unbounded domains and, in particular, infinite energy solutions\index{infinite energy solutions}.

In a different direction, one could hope to use the mechanism of inviscid damping\index{inviscid damping} to construct families of global-in-time solutions around explicit stationary solutions of finite energy\index{finite energy}, such as smooth shear flows\index{shear flows} or vortices. This has been successfully implemented in recent years in the 2D Euler case $\alpha=0$, for perturbations of the Couette flow\index{Couette flow} (by Bedrossian--Masmoudi \cite{Bedrossian-Masmoudi:inviscid-damping-2d-euler} and Ionescu--Jia \cite{Ionescu-Jia:inviscid-damping-couette-channel}) and then general monotonic shear flows\index{monotonic shear flows} \cite{Ionescu-Jia:nonlinear-inviscid-damping-monotonic-shear,Masmoudi-Zhao:nonlinear-inviscid-damping-monotone-shear-finite-channel}. It is tempting to try to adapt the mechanism of inviscid damping to construct families of nontrivial global solutions of the gSQG equations, at least for some parameters $\alpha>0$ small. The easiest would be to perturb around the Couette flow corresponding to $\theta(t,x,y)=-1$ on the bounded channel $\mathcal{D}=\mathbb{T}\times[0,1]$. Unfortunately and surprisingly, recent work of G\'omez-Serrano--Ionescu--Jia (discussed in \cite{Ionescu-Jia:icm-survey}) shows that this fails to produce global solutions for any parameters $\alpha>0$, due to a forward cascade that leads to loss of regularity in  finite time.

\subsection{Quasiperiodic solutions in PDE}

Our main goal in this monograph is to demonstrate the existence of large families of global solutions of the generalized SQG equations. We do this using KAM theory\index{KAM theory}, by constructing quasiperiodic solutions\index{quasiperiodic solutions} for almost all initial data\index{initial data} in a neighborhood of the unit disk\index{unit disk} (the simplest stationary patch solution with finite energy).

The first application of KAM theory \cite{Kolmogorov:kam-theory,Arnold:kam-theory,Moser:kam-theory} was to prove the existence of invariant tori that were small perturbations of finite dimensional nearly integrable Hamiltonian systems\index{integrable Hamiltonian system}. In order to upgrade it to the infinite dimensional (PDE) case, the first results are due to Kuksin \cite{Kuksin:hamiltonian-perturbations-kam}, Wayne \cite{Wayne:periodic-quasiperiodic-nlw-kam}, P\"oschel \cite{Poschel:quasiperiodic-solutions-nlw} for 1-d semilinear wave and Schr\"odinger equations with Dirichlet boundary conditions and Craig--Wayne 
\cite{Craig-Wayne:newton-periodic-nlw}, Bourgain 
\cite{Bourgain:green-function-estimates-schrodinger-operators}, Gr\'ebert--Kappeler \cite{Grebert-Kappeler:defocusing-nls-book} and Chierchia--You 
\cite{Chierchia-You:kam-1d-nlw-periodic} with periodic boundary conditions. See also \cite{Kuksin:analysis-hamiltonian-pde}. In the semilinear multidimensional case, we refer to the works of Bourgain \cite{Bourgain:quasiperiodic-solutions-hamiltonian-nonlinear-pde}, Eliasson--Kuksin \cite{Eliasson-Kuksin:kam-nls}, Gr\'ebert--Paturel \cite{Grebert-Paturel:kam-klein-gordon-Sd}, Wang \cite{Wang:quasiperiodic-nls} and Berti--Bolle 
\cite{Berti-Bolle:quasiperiodic-nlw-d-dim} and references therein. See also De la Llave--Sire
\cite{DeLaLlave-Sire:a-posteriori-kam-ill-posed-pde}. Note that all the previous results only were able to deal with semilinear problems\index{semilinear problems}.

In the last decade there has been an emergence of results of quasiperiodic solutions for quasilinear PDE\index{quasilinear PDE}, motivated by applications to the dynamics of confined fluids, building up and polishing the techniques and the methods and culminating with excellent theorems. Baldi--Berti--Montalto constructed quasiperiodic solutions to the Airy equation \cite{Baldi-Berti-Montalto:KAM-quasilinear-airy} and KdV and mKdV \cite{Baldi-Berti-Montalto:KAM-quasilinear-kdv-announcement,Baldi-Berti-Montalto:KAM-quasilinear-kdv,Baldi-Berti-Montalto:KAM-quasilinear-mkdv}. See also the results of Giuliani for gKdV \cite{Giuliani:quasiperiodic-generalized-kdv}, and
\cite{Feola:kam-quasilinear-nls,Feola-Procesi:quasiperiodic-schrodinger,Montalto:quasi-periodic-kirchhoff,Feola-Giuliani-Procesi:kam-tori-degasperis-procesi} and references therein for other relevant models. In the context of water waves, Baldi--Berti--Haus--Montalto \cite{Baldi-Berti-Haus-Montalto:quasiperiodic-gravity-waves-finite-depth} (gravity case), Berti--Montalto \cite{Berti-Montalto:quasiperiodic-standing-gravity-capillary} (gravity-capillary case), Feola--Giuliani \cite{Feola-Giuliani:quasiperiodic-water-waves,Feola-Giuliani:quasiperiodic-water-waves-announcement} (infinite depth) and Berti--Franzoi--Maspero \cite{Berti-Franzoi-Maspero:traveling-quasiperiodic-water-waves-vorticity} (constant non-zero vorticity) constructed quasiperiodic solutions. Numerically, Wilkening--Zhao \cite{Wilkening-Zhao:quasiperiodic-gravity-capillary,Wilkening-Zhao:spatially-quasiperiodic-infinite-depth} computed quasiperiodic gravity-capillary water waves\index{gravity-capillary water waves} in the infinite depth case.


Berti--Hassainia--Masmoudi \cite{Berti-Hassainia-Masmoudi:time-quasiperiodic-vortex-patches} constructed quasiperiodic solutions close to elliptical vortex\index{elliptical vortex} patches\index{vortex patches}, introducing the angular momentum as a symplectic variable. Hassainia--Roulley \cite{Hassainia-Roulley:quasiperiodic-euler-boundary} constructed quasiperiodic solutions of the 2D Euler equations in a bounded domain, Roulley \cite{Roulley:periodic-quasiperiodic-euler-alpha} proved its existence for the Euler-$\alpha$ equation and Hmidi--Roulley \cite{Hmidi-Roulley:quasiperiodic-qgsw} for the QGSW equations.

Other examples of quasiperiodic solutions in the context of the incompressible Euler and Navier-Stokes equations, even in high dimensions, were obtained by Crouseilles--Faou, {Elgindi--Jeong}, Enciso--Peralta-Salas--Torres de Lizaur \cite{Crouseilles-Faou:quasiperiodic-2d-euler,Elgindi-Jeong:symmetries-fluids,
Enciso-PeraltaSalas-TorresdeLizaur:quasiperiodic-solutions-euler} for Euler, using non-KAM constructions, Baldi--Montalto \cite{Baldi-Montalto:quasiperiodic-euler-3d} for forced Euler, using a KAM construction  and Franzoi--Montalto, Montalto \cite{Franzoi-Montalto:kam-inviscid-limit-navier-stokes,Montalto:quasiperiodic-navier-stokes} for forced Navier-Stokes, using a KAM construction. Finally, we would like to draw the attention to the recent results by Hassainia--Hmidi--Masmoudi \cite{Hassainia-Hmidi-Masmoudi:kam-gsqg} who proved the existence of global quasiperiodic solutions for the generalized SQG equations, for a set of parameters $\alpha\in\left(0,1/2\right)$. The set of acceptable parameters $\alpha$ is unknown, but of full measure in $(0,1/2)$.

We emphasize that most of these recent results in the quasilinear case (with the notable exception of the papers \cite{Feola-Giuliani:quasiperiodic-water-waves}, \cite{Giuliani:quasiperiodic-generalized-kdv} and \cite{Baldi-Berti-Montalto:KAM-quasilinear-kdv}) rely on the use of {\it{external parameters}}. Quasiperiodic solutions are then constructed for all initial data, but for an unknown set of parameters, usually generic\index{generic} of full measure. The point is that the presence of external parameters\index{external parameter} improves significantly the structure of the resonances of the system, which plays a key role in the analysis. 

The drawback is that the family of acceptable parameters is not explicit, and one cannot guarantee that quasiperiodic solutions exist for a specific given equation. Our broad goal in this monograph is to develop a robust and flexible method to construct quasiperiodic solutions for certain fluid models, without requiring the presence of external parameters. The basic idea is to replace the genericity of the external parameters with genericity of the initial data. This leads however to very significant difficulties at the implementation level; see below for a more detailed discussion.

 \subsection{Weak solutions and finite time singularities}

The generalized SQG equations have been studied extensively, by many authors. In this subsection we discuss two other areas of active research, and provide some references.

In his thesis \cite{Resnick:phd-thesis-sqg-chicago}, Resnick demonstrated the global existence of weak solutions\index{weak solutions} in $L^2$ through the use of the oddness of the Riesz transform\index{Riesz transform} to achieve additional cancellation. Marchand \cite{Marchand:existence-regularity-weak-solutions-sqg} later extended this result to include initial data belonging to $L^p$ with $p$ greater than $\frac43$. See also \cite{Nahmod-Pavlovic-Staffilani-Totz:global-invariant-measures-gsqg} for other existence results concerning weak solutions. Non-uniqueness of weak solutions of SQG remains a difficult problem, with progress being made through works such as Azzam--Bedrossian \cite{Azzam-Bedrossian:bmo-uniqueness-active-scalar-equations} or Isett--Vicol \cite{Isett-Vicol:holder-continuous-active-scalar}, and most importantly, Buckmaster--Shkoller--Vicol \cite{Buckmaster-Shkoller-Vicol:nonuniqueness-sqg}, as well as alternative proofs by Isett--Ma \cite{Isett-Ma:non-uniqueness-sqg} and the investigation of the stationary problem by Cheng--Kwon--Li \cite{Cheng-Kwon-Li:non-uniqueness-steady-state-sqg}.

 One of the most significant questions in mathematical fluid mechanics is whether the SQG and gSQG system exhibits finite time singularities or has global existence. Kiselev--Nazarov \cite{Kiselev-Nazarov:simple-energy-pump-sqg} created solutions that exhibited norm inflation, and Friedlander--Shvydkoy \cite{Friedlander-Shvydkoy:unstable-spectrum-sqg} demonstrated the presence of unstable eigenvalues\index{unstable eigenvalues} in the spectrum. He--Kiselev \cite{He-Kiselev:small-scale-creation-sqg} proved an exponential in time growth of the $C^2$-norm. See also the construction of singular solutions with infinite energy by Castro--C\'ordoba \cite{Castro-Cordoba:infinite-energy-sqg} and ill-posedness\index{ill-posedness} results by C\'ordoba--Mart\'inez-Zoroa and Jeong--Kim \cite{Cordoba-MartinezZoroa:ill-posedness-sqg,Jeong-Kim:illposedness-sqg}.

In order to understand the possibility of a finite time blow-up scenario, numerical studies have been conducted. Constantin--Majda--Tabak \cite{Constantin-Majda-Tabak:formation-fronts-qg} suggested that a singularity in the form of a hyperbolic saddle\index{hyperbolic saddle} may occur, closing in a finite amount of time. However, Ohkitani--Yamada \cite{Ohkitani-Yamada:inviscid-limit-sqg} and Constantin--Nie--Sch\"orghofer \cite{Constantin-Nie-Schorghofer:nonsingular-sqg-flow} proposed that the growth was actually double exponential. C\'ordoba \cite{Cordoba:nonexistence-hyperbolic-blowup-qg} bounded the growth at quadruple exponential, and later C\'ordoba and Fefferman \cite{Cordoba-Fefferman:growth-solutions-qg-2d-euler} proposed a double exponential bound, which was supported by numerical simulations from Deng--Hou--Li--Yu \cite{Deng-Hou-Li-Yu:non-blowup-2d-sqg}. Constantin--Lai--Sharma--Tseng--Wu \cite{Constantin-Lai-Sharma-Tseng-Wu:new-numerics-sqg} later reexamined the hyperbolic saddle scenario using improved algorithms and found no evidence of blowup. Scott \cite{Scott:scenario-singularity-quasigeostrophic} proposed a scenario in which filamentation\index{filamentation} occurs and blowup\index{blowup} of $\nabla \theta$ occurs after several cascades, starting from elliptical configurations. This is currently the only scenario that remains valid in the smooth setting. In \cite{Garcia-GomezSerrano:self-similar-spirals-gsqg}, very recently, Garc\'ia--G\'omez-Serrano constructed a big class of non-trivial self-similar spiral solutions close to radial ones with a mild singularity at the origin.

Even though the finite time singularity problem seems elusive, there exist several numerical scenarios suggesting such a singularity. The first one, proposed by C\'ordoba--Fontelos--Mancho--Rodrigo \cite{Cordoba-Fontelos-Mancho-Rodrigo:evidence-singularities-contour-dynamics} initially starts as two patches rolling onto each other and finally collapsing. At the intersection point the curvature\index{curvature} blows up (the curve should lose regularity due to the results by Gancedo and Strain \cite{Gancedo-Strain:absence-splash-muskat-SQG}, see also \cite{Kiselev-Luo:nonexistence-splash-sqg,Jeon-Zlatos:no-splash-gsqg}) and the collapse is suggested to be asymptotically self-similar. The second scenario was proposed by Scott--Dritschel \cite{Scott-Dritschel:self-similar-sqg}, taking ellipses as initial condition; starting with an aspect ratio of $0.16$, they report a self-similar cascade of filamentation. In \cite{Scott-Dritschel:self-similar-sqg-patch}, again taking ellipses as initial condition and combining numerical analysis with asymptotic calculations\index{asymptotic calculations}, they conjecture a scenario where the patch develops a corner in finite time, together with a self-similar spiral\index{self-similar spiral}. Finally, Kiselev--Ryzhik--Yao--Zlato\v{s} \cite{Kiselev-Ryzhik-Yao-Zlatos:singularity-alpha-patch-boundary} (for $ 0 < \alpha < \frac{1}{12}$) and later Gancedo--Patel \cite{Gancedo-Patel:local-existence-blowup-gsqg} (for $ 0 < \alpha < \frac{1}{3}$) construct finite time singularities in the presence of a boundary.

\section{Main result}\label{main_part}
Before we state the main result, let us first recall the definition of a quasiperiodic function\index{quasiperiodic function}:
\begin{definition}\label{def_quasi_per}
Let $X$ be a Hilbert space\index{Hilbert space} and $\nu\in\mathbb{N}$ be a fixed natural number\index{natural number}. A function\index{function} $f:\mathbb{R}\mapsto  X$ is said to be quasiperiodic\index{quasiperiodic} with frequency\index{frequency} ${\omega}\in\mathbb{R}^\nu$, if there exists  $i:\mathbb{T}^\nu\mapsto X$ such that $f(t)=i({\omega} t).$
\end{definition} 

In this monograph, we consider a  patch solution\index{patch solution} to \eqref{sqg-patch} of the form:
\begin{align}\label{parametrization_z}
z(x,t):= \sqrt{1+f(x,t)}(\cos x, \sin x),\text{ for some $f(\cdot,t):\mathbb{T}\mapsto (-1,\infty)$.}
\end{align}
Note that one of the advantages of the use of the variable $f$, instead of a more natural parametrization\index{parametrization}  $z(x,t)=R(x,t)(\cos x,\sin x)$ relies on the conservation of the area\index{conservation of area} of the patch in the dynamics in \eqref{gSQG_1}; if the patch initially has area  $|D(0)|=\pi$, then $|D(t)|=\pi$ for all $t\ge0$, therefore
\begin{align}\label{total_vorticity_preserve_intro}
\pi = |D(t)| =\frac{1}{2}\int_{\mathbb{T}}R(x,t)^2dx = \pi + \int_{\mathbb{T}}f(x,t)dx.
\end{align}
Thus, we can assume that $f$ has zero average in the variable $x$.
 
Plugging \eqref{parametrization_z} into \eqref{sqg-patch}, one can  find that the evolution of $f$ can be expressed as (we refer to Section~\ref{11Hamitons1} for more detailed computations)
\begin{align}\label{evol_Req}
\partial_t f(x,t) = \frac{2}{2-\alpha}\partial_\theta\left(\int_{\mathbb{T}}\frac{(z(x,t)-z(y,t))\cdot\partial_x z(y,t)^{\perp}}{|z(x,t)-z(y,t)|^{\alpha}}dy\right)=:X_{\text{gSQG}}(f(x,t)).
\end{align}
As noted in \cite{Hassainia-Hmidi:v-states-generalized-sqg,Resnick:phd-thesis-sqg-chicago,Marsden-Weinstein:coadjoint-orbits-vortex-patch,Saffman:book-vortex-dynamics}, the equation \eqref{evol_Req} can be seen as a Hamiltonian system\index{Hamiltonian system} with the associated Hamiltonian 
\begin{align}\label{hamiltonian_intro11}
\mathcal{H}(f):=\int_{D} 1_D * \frac{1}{|\cdot|^\alpha}(x)dx,
\end{align} where $D$ is the patch determined by the parametrization $f$ as in \eqref{parametrization_z} (see Chapter~\ref{11jshsdsss2wxowd2}). More precisely, the vector field $X_{\text{gSQG}}(f)$ is given by
\begin{align}\label{Hamil_vec_intro}
X_{\text{gSQG}}(f)=\partial_x \left( \nabla_{L^2}\mathcal{H}(f)\right),
\end{align}
where $\nabla_{L^2}\mathcal{H}(f)$ denotes the gradient vector field\index{gradient vector field} of $\mathcal{H}$ at $f$ in the space $L^2(\mathbb{T})$.
 
The linearized equation\index{linearized equation} of \eqref{evol_Req} at the unit disk ($f=0$) can be written as (see Proposition~\ref{expansion_1})
\begin{align}\label{linear_intro_operator_language1}
f_t = \frac{d}{dt}X_{\text{gSQG}}(tf)\bigg|_{t=0} = \partial_x \left( -\frac{1}{2}\Lambda^{\alpha-1}f + \frac{T_\alpha}{4}f\right), 
\end{align}
where
\begin{align}\label{defoflamba_intro}
\Lambda^{\alpha-1}f(x)& :=\int_\mathbb{T}(2-2\cos(x-y))^{-\frac{\alpha}{2}}(f(x)-f(y))dy, \\
& \text{ and } T_\alpha:=  \frac{2\pi\Gamma(3-\alpha)}{\Gamma(2-\frac{\alpha}{2})\Gamma(2-\frac{\alpha}{2})}. \nonumber
\end{align}
One can also rewrite the linearized equation \eqref{linear_intro_operator_language1} as 
\begin{align}\label{linear_intro}
\partial_tf =Op(\ii W(j))[f](x,t),
\end{align}
where $Op(W(j))$ denotes the pseudo differential operator\index{Pseudo differential operator} associated to the the Fourier\index{multipliers} multiplier\index{Fourier multiplier} $W(j)$, defined as
\begin{align}\label{intro_linear_rfs1}
W(j):=  j\left(-\frac{1}{2} \mathcal{C}_\alpha\left(\frac{\Gamma(|j|+\frac{\alpha}{2})}{\Gamma(1+|j|-\frac{\alpha}2)} - \frac{\Gamma(\frac{\alpha}{2})}{\Gamma(1-\frac{\alpha}2)}\right) + \frac{\pi(-1)^{j}\Gamma(3-\alpha)}{2\Gamma(2 - \frac{\alpha}2)\Gamma(2 - \frac{\alpha}2)}\right), 
\end{align}
where $\mathcal{C}_\alpha:=-\frac{2\pi\Gamma(1-\alpha)}{\Gamma\left(\frac{\alpha}{2}\right)\Gamma\left(1-\frac{\alpha}{2}\right)}$.
A classical asymptotic analysis for the Gamma function\index{Gamma function} tells us that $W(j)$ exhibits an asymptotic behavior like $j|j|^{\alpha-1}$ (e.g.  \cite[Theorem 2.1]{Laforgia-Natalini:asymptotic-expansion-ratio-gamma}), more precisely,
\[
W(j) = C(\alpha) j |j|^{\alpha-1} + O(1),\text{ for some constant $C(\alpha)\in\mathbb{R}$ for  $\alpha\in (0,2)\backslash \left\{1\right\}$}.
\]
Given a set of natural numbers $S^+:=\left\{ j_1,\ldots, j_\nu \right\}\subset \mathbb{N}$ (also denoting $S:=\left\{ \pm j: j\in S^+\right\}$), the linear equation \eqref{linear_intro} possesses time-quasiperiodic\index{time-quasiperiodic} solutions of the form
\begin{align}\label{linear_sol_intro}
f(x,t) = \sum_{j_k\in {S}}\sqrt{|j_k|\zeta_{k}}e^{\ii (W(j_k) t + j_k x)} = \sum_{j_k\in S^+} 2\sqrt{|j_k|\zeta_{k}}\cos(W(j_k)t+ j_k x), 
\end{align}
for some $\zeta_1,\ldots,\zeta_\nu>0$, for which the $j$-th Fourier coefficient\index{Fourier coefficient} is oscillating in time with frequency $W(j)$. Indeed, according to Definition~\ref{def_quasi_per}, the solution \eqref{linear_sol_intro} to the linearized equation can be expressed as
\begin{align}\nonumber
f(t,\cdot)&=f(t)=i_{\text{linear}}(\overline{\omega}t), \\
\text{ where }i_{\text{linear}}(\varphi)&:=\sum_{j_k\in S^+}2\sqrt{\zeta_k}\cos(\varphi_k + j_kx) \text{ and  $\overline{\omega}_{k}:=W(j_k)$.}\label{linear_embedding}
\end{align}
This naturally leads to the question whether there exists such a time-quasiperiodic solution to the full nonlinear problem \eqref{evol_Req} around the steady state $f=0$.

 In our analysis, we make use of several invariance properties of the equation~\eqref{parametrization_z}. One is the so-called time-reversibility with respect to the involution $\rho: f(x)\mapsto f(-x)$, namely, 
 \[
 \mathcal{H}(\rho (f)) = \mathcal{H}(f), \text{ where $\mathcal{H}$ is defined as in \eqref{hamiltonian_intro11}.}
 \]
 We say that a solution $f(x,t)$ to \eqref{evol_Req} is reversible\index{reversible}, if $\rho(f)(x,-t)$ is also a solution. Another invariance property is the rotational invariance\index{rotational invariance} of solutions. More precisely, given an integer $\mathtt{M}\in\mathbb{N}$, if the initial data of the gSQG equation is invariant under a $\frac{2\pi}{\mathtt{M}}$-rotation, then the solution at any time is also invariant under a $\frac{2\pi}{\mathtt{M}}$-rotation. Such an $\mathtt{M}$-fold symmetric patch can be associated to a parametrization $f$ in \eqref{parametrization_z} being invariant under a  $\frac{2\pi}{\mathtt{M}}$-translation of the variable $x$:
 \begin{align}\label{mforl_intro_1}
 f(t,x+\frac{2\pi}{\mathtt{M}}) = f(t,x).
 \end{align}
  Those properties of the gSQG equation will be studied in detail in Chapter~\ref{11jshsdsss2wxowd2}.

The main theorem we prove in this monograph is the following (stated informally, we refer to Theorem~\ref{main1} for a precise statement):

\begin{theorem}\label{main_2}
[\textbf{= Theorem \ref{main1}}]
Let $\alpha\in (1,2)$, $S^+ \subset \mathbb{N}$ and a symmetry class\index{symmetry class} $\mathtt{M}$\index{$\mathtt{M}$} satisfying some non-resonance conditions (cf. Section~\ref{posdsdtangeinal} and \ref{rpoisdsd1sd}) be fixed. Then, for all sufficiently small $\epsilon>0$, there exists a set of amplitudes $A_\epsilon\subset [1,2]^{|S^+|}\subset\mathbb{R}^{|S^+|}$ such that for each $\vec{\zeta}\in A_\epsilon$, there exist a frequency vector $\omega=\omega(\vec{\zeta})$ and a time-quasiperiodic solution to \eqref{evol_Req} of the form
\begin{align}\label{sol_expansion_intro}
f(\theta,t) = 2\epsilon\sum_{j_k\in S^+}\sqrt{|j_k|\vec{\zeta}_k}\cos(\omega_kt + j_k x) + o (\epsilon).
\end{align}
The set $A_\epsilon$ is a Cantor-like set\index{Cantor-like set} of asymptotically full measure\index{full measure}, in the sense that $\lim_{\epsilon\to 0}\frac{|A_\epsilon|}{\left|[1,2]^{|S^+|}\right|}=1$. The solution \eqref{sol_expansion_intro} is in some Sobolev space $H^{s_0}$ for some $s_0\gg 1$, and it is reversible and invariant under $\frac{2\pi}{\mathtt{M}}$-translation\index{$\frac{2\pi}{\mathtt{M}}-translation invariance$} in the variable $\theta$. Lastly, the solution is linearly stable\index{linearly stable} under $\frac{2\pi}{\mathtt{M}}$-translation invariant perturbations.
\end{theorem}

Some remarks are in order:

\begin{remark}\label{rmk_allalpha}
As stated above, our proof does not make use of any external parameters\index{external parameter} ($\alpha$ would be the natural candidate) as opposed to \cite{Hassainia-Hmidi-Masmoudi:kam-gsqg} and indeed this results in needing the Diophantine constant\index{Diophantine constant} $\gamma$\index{$\gamma$} to be $\gamma = o(\epsilon^{2})$, which in turn requires Normal Form expansions (cf. Chapter \ref{skpssisodssuw}), and also the computation of the explicit terms of size $O(\epsilon)$ and $O(\epsilon^2)$. Relaxing this constraint  would significantly shorten the length and the complexity of this monograph.
\end{remark}

\begin{remark}
It is conceivable that our proof of Theorem \ref{main1} would also work in the case $\alpha < 1$, changing the relevant sections and estimates. In the case of the SQG equation $(\alpha = 1)$, the analysis in Section~\ref{sjdwjjjwdsegosdeoeo} breaks down since the sum of pseudo differential symbols\index{pseudo differential symbols} is not finite anymore and the regularity losses coming from the Egorov method are not finite. Most of the other parts of the reduction\index{reduction} also hold for all cases of $\alpha$, possibly with minimal changes.
\end{remark}

\begin{remark}
The closer $\alpha$ is to $1$, the more conjugations are required in the reduction of the linearized operator to a constant coefficients operator. In the adaptation of the Egorov method\index{Egorov method}, inspired by \cite{Berti-Montalto:quasiperiodic-standing-gravity-capillary}, we use a slightly more general flow, compared to the ones in \cite{Berti-Montalto:quasiperiodic-standing-gravity-capillary,Hassainia-Hmidi:v-states-generalized-sqg,Feola-Giuliani:quasiperiodic-water-waves}, to avoid a large number of iterations that might cause potential complexity. See Remark~\ref{Egorog_32intsrrqqw1}. 
\end{remark}

\begin{remark}\label{ampli_intro_3}All the amplitudes in the set of  $A_\epsilon$\index{$A_\epsilon$} in Theorem~\ref{main_2}, which takes asymptotically full measure, can possess quasiperiodic solutions, if the choice of tangential sites\index{tangential sites} $S^+$ can be made properly so that some non-resonance conditions are satisfied. The precise conditions are stated in Section~\ref{rpoisdsd1sd}.  Then a natural question is whether such a set $S^+=\left\{ j_1,\ldots,j_\nu\right\}\subset\mathbb{N}$ is generic or not. The non-resonance conditions that we require can be roughly expressed as
\[
P(j_1,j_2,\ldots,j_\nu)\ne 0,\text{ for some function $P:\mathbb{Z}^{\nu}\mapsto \mathbb{R}$}.
\]
Compared to  previous works (e.g. \cite{Giuliani:quasiperiodic-generalized-kdv,Baldi-Berti-Montalto:KAM-quasilinear-kdv,Baldi-Berti-Montalto:KAM-quasilinear-mkdv}), our $P$ involves Gamma functions and the verification of the non-resonance condition is much more complicated. While we expect that such conditions can be satisfied by "generic" choices of $S$ with small $\mathtt{M}$, we will give a rigorous proof  only for the case where the frequencies are supported on  multiples of  sufficiently large  $\mathtt{M}$. This allows us to focus on the asymptotic behavior of the Gamma function in the analysis.
\end{remark}

Before describing the idea of the proof in more detail, we give more explanation about the internal parameter, which we think of the most crucial part of the proof.

 We consider a finite number of Fourier modes $S^+\subset \mathbb{N}$ such that $|S^+| =\nu$. Setting $S:=\left\{ \pm j: j\in S^+\right\}$  and decomposing
\[
L^2(\mathbb{T}) = H_S \oplus \left(H_{S}\right)^{\perp},\quad H_S:=\left\{ f\in{L}^2: f_j = 0,\text{ if $j\notin S$}\right\},\quad H_{S^\perp}=\left( H_S\right)^\perp,
\]
 one can see the SQG dynamics of the tangential component\index{tangential component} in $H_S$\index{$H_S$}, and the normal component\index{normal component} in $H_{S^\perp}$\index{$H_{S^\perp}$}. While the dynamics of the tangential component is finite dimensional, the dynamics of the normal component will be given as a slight variation of the SQG equation~\ref{evol_Req}, due to the influence of the tangential component.  Then under a suitable symplectic transformation\index{symplectic transformation} $\Phi^{WB}$ (See Proposition~\ref{normal_form_prop11231}) and the use of angle-action variables to reparametrize the tangential component (See \eqref{actionssx}), the SQG Hamiltonian $\mathcal{H}$ in \eqref{hamiltonian_intro11} on $L^2(\mathbb{T})$ can be reformulated in terms of another Hamiltonian $H_{\zeta}$ in \eqref{rescaled_hamiltonian} on $\mathbb{T}^\nu\times \mathbb{R}^\nu\times H_{S^\perp}$ taking the form (see Proposition~\ref{aavariable}),
 \[
H_{\zeta}(\theta,y,z) = C_\epsilon(\zeta) + 2\pi{\omega(\zeta)}\cdot y + \mathcal{N}(\theta)(z,z) + P(\theta,y,z),
 \]
 for some $\theta$-dependent bilinear map\index{bilinear map} $\mathcal{N}(\theta)$, and a perturbative term $P(\theta,y,z)$, which must sufficiently small. Note that $\zeta\in \mathbb{R}^{\nu}$ in the change of variable in \eqref{actionssx} denotes a rescaled amplitude vector of the tangential component.
 Without the perturbative term $P$, the corresponding Hamiltonian system give a quasiperiodic motion with the $\zeta$-dependent frequency $\omega(\zeta) = \overline{\omega} + \epsilon^2 \mathbb{A}\zeta$ for some constant frequency $\overline{\omega}$ and $\mathbb{A}\in\mathbb{R}^{\nu\times \nu}$. The key point is that under such process, we see that the amplitude $\zeta$ can modulate the linear frequency $\omega(\zeta)$, which will serve as an internal parameter to perform the KAM theory. Furthermore, since the size of the modulation at the linear level is $O(\epsilon^2)$, the perturbation $P$ must be $o(\epsilon^2)$. 
 
  Such a derivation of the Hamiltonian $H_\zeta$\index{$H_\zeta$} exhibiting the dependence on the amplitude, while ensuring $P=o(\epsilon^2)$ has been successfully implemented in earlier works in the literature, especially for the first time in \cite{Baldi-Berti-Montalto:KAM-quasilinear-kdv-announcement} to the best of our knowledge. We emphasize that the technique that we adapt in this paper  requires a careful choice of the tangential modes $S$ to exclude possible resonances. All the sufficient conditions on the choice of $S$ is described in Section~\ref{rpoisdsd1sd}.

\color{black}
\section{Strategy of the proof and the structure of the monograph}\label{idea_of_the_proof_12}

We first outline the main ideas of the proof of Theorem~\ref{main_2} and  link them with the sections of this monograph afterwards.

In order to describe the strategy of the proof, let us fix
\begin{align}\label{intro_setandnu}
\nu\in \mathbb{N},\quad S^+:=\left\{ j_1<\cdots < j_\nu\right\}\subset\mathbb{N},
\end{align}
and denote 
\begin{align}\label{tangential_sites_normal_intro}
S:=\left\{ \pm j: j\in S^+\right\},\quad S^\perp:=\mathbb{Z}\backslash(S\cup\left\{ 0 \right\}).
\end{align}
In view of \eqref{linear_sol_intro}, one can think of $S$ as a set of  Fourier modes of the solution at the linear level, and $S^\perp$ as the  support of the orthogonal correction term for the solution to the nonlinear problem, while the $0$-th mode is excluded, since we look for a solution with zero average (see \eqref{total_vorticity_preserve_intro}).
We also denote the linear frequency by
\begin{align}\label{linear_fre_intro2}
{\overline{\omega}}\in \mathbb{R}^\nu,\quad (\overline{\omega})_k = W(j_k),\text{ for $k=1,\ldots,\nu$.}
\end{align}
Using the notations in \eqref{evol_Req} and \eqref{linear_intro}, we can rewrite \eqref{evol_Req}  as
 \begin{align}\label{gSQG_full_intro}
 f_t = X_{\text{gSQG}}(f) = Op(\ii W(j))[f] + P_{\text{gSQG}}(f),
 \end{align}
 where $P_{\text{gSQG}}(f)$ collects the nonlinear contribution of the vector field\index{vector field} $X_{\text{gSQG}}(f)$. Since we are interested in the solutions near $f=0$, replacing $f$ by $\epsilon f$ for small $\epsilon>0$, we are led to study the equation of the form:
 \begin{align}\label{gSQG_full_intro2}
 f_t = Op(\ii W(j))[f] + P_{\epsilon,\text{gSQG}}( f),\text{ where $P_{\epsilon,\text{gSQG}}(f):=\frac{1}{\epsilon}P_{\text{gSQG}}(\epsilon f)$}.
 \end{align}
As we observed in Section~\ref{main_part}, we have an embedding $i_{\text{linear}}$ (see \eqref{linear_embedding}) for which $f_{\text{lin}}(t,\theta):=i_{\text{linear}}(\overline{\omega}t)$ solves the linear equation $\partial_tf_{\text{lin}}=Op(\ii W(j))[f_{\text{lin}}]$. The question is whether such an embedding can persist under the nonlinear perturbation as in \eqref{gSQG_full_intro}. Perhaps, one of the most naive attempts  could be plugging the ansatz,
\[
f(t,x)=i(\overline{\omega} t,x),\text{ for some $i:\mathbb{T}^\nu\times \mathbb{T}\mapsto \mathbb{R}$},
\]
into \eqref{gSQG_full_intro}, which leads us to find  $i$ such that
\begin{align}\label{functional_intro_1}
\mathcal{F}_{\text{gSQG}}(i):=\overline{\omega}\cdot \partial_{\varphi}i(\varphi) -Op(\ii W(j))[i(\varphi)] - P_{\epsilon,\text{gSQG}}(i(\varphi)) = 0,\text{ for $\varphi\in\mathbb{T}^\nu$.}
\end{align}
We can think of  $\mathcal{F}_{\text{gSQG}}$ as a map between spaces of functions of $(\varphi,\theta)$. Having the explicit solution $i_{\text{linear}}$ at the linear level given in \eqref{linear_embedding} and noting
 that $P_{\epsilon,\text{gSQG}}$ is "small" depending on $\epsilon>0$, we might expect the sequence of embeddings $i_n$, formally defined in the spirit of Newton's method\index{Newton's method},
 \begin{align}\label{newton_intro}
 i_0=i_{\text{linear}},\quad  i_{n+1}:=i_n -\left(d_i\mathcal{F}_{\text{gSQG}}(i_n)\right)^{-1}[\mathcal{F}_{\text{gSQG}}(i_n)],  \text{ for $n\ge 0$},
 \end{align}
 where
 \begin{align}\nonumber
 d_i\mathcal{F}_{\text{gSQG}}(i)[\ihat]& :=\frac{d}{dt}\mathcal{F}_{\text{gSQG}}(i+t\ihat)\bigg|_{t=0} \\
 & \overset{\eqref{functional_intro_1}}=\overline{\omega}\cdot \partial_\varphi \ihat - Op(\ii W(j))[\ihat] - d_iP_{\epsilon,\text{gSQG}}(i)[\ihat]\label{newton_intro2}\end{align}
 to converge to a solution for \eqref{functional_intro_1}. Clearly, the above argument is far less rigorous, and we will investigate how to modify the strategy.

\subsection{Sketch of the proof, Part 1: Internal parameter and the weak Birkhoff normal form\index{Birkhoff normal form}}

  \subsubsection{ {Loss of derivatives and the Nash-Moser scheme\index{Nash-Moser scheme}}.}
 As a rule of thumb in  usual perturbative problems, the invertibility of the linearized operator\index{linearized operator} $ d_i\mathcal{F}_{\text{gSQG}}(i)$ in \eqref{newton_intro} would rely on the invertibility of the linear part:
 \begin{align}\label{linearpart_1}
 L[\ihat]:=\overline{\omega}\cdot \partial_\varphi \ihat - Op(\ii W(j))[\ihat], 
 \end{align} assuming that the contribution of the perturbative part is negligible\index{negligible}. While looking for an embedding $i$ in Sobolev spaces\index{Sobolev spaces} $H^{s}(\mathbb{T}^{\nu+1})$\index{$H^s(\mathbb{T}^{\nu+1})$} (for large $s\gg 1$), it is not trivial  whether the operator $L$ can be invertible between two fixed Sobolev spaces. However, the classical KAM theory tells us that the invertibility of $L$ can be achieved depending on the frequency vector $\overline{\omega}$. Indeed, for the frequency vectors that satisfy  the so-called "Melnikov condition"\index{Melnikov condition} with some $\gamma,\tau>0$:
 \begin{align}\label{first_menikov_intro1}
 \left\{ \omega: |\ii \omega \cdot l - \ii W(j)|\ge \gamma |l|^{\tau} \text{ for all $(l,j)\in\mathbb{Z}^\nu\times\mathbb{Z}$}\right\},
 \end{align}
 one can formally invert $L$ using the Fourier series, that is,
 \begin{align}\label{formal_inverse1_intro}
 L[\ihat](\varphi,x)=g(\varphi,x) \overset{\eqref{linearpart_1}}\iff \ihat(\varphi,x)=\sum_{(l,j)\in \mathbb{Z}^\nu\times\mathbb{Z}}\frac{\widehat{g}(l,j)}{\ii (\overline{\omega}\cdot l - W(j))}e^{\ii (\varphi\cdot l + j x)},  
  \end{align}
  where $\hat{g}(l,j):=\frac{1}{(2\pi)^{\nu+1}}\int_{\mathbb{T}^{\nu+1}}g(\varphi,x)e^{\ii (\varphi\cdot l + j x)}d\varphi dx$.
Note that the expression of the inverse in \eqref{formal_inverse1_intro} tells us that there is a  regularity mismatch between the image and the domain spaces. If $\ihat\in H^{s}$ for some $s>0$, we see that there is a loss of derivatives\index{loss of derivatives} due to the differential operators  $\partial_\varphi$ and $Op(W(j))$, while inverting $L$ as in \eqref{formal_inverse1_intro} does not gain the same amount of the regularity, and it actually causes another loss of derivatives by $\tau$; the best estimate one can expect under the condition is that
\begin{align}\label{inverse_estimate_1_intro}\rVert L^{-1}[g]\rVert_{H^{s-\tau}}\lesssim \gamma^{-1} \rVert g \rVert_{H^{s}} \text{ for $g\in H^s(\mathbb{T}^{\nu+1})$}.
\end{align}
Therefore, the formal sequence of $i_n$ in \eqref{newton_intro} does not seem to be closed, since $i_{n+1}$ must be less regular than $i_n$ at each iteration. Hence, the crude iteration procedure in \eqref{newton_intro} needs to be replaced by the Nash-Moser scheme\index{Nash-Moser scheme}, projecting each approximate solution $i_n$ into a finite dimensional space so that $i_n$ remains in $C^\infty(\mathbb{T}^{\nu+1})$ for each $n\ge0$.
\subsubsection{Internal parameter.}
We have observed above that the condition on the frequency vector in \eqref{first_menikov_intro1} is one of the necessary conditions to perform the iteration \eqref{newton_intro}. Then a very natural question is how to check whether the linear frequency $\overline{\omega}$ satisfies such a condition. In general (for fixed $\gamma,\tau>0$), it is very hard to determine whether a given vector $\omega\in\mathbb{R}^\nu$ satisfies even a more relaxed condition (so-called Diophantine condition):
 \begin{align}\label{diophaintime_intro}
 |\omega \cdot l| > \gamma |l|^{\tau} ,\text{ for all $l\in\mathbb{Z}^\nu$.}
 \end{align}
However, it is well-known that given an open set $\Omega\subset \mathbb{R}^\nu$, "almost all" $\omega\in\Omega$ satisfy \eqref{diophaintime_intro}, more precisely, such  non-resonance frequency vectors take asymptotically full measure in $\Omega$ as $\gamma\to 0$. Indeed, the KAM theory does not tell us exactly which frequency vector can possess a quasiperiodic solution, but rather it tells us that the set of frequency vectors that possess a quasiperiodic solution has non-zero measure in a given set of frequencies. This is why we need  parameter-dependent equations to perform the KAM theory; if our equation depends on a parameter, and the parameter can "properly"modulate the linear frequency $\overline{\omega}$, then for almost all parameter values, we might expect to invert the operator $L$. The gSQG equation \eqref{gSQG_1} certainly involves a parameter $\alpha\in (0,2)$, therefore one might be tempted to use $\alpha$ to modulate the linear frequency by looking at $\overline{\omega}$ in \eqref{linear_fre_intro2} and \eqref{intro_linear_rfs1} as a function of $\alpha$, that is, $\overline{\omega}=\overline{\omega}(\alpha)$. This attempt would enable us to obtain quasiperiodic solutions for "almost every" $\alpha$ (without knowing precisely which $\alpha$ satisfies the condition), while such a result cannot be, in principle, obtained for \emph{every} $\alpha$. For this reason, we follow the strategy in \cite{Baldi-Berti-Montalto:KAM-quasilinear-kdv,Baldi-Berti-Montalto:KAM-quasilinear-mkdv,Feola-Giuliani:quasiperiodic-water-waves} and derive a weak Birkhoff normal form of the Hamiltonian $\mathcal{H}$, from which we can see a modulation effect of the linear frequency by the amplitude. In this regard, a bit more precise explanation will follow.
 
\subsubsection{Action-angle variables.\index{Action-angle variables}} 
  According to the decomposition in \eqref{tangential_sites_normal_intro}, we denote
  \begin{align}\label{tangential_normal_intro_3}
  H_S:=\left\{ f : f_j =0, \text{ if $j\notin S$}\right\},\quad H_{S^\perp}:=\left\{ f: f_j=0,\text{ if $j\notin S^\perp$}\right\}, \quad f_j:=\int_{\mathbb{T}}e^{\ii j  x}dx,
  \end{align}
  and we will refer to $H_S$ and $H_{S^\perp}$ as \emph{tangential space} and \emph{normal space}, respectively.
 We introduce the amplitude variable ${\zeta}$:
 \begin{align}\label{amplitude_intro_set}
 \zeta \in [1,2]^{\nu}\subset  \left(\mathbb{R}^+\right)^\nu,
 \end{align}
  and consider a $\zeta$-dependent change of variables\index{$U_\zeta$}, $U_\zeta:\mathbb{T}^\nu\times \mathbb{R}^\nu\times H_{S^\perp}\mapsto L^2(\mathbb{T})$:
 \begin{equation}
  \begin{aligned}\label{intro_change_val}
  &U_\zeta(\theta,y,z):=\epsilon\left( \sum_{j_k\in S}\sqrt{|j_k|(\zeta_k+\epsilon^{2(b-1)}y_k)}e^{\ii(\theta_k + j_k x)} + \epsilon^{(b-1)}z\right), \\& \text{ for some $b\in (1,1+1/12)$},\text{ with }\zeta_{-k}:=\zeta_k,\quad \theta_{-k}:=-\theta_k,\quad y_{-k}:=y_k.
  \end{aligned}
  \end{equation}
 The variables $(\theta,y)$ are the so-called action-angle variables and the above change of variables can be thought of as a reparametrization of functions in $H_{S}$. Also, in order to see the motivation of the constant $b$ in \eqref{intro_change_val}, we note that in  view of \eqref{linear_embedding}, the map 
 \begin{align}\label{linear_intro_embedd2}
 \mathbb{T}^\nu \ni \varphi \mapsto U_\zeta(\varphi,0,0)\in L^2(\mathbb{T})
 \end{align} corresponds to $i_{\text{linear}}$ up to the rescaling factor $\epsilon$, therefore the terms $\epsilon^{2(b-1)}y$ and $\epsilon^{b-1}z$ in \eqref{intro_change_val} can be thought of as  correction terms to solve the nonlinear problem, which justifies the requirement $b>1$. 
 
 Now we define a $\zeta$-dependent Hamiltonian $\mathcal{H}_{\zeta}$ on $\mathbb{T}^\nu\times\mathbb{R}^\nu\times H_{S^\perp}$ as
  \begin{align}\label{newHamiltonian_intro}
 \mathcal{H}_{\zeta}(\theta,y,z):=\epsilon^{-2b}\mathcal{H}\circ U_\zeta(\theta,y,z).
 \end{align}
 We note that the factor $\epsilon^{-2b}$ in \eqref{newHamiltonian_intro} naturally arises in the rescaling of the change of variables to describe the Hamiltonian equation in the new phase space $\mathbb{T}^\nu\times\mathbb{R}^\nu\times H_{S^\perp}$. Indeed, the gSQG dynamics in $L^2(\mathbb{T})$ in \eqref{evol_Req} can be  easily rewritten as an evolution equation in the new phase space by pulling back the vector field by $U_{\zeta}$, and one can obtain the equivalent equation is given by (one can follow the same  computations  given in Chapter~\ref{aavariab} after the proof of Proposition~\ref{aavariable}. Note that the factor $\frac{1}{2\pi}$ is due to our definition for the symplectic form $\sigma$ in \eqref{symplectic} and the gradient in \eqref{sjdsdsdxcxc12s}, but does not play a crucial role throughout the proof)
 \begin{align}\label{new_Hamiltonian_intro12}
\frac{d}{dt} \colvec{\theta(t) \\ y(t)\\ z(t)} = \colvec{\frac{1}{2\pi} \partial_y \mathcal{H}_\zeta(\theta,y,z) \\ - \frac{1}{2\pi} \partial_\theta \mathcal{H}_\zeta(\theta,y,z) \\ \partial_x \left(\nabla_z \mathcal{H}_\zeta(\theta,y,z)\right) },
 \end{align}
 where $\nabla_z \mathcal{H}_\zeta$ is the gradient vector field of $\mathcal{H}_\zeta$ restricted to the subspace $H_{S^\perp}$. Indeed, one can show that if $(\theta(t),y(t),z(t))$ is a solution to \eqref{new_Hamiltonian_intro12}, then $f(t):=U_\zeta(\theta(t),y(t),z(t))$ is a solution to the gSQG equation \eqref{hamiltonian_intro11} (again, see Chapter~\ref{aavariab}). 
  
   \begin{itemize}
  
  \item[$\bullet$] \textbf{A toy model 1: A perturbed Airy equation\index{perturbed Airy equation}.}  
  
   Recall that our goal at this moment is to see whether we can extract a dependence of the linear frequency on $\zeta$. Therefore the question becomes whether the new system \eqref{new_Hamiltonian_intro12} has a linear frequency that can be modulated by $\zeta$. Of course, the answer depends on the structure of $\mathcal{H}$. To this end, let us consider a simpler example, where we can see that the modulation of the linear frequency relies on the quartic homogeneous term of the Hamiltonian. 
   
    As a toy model, let us consider a perturbed Airy equation (see \cite{Baldi-Berti-Montalto:KAM-quasilinear-airy,Baldi-Berti-Montalto:KAM-quasilinear-kdv,Baldi-Berti-Montalto:KAM-quasilinear-mkdv,Giuliani:quasiperiodic-generalized-kdv} for the results of relevant but more complicated models). We define a  Hamiltonian $\mathcal{G}:L^2(\mathbb{T})\mapsto \mathbb{R}\cup\left\{ \infty\right\}$ of the form:    \begin{equation} 
       \begin{aligned}
    &\mathcal{G}(f):=\mathcal{G}_2(f)+\mathcal{G}_{4}(f),\text{ where }\\
    &\mathcal{G}_2(f):=\frac{1}{2}\int_{\mathbb{T}}(\partial_{x}f)^2(x)dx,\\
    & \mathcal{G}_{4}(f):=\sum_{\substack{j_1,j_2,j_3,j_4\in \mathbb{Z}\backslash\left\{ 0 \right\}}}G({j_1,j_2,j_3,j_4})f_{j_1}f_{j_2}f_{j_3}f_{j_4} \text{ for some $G:\mathbb{Z}^4\mapsto \mathbb{C}$}.\label{g4_intro_1}
    \end{aligned}
    \end{equation}  The associated nonlinear Hamiltonian equation to $\mathcal{G}$ is written as (see the comparison with \eqref{Hamil_vec_intro})
    \begin{align}\label{nonliear_Hamiltonian_Gtoy}
    f_t = \partial_x(\nabla_{L^2}\mathcal{G}(f)) = -\partial_{xxx}f + \partial_x(\nabla_{L^2}\mathcal{G}_4(f)).
    \end{align}

       The linearized Hamiltonian equation associated to $\mathcal{G}$ depends on only the quadratic term $\mathcal{G}_2$ and it corresponds to the Airy equation, 
       \begin{align}\label{Airy_toy_1}
       f_t=-\partial_{xxx}f.
       \end{align}
        It is trivial to see that the Airy equation  possesses quasiperiodic solutions with the linear frequency $\omega^{\text{Airy}}\in \mathbb{R}^\nu$ given by
    \begin{align}\label{airy_intro}
  (\overline\omega^{\text{Airy}})_{k}:=j_k^3,\text{ for $S^+:=\left\{ j_1<\cdots < j_\nu\right\}$.}
    \end{align}
    At this point, the linearized equation does not reveal the modulation of the linear frequency by the amplitude. To make the computations easier, let us make the following assumptions on the quartic term $\mathcal{G}_4$:
    \begin{enumerate}[label=(A\arabic*)]
    \item \label{Assumpp1_intro1}${G}(j_1,j_2,j_3,j_4)$ is invariant under any permutation on $\left\{ j_1,j_2,j_3,j_4\right\}$. This assumption is simply to make the computations easier and can be assumed for general quartic Hamiltonian, since we can take the average of the summation in \eqref{g4_intro_1} over all the permutations.
    \item \label{Assumpp1_intro2}${G}$ is supported  only on modes $S$, that is, ${G}(j_1,j_2,j_3,j_4)=0$, if $j_i\notin S$ for some $=1,2,3,4$. This is to focus on the contribution of $\zeta$ through this example, since $\zeta$ presents only in the modes in the set $S$ (see \eqref{intro_change_val}).
    \item \label{Assumpp1_intro3} If $G(j_1,j_2,j_3,j_4)\ne 0$, then $j_a=-j_b$ and $j_c=-j_d$  for a permutation $(a,b,c,d)$ of $(1,2,3,4)$. In other words, there is no nontrivial resonance in $\mathcal{G}_4$. Therefore, using the assumption \ref{Assumpp1_intro1} and \ref{Assumpp1_intro2}, we can define a $\nu\times \nu$ matrix $\mathbb{G}$ as 
    \begin{align}\label{G_symple}
    \mathbb{G}^i_k:=G(j_i,-j_i,j_k,-j_k)\text{ for $i,k=1,\ldots,\nu$ and $j_i,j_k\in S^+$.}
    \end{align}
    \end{enumerate}

 To see the modulation of the frequency by $\zeta$, we compute the composition with the change of variables $U_\zeta$:
  \begin{align}\label{toy_normal_from1}
  \mathcal{G}_\zeta(\theta,y,z):=\epsilon^{-2b}\mathcal{G}\circ U_\zeta(\theta,y,z).
  \end{align}
  For the quadratic term, we see that 
  \begin{align}
  &\epsilon^{-2b}\mathcal{G}_2\circ U_\zeta(\theta,y,z)\nonumber\\ &\overset{\eqref{intro_change_val}}=\epsilon^{-2b}\left(2\pi \sum_{j_k\in S^+}j_k^2\epsilon^2(j_k\zeta_k + \epsilon^{2(b-1)}j_ky_k) \right. +\left.\frac{1}{2}\int_{\mathbb{T}}\epsilon^{2b}(\partial_x z)^2(x)dx\right)\nonumber\\
  &=C_{\epsilon,\zeta} + 2\pi \sum_{j_k\in S^+}(j_k)^3y_k + \frac{1}{2}\int_{\mathbb{T}}(\partial_x z)^2(x)dx\nonumber\\
  &\overset{\eqref{airy_intro}}=C_{\epsilon,\zeta} + 2\pi\left( \overline\omega^{\text{Airy}}\cdot y \right)+  \frac{1}{2}\int_{\mathbb{T}}(\partial_x z)^2(x)dx.\label{quadratic_toy_model1}
  \end{align}
  for some $C_{\epsilon,\zeta}$ that does not depend on $(\theta,y,z)$. For the quartic term, we have
  \begin{align*}
  \mathcal{G}_{4}(f)&\overset{\eqref{g4_intro_1},\ref{Assumpp1_intro2}}= \sum_{\substack{j_1,j_2,j_3,j_4\in S}}G(j_1,j_2,j_3,j_4)f_{j_1}f_{j_2}f_{j_3}f_{j_4} \nonumber\\
  &\overset{\ref{Assumpp1_intro1},\ref{Assumpp1_intro3}}=6\sum_{j_1=j_2\in S^+}G(j_1,-j_1,j_2,-j_2)|f_{j_1}|^2|f_{j_2}|^2\nonumber\\
  &+12\sum_{j_1,j_2\in S^+,\ j_1\ne j_2}G(j_1,j_2,-j_1,-j_2)|f_{j_1}|^2|f_{j_2}|^2\nonumber\\
  &\overset{\eqref{G_symple}}=6\sum_{i=1}^\nu \mathbb{G}^i_i|f_{j_i}|^4 + 12\sum_{i,k=1,\ i\ne k}^{\nu}\mathbb{G}^i_k|f_{j_i}|^2|f_{j_k}|^2.
  \end{align*}
Therefore,  using \eqref{intro_change_val}, we can see that
  \begin{align}\label{quartic_toy_model2}
   & \epsilon^{-2b}\mathcal{G}_4\circ U_\zeta(\theta,y,z)\nonumber \\
   &\overset{\eqref{intro_change_val}}=6\epsilon^{-2b}\sum_{i=1}^\nu \mathbb{G}^i_i\epsilon^4(j_i\zeta_i + \epsilon^{2(b-1)}j_iy_i)^2 \nonumber\\
   &+ 12\epsilon^{-2b}\sum_{i,k=1,\ i\ne k}^\nu\mathbb{G}^i_k \epsilon^4(j_i\zeta_i + \epsilon^{2(b-1)}j_iy_i)(j_k\zeta_k + \epsilon^{2(b-1)}j_ky_k)\nonumber\\
   &=C_{\epsilon,\zeta} + 12\epsilon^2\left(\sum_{i=1}^\nu\mathbb{G}^i_i j_i^2\zeta_iy_i + 2\sum_{i,k=1,\ i\ne k}^\nu\mathbb{G}^i_k j_ij_k\zeta_i y_k\right)\nonumber\\
   & \ + 6\epsilon^{2b}\left(\sum_{i=1}^\nu\mathbb{G}^i_i j_i^2y_i^2 + 2\sum_{i,k=1,\ i\ne k}^\nu\mathbb{G}^i_k j_ij_k y_i y_k\right)\nonumber\\
   &=C_{\epsilon,\zeta} + 12\epsilon^2 G^{\text{mod}}\zeta\cdot y + 6\epsilon^{2b}G^{\text{mod}}y\cdot y,
  \end{align}
  where $G^{\text{mod}}$ is a  $\nu\times\nu$ symmetric matrix defined as
  \begin{align}\label{G_mod_def_intro}
  (G^{\text{mod}})^i_k :=\begin{cases}
   j_i^2\mathbb{G}^{i}_i &\text{ if $i=k$}\\
   2j_ij_k\mathbb{G}^{i}_k & \text{ if $i\ne k$}
   \end{cases} , \text{ for $j_i,j_k\in S^+=\left\{ j_1<\cdots<j_\nu\right\}$.}
  \end{align}
  Plugging  \eqref{quartic_toy_model2} and \eqref{quadratic_toy_model1} into \eqref{toy_normal_from1}, we see a normal form of the nonlinear Hamiltonian $\mathcal{G}_\zeta$:
  \begin{align}\label{normal_form_42}
  \mathcal{G}_\zeta(\theta,y,z) &= C_{\epsilon,\zeta} + \left(2\pi \overline\omega^{\text{Airy}} + 12\epsilon^2G^{\text{mod}}\zeta\right)\cdot y \nonumber\\
  &+ 6\epsilon^{2b}G^{\text{mod}}y\cdot y + \frac{1}{2}\int_{\mathbb{T}}(\partial_x z)^2(x)dx\nonumber\\
  &=C_{\epsilon,\zeta} + 2\pi\omega^{\text{Airy}}(\zeta)\cdot y + 6\epsilon^{2b}G^{\text{mod}}y\cdot y + \frac{1}{2}\int_{\mathbb{T}}(\partial_x z)^2(x)dx,
  \end{align}
  where
  \begin{align}\label{modulation_toymodel}
  \omega^{\text{Airy}}(\zeta):= \overline{\omega}^{\text{Airy}} + \frac{6}{\pi}\epsilon^2 G^{\text{mod}}\zeta \in \mathbb{R}^{\nu}.
  \end{align}
 Recall that we can write the Hamiltonian equation \eqref{nonliear_Hamiltonian_Gtoy} in the phase space $L^2(\mathbb{T})$ as another Hamiltonian equation in the new phase space $\mathbb{T}^\nu\times\mathbb{R}^\nu\times H_{S^\perp}$ by pulling back the vector field by $U_\zeta$, which should be written as (compare to \eqref{new_Hamiltonian_intro12})
 \begin{align}\label{Gtoy_mode_new1}
 \frac{d}{dt}\colvec{\theta(t) \\ y(t) \\ z(t)} = \colvec{\frac{1}{2\pi} \partial_y \mathcal{G}_\zeta(\theta,y,z) \\ - \frac{1}{2\pi} \partial_\theta \mathcal{G}_\zeta(\theta,y,z) \\ \partial_x \left(\nabla_z \mathcal{G}_\zeta(\theta,y,z)\right) }\overset{\eqref{normal_form_42}}=\colvec{\omega^{\text{Airy}}(\zeta) \\ 0 \\ 0 } + \colvec{\frac{6}{\pi}\epsilon^{2b}G^{\text{mod}}y(t) \\ 0 \\ -\partial_{xxx}(z(t))}.
 \end{align}
  From the above equation, we can easily see that the trivial embedding $i_{triv}:\mathbb{T}^\nu\mapsto \mathbb{T}^\nu\times \mathbb{R}^\nu \times H_{S^\perp}$, defined as
  \begin{align}\label{trivial_embedding_intro_toy}
  i_{triv}(\varphi):=(\varphi,0,0).
  \end{align}
  is invariant under the vector field in \eqref{Gtoy_mode_new1},  and the solution to \eqref{Gtoy_mode_new1} can be written as
  \begin{align}\label{trivial_toy_sol}
  (\theta(t),y(t),z(t))=i_{triv}(\omega^{\text{Airy}}(\zeta)t) = (\omega^{\text{Airy}}(\zeta)t,0,0).
  \end{align}
  Clearly, the solution \eqref{trivial_toy_sol} is time-quasiperiodic because each "angular component" $\theta_i(t)$ oscillates with frequency $(\omega^{\text{Airy}}(\zeta))_i$ for each $i=1,\ldots,\nu$ and the frequency vector $\omega^{\text{Airy}}(\zeta)$ is modulated by the amplitude $\zeta$ through the relation in \eqref{modulation_toymodel}, which cannot be observed by just looking at the linear Airy equation \eqref{Airy_toy_1}. Therefore, a quasiperiodic solution to \eqref{nonliear_Hamiltonian_Gtoy} can be obtained as $f(t):=U_\zeta(i_{triv}(\omega^{\text{Airy}}(\zeta)t))$. As shown in this example, our strategy to study the equation \eqref{new_Hamiltonian_intro12} is to derive a "normal form" where we can see a modulation of the frequency by the amplitude $\zeta$ at the linear level of the equation (that is, quadratic level of the Hamiltonian).  
   \end{itemize}
   
  \subsubsection{Weak Birkhoff normal form\index{Weak Birkhoff normal form}.} Our Hamiltonian $\mathcal{H}$ in \eqref{hamiltonian_intro11} does not possess a simple structure as in the toy-model example, therefore it is hopeless to expect $\mathcal{H}_\zeta$ in \eqref{newHamiltonian_intro} to have a simple form as in \eqref{normal_form_42}. However, we will construct a symplectic transformation $\Phi^{WB}:L^2(\mathbb{T})\mapsto L^2(\mathbb{T})$ so that defining another Hamiltonian $H$ as
  \begin{align}\label{H_intro_def_1}
  H(f):=\mathcal{H}\circ \Phi^{WB}(f),
  \end{align}
  we can rewrite $H$, under the composition with $U_\zeta$,  as (compare  below with \eqref{normal_form_42})
  \begin{align}\label{newHamiltonian_int11ro}
  H_\zeta(\theta,y,z)& :=\epsilon^{-2b}H\circ(U_\zeta(\theta,y,z)) \nonumber \\ & = C_{\epsilon,\zeta} + 2\pi \omega^{\text{gSQG}}(\zeta)\cdot y + 6\epsilon^{2b}\mathbb{A}y\cdot y + \mathcal{N}(\theta)(z,z) + P(\theta,y,z),
  \end{align}
  where
  \begin{align}\label{omegagSQG_def_intro}
  \omega^{\text{gSQG}}(\zeta) =\overline{\omega} + \frac{6}{\pi}\epsilon^2 \mathbb{A}\zeta, \text{ (see \eqref{linear_fre_intro2} for the definition of $\overline{\omega}$)},
  \end{align}
  for some $\nu\times \nu$ symmetric matrix $\mathbb{A}$, for some $\theta$-dependent bilinear form $\mathcal{N}(\theta)$ on $H_{S^\perp}$, and for some  perturbation $P$ satisfying some smallness condition.  Note that in \eqref{newHamiltonian_int11ro}, the bilinear form $\mathcal{N}$ and the perturbation $P$ depend on $\zeta$, while we do not denote its dependence to avoid notational complication. 
  
   Certainly, the following concerns need to be taken into account:
  \begin{enumerate}[label=(C\arabic*)]
  \item \label{concern1} What structure of $\omega^{\text{gSQG}}(\zeta)$ do we need? More precisely, what do we require on the matrix $\mathbb{A}$?.
  \item \label{concern2}How to construct $\Phi^{WB}$ so that $H$ defined by \eqref{H_intro_def_1} has the structure in \eqref{newHamiltonian_int11ro} with sufficiently small $P$.
    \end{enumerate}  
  We postpone more detailed comments on the above concerns to the next paragraph  but for now focus on how to transform the functional equation  \eqref{functional_intro_1}, according to  the new Hamiltonian \eqref{newHamiltonian_int11ro}. By requiring $\Phi^{WB}$ to be symplectic, we have that the Hamiltonian equation associated to $\mathcal{H}$ is equivalent to the Hamiltonian equation associated to $H$, therefore, again pulling back the Hamiltonian vector field of $H$  by $U_\zeta$. Thus we are led to study the equation:
  \begin{equation}\label{Hamil_final_intro1}
  \begin{aligned}
    \frac{d}{dt}\colvec{\theta(t)\\y(t)\\z(t)} &=\colvec{\frac{1}{2\pi} \partial_y H_\zeta(\theta,y,z) \\ - \frac{1}{2\pi} \partial_\theta H_\zeta(\theta,y,z) \\ \partial_x \left(\nabla_z H_\zeta(\theta,y,z)\right) }\\
    & \overset{\eqref{newHamiltonian_int11ro}}=\colvec{\omega^{\text{gSQG}}(\zeta) \\ 0 \\ 0 }+\colvec{\frac{6}{\pi}\epsilon^{2b}\mathbb{A}y \\ -\frac{1}{2\pi}\partial_\theta\mathcal{N}(\theta)(z,z) \\ \partial_x((\mathcal{N}(\theta))^{T}[z])} + \colvec{\frac{1}{2\pi} \partial_y P(\theta,y,z) \\ - \frac{1}{2\pi} \partial_\theta P(\theta,y,z) \\ \partial_x \left(\nabla_z P(\theta,y,z)\right)}\\
    &=:X_{H_\zeta}(\theta(t),y(t),z(t)),
  \end{aligned}
  \end{equation}
  where $z\mapsto \mathcal{N}(\theta)^{T}[z]$ is the linear map on $H_{S^\perp}$ such that
  \[
  \int_{\mathbb{T}}\mathcal{N}(\theta)^{T}[z](x)h(x)dx =\nabla_z(\mathcal{N}(\theta)(z,z))[h],\text{ for all $h\in S^\perp$.}
  \]
  If $t\mapsto i_{\infty}(\omega^{\text{gSQG}}(\zeta)t)\in\mathbb{T}^\nu\times\mathbb{R}^\nu\times H_{S^\perp}$ is a quasiperiodic solution to \eqref{Hamil_final_intro1}, then the  quasiperiodic solution to the gSQG equation \eqref{evol_Req} will be recovered by
\[
f(t)=\Phi^{WB}(U_\zeta(i_\infty(\omega^{\text{gSQG}}(\zeta)t))).
\]
Assuming $P\equiv0$, the equation \eqref{Hamil_final_intro1} indeed possesses the trivial embedding\index{trivial embedding} \eqref{trivial_embedding_intro_toy} as a quasiperiodic solution with the frequency vector $\omega^{\text{gSQG}}(\zeta)$. Therefore, our main goal of this monograph becomes to study whether such an embedding can persist under a perturbation $P$ in the system \eqref{Hamil_final_intro1}. Note that taking into account  the dependence of the frequency on $\zeta$, more precise statements to be proved are that "for almost every $\zeta$" in \eqref{amplitude_intro_set}, the quasiperiodic solution with frequency $\omega^{\text{gSQG}}(\zeta)$ can survive under the perturbation.  Making an ansatz,
\begin{align}\label{ansatz_intro_final1}
\text{$t\mapsto i(\omega^{\text{gSQG}}(\zeta)t)$} &  \text{ solves \eqref{Hamil_final_intro1}},\nonumber \\
& \text{ for some $i:\mathbb{T}^\nu\mapsto \mathbb{T}^\nu\times \mathbb{R}^\nu\times H_{S^\perp}$, $i(\varphi)=:(\theta(\varphi),y(\varphi),z(\varphi))$},
\end{align}
we are led to find $i$  such that
\begin{align}\label{functional_intro_4}
\mathcal{F}(i):=\omega^{\text{gSQG}}(\zeta)\cdot\partial_\varphi i(\varphi) - X_{H_\zeta}(i(\varphi)) = 0.
\end{align}
 Now,  let us discuss the concerns \ref{concern1}-\ref{concern2}. 

 \paragraph*{Item \ref{concern1}: Frequency vector $\omega^{\text{gSQG}}(\zeta)$ and  use of $\omega$ as a parameter.} 
Let us first make it clear why we need to care about the structure of $\omega^{\text{gSQG}}(\zeta)$. We recall from \eqref{amplitude_intro_set} that our parameter $\zeta$ lies in a fixed subset $[1,2]^\nu$, and we want to perform the iteration \eqref{refined_sequence_intro2} for sufficiently many $\zeta$ in $[1,2]^\nu$. Denoting
\begin{align}\label{set_fre_intro1}
\Omega:=\left\{ \omega\in \mathbb{R}^\nu: \omega = \omega^{\text{gSQG}}(\zeta),\ \zeta\in [1,2]^\nu\right\},
\end{align}
we can only hope that "almost every" $\omega$ in $\Omega$ satisfy all the necessary non-resonance conditions, such as  \eqref{diophaintime_intro}. This  does not necessarily imply that for "almost every" $\zeta\in[1,2]^\nu$, $\omega^{\text{gSQG}}(\zeta)$ satisfies the necessary non-resonance conditions, especially in case $\zeta\mapsto \omega^{\text{gSQG}}(\zeta)$ is not one-to-one.  Thus, we require that the matrix $\mathbb{A}$ in \eqref{omegagSQG_def_intro} is invertible. The explicit form of $\mathbb{A}$ is not important at this point, but it is important to note that $\mathbb{A}$ is completely determined by the choice of the set $S^+$ in \eqref{intro_setandnu}. The invertibility of $\mathbb{A}$ is one of the "non-resonance conditions\index{non-resonance conditions}" described in the statement of Theorem~\ref{main_2}. In the proof, this condition will be verified (See Section~\ref{rpoisdsd1sd} and Proposition~\ref{freq_amp}).

  Once we have the invertibility of the map $\zeta\mapsto \omega^{\text{gSQG}}(\zeta)$, we will use the frequency $\omega$ as a parameter of the system \eqref{Hamil_final_intro1} and think of $\zeta$ as a quantity determined by $\omega$. More precisely, defining
  \begin{align}\label{omega_set_intro1}
  \Omega_\epsilon:=\left\{ \omega\in \mathbb{R}^\nu: \omega = \overline{\omega} +\frac{6}{\pi}\epsilon^2\mathbb{A}\zeta,\text{ for }\zeta\in [1,2]^{\nu}\right\},
  \end{align}
  we think of the amplitude $\zeta$ to be a function of $\omega$ for $\omega\in \Omega_\epsilon$. With a slight abuse of notation, we will still denote $\zeta$ as if it is an independent variable, but it is actually a function of $\omega$, determined by
  \begin{align}\label{zeta_def_intro12}
  \zeta=\zeta(\omega):=\left(\frac{6}{\pi}\epsilon^2\mathbb{A}\right)^{-1}(\omega-\overline{\omega}),\text{ for $\omega\in \Omega_\epsilon$.}
  \end{align} Then the Hamiltonian functional \eqref{newHamiltonian_int11ro}, the system \eqref{Hamil_final_intro1} and the functional equation \eqref{functional_intro_4} become
 \begin{align}
 H_\zeta(\theta,y,z) &= C_{\epsilon,\zeta} + 2\pi \omega\cdot y + 6\epsilon^{2b}\mathbb{A}y\cdot y + \mathcal{N}(\theta)(z,z) + P(\theta,y,z), \label{actual_system_introHamil_1}\\
 \frac{d}{dt}\colvec{\theta(t)\\y(t)\\z(t)} & =\colvec{\omega \\ 0 \\ 0 }+{\colvec{\frac{6}{\pi}\epsilon^{2b}\mathbb{A}y  \\ -\frac{1}{2\pi}\partial_\theta\mathcal{N}(\theta)(z,z) \\ \partial_x((\mathcal{N}(\theta))^{T}[z])}}+ {\colvec{\frac{1}{2\pi} \partial_y P(\theta,y,z) \\ - \frac{1}{2\pi} \partial_\theta P(\theta,y,z) \\ \partial_x \left(\nabla_z P(\theta,y,z)\right)}}\nonumber\\& =:X_{H_\zeta}(\theta(t),y(t),z(t)),\label{actual_system_intro}
 \end{align}
 and
 \begin{align}
 \mathcal{F}_\omega(i)&:=\omega\cdot\partial_\varphi i(\varphi)- X_{H_{\zeta}}(i(\varphi))=\colvec{\omega\cdot \partial_\varphi \theta(\varphi) \\\omega\cdot \partial_\varphi y(\varphi)\\\omega\cdot  \partial_\varphi z(\varphi) } - X_{H_\zeta}(i(\varphi)) = 0.\label{actual_functioanleq}
 \end{align}
 Here $X_{H_\zeta} = X_{H_{\zeta(\omega)}}$ is now thought of as an $\omega$-dependent vector field (instead of a $\zeta$-dependent vector field) and  the Hamiltonian $H_\zeta$, as well as $\mathcal{N}$ and $P$,  depends on $\omega$ implicitly through \eqref{zeta_def_intro12}. Clearly, if an embedding $i_\infty$ solves \eqref{actual_functioanleq} with some $\omega\in \Omega_\epsilon$, then $i_\infty$ solves \eqref{functional_intro_4} with $\zeta(\omega)$. The reason why we use  $\omega$ as a parameter is that by doing so, it is easier to check the non-resonance conditions such as \eqref{first_menikov_intro1} or \eqref{diophaintime_intro} throughout the proof. 
 
   Now, we transform the initial iteration scheme given in \eqref{newton_intro}, since our new system \eqref{actual_system_intro} has a  slightly different form, compared to \eqref{gSQG_full_intro2}.  The main scheme is quite same as described in \eqref{newton_intro} based on Newton's method (more precisely, Nash-Moser scheme as we discussed before): Noting that the trivial embedding \eqref{trivial_embedding_intro_toy} solves \eqref{actual_system_intro} without the perturbation $P$, we set up a formal sequence of approximate solutions $i_n:\mathbb{T}^\nu\mapsto \mathbb{T}^\nu\times \mathbb{R}^\nu\times H_{S^\perp}$,
 \begin{align}\label{refined_sequence_intro2}
 \begin{cases}
 i_0:=i_{triv},\quad i_{n+1}:=i_{n}-(d_i\mathcal{F}_\omega(i_n))^{-1}[\mathcal{F}_\omega(i_n)],\text{ for $n\ge 0$},\\
 d_i\mathcal{F}_\omega(i)[\ihat]:=\frac{d}{dt}\mathcal{F}_\omega(i+t\ihat)\bigg|_{t=0}\overset{\eqref{actual_functioanleq}}=\omega\cdot\partial_\varphi \ihat - d_iX_{H_\zeta}(i)[\ihat],
 \end{cases}
 \end{align}
 and study the convergence of $i_n$ to a solution to \eqref{actual_functioanleq}. In the iteration scheme\index{iteration scheme} \eqref{refined_sequence_intro2}, we do not expect the inverse of the linearized operator $d_i\mathcal{F}_\omega(i)[\ihat]$ to be obtained for every $\omega\in \Omega_\epsilon$, but we select $\omega$ such that $\omega$ satisfies all the necessary non-resonance conditions to obtain an inverse of the linearized operator. 
 
    Before we close our discussion concerning the use of $\omega$ as a parameter, we emphasize that the constant $\gamma$ arising in the Diophantine condition \eqref{diophaintime_intro} needs to be small depending on $\epsilon$. We wish to select non-resonant frequencies from the set $\Omega_\epsilon$, that is, we wish the set (for some fixed $\tau>0$)
  \begin{align}\label{Resonance_introvector1}
  \Omega_0:=\left\{ \omega\in \Omega_\epsilon: |\omega\cdot l| > \gamma|l|^{\tau},\text{ for all $l\in \mathbb{Z}^\nu$}\right\}
  \end{align}
  to take asymptotically  full measure in $\Omega_\epsilon$.
  However, $\Omega_\epsilon$ in \eqref{omega_set_intro1} is an $\epsilon^2$-neighborhood of $\overline{\omega}$. Therefore, in case $\overline{\omega}$ is resonant, that is, $\overline{\omega}\cdot l_* = 0$ for some $l_*\in \mathbb{Z}^\nu$, we  have
  \[
  |\omega\cdot l_*| \le |(\omega-\overline{\omega})\cdot l_*| + |\overline{\omega}\cdot l_*| \le |\omega-\overline{\omega}||l_*|\lesssim \epsilon^2|l_*|, \text{ for all $\omega\in \Omega_\epsilon$}.
  \]
 Thus, it is not, in general, possible for $\Omega_0$ to obtain asymptotically full measure in $\Omega_\epsilon$, unless $\gamma=o(\epsilon^2)$. For this reason, we will fix $\gamma$ to be
 \begin{align}\label{gamma_impose_intro1}
 \gamma:=\epsilon^{2b},\text{ where $b>1$ is chosen as in \eqref{intro_change_val}.}
 \end{align}
  \paragraph*{Item \ref{concern2}: Construction of $\Phi^{WB}$.} Now, we discuss how to construct the transformation $\Phi^{WB}$ so that we can obtain \eqref{newHamiltonian_int11ro} through \eqref{H_intro_def_1}. 
  
   Before we start, we first fix some notation.   In view of \eqref{tangential_normal_intro_3}, we denote by $v$\index{$v$} and $z$\index{$z$} the variables in spaces $H_S$ and $H_{S^\perp}$ respectively, so that a given $f\in L^2(\mathbb{T})$ can be written as
   \begin{align}\label{f_decomp_intro1}
   f= v+ z,\text{ for some $v\in H_S$ and $z\in H_{S^\perp}$ in a unique way.}
   \end{align}
  The variables  $v$ and $z$ will be called a "tangential variable" and a "normal variable".
   We  define $v_\zeta:\mathbb{T}^\nu\times \mathbb{R}^\nu\mapsto H_S$ by
 \begin{align}\label{v_def_intro_21}
 v_\zeta(\theta,y):=\sum_{j_k\in S}\sqrt{|j_k|(\zeta_k +y_k)}e^{\ii( \theta_k +j_k x)},
 \end{align}
 so that  $U_\zeta$ in \eqref{intro_change_val} can be written as
 \begin{align}\label{Uzeta_rewrite_intro}
 U_\zeta(\theta,y,z) = \epsilon v_{\zeta}(\theta,\epsilon^{2(b-1)}y) + \epsilon^b z=: \epsilon{v_{\epsilon}}(\theta,y) + \epsilon^bz.
 \end{align}

Let us consider a homogeneous expansion\index{homogeneous expansion} of $H$ defined by \eqref{H_intro_def_1}:
  \begin{align}\label{homo_H_intro_2}
 H(f)=H_2(f)+H_3(f)+H_4(f)+H_{5}(f) + H_{\ge 6}(f),
 \end{align}
 where $H_i$ is homogeneous of degree\index{homogeneous of degree $i$} $i$ for $i=1,\ldots 5$, and $H_{\ge6}$ collects all the terms that are homogenous of degree at least $6$.  Also, for each $H_i$ and $0\le k\le i$, we denote by $H_{i,k}$ the term in $H_{i}$ that is homogeneous of degree $k$ in the variable $z$ . For example, recalling the Airy equation in \eqref{g4_intro_1}, we can write 
 \begin{align*}
 \mathcal{G}_2(f)&\overset{\eqref{g4_intro_1}}=\frac{1}{2}\int_{\mathbb{T}}(\partial_x f)^2(x) dx =\frac{1}{2}\int_{\mathbb{T}}(\partial_x v)^2(x) dx+\frac{1}{2}\int_{\mathbb{T}}(\partial_x z)^2(x) dx \\
 & = \mathcal{G}_{2,0}(f) + \mathcal{G}_{2,2}(f).
 \end{align*}
  The reason of introducing the above notation is to see the contribution of $v$ and $z$ in each homogeneous term $H_i$ separately.
  
   Using the above notation, we can rearrange the expansion in \eqref{homo_H_intro_2} as 
  \begin{equation}\label{reaarange_expansion_intro}
  \begin{aligned}  H(f) &=\underbrace{(H_{2,0} + H_{4,0})}_{=:Z_0(f)} + \underbrace{(H_{2,1}+ H_{3,0}+H_{3,1}+H_{4,1} + H_{5,0}+H_{5,1})}_{=:Z_1(f)} \\
  &+ \underbrace{(H_{2,2} + H_{3,2} +H_{4,2})}_{=:Z_2(f)}  + \underbrace{(H_{3,3} + H_{4,3}+H_{4,4} + \sum_{k=2}^5H_{5,k} + H_{\ge 6})}_{=:Z_3(f)}.
  \end{aligned}
  \end{equation}
  Using the change of variables from \eqref{f_decomp_intro1} and \eqref{v_def_intro_21}, and splitting $Z_2(f)=Z_2(v_\zeta(\theta,y)+z)$ into $Z_2(v_\zeta(\theta,0)+z) + (Z_2(v_\zeta(\theta,y)+z)-Z_2(v_\zeta(\theta,0)+z))$, we can rewrite the above as
  \begin{align}\label{invandz_intro1}
  H(v_\zeta(\theta,y)+z) &= Z_0(v_\zeta(\theta,y)) + \mathcal{N}(\theta)(z,z)+ \tilde{P}(v_\zeta(\theta,y)+z),
  \end{align}
  where
  \begin{equation}\label{NandP_intro}
  \begin{aligned}
   \mathcal{N}(\theta)(z,z)&:=Z_2(v_\zeta(\theta,0) + z),\\
    \tilde{P}(v_\zeta(\theta,y)+z)&:=Z_1(v_\zeta(\theta,y)+z) + (Z_2(v_\zeta(\theta,y)+z)-Z_2(v_\zeta(\theta,0)+z)) \\
    &+ Z_3(v_\zeta(\theta,y)+z).
  \end{aligned}
  \end{equation}
  Therefore, using \eqref{Uzeta_rewrite_intro}, and the homogeneity of each $H_{i,k}$ in \eqref{reaarange_expansion_intro}, we obtain
  \begin{align}\label{Hamilton_221introsd1srfotj}
  H_\zeta(\theta,y,z)&\overset{\eqref{newHamiltonian_int11ro}}=\epsilon^{-2b}H(U_\zeta(\theta,y,z))\nonumber\\
  & \overset{\eqref{Uzeta_rewrite_intro}}= \epsilon^{-2b}H(\epsilon v_\zeta(\theta,\epsilon^{2(b-1)}y) +\epsilon^b z)\nonumber\\
  & \overset{\eqref{invandz_intro1}}=\epsilon^{-2b}Z_0(\epsilon v_\zeta(\theta,\epsilon^{2(b-1)}y)) +\epsilon^{-2b}\mathcal{N}(\theta)(\epsilon^bz,\epsilon^bz)  \\
  &+\epsilon^{-2b}\tilde{P}(\epsilon v_\zeta(\theta,\epsilon^{2(b-1)}y)+\epsilon^bz)\nonumber\\
  & = \epsilon^{-2b}Z_0(\epsilon v_\zeta(\theta,\epsilon^{2(b-1)}y))  + \mathcal{N}(\theta)(z,z) \\
  &+\epsilon^{-2b}\tilde{P}(\epsilon v_\zeta(\theta,\epsilon^{2(b-1)}y)+\epsilon^bz),
  \end{align}
  where the last equality follows from the homogeneity of $z\mapsto \mathcal{N}(\theta)(z,z)$, which can be seen from its definition in \eqref{NandP_intro} and the definition of $Z_2$ in \eqref{reaarange_expansion_intro} noticing that $H_{2,2},H_{3,2}$ and $H_{4,2}$ are homogeneous of degree $2$ in the normal variable\index{normal variable} $z$. Comparing the above formula to  \eqref{newHamiltonian_int11ro}, which we aim to obtain, we now see that in order to have a normal form as in \eqref{newHamiltonian_int11ro}, we need to choose $\Phi^{WB}$ in \eqref{H_intro_def_1} so that 
  \begin{align}
   \epsilon^{-2b}Z_0(\epsilon v_\zeta(\theta,\epsilon^{2(b-1)}y)) & = C_{\epsilon,\zeta} + 2\pi \omega^{\text{gSQG}}(\zeta)\cdot y + 6\epsilon^{2b}\mathbb{A}y\cdot y.\label{lieanrfreq1}
    \end{align}
Once $\Phi^{WB}$ is chosen,  the perturbation $P$ in \eqref{newHamiltonian_int11ro} is given by
  \begin{align}\label{Perturbation_afternormal_form_intro}
    P(\theta,y,z) =  \epsilon^{-2b}\tilde{P}(\epsilon v_\zeta(\theta,\epsilon^{2(b-1)}) +\epsilon^bz),
  \end{align}
  which can be seen from the comparison between \eqref{Hamilton_221introsd1srfotj} and \eqref{newHamiltonian_int11ro}.
  In our proof, we will construct $\Phi^{WB}$ so that (see Proposition~\ref{normal_form_prop11231}):
  \begin{enumerate}[label=(WB\arabic*)]
\item \label{WB22_intro123} \eqref{lieanrfreq1} holds. Note that $Z_0$, that is $H_{2,0} + H_{4,0}$ (see \eqref{reaarange_expansion_intro}), is restricted to the tangential space $H_{S}$. As we saw in the example of the Airy equation \eqref{g4_intro_1}, whether \eqref{lieanrfreq1} holds or not depends on the non-existence of nontrivial resonances (see assumption \ref{Assumpp1_intro3} above the equation \eqref{G_symple}).  
\item \label{WB22_intro1} $Z_1(f) \equiv 0$ in \eqref{reaarange_expansion_intro}, therefore, \eqref{reaarange_expansion_intro} reduces to 
\begin{align}
H & = (H_{2,0} + H_{4,0}) + (H_{2,2} + H_{3,2} + H_{4,2}) \nonumber \\
&+ (H_{3,3} + H_{4,3}+H_{4,4} + \sum_{k=2}^5H_{5,k} + H_{\ge 6}). \label{BKIN_intro_32}
\end{align}
We do not write the exact formula here, but compare the structure above to the formulae in Proposition~\ref{normal_form_prop11231}. From \eqref{Perturbation_afternormal_form_intro}, we see that the higher the degree of homogeneity\index{degree of homogeneity} of $\tilde{P}$ is, the smaller $P$ becomes, thanks to the rescaling factor $\epsilon$. Compared to $Z_3$, $Z_1$ has the low homogeneous contribution to $\tilde{P}$ in \eqref{NandP_intro} and the elimination of such low homogeneous contribution yields a sufficient smallness  of $P_\zeta$ to perform  the iteration \eqref{refined_sequence_intro2} (See Lemma~\ref{X_estimate}).
\end{enumerate}
In regards to \ref{WB22_intro1},  indeed, one can formally check whether the perturbative term $P$, defined in \eqref{Perturbation_afternormal_form_intro} and \eqref{NandP_intro} assuming $Z_1\equiv0$, is small enough. To see this,
counting the degree of homogeneity of each term in \eqref{reaarange_expansion_intro}, one can find formally (see \eqref{tame92} for the precise estimates) that:
\begin{align}
|\mathcal{F}(i_{triv})|&=O(\epsilon^{6-2b}),\label{formal_estimate_1intro1}\\
 |d_i^2\mathcal{F}_\omega(i_{triv})[\ihat,\ihat]|&=O(\epsilon|\ihat|^2).\label{formal_estimate_1intro2}
\end{align}
Furthermore, for $\omega\in \Omega_0$ in \eqref{Resonance_introvector1}, the best possible estimate we can expect for $(D_i\mathcal{F}_\omega(i))^{-1}$ would be (similar to \eqref{inverse_estimate_1_intro})
\begin{align}\label{formal_inverse_1intro}
|d_i\mathcal{F}_\omega(i))^{-1}[\ihat]| = O(\gamma^{-1}|\ihat|).
\end{align}
Therefore, the first iteration in \eqref{refined_sequence_intro2} formally gives us that
\begin{align}\label{first_intro_large11}
|i_1 - i_{0}|\overset{\eqref{refined_sequence_intro2}}=|d_i\mathcal{F}_\omega(i_{triv}))^{-1}[\mathcal{F}_\omega(i_{triv})]| \overset{\eqref{formal_inverse_1intro}}\lesssim \gamma^{-1}|\mathcal{F}(i_{triv})|\overset{\eqref{formal_estimate_1intro1}}\lesssim \epsilon^{6-2b}\gamma^{-1},
\end{align}
while
\begin{align}\label{second_value_intro}
|\mathcal{F}_\omega(i_1)|&\le |\mathcal{F}_\omega(i_1)-\mathcal{F}_\omega(i_0)-d_i\mathcal{F}_\omega(i_0)[i_1-i_0]| +|\mathcal{F}_\omega(i_0)+d_i\mathcal{F}_\omega(i_0)[i_1-i_0]|\nonumber\\
&\overset{\eqref{refined_sequence_intro2}}=|\mathcal{F}_\omega(i_1)-\mathcal{F}_\omega(i_0)-d_i\mathcal{F}_\omega(i_0)[i_1-i_0]| \nonumber \\
& + \underbrace{|\mathcal{F}_\omega(i_0)-d_i\mathcal{F}_\omega(i_0)[(d_i\mathcal{F}_\omega(i_0))^{-1}[\mathcal{F}_\omega(i_0)]]}_{=0}\nonumber\\
&\overset{\eqref{formal_estimate_1intro2}}\lesssim \epsilon|i_1-i_0|^2\overset{\eqref{first_intro_large11}}\lesssim \epsilon(\epsilon^{6-2b}\gamma^{-1})^{2}\overset{\eqref{gamma_impose_intro1}} = \epsilon^{13-8b}.
\end{align}
Since $\epsilon^{13-8b}\le \epsilon^{6-2b}$ for $b\in (1,1+1/6)$ (see \eqref{intro_change_val} for the range of $b$), we see that the bound of  $\mathcal{F}_\omega(i_1)$ in \eqref{second_value_intro} has indeed improved compared to $\mathcal{F}_\omega(i_{triv})$ in \eqref{formal_estimate_1intro1}. This formal computation can be thought of as  evidence that $P$ in \eqref{Perturbation_afternormal_form_intro} is small enough for the iteration to close. Furthermore, from the  usual super-exponential convergence rate\index{super-exponential convergence rate} of the Nash-Moser scheme, we expect that each $i_n$ stays close to the trivial embedding, 
\begin{align}
|\mathfrak{I}_n| = O(\epsilon^{6-2b}\gamma^{-1}) \text{ and }|\mathcal{F}_\omega(i_n)| = O(\epsilon^{6-2b}),\text{ for each $n\ge 0$}, \nonumber \\
\text{ where }\mathfrak{I}_n:=i_n - i_{triv}. \label{smallness_intro_4}
\end{align}
Naturally, the above computations are not rigorous at all. However in our proof, we will see that \eqref{smallness_intro_4} actually holds true (see \eqref{rosisidsd22sd}). Therefore, in the further analysis of the invertibility of the linearized operator $\ihat\mapsto d_i\mathcal{F}(i)[\ihat]$, which we will discuss in more details later, we will focus on the embeddings $i$ which are close enough to the trivial embedding $i_{triv}$.
 
Lastly, we note that such a transformation $\Phi^{WB}$ will be constructed (in Chapter~\ref{skpssisodssuw}) by means of  time-$1$ flow maps of auxiliary Hamiltonian\index{auxiliary Hamiltonian} equations, which guarantee that the resulting transformation is  symplectic. We also note that such a procedure only modifies a finite number of Fourier modes, which can be deduced from the fact that $Z_0$ and $Z_1$ in \eqref{reaarange_expansion_intro} involve at most one normal variable. Therefore, the flow maps that will be used to construct $\Phi^{WB}$ can be constructed safely in terms of the well-posedness of the auxiliary Hamiltonian equations.

 \subsubsection{Summary of part 1.} In the above discussion, we explored how to extract an internal parameter yielding a parameter-dependent Hamiltonian $H_\zeta$ in \eqref{actual_system_introHamil_1} with sufficiently small perturbative term $P$. Afterwards, we transformed the search for a quasiperiodic solution to the search for a zero of the functional equation in \eqref{actual_functioanleq}, which can be solved by means of the iterative scheme in \eqref{refined_sequence_intro2}.

\subsection{Sketch of the proof, Part 2: Analysis of the linearized operator}

 In the second part of the sketch of the proof, we discuss in more details how to "invert" the linearized operator arising in the iterative scheme \eqref{refined_sequence_intro2}:
  \begin{align}\label{linearized_operatorintro_iv}
  \ihat\mapsto d_i \mathcal{F}_\omega(i)[\ihat] =  (\omega\cdot\partial_{\varphi} -d_i X_{H_\zeta}(i))[\ihat],
  \end{align} for a fixed embedding  $i:\mathbb{T}^\nu\mapsto \mathbb{T}^\nu\times\mathbb{R}^\nu\times H_{S^\perp}$. Denoting the embedding $i$ as  
  \begin{align}\label{given_embedding}
  i:\varphi\mapsto (\theta(\varphi),y(\varphi),z(\varphi)),
  \end{align} we think of $\theta,y$ as  smooth functions  on $\mathbb{T}^{\nu}$  and $z$ as a smooth function on $\mathbb{T}^\nu\times \mathbb{T}$, such that $z(\varphi)$ for each $\varphi\in\mathbb{T}^\nu$ is restricted to the normal space $H_{S^\perp}$ (see \eqref{tangential_normal_intro_3} for $H_{S^\perp}$). Such functions will be measured in the Sobolev spaces $H^{s}(\mathbb{T}^\nu)\times H^{s}(\mathbb{T}^\nu)\times H^{s}(\mathbb{T}^{\nu+1})$.  In the analysis of the linearized operator $\ihat\mapsto d_i \mathcal{F}_\omega(i)[\ihat]$, we will assume on $i$ that 
  \begin{align}\label{size_assumption_intro1}
  \rVert \mathfrak{I}\rVert^{\Lip(\gamma,\Omega_1)}_{s_0+\mu}\le \mathtt{C}\epsilon^{6-2b}\gamma^{-1},\text{ for some $\Omega_1\subset \Omega_0$, $s_0,\mu,\mathtt{C}>0$}, \text{ where $\mathfrak{I}:=i-i_{triv}$.}
  \end{align}
   Before we discuss the linearized operator, some comments to understand the above assumption regarding the norm $ \rVert \cdot \rVert^{\Lip(\gamma,\Omega_1)}_{s_0+\mu}$,  the set $\Omega_1$, the constants $s_0,\mu,\mathtt{C}>0$ and some motivation of such an assumption will follow: 

\paragraph*{Constants arising in the proof}. We first clarify some constants arising in the proof. Let us  denote\index{$\mathtt{p}$}
 \begin{align}\label{free_parameters_intro}
 \mathtt{p}:=\left\{ (\alpha,\nu,s_0) : \alpha\in (1,2),\ 2\le \nu\in\mathbb{N}, \ s_0\ge \frac{\nu+2}{2}\right\}.
  \end{align}
 The $3$-tuple $\mathtt{p}$ consists of the quantities that we can  freely choose. For instance, if we want to find a quasiperiodic solution to the gSQG equation for $\alpha=3/2$ that can be represented by a $3$-dimensional embedded torus $i_\infty$ (see Definition~\ref{def_quasi_per}) satisfying the Sobolev regularity 
 \[
 \rVert i_{\infty}\rVert_{H^{10}(\mathbb{T}^\nu)\times H^{10}(\mathbb{T}^\nu)\times H^{10}(\mathbb{T}^{\nu+1})}<\infty,
 \]
 then we can just choose $\mathtt{p}=(3/2,3,10)$.  The requirement of $s_0\ge\frac{\nu+2}{2}$ is to guarantee that the corresponding Sobolev space becomes an algebra. Once $\mathtt{p}$ is fixed, then our previous discussion can be summarized as follows: We choose the tangential sites $S^+$ as in \eqref{intro_setandnu} (thus the sets of modes $S$ and $S^\perp$ are fixed accordingly) so that we can derive a weak Birkhoff normal form as in \eqref{newHamiltonian_int11ro}, and we obtain an $\omega$-dependent functional equation \eqref{actual_functioanleq} for $\omega\in \Omega_{\epsilon}$ in \eqref{omega_set_intro1}. Therefore, we are led to study the convergence of the iteration \eqref{refined_sequence_intro2}.  We will also consider the following constants fixed throughout the proof:
 \begin{align}\label{parameter_fixed1}
 b\in (1, 1+1/12),\quad \mathtt{M}\in \mathbb{N}, \quad \tau:=\frac{2}{\alpha-1} + \nu + 2.
 \end{align}
 \begin{enumerate}[label=$\circ$]
 \item Constant $b\in (1,1+1/12)$\index{$b$}: This appears in the introduction of the  action-angle variables in \eqref{intro_change_val}, where we already discussed why $b$ needs to be larger than $1$. For some technical reasons, we will consider $b$ strictly less than $1+1/12$. Note that the constant $b$ also appears in the Diophantine condition to measure the non-resonance of frequency vectors (see \eqref{Resonance_introvector1} and \eqref{gamma_impose_intro1}). 
 
 \item Constant $\mathtt{M}\in\mathbb{N}$\index{$\mathtt{M}$}: As stated in our main theorem (Theorem~\ref{main_2}), we will find solutions that are invariant under  a $\frac{2\pi}{\mathtt{M}}$-translation, as in \eqref{mforl_intro_1} (at the patch level, this corresponds to the invariance under a $\frac{2\pi}{\mathtt{M}}$- rotation). For this purpose, we will choose $S^+$ in \eqref{intro_setandnu} to be multiples of $\mathtt{M}$ (see Remark~\ref{ampli_intro_3} and Proposition~\ref{freq_amp}). We note that $\mathtt{M}$ does  play a crucial role in the convergence of the iteration scheme (its role mainly lies in the rigorous analysis of the non-resonance conditions, see Appendix~\ref{applanxxx}).
 \item Constant $\tau$\index{$\tau$}: The constant $\tau$ arises in  non-resonance conditions, for example in \eqref{Resonance_introvector1}. The motivation is that $\tau$ needs to be large enough to obtain the measure estimate for the non-resonance frequencies (See Proposition~\ref{set_measure_last}).
 \end{enumerate}  
\paragraph*{Size assumption~\eqref{size_assumption_intro1} on the embedding $i$.} Firstly, for the definition of the norm $\rVert \cdot\rVert_s^{\Lip(\gamma,\Omega)}$, we refer to \eqref{gbxx11ssx} and \eqref{omega_dep_norm1}. To understand the assumption \eqref{size_assumption_intro1},  note that even though the iteration \eqref{refined_sequence_intro2} starts with $i_{triv}$\index{$i_{triv}$}, which is independent of the parameter $\omega$, all the other approximate solutions $i_n$\index{$i_n$} for $n>0$ will implicitly depend on $\omega$ since the functional $\mathcal{F}_\omega$ and the linearized operator $\ihat \mapsto D_i\mathcal{F}_\omega(i)[\ihat]$ depend on $\omega$. Hence, the dependence of the embeddings on $\omega$ must be considered as well, and we will measure it in a  Lipschitz way \eqref{omega_dep_norm1}. Also, as we noted, we do not expect that the iteration \eqref{refined_sequence_intro2} runs for every $\omega$. Therefore, at each iteration, we will exclude some resonant frequencies for which the next approximate solution may not be well-defined. Indeed, in our proof, we will have a sequence of nested sets $\mathcal{G}_{n+1}\subset \mathcal{G}_n\subset \Omega_0$ such that $n$-th approximate solution $i_n$ is well-defined for $\omega\in \mathcal{G}_n$ (see item \ref{embeddingsdsd2sd1} of Proposition~\ref{nashmoser2d}). Therefore, to analyze the linearized operator at a given embedding $i$,  we assume that $i$ is defined for $\omega\in \Omega_1$ for some set $\Omega_1\subset \Omega_0$.

Another motivation that we need to keep in mind is  that the approximate solutions in our iteration are expected to stay close to the trivial embedding, in view of our formal computations in \eqref{smallness_intro_4}. For this reason, we introduce the constant $\mathtt{C}$  in the assumption \eqref{size_assumption_intro1}, and presumably, we expect such a constant to depend on only $\mathtt{p},{S^+}$ and the constant $b$ (in \eqref{free_parameters_intro},  \eqref{intro_setandnu} and \eqref{parameter_fixed1} respectively), since these quantities completely determine the functional equation \eqref{actual_functioanleq} (see \eqref{rosisidsd22sd}).

Lastly, $\mu$ in \eqref{size_assumption_intro1} denotes an additional necessary regularity to obtain numerous estimates arising in the analysis of the linearized operator. Roughly speaking, in order to obtain the invertibility or desired estimates in the  analysis of the linearized operator at $i$, we need the embedding $i$ to be more regular than $H^{s_0}$, in which we wish to see a convergence of the approximate solutions. We emphasize that as long as such an additional regularity $\mu$ depends on only $\mathtt{p}$ in \eqref{free_parameters_intro}, which is fixed from the beginning, the usual Nash-Moser scheme tells us that the sequence of approximate solutions in the iteration can be constructed in such a way that they converge in $H^{s_0}$, diverge in a finer space ($H^{\mathtt{S}}$ for some very large $\mathtt{S}\gg s_0$)\index{$\mathtt{S}$}, while stay small in the intermediate space $H^{s_0+\mu}$ (interpolating the low-high norms. See item~\ref{embeddingsdsd2sd1} of  Proposition~\ref{nashmoser2d}). Thus our analysis (especially Chapter~\ref{aprisdinv}-\ref{reduction}) should read as "there exists $\mu(\mathtt{p})>0$ such that if \eqref{size_assumption_intro1} holds for $\mu=\mu(\mathtt{p})$\index{$\mu$}, then all the necessary properties of the linearized operator hold". Afterwards, we will see every approximate solution satisfies such an assumption. See \eqref{rosisidsd22sd}.   The upper bound of $\mu$  is  computable in principle, but we do not do so to avoid additional complexity of the proof. Instead, we will only focus on the non-dependence of $\mu$ on the parameters other than $\mathtt{p}$. We will always assume $\mu$ to be as large as we want, as long as it depends only on $\mathtt{p}$.

\subsubsection{Approximate inverse\index{Approximate inverse}.} In practice, we will not look for the actual inverse, $(d_i\mathcal{F}_\omega(i))^{-1}$. Instead, we will look for an \emph{approximate inverse}, $T(i)$, such that 
\begin{align}\label{approximate_inverse_intro11}
|(d_i\mathcal{F}_\omega(i)\circ T(i) - I)[g]| = O(|\mathcal{F}_\omega(i)|), \text{ where $I$ is the identity operator}.
\end{align}
(See \eqref{inverse_estimate1} for the actual estimate we obtain in our proof). The estimate in \eqref{inverse_estimate1} is more complicated due to the consideration of the dependence on $\omega$\index{$\omega$} and the contribution of the high modes  when $s\gg s_0$, which is split into the ones from from $i_0, Z:=\mathcal{F}_\omega(i_0)$ and $g$).
Although the operator $T(i)$ above is not the actual inverse, we expect that as the approximate solutions in \eqref{refined_sequence_intro2}  approach to the solution $i_\infty$ (that is $\mathcal{F}_\omega(i_n)$ gets smaller in the iteration), the error becomes smaller and smaller, therefore eventually, $T(i_\infty)$ is indeed equal to the actual inverse. Thus, our iteration scheme \eqref{refined_sequence_intro2} is now replaced by
 \begin{align}\label{refined_sequence_intro3}
i_0:=i_{triv},\quad  i_{n+1}:=i_{n}-T(i_n)[\mathcal{F}_\omega(i_n)].
 \end{align}
The error of size $O(|\mathcal{F}_\omega(i)|)$ in \eqref{approximate_inverse_intro11} can be managed in the usual Nash-Moser scheme. We also note that the operator $T(i)$  depends on $\omega$ as well, since so do $\mathcal{F}_\omega$ and $i$. However, we omit the notation of such dependence to avoid notational complexity.

 Now the question is how we can find such an approximate inverse $T(i)$ to run the iteration \eqref{refined_sequence_intro3}. To this end, we follow the theory of Berti and Bolle developed in \cite{Berti-Bolle:nash-moser-kam}, which tells us that an approximate inverse of $\ihat\mapsto d_i \mathcal{F}_\omega(i)[\ihat]$ can be obtained by an approximate inverse in the normal direction.  In order to see this sentence in more detail, let us suppose $i_\infty$ is a solution, $\mathcal{F}_\omega(i_\infty)=0$, and investigate what the linearized operator $d_i\mathcal{F}_\omega(i_\infty)$ looks like, because we can expect that the  different nature between $d_i\mathcal{F}_\omega(i)$ and $d_i\mathcal{F}_\omega(i_\infty)$ can be measured quantitatively by $|\mathcal{F}_\omega(i)|$, which is admissible in view of \eqref{approximate_inverse_intro11}.

The key implication in \cite{Berti-Bolle:nash-moser-kam} is the derivation of a normal form near invariant tori\index{invariant tori}. Indeed \cite[Theorem 1]{Berti-Bolle:nash-moser-kam} tells us that (roughly speaking) if 
\begin{align}\label{assumption_solve_intro2}
\text{$i_\infty$ solves $\mathcal{F}_\omega(i_\infty)=0$},
\end{align} there exists a symplectic diffeomorphism $G:\mathbb{T}^\nu\times \mathbb{R}^\nu\times H_{S^\perp}\mapsto \mathbb{T}^\nu\times \mathbb{R}^\nu\times H_{S^\perp}$ such that $G(i_{triv}(\varphi)) = i_\infty(\varphi)$. Moreover the Hamiltonian $K_\infty$, defined as
 \begin{align}\label{K_infty3def_intro}
 K_\infty(\psi,\eta,w):=H_\zeta\circ G(\psi,\eta,w),\text{ for $(\psi,\eta,w)\in \mathbb{T}^\nu\times \mathbb{R}^\nu\times H_{S^\perp}$},
 \end{align}
 possesses a normal form such that
 \begin{equation}\label{normal_K123}
 \begin{aligned}
  K_\infty(\psi,\eta,w) &=C + 2\pi \omega\cdot \eta \\
  &\ + \frac{1}{2}\partial_{\eta}^2(K_\infty)(\psi,0,0)[\eta,\eta] + \frac{1}{2}\nabla_w^2(K_\infty)(\psi,0,0)[w,w] \\
 & \ + \partial_{\eta} \nabla_w(K_\infty)(\psi,0,0)[\eta,w] + O_3(\eta,w),
 \end{aligned}
 \end{equation} 
 where $\nabla_w^2$ denotes the Hessian operator in the variable $w$ and   $O_3(\eta,w)$ collects all the terms that are homogeneous of degree at least $3$.  This  result implies that the linearized Hamiltonian vector field of $K_\infty$ on the  torus $$\left\{(\psi,0,0)\in \mathbb{T}^\nu\times \mathbb{R}^\nu\times H_{S^\perp}: \psi\in \mathbb{T}^\nu\right\}$$ can be easily computed from \eqref{normal_K123} (since such a linearized vector field does not depend on the cubic contribution of $\eta,w$).  With simple computations, one can obtain that (in the direction $\ihat:=(\hat{\psi},\hat{\eta},\hat{w})$) 
\begin{align}\label{linearized_intro_KHamil}
d_iX_{K_\infty}(i_{triv}(\psi))[\ihat]& :=\frac{d}{dt}\colvec{\frac{1}{2\pi} \partial_\eta K_\infty(i_{triv}(\psi)+t\ihat) \\ - \frac{1}{2\pi} \partial_\psi K_\infty(i_{triv}(\psi)+t\ihat) \\ \partial_x \left(\nabla_w K_\infty(i_{triv}(\psi)+t\ihat)\right)} \bigg|_{t=0} \nonumber\\
&= \colvec{\frac{1}{2\pi}\left( \partial_{\eta}^2(K_\infty)(i_{triv}(\psi))[\hat{\eta}] + (\partial_{\eta}(\nabla_wK_\infty)(i_{triv}(\psi)))^{T}[\hat{w}]\right) \\ 0 \\ \partial_x (\nabla_w^2(K_\infty)(i_{triv}(\psi))[\hat{w}]) + \partial_{\eta}(\nabla_wK_\infty)(i_{triv}(\psi))[\hat{\eta}]}
\end{align}
where $(\partial_{\eta}(\nabla_wK_\infty)(i_{triv}(\psi)))^{T}$ denotes the linear map $H_{S^\perp}\mapsto \mathbb{R}^\nu$ such that for all $\hat{w}\in H_{S^\perp},\hat{\eta}\in\mathbb{R}^\nu$,
\[
(\partial_{\eta}(\nabla_wK_\infty)(i_{triv}(\psi)))^{T}[\hat{w}]\cdot \hat{\eta} = \frac{d}{dt}(\hat{w},\nabla_wK_\infty(\psi,t\hat{\eta},0))_{L^2(\mathbb{T})}|_{t=0}.
\]
Therefore, the linearization\index{linearization} of $i\mapsto \omega\cdot\partial_\varphi i- X_{K_\infty}(i)$ at $i_{triv}$  is given by
\small
\begin{align}\label{Dinfty_intro}
\mathbb{D}_\infty[\ihat]:= \colvec{\omega\cdot\partial_\varphi \hat{\psi}(\varphi) \\ \omega\cdot\partial_\varphi \hat{\eta}(\varphi) \\ \omega\cdot\partial_\varphi \hat{w}(\varphi)} - \colvec{\frac{1}{2\pi}\left( \partial_{\eta}^2(K_\infty)(i_{triv}(\varphi))[\hat{\eta}] + (\partial_{\eta}\nabla(K_\infty)(i_{triv}(\varphi)))^{T}[\hat{w}]\right) \\ 0 \\ \partial_x (\nabla_w^2(K_\infty)(i_{triv}(\varphi))[\hat{w}]) + \partial_{\eta}(\nabla_wK_\infty)(i_{triv}(\varphi))[\hat{\eta}]}.
\end{align}
\normalsize
As one might have already noticed from \eqref{K_infty3def_intro}, the linearized operator, $\ihat\mapsto \mathbb{D}_\infty[\ihat]$, is indeed the "pullback"  of the linearized operator of $i\mapsto \mathcal{F}_\omega(i)$ at $i_\infty$. More precisely, one can obtain (we omit the computations here but one can easily follow the proof of Proposition~\ref{approx_inverse})
\begin{align}\label{Inver_intro31}
d_i\mathcal{F}_\omega(i_\infty)[\ihat] = DG(i_{triv})\circ \mathbb{D}_\infty \circ (DG(i_{triv}))^{-1}[\ihat].
\end{align}
Therefore, the invertibility of $\ihat\mapsto d_i\mathcal{F}_\omega(i_\infty)[\ihat] $ reduces to the invertibility of $\ihat\mapsto \mathbb{D}_\infty[\ihat]$.

From \eqref{refined_sequence_intro3}, recall that we need to find an approximate inverse\index{Approximate inverse} at each $i_n$, which is not necessarily a solution as in \eqref{assumption_solve_intro2}. However the consequence of \cite{Berti-Bolle:nash-moser-kam} even tells us that given $i$, which is not necessarily a solution, there exists a symplectic diffeomorphism $G_\delta:\mathbb{T}^\nu\times \mathbb{R}^\nu\times H_{S^\perp}\mapsto \mathbb{T}^\nu\times \mathbb{R}^\nu\times H_{S^\perp}$ such that defining
\begin{align}
K(\psi,\eta,w)&:=H_\zeta\circ G_\delta(\psi,\eta,w),\label{K_another_Hamiltonian_haaaaaa}\\
\mathbb{D}[\ihat]&:= \colvec{\omega\cdot\partial_\varphi \hat{\psi}(\varphi) \\ \omega\cdot\partial_\varphi \hat{\eta}(\varphi) \\ \omega\cdot\partial_\varphi \hat{w}(\varphi)} \nonumber \\
& - \colvec{\frac{1}{2\pi}\left( \partial_{\eta}^2(K)(i_{triv}(\varphi))[\hat{\eta}] + (\partial_{\eta}(\nabla_wK)(i_{triv}(\varphi)))^{T}[\hat{w}]\right) \\ 0 \\ \partial_x (\nabla_w^2(K)(i_{triv}(\varphi))[\hat{w}]) + \partial_{\eta}(\nabla_wK)(i_{triv}(\varphi))[\hat{\eta}]},\label{K_another_Hamiltonian_haaaaaa2}
\end{align}
it holds that (compare  below with \eqref{Inver_intro31})
\begin{align}\label{difkoosdsdw1_intro_tlqkf}
d_i\mathcal{F}_\omega(i) = DG_\delta(i_{triv})\circ \mathbb{D} \circ (DG_\delta(i_{triv}))^{-1} + O(|\mathcal{F}_\omega(i)|).
\end{align}
 (See \eqref{e3ishere11} and \eqref{esti2m2ate_again_freee1ints1} for the precise result of \eqref{difkoosdsdw1_intro_tlqkf} in our proof).
Hence, if $\ihat\mapsto \mathbb{D}[\ihat]$ is invertible, then it is expected to be an approximate inverse, since the difference from the actual inverse can be quantitatively estimated by the size of $\mathcal{F}_\omega(i)$. This leads us to study the invertibility of $\mathbb{D}$.

 Now, let us see how to achieve the invertibility of the linear operator $\mathbb{D}$ defined in \eqref{K_another_Hamiltonian_haaaaaa2}. To simplify the notation, we denote $K_{ij}(\varphi):=(\partial_\eta)^i(\nabla_w)^jK(\varphi,0,0)$, that is,
 \begin{align}
 \partial_{\eta}^2K(i_{triv}(\varphi))=:K_{20}(\varphi), \quad \partial_{\eta}(\nabla_wK)(i_{triv}(\varphi))=:K_{11}(\varphi), \nonumber \\
 \nabla_{w}^2K(i_{triv}(\varphi))=:K_{02}(\varphi). \label{K_02def_intro_almosttheend}
 \end{align}
Recalling that $K$ maps $\mathbb{T}^\nu\times\mathbb{R}^\nu\times H_{S^\perp}$ to  $\mathbb{R}\cup\left\{ \infty\right\}$, we can think of $K_{20}(\varphi)$ as a $\varphi$-dependent $\nu\times \nu$ symmetric matrix, $K_{11}(\varphi)$ as a $\varphi$-dependent linear operator from $\mathbb{R}^\nu$ to $H_{S^\perp}$, and $K_{02}(\varphi)$ as a $\varphi$-dependent linear operator from $H_{S^\perp}$ to $H_{S^\perp}$. Therefore, given $g$, denoted as 
\begin{align*}\mathbb{T}^\nu\ni \varphi\mapsto g(\varphi)=(g_1(\varphi),g_2(\varphi),g_3(\varphi))\in \mathbb{T}^\nu\times \mathbb{R}^\nu\times H_{S^\perp},\end{align*}
we aim to find $\ihat(\varphi)=(\hat{\psi}(\varphi),\hat{\eta}(\varphi),\hat{w}(\varphi))$ such that $\mathbb{D}[\ihat]=g$, that is,
\begin{subequations}\label{Threeequations_intro3binded}
\begin{align}
\omega\cdot \partial_\varphi\hat{\psi} - \frac{1}{2\pi}\left( K_{20}[\hat{\eta}] +K_{11}^{T}[\hat{w}]\right) &= g_1, \label{Threeequations_intro}\\
\omega\cdot\partial_\varphi\hat{\eta} &= g_2,\label{Threeequations_intro2}\\
\omega \cdot \partial_\varphi\hat{w} - \partial_x (K_{02}[\hat{w}] )-K_{11}[\hat{\eta}] &= g_3.\label{Threeequations_intro3}
\end{align}
\end{subequations}

 We  start with the simplest equation, the second one. One necessary condition for  \eqref{Threeequations_intro2} to have a solution $\hat{\eta}$ is that $g_2$ has zero-average in $\varphi\in \mathbb{T}^\nu$, which can be easily seen from  the elementary Stokes theorem\index{Stokes theorem}. In order to overcome this issue, the authors in \cite{Baldi-Berti-Montalto:KAM-quasilinear-kdv,Berti-Bolle:nash-moser-kam} introduced an extra factor in the Hamiltonian. In our proof, we restrict  the search for an invariant embedding to the search for an invariant \emph{reversible embedding} which ensures that we only need to solve \eqref{Threeequations_intro2} for $g_2$ such that  $g_2(\varphi)=-g_2(-\varphi)$ \color{black} (In the statement of Proposition~\ref{approx_inverse}, we consider the approximate inverse $T$ as a map restricted to the space $Y_i$, which is defined in \eqref{space_reversibles}). In such  a case, the average of $g_2$ is always zero. Therefore \eqref{Threeequations_intro2} has a  solution when $\omega$ is Diophantine. More precisely, if $\omega\in \Omega_0$ (see \eqref{Resonance_introvector1}) and $g_2 = \sum_{l\in\mathbb{Z}^\nu\backslash\left\{ 0 \right\}} \widehat{g_2}(l)e^{\ii l\cdot \varphi}$, then $\hat{\eta}$ is determined by
 \begin{align}\label{g_2inverted_intro}
 \hat{\eta}(\varphi) = \langle \hat{\eta}\rangle + \sum_{l\in\mathbb{Z}^\nu\backslash\left\{ 0 \right\}}\frac{\widehat{g_2}(l)}{\omega\cdot l }e^{\ii l\cdot \varphi} =: \langle \hat{\eta}\rangle + (\omega\cdot\partial_\varphi)^{-1} g_2,
 \end{align}
 where $\langle \hat{\eta}\rangle$ is a constant that does not depend on $\varphi$, which can be freely chosen in view of \eqref{Threeequations_intro2}. 
  Now we move to the third equation \eqref{Threeequations_intro3} and let us denote 
  \begin{align}\label{normal_direction_inverse}
  \mathcal{L}_\omega:=\omega\cdot \partial_{\varphi} - \partial_x K_{02},
  \end{align}
 so that we are led to solve $\mathcal{L}_\omega[\hat{w}] = g_3 + K_{11}[\hat{\eta}]$. Unlike the first two equations \eqref{Threeequations_intro} and \eqref{Threeequations_intro2}, the invertibility of the operator $\mathcal{L}_\omega$ is not simply guaranteed by imposing the Diophantine condition on $\omega$, and its solvability is much more complicated, since $H_{S^\perp}$ is an infinite dimensional space. The invertibility of $\mathcal{L}_\omega$ is our main goal in Chapter~\ref{Linear_op}-\ref{reduction} and we will assume for the moment that (see Proposition~\ref{normal_inversion} for the precise statement with necessary estimates)
\medskip

\begin{minipage}{0.95\textwidth}
\emph{Invertibility in the normal direction:} There exists a set of frequencies $\Omega_\infty(i) \subset \Omega_0$ such that for $\omega \in \Omega_\infty(i)$, the operator $\mathcal{L}_\omega$ is invertible.
\end{minipage}

\medskip
\noindent With the above assumption, we can write a solution $\hat{w}$ in \eqref{Threeequations_intro3} as
\begin{align}\label{solution_normal_31intro}
\hat{w} = \mathcal{L}_\omega^{-1}(g_3 + K_{11}[\hat{\eta}]).
\end{align}
Lastly, having $\hat{\eta}$ and $\hat{w}$ as in \eqref{g_2inverted_intro} and \eqref{solution_normal_31intro}, the first equation \eqref{Threeequations_intro} can be solved as
\begin{align}
\hat{\psi} = (\omega\cdot\partial_\varphi)^{-1}\left(g_1 + \frac{1}{2\pi}\left( K_{20}[\hat{\eta}] +K_{11}^{T}[\hat{w}]\right) \right),
\end{align}
where $\langle \hat{\eta}\rangle$ in \eqref{g_2inverted_intro} must be chosen to guarantee that the $g_1 + \frac{1}{2\pi}\left( K_{20}[\hat{\eta}] +K_{11}^{T}[\hat{w}]\right)$ has zero-average in $\varphi$. Such a choice of $(\hat{\psi},\hat{\eta},\hat{w})$ is  a solution to \eqref{Threeequations_intro3binded}.

The summary of the above discussion  is that the search for an approximate inverse leads us to the investigation of the invertibility of $\mathcal{L}_\omega$ in \eqref{normal_direction_inverse}, which is restricted to the normal direction.  In view of the invertibility assumption, we will find for each $i_n$ in the sequence \eqref{refined_sequence_intro3}, a nested sequence of frequency sets $\mathcal{G}_{n+1}\subset\mathcal{G}_n\subset \Omega_\infty(i_n) \subset \Omega_0$ such that for $\omega\in \mathcal{G}_{n+1}$, the operator $T(i_n)$ is well-defined, therefore so is $i_{n+1}$  (over $\omega \in \mathcal{G}_{n+1}$). This implies that we can  run the iteration for $\omega\in \cap_{n=0}^\infty\mathcal{G}_n$ towards the solution $i_\infty$, which will be rigorously proved in Proposition~\ref{nashmoser2d}. In the rest, we will focus on how to invert $\mathcal{L}_\omega$.

\subsubsection{Structure of $\mathcal{L}_\omega$ at a fixed embedding $i_0$.} Let us fix an embedding $i_0$\index{$i_0$} so that
\begin{align}\label{i_0notation_introlevel1}
i_0(\varphi):=(\theta_0(\varphi), y_0(\varphi),z_0(\varphi)) \text{  satisfies assumption \eqref{size_assumption_intro1} with $\mathfrak{I}_0:=i_0-i_{triv}$,\index{$\mathfrak{I}_0$}}
\end{align} and investigate the invertibility of $\mathcal{L}_\omega$\index{$\mathcal{L}_\omega$} determined by $i_0$ (we use a slight abuse of notation that from now on, $i_0$ is not necessarily the trivial embedding).  A very natural question at this point is {what   the operator $\mathcal{L}_\omega$ does really look like}. We started from the Hamiltonian $\mathcal{H}$ in \eqref{hamiltonian_intro11}, which is already complicated but we have composed it with multiple changes of variables.
To understand the operator $\mathcal{L}_\omega$, let us trace back the compositions. From \eqref{H_intro_def_1}, \eqref{newHamiltonian_int11ro} and \eqref{K_another_Hamiltonian_haaaaaa}, we see that the Hamiltonian $K$ is given by
  \begin{align}\label{def_K_intro_traceback}
  K(\psi,\eta,w) = \epsilon^{-2b}\mathcal{H}\circ \Phi^{WB}\circ U_\zeta \circ G_\delta (\psi,\eta,w),\quad (\psi,\eta,w)\in\mathbb{T}^\nu\times\mathbb{R}^\nu\times H_{S^\perp},
  \end{align} 
  where $\Phi^{WB}$\index{$\Phi^{WB}$} is a  transformation designed to eliminate some homogeneous contributions of the tangential variables (see \ref{WB22_intro123} and \ref{WB22_intro1} in Part 1), $U_\zeta$ is the introduction of the action-angle variables as well as the amplitude variable (see \eqref{intro_change_val}), and $G_\delta$ is introduced in the study of the approximate inverse of the linearized operator.  The linear operator $\mathcal{L}_\omega$ is defined as 
  \begin{align}
  \mathcal{L}_\omega[\hat{w}]:=\omega\cdot \partial_\varphi \hat{w} - \partial_x(\nabla_w^2 K(\varphi,0,0))[\hat{w}]), \nonumber \\
  \text{ for $\hat{w}\in C^\infty(\mathbb{T}^\nu\times \mathbb{T})$ such that $\hat{w}(\varphi)\in H_{S^\perp}$ for each $\varphi\in\mathbb{T}^\nu$,}\label{def_omegaKLintro}
  \end{align}
  which follows from \eqref{normal_direction_inverse} and \eqref{K_02def_intro_almosttheend}. Notice that the operator $\hat{w}\mapsto \nabla_w^2 K(\varphi,0,0))[\hat{w}]$ also depends on the frequency $\omega$ since $K$ in \eqref{def_K_intro_traceback} depends on $\zeta$, while we have the one-to-one correspondence between $\omega$ and $\zeta$ through \eqref{zeta_def_intro12}. 
  
   The earlier discussion tells us that the existence of an approximate inverse at $i_0$, $T(i_0)$, is determined by the invertibility of the operator $\mathcal{L}_\omega$. Then, where is the information of $i_0$  encoded in the expression \eqref{def_omegaKLintro}? Indeed, the symplectic diffeomorphism  $G_\delta$ in \eqref{def_K_intro_traceback} is determined by $i_0$, while we do not explicitly write it in this sketch level discussion. We only note that the image of $i_{triv}$ under $G_\delta$ corresponds to an \emph{isotropic embedding}\index{isotropic embedding} $i_\delta$\index{$i_\delta$} induced from $i_0$ (see \eqref{iso_tropic}, \eqref{def_G_delta} and Lemma~\ref{i_delta_property}), for which only the "action" coordinate differs from $i_0$, that is,
   \begin{align}\label{idlsotrophi_intro_3}
   i_\delta(\varphi):=G_\delta(i_{triv}(\varphi))\overset{\eqref{iso_tropic}}=(\theta_0(\varphi),y_\delta(\varphi),z_0(\varphi)),
   \end{align}
  and $i_\delta$ is also close to the trivial embedding (see \eqref{size_assumption_3}); 
   \begin{align}\label{isotropic_123intro231}
   i_\delta(\varphi)\text{  satisfies the assumption \eqref{size_assumption_intro1} with $\mathfrak{I}_\delta:=i_\delta-i_{triv}$,\index{$\mathfrak{I}_\delta$}}
   \end{align}
 In view of \eqref{def_K_intro_traceback} and \eqref{def_omegaKLintro}, let us take a simple example to see what a Hessian operator looks like on the embedded torus $G_\delta(i_{triv})=i_\delta$.
  \begin{enumerate}[label=$\bullet$]
  \item \textbf{A toy model 2: A perturbed Airy equation.\index{perturbed Airy equation}} Let us denote the orthogonal projections to $H_{S}$ and $H_{S^\perp}$ by (see \eqref{tangential_normal_intro_3}  and \eqref{tangential_sites_normal_intro} for the notations)
  \begin{align}\label{projections}
  \Pi_{S}f = \sum_{j\in S}f_j e^{\ii j x},\quad \Pi_{S^\perp}f:= \sum_{j\in S^\perp}f_j e^{\ii j x},\text{ for $f\in L^2(\mathbb{T})$ with $\int_{\mathbb{T}}f(x)dx = 0$.}
  \end{align}
  We also denote the variables in $H_S$ and $H_{S^\perp}$ by $v$ and $z$ respectively, so that given $f$, we have a unique decomposition as 
  \begin{align}\label{decomp_intro_32123}
  f=v+z,\text{ where $v\in H_S$ and $z\in H_{S^\perp}$ for $f\in L^2(\mathbb{T})$ with $\int_{\mathbb{T}}f(x)dx = 0$.}
  \end{align} 
     Let us consider a Hamiltonian $\mathcal{G}:L^2(\mathbb{T})\mapsto\mathbb{R}\cup\left\{ \infty\right\}$, 
\begin{align}\label{Airy_Hamiltonain_intro2}
\mathcal{G}(f):=\mathcal{G}_{2,2}(f) + \mathcal{G}_{3,2}(f) + \mathcal{G}_{4,2}(f) +\mathcal{G}_{5,2}(f),
\end{align}
where
\begin{equation}\label{Hamiltonian_21231}
\begin{aligned}
&\mathcal{G}_{2,2}(f):=\frac{1}{2}\int_{\mathbb{T}}(\partial_xz)^2dx,\quad \mathcal{G}_{3,2}(f):=\frac{1}{2}\int_{\mathbb{T}}v z^2dx,\\
&\mathcal{G}_{4,2}(f):=\frac{1}{2}\int_{\mathbb{T}}vz\Pi_{S}(vz) + v^2z^2dx,\quad \mathcal{G}_{5,2}(f):=\frac{1}{2}\int_{\mathbb{T}}v^3z^2dx.
\end{aligned}
\end{equation}
 Notice that $\mathcal{G}_{i,k}$ is homogeneous of degree $i$ in the variable $f$ and homogeneous of degree $k$ in the variable $z$, for a fixed $v$. The motivation for the above example is that in our case, we constructed the symplectic transformation $\Phi^{WB}$ in \eqref{H_intro_def_1} so that $H=\mathcal{H}\circ\Phi^{WB}$ reduces to the form of \eqref{BKIN_intro_32}. We designed our toy model \eqref{Airy_Hamiltonain_intro2} to include quadratic contributions of the normal variable $z$. The Hessian operator in the normal direction of $\mathcal{G}$ at $f=v+z$ can be easily computed as
 \begin{align*}
 \nabla_{z}^2\mathcal{G}(f)[\hat{z}] & = -\partial_{xx}\hat{z} +\Pi_{S^\perp}(v\hat{z}) + \Pi_{S^\perp}(v\Pi_{S}(v\hat{z}) + v^2\hat{z})+\Pi_{S^\perp}(v^3\hat{z})\\
 & = \Pi_{S^\perp}\left(-\partial_{xx}\hat{z} + (v + v^2 + v^3)\hat{z}\right) + \Pi_{S^\perp}(v\Pi_S(v\hat{z})).
 \end{align*}
  Therefore, defining
  \begin{align}\label{airy_intro_whatareyoudoing}
  \mathcal{G}_\zeta(\theta,y,z):=\epsilon^{-2b}\mathcal{G}\circ U_\zeta(\theta,y,z),
  \end{align}
  and using the rescaled variables in \eqref{Uzeta_rewrite_intro}, we obtain the Hessian operator of $\mathcal{G}_\zeta(\theta,y,z)$ in the normal direction, 
  \begin{align}
  \nabla_z^2\mathcal{G}_\zeta(\theta,y,z)[\hat{z}]&= \Pi_{S^\perp}\left(-\partial_{xx}\hat{z} + (\epsilon v_\epsilon(\theta,y) + \epsilon^2 v_\epsilon(\theta,y)^2 +\epsilon^3  v_\epsilon(\theta,y)^3)\hat{z}\right) \nonumber \\
  &+\epsilon^2 \Pi_{S^\perp}(v_\epsilon(\theta,y)\Pi_S(v_\epsilon(\theta,y)\hat{z}))\label{normal_driec1nont1}
  \end{align}
  On the embedded torus $i_\delta(\mathbb{T}^\nu)\subset \mathbb{T}^\nu\times \mathbb{R}^\nu\times H_{S^\perp}$, we obtain  a $\varphi$-dependent operator in $H_{S^\perp}$:
  \begin{equation}\label{varphi_dependent_op_intro1}
  \begin{aligned}
  \nabla_z^2\mathcal{G}_\zeta(i_\delta(\varphi))[\hat{z}]&\overset{\eqref{idlsotrophi_intro_3}}= \nabla_z^2\mathcal{G}_\zeta(\theta_0(\varphi),y_\delta(\varphi),z_0(\varphi))[\hat{z}]\\
  & \overset{\eqref{normal_driec1nont1}}= \Pi_{S^\perp}\left( a^{\text{Airy}}_2(\varphi,x)\partial_{xx} \hat{z} +a^{\text{Airy}}_0(\varphi,x)\hat{z} \right) +  R^{\text{Airy}}(\varphi)[\hat{z}] \\& =:\Pi_{S^\perp}M^{\text{Airy}}(\varphi)+R^{\text{Airy}}(\varphi)[\hat{z}],
  \end{aligned}
  \end{equation}
  where the coefficients $a_2,a_0$ are given by
  \begin{align}\label{coefficients_intro}
  a^{\text{Airy}}_2(\varphi,x):= -1,\quad a^{\text{Airy}}_0(\varphi,x):= \sum_{i=1}^{3}\epsilon^i (v_\epsilon(\theta_0(\varphi),y_\delta(\varphi)))^i,
  \end{align}
  and the operator $R^{\text{Airy}}(\varphi)$ is defined as
  \begin{align}\label{finite_dim_op_intro}
  R^{\text{Airy}}(\varphi)[\hat{z}]:= \epsilon^2\Pi_{S^\perp}(v_\epsilon(\theta_0(\varphi),y_\delta(\varphi))\Pi_S(v_\epsilon(\theta_0(\varphi),y_\delta(\varphi))\hat{z})).
  \end{align}
 Note that for each fixed $\varphi$, the operator $\nabla_z^2\mathcal{G}_\zeta(i_\delta(\varphi))$ consists of two types of operators;  the operator $\hat{z}\mapsto  M^{\text{Airy}}(\varphi)[\hat{z}]$ is a pseudo differential operator of order $2$, and $\hat{z}\mapsto R^{\text{Airy}}(\varphi)[\hat{z}]$ is a smoothing operator. Indeed, since $v_\epsilon(\theta_0(\varphi),y_\delta(\varphi))$ is supported on a finite number of Fourier modes, $S$ (see \eqref{v_def_intro_21} and \eqref{Uzeta_rewrite_intro}), only a finite number of modes of $\hat{z}(\varphi)$ is involved in \eqref{finite_dim_op_intro}, therefore $R^{\text{Airy}}(\varphi):H_{S^\perp}\mapsto H_{S^\perp}$ is a smoothing operator. 
 
  We first analyze the coefficients  $a_2^{\text{Airy}},a_0^{\text{Airy}}$, and then study the operator $R^{\text{Airy}}$. From the definition of $v_\epsilon$ in \eqref{Uzeta_rewrite_intro} and \eqref{v_def_intro_21},  we have
  \begin{align}\label{totalloser_inter}
  v_\epsilon(\theta_0(\varphi),y_\delta(\varphi))(x) &= \sum_{j_k\in S}\sqrt{|j_k|(\zeta_k + \epsilon^{2(b-1)}(y_\delta(\varphi))_k)}e^{\ii (\theta_0(\varphi)_k + j_k x)},
  \end{align}
  while each summand can be split as
  \begin{align}\label{monomial_intro1}
  \sqrt{|j_k|(\zeta_k + \epsilon^{2(b-1)}(y_\delta(\varphi))_k)}e^{\ii (\theta_0(\varphi)_k + j_k x)} = \underbrace{\sqrt{|j_k|\zeta_k}e^{\ii (\varphi \cdot \mathtt{l}(j_k) + j_k x)}}_{=:\overline{v}_{j_k}(\varphi,x)} + O(|\mathfrak{I}_\delta|).
    \end{align}
  where $\mathtt{l}(j_k)$ is defind as
 \begin{align}\label{def_lll_intro_lo}
\mathtt{l}(j_k) := \mathtt{e}_{k}, \text{ and }\mathtt{l}(-j_k):=-\mathtt{l}(j_k) \text { for $j_k\in S^{+}=\left\{ j_1,\ldots,j_{\nu}\right\}$}.
\end{align}
Writing 
  \begin{align}\label{simple_monomials_intro32}
  v_\epsilon := v_\epsilon(\theta_0(\varphi),y_\delta(\varphi)),\quad \overline{v} := \sum_{j_k\in S}\overline{v}_{j_k}(\varphi,x),
  \end{align}
  for simplicity, the decomposition \eqref{monomial_intro1} shows that the tangential variable $v_\epsilon$ can be decomposed into a sum of the $i_\delta$-independent monomial ($\overline{v}_{j_k}$) and a small contribution of the size of the embedding $O(|\mathfrak{I}_\delta|) \overset{\eqref{isotropic_123intro231},\eqref{size_assumption_intro1}}=O(\epsilon^{6-2b}\gamma^{-1})$. Hence, from \eqref{coefficients_intro}, we see that the coefficient $a_0^{\text{Airy}}$ consists of at most cubic interactions of the monomials $\overline{v}_{j_k}$ and the small contribution of size $O(\epsilon |\mathfrak{I}_0|)$, that is,
\begin{align}\label{fire_rossgaller_intro}
a_0^{\text{Airy}}(\varphi,x) = \sum_{i=1}^{3}\epsilon^i \left( \sum_{j_{k_1},\ldots,j_{k_i}\in S}\overline{v}_{j_{k_1}}(\varphi,x)\cdots\overline{v}_{j_{k_i}}(\varphi,x) \right) + O(\epsilon |\mathfrak{I}_\delta|).
\end{align}
  Now, we consider the operator $R^{\text{Airy}}$. As above, we can decompose it into 
  \begin{equation}\label{R_airy_finite_intro3}
  \begin{aligned}
  R^{\text{Airy}}(\varphi)[\hat{z}] &= \epsilon^2\Pi_{S^\perp} \left( \overline{v} \Pi_{S}(\overline{v}\hat{z})\right) + R_2^{\text{Airy}},\\
  R^{\text{Airy}}_2(\varphi)[\hat{z}]&:=\epsilon^2\left(\Pi_{S^\perp}(v_\epsilon -\overline{v})\Pi_{S}( v_\epsilon\hat{z})\right) + \epsilon^2 \left(\overline{v}\Pi_S((v_\epsilon - \overline{v})\hat{z})\right).
  \end{aligned}
  \end{equation}
  This shows that $  R^{\text{Airy}}$ can be also decomposed into a sum of an operator that does not depend on the embedding $i_\delta$ and a operator $R^{\text{Airy}}_2$ whose size is $O(\epsilon^2|v_\epsilon -\overline{v}|) = O(\epsilon^2|\mathfrak{I}_\delta|)$, which is much smaller than $O(\epsilon^2)$.
   
   The summary of this example is  that while the Hessian operator of $\mathcal{G}$ at $f=0$  is simply $$\nabla_z^2\mathcal{G}(0)[\hat{z}]=-\partial_{xx}\hat{z},$$ the Hessian operator of rescaled Hamiltonian $\mathcal{G}_\zeta$ on the embedded torus $i_\delta$, $\nabla_z^2\mathcal{G}_\zeta(i_\delta(\varphi))[\hat{z}]$ in \eqref{varphi_dependent_op_intro1} reveals  perturbations in the coefficients and a smoothing operator $R^{\text{Airy}}$. Moreover, the deviation between these two Hessian operators is composed of the contribution of simple monomials $\overline{v}$ in \eqref{simple_monomials_intro32}, and smaller contributions from the embedding $i_\delta$.
  \end{enumerate}
Comparing $K$ in \eqref{def_K_intro_traceback} and $\mathcal{G}_\zeta$ in  \eqref{airy_intro_whatareyoudoing},  we can deduce from the above example that in the case of gSQG the Hessian operator, $\nabla_{w}^2 K(i_{triv}(\varphi))$, will be expressed as a perturbation from the  $\nabla_z^2\mathcal{H}(0)$,
\begin{align}\label{Hsssss1zsintro_ll221}
\nabla_z^2\mathcal{H}(0)[\hat{z}] \overset{\eqref{linear_intro_operator_language1}}= -\frac{1}{2}\Lambda^{\alpha-1}\hat{z} + \frac{T_\alpha}{4}\hat{z}.
\end{align}
Indeed, in Chapter~\ref{Linear_op}, we will find that  $\nabla_{w}^2 K(i_{triv}(\varphi))$ has the form (compare  below to \eqref{varphi_dependent_op_intro1}):
\begin{align}\label{structureof_MitntroM1}
\nabla_w^2K(i_{triv}(\varphi))[\hat{w}]= \Pi_{S^\perp}M(\varphi)[\hat{w}] + R(\varphi)[\hat{w}],
\end{align}
where $M(\varphi)$ is an operator with perturbed coefficients and and $R(\varphi)$ is a smoothing operator. The operator $M(\varphi)$ is of the form:
\begin{align}\label{Mdefinitiointrosn1}
M(\varphi)[\hat{w}]:=\Lambda^{\alpha-1}_{-\frac{1}{2} + a_1(\varphi)}\hat{w} + \left(\frac{T_\alpha}{4}+a_2(\varphi) \right)\hat{w} + \Upsilon^{\alpha-3}_{a_3(\varphi)}\hat{w},
\end{align}
where 
\begin{equation}\label{Lambda_alpha_intro_def}
\begin{aligned}
&\Lambda^{\alpha-1}_{a}h:=\int_{\mathbb{T}}(2-2\cos(x-y))^{-\frac{\alpha}{2}}a(x,y)(h(x)-h(y))dy,\\
&\Upsilon^{\alpha-3}_{a}h := \int_{\mathbb{T}}(2-2\cos(x-y))^{1-\frac{\alpha}{2}}a(x,y)h(y)dy.
\end{aligned}
\end{equation}
The $\varphi$-dependent "coefficients" $a_1(\varphi)=a_1(\varphi,x,y),a_2(\varphi)=a_2(\varphi,x),a_3(\varphi)=a_3(\varphi,x,y)$ and the  operator $R(\varphi)$ are small;
\begin{align}\label{superrough}
a_i = O(\epsilon)\text{ for $i=1,2,3$},\quad R=O(\epsilon).
\end{align}
With \eqref{structureof_MitntroM1}, the definition of $\mathcal{L}_\omega$ in \eqref{def_omegaKLintro} tells us that the linear operator $\mathcal{L}_\omega$ is of the form :
\begin{align}\label{such_intro_inverse11}
\mathcal{L}_\omega = \omega\cdot \partial_\varphi -\Pi_{S^\perp} \partial_x M(\varphi) +\mathcal{R}(\varphi),\text{ where $\mathcal{R}:=-\partial_x R$.}
\end{align}
The derivation of the from of $\mathcal{L}_\omega$ is the main result of Chapter~\ref{Linear_op} and Proposition~\ref{linearized_opspdosodw}. 

 Actually, Proposition~\ref{linearized_opspdosodw} tells us something more. As we observed in the toy example (see \eqref{fire_rossgaller_intro} and \eqref{R_airy_finite_intro3}), we can decompose the contribution to the coefficients and $R$ into the ones from the monomials $\overline{v}$ and the ones from $i_\delta$. This is the motivation for the definition of the class of coefficients $\mathfrak{C}_1(i_0)$ and $\mathfrak{C}_2(i_0)$ in Definition~\ref{homogeneous_expansion1} and the class of a "finite dimensional operator" $\mathfrak{R}(i_0)$ in Definition~\ref{remiander_class12}. In Proposition~\ref{linearized_opspdosodw}, we prove that 
\begin{align}\label{dktlqkfgkrl1p0intro}
a_1,a_3\in\mathfrak{C}_2(i_0),\quad a_2\in\mathfrak{C}_1(i_0),\quad \mathcal{R}\in\mathfrak{R}(i_0).
\end{align} 
Indeed, we will see that the coefficients/operators of size $O(\epsilon^2)$ are not perturbative (see \eqref{sizeofremainder_321}). Therefore, in the reduction procedure, which we will discuss later, we need more precise information about coefficients/operators of sizes $O(\epsilon), O(\epsilon^2)$.  The inclusions \eqref{dktlqkfgkrl1p0intro} tell us that the eigenvalues of the Hessian $\nabla_w^2K(i_{triv}(\varphi))$ are independent of $i_0,i_\delta$ up to $o(\epsilon^2)$, since the contributions  of size $O(\epsilon),O(\epsilon^2)$ are solely determined by the monomials $\overline{v}$ in \eqref{simple_monomials_intro32}.

\subsubsection{KAM reducibility\index{KAM reducibility} and homological equation\index{homological equation}}
The invertibility of $\mathcal{L}_\omega$ is achieved by diagonalizing the operator. Roughly speaking, we aim to find $\varphi$-dependent invertible transformations $\Phi_{1,\infty}(\varphi),\Phi_{2,\infty}(\varphi):H_{S^\perp}\mapsto H_{S^\perp}$\index{$\Phi_{1,\infty}$}\index{$\Phi_{2,\infty}$}  and Fourier multipliers $j\mapsto d_\infty(j)$ such that
\begin{align}\label{diagona_intro_23123}
\Phi_{1,\infty}\circ \mathcal{L}_\omega\circ \Phi_{2,\infty} =\mathcal{L}_\infty:=\omega\cdot \partial_\varphi - \text{diag}_{j\in S^\perp}d_\infty(j).
\end{align}
Note that since the multipliers $d_\infty(j)$ are the eigenvalues of $\hat{w}\mapsto \partial_x(\nabla_w^2 K(\varphi,0,0))[\hat{w}])$\index{$d_{\infty}(j)$} (we remove the dependence on $\varphi$ by reparametrizing in the variable $\varphi$) and  $K$ depends on $\omega$, so does $d_\infty$ as well, therefore, we denote
\begin{align}\label{eigen_dep_omega_intro}
d_\infty(j) = d_\infty(\omega,j),\text{ for $j\in S^\perp$}.
\end{align}
Then,  for the frequencies satisfying the so-called  first order Melnikov condition\index{first order Melnikov condition} (as we saw in \eqref{formal_inverse1_intro}):
\begin{align}\label{first_order_melinintro_3212}
\left\{ \omega : |\ii \omega\cdot l - d_\infty(\omega,j)| \ge \gamma |l|^{-\tau}|j|^{\alpha},\text{ for all $(l,j)\in\mathbb{Z}^\nu\times S^\perp$}\right\},
\end{align}
we can invert the operator in \eqref{diagona_intro_23123}. 

  In order to achieve \eqref{diagona_intro_23123}, let us suppose we have an operator $L$ of the form:
\begin{align}\label{l_0inftronotext123}
L:=\omega\cdot\partial_\varphi - \text{diag}_{j\in S^\perp}d(\omega,j) + \mathcal{R}(\varphi) =:\omega\cdot\partial_\varphi - D + \mathcal{R}(\varphi),
\end{align}
for some Fourier multipliers $j\mapsto d(\omega,j)$ and for some $\varphi$-dependent operator $\mathcal{R}(\varphi):H_{S^\perp}\mapsto H_{S^\perp}$ and investigate how to obtain, from $L$, a fully diagonalized operator such as \eqref{diagona_intro_23123}. We consider a transformation $\phi(\varphi)$ such that
\begin{align}
\phi(\varphi) = I + \psi(\varphi),\text{ for some small operator $\psi(\varphi):H_{S^\perp}\mapsto H_{S^\perp}$,} \nonumber \\ \text{so that $\phi(\varphi)^{-1}=:I+\tilde{\psi}(\varphi)$ exists.} \label{phi_ansats1231_intr2}
\end{align}
We compute the conjugation of $L$ with $\phi$, that is $\phi^{-1}\circ L\circ \phi$. First, we compute
\begin{align}\label{computethehomoho_intro1}
L\circ \phi & = L \circ(I + \psi) \overset{\eqref{l_0inftronotext123}}= L + (\omega\cdot\partial_\varphi \circ \psi- D\circ\psi + \mathcal{R}\circ \psi)\nonumber\\
& = L + \omega\cdot\partial_\varphi(\psi) + \psi\circ (\omega\cdot\partial_\varphi)- D\circ\psi + \mathcal{R}\circ \psi\nonumber\\
& = (\omega\cdot\partial_\varphi - D) + \mathcal{R} + \psi\circ (\omega\cdot\partial_\varphi) - \psi\circ D + \omega\cdot\partial_\varphi(\psi) \nonumber \\ & +( \psi\circ D - D\circ\psi) + \mathcal{R}\circ \psi\nonumber\\
& = (I+\psi)\circ (\omega\cdot\partial_\varphi - D) + (\omega\cdot\partial_\varphi(\psi) + [\psi, D] + \mathcal{R}) + \mathcal{R}\circ \psi.
\end{align}
Denoting 
\begin{align}\label{newdiago_intro13}
r(j):=\mathcal{R}^j_j(0),\quad D_{r}:=\text{diag}_{j\in S^\perp}r(j),
\end{align}where $\mathcal{R}^{j_1}_{j_2}(l)$ is the matrix representation of the operator $\mathcal{R}(\varphi)$ such that
\[
\mathcal{R}(\varphi)[f] = \sum_{(l,j_1)\in\mathbb{Z}^\nu\times\mathbb{Z}}\mathcal{R}^{j_1}_{j_2}(l)f_{j_1}e^{\ii (\varphi\cdot l + j_2 x)},
\]
we solve a homological equation for the operator $\psi$:
\begin{align}\label{homogolog_intro_32}
\omega\cdot\partial_\varphi(\psi) + [\psi, D] + \mathcal{R} = \text{diag}_{j\in S^\perp}r(j)=:D_{r}.
\end{align}
For the solution $\psi$, it follows from \eqref{computethehomoho_intro1} that
\[
L\circ \phi = \phi\circ (\omega\cdot\partial_\varphi - D)  + D_{r} + \mathcal{R}\circ \psi = \phi\circ (\omega\cdot\partial_\varphi - (D+D_r)) -\psi\circ D_r + \mathcal{R}\circ\psi.
\]
Therefore, we see that 
\begin{align}\label{on22222222eite221ration1}
\phi^{-1}\circ L \circ \phi = \omega\cdot\partial_\varphi - \underbrace{(D-D_r)}_{=:D_{\text{new}}}  +\underbrace{ \phi^{-1}\circ\left( -\psi\circ D_r + \mathcal{R}\circ\psi \right)}_{=:\mathcal{R}_{\text{new}}}.
\end{align}
In the above formal computations, we observe that given a linear operator of the form \eqref{l_0inftronotext123}, a solution to the homological equation \eqref{homogolog_intro_32} leads us to an updated linear operator with a new diagonal part and a new remainder part.  If the size of the new remainder $\mathcal{R}_{new}$ is shown to be significantly smaller than the previous remainder $\mathcal{R}$, we can deduce that an infinite number of such conjugations converge to a fully diagonalized operator \eqref{diagona_intro_23123}. The first thing we can observe from the homological equation \eqref{homogolog_intro_32} is that the size of $\psi$ will depend on $\mathcal{R}$, while  we do not expect the solution to exist for every frequency $\omega$ and we need to select $\omega$ that satisfies sufficient non-resonance conditions. From the Diophantine condition~\eqref{diophaintime_intro}, heuristically we can deduce that the size of $\psi$ must be
  \begin{align}\label{size_psi_intro3}
  \psi = O(\gamma^{-1}|\mathcal{R}|).
  \end{align}
Therefore, the new remainder \eqref{on22222222eite221ration1} can be  inferred to be
\begin{align}\label{size_psi_intro3}
  \mathcal{R}_{\text{new}}= O(\gamma^{-1}|\mathcal{R}|^2).
  \end{align}
In order for this new remainder to be smaller than the previous one, $|\mathcal{R}|\gtrsim \gamma^{-1}|\mathcal{R}|^2 $, we must have
\begin{align}\label{sizeofremainder_321}
\mathcal{R}=o(\gamma^{-1})\overset{\eqref{gamma_impose_intro1}}=o(\epsilon^{2b}).
\end{align}
Furthermore, since $\mathcal{R}_{\text{new}}$ is  quadratic in $\mathcal{R}$, we need $\mathcal{R}$ to be bounded in the sense that the composition $\mathcal{R}\circ \mathcal{R}$ does not accumulate the derivatives (for example, compositions of an operator such as $\epsilon^{3}\partial_x$ do not get smaller after the compositions, since its contribution to high modes are  increasing because of the accumulation of the derivatives). 

 In our proof, we will measure the size of the remainder using the notion of "modulo-tame constant, Lip-modulo-tame constant" (see Definition~\ref{def_tame1}, \ref{lip_tamemodulo}). Such notions were introduced in \cite{Berti-Montalto:quasiperiodic-standing-gravity-capillary}. In regard to the above discussion, compositions of $0$-modulo-tame  operators remain as a $0$-modulo-tame operator, which implies that compositions of such operators do not accumulate the derivatives (Lemma~\ref{consdsdpcosdldsx}). We also note that in order to rigorously construct the infinite sequence of conjugations \eqref{on22222222eite221ration1} and  \eqref{phi_ansats1231_intr2}, we need to regularize the remainder $\mathcal{R}(\varphi)$ at each step (regularize with respect to the variable $\varphi$), which requires the estimates of not only the remainder $\mathcal{R}$ but also the derivatives of it with respect to $\varphi$, that is $\partial_{\varphi}^{\mathtt{b}_0}(\mathcal{R}_0)$ for sufficiently large $\mathtt{b}_0>0$. A question how large $\mathtt{b}_0$ needs to be could  be answered by solving \eqref{homogolog_intro_32} rigorously (in our paper, our choice is given in \eqref{nsjdjwdsdnsde}).

 Our operator $\mathcal{L}_\omega$ in \eqref{sizeofremainder_321}  does not satisfy such conditions. Therefore, inspired by \cite{Baldi-Berti-Montalto:KAM-quasilinear-airy}, we will conjugate $\mathcal{L}_\omega$ a finite number of times into the form with a sufficiently small remainder. That is, we will construct invertible  transformations $\Phi_{1-6,1}, \Phi_{1-6,2}$ such that the conjugated operator has the form:
 \begin{align}\label{result_of_resppsdwintro}
 L_0:=\Phi_{1-6,1}\circ \mathcal{L}_\omega \circ \Phi_{1-6,2} =\omega\cdot \partial_\varphi - D_0 + \mathcal{R}_0, 
 \end{align}
 with
 \begin{enumerate}[label=(R\arabic*)]
 \item \label{itntro_poep1} $D_0$ is a Fourier multiplier, that is, $D_0$ does not depend on $\varphi,x$.
 \item \label{itntro_poep12}$\partial_\varphi^{\mathtt{b}_0}\mathcal{R}_0(\varphi),\ \mathcal{R}_0$ satisfies the smallness condition $\mathcal{R}_0 = o(\epsilon^{2b})$ (as a $0$-Lip-modulo-tame operator) for sufficiently large $\mathtt{b}_0>0$. 
 \end{enumerate}
 This is the main task in Subsections~\ref{change_of_the_space_variables}-\ref{taksd2tamesd} yielding  Proposition~\ref{modulut2sosdtame}. After achieving \eqref{result_of_resppsdwintro}, the rigorous iteration to construct $\phi$ in \eqref{on22222222eite221ration1} towards the full diagonal will be proved rigorously throughout Section~\ref{rkppsdsdwdwdx1kaksd}.

\subsubsection{Reduction to a constant coefficient operator: Sections~\ref{change_of_the_space_variables}-\ref{taksd2tamesd}} In view of the above discussion, we are left to study how to construct transformations $\Phi_{1-6,1},\Phi_{1-6,2}$ with which the conjugation of $\mathcal{L}_\omega$ is reduced to the form of \eqref{result_of_resppsdwintro} satisfying \ref{itntro_poep1} and \ref{itntro_poep12}.

 Recall from \eqref{such_intro_inverse11} and \eqref{Mdefinitiointrosn1} that $\mathcal{L}_\omega$ has the form
 \begin{align}\label{first_opertor_intro_harrysocererstone}
 \mathcal{L}_\omega = \omega\cdot \partial_\varphi - \Pi_{S^\perp}\partial_x\left(\Lambda^{\alpha-1}_{-\frac{1}{2} + a_1(\varphi)} + \left(\frac{T_\alpha}{4}+a_2(\varphi) \right)+ \Upsilon^{\alpha-3}_{a_3(\varphi)}\right) + \mathcal{R}(\varphi),
 \end{align}
 for some $a_1(\varphi)=a_1(\varphi,x,y),\ a_2(\varphi)=a_2(\varphi,x)$ and $a_3(\varphi)=a_3(\varphi,x,y)$ such that $a_1,a_2,a_3 = O(\epsilon)$ and for some smoothing operator $\mathcal{R}$ of size $O(\epsilon^2)$. 
 We reduce the coefficients to constants (remove the dependence on $x$ and $\varphi$) from the highest order term to the lower ones, using the conjugation techniques motivated from \cite{Baldi-Berti-Montalto:KAM-quasilinear-airy,Baldi-Berti-Montalto:KAM-quasilinear-kdv,Berti-Montalto:quasiperiodic-standing-gravity-capillary}.
 \begin{enumerate}[label=$\bullet$]
 \item Change of the space variables (Section~\ref{change_of_the_space_variables}):  We conjugate $\mathcal{L}_\omega$ with a $\varphi$-dependent time-$1$ flow map $\Phi_1$ determined by the PDE,
$ \partial_t u = \Pi_{S^\perp}\partial_x\left(\frac{\beta(\varphi,x)}{1+t\beta_x(\varphi,x)}u(t,\varphi,x)\right),\text{ for some function $\beta(\varphi,x)$}$.
  The function $\beta$ is chosen so that the conjugated operator is of the form:
 \[
 \mathcal{L}^1:=(\Phi_1)^{-1}\circ \mathcal{L}_\omega\circ \Phi_1 =\omega\cdot \partial_\varphi - \Pi_{S^\perp}\partial_x \left( b_1(\varphi)\Lambda^{\alpha-1} + b_2(\varphi,x) + \Upsilon^{\alpha-3}_{b_3(\varphi)} \right) + R_1(\varphi),
 \]
 for some  $b_1(\varphi), b_2(\varphi,x), b_3(\varphi,x,y)$, where $b_1$, does not depend on $x$ anymore and $b_2-\frac{T_\alpha}4, b_3=O(\epsilon)$.  ${R}_1(\varphi)$ is a smoothing operator of size $O(\epsilon^2)$. A more precise statement can be found in Proposition~\ref{toohard_2_3}.
 
  \item Reparametrization of time (Section~\ref{reparmedumbple}):  We conjugate $\mathcal{L}^1$ with a transformation $\Phi_2:[h](\varphi,x):=h(\varphi + \omega p_1(\varphi),x)$ for some function $p_1(\varphi)$. The function $p_1$ is chosen so that the conjugated operator is of the form:
 \begin{align}
 \mathcal{L}^2&:=\frac{1}{\rho(\varphi)}(\Phi_2)^{-1}\circ \mathcal{L}^1\circ \Phi_2 =\omega\cdot \partial_\varphi - \Pi_{S^\perp}\partial_x \left( \mathtt{m}_\alpha\Lambda^{\alpha-1} + d_2(\varphi,x) + \Upsilon^{\alpha-3}_{d_3(\varphi)} \right) \nonumber \\
 &+ R_2(\varphi),\label{L2nevil_intro3}
 \end{align}
 for some functions $\rho(\varphi), d_2(\varphi,x), d_3(\varphi,x,y)$ and for a constant $\mathtt{m}_\alpha$ that does not depend on any of $\varphi,x$, and $d_2-\frac{T_\alpha}{4},d_3=O(\epsilon)$.  $R_2(\varphi)$ is  a smoothing operator of size $O(\epsilon^2)$. A more precise statement can be found in Proposition~\ref{prop_time_rep_92}.
  \item Symbolic representation of the operator (Subsection~\ref{symbolic_1231_deulores}):  We rewrite the operator $ \mathtt{m}_\alpha\Lambda^{\alpha-1} + d_2(\varphi,x) + \Upsilon^{\alpha-3}_{d_3(\varphi)}$ in \eqref{L2nevil_intro3} as a pseudo differential operator with associated symbols. Using the Fourier multiplier $m_{1,\alpha}(\xi)$ associated to the operator $\Lambda^{\alpha-1}$, that is, $Op^W(m_{1,\alpha}(\xi)) = \Lambda^{\alpha-1}$ (see {Proposition}~\ref{nichts2}), we rewrite $\mathcal{L}^2$ as 
  \begin{align}
\mathcal{L}^2 = \omega\cdot\partial_\varphi - \Pi_{S^\perp}\partial_x Op^W\left(\mathtt{m}_\alpha m_{1,\alpha}(\xi) +\frac{T_\alpha}4 + \mathfrak{d}_0(\varphi)\right) + R_2(\varphi),
  \end{align}
  for some $\varphi$-dependent classical symbol $\mathfrak{d}_0(\varphi,x,\xi)$ of order $0$ such that $\mathfrak{d}_0=O(\epsilon)$, where $Op^W$ denotes the Weyl quantization of the symbol \eqref{weyl_quant}. A more detailed statement can be found in \eqref{rewrite}.
  
 
  \item Egorov method (step 1) (Subsection~\ref{subsuhamgqiwkinsl1}):  We conjugate $\mathcal{L}^2$ with  $\varphi$-dependent time-$1$ flow map $\Phi_3$ determined by the PDE,
\begin{align}\label{egorov_pdeintro_3}
 \partial_t u =  \Pi_{S^\perp}\partial_x Op^W(\mathfrak{a}(\varphi))[u],
 \end{align}
 for some $\varphi$-dependent symbol function $\mathfrak{a}(\varphi)=\mathfrak{a}(\varphi,x,\xi)$ of order $1-\alpha$.
 The symbol $\mathfrak{a}$ chosen so that the conjugated operator is of the form:
 \begin{align*}
 \mathcal{L}^3&:=(\Phi_3)^{-1}\circ \mathcal{L}^2\circ \Phi_3 \\
 &=\omega\cdot \partial_\varphi - \Pi_{S^\perp}\partial_x Op^W\left(\mathtt{m}_\alpha m_{1,\alpha}(\xi) +\frac{T_\alpha}4 +  \sum_{k=0}^{\mathtt{N}_\alpha}M_x(\mathfrak{d}_k(\varphi)) + \mathfrak{r}_{-2}(\varphi) \right)\\
 &+ \partial_x\Pi_{S^\perp} W_0(\varphi) + R_3(\varphi).
 \end{align*}
 where
 \begin{enumerate}[label=$\circ$]
 \item $\mathfrak{d}_k(\varphi)=\mathfrak{d}_k(\varphi,x,\xi)$are  $\varphi$-dependent symbols of order $k(1-\alpha)$ for $k=1\ldots, \mathtt{N}_\alpha$\index{$\mathtt{N}_\alpha$}. $\mathtt{N}_\alpha$ is a constant that depends on only $\mathtt{p}$ (see \eqref{free_parameters_intro} and \eqref{nsjdjwdsdnsde}). The motivation of $\mathtt{N}_\alpha$ is that this constant is arising in the Taylor expansion of the flow map of \eqref{egorov_pdeintro_3} in the time variable $t$ to obtain necessary estimates.
 \item $M_x(\mathfrak{d}_k(\varphi))$ is the average of the symbol $\mathfrak{d}_k(\varphi,x,\xi)$, therefore the symbol $\sum_{k=0}^{\mathtt{N}_\alpha}M_x(\mathfrak{d}_k(\varphi))$ in the expression of $\mathcal{L}^3$ does not depend on the variable $x$.
 \item $\mathfrak{r}_{-2}(\varphi)=\mathfrak{r}_{-2}(\varphi,x,\xi)$ is a $\varphi$-dependent symbol of order $-2$ such that $\mathfrak{r}_{-2}= O(\epsilon)$.  
 \item $\partial_x \Pi_{S^\perp}W_0(\varphi)$ is a remainder satisfying the requirement \ref{itntro_poep12} for \eqref{result_of_resppsdwintro}.
 \item $R_3$ is a smoothing operator of size $O(\epsilon^2)$.
  \end{enumerate}
  A more detailed statement can be found in Proposition~\ref{induction_egorov}.
\begin{remark}\label{Egorog_32intsrrqqw1} As opposed to the earlier works \cite{Berti-Montalto:quasiperiodic-standing-gravity-capillary,Hassainia-Hmidi-Masmoudi:kam-gsqg}, we choose a slightly more general symbol $\mathfrak{a}(\varphi)$ in \eqref{egorov_pdeintro_3} to eliminate the dependence of the lower order symbols on the variable $x$. This is why we use the Weyl quantization in the reduction procedure. The PDE in \eqref{egorov_pdeintro_3}  may not be well-posed in general, unless the operator $Op^W(\mathfrak{a}(\varphi))$ is a  symmetric operator. The use of the Weyl quantization makes it easier to verify the symmetry of the operator, since  $Op^W(\mathfrak{a}(\varphi))$ is symmetric if and only if $\mathfrak{a}(\varphi)$ is a real-valued symbol. Furthermore, compared to \cite{Berti-Montalto:quasiperiodic-standing-gravity-capillary,Hassainia-Hmidi-Masmoudi:kam-gsqg}, the choice of the symbol $\mathfrak{a}(\varphi)$ is more delicate, since the highest order of our equation is "weak". We overcome this issue by writing $\mathfrak{a}(\varphi)$ as a finite sum of lower order symbols, $\mathfrak{a}(\varphi)=\sum_{k=0}^{\mathtt{N}_\alpha}\mathfrak{a}_{k}(\varphi)$ for some $\mathfrak{a}_k\in\mathcal{S}^{(k+1)(1-\alpha)}$ and search for $\mathfrak{a}_k$ to obtain the desired conjugated operator of the form \eqref{egorov_pdeintro_3}.
  \end{remark}

  \item Egorov method (step 2) (Subsection~\ref{subsuhamgqiwkinsl2}):  We conjugate $\mathcal{L}^3$ with a $\varphi$-dependent time-$1$ flow map $\Phi_4$ determined by the PDEs $ \partial_t u =\Pi_{S^\perp} \partial_x Op^W(\mathfrak{p}_2(\varphi))[u],$ for some $\varphi$-dependent Fourier multiplier $\mathfrak{p}_2(\varphi) = \mathfrak{p}_2(\varphi,\xi)$ of order $1-\alpha$ that does not depend on $x$. 
The function $\mathfrak{p}_2$ is chosen so that the conjugated operator is of the form:
\begin{equation}\label{pppsowpintro3222123}
 \begin{aligned}
 \mathcal{L}^4&:=(\Phi_4)^{-1}\circ \mathcal{L}^3\circ \Phi_4 \\
 &=\omega\cdot \partial_\varphi - \Pi_{S^\perp}\partial_x  Op^W\left(\mathtt{m}_\alpha m_{1,\alpha}(\xi) +\frac{T_\alpha}4  + \mathfrak{m}_{\le0}(\xi) + \mathfrak{r}_{-2,1}(\varphi)\right)\\
 &+ \partial_x\Pi_{S^\perp} W_1(\varphi) + R_4(\varphi),
 \end{aligned}
 \end{equation}
  where
 \begin{enumerate}[label=$\circ$]
 \item $\mathfrak{m}_{\le0}(\xi)$ is a Fourier multiplier of order $0$ that is independent of  $\varphi,x$. 
 \item $\mathfrak{r}_{-2,1}(\varphi)=\mathfrak{r}_{-2,1}(\varphi,x,\xi)$ is a $\varphi$-dependent symbol of order $-2$ such that $\mathfrak{r}_{-2,1}=O(\epsilon)$.
 \item $\partial_x \Pi_{S^\perp}W_1(\varphi)$ is a remainder satisfying the requirement \ref{itntro_poep12} for \eqref{result_of_resppsdwintro}.
 \item $R_4$ is a smoothing operator of size $O(\epsilon^2)$.
  \end{enumerate}
A more detailed statement can be found in Proposition~\ref{rlaqkqdpsanfsdfmf}.
 
  \item Linear Birkhoff normal form (step 1) (Subsection~\ref{firststeplbnf}):  Comparing \eqref{pppsowpintro3222123} to \eqref{result_of_resppsdwintro} and recalling the requirement for the remainder $\mathcal{R}_0=o(\epsilon^{2b})$ from \ref{itntro_poep12},  the operators $Op^W(\mathfrak{r}_{-2,1}(\varphi))$ and $R_4(\varphi)$ are not perturbative, that is, their size is not $o(\epsilon^{2b})$ since $b >1$ (see \eqref{parameter_fixed1}). In order to remove  terms of size $O(\epsilon)$, we conjugate $\mathcal{L}^4$ with a $\varphi$-dependent time-$1$ flow map $\Phi_5$ determined by the PDEs $\frac{d}{dt}u = \partial_x \Pi_{S^\perp}\epsilon Op^W(\rho_1(\varphi))[u]$ for some $\varphi$-dependent symbol $\rho_1(\varphi)=\rho_1(\varphi,x,\xi)$ of order $-1-\alpha$. The function $\rho_1(\varphi,x,\xi)$ is chosen so that the conjugated operator has the form:
 \begin{align*}
 \mathcal{L}^5&=(\Phi_5)^{-1}\circ \mathcal{L}^4\circ \Phi_5\\
 & =\omega\cdot \partial_\varphi  \\
 & -\Pi_{S^\perp}\partial_x Op^W\left(\mathtt{m}_\alpha m_{1,\alpha}(\xi) +\left(\frac{T_\alpha}4  + \mathfrak{m}_{\le0}(\xi)\right) +\epsilon^2\tilde{\mathfrak{b}}_2(\varphi) + {\mathfrak{r}}_{-2,\le 3,*}(\varphi)\right) \\
 & + \partial_x\Pi_{S^\perp}W_2(\varphi) + R_5(\varphi),
 \end{align*} where
 \begin{enumerate}[label=$\circ$]
 \item $\epsilon^2\tilde{\mathfrak{b}}_2(\varphi) = \epsilon^2\tilde{\mathfrak{b}}_2(\varphi,x,\xi)$ is a $\varphi$-dependent symbol of order $-2$ such that $\epsilon^2\tilde{\mathfrak{b}}_2=O(\epsilon^2)$.
 \item $\mathfrak{r}_{-2,\le 3, *}(\varphi)=\mathfrak{r}_{-2,\le 3, *}(\varphi,x,\xi)$ is a $\varphi$-dependent symbol of order $-2$ such that $\mathfrak{r}_{-2,\le 3, *}=o(\epsilon^{2b})$.
 \item $\partial_x \Pi_{S^\perp}W_2(\varphi)$ is a remainder satisfying the requirement \ref{itntro_poep12} for \eqref{result_of_resppsdwintro}.
 \item $R_5$ is a smoothing operator of size $O(\epsilon^2)$.
  \end{enumerate}  
  A more detailed statement can be found in Proposition~\ref{linearstep17}.
  \item Linear Birkhoff normal form (step 2) (Subsection~\ref{firststeplbnf2}): Comparing \eqref{pppsowpintro3222123} to \eqref{result_of_resppsdwintro} and recalling the requirement for the remainder $\mathcal{R}_0=o(\epsilon^{2b})$ from \ref{itntro_poep12},  the operators $Op^W(\epsilon^2\tilde{\mathfrak{b}}_{2}(\varphi))$ and $R_5(\varphi)$ are not perturbative, that is, their size is not $o(\epsilon^{2b})$ since $b >1$ (see \eqref{parameter_fixed1}).  To eliminate the terms of size $O(\epsilon^2)$, we conjugate $\mathcal{L}^5$ with a $\varphi$-dependent time-$1$ flow map $\Phi_6$ determined by the PDEs $\frac{d}{dt}u = \partial_x \Pi_{S^\perp}\epsilon Op^W(\rho_2(\varphi))[u]$ for some $\varphi$-dependent symbol $\rho_2(\varphi)=\rho_2(\varphi,x,\xi)$ of order $-1-\alpha$. The function $\rho_2(\varphi,x,\xi)$ is chosen so that the conjugated operator has the form:
\begin{equation}\label{ubtr_ju_fn11}
 \begin{aligned}
 \mathcal{L}^6&:=(\Phi_6)^{-1}\circ \mathcal{L}_\omega\circ \Phi_6\\
 & =\omega\cdot \partial_\varphi - \Pi_{S^\perp}\partial_x Op^W\left(\mathtt{m}_\alpha m_{1,\alpha}(\xi) + \left(\frac{T_\alpha}4 + \mathfrak{m}_{\le 0}(\xi) + \epsilon^2\mathfrak{m}_{\mathfrak{b}}(\xi) \right) \right)\\
 & \ +\Pi_{S^\perp}\partial_x Op^W(\mathfrak{r}_{-2,\le3,\sharp}(\varphi)) + \partial_x\Pi_{S^\perp} W_3(\varphi) + R_6(\varphi),
 \end{aligned}
 \end{equation}
where
 \begin{enumerate}[label=$\circ$]
 \item $\mathfrak{m}_{\mathfrak{b}}(\xi)$ is a Fourier multiplier of order $-2$ that is independent of $\varphi,x$.
 \item $\mathfrak{r}_{-2,\le3,\sharp}(\varphi)=\mathfrak{r}_{-2,\le3,\sharp}(\varphi,x,\xi)$ is a $\varphi$-dependent symbol of order $-2$ such that $\mathfrak{r}_{-2,\le 3, \sharp}=o(\epsilon^{2b})$.
 \item $\partial_x \Pi_{S^\perp}W_3(\varphi)$ is a remainder satisfying the requirement \ref{itntro_poep12} for \eqref{result_of_resppsdwintro}.
 \item $R_6$ is a smoothing operator of size $o(\epsilon^{2b})$.
  \end{enumerate} 
   A more detailed statement can be found in Proposition~\ref{linearstep172sd}.
   \end{enumerate}
  Finally, denoting
  \begin{align*}
  \mathcal{R}_0&:=\Pi_{S^\perp}\partial_x Op^W(\mathfrak{r}_{-2,\le3,\sharp}(\varphi)) + \partial_x\Pi_{S^\perp} W_3(\varphi) + R_6(\varphi),\\
  D_0&:=\Pi_{S^\perp}\partial_x Op^W\left(\mathtt{m}_\alpha m_{1,\alpha}(\xi) + \left(\frac{T_\alpha}4 + \mathfrak{m}_{\le 0}(\xi) + \epsilon^2\mathfrak{m}_{\mathfrak{b}}(\xi) \right) \right),\\
  L_0&:=\mathcal{L}^6\overset{\eqref{ubtr_ju_fn11}}= \omega\cdot\partial_\varphi - D_0 + \mathcal{R}_0,
  \end{align*}
  we prove in Section~\ref{taksd2tamesd} that $ \mathcal{R}_0$ indeed satisfies the smallness condition \ref{itntro_poep12} and thus achieve the desired reduction to \eqref{result_of_resppsdwintro}.

Lastly, we note that  all transformations, constructed in the reduction  procedure from $\mathcal{L}_\omega$ to $\mathcal{L}^6$, are required to be  1) reversibility preserving, 2) $\frac{2\pi}{\mathtt{M}}$-translation invariance preserving and 3) real (mapping a real-valued function to a real-valued function), to guarantee that the quasiperiodic solution $f(t,x)$ in Theorem~\ref{main_2} is a reversible, $\frac{2\pi}{\mathtt{M}}$-translation invariant and real-valued solution.

\subsection{Summary of the sketch  and the structure of the monograph}\label{roopow02923}
 In summary, we derive a weak Birkhoff normal form from the Hamiltonian of the gSQG equation with two purposes: 1) extraction of an internal parameter and 2) reduction of the size of the perturbative term, inspired by \cite{Baldi-Berti-Montalto:KAM-quasilinear-kdv,Baldi-Berti-Montalto:KAM-quasilinear-mkdv,Feola-Giuliani:quasiperiodic-water-waves}. Once the weak Birkhoff normal from is derived, we aim to perform a Nash-Moser scheme, which leads us to study an approximate inverse of the linearized operator. The search for an approximate inverse reduces to the invertibility of the linearized operator restricted to the normal space in the spirit of  the derivation of the normal form derivation near an invariant torus in \cite{Berti-Bolle:nash-moser-kam}. The invertibility of the linearized operator in the normal direction will be achieved by reduction of the operator to a constant coefficients operator (a diagonalized operator) by means of symplectic transformations studied in \cite{Baldi-Berti-Montalto:KAM-quasilinear-airy,Berti-Montalto:quasiperiodic-standing-gravity-capillary}. Once the reduction is complete, then we finally prove that the approximate solutions constructed by the Nash-Moser iteration converge to the desired quasiperiodic solution. Meanwhile, necessary non-resonance conditions on the frequencies will be taken into account following the strategy in \cite{Giuliani:quasiperiodic-generalized-kdv}. 
 
In this monograph, the above scheme has been structured as follows:
 
Chapter~\ref{prelinsds111} is devoted to the basic definitions of the spaces and operators, and also to the spaces we will work in. In particular, we also analyze the different multipliers that will appear throughout the proof and the properties of the special functions involved.

Chapter~\ref{11jshsdsss2wxowd2} explains the Hamiltonian character of the gSQG equation and recasts it in a way which is compatible with its associated Poisson bracket. We also perform expansions of the Hamiltonian that will later prove useful in the upcoming chapters. Some useful invariance properties of the gSQG equation will be considered as well.

Chapter~\ref{skpssisodssuw} performs the calculation of the weak Birkhoff normal form, reducing the nonlinear interactions in a way that for homogeneous terms of degree $n$, $n \leq 5$, which ensures the sufficient smallness condition of the perturbative term. The nonexistence of  non-trivial resonances in the quartic Hamiltonian will be investigated, which serves as a key ingredient to derive a desired normal form. 

Chapters~\ref{aavariab} and \ref{nonlads} set the problem up into action-angle coordinates and split the nonlinear operator into three different components. We summarize the necessary conditions on the choice of the tangential sites $S$, which will be mainly used in the measure estimate of frequency set. We also state our main theorem there in terms of finding a zero of a nonlinear operator $\mathcal{F}_\omega$.

Chapter~\ref{aprisdinv} reduces the problem of understanding
the linearization of $\mathcal{F}_\omega$ at a given embedding into a linearization of a normal form around the trivial embedding and proving it is invertible as an operator between the corresponding spaces. Furthermore, the system is diagonal and the condition can be further reduced to prove the invertibility of the \textit{normal part} of the operator.

Chapter~\ref{Linear_op} computes a more explicit representation of the aforementioned operator, emphasizing that the calculation is not restricted to the $O(1)$ terms but also to the $O(\epsilon)$ and $O(\epsilon^2)$ as well, since these terms are not perturbative. For this purpose, we will classify the coefficients/smoothing operators arising from the linearized operator as special classes $\mathfrak{C}_1, \mathfrak{C}_2/ \mathfrak{R}$, which reveals that the non-perturbative terms (terms of size $O(\epsilon), O(\epsilon^2)$) are independent of the approximate solutions.

Chapter~\ref{proper_transform_9chapter} collects some useful properties of symplectic transformations that will be used in the reduction process towards a constant coefficients operator. 

In Chapter~\ref{reduction}, we will conjugate the operator via reversible transformations in order to reduce it to a constant-coefficient operator, modulo semilinear and $O(\epsilon^3)$ parts. This is the most important piece of the manuscript and the most demanding one. These transformations are discussed in Chapter~\ref{proper_transform_9chapter}, along with their corresponding spaces. 

Finally, in Chapters~\ref{NSmoes} and \ref{Final_chapter_3} we finalize the proof of the Theorem by using a Nash-Moser argument and computing the estimates on the measure of set of  the non-resonant frequencies.

Appendix \ref{Appeo100021} contains technical lemmas related to the computations of the different pseudo-differential operators. 

 Appendix~\ref{applanxxx} contains a rigorous construction of tangential sites satisfying all the non-resonance conditions described in Section~\ref{rpoisdsd1sd}.

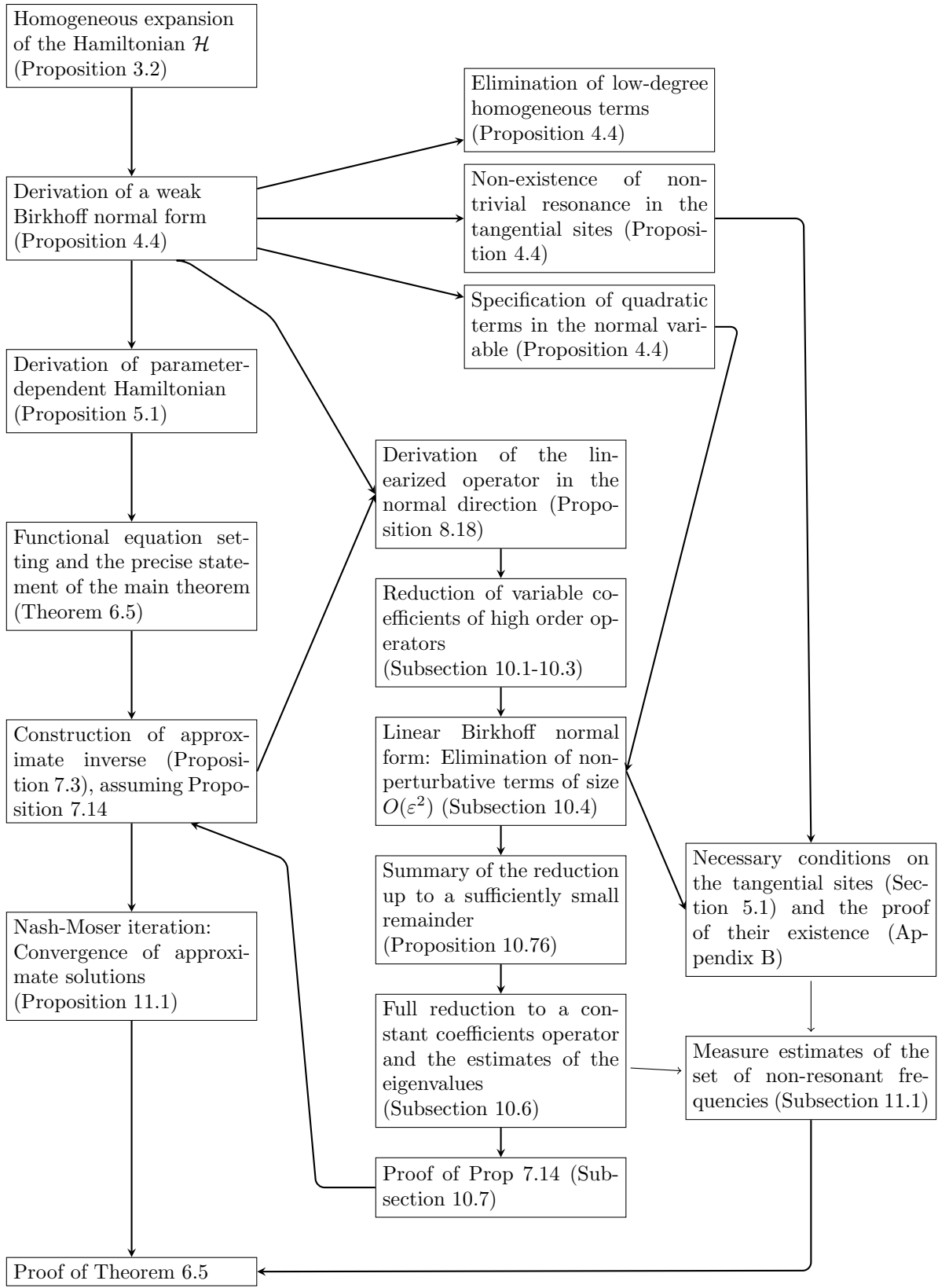
\begin{figure}
    \centering

    \usetikzlibrary {calc,positioning,shapes.misc}
\begin{tikzpicture}[node distance=2cm]
    \node (Hom_gSQG) [rectangle, draw] {\begin{minipage}{4cm} Homogeneous expansion \\ of the Hamiltonian $\mathcal{H}$ \\(Proposition \ref{expansion_1}) \end{minipage}};
    \node (Weak_B) [rectangle, draw, below=1.5cm  of Hom_gSQG] {\begin{minipage}{4cm}Derivation of  a weak \\ Birkhoff normal form \\ (Proposition~\ref{normal_form_prop11231})\end{minipage}};   
    \node (AA_var) [rectangle, draw, below=1.5cm of Weak_B] {\begin{minipage}{4cm} Derivation of parameter-dependent Hamiltonian  \\(Proposition~\ref{aavariable}) \end{minipage}};
      \node (AA_var_1) [rectangle, draw, below=1.5cm of AA_var] {\begin{minipage}{4cm} Functional equation setting and the precise statement of the main theorem (Theorem~\ref{main1})\end{minipage}};
  \node (Approx_inv) [rectangle, draw, below=1.5cm of AA_var_1] {\begin{minipage}{4cm} Construction of approximate inverse (Proposition~\ref{approx_inverse}), assuming Proposition~\ref{normal_inversion} \end{minipage}};
    \node (Nash_moser) [rectangle, draw, below=1.5cm of Approx_inv] {\begin{minipage}{4cm}Nash-Moser iteration:\\ Convergence of approximate solutions \\(Proposition~\ref{nashmoser2d})\end{minipage}};
    \node (Proof_1) [rectangle, draw, below=4cm of Nash_moser] {\begin{minipage}{4cm} Proof of Theorem~\ref{main1} \end{minipage}};
    \node (Weak_B2) [rectangle, draw, right=3.5cm of Weak_B]{\begin{minipage}{4cm} Non-existence of nontrivial  resonance in the tangential sites (Proposition~\ref{normal_form_prop11231})\end{minipage}};
       \node (Weak_B1) [rectangle, draw, above=0.2cm of Weak_B2]{\begin{minipage}{4cm} Elimination of low-degree homogeneous terms \\   (Proposition~\ref{normal_form_prop11231})\end{minipage}}; 
      \node (Weak_B3) [rectangle, draw, below=0.2cm of Weak_B2]{\begin{minipage}{4cm} Specification of quadratic terms in the normal variable (Proposition~\ref{normal_form_prop11231})\end{minipage}};

        \node (Reduction_2) [rectangle, draw, right=2cm of Approx_inv]{\begin{minipage}{4cm}  Linear Birkhoff normal form: Elimination of non-perturbative terms of size $O(\epsilon^2)$ (Section~\ref{linearbbss22}) \end{minipage}};
    \node (Reduction_1) [rectangle, draw, above=0.5cm of Reduction_2]{\begin{minipage}{4cm} Reduction of variable coefficients of high order operators\\ (Section~\ref{change_of_the_space_variables}-\ref{sjdwjjjwdsegosdeoeo})\end{minipage}};
     \node (Reduction_6) [rectangle, draw, above=0.5cm of Reduction_1]{\begin{minipage}{4cm} Derivation of the linearized operator in the normal direction (Proposition~\ref{linearized_opspdosodw})\end{minipage}};
       \node (Reduction_22) [rectangle, draw, below=0.5cm of Reduction_2]{\begin{minipage}{4cm}  Summary of the reduction\\ up to a sufficiently small remainder\\ (Proposition~\ref{modulut2sosdtame}) \end{minipage}};
    \node (Reduction_3) [rectangle, draw, below=0.5cm of Reduction_22]{\begin{minipage}{4cm} Full reduction to a constant coefficients operator and the estimates of the eigenvalues\\ (Section~\ref{rkppsdsdwdwdx1kaksd}) \end{minipage}};

     \node (Reduction_5) [rectangle, draw, below=0.5cm of Reduction_3]{\begin{minipage}{4cm} Proof of Prop 7.14 (Subsection~\ref{loidneks1s2})\end{minipage}};

        \node (Hypo_onS) [rectangle, draw, right=1cm of Reduction_22]{\begin{minipage}{4cm} Necessary conditions on the tangential sites (Section~\ref{rpoisdsd1sd}) and the proof of their existence (Appendix~\ref{applanxxx}) \end{minipage}};
 \node (Measure_estimate) [rectangle, draw, below=1cm of Hypo_onS] {\begin{minipage}{4cm} Measure estimates of the set of non-resonant frequencies (Section~\ref{measure_estimate2ssd2})\end{minipage}};


   \draw[-stealth,thick,rounded corners](Hom_gSQG) -- (Weak_B);
   \draw[-stealth,thick,rounded corners](Weak_B) -- (AA_var);
    \draw[-stealth,thick,rounded corners] (AA_var) -- (AA_var_1);
     \draw[-stealth,thick,rounded corners](AA_var_1) -- (Approx_inv);
      \draw[-stealth,thick,rounded corners] (Approx_inv) -- (Nash_moser);
      \draw[-stealth,thick,rounded corners] (Nash_moser) -- (Proof_1);

  \draw[-stealth,thick,rounded corners](Reduction_6) -- (Reduction_1);
   \draw[-stealth,thick,rounded corners](Reduction_1) -- (Reduction_2);
    \draw[-stealth,thick,rounded corners](Reduction_2) -- (Reduction_22);
     \draw[-stealth,thick,rounded corners] (Reduction_22) -- (Reduction_3);
      \draw[-stealth,thick,rounded corners](Reduction_3) -- (Reduction_5);

      \draw[-stealth,thick,rounded corners](Weak_B) -- (Weak_B2);
      \draw[-stealth,thick,rounded corners](Weak_B) -- (Weak_B1);
      \draw[-stealth,thick,rounded corners](Weak_B) -- (Weak_B3);

      \draw[-stealth,thick,rounded corners](Approx_inv.east) -- (Reduction_6.west);   
      \draw[-stealth,thick,rounded corners] (Weak_B) -- ($(Weak_B.south)+(.9,0)$) -- ($(Reduction_6.west)+(-1.75,3)$) -- (Reduction_6.west);
      
      \draw[-stealth,thick,rounded corners] (Weak_B3.east) -- ($(Weak_B3.east)+(.4,0)$)    --  (Reduction_2.east);

     \draw[-stealth,thick,rounded corners] (Weak_B2.east) --  ($(Weak_B2.east)+(1.5,0)$) -- (Hypo_onS.north);
     
      \draw[-stealth,thick,rounded corners] (Reduction_2.east) -- (Hypo_onS.west);
        \draw[->, shorten >=2pt, shorten <=2pt] (Hypo_onS) -- (Measure_estimate);
         \draw[->, shorten >=2pt, shorten <=2pt] (Reduction_3) -- (Measure_estimate);

     \draw[-stealth,thick,rounded corners] (Measure_estimate.south) --  ($(Measure_estimate.south)+(0,-2.55)$) -- (Proof_1.east);
     
      \draw[-stealth,thick,rounded corners] (Reduction_5.west) -- ($(Reduction_5.west)+(-1,0)$) -- ($(Approx_inv.east)+(0.5,-1.5)$) -- ($(Approx_inv.south) + (1,0)$);
 
    \end{tikzpicture}
    \caption{Connections between  relevant sections and propositions}
    \label{diag2ram_1}
\end{figure}

The relevant connections between sections/propositions in the proof are summarized in Figure~\ref{diag2ram_1}.

\chapter{Preliminaries and notations}\label{prelinsds111}
In this chapter, we recall some classical concepts and specify notations that are used throughout the paper. 
\section{Basic notations}
We denote by $\mathbb{N},\mathbb{Z}$\index{$\mathbb{N}$}\index{$\mathbb{Z}$}, the sets of natural numbers and integers, respectively, and denote $\mathbb{N}_0:=\mathbb{N}\cup \left\{ 0 \right\}$.

 For an object $O$, which possibly is a number, operator or a set, and two real numbers $A,B$, we denote $A \le_O B,$ if there exists a constant $C$ that depends on $O$ such that $A \le CB$. Furthermore, we denote by $C_O,c_O$ or $C(O),c(o)$,  positive constants that depend on $O$, which might vary from line to line in  the computations.
 
 For a linear operator $A$ acting on a vector space $X$, we denote the action of $A$ on $h\in X$ by $Ah$ or $A[h]$, depending on whether a   a clarification is necessary. In case where $X$ is a Hilbert space, we denote by $A^{T}$ be the adjoint operator of $A$.

\section{Function spaces and norms}\label{functionspaces}
\subsection{Usual Sobolev spaces}
We denote by $L^2(\mathbb{T})$ the space of square-integrable $2\pi$-periodic real-valued functions. We denote for $s\in \mathbb{N}$, 
\begin{align}
 L^2_0:=\left\{ u\in L^2(\mathbb{T}): \int_{\mathbb{T}}u(x)dx = 0 \right\}, \nonumber \\ H^s_0:=\left\{ u\in H^s(\mathbb{T}) : \int_{\mathbb{T}}u(x)dx = 0\right\}, \quad C^\infty_0:=\cap_{s\ge 0}H^s_0. 
    \label{rkasd2sdadazzxc}
\end{align}
For $\varphi\in \mathbb{T}^\nu$, we think of a $\varphi$-dependent function $f(\varphi)\in C^\infty_0$ as a function  $f=f(\varphi,x)$ on $\mathbb{T}^{\nu+1}$. 
We denote the Sobolev spaces:\index{$H^s_{\varphi,x}$}\index{$L^2_{\varphi,x}$}\index{$C^\infty_{\varphi,x}$}
\begin{equation}\label{norms22}
\begin{aligned}
&H^s_{\varphi,x}:=H^s(\mathbb{T}^{\nu+1}),\quad H^s_{\varphi}:=H^s(\mathbb{T}^\nu),\quad H^s_x:=H^s(\mathbb{T}),\\
&L^2_{\varphi,x}:=L^2(\mathbb{T}^{\nu+1}),\quad L^2_{\varphi}:=L^2(\mathbb{T}^{\nu}),\quad L^2_{x}:=L^2(\mathbb{T}),\\
& C^\infty_{\varphi,x}:=\cap_{s\ge 0 }H^s_{\varphi,x}, \quad  C^\infty_\varphi := \cap_{s\ge 0}H^s_{\varphi},\quad C^\infty_x:=\cap_{s\ge 0}H^s_x.
\end{aligned}
\end{equation}
For $u\in C^\infty_0$ and $f\in C^\infty_{\varphi,x}$, we denote its Fourier modes by
\begin{align*}
\widehat{u}_j :=\frac{1}{2\pi}\int_{\mathbb{T}}u(x)e^{- \ii j x}dx, \, \text{ and }\, \widehat{f}_j(l) := \frac{1}{(2\pi)^{\nu + 1}}\int_{\mathbb{T}^{\nu}}\int_{\mathbb{T}} f(\varphi,x)e^{-\ii (l\cdot \varphi +  j x)}dxd\varphi.
\end{align*}
so that the Fourier inversion formula reads
\begin{align}\label{sjd2sjjjsd}
u(x)=\sum_{j\in \mathbb{Z}}\widehat{u}_je^{\ii j x}, \quad \text{ and } \, f(\varphi,x) =\sum_{l\in \mathbb{Z}^\nu, j\in \mathbb{Z}}\widehat{f}_j(l)e^{\ii (l\cdot \varphi + j x)}.
\end{align}
We use the norm in  $H^{s}_{\varphi,x}$, defined as \index{$\rVert \cdot\rVert_s$}
\begin{align}\label{gbxx11ssx}
\rVert f\rVert_{s}:=\rVert f \rVert_{H^s_{\varphi,x}}:= \sqrt{\sum_{(l,j)\in \mathbb{Z}^{\nu + 1}}\langle l, j \rangle^{2s}|\widehat{f}_i(l)|^2},\text{ where $\langle l,j \rangle := \sqrt{1+|l|^2+|j|^2)}$.}
\end{align}
Throughout the paper, we consider $H^s_{\varphi,x}$, where  the regularity index $s$ is assumed to be $s\ge s_0$, unless specified, and $s_0$ is fixed so that 
\begin{align}\label{range_of_s}
s\ge s_0 > \frac{\nu+2}{2}.
\end{align}
For such $s,s_0$, the usual Sobolev embedding theorem\index{Sobolev embedding theorem} gives us that
\begin{align}\label{bncal}
\rVert fg\rVert_{s} \le_{s}\rVert f \rVert_{s}\rVert g \rVert_{s_0} + \rVert f \rVert_{s_0}\rVert g \rVert_{s},\text{ for all $f,g\in C^\infty_{\varphi,x}$}.
\end{align}

\begin{lemma}\label{intesdsdsd}\cite[Lemma 2.2]{Berti-Montalto:quasiperiodic-standing-gravity-capillary}
Let $a_0,b_0\ge 0$ and $p,q>0$. For all $\epsilon>0$, there exists a constant $C=C(\epsilon,p,q)>0$ such that 
\[
\rVert u\rVert_{{a_0+p}}\rVert v\rVert_{{b_0+q}} \le \epsilon \rVert u\rVert_{{a_0+p+q}}\rVert v\rVert_{{b_0}} + C\rVert u\rVert_{a_0}\rVert v \rVert_{{b_0+p + q}}, \text{ for all $u,v\in C^\infty_{\varphi,x}$.}
\]
\end{lemma}

\subsubsection{$\omega$-dependent functions}

For a set $\Omega\subset \mathbb{R}^\nu$ and a Hilbert space $(E,\rVert \cdot \rVert_E)$, we consider $E$-valued functions $u=u(\omega)$ that depend on $\omega\in \Omega$. For such functions, we define the norms
\begin{align*}
\rVert u \rVert_{E,\Omega}^{\sup}:=\sup_{\omega\in \Omega} \rVert u(\omega)\rVert_E, \quad \text{ and } \quad \rVert u \rVert_{E,\Omega}^{\lip}:=\sup_{\omega_1,\omega_2\in \Omega, \omega_1\ne \omega_2}\frac{\rVert u(\omega_1)-u(\omega_2)\rVert_{E}}{|\omega_1-\omega_2|}.
\end{align*}
Given $0<\gamma<1$, we denote\index{$\rVert \cdot \rVert_E^{Lip(\gamma,\Omega)}$}
\begin{align}\label{omega_dep_norm1}
\rVert u \rVert_E^{\Lip(\gamma,\Omega)}:= \rVert  u \rVert_{E,\Omega}^{\sup} + \gamma \rVert u \rVert_{E,\Omega}^\lip,\, \Lip(\Omega,E):=\left\{ f:\Omega\mapsto E  \ | \ \rVert f \rVert_{E}^{\Lip(\gamma,\Omega)}< \infty\right\}.
\end{align}
In case $E=\mathbb{R}^n$ for some $n\in\mathbb{N}$, we simply denote
\begin{align}\label{constnas1}
|u|^{\Lip(\gamma,\Omega)} := \sup_{\omega\in \Omega} |u(\omega)| + \gamma \sup_{\omega_1,\omega_2\in \Omega}\frac{|u(\omega_1) - u(\omega_2)|}{|\omega_1-\omega_2|}.
\end{align}
When $E=H^s_{\varphi,x}$, we simply denote\index{$\rVert \cdot \rVert_s^{Lip(\gamma,\Omega)}$}
\[
\rVert u \rVert^{\Lip(\gamma,\Omega)}_s=:\rVert u \rVert^{\Lip(\gamma,\Omega)}_{H^s_{\varphi,x}}.
\]
We have useful estimates for such $\omega$-dependent functions:
\begin{lemma}[Multiplication]\label{interpolation_2s}\cite[Lemma 2.3]{Baldi-Berti-Montalto:KAM-quasilinear-airy}
\[
\rVert uv\rVert_s^{\Lip(\gamma,\Omega)} \le_{s}\rVert u \rVert_s^{\Lip(\gamma,\Omega)}\rVert v \rVert_{s_0}^{\Lip(\gamma,\Omega)}+ \rVert u \rVert_{s_0}^{\Lip(\gamma,\Omega)}\rVert v \rVert_s^{\Lip(\gamma,\Omega)},
\]
for all $f,g\in \Lip(\Omega,H^s_{\varphi,x})$.
\end{lemma}
\begin{lemma}[Interpolation]\label{GNint}
Let $p\ge q \ge 0$. Then, there exists a constant $C=C(s_0,q,p)$ such that
\begin{align*}
\rVert u \rVert_{s_0 + q }^{\Lip(\gamma,\Omega)}\le C\left(\rVert u \rVert^{\Lip(\gamma,\Omega)}_{s_0}\right)^{1-\frac{q}p}\left(\rVert u \rVert_{s_0+p}^{\Lip(\gamma,\Omega)}\right)^{\frac{q}{p}}.
\end{align*}
\end{lemma}
\begin{proof}
Apply the Gagliardo-Nirenberg interpolation inequality to $u(\omega)$ and  $ \displaystyle \gamma\frac{u(\omega_1)-u(\omega_2)}{|\omega_1-\omega_2|}$ for $\omega,\omega_1,\omega_2\in \Omega$ and take the supremum over $\omega,\omega_1,\omega_2$. 
\end{proof}

\begin{lemma}[Change of variables]\label{change_ofsd2}\cite[Lemma 2.4]{Baldi-Berti-Montalto:KAM-quasilinear-kdv}
Let $\beta:\mathbb{T}\mapsto \mathbb{T}$ be such that $\rVert \beta \rVert_{W^{1,\infty}(\mathbb{T})} < \frac{1}2$. Then,
\begin{enumerate}[label=(\arabic*)]
\item $f_{\beta}: x\mapsto x+\beta(x)$ is invertible, and there exists $\tilde{\beta}=\tilde{\beta}(x)$ such that $f_\beta^{-1}(x)=x+\tilde{\beta}(x)$.
\item If $\beta$ depends on $\omega\in \Omega$ and $\varphi\in \mathbb{T}^\nu$, it holds that $\rVert \tilde{\beta} \rVert^{\Lip(\gamma,\Omega)}_{W^{s,\infty}(\mathbb{T}^{\nu+1})}\le_s \rVert \beta \rVert^{\Lip(\gamma,\Omega_1)}_{W^{s+1,\infty}(\mathbb{T}^{\nu+1})}$.
\item \label{sd2sd} Given $u\in C^\infty_{\varphi,x}$,  it holds that  denoting $(u\circ f_\beta)(\varphi,x) :=u(\varphi,x + \beta(\varphi,x))$,
\[
\rVert u\circ f_\beta \rVert_{s}^{\Lip(\gamma,\Omega)}\le_s \left(\rVert u \rVert^{\Lip(\gamma,\Omega)}_{s+1}  + \rVert \beta \rVert^{\Lip(\gamma,\Omega)}_{W^{s,\infty}(\mathbb{T}^{\nu+1})}\rVert u \rVert_{s_0}^{\Lip(\gamma,\Omega)}\right).
\]
\end{enumerate}
\end{lemma}

\color{black}
\section{Linear operators}\label{linsesdjsdjsd}
Let $\mathcal{A}:C^\infty_x\mapsto C^\infty_x$ be a linear operator.
We denote $\mathcal{A}^T$ to be its adjoint operator, that is,
\[
(\mathcal{A}[u],v)_{L^2_x}:=\int_{\mathbb{T}}\mathcal{A}u(x)\overline{v}(x)dx = \int_{\mathbb{T}}u(x)\overline{\mathcal{A}^Tv}(x)dx=(u,\mathcal{A}^T[v])_{L^2_x},
\]
where $\overline{u}(x)$ is the complex conjugation of $u(x)$. We also consider linear operators $\mathcal{A}:C^\infty_{\varphi,x}\mapsto C^\infty_{\varphi,x}$. Clearly, we can think of $\mathcal{A}=\mathcal{A}(\varphi)$ as a $\varphi$-dependent linear operator in $C^\infty_x$. Thanks to the following elementary relation:
\begin{align*}
(\mathcal{A}[f],g)_{L^2_{\varphi,x}} & := \int_{\mathbb{T}^{\nu}}\int_{\mathbb{T}}\mathcal{A}(\varphi)\left(f(\varphi,x)\right)\overline{g}(\varphi,x)dxd\varphi \\
& = \int_{\mathbb{T}^\nu}\int_\mathbb{T}f(\varphi,x)\overline{(\mathcal{A}(\varphi))^Tg(\varphi,x)}dxd\varphi = (f,\mathcal{A}^T[g])_{L^2_{\varphi,x}},
\end{align*}
we have that $\mathcal{A}^{T}(\varphi) = \mathcal{A}(\varphi)^T$, we simply use the notation $\mathcal{A}^T$ to denote the adjoint of $\mathcal{A}$\index{adjoint} as a linear operator in $C^\infty_{\varphi,x}$ or $C^\infty_x$, without any confusion.

For a  linear operator $\mathcal{A}:C^\infty_{\varphi,x}\mapsto C^\infty_{\varphi,x}$, we have its matrix representation as (recalling our convention of the Fourier series expansion in \eqref{sjd2sjjjsd})
\begin{align}\label{matrix_rep_2}
\mathcal{A}f(\varphi,x) = \sum_{j_1,j_2\in \mathbb{Z},\ l,l'\in \mathbb{Z}^\nu}\mathcal{A}^{j_1}_{j_2}(l-l')\widehat{f}_{j_1}(l')e^{\ii (l\cdot\varphi +j_2x)}.
\end{align}

\begin{definition}\label{realtorealop}
We say that an operator (not necessarily linear) $\mathcal{A}$ is real if $\mathcal{A}$ maps a real-valued function to a real-valued function.
\end{definition}

\subsection{Pseudo differential operators}\label{pdo_section}\index{Pseudo differential operator}
We consider symbols\index{symbol} $\mathfrak{a}=\mathfrak{a}(x,\xi)\in C^\infty(\mathbb{T}\times \mathbb{R})$.  We say $\mathfrak{a}$ is a classical symbol\index{classical symbol} of order $m$ and denote $\mathfrak{a}\in \mathcal{S}^m$ for some $m\in \mathbb{R}$, if 
\begin{align}\label{symbol_norm}\index{$\mathcal{S}^m$}
\sup_{x\in\mathbb{T}}|\partial_x^s\partial_\xi^\eta \mathfrak{a}(x,\xi)| \langle \xi \rangle^{-(m-\eta)} < \infty, \text{ for all $s,\eta\ge0$}, \text{ where $\langle \xi \rangle := \sqrt{1+ |\xi|^2}$.}
\end{align}
We denote 
\begin{align}\label{norm_sd2sdsymbol_xsd}
|\mathfrak{a}|_{m,H^s_x,\eta_0} := \max_{0\le\eta\le \eta_0}\sup_{\xi\in \mathbb{R}} \rVert \partial_\xi^{\eta}\mathfrak{a}(\cdot,\xi)\rVert_{H^s_x}\langle \xi \rangle^{\eta-m}.
\end{align}

\subsubsection{Standard quantization and Weyl quantization}\index{Standard quantization}\index{Weyl quantization}
For a given symbol $\mathfrak{a}$, we denote its standard quantization by
\begin{align}\label{usualquati}
Op(\mathfrak{a})[u]=\sum_{j\in \mathbb{Z}}\mathfrak{a}(x,j)\widehat{u}_je^{\ii j x}.
\end{align}\index{$Op$}
We also define the Weyl quantization of $\mathfrak{a}$ by 
\begin{align}\nonumber
Op^W(\mathfrak{a})[u]& :=\sum_{{k,j}\in \mathbb{Z}}\widehat{\mathfrak{a}}\left(k-j,\frac{k+j}{2}\right)\widehat{u}_je^{\ii k x},\\
\text{ where }\quad \widehat{\mathfrak{a}}(k,\xi)& :=\frac{1}{2\pi}\int_{\mathbb{T}}\mathfrak{a}(x,\xi)e^{-\ii k x}dx. \label{weyl_quant}
\end{align}\index{$Op^W$}
Compared to the usual quantization in \eqref{usualquati}, one can easily check that\index{$\mathcal{S}^m$}
\begin{align}\label{usual_weyl}
Op^W(\mathfrak{a}) = Op(\mathfrak{b})\iff  \widehat{\mathfrak{a}}\left(k,\frac{k}2 +\xi\right) =\widehat{\mathfrak{b}}\left(k,\xi\right) \text{ for all $k,\xi \in \mathbb{Z}$}.
\end{align}

For a linear operator $A$, we say a linear operator $A\in OPS^m$\index{$OPS^m$}, if there exists a classical symbol $\mathfrak{a}\in \mathcal{S}^m$ such that $A=Op(\mathfrak{a})$ or $A=Op^W(\mathfrak{a})$. Note that thanks to Lemma~\ref{sdtosd}, the statement $A\in OPS^m$ is well-defined, whether we use the standard quantization or the Weyl quantization.

With the above definition of $Op^W(\mathfrak{a})$, one can easily check, using \eqref{weyl_quant},  that
\begin{align}\label{transfposdsdsd}
Op^W(\mathfrak{a})^T = Op^W(\overline{\mathfrak{a}}),
\end{align}
therefore,
\begin{align}\label{adjoint_weyl}
\text{$Op^W(\mathfrak{a})$ is a  symmetric operator } \iff \text{ $\mathfrak{a}$ is a real-valued symbol.}
\end{align}

The following results are classical:
\begin{lemma}\label{sdtosd}
We have that:
\begin{enumerate}[label=(\arabic*)]
\item \label{sjdsdjclasiss1} (Equivalence of the quantizations) For $\mathfrak{a},\mathfrak{b}\in \mathcal{S}^m$ such that $Op^W(\mathfrak{a})=Op(\mathfrak{b})$, there exists $\mu=\mu(m,\eta)$ such that  for all $s\ge0$ and $\eta\ge0$,
\begin{align*}
|\mathfrak{b}|_{m,H^s_x,\eta}\le_{\eta,s}|\mathfrak{a}|_{m,H^{s+\mu}_x,\eta},\quad |\mathfrak{a}|_{m,H^s_x,\eta}\le_{\eta,s}|\mathfrak{b}|_{m,H^{s+\mu}_x,\eta}.
\end{align*}
\item \label{sjdsdjclasiss2} (Action as a linear operator)  Let $\mathfrak{a}\in \mathcal{S}^m$. There exists $\mu=\mu(m)$ such that  for all $s\ge 0$, 
\begin{align}\label{j2sdsdpsd01}
\rVert Op^W(\mathfrak{a})h\rVert_{H^s_x} \le |\mathfrak{a}|_{m,H^{s+\mu}_x,0}\rVert h \rVert_{H^{m}_x} + |\mathfrak{a}|_{m,H^{\mu}_x,0}\rVert h \rVert_{H^{s+m}_x}.
\end{align}
\end{enumerate}
\end{lemma}
\begin{proof}
\vspace{0.5\baselineskip}
\noindent\textit{Proof of \ref{sjdsdjclasiss1}.} We prove only
\begin{align}\label{onenightonly}
|\mathfrak{b}|_{m,H^s_x,\eta}\le_{\eta,s,m}|\mathfrak{a}|_{m,H^{s+\mu}_x,\eta},
\end{align}
since the other estimate follows in the same way.
By definition of the norm in \eqref{norm_sd2sdsymbol_xsd}, we have that for all $k\in \mathbb{Z},\xi\in \mathbb{R}$, 
\begin{align}\label{jsd2sdszzz}
|\widehat{\partial_\xi^\eta\mathfrak{a}}(k,\xi)|^2\langle k\rangle^{2s} \le \sum_{j\in \mathbb{Z}}|\widehat{\partial_\xi^\eta\mathfrak{a}}(j,\xi)|^2\langle j\rangle^{2s}\le_{s,\eta} |\mathfrak{a}|_{m,s,\eta}^2\langle \xi \rangle^{2(m-\eta)}.
\end{align}
Therefore,  using \eqref{usual_weyl}, we compute
\begin{align}
\rVert \partial_\xi^\eta\mathfrak{b}(\cdot,\xi)\rVert_{H^s_{x}}^2 &= \sum_{k\in \mathbb{Z}}\langle k \rangle^{2s}|\widehat{\partial_\xi^\eta\mathfrak{b}}(k,\xi)|^2 = \sum_{k\in \mathbb{Z}}\langle k\rangle^{2s}\left|\widehat{\partial_\xi^\eta\mathfrak{a}}(k,\xi+\frac{k}2)\right|^2 \nonumber\\
&\le  \sum_{k\in \mathbb{Z}} \langle k \rangle^{2s+\mu}\left| \widehat{\partial_\xi^\eta\mathfrak{a}}(k,\xi+\frac{k}2)\right|^2 \langle k \rangle^{-\mu}\nonumber \\
&\overset{\eqref{jsd2sdszzz}}{\le_{s,\eta}} |\mathfrak{a}|_{m,s+\mu,\eta}^2 \sum_{k\in \mathbb{Z}}\langle \xi + \frac{k}{2}\rangle^{2(m-\eta)}\langle k \rangle^{-\mu},\label{skd2rlaqkq}
\end{align}
for any $\mu\ge 0$. By choosing $\mu:=4|m-\eta|+2$, it is clear that
\[
\langle \xi + \frac{k}2\rangle^{2(m-\eta)}\langle k \rangle^{-\mu}\le_{m,\eta} 
\begin{cases}
\langle k\rangle^{2(m-\eta)}\langle k \rangle^{-\mu} \le_{m,\eta} \langle k\rangle^{-2|m-\eta|-2}, & \text{ if $\frac{2}{3}|k| \ge  |\xi|$} \\
\langle \xi \rangle^{2(m-\eta)}\langle k \rangle^{-\mu}\le_{m,\eta} \langle \xi \rangle^{2(m-\eta)}\langle k \rangle^{-2},&\text{ if $\frac{2}{3}|k|\le |\xi|$}
\end{cases}
\]
Hence, \eqref{skd2rlaqkq} yields that
\begin{align*}
& \rVert \partial_\xi^\eta\mathfrak{b}(\cdot,\xi)\rVert^2_{H^s_x} \\&\le_{s,m,\eta} |\mathfrak{a}|^2_{m,H^{s+\mu}_x,\eta}\left( \sum_{k\in \mathbb{Z},\ \frac{2}{3}|k|\ge |\xi|} \langle k \rangle^{-2|m-\eta| -2} +\langle \xi \rangle^{2(m-\eta)} \sum_{k\in \mathbb{Z},\ \frac{2}{3}|k|\le4|\xi|}\langle k \rangle^{-2}\right)\\
& \le_{m,\eta} |\mathfrak{a}|^2_{m,H^{s+\mu}_x,\eta}\langle \xi \rangle^{2(m-\eta)},
\end{align*}
which proves \eqref{onenightonly}.

\vspace{0.5\baselineskip}
\noindent\textit{Proof of \ref{sjdsdjclasiss2}.} Let us choose $\mathfrak{b}$ so that $Op(\mathfrak{b})=Op^W(\mathfrak{a})$ using \eqref{usual_weyl}. We compute
\[
\rVert Op^W(\mathfrak{a})h\rVert_{H^s_x}^2 = \rVert Op(\mathfrak{b})h\rVert_{H^s_x}^2 = \sum_{k\in \mathbb{Z}}\left(\sum_{j\in \mathbb{Z}}\langle k \rangle^{s}\widehat{\mathfrak{b}}(k-j,j) h_j \right)^2.
\]
Using $\langle k\rangle^s \le_s \langle k-j \rangle^s + \langle j \rangle^s$, we have
\begin{align}
\rVert Op^W(\mathfrak{a})h\rVert_{H^s_x}^2\le_s \sum_{k\in \mathbb{Z}}&\left(\left(\sum_{j\in \mathbb{Z}}\langle k-j \rangle^{s}|\widehat{\mathfrak{b}}(k-j,j)|| h_j |\right)^2\right. \nonumber \\
& \left.+\left(\sum_{j\in \mathbb{Z}}|\widehat{\mathfrak{b}}(k-j,j)|\langle j \rangle^s|h_j|\right)^2\right) \label{sd2sdrlaqkqsd2}
\end{align}
Using \eqref{jsd2sdszzz} with $\eta=0$, we have
\begin{align*}
|\widehat{\mathfrak{b}}(k-j,j)|& \le |\mathfrak{b}|_{m,H^{s+2}_x,0}\langle k-j\rangle^{-s-2}\langle j \rangle^m, \\  |\widehat{\mathfrak{b}}(k-j,j)|& \le |\mathfrak{b}|_{m,H^{2}_x,0}\langle k-j\rangle^{-2}\langle j \rangle^m.
\end{align*}
Plugging this into \eqref{sd2sdrlaqkqsd2}, we get
\begin{align}\nonumber
\rVert Op^W(\mathfrak{a})h\rVert_{H^s_x}^2 & \le_s |\mathfrak{b}|_{m,H^{s+2}_x,0}^2 \sum_{k\in \mathbb{Z}}\left(\sum_{j\in \mathbb{Z}} \langle k-j \rangle^{-2}\langle j\rangle^m|h_j| \right)^2 \\
& + |\mathfrak{b}|^2_{m,H^2_x,0} \sum_{k\in \mathbb{Z}}\left( \sum_{j\in \mathbb{Z}} \langle k-j \rangle^{-2}\langle j \rangle^{s+m}|h_j|\right)^2 \label{2sd23skkk}
\end{align}
Let $\underline{g},\underline{h}_{m},\underline{h}_{s+m}$, be such that
\[
\underline{g}(x):=\sum_{k}\langle k \rangle^{-2}e^{\ii k x},\quad \underline{h}_m(x):=\sum_j \langle j \rangle^{m}|h_j|e^{\ii j x},
\]
\[ \underline{h}_{s+m}(x):=\sum_j \langle j \rangle^{s+m}|h_j|e^{\ii j x}.
\]
Clearly, 
\[
\rVert \underline{g} \rVert_{L^\infty_x}\le C,\quad \rVert \underline{h}_m\rVert_{L^2_x}\le \rVert h \rVert_{H^m_x},\quad \rVert \underline{h}_{s+m}\rVert_{L^2_x}\le \rVert h \rVert_{H^{s+m}_x},
\]
for some $C>0$.
Hence, it follows from \eqref{2sd23skkk} that
\begin{align*}
\rVert Op^W(\mathfrak{a})h\rVert_{H^s_x}^2& \le_s |\mathfrak{b}|_{m,H^{s+2}_x,0}^2\int_{\mathbb{T}}(\underline{g}(x)\underline{h}_{m}(x))^2dx + |\mathfrak{b}|_{m,H^{2}_x,0}^2\int_{\mathbb{T}}(\underline{g}(x)\underline{h}_{s+m}(x))^2dx\\
&\le C\left( |\mathfrak{b}|_{m,H^{s+2}_x,0}^2\rVert h\rVert_{H^m_x}^2 + |\mathfrak{b}|_{m,H^{2}_x,0}^2\rVert h\rVert_{H^{s+m}_x}^2\right).
\end{align*}
Using \ref{sjdsdjclasiss1}, we obtain \eqref{j2sdsdpsd01}.
\end{proof}

\paragraph{Compositions:}
For two symbols $\mathfrak{a}\in\mathcal{S}^m$, $\mathfrak{b}\in\mathcal{S}^{m'}$, it follows straightforwardly from the definitions in \eqref{usualquati} and \eqref{weyl_quant} that
\begin{equation}\label{compose_289}
\begin{aligned}
&Op(\mathfrak{a})\circ Op(\mathfrak{b}) = Op(\mathfrak{c}),\quad \mathfrak{c}\in \mathcal{S}^{m+m'},\quad \widehat{\mathfrak{c}}(k,\xi) = \sum_{j_1+j_2=k}\widehat{\mathfrak{a}}(j_1,\xi +j_2)\widehat{\mathfrak{b}}(j_2,\xi),\\
&Op^{W}(\mathfrak{a})\circ Op^W(\mathfrak{b}) = Op^W(\mathfrak{c}),\quad \mathfrak{c}\in \mathcal{S}^{m+m'}, \\
& \widehat{\mathfrak{c}}(k,\xi) = \sum_{j_1+j_2 = k}\widehat{\mathfrak{a}}\left(j_1, \xi + \frac{j_2}{2} \right)\widehat{\mathfrak{b}}\left(j_2,\xi - \frac{j_1}{2}\right).
\end{aligned}
\end{equation}

\begin{lemma}\label{composition_lemmaaa}
For two symbols $\mathfrak{a}\in \mathcal{S}^m$, $\mathfrak{b}\in \mathcal{S}^{m'}$, the composition of $Op^W(\mathfrak{a}),Op^W(\mathfrak{b})$ can be expressed as
\begin{equation}\label{composition_weyl}
\begin{aligned}
Op^W(\mathfrak{a})\circ Op^W(\mathfrak{b}) &= Op^W(\mathfrak{c}),\\
{\mathfrak{c}}(x,\xi) &=\mathfrak{a}\mathfrak{b}  + \frac{1}{2\ii}\left\{ \mathfrak{a},\mathfrak{b}\right\} + \sigma_{\mathfrak{a},\mathfrak{b}}^1 +\sigma^2_{\mathfrak{a},\mathfrak{b}},
\end{aligned}
\end{equation}
where $\left\{\mathfrak{a},\mathfrak{b} \right\} := \partial_\xi \mathfrak{a} \partial_x \mathfrak{b} - \partial_\xi \mathfrak{b} \partial_x \mathfrak{a}$, and $\sigma_{\mathfrak{a},\mathfrak{b}}^1\in \mathcal{S}^{m+m'-2},\ \sigma^2_{\mathfrak{a},\mathfrak{b}}\in \mathcal{S}^{m+m' - 3}$  with
\begin{align*}
\sigma_{\mathfrak{a},\mathfrak{b}}^1(x,\xi) & =  \frac{1}{8}\left(-\partial_{xx}\mathfrak{a}\partial_{\xi\xi}\mathfrak{b} +2\partial_{x\xi}\mathfrak{a}\partial_{x\xi}\mathfrak{b}-\partial_{\xi\xi}\mathfrak{a}\partial_{xx}\mathfrak{b} \right),\\
\widehat{\sigma^2_{\mathfrak{a},\mathfrak{b}}}(k,\xi) & =(a_1b_4 + a_2(b_3+b_4)+a_3(b_2+b_3+b_4)+a_4(b_1+b_2+b_3+b_4)), 
\end{align*}
where $a_i,b_i$ are given in \eqref{fsdsdnotsds}.
\end{lemma}
\begin{proof}
Recall the Fourier series representation of the composed symbol from \eqref{compose_289}:
\begin{align}\label{symbol_compositions26}
\widehat{\mathfrak{c}}(k,\xi) = \sum_{j_1+j_2=k}\widehat{\mathfrak{a}}\left(j_1,\xi+\frac{j_2}2\right)\widehat{\mathfrak{b}}\left(j_2,\xi-\frac{j_1}{2}\right).
\end{align}
Using the Taylor expansion of smooth functions,
\begin{align}\label{elementary_taylor}
f(x)=\sum_{n=0}^{N}f^{(n)}(0)\frac{x^n}{n!} + x^{N+1}\int_0^{1}f^{(N+1)}(x\tau)\frac{(1-\tau)^{N}}{N!}d\tau \text{ for $x\in \mathbb{R}$},
\end{align}
we consider the expansion of $\widehat{\mathfrak{a}},\widehat{\mathfrak{b}}$ in the variable $\xi$:\begin{equation}\label{fsdsdnotsds}
\begin{aligned}
\widehat{\mathfrak{a}}\left(j_1,\xi+\frac{j_2}2\right) &=  \widehat{\mathfrak{a}}(j_1,\xi) + \widehat{\partial_\xi\mathfrak{a}}(j_1,\xi)\frac{\ii j_2}{2\ii} + \frac{1}{2!}\widehat{\partial_{\xi \xi}\mathfrak{a}}(j_1,\xi)\left(\frac{\ii j_2 }{2\ii}\right)^2 \\
&+ \left(\frac{j_2}{2}\right)^3\int_{0}^1\widehat{\partial_{\xi\xi\xi}\mathfrak{a}}(j_1,\xi + \tau\frac{j_2}2)\frac{(1-\tau)^2}{2!}d\tau\\
& =: a_1+a_2 + a_3 + a_4,\\
\widehat{\mathfrak{b}}\left(j_2,\xi-\frac{j_1}2\right) &=  \widehat{\mathfrak{b}}(j_2,\xi) - \widehat{\partial_\xi\mathfrak{b}}(j_2,\xi)\frac{\ii j_1}{2\ii } + \frac{1}{2!}\widehat{\partial_{\xi \xi}\mathfrak{b}}(j_2,\xi)\left(\frac{\ii j_1}{2\ii }\right)^2 \\
&- \left(\frac{j_1}{2}\right)^3\int_{0}^1\widehat{\partial_{\xi\xi\xi}\mathfrak{b}}(j_2,\xi - \tau\frac{j_1}2)\frac{(1-\tau)^2}{2!}d\tau\\
& =: b_1+b_2+b_3+b_4.
\end{aligned}
\end{equation}
Then, we can collect the terms as 
\begin{align*}
\widehat{\mathfrak{a}}\left(j_1,\xi+\frac{j_2}2\right) \widehat{\mathfrak{b}}\left(j_2,\xi-\frac{j_1}2\right) &= (a_1+a_2+a_3+a_4)(b_1+b_2+b_3+b_4)\\
&  = a_1b_1 + (a_1b_2+a_2b_1)+(a_1b_3+a_2b_2+a_3b_1)\\
& \ + (a_1b_4 + a_2(b_3+b_4)+a_3(b_2+b_3+b_4)\\
&\ +a_4(b_1+b_2+b_3+b_4)).
\end{align*}
Summing up all $j_1,j_2$, we have
\begin{align*}
\sum_{j_1+j_2=k}a_1b_1 &= \mathfrak{a}\mathfrak{b},\\
\sum_{j_1+j_2=k} (a_1b_2+a_2b_1) &= -\frac{1}{2\ii}\left( \partial_x\mathfrak{a}\partial_\xi\mathfrak{b}  - \partial_\xi\mathfrak{a}\partial_x\mathfrak{b}\right)\\
\sum_{j_1+j_2=k}(a_1b_3+a_2b_2 + a_3b_1) & =\frac{1}{4}\left(-\frac{1}{2}\partial_{xx}\mathfrak{a}\partial_{\xi\xi}\mathfrak{b} +\partial_{x\xi}\mathfrak{a}\partial_{x\xi}\mathfrak{b} -\frac{1}{2} \partial_{\xi\xi}\mathfrak{a}\partial_{xx}\mathfrak{b} \right)\\
&=: \sigma^1_{\mathfrak{a},\mathfrak{b}}.
\end{align*} 
Denoting the rest by 
\[
\sigma^2_{\mathfrak{a},\mathfrak{b}}:= \sum_{k}\sum_{j_1+j_2=k} (a_1b_4 + a_2(b_3+b_4)+a_3(b_2+b_3+b_4)+a_4(b_1+b_2+b_3+b_4)),
\]
we obtain the expression for $\mathfrak{c}$ in \eqref{composition_weyl}. Recalling the definition of $\mathcal{S}^{m}$ in \eqref{symbol_norm} and noting that $\mathfrak{a},\mathfrak{b}$ are smooth in the variable $x$, it follows straightforwardly that $\sigma^1_{\mathfrak{a},\mathfrak{b}}\in \mathcal{S}^{m+m'-2},\ \sigma^2_{\mathfrak{a},\mathfrak{b}}\in \mathcal{S}^{m+m'-3}$ (we do not prove quantitative estimates for $\sigma^1_{\mathfrak{a},\mathfrak{b}},\sigma^2_{\mathfrak{a},\mathfrak{b}}$ here, but we  postpone them to Lemma~\ref{compositionlemma25}). 
\end{proof}


\begin{remark}\label{indepofxxxx}In case $\mathfrak{a}$ or $\mathfrak{b}$ is independent of $x$, we have cleaner expression for the composition. Indeed, if $\mathfrak{a}(\xi)$ is independent of $x$, then it immediately follows from \eqref{compose_289} that
\begin{align}\label{compose_28926}
Op^W(\mathfrak{a})\circ Op^W(\mathfrak{b}) = \mathfrak{a}\left(\xi+\frac{k}{2}\right)\widehat{\mathfrak{b}}(k,\xi),\, Op^W(\mathfrak{b})\circ Op^W(\mathfrak{a}) = \widehat{\mathfrak{b}}(k,\xi)\mathfrak{a}\left( \xi - \frac{k}{2}\right).
\end{align}
This immediately implies that 
\begin{align}\label{compsewithpartialx}
\partial_xOp^W(\mathfrak{a}) = Op^W\left(\frac{1}{2}\partial_x \mathfrak{a} + i\xi \mathfrak{a}\right),\quad Op^W(\mathfrak{a})\partial_x = Op^W\left(-\frac{1}{2}\partial_x \mathfrak{a} +\ii \xi \mathfrak{a}\right).
\end{align}
\end{remark}

\color{black}
\subsubsection{Symmetric commutator}\index{Symmetric commutator}
The usual commutator of two linear operators is defined as
\begin{align}\label{comusd2sd}
[A,B] := AB - BA.
\end{align}
We also consider a slight variation:
\begin{align}\label{commutator}
[A,B]_x:= A\partial_x B - B\partial_x A.
\end{align}
Note that if $A,B$ are symmetric operators, then $[A,B]_x$ is also a symmetric operator, since $\partial_x$ is asymmetric\index{asymmetric}.

\begin{lemma}\label{symsdcodmyysd}
Let $A=Op^W(\mathfrak{a})$, $B=Op^W(\mathfrak{b})$ for some $\mathfrak{a}\in \mathcal{S}^m$, $\mathfrak{b}\in \mathcal{S}^{m'}$. Then $[A,B]_x\in OPS^{m+m'}$. That is, there exists a symbol $\mathfrak{a}\star \mathfrak{b}\in \mathcal{S}^{m+m'}$ such that  \index{$\star$}
\begin{align}\label{syysdsdsd}
[A,B]_x =  Op^W(\mathfrak{a}\star \mathfrak{b}).
\end{align}
The symbol $\mathfrak{a}\star \mathfrak{b}$ admits an expansion
\begin{equation}\label{commutator_pseudo}
\begin{aligned}
(\mathfrak{a}\star \mathfrak{b})(x,\xi) &= \left\{\mathfrak{a},\mathfrak{b} \right\}_x + r_{\mathfrak{a},\mathfrak{b}},\\
\left\{ \mathfrak{a},\mathfrak{b}\right\}_{x} &:= \mathfrak{a}\partial_x\mathfrak{b} - \mathfrak{b}\partial_x \mathfrak{a} +\xi\left\{\mathfrak{a},\mathfrak{b}\right\} = \partial_\xi(\xi \mathfrak{a})\partial_x \mathfrak{b} - \partial_x \mathfrak{a} \partial_{\xi}(\xi \mathfrak{b}) \in \mathcal{S}^{m+m'},\\
r_{\mathfrak{a},\mathfrak{b}} & := \left(\frac{1}{2}(\sigma^1_{\mathfrak{a},\partial_x\mathfrak{b}} + \sigma^2_{\mathfrak{a},\partial_x\mathfrak{b}}) + \ii \sigma^2_{\mathfrak{a},\xi\mathfrak{b}} \right) \\
& \qquad - \left(\frac{1}{2}(\sigma^1_{\mathfrak{b},\partial_x\mathfrak{a}} + \sigma^2_{\mathfrak{b},\partial_x\mathfrak{a}}) + \ii \sigma^2_{\mathfrak{b},\xi\mathfrak{a}} \right)\in \mathcal{S}^{m+m' - 2},
\end{aligned}
\end{equation}
where $\sigma_{\mathfrak{a},\mathfrak{b}}^{i}$ for $i=1,2$ are as in Lemma~\ref{composition_lemmaaa}.
\end{lemma}
\begin{proof}
From Lemma~\ref{composition_lemmaaa} and \eqref{compsewithpartialx}, we have that
\begin{align*}
A\partial_x B &= Op^W(\mathfrak{a})\circ Op^W(\frac{1}{2}\partial_x\mathfrak{b} + \ii \xi \mathfrak{b}) \\
& = \frac{1}{2}Op^W(\mathfrak{a})\circ Op^W(\partial_x \mathfrak{b}) + \ii Op^W(\mathfrak{a})\circ Op^W( \mathfrak{\xi \mathfrak{b}})\\
& = \frac{1}{2}Op^W\left(\mathfrak{a}\partial_x\mathfrak{b} + \frac{1}{2\ii }\left\{ \mathfrak{a},\partial_x\mathfrak{b}\right\} + (\sigma^1_{\mathfrak{a},\partial_x\mathfrak{b}} + \sigma^2_{\mathfrak{a},\partial_x\mathfrak{b}}) \right) \\
& \qquad+ \ii Op^W\left(\left(\mathfrak{a}\xi\mathfrak{b} + \frac{1}{2\ii}\left\{ \mathfrak{a},\xi\mathfrak{b}\right\} + \sigma^1_{\mathfrak{a},\xi\mathfrak{b}}\right) + \sigma^2_{\mathfrak{a},\xi\mathfrak{b}} \right)\\
& = Op^W\left( \frac{1}{2}\mathfrak{a}\partial_x \mathfrak{b} + \frac{1}{4\ii}\left\{ \mathfrak{a},\partial_x\mathfrak{b}\right\} + \ii \mathfrak{a}\xi\mathfrak{b} + \frac{1}{2}\left\{ \mathfrak{a},\xi\mathfrak{b}\right\} + \ii \sigma_{\mathfrak{a},\xi\mathfrak{b}}^1 \right) \\
& \qquad + Op^W\left(\frac{1}{2}(\sigma^1_{\mathfrak{a},\partial_x\mathfrak{b}} + \sigma^2_{\mathfrak{a},\partial_x\mathfrak{b}}) + \ii \sigma^2_{\mathfrak{a},\xi\mathfrak{b}} \right).
\end{align*}
Similarly, we have
\begin{align*}
B\partial_x A & = Op^W\left(  \frac{1}{2}\mathfrak{b}\partial_x \mathfrak{a} + \frac{1}{4\ii}\left\{ \mathfrak{b},\partial_x\mathfrak{a}\right\} + \ii \mathfrak{b}\xi\mathfrak{a} + \frac{1}{2}\left\{ \mathfrak{b},\xi\mathfrak{a}\right\} + \ii \sigma_{\mathfrak{b},\xi\mathfrak{a}}^1\right) \\
& + Op^W\left(\frac{1}{2}(\sigma^1_{\mathfrak{b},\partial_x\mathfrak{a}} + \sigma^2_{\mathfrak{b},\partial_x\mathfrak{a}}) + \ii \sigma^2_{\mathfrak{b},\xi\mathfrak{a}} \right)
\end{align*}
Using the definitions of $\left\{ \mathfrak{a},\mathfrak{b}\right\}$ and $\sigma_{\mathfrak{a},\mathfrak{b}}^1$ in Lemma~\ref{composition_lemmaaa}, we have
\begin{align*}
& A\partial_x B - B\partial_x A = Op^W \left(\left\{ \mathfrak{a},\mathfrak{b}\right\}_x \right)  \\
&+  Op^W\left(\left(\frac{1}{2}(\sigma^1_{\mathfrak{a},\partial_x\mathfrak{b}} + \sigma^2_{\mathfrak{a},\partial_x\mathfrak{b}}) + \ii \sigma^2_{\mathfrak{a},\xi\mathfrak{b}} \right) - \left(\frac{1}{2}(\sigma^1_{\mathfrak{b},\partial_x\mathfrak{a}} + \sigma^2_{\mathfrak{b},\partial_x\mathfrak{a}}) + \ii \sigma^2_{\mathfrak{b},\xi\mathfrak{a}} \right) \right),
\end{align*}
which gives us \eqref{commutator_pseudo}.
\end{proof}

 \subsubsection{Norms of  symbols}
As in Section~\ref{functionspaces}, we will consider symbols that depend on $\omega\in \Omega\subset \mathbb{R}^\nu$ and $\varphi\in \mathbb{T}^\nu$, that is, $\mathfrak{a}(\omega)=\mathfrak{a}(\omega,\varphi,x,\xi)\in \mathcal{S}^m$. We define the norms:\begin{equation}\label{symbol_norms6223}
\begin{aligned}
|\mathfrak{a}(\omega)|_{m,s,\eta_0} &:= \max_{0\le\eta\le \eta_0}\sup_{\xi\in \mathbb{R}} \rVert \partial_\xi^{\eta}\mathfrak{a}(\omega,\cdot,\cdot,\xi)\rVert_{H^s_{\varphi,x}}\langle \xi \rangle^{\eta-m},\\
|\mathfrak{a}|_{m,s,\eta_0}^{\sup} &:= \sup_{\omega\in \Omega}|\mathfrak{a}(\omega)|_{m,s,\eta_0},\\
|\mathfrak{a}|_{m,s,\eta_0}^{\lip}& := \sup_{\omega_1,\omega_2\in \Omega}\frac{|\mathfrak{a}(\omega_1) - \mathfrak{a}(\omega_2)|_{m,s,\eta_0}}{|\omega_1-\omega_2|},\\
|\mathfrak{a}|_{m,s,\eta_0}^{\Lip(\gamma,\Omega)}&:= |\mathfrak{a}|_{m,s,\eta_0}^{\sup} + \gamma |\mathfrak{a}|_{m,s,\eta_0}^{\lip},\text{ for some $\gamma\in (0,1)$.}
\end{aligned}
\end{equation}

The following lemma can be easily proved by following the proof of Lemma~\ref{sdtosd}, word by word. For simplicity, we omit the proof:

\begin{lemma}\label{sdt2osd1}
We have that:
\begin{enumerate}[label=(\arabic*)]
\item \label{sjdsdjclasiss21} (Equivalence of the quantizations) For $\mathfrak{a},\mathfrak{b}\in \mathcal{S}^m$ such that $Op^W(\mathfrak{a})=Op(\mathfrak{b})$, there exists $\mu=\mu(m,\eta)$ such that  for all $s\ge s_0$ and $\eta\ge0$,
\begin{align*}
|\mathfrak{b}|_{m,s,\eta}^{\Lip(\gamma,\Omega)}\le_{\eta,s}|\mathfrak{a}|_{m,{s+\mu},\eta}^{\Lip(\gamma,\Omega)},\quad |\mathfrak{a}|_{m,s,\eta}^{\Lip(\gamma,\Omega)}\le_{\eta,s}|\mathfrak{b}|_{m,{s+\mu},\eta}^{\Lip(\gamma,\Omega)}.
\end{align*}
\item \label{sjdsdjclasiss22} (Action as a linear operator)  Let $\mathfrak{a}\in \mathcal{S}^m$. There exists $\mu=\mu(m)$ such that  for all $s\ge s_0$, 
\begin{align}\label{j2sdsxxdpsd01}
\rVert Op^W(\mathfrak{a})h\rVert^{\Lip(\gamma,\Omega)}_{s} \le_s |\mathfrak{a}|_{m,{s+\mu},0}^{\Lip(\gamma,\Omega)}\rVert h \rVert^{\Lip(\gamma,\Omega)}_{s_0+m} + |\mathfrak{a}|^{\Lip(\gamma,\Omega)}_{m,{s_0+\mu},0}\rVert h \rVert^{\Lip(\gamma,\Omega)}_{{s+m}}.
\end{align}
\end{enumerate}
\end{lemma}

In what follows, we collect useful properties of the norms for the symbols depending on $\omega,\varphi$. We denote the Fourier modes of a symbol $\mathfrak{a}$ in $\varphi,x$ by
\begin{equation}\label{fourier_ssd2s}
\begin{aligned}\widehat{\mathfrak{a}}(\omega,\varphi,j,\xi)&:= \frac{1}{2\pi}\int_{\mathbb{T}}\mathfrak{a}(\omega,\varphi,x,\xi)e^{-\ii j x}dx,\\ 
\widehat{\mathfrak{a}}^{\varphi,x}(\omega,l,j,\xi)&:=\frac{1}{(2\pi)^{\nu+1}}\int_{\mathbb{T}^\nu\times \mathbb{T}}\mathfrak{a}(\omega,\varphi,x,\xi)e^{-\ii (l\cdot \varphi + jx)}d\varphi dx.
\end{aligned}
\end{equation}
We first prove the estimates for a product of two symbols:
\begin{lemma}\label{product_symbols}
Given $\mathfrak{a}\in \mathcal{S}^{m}$, $\mathfrak{b}\in \mathcal{S}^{m'}$ for some $m,m'\in \mathbb{R}$, it holds that
\begin{align}
|\mathfrak{a}\mathfrak{b}|_{m+m',s,\eta}^{\Lip(\gamma,\Omega)}&\le_{m,m',s,\eta} |\mathfrak{a}|^{\Lip(\gamma,\Omega)}_{m,s,\eta} |\mathfrak{b}|^{\Lip(\gamma,\Omega)}_{m',s_0,\eta} +  |\mathfrak{a}|^{\Lip(\gamma,\Omega)}_{m,s_0,\eta} |\mathfrak{b}|^{\Lip(\gamma,\Omega)}_{m,s,\eta}  \text{ for $\gamma\in (0,1)$.} \label{lipnorm_symbol2}
\end{align}
\end{lemma}
\begin{proof}
We first prove
\begin{align}
|\mathfrak{a}\mathfrak{b}|_{m+m',s,\eta_0}^{\sup}&\le_{m,m',s,\eta}|\mathfrak{a}|^{\sup}_{m,s,\eta}|\mathfrak{b}|^{\sup}_{m',s_0,\eta} +  |\mathfrak{a}|^{\sup}_{m,s_0,\eta}|\mathfrak{b}|^{\sup}_{m',s,\eta}. \label{supnorm_symbols}
\end{align}
Let $\omega$ be fixed. For simplicity, we omit the dependence on $\varphi,x,\xi$ in the notation. For a fixed $0\le \eta\in\mathbb{N}_0$, it is clear that 
\[
\left(\partial_\xi\right)^{\eta}(\mathfrak{a}(\omega)\mathfrak{b}(\omega)) =   \sum_{n=0}^{\eta} \binom{\eta}{n}((\partial_\xi)^{\eta-n}\mathfrak{a}(\omega))((\partial_\xi)^{n}\mathfrak{b}(\omega)).
\] Therefore using  \eqref{bncal} and splitting $\langle \xi\rangle^{\eta -(m+m')}\le_{\eta,m,m'} \langle \xi \rangle^{\eta-n-m}\langle \xi \rangle^{n-m'}$, we get
\begin{align*}
& \rVert \partial_\xi^{\eta}(\mathfrak{a}\mathfrak{b})(\omega,\cdot,\cdot,\xi)\rVert_{H^s_{\varphi,x}}\langle \xi  \rangle^{\eta-(m+m')}\\
&\le_{\eta,s,m,m'} |\mathfrak{a}(\omega)|_{m,s,\eta}|\mathfrak{b}(\omega)|_{m',s_0,\eta} +  |\mathfrak{a}(\omega)|_{m,s_0,\eta}|\mathfrak{b}(\omega)|_{m',s,\eta}.
\end{align*}
Therefore,
\begin{align}\label{mms2s}
|\mathfrak{a}(\omega)\mathfrak{b}(\omega)|_{m+m',s,\eta}\le_{\eta,s,m,m'} |\mathfrak{a}(\omega)|_{m,s,\eta}|\mathfrak{b}(\omega)|_{m',s_0,\eta} +  |\mathfrak{a}(\omega)|_{m,s_0,\eta}|\mathfrak{b}(\omega)|_{m',s,\eta},
\end{align}
This immediately implies \eqref{supnorm_symbols} by taking the supremum in $\omega$. For the Lipschitz dependence in \eqref{lipnorm_symbol2}, we have (denoting $\Delta_{12} \mathfrak{a}:=\mathfrak{a}(\omega_1)-\mathfrak{a}(\omega_2)$), that 
$
\mathfrak{a}(\omega_1)\mathfrak{b}(\omega_1) - \mathfrak{a}(\omega_2)\mathfrak{b}(\omega_2) = (\Delta_{12}\mathfrak{a}) \mathfrak{b}(\omega_1) + \mathfrak{a}(\omega_2)\Delta_{12} \mathfrak{b}.$
Hence,
\[
\frac{\mathfrak{a}(\omega_1)\mathfrak{b}(\omega_1) - \mathfrak{a}(\omega_2)\mathfrak{b}(\omega_2)}{|\omega_1-\omega_2|} = \frac{\Delta_{12} \mathfrak{a}}{|\omega_1-\omega_2|}\mathfrak{b}(\omega_1) + \mathfrak{a}_2(\omega_2)\frac{\Delta_{12} \mathfrak{b}}{|\omega_1-\omega_2|}.
\]
Again using \eqref{mms2s},  we have
\begin{align*}
&\frac{\left|\mathfrak{a}(\omega_1)\mathfrak{b}(\omega_1)-\mathfrak{a}(\omega_2)\mathfrak{b}(\omega_2) \right|_{m+m',s,\eta}}{|\omega_1-\omega_2|} \\
& \quad \le_{\eta,s,m,m'} |\mathfrak{a}|^{\lip}_{m,s,\eta}|\mathfrak{b}|^{\sup}_{m',s_0,\eta} + |\mathfrak{a}|^{\lip}_{m,s_0,\eta}|\mathfrak{b}|^{\sup}_{m',s,\eta} \\
& \quad +  |\mathfrak{b}|^{\lip}_{m',s,\eta}|\mathfrak{a}|^{\sup}_{m,s_0,\eta} + |\mathfrak{b}|^{\lip}_{m',s_0,\eta}|\mathfrak{a}|^{\sup}_{m,s,\eta}\\
&\quad  \le (|\mathfrak{a}|^{\lip}_{m,s,\eta} + |\mathfrak{a}|^{\sup}_{m,s,\eta})(|\mathfrak{b}|^{\lip}_{m',s_0,\eta} +|\mathfrak{b}|^{\sup}_{m',s_0,\eta}) \\
& \quad + (|\mathfrak{a}|^{\lip}_{m,s_0,\eta} + |\mathfrak{a}|^{\sup}_{m,s_0,\eta})(|\mathfrak{b}|^{\lip}_{m',s,\eta} +|\mathfrak{b}|^{\sup}_{m',s,\eta}).
\end{align*}
Then, the result follows by multiplying by $\gamma$ on both sides and  taking the supremum in $\omega_1,\omega_2$.
\end{proof}

\begin{remark}\label{remarksd2dsds}
In case $\mathfrak{a}=\mathfrak{a}(\omega,\xi)\in \mathcal{S}^m$ and $\mathfrak{b}=\mathfrak{b}(\omega,\xi)\in \mathcal{S}^{m'}$ are  Fourier multipliers, that is, they do not depend on $x,\varphi$, then \eqref{compose_28926} shows that
\[
Op^W(\mathfrak{a})\circ Op^W(\mathfrak{b}) = Op^W(\mathfrak{a}\mathfrak{b}). 
\]Using the norms in \eqref{norm_sd2sdsymbol_xsd} and \eqref{omega_dep_norm1}, it follows straightforwardly from Lemma~\ref{product_symbols} that
\[
|\mathfrak{a}\mathfrak{b}|^{\Lip(\gamma,\Omega)}_{m,H^0_x,\eta}\le_{\eta} |\mathfrak{a}|^{\Lip(\gamma,\Omega)}_{m,H^0_x,\eta}|\mathfrak{b}|^{\Lip(\gamma,\Omega)}_{m,H^0_x,\eta}.
\]
\end{remark}

Now, we estimate symbols obtained from a composition of two symbols. In view of Lemma~\ref{composition_lemmaaa}, we prove the following lemma:

\begin{lemma}\label{compositionlemma2502}
Given $\omega,\varphi$ dependent symbols $\mathfrak{a}\in \mathcal{S}^{m}$, $\mathfrak{b}\in \mathcal{S}^{m'}$ for some $m,m'\in \mathbb{R}$, we define $\mathfrak{c}_{\tau_1,\tau_2}$ be such that
\[
\widehat{\mathfrak{c}_{\tau_1,\tau_2}}(\omega,\varphi,k,\xi):=\sum_{j_1+j_2=k}\widehat{\mathfrak{a}}(\omega,\varphi,j_1,\xi+\tau_1j_2)\widehat{\mathfrak{b}}(\omega,\varphi,j_2,\xi+\tau_2j_1).
\]
Then, for each $\eta\in \mathbb{N}_0$,  there exists $\mu=\mu(m,m',\eta)$ such that 
\begin{align}\label{sup_norm21s}
\sup_{|\tau_1|,|\tau_2|\le 1}|\mathfrak{c}_{\tau_1,\tau_2}|^{\sup}_{m+m',s,\eta} \le_{m,m',s,\eta}|\mathfrak{a}|^{\sup}_{m,s+\mu,\eta}|\mathfrak{b}|^{\sup}_{m',s_0+\mu,\eta} + |\mathfrak{a}|^{\sup}_{m,s_0+\mu,\eta}|\mathfrak{b}|^{\sup}_{m',s+\mu,\eta}.
\end{align}
And, 
\begin{align}\label{lippss2}
\sup_{|\tau_1|,|\tau_2|\le 1}&|\mathfrak{c}_{\tau_1,\tau_2}|^{\Lip(\gamma,\Omega)}_{m+m',s,\eta} \le_{\nu,m,m',s,\eta}|\mathfrak{a}|^{\Lip(\gamma,\Omega)}_{m,s+\mu,\eta}|\mathfrak{b}|^{\Lip(\gamma,\Omega)}_{m',s_0+\mu,\eta} + |\mathfrak{a}|^{\Lip(\gamma,\Omega)}_{m,s_0+\mu,\eta}|\mathfrak{b}|^{\Lip(\gamma,\Omega)}_{m',s+\mu,\eta}
\end{align}
for $|\gamma|\le 1$.
\end{lemma}

\begin{proof}

Using the notation in \eqref{fourier_ssd2s} and recalling the norm $|\cdot|_{m,s,\eta}$ in \eqref{symbol_norms6223}, it is clear that  for each fixed $\omega\in \Omega$,
\[
 \sup_{\xi\in \mathbb{R}}\sum_{j\in \mathbb{Z},\ l\in \mathbb{Z}^{\nu}}(\langle j \rangle + \langle l \rangle)^{2s}|\partial_{\xi}^{\eta}\widehat{\mathfrak{a}}^{\varphi,x}(\omega,l,j,\xi)|^{2}\langle \xi \rangle^{2(\eta-m)}\le_{m,\eta,s} |\mathfrak{a}(\omega)|^2_{m,s,\eta}.
\]
Therefore for each fixed $(l,j,\xi)$,
\begin{align}\label{as2ss}
|\partial_\xi^{\eta}\widehat{\mathfrak{a}}^{\varphi,x}(\omega,l,j,\xi)|\le_{m,s,\eta} |\mathfrak{a}(\omega)|_{m,s,\eta} \frac{\langle \xi \rangle^{m-\eta}}{(\langle j \rangle + \langle l \rangle)^s},\text{ for each $s\ge s_0$}.
\end{align}
Similarly, we have
\begin{align}\label{as2ss2}
|\partial_\xi^{\eta}\widehat{\mathfrak{b}}^{\varphi,x}(\omega,l,j,\xi)|\le_{m',s,\eta} |\mathfrak{b}(\omega)|_{m',s,\eta_0}\frac{\langle\xi\rangle^{m'-\eta}}{(\langle j \rangle + \langle l \rangle)^s},\text{ for each $s\ge s_0$} .
\end{align}
Then we can write each Fourier mode of $\partial_\xi^{\eta}\mathfrak{c}_{\tau_1,\tau_2}$ as
\begin{align*}
&\partial_\xi^{\eta}\widehat{\mathfrak{c}_{\tau_1,\tau_2}}^{\varphi,x}(\omega,l,k,\xi)\\
&=\sum_{n=0}^{\eta}\binom{\eta}{n}\sum_{\substack{j_1+j_2 = k,\\ l_1+l_2=l\in \mathbb{Z}^{\nu}}}\partial_{\xi}^{n}\widehat{\mathfrak{a}}^{\varphi,x}(\omega,l_1,j_1,\xi+\tau_1j_2)\partial_{\xi}^{\eta-n}\widehat{b}^{\varphi,x}(\omega,l_2,j_2,\xi+\tau_2j_1).
\end{align*}
Using that $(\langle k \rangle + \langle l\rangle)^{2s} \le_s (\langle j_1 \rangle + \langle l_1 \rangle)^{2s} + (\langle j_2 \rangle + \langle l_2 \rangle)^{2s}$ for $j_1+j_2=k,\ l_1+l_2=l$, we have

\begin{equation}\label{summm2}
\begin{aligned}
(\langle k\rangle& + \langle l \rangle )^{2s}|\partial_\xi^{\eta}\widehat{\mathfrak{c}_{\tau_1,\tau_2}}^{\varphi,x}(\omega,l,k,\xi)|^2\\
& \le_{s,\eta,n} \sum_{n=0}^{\eta}\left( \left(\sum_{\substack{j_1+j_2=k,\\ l_1+l_2 = l}}(\langle j_1 \rangle + \langle l_1 \rangle)^{s}|(\partial_\xi)^{\eta-n}\widehat{\mathfrak{a}}^{\varphi,x}(\omega,l_1,j_1,\xi+\tau_1j_2)|\right.\right. \\
& \left.\left.\qquad \qquad \qquad \qquad\times|(\partial_\xi)^{n}\widehat{\mathfrak{b}}^{\varphi,x}(\omega,l_2,j_2,\xi+\tau_2j_1)|\right)^2 \right. \\
&\left.  \ + \left(\sum_{\substack{j_1+j_2=k,\\ l_1+l_2 = l}}(\langle j_2 \rangle + \langle l_2 \rangle)^{s}|(\partial_\xi)^{\eta-n}\widehat{\mathfrak{a}}^{\varphi,x}(\omega,l_1,j_1,\xi+\tau_1j_2)|\right.\right. \\
& \left.\left.\qquad \qquad \qquad \qquad\times|(\partial_\xi)^{n}\widehat{\mathfrak{b}}^{\varphi,x}(\omega,l_2,j_2,\xi+\tau_2j_1)|\right)^2\right).
\end{aligned}
\end{equation}
Using \eqref{as2ss} and \eqref{as2ss2}, we have that for some $\mu=\mu(m,m',\eta)\ge 0$, which will be chosen later,
\begin{align}
&|(\partial_\xi)^{\eta-n}\widehat{\mathfrak{a}}^{\varphi,x}(\omega,l_1,j_1,\xi+\tau_1j_2)||(\partial_\xi)^{n}\widehat{\mathfrak{b}}^{\varphi,x}(\omega,l_2,j_2,\xi+\tau_2j_1)| \nonumber\\
&\le_{m,m',\eta,s} |\mathfrak{a}(\omega)|_{m,s+\mu,\eta}\frac{\langle \xi+\tau_1j_2 \rangle^{m-(\eta-n)}}{(\langle j_1 \rangle + \langle l_1 \rangle)^{s+\mu}}|\mathfrak{b}(\omega)|_{m',s_0+\mu,\eta_0}\frac{\langle \xi +\tau_2j_1\rangle^{m'-n}}{(\langle j_2 \rangle + \langle l_2 \rangle)^{s_0+\mu}}\nonumber\\
&\le |\mathfrak{a}(\omega)|_{m,s+\mu,\eta}|\mathfrak{b}(\omega)|_{m',s_0+\mu,\eta}\frac{\langle \xi +\tau_1j_2\rangle^{m-(\eta-n)} \langle \xi + \tau_2j_1\rangle^{m'-n}}{(\langle j_1 \rangle + \langle l_1 \rangle)^{s+\mu}(\langle j_2 \rangle + \langle l_2 \rangle)^{s_0+\mu}}. \label{firstsum1}
\end{align}
Similarly, the same quantity can be estimated as
\begin{align}\label{firstsum2}
&|(\partial_\xi)^{\eta-n}\widehat{\mathfrak{a}}^{\varphi,x}(\omega,l_1,j_1,\xi+\tau_1j_2)||(\partial_\xi)^{n}\widehat{\mathfrak{b}}^{\varphi,x}(\omega,l_2,j_2,\xi+\tau_2j_1)| \nonumber\\
&\le _{m,m',\eta,s}|\mathfrak{a}(\omega)|_{m,s_0+\mu,\eta}|\mathfrak{b}(\omega)|_{m',s+\mu,\eta}\frac{\langle \xi +\tau_1j_2\rangle^{m-(\eta-n)}\langle \xi + \tau_2j_1\rangle^{m'-n}}{(\langle j_1 \rangle + \langle l_1 \rangle)^{s_0+\mu}(\langle j_2 \rangle + \langle l_2 \rangle)^{s+\mu}}.
\end{align}
Plugging these two estimates into  \eqref{summm2} (Plugging \eqref{firstsum1} into the first summation and plugging \eqref{firstsum2} into the second summation into \eqref{summm2}), we obtain
\begin{align*}
&(\langle k\rangle + \langle l \rangle )^{2s}|(\partial_\xi)^{\eta}\widehat{\mathfrak{c}_{\tau_1,\tau_2}}^{\varphi,x}(\omega,l,k,\xi)|^2 \\
& \le_{m,m',s,\eta}\left(\sum_{\substack{j_1+j_2=k \\ l_1+l_2=l}}|\mathfrak{a}(\omega)|_{m,s+\mu,\eta}|\mathfrak{b}(\omega)|_{m',s_0+\mu,\eta}\frac{\langle \xi +\tau_1j_2\rangle^{m-(\eta-n)} \langle \xi + \tau_2j_1\rangle^{m'-n}}{(\langle j_1 \rangle + \langle l_1 \rangle)^{\mu}(\langle j_2 \rangle + \langle l_2 \rangle)^{s_0+\mu}}\right)^2\\
& \ + \left(\sum_{\substack{j_1+j_2=k\\ l_1+l_2=l}}|\mathfrak{a}(\omega)|_{m,s_0+\mu,\eta}|\mathfrak{b}(\omega)|_{m',s+\mu,\eta}\frac{\langle \xi +\tau_1j_2\rangle^{m-(\eta-n)} \langle \xi + \tau_2j_1\rangle^{m'-n}}{(\langle j_1 \rangle + \langle l_1 \rangle)^{s_0+\mu}(\langle j_2 \rangle + \langle l_2 \rangle)^{\mu}}\right)^2
\end{align*}
Using the Cauchy-Schwarz inequality, we get  (using that $\sum_{j,l} \left(\frac{1}{\langle j \rangle + \langle l \rangle}\right)^{s_0} <_{s_0} 1 $, thanks to \eqref{range_of_s}),
\begin{equation}\label{hhssx2}
\begin{aligned}
& \sum_{k\in \mathbb{Z},\ l \in \mathbb{Z}^\nu}(\langle k\rangle + \langle l \rangle )^{2s}|(\partial_\xi)^{\eta}\widehat{\mathfrak{c}_{\tau_1,\tau_2}}^{\varphi,x}(\omega,l,k,\xi)|^2 \\
& \le_{m,m',\eta,s}(|\mathfrak{a}(\omega)|_{m,s+\mu,\eta_0}|\mathfrak{b}(\omega)|_{m',s_0+\mu,\eta})^2 \\
& \ \times \sum_{j_1,j_2\in\mathbb{Z},\ l_1,l_2\in \mathbb{Z}^\nu}\left(\frac{\langle \xi +\tau_1j_2\rangle^{m-(\eta-n)} \langle \xi + \tau_2j_1\rangle^{m'-n}}{(\langle j_1 \rangle + \langle l_1 \rangle)^{\mu}(\langle j_2 \rangle + \langle l_2 \rangle)^{\mu}}\right)^2\\
& \ + \left(|\mathfrak{a}(\omega)|_{m,s_0+\mu,\eta}|\mathfrak{b}(\omega)|_{m',s+\mu,\eta} \right)^2 \\
& \ \times \sum_{j_1,j_2\in\mathbb{Z},\ l_1,l_2\in \mathbb{Z}^\nu}\left(\frac{\langle \xi +\tau_1j_2\rangle^{m-(\eta-n)} \langle \xi + \tau_2j_1\rangle^{m'-n}}{(\langle j_1 \rangle + \langle l_1 \rangle)^{\mu}(\langle j_2 \rangle + \langle l_2 \rangle)^{\mu}}\right)^2.
\end{aligned}
\end{equation}
Now, we estimate the summations in the above estimate. Since $|\tau_1|,|\tau_2|\le 1$, we have that 
\begin{align}\label{xi_des}
\langle \xi + \tau_1 j_2\rangle^{m-(\eta-n)}\langle \xi +\tau_2j_1\rangle^{m'-n}\le_{m,m',\eta_0} \langle \xi \rangle^{m+m' - \eta},\text{ for $|j_1|,|j_2|\le \frac{|\xi|}{2}$}.
\end{align}
If $|j_1|>\frac{|\xi|}{2}$ or $|j_2|> \frac{|\xi|}{2}$, then one can choose large enough $\mu$ depending only on $m,m',\eta$ such that
\begin{align*}
&\sum_{\substack{j_1,j_2\in\mathbb{Z},\\ \text{$|j_1|\ge \frac{|\xi|}{2}$ or $|j_2|\ge \frac{|\xi|}{2}$},\\ l_1,l_2\in \mathbb{Z}^\nu}}\left(\frac{\langle \xi +\tau_1j_2\rangle^{m-(\eta-n)} \langle \xi + \tau_2j_1\rangle^{m'-n}}{(\langle j_1 \rangle + \langle l_1 \rangle)^{\mu}(\langle j_2 \rangle + \langle l_2 \rangle)^{\mu}}\right)^2 \le_{m,m',\eta}\langle \xi \rangle^{2(m+m' - \eta)}.
\end{align*}
Combining this with \eqref{xi_des} and taking the supremum in $\omega$, the estimate in \eqref{hhssx2} yields \eqref{sup_norm21s}, 
Now we turn to  \eqref{lippss2}.  it is clear that (denoting $\Delta_{12}\mathfrak{a}:=\mathfrak{a}(\omega_1)-\mathfrak{a}(\omega_2)$),
\[
\frac{\Delta_{12} \mathfrak{c}_{\tau_1,\tau_2}}{|\omega_1-\omega_2|} = \frac{\Delta_{12}\mathfrak{a}}{|\omega_1-\omega_2|}\mathfrak{b}(\omega_1) + \mathfrak{a}(\omega_2)\frac{\Delta_{12} \mathfrak{b}}{|\omega_1-\omega_2|}.
\]
Applying \eqref{sup_norm21s} to $\frac{\Delta_{12}\mathfrak{a}}{|\omega_1-\omega_2|}\mathfrak{b}(\omega_1)$ and $\mathfrak{a}(\omega_2)\frac{\Delta_{12} \mathfrak{b}}{|\omega_1-\omega_2|}$, we obtain that for each $\omega_1,\omega_2$,
\begin{equation}\label{lisd2sdzxc}
\begin{aligned}
|&\frac{\Delta_{12} \mathfrak{c}_{\tau_1,\tau_2}}{|\omega_1-\omega_2|}|_{m+m',s,\eta_0}\\
&\quad \le_{\nu,m,m',s,\eta}|\mathfrak{a}|^{\lip}_{m,s+\mu,\eta}|\mathfrak{b}|^{\sup}_{m',s_0+\mu,\eta}+|\mathfrak{a}|^{\lip}_{m,s_0+\mu,\eta}|\mathfrak{b}|^{\sup}_{m,s+\mu,\eta}\\
& \quad \ +|\mathfrak{a}|^{\sup}_{m,s+\mu,\eta}|\mathfrak{b}|^{\lip}_{m,s_0+\mu,\eta}+|\mathfrak{a}|^{\sup}_{m,s_0+\mu,\eta}|\mathfrak{b}|^{\lip}_{m,s+\mu,\eta}\\
& \quad \le (|\mathfrak{a}|^{\lip}_{m,s+\mu,\eta} + |\mathfrak{a}|^{\sup}_{m,s+\mu,\eta})(|\mathfrak{b}|^{\lip}_{m',s_0+\mu,\eta} +|\mathfrak{b}|^{\sup}_{m',s_0+\mu,\eta}) \\
& \ + (|\mathfrak{a}|^{\lip}_{m,s_0+\mu,\eta} + |\mathfrak{a}|^{\sup}_{m,s_0+\mu,\eta})(|\mathfrak{b}|^{\lip}_{m',s+\mu,\eta} +|\mathfrak{b}|^{\sup}_{m',s+\mu,\eta}).
\end{aligned}
\end{equation}
Noting that $|\gamma|<1$ and taking the supremum in $\omega_1,\omega_2$, we get
\[
\gamma |\mathfrak{c}_{\tau_1,\tau_2}|^{\lip}_{m+m',s,\eta}\le_{m,m',\eta} |\mathfrak{a}|^{\Lip(\gamma,\Omega)}_{m,s+\mu,\eta}|\mathfrak{b}|^{\Lip(\gamma,\Omega)}_{m',s_0+\mu,\eta} + |\mathfrak{a}|^{\Lip(\gamma,\Omega)}_{m,s_0+\mu,\eta}|\mathfrak{b}|^{\Lip(\gamma,\Omega)}_{m',s+\mu,\eta}.
\]
Combining this with \eqref{sup_norm21s} and recalling the definition of the norm in \eqref{omega_dep_norm1}, the desired estimate \eqref{lippss2} follows.
\end{proof}

\begin{lemma}\label{compositionlemma25}
Given $\mathfrak{a}\in \mathcal{S}^{m}$, $\mathfrak{b}\in \mathcal{S}^{m'}$ for some $m,m'\in \mathbb{R}$, let $\sigma^1_{\mathfrak{a},\mathfrak{b}},\sigma^2_{\mathfrak{a},\mathfrak{b}}$ be the symbols given in \eqref{composition_weyl}.
Then, we have $\sigma_{\mathfrak{a},\mathfrak{b}}^1\in \mathcal{S}^{m+m'-2}, \sigma^2_{\mathfrak{a},\mathfrak{b}}\in \mathcal{S}^{m+m'-3}$. Furthermore, for each $\eta\in \mathbb{N}_0$, there exist $\mu=\mu({m,m',\eta})$ and an absolute constant $N>0$, such that
\begin{align*}
 |\sigma^1_{\mathfrak{a},\mathfrak{b}}|_{m+m'-2,s,\eta}^{\Lip(\gamma,\Omega)} &\le_{m,m',s,\eta} |\mathfrak{a}|_{m,s_0+\mu,\eta+N}^{\Lip(\gamma,\Omega)} |\mathfrak{b}|_{m',s+\mu,\eta+N}^{\Lip(\gamma,\Omega)} + |\mathfrak{a}|_{m,s+\mu,\eta+N}^{\Lip(\gamma,\Omega)} |\mathfrak{b}|_{m',s_0+\mu, \eta+N}^{\Lip(\gamma,\Omega)}\\
 |\sigma^2_{\mathfrak{a},\mathfrak{b}}|_{m+m'-3,s,\eta}^{\Lip(\gamma,\Omega)} &\le_{m,m',s,\eta} |\mathfrak{a}|_{m,s_0+\mu,\eta+N}^{\Lip(\gamma,\Omega)} |\mathfrak{b}|_{m',s+\mu,\eta+N}^{\Lip(\gamma,\Omega)} + |\mathfrak{a}|_{m,s+\mu,\eta+N}^{\Lip(\gamma,\Omega)} |\mathfrak{b}|_{m',s_0+\mu, \eta+N}^{\Lip(\gamma,\Omega)}\end{align*}
 Consequently, $r_{\mathfrak{a},\mathfrak{b}}\in \mathcal{S}^{m+m'-2}$ in \eqref{commutator_pseudo} satisfies
 \[
 |r_{\mathfrak{a},\mathfrak{b}}|_{m+m'-2,s,\eta}^{\Lip(\gamma,\Omega)} \le_{m,m',s,\eta} |\mathfrak{a}|_{m,s_0+\mu,\eta+N}^{\Lip(\gamma,\Omega)} |\mathfrak{b}|_{m',s+\mu,\eta+N}^{\Lip(\gamma,\Omega)} + |\mathfrak{a}|_{m,s+\mu,\eta+N}^{\Lip(\gamma,\Omega)} |\mathfrak{b}|_{m',s_0+\mu, \eta+N}^{\Lip(\gamma,\Omega)}
 \] 
\end{lemma}
\begin{proof}
It follows straightforwardly from expressions for $\sigma^1_{\mathfrak{a},\mathfrak{b}},\sigma^2_{\mathfrak{a},\mathfrak{b}}$ in Lemma~\ref{composition_lemmaaa} and the estimates in Lemma~\ref{product_symbols}, and Lemma~\ref{compositionlemma2502}.
\end{proof}

\begin{lemma}\label{compandkskd2sd}
Let $\mathfrak{a}\in \mathcal{S}^m$ and $\mathfrak{b}\in \mathcal{S}^{m'}$ for some $m,m'\in \mathbb{R}$. For each $\eta\in \mathbb{N}_0$, there exist $\mu=\mu(m,m',\eta)$ and an absolute constant $N>0$ such that the followings hold true: 
\begin{enumerate}[label=(\arabic*)]
\item \label{compositsd} For $\mathfrak{c}$ such that $Op^W(\mathfrak{a})\circ Op^W(\mathfrak{b})=Op^W(\mathfrak{c})$, we have  $\mathfrak{c}\in S^{m+m'}$ and for all $s\ge s_0$,
\begin{align*}
|\mathfrak{c}|^{\Lip(\gamma,\Omega)}_{m+m',s,\eta}\le_{m,m',\eta,s} |\mathfrak{a}|^{\Lip(\gamma,\Omega)}_{m,s+\mu,\eta + N}|\mathfrak{b}|^{\Lip(\gamma,\Omega)}_{m',s_0+\mu,\eta + N} +  |\mathfrak{a}|^{\Lip(\gamma,\Omega)}_{m,s_0+\mu,\eta + N}|\mathfrak{b}|^{\Lip(\gamma,\Omega)}_{m',s+\mu,\eta + N }.
\end{align*}
\item \label{symmed} For $\mathfrak{c}$ such that $[Op^W(\mathfrak{a}),Op^{W}(\mathfrak{b})] = Op^W(\mathfrak{c})$, we have $\mathfrak{c}\in \mathcal{S}^{m+m' -1}$ and for all $s\ge s_0$,
 \begin{align*}
|\mathfrak{c}|^{\Lip(\gamma,\Omega)}_{m+m'-1,s,\eta}\le_{m,m',\eta,s}  |\mathfrak{a}|^{\Lip(\gamma,\Omega)}_{m,s+\mu,\eta+N}|\mathfrak{b}|^{\Lip(\gamma,\Omega)}_{m',s_0+\mu,\eta+N} +  |\mathfrak{a}|^{\Lip(\gamma,\Omega)}_{m,s_0+\mu,\eta+N}|\mathfrak{b}|^{\Lip(\gamma,\Omega)}_{m',s+\mu,\eta+N}
\end{align*}
\item \label{symocsdwd} 
For $\mathfrak{c}$ such that $[Op^W(\mathfrak{a}),Op^{W}(\mathfrak{b})]_x = Op^W(\mathfrak{c})$, we have $\mathfrak{c}\in \mathcal{S}^{m+m' }$ and for all $s\ge s_0$,
 \begin{align*}
|\mathfrak{c}|^{\Lip(\gamma,\Omega)}_{m+m',s,\eta}\le_{m,m',\eta,s}  |\mathfrak{a}|^{\Lip(\gamma,\Omega)}_{m,s+\mu,\eta+N}|\mathfrak{b}|^{\Lip(\gamma,\Omega)}_{m',s_0+\mu,\eta+N} +  |\mathfrak{a}|^{\Lip(\gamma,\Omega)}_{m,s_0+\mu,\eta+N}|\mathfrak{b}|^{\Lip(\gamma,\Omega)}_{m',s+\mu,\eta+N}.
\end{align*}
\end{enumerate}
\end{lemma}
\begin{proof}
The result follows from Lemma~\ref{compositionlemma2502}, Lemma~\ref{composition_lemmaaa}, Lemma~\ref{product_symbols}, Lemma~\ref{symsdcodmyysd} and Lemma~\ref{compositionlemma25}.
\end{proof}
For \ref{compositsd} of Lemma~\ref{compandkskd2sd}, we will obtain a finer estimate in case that $\mathfrak{a}$ or $\mathfrak{b}$ is supported on a finite number of Fourier modes:
\begin{lemma}\label{finitemodes_foudsd}
Let $\mathfrak{a} \in \mathcal{S}^m$ and $\mathfrak{b}\in \mathcal{S}^{m'}$. Assume that $\mathfrak{b}=\mathfrak{b}(\omega,\varphi,x,\xi)$ is supported on a finite number of Fourier modes, then we can estimate the symbol of the composition of $Op(\mathfrak{a})$ and $Op(\mathfrak{b})$ without any loss of derivatives. More precisely, if
\begin{align}\label{xxxxx12323sd}
\widehat{\mathfrak{b}}(\omega,l,k,\xi) = 0,\text{ if $|l| + |k| \ge N$ for some $N\ge0$.}
\end{align}
then, for a symbol $\mathfrak{c}$ such that $Op(\mathfrak{a})\circ Op(\mathfrak{b}) = Op(\mathfrak{c})$, we have that
\begin{align}\label{sliisdpwpwpersx}
|\mathfrak{c}|^{\Lip(\gamma,\Omega)}_{m+m',s,\eta}\le_{m,m',s,\eta,N}  |\mathfrak{a}|^{\Lip(\gamma,\Omega)}_{m,s,\eta}|\mathfrak{b}|^{\Lip(\gamma,\Omega)}_{m',s_0,\eta} + |\mathfrak{a}|^{\Lip(\gamma,\Omega)}_{m,s_0,\eta}|\mathfrak{b}|^{\Lip(\gamma,\Omega)}_{m',s,\eta},\text{ for all $s\ge s_0$.}
\end{align}
\end{lemma}
\begin{proof}
The proof is similar to the one of Lemma~\ref{compositionlemma2502}. Using \eqref{compose_289}, we have (omitting $\omega$ dependence for simplicity),
\[
\partial_\xi^\eta\widehat{\mathfrak{c}}^{\varphi,x}(l,k,\xi) = \sum_{n=0}^{\eta}C_{n}\sum_{j_1+j_2 = k,\ l_1+l_2 = l}\partial_\xi^{\eta -n}\widehat{\mathfrak{a}}^{\varphi,x}(l_1,j_1,\xi+{j_2})\partial_\xi^n\widehat{\mathfrak{b}}^{\varphi,x}(l_2,j_2,\xi).
\]
Hence,
\begin{align}
& |\partial_\xi^\eta\widehat{\mathfrak{c}}^{\varphi,x}(l,k,\xi)|^2(\langle l \rangle + \langle k \rangle)^{2s}\\
&\le\sum_{n=0}^{\eta}C_{n} \left| \sum_{j_1+j_2 = k,\ l_1+l_2 = l}\partial_\xi^{\eta-n}\widehat{\mathfrak{a}}^{\varphi,x}(l_1,j_1,\xi+{j_2})\partial_\xi^{n}\widehat{\mathfrak{b}}^{\varphi,x}(l_2,j_2,\xi )\right|^2(\langle l \rangle + \langle k \rangle)^{2s}\nonumber\\
& \le_{s,\eta}\sum_{n=0}^\eta C_n \left| \sum_{\substack{j_1+j_2 = k,\\ l_1+l_2 = l}}(\langle l_1 \rangle + \langle j_1 \rangle)^s|\partial_\xi^{\eta-n}\widehat{\mathfrak{a}}^{\varphi,x}(l_1,j_1,\xi+{j_2})||\partial_\xi^n\widehat{\mathfrak{b}}^{\varphi,x}(l_2,j_2,\xi )|\right|^2\nonumber\\
& \ +\sum_{n=0}^{\eta}C_{n} \left| \sum_{j_1+j_2 = k,\ l_1+l_2 = l}|\partial_\xi^{\eta-n}\widehat{\mathfrak{a}}^{\varphi,x}(l_1,j_1,\xi+{j_2})|(\langle l_2 \rangle + \langle j_2 \rangle)^s|\partial_\xi^n\widehat{\mathfrak{b}}^{\varphi,x}(l_2,j_2,\xi )|\right|^2\nonumber \\
& = A_1 + A_2.\label{rsddsdsdsds222sd}
\end{align}
For $A_1$, we have that
\begin{align}
& \sum_{l\in \mathbb{Z}^\nu,k\in \mathbb{Z}} A_{1} \nonumber \\
& \le \sum_{n=0}^\eta C_n\sum_{k,l}\left| \sum_{j_1+j_2 = k,\ l_1+l_2 = l}(\langle l_1 \rangle + \langle j_1 \rangle)^s|\partial_\xi^{\eta-n}\widehat{\mathfrak{a}}^{\varphi,x}(l_1,j_1,\xi+{j_2})|\right. \nonumber \\
& \qquad \qquad \qquad \qquad \times \left.(\langle l_2\rangle + \langle j_2 \rangle)^{s_0}|\partial_\xi^n\widehat{\mathfrak{b}}^{\varphi,x}(l_2,j_2,\xi)|(\langle l_2\rangle + \langle j_2 \rangle)^{-s_0}\right|^2\nonumber\\
&  \le_{s_0} \sum_{n=0}^{\eta}C_{n}\sum_{j_1,j_2\in \mathbb{Z},\ l_1,l_2\in \mathbb{Z}^{\nu}}\left((\langle l_1 \rangle + \langle j_1 \rangle)^s|\widehat{\partial_\xi^{\eta-n}\mathfrak{a}}(l_1,j_1,\xi+{j_2})|\right. \nonumber  \\
& \qquad \qquad \qquad \qquad \qquad \qquad \times \left.(\langle l_2\rangle + \langle j_2 \rangle)^{s_0}|\widehat{\partial_\xi^n\mathfrak{b}}(l_2,j_2,\xi )|\right)^2\nonumber\\
& \le \sum_{n=0}^{\eta}C_n  \sum_{l_2\in\mathbb{Z}^\nu,j_2\in\mathbb{Z}}\rVert \partial_\xi^{\eta-n}\mathfrak{a}(\cdot,\cdot,\xi+j_2)\rVert_{s}^2\left((\langle l_2\rangle + \langle j_2 \rangle)^{s_0}|\widehat{\partial_\xi^n\mathfrak{b}}(l_2,j_2,\xi )|\right)^2\nonumber\\
&\overset{\eqref{xxxxx12323sd}}{\le} \sum_{n=0}^\eta C_n \sum_{|l_2| + |j_2|\le N} |\mathfrak{a}|_{m,s,\eta}^2\langle \xi + j_2 \rangle^{2(m-(\eta-n))}\left((\langle l_2\rangle + \langle j_2 \rangle)^{s_0}|\widehat{\partial_\xi^n\mathfrak{b}}(l_2,j_2,\xi )|\right)^2\nonumber\\
&\le_{N,m,\eta} \sum_{n=0}^\eta C_n|\mathfrak{a}|_{m,s,\eta}^{2}\langle \xi \rangle^{2(m-(\eta-n))}|\mathfrak{b}|_{m',s_0,\eta}^2\langle \xi \rangle^{2(m' - n)}\nonumber\\
&\le_\eta |\mathfrak{a}|^2_{m,s,\eta_0}|\mathfrak{b}|^2_{m',s_0,\eta}\langle \xi \rangle^{2(m+m' - \eta)}.\label{rlaplsud2ssd}
\end{align}
where the second inequality follows from the Cauchy-Schwarz inequality, and the fourth and fifth inequalities follow from the definition of the norm $|\cdot|_{m,s,\eta}$.
The same computations give us that 
\[
\sum_{l\in \mathbb{Z}^\nu,k\in \mathbb{Z}} A_{2}\le_{s,N,\eta}|\mathfrak{a}|^2_{m,s_0,\eta}|\mathfrak{b}|^2_{m',s,\eta}\langle \xi \rangle^{2(m+m' - \eta)}.
\]
Plugging this and \eqref{rlaplsud2ssd} into \eqref{rsddsdsdsds222sd}, 
\[
\rVert \partial_\xi^{\eta}\mathfrak{c}(\cdot,\cdot,\xi)\rVert_{s}^2\le_{s,N,m,m',\eta} \left( |\mathfrak{a}|_{m,s,\eta}|\mathfrak{b}|_{m',s_0,\eta} + |\mathfrak{a}|_{m,s_0,\eta}|\mathfrak{b}|_{m',s,\eta} \right)^2\langle \xi  \rangle^{2(m+m'-\eta)},
\]
which implies
\[
|\mathfrak{c}|_{m+m',s,\eta} \le_{s,N,m,m',\eta}  |\mathfrak{a}|_{m,s,\eta}|\mathfrak{b}|_{m',s_0,\eta} + |\mathfrak{a}|_{m,s_0,\eta}|\mathfrak{b}|_{m',s,\eta}.
\]
The Lipschitz dependence on $\omega$ can be proved following the same argument as in the proof of Lemma~\ref{compositionlemma2502} (especially following the computations in \eqref{lisd2sdzxc}). Then, \eqref{sliisdpwpwpersx} follows immediately.
\end{proof}

\subsection{Lipschitz-tame operators}\label{s_decay_norm}\index{Lipschitz-tame operators}
We recall from \cite[Section 2]{Berti-Montalto:quasiperiodic-standing-gravity-capillary}, \cite[Section 2]{Feola-Giuliani:quasiperiodic-water-waves} the notion of the Lipschitz-tame estimates of linear operators.\index{Lipschitz-tame estimates}

\begin{definition}[$\mu$-tame operators]\label{def_tame}
For $\mu\ge0$ and $\mathtt{S}\gg s_0$, a linear operator $\mathcal{A}:C^\infty_{\varphi,x}\mapsto C^\infty_{\varphi,x}$ is said to be $\mu$-tame up to $\mathtt{S}$, if there exists a non-decreasing function $[s_0,\mathtt{S}]\ni s\mapsto \mathfrak{M}_{\mathcal{A}}(\mu,s)\ge 0$, such\index{$\mu$-tame up to $\mathtt{S}$} that\index{$\mathfrak{M}_{\mathcal{A}}$}
\begin{align*}
\rVert \mathcal{A}h \rVert_{s}\le \mathfrak{M}_{\mathcal{A}}(\mu,s)\rVert h \rVert_{s_0+\mu} + \mathfrak{M}_{\mathcal{A}}(\mu,s_0)\rVert h \rVert_{s+\mu}, \text{ for all $h\in C^\infty_{\varphi,x}$.}
\end{align*}
$\mathtt{S}$ can be possibly chosen to be $\infty$, in which case, we simply say that the operator $\mathcal{A}$ is $\mu$-tame.
\end{definition}
We also consider Lipschitz dependence of the operators on the frequency $\omega\in \Omega\subset\mathbb{R}^{\nu}$.  To estimate such operators, we recall the following definition:
\begin{definition}[Lip-$\mu$-tame operators]\label{lip_tame}\index{Lip-$\mu$-tame operators}
Let $\mathcal{A}=\mathcal{A}(\omega):C^\infty_{\varphi,x}\mapsto C^\infty_{\varphi,x}$ be a $\omega$-dependent linear operator, for $\omega\in \Omega\subset \mathbb{R}^\nu$. We\index{$\mathfrak{M}^\gamma_{\mathcal{A}}$} denote
\[
\Delta_{12}\mathcal{A}:=\frac{\mathcal{A}(\omega_1)-\mathcal{A}(\omega_2)}{|\omega_1-\omega_2|}.
\]
 For $\mu\ge0$, $\mathtt{S}\gg s_0$ and $\gamma\in (0,1)$, the operator $\mathcal{A}$ is said to be Lip-$\mu$-tame up to $\mathtt{S}$, if there exists a non-decreasing function $[s_0,\mathtt{S}]\ni s\mapsto \mathfrak{M}^{\gamma}_{\mathcal{A}}(\mu,s)\ge 0$, such that
\begin{align*}
\sup_{\omega\in \Omega}\rVert \mathcal{A}(\omega)h \rVert_{s} + \gamma \sup_{\omega_1\ne \omega_2\in \Omega}\rVert \Delta_{12}\mathcal{A}h\rVert_{s}\le \mathfrak{M}^\gamma_{\mathcal{A}}(\mu,s)\rVert h \rVert_{s_0+\mu} + \mathfrak{M}^\gamma_{\mathcal{A}}(\mu,s_0)\rVert h \rVert_{s+\mu}, 
\end{align*}
for all $h\in C^\infty_{\varphi,x}$.
$\mathtt{S}$ can be possibly chosen to be $\infty$, in which case, we simply say that the operator $\mathcal{A}$ is Lip-$\mu$-tame\index{Lip-$\mu$-tame}.
\end{definition}

\subsubsection{Majorant operators and modulo-tame operators} We recall the notion of majorant functions and majorant operators\index{Majorant operators}\index{modulo-tame operators}:
\begin{definition}[Majorant functions/operators]\label{Majorant}
Given a function $u\in C^\infty_{\varphi,x}$, we define its majorant\index{Majorant functions} as
\begin{align}\label{maojisd2sd1}
\underline{u}(\varphi,x):= \sum_{l\in \mathbb{Z}^\nu,\ j\in \mathbb{Z}}|u_j(l)|e^{\ii (l\cdot \varphi + j x)}.
\end{align}
Given a linear operator $\mathcal{A}$, we define its majorant as
\begin{align}\label{maojisd2sd}
\underline{\mathcal{A}}h := \sum_{j_1,j_2\in \mathbb{Z},\ l,l'\in \mathbb{Z}^\nu} |\mathcal{A}_{j_2}^{j_1}(l-l')|h_{j_1}(l')e^{\ii(l\cdot \varphi + j_2 x)}.
\end{align}
\end{definition}
Using the notion of majorant operators, we define modulo-tame operators\index{$\mu$-modulo tame operators}:
\begin{definition}[$\mu$-modulo tame operators]\label{def_tame1}
For $\mu\ge0$ and $\mathtt{S}\gg s_0$, a linear operator $\mathcal{A}:C^\infty_{\varphi,x}\mapsto C^\infty_{\varphi,x}$ is said to be $\mu$-modulo tame up\index{$\mu$-modulo tame up to $\mathtt{S}$} to $\mathtt{S}$, if there exists a non-decreasing function $[s_0,\mathtt{S}]\ni s\mapsto \mathfrak{M}^\sharp_{\mathcal{A}}(\mu,s)\ge 0$, such that
\begin{align*}
\rVert \underline{\mathcal{A}}h \rVert_{s}\le \mathfrak{M}^\sharp_{\mathcal{A}}(\mu,s)\rVert h \rVert_{s_0+\mu} + \mathfrak{M}^\sharp_{\mathcal{A}}(\mu,s_0)\rVert h \rVert_{s+\mu}, \text{ for all $h\in C^\infty_{\varphi,x}$.}
\end{align*}
$\mathtt{S}$ can be possibly chosen to be $\infty$, in which case, we simply say that the operator $\mathcal{A}$ is $\mu$-tame\index{$\mu$-tame}.
\end{definition}

\begin{remark}
Clearly, we have that
\begin{align}\label{endedupspending}
\rVert \underline{u}\rVert_s = \rVert u \rVert_{s},
\end{align}
while we cannot have $\mathfrak{M}^\sharp_{{\mathcal{A}}}(\mu,s) \le \mathfrak{M}_{\mathcal{A}}(\mu,s)$ in general. In Chapter~\ref{reduction},  we will estimate tame constants of $[\partial_{\varphi}^{\vec{b}}(\mathcal{A}),\partial_x]$ for sufficiently large vector $\vec{b}$,  as well as $\mathcal{A}$, to bound the modulo tame constant by tame constant (see Lemma~\ref{modulsdo2ds}).\index{$\mathfrak{M}^\sharp_{\mathcal{A}}$}
\end{remark}
\begin{definition}[Lip-$\mu$-modulo tame operators]\label{lip_tamemodulo}
Let $\mathcal{A}=\mathcal{A}(\omega):C^\infty_{\varphi,x}\mapsto C^\infty_{\varphi,x}$ be\index{Lip-$\mu$-modulo tame operators} a $\omega$-dependent linear operator, for $\omega\in \Omega\subset \mathbb{R}^\nu$. For $\mu\ge0$, $\mathtt{S}\gg s_0$ and $\gamma\in (0,1)$, the operator $\mathcal{A}$ is said to be Lip-$\mu$-modulo tame up to $\mathtt{S}$, if there exists a non-decreasing function $[s_0,\mathtt{S}]\ni s\mapsto \mathfrak{M}^{\sharp,\gamma}_{\mathcal{A}}(\mu,s)\ge 0$, such that\index{$\mathfrak{M}^{\sharp,\gamma}_{\mathcal{A}}$}
\begin{align*}
\sup_{\omega\in \Omega}\rVert \underline{\mathcal{A}(\omega)}h \rVert_{s} + \gamma \sup_{\omega_1\ne \omega_2\in \Omega}\rVert \underline{\Delta_{12}\mathcal{A}}h\rVert_{s}\le \mathfrak{M}^{^\sharp,\gamma}_{\mathcal{A}}(\mu,s)\rVert h \rVert_{s_0+\mu} + \mathfrak{M}^{^\sharp,\gamma}_{\mathcal{A}}(\mu,s_0)\rVert h \rVert_{s+\mu},
\end{align*} 
for all $h\in C^\infty_{\varphi,x}$.
$\mathtt{S}$ can be possibly chosen to be $\infty$, in which case, we simply say that the operator $\mathcal{A}$ is Lip-$\mu$-modulo tame.
\end{definition}
\subsubsection{Properties of the tame operators}
In what follows, we list useful properties of the tame operators. The proofs can be found in \cite[Section 2]{Berti-Montalto:quasiperiodic-standing-gravity-capillary} and \cite[Section 2]{Feola-Giuliani:quasiperiodic-water-waves}.

We use the following notations: For $\varphi\in \mathbb{T}^\nu$, and $\vec{b}\in \mathbb{N}_0^{\nu}$ and $l\in \mathbb{Z}^\nu$, we denote
\[
\partial_{\varphi}^{\vec{b}}:=\partial_{\varphi_1}^{\vec{b}_1}\ldots\partial_{\varphi_\nu}^{\vec{b}_\nu},\quad l^{\vec{b}}:=\Pi_{i=1}^{\nu}(l_i)^{\vec{b}_i},\text{ where $l=(l_1,\ldots,l_\nu),\ \vec{b}=(\vec{b}_1,\ldots,\vec{b}_\nu)$.}
\]

Given $N\in \mathbb{R}^+$, and a linear operator $\mathcal{A}$,  let us denote by $\Pi_{\le N} \mathcal{A}$, the operator\index{$\Pi_{\le N}$}\index{$\Pi_N$} defined as
\begin{align}\label{truncation_1}
(\Pi_N \mathcal{A})^{j_1}_{j_2}(l) := \begin{cases}
 \mathcal{A}^{j_1}_{j_2}(l), & \text{ if $|l| \le N$},\\
 0, & \text{ otherwise},
\end{cases}
\end{align}
where $\mathcal{A}^{j_1}_{j_2}(l)$ is as in \eqref{matrix_rep_2},  and we denote $\Pi_{> N}:= I-\Pi_{\le N }$.

\begin{lemma}\cite[Lemma A.1]{Feola-Giuliani-Procesi:reducibility-weakly-dispersive-degasperis-procesi}\cite[Lemma 2.25]{Berti-Montalto:quasiperiodic-standing-gravity-capillary}\label{consdsdpcosdldsx} Let $\mathcal{A},\mathcal{B}$ be Lip-$\mu_{\mathcal{A}}$-tame and Lip-$\mu_{\mathcal{B}}$-tame operators up to  $\mathtt{S}$. Then, $\mathcal{A}\circ \mathcal{B}$ is Lip-$\left(\mu_{A}+\mu_{B}\right)$-tame operator up to $\mathtt{S}$ with tame constant
\[
\mathfrak{M}_{\mathcal{A}\circ\mathcal{B}}^{\gamma}(\mu_{A}+\mu_{B},s)\le \mathfrak{M}_{\mathcal{A}}^\gamma(\mu_A,s)\mathfrak{M}_{\mathcal{B}}^{\gamma}(\mu_\mathcal{B},s_0+\mu_{A}) +  \mathfrak{M}_{\mathcal{A}}^\gamma(\mu_A,s_0)\mathfrak{M}_{\mathcal{B}}^{\gamma}(\mu_B,s_0+\mu_{A})
\]
If $\mu_{\mathcal{A}}=\mu_{\mathcal{B}}=0$, then we have that  for $\vec{b}\in \mathbb{N}_0^\nu$, 
\begin{equation}\label{chanrul2}
\begin{aligned}
\mathfrak{M}^{\sharp,\gamma}_{\partial_{\varphi}^{\vec{b}}(\mathcal{A}\mathcal{B})}(0,s) &\le C(\vec{b})\left(\mathfrak{M}^{\sharp,\gamma}_{\partial_{\varphi}^{\vec{b}}(\mathcal{A})}(0,s)\mathfrak{M}^{\sharp,\gamma}_{\mathcal{B}}(0,s_0)+\mathfrak{M}^{\sharp,\gamma}_{\partial_{\varphi}^{\vec{b}}(\mathcal{A})}(0,s_0)\mathfrak{M}^{\sharp,\gamma}_{\mathcal{B}}(0,s) \right.\\
&\left. \ +\mathfrak{M}^{\sharp,\gamma}_{\mathcal{A}}(0,s)\mathfrak{M}^{\sharp,\gamma}_{\partial_{\varphi}^{\vec{b}}(\mathcal{B})}(0,s_0)+\mathfrak{M}^{\sharp,\gamma}_{\mathcal{A}}(0,s_0)\mathfrak{M}^{\sharp,\gamma}_{\partial_{\varphi}^{\vec{b}}(\mathcal{B})}(0,s)\right)
\end{aligned}
\end{equation}
\end{lemma}

\begin{lemma}\cite[Lemma A.2]{Feola-Giuliani-Procesi:reducibility-weakly-dispersive-degasperis-procesi}\label{lipdeopsonh}
Let $\mathcal{A}$ be a Lip-$\mu$-tame operator up to $\mathtt{S}$ for some $\mu\ge 0$ and let $h=h(\omega)$ be an $\omega$-dependent function in $C^\infty_{\varphi,x}$. Then,
\begin{align*}
\rVert \mathcal{A}h \rVert_{s}^{\Lip(\gamma,\Omega)} \le_s \mathfrak{M}_{\mathcal{A}}^{\gamma}(\mu,s)\rVert   h \rVert_{s_0+\mu}^{\Lip(\gamma,\Omega)} + \mathfrak{M}_{\mathcal{A}}^{\gamma}(\mu,s_0)\rVert h \rVert_{s+\mu}^{\Lip(\gamma,\Omega)}.
\end{align*} 
\end{lemma}

\begin{lemma}[Pseudo differential operator]\cite[Lemma 2.21]{Berti-Montalto:quasiperiodic-standing-gravity-capillary}\label{rjsdj2jsdsd}\label{symbosldtame} Let $\mathfrak{a}=\mathfrak{a}(\omega,\varphi,x,\xi) \in \mathcal{S}^0$ be a symbol that depends on $\omega,\varphi$. For $s\ge s_0$, $\mathcal{A}:=Op^W(\mathfrak{a})$ is a Lip-$0$-tame operator with tame constant
\[
\mathfrak{M}^{\sharp,\gamma}_{\mathcal{A}}(0,s)\le_s |\mathfrak{a}|^{\Lip(\gamma,\Omega)}_{0,s,0}.
\]
\end{lemma}

\begin{lemma}\cite[Lemma 2.27, 2.28]{Berti-Montalto:quasiperiodic-standing-gravity-capillary}\label{smoothing1123}
For $\vec{\mathtt{b}}\in \mathbb{N}_0^\nu$, it holds that
\begin{align}
\mathfrak{M}^{\sharp}_{\Pi_{> N}\mathcal{A}}(0,s)\le N^{-|\mathtt{b}|}\mathfrak{M}^{\sharp}_{\partial_\varphi^{\vec{b}}\mathcal{A}}(0,s),\quad \mathfrak{M}^{\sharp}_{\Pi_{> N}\mathcal{A}}(0,s)\le \mathfrak{M}^{\sharp}_{\mathcal{A}}(0,s),\label{modulo_tame22323}\\
\mathfrak{M}^{\sharp,\gamma}_{\Pi_{> N}\mathcal{A}}(0,s)\le N^{-|\mathtt{b}|}\mathfrak{M}^{\sharp,\gamma}_{\partial_\varphi^{\vec{b}}\mathcal{A}}(0,s),\quad \mathfrak{M}^{\sharp,\gamma}_{\Pi_{> N}\mathcal{A}}(0,s)\le \mathfrak{M}^{\sharp,\gamma}_{\mathcal{A}}(0,s).\label{lipmodulo_tame22323}
\end{align}
\end{lemma}

\begin{lemma}(Lipschitz-tame to Lipschitz-modulo-tame)\label{modulsdo2ds}
Let $\mathtt{b}_0\in \mathbb{N}$ be fixed. For each $\vec{\mathtt{b}}\in \mathbb{N}^\nu$ such that $|\vec{\mathtt{b}}|\le s_0+\mathtt{b}_0$, assume that  $\partial_{\varphi}^{\vec{\mathtt{b}}}(\mathcal{A})$ and $ [\partial_{\varphi}^{\vec{\mathtt{b}}}(\mathcal{A}),\partial_x]$ are Lip-$0$-tame operators up to $\mathtt{S}$ with tame constants $\mathfrak{M}^\gamma_{\partial_{\varphi}^{\vec{\mathtt{b}}}(\mathcal{A})}(0,s),\mathfrak{M}^\gamma_{[\partial_{\varphi}^{\vec{\mathtt{b}}}(\mathcal{A}),\partial_x]}(0,s)$ for $s\in [s_0,\mathtt{S}]$. Then, for all $\vec{b}\in \mathbb{N}^\nu$  such that $|\vec{b}|\le \mathtt{b}_0$, $\partial_{\varphi}^{\vec{b}}(\mathcal{A})$ is  Lipschitz-0-modulo tame up to $\mathtt{S}$ and
\begin{align}\label{moduks222sdsdd2sd}
\mathfrak{M}^{\sharp,\gamma}_{\partial_{\varphi}^{\vec{b}}(\mathcal{A})}(0,s)\le_{s_0,\mathtt{b}_0} \sum_{\substack{\vec{\mathtt{b}}\in \mathbb{N}^\nu, \\ |\vec{\mathtt{b}}|\le s_0+\mathtt{b}_0}}\left( \mathfrak{M}^\gamma_{\partial_{\varphi}^{\vec{\mathtt{b}}}(\mathcal{A})}(0,s) + \mathfrak{M}^\gamma_{[\partial_{\varphi}^{\vec{\mathtt{b}}}(\mathcal{A}),\partial_x]}(0,s)\right).
\end{align}
\end{lemma}
\begin{proof}
To simplify the notation, let us denote for fixed $\omega,\omega_1,\omega_2\in \Omega$, 
\begin{align}\label{reminiscings1}
\Delta_{12}^{k}\mathcal{A}:=\begin{cases}
\mathcal{A}(\omega), & \text{ if  $k=0$,}\\
\gamma \frac{\mathcal{A}(\omega_1)-\mathcal{A}(\omega_2)}{|\omega_1-\omega_2|}, & \text{ if $k=1$.}
\end{cases}
\end{align}

First, we test the operator $\mathcal{A}$ with a monomial\index{monomial} $ h(\varphi,x):=e^{\ii (l' \cdot \varphi + j' x)}$, then it follows from the definition of $\mathfrak{M}^\gamma_{\partial_{\varphi}^{\vec{\mathtt{b}}}(\mathcal{A})}(0,s)$ that\index{tame constant}
\begin{align*}
\rVert \Delta_{12}^{k}\partial_{\varphi}^{\vec{\mathtt{b}}}(\mathcal{A})h\rVert_s^2 &\le C\left(\left(\mathfrak{M}^\gamma_{\partial_{\varphi}^{\vec{\mathtt{b}}}(\mathcal{A})}(0,s)  \right)^2\rVert h\rVert_{s_0}^2 + \left(\mathfrak{M}^\gamma_{\partial_{\varphi}^{\vec{\mathtt{b}}}(\mathcal{A})}(0,s_0)  \right)^2\rVert h\rVert_{s}^2\right)\\
& = C\left(\left(\mathfrak{M}^\gamma_{\partial_{\varphi}^{\vec{\mathtt{b}}}(\mathcal{A})}(0,s)  \right)^2 \langle l', j'\rangle^{2s_0}  + \left(\mathfrak{M}^\gamma_{\partial_{\varphi}^{\vec{\mathtt{b}}}(\mathcal{A})}(0,s_0)  \right)^2\langle l', j' \rangle^{2s}\right),
\end{align*}
Using the matrix representation of a linear  operator in \eqref{matrix_rep_2},  we have
\[
\left(\Delta_{12}^k\partial_{\varphi}^{\vec{\mathtt{b}}}(\mathcal{A})\right)^{j'}_j(l-l')=\ii^{|\vec{\mathtt{b}}|}(l-l')^{\vec{\mathtt{b}}}(\Delta_{12}^k\mathcal{A})^{j'}_j(l-l').
\]
Therefore,  the above inequality for $\rVert \Delta_{12}^{k}\partial_{\varphi}^{\vec{\mathtt{b}}}(\mathcal{A})h\rVert_s^2 $ gives us
\begin{align}
& \sum_{(l,j)\in \mathbb{Z}^{\nu+1}}\langle l,j\rangle^{2s}|(l-l')^{\vec{\mathtt{b}}}|^2 \left|\left( \Delta_{12}^{k}\mathcal{A}\right)^{j'}_{j}(l-l')\right|^2 \nonumber \\
& \le C\left(\left(\mathfrak{M}^\gamma_{\partial_\varphi^{\vec{\mathtt{b}}}(\mathcal{A})}(0,s)  \right)^2 \langle l', j'\rangle^{2s_0}  + \left(\mathfrak{M}^\gamma_{\partial_\varphi^{\vec{\mathtt{b}}}(\mathcal{A})}(0,s_0)  \right)^2\langle l', j' \rangle^{2s}\right). \label{fixsdsd2sd}
\end{align}
Similarly, using the matrix representation, 
\[
([\Delta_{12}^k\partial_\varphi^{\vec{\mathtt{b}}}(\mathcal{A}),\partial_x])^{j'}_j(l-l') = \ii^{|\vec{\mathtt{b}}| + 1}(l-l')^{\vec{\mathtt{b}}}(j-j')(\Delta_{12}^k\mathcal{A})^{j'}_j(l-l'),
\]
 and using that $[\partial_{\varphi}^{\vec{\mathtt{b}}}(\mathcal{A}),\partial_x]$ is also Lip-$0$-tame, we have that  
\begin{align}
\sum_{(l,j)\in \mathbb{Z}^{\nu+1}}&\langle l,j\rangle^{2s}|(l-l')^{\vec{\mathtt{b}}}|^2|j-j'|^2 \left|\left( \Delta_{12}^{k}\mathcal{A}\right)^{j'}_{j}(l-l')\right|^2\nonumber\\
& \le C\left(\left(\mathfrak{M}^\gamma_{[\partial_{\varphi}^{\vec{\mathtt{b}}}(\mathcal{A}),\partial_x]}(0,s)  \right)^2 \langle l', j'\rangle^{2s_0}  + \left(\mathfrak{M}^\gamma_{[\partial_{\varphi}^{\vec{\mathtt{b}}}(\mathcal{A}),\partial_x]}(0,s_0)  \right)^2\langle l', j' \rangle^{2s}\right).\label{commtsdsd12}
\end{align}
Now, we recall the majorant operator from \eqref{maojisd2sd} and see that for all $h\in C^\infty_{\varphi,x}$,
\begin{align}
& \rVert \underline{\Delta_{12}^k\partial_{\varphi}^{\vec{b}}(\mathcal{A}})h\rVert_s^2 =  \sum_{(l,j)\in \mathbb{Z}^{\nu+1}}\langle l,j\rangle^{2s}\left(\sum_{(l',j')\in \mathbb{T}^{\nu+1}}|(\Delta_{12}^k\partial_{\varphi}^{\vec{b}}(\mathcal{A}))^{j'}_j(l-l')|h_{j'}(l') \right)^2\nonumber\\
&\le \sum_{(l,j)\in \mathbb{Z}^{\nu+1}}\langle l,j\rangle^{2s}\left(\sum_{(l',j')\in \mathbb{T}^{\nu+1}}|(\Delta_{12}^k\mathcal{A})^{j'}_j(l-l')| |(l-l')^{\vec{b}}|\right. \nonumber\\
& \quad \left.\times \langle l-l'\rangle^{s_0}\langle j-j'\rangle|h_{j'}(l')|\frac{1}{\langle l-l'\rangle^{s_0}\langle j-j'\rangle} \right)^2\nonumber\\
&\le _{s_0}\sum_{\substack{(l,j)\in \mathbb{Z}^{\nu+1} \\ (l',j')\in \mathbb{T}^{\nu+1}}}\langle l,j\rangle^{2s}|(\Delta_{12}^k\mathcal{A})^{j'}_j(l-l')|^2\langle l-l'\rangle^{2(s_0+|\vec{b}|)}\langle j-j'\rangle^2|h_{j'}(l')|^2\nonumber\\
&=\sum_{(l',j')\in \mathbb{Z}^{\nu+1}}|h_{j'}(l')|^2\sum_{(l,j)\in \mathbb{T}^{\nu+1}}|(\Delta_{12}^k\mathcal{A})^{j'}_j(l-l')|^2\langle l-l'\rangle^{2(s_0+|\vec{b}|)}\langle j-j'\rangle^2\langle l,j\rangle^{2s},\label{rjsdpo02s2}
\end{align}
where the first inequality follows from the Cauchy-Schwarz inequality with \eqref{range_of_s}. Since $|\vec{b}|\le \mathtt{b}_0$, we have
\begin{align*}
 \langle l-l'\rangle^{2(s_0+|\vec{b}|)}\langle j-j'\rangle^2 &\le_{s_0} \sum_{\substack{\vec{\mathtt{b}}\in \mathbb{N}^\nu, \\ |\vec{\mathtt{b}}|\le s_0+\mathtt{b}_0}}\left(|(l-l')^{\vec{\mathtt{b}}}|^2 + |(l-l')^{\vec{\mathtt{b}}}|^2|j-j'|^2\right),
 \end{align*}
  therefore, \eqref{fixsdsd2sd}, \eqref{commtsdsd12} and \eqref{rjsdpo02s2} give us that
  \begin{align*}
  & \rVert \underline{\Delta_{12}^k\partial_{\varphi}^{\vec{b}}(\mathcal{A}})h\rVert_s^2 \\
  &\le_{s_0,\mathtt{b}_0}\sum_{\substack{\vec{\mathtt{b}}\in \mathbb{N}^\nu, \\ |\vec{\mathtt{b}}|\le s_0+\mathtt{b}_0}}\left( \sum_{(l',j')\in \mathbb{T}^{\nu+1}}|h_{j'}(l')|^2\right. \\
  & \times \left.\left( \langle l',j'\rangle^{2s_0} \left(\left(\mathfrak{M}^\gamma_{\partial_\varphi^{\vec{\mathtt{b}}}(\mathcal{A})}(0,s)\right)^2 + \left(\mathfrak{M}^\gamma_{[\partial_{\varphi}^{\vec{\mathtt{b}}}(\mathcal{A}),\partial_x]}(0,s)  \right)^2 \right)\right. \right. \\
  & \left. \left. +  \langle l',j'\rangle^{2s} \left(\left(\mathfrak{M}^\gamma_{\partial_\varphi^{\vec{\mathtt{b}}}(\mathcal{A})}(0,s_0)\right)^2 + \left(\mathfrak{M}^\gamma_{[\partial_{\varphi}^{\vec{\mathtt{b}}}(\mathcal{A}),\partial_x]}(0,s_0)  \right)^2 \right)  \right)\right)\\
  & \le \sum_{\substack{\vec{\mathtt{b}}\in \mathbb{N}^\nu, \\ |\vec{\mathtt{b}}|\le s_0+\mathtt{b}_0}}\left(\left(\mathfrak{M}^\gamma_{\partial_\varphi^{\vec{\mathtt{b}}}(\mathcal{A})}(0,s)\right)^2 + \left(\mathfrak{M}^\gamma_{[\partial_{\varphi}^{\vec{\mathtt{b}}}(\mathcal{A}),\partial_x]}(0,s)  \right)^2 \right)\rVert h\rVert_{s_0} \\
  & \ + \sum_{\substack{\vec{\mathtt{b}}\in \mathbb{N}^\nu, \\ |\vec{\mathtt{b}}|\le s_0+\mathtt{b}_0}} \left(\left(\mathfrak{M}^\gamma_{\partial_\varphi^{\vec{\mathtt{b}}}(\mathcal{A})}(0,s_0)\right)^2 + \left(\mathfrak{M}^\gamma_{[\partial_{\varphi}^{\vec{\mathtt{b}}}(\mathcal{A}),\partial_x]}(0,s_0)  \right)^2 \right)\rVert h\rVert_{s}.
  \end{align*}
  Recalling the definition of the Lip-$0$-modulo tame operator in Definition~\ref{lip_tamemodulo}, we obtain \eqref{moduks222sdsdd2sd}.
\end{proof}

\begin{lemma}\cite[Lemma 2.26]{Berti-Montalto:quasiperiodic-standing-gravity-capillary}\label{kinverson_modulsd}
Let $\vec{b}\in \mathbb{N}^\nu_0$ and  $\Phi:= I + \mathcal{A}$ where $\mathcal{A},\partial_{\varphi}^{\vec{b}}(\mathcal{A})$ are   Lip-$0$-modulo tame operators up to $\mathtt{S}$ with modulo-tame constants $\mathfrak{M}^{\sharp,\gamma}_{\mathcal{A}}(0,s),\mathfrak{M}^{\sharp,\gamma}_{\partial_\varphi^{\vec{b}}(\mathcal{A})}(0,s)$. Then, there exist  constants $C=C(|\vec{b}|)$ and  $\delta=\delta(|\vec{b}|)>0$ such that if $\mathfrak{M}^{\sharp,\gamma}_{\mathcal{A}}(0,s_0)\le \delta$, then $\Phi$ is invertible and $\tilde{\mathcal{A}}:=\Phi^{-1}-I$ satisfies
\begin{equation}\label{moudksd2sdkspcsd}
\begin{aligned}
\mathfrak{M}^{\sharp,\gamma}_{\tilde{\mathcal{A}}}(0,s)&\le C\mathfrak{M}^{\sharp,\gamma}_{\mathcal{A}}(0,s),\\
\mathfrak{M}^{\sharp,\gamma}_{\partial_\varphi^{\vec{b}}(\tilde{\mathcal{A}})}(0,s)&\le C\left(\mathfrak{M}^{\sharp,\gamma}_{\partial_\varphi^{\vec{b}}(\mathcal{A})}(0,s) +\mathfrak{M}^{\sharp,\gamma}_{\partial_\varphi^{\vec{b}}(\mathcal{A})}(0,s_0)\mathfrak{M}^{\sharp,\gamma}_{\mathcal{A}}(0,s) \right),\text{ for $s\in [s_0,\mathtt{S}]$.}
\end{aligned}
\end{equation}
Furthermore, the same statement holds true, replacing the Lip-$0$-modulo-tame constant by the $0$-modulo tame constant.
\end{lemma}
\begin{proof}
The estimates \eqref{moudksd2sdkspcsd} can be proved in the same way as in \cite[Lemma 2.26]{Berti-Montalto:quasiperiodic-standing-gravity-capillary}, just replacing the $\mathcal{D}^{k_0}$-$0$-modulo-tame constant by the Lip-$0$-tame constant. To obtain the result for the $0$-modulo-tame estimate, we can simply apply the estimates \eqref{moudksd2sdkspcsd} for an operator that does not depend on $\omega$, noting that by definition, $\mathfrak{M}^{\sharp}_{\mathcal{A}}(0,s)=\mathfrak{M}^{\sharp,\gamma}_{\mathcal{A}}(0,s)$, if $\mathcal{A}$ is independent of $\omega$.
\end{proof}

\begin{lemma}\label{stsdjwdsymbodlsd}
Let $\mathcal{A}$ be a Lip-$0$-modulo-tame operator up to $\mathtt{S}$, with tame constant $\mathfrak{M}^{\sharp,\gamma}_{\mathcal{A}}(0,s)$. We define
\[
r(\omega,j):=
\begin{cases}\mathcal{A}^j_j(0),&\text{ for $j\in \mathbb{Z}\backslash\left\{ 0 \right\}$,}\\
0,&\text{ otherwise}.
\end{cases}
\]
Then, it holds that
\[
|r|^{\Lip(\gamma,\Omega)}_{0,0,0}:=\sup_{\omega,\omega_1,\omega_2\in \Omega,\ j\in \mathbb{Z}\backslash \left\{ 0 \right\}} |r(\omega,j)| +  \gamma \frac{|r(\omega_1,j)-r(\omega_2,j)|}{|\omega_1-\omega_2|}  \le_{s_0} \mathfrak{M}^{\sharp,\gamma}_{\mathcal{A}}(0,s_0).
\]
\end{lemma}
\begin{proof}
Using the notation \eqref{reminiscings1}, let us denote
\[
\Delta^{k}_{12}r(j):=\begin{cases}
r(\omega,j), & \text{ if $k=0$,}\\
\gamma \frac{r(\omega_1,j)-r(\omega_2,j)}{|\omega_1-\omega_2|},& \text{ if $k=1$},
\end{cases}
\]
so that $\Delta^k_{12}r(j) = \Delta_{12}^k\mathcal{A}^j_j$. Plugging $l'=0,\vec{\mathtt{b}} = 0,s=s_0$ in \eqref{fixsdsd2sd}, we have
\[
\sum_{(l,j)\in \mathbb{Z}^{\nu+1}}\langle l,j\rangle^{2s_0} \left|\left( \Delta_{12}^{k}\mathcal{A}\right)^{j'}_{j}(l)\right|^2 \le_{s_0}\left(\mathfrak{M}^\gamma_{\mathcal{A}}(0,s_0)  \right)^2 \langle  j'\rangle^{2s_0},
\]
while we obviously have 
\[
\langle j'\rangle^{2s_0}\left|\left(\Delta^{k}_{12}\mathcal{A}\right)^{j'}_{j'}(0)\right|^2 \le \sum_{(l,j)\in \mathbb{Z}^{\nu+1}}\langle l,j\rangle^{2s_0} \left|\left( \Delta_{12}^{k}\mathcal{A}\right)^{j'}_{j}(l)\right|^2.
\]
Therefore, combining the above two inequalities, we obtain $|\Delta^k_{12}r(j) | \le_{s_0} \mathfrak{M}^\gamma_{\mathcal{A}}(0,s_0)\le  \mathfrak{M}^{\sharp,\gamma}_{\mathcal{A}}(0,s_0)  $, which proves the desired estimate, taking the supremum in $\omega,\omega_1,\omega_2$.
\end{proof}

\section{Modified fractional Laplacians: $\Lambda^{\alpha-1}$ and $\Upsilon^{\alpha-3}$.}
For $f\in C^\infty_0$, we define\index{Modified fractional Laplacians}\index{$\Lambda^{\alpha-1}$}\index{$\Upsilon^{\alpha-3}$}
\begin{align}
\Lambda^{\alpha-1}f(x)& := \int_\mathbb{T}(2-2\cos(x-y))^{-\frac{\alpha}{2}}(f(x)-f(y))dy,\nonumber \\\Upsilon^{\alpha-3}f(x) & :=\int_{\mathbb{T}}(2-2\cos(x-y))^{1-\frac{\alpha}{2}}f(y)dy.\label{fractional_1}
\end{align}
We denote the multipliers\index{multipliers} of $\Lambda^{\alpha-1}$ and $\Upsilon^{\alpha-3}$ by 
\begin{align}
m^\circ_{1,\alpha}(j)& := \int_{\mathbb{T}}(2-2\cos y)^{-\frac{\alpha}{2}}(1-e^{-\ii j y})dy, \nonumber \\\quad m^\circ_{2,\alpha}(j) & :=\int_{\mathbb{T}}(2-2\cos y)^{1-\frac{\alpha}{2}}e^{-\ii j y}dy.
\label{multiplier_2}
\end{align}\index{$m^{\circ}_{1,\alpha}$}\index{$m^\circ_{2,\alpha}$}
More explicitly, We have that {\cite[Lemma 2.6 and Eq. (2.1)]{Castro-Cordoba-GomezSerrano:existence-regularity-vstates-gsqg}}
\begin{equation}\label{explicit_multiplier}
\begin{aligned}
\displaystyle m^\circ_{1,\alpha}(j):&=\mathcal{C}_\alpha\left(\frac{\Gamma\left(|j|+\frac{\alpha}{2}\right)}{\Gamma\left(1+|j|-\frac{\alpha}{2}\right)} - \frac{\Gamma\left(\frac{\alpha}{2}\right)}{\Gamma\left(1-\frac{\alpha}{2}\right)} \right),\\ \displaystyle 
 m_{2,\alpha}^\circ(j)& := 2\pi \frac{(-1)^{j}\Gamma(3-\alpha)}{\Gamma(2+|j|-\frac{\alpha}{2})\Gamma(2-|j|-\frac{\alpha}{2})},
\end{aligned}
\end{equation}
where \index{$\mathcal{C}_\alpha$}
\begin{align}\label{Caohsdsdc}
\mathcal{C}_{\alpha}:=-\frac{2\pi\Gamma(1-\alpha)}{\Gamma\left(\frac{\alpha}{2}\right)\Gamma\left(1-\frac{\alpha}{2}\right)}>0.
\end{align}
With $m^\circ_{1,\alpha},m^\circ_{2,\alpha}$ above, we have 
\begin{align}\label{rossandrachel}
\Lambda^{\alpha-1}f(x) = \sum_{j\in \mathbb{Z}}{m^\circ_{1,\alpha}}(j)\widehat{f}_je^{\ii j x},\quad \Upsilon^{\alpha-3}f(x)=\sum_{j\in \mathbb{Z}}m^\circ_{2,\alpha}(j)\widehat{f}_je^{\ii j x}.
\end{align}
One can show that using the asymptotic expansion of Gamma function\index{Gamma function} and Euler's reflection formula\index{Euler's reflection formula} \cite[Theorem 2.1]{Laforgia-Natalini:asymptotic-expansion-ratio-gamma} that for $\alpha\in (0,2)\backslash \left\{ 1\right\}$,
 \begin{align}\label{oosdsd1asmod2s}
 \lim_{|j|\to \infty}m^\circ_{1,\alpha}(j)|j|^{-(\alpha-1)} =C_{1,\alpha} ,\quad  \lim_{|j|\to \infty}m^\circ_{2,\alpha}(j) |j|^{-(\alpha-3)} =C_{2,\alpha},
 \end{align}
 for some constants $C_{1,\alpha}, C_{2,\alpha}$.
 
We will also consider the  weighted operators:
\begin{equation}\label{weighted_operators_1}
\begin{aligned}
&\Lambda^{\alpha-1}_{a}h:=\int_{\mathbb{T}}(2-2\cos(x-y))^{-\frac{\alpha}{2}}a(x,y)(h(x)-h(y))dy,\\
&\Upsilon^{\alpha-3}_{a}h := \int_{\mathbb{T}}(2-2\cos(x-y))^{1-\frac{\alpha}{2}}a(x,y)h(y)dy.
\end{aligned}
\end{equation}
Therefore defining for each $j\in \mathbb{Z}$,
\begin{equation}\label{symbolsoffractions}
\begin{aligned}
\mathfrak{a}_{\Lambda}(x,j)&:= \int_{\mathbb{T}}(2-2\cos y)^{-\frac{\alpha}2}a(x,x-y)(1-e^{-\ii j y})dy,\\
\mathfrak{a}_{\Upsilon}(x,j)&:= \int_{\mathbb{T}}(2-2\cos y)^{1-\frac{\alpha}2}a(x,x-y)e^{-\ii j y}dy,
\end{aligned}
\end{equation}
we have that 
\begin{align}\label{symbolsoffractions2}
\Lambda_{a}^{\alpha-1}h(x) = \sum_{j\in \mathbb{Z}}\mathfrak{a}_\Lambda(x,j)\widehat{h}_je^{\ii j x},\quad  \Upsilon_{a}^{\alpha-3}h(x) = \sum_{j\in \mathbb{Z}}\mathfrak{a}_\Upsilon(x,j)\widehat{h}_je^{\ii j x}.
\end{align}\index{weighted operators}
In the rest of this subsection, we aim to find  smooth symbols to represent the operators $\Lambda^{\alpha-1},\Upsilon^{\alpha-3}$, In order to use the lemmas introduced in Section~\ref{linsesdjsdjsd}. More precisely, we will extend the above symbols to be defined on $\mathbb{T}\times \mathbb{R}$. 

\subsubsection{Analysis of $\Lambda^{\alpha-1}$:}
We study the operator $\Lambda^{\alpha-1}$. We mainly focus on the properties of the multiplier:
\begin{align}\label{lambdadeffsd1}
\lambda^\circ_{\alpha}(j):=j m^\circ_{1,\alpha}(j),\text{ for $j \in \mathbb{Z}$.}
\end{align}
Indeed, \eqref{rossandrachel} tells us that the operator $\partial_x\Lambda^{\alpha-1}$ can be characterized by the symbol $i\lambda^\circ_\alpha(j)$.

Let us denote\index{$\lambda^\circ_{\alpha-1}$}
\begin{align}\label{alpha_1aslongas}
\lambda^\circ_{\alpha-1}(j):=\frac{\Gamma(|j|+\frac{\alpha}{2})}{\Gamma(1+|j|-\frac{\alpha}2)},
\end{align}
so that we can write $\lambda^\circ_\alpha$ as, recalling $m^\circ_{1,\alpha}$ from \eqref{explicit_multiplier}, 
\begin{align}\label{rkaksd2sdsdmental}
\lambda^\circ_\alpha(j) = \mathcal{C}_\alpha j(\lambda^\circ_{\alpha-1}(j)-\lambda^\circ_{\alpha-1}(0)) =  \mathcal{C}_\alpha j \left(\frac{\Gamma(|j|+\frac{\alpha}{2})}{\Gamma(1+|j|-\frac{\alpha}2)} - \frac{\Gamma(\frac{\alpha}{2})}{\Gamma(1-\frac{\alpha}2)}\right).
\end{align}
Note that $\lambda^\circ_{\alpha-1}$ is  defined only on the integers. Recall that the function $z\mapsto \Gamma(z)$ is holomorphic\index{holomorphic} on $\left\{ z\in \mathbb{C}: \text{Re}(z)>0\right\}$. Since we always assume that $\alpha\in (1,2)$,  $\lambda^\circ_{\alpha-1}$ can be naturally extended to a smooth function on $\mathbb{R}\backslash \left\{ 0 \right\}$ as\index{$\lambda_{\alpha-1}$}:
\begin{align}\label{def_inmu}
\lambda_{\alpha-1}(\xi) := \frac{\Gamma(|\xi|+\frac{\alpha}{2})}{\Gamma(1+|\xi|-\frac{\alpha}2)} \text{ for $\xi\in \mathbb{R}$.}
\end{align}

\begin{lemma}\label{secretfreidns}
 For each $\eta\in \mathbb{N}_0$, it holds that
\begin{align}\label{kyleandcartman12}
\sup_{\xi\in \mathbb{R}, |\xi|\ge \frac{1}{2}}\partial_\xi^\eta \lambda_{\alpha-1}(\xi)\langle \xi \rangle^{-(\alpha-1)+\eta} \le_{\alpha,\eta} 1.
\end{align}
\end{lemma}
\begin{proof}
We prove it by usual induction. For $\eta =0$, we recall Gautschi's inequality\index{Gautschi's inequality} \cite[Remark 2.1.1]{Qi-Luo:bounds-ratio-gamma-functions}:
\begin{align}\label{jjsd1sdsd}
x^{1-s} < \frac{\Gamma(x+1)}{\Gamma(x+s)} < (x+s)^{1-s},\text{ for all $x > 0$ and $s\in (0,1)$.}
\end{align}
Plugging $x=\frac{\alpha}{2}-1 + |\xi| >0$ for $|\xi|\ge \frac{1}{2}$ and $s=(2-\alpha)\in(0,1)$, we have that
\begin{align}\label{jsdguasd1sd}
\left( |\xi| - \frac{2-\alpha}2\right)^{\alpha-1} < \lambda_{\alpha-1}(\xi) < \left(|\xi|+\frac{2-\alpha}{2}\right)^{\alpha-1},\text{ for all $\xi\in \mathbb{R}$},
\end{align} This proves \eqref{kyleandcartman12} when $\eta=0$.\index{$\lambda^\circ_\alpha$}

 Now we assume that \eqref{kyleandcartman12} holds true for some $\eta\ge 0$ and aim to prove it for $\partial_\xi^{\eta+1}\lambda_{\alpha-1}$. Since $\xi\mapsto \lambda_{\alpha-1}$ is even, let us assume, without loss of generality, that $\xi > 0$ to obtain \eqref{kyleandcartman12}.  We denote by $\psi^{(\eta)}(\xi)$ be the $\eta$-th polygamma function\index{polygamma function}:
 \begin{equation}\label{polygamma_1231}
 \begin{aligned}
 \psi^{(\eta)}(\xi)& :=-\int_0^1 \frac{t^{\xi-1}}{1-t}(\log{t})^\eta dt, \\ \psi^{(0)}(\xi)& :=\log \xi - \frac{1}{2\xi} - 2\int_0^{\infty}\frac{t}{(t^2 + \xi^2)(e^{2\pi t}-1)}dt, \text{ for $\eta\in \mathbb{N}$, $0<\xi\in \mathbb{R}$}.
 \end{aligned}
 \end{equation}
We collect necessary properties of the polygamma functions \cite[Lemma 1]{Guo-Qi:refinements-lower-bounds-polygamma}: 
\begin{align}
\label{poly1}\psi^{(0)}(\xi) &= \frac{\Gamma'(\xi)}{\Gamma(\xi)}\text{ for all $\xi>0$.}\\
 \label{poly3}  \psi^{(\eta)}(\xi) &= \left(\frac{d}{d\xi}\right)^{\eta}\psi^{(0)}(\xi),\text{ for $\eta,\xi \ge0$.}\\
\label{poly2} \frac{(\eta-1)!}{\xi^\eta} + \frac{\eta!}{2\xi^{\eta+1}} &\le (-1)^{\eta+1}\psi^{(\eta)}(\xi)\le \frac{(\eta-1)!}{\xi^\eta} + \frac{\eta!}{\xi^{\eta+1}},\text{ for $\eta \ge 1$, $\xi>0$.}
\end{align}
From \eqref{poly2} and the definition of $\psi^{(0)}(\xi)$ in \eqref{polygamma_1231}, it is clear that
\begin{align}\label{cartmanandbunny}
|\psi^{(\eta)}(\xi +\alpha/2) - \psi^{(\eta)}(\xi +1 -\alpha/2)|\le_\eta \langle \xi \rangle^{-\eta-1},\text{ for $\xi\ge\frac{1}{2}$ and $\eta \in \mathbb{N}_0$.}
\end{align}
Also, using \eqref{poly1} and $\lambda_{\alpha-1}$ in \eqref{def_inmu}, we have that
\begin{align}\label{increasingalphs}
\frac{d}{d\xi}\lambda_{\alpha-1}(\xi) = \lambda_{\alpha-1}(\xi)(\psi^{(0)}(\xi + \frac{\alpha}2) - \psi^{(0)}(\xi +1 -\frac{\alpha}2)),\text{ for $\xi >0$}.
\end{align}
Therefore, it is straightforward that
\begin{align*}
& \left| \left(\frac{d}{d\xi}\right)^{\eta+1}\lambda_{\alpha-1}(\xi)\right| =\left|\left(\frac{d}{d\xi}\right)^{\eta}\left(\lambda_{\alpha-1}(\xi)(\psi^{(0)}(\xi + \frac{\alpha}2) - \psi^{(0)}(\xi +1 -\frac{\alpha}2)) \right)\right|\\
& \overset{\eqref{poly3}}= \left| \sum_{k=0}^{\eta}\binom{\eta}{k}\left(\left(\frac{d}{d\xi}\right)^{k}(\lambda_{\alpha-1})(\xi) (\psi^{(\eta-k)}(\xi + \alpha/2) - \psi^{(\eta-k)}(\xi+1-\alpha/2)) \right)\right|\\
&\overset{\eqref{cartmanandbunny}}\le_{\eta}\sum_{k=0}^{\eta}\binom{\eta}{k}\left| \left(\frac{d}{d\xi}\right)^{k}(\lambda_{\alpha-1})(\xi)\right| \langle \xi \rangle^{-(\eta-k+1)}\\
&\overset{\eqref{kyleandcartman12}}\le_\eta \sum_{k=0}^\eta\binom{\eta}{k} \langle \xi \rangle^{(\alpha-1)-k}\langle \xi \rangle^{-(\eta-k+1)}\\
&\le_\eta \langle \xi \rangle^{(\alpha-1)-(\eta+1)},
\end{align*}
which gives us \eqref{kyleandcartman12} for $\partial_\xi^{(\eta+1)}\lambda_{\alpha-1}$. This finishes the proof.
\end{proof}

\begin{lemma}\label{jrhsjds2sdswehave}
For $\xi \ge\frac{1}{2}$, $\xi\mapsto \lambda_{\alpha-1}(\xi)$ is strictly monotone-increasing and
\begin{align}\label{sdsds11taxdsd}
\partial_\xi \lambda_{\alpha-1}(\xi) \ge_\alpha \xi^{\alpha-2},\quad \text{ for $\xi \ge\frac{1}{2}$}.
\end{align}
\end{lemma}
\begin{proof}
From \eqref{jsdguasd1sd}, we have that $\lambda_{\alpha-1}(\xi) \ge_\alpha \xi^{\alpha-1}$ for $\xi>\frac{1}{2}$. Also, from \eqref{poly2}, we have $\psi^{(1)}(\xi) \ge \frac{1}{\xi} + \frac{1}{2\xi^2}$ for $\xi > \frac{1}{2}$, especially (since $\alpha>1$) \eqref{poly3} and the mean-value theorem tells us that
\begin{align*}
\psi^{(0)}(\xi + \frac{\alpha}2) - \psi^{(0)}(\xi+1-\frac{\alpha}2) \ge_{\alpha} \frac{1}{\xi}.
\end{align*}
Therefore, the result follows from \eqref{increasingalphs}.
\end{proof}

Now, we study the multipliers $\lambda^\circ_{\alpha}$ and $\lambda_{\alpha-1}$ at  integer points.

\begin{lemma}\label{melnikove_23s}
For all $j,k\in \mathbb{Z}\backslash\left\{ 0 \right\}$ such that $j\ne k$, it holds that
\begin{align}\label{pheobe}
|\lambda^\circ_{\alpha}(j)-\lambda^\circ_\alpha(k)|\ge_\alpha |j-k|\left(|j|^{\alpha-1} + |k|^{\alpha-1}\right).
\end{align}
\end{lemma}
\begin{proof}
Note that $j\mapsto \lambda_\alpha^\circ(j)$ is odd. Therefore, without loss of generality, we only need to prove \eqref{pheobe} for the cases: $j>0>k$ and $ j > k > 0$. 

\vspace{0.5\baselineskip}
\noindent\textit{Proof for  $j>0>k$.} 
We first observe from $\psi^{(0)}$ in \eqref{polygamma_1231} that $\xi\mapsto \psi^{(0)}(\xi)$ is strictly increasing for $\xi >0$. Since  $\Gamma(\xi)$ is strictly positive for $\xi>0$, we have $\lambda_{\alpha-1}(\xi)>0$. Especially, \eqref{increasingalphs} tells us that $\lambda_{\alpha-1}$ is strictly increasing function for $\xi>0$. Therefore, we have  
\begin{align}\label{freidnsd}
\lambda_{\alpha-1}(\xi) > \lambda_{\alpha-1}(0)>0, \text{ for all $\xi \ge 0$.}
\end{align}
 Thus, for $j>0>k$, it follows from \eqref{rkaksd2sdsdmental}  that
\begin{align*}
\lambda_\alpha^\circ(j)  -\lambda_\alpha^\circ(k) & \ge\mathcal{C}_\alpha\left( j\left(\lambda_{\alpha-1}(j)-\lambda_{\alpha-1}(0) \right) + |k|\left(\lambda_{\alpha-1}(k)-\lambda_{\alpha-1}(0) \right)\right) \\
& \ge_{\alpha} j^{\alpha}+|k|^{\alpha},
\end{align*}
where the last inequality follows from Lemma~\ref{jrhsjds2sdswehave}. The above inequality certainly implies \eqref{pheobe}.

\vspace{0.5\baselineskip}
\noindent\textit{Proof for  $j>k>0$.} 
It follows from \eqref{freidnsd} and 
\begin{align*}
\lambda_\alpha^\circ(j) - \lambda_\alpha^\circ(k) & =\mathcal{C}_\alpha \int_k^j \left(\lambda_{\alpha-1}(\xi) - \lambda_{\alpha-1}(0)\right) + \xi \partial_\xi\lambda_{\alpha-1}(\xi)d\xi \\& \ge \mathcal{C}_\alpha \int_k^j \xi \lambda_{\xi}\lambda_{\alpha-1}(\xi)d\xi.
\end{align*}
Using Lemma~\ref{jrhsjds2sdswehave}, we have
\[
\lambda_\alpha^\circ(j) - \lambda_\alpha^\circ(k) \ge_\alpha j^{\alpha}-k^{\alpha},
\]
which implies \eqref{pheobe}.
\end{proof}

\begin{lemma}\label{integralsd2sd}
For all $j,k\in \mathbb{Z}\backslash \left\{ 0 \right\}$ such that $j+k\ne 0$, it holds that
\[
|\lambda^\circ_{\alpha}(j+k) - \lambda^\circ_{\alpha}(j)-\lambda^\circ_\alpha(k)|\ge_{\alpha} (\max\left\{ |j|,|k|\right\})^{\alpha-1}\min\left\{ |j|,|k|\right\}.
\]
\end{lemma}
\begin{proof}
In view of  \eqref{rkaksd2sdsdmental}, it suffices to show that
\begin{align}\nonumber
& |(j+k)\lambda_{\alpha-1}(j+k)-j\lambda_{\alpha-1}(j) -k \lambda_{\alpha-1}(k)|  \\
& \ge_{\alpha}(\max\left\{ |j|,|k|\right\})^{\alpha-1}\min\left\{ |j|,|k|\right\}.\label{joyandrache}
\end{align}
Furthermore, we will assume, without loss of generality, that $j>0$.  Under this assumption, we consider three cases: $j=k>0,\ j> k >0,\ k>j>0$,  and $ j>0>k$. 

\textbf{Case $j= k>0$:} In this case, we have
\begin{align}\label{22xcxc11}
|(j+k)\lambda_{\alpha-1}(j+k)-j\lambda_{\alpha-1}(j) -k \lambda_{\alpha-1}(k)| = 2k(\lambda_{\alpha-1}(2k)-\lambda_{\alpha-1}(k)).
\end{align}
If $k=1$, then it follows from Lemma~\ref{jrhsjds2sdswehave} that
\begin{align}\label{lllxcxcxcxxc1}
2k(\lambda_{\alpha-1}(2k)-\lambda_{\alpha-1}(k)) \ge_\alpha 1 = k^\alpha.
\end{align}
If $k\ge 2$, then, using \eqref{jsdguasd1sd},  we have
\begin{align*}
\lambda_{\alpha-1}(2k) &> \left(2k - \frac{2-\alpha}{2}\right)^{\alpha-1}> k^{\alpha-1}\left(2-\frac{2-\alpha}{2k} \right)^{\alpha-1}\overset{k\ge 2}>k^{\alpha-1}\left(\frac{7}{4}\right)^{\alpha-1},\\
\quad \lambda_{\alpha-1}(k) &< \left(k + \frac{2-\alpha}2\right)^{\alpha-1}\le k^{\alpha-1}\left( 1+ \frac{(2-\alpha)}{2k}\right)^{\alpha-1}\overset{k\ge 2}\le k^{\alpha-1}\left(\frac{5}{4}\right)^{\alpha-1}.
\end{align*}
Hence, it follows that
\[
\lambda_{\alpha-1}(2k)-\lambda_{\alpha-1}(k) >_\alpha k^{\alpha-1}.
\]
Plugging this and \eqref{lllxcxcxcxxc1} into \eqref{22xcxc11}, we get
\begin{align}\label{kequaltojand11}
|(j+k)\lambda_{\alpha-1}(j+k)-j\lambda_{\alpha-1}(j) -k \lambda_{\alpha-1}(k)| >_\alpha k^\alpha,\text{ if $k=j > 0$},
\end{align}
which gives \eqref{joyandrache}.

\textbf{Case $j> k>0$:} Using Lemma~\ref{jrhsjds2sdswehave}, we compute
\begin{align}
& |(j+k)\lambda_{\alpha-1}(j+k)-j\lambda_{\alpha-1}(j) -k \lambda_{\alpha-1}(k)|\nonumber\\
&  = j\left(\lambda_{\alpha-1}(j+k)-\lambda_\alpha(j)\right) + k\left( \lambda_{\alpha-1}(j+k) - \lambda_{\alpha-1}(k)\right)\nonumber\\
& >_\alpha j \int_{j}^{j+k}\partial_\xi \lambda_{\alpha-1}(\xi)d\xi\nonumber\\
& >_\alpha j|j+k|^{\alpha-2}k.
\label{kislex}
\end{align}
Then, we can further compute
\begin{align*}
j|j+k|^{\alpha-2}k\ge_{\alpha}\begin{cases}
jk^{\alpha-1}\ge_\alpha j^{\alpha-1}k,&\text{ if $k < j \le 2k$,}\\
j j^{\alpha-2}k \ge_\alpha j^{\alpha-1}k, & \text{ if $j > 2k$.}
\end{cases}
\end{align*}
Plugging this into \eqref{kislex}, we obtain \eqref{joyandrache}.

\textbf{Case $k> j > 0$:} Changing the role of $j,k$ in the case where $j>k>0$ above, we get \eqref{joyandrache}.

\textbf{Case $j >0 >  k$:} In this case, we replace $k$ by $-k$, then \eqref{joyandrache} follows from the above cases.
\end{proof}
 \color{black}
Finally we summarize the results for the multipliers $m_{1,\alpha}$ and $\lambda_{\alpha}$, which are proper extensions of $m_{1,\alpha}^\circ$ and $\lambda_\alpha^\circ$ to $\mathbb{R}$.

\begin{proposition}\label{nichts2}
There exist real-valued classical symbols  $\lambda_{\alpha}\in \mathcal{S}^{\alpha}, m_{1,\alpha}\in \mathcal{S}^{\alpha-1}$ such that
\begin{enumerate}[label=(\arabic*)]
\item  \label{extension_ofoperators} It holds that 
and  $\Lambda^{\alpha-1} = Op^W(m_{1,\alpha})$ and $\partial_x \Lambda^{\alpha-1} =  Op^W(\ii  \lambda_{\alpha})$ as linear operators on $C^\infty_{0}$ (set of smooth periodic functions with zero-average. See \eqref{rkasd2sdadazzxc}). That is, it holds that \index{$m_{1,\alpha}$} 
\begin{align}\label{def_of23sm}
 m_{1,\alpha}(j) = m_{1,\alpha}^\circ(j)=\mathcal{C}_\alpha\left( \frac{\Gamma\left(|j|+\frac{\alpha}{2}\right)}{\Gamma\left(1+|j|-\frac{\alpha}{2}\right)}-\frac{\Gamma\left(\frac{\alpha}{2}\right)}{\Gamma\left(1-\frac{\alpha}{2}\right)}\right), \text{ for all $j\in \mathbb{Z}\backslash\left\{ 0 \right\}$,}
\end{align}
and
\begin{align}\label{lambdadeffsd}
\lambda_{\alpha}(\xi) = \xi m_{1,\alpha}(\xi),\text{ for all $\xi\in \mathbb{R}$},
\end{align}

\item  \label{extension_ofoperators2} $\xi \mapsto m_{1,\alpha}(\xi)$ is strictly positive, even and monotone increasing.
\item \label{extension_ofoperators3} For all $\xi\in \mathbb{R}$, $ \langle \xi \rangle^{\alpha-1} \le_\alpha  |\partial_\xi \lambda_{\alpha}(\xi)| \le_\alpha \langle \xi \rangle^{\alpha-1}$, therefore $\frac{1}{\partial_{\xi}\lambda_{\alpha}}\in \mathcal{S}^{1-\alpha}$.

\item \label{extension_ofoperators4} For all $j,k\in \mathbb{Z}\backslash \left\{ 0 \right\}$ such that $j+k\ne 0$, it holds that
\[
|\lambda_{\alpha}(j+k) - \lambda_{\alpha}(j)-\lambda_\alpha(k)|\ge_{\alpha,k}|j|^{\alpha-1}.
\]

\item \label{extension_ofoperators5} For all $j,k\in \mathbb{Z}\backslash\left\{ 0 \right\}$ such that $j\ne k$, it holds that
\begin{align}\label{pheobe1}
|j-k|\left(|j|^{\alpha-1} + |k|^{\alpha-1}\right)\le_\alpha |\lambda_{\alpha}(j)-\lambda_\alpha(k)|\le_\alpha |j-k|\left(|j|^{\alpha-1} + |k|^{\alpha-1}\right).
\end{align}
\end{enumerate}
\end{proposition}

\begin{proof}
We recall from \eqref{def_inmu} that $\lambda_{\alpha-1}$ is well-defined for all real numbers and smooth everywhere except for $0$. In view of \eqref{explicit_multiplier}, we define $m_{1,\alpha}$ on $\mathbb{R}$ so that  $m_{1,\alpha}(\xi)$ is strictly positive, even, smooth on $\mathbb{R}$, and monotone increasing for $0 < \xi < \frac{1}{2}$ and 
\[
m_{1,\alpha}(\xi):= 
\mathcal{C}_\alpha\left(\lambda_{\alpha-1}(\xi)-\lambda_{\alpha-1}(0) \right),  \text{ if $|\xi| \ge \frac{1}{2}$,}
\]
and even and smooth in $\mathbb{R}$. Such an extension is possible since $\lambda_{\alpha-1}$ is even, smooth,  strictly increasing and $\lambda_{\alpha-1}(\xi) -\lambda_{\alpha-1}(0)>0$ for $\xi\ge \frac{1}2$, and  which follows from \eqref{freidnsd},  Lemma~\ref{jrhsjds2sdswehave}, Lemma~\ref{secretfreidns}, and $\mathcal{C}_\alpha$ defined in \eqref{Caohsdsdc} is a strictly positive constant  for each $\alpha\in (1,2)$. Once $m_{1,\alpha}(\xi)$ is defined for all $\xi\in \mathbb{R}$, we define $\lambda_{\alpha}:=\xi m_{1,\alpha}$ for all $\xi\in \mathbb{R}$.

By its definition, we have $m_{1,\alpha}=m_{1,\alpha}^\circ$ and $\lambda_\alpha=\lambda_\alpha^\circ$  on $\mathbb{Z}\backslash\left\{ 0 \right\}$, therefore we have $Op^W(m_{1,\alpha})=\Lambda^{\alpha-1}$ and $Op^W(\ii \xi \lambda_\alpha)=\partial_x\Lambda^{\alpha-1}$ on $C^\infty_0$. Thanks to Lemma~\ref{secretfreidns}, we have that  $m_{1,\alpha}\in \mathcal{S}^{\alpha-1}$ and $\lambda_{\alpha}\in \mathcal{S}^\alpha$ are classical symbols. Threfore  \ref{extension_ofoperators} and \ref{extension_ofoperators2} follow.  Since $m_{1,\alpha}$ is strictly positive and monotone increasing, $\partial_\xi \lambda_\alpha(\xi) = m_{1,\alpha}(\xi) + \xi \partial_\xi m_{1,\alpha}(\xi) >0$ for all $\xi\in \mathbb{R}$. Therefore, \ref{extension_ofoperators3} follows from Lemma~\ref{jrhsjds2sdswehave}. Again, since $\lambda_{\alpha}=\lambda^\circ_\alpha$ on $\mathbb{Z}\backslash\left\{ 0 \right\}$, \ref{extension_ofoperators4} and the lower bound in \ref{extension_ofoperators5} follow from Lemma~\ref{melnikove_23s}  and Lemma~\ref{integralsd2sd}, respectively. 

 Now, it suffices to prove the upper bound in \eqref{pheobe1}. Using $\lambda_\alpha\in \mathcal{S}^{\alpha}$ and assuming $j>k$, without loss of generality, we have 
 \begin{align*}
 |\lambda_\alpha(j)-\lambda_\alpha(k)| & =\left| \int_k^j \partial_\xi(\lambda_\alpha(\xi)) d\xi \right| \\
 & \le_\alpha \left|\int_k^j |\xi|^{\alpha-1}d\xi\right| \le_\alpha |j-k| (|j|^{\alpha-1} + |k|^{\alpha-1}),
 \end{align*}
 which gives the upper bound in \eqref{pheobe1}.
\end{proof}

\begin{lemma}\label{kaoxixcc2xc}
 For each $k\in \mathbb{Z}\backslash\left\{ 0 \right\}$, there exists a classical symbol  $\kappa_{k,1-\alpha}(\xi)\in \mathcal{S}^{1-\alpha}$ such that\index{$\kappa_{k,1-\alpha}$}
\begin{align}\label{rjjsdsdjxcxcd}
\kappa_{k,1-\alpha}(\xi) = \frac{1}{\lambda_{\alpha}(\xi)-\lambda_{\alpha}(\xi - k)  -\lambda_{\alpha}(k)},\text{ for all $\xi\in \mathbb{Z}\backslash \left\{ 0,k \right\}$ }.
\end{align}
\end{lemma}
\begin{proof}
Thanks to \ref{extension_ofoperators4} of Proposition~\ref{nichts2} (replacing $j$ by $-j$), we see that $\kappa_{k,1-\alpha}(\xi)$ is well defined for all $\xi\in \mathbb{R}$ except for $\xi =0$ and $\xi = k$ for which the denominator vanishes. Therefore, we can construct $\kappa_{k,1-\alpha}(\xi)$ for $\xi\in \mathbb{R}$ in such a way that
\begin{align*}
\kappa_{k,1-\alpha}(\xi) =
\begin{cases}
 \frac{1}{\lambda_{\alpha}(\xi)-\lambda_{\alpha}(\xi - k)  -\lambda_{\alpha}(k)}, &\text{ if $|\xi-0|\ge \frac{1}{2}$ or $|\xi - k| \ge \frac{1}{2}$ },\\
 0, & \text{ if $\xi =0$ or $\xi = k$}, 
\end{cases}
\end{align*}
and it is smooth in $\mathbb{R}$. In order to see $\kappa_{k,1-\alpha}(\xi)\in \mathcal{S}^{1-\alpha}$, we only need to focus on large $\xi$, that is, it is sufficient to show that
\begin{align}\label{alinasd1xc}
\lim_{|\xi|\to \infty}\partial_\xi^{\eta}\kappa_{k,1-\alpha}(\xi) \langle \xi \rangle^{-(1-\alpha)+\eta}\le_{k,\eta,\alpha} 1.
\end{align}
This follows straightforwardly from the usual induction argument, using that $\lambda_\alpha\in \mathcal{S}^\alpha$. Without loss of generality, let us assume that $k >0$. We will show
\begin{align}\label{seasxcxc}
\sup_{\xi \to \infty }\partial_\xi^{\eta}\kappa_{k,1-\alpha}(\xi) \langle \xi \rangle^{-(1-\alpha)+\eta}\le_{k,\eta,\alpha} 1.
\end{align}
Note that the other necessary estimate $\sup_{\xi \to -\infty }\partial_\xi^{\eta}\kappa_{k,1-\alpha}(\xi) \langle \xi \rangle^{-(1-\alpha)+\eta}\le_{k,\eta} 1$ follows in the same way, hence we will omit the proof.

Towards \eqref{seasxcxc}, we see that for $\xi\gg 0$ and $\eta =0$,
\begin{align}\nonumber
|\kappa_{k,1-\alpha}(\xi)\langle \xi \rangle^{-(1-\alpha)}| & = \left|\frac{\langle \xi \rangle^{\alpha-1}}{\lambda_{\alpha}(\xi) - \lambda_\alpha({\xi - k}) - \lambda_{\alpha}(k)}\right|\\
& \le \frac{1}{\left|\frac{\int_{0}^1\partial_\xi\lambda_\alpha(\xi + (t-1)k)dt}{\langle \xi \rangle^{\alpha-1}}\right| - \frac{\lambda_\alpha(k)}{\langle \xi \rangle^{\alpha-1}}}.\label{sdsxx11xcxc1gcxcwd}
\end{align}
Since $\xi\mapsto m_{1,\alpha}$ is monotone increasing for $\xi\ge 0$ (\ref{extension_ofoperators2} of Proposition~\ref{nichts2}), we see from \eqref{lambdadeffsd} that 
\[
\lim_{\xi \to \infty}\left|\frac{\int_{0}^1\partial_\xi\lambda_\alpha(\xi + (t-1)k)dt}{\langle \xi \rangle^{\alpha-1}}\right| =\lim_{\xi \to \infty}  \frac{\int_{0}^1|\partial_\xi\lambda_\alpha(\xi + (t-1)k)|dt}{\langle \xi \rangle^{\alpha-1}} \ge_\alpha 1,
\]
where the last inequality follows from \ref{extension_ofoperators3} of Proposition~\ref{nichts2}. Therefore, taking $\xi$ to $\infty$ in \eqref{sdsxx11xcxc1gcxcwd}, we obtain \eqref{seasxcxc}, when $\eta=0$. For $\eta+1\ge 1$, we have that for $\xi \ge k+1$,
\begin{align*}
|\partial_\xi^{\eta+1}\kappa_{k,1-\alpha}(\xi)| &= |\partial_\xi^{\eta}(|\kappa_{k,1-\alpha}(\xi)|^2(\partial_\xi \lambda_\alpha(\xi) - \partial_\xi\lambda_\alpha(\xi - k)))|\\
&=\left| \partial_\xi^{\eta}\left(|\kappa_{k,1-\alpha}(\xi)|^2 \int_0^1\partial_\xi^{2}\lambda_{\alpha}(\xi + (t-1)k)dt\right)\right|\\
& = \left|\sum_{\eta_1+\eta_2=\eta}C_{\eta_1,\eta_2}  \partial_\xi^{\eta_1}(\kappa_{k,1-\alpha}(\xi)^2)\int_0^{1}\partial_\xi^{\eta_2+2}\lambda_{\alpha}(\xi + (t-1)k)dt\right|.
\end{align*}
Using the induction hypothesis \eqref{seasxcxc}, we have $|\partial_\xi^{\eta_1}(\kappa_{k,1-\alpha}(\xi)^2)|\le_{\eta,\alpha,k}\langle \xi \rangle^{2(1-\alpha)-\eta_1}$, while \ref{extension_ofoperators3} of  Proposition~\ref{nichts2} gives us that $|\partial_\xi^{\eta_2+2}\lambda_{\alpha}(\xi + (t-1)k)|\le_{\eta,k,\alpha} \langle \xi \rangle^{\alpha-\eta_2-2}$. Therefore, we have that
\[
|\partial_\xi^{\eta+1}\kappa_{k,1-\alpha}(\xi)|\le_{\eta,\alpha,k} \langle \xi \rangle^{-\alpha-\eta}=\langle \xi \rangle^{(1-\alpha)-(\eta+1)}.
\]
This proves \eqref{seasxcxc}  for $\eta+1$.
\end{proof}

\subsubsection{Analysis of $\Upsilon_a^{\alpha-3}$}
We study the operator $\Upsilon_a^{\alpha-3}$ defined in \eqref{weighted_operators_1}, for some $a=a(\omega,\varphi,x,y)$ and $\omega\in \Omega$ and $ (\varphi,x)\in \mathbb{T}^{\nu}\times \mathbb{T}$. Especially, we aim to express the operator $\Upsilon^{\alpha-3}_a$ as a pseudo differential operator with a classical symbol. 

Using \eqref{usualquati} and \eqref{weighted_operators_1}, it is easy to see that
\begin{align}\nonumber
\Upsilon^{\alpha-3}_a & = Op(\mathfrak{a}),\text{ where } \\
\mathfrak{a}(\omega,\varphi,x,j)& := \int_{\mathbb{T}}(2-2\cos y)^{1-\frac{\alpha}2}a(\omega,\varphi,x,x-y)e^{-\ii j y}dy\text{ for $j\in \mathbb{Z}$}. \label{usualsymboldoa3}
\end{align}
Here, the "symbol", $\mathfrak{a}(\omega,\varphi,x,j),$ is defined only for $j\in \mathbb{Z}$. In order to make a use of the lemmas studied in Section~\ref{linsesdjsdjsd}, we aim to find an extension $\mathfrak{a}(\omega,\varphi,x,\xi)$, which is well-defined for all $\xi\in \mathbb{R}$, so that the symbol can be measured by the norms in \eqref{symbol_norms6223}. Of course, a particular choice of an extension does not change the operator, since our pseudo differential operators that we consider in this paper always act on periodic functions. 

 To this end, we pick $\psi_1,\psi_2$ to be  smooth non-negative bump functions\index{bump functions} on $\mathbb{R}$ such that $\text{supp}(\psi_1)\in (-1/4,1/4)$ with $\int_{\mathbb{R}}\psi_1(\xi)d\xi = 1$, and $\text{supp}(\psi_2)\in (-\pi/2,\pi/2)$ with $\psi_2(y)=1$ for $y\in (-\pi/3,\pi/3)$. For such $\psi_1,\psi_2$, we define\index{$S_{\Upsilon}$}\index{$S_{1,\Upsilon}$}\index{$S_{\text{step},\Upsilon}$}\index{$S_{2,\Upsilon}$} $S_{\Upsilon}(a)(\omega,\varphi,x,\xi)$ as
 \begin{align}
S_{1,\Upsilon}(a)(\omega,\varphi,x,\xi) &:= \int_{-\pi}^{\pi}(2-2\cos y)^{1-\frac{\alpha}2}a(\omega,\varphi,x,x-y)\psi_2(y)e^{- \ii \xi y}dy,\label{keypart}\\
S_{\text{step},\Upsilon}(a)(\omega,\varphi,x,\xi)&:=\int_{-\pi}^{\pi}(2-2\cos y)^{1-\frac{\alpha}2}a(\omega,\varphi,x,x-y)(1-\psi_2(y))e^{- \ii j y}dy,\label{stepak1sd}\\
&\text{ for a unique $j\in \mathbb{Z}$ such that $\xi \in (j-1/2,j+1/2]$},\nonumber\\
S_{2,\Upsilon}(a)(\omega,\varphi,x,\xi)&:=\psi_1*S_{\text{step},\Upsilon}(a):= \int_{\mathbb{R}}\psi_1(\xi - \xi_0)  S_{\text{step},\Upsilon}(a)(\omega,\varphi,x,\xi_0)d\xi_0,\label{sd2sdesx}\\
S_{\Upsilon}(a)(\omega,\varphi,x,\xi)&:=S_{1,\Upsilon}(a)(\omega,\varphi,x,\xi) + S_{2,\Upsilon}(a)(\omega,\varphi,x,\xi).\label{def_ofsyma}
 \end{align}
 Note that $S_{\Upsilon}(a)$ is smooth in $\varphi,x,\xi$, if so is $a$, and 
 \[
 S_{\Upsilon}(a)(\omega,\varphi,x,j) = \int_{\mathbb{T}}(2-2\cos y)^{1-\frac{\alpha}2}a(\omega,\varphi,x,x-y)e^{-\ii j y}dy,
 \]
 when $\xi = j$ for some $j\in \mathbb{Z}$. Therefore, \eqref{usualsymboldoa3} tells us that 
 \begin{align}\label{rkpppsdsd1p0p2}
 Op(S_{\Upsilon}(a)) = \Upsilon^{\alpha-3}_a.
 \end{align}
Now, we aim to prove that $S_\Upsilon(a)$ is a classical symbol in $\mathcal{S}^{\alpha-3}$, that is, $\partial_\xi^\eta(S_{\Upsilon}(a))\langle \xi \rangle^{-(\alpha-3)+\eta}  $ is bounded for each $\eta\in \mathbb{N}_0$. We will consider $S_{1,\Upsilon}$ and $S_{2,\Upsilon}$ in \eqref{def_ofsyma} separately. 

We first analyze the symbol $S_{1,\Upsilon}(a)$.
\begin{lemma}\label{s1treat}
For each $\eta\in \mathbb{N}_0$, there exists $\mu=\mu(\eta,\nu)\ge 0$ such that  
\begin{align}\label{a1park1sd}
|S_{1,\Upsilon}(a)|^{\Lip(\gamma,\Omega)}_{\alpha-3,s,\eta} \le_{s,\eta,\alpha,\nu}\rVert a \rVert^{\Lip(\gamma,\Omega)}_{H^{s+\mu}(\mathbb{T}^{\nu + 2})},\text{ for all $s\ge s_0$.}
\end{align}
\end{lemma}
\begin{proof}
Let $\eta\in \mathbb{N}_0$ and $s\ge s_0$ be fixed, and let us choose $s_1,s_2\in \mathbb{N}_0$ such that $s_1+s_2 = s$. For $\omega,\omega_1,\omega_2\in \Omega$, we denote
\begin{align}
\Delta_{12}^k a(\varphi,x,y)&:=\begin{cases}
a(\omega,\varphi,x,y), & \text{ if $k=0$,}\\
\gamma\frac{a(\omega_1,\varphi,x,y)-a(\omega_2,\varphi,x,y)}{|\omega_1-\omega_2|},& \text{ if $k=1$},
\end{cases}\label{deltanotationfora}\\
\Delta_{12}^k S_{1,\Upsilon}(a)(\varphi,x,\xi)&:=\begin{cases}
S_{1,\Upsilon}(a)(\omega,\varphi,x,y), & \text{ if $k=0$,}\\
\gamma\frac{S_{1,\Upsilon}(a)(\omega_1,\varphi,x,\xi)-S_{1,\Upsilon}(a)(\omega_2,\varphi,x,\xi)}{|\omega_1-\omega_2|},& \text{ if $k=1$}.
\end{cases}
\end{align}
Since $a\mapsto S_{1,\Upsilon}(a)$ is linear, we have that 
\begin{align}\label{rjjsdsd1sdsd}
\Delta_{12}^k(\partial_\varphi^{s_1}\partial_x^{s_2}S_{1,\Upsilon}(a))(\varphi,x,\xi) = S_{1,\Upsilon}(\Delta_{12}^k(\partial_\varphi^{s_1}\partial_x^{s_2}a))(\varphi,x,\xi),
\end{align}
while \eqref{keypart} tells us that
\begin{align}\nonumber
& S_{1,\Upsilon}(\Delta_{12}^k(\partial_\varphi^{s_1}\partial_x^{s_2}a))(\varphi,x,\xi) \\
& = \int_{-\pi}^{\pi}|y|^{2-\alpha}\underbrace{\left(\frac{2-2\cos y}{|y|^2}\right)^{1-\frac{\alpha}2}\psi_2(y)\Delta^{k}_{12}\partial_\varphi^{s_1}\partial_x^{s_2}(a(\varphi,x,x-y))}_{=:f(\varphi,x,y)}e^{-\ii \xi y}dy. \label{defsd1sdx1}
\end{align}
Since $\text{supp}(\psi_2)\subset (-\pi / 2, \pi/2)$, it holds that $\text{supp}(f(\varphi,x,\cdot))\subset (-\pi/2,\pi/2)$. Therefore, applying Lemma~\ref{immasd21sd}, which will be proved after the proof of this Lemma, we get
\begin{align}\label{whattobys}
\left|  \partial_\xi^\eta (S_{1,\Upsilon}(\Delta_{12}^k(\partial_\varphi^{s_1}\partial_x^{s_2}a)))(\varphi,x,\xi)\right| \le_{\alpha,\eta} \rVert 
f(\varphi,x,\cdot)\rVert_{H^{\eta + 3}(\mathbb{T})}\langle \xi \rangle^{\alpha-3  - \eta}.
\end{align}
From the definition of $f$ in \eqref{defsd1sdx1}, it is clear that (using the usual Sobolev embedding theorem) we can choose $\mu(\nu)>0$ such that
\begin{align}\label{fsjsdpwosmds}
\rVert f(\varphi,x,\cdot)\rVert_{H^{\eta+3}(\mathbb{T})}\le_{s,\nu}  \rVert a \rVert_{H^{s+\mu(\nu) + \eta +3}(\mathbb{T}^{\nu+2})}^{\Lip(\gamma,\Omega)}.
\end{align}
Therefore, 
\begin{align*}
\rVert  \partial_\xi^{\eta} \Delta^{k}_{12}(S_{1,\Upsilon}(a))(\cdot,\cdot,\xi)\rVert_{s}& \le_s \sum_{\substack{s_1+s_2=s,\\ s_1,s_2\in \mathbb{N}_0}}\rVert \partial_\xi^\eta \Delta_{12}^k\partial_{\varphi}^{s_1}\partial_x^{s_2}S_{1,\Upsilon}(a)\rVert_{L^2_{\varphi,x}}\\
&\overset{\eqref{rjjsdsd1sdsd}}= \sum_{\substack{s_1+s_2=s,\\ s_1,s_2\in \mathbb{N}_0}}\rVert \partial_\xi^\eta S_{1,\Upsilon}(\Delta_{12}^k\partial_{\varphi}^{s_1}\partial_x^{s_2}a)\rVert_{L^2_{\varphi,x}}\\
&\overset{\eqref{whattobys},\eqref{fsjsdpwosmds}}{\le_{s,\alpha,\eta,\nu}} \rVert a \rVert_{H^{s+\mu(\nu) + \eta +3}(\mathbb{T}^{\nu+2})}^{\Lip(\gamma,\Omega)}\langle \xi \rangle^{\alpha-3 - \eta}.
\end{align*}
Therefore, replacing $\mu(\nu) + \eta+3$ by $\mu(\nu,\eta)$, we obtain \eqref{a1park1sd}.
\end{proof}

\begin{lemma}\label{immasd21sd}
Let $f\in C^\infty_c(\mathbb{-\pi,\pi})$ be such that $\text{supp}(f)\subset (-\frac{\pi}2,\frac{\pi}2)$.  We denote
\[
\mathfrak{b}(\xi):= \int_{-\pi}^\pi |y|^{2-\alpha}f(y)e^{-\ii \xi y}dy,\text{ for $(x,\xi)\in \mathbb{T}\times \mathbb{R}$.}
\]
Then, for each $\eta\in \mathbb{N}_0$, $\mathfrak{b}$ satisfies
\begin{align}\label{rkkasd1s22d}
\sup_{\xi\in \mathbb{R}}|\partial_\xi^\eta{\mathfrak{b}}(\xi)|\langle \xi \rangle^{-(\alpha-3)+\eta}\le_{\alpha,\eta} \rVert f \rVert_{H^{\eta + 3}(\mathbb{T})},\text{ for all $\eta\in \mathbb{N}_0$.}
\end{align}
\end{lemma}

\begin{proof}
Since \eqref{rkkasd1s22d} concerns only large $|\xi|$, we will assume that $|\xi|\ge 3$. We also assume that $\xi>0$, since the case where $\xi<0$ can be done in the same manner. We argue by induction. 

\vspace{0.5\baselineskip}
\noindent\textit{Proof for $\eta=0$.} For $\eta=0$, the integration by parts gives us that
\begin{align*}
\mathfrak{b}(\xi) &= \frac{1}{\ii\xi}\int_{-\pi}^{\pi}((2-\alpha)y|y|^{-\alpha}f(y) + |y|^{2-\alpha}f'(y)e^{-\ii \xi y})dy\\
&=-\frac{\ii}{\xi}\int_{-\pi}^{\pi}(2-\alpha)y|y|^{-\alpha}f(y)e^{-\ii \xi y}dy - \frac{1}{\xi^2}\int_{-\pi}^\pi \frac{d}{dy}(|y|^{2-\alpha}f'(y))e^{-\ii \xi y}dy\\
&=: -\frac{\ii}{\xi}\int_{-\pi}^{\pi}(2-\alpha)y|y|^{-\alpha}f(y)e^{-\ii \xi y}dy  + \mathfrak{b}_1(\xi).
\end{align*}
Since $\alpha\in (1,2)$, the second integral already satisfies \eqref{rkkasd1s22d}, that is
\begin{align*}
\sup_{\xi\in \mathbb{R}}|\mathfrak{b}_1(\xi)| \langle \xi \rangle^{-(\alpha-3)}\le \sup_{\xi\in \mathbb{R}}\langle \xi \rangle^{1-\alpha}\rVert f \rVert_{W^{2,\infty}(\mathbb{T})}\le \rVert f \rVert_{H^3(\mathbb{T})}.
\end{align*}
For the first integral, using the change of variables ($\xi y\to y$), we see that 
\begin{align}\label{rlajsdsdxc}
 \frac{\ii}{\xi}\int_{-\pi}^{\pi}y|y|^{-\alpha}f(y)e^{-\ii \xi y}dy= \ii |\xi|^{\alpha}\xi^{-3}\int_{-\xi \pi}^{\xi \pi} y |y|^{-\alpha}f(y/\xi)e^{-\ii y} dy.
\end{align}
We choose $j\in \mathbb{\mathbb{Z}}$ such that $\xi\in (j-1/2,j+1/2]$.  Then, we can write the integral above as
\begin{align*}
\int_{-\xi \pi}^{\xi \pi} y |y|^{-\alpha}f(y/\xi)e^{-\ii y} dy & = \int_{\pi}^{j\pi}y |y|^{-\alpha}f(y/\xi)e^{-\ii y} dy\\
& +  \int_{-j \pi}^{-\pi}y |y|^{-\alpha}f(y/\xi)e^{-\ii y} dy + \mathfrak{b}_2(\xi),
\end{align*}
where $\sup_{\xi\in \mathbb{R}}|\mathfrak{b}_2(\xi)|\le C\rVert f\rVert_{L^\infty(\mathbb{T})}$ for some $C>0$. In view of \eqref{rlajsdsdxc}, it suffices to show that
\begin{align}\label{itsuffices}
\left|\int_{\pi}^{j\pi}y |y|^{-\alpha}f(y/\xi)e^{-\ii y} dy+  \int_{-j \pi}^{-\pi}y |y|^{-\alpha}f(y/\xi)e^{-\ii y} dy \right| \le_\alpha \rVert f \rVert_{H^3(\mathbb{T})}.
\end{align}
 To this end, we write the first integral as
\begin{align*}
& \int_{\pi}^{j\pi}y |y|^{-\alpha}f(y/\xi)e^{-\ii y}dy  \\
&= \int_{ \pi }^{j\pi}y |y|^{-\alpha}f(y/\xi)\frac{1}{-\ii}\frac{d}{dy}(e^{-\ii y} )dy\\
& = \frac{1}{\ii }\int_\pi^{j\pi}  \left( (1-\alpha)y^{-\alpha}f(y/\xi) +\frac{1}{\xi} y^{1-\alpha}f'(y/\xi)\right) e^{- \ii y} dy  +(y^{1-\alpha}f(y/\xi)\frac{1}{-\ii}\bigg|_{y=\pi}^{y=j\pi},
\end{align*}
therefore,
\[
\left | \int_{\pi}^{j\pi}y |y|^{-\alpha}f(y/\xi)e^{-\ii y}dy \right|\le_\alpha \left( \rVert f \rVert_{W^{1,\infty}(\mathbb{T})}\frac{j}{|\xi|}\right) \le C  \rVert f \rVert_{H^{3}(\mathbb{T})},
\]
where the last inequality follows from $\xi\in (j-1/2,j+1/2]$ and the Sobolev embedding theorem.
Similarly, we obtain  $\left| \int_{-j \pi}^{-\pi}y |y|^{-\alpha}f(y/\xi)e^{-\ii y} dy\right| \le_\alpha \rVert f \rVert_{H^{3}(\mathbb{T})}$ which gives \eqref{itsuffices}.

\vspace{0.5\baselineskip}
\noindent\textit{Proof for $\eta \ge 1$.} For $\eta\ge 1$,  integration by parts gives us that
\begin{align}\label{javisdsd1sdstlqkftoasd}
\partial_\xi^\eta\mathfrak{b}(\xi) =- \frac{\ii}{\xi}\int_{-\pi}^{\pi} (2-\alpha+\eta)|y|^{2-\alpha}y^{\eta-1}f(y) + |y|^{2-\alpha}y^{\eta-1}(yf'(y))e^{-\ii \xi y}dy.
\end{align}
Using the induction hypothesis, we have that
\begin{align*}
& \sup_{\xi\in \mathbb{R}}\frac{\left| \int_{-\pi}^{\pi} (2-\alpha+\eta)|y|^{2-\alpha}y^{\eta-1}f(y) + |y|^{2-\alpha}y^{\eta-1}(yf'(y))e^{-\ii \xi y}dy\right|}{ \langle \xi \rangle^{(\alpha-3)-\eta+1}}\\
& \le_{\alpha,\eta} \rVert f \rVert_{H^{\eta + 3}(\mathbb{T})}.
\end{align*}
Taking account of the factor $\xi^{-1}$ in \eqref{javisdsd1sdstlqkftoasd}, this proves the desired estimate \eqref{rkkasd1s22d}.
\end{proof}

Now, we turn to $S_{2,\Upsilon}(a)$ in \eqref{sd2sdesx}. 
\begin{lemma}\label{s1treat2}
For each $\eta\in \mathbb{N}_0$, there exists $\mu=\mu(\eta)\ge 0$ such that
\begin{align*}
|S_{2,\Upsilon}|^{\Lip(\gamma,\Omega)}_{\alpha-3,s,\eta} \le_{s,\eta,\alpha,\nu}\rVert a \rVert^{\Lip(\gamma,\Omega)}_{H^{s+\mu}(\mathbb{T}^{\nu+2})},\text{ for all $s\ge s_0$}.
\end{align*} 
\end{lemma}
\begin{proof}
Let us fix $\xi$ and let $j$ be the unique integer such that $\xi \in (j-1/2,j+1/2]$. From \eqref{stepak1sd},  \eqref{sd2sdesx} and noting that $\xi\mapsto \psi_1(\xi)$ is supported on  $\xi\in (-1/4,1/4)$, it suffices to show that
\begin{align}\label{extensiobsd1ksd}
\rVert S_{\text{step},\Upsilon}(a)(\cdot,\cdot,\cdot,j)\rVert_{s}^{\Lip(\gamma,\Omega)}\le_{\alpha,\nu,s} \rVert a \rVert^{\Lip(\gamma,\Omega)}_{H^{s+\mu}(\mathbb{T}^{\nu+2})}\langle j \rangle^{-\eta},\text{ for all $s\ge s_0$ and $\eta\in \mathbb{N}_0$},
\end{align}
for some $\mu=\mu(\eta)>0$. To this end, we use the notation in \eqref{deltanotationfora}. For $s_1,s_2\in \mathbb{N}_0$ such that $s_1+s_2=s$, we have
\begin{align}
&\Delta^k_{12}\partial_{x}^{s_1}\partial_{\varphi}^{s_2}S_{\text{step},\Upsilon}(a)(\varphi,x,j) =S_{\text{step},\Upsilon}(\Delta^k_{12}\partial_{x}^{s_1}\partial_{\varphi}^{s_2}a) = \int_{-\pi}^{\pi} \Delta^{k}_{12}f(\varphi,x,y)e^{-\ii j y}dy,\nonumber\\
& \Delta^k_{12}f(\varphi,x,y):=(2-2\cos y)^{1-\frac{\alpha}2}(1-\psi_1(y))\partial_{x}^{s_1}\partial_{\varphi}^{s_2}(\Delta^k_{12}a(\varphi,x,x-y)).\label{rasdasd}
\end{align}
From this, we note that $\Delta^k_{12}\partial_{x}^{s_1}\partial_{\varphi}^{s_2}S_{\text{step},\Upsilon}(a)$ is the $j$-th Fourier coefficient of the smooth map $y\mapsto \Delta^k_{12}f(\varphi,x,y)$ for each fixed $\varphi,x$, since $\psi_2(y)=1$ for $y\in (-\pi/3,\pi/3)$ and $\text{supp}(\psi_2)\in (-\pi/2,\pi/2)$ . Therefore, 
\begin{align}\nonumber
|\Delta^k_{12}\partial_{x}^{s_1}\partial_{\varphi}^{s_2}S_{\text{step},\Upsilon}(a)(\varphi,x,j)| &\le \rVert \Delta^k_{12}f(\varphi,x,\cdot)\rVert_{H^{\eta}(\mathbb{T})}\langle j \rangle^{-\eta} \\
& \le \rVert \Delta^{k}_{12}f\rVert_{H^{s+\mu(\nu)+\eta}(\mathbb{T}^{\nu + 2})}\langle j \rangle^{-\eta}, \label{11xxcxcxc2sdsxc}
\end{align}
for some $\mu(\nu)>0$, where the last inequality follows from the usual Sobolev embedding theorem. 
From \eqref{rasdasd}, it is clear that
\[
\rVert \Delta^k_{12}f\rVert_{H^{s+\mu(\nu) + \eta}(\mathbb{T}^{\nu+2})}\le_{\nu} \rVert \Delta^{k}_{12}a\rVert_{H^{s+\mu(\nu) + \eta}(\mathbb{T}^{\nu+2})} \le \rVert a \rVert^{\Lip(\gamma,\Omega)}_{H^{s+\mu(\nu) + \eta}(\mathbb{T}^{\nu+2})}.
\]
Therefore, combining this with \eqref{11xxcxcxc2sdsxc}, we get
\begin{align*}
& \rVert S_{\text{step},\Upsilon}(a)(\cdot,\cdot,\cdot,j)\rVert^{\Lip(\gamma,\Omega)}_s \\
&\le \sup_{\omega,\omega_1,\omega_2\in \Omega}\sum_{s_1+s_2=s}  |\Delta^k_{12}\partial_x^{s_1}\partial_{\varphi}^{s_2}S_{\text{step},\Upsilon}(a)(\cdot,\cdot,j)|_{L^\infty(\mathbb{T}^{\nu+1})}\\
& \le_s \rVert a \rVert^{\Lip(\gamma,\Omega)}_{H^{s+\mu(\nu) + \eta}(\mathbb{T}^{\nu+2})}\langle j \rangle^{-\eta},
\end{align*}
which proves \eqref{extensiobsd1ksd}.
\end{proof}

The following proposition follows straightforwardly from Lemma~\ref{s1treat}, Lemma~\ref{s1treat2}, \eqref{def_ofsyma} and \eqref{rkpppsdsd1p0p2}:

\begin{lemma}\label{lemmasd1sdxcxc}
Given a function $a=a(\omega,\varphi,x,y)$, $S_{\Upsilon}(a)$ defined in \eqref{def_ofsyma} is a classical symbol in  $S^{\alpha-3}$ and it satisfies 
\[
Op(S_{\Upsilon}(a)) = \Upsilon^{\alpha-3}_{a},
\]
with the estimates
\begin{align}\label{estsxzsdxc}
|S_{\Upsilon}(a)|^{\Lip(\gamma,\Omega)}_{\alpha-3,s,\eta} \le_{s,\eta,\alpha,\nu}\rVert a \rVert^{\Lip(\gamma,\Omega)}_{H^{s+\mu}(\mathbb{T}^{\nu + 2})},\text{ for all $s\ge s_0$.}
\end{align}
\end{lemma}

\color{black}
\section{Hamiltonian structure  in $L^2_0$}
Given a Hamiltonian\index{Hamiltonian} $H:L^2_0\mapsto \mathbb{R}\cup\left\{ \infty\right\}$, we define  the $L_0^2$-gradient of $H$ at $f\in L^2_0$\index{$L^2_0$-gradient}, $\nabla H(f)$, as the unique vector field such that 
\begin{align}\label{sjdsdsdxcxc12s}
\int_{\mathbb{T}}\nabla_{L^2_0} H(f)(x)g(x)dx = \frac{d}{dt}H(f+tg)\bigg|_{t=0}.
\end{align}
Similarly, we define the $L^2_0$-Hessian of $H$ at $f$, $\nabla^2H(f)$, as the unique linear operator such that
\begin{align}\label{hessisd2sddefsx}
\int_{\mathbb{T}}\nabla_{L^2_0}^{2}H(f)[h](x)g(x)dx = \frac{d^2}{dsdt}H(f+th+sg)\bigg|_{s=t=0}.
\end{align}
Clearly, the gradient and the Hessian\index{Hessian} of $H$ above are well-defined thanks to the classical Riesz representation theorem\index{Riesz representation theorem}.

We consider the symplectic form\index{symplectic form} $\sigma$ on $L^2_0$, given by
\begin{align}\label{symplectic}
\sigma(f,g) :=\int_{\mathbb{T}}\left(\partial_{x}^{-1}f\right)(x)g(x)dx.
\end{align}

We define the Hamiltonian vector field\index{Hamiltonian vector field} $X_{H}:L^2_0\mapsto L^2_0$ as the unique vector field such that
\begin{align*}
(\nabla_{L^2_0} H(f),h)_{L^2} = \sigma(X_H(f),h),
\end{align*}
Therefore, from \eqref{symplectic}, we see that the Hamiltonian vector field can be written as 
\[
X_{H}(f) =\partial_x\nabla_{L^2_0} H(f).
\]
Let $\Phi:L^2_0\mapsto L^2_0$ be a symplectic transformation\index{symplectic transformation}, that is,
\[
\sigma(\Phi(f),\Phi(g)) = \sigma(f,g), \quad \text{ for all $f,g\in L^2_0$.}
\]
Using the definition of $\sigma$ in \eqref{symplectic}, it is clear that a linear operator $\Phi$ is a symplectic transformation if and only if
\begin{align}\label{symplectic73}
\Phi^{T} \partial_x^{-1}\Phi = \partial_x^{-1}.
\end{align}

\subsection{Hamiltonian operators}

We define a class of linear operators on $C^\infty_{\varphi,x}$, generated by a Hamiltonian\index{Hamiltonian operators}.
\begin{definition}\label{hamiltonian_operatior230}
Let $M=M(\varphi):L^2_0\mapsto L_0^2$ be a $\varphi$-dependent linear operator.  We say that a linear operator $\mathcal{L}:C^\infty_{\varphi,x}\mapsto C^\infty_{\varphi,x}$ of the form
\begin{align}\label{hmtoperator}
\mathcal{L}[h]  = \omega\cdot \partial_\varphi h - \partial_x M[h],
\end{align}
 is Hamiltonian, if $M(\varphi):H_{S^\perp}\mapsto H_{S^\perp}$ is symmetric, that is $M(\varphi)=M(\varphi)^{T}$ for each $\varphi\in \mathbb{T}^\nu$.
\end{definition}

Let $\Phi=\Phi(\varphi)$ be a  symplectic transformation for each $\varphi$. The conjugation\index{conjugation} of the linear operator of the form in \eqref{hmtoperator} with $\Phi$ is also Hamiltonian, that is,
\begin{align*}
\Phi^{-1}\mathcal{L}\Phi = \omega\partial_\varphi - \partial_x{N},\quad N =\omega\cdot\partial_\varphi(\Phi^{T})\partial_x^{-1}\Phi + \Phi^{T}M\Phi.
\end{align*}
Indeed, $N$ is symmetric, since differentiating \eqref{symplectic73}, we have $
\omega\cdot \partial_\varphi(\Phi^{T})\partial_x^{-1}\Phi + \Phi^{T}\partial_x^{-1}\omega\cdot \partial_\varphi(\Phi) = 0$. Therefore we have the following:
\begin{lemma}\label{Hamiltonian_operators_under_symplectic}
If $\mathcal{L}$ is a Hamiltonian operator and $\Phi$ is a symplectic transformation, then $\Phi^{-1}\mathcal{L}\Phi$ is also a Hamiltonian operator.
\end{lemma}
\color{black}

\subsection{Homogeneous Hamiltonians}
 Typically\index{Homogeneous Hamiltonians}, a Hamiltonian $H$ in $L^2_0$ is not well-defined everywhere in $L^2_0$ but its $L^2_0$-gradient, $\nabla_{L^2_0} H:f\mapsto \nabla H(f)$, can be defined as a nonlinear operator between two Sobolev spaces. We specify an elementary notion arising from the $L^2_0$-gradient of a Hamiltonian.
\begin{definition}\label{musmooth}
We say a $k$-linear map $A:(C^\infty_x)^k\mapsto C^\infty_x$ is $\mu$-regular\index{$\mu$-regular},  if 
\begin{align*}
\sup_{f=(f_1,\ldots,f_k)\in (C^\infty_x)^k}\frac{\rVert A(f) \rVert_{H^s_x}}{\sum_{i=1}^{k}\left(\rVert f_i \rVert_{H^{s+\mu}_x}\Pi_{j=1,\ j\ne i}^k\rVert f_j \rVert_{H^{\mu}_x}\right)} <\infty,\text{ for all $s\ge 0$.}
\end{align*}
\end{definition}
For example, a trilinear map $A(f):=\partial_x^2(f)f^2$, is $\mu$-regular, for some $\mu\ge 0$, which can be seen from \eqref{bncal}.
We also specify the notion of a homogeneous Hamiltonian in $L^2_0$:
\begin{definition}\label{mureft}
 We say a Hamiltonian is $\mu$-regular if  $H$ admits a homogeneous expansion,
\[
H:=\sum_{k=2}^{\infty} H_k, 
\]
where $H_k$ is homogeneous of degree $k$ such  that each $\nabla_{L^2_0}H_k$, which is a $(k-1)$-linear map, is $\mu$-regular.
\end{definition}

\color{black}

\section{Time-reversible Hamiltonians}
In this section, we briefly recall some basic properties of a reversible Hamiltonian system\index{reversible Hamiltonians}\index{Time-reversible Hamiltonians}.

\begin{definition}\label{shdj2sdsdsd}
 Let $X$ be a Hilbert space\index{Hilbert space} and $\rho:X\mapsto X$ be an involution\index{involution}, that is, $\rho\circ \rho(f) =f,$ for all $ f\in X$. We say that a Hamiltonian $H$ on $X$ is reversible (or time-reversible) with respect to $\rho$,  if $H\circ \rho = H$. 
 We also say that  $H$ is reversibility preserving\index{reversibility preserving Hamiltonian} if $H\circ \rho = - H$.  Furthermore, we say that a transformation\index{reversibility preserving transformation} $\Phi:X\mapsto X$ is reversibility preserving if $\rho\circ \Phi \circ \rho = \Phi$.
 \end{definition}

\subsection{Reversible Hamiltonian on $L^2_0$.}

We denote by $\rho:L^2_0\mapsto L^2_0$, the involution given as
\begin{align}\label{definiteionofinv}
\rho(f)(x) := f(-x).
\end{align}

\begin{lemma}\label{time_reversibility_map}
Let $G$ be a reversibility preserving  Hamiltonian in $L^2_0$ with respect to an involution $\rho$, and let $\Phi_G^t$ be the unique flow map determined by the Hamiltonian PDE\index{Hamiltonian PDE} at time $t$:
\begin{align}\label{hamiltosnains}
f_t = X_{G}(f).
\end{align}
Then, $\Phi^t_G$ is a reversibility preserving map.
\end{lemma}
\begin{proof}
Since $G\circ \rho=-G$, differentiating $G(\rho(f))=-G(f)$ in $f$, we have $
\rho\left(\nabla_{L^2_0} G(\rho(f))\right) = -\nabla_{L^2_0} G(f).$
Hence, taking $\partial_x$ on both sides, we see that 
\begin{align}\label{Hamitlon_rev}
\rho X_G(\rho f) = X_G(f), \text{ for all $f\in L^2_0$}.
\end{align}
In order to show that $\Phi_G^t$ is reversibility preserving, we need to show that \begin{align}\label{rrss}
\rho \Phi^t_G(\rho(f)) =\Phi^t_G(f).
\end{align}
By its definition, $\Phi_G^t$ solves $\partial_t \Phi^t_G(f) = X_G(\Phi^t_G(f))$ with $\Phi_G^0(f)=f$. Furthermore, $\rho\Phi_G^t(\rho f)$ also solves
\begin{align*}
\partial_t \left( \rho\Phi_G^t(\rho f)  \right) & = \rho X_G(\Phi_G^t(\rho f)) \overset{\eqref{Hamitlon_rev}}=X_G(\rho\Phi_G^t(\rho f)),\\
\text{ with }\rho\Phi^t_G(\rho f)|_{t=0} & = \rho\circ \rho (f) = f. 
\end{align*}
Since the flow map $\Phi^t_G$ is unique, we have $\rho \Phi_G^t (\rho f) = f$ and this implies \eqref{rrss}.\end{proof}

\begin{definition}\cite[Definition 2.2]{Baldi-Berti-Montalto:KAM-quasilinear-airy}\label{reversible_operators2s}
Let
\begin{align}\nonumber
X& :=\left\{ f\in C^\infty_{\varphi,x}: f(-\varphi,-x) = f(\varphi,x)\right\},\\  Y & :=\left\{ f\in C^\infty_{\varphi,x}: f(-\varphi,-x) = - f(\varphi,x)\right\} \label{sdssXandYsax}
\end{align}
We say that a  linear operator $\mathcal{A}:C^\infty_{\varphi,x}\mapsto C^\infty_{\varphi,x}$ is 
\begin{enumerate}
\item reversible if $X\mapsto Y$ and $Y\mapsto X$.
\item reversibility preserving $X\mapsto X,\ Y\mapsto Y$.
\end{enumerate}
\end{definition}
\begin{definition}\label{sjdwdsd}
We say a symbol $\mathfrak{a}$ is a reversible symbol\index{reversible symbol}, if $\partial_x Op^W(\mathfrak{a})$ is a reversible operator\index{reversible operator}. We say $\mathfrak{a}$ is a reversibility preserving symbol\index{reversibility preserving symbol}, if $\partial_xOp^W(\mathfrak{a})$ is  a reversibility preserving operator\index{reversibility preserving operator}. 
\end{definition}

In view of the real operators (Definition~\ref{realtorealop}), symmetric operators, and reversible/reversibility preserving operators (Definition~\ref{sjdwdsd}), one can easily show from the definition of Weyl quantization in \eqref{weyl_quant} that 
\begin{equation}\label{llsd1sdxcxcxcmd19sd2}
\begin{aligned}
&\text{$\partial_xOp^W(\mathfrak{a})$ is  a real operator if $\mathfrak{a}(\varphi,x,\xi) = \overline{a}(\varphi,x,-\xi)$},\\
&\text{$Op^W(\mathfrak{a})$ is a symmetric operator if $\mathfrak{a}(\varphi,x,\xi) = \overline{\mathfrak{a}}(\varphi,x,\xi)$},\\
&\text{$\partial_xOp^W(\mathfrak{a})$ is a reversible operator if $\mathfrak{a}(\varphi,x,\xi) = \mathfrak{a}(-\varphi,-x,-\xi)$}.
\end{aligned}
\end{equation}
For a general linear operator $\mathcal{A}$, we have, using the notation in \eqref{matrix_rep_2}\index{real operator}\index{symmetric operator}\index{reversible operator}, 
\begin{equation}\label{lsdrjdmlRmxdlqhdlsekx}
\begin{aligned}
&\text{$\mathcal{A}$ is a real operator if $\mathcal{A}^{j_1}_{j_2}(l) = \overline{\mathcal{A}^{-j_1}_{-j_2}(-l)}$},\\
&\text{$\mathcal{A}$ is a symmetric operator if $(\mathcal{A})^{j_1}_{j_2}(l) = (\overline{\mathcal{A})^{j_2}_{j_1}(l)}$},\\
&\text{$\mathcal{A}$ is a reversible operator if $(\mathcal{A})^{j_1}_{j_2}(l) = \mathcal{A}^{-j_1}_{-j_2}(-l)$}.
\end{aligned}
\end{equation}

\begin{lemma}\label{reversible_reversbsd1}
If $\mathfrak{a}$ is a reversible symbol, and $\mathfrak{b}$ is a reversibility preserving symbol, then $\mathfrak{a}\star \mathfrak{b}$ is a reversible symbol.
\end{lemma}
\begin{proof}
Recalling the definition of $\star$ from \eqref{syysdsdsd}, we have that
\[
\partial_xOp^W(\mathfrak{a}\star \mathfrak{b}) = \partial_xOp^W(\mathfrak{a})\circ \partial_x Op^W(\mathfrak{b}) -  \partial_xOp^W(\mathfrak{b})\circ \partial_x Op^W(\mathfrak{a}),
\]
where each term on the right-hand side is a composition of a reversible operator and a reversibility preserving operator. Therefore, $\partial_xOp^W(\mathfrak{a}\star \mathfrak{b}) $ is a reversible operator. 
\end{proof}

\subsection{Reversible Hamiltonian on $\mathbb{T}^\nu\times \mathbb{R}^\nu\times H_{S^\perp}$.}
Another important phase space that we will work on is $\mathbb{T}^\nu\times \mathbb{R}^\nu\times H_{S^\perp}$, where $H_{S^\perp}$, a subspace of $L^2_0$ is defined as follows:
\begin{align*}
&H_{S^\perp}:=\left\{ f\in L^2_0: f_j = 0 \text{ if  $j\in S$},\right\},\\
&S^+=\left\{ 0< j_1,\ldots,<j_\nu \in \mathbb{N}\right\},\ S:=S\cup(-S),\quad  S^{\perp}:=\mathbb{Z}\backslash\left( S\cup \left\{ 0\right\}\right),
\end{align*}
for a finite subset $S\subset \mathbb{Z}$. Let $H:\mathbb{T}^\nu\times  \mathbb{R}^\nu\times H_{S^\perp}\mapsto \mathbb{R}\cup \left\{ \infty\right\}$ be a Hamiltonian, and let a symplectic two form $\mathcal{W}$ be such that 
\[
\mathcal{W}:=\sum_{i=1}^\nu d_{\theta_i}\wedge d_{y_i}\oplus \sigma_{S^{\perp}},
\]
where $\sigma_{S^\perp}$ is the restriction of $\sigma $ to $H_{S^\perp}$, that is,
$\sigma_{S^{\perp}}(z,z) = \int_{\mathbb{T}} \left( \partial_x^{-1}z\right)(x)z(x)dx$. Furthermore, let us assume that $H$ is reversible with respect to $\rho_*$, defined as
\begin{align}\label{definiteionofinv2}
\rho_*(\theta,y,z):=(-\theta,y,\rho(z)),
\end{align} where $\rho$ is as in \eqref{definiteionofinv}. In the rest of this section, we recall some basic properties of a dynamical system induced by a reversible Hamiltonian $H$ with respect to the symplectic form $\mathcal{W}$ and the involution $\rho_*$.

\begin{definition}\label{reversubsd1nsdsd}
We say a smooth embedding $i:\mathbb{T}^\nu\mapsto \mathbb{T}^\nu\times \mathbb{R}^\nu\times H_{S^\perp}$, $i(\varphi):=(\theta(\varphi),y(\varphi),z(\varphi))$ is reversible\index{reversible embedding} if  $i(-\varphi)=\rho_* i(\varphi)$ for all $\varphi\in \mathbb{T}^\nu$, that is,
\[
\theta(-\varphi)=-\theta(\varphi),\ y(-\varphi)=y(\varphi), \ z(-\varphi)=\rho z(\varphi).
\]
\end{definition}
We denote 
\begin{equation}\label{space_reversibles}
\begin{aligned}
X_i&:=\left\{i:\mathbb{T}^\nu\mapsto \mathbb{T}^\nu\times \mathbb{R}^{\nu}\times H_{S^\perp}: \rho_*i(\varphi) = i(-\varphi))\right\},\\
 Y_i& :=\left\{i:\mathbb{T}^\varphi\mapsto \mathbb{T}^\nu\times \mathbb{R}^{\nu}\times H_{S^\perp}: \rho_*i(\varphi) = -i(-\varphi))\right\},\\
X_{e}&:=\left\{ f:\mathbb{R}^\nu\mapsto \mathbb{R}^\nu: f(\varphi)=f(-\varphi)\right\}, \, Y_{o}:=\left\{ f:\mathbb{R}^\nu\mapsto \mathbb{R}^\nu: f(\varphi)=-f(-\varphi)\right\},\\
X_\perp&:=\left\{ z\in H_{S^\perp}: z(\varphi,x) = z(-\varphi,-x)\right\},\, Y_\perp:=\left\{ z\in H_{S^\perp}: z(\varphi,x) = -z(-\varphi,-x)\right\}.
\end{aligned}
\end{equation}
It is clear that if $i=(\theta,y,z)\in X_i$ then $\theta\in Y_o,y\in X_e$ and $z\in X_\perp$.

In view of \eqref{sjdsdsdxcxc12s} and \eqref{hessisd2sddefsx}, we define $\nabla_zH(i),\ \nabla_z^2H(i)[z_0]\in H_{S^\perp}$ for $z_0\in H_{S^\perp}$ be such that
\begin{align*}
\int_{\mathbb{T}}\nabla_zH(i)(x)z(x)dx &= \frac{d}{dt}H(i+t(0,0,z(x)))\bigg|_{t=0},\\ \int_{\mathbb{T}}\nabla_z^2H(i)[z_0](x)z(x)dx&=\frac{d^2}{dtds}H(i+(0,0,tz_0+sz))\bigg|_{t=s=0},
\end{align*}
 for all $z\in H_{S^\perp}$. 
\begin{lemma}\label{reversible_sdesd}
Let $\omega\in \mathbb{R}^\nu$ and let  $H$ be a reversible Hamiltonian on $\mathbb{T}^\nu\times \mathbb{R}^\nu\times H_{S^\perp}$. Then, the linear map $L: h\mapsto \omega\cdot \partial_\varphi h - \partial_x \nabla_z^2H(\bar{i}(\varphi))[h]$ is reversible.
\end{lemma}
\begin{proof}
In view of Definition~\ref{reversible_operators2s}, we need to show that $L$ maps $X$ to $Y$, and $Y$ to $X$. We will show that $L$ maps $X$ to $Y$ only, because the other case can be done in the same way.
Since $H$ is reversible, we have that for all $\varphi\in\mathbb{T}^\nu,t\in \mathbb{R}$ and $h\in X$,
\[
H(\varphi,0,th(\varphi)) = H(\rho_*(\varphi,0,th(\varphi))) = H(-\varphi,0,t\rho(h)(-\varphi)).
\]
Differentiating in $t$ at $t=0$, we have
\[
\nabla_z^2H(\varphi,0,0)[h(\varphi)] = \rho(\nabla_z^2H(-\varphi,0,0)[\rho(h)(-\varphi)]).
\]
Using $\rho(h)(\varphi) = h(-\varphi)$, and $\partial_x \rho(h)(\varphi) = - \rho(\partial_x h)(\varphi)$, we see that 
\begin{align}\label{revsd22323s}
\partial_x\nabla_z^2H(\bar{i}(\varphi))[h(\varphi)]=\partial_x\rho( \nabla_z^2H(\bar{i}(-\varphi))[h(-\varphi)]) = -\rho(\partial_x\nabla_z^2H(\bar{i}(-\varphi))[h(-\varphi)]).
\end{align}
Furthermore, for $h\in X$, that is, $h(\varphi,x)=h(-\varphi,-x)$, we have that $\omega\cdot\partial_\varphi h(-\varphi,-x) = - \omega\cdot \partial_\varphi h(\varphi,x)$. Thus, it follows that
\begin{align*}
Lh(-\varphi,-x)& = \omega\cdot\partial_\varphi h(-\varphi,-x) - \rho(\partial_x\nabla_z^2H(\bar{i}(-\varphi))[h(-\varphi)])\\
&  \overset{\eqref{revsd22323s}}{=}-\omega\cdot\partial_\varphi h(\varphi,x)  +\partial_x\nabla_z^2H(\bar{i}(\varphi))[h(\varphi)]\\
& =-Lh(\varphi,x).
\end{align*}
Hence, $Lh\in Y$, which proves $L$ maps $X$ to $Y$.
\end{proof}


\section{Translation invariance}
In this section, we collect properties of  a Hamiltonian in $L^2_0$ that is invariant under\index{$\frac{2\pi}{\mathtt{M}}$-translation invariance preserving transformation} $\frac{2\pi}{\mathtt{M}}$-translation for some $\mathtt{M}\in \mathbb{N}$.

 For $\mathtt{M}\in \mathbb{N}$, let us denote
\begin{align}\label{invariant_m}
\rho_{\mathtt{M}}(f)(x):= f(x + \frac{2\pi}{\mathtt{M}} ), \text{ for $f\in C^\infty_x$,}\quad X_{\mathtt{M}}:=\left\{ f\in L^2_0: \rho_{\mathtt{M}}(f) =f\right\}.
\end{align}
It is clear that 
\begin{align}\label{mfosd2sd}
f\in X_{\mathtt{M}} \iff f(x) = \sum_{j \in \mathbb{Z}\backslash \left\{ 0 \right\}} \widehat{f}_{\mathtt{M}j}e^{\ii \mathtt{M} j x},
\end{align} 
that is, the Fourier coefficients of $f$ are supported only in the modes of multiples of $\mathtt{M}$.\begin{definition}\label{mfold_symm_def}
For $\mathtt{M}\in \mathbb{N}$, we say that a Hamiltonian $H$ in $L^2_0$  is invariant under a $\frac{2\pi}{\mathtt{M}}$-translation if  $
H(\rho_{\mathtt{M}}(f)) = H(f).$ Furthermore, we say that a transformation $\Phi:C^\infty_x\mapsto C^\infty_x$ is $\frac{2\pi}{\mathtt{M}}$-translation invariance preserving, if $\rho_{\mathtt{M}}^{-1}\circ\Phi\circ\rho_{\mathtt{M}}(f) = \Phi(f)$. 
\end{definition}

\begin{lemma}
If a linear map $\mathcal{A}:L^2_0\mapsto L^2_0$ is $\frac{2\pi}{\mathtt{M}}$-translation invariance preserving, then $\mathcal{A}:X_{\mathtt{M}}\mapsto X_{\mathtt{M}}$.
\end{lemma}
\begin{proof}
  For $f\in X_{\mathtt{M}}$, we have $\rho_{\mathtt{M}}\mathcal{A}(f) = \rho_{\mathtt{M}}^2 \rho_{\mathtt{M}}^{-1}\mathcal{A}(\rho_\mathtt{M}f) = \rho_{\mathtt{M}}^2 \mathcal{A}(f)$. Therefore, we have $\mathcal{A}=\rho_{\mathtt{M}}\mathcal{A}(f)$, which implies $\mathcal{A}f\in X_{\mathtt{M}}$.
\end{proof}

\begin{lemma}\label{minvariantce}
Let $G$ be a  $\frac{2\pi}{\mathtt{M}}$-translation invariant Hamiltonian\index{$\frac{2\pi}{\mathtt{M}}$-translation invariant Hamiltonian}, and let $\Phi^t_G$ be the unique flow map determined by the Hamiltonian PDE at time $t$:
\begin{align}\label{msjd2ds}
f_t = X_G(f).
\end{align}
Then, $\Phi^t_G$ is $\frac{2\pi}{\mathtt{M}}$-translation invariance preserving.
\end{lemma}
\begin{proof}
Since $G$ is $\frac{2\pi}{\mathtt{M}}$-translation invariant, it follows that $\rho_{\mathtt{M}}^{-1}\nabla G(\rho_{\mathtt{M}}(f))$, therefore,
\begin{align}\label{wearenoty}
\rho_{\mathtt{M}}^{-1}X_{G}(\rho_{\mathtt{M}}(f)) = X_{G}(f).
\end{align}
Hence, we have that
\begin{align*}
\frac{d}{dt}(\rho_{\mathtt{M}}^{-1}\circ \Phi_G^t \circ \rho_\mathtt{M}(f)) & = \rho_{\mathtt{M}}^{-1}X_G(\Phi_G^t\circ \rho_\mathtt{M} (f)) = \rho_\mathtt{M}^{-1}X_G(\rho_\mathtt{M}\circ \rho_\mathtt{M}^{-1}\circ\Phi_G^t\circ \rho_\mathtt{M}(f)) \\& \overset{\eqref{wearenoty}}=X_G(\rho_{\mathtt{M}}^{-1}\circ \Phi_G^t \circ \rho_\mathtt{M}(f)),
\end{align*}
which solves \eqref{msjd2ds}. Since $\rho_{\mathtt{M}}^{-1}\circ \Phi_G^t \circ \rho_\mathtt{M}(f)\bigg|_{t=0}= f$, we have $\rho_{\mathtt{M}}^{-1}\circ \Phi_G^t \circ \rho_\mathtt{M}(f) = f$ for all $t\ge 0$.
\end{proof}

\subsection{$\frac{2\pi}{\mathtt{M}}$-translation invariant Hamiltonian on $\mathbb{T}^\nu\times \mathbb{R}^\nu\times H_{S^\perp}$.}
We denote
\begin{align}\nonumber
X_{i,\mathtt{M}}& := \left\{ i=(\theta,y,z):\mathbb{T}^\nu\mapsto \mathbb{T}^\nu \times \mathbb{R}^\nu \times H_{S^\perp}: \rho_\mathtt{M}(z(\varphi)) = z(\varphi)\right\},\\
\rho_{\mathtt{M},*}(i)& :=(\theta,y,\rho_{\mathtt{M}}(z)).\label{mforlsd1sd}
\end{align}

\begin{definition}\label{dekfhsspdiw11}
We say that a Hamiltonian $H$ on $\mathbb{T}^\nu\times \mathbb{R}^\nu\times H_{S^\perp}$ is $\frac{2\pi}{\mathtt{M}}$-translation invariant if $H(\theta,y,\rho_{\mathtt{M}}z) = H(\theta,y,z)$. We say an embedding $i:\mathbb{T}^\nu\mapsto \mathbb{T}^\nu\times \mathbb{R}^\nu\times H_{S^\perp}$ is $\frac{2\pi}{\mathtt{M}}$-translation invariant, if $\rho_{\mathtt{M},*}(i) = i$. Lastly, we say that a symbol $\mathfrak{a}(x,\xi)$ is $\frac{2\pi}{\mathtt{M}}$-translation invariant if $x\mapsto \mathfrak{a}(x,\xi)$ is  $\frac{2\pi}{\mathtt{M}}$-translation invariant\index{$\frac{2\pi}{\mathtt{M}}$-translation invariant symbol}\index{$\frac{2\pi}{\mathtt{M}}$-translation invariant embedding}.
\end{definition}

\color{black}


\chapter{Hamiltonian structure of the gSQG equations}\label{11jshsdsss2wxowd2}
\section{Hamiltonian equation in the patch setting}\label{11Hamitons1}
We consider the Hamiltonian structure of the generalized SQG equations in the patch setting.  We recall the symplectic structure  $(L^2_0,\sigma)$, where the symplectic $2$-form $\sigma$ is defined as in \eqref{symplectic}.

Throughout the paper, we consider a star-shaped domain\index{star-shaped domain} $D$ defined as
\begin{align*}
D:=\left\{ r(\cos\theta,\sin\theta)\in \mathbb{R}^2 : 0\le r  < R(\theta) \right\}.
\end{align*}
We denote by $z(\theta)$ and $\Psi$ the boundary of  $D$ and  stream function\index{stream function} determined by $1_D$ respectively, that is,
\begin{align}\label{stream_function}
&z(\theta):=R(\theta)(\cos(\theta),\sin(\theta)),\quad \text{ and }\quad \Psi(x):= 1_D * \frac{1}{|x|^\alpha}.
\end{align}
Recall that the patch solutions\index{patch solution} for the gSQG equations in \eqref{gSQG_1} can be written in terms of $R(t,\theta)$ as
\begin{align}\label{gSQG}
R(t,\theta)\partial_tR(t,\theta) = \nabla\Psi(z(t,\theta))\cdot \partial_\theta z(t,\theta).
\end{align}
In order to rewrite \eqref{gSQG} in terms of the Hamiltonian vector field, we define ${H}:L^2\mapsto \mathbb{R}\cup\left\{ \infty\right\}$  as\index{$H$}
\begin{align}\label{Hamiltonian_2}
H(R) &:= \int_{D}\Psi(x)dx\nonumber\\
& = \int_{D}\int_{D}\frac{1}{|x-y|^{\alpha}}dxdy \nonumber\\
&= \int_{\mathbb{T}}\int_{\mathbb{T}}\int_0^{R(\theta)}\int_{0}^{R(\eta)}\frac{1}{|r(\cos \theta,\sin\theta)-\rho(\cos\eta,\sin\eta)|^{\alpha}}r\rho drd\rho d\theta d\eta.
\end{align} 
Computing $\frac{d}{ds}H(R+sh)\bigg|_{s=0}$, one can immediately find that 
\[
\frac{d}{ds}H(R+sh)\bigg|_{s=0} = 2\int_{\mathbb{T}}\Psi(z(\theta))R(\theta)h(\theta)d\theta.
\] In other words, we have
\begin{align}\label{L2gradient}
\nabla_{L^2}H(R)(\theta)  = 2R(\theta)\Psi(z(\theta)),
\end{align}
where $\nabla_{L^2}$ denotes the gradient in $L^2(\mathbb{T})$.
Now, we define a Hamiltonian $\mathcal{H}: L^2_0\mapsto \mathbb{R}\cup \left\{\infty \right\} $ as\index{$\mathcal{H}$}
\begin{align}\label{Hamiltonian}
\mathcal{H}(f) :=H(\sqrt{1+f}).
\end{align}
Using this and the chain rule, one can easily see from \eqref{L2gradient} and \eqref{Hamiltonian} that
\begin{align}\label{l20gradient}
\nabla_{L^2_0}\mathcal{H}(f)(\theta) = \frac{1}{2\sqrt{1+f(\theta)}}\nabla_{L^2}H(\sqrt{1+f})(\theta) = \Psi(z(\theta)).
\end{align}
Also, for $f(t,\cdot)\in L^2_0 $ such that $\sqrt{1+f(t,\theta)} = R(t,\theta)$, it follows immediately that $R(t,\theta)\partial_t R(t,\theta) = \frac{1}{2}\partial_t f(t,\theta)$. Therefore, we obtain from \eqref{l20gradient} and \eqref{gSQG} that
\begin{align}\label{gsqg_2}
\partial_tf(t,\theta) =2 \partial_\theta \left(\Psi(z(t,\theta))\right) = 2\partial_{\theta}\nabla_{L^2_0}\mathcal{H}(f).
\end{align}
By reparametrizing the time as $g(t,\theta):=f\left(\frac{1}{2}t,\theta\right)$, we obtain from  \eqref{gsqg_2} that
\begin{align*}
\partial_{t}g = 2\partial_\theta \nabla_{L^2_0}\mathcal{H}(g) = X_{\mathcal{H}}(g).
\end{align*}
We summarize the above considerations in the following proposition:
\begin{proposition}\label{Hamiltonian_prop}
Let $\mathcal{H}$ be the Hamiltonian\index{Hamiltonian} given as in \eqref{Hamiltonian} and \eqref{Hamiltonian_2} and $\sigma$ be the symplectic $2$-form given in \eqref{symplectic}. Let $f(t,\cdot) \in L^2_0$\index{$L^2_0$} be a solution to the Hamiltonian equation
\begin{align}\label{gSQG_hamiltonian}
\partial_t f(t,\theta) = X_{\mathcal{H}}(f(t,\cdot)).
\end{align}
Then $\omega(t,x):=1_{D_t}(x)$, where $
D_t:=\left\{ r(\cos \theta, \sin \theta)\in \mathbb{R}^2 : 0 \le r < \sqrt{1+f(2t,\theta)}\right\}$ is a weak solution to \eqref{gSQG_1}.
\end{proposition}
 In the rest of the paper, we denote the $L^2_0$-gradient by $\nabla$, instead of $\nabla_{L^2_0}$.\section{Expansion of $\mathcal{H}$}
The main result of this section can be summarized in the following proposition\index{homogeneous expansion}:
\begin{proposition}\label{expansion_1}
$\mathcal{H}(f)$, defined in \eqref{Hamiltonian} and \eqref{Hamiltonian_2}, admits the following homogeneous expansion:
\begin{align*}
\mathcal{H}(f) = \mathcal{H}_2(f)+\mathcal{H}_3(f)+\mathcal{H}_4(f)+\mathcal{H}_{\ge 5}(f),
\end{align*}
where 
\begin{equation}\label{expansion_2}
\begin{aligned}
\mathcal{H}_2(f)&:= -\frac{1}{4}\int_{\mathbb{T}} f\Lambda^{\alpha-1}f d\theta + \frac{1}{8}T_\alpha\int_{\mathbb{T}} f^2 d\theta,\\
 \mathcal{H}_3(f)&:= \frac{\alpha}{16}\int_{\mathbb{T}} f^2 \Lambda^{\alpha-1}fd\theta   -\frac{\alpha}{48}T_\alpha\int_{\mathbb{T}} f^3d\theta,\\
\mathcal{H}_4(f) &:=-\frac{\alpha(\alpha+4)}{192}\int_{\mathbb{T}} f^3\Lambda^{\alpha-1}fd\theta -\frac{\alpha^2}{256}\int_{\mathbb{T}} f^2 \Lambda^{\alpha-1}(f^2)d\theta \\
& + \frac{\alpha}{384}\int_{\mathbb{T}} M_{4}(f)(\theta)d\theta +  \frac{\alpha(2+\alpha)T_\alpha}{384}\int_{\mathbb{T}} f(\theta)^4 d\theta,
\end{aligned}
\end{equation}
where 
\begin{align}\nonumber
T_\alpha & = m_{2,\alpha}(0) \text{ (see \eqref{multiplier_2})},\\
M_4(f)(\theta)& :=\int_{\mathbb{T}}(2-2\cos(\theta-\eta))^{-1-\frac{\alpha}{2}}(f(\theta)-f(\eta))^4d\eta,  \label{def_talpha}
\end{align}
and $\Lambda^{\alpha-1}$ is as in \eqref{fractional_1}. $\mathcal{H}_{\ge 5}$ collects all the terms homogeneous of degree at least $5$ in $f$ and its gradient is of the form:
\begin{equation}\label{gradh5}
\begin{aligned}
\nabla \mathcal{H}_{\ge 5}(f) &= \int (2-2\cos(\theta-\eta))^{-\frac{\alpha}{2}}F_{1,\ge3}(f(\theta),f(\eta),J(\theta,\eta))(f(\theta)-f(\eta))d\eta \\
& \ +\int (2-2\cos(\theta-\eta))^{1-\frac{\alpha}{2}}F_{2,\ge4}(f(\theta),f(\eta))d\eta,
\end{aligned}
\end{equation}
where $J(\theta,\eta):=\frac{(f(\theta)-f(\eta))^2}{2-2\cos(\theta-\eta)}$ for some  functions $F_{1,\ge 3}:\mathbb{R}^3\mapsto \mathbb{R}$ and $F_{2,\ge4}:\mathbb{R}^2\mapsto \mathbb{R}$ that are analytic functions\index{analytic functions} near the origin and homogeneous of degree at least $3$ and $4$ respectively.
\end{proposition}

\begin{proof}
We split the proof into three steps. In the first step, we rewrite the Hamiltonian $\mathcal{H}$ in terms of $f$. In the second step, we specify the linear, quadratic and cubic terms of $\nabla\mathcal{H}(f)$. Afterwards, we integrate them to recover the quadratic, cubic and quartic terms of $\mathcal{H}$. In the last step, we show that the terms in $\nabla \mathcal{H}$ that are homogeneous of degree at least $4$ in $f$ are of the form in \eqref{gradh5}. 

\textbf{Expression for $\nabla\mathcal{H}$ in terms of $f$.}
 We first rewrite the Hamiltonian $\mathcal{H}$ more explicitly. Given $f\in L^2_0$, we set 
 \begin{align}\label{rtof}
 R(\theta):=\sqrt{1+f(\theta)}, \quad \text{ and }z(\theta):=R(\theta)(\cos \theta, \sin \theta),
\end{align}
so that $z(\theta)$ represents a parametrization of $\partial D$, determined by $f$. Using the identity:
\begin{align*}
\int_{D}\frac{1}{|x-y|^{\alpha}}dy & = -\frac{1}{2-\alpha}\int_{\partial_D}\frac{(x-y)}{|x-y|^{\alpha}}\cdot \vec{n}(y)d\sigma \\
& = \frac{1}{2-\alpha}\int_{\mathbb{T}}\frac{(x-z(\eta))\cdot \partial_{\eta}z(\eta)^{\perp}}{|x-z(\eta)|^{\alpha}}d\eta,
\end{align*}
which simply follows from  integration by parts, we find that the corresponding stream function $\Psi(z(\theta))$ in \eqref{stream_function} can be written as
\begin{align} \nonumber
\Psi(z(\theta)) = \frac{1}{2-\alpha}\int_{\mathbb{T}} (2-2\cos(\theta-\eta))^{-\frac{\alpha}{2}}\left(R(\theta)R(\eta) + \frac{(R(\theta)-R(\eta))^2}{2-2\cos(\theta-\eta)}\right)^{-\frac{\alpha}{2}} \\
\times (A_1+A_2+A_3)d\eta, \label{stream_2}
\end{align}
where 
\begin{align*}
A_1 &:=R(\theta)R'(\eta)\sin(\theta-\eta), \quad A_2 := R(\eta)R(\theta)(1-\cos(\theta-\eta)), \\
\quad A_3& :=-R(\eta)(R(\theta)-R(\eta)).
\end{align*}
Recalling \eqref{l20gradient} and replacing $R(\theta)$ in \eqref{stream_2} by $\sqrt{1+f(\theta)}$, we obtain that 
\begin{equation}\label{grad_hamiltonian1}
\begin{aligned}
& \nabla\mathcal{H}(f) = \Psi(z(\theta))  \\
&=\frac{1}{2-\alpha} \int(2-2\cos(\theta-\eta))^{-\frac{\alpha}{2}}G_1(f(\theta),f(\eta),{J})f'(\eta)\sin(\theta-\eta)d\eta\\
& \ +\frac{1}{2(2-\alpha)} \int(2-2\cos(\theta-\eta))^{1-\frac{\alpha}{2}}G_2(f(\theta),f(\eta),{J})d\eta \\
& \ - \frac{1}{2-\alpha}\int(2-2\cos(\theta-\eta))^{-\frac{\alpha}{2}}G_3(f(\theta),f(\eta),{J})(f(\theta)-f(\eta))d\eta,
\end{aligned}
\end{equation}
where
\begin{equation}\label{grad_hamiltonian2}
\begin{aligned}
{J}=J(\theta,\eta)&:=\frac{(f(\theta)-f(\eta))^2}{2-2\cos(\theta-\eta)},\\
G_1(x,y,z) & = \left(\sqrt{1+x}\sqrt{1+y}+z\frac{\left(\sqrt{1+x}-\sqrt{1+y}\right)^2}{(x-y)^2}\right)^{-\frac{\alpha}{2}}\frac{\sqrt{1+x}}{2\sqrt{1+y}}\\
G_2(x,y,z) & = \left(\sqrt{1+x}\sqrt{1+y}+z\frac{\left(\sqrt{1+x}-\sqrt{1+y}\right)^2}{(x-y)^2}\right)^{-\frac{\alpha}{2}}\sqrt{1+x}\sqrt{1+y}\\
G_3(x,y,z) & = \left( \sqrt{1+x}\sqrt{1+y}+z\frac{\left(\sqrt{1+x}-\sqrt{1+y}\right)^2}{(x-y)^2} \right)^{-\frac{\alpha}{2}}
\\&\times\left( \frac{\sqrt{1+y}(\sqrt{1+x}-\sqrt{1+y})}{x-y}\right).
\end{aligned}
\end{equation}

\textbf{Explicit formulae for $\nabla\mathcal{H}_2$, $\nabla\mathcal{H}_3$ and $\nabla\mathcal{H}_4$.} We note that the functions $G_1,G_2,G_3$ are analytic near the origin in $\mathbb{R}^3$ (more precisely, they can be extended so that they are analytic near the origin). In order to find $\mathcal{H}_2,\mathcal{H}_3,\mathcal{H}_4$, which will correspond to the quadratic, cubic and quartic parts of the Hamiltonian $\mathcal{H}$,  we will compute the Taylor series coefficients of $G_i$ near the origin. For $\mathcal{H}_2$, we find that
\[
G_1(0,0,0) = \frac{1}{2},\quad \partial_xG_2(0,0,0)=\partial_yG_2(0,0,0)=\frac{2-\alpha}{4},\quad G_3(0,0,0) = \frac{1}{2}.
\]
This gives the linear term of $\nabla\mathcal{H}$, that is, from \eqref{grad_hamiltonian1}, it follows that
\[
\nabla\mathcal{H}_2(f) = -\frac{1}{2}\Lambda^{\alpha-1}f(\theta) + \frac{1}{4}T_\alpha f(\theta).
\]
Integrating this with respect to $f$, we find that
\begin{align}\label{H2}
\mathcal{H}_2 = -\frac{1}{4}\int_{\mathbb{T}}f\Lambda^{\alpha-1}fd\theta + \frac{T_\alpha}{8}f^2d\theta.
\end{align}
Similarly, we compute $\nabla\mathcal{H}_3,\nabla\mathcal{H}_4$ and find that
\begin{align*}
\nabla\mathcal{H}_3(f) & = \frac{\alpha}{16}\int (2-2\cos(\theta-\eta))^{-\frac{\alpha}{2}}(3f(\theta)^2 - 2f(\theta)f(\eta)-f(\eta)^2)d\eta  \\
& - \frac{\alpha}{16}\int (2-2\cos(\theta-\eta))^{1-\frac{\alpha}{2}}d\eta f(\theta)^2.
\end{align*}
and
\begin{align*}
&\nabla \mathcal{H}_4(f)  = -\frac{\alpha}{192}\int (2-2\cos(\theta-\eta))^{-\frac{\alpha}{2}}\\
&\times((7\alpha+16)f(\theta)^3 - 3(\alpha+4)f(\theta)^2f(\eta)-3\alpha f(\theta)f(\eta)^2 - (\alpha+4)f(\eta)^3)d\eta \nonumber\\
& \ + \frac{\alpha}{48}\int (2-2\cos(\theta-\eta))^{-1-\frac{\alpha}{2}} (f(\theta)-f(\eta))^3 d\eta + \frac{\alpha}{96}(2+\alpha)T_\alpha f(\theta)^3.
\end{align*}
Integrating them with respect to $f$, one can find that
\begin{align*}
\mathcal{H}_3(f) &=  \frac{\alpha}{16}\int f^2 \Lambda^{\alpha-1}fd\theta   -\frac{\alpha}{48}T_\alpha\int f^3d\theta \\
\mathcal{H}_4(f) &= -\frac{\alpha(\alpha+4)}{192}\int f^3\Lambda^{\alpha-1}fd\theta -\frac{\alpha^2}{256}\int f^2 \Lambda^{\alpha-1}(f^2)d\theta \\
& + \frac{\alpha}{384}\int M_{4}(f)d\theta + \frac{\alpha}{384}(2+\alpha)T_\alpha \int f(\theta)^4 d\theta,
\end{align*}
where $T_\alpha$ and  $M_4$ are as in \eqref{def_talpha}. With \eqref{H2}, we obtain \eqref{expansion_2}.

\textbf{Expression for $\mathcal{H}_5$.}
Let $\nabla\mathcal{H}_5:=\nabla\mathcal{H}-\sum_{i=2}^4\nabla\mathcal{H}_i$. Since all linear, quadratic and cubic contributions of $f$ on $\nabla\mathcal{H}$ are contained in $\sum_{i=2}^4\nabla\mathcal{H}_i$, we see from \eqref{grad_hamiltonian1} that
\begin{equation}\label{H_52}
\begin{aligned}
\nabla \mathcal{H}_5(f)  & =\frac{1}{2-\alpha} \int(2-2\cos(\theta-\eta))^{-\frac{\alpha}{2}}G_{1,\ge 3}(f(\theta),f(\eta),{J})f'(\eta)\sin(\theta-\eta)d\eta\\
& \ +\frac{1}{2(2-\alpha)} \int(2-2\cos(\theta-\eta))^{1-\frac{\alpha}{2}}G_{2,\ge 4}(f(\theta),f(\eta),{J})d\eta \\
& \ - \frac{1}{2-\alpha}\int(2-2\cos(\theta-\eta))^{-\frac{\alpha}{2}}G_{3,\ge 3}(f(\theta),f(\eta),{J})(f(\theta)-f(\eta))d\eta,
\end{aligned}
\end{equation}
where $G_{1,\ge3}$ and $G_{2,\ge3}$ collect all the terms from $G_1$ and $G_2$ that are homogeneous of degree at least $3$ and $G_{2\ge 4}$ collects the terms from $G_2$ that are homogeneous of degree at least $4$. We claim that there exist $F_{1,\ge3}(x,y,z),F_{2,\ge4}(x,y,z)$, which are analytic near $(0,0,0)$, homogeneous of degree at least $3$ and $4$ respectively, and
\begin{align}\label{H_53}
\nabla \mathcal{H}_5(f) &= \int (2-2\cos(\theta-\eta))^{-\frac{\alpha}{2}}F_{1,\ge3}(f(\theta),f(\eta),{J})(f(\theta)-f(\eta))d\eta\nonumber \\
& \ +\int (2-2\cos(\theta-\eta))^{1-\frac{\alpha}{2}}F_{2,\ge4}(f(\theta),f(\eta),{J})d\eta.
\end{align}
Clearly, the second and the third integrals in \eqref{H_52} are already of the form in \eqref{H_53}, thus once we prove that the first integral in \eqref{H_52} is of the form \eqref{H_53}, we will finish the proof of the proposition.

 For the first integral in \eqref{H_52}, let us write the integral, using the Taylor expansion of $G_{1,\ge3}$, as 
\begin{align}
& I_1:=\int(2-2\cos(\theta-\eta))^{-\frac{\alpha}{2}}G_{1,\ge 3}(f(\theta),f(\eta),{J})f'(\eta)\sin(\theta-\eta)d\eta\nonumber\\
&  = \sum_{n=4}^{\infty}\sum_{k+l+2m=n}C_{k,l,m}f(\theta)^{k}\underbrace{\int (2-2\cos(\theta-\eta))^{-\frac{\alpha}{2}}f(\eta)^lf'(\eta){J}^{m}\sin(\theta-\eta)d\eta}_{=:a_{l,m}}.\label{I_def}
\end{align}
 Using integration by parts, it follows immediately that
\begin{align}
&a_{l,m} = f(\theta)d_{l,m}a_{l-1,m} + c_{l,m},\label{a_sequence}\\
&a_{0,m} = \frac{2-\alpha-2m}{4(2m+1)}\int (2-2\cos(\theta-\eta))^{1-\frac{\alpha}{2}}(f(\theta)-f(\eta)){J}^m d\eta,\nonumber\\
 &\quad \quad + \frac{2m+\alpha-1}{2m+1}\int (2-2\cos(\theta-\eta))^{-\frac{\alpha}{2}}(f(\theta)-f(\eta)){J}^m d\eta, \label{a_sequence_0}
\end{align}
where
\begin{align*}
&d_{l,m}:=\frac{l}{2m+1+l}\\
&c_{l,m}:=\frac{2-\alpha-2m}{4(2m+1+l)}\left(\int (2-2\cos(\theta-\eta))^{1-\frac{\alpha}{2}}f(\eta)^l (f(\theta)-f(\eta)){J}^md\eta\right),\\& \quad \quad +\frac{2m+\alpha-1}{2m+1+l}(\int (2-2\cos(\theta-\eta))^{-\frac{\alpha}{2}}f(\eta)^l(f(\theta)-f(\eta)) {J}^md\eta).\\
\end{align*}
Hence, it follows from \eqref{a_sequence} that
\[
a_{l,m} = f(\theta)^{l}\left(\Pi_{j=1}^{l}d_{l,m}\right)a_{0,m} + \sum_{i=0}^{l-1}\left(\Pi_{j=l-i}^{l}d_{j,m}\right)f(\theta)^{i+1}c_{l-i-1,m} + c_{l,m}.
\]
Plugging this into \eqref{I_def}, it is clear that $I_1$ in \eqref{I_def} is of the form in \eqref{H_53}. The convergence of the series is guaranteed from the fact that $G_{1}$ is analytic near $(0,0,0)$. Therefore we have \eqref{H_53}.

 Lastly, since any term in the second integral in \eqref{H_53} that is homogeneous of degree at least order $1$ in ${J}$ is of the form as the first integral in \eqref{H_53}, we obtain \eqref{gradh5}.
\end{proof}

\begin{corollary}\label{structure_hanams}
The Hamiltonian $\mathcal{H}$ is $\mu_0$-regular for some $\mu_0=\mu_0(\alpha)>0$.
\end{corollary}
\begin{proof}
The result follows straightforwardly from the expansion of $\mathcal{H}$ in Proposition~\ref{expansion_1} and successive applications of \eqref{bncal}. \end{proof}

\section{Conservation of momentum, time-reversibility and $\mathtt{M}$-fold symmetry}\label{momentum}
Before we close this chapter, we point out three useful properties of the gSQG equations: Conservation of momentum\index{Conservation of momentum}, time-reversibility\index{time-reversibility} and invariance under rotations\index{invariance under rotations}.
\subsection{Conservation of momentum}
We recall that the gSQG dynamics preserves the angular momentum\index{angular momentum}, $\int_{\mathbb{R}^2}\omega(t,x)|x|^2dx$.  In the patch setting, this corresponds to (using the notations in \eqref{rtof}), 
\begin{align*}
M(f):=\int_{\mathbb{T}}|f(x)|^2dx \text{ is a preserved quantity.}
\end{align*}
In other words, we have $\left\{ \mathcal{H},M\right\} = 0$, where the Poisson bracket\index{Poisson bracket} $\left\{ \cdot, \cdot \right\}$ is defined as 
\begin{align}\label{poisson_bracket}
\left\{ H,F \right\}(f) :=  \int_{\mathbb{T}} \nabla H(f) \partial_x \nabla F(f) dx = -\frac{1}{2\pi}\sum_{j\ne  0} \ii j\left( \partial_{f_{-j}}H\right)(f)\left(\partial_{f_j}F\right)(f) .
\end{align} 
Note that for a homogeneous Hamiltonian $H_{n}$ of degree $n$ such that
$
H_n(f) = \sum_{j_1,\ldots,j_n}H_{n,j_1,\ldots,j_n}f_{j_1}\cdots f_{j_n},
$
$H_n$ preserves the momentum, that is, $\left\{ H_n,M\right\} = 0$ if 
\begin{align}\label{conserve_momentum}
H_{n,j_1,\ldots, j_n} = 0, \text{ for } j_1+\cdots + j_n \ne 0.
\end{align} 
It is straightforward to see from Proposition~\ref{expansion_1} that each homogeneous term $\mathcal{H}_n$ satisfies \eqref{conserve_momentum}, thus the Hamiltonian $\mathcal{H}$ preserves the momentum.

\subsection{Time-reversibility}\label{reversibility22}

From \eqref{Hamiltonian_2} and \eqref{Hamiltonian}, one can easily check that the gSQG Hamiltonian $\mathcal{H}$ in \eqref{Hamiltonian} is time-reversible\index{time-reversibility} with respect to the involution $\rho$ in \eqref{definiteionofinv}. 
 Moreover, one can easily see from Proposition~\ref{expansion_1} that $\mathcal{H}_{i}$ for  $i=1,\ldots,4$ are time-reversible as well.
Note that for a homogeneous Hamiltonian $H_n$ of degree $n$ such that $
H_{n}(f) = \sum_{j_1,\ldots,j_n}H_{n,j_1,\ldots,j_n}f_{j_1}\cdots f_{j_n}$, the time-reversibility of $H_n$ with respect to $\rho$ implies that
\[
H_{n,-j_1,\ldots,-j_n} = H_{n,j_1,\ldots,j_n}.
\]
On the other hand, if $H_{n}\circ \rho = - H_n$, that is, $H_n$ is reversibility preserving (see Definition~\ref{shdj2sdsdsd}), then  it holds that
\begin{align}\label{2jsdsdsdsd2dp0s}
H_{n,-j_1,\ldots,-j_n} =  -H_{n,j_1,\ldots,j_n}.
\end{align}
In such case, Lemma~\ref{time_reversibility_map} implies that the time-1 flow map of the Hamiltonian PDE associated to $H_{n}$ is a reversibility preserving map. In the next chapter, our symplectic transformations to obtain the weak Birkhoff normal form will be obtained as flow maps generated by  Hamiltonians of the form in \eqref{2jsdsdsdsd2dp0s}. Thus, they are not only symplectic, but also reversibility preserving.

\begin{remark}\label{real_jgasddsd}
It is clear from \eqref{Hamiltonian_2} that $H$ is a real-valued Hamiltonian. In general, a homogeneous Hamiltonian $H_{n}(f)\mapsto \sum_{j_1,\ldots,j_n}H_{n,j_1,\ldots,j_n}f_{j_1},\ldots,f_{j_n}$ is real-valued if and only if
\begin{align}\label{real_bonsdsd}
\overline{H_{n,j_1,\ldots j_n}} = H_{n,-j_1,\ldots,j_n}.
\end{align}
Since our phase space is a set of real-valued functions, it follows that its Hessian operator is a real linear operator (see Definition~\ref{realtorealop}).
\end{remark}

\subsection{$\mathtt{M}$-fold symmetric patches}
Another key property of the gSQG equations\index{$\mathtt{M}$-fold symmetric patches} is that the solutions of \eqref{gSQG_1} enjoy certain rotational symmetry. More precisely, let us denote by $R_\theta$ the rotation matrix:
\[
R_{\theta}:=\colvec{\cos \theta &  -\sin \theta \\ \sin\theta & \cos\theta},\text{ for $\theta\in \mathbb{T}$.}
\]
One can easily check that if $\omega(t,x)$ is a solution to \eqref{gSQG_1} with an initial datum $\omega_0(x)$, then $\omega_{\theta}(t,x):=\omega(t,R_\theta x)$ solves \eqref{gSQG_1} with the initial datum $\omega_0(t,R_\theta x)$. Especially, if the initial datum $\omega_0$ is invariant under a $\frac{2\pi}{\mathtt{M}}$-rotation, that is, $\omega_0(x) = \omega_0(R_{2\pi/\mathtt{M}}x)$, then the solution $\omega(t,x)$ is also invariant under  a $\frac{2\pi}{\mathtt{M}}$-rotation. In the phase space $L^2_0$, such invariance can be rephrased as $\frac{2\pi}{\mathtt{M}}$-translation invariance. Indeed, we have the following:

\begin{lemma}
The Hamiltonian $\mathcal{H}$ in \eqref{Hamiltonian} is $\frac{2\pi}{\mathtt{M}}$-translation invariant.
\end{lemma}
\begin{proof}
In view of \eqref{Hamiltonian}, it suffices to show that $H$ in \eqref{Hamiltonian_2} satisfies $H(\rho_{\mathtt{M}}(R)) = H(R)$, that is,
\begin{align}\label{rlsdsdclaim}
\int_{R_{\frac{2\pi}{\mathtt{M}}}D}\int_{R_{\frac{2\pi}{\mathtt{M}}}D} \frac{1}{|x-y|^{\alpha}}dxdy = \int_D{\int_{D}} \frac{1}{|x-y|^{\alpha}}dxdy,
\end{align}
which follows immediately from   a change of variables.
\end{proof}

 \begin{remark}\label{krsds2221sd}One can easily see from Proposition~\ref{expansion_1} that $\mathcal{H}_{i}$ for  $i=1,\ldots,4$ are  $\frac{2\pi}{\mathtt{M}}$-translation invariant  as well. Note that for a homogeneous Hamiltonian $H_n$ of degree $n$ such that $
H_{n}(f) = \sum_{j_1,\ldots,j_n}H_{n,j_1,\ldots,j_n}f_{j_1}\cdots f_{j_n}$, the condition \eqref{conserve_momentum} implies  $\frac{2\pi}{\mathtt{M}}$-translation invariance of $H_n$ as well. Indeed, 
\begin{align*}
H_n(\rho_{\mathtt{M}}f) & = \sum_{j_1,\ldots,j_n}H_{j_1,\ldots,j_n}\widehat{\rho_{\mathtt{M}}(f)}_{j_1}\cdots \widehat{\rho_{\mathtt{M}}(f)}_{j_n} \\ & =  \sum_{j_1,\ldots,j_n} H_{j_1,\ldots,j_n} e^{\frac{2\pi}{\mathtt{M}}\ii (j_1 + \cdots + j_n)}\widehat{f}_{j_1}\cdots \widehat{f}_{j_n}
\end{align*}
which is equal to $H_n(f)$ if \eqref{conserve_momentum} holds. 
\end{remark}

\color{black}

\chapter{Weak Birkhoff normal form}\label{skpssisodssuw}
In this chapter, we aim to normalize the Hamiltonian $\mathcal{H}$, up to the quintic term. 
\section{Tangential sites and normal sites}\label{posdsdtangeinal}
Given $f\in L^2_0$, we recall the Fourier series expression (compared to \eqref{sjd2sjjjsd}, we replace $\widehat{f}_j$ by $f_j$ for simplicitiy).
\begin{align}\label{fourier}
f(x) = \sum_{j\ne 0}f_je^{\ii j x},  \text{ where } f_j:=\frac{1}{2\pi}\int_{\mathbb{T}}f(x)e^{-\ii j x}dx.  
\end{align} 
Plugging the series expansion in \eqref{fourier} into \eqref{expansion_2}, one can rewrite $\mathcal{H}_n$, for $n=1,...4$,  as
\begin{equation}\label{fourier_hamiltonian}
\begin{aligned}
&\mathcal{H}_2(f) =\sum_{j\ne 0}\kappa_{j}|f_j|^2, \quad \mathcal{H}_3(f) = \sum_{j_1+j_2+j_3=0}H_{3,j_1,j_2,j_3}f_{j_1}f_{j_2}f_{j_3},\\& \mathcal{H}_4(f) = \sum_{j_1+j_2+j_3+j_4=0}H_{4,j_1,j_2,j_3,j_4}f_{j_1}f_{j_2}f_{j_3}f_{j_4}.
\end{aligned}
\end{equation}
where 
\begin{equation}\label{moreexplicits2}
\begin{aligned}
\kappa_j&:={\pi}\left(-\frac{1}{2}m^\circ_{1,\alpha}(j)  + \frac{T_\alpha}{4} \right),\\
  H_{3,j_1,j_2,j_3}&:=\frac{\pi\alpha}{24}\left(m^\circ_{1,\alpha}(j_1)+m^\circ_{1,\alpha}(j_2)+m^\circ_{1,\alpha}(j_3) - T_\alpha \right),\\
H_{4,j_1,j_2,j_3,j_4}&:= -\frac{\alpha(\alpha+4)\pi}{384}(m_{1,\alpha}(j_1)+m_{1,\alpha}(j_2)+m_{1,\alpha}(j_3)+m_{1,\alpha}(j_4))\\
& \ -\frac{\alpha^2\pi}{768}(m^\circ_{1,\alpha}(j_1+j_2)+m^\circ_{1,\alpha}(j_1+j_3)\\
&+m^\circ_{1,\alpha}(j_1+j_4)+m^\circ_{1,\alpha}(j_2+j_3)+m^\circ_{1,\alpha}(j_2+j_4)+m^\circ_{1,\alpha}(j_3+j_4))\\
& \ + \frac{\alpha\pi}{192}m^\circ_{4,\alpha}(j_1,j_2,j_3,j_4) + \frac{\alpha(\alpha+2)T_\alpha\pi}{192}.
\end{aligned}
\end{equation}
and\index{$T_\alpha$}
\begin{equation}\label{mutsd2xcxcxc}
\begin{aligned}
m^\circ_{1,\alpha}(j)& :=\int_{\mathbb{T}}(2-2\cos(\theta))^{-\frac{\alpha}2}(1-e^{-\ii j \theta})d\theta,\\
T_{\alpha}&:= \int_{\mathbb{T}}(2-2\cos(\theta))^{1-\frac{\alpha}2}d\theta,\\
m^\circ_{4,\alpha}(j_1,j_2,j_3,j_4)& =\int_{\mathbb{T}}(2-2\cos(\theta))^{-1-\frac{\alpha}2}(1-e^{-\ii j_1 \theta})(1-e^{-\ii j_2 \theta}) \\
&\times (1-e^{-\ii j_3 \theta})(1-e^{-\ii j_4 \theta})d\theta.
\end{aligned}
\end{equation}
Note that $H_{3,j_1,j_2,j_3}$ and $H_{4,j_1,j_2,j_3,j_4}$ are expressed in such a way that they are invariant under any permutation of the sub-indices.

Now we pick the tangential sites of the Fourier modes, for which we normalize the Hamiltonian. For  $\nu \in \mathbb{N}$, we pick natural numbers $0 <\mathtt{j_1} < \cdots < \mathtt{j}_\nu$ and set 
\begin{align}\label{s00sd1sd}
S_0:=\left\{ \pm \mathtt{j} : \mathtt{j}\in S_0^+\right\},\quad S_0^+:=\left\{ 0 < \mathtt{j}_1 < \cdots < \mathtt{j}_\nu,\text{ for $i=1,\ldots,\nu$}\right\}.
\end{align}
The tangential sites $S$\index{tangential sites} and the normal sites\index{normal sites} $S^\perp$ are chosen to be  
\begin{align}\label{tan_site}
S^+:=\mathtt{M}S_0^+:=\left\{ \mathtt{M}j: j\in S_0^{+}\right\},\quad S:=\left\{ \pm j : |j|\in S^+\right\}, \quad S^\perp := \mathbb{Z}\backslash \left( S \cup \left\{ 0\right\}\right).
\end{align}
In order to perform derive the weak Birkhoff normal form\index{Weak Birkhoff normal form}, we impose the following conditions on the choice of $S$:
\begin{enumerate}[label=\textbf{S.\arabic*}]
\item \label{tangent_2} If $j_1,j_2\in S$ and $j_1\ne -j_2$, then $j_1+j_2\in S^\perp$.
\item  If more than three of  $j_1,j_2,j_3,j_4,j_5$ are in $S$ and none of them is $0$, then $j_1\kappa_{j_1}+j_2\kappa_{j_2}+j_3\kappa_{j_3}+j_4\kappa_{j_4}+j_5\kappa_{j_5} \ne 0. 
$\label{tangent_1}
\end{enumerate}
Let $S,S^\perp$ be fixed as above. Before performing the weak Birkhoff normal form, we specify some useful notations:
\subsection{Notations} We denote
\begin{align}\nonumber
H_S & := \left\{ f\in L^2_0: f_j = 0, \text{ for }j\in S^\perp\right\}, \\H_{S^\perp}& :=  (H_{S})^{\perp} \ (\text{orthogonal complement of $H_S$ in $L^2_0$}). \label{tangential_function}
\end{align}
We denote by $\Pi_{S}$ and $\Pi_{S^\perp}$ the projections onto the 
subspaces $H_{S}$ and $H_{S^\perp}$ respectively. Given $f\in L^2_0$, we can write it as
\begin{align}\label{decom_tandn}
f := v+ z = \sum_{j\in S}v_je^{\ii  j x} + \sum_{j\in S^\perp}z_je^{\ii  j x},  \text{ where }  v:=\Pi_S f\in H_{S}, z := \Pi_{S^\perp}f \in H_{S^\perp}.
\end{align}

Let $E$ be a finite dimensional subspace\index{finite dimensional subspace} in $L^2_0$ such that
\begin{align}\label{finite_dim}
E:= \text{span}\left\{ e^{\ii jx}: |j| \le C \right\},
\end{align}
for some constant $C>0$, and denote by $\Pi_{E}$ the orthogonal projection\index{orthogonal projection} to $E$.

In the weak Birkhoff normal form procedure, we do not need an explicit expression for some multilinear  maps\index{multilinear maps}, but we  only need how such multilinear maps depend on the tangential/normal component.  For simplicity, we denote a $k$-linear map\index{$k$-linear map} in $(v,z)$ by
\begin{align}\label{multilinear52}R(v^{k-q}z^q) = M[\underbrace{v,\ldots,v}_{k-q \text{ times}},\underbrace{z,\ldots,z}_{q \text{ times}}], \quad M=k\text{-linear}.\end{align}
For a homogeneous Hamiltonian $H_n$ of degree $n$, we write it as
\begin{align}\label{degree_on_normal}
H_{n} = \sum_{i=0}^nH_{n,i}, \text{ where  each of $H_{n,i}$ is of the form $R(v^{n-i}z^i)$.}
\end{align}

\section{Composition with a time-1 flow}\index{time-1 flow}
In order to find the weak Birkhoff normal form of the Hamiltonian $\mathcal{H}$, we construct a sequence of symplectic transformations to remove the trivial resonances arising from the tangential sites $S$. We first recall the following classical lemma to find an expansion of a composition of a Hamiltonian with a symplectic transformation:
\begin{lemma}\label{symplectic_1} Let $H,F:L^2_0 :\mapsto \R\cup\left\{ \infty\right\}$ be  Hamiltonians in $L^2_0$ such that $F$ is supported on a finite number of Fourier modes, more precisely, $F(\Pi_{E^c}u) = 0$ for a finite dimensional space of the form in \eqref{finite_dim}. Also, let $\Phi_t: L^2_0\mapsto L^2_0$ be the associated flow map at time $t$, that is, $\frac{d}{dt}\Phi_t(u) = X_{F}(\Phi_t(u))$ for $t\ge 0$. Then for each $t\ge 0$, $\Phi_t$ is a symplectic transformation and admits the following Taylor expansion:
\begin{equation}\label{expansion2s3dsss}
\begin{aligned}H\circ\Phi_1(u) &=H(u) + \left\{ H,F\right\} + \left\{\left\{ H,F\right\},F\right\} +\frac{1}{2!}\left\{ \left\{\left\{ H,F\right\},F\right\},F\right\} + \ldots \\
& = \frac{1}{n!}\sum_{i=0}^{n}A_n(u) +\frac{1}{n!}\int_0^1 A_{n+1}(\Phi_t(u))(1-t)^ndt,\\
&A_{0}:=H,\quad A_{n+1} :=\left\{ A_n,F\right\}, \text{ for $n\ge 0$.}
\end{aligned}
\end{equation}
where $\left\{ \cdot, \cdot \right\}$ denotes the Poisson Bracket given in \eqref{poisson_bracket}. 
\end{lemma}
\begin{proof}
Since $F$ is supported on a finite dimensional space $E$, the map $\Phi^t$ is well-defined for all $t>0$ which depends on $t$ smoothly. That $\Phi_t$ is symplectic transformation is a classical result  (see \cite{Hubbard-Ilyashenko:proof-kolmogorov} for example). To see the expansion in \eqref{expansion2s3dsss}, let $f(t):=H(\Phi_t(u))$. From \eqref{poisson_bracket} and the fact that $\Phi_t$ is the flow map of the Hamiltonian $F$, it follows that
\begin{align*}
&f^{(n+1)}(t)=\left\{ A^n,F\right\}(\Phi_t(u)).
\end{align*}
Therefore \eqref{expansion2s3dsss} follows from \eqref{elementary_taylor}.
\end{proof}

We already observed in Section~\ref{momentum} that each $\mathcal{H}_n$ in Proposition~\ref{expansion_1} preserves the momentum, that is, it  satisfies \eqref{conserve_momentum}. In the sequel, we will construct a symplectic transformation as a time-1 flow of a homogenous Hamiltonian of degree $3,4,5$ satisfying \eqref{conserve_momentum}. The next lemma shows that the composition with such transformation preserves the property \eqref{conserve_momentum}.

\begin{lemma}\label{momentum1}
Let $H_n,F_m:L^2_0\mapsto \mathbb{R}\cup \left\{ \infty\right\}$ be homogeneous Hamiltonians of degree $n$ and $m$ respectively, such that
\begin{align*}
H_{n,j_1,\ldots,j_n}= 0, \text{ if }j_1+\cdots j_n \ne 0,\quad F_{m,j_1,\ldots,j_m} =0, \text{ if }j_1+\cdots j_m \ne 0.
\end{align*}
Then,  $G:=\left\{ H_n,F_m\right\}$ is a homogeneous Hamiltonian of degree $n+m-1$ such that
\[
G_{n+m,j_1,\ldots,j_{n+m-1}} = 0, \text{ if }j_1 + \cdots + j_{n+m-1} \ne 0.
\]
\end{lemma}
\begin{proof}
The lemma follows straightforwardly from the definition of the Poisson bracket \eqref{poisson_bracket}.
\end{proof}

From Corollary~\ref{structure_hanams}, we also know that $\mathcal{H}$ is $\mu$-regular for some $\mu=\mu(\alpha)>0$. The following lemma shows that the weak Birkhoff normal form that we will obtain in the next section satisfies the same property:
\begin{lemma}\label{muregular_kep}
Let $H$ be $\mu$-regular and $F:L^2_0\mapsto \mathbb{R}\cup\left\{ \infty\right\}$ be a Hamiltonian supported on a finite number of Fourier modes and $\Phi:C^\infty(\mathbb{T})\mapsto C^\infty(\mathbb{T})$ be the time-1 flow map of its Hamiltonian PDE  as in Lemma~\ref{symplectic_1}. Then $H$ is $\mu$-regular.
\end{lemma}
\begin{proof}
Note that the Hamiltonian PDE of $F$, $f_t:=\partial_x\nabla F(f)$, is actually a  finite dimensional ODE, where with a smooth vector field. This gives us a smooth flow map $\Phi_t$. Therefore the result follows straightforwardly from the Talyor expansion of $t\mapsto \Phi_t$ at $t=0$.
\end{proof}

\section{Weak Birkhoff normal form}
The main result of this section is as follows\index{Weak Birkhoff normal form}:
\begin{proposition}\label{normal_form_prop11231}
Let $S,S^{\perp}$ be fixed as in \eqref{tan_site} satisfying \ref{tangent_2} and \ref{tangent_1}. Then, there exists a symplectic transformation\index{$\Phi^{WB}$} $\Phi^{WB}:L^2_0\mapsto L^2_0$ such that 
\begin{align}\label{normal_form_f}
H :=\mathcal{H}\circ \Phi^{WB}= H_2+H_3+H_4+H_5+H_{\ge 6},
\end{align}
where
\begin{enumerate}[label=\textbf{\arabic*}]
\item (Quadratic and cubic terms) we have
\begin{equation}\label{normal_formsf}
\begin{aligned}
H_2(f)& = \int_{\mathbb{T}}-\frac{1}{4} v\Lambda^{\alpha-1}v +\frac{T_\alpha}{8}v^2dx + \int_{\mathbb{T}}-\frac{1}{4} z\Lambda^{\alpha-1}z +\frac{T_\alpha}{8}z^2dx =:H_{2,0} + H_{2,2} \\
H_3(f) & = \frac{\alpha}{16}\int_\mathbb{T} v\left( 2z \Lambda^{\alpha-1}z + \Lambda^{\alpha-1}(z^2)  -T_\alpha z^2\right) dx  + \frac{\alpha}{16} \int_\mathbb{T} z^2 \Lambda^{\alpha-1}z -\frac{T_\alpha}{3}z^3dx \\
& =:H_{3,2} + H_{3,3}.
\end{aligned}
\end{equation}
\item \label{sdsd22sd}(Quartic and quintic terms) We have $H_{4}(f) = H_{4,0} + H_{4,2} + H_{4,3}+H_{4,4}$ where
\begin{align}
& H_{4,0}(f) = 6\sum_{j_1\in S^{+}}H^{(3)}_{4,j_1,-j_1,j_1,-j_1}|v_{j_1}|^4 + 12\sum_{\substack{j_1,j_2\in S^+,\\  j_1\ne j_2}}H^{(3)}_{4,j_1,-j_1,j_2,-j_2}|v_{j_1}|^2|v_{j_2}|^2, \label{h4h5}\\
& H^{(3)}_{4,j_1,-j_1,j_2,-j_2} = H_{4,j_1,-j_1,j_2,-j_2}  \nonumber\\
&-\begin{cases}
\frac{3}{4}\frac{(j_1+j_2)(H_{3,j_1,j_2,-(j_1+j_2)})^2}{(j_1+j_2)\kappa_{j_1+j_2}-j_1\kappa_{j_1}-j_2\kappa_{j_2}} & \text{ if $j_1=j_2$,}\\
\frac{3}{4}\left( \frac{(j_1+j_2)(H_{3,j_1,j_2,-(j_1+j_2)})^2}{(j_1+j_2)\kappa_{j_1+j_2}-j_1\kappa_{j_1}-j_2\kappa_{j_2}} +\frac{(j_1-j_2)(H_{3,j_1,-j_2,-(j_1-j_2)})^2}{(j_1-j_2)\kappa_{j_1-j_2} - j_1\kappa_{j_1}+j_2\kappa_{j_2}}\right)& \text{ if $j_1\ne j_2$,}
\end{cases}\label{h4h52s}\\
& H_{4,2}(f)  :=  \mathcal{H}_{4,2}(f) + \mathfrak{H}_{1}(f) + \mathfrak{H}_2(f)\label{sjdjsjs93jx} \\
& \mathfrak{H}_1(f) := \frac{3\alpha}{8}\int \left(v\Lambda^{\alpha-1}z + z\Lambda^{\alpha-1} v + \Lambda^{\alpha-1}(vz) - T_\alpha vz \right)\partial_x \Pi_{S}K_1(v,z)dx \label{sjdjsjs93jx32x}\\
& \mathfrak{H}_2(f) :=\frac{3\alpha}{16} \int \left( 2z\Lambda^{\alpha-1}z + \Lambda^{\alpha-1}(z^2) - T_\alpha z^2\right)\partial_x \Pi_{S^\perp}K_2(v,v)dx, \nonumber \\
& \text{ where $K_1,K_2$ are as in \eqref{K12def},}\nonumber\\
 &H_{4,3}=R(vz^3),\quad H_{4,4}=R(z^4),\nonumber\\
&H_{5}(f) =\sum_{i=2}^5H_{5,i}= \sum_{i=2}^5R(v^{5-i}z^i),\nonumber
\end{align}
and $H_{\ge 6}$ collects all the terms homogeneous of degree at least $6$.\index{$\mathfrak{H}_1$}\index{$\mathfrak{H}_2$}
\item \label{item3_normal} There exists a finite dimensional space $E$ of the form in \eqref{finite_dim} such that $\Phi^{WB}(f) = f + \Psi(f) $ where $\Psi$ satisfies $\Psi(f) = \Pi_E \Psi(\Pi_E f)$. 
\item \label{item3_normal_2} $\Phi^{WB}$ is real, reversibility preserving and $\frac{2\pi}{\mathtt{M}}$-translation invariance preserving. 
\item \label{item3_normal_3} $H$ is $\mu_1$-regular for some $\mu_1=\mu_1(\alpha)>0$.
\end{enumerate}
\end{proposition}

The proof of the proposition will be given at the end of the section, using several lemmas. 

\begin{lemma}\label{ssslem} (Normalization of the cubic Hamiltonian) There exists a symplectic transformation $\Phi^{(3)}$ such that 
\begin{align*}
H^{(3)}(f):=\mathcal{H}\circ \Phi^{(3)}(f) = H^{(3)}_2(f) + H^{(3)}_3(f) +H^{(3)}_4(f) + H^{(3)}_{\ge 5}(f),
\end{align*}
where
\begin{itemize}
\item[(1)] (Quadratic and cubic terms) We have
\begin{equation}\label{m_H3}
\begin{aligned}
H^{(3)}_2(f) &= \mathcal{H}_2(f), \\
H^{(3)}_3(f) &= \mathcal{H}_{3,\ge2}=  \frac{\alpha}{16}\int_\mathbb{T} v\left( 2z \Lambda^{\alpha-1}z + \Lambda^{\alpha-1}(z^2)  -T_\alpha z^2\right) dx  \\
& + \frac{\alpha}{16} \int_\mathbb{T} z^2 \Lambda^{\alpha-1}z -\frac{T_\alpha}{3}z^3dx.
\end{aligned}
\end{equation}
\item[(2)] (Quartic and higher degree terms) We have
\begin{equation}\label{quartic_11122}
\begin{aligned}
&H^{(3)}_4(f) = \sum_{i=0}^4\mathcal{H}^{(3)}_{4,i}(f),\\
&H^{(3)}_{4,0} = \mathcal{H}_{4,0} + \frac{3\alpha}{32}\int \left( \Lambda^{\alpha-1}(v^2) + 2v \Lambda^{\alpha-1} v - T_\alpha v^2 \right)\partial_x\Pi_{S^\perp}K_2(v,v)dx,\\
& H^{(3)}_{4,2} = \mathcal{H}_{4,2} + \frac{3\alpha}{8}\int \left(v\Lambda^{\alpha-1}z + z\Lambda^{\alpha-1} v + \Lambda^{\alpha-1}(vz) - T_\alpha vz \right)\partial_x \Pi_{S}K_1(v,z)dx \\
& \quad \quad \  + \frac{3\alpha}{16} \int \left( 2z\Lambda^{\alpha-1}z + \Lambda^{\alpha-1}(z^2) - T_\alpha z^2\right)\partial_x \Pi_{S^\perp}K_2(v,v)dx.
\end{aligned}
\end{equation}
and $H^{(3)}_{\ge 5}(f)$ collects all the terms homogeneous of degree at least $5$ and $H^{(3)}_5$ satisfies \eqref{conserve_momentum}.
\item[(3)] There exists a finite dimensional space $E$ of the form in \eqref{finite_dim} such that $\Phi^{(3)}(f) = f + \Psi^{(3)}(f) $ where $\Psi^{(3)}$ satisfies $\Psi^{(3)}(f) = \Pi_E \Psi^{(3)}(\Pi_E f)$.
\item[(4)] $\Phi^{(3)}$ is real and reversibility preserving.
\end{itemize}
\end{lemma}
\begin{proof}

\vspace{0.5\baselineskip}
\noindent\textit{Proof of (1).} Let us consider a cubic Hamiltonian:
\begin{align}\label{F3def}
F^{(3)}(f) := \sum_{j_1+j_2+j_3 = 0}F^{(3)}_{j_1,j_2,j_3}f_{j_1}f_{j_2}f_{j_3},
\end{align}
which will be determined later, but we assume that $F^{(3)}_{j_1,j_2,j_3}$ is invariant under any permutation of the sub-indices (see \eqref{def_F3} for the choice of $F^{(3)}$). We denote its time-$1$ flow map by $\Phi^{(3)}$.  We compute using Lemma~\ref{symplectic_1}, 
\begin{align}\nonumber
{H}^{(3)} & := \mathcal{H}\circ \Phi^{F^{(3)}}_{1} = \underbrace{\mathcal{H}_2}_{=:{H}^{(3)}_2} + \underbrace{\left\{\mathcal{H}_2, F^{(3)} \right\} + \mathcal{H}_3}_{=:{H}^{(3)}_3} \\
& + \underbrace{\mathcal{H}_4 + \left\{ \mathcal{H}_3, F^{(3)} \right\} + \frac{1}{2} \left\{ \left\{ \mathcal{H}_2, F^{(3)}\right\}, F^{(3)}\right\}}_{=:{H}^{(3)}_4} + {H}^{(3)}_{\ge 5}.
\label{H33expansion}
\end{align}
Using \eqref{moreexplicits2}, we have  $\left(\partial_{f_{-j}}\mathcal{H}_2\right)(f) = 2\kappa_{j}f_j$, $\left(\partial_{f_j}F^{(3)}\right)(f) = 3\sum_{j_1+j_2 = -j}F^{(3)}_{j_1,j_2,j}f_{j_1}f_{j_2}$. Thus, using \eqref{poisson_bracket}, we obtain (relabeling $j$ by $j_3$),
\begin{align*}
\left\{\mathcal{H}_2, F^{(3)} \right\}  & =-\frac{6\ii}{2\pi} \sum_{j_3\ne 0 } \left(j_3\kappa_{j_3}f_{j_3}\sum_{j_1+j_2 +j_3=0}F^{(3)}_{j_1,j_2,j_3}f_{j_1}f_{j_2}\right) \\
& = -\frac{6\ii}{2\pi}\sum_{j_1+j_2+j_3 = 0}j_3\kappa_{j_3}F^{(3)}_{j_1,j_2,j_3}f_{j_1}f_{j_2}f_{j_3}\\
& = -\frac{\ii}{\pi} \sum_{j_1 + j_2 + j_3 = 0}\left(j_1\kappa_{j_1} +j_2\kappa_{j_2} +j_3\kappa_{j_3}  \right)F^{(3)}_{j_1,j_2,j_3}f_{j_1}f_{j_2}f_{j_3},
\end{align*}
where the last equality follows from the assumption that $F^{(3)}_{j_1,j_2,j_3}$ is invariant under a permutation on the sub-indices $j_1,j_2,j_3$.
Therefore, it follows from \eqref{fourier_hamiltonian} and \eqref{H33expansion} that 
\begin{align}
{H}^{(3)}_3 &=\sum_{j_1+j_2+j_3 = 0} \left( H_{3,j_1,j_2,j_3}- \frac{\ii}{\pi} \left(j_1\kappa_{j_1} +j_2\kappa_{j_2} +j_3\kappa_{j_3}  \right)F^{(3)}_{j_1,j_2,j_3}\right)f_{j_1}f_{j_2} f_{j_3}\nonumber\\
& = \sum_{\substack{j_1+j_2+j_3 = 0 \\ (j_1,j_2,j_3 )\in \mathcal{A}_3}} \left( H_{3,j_1,j_2,j_3}- \frac{\ii}{\pi} \left(j_1\kappa_{j_1} +j_2\kappa_{j_2} +j_3\kappa_{j_3}  \right)F^{(3)}_{j_1,j_2,j_3}\right)f_{j_1}f_{j_2} f_{j_3}\label{sssperp}\\
& \ + \sum_{\substack{j_1+j_2+j_3 = 0 \\ (j_1,j_2,j_3 ) \in \mathcal{A}_3^c}} \left( H_{3,j_1,j_2,j_3}- \frac{\ii}{\pi} \left(j_1\kappa_{j_1} +j_2\kappa_{j_2} +j_3\kappa_{j_3}  \right)F^{(3)}_{j_1,j_2,j_3}\right)f_{j_1}f_{j_2} f_{j_3},\label{sssperp2}
\end{align}
where 
\begin{align}\label{defi_1sd3a}
\mathcal{A}_3 := \left\{ (j_1,j_2,j_3)\in  \mathcal \Z^3:\ \text{at least two of $\left\{j_1,j_2,j_3 \right\}$ are in $S$}\right\}.
\end{align}  We claim that
\begin{align}\label{sssclaim}
j_1\kappa_{j_1} +j_2\kappa_{j_2} +j_3\kappa_{j_3}  \ne 0, \quad \text{if $j_1+j_2+j_3 = 0$.}
\end{align}
To see this, let us suppose to the contrary that there exist $j_1,j_2,j_3$ such that
\[
j_1+j_2+j_3=0, \quad j_1\kappa_{j_1} +j_2\kappa_{j_2} +j_3\kappa_{j_3}  = 0, \quad j_1,j_2,j_3\ne 0.
\]
Without loss of generality, we assume that $j_1>0$, $j_2>0$, hence $|j_3|>|j_1|,|j_2|$. From $\kappa_j$ in \eqref{moreexplicits2} and $j_1+j_2+j_3=0$, it holds that
\[
j_1\kappa_{j_1} +j_2\kappa_{j_2} +j_3\kappa_{j_3} = 0 \implies j_1m^\circ_{1,\alpha}(j_1)+j_2m^\circ_{1,\alpha}(j_2)+j_3m^\circ_{1,\alpha}(j_3) = 0.
\]
Therefore, we have $j_1\left( m^\circ_{1,\alpha}(j_1)-m^\circ_{1,\alpha}(j_3) \right) + j_2\left( m^\circ_{1,\alpha}(j_2)-m^\circ_{1,\alpha}(j_3) \right) = 0$. Thanks to \ref{extension_ofoperators2} of Proposition~\ref{nichts2}, we have that 
\begin{align}\label{evensdsysd}
j\mapsto m^\circ_{1,\alpha}(j)\text{ is monotone increasing and even,}
\end{align} which yields a contradiction. This proves the claim~\eqref{sssclaim}. Therefore,  we can define $F^{(3)}$ to be 
\begin{align}\label{def_F3}
F^{(3)}_{j_1,j_2,j_3} = 
\begin{cases}
\frac{2\pi H_{3,j_1,j_2,j_3}}{2\ii \left(j_1\kappa_{j_1} +j_2\kappa_{j_2} +j_3\kappa_{j_3} \right)}  & \text{ if $(j_1,j_2,j_3) \in \mathcal{A}_3$},\\
0 & \text{ otherwise.}
\end{cases}
\end{align}
From \eqref{H33expansion}, \eqref{sssperp2} and \eqref{def_F3}, we obtain \eqref{m_H3}. Indeed,  the modified Hamiltonian does not change the terms involving more than $2$ components in the normal site, that is 
\begin{align}\label{sd11sdsdx2x}
H^{(3)}_3(f) = \mathcal{H}_{3,2}+\mathcal{H}_{3,3},
\end{align}
which gives \eqref{m_H3}.

\vspace{0.5\baselineskip}
\noindent\textit{Proof of (2).} From \eqref{m_H3}, and $\mathcal{H}^{(3)}_3$ in \eqref{H33expansion}, we have
\begin{align}\label{h2andf3}
\left\{ \mathcal{H}_2, F^{(3)}\right\} = -\mathcal{H}_{3,\le 1}:=-(\mathcal{H}_{3,0} + \mathcal{H}_{3,1}),
\end{align}
and 
\begin{align}\label{H_34}
{H}^{(3)}_4 = \mathcal{H}_4 + \frac{1}{2}\left\{ \mathcal{H}_{3,\le 1},F^{(3)}\right\} + \left\{ {H}^{(3)}_3,F^{(3)}\right\}.
\end{align}
To compute the Poisson bracket, we write $F^{(3)}$ as
\begin{align}\label{def_F3333333}
F^{(3)}(f) = 3\sum_{\substack{j_1+j_2+j_3 = 0 \\ j_1,j_2\in S\\ j_3\in S^\perp}}F^{(3)}_{j_1,j_2,j_3}f_{j_1}f_{j_2}f_{j_3},
\end{align}
which follows from \eqref{def_F3} and that there is no $j_1,j_2,j_3\in S$ such that $j_1+j_2+j_3 =0$ (see \ref{tangent_2}).  Therefore, we have that for $h\in L^2(\mathbb{T})$, 
\begin{align}\nonumber
\left( \nabla F^{(3)}(f), h \right)_{L^2} &= 2\pi \cdot 6\sum_{j_1 \in S}\left(\sum_{\substack{j_2+j_3 = j_1 \\ j_2\in S, \\ j_3\in S^\perp}}F^{(3)}_{-j_1,j_2,j_3}f_{j_2}f_{j_3} \right)\overline{h}_{j_1} \\
& + 2\pi \cdot 3\sum_{j_3 \in S^\perp}\left(\sum_{\substack{j_1+j_2 = j_3 \\j_1, j_2\in S}}F^{(3)}_{j_1,j_2,-j_3}f_{j_1}f_{j_2} \right)\overline{h}_{j_3}\nonumber\\
& =:  6\int  \Pi_{S}K_1(v,z)h dx + 3\int \Pi_{S^\perp}K_2(v,v)hdx,\label{rlsd1sdfff}
\end{align}
where 
\begin{align}\label{K12def}
K_1(v,z) =\sum_{j_2+j_3 = j}F^{(3)}_{-j,j_2,j_3}v_{j_2}z_{j_3}e^{\ii j x}, \, \text{ and }\, K_2(v,v) =\sum_{j_1+j_2 = j}F^{(3)}_{j_1,j_2,-j}v_{j_1}v_{j_2}e^{\ii j x}.
\end{align}
From \eqref{K12def}, one can easily see that
\begin{align}\label{K1K2relation}
\int v\Pi_{S}K_1(w,z)dx = \int z\Pi_{S^\perp}K_2(w,v)dx, \quad \text{ for $v,w\in H_S$ and $z\in H_{S^\perp}$}.
\end{align}
Hence \eqref{rlsd1sdfff}  tells us that \begin{align}\label{gradf2s2s}
\nabla F^{(3)}(f) = 6\Pi_SK_1(v,z) + 3\Pi_{S^\perp}K_2(v,v).
\end{align}
Furthermore, from $\mathcal{H}_3$ and $\mathcal{H}_{3,\le 1}$ in \eqref{expansion_2} and \eqref{h2andf3}, we see that 
\begin{align}
\nabla \mathcal{H}_{3,\le 1}(f) & =  \frac{\alpha}{16}\Pi_S\left(2v\Lambda^{\alpha-1} v + \Lambda^{\alpha-1}(v^2) + 2v\Lambda^{\alpha-1}z + 2z\Lambda^{\alpha-1} v + 2\Lambda^{\alpha-1}(vz) \right) \nonumber \\
& - \frac{\alpha}{16}T_\alpha\Pi_S \left(v^2+2vz \right) +\frac{\alpha}{16}\Pi_{S^\perp}\left(\Lambda^{\alpha-1}(v^2) + 2v\Lambda^{\alpha-1} v \right)  \nonumber \\&-\frac{\alpha}{16}T_\alpha\Pi_{S^\perp}(v^2).\label{131hsd1}
\end{align}
and from \eqref{m_H3}, we have
\begin{align}\nonumber
\nabla {H}^{(3)}_3(f) &= \frac{\alpha}{16}\Pi_S\left(2z\Lambda^{\alpha-1}z + \Lambda^{\alpha-1}(z^2) \right) - \frac{\alpha}{16}T_\alpha \Pi_S z^2 \nonumber\\
& \ + \frac{\alpha}{16}\Pi_{S^\perp}\left( 2v\Lambda^{\alpha-1}z + 2\Lambda^{\alpha-1}(vz) + 2z \Lambda^{\alpha-1}v + 2z \Lambda^{\alpha-1}z + \Lambda^{\alpha-1}(z^2)\right) \nonumber \\& - \frac{\alpha}{16}T_\alpha \Pi_{S^\perp}\left( 2vz + z^2\right).
\label{h2msd1x}
\end{align}
Using \eqref{gradf2s2s}, \eqref{131hsd1} and the definition of the Poisson bracket in \eqref{poisson_bracket}, we have
\begin{align}\label{eq10}
\left\{ \mathcal{H}_{3,\le 1} , F^{(3)}\right\} & = \frac{3\alpha}{16}\int \left( \Lambda^{\alpha-1}(v^2) + 2v \Lambda^{\alpha-1} v - T_\alpha v^2 \right)\partial_x\Pi_{S^\perp}K_2(v,v)dx\nonumber \\
& \ +\underbrace{ \frac{3\alpha}{8}\int \left(2v\Lambda^{\alpha-1} v + \Lambda^{\alpha-1}(v^2) - T_\alpha v^2\right)\partial_x \Pi_SK_1(v,z) dx}_{=0 \text{ thanks to \ref{tangent_2}}} \nonumber\\
& \ + \frac{3\alpha}{4}\int \left(v\Lambda^{\alpha-1}z + z\Lambda^{\alpha-1} v + \Lambda^{\alpha-1}(vz) - T_\alpha vz \right)\partial_x \Pi_{S}K_1(v,z)dx.
\end{align}
Similarly, using \eqref{h2msd1x} and \eqref{gradf2s2s}, we have
\begin{align}\label{eq11}
\left\{ {H}^{(3)}_3,F^{(3)}\right\} & =  \frac{3\alpha}{8}\int \left( v\Lambda^{\alpha-1}z + \Lambda^{\alpha-1}(vz) + z\Lambda^{\alpha-1}v - T_\alpha vz \right)\partial_x \Pi_{S^\perp}K_2(v,v)dx \nonumber\\
& \ + \frac{3\alpha}{16} \int \left( 2z\Lambda^{\alpha-1}z + \Lambda^{\alpha-1}(z^2) - T_\alpha z^2\right)\partial_x \Pi_{S^\perp}K_2(v,v)dx\\
& \ +\frac{3\alpha}{8}\int \left(2z\Lambda^{\alpha-1}z + \Lambda^{\alpha-1}(z^2) - T_\alpha z^2 \right)\partial_x \Pi_S K_1(v,z)dx \nonumber
\end{align}
Collecting the terms involving only tangential sites or  two normal sites, we find that
\begin{align*}
\frac{1}{2}&\left\{ \mathcal{H}_{3,\le 1},F^{(3)}\right\} + \left\{ {H}^{(3)}_3,F^{(3)}\right\}\\
& = \frac{3\alpha}{32}\int \left( \Lambda^{\alpha-1}(v^2) + 2v \Lambda^{\alpha-1} v - T_\alpha v^2 \right)\partial_x\Pi_{S^\perp}K_2(v,v)dx\\
& \ + \frac{3\alpha}{8}\int \left(v\Lambda^{\alpha-1}z + z\Lambda^{\alpha-1} v + \Lambda^{\alpha-1}(vz) - T_\alpha vz \right)\partial_x \Pi_{S}K_1(v,z)dx \\
& \ + \frac{3\alpha}{16} \int \left( 2z\Lambda^{\alpha-1}z + \Lambda^{\alpha-1}(z^2) - T_\alpha z^2\right)\partial_x \Pi_{S^\perp}K_2(v,v)dx\\
& \ + R(v^3z) + R(vz^3)+R(z^4),
\end{align*}
for some multilinear map $R$. With \eqref{H_34}, this proves \eqref{quartic_11122}. Also, Lemma~\ref{momentum1} implies that $H^{(3)}_5$ satisfies \eqref{conserve_momentum}. 

\vspace{0.5\baselineskip}
\noindent\textit{Proof of (3).} $\mathcal{A}_3$ is a finite set since it is contained in $[-2K,2K]^3$, where $K = \max_{i=1,\ldots \nu, \quad j_i \in S} |j_{i}|$

\vspace{0.5\baselineskip}
\noindent\textit{Proof of (4).} Since $\mathcal{H}$ is real-valued Hamiltonian\index{real-valued Hamiltonian}, $H_{3,j_1,j_2,j_3}$ satisfies \eqref{real_bonsdsd}. Using that $j\mapsto \kappa_j$ is even, we see that $F^{(3)}$ in \eqref{def_F3} also satisfies \eqref{real_bonsdsd}, thus $F^{(3)}$ is a real-valued Hamiltonian as well. Since $\Phi^{(3)}$ is the time-$1$ flow map of Hamiltonian PDE, generated by $F^{(3)}$, $\Phi^{(3)}$ maps a real-valued function to a real-valued function. For the reversibility preserving property, it follows immediately from the definition of $F^{(3)}$ in \eqref{def_F3}, that $F^{(3)}\circ \rho = - F^{(3)}$ with the involution $\rho$ in  \eqref{definiteionofinv}  (see Section~\ref{momentum} and the definition of $\kappa_j$ in \eqref{fourier_hamiltonian}). Then the result follows from Lemma~\ref{time_reversibility_map}.  
\end{proof}

We specify the Hamiltonian $H^{(3)}$ in terms of the Fourier modes:
\begin{lemma}\label{fourier_H3333s2s}
We have
\begin{align}
H_3^{(3)}(f)&=\sum_{\substack{j_1+j_2+j_3=0,\\ (j_1,j_2,j_3)\in \mathcal{A}_3^c}}H_{3,j_1,j_2,j_3}f_{j_1}f_{j_2}f_{j_3},\label{cucucucuc}\\
H^{(3)}_{4,0}(f) &= \sum_{\substack{j_1+j_2+j_3+j_4=0,\\ j_1,j_2,j_3,j_4\in S} }H_{4,j_1,j_2,j_3,j_4}f_{j_1}f_{j_2}f_{j_3}f_{j_4} \nonumber\\
& \  -\frac{9}{4}\sum_{\substack{j_1+j_2+j_3+j_4=0,\\ j_1,j_2,j_3,j_4\in S,\\ j_1+j_2\in S^\perp}}(j_1+j_2)\frac{H_{3,j_1,j_2,-(j_1+j_2)}H_{3,j_3,j_4,-(j_3+j_4)}}{j_3\kappa_{j_3}+j_4\kappa_{j_4}-(j_3+j_4)\kappa_{j_3+j_4}}f_{j_1}f_{j_2}f_{j_3}f_{j_4}\label{quququququ}\\
H_{4,2}^{(3)}(f) &= \sum_{\substack{j_1+j_2+j_3+j_4=0,\\\text{exactly two of $j_1,j_2,j_3,j_4$ are in $S^\perp$}}}H_{4,j_1,j_2,j_3,j_4}f_{j_1}f_{j_2}f_{j_3}f_{j_4}\nonumber\\
& \ -9\sum_{\substack{j_1+j_2+j_3+j_4=0,\\ j_1,j_2\in S,\ j_3,j_4\in S^\perp,\\ j_2+j_3\in S}}(j_2+j_3)\frac{H_{-(j_2+j_3),j_2,j_3}H_{3,j_1,j_4,-(j_1+j_4)}}{j_1\kappa_{j_1}+j_4\kappa_{j_4}-(j_1+j_4)\kappa_{j_1+j_4}}f_{j_1}f_{j_2}f_{j_3}f_{j_4}\nonumber\\
& \ -\frac{9}{2}\sum_{\substack{j_1+j_2+j_3+j_4=0,\\ j_1,j_2\in S,\ j_3,j_4\in S^\perp,\\ j_1+j_2\in S^\perp}}(j_3+j_4)\frac{H_{3,j_1,j_2,-(j_1+j_2)}H_{3,j_3,j_4,-(j_3+j_4)}}{j_1\kappa_{j_1}+j_2\kappa_{j_2}-(j_1+j_2)\kappa_{j_1+j_2}}f_{j_1}f_{j_2}f_{j_3}f_{j_4} \label{sswsw2s}
\end{align}
\end{lemma}
\begin{proof}
The expression \eqref{cucucucuc} follows immediately from \eqref{sssperp}, \eqref{sssperp2} and \eqref{def_F3}.

Now, let us derive \eqref{quququququ} and \eqref{sswsw2s}. In view of \eqref{H_34}, we compute $\nabla \mathcal{H}_{3,\le 1}$ and $\nabla{H}^{(3)}_3$ in terms of the Fourier modes. For $\nabla \mathcal{H}_{3,\le 1}$, it follows from \eqref{h2andf3} and \eqref{defi_1sd3a} that  (recalling \ref{tangent_2} in the condition of the tangential sites)
\[\mathcal{H}_{3,\le 1}(f)=\sum_{\substack{j_1+j_2+j_3=0,\\(j_1,j_2,j_3)\in\mathcal{A}_3}}H_{3,j_1,j_2,j_3}f_{j_1}f_{j_2}f_{j_3} =  3\sum_{\substack{j_1+j_2+j_3=0,\\ j_1,j_2\in S, j_3\in S^\perp}}H_{3,j_1,j_2,j_3}f_{j_1}f_{j_2}f_{j_3}.
\]
Hence, we have
\begin{align}\label{rjsdsd2sdsdsdsdssds}
\partial_{f_{-j}} \mathcal{H}_{3,\le 1}(f)=\begin{cases}
 6\sum_{\substack{j_2+j_3=j,\\ j_2\in S,\ j_3\in S^\perp}}H_{3,-j,j_2,j_3}f_{j_2}f_{j_3}& \text{ if $j\in S$,}\\
3\sum_{\substack{j_1+j_2=j,\\ j_1,j_2\in S}}H_{3,j_1,j_2,-j}f_{j_1}f_{j_2}& \text{ if $j\in S^\perp$}.
\end{cases}
\end{align}
Similarly, we have for $H^{(3)}_3$ (see \eqref{cucucucuc}), 
\begin{align*}
H^{(3)}_3(f)& =\sum_{\substack{j_1+j_2+j_3=0,\\(j_1,j_2,j_3)\in\mathcal{A}_3^c}}H_{3,j_1,j_2,j_3}f_{j_1}f_{j_2}f_{j_3}\\
& = 3\sum_{\substack{j_1+j_2+j_3=0,\\ j_1\in S, j_2, j_3\in S^\perp}}H_{3,j_1,j_2,j_3}f_{j_1}f_{j_2}f_{j_3} + \sum_{\substack{j_1+j_2+j_3=0,\\ j_1,j_2, j_3\in S^\perp}}H_{3,j_1,j_2,j_3}f_{j_1}f_{j_2}f_{j_3},
\end{align*}
therefore
\begin{align}\nonumber
& \partial_{f_{-j}}H_3^{(3)}(f)\\
& =\begin{cases}
 3\sum_{\substack{j_2+j_3=j,\\ j_2,j_3\in S^\perp}}H_{3,-j,j_2,j_3}f_{j_2}f_{j_3} & \text{ if $j\in S$,}\\
  6\sum_{\substack{j_1+j_2 = j,\\ j_1\in S,\ j_2\in S^\perp}}H_{3,j_1,j_2,-j}f_{j_1}f_{j_2} + 3\sum_{\substack{j_1+j_2 = j,\\ j_1,j_2\in S^\perp}}H_{3,j_1,j_2,-j}f_{j_1}f_{j_2} & \text{ if $j\in S^\perp$}.
\end{cases}\label{threelosers1}
\end{align}
Also, from \eqref{F3def} and \eqref{def_F3}, it follows that
\[
F^{(3)}(f) = \sum_{\substack{j_1+j_2+j_3=0,\\ (j_1,j_2,j_3)\in \mathcal{A}_3}}F^{(3)}_{j_1,j_2,j_3}f_{j_1}f_{j_2}f_{j_3},
\]
thus (again using \ref{tangent_2} in the tangential site condition),
\begin{align}\label{ffoursd2ssxzz}
\partial_{f_j}F^{(3)}(f)=\begin{cases}
6\sum_{\substack{j_1+j_2=-j,\\ j_1\in S,\ j_2\in S^\perp}}F^{(3)}_{j_1,j_2,j}f_{j_1}f_{j_2}& \text{ if $j\in S$},\\
3\sum_{\substack{j_1+j_2=-j,\\ j_1,j_2\in S}}F^{(3)}_{j_1,j_2,j}f_{j_1}f_{j_2}& \text{ if $j\in S^\perp$}.
\end{cases}
\end{align}
Therefore, using \eqref{poisson_bracket}, \eqref{rjsdsd2sdsdsdsdssds} and \eqref{ffoursd2ssxzz}, we can find
\begin{align*}
& \left\{\mathcal{H}_{3,\le 1},F^{(3)} \right\}\\
& = -\frac{9}{2\pi}\sum_{\substack{j_1+j_2+j_3+j_4=0,\\ j_1,j_2,j_3,j_4\in S,\\ j_1+j_2\in S^\perp}}\ii (j_1+j_2)H_{3,j_1,j_2,-(j_1+j_2)}F^{(3)}_{j_3,j_4,-(j_3+j_4)}f_{j_1}f_{j_2}f_{j_3}f_{j_4}\\
& \ -\frac{36}{2\pi}\sum_{\substack{j_1+j_2+j_3+j_4=0,\\ j_1,j_2\in S,\ j_3,j_4\in S^\perp,\\ j_2+j_3\in S}}\ii(j_2+j_3)H_{-(j_2+j_3),j_2,j_3}F^{(3)}_{j_1,j_4,-(j_1+j_4)}f_{j_1}f_{j_2}f_{j_3}f_{j_4},
\end{align*}
and using \eqref{threelosers1} and \eqref{ffoursd2ssxzz}, 
\begin{align*}
& \left\{ H_3^{(3)},F^{(3)}\right\}\\
& = -\frac{9}{2\pi}\sum_{\substack{j_1+j_2+j_3+j_4=0,\\ j_1,j_2\in S,\ j_3,j_4\in S^\perp,\\ j_1+j_2\in S^\perp}}\ii(j_3+j_4)F^{(3)}_{j_1,j_2,-(j_1+j_2)}H_{3,j_3,j_4,-(j_3+j_4)}f_{j_1}f_{j_2}f_{j_3}f_{j_4} \\
& + R(vz^3)+R(v^3z).
\end{align*}
Therefore, it follows from the above expressions and \eqref{def_F3} that
\begin{align*}
& \frac{1}{2}\left\{ \mathcal{H}_{3,\le 1},F^{(3)}\right\} + \left\{H^{(3)}_3, F^{(3)}\right\} \\
&= -\frac{9}{4}\sum_{\substack{j_1+j_2+j_3+j_4=0,\\ j_1,j_2,j_3,j_4\in S,\\ j_1+j_2\in S^\perp}}(j_1+j_2)\frac{H_{3,j_1,j_2,-(j_1+j_2)}H_{3,j_3,j_4,-(j_3+j_4)}}{j_3\kappa_{j_3}+j_4\kappa_{j_4}-(j_3+j_4)\kappa_{j_3+j_4}}f_{j_1}f_{j_2}f_{j_3}f_{j_4}\\
& \ -9\sum_{\substack{j_1+j_2+j_3+j_4=0,\\ j_1,j_2\in S,\ j_3,j_4\in S^\perp,\\ j_2+j_3\in S}}(j_2+j_3)\frac{H_{-(j_2+j_3),j_2,j_3}H_{3,j_1,j_4,-(j_1+j_4)}}{j_1\kappa_{j_1}+j_4\kappa_{j_4}-(j_1+j_4)\kappa_{j_1+j_4}}f_{j_1}f_{j_2}f_{j_3}f_{j_4}\\
& \ -\frac{9}{2}\sum_{\substack{j_1+j_2+j_3+j_4=0,\\ j_1,j_2\in S,\ j_3,j_4\in S^\perp,\\ j_1+j_2\in S^\perp}}(j_3+j_4)\frac{H_{3,j_1,j_2,-(j_1+j_2)}H_{3,j_3,j_4,-(j_3+j_4)}}{j_1\kappa_{j_1}+j_2\kappa_{j_2}-(j_1+j_2)\kappa_{j_1+j_2}}f_{j_1}f_{j_2}f_{j_3}f_{j_4} \\
& + R(vz^3) + R(v^3 z).
\end{align*}
 Plugging this into \eqref{H_34}, we obtain \eqref{quququququ} and \eqref{sswsw2s}.
\end{proof}

Now we turn to the quartic Hamiltonian. The normalization of the quartic Hamiltonian relies on the non-existence of nontrivial resonances, which we will prove in Lemma~\ref{sssslemma}.

\begin{lemma}\label{quartic}(Normalization\index{Normalization} of the quartic Hamiltonian) There exists a symplectic transformation $\Phi^{(4)}$ such that 
\begin{align*}
H^{(4)}(f):=H^{(3)}\circ \Phi^{(4)}(f) = H^{(4)}_2(f) + H^{(4)}_3(f) +H^{(4)}_4(f) + H^{(4)}_{\ge 5}(f),
\end{align*}
where
\begin{enumerate}[label=\text{(\arabic*)}]
\item \label{232xdsd22}(Quadratic and cubic terms) We have
\begin{equation}\label{m_H4222}
\begin{aligned}
H^{(4)}_2(f) = \mathcal{H}_2(f), \quad H^{(4)}_3(f) = {H}^{(3)}_3(f)
\end{aligned}
\end{equation} 
\item \label{item2_lemsss}(Quartic and higher degree terms) Using the Fourier series expansion in  \eqref{fourier}, let us rewrite $H^{(3)}_{4,0}$ in \eqref{quartic_11122} as
\begin{align}\label{mutliplier1321}
H^{(3)}_{4,0} =: \sum_{\substack{j_1,j_2,j_3,j_4\in S \\ j_1+j_2+j_3+j_4=0}}H^{(3)}_{4,j_1,j_2,j_3,j_4}v_{j_1}v_{j_2}v_{j_3}v_{j_4},
\end{align} where $H^{(3)}_{4,j_1,j_2,j_3,j_4}$ is invariant under any permutation on $\left\{ j_1,j_2,j_3,j_4\right\}$. Then we have
\begin{equation}\label{quartic_1112}
\begin{aligned}
H^{(4)}_4(f) &= \sum_{i=0}^4H^{(4)}_{4,i}(f), \text{ where }\\
H^{(4)}_{4,0} &= 6\sum_{j_1\in S^{+}}H^{(3)}_{4,j_1,-j_1,j_1,-j_1}|v_{j_1}|^4 + 12\sum_{\substack{j_1,j_2\in S^+,\\  j_1\ne j_2}}H^{(3)}_{4,j_1,-j_1,j_2,-j_2}|v_{j_1}|^2|v_{j_2}|^2,\\
 \quad H^{(4)}_{4,1} &= 0, \quad H^{(4)}_{4,i}=H^{(3)}_{4,i}, \text{ for }i=2,3,4,
\end{aligned}
\end{equation}
and $H^{(4)}_{\ge 5}(f)$ collects all the terms homogeneous of degree at least $5$ and $H^{(4)}_5$ satisfies \eqref{conserve_momentum}. 
\item \label{232xdsd222}There exists a finite dimensional space $E$ of the form in \eqref{finite_dim} such that $\Phi^{(4)}(f) = f + \Psi^{(4)}(f) $ where $\Psi^{(4)}$ satisfies $\Psi^{(4)}(f) = \Pi_E \Psi^{(4)}(\Pi_E f)$.
\item\label{232xdsd2223} $\Phi^{(4)}$ is  real and reversibility preserving.
\end{enumerate}
\end{lemma}

\begin{proof}
We will prove \ref{item2_lemsss} first and then prove  \ref{232xdsd22}.

\vspace{0.5\baselineskip}
\noindent\textit{Proof of \ref{item2_lemsss}.} Let us consider a Hamiltonian:
\begin{align}\label{F4def}
F^{(4)}(f) := \sum_{j_1+j_2+j_3+j_4 = 0}F^{(4)}_{j_1,j_2,j_3,j_4}f_{j_1}f_{j_2}f_{j_3}f_{j_4},
\end{align}
which will be determined later and let us denote its time-1 flow map by $\Phi^{(4)}$.  We compute using Lemma~\ref{symplectic_1},
\begin{align}\nonumber
{H}^{(4)} & := {H}^{(3)}\circ \Phi^{(4)} = H^{(3)}_2 + {H}^{(3)}_3 + {H}^{(4)}_4 + {H}^{(4)}_{\ge 5}, \quad \text{ where }\\
\quad {H}^{(4)}_4 & = {H}^{(3)}_4 + \left\{ {H}_2^{(3)},F^{(4)}\right\}.
\label{H44expansion}
\end{align}
From \eqref{m_H3}, we have $H^{(3)}_2=\mathcal{H}_2$. Also, it follows from \eqref{poisson_bracket}, \eqref{F4def} and \eqref{fourier_hamiltonian} that
\begin{align*}
 \left\{ \mathcal{H}_2,F^{(4)}\right\} &= -\sum_{j_1+j_2+j_3+j_4=0}\frac{8\ii}{2\pi} j_1\kappa_{j_1}F^{(4)}_{j_1,j_2,j_3,j_4}f_{j_1}f_{j_2}f_{j_3}f_{j_4} \\
&=-\frac{2\ii}{2\pi}\sum_{j_1+j_2+j_3+j_4=0}\left(j_1\kappa_{j_1} +j_2\kappa_{j_2} +j_3\kappa_{j_3}  + j_4\kappa_{j_4}\right)F^{(4)}_{j_1,j_2,j_3,j_4}f_{j_1}f_{j_2}f_{j_3}f_{j_4}
\end{align*}
Therefore, we have
\begin{align}
{H}^{(4)}_4  &=\sum_{j_1+j_2+j_3+j_4 = 0} \left( H^{(3)}_{4,j_1,j_2,j_3,j_4}- \frac{2\ii}{2\pi} \left(j_1\kappa_{j_1} +j_2\kappa_{j_2} +j_3\kappa_{j_3}  + j_4\kappa_{j_4}\right)F^{(4)}_{j_1,j_2,j_3,j_4}\right) \nonumber\\
& \qquad \qquad \qquad \qquad \times f_{j_1}f_{j_2} f_{j_3}f_{j_4}\nonumber\\
& = \sum_{\substack{j_1+j_2+j_3+j_4 = 0 \\ (j_1,j_2,j_3,j_4 )\in \mathcal{A}_4}} \left( H^{(3)}_{4,j_1,j_2,j_3,j_4}- \frac{2\ii}{2\pi} \left(j_1\kappa_{j_1} +j_2\kappa_{j_2} +j_3\kappa_{j_3}  + j_4\kappa_{j_4} \right)F^{(4)}_{j_1,j_2,j_3,j_4}\right)\nonumber\\
& \qquad \qquad \qquad \qquad \times f_{j_1}f_{j_2} f_{j_3}f_{j_4}\nonumber\\
& \ + \sum_{\substack{j_1+j_2+j_3+j_4 = 0 \\ (j_1,j_2,j_3,j_4 )\notin \mathcal{A}_4}} \left( H^{(3)}_{4,j_1,j_2,j_3,j_4}- \frac{2\ii}{2\pi} \left(j_1\kappa_{j_1} +j_2\kappa_{j_2} +j_3\kappa_{j_3}  + j_4\kappa_{j_4} \right)F^{(4)}_{j_1,j_2,j_3,j_4}\right)\nonumber\\
& \qquad \qquad \qquad \qquad \times f_{j_1}f_{j_2} f_{j_3}f_{j_4},\label{ssssperp4}
\end{align}
where \begin{align*}
\mathcal{A}_4 := \left\{(j_1,j_2,j_3,j_4)\in \mathbb{Z}^4 \right.& : j_1\kappa_{j_1} +j_2\kappa_{j_2} +j_3\kappa_{j_3}  + j_4\kappa_{j_4} \ne 0  \\
& \left. \text{  and at least three of $\left\{j_1,j_2,j_3,j_4 \right\}$  are in $S$ }\right\}.\end{align*}
Hence we can define $F^{(4)}$ to be
\begin{align}\label{def_F4}
F^{(4)}_{j_1,j_2,j_3,j_4} = 
\begin{cases}
\frac{2\pi H^{(3)}_{4,j_1,j_2,j_3,j_4}}{2\ii \left(j_1\kappa_{j_1} +j_2\kappa_{j_2} +j_3\kappa_{j_3}  + j_4\kappa_{j_4}  \right)}  & \text{ if $(j_1,j_2,j_3,j_4) \in \mathcal{A}_4$},\\
0 & \text{ otherwise},
\end{cases}
\end{align}
so that \eqref{ssssperp4} yields that
\begin{align}\label{ssssperp5}
{H}^{(4)}_4 = \sum_{\substack{j_1+j_2+j_3+j_4 = 0 \\ (j_1,j_2,j_3,j_4)\notin \mathcal{A}_4}} H^{(3)}_{4,j_1,j_2,j_3,j_4}f_{j_1}f_{j_2}f_{j_3}f_{j_4}.
\end{align}
Recalling the notation in \eqref{degree_on_normal},  we see that
\[
{H}^{(4)}_{4,i} = {H}^{(3)}_{4,i} \quad \text{ for }i=2,3,4,
\]
thanks to \eqref{ssssperp5} and the definition of $\mathcal{A}_4$, which show that if  at least two of $j_1,\ldots j_4$ are in $S^\perp$, then $(j_1,\ldots,j_4)\notin \mathcal{A}_4$. Furthermore, it follows from  Lemma~\ref{sssslemma} and \eqref{tan_site} that there is no $(j_1,j_2,j_3,j_4)\in \mathcal{A}_4^c$ such that $j_1,j_2,j_3\in S$ and $j_4\in S^\perp$, therefore
\begin{align}\label{Sperp_nochange}
{H}^{(4)}_{4,1} = 0.
\end{align}
Thus, it follows from  \eqref{ssssperp5} and Lemma~\ref{sssslemma} that
\begin{equation}\label{ssss}
\begin{aligned}
{H}^{(4)}_{4,0} &= \sum_{\substack{j_1+j_2+j_3+j_4 = 0, \\ j_1\kappa_{j_1} +j_2\kappa_{j_2} +j_3\kappa_{j_3}  + j_4\kappa_{j_4}= 0 \\ j_1,j_2,j_3,j_4 \in S \\}}H^{(3)}_{4,j_1,j_2,j_3,j_4}v_{j_1}v_{j_2}v_{j_3}v_{j_4}\\
& = 6\sum_{\substack{j_1+j_2+j_3+j_4=0,\\j_1\kappa_{j_1} +j_2\kappa_{j_2} +j_3\kappa_{j_3}  + j_4\kappa_{j_4}= 0,\\  j_1,j_2\in S^+,\quad j_3,j_4\in S\backslash S^{+}}}H^{(3)}_{4,j_1,j_2,j_3,j_4}v_{j_1}v_{j_2}v_{j_3}v_{j_4}\\
& = 6\sum_{\substack{j_1+j_2+j_3+j_4=0,\\j_1\kappa_{j_1} +j_2\kappa_{j_2} +j_3\kappa_{j_3}  + j_4\kappa_{j_4}= 0,\\  j_1,j_2\in S^+,\quad j_3,j_4\in S\backslash S^{+},\\ j_1=j_2}}H^{(3)}_{4,j_1,j_2,j_3,j_4}v_{j_1}v_{j_2}v_{j_3}v_{j_4} \\
& + 6\sum_{\substack{j_1+j_2+j_3+j_4=0,\\j_1\kappa_{j_1} +j_2\kappa_{j_2} +j_3\kappa_{j_3}  + j_4\kappa_{j_4}= 0,\\  j_1,j_2\in S^+,\quad j_3,j_4\in S\backslash S^{+},\\ j_1\ne j_2}}H^{(3)}_{4,j_1,j_2,j_3,j_4}v_{j_1}v_{j_2}v_{j_3}v_{j_4}\\
& = 6\sum_{j_1\in S^{+}}H^{(3)}_{4,j_1,-j_1,j_1,-j_1}|v_{j_1}|^4 + 12\sum_{\substack{j_1,j_2\in S^+,\\ j_1\ne j_2, \\ j_1=-j_3,\ j_2=-j_4}}H^{(3)}_{4,j_1,j_2,j_3,j_4}v_{j_1}v_{j_2}v_{j_3}v_{j_4}\\
& =  6\sum_{j_1\in S^{+}}H^{(3)}_{4,j_1,-j_1,j_1,-j_1}|v_{j_1}|^4 + 12\sum_{j_1,j_2\in S^+,\ j_1\ne j_2}H^{(3)}_{4,j_1,-j_1,j_2,-j_2}|v_{j_1}|^2|v_{j_2}|^2.
\end{aligned}
\end{equation}
  This proves \eqref{quartic_1112}. Also, Lemma~\ref{momentum1} implies that $H^{(4)}_5$ satisfies \eqref{conserve_momentum}.

\vspace{0.5\baselineskip}
\noindent\textit{Proof of \ref{232xdsd22}.} This immediately follows from \eqref{H44expansion} since the symplectic transformation $\Phi^{(4)}$ does not change the quadratic and cubic terms (see \eqref{H44expansion} and \eqref{m_H3}).

\vspace{0.5\baselineskip}
\noindent\textit{Proof of \ref{232xdsd222}.}: $\mathcal{A}_4$ is a finite set since it is contained in $[-3K,3K]^4$, where $K = \max_{i=1,\ldots \nu, \quad j_i \in S} |j_{i}|$

\vspace{0.5\baselineskip}
\noindent\textit{Proof of \ref{232xdsd2223}.} The proof is identical to (4) in Lemma~\ref{ssslem}.
\end{proof}

\begin{lemma}\label{sssslemma}
Let $(j_1,j_2,j_3,j_4) \in \mathbb{Z}^{4}$, $j_i \neq 0$ be solutions of

\begin{align*}
j_1 + j_2 + j_3 + j_4 = 0, \quad j_1\kappa_{j_1} +j_2\kappa_{j_2} +j_3\kappa_{j_3}  + j_4\kappa_{j_4}= 0
\end{align*}
Then the only solutions are given by $j_a = -j_b, j_c = -j_d$ where $(a,b,c,d)$ is a permutation of $(1,2,3,4)$.
\end{lemma}
\begin{proof}
It is enough to consider the case when two of the $j$'s are positive and two are negative since if sign$(j_1)$ = sign$(j_2)$ = sign$(j_3)$ we can rewrite the equation as
\begin{align*}
0 = j_1 (\kappa_{j_1} - \kappa_{j_1+j_2+j_3})
+ j_2 (\kappa_{j_2} -\kappa_{j_1+j_2+j_3})
+ j_3 (\kappa_{j_3} - \kappa_{j_1+j_2+j_3})
\end{align*}
and we get a contradiction by the monotonicity of $\kappa_{j}$. Without loss of generality we may assume that $0 < j_1, j_4$ and $0 > j_2, j_3$. We claim that the following equation (in $j_1$) has exactly two solutions ($j_1 = -j_2$, $j_1 = -j_3$):

\begin{align}
\label{eq_reduced_j1}
F(j_1) := j_1 \kappa_{j_1} + j_2 \kappa_{j_2}  + j_3 \kappa_{j_3} - (j_1 + j_2 + j_3)\kappa_{j_1+j_2+j_3} = 0
\end{align}

To do so, we will prove that for $\alpha \geq 1$, $F$ is increasing for $j_1 \geq -j_2 - j_3$ and that $F$ is convex for $0 \leq j_1 \leq -j_2 - j_3$ and for $\alpha < 1$, $F$ is decreasing for $j_1 \geq -j_2 - j_3$ and $F$ is concave for $0 \leq j_1 \leq -j_2 - j_3$. Together with $\lim_{\pm \infty} F(j) = +\infty$ for $\alpha \geq 1$ and $\lim_{\pm \infty} F(j) = -\infty$ for $\alpha < 1$ this ensures that there are at most 2 solutions, which are the ones given above.

Under the hypotheses $j_1 > 0, j_1 + j_2 + j_3 \geq 0$ and using the definition of $\kappa_{j}$ in \eqref{fourier_hamiltonian} and \eqref{explicit_multiplier}, we obtain, for a positive constant $C$:

\begin{align*}
 F(j_1+1) - F(j_1) & = C\alpha\left(\left(j_1 + \frac12\right)\frac{\Gamma\left(\frac{\alpha}{2}+j_1\right)}{\Gamma\left(2-\frac{\alpha}{2}+j_1\right)} \right.\\& \left.- \left(j_1 + j_2 + j_3 + \frac12\right)\frac{\Gamma\left(\frac{\alpha}{2}+j_1 + j_2 + j_3\right)}{\Gamma\left(2-\frac{\alpha}{2}+j_1 + j_2 + j_3\right)}\right),
\end{align*}

and the monotonicity result follows from the fact that
\begin{align*}
\left(j_1 + \frac12\right)\frac{\Gamma\left(\frac{\alpha}{2}+j_1\right)}{\Gamma\left(2-\frac{\alpha}{2}+j_1\right)} & - \left(j_1 + 1 + \frac12\right)\frac{\Gamma\left(\frac{\alpha}{2}+j_1 + 1\right)}{\Gamma\left(2-\frac{\alpha}{2}+j_1 + 1\right)}
\\
& = -(\alpha-1)(1+j_1) \frac{\Gamma\left(\frac{\alpha}{2} + j_1\right)}{\Gamma\left(3-\frac{\alpha}{2} + j_1\right)}.
\end{align*}

To prove the convexity part of the statement, we calculate, for $j_1 \geq 1, j_1 + j_2 + j_3 \leq -1$:

\begin{align*}
& (F(j_1+1) - F(j_1)) - (F(j_1) - F(j_1 - 1))
\\
&= C\alpha(\alpha-1)\left(j_1 \frac{\Gamma\left(-1+\frac{\alpha}{2}+j_1\right)}{\Gamma\left(2-\frac{\alpha}{2}+j_1\right)} - (j_1 + j_2 + j_3) \frac{\Gamma\left(-1+\frac{\alpha}{2}+j_1+j_2+j_3\right)}{\Gamma\left(2-\frac{\alpha}{2}+j_1+j_2+j_3\right)}\right),
\end{align*}

and the convexity result follows from

\begin{align*}
j_1 \frac{\Gamma\left(-1+\frac{\alpha}{2}+j_1\right)}{\Gamma\left(2-\frac{\alpha}{2}+j_1\right)}
& - (j_1+1) \frac{\Gamma\left(-1+\frac{\alpha}{2}+j_1+1\right)}{\Gamma\left(2-\frac{\alpha}{2}+j_1+1\right)} \\
& = (2-\alpha)\left(j_1 + \frac12\right)\frac{\Gamma\left(-1+\frac{\alpha}{2}+j_1\right)}{\Gamma\left(3-\frac{\alpha}{2}+j_1\right)}.
\end{align*}
\end{proof}

Before normalizing the quintic Hamiltonian, we need to have a symmetrized expression of $H^{(3)}_{4,j_1,j_2,j_3,j_4}$ and compute $H_{4,j_1,-j_1,j_2,-j_2}$ for $j_1,j_2\in S^+$ in view of $H^{(4)}_{4,0}$ in \eqref{quartic_1112}.

\begin{lemma}\label{ammsmw}
For $j_1,j_2\in S^+$, $H^{(3)}_{4,j_1,-j_1,j_2,-j_2}$ in \eqref{mutliplier1321} and \eqref{quartic_1112} can be written as
\begin{align*}
& H^{(3)}_{4,j_1,-j_1,j_2,-j_2}=H_{4,j_1,-j_1,j_2,-j_2} \\
& - \begin{cases}
 \frac{3}{4}\frac{(j_1+j_2)(H_{3,j_1,j_2,-(j_1+j_2)})^2}{(j_1+j_2)\kappa_{j_1+j_2}-j_1\kappa_{j_1}-j_2\kappa_{j_2}} & \text{ if $j_1=j_2$,}\\
\frac{3}{4}\left( \frac{(j_1+j_2)(H_{3,j_1,j_2,-(j_1+j_2)})^2}{(j_1+j_2)\kappa_{j_1+j_2}-j_1\kappa_{j_1}-j_2\kappa_{j_2}} +\frac{(j_1-j_2)(H_{3,j_1,-j_2,-(j_1-j_2)})^2}{(j_1-j_2)\kappa_{j_1-j_2} - j_1\kappa_{j_1}+j_2\kappa_{j_2}}\right)& \text{ if $j_1\ne j_2$.}
\end{cases}
\end{align*}
\end{lemma}
\begin{proof}
From \eqref{quququququ} we see that
\begin{align}\label{hhhdeds}
H^{(3)}_{4,0}&=\mathbb{H}_1 - \frac{9}{4}\mathbb{H}_2,\text{ where}\nonumber\\
\mathbb{H}_1&:=\sum_{\substack{j_1+j_2+j_3+j_4=0,\\ j_1,j_2,j_3,j_4\in S} }H_{4,j_1,j_2,j_3,j_4}f_{j_1}f_{j_2}f_{j_3}f_{j_4}\nonumber\\
\mathbb{H}_2&:=\sum_{\substack{j_1+j_2+j_3+j_4=0,\\ j_1,j_2,j_3,j_4\in S,\\ j_1+j_2\in S^\perp}}\underbrace{1_{\left(j_1+j_2\ne 0\right\}}(j_1+j_2)\frac{H_{3,j_1,j_2,-(j_1+j_2)}H_{3,j_3,j_4,-(j_3+j_4)}}{j_3\kappa_{j_3}+j_4\kappa_{j_4}-(j_3+j_4)\kappa_{j_3+j_4}}}_{=:\mathbb{H}_{2,j_1,j_2,j_3,j_4}}f_{j_1}f_{j_2}f_{j_3}f_{j_4} \\
&=\sum_{\substack{j_1+j_2+j_3+j_4=0,\\ j_1,j_2,j_3,j_4\in S}}\mathbb{H}_{2,j_1,j_2,j_3,j_4}f_{j_1}f_{j_2}f_{j_3}f_{j_4}.\nonumber
\end{align}
Note that in the very last equality, we removed the restriction on the modes $j_1+j_2\in S^{\perp}$, which is allowed thanks to the condition on the tangential sites in \ref{tangent_2}.
Hence, we can write $H^{(3)}_{4,j_1,-j_1,j_2,-j_2}$ for $j_1,j_2\in S^+$ in \eqref{mutliplier1321}, in case where $(j_1,j_2,j_3,j_4)=(j_1,-j_1,j_2,-j_2)$, as (noting that $H_{4,j_1,j_2,j_3,j_4}$ in \eqref{moreexplicits2} is already invariant under permutations on the sub-indices) 
\begin{align*}
&H^{(3)}_{4,j_1,-j_1,j_2,-j_2}= H_{4,j_1,-j_1,j_2,-j_2} \\
& -\frac{9}{4}\times \frac{1}{24}\times 4\underbrace{\left(\mathbb{H}_{2,j_1,j_2,-j_1,-j_2}+\mathbb{H}_{2,j_1,-j_2,-j_1,j_2} +\mathbb{H}_{2,-j_1,j_2,j_1,-j_2} + \mathbb{H}_{2,-j_1,-j_2,j_1,j_2}\right)}_{=:K_{j_1,j_2}}.
\end{align*}
Using the definition of $\mathbb{H}_{2,j_1,j_2,j_3,j_4}$  in \eqref{hhhdeds}, we have that for $j_1,j_2\in S^+$,
\[
K_{j_1,j_2} = \begin{cases}
\frac{2(j_1+j_2)(H_{3,j_1,j_2,-(j_1+j_2)})^2}{(j_1+j_2)\kappa_{j_1+j_2}-j_1\kappa_{j_1}-j_2\kappa_{j_2}} & \text{ if $j_1=j_2$,}\\
\frac{2(j_1+j_2)(H_{3,j_1,j_2,-(j_1+j_2)})^2}{(j_1+j_2)\kappa_{j_1+j_2}-j_1\kappa_{j_1}-j_2\kappa_{j_2}} + \frac{2(j_1-j_2)(H_{3,j_1,-j_2,-(j_1-j_2)})^2}{(j_1-j_2)\kappa_{j_1-j_2} - j_1\kappa_{j_1}+j_2\kappa_{j_2}} & \text{ if $j_1\ne j_2$}.
\end{cases}
\]
Plugging this into the above expression for $H^{(3)}_{4,j_1,-j_1,j_2,-j_2}$, the desired result follows.
\end{proof}

Lastly, we deal with the quintic Hamiltonian.

\begin{lemma}\label{quintic} (Normalization of the quintic Hamiltonian) There exists a symplectic transformation $\Phi^{(4)}$ such that 
\begin{align*}
H^{(5)}(f):=H^{(4)}\circ \Phi^{(5)}(f) = H^{(5)}_2(f) + H^{(5)}_3(f) +H^{(5)}_4(f) + H^{(5)}_5(f) + H^{(5)}_{\ge6}(f)
\end{align*}
where
\begin{enumerate}[label=\text{(\arabic*)}]
\item (Quadratic, cubic and quartic terms) We have
\begin{equation}\label{m_H5222}
\begin{aligned}
H^{(5)}_2(f) & = \mathcal{H}_2(f), \quad H^{(5)}_3(f) = {H}^{(3)}_3(f),\\
\quad H^{(5)}_4(f) & = H^{(4)}_4(f), \quad H^{(5)}_5 = \sum_{i=2}^{5}R(v^{5-i}z^i), 
\end{aligned}
\end{equation} \label{item1_lemsss}
and $H^{(5)}_{\ge 6}$ collects all the terms homogeneous of degree at least $6$.
\item There exists a finite dimensional space $E$ of the form in \eqref{finite_dim} such that $\Phi^{(5)}(f) = f + \Psi^{(4)}(f) $ where $\Psi^{(5)}$ satisfies $\Psi^{(5)}(f) = \Pi_E \Psi^{(5)}(\Pi_E f)$.
\item[(3)] $\Phi^{(5)}$ is real and reversibility preserving.
\end{enumerate}
\end{lemma}
\begin{proof}
\vspace{0.5\baselineskip}
\noindent\textit{Proof of (1).}
Let us write $\mathcal{H}^{(4)}_5$ obtained in Lemma~\ref{quartic} as
\begin{align}\label{1sdcxxcx}
\mathcal{H}^{(4)}_5(f) = \sum_{j_1+j_2+j_3+j_4+j_5 = 0}H^{(4)}_{5,j_1,j_2,j_3,j_4,j_5}f_{j_1}f_{j_2}f_{j_3}f_{j_4}f_{j_5},
\end{align}
where $H^{(4)}_{5,j_1,j_2,j_3,j_4,j_5}$ is assumed to be invariant under any permutations\index{permutations} on $\left\{ j_1,j_2.j_3,j_4,j_5\right\}$ without loss of generality. Indeed, we already prove in \ref{item2_lemsss} of Lemma~\ref{quartic} that $H^{(4)}_{5,j_1,j_2,j_3,j_4,j_5}=0$, if $j_{1}+\cdots,j_5\ne 0$, hence we have the expression in \eqref{1sdcxxcx}.  As before, we consider a Hamiltonian:
\begin{align}\label{F5def}
F^{(5)}(f) := \sum_{j_1+j_2+j_3 + j_4 + j_5 = 0}F^{(5)}_{j_1,j_2,j_3,j_4,j_5}f_{j_1}f_{j_2}f_{j_3}f_{j_4}f_{j_5},
\end{align}
which will be determined later and let us denote its time-1 flow map by $\Phi^{(5)}$.  We compute using Lemma~\ref{symplectic_1}, 
\begin{align}\nonumber
{H}^{(5)} := {H}^{(4)}\circ \Phi^{F^{(5)}}_{1} = {H}^{(4)}_2 + {H}^{(4)}_3 + {H}^{(4)}_4 +{H}^{(5)}_5+ {H}^{(5)}_{\ge 6}, \\
\quad \text{ where }\quad  {H}^{(5)}_5 + \left\{{H}^{(4)}_2,F^{(5)} \right\}.\label{H55expansion}
\end{align}
Again, it follows from \eqref{m_H4222} that $H^{(4)}_2=\mathcal{H}_2$, thus using \eqref{poisson_bracket}, we have
\begin{align*}
 \left\{ \mathcal{H}_2,F^{(5)}\right\}&=-\frac{2\ii}{2\pi}\sum_{j_1+j_2+j_3+j_4+j_5=0}\left( j_1\kappa_{j_1}+j_2\kappa_{j_2}+j_3\kappa_{j_3}+j_4\kappa_{j_4} + j_5\kappa_{j_5}\right)\\
 & \times F^{(5)}_{j_1,j_2,j_3,j_4,j_5}f_{j_1}f_{j_2}f_{j_3}f_{j_4}f_{j_5}.
\end{align*}
With  \ref{item2_lemsss} in Lemma~\ref{quartic}, which implies that $H^{(4)}_5$ satisfies \eqref{conserve_momentum}, we have
\begin{align}
{H}^{(5)}_5  &=\sum_{j_1+j_2+j_3+j_4+j_5 = 0} \left( H^{(4)}_{5,j_1,j_2,j_3,j_4,j_5}\right. \nonumber \\& \qquad \left. -\frac{ \ii}{\pi} \left(j_1\kappa_{j_1}+j_2\kappa_{j_2}+j_3\kappa_{j_3}+j_4\kappa_{j_4} + j_5\kappa_{j_5}\right) F^{(5)}_{j_1,j_2,j_3,j_4,j_5}\right)f_{j_1}f_{j_2} f_{j_3}f_{j_4}f_{j_5}\nonumber\\
& = \sum_{\substack{j_1+j_2+j_3+j_4 + j_5 = 0 \\ (j_1,j_2,j_3,j_4 ,j_5)\in \mathcal{A}_5}} \left( H^{(4)}_{5,j_1,j_2,j_3,j_4,j_5}\right. \nonumber \\& \qquad\left.- \frac{\ii}{\pi} \left(j_1\kappa_{j_1}+j_2\kappa_{j_2}+j_3\kappa_{j_3}+j_4\kappa_{j_4} + j_5\kappa_{j_5}\right)F^{(5)}_{j_1,j_2,j_3,j_4,j_5}\right)f_{j_1}f_{j_2} f_{j_3}f_{j_4}f_{j_5}\nonumber\\
& \ +  \sum_{\substack{j_1+j_2+j_3+j_4 + j_5 = 0 \\ (j_1,j_2,j_3,j_4 ,j_5)\notin \mathcal{A}_5}} \left( H^{(4)}_{5,j_1,j_2,j_3,j_4,j_5}\right. \nonumber \\& \qquad\left.- \frac{\ii}{\pi} \left(j_1\kappa_{j_1}+j_2\kappa_{j_2}+j_3\kappa_{j_3}+j_4\kappa_{j_4} + j_5\kappa_{j_5}\right)F^{(5)}_{j_1,j_2,j_3,j_4,j_5}\right)f_{j_1}f_{j_2} f_{j_3}f_{j_4}f_{j_5},\label{sssssperp6}
\end{align}
where \begin{align*}
\mathcal{A}_5 := & \left\{(j_1,j_2,j_3,j_4,j_5)\in \mathbb{Z}^5 : \sum_{i=1}^5 j_{i}\kappa_{j_i}\ne 0\right. \\  &
\text{  and at least four of $\left\{j_1,j_2,j_3,j_4,j_5 \right\}$  are in $S$ }\Big\}.\end{align*}
Hence we can define $F^{(5)}$ to be
\begin{align}\label{def_F5}
F^{(5)}_{j_1,j_2,j_3,j_4,j_5} = 
\begin{cases}
\frac{\pi H^{(4)}_{5,j_1,j_2,j_3,j_4,j_5}}{\ii \sum_{i=1}^5 j_{i}\kappa_{j_i}}  & \text{ if $(j_1,j_2,j_3,j_4,j_5) \in \mathcal{A}_5$ },\\
0 & \text{ otherwise,}
\end{cases}
\end{align}
so that \eqref{sssssperp6} yields that
\[
{H}^{(5)}_5 = \sum_{\substack{j_1+j_2+j_3+j_4+j_5=0 \\ (j_1,j_2,j_3,j_4,j_5)\notin \mathcal{A}_5}}H^{(4)}_{5,j_1,j_2,j_3,j_4,j_5}f_{j_1}f_{j_2}f_{j_3}f_{j_4}f_{j_5}.
\]
Indeed, \ref{tangent_1} tells us that \eqref{def_F5} is well-defined.

Now, let us write ${H}^{(5)}_5 = \sum_{i=0}^5 {H}^{(5)}_{5,i}$, where ${H}^{(5)}_{5,i}$ is of the form $R(v^iz^{5-i})$. By the definition of $\mathcal{A}_5$, it holds that 
\begin{align}\label{tangential_site_condition4}
{H}^{(5)}_{5,0} = {H}^{(5)}_{5,1} = 0.
\end{align}
This gives \eqref{m_H5222}.

\vspace{0.5\baselineskip}
\noindent\textit{Proof of (2)} $\mathcal{A}_5$ is a finite set since it is contained in $[-4K,4K]^5$, where $K = \max_{i=1,\ldots \nu, \quad j_i \in S} |j_{i}|$

\vspace{0.5\baselineskip}
\noindent\textit{Proof of (3).} The proof is identical to (4) in Lemma~\ref{ssslem}.
\end{proof}

Now we are ready to prove Proposition~\ref{normal_form_prop11231}.
\begin{proofprop}{normal_form_prop11231}
Let 
\begin{align}\label{phiwbdefsx}
\Phi^{WB} := \Phi^{(3)}\circ \Phi^{(4)} \circ \Phi^{(5)},
\end{align} where $\Phi^{(3)},\Phi^{(4)},\Phi^{(5)}$ are the symplectic transformations obtained in Lemmas~\ref{ssslem}, \ref{quartic} and \ref{quintic}. Then it follows from the lemmas that
\[
H:=\mathcal{H}\circ \Phi^{WB} = H^{(5)}_2(f) + H^{(5)}_3(f) +H^{(5)}_4(f) + H^{(5)}_5(f) + H^{(5)}_{\ge6}(f).
\] 

Then, \eqref{normal_formsf} follows from \eqref{m_H5222} and \eqref{m_H3}. For the quartic and quintic terms, \eqref{h4h5} follows from \eqref{quartic_1112}, \eqref{quartic_11122} and Lemma~\ref{ammsmw}.  Also, \ref{item3_normal} follows from Lemma \ref{ssslem}, \ref{quartic} and \ref{quintic} as well. $\Phi^{WB}$ is real and reversibility preserving because it is a composition of real, reversibility preserving transformations. Furthermore, from \eqref{F3def}, \eqref{F4def} and \eqref{F5def}, it follows that $\Phi^{WB}$ is a composition of $\frac{2\pi}{\mathtt{M}}$-translation invariant preserving transformations (see Remark~\ref{krsds2221sd} and Lemma~\ref{minvariantce}), there for $\Phi^{WB}$ is $\frac{2\pi}{\mathtt{M}}$-translation invariant preserving as well.  \ref{item3_normal_3} follows immediately from Lemma~\ref{muregular_kep} since $F^{(3)},F^{(4)},F^{(5)}$ are supported on a finite number of Fourier modes.
\end{proofprop}

We specify $H_{3,2},H_{4,2}$ in the following lemma:
\begin{lemma}\label{h24sdx2s}
$H_{3,2}$ in \eqref{normal_formsf} and $H_{4,2}= \mathcal{H}_{4,2} + \mathfrak{H}_1 + \mathfrak{H}_2(f)$ in \eqref{sjdjsjs93jx} has the following formula in terms of the Fourier modes:
\begin{align}
H_{3,2}(f) &=\sum_{\substack{j_1+j_2+j_3 = 0,\\ \text{exactly two of $_{j_1,j_2,j_3}$ are in $S^{\perp}$} }} H_{3,j_1,j_2,j_3}f_{j_1}f_{j_2}f_{j_3}, \label{rlaksd2sdsd}\\
\mathcal{H}_{4,2}(f)&=\sum_{\substack{j_1+j_2+j_3+j_4=0,\\\text{exactly two of $j_1,j_2,j_3,j_4$ are in $S^\perp$}}}H_{4,j_1,j_2,j_3,j_4}f_{j_1}f_{j_2}f_{j_3}f_{j_4},\label{rlaksd2sdsd1}\\
\mathfrak{H}_1(f)& = -9\sum_{\substack{j_1+j_2+j_3+j_4=0,\\ j_1,j_2\in S,\ j_3,j_4\in S^\perp, \\ j_2+j_3\in S}}(j_2+j_3)\frac{H_{3,-(j_2+j_3),j_2,j_3}H_{3,j_1,j_4,-(j_1+j_4)}}{j_1\kappa_{j_1}+j_4\kappa_{j_4}-(j_1+j_4)\kappa_{j_1+j_4}}f_{j_1}f_{j_2}f_{j_3}f_{j_4},\label{rlaksd2sdsd3}\\
\mathfrak{H}_2(f)&=-\frac{9}{2}\sum_{\substack{j_1+j_2+j_3+j_4=0,\\ j_1,j_2\in S,\ j_3,j_4\in S^\perp,\\ j_1+j_2\in S^\perp}}(j_3+j_4)\frac{H_{3,j_1,j_2,-(j_1+j_2)}H_{3,j_3,j_4,-(j_3+j_4)}}{j_1\kappa_{j_1}+j_2\kappa_{j_2}-(j_1+j_2)\kappa_{j_1+j_2}}f_{j_1}f_{j_2}f_{j_3}f_{j_4} .\label{rlaksd2sdsd4}
\end{align}
\end{lemma}
\begin{proof}
For $H_{3,2}$, the expansion follows from  \eqref{normal_formsf}. From Lemma~\ref{quartic} and Lemma~\ref{quintic}, it follows that $H_{4,2} = H^{(3)}_{4,2}$. Therefore, the result for $H_{4,2}$ follows from \eqref{sswsw2s} in Lemma~\ref{fourier_H3333s2s}. 
\end{proof}

Before we close this chapter, we prove the expansion of $\Phi^{WB}$. 
\begin{lemma}\label{WB_expansion11}
For $f\in L^2_0$, we have the expansion,
\[\Phi^{WB}(f) = f + 6\partial_{x}\Pi_S K_1(v,z) + 3\partial_x\Pi_{S^\perp}K_2(v,v) + \Phi^{WB}_{\ge 3}(f),\]
where  $K_1,K_2$ are as in \eqref{K12def} and $\Phi^{WB}_{\ge 3}:L^2_0\mapsto L^2_0$ collects all the homogeneous terms of degree at least $3$.
\end{lemma}

\begin{proof}
Recall the definition of $\Phi^{WB}$ \eqref{phiwbdefsx} and that each $\Phi^{(i)},\ i=3,4,5$ are time-1 flow map of Hamiltonians $F^{(i)}$. Since $F^{(4)},F^{(5)}$ are quartic and quintic Hamiltonian respectively, the linear and the quadratic terms of $\Phi^{WB}$ coincide with those of $\Phi^{(3)}$. Denoting $\Phi^{(3)}_t$ be the Hamiltonian flow of $F^{(3)}$ in the proof of Lemma~\ref{ssslem}, we have
\[
\frac{d}{dt}\Phi^{(3)}_t = \partial_x\nabla F^{(3)}(\Phi_t^{(3)}(f)),
\]
hence, the Taylor expansion as in \eqref{elementary_taylor} of $t\mapsto \Phi_t^{(3)}$ gives us that
\begin{align*}
\Phi^{(3)}(f) = \Phi_1^{(3)}(f) & = f + \partial_x\nabla F^{(3)}(f) \\
& + \int_0^{1}\partial_x\nabla^2F^{(3)}(\Phi_t^{(3)}(f))[\partial_x\nabla F^{(3)}(\Phi_t^{(3)}(f))]dt.
\end{align*}
Since $F^{(3)}$ is cubic, the last term in integral is homogeneous of degree $\ge 3$.
Recalling \eqref{gradf2s2s} for $\nabla F^{(3)}(f)$, the desired result follows.
\end{proof}

\chapter{Action-angle variables}\label{aavariab}
In this chapter, we introduce action-angle variables\index{Action-angle variables} and rewrite the Hamiltonian $H$ that was obtained in Proposition~\ref{normal_form_prop11231} in terms of those variables. First, we denote the amplitude vector\index{amplitude vector} by $\zeta$:
 \begin{align}\label{amskoisdwdw21s}
 \zeta\in [1,2]^{\nu}\subset (\mathbb{R}^+)^\nu,
 \end{align} 
 For each $\zeta$, we define a change of variables, $U_{\zeta}:\mathbb{T}^\nu \times \mathbb{R}^\nu\times H_{S^\perp}\mapsto L^2_0$, as follows
\begin{align}\label{actionssx}
&U_\zeta(\theta,y,z):=\epsilon{\sum_{j_k\in S}\sqrt{|j_k|(\zeta_k + \epsilon^{2(b-1)}y_k)}e^{\ii \theta_k}}e^{\ii j_k x} + \epsilon^b z{=:\epsilon v_{\epsilon}(\theta,y)} + \epsilon^bz,\\
&\zeta_{-k} :=\zeta_k,\quad y_{-k} = y_k,\quad \theta_{-k}=-\theta_k,\quad k=1,\ldots,\nu,\quad \epsilon>0\quad b-1\in (0,\frac{1}{12}).\nonumber
\end{align}
Note that $v_\epsilon$ depends on $\zeta$ as well, but we omit it in the notations for simplicity. Clearly,  we have $\overline{v}_j=v_{-j}$ hence $v_\epsilon$ is real-valued and its Fourier modes\index{Fourier mode} are supported on the tangential sites. For sufficiently small $\epsilon>0$, we can easily see that $\zeta$ determines the size of the amplitude of the function determined by $(\theta,y,z)$. Using the new variables $(\theta,y,z)$, we define a new $\zeta$-dependent Hamiltonian\index{$\zeta$-dependent Hamiltonian} in $\mathbb{T}^\nu\times \mathbb{R}^\nu \times H_{S^\perp}$:
\begin{align}\label{rescaled_hamiltonian}
H_{\zeta}(\theta,y,z):=\frac{1}{\epsilon^{2b}}H(U_\zeta(\theta,y,z))=\frac{1}{\epsilon^{2b}}\mathcal{H}\circ \Phi^{WB}\circ U_\zeta,
\end{align}where $H$ is as in Proposition~\ref{normal_form_prop11231}.

\begin{proposition}\label{aavariable}
Let $H$ be the Hamiltonian obtained in Proposition~\ref{normal_form_prop11231} and let $v_\epsilon,U_\zeta$ and $H_{\zeta}$  be as in  \eqref{actionssx} and \eqref{rescaled_hamiltonian}. Then we have that for some constant $C_\epsilon(\zeta)$, 
\begin{align}\label{normal_formaa2}
H_{\zeta}(\theta,y,z) = C_\epsilon(\zeta) + 2\pi{\omega(\zeta)}\cdot y + \mathcal{N}(\theta)(z,z) + P(\theta,y,z),
\end{align} 
where
\begin{align}
&\omega(\zeta) := \left( \overline{\omega} + \frac{6\epsilon^2}{\pi}\mathbb{A}\zeta  \right), \text{ where $\overline{\omega}\in\mathbb{R}^\nu$ and $\mathbb{A}\in \mathbb{R}^{\nu\times\nu}$ are in \eqref{linear_frequency_aa} and \eqref{amplitude_modulation1}},\label{normal_formaa3}\\
& \mathcal{N}(\theta)(z,z) :=H_{2,2}(z) +\epsilon H_{3,2}(v_\epsilon(\theta,0)+z) + \epsilon^2H_{4,2}(v_\epsilon(\theta,0)+z),\label{normal_formaa4}\\
&P(\theta,y,z) := \epsilon\left( H_{3,2}( v_\epsilon(\theta,y)+z) - H_{3,2}(v_\epsilon(\theta,0)+z)\right) \nonumber\\
& \quad \quad \quad \quad+ \epsilon^2 \left( H_{4,2}(v_\epsilon(\theta,y)+z) - H_{4,2}(v_\epsilon(\theta,0) + z)\right)\nonumber\\
& \quad \quad \quad \quad \  +6\epsilon^{2b}y^T \mathbb{A}y + \epsilon^{-2b}\left( H_{3,3} +  \sum_{i=3}^4H_{4,i}+\sum_{i=2}^{5}{H_{5,i}} + H_{\ge6}\right)\circ U_\zeta.\label{normal_formaa5}
\end{align}
\end{proposition}

\begin{proof}
We recall $H$ in Proposition~\ref{normal_form_prop11231} and see that
\begin{align}\label{H_decomposition}
H =\left(H_{2,0} + H_{4,0} \right) + \left(H_{2,2} + H_{3,2}+H_{4,2} \right) + \left(H_{3,3} + \sum_{i=3}^4H_{4,i}+\sum_{i=2}^{5}{H_{5,i}} + H_{\ge6} \right).
\end{align}
Let us compute the terms in each parenthesis separately:

\textbf{Computation for $H_{2,0} + H_{4,0}$.}
Using $H_{2,0}$ in \eqref{normal_formsf}, we compuate
\begin{align*}
H_{2,0}(U_{\zeta}(\theta,y,z))& \overset{\eqref{expansion_2},\eqref{normal_formsf}}=\mathcal{H}_2(\epsilon v_\epsilon(\theta,y)) \\
& \overset{\eqref{fourier_hamiltonian}}= 2\epsilon^2\sum_{j_k\in S^+}j_k\kappa_{j_k}\left(\zeta_k + \epsilon^{2(b-1)}y_k \right)\\
& = \epsilon^2C_1(\zeta) + 2\epsilon^{2b}\sum_{j_{k}\in S^+}j_k\kappa_{j_k}y_k \\
&=:\epsilon^2C_1(\zeta) +2\pi \epsilon^{2b}\overline{\omega}\cdot y,
\end{align*}
where $C_1(\zeta)$ is a constant that does not depend on $(\theta,y,z)$ and  $\overline{\omega}\in \mathbb{R}^\nu$\index{$\overline{\omega}$} is defined as 
\begin{align}\label{linear_frequency_aa}
(\overline{\omega})_k :=\frac{1}{\pi} \kappa_{j_k}\overset{\eqref{moreexplicits2}}= j_k\left(-\frac{1}{2}m_{1,\alpha}(j_k) + \frac{T_\alpha}{4} \right) \quad \text{ for $j_1,\ldots,j_\nu \in S^{+}$}.
\end{align}
Using $H_{4,0}$ in  \eqref{h4h5}, we have
\begin{align*}
& H_{4,0}(U_\zeta)(\theta,y,z) \\
& =6\sum_{j_k\in S^+}H^{(3)}_{4,j_k,-j_k,j_k,-j_k}\epsilon^4j_k^2(\zeta_k+\epsilon^{2(b-1)}y_k)^2 \\
& \ + 12\sum_{\substack{j_{k_1},j_{k_2}\in S^+,\\ j_{k_1}\ne j_{k_2}}}H^{(3)}_{4,j_{k_1},-j_{k_1},j_{k_2},-j_{k_2}}\epsilon^4j_{k_1}j_{k_2}(\zeta_{k_1}+\epsilon^{2(b-1)}y_{k_1})(\zeta_{k_2}+\epsilon^{2(b-1)}y_{k_2})\\
& = C_2(\zeta) + 12\epsilon^{2b+2}\left( \sum_{j_k\in S^{+}}H_{4,j_k,-j_k,j_k,-j_k}^{(3)}j_k^2\zeta_ky_k \right. \\
& \left. \qquad \qquad \qquad \qquad + 2\sum_{\substack{j_{k_1},j_{k_2}\in S^+,\\ j_{k_1}\ne j_{k_2}}}H^{(3)}_{j_{k_1},-j_{k_1},j_{k_2},-j_{k_2}}j_{k_1}j_{k_2}\zeta_{k_1}y_{k_2} \right)\\
& \ +6\epsilon^{4b}\left( \sum_{j_k\in S^{+}}H_{4,j_k,-j_k,j_k,-j_k}^{(3)}j_k^2y_k^2 + 2\sum_{\substack{j_{k_1},j_{k_2}\in S^+,\\ j_{k_1}\ne j_{k_2}}}H^{(3)}_{j_{k_1},-j_{k_1},j_{k_2},-j_{k_2}}j_{k_1}j_{k_2}y_{k_1}y_{k_2} \right)\\
& = C_2(\zeta) + 12\epsilon^{2b+2}\mathbb{A}\zeta\cdot y + 6\epsilon^{4b}\mathbb{A}y\cdot y,
 \end{align*}
 where  $C_2(\zeta)$ is a constant that does not depend on $(\theta,y,z)$ and \index{$\mathbb{A}$}$\mathbb{A}\in \mathbb{R}^{\nu\times \nu}$ is given by
 \begin{align}\label{amplitude_modulation1}
 \quad \mathbb{A}^{k_1}_{k_2}:=\begin{cases}
j_{k_1}^2H^{(3)}_{4,j_{k_1},-j_{k_1},j_{k_2},-j_{k_2}} &\text{ if $k_1=k_2$,}\\
2j_{k_1}j_{k_2}H^{(3)}_{4,j_{k_1},-j_{k_1},j_{k_2},-j_{k_2}} &\text{ if $k_1\ne k_2$},
 \end{cases}  \end{align}
  where $k_1,k_2=1,\ldots,\nu$, and $j_{k_1},j_{k_2}\in S^+$.
Therefore, we have that for some constant $C_\epsilon(\zeta)$,
\begin{align}\label{zero_normal}
H_{2,0}(U_\zeta)+H_{4,0}(U_\zeta) = C_\epsilon(\zeta) + 2\pi\epsilon^{2b}\left(\overline{\omega}+ \frac{6}{\pi}\epsilon^2\mathbb{A}\zeta \right)\cdot y + 6\epsilon^{4b} \mathbb{A}y\cdot y.
\end{align}

\textbf{Computation for $H_{2,2}+H_{3,2}+H_{4,2}$.}
For $H_{2,2}$, we have
\begin{align}\label{2233}
H_{2,2}(U_\zeta) =\epsilon^{2b} H_{2,2}(z) .
\end{align}

For $H_{3,2} + H_{4,2}$, recall from \eqref{m_H3} and \eqref{quartic_11122} that they are of the form $R(vz^2)$ and $R(v^2z^2)$ respectively, therefore we have (recalling $U_\zeta(\theta,y,z) = \epsilon v_\epsilon(\theta,y) + \epsilon^bz$ as in \eqref{actionssx}),
\begin{align*}
H_{3,2}(U_\zeta) + H_{4,2}(U_\zeta) &= H_{3,2}(\epsilon v_\epsilon(\theta,y) + \epsilon^b z) + H_{4,2}(\epsilon v_\epsilon(\theta,y) + \epsilon^bz)\\
& = \epsilon^{1+2b}H_{3,2}(v_\epsilon(\theta,y)+z) + \epsilon^{2+2b}H_{4,2}(v_\epsilon(\theta,y) + z)\\
& = \epsilon^{1+2b}H_{3,2}(v_\epsilon(\theta,0) + z) + \epsilon^{2+2b}H_{4,2}(v_\epsilon(\theta,0) + z) \\
& + \epsilon^{1+2b}P_1(\theta,y,z),
\end{align*}
where
\begin{align}\nonumber 
P_1(\theta,y,z) & = \left( H_{3,2}( v_\epsilon(\theta,y)+z) - H_{3,2}(v_\epsilon(\theta,0)+z)\right) \\
& + \epsilon \left( H_{4,2}(v_\epsilon(\theta,y)+z) - H_{4,2}(v_\epsilon(\theta,0) + z)\right).\label{p1sdsdright}
\end{align}
 With \eqref{2233}, we denote  by $\mathcal{N}(\theta)(z,z)$ the bilinear form in $z$, that is,
 \begin{equation}\label{def_n1andp1}
 \begin{aligned}
 &\mathcal{N}(\theta)(z,z):= H_{2,2}(z) +\epsilon H_{3,2}(v_\epsilon(\theta,0)+z) + \epsilon^2H_{4,2}(v_\epsilon(\theta,0)+z),
  \end{aligned}
  \end{equation}
then, we obtain
\begin{align}\label{22333242}
H_{2,2}(U_\zeta)+H_{3,2}(U_\zeta)+H_{4,2}(U_\zeta) = \epsilon^{2b} \mathcal{N}(\theta)(z,z) + \epsilon^{1+2b} P_1(\theta,y,z).
\end{align}

\textbf{Computation for $H_{3,3} + \sum_{i=3}^4H_{4,i}+\sum_{i=2}^{5}{H_{5,i}} + H_{\ge6}$.} Since all of these terms are small enough, we simply denote
\begin{align}\label{p2_def}
P_2(\theta,y,z) := \epsilon^{-2b}\left( H_{3,3} +  \sum_{i=3}^4H_{4,i}+\sum_{i=2}^{5}{H_{5,i}} + H_{\ge6}\right),
\end{align}
so that 
\begin{align}\label{perturbative_term}
H_{3,3}(U_\zeta) + \sum_{i=3}^4H_{4,i}(U_\zeta)+\sum_{i=2}^{5}{H_{5,i}}(U_\zeta) + H_{\ge6}(U_\zeta) =: \epsilon^{2b} P_2(\theta,y,z).
\end{align}
Therefore, plugging \eqref{zero_normal}, \eqref{22333242} and \eqref{perturbative_term} into \eqref{H_decomposition} and \eqref{rescaled_hamiltonian}, we obtain
\begin{align}\nonumber
H_{\zeta}(\theta,y,z) & =\frac{C_\epsilon(\zeta)}{\epsilon^{-2b}} + 2\pi\left(\overline{\omega} + \frac{6}{\pi}\epsilon^2\mathbb{A}\zeta \right)\cdot y +6\epsilon^{2b} \mathbb{A}y\cdot y  \\
& + \mathcal{N}(\theta)(z,z) + \epsilon P_1(\theta,y,z) + P_2(\theta,y,z).\label{normal_form_aa}
\end{align}
This is equivalent to  \eqref{normal_formaa2}.
\end{proof}

Now, we consider the corresponding symplectic $2$-form $\mathcal{W}$ on $\mathbb{T}^\nu\times \mathbb{R}^\nu \times H_{S^\perp}$, which is obtained as the pullback of $\sigma$ in \eqref{symplectic} by the transformations that we have performed. From Proposition~\ref{normal_form_prop11231} ($\Phi^{WB}$ is symplectic), and \eqref{actionssx}, we have
\begin{align}\label{symplectic_form1}
\widehat{\mathcal{W}} := (\Phi^{WB}\circ U_\zeta)^{*}\sigma = U_\zeta^*\sigma = \epsilon^{2b}\left( \frac{1}{2\pi}\sum_{k=1}^{\nu}d\theta_k \wedge dy_k \right)\oplus  \epsilon^{2b}\sigma_{S^\perp}=:\epsilon^{2b}{\mathcal{W}},
\end{align}
where $\sigma_{S^\perp}$ is the restriction of $\sigma $ to $H_{S^\perp}$, that is,
$\sigma_{S^{\perp}}(z,z) = \int \left( \partial_x^{-1}z\right)zdx$. From \eqref{rescaled_hamiltonian}, it follows that the gSQG dynamics are equivalent to
\begin{align}\label{dynamics1}
\colvec{\dot{\theta}\\ \dot{y} \\ \dot{z}} = \widehat{X}_{H\circ U_\zeta}(\theta,y,z) = \epsilon^{2b}\widehat{X}_{H_\zeta}(\theta,y,z), 
\end{align}
where $\widehat{X}_{H\circ U_\zeta},\widehat{X}_{\zeta}$ are the Hamiltonian vector fields induced by the symplectic form $\widehat{\mathcal{W}}$.   From \eqref{rescaled_hamiltonian}, we have that for any vector field $h$ on $\mathbb{T}^\nu\times \mathbb{R}^\nu\times H_{S^\perp}$,
\begin{align*}
\nabla_{\theta,y,z}H_\zeta[h] &= \epsilon^{-2b}(\nabla_{\theta,y,z}(H\circ U_\zeta),h) =\epsilon^{-2b}\widehat{\mathcal{W}}(\widehat{X}_{H\circ U_\zeta},h) = \widehat{\mathcal{W}}(\widehat{X}_{H_\zeta},h) \\
&=\epsilon^{2b}\mathcal{W}(\widehat{X}_{H_\zeta},h),
\end{align*}
where the last equality follows from \eqref{symplectic_form1}.
This implies (using the definition of $\mathcal{W}$ in \eqref{symplectic_form1}),
\[
\epsilon^{2b}\widehat{X}_{H_\zeta}(\theta,y,z) = \colvec{\frac{1}{2\pi}\partial_y H_\zeta(\theta,y,z) \\ -\frac{1}{2\pi}\partial_\theta H_\zeta(\theta,y,z) \\ \partial_x\nabla_z H_\zeta(\theta,y,z)}.
\]
Therefore, \eqref{dynamics1} is equivalent to 
\begin{align}\label{dynamics2}
\colvec{\dot{\theta}\\ \dot{y} \\ \dot{z}} = \colvec{\frac{1}{2\pi}\partial_y H_\zeta(\theta,y,z) \\ -\frac{1}{2\pi}\partial_\theta H_\zeta(\theta,y,z) \\ \partial_x\nabla_z H_\zeta(\theta,y,z)} =: X_{H_\zeta}.
\end{align}
 Note that $X_{H_\zeta}$ is the Hamiltonian vector field of $H_\zeta$ induced by $\mathcal{W}$. We summarize these considerations in the following proposition:
 \begin{proposition}\label{newhamaiam}
Let $H_\zeta$\index{$H_\zeta$} be a $\zeta$-dependent Hamiltonians given in Proposition~\ref{aavariable}. Then a quasiperiodic solution $t\mapsto (\theta(t),y(t),z(t))$ to \eqref{dynamics2} determines a quasiperiodic solution for the gSQG, that is, $
f(t):=\Phi^{WB}\circ U_\zeta(\theta(t),y(t),z(t))$ solves \eqref{gSQG_hamiltonian}.
 \end{proposition}
 
 Our new Hamiltonian $H_\zeta$ on $\mathbb{T}^\nu\times \mathbb{R}^\nu\times H_{S^\perp}$ satisfies the time-reversibility and the $\frac{2\pi}{\mathtt{M}}$-translation invariance:

 \begin{lemma}\label{rsksdtimerevsd}
  $H_\zeta$ is  time-reversible with respect to $\rho_*$, defined in \eqref{definiteionofinv2}.
 \end{lemma}
\begin{proof}
It is easy to see from \eqref{actionssx} that $
\rho(f) = \rho_* (U_\zeta)$, where $\rho$ is defined in \eqref{definiteionofinv}. Since $H$ is reversible (see \ref{item3_normal_2} in Proposition~\ref{normal_form_prop11231}) with respect to $\rho$, we have
\[
H_\zeta\circ \rho_* = \epsilon^{-2b}H(\rho_*(U_\zeta)) = \epsilon^{-2b}H(\rho(f))=\epsilon^{-2b}H(f)=\epsilon^{-2b}H(U_\zeta)=H_\zeta,
\]
which implies that $H_\zeta$ is a time-reversible Hamiltonian with respect to the involution $\rho_*$.
\end{proof}

\begin{lemma}\label{msdjsd11111sd}
$H_\zeta$ is $\frac{2\pi}{\mathtt{M}}$-translation invariant.
\end{lemma}
\begin{proof}
Recalling the property of $S$ in \eqref{tan_site}, we see that 
\begin{align}\label{uusd11}
U_\zeta(\theta,y,\rho_{\mathtt{M}}z) = \rho_{\mathtt{M}}\circ U_\zeta(\theta,y,z).
\end{align} Since $H$ is $\frac{2\pi}{\mathtt{M}}$-translation invariant, which follows from \ref{item3_normal_2} of Proposition~\ref{normal_form_prop11231},  the definition of $H_\zeta$ in \eqref{rescaled_hamiltonian} tells us that
\begin{align}\label{rmsdsdsdsdx}
H_{\zeta}(\theta,y,\rho_{\mathtt{M}}z) = \frac{1}{\epsilon^{2b}}H(U_\zeta(\theta,y,\rho_{\mathtt{M}}z)) = \frac{1}{\epsilon^{2b}}H(\rho_{\mathtt{M}}\circ U_\zeta(\theta,y,z)) = H_{\zeta}(\theta,y,z).
\end{align}
This tells us that $H_\zeta$ is $\frac{2\pi}{\mathtt{M}}$-translation invariant (see Definition~\ref{dekfhsspdiw11}).
\end{proof}
\color{black}
\section{Hypotheses on the tangential sites}\label{rpoisdsd1sd}
In this section, we will specify precise conditions on the choice of our tangential sites $S$ in \eqref{tan_site}. As seen in Chapter~\ref{skpssisodssuw}, the conditions, \ref{tangent_2} and \ref{tangent_1} are imposed to derive the weak Birkhoff normal from in Proposition~\ref{normal_form_prop11231}. In order to estimate the measure of frequencies, which will be explained below, we require more conditions on the choice of $S$. Those conditions are mainly due to two different reasons, 1)  Dependence of the Hamiltonian on a parameter to perform the KAM theory and 2) Measure estimate for the set of non-resonant frequencies. We explain more details separate in what follows:

\subsection{Parameter dependent Hamiltonian}
A crucial consequence of Proposition~\ref{newhamaiam} (and Proposition~\ref{aavariable}, \ref{normal_form_prop11231}) is that if we choose the tangential sites $S$ satisfying \ref{tangent_2} and \ref{tangent_1}, we obtain a $\zeta$-dependent Hamiltonian $H_\zeta$ on $\mathbb{T}^\nu\times \mathbb{R}^\nu\times H_{S^\perp}$, which gives us  equivalent dynamics of the gSQG equations. Compared to the Hamiltonian $\mathcal{H}$ in \eqref{Hamiltonian}, this parameter-dependent Hamiltonians is necessary for the proof of the persistence of quasiperiodic solutions, because the KAM theory does not tell us exactly which frequency can produce a persistent quasiperiodic solutions, but rather tell us the measure of the set of frequencies which produce persistent quasiperiodic solutions. More precisely, in view of \eqref{normal_formaa2} and \eqref{dynamics2}, we are led to find a solution $(\theta(t),y(t),z(t))$ satisfying
\begin{align}\label{shoesxbue}
\colvec{\dot{\theta}\\ \dot{y} \\ \dot{z}} =X_{H_\zeta} =   \colvec{\omega(\zeta) \\ 0 \\ 0  } + \colvec{0 \\ -\partial_\theta\mathcal{N}(\theta)(z,z)\\ \partial_x\nabla_z(\mathcal{N}(\theta)(z,z))} + X_{P}(\theta,y,z).
\end{align}
Neglecting the perturbed term $X_P$, one can easily see that $(\theta(t),y(t),z(t)):=(\omega(\zeta)t,0,0)$ solves \eqref{shoesxbue}, in other words, the linearized equation for \eqref{shoesxbue} possesses a quasiperiodic solution with frequency $\omega(\zeta)$ for each $\zeta\in [1,2]^\nu$ (see \eqref{amskoisdwdw21s}). Therefore the aim of the rest of the paper is to find out the set of $\zeta$ for which such quasiperiodic solutions can survive under the presence of the nonlinear perturbation $X_P$, which has been already made small enough in the weak Birkhoff normal form procedure (Proposition~\ref{normal_form_prop11231}).

 Following the convention in the literature, we will use \index{$\omega(\zeta)$}$\omega(\zeta)\in \mathbb{R}^\nu$ as a parameter instead of $\zeta$. That is, in view of \eqref{normal_formaa3} and \eqref{amskoisdwdw21s}, we define $\Omega_\epsilon$\index{$\Omega_\epsilon$}, the set of frequencies, by
  \begin{align}\label{frequency_set1}
\Omega_\epsilon:=\left\{\omega \in \mathbb{R}^\nu : \omega = \omega(\zeta)=\frac{6\epsilon^2}{\pi}\mathbb{A}\zeta + \overline{\omega} \ \text{ for some }\  \zeta\in [1,2]^{\nu} \right\}.
\end{align}
We will investigate for which $\omega\in \Omega_\epsilon$, the following $\omega$-dependent equation has a quasiperiodic solution:
\begin{align}\label{shoesxbue2}
\colvec{\dot{\theta}\\ \dot{y} \\ \dot{z}}  - \colvec{\omega \\ 0 \\ 0  } - \colvec{0 \\ -\partial_\theta\mathcal{N}(\theta)(z,z)\\ \partial_x\nabla_z(\mathcal{N}(\theta)(z,z))} - X_{P}(\theta,y,z) = 0.
\end{align}
We note that $\mathcal{N}$ and $P$ in \eqref{shoesxbue} depend on $\zeta$ (this follows from their definitions in \eqref{normal_formaa4} and \eqref{normal_formaa5}, while $v_\epsilon$ depends on $\zeta$ as in \eqref{actionssx}). Therefore, $\mathcal{N}$ and $P$ in \eqref{shoesxbue2} must depend on $\omega$ (however, we omit their dependence in the notation for simplicity).

 A necessary requirement for the use of $\omega$ as a parameter as described above is a one-to-one correspondence between $\omega$ and $\zeta$, that is, the invertibility of the matrix $\mathbb{A}$ in \eqref{frequency_set1}. Recalling the definition of $\mathbb{A}$ from \eqref{amplitude_modulation1}, we see that the invertibility of $\mathbb{A}$ completely depends on the choice of $S^+$ in \eqref{tan_site}. Thus, we are led to make the following hypothesis on $S^+$:
 \begin{enumerate}[label=$(\mathtt{H}1)$]
 \item \label{hypos1s1}$\mathbb{A}$ defined in \eqref{amplitude_modulation1} is invertible, thus, it satisfies
 \begin{align}\label{rpishsnsdiwsxcsd}
 |\vec{w}| \le_{\alpha,\nu,S} |\mathbb{A}\vec{w} |\le_{\alpha,\nu,S}|\vec{w}|,\text{ for all $\vec{w}\in \mathbb{R}^\nu$}.
 \end{align}
 \end{enumerate}

 \subsection{Non-resonant frequencies} Other\index{Non-resonant frequencies} requirements for the choice of $S$, that is the choice of $S^+$, arise in the measure estimate of the non-resonant frequency set in Section~\ref{measure_estimate2ssd2}. In order to state the conditions, we recall $S^+=\left\{ 0<j_1<\cdots<j_\nu\right\}$ from \eqref{tan_site} and define\index{$\vec{D}(\xi)$} $\vec{D}(\xi)\in \mathbb{R}^\nu$ for each $\xi\in \mathbb{Z}\backslash S$ by 
 \begin{align}\label{eigenvalusdsd1}
 (\vec{D}(\xi))_k&:=j_k\left( \frac{12}{\pi}H_{4,j_k,-j_k,\xi,-\xi}\right. \nonumber\\
& \left. \  +\frac{18}{\pi^2}\left(\frac{(\xi - j_k)H_{3,-j_k,\xi,-(\xi - j_k)}H_{3,j_k,\xi-j_k,-\xi}}{\lambda^\circ_\alpha(j_k)+\lambda^\circ_\alpha(\xi-j_k)-\lambda^\circ_\alpha(\xi)} \right.\right. \\
& \left.\left.\qquad \qquad \qquad +\frac{(\xi + j_k)H_{3,j_k,\xi,-(\xi + j_k)}H_{3,-j_k,\xi + j_k,-\xi}}{-\lambda^\circ_\alpha(j_k)+\lambda^\circ_\alpha(\xi+j_k)-\lambda^\circ_\alpha(\xi)} \right)\right),
 \end{align}
 if $\xi\ne 0$ and
 \begin{align}\label{eigenvalusdsd2}
 (\vec{D}(\xi))_k = 0 , \text{ if $\xi = 0$.}
 \end{align}
 Note that the explicit expressions for  $H_{4,j_1,j_2,j_3,j_4}, H_{3,j_1,j_2,j_3}$ and $\lambda_\alpha^\circ(j_1)$ for integers $j_1,\ldots,j_4$ can be find in \eqref{moreexplicits2} and \eqref{rkaksd2sdsdmental}, therefore $\vec{D}(\xi)$ is also completely determined by the choice of $S^+$. The motivation of the definition of $\vec{D}(\xi)$ is that $\nabla_{z}^2H_{\zeta}(\theta,0,0)$, the Hessian of $H_{\zeta}$ in the normal direction at $(\theta,0,0)$, can be conjugated by a symplectic transformation to a diagonal operator with the Fourier multiplier $ -\frac{1}{2}m^\circ_{1,\alpha}(\xi) +\frac{T_\alpha}{4} + \epsilon^2 \vec{D}(\xi)\cdot \zeta + o(\epsilon^2)$ (see \eqref{mathfrakFxe} and \eqref{seamxxsxcondresauesc1}).  We observe that the Fourier multiplier acts on the normal sites $S^\perp$, while we can infer  from our choice of $S$ in \eqref{tan_site} and Lemma~\ref{msdjsd11111sd}  that we can even restrict the action of the multiplier to the $\frac{2\pi}{\mathtt{M}}$-translation invariant functions. For this reason, we denote\index{$S^\perp_{\mathtt{M}}$} 
 \begin{align}\label{mnorlam}
 S^\perp_{\mathtt{M}}:=\left\{ \mathtt{M}j\in \mathbb{Z} : j\in S^\perp \right\}.
 \end{align}
With $S$ in \eqref{tan_site}, we see that $S\cup S^\perp_{\mathtt{M}}$ is the set of all $\mathtt{M}$ multiples of non-zero integers. With $\vec{D}(\xi)$ in \eqref{eigenvalusdsd1} and \eqref{eigenvalusdsd2}, and the set $S^\perp_\mathtt{M}$ in \eqref{mnorlam}, we state the conditions on the choice of $S$ as follows: Denoting\index{$W(j)$}
\begin{align}\label{defoflsksdjwwww}
W(j):=\frac{1}{\pi}j \kappa_k\overset{\eqref{moreexplicits2},\eqref{lambdadeffsd}}=-\frac{1}{2}\lambda_\alpha(j) + \frac{T_\alpha}{4}j, 
\end{align} 
\begin{enumerate}[label=$(\mathtt{H}\arabic*)$]
  \setcounter{enumi}{1}
  \item \label{hypothsdj2} There exists a constant $C_{\mathtt{H}2}(\alpha,\nu,S)>0$ such that for $\xi_1,\xi_2\in S_{\mathtt{M}}^\perp$ satisfying $\max\left\{ \xi_1,\xi_2\right\} \ge C_{\mathtt{H2}}$, it holds that
  \begin{align}
  & |\xi_1 - \xi_2| (|\xi_1|^{\alpha-1} + |\xi_2|^{\alpha-1})\nonumber \\
  & \le_{\alpha,\nu,S} \left|\left(W(\xi_1) - \frac{\pi}6\xi_1\vec{D}(\xi_1)\cdot\mathbb{A}^{-1}\overline{\omega}\right)-\left(W(\xi_2) - \frac{\pi}6\xi_2\vec{D}(\xi_2)\cdot\mathbb{A}^{-1}\overline{\omega}\right)\right|,\label{poisosuduwoqs}
  \end{align}
  where $\kappa_{\xi}$ and $\overline{\omega}$ are as in \eqref{moreexplicits2} and \eqref{linear_frequency_aa}, respectively.
  \item \label{hypothsdj22} Define\index{$\mathbb{C}_{\xi_1,\xi_2}$}\index{$\mathbb{B}_{\xi_1,\xi_2}$} $\mathbb{C}_{\xi_1,\xi_1},\mathbb{B}_{\xi_1,\xi_2}\in\mathbb{R}^{\nu\times\nu}$ by
\begin{align}\label{mathbsk2sdx}
\left(\mathbb{B}_{\xi_1,\xi_2}\right)^{i}_{k}:= \frac{(\xi_1\vec{D}(\xi_1) - \xi_2\vec{D}(\xi_2))_{i}}{W(\xi_1) - W(\xi_2)}(\overline{\omega})_{k},\quad  \mathbb{C}_{\xi_1,\xi_2}:= \mathbb{A}^T - \frac{\pi}6\mathbb{B}_{\xi_1,\xi_2}.
\end{align}
Then, 
\begin{align}\label{oosdhwwdwd1ssd}\text{det}(\mathbb{C}_{\xi_1,\xi_2})\ne 0,\text{ for each $\xi_1,\xi_2\in S_\mathtt{M}^\perp\cup \left\{ 0 \right\}$ such that $\xi_1\ne \xi_2$.}
\end{align}
\end{enumerate}

 The hypotheses \ref{hypothsdj2} and \ref{hypothsdj22} are sufficient conditions for us to prove that the set of $\omega\in \Omega_\epsilon$ that produces persistent quasiperiodic solutions have positive measure.  We will not state the precise motivation for the above hypotheses here (see Lemma~\ref{rkhxcxc11}), however, we emphasize that \ref{hypothsdj2} and \ref{hypothsdj22} solely depend on the choice of $S$, that is, the choice of $S_0$ and $\mathtt{M}$ in \eqref{tan_site}.

Now the question is whether one can pick the tangential sites $S$ so that $S$ satisfies \ref{tangent_2}-\ref{tangent_1} as well as \ref{hypos1s1}-\ref{hypothsdj22}. The following proposition ensures that there are infinitely many choices of $S$ for which all the properties are satisfied. To avoid the complexity of the paper, we postpone the proof to Appendix~\ref{applanxxx}. See Proposition~\ref{vefiisdkswjsjdsdwsd}.
\begin{proposition}\label{freq_amp}
Let $\alpha\in (1,2)$ and $2\le \nu\in\mathbb{N}$. There exist infinitely many $S_0^+$ and $\mathtt{M}$ such that $S$ constructed in \eqref{tan_site} satisfies \ref{tangent_2},\ref{tangent_1} and the hypotheses \ref{hypos1s1}-\ref{hypothsdj22}.
\end{proposition}

\color{black}

\chapter{The nonlinear functional setting}\label{nonlads}
We look for $\omega$ and an embedding $i:\mathbb{T}^\nu \ni\varphi \mapsto (\theta(\varphi),y(\varphi),z(\varphi))\in  \mathbb{T}^\nu\times \mathbb{R}^\nu \times H_{S^\perp}$ such that\index{$\mathcal{F}_\omega$} 
 \begin{align}\label{invariant_torus3}
\mathcal{F}_\omega(i)(\varphi):=\omega\cdot\partial_\varphi i(\varphi) - X_{H_{\zeta}}(i(\varphi)) =:\D_\omega i (\varphi) - X_{H_\zeta}(i(\varphi))= 0.
 \end{align}
 Clearly, a solution $(\omega,i)$ to the above equation determines a quasiperiodic solution\index{quasiperiodic solutions} to the Hamiltonian equation \eqref{dynamics2} as
 \[
 t\mapsto i(\omega t).
 \] 
 As described in Section~\ref{rpoisdsd1sd}, we will use $\omega$ as a parameter,  instead of $\zeta$, more precisely, $\zeta$ in \eqref{invariant_torus3} will be thought of as a function of $\omega$ in such a way that (see $\omega(\zeta)$ in \eqref{normal_formaa3})
 \begin{align}\label{xi_omega_dependent}
 \zeta(\omega):= \frac{\pi}{6\epsilon^2}\mathbb{A}^{-1}(\omega -\overline{\omega}).
 \end{align}
On the frequencies in $\Omega_\epsilon$ (see \eqref{frequency_set1}), we impose the Diophantine condition\index{Diophantine condition}\index{$\Omega_0$}\index{$\gamma$}\index{$\tau$}:
\begin{align}
\Omega_0 &:=\left\{ \omega\in \Omega_\epsilon : |\omega\cdot l | > \gamma|l|^{\tau}, \quad \forall l\in \mathbb{Z}^\nu\backslash \left\{0 \right\}\right\}, \nonumber \\ \gamma &:=\epsilon^{2b} \text{ (recall $b>1$ from \eqref{actionssx})},\quad \tau:=\frac{2}{\alpha-1}+\nu + 2.
\label{frequency_set2}
\end{align}
   Now we consider the functional spaces for the problem \eqref{invariant_torus3}. Recall the function spaces in Section~\ref{functionspaces} with $s> s_0$, where $s_0$ is fixed so that 
\begin{align}\label{s0_fixed}
s_0 > \frac{\nu + 2 }{2}.
\end{align}
Given an embedding $i(\varphi):=(\theta(\varphi),y(\varphi),z(\varphi))$, we think of $i$ as an element in $H^s_\varphi\times H^s_\varphi\times H^s_{\varphi,x}$. 
We note that all the embeddings that will be considered throughout the paper are allowed to depend on the parameter  $\omega\in \Omega_\epsilon$.  For a function on $\Omega_\epsilon$, 
\[
F:\Omega_\epsilon \ni \omega \mapsto (F_1(\omega),F_2(\omega),F_3(\omega))\in H^s_\varphi\times H^s_\varphi\times H^s_{\varphi,x},\] we denote its norm by (for a set $\Omega\subset \Omega_\epsilon$)\index{$\rVert\cdot\rVert_s$}\index{$\rVert\cdot \rVert_s^{\Lip(\gamma,\Omega)}$},
\begin{equation}\label{norms_def_lip}
 \begin{aligned}
 \rVert F(\omega) \rVert_{s} &:= \rVert F_1(\omega)\rVert_{H^s_{\varphi}} +\rVert F_2(\omega) \rVert_{H^s_{\varphi}} + \rVert F_3(\omega)\rVert_{H^s_{\varphi,x}}, \\ & \text{ (see \eqref{norms22} for $\rVert \cdot\rVert_{H^s_{\varphi}},\ \rVert \cdot \rVert_{H^s_{\varphi,x}}$)},\\
  \rVert F \rVert_{s}^{\Lip(\gamma,\Omega)} & := \rVert F_1\rVert_{H^s_{\varphi}}^{\Lip(\gamma,\Omega)} +\rVert F_2\rVert_{H^s_{\varphi}}^{\Lip(\gamma,\Omega)} + \rVert F_3\rVert_{H^s_{\varphi,x}}^{\Lip(\gamma,\Omega)}, \\ & \text{ (see \eqref{omega_dep_norm1} for  $\rVert \cdot\rVert_{H^s_{\varphi}}^{\Lip(\gamma,\Omega)},\ \rVert \cdot \rVert_{H^s_{\varphi,x}}^{\Lip(\gamma,\Omega)}$)},
 \end{aligned}
 \end{equation}
 and we denote
 \[
 \Lip(\Omega,C^\infty):=\cap_{s\ge 0}\Lip\left(\Omega,H_{\varphi}^s\times H^s_{\varphi}\times H^s_{\varphi,x}\right), 
 \]
 see \eqref{omega_dep_norm1} for the definition of $\Lip(\Omega,E)$.
 \begin{remark}\label{lipdeps}
$v_\epsilon$ depends on $\zeta$ smoothly for small enough $y$ and $\epsilon$, with a Lipschitz constant $O(1)$  (see \eqref{actionssx}). This yields that $v_\epsilon$ depends on $\omega$ smoothly with a Lipschitz constant $O(\epsilon^{-2})$, because of \eqref{xi_omega_dependent}. 
 \end{remark}

Regarding the Diophantine condition in \eqref{frequency_set2}, for $\omega\in \Omega_0$ and  a function $u$ with zero-average in $\varphi$, that is, $\int u(\varphi) d\vphi = 0$, one can easily check (see \cite[eq. (5.21)]{Baldi-Berti-Montalto:KAM-quasilinear-kdv})
\begin{align}
&\rVert \D_\omega^{-1}u\rVert_s \le_{s,\tau} \gamma^{-1}\rVert u \rVert_{s+\tau} \quad \text{ for $u\in C_{\varphi,x}^\infty$}\label{domega_inverse1}
\\
&\rVert \D_\omega^{-1}u\rVert_s^{\Lip(\gamma,\Omega_1)} \le_{s,\tau} \gamma^{-1}\rVert u \rVert_{s+2\tau+1}^{\Lip(\gamma,\Omega_1)} \quad \text{ for $u\in \Lip(\Omega_0,C^\infty)$}, \label{domega_inverse}
\end{align}
where $\mathcal{D}_\omega:=\omega\cdot\partial_\varphi$\index{$\mathcal{D}_\omega$}.
\section{Regularity of the functional $\mathcal{F}_\omega$}\label{regularity_functionss}
In this section, we study the regularity of the domain/image spaces of the functional $\mathcal{F}_\omega$.   Throughout the paper, we will focus on the embeddings that are close to the trivial embedding. More precisely, we consider the embeddings $i$ such that 
 \begin{align}\label{small_I}
\rVert \mathfrak{I}\rVert^{\Lip(\gamma,\Omega_1)}_{s_0+\mu} \le C \epsilon^{6-2b}\gamma^{-1}= C\epsilon^{6-4b},  \text{ for some $\Omega_1\subset \Omega_0$, $\mu,C>0$}, 
\end{align}
where $\mathfrak{I}(\varphi):=i(\varphi)-i_{triv}(\varphi).$
Recalling the norms in \eqref{norms_def_lip} and the smallness assumption in \eqref{small_I}, we give the estimates on  the functional $\mathcal{F}_\omega$:
\begin{lemma}\label{X_estimate}
There exists  $\mu_2=\mu_2(\alpha, \nu)>0$ such that if an embedding $i\in \Lip(\Omega_1,C^\infty)$ satisfies \eqref{small_I}  for $\mu=\mu_2$ and $C>0$, then (denoting $\rVert \cdot \rVert^{\Lip(\gamma,\Omega_1)}$ by $\rVert \cdot \rVert^{\Lip}$ and $\le_{s,C}$ by $\le$ for simplicity), 
\begin{align}
&\rVert\partial_y P(i)\rVert_s^{\Lip} \le\epsilon^{4} + \epsilon^{2b}\rVert \mathfrak{I}\rVert_{s+\mu_2}^{\Lip}, \quad \quad \rVert\partial_\theta P(i)\rVert_s^{\Lip} \le \epsilon^{6-2b}+\epsilon^{11-6b}\rVert \mathfrak{I}\rVert_{s+\mu_2}^{\Lip}\label{tame1}\\
&\rVert \nabla_z P(i)\rVert_s^{\Lip} \le \epsilon^{5-b} + \epsilon^{6-3b}\rVert \mathfrak{I}\rVert_{s+\mu_2}^{\Lip}, \quad \quad \rVert X_P(i)\rVert_s^{\Lip} \le \epsilon^{6-2b}+\epsilon^{2b}\rVert \mathfrak{I}\rVert_{s+\mu_2}^{\Lip}\label{tame2}\\
&\rVert \partial_\theta\partial_y P(i)\rVert_s^{\Lip} \le\epsilon^{4} + \epsilon^{5-2b}\rVert \mathfrak{I}\rVert_{s+\mu_2}^{\Lip}, \quad \quad \rVert \partial_y\nabla_zP(i)\rVert_s^{\Lip} \le \epsilon^{3+b}+\epsilon^{2b-1}\rVert \mathfrak{I}\rVert_{s+\mu_2}^{\Lip}\label{tame3}\\
&\rVert \partial_{yy}P(i) - 6\epsilon^{2b}\mathbb{A}\rVert_{s}^{\Lip}\le \epsilon^{2+2b}+\epsilon^{3}\rVert \mathfrak{I}\rVert_{s+\mu_2}^{\Lip}, \quad \text{ where $\mathbb{A}$ is as in \eqref{amplitude_modulation1}}.\label{tame4}
\end{align}
Furthermore, for $\ihat(\omega)\in \Lip(\Omega_1,C^\infty)$, it holds that
\begin{align}
&\rVert \partial_y d_iX_P(i)[\ihat] \rVert_{s}^{\Lip} \le  \epsilon^{2b-1}(\rVert \ihat \rVert_{s+\mu_2}^{\Lip} + \rVert \mathfrak{I} \rVert_{s+\mu_2}^{\Lip}\rVert \ihat \rVert_{s_0+\mu_2}^{\Lip})\label{tame5}\\
&\rVert d_iX_{H_{\zeta}}(i)[\ihat] - (0,0, (-\frac{1}{2}\Lambda^{\alpha-1}\partial_x+\frac{T_\alpha}{4}\partial_x)\hat{z}) \rVert_{s}^{\Lip} \le  \epsilon(\rVert \ihat \rVert_{s+\mu_2}^{\Lip} + \rVert \mathfrak{I} \rVert_{s+\mu_2}^{\Lip}\rVert \ihat \rVert_{s_0+\mu_2}^{\Lip})\label{tame6}\\
&\rVert d_i^2X_{H_{\zeta}}(i)[\ihat,\ihat] \rVert_{s}^{\Lip} \le  \epsilon(\rVert \ihat \rVert_{s_0+\mu_2}^{\Lip}\rVert \ihat \rVert_{s+\mu_2}^{\Lip} + \rVert \mathfrak{I} \rVert_{s+\mu_2}^{\Lip}\left(\rVert \ihat \rVert_{s_0+\mu_2}^{\Lip}\right)^2),\label{tame7}\\
&\rVert \mathcal{F}_\omega(i_{triv})\rVert^{\Lip}_s\le \epsilon^{6-2b},\nonumber \\
& \rVert d_i^2\mathcal{F}_\omega(i)[\ihat,\ihat]\rVert_s^{\Lip}\le \epsilon(\rVert \ihat \rVert_{s_0+\mu_2}^{\Lip}\rVert \ihat \rVert_{s+\mu_2}^{\Lip} + \rVert \mathfrak{I} \rVert_{s+\mu_2}^{\Lip}\left(\rVert \ihat \rVert_{s_0+\mu_2}^{\Lip}\right)^2).\label{tame92}
\end{align}
Especially, it holds that for $i,i_1,i_2\in \Lip(\Omega_1,C^\infty)$ satisfying \eqref{small_I},
\begin{align}\label{fffs2x}
&\rVert \mathcal{F}_\omega(i_1)-\mathcal{F}_\omega(i_2) - d_i\mathcal{F}_\omega(i_2)[i_1-i_2]\rVert^{\Lip}_{s}\nonumber\\
& \quad \le \epsilon\left(\rVert i_1-i_2\rVert^{\Lip}_{s+\mu_2}\rVert i_1-i_2\rVert^{\Lip}_{s_0+\mu_2} +(\rVert \mathfrak{I}_1\rVert^{\Lip}_{s+\mu_2} + \rVert \mathfrak{I}_2\rVert^{\Lip}_{s+\mu_2}) \left(\rVert i_1-i_2\rVert^{\Lip}_{s_0+\mu_2}\right)^2\right),\end{align}
where $\mathfrak{I}_2:=i_2 - i_{triv}$.
\end{lemma}
\begin{proof}
We give a proof for $\partial_yP$ in \eqref{tame1} only, since all the other estimates can be obtained in a similar manner.

Recall the definition of $P$ from Proposition~\ref{normal_formaa5}:
\begin{align}
P(\theta,y,z) &:= \underbrace{\epsilon\left( H_{3,2}( v_\epsilon(\theta,y)+z) - H_{3,2}(v_\epsilon(\theta,0)+z)\right)}_{=:P_1} \nonumber \\& + \underbrace{\epsilon^2 \left( H_{4,2}(v_\epsilon(\theta,y)+z) - H_{4,2}(v_\epsilon(\theta,0) + z)\right)}_{=:P_2}\nonumber\\
&\  +6\underbrace{\epsilon^{2b}y^T \mathbb{A}y}_{=:P_3} + \underbrace{\epsilon^{-2b}\left( H_{3,3} +  \sum_{i=3}^4H_{4,i}+\sum_{i=2}^{5}{H_{5,i}} + H_{\ge6}\right)\circ U_\zeta(\theta,y,z)}_{=:P_4}.\label{where4sd}
\end{align}
We prove  the estimate for $P_1$ only, that is,
\begin{align}\label{main_contriris}
\rVert\partial_y P_1(i)\rVert_s^{\Lip} \le_{s,C}  \epsilon^{2b}\rVert \mathfrak{I}\rVert_{s+\mu_2}^{\Lip},\text{ for some $\mu_2\ge0$,}
\end{align}
since the other terms can be treated in the same way. For $P_1$, 
we recall from \ref{item3_normal_3} in Proposition~\ref{normal_form_prop11231} that $H$ is $\mu_1$-regular for some $\mu_1(\alpha)>0$. From its definition in Definition~\ref{mureft}, we can easily see that there exists a bilinear map $R:(H^{s+\mu_1}(\mathbb{T}))^2\mapsto H^{s}(\mathbb{T})$ such that $\nabla_v H_{3,2}(f)=\Pi_S R(z,z)$ (denoting $f=v+z$ and $\nabla_{v}:=\Pi_{S}\nabla_{L^2_0}$) and
\begin{align*}
\rVert R(z,z) \rVert_{H^{s}(\mathbb{T})}&\le_{s,C} \rVert z \rVert_{H^{s+\mu_1}(\mathbb{T})}\rVert z\rVert_{H^{s_0+\mu_1}(\mathbb{T})},\end{align*}
which implies (denoting $i(\varphi)=(\theta_0(\varphi),y_0(\varphi),z_0(\varphi))$), 
\begin{align}\label{2mjs2}
\rVert R(z_0(\cdot),z_0(\cdot))\rVert_{s}^{\Lip} \le_s \rVert z_0 \rVert^{\Lip}_{s+\mu_1}\rVert z_0\rVert_{s_0+\mu_1}^\Lip\le_{s,C} \epsilon^{6-4b}\rVert \mathfrak{I} \rVert^{\Lip}_{s+\mu_1},
\end{align}
where the last inequality follows from \eqref{small_I}.
Note that the elementary chain rule gives us that
 \begin{align}\nonumber
\partial_yP_1(\theta,y,z)& =\epsilon(\nabla_v( H_{3,2}( v_\epsilon(\theta,y)+z) - H_{3,2}(v_\epsilon(\theta,0)+z)),\partial_yv_\epsilon(\theta,y))_{L_x^2} \\
& = \epsilon(R(z,z),\partial_yv_\epsilon(\theta,y))_{L^2_x}.\label{parts_2}
\end{align}
Therefore, $\partial_y P_1(i)(\varphi)= \epsilon(R(z_0(\varphi),z_0(\varphi)),\partial_yv_\epsilon(\theta_0(\varphi),y_0(\varphi)))_{L^2_x}$, while it follows from \eqref{actionssx} that
\[
\partial_y v_\epsilon(\theta_0,y_0)=\sum_{j_k\in S}\frac{\epsilon^{2(b-1)}|j_k|}{2\sqrt{|j_k|(\zeta_k + \epsilon^{2(b-1)}(y_0)_k)}}e^{\ii (\theta_0)_k }e^{\ii j_k x}.
\]
Since $\zeta_k\in [1,2]$ (see \eqref{frequency_set1}), the denominator is strictly positive for all sufficiently small $\epsilon>0$. Therefore  (see Remark~\ref{lipdeps})
\[
\rVert\partial_y v_\epsilon(\theta_0(\cdot),y_0(\cdot))\rVert_s^{\Lip}\le_{s,C} \epsilon^{2(b-1)}(1+\rVert \mathfrak{I}\rVert^{\Lip}_{s}).
\]
  Combining this with \eqref{2mjs2} and applying Lemma~\ref{interpolation_2s} to \eqref{parts_2}, we obtain 
 \begin{align*}
 \rVert \partial_y P_1(i) \rVert_s^{\Lip}&\le_{s,C} \epsilon\left(\rVert \partial_y v_\epsilon\rVert_{s}^{\Lip}\rVert R(z_0,z_0)\rVert^{\Lip}_{s_0} +\rVert \partial_y v_\epsilon\rVert_{s_0}^{\Lip}\rVert R(z_0,z_0)\rVert^{\Lip}_{s}\right)\\
 & \le_{s,C} \epsilon^{2b-1}\left( \epsilon^{6-4b} (1+\rVert \mathfrak{I}\rVert^{\Lip}_s)\rVert \mathfrak{I}\rVert_{s_0+\mu_1}^{\Lip} + \epsilon^{6-4b}(1+\rVert\mathfrak{I}\rVert^{\Lip}_{s_0})\rVert \mathfrak{I}\rVert_{s+\mu_1}\right).
 \end{align*}
  Assuming \eqref{small_I} for some $\mu=\mu_2$ and $\mu_2$ is large enough depending in $\mu_1$, the above inequality implies 
  \[
  \rVert \partial_y P_1(i) \rVert_s^{\Lip}\le_{s,C} \epsilon^{5-2b}\rVert \mathfrak{I}\rVert_{s+\mu_2}^{\Lip}.
  \]
  Since $b-1 \in (0,1/12)$ as fixed in \eqref{actionssx}, we have $\epsilon^{5-2b}<\epsilon^{2b}$, thus the estimate \eqref{main_contriris} follows.
  
  We point out that $\epsilon^4$ in the estimate for $\rVert \partial_y P \rVert^{\Lip}_s$ in \eqref{tame1} is due to the term $P_4$ in \eqref{where4sd}, especially $H_{\ge 6}$. 
\end{proof}

\section{Reversible and $\frac{2\pi}{\mathtt{M}}$-translation invariant  solutions}
\subsubsection{Reversible solutions}
From Lemma~\ref{rsksdtimerevsd}, we already know that $H_\zeta$ is reversible with respect to $\rho_*$ in \eqref{definiteionofinv2}. Therefore we are led to look for a reversible embedding $i$ such that
\[
\rho_*(i(\varphi)) = i(-\varphi),\text{ that is, } (\theta(\varphi),y(\varphi),z(\varphi)) =  (-\theta(-\varphi),y(-\varphi),\rho(z(-\varphi))).
\]
Recalling the function spaces in \eqref{space_reversibles}, we have the following lemma:
\begin{lemma}\label{function_spacess}
$\mathcal{F}_\omega(\cdot):X_{i}\mapsto Y_{i}$.
\end{lemma}
\begin{proof}
It is clear that, $\partial_\varphi:X_i\mapsto Y_i$. For $X_{H_\zeta}$,  with a simple computation, the reversibility of $H_\zeta$ implies that for each $i=(\theta,y,z)$,
\[
-\rho_*(X_{H_\zeta}\circ \rho_*(i)) = X_{H_\zeta}(i).
\]
Note that for $i\in X_{i}$, it holds that $\rho_*(i(-\varphi)) =i(\varphi)$, therefore,
\[
X_{H_{\zeta}}(i(-\varphi)) = -\rho_*(X_{H_\zeta}(\rho_*i(-\varphi))) = -\rho_*(X_{H_\zeta}(i(\varphi))),
\]
which implies $X_{H_\zeta}(i)\in Y_{i}$.
\end{proof}

\subsubsection{$\frac{2\pi}{\mathtt{M}}$-translation invariant  solutions}
From Lemma~\ref{msdjsd11111sd}, we already know that $H_\zeta$ is $\frac{2\pi}{\mathtt{M}}$-translation invariant. Recalling the space $X_{i,\mathtt{M}}$ from \eqref{mforlsd1sd}, we have the following lemma:
\begin{lemma}\label{sjdsd1sdsd}
$
\mathcal{F}_\omega(\cdot):X_{i,\mathtt{M}}\mapsto X_{i,\mathtt{M}}.$
\end{lemma}
\begin{proof}
Recalling  $\rho_{\mathtt{M},*}$\index{$\rho_{\mathtt{M},*}$} in \eqref{mforlsd1sd}, Lemma~\ref{msdjsd11111sd} tells us that $H_\zeta\circ \rho_{\mathtt{M},*} = H_\zeta$, which implies
\[
X_{H_\zeta}(\rho_{\mathtt{M},*}(i)) = \rho_{\mathtt{M},*}( X_{H_\zeta}(i)).
\]
Since $i\in X_{i,\mathtt{M}}$ implies $\rho_{\mathtt{M},*}(i) = i$, the above equality gives us that 
\begin{align}\label{jsd11sdsd}
\rho_{\mathtt{M},*}(X_{H_\zeta}(i)) = X_{H_{\zeta}(i)}.
\end{align}
Furthermore, it is trivial that $\rho_{\mathtt{M},*}\left( \D_\omega (i)\right) = \D_\omega(\rho_{\mathtt{M},*}(i)) = \D_\omega(i)$. Combining this with  \eqref{jsd11sdsd}, we see that
\begin{align}\label{ksd1sdsds}
\rho_{\mathtt{M},*}(\mathcal{F}_\omega(i)) = \rho_{\mathtt{M},*}(\D_\omega (i) - X_{H_\zeta}(i))  = \mathcal{F}_\omega(i),
\end{align}
hence, $\mathcal{F}_\omega(i)\in X_{i,\mathtt{M}}$.
\end{proof}

 \color{black}
  \section{Statement of the main theorem}  
 We fix\index{$\mathtt{p}$}\index{$b$}
  \begin{equation}\label{parametersets1231}
\begin{aligned}
\mathtt{p}&:=\left\{(\alpha,\nu,s_0): \ \alpha\in (1,2),\ 2\le \nu\in \mathbb{N},\ s_0\ge \frac{\nu+2}2\right\},\quad b\in (1,1+1/12).
\end{aligned}
\end{equation}
We are ready to state our main theorem of this paper: 
\begin{theorem}\label{main1}
Given $\mathtt{p}$ in \eqref{parametersets1231}, let us fix the tangential sites $S$\index{tangential sites} as in \eqref{s00sd1sd}, that is, fix $S_0^+$ and $\mathtt{M}$ so that \ref{tangent_2},\ref{tangent_1} and the hypotheses \ref{hypos1s1}-\ref{hypothsdj22} hold. Then there exists $\epsilon_0=\epsilon_0(\mathtt{p},S^+_0,\mathtt{M},b)>0$ such that 
 for all $\epsilon \in (0,\epsilon_0)$,  there exist a Cantor-like  set\index{Cantor-like set}\index{$\mathcal{C}_\epsilon$} $\mathcal{C}_\epsilon\subset \Omega_\epsilon$ (see \eqref{frequency_set1}) such that
\begin{align}\label{rpsdsetfsd}
\lim_{\epsilon\mapsto 0}\frac{\left| \mathcal{C}_\epsilon \right|}{\left| \Omega_\epsilon \right|} = 1,
\end{align}
and for each $\omega\in \mathcal{C}_\epsilon$, there exists a reversible, $\frac{2\pi}{\mathtt{M}}$-translation invariant embedding $i_\infty(\omega)\in H_{\varphi}^{s_0}\times H_{\varphi}^{s_0}\times H_{\varphi,x}^{s_0}$ which solves \eqref{invariant_torus3}. Furthermore,  $i_\infty(\omega)$ is linearly stable under $\frac{2\pi}{\mathtt{M}}$-translation invariant perturbations.
\end{theorem}
\begin{remark}\label{jrjsdpppwwdds}
By the linear stability\index{linear stability} under $\frac{2\pi}{\mathtt{M}}$-translation invariant perturbations, we mean the following: If $I(0)=(\theta(0),y(0),z(0))$ is $\frac{2\pi}{\mathtt{M}}$-translation invariant, and $I(t)$ solves the linearized Hamiltonian equation,
\[
\dot{I}(t) = d_i X_{H_\zeta}(i_\infty(\omega)(\omega t))[I(t)],
\]
then the norm,  $\sup_{t\ge 0}\left( |\theta(t)| +|y(t)| +\rVert z(t)\rVert_{H^{s_0}_{x}} \right) < C$, for some constant $C>0$.  
\end{remark}
\color{black}



\color{black}

\chapter{Approximate inverse}\label{aprisdinv}

In this chapter, we will study the linearized operator of $i\mapsto  \mathcal{F}_{\omega}(i)$ in \eqref{invariant_torus3}, in order to implement the Nash-Moser scheme to find a solution to \eqref{invariant_torus3}. For each fixed $\omega$, the linearized operator\index{linearized operator} of $\mathcal{F}_\omega$ at $(i_0(\omega)$ is given by
\begin{align}\label{tame200zwanzig}
d_{i}\mathcal{F}_\omega(i_0)[\ihat] = \omega\cdot\partial_\varphi \ihat - d_iX_{H_{\zeta}}(i_0)[\ihat].
\end{align}
As mentioned in the previous chapter, we will focus on the reversible, $\frac{2\pi}{\mathtt{M}}$-translation invariant embeddings that are close to the trivial embedding\index{trivial embedding}. More precisely, we assume the following:
The embedding $i_0$ depends on $\omega$ for some $\Omega_1\subset \Omega_0$ and denoting\index{$Z$}
\begin{align}\label{IZdef}
\mathfrak{I}_0:=i_0- i_{triv},\quad Z := \mathcal{F}_\omega(i_0),
\end{align}
we assume that
\begin{align}
&\rVert \mathfrak{I}_0\rVert^{\Lip(\gamma, \Omega_1)}_{s_0+\mu} \le \mathtt{C}\epsilon^{6-2b}\gamma^{-1},\quad \rVert Z \rVert_{s_0+\mu}^{\Lip(\gamma,\Omega_1)}\le \mathtt{C}\epsilon^{6-2b},\quad \text{ for some $\mu,\mathtt{C}>0$},\label{size_assumption_2}\\
&\text{$i_0(\omega)$ is a reversible and $\frac{2\pi}{\mathtt{M}}$-translation invariant.}\label{sjwojsdsdiwosdsd}
\end{align}

We will look for an approximate inverse of the linear operator $d_i\mathcal{F}_\omega(i_0)$.  
Towards the construction of an approximate inverse, we first construct an isotropic embedding. Given $ i_0 = i_0(\vphi)=(\theta(\vphi),y(\vphi),z(\vphi))$, we define \index{$i_\delta$}
\begin{align}\label{iso_tropic}
i_\delta(\vphi) = (\theta_0(\vphi),y_\delta(\vphi),z_0(\vphi)),
\end{align}
where
\begin{align}
y_\delta(\vphi)&:=y_0(\vphi) + (\partial_\vphi\theta_0(\vphi))^{-T}[\rho(\vphi)], \quad \rho_j(\vphi):= \Delta_{\vphi}^{-1}\left( \sum_{k=1}^\nu \partial_{\vphi_j}A_{k,j}(\vphi)\right),\label{Y_def1}\\
A_{k,j}(\vphi)&:=\partial_{\vphi_k}a_j(\vphi)-\partial_{\vphi_j} a_k(\vphi),\nonumber \\
a_k(\vphi)& :=-((\partial_\vphi \theta_0(\vphi))^{T}[y_0(\vphi)])_k + \frac{1}{2}(\partial_{\vphi_k}z_0(\vphi),\partial_x^{-1}z_0(\vphi))_{L^2}.\label{Y_def2}
\end{align}

\begin{lemma}\label{isotropicisrever}
$i_\delta$ is a reversible embedding, that is, $\rho_*\circ i_\delta(\varphi) =  i_\delta(\varphi)$, where $\rho_*$ is the involution given in \eqref{definiteionofinv2}. Furthermore, $i_\delta$ is $\frac{2\pi}{\mathtt{M}}$-translation invariant.
\end{lemma}
\begin{proof}
For the reversibility, it is sufficient to prove that $y_\delta(-\varphi)=y_\delta(\varphi)$. Using that $i_0$ is reversible ($\rho_*\circ i_0 = i_0$), it is straightforward that $\varphi\mapsto a_k(\varphi)$ in \eqref{Y_def2} is even. Using \eqref{Y_def1}, one can easily show that $y_\delta$ is even as well. Since $i_\delta$ does not change $z_0$, it is trivial that $i_\delta$ is $\frac{2\pi}{\mathtt{M}}$-translation invariant.
\end{proof}
\begin{lemma}\cite[Lemma 6.2, Lemma 6.3]{Baldi-Berti-Montalto:KAM-quasilinear-kdv}\label{i_delta_property}
 $i_\delta$ is isotropic\index{isotropic}, that is, $i_\delta(\omega)^*\mathcal{W}=0$ on $\mathbb{T}^\nu$, where $i_\delta(\omega)^*\mathcal{W}$ is the pullback of $\mathcal{W}$ by $i_\delta(\omega)$ and  $\mathcal{W}$ is as in \eqref{symplectic_form1}. Furthermore, there exists $\mu_3(\tau)>0$ such that   for all $s\ge s_0$,
\begin{align}
\rVert y_\delta - y_0 \rVert_s^{\Lip(\gamma,\Omega_1)} &\le_s \rVert\mathfrak{I}_0 \rVert_{s+\mu_3}^{\Lip(\gamma,\Omega_1)},\label{krjsxcxcxc1} \\
\rVert y_\delta - y_0 \rVert_s^{\Lip(\gamma,\Omega_1)}& \le_s \gamma^{-1}\left( \rVert  Z\rVert_{s+\mu_3}^{\Lip(\gamma,\Omega_1)}+ \rVert  Z\rVert_{s_0+\mu_3}^{\Lip(\gamma,\Omega_1)}\rVert \mathfrak{I}_0\rVert_{s+\mu_3}^{\Lip(\gamma,\Omega^1)} \right)\label{y_estimate}\\
 \rVert \mathcal{F}_\omega(i_\delta) \rVert_{s}^{\Lip(\gamma,\Omega_1)} &\le_s \rVert Z \rVert_{s+\mu_3}^{\Lip(\gamma,\Omega_1)} + \epsilon^{-1}\rVert Z \rVert_{s_0+\mu_3}^{\Lip(\gamma,\Omega_1)}\rVert \mathfrak{I}_0 \rVert_{s+\mu_3}^{\Lip(\gamma,\Omega_1)},\label{F_estimate}\\
\rVert d_ii_\delta(i_0)[\ihat] \rVert_s^{\Lip(\gamma,\Omega_1)}  &\le_s \rVert \ihat \rVert_{s+\mu_3}^{\Lip(\gamma,\Omega_1)} + \rVert \mathfrak{I}_0 \rVert_{s+\mu_3}^{\Lip(\gamma,\Omega_1)}\rVert \ihat \rVert_{s_0+\mu_3}^{\Lip(\gamma,\Omega_1)},\label{di_estimate}
\end{align}
for $\ihat\in \Lip(\Omega_1,C^\infty)$.
\end{lemma}
In view of \eqref{iso_tropic}, the estimate \eqref{krjsxcxcxc1} tells us that
\begin{align}\label{isjkxcoksd1sd}
\rVert \mathfrak{I}_\delta\rVert^{\Lip(\gamma, \Omega_\epsilon)}_{s}\le_s \rVert \mathfrak{I}_0\rVert^{\Lip(\gamma, \Omega_\epsilon)}_{s+\mu_3}
\end{align}

The goal of this chapter is summarized in the next proposition:
\begin{proposition}\label{approx_inverse}
Let $\mathtt{S}\gg s_0$, $\mathtt{C}>0$ be fixed. There exist $\mu_\mathtt{p}=\mu_{\mathtt{p}}(\mathtt{p})>0$ and $\epsilon_0=\epsilon_0(\mathtt{p},\mathtt{M},S_0,\mathtt{S},b,\mathtt{C})>0$ such that if a reversible, $\frac{2\pi}{\mathtt{M}}$-translation invariant embedding $i_0$ satisfies \eqref{size_assumption_2} for $\mu=\mu_\mathtt{p}$ and $\epsilon\in (0,\epsilon_0)$ for some $\Omega_1\subset \Omega_0$, then the following holds:
There exist a set of frequencies $\Omega_\infty=\Omega_\infty(i_0)\subset \Omega_1$  and a linear map $T(i_0)(\omega):Y_{i}\cap X_{i,\mathtt{M}}\mapsto X_{i}\cap X_{i,\mathtt{M}}$ for $\omega\in \Omega_{\infty}$ such that $T(i_0)(\omega)$ is an approximate inverse of $d_{i}X_{H_{\zeta}}(i_0(\omega))$. More precisely, for all $g(\omega)\in Y_{i}\cap X_{i,\mathtt{M}}$, it holds that $T(i_0)(\omega)[g(\omega)]\in X_i\cap X_{i,\mathtt{M}}$ and  (denoting $\le_{\mathtt{p},b,\mathtt{M},S_0,\mathtt{S},\mathtt{C}}$ by $\le$)
\begin{align}\label{inverse_estimate1}
\rVert &\left(d_{i}\mathcal{F}_\omega(i_0)\circ T(i_0) - I\right)[g]\rVert_s^{\Lip(\gamma,\Omega_\infty)}\nonumber\\
& \le \epsilon^{2b-1}\gamma^{-2}\left(\rVert Z\rVert_{s_0+\mu_\mathtt{p}}^{\Lip(\gamma,\Omega_\infty)}\rVert g\rVert_{s+\mu_\mathtt{p}}^{\Lip(\gamma,\Omega_\infty)}\right. \nonumber \\
& \qquad \left.+ (\rVert Z\rVert_{s+\mu_\mathtt{p}}^{\Lip(\gamma,\Omega_\infty)} + \epsilon^7\gamma^{-4}\rVert Z\rVert_{s_0+\mu_\mathtt{p}}^{\Lip(\gamma,\Omega_\infty)}\rVert \mathfrak{I}_0\rVert_{s+\mu_\mathtt{p}}^{\Lip(\gamma,\Omega_\infty)})\rVert g\rVert_{s_0+\mu_\mathtt{p}}^{\Lip(\gamma,\Omega_\infty)}\right),
\end{align}
and
\begin{align}\label{inverse_estimate2}
\rVert T(i_0)[g] \rVert_s^{\Lip(\gamma,\Omega_\infty)} \le  \gamma^{-1}\left( \rVert g \rVert_{s+\mu_\mathtt{p}}^{\Lip(\gamma,\Omega_\infty)} + \epsilon^{7}\gamma^{-4}\rVert \mathfrak{I}_0\rVert_{s+\mu_\mathtt{p}}^{\Lip(\gamma,\Omega_\infty)} \rVert g \rVert_{s_0+\mu_\mathtt{p}}^{\Lip(\gamma,\Omega_\infty)} \right),
\end{align}
for all $s\in [s_0,\mathtt{S}]$.
\end{proposition}
\begin{remark}\label{derv_loss}
 In Chapter~\ref{aprisdinv}-Chapter~\ref{reduction}, the loss of derivatives will be denoted by  $\mu_0$, which might vary from line to line but depend on only the fixed parameters $\mathtt{p}$ in \eqref{parametersets1231}. The loss of derivatives $\mu_{\mathtt{p}}$ in Proposition~\ref{approx_inverse} will be chosen large enough so that all the estimates in what follows hold. Especially $\mu_{\mathtt{p}}$ is assumed to be larger than $\mu_1,\mu_2,\mu_3$ that appear in \ref{item3_normal_3} of Proposition~\ref{normal_form_f}, Proposition~\ref{X_estimate} and Lemma~\ref{i_delta_property}. 
 
  Also, in view of the dependence of $\epsilon_0$, let us denote\index{$\mathtt{pe}$}
  \begin{align}\label{psd1sdsdsd}
  \mathtt{pe}:=\left\{ \mathtt{p},\mathtt{M},S_0,b,\mathtt{C}\right\}.
  \end{align} 
 Here, $\mathtt{p}$ and $b$ denote the parameters that we fixed in \eqref{parametersets1231}, $\mathtt{M}$ and $S_0$ correspond to the choice of the tangential sites $S$ as in \eqref{tan_site}, and $\mathtt{C}$ is a constant that appears in the size assumption of the embedding in \eqref{size_assumption_2}. Note that $\epsilon$ needs to be small enough depending on $\mathtt{pe}$ and the range of $s$, that is $[s_0,\mathtt{S}]$.  $\mathtt{S}$  will be fixed in the Nash-Moser iteration in Chapter~\ref{NSmoes} (see \eqref{parambesd2sdjj22}). We remark that $\mathtt{S}$ implicitly depends on $\mathtt{p}$, however we will not explicitly trace its size, for simplicity of the proof. Instead, we will assume that $\epsilon$ is always small enough, depending on $s$ as well as $\mathtt{pe}$.
 \end{remark}
 
\begin{remark}\label{Insteadofimbd}  It is an immediate consequence of \eqref{isjkxcoksd1sd} that
given $\mu_0\ge0$, we can choose $\mu=\mu(\mu_0,\mathtt{p})$ in \eqref{size_assumption_2} so that
\begin{align}\label{size_assumption_3}
\rVert \mathfrak{I}_\delta\rVert^{\Lip(\gamma, \Omega_\epsilon)}_{s_0+\mu_0} \le_{\mathtt{pe}} \epsilon^{6-2b}\gamma^{-1}.
\end{align}
We will frequently use  the smallness of $i_\delta$ in \eqref{size_assumption_3}, instead of the smallness of $i_0$, assuming that $\mu_0$ in \eqref{size_assumption_3} is large enough depending on $\mathtt{p}$.
 \end{remark}

The proof of Proposition~\ref{approx_inverse} is identical to \cite[Theorem 6.10]{Baldi-Berti-Montalto:KAM-quasilinear-kdv}, since the estimates in Lemma~\ref{X_estimate} are same as \cite[Lemma 5.3]{Baldi-Berti-Montalto:KAM-quasilinear-kdv}.  We will describe the general idea of the proof for the sake of completeness. The proof of Proposition~\ref{approx_inverse} will be given at the end of this chapter by using Proposition~\ref{normal_inversion}, which will be proved in Chapter~\ref{Linear_op}-\ref{reduction}.

 We define a change of variables $G_\delta:\mathbb{T}^\nu \times \mathbb{R}^\nu \times H_{S^\perp}\mapsto \mathbb{T}^\nu \times \mathbb{R}^\nu \times H_{S^\perp}$ as
  \begin{align}\label{def_G_delta}
  G_\delta \colvec{\psi \\ \eta \\ w} := \colvec{\theta \\ y \\ z }:= \colvec{\theta_0(\psi) \\ y_\delta(\psi) + (\partial_\psi\theta_0(\psi))^{-T}\eta + ((\partial_\theta\tilde{z}_0)(\theta_0(\psi)))^{T}\partial_xw \\ z_0(\psi) + w},
  \end{align}
  where $\tilde{z}_0(\theta):= z_0(\theta_0^{-1}(\theta))$ for $\theta\in \mathbb{T}^\nu$.
  Clearly, the image of the trivial embedding $\vphi\mapsto (\vphi,0,0)$ is mapped to the embedded torus $i_\delta$ by the map $G_\delta$.   We also list some properties of $G$ in the following lemmas.
  \begin{lemma}\label{Gsers}
  $G_\delta$ is reversibility preserving, that is $G_\delta:X_i\mapsto X_i$ and $G_\delta:Y_i\mapsto Y_i$. Furthermore, for the $\frac{2\pi}{\mathtt{M}}$-translation invarance,  we have that, $G_\delta:X_{i,\mathtt{M}}\mapsto X_{i,\mathtt{M}}$.
  \end{lemma}
  \begin{proof}
 We first prove the reversibility preserving of $G_\delta$. It suffices to show that $\rho_*\circ G_\delta\circ \rho_* = G_\delta$, where $\rho_*$ is given in \eqref{definiteionofinv2}. Since $i_\delta$ is a  reversible embedding (see Lemma~\ref{isotropicisrever}),   we have that $\theta_0$ is odd, $y_\delta$ is even and $z_0$ satisfies $\rho z_0(\varphi) = z_0(-\varphi)$. Especially it holds that $
\rho\tilde{z}_0(\theta) = \tilde{z}_0(-\theta)$.
 This implies that $\rho\partial_\theta \tilde{z}_0(\theta) =-\partial_\theta\tilde{z}_0(-\theta)$. Therefore, using $\rho^{T}=\rho$ and $\partial_x(\rho w) = -\rho(\partial_x w)$, which easily follows from its definition in \eqref{definiteionofinv}, we see that
 \begin{align*}
(\partial_\theta \tilde{z}_0(-\theta))^{T}\partial_x(\rho w) & = -( \partial_\theta\tilde{z}_0(-\theta))^{T}\rho(\partial_x w) \\
& = -(\rho\partial_\theta\tilde{z}_0(-\theta))^{T}\partial_x w = (\partial_\theta\tilde{z}_0(\theta))^T\partial_xw,
 \end{align*}
 for all $\theta\in \mathbb{T}^\nu,\ w\in H_{S^\perp}$.
 Hence, we have 
\begin{align*}
G_\delta\circ \rho_*\colvec{\psi \\ \eta \\ w} &= \colvec{\theta_0(-\psi) \\ y_\delta(-\psi) + (\partial_\psi\theta_0(-\psi))^{-T}\eta +(\partial_\theta\tilde{z}(\theta_0(-\psi)))^{T}\partial_x(\rho w),\\  z_0(-\psi) +\rho w } \\
&= \colvec{-\theta_0(\psi) \\ y_\delta(\psi) + (\partial_\psi \theta_0(\psi))^T\eta + (\partial_\theta\tilde{z}(\theta_0(\psi)))^T\partial_xw \\ \rho (z_0(\psi) +w)},
\end{align*}
which gives $\rho_*(G_\delta\circ \rho_*) = G_\delta$. 

 For the preservation of the $\frac{2\pi}{\mathtt{M}}$-translation invariance, the result follows immediately from the definition of $X_{i,\mathtt{M}}$ in \eqref{mforlsd1sd} and the assumption on $i_0$ in \eqref{sjwojsdsdiwosdsd}.
  \end{proof}
  
  \begin{lemma}\cite[Lemma 2]{Berti-Bolle:nash-moser-kam}\label{G_delta_simplectic}
  $G_\delta$ is symplectic.
  \end{lemma}
   \begin{lemma}\cite[Lemma 6.7]{Baldi-Berti-Montalto:KAM-quasilinear-kdv}\label{G_delta_estimate}
  We have that 
  \begin{align*}
  & \rVert DG_\delta(\varphi,0,0)[\ihat]\rVert_s^{\Lip(\gamma,\Omega_1)} + \rVert (DG_\delta(\varphi,0,0))^{-1}[\ihat]\rVert_s^{\Lip(\gamma,\Omega_1)} \\
  &\le_{\mathtt{pe},s}  \rVert \ihat\rVert_{s+\mu_0}^{\Lip(\gamma,\Omega_1)} + \rVert \ihat\rVert_{s_0+\mu_0}^{\Lip(\gamma,\Omega_1)}\rVert \mathfrak{I}_0\rVert_{s+\mu_0}^{\Lip(\gamma,\Omega_1)},  \end{align*}
  and
   \begin{align*}
   \rVert D^2G_\delta(\bar{i})[\ihat_1,\ihat_2])\rVert_s^{\Lip(\gamma,\Omega_1)}
   & \le_{\mathtt{pe},s} \rVert \ihat_1 \rVert_{s}^{\Lip(\gamma,\Omega_1)} \rVert \ihat_2 \rVert_{s_0}^{\Lip(\gamma,\Omega_1)}\\
   & + \rVert \ihat_1 \rVert_{s_0}^{\Lip(\gamma,\Omega_1)} \rVert \ihat_2 \rVert_{s}^{\Lip(\gamma,\Omega_1)}\\
   & \qquad + \rVert \mathfrak{I}_0 \rVert_{s+\mu_0}^{\Lip(\gamma,\Omega_1)} \rVert \ihat_1 \rVert_{s_0}^{\Lip(\gamma,\Omega_1)} \rVert \ihat_1 \rVert_{s_0}^{\Lip(\gamma,\Omega_1)}.
 \end{align*}
  for $\ihat,\ihat_1,\ihat_2\in \Lip(\Omega_1,C^\infty)$.
  \end{lemma}
  
 Composing $H_{\zeta}$ in \eqref{normal_formaa2} with the transformation $G_\delta$, we define a Hamiltonian $K$\index{$K$} to be
 \begin{align}\label{modified_Hamiltonian_2}
 K(\psi,\eta,w):=H_{\zeta}(G_\delta(\psi,\eta,w)).
 \end{align}
 
For each $\psi$, we consider the Taylor expansion of $K$ in $\eta,w$ (in the space $\mathbb{R}^\nu \times H_{S^\perp}$) at $(\psi,0,0)$:
 \begin{equation}\label{taylor_K}
 \begin{aligned}
 K(\psi,\eta,w) & = K_{00}(\psi) +K_{10}(\psi)[\eta] +K_{01}(\psi)[w]+ \frac{1}{2}K_{20}(\psi)[\eta,\eta] \\
 & \ +\frac{1}{2}K_{02}(\psi)[w,w] +K_{11}(\psi)[\eta,w]  +K_{\ge 3}(\psi,\eta,w),  
 \end{aligned}
 \end{equation}
 where 
 \begin{align}\label{rksksd1kpsxx}
 K_{ij}(\psi):=\partial_\eta^i\nabla_w^jK(\psi,0,0),
 \end{align} and  $K_{\ge 3}$ collects all the terms homogeneous in $(\eta,w)$ of order at least $3$. Note that the Hamiltonian vector field $X_{K}$ generated by $K$\index{$X_K$} with the symplectic form $\mathcal{W}$ in \eqref{symplectic_form1} at $(\psi,0,0)$ is
 \[
 X_{K}(\psi,0,0) = \colvec{\frac{1}{2\pi}K_{10}(\psi) \\- \frac{1}{2\pi}\partial_\psi K_{00}(\psi) \\ \partial_x (K_{01}(\psi))}.
 \]

 \begin{lemma}\label{lkksred}
 $X_K$ is reversible, that is, $X_K:X_{i}\mapsto Y_{i}$.
 \end{lemma}
 \begin{proof}
 The proof is identical to Lemma~\ref{function_spacess}. Indeed,  $K$ is also reversible with respect to the involution $\rho_*$, thanks to Lemma~\ref{rsksdtimerevsd}, Lemma~\ref{Gsers} and \eqref{modified_Hamiltonian_2}.
 \end{proof}
  Let us denote the trivial embedding by 
 \begin{align}\label{trivial_embedding}
 \bar{i}(\vphi):=i_{triv}(\vphi)=(\varphi,0,0).
 \end{align}
 Then the linearization of $X_K$ at $\bar{i}$ in the direction $ \ihat(\vphi)=(\hat{\psi}(\vphi),\hat{\eta}(\vphi),\hat{w}(\vphi))$ is given by (using \eqref{taylor_K}), 
 \begin{align}\label{rjoosdwiwnds1s2s2}
 & d_iX_{K}(\bar{i})[\ihat](\vphi) = \frac{d}{dt}\colvec{\frac{1}{2\pi}\partial_\eta K(\bar{i}+t\ihat) \\ -\frac{1}{2\pi} \partial_\psi K(\bar{i}+t\ihat)  \\ \partial_x \nabla_w K(\bar{i}+t\ihat)}\bigg|_{t=0}\nonumber\\
 & = \colvec{\frac{1}{2\pi}\left(\partial_\psi K_{10}(\vphi)[\hat{\psi}(\vphi)] + K_{20}(\vphi)[\hat{\eta}(\vphi)] + K_{11}(\vphi)^T[\hat{w}(\vphi)]\right) \\ \frac{1}{2\pi}\left(-\partial_{\psi\psi}K_{00}(\vphi)[\hat{\psi}(\vphi)] - \partial_\psi K_{10}(\vphi)^T[\hat{\eta}(\vphi)] - \partial_\psi K_{01}(\vphi)^T[\hat{w}(\vphi)]\right) \\ \partial_x(\partial_\psi K_{01}(\vphi)[\hat{\psi(\vphi)}] + K_{02}(\vphi)[\hat{w}(\vphi)]+K_{11}(\vphi)[\hat{\eta}(\vphi)])}. 
 \end{align}
  Now we consider the linear operator  $L$:
  \begin{align}\label{linearized_XK}
 L[\ihat](\vphi):= \D_{\omega}\ihat(\vphi) - d_iX_{K}(\bar{i})[\ihat](\vphi),
  \end{align}
  which is obtained by linearizing $i\mapsto \D_\omega i - X_{K}(i)$ at $\bar{i}$. In order to find an approximate inverse of $L$, we decompose it as
  \begin{equation}\label{decomp_L}
  \begin{aligned}
  L[\ihat] &= \underbrace{\colvec{\D_\omega \hat{\psi}(\vphi) \\ \D_\omega \hat{\eta}(\vphi) \\ \D_\omega \hat{w}(\vphi)} -\colvec{ \frac{1}{2\pi}\left( K_{20}(\vphi)[\hat{\eta}(\vphi)] + K_{11}(\vphi)^T[\hat{w}(\vphi)]\right) \\ 0 \\ \partial_x (K_{02}(\vphi)[\hat{w}(\vphi)]+K_{11}(\vphi)[\hat{\eta}(\vphi)])}}_{=:\mathbb{D}[\ihat]}\\
  & \  - \colvec{\frac{1}{2\pi}\partial_\psi K_{10}(\vphi)[\hat{\psi}(\vphi)]  \\ \frac{1}{2\pi}\left(-\partial_{\psi\psi}K_{00}(\vphi)[\hat{\psi}(\vphi)] - \partial_\psi K_{10}(\vphi)^T[\hat{\eta}(\vphi)] - \partial_\psi K_{01}(\vphi)^T[\hat{w}(\vphi)]\right) \\ \partial_x(\partial_\psi K_{01}(\vphi)[\hat{\psi(\vphi)}])}.
  \end{aligned}
  \end{equation}
  We check that $\mathbb{D}$ is reversible, that is, $\mathbb{D}:X_i\mapsto Y_i$: \index{$\mathbb{D}$}
\begin{lemma}\label{ress2} For 
\[
K_{20}(\varphi):=\partial_{yy}K(\bar{i}(\varphi)),\quad K_{11}(\varphi):=\partial_y\nabla_zK(\bar{i}(\varphi)),\quad K_{02}(\varphi):=\nabla_z^2K(\bar{i}(\varphi)),
\]
 it holds that if $(\hat{\psi},\hat{\eta},\hat{w})\in Y_o\times X_e\times X_{\perp} = X_i$, 
\[
K_{20}[\hat{\eta}]\in X_e,\quad K_{11}[\hat{\eta}]\in Y_{\perp},\quad K_{02}[\hat{w}]\in X_{\perp},\quad K_{11}^{T}[\hat{w}]\in X_e.
\]
Therefore, $\mathbb{D}[\ihat] \in Y_i$.
\end{lemma}
\begin{proof}
We prove $K_{20}$ only since the other operators can be managed in the same way. 

 Since $K$ is reversible (see Lemma~\ref{Gsers} and \eqref{modified_Hamiltonian_2} and Lemma~\ref{rsksdtimerevsd}, which say that $H_\zeta$ is reversible and $G_\delta$ is reversibility preserving) and $\hat{\eta}$ is even, we have that for all $\varphi,t$,
 \[
 K(\varphi,t\hat{\eta}(\varphi),0) = K(-\varphi,t\hat{\eta}(\varphi),0).
 \] 
 Differentiating in $t$ at $t=0$, we have $\partial_{yy}K(\bar{i})[\hat{\eta}(\varphi)] = \partial_{yy}K(\bar{i}(-\varphi))[\hat{\eta}(\varphi)]$. Therefore,
 \begin{align*}
 K_{20}(-\varphi)[\hat{\eta}(-\varphi)] &\overset{\eqref{rksksd1kpsxx}}= \partial_{yy}K(-\bar{i}(\varphi))[\hat{\eta}(-\varphi)] \overset{\hat{y}\in X_e}= \partial_{yy}K(-\bar{i}(\varphi))[\hat{\eta}(\varphi)]\\
 & =\partial_{yy}K(\bar{i})[\hat{\eta}(\varphi)] \overset{\eqref{rksksd1kpsxx}}= K_{20}(\varphi)[\hat{\eta}(\varphi)].
 \end{align*}
 Hence, $\varphi\mapsto K_{20}[\hat{\eta}]$ is even.
\end{proof}

We check $\mathbb{D}$ is $\frac{2\pi}{\mathtt{M}}$-translation invariant, that is, $\mathbb{D}:X_{i,\mathtt{M}}\mapsto X_{i,\mathtt{M}}$.
\begin{lemma}\label{mfisodswwdsd2s2}
For $(\hat{\psi},\hat{\eta},\hat{w})\in X_{i,\mathtt{M}}$, we have that 
\[
K_{02}(\varphi)[\hat{w}], \ K_{11}(\varphi)[\hat{\eta}]\in X_{\mathtt{M}}.
\]
\end{lemma}
  \begin{proof}
  We prove $K_{02}$ only since the other operator can be managed in the same way. 
  
 From Lemma~\ref{msdjsd11111sd}, and Lemma~\ref{Gsers},  we have that  $K$ defined in \eqref{modified_Hamiltonian_2} is $\frac{2\pi}{\mathtt{M}}$-translation invariant.  
 Therefore, for $\hat{w}\in X_{\mathtt{M}}$, we have that for real numbers $t,s$ close to $0$ and $h\in H_{S^\perp}$, we have
 \[
K(\varphi,0,t\hat{w}(\varphi) + s h) = K(\varphi,0,\rho_{\mathtt{M}}(t\hat{w}(\varphi)+sh) = K({\varphi,0,t\hat{w}(\varphi)+s\rho_{\mathtt{M}}h}).
 \]
 Differentiating in $t,s$ at $t=s=0$, we have $(K_{02}(\varphi)[\hat{w}],h)_{L^2_x} = (K_{02}(\varphi)[\hat{w}],\rho_{\mathtt{M}}h)_{L^2_x}$.
 This implies that $\rho_{\mathtt{M}}K_{02}(\varphi)[\hat{w}]=K_{02}(\varphi)[\hat{w}]$, that is, $K_{02}(\varphi)[\hat{w}]$ is invariant under the $\frac{2\pi}{\mathtt{M}}$-translation. Therefore, $K_{02}(\varphi)[\hat{w}]\in X_{\mathtt{M}}$.
  \end{proof}
  The next lemma shows that the approximate inverse of $L$ can be obtained by inverting $\mathbb{D}$. 
  \begin{lemma}\cite[Lemma 6.4]{Baldi-Berti-Montalto:KAM-quasilinear-kdv}\label{vanishing_terms}
  If $i_0$ is an invariant torus, then $L=\mathbb{D}$. More precisely, we have 
    \begin{align}\nonumber
  &\rVert \partial_\psi K_{00}\rVert_s^{\Lip(\gamma,\Omega_1)} + \rVert  K_{10} - \omega \rVert_s^{\Lip(\gamma,\Omega_1)} + \rVert \partial_\psi K_{01} \rVert_s^{\Lip(\gamma,\Omega_1)} \\
  & \le_{\mathtt{pe},s} \rVert  Z \rVert_{s+\mu_0}^{\Lip(\gamma,\Omega_1)} + \epsilon^{-1}\rVert Z \rVert_{s_0+\mu_0}^{\Lip(\gamma,\Omega_1)}\rVert \mathfrak{I}_0\rVert_{s+\mu_0}^{\Lip(\gamma,\Omega_1)}.\label{difference_1}
  \end{align}
  Therefore, we have that,
  \begin{align}\label{difference_2}
  \rVert (L-\mathbb{D})[\ihat]\rVert_s^{\Lip(\gamma,\Omega_1)}
  &  \le_{\mathtt{pe},s}  \rVert Z \rVert_{s_0+\mu_0}^{\Lip(\gamma,\Omega^1)}\rVert \ihat \rVert_{s+\mu_0}^{\Lip(\gamma,\Omega_1)} \nonumber \\
  & +  (\rVert Z \rVert_{s+\mu_0}^{\Lip(\gamma,\Omega_1)} + \epsilon^{-1} \rVert Z\rVert_{s_0+\mu_0}^{\Lip(\gamma,\Omega_1)} \rVert \mathfrak{I}_0\rVert_{s+\mu_0}^{\Lip(\gamma,\Omega_1)}) \rVert \ihat \rVert_{s_0+\mu_0}^{\Lip(\gamma,\Omega^1)}.
  \end{align}
  \end{lemma}

In order to prove the invertibility of $\mathbb{D}$, we recall the estimates for $K_{20}$ and $K_{11}$:
\begin{lemma}\label{Kestimate1}\cite[Lemma 6.6]{Baldi-Berti-Montalto:KAM-quasilinear-kdv}
For $\ihat\in \Lip(\Omega_1,C^\infty)$, we have
\begin{equation}\label{rlaqkqdpqkqakfdjsd}
\begin{aligned}
&\rVert K_{20} - 6\epsilon^{2b}\mathbb{A} \rVert_s^{\Lip(\gamma,\Omega_1)} \le_{\mathtt{pe},s} \epsilon^{2b+2} + \epsilon^{2b}\rVert \mathfrak{I}_0 \rVert_{s+\mu_0}^{\Lip(\gamma,\Omega_1)} \\
&\rVert K_{11}[\hat{\eta}] \rVert_s^{\Lip(\gamma,\Omega_1)} \le_{\mathtt{pe},s} \epsilon^{5-2b}\rVert \hat{\eta} \rVert_{s+\mu_0}^{\Lip(\gamma,\Omega_1)} + \epsilon^{2b-1}\rVert \mathfrak{I}_0 \rVert_{s+\mu_0}^{\Lip(\gamma,\Omega_1)}\rVert \hat{\eta} \rVert_{s_0+\mu_0}^{\Lip(\gamma,\Omega_1)}\\
& \rVert K_{11}^{T}[\hat{w}] \rVert_s^{\Lip(\gamma,\Omega^1)} \le_{\mathtt{pe},s} \epsilon^{5-2b}\rVert \hat{w} \rVert_{s+\mu_0}^{\Lip(\gamma,\Omega_1)} + \epsilon^{2b-1}\rVert \mathfrak{I}_0 \rVert_{s+\mu_0}^{\Lip(\gamma,\Omega_1)}\rVert \hat{w}\rVert_{s_0+\mu_0}^{\Lip(\gamma,\Omega_1)}.
\end{aligned}
\end{equation}
In particular, $\rVert K_{20} - 6\epsilon^{2b}\mathbb{A} \rVert_{s_0}^{\Lip(\gamma,\Omega_1)} \lepe \epsilon^{6-2b}$ and 
\[
\rVert K_{11}[\hat{\eta}] \rVert_{s_0}^{\Lip(\gamma,\Omega_1)} \lepe  \epsilon^{5-2b}\rVert \hat{\eta} \rVert_{s_0+\mu_0}^{\Lip(\gamma,\Omega_1)},  \rVert K_{11}[\hat{w}] \rVert_{s_0}^{\Lip(\gamma,\Omega_1)} \lepe \epsilon^{5-2b}\rVert \hat{w} \rVert_{s_0+\mu_0}^{\Lip(\gamma,\Omega_1)}
\]
\end{lemma}
Now, we study how to invert $\mathbb{D}$.  Given $g:=(g_1,g_2,g_3)\in Y_{i}\cap X_{i,\mathtt{M}}$, we will find $\ihat = (\widehat{\psi},\widehat{\eta},\widehat{w}) \in X_{i}\cap X_{i,\mathtt{M}}$ such that  (recalling $\mathbb{D}$ from \eqref{decomp_L})
\begin{align}\label{inverse_image}\mathbb{D}[\ihat] =\colvec{\D_\omega \hat{\psi}(\vphi) \\ \D_\omega \hat{\eta}(\vphi) \\ \D_\omega \hat{w}(\vphi)} -\colvec{ \frac{1}{2\pi} \left( K_{20}(\vphi)[\hat{\eta}(\vphi)] + K_{11}(\vphi)^T[\hat{w}(\vphi)] \right)\\ 0 \\ \partial_x (K_{02}(\vphi)[\hat{w}(\vphi)]+K_{11}(\vphi)[\hat{\eta}(\vphi)])} =  g\end{align} 
Let us begin with the second equation, $\D_\omega\hat{\eta} = g_2$. Since the right-hand side has zero mean ($g\in Y_{i}$, thus  $\varphi\mapsto g_2(\varphi)$ is odd), we can invert $\D_\omega$, thus,
\begin{align}\label{eta_hat}
\hat{\eta} := \D_\omega^{-1}g_2 + \langle \hat{\eta}\rangle \in X_{e},
\end{align}
where $\langle \hat{\eta}\rangle$ is the average of $\hat{\eta}$ in $\varphi$, which can be freely chosen.  This will be fixed later soon (see \eqref{mean_hateta}). Therefore, we have (using \eqref{domega_inverse} and recalling $|\cdot|^{\Lip}$ from \eqref{constnas1}),
\begin{align}\label{sizeofeta}
\rVert \hat{\eta}\rVert^{\Lip(\gamma,\Omega_1)}_s\le_{\mathtt{pe},s} \gamma^{-1}\rVert g_2 \rVert^{\Lip(\gamma,\Omega_1)}_{s+\mu_0} +|\langle\hat{\eta}\rangle|^{\Lip(\gamma,\Omega_1)}.
\end{align}
 Now we turn to the third equation from \eqref{inverse_image}\index{$\mathcal{L}_\omega$},
 \begin{align}\label{third_equation}
 \mathcal{L}_\omega \hat{w}:=\D_\omega \hat{w} - \partial_x(K_{02}[\hat{w}]) = g_3 + \partial_x(K_{11}[\hat{\eta}]).
 \end{align} 
 Note that $g\in Y_{i}\cap X_{i,\mathtt{M}}$ and Lemma~\ref{mfisodswwdsd2s2},\ \ref{ress2} imply that 
 \begin{align}\label{righsdpwdsdsd}
 g_3 + \partial_x(K_{11}[\hat{\eta}])\in Y_{\perp}\cap X_{\mathtt{M}}.
 \end{align}

 The proof of the following proposition will be given through Chapter~\ref{Linear_op}-\ref{reduction}. 
 
\begin{proposition}\label{normal_inversion}
Given $\mathtt{S}\gg s_0$, there exist constants\index{$\Omega_\infty$} $\mu_{\mathtt{p},0}=\mu_{\mathtt{p},0}(\mathtt{p})\ge 0$, $\epsilon_0=\epsilon_0(\mathtt{pe},\mathtt{S})>0$ and  a set $\Omega_{\infty}=\Omega_\infty(i_0)\subset\Omega_1$  such that if \eqref{size_assumption_2} holds for some $\mu=\mu_{\mathtt{p},0},\mathtt{C}>0$ and $\epsilon\in (0,\epsilon_0)$, then for all $\omega\in \Omega_\infty$, $\mathcal{L}_\omega:X_{\perp}\cap X_{\mathtt{M}}\mapsto Y_{\perp}\cap X_{\mathtt{M}}$ is invertible. More precisely,  given $f(\omega)\in Y_{\perp}\cap X_{\mathtt{M}} $, there exists $w(\omega)\in X_{\perp}\cap X_{\mathtt{M}}$ such that $\mathcal{L}_\omega[w(\omega)] = f(\omega)$ for each $\omega\in \Omega_\infty$ and 
\begin{align}\nonumber
\rVert w \rVert_s^{\Lip(\gamma,\Omega_\infty)} & = \rVert \mathcal{L}_\omega^{-1}[g]\rVert_s^{\Lip(\gamma,\Omega_\infty)} \\
& \le_{\mathtt{pe},\mathtt{S}} \gamma^{-1}\left(\rVert f \rVert^{\Lip(\gamma,\Omega_\infty)}_{s+\mu_{\mathtt{p},0}} +\epsilon^{7}\gamma^{-4}\rVert \mathfrak{I}_0 \rVert^{\Lip(\gamma,\Omega_\infty)}_{s+\mu_{\mathtt{p},0}}\rVert f \rVert^{\Lip(\gamma,\Omega_\infty)}_{s_0+\mu_{\mathtt{p},0}}\right),
\label{rjjsxcxcxcxc2sdpp}
\end{align}
for all $s\in [s_0,\mathtt{S}]$. 
\end{proposition}
Since $\mu_{\mathtt{p},0}$ in the above proposition depends only on $\mathtt{p}$, we still denote it by $\mu_{0}$, which always varies from line to line.

Using the above proposition, we can solve \eqref{third_equation} by
\begin{align}\label{what}
\hat{w} = \mathcal{L}_\omega^{-1}[g_3 + \partial_x(K_{11}[\hat{\eta}])] \in X_{\perp}\cap X_{\mathtt{M}},
\end{align}
and obtain the estimates,
\begin{equation}\label{whatsize3d}
\begin{aligned}
\rVert \hat{w}\rVert^{\Lip(\gamma,\Omega_\infty)}_s& \le_{\mathtt{pe},s} \gamma^{-1}\left(\rVert g_3 \rVert^{\Lip(\gamma,\Omega_\infty)}_{s+\mu_0} +\epsilon^{7}\gamma^{-4}\rVert \mathfrak{I}_0 \rVert^{\Lip(\gamma,\Omega_\infty)}_{s+\mu_0}\rVert g_3 \rVert^{\Lip(\gamma,\Omega_\infty)}_{s_0+\mu_0}\right)\\
&\ +  \gamma^{-1}\left(\rVert K_{11}[\hat{\eta}] \rVert^{\Lip(\gamma,\Omega_\infty)}_{s+\mu_0} +\epsilon^{7}\gamma^{-4}\rVert \mathfrak{I}_0 \rVert^{\Lip(\gamma,\Omega_\infty)}_{s+\mu_0}\rVert  K_{11}[\hat{\eta}]\rVert^{\Lip(\gamma,\Omega_\infty)}_{s_0+\mu_0}\right),
\end{aligned}
\end{equation}
while Lemma~\ref{Kestimate1}  implies that
\begin{align*}
\rVert K_{11}&[\hat{\eta}] \rVert^{\Lip(\gamma,\Omega_\infty)}_{s+\mu_0} +\epsilon^{7}\gamma^{-4}\rVert \mathfrak{I}_0 \rVert^{\Lip(\gamma,\Omega_\infty)}_{s+\mu_0}\rVert  K_{11}[\hat{\eta}]\rVert^{\Lip(\gamma,\Omega_\infty)}_{s_0+\mu_0}\\
&  \le_{\mathtt{pe},s}  \epsilon^{5-2b}\rVert \hat{\eta}\rVert^{\Lip(\gamma,\Omega_\infty)}_{s+\mu_0} + \epsilon^{2b-1}\rVert \mathfrak{I}_0\rVert^{\Lip(\gamma,\Omega_\infty)}_{s+\mu_0} \rVert\hat{\eta}\rVert^{\Lip(\gamma,\Omega_\infty)}_{s_0+\mu_0}\\
& \overset{\eqref{sizeofeta}}\le_{\mathtt{pe},s}  \rVert g \rVert^{\Lip(\gamma,\Omega_\infty)}_{s+\mu_0} + \epsilon^{-1}\rVert \mathfrak{I}_0 \rVert^{\Lip(\gamma,\Omega_\infty)}_{s+\mu_0} \rVert g \rVert^{\Lip(\gamma,\Omega_\infty)}_{s_0+\mu_0} \\
& \qquad + \left(\epsilon^{5-2b}+\epsilon^{2b-1}\rVert \mathfrak{I}_0 \rVert^{\Lip(\gamma,\Omega_\infty)}_{s+\mu_0}\right)|\langle \hat{\eta}\rangle|^{\Lip(\gamma,\Omega_1)}.
\end{align*}
 Therefore,
\begin{equation}\label{wsize2s}
\begin{aligned}
\rVert \hat{w}\rVert^{\Lip(\gamma,\Omega_\infty)}_s& \le_{\mathtt{pe},s}  \gamma^{-1}\left(\rVert g \rVert^{\Lip(\gamma,\Omega_\infty)}_{s+\mu_0} +\epsilon^{7}\gamma^{-4}\rVert \mathfrak{I}_0 \rVert^{\Lip(\gamma,\Omega_\infty)}_{s+\mu_0}\rVert g \rVert^{\Lip(\gamma,\Omega_\infty)}_{s_0+\mu_0}\right)\\
& \ + \gamma^{-1}\left(\epsilon^{5-2b}+\epsilon^{2b-1}\rVert \mathfrak{I}_0 \rVert^{\Lip(\gamma,\Omega_\infty)}_{s+\mu_0} \right)|\langle \hat{\eta}\rangle|^{\Lip(\gamma,\Omega_1)}.
\end{aligned}
\end{equation}

Lastly, the first equation in \eqref{inverse_image} can be written as
\begin{align}\label{rossshsdwdsd}
\D_\omega \hat{\psi} = g_1 + \frac{1}{2\pi}(K_{20}[\hat{\eta}] + K_{11}^T[\hat{w}])\overset{Lemma~\ref{ress2}}\in X_e.
\end{align}
Using \eqref{eta_hat} and \eqref{what}, this is equivalent to
\begin{align}
2\pi \D_\omega \hat{\psi} &=2\pi g_1 + \underbrace{(K_{20} +K_{11}^{T}\mathcal{L}_\omega^{-1}\partial_x(K_{11}))[\langle \hat{\eta}\rangle]}_{=:M_1(\vphi)[\langle \hat{\eta} \rangle]} \nonumber \\
& + \underbrace{(K_{20} +K_{11}^{T}\mathcal{L}_\omega^{-1}\partial_x(K_{11}))[\D_\omega^{-1}g_2 ] }_{=M_1(\vphi)[\D_\omega^{-1} g_2]}+ \underbrace{K_{11}^T\mathcal{L}_\omega^{-1}[g_3]}_{=:M_2(\vphi)[g_3]}.\label{lastequas2}
\end{align}
To invert $\D_\omega$, we will choose $\langle \hat{\eta} \rangle$ so that the right-hand side has zero average. To do so, let us write $M_1$ in \eqref{lastequas2} as
\begin{align}\label{m1qejsd1x}
M_1(\vphi)[\langle \hat{\eta}\rangle] = (6\epsilon^{2b}\mathbb{A} +\underbrace{(K_{20}-6\epsilon^{2b}\mathbb{A} + K_{11}^{T}\mathcal{L}_\omega^{-1}\partial_x(K_{11})}_{=:M_3(\vphi)}))[\langle \hat{\eta}\rangle].
\end{align}
Using Lemma~\ref{Kestimate1}, Proposition~\ref{normal_inversion} and \eqref{size_assumption_2}, we have that 
\begin{align*}
&\rVert K_{20} - 6\epsilon^{2b}\mathbb{A} \rVert_{s}^{\Lip(\gamma,\Omega_\infty)} \le_{\mathtt{pe},s}  \epsilon^{2b+2} + \epsilon^{2b}\rVert \mathfrak{I}_0 \rVert_{s+\mu_0}^{\Lip(\gamma,\Omega_1)}\\
&\rVert K_{11}^{T}\mathcal{L}_\omega^{-1}\partial_x(K_{11})  \rVert_{s}^{\Lip(\gamma,\Omega_\infty)} \le_{\mathtt{pe},s}  \epsilon^{10-6b} + \epsilon^{4-2b}\rVert \mathfrak{I}_0  \rVert_{s+\mu_0}^{\Lip(\gamma,\Omega_\infty)} ,
\end{align*}
hence  $M_3$ in \eqref{m1qejsd1x} can be estimated as 
\begin{align}\label{m33s}
\rVert M_3 \rVert_{s}^{\Lip(\gamma,\Omega_\infty)} \le_{\mathtt{pe},s}  \epsilon^{10-6b} + \epsilon^{4-2b}\rVert \mathfrak{I}_0  \rVert_{s+\mu_0}^{\Lip(\gamma,\Omega_\infty)} .\end{align}
Plugging $s=s_0$, we have $\rVert M_3 \rVert_{s_0}^{\Lip(\gamma,\Omega_\infty)} \le\epsilon^{10-6b}$ (see \eqref{size_assumption_2}, \eqref{size_assumption_3}). 
Hence,  \eqref{rpishsnsdiwsxcsd} tells us that $
M_{1} = 6\epsilon^{2b}\left(\mathbb{A} + \frac{M_3}{\epsilon^{2b}} \right),$  is invertible and satisfies
\begin{equation}\label{m2sd211}
\begin{aligned}
&\rVert M_1h\rVert^{\Lip(\gamma,\Omega_\infty)}_{s}\le_{\mathtt{pe},s}  \epsilon^{2b}\rVert h\rVert_{s+\mu_0}^{\Lip(\gamma,\Omega_\infty)} +\epsilon^{4-2b} \rVert \mathfrak{I}_0\rVert^{\Lip(\gamma,\Omega_\infty)}_{s+\mu_0}\rVert h\rVert^{\Lip(\gamma,\Omega_\infty)}_{s_0+\mu_0} ,\\
 &\rVert M_1^{-1}\rVert^{\Lip(\gamma,\Omega_\infty)}_{s_0}\le_{\mathtt{pe}} \epsilon^{-2b}= \gamma^{-1}.
\end{aligned}
\end{equation}
Especially, $\langle M_1\rangle$, the average of $M_1$ in $\varphi$ has an inverse and 
\begin{align}\label{sds112dk2sdsd}
|\langle M_1\rangle^{-1}|^{\Lip(\gamma,\Omega_\infty)}\lepe \epsilon^{-2b}.
\end{align}
Again using Lemma~\ref{Kestimate1},  and Proposition~\ref{normal_inversion}, we have for $M_2$ in \eqref{lastequas2} that
\begin{align}\label{M23s3}
\rVert M_2[g_3]\rVert^{\Lip(\gamma,\Omega_\infty)}_{s}\le_{\mathtt{pe},s}  \epsilon^{5-4b}\rVert g\rVert^{\Lip(\gamma,\Omega_\infty)}_{s+\mu_0} + \epsilon^{-1}\rVert \mathfrak{I}_0\rVert^{\Lip(\gamma,\Omega_\infty)}_{s+\mu_0} \rVert g\rVert^{\Lip(\gamma,\Omega_\infty)}_{s_0+\mu_0},
\end{align}
especially when $s=s_0$ (thanks to \eqref{size_assumption_2}), 
\begin{align}\label{sjsd00sdsdx}
\rVert M_2[g_3]\rVert^{\Lip(\gamma,\Omega_\infty)}_{s_0}\le_{\mathtt{pe}}  \epsilon^{5-4b}\rVert g\rVert^{\Lip(\gamma,\Omega_\infty)}_{s_0+\mu_0}.
\end{align}
 Hence, we choose $\langle \hat{\eta}\rangle$ so that the right-hand side of  \eqref{lastequas2} is zero, in other words,
 \begin{align}\label{mean_hateta}
\langle \hat{\eta} \rangle := -\langle M_1\rangle^{-1} \left(\langle 2\pi g_1 \rangle + \langle M_1[\D_\omega^{-1}g_2] \rangle + \langle M_2[g_3]\rangle \right).
\end{align}
Then, it follows from  \eqref{m2sd211}, \eqref{sds112dk2sdsd}, \eqref{sjsd00sdsdx} and \eqref{domega_inverse} that
\begin{align}\label{averageparts2}
|\langle \hat{\eta}\rangle|^{\Lip(\gamma,\Omega_\infty)}\le_{\mathtt{pe}}  \gamma^{-1}\rVert g\rVert_{s_0+\mu_0}.
\end{align}
Plugging this into \eqref{wsize2s} and \eqref{sizeofeta}, we have
\begin{align}
\rVert \hat{w}\rVert^{\Lip(\gamma,\Omega_\infty)}_s& \le_{\mathtt{pe},s}  \gamma^{-1}\left(\rVert g \rVert^{\Lip(\gamma,\Omega_\infty)}_{s+\mu_0} +\epsilon^{7}\gamma^{-4}\rVert \mathfrak{I}_0 \rVert^{\Lip(\gamma,\Omega_\infty)}_{s+\mu_0}\rVert g \rVert^{\Lip(\gamma,\Omega_\infty)}_{s_0+\mu_0}\right),\label{wwsdsxxc2}\\
\rVert \hat{\eta}\rVert^{\Lip(\gamma,\Omega_1)}_s&\le_{\mathtt{pe},s} \gamma^{-1}\rVert g_2 \rVert^{\Lip(\gamma,\Omega_1)}_{s+\mu_0}.\label{wwsdsxxc22}
\end{align}

Finally, we can choose (see \eqref{lastequas2}) 
\begin{align}\label{psihat}
\hat{\psi} :=\D_\omega^{-1} \left(g_1 +\frac{1}{2\pi}\left( M_1[\langle \hat{\eta} \rangle] +M_1(\vphi)[\D_\omega^{-1}g_2] + M_2[g_3]\right)\right),
\end{align}
then it follows from \eqref{m2sd211}, \eqref{M23s3} and \eqref{domega_inverse} that
\begin{align}\label{psihasd2s}
\rVert \hat{\psi}\rVert^{\Lip(\gamma,\Omega_\infty)}_{s}\le_{\mathtt{pe},s}  \gamma^{-1}\left( \rVert g\rVert^{\Lip(\gamma,\Omega_\infty)}_{s+\mu_0} + \epsilon^{-1}\rVert \mathfrak{I}_0\rVert^{\Lip(\gamma,\Omega_\infty)}_{s+\mu_0} \rVert g\rVert^{\Lip(\gamma,\Omega_\infty)}_{s_0+\mu_0}\right),
\end{align}
Thus the inverse image of $g$ in \eqref{inverse_image} can be chosen as in  \eqref{eta_hat}, \eqref{mean_hateta}, \eqref{what} and \eqref{psihat}. The estimates for $\mathbb{D}^{-1}$ is given in the next lemma.
\begin{lemma}\label{Dinverse}\cite[Proposition 6.9]{Baldi-Berti-Montalto:KAM-quasilinear-kdv}
For all $\omega \in \Omega_\infty$ where $\Omega_\infty$ is as defined in Proposition~\ref{normal_inversion}, $\mathbb{D}(\omega)^{-1}:Y_{i}\cap X_{i,\mathtt{M}}\mapsto X_{i}\cap X_{i,\mathtt{M}}$ is well-defined. More precisely, given $g\in \Lip(\Omega_\infty,Y_{i}\cap X_{i,\mathtt{M}})$, it holds that
\begin{align*}
\rVert \mathbb{D}^{-1}[g]\rVert_{s}^{\Lip(\gamma,\Omega_\infty)} \le_{\mathtt{pe},s}  \gamma^{-1}(\rVert g \rVert_{s+\mu_0}^{\Lip(\gamma,\Omega_\infty)}+\epsilon^{7}\gamma^{-4}\rVert \mathfrak{I}_0 \rVert_{s+\mu_0}^{\Lip(\gamma,\Omega_\infty)} \rVert g\rVert_{s_0+\mu_0}^{\Lip(\gamma,\Omega_\infty)}).
\end{align*} 
\end{lemma}
\begin{proof}
The estimates for $\rVert \hat{\psi} \rVert^{\Lip(\gamma,\Omega_\infty)}_s+\rVert \hat{\eta} \rVert^{\Lip(\gamma,\Omega_\infty)}_s+\rVert \hat{w} \rVert^{\Lip(\gamma,\Omega_\infty)}_s$ follow from \eqref{wwsdsxxc2}, \eqref{wwsdsxxc22} and \eqref{psihasd2s}. To show that $\mathbb{D}^{-1}[g]\in X_i$,  we see from \eqref{eta_hat}, \eqref{what} and \eqref{rossshsdwdsd}, which give us that $\hat{\psi}$ is odd, $\hat{\eta}$ is even, and $\hat{w}\in X_{\perp}\cap X_{\mathtt{M}}$. This implies $(\hat{\psi},\hat{\eta},\hat{w})\in X_i\cap X_{i,\mathtt{M}}.$
\end{proof}
Now, we are ready to prove Proposition~\ref{approx_inverse}.
\begin{proofprop}{approx_inverse}
For $\omega\in \Omega_\infty$, which is given in Proposition~\ref{normal_inversion}, let us define an $\omega$-dependent linear map $T(i_0)=T(i_0)(\omega)$,
\begin{align}\label{tes3s}
T(i_0)[g]:=D{G}_\delta(\bar{i})\circ \mathbb{D}^{-1}\circ DG_\delta(\bar{i})^{-1}[g], \quad g\in Y_{i}.
\end{align}
Thanks to Lemma~\ref{Gsers} and Lemma~\ref{Dinverse}, we have $T(i_0)[g]\in X_i\cap X_{i,\mathtt{M}}$.  

Therefore it suffices to prove the estimates \eqref{inverse_estimate1} and \eqref{inverse_estimate2}.
Let us first derive the estimate  \ref{inverse_estimate2}. From Lemma~\ref{G_delta_estimate} and \ref{Dinverse}, direct computations give us that
\begin{align}
\rVert T(i_0)[g]\rVert^{\Lip(\gamma,\Omega_\infty)}_{s}&\le_{\mathtt{pe},s}  \gamma^{-1}\left(\rVert g \rVert^{\Lip(\gamma,\Omega_{\infty})}_{s+\mu_{0}} + \epsilon^7\gamma^{-4}\rVert \mathfrak{I}_0\rVert^{\Lip(\gamma,\Omega_\infty)}_{s+\mu_0} \rVert g\rVert^{\Lip(\gamma,\Omega_\infty)}_{s+\mu_0} \right),
\end{align}
which gives  \eqref{inverse_estimate2}.

In order to prove \eqref{inverse_estimate1}, we write
\begin{align}\label{e1}
d_{i}\mathcal{F}_\omega(i_0) = d_{i}\mathcal{F}_\omega(i_\delta) +\underbrace{d_{i}\mathcal{F}_\omega(i_0)-d_{i}\mathcal{F}_\omega(i_\delta)}_{=:\mathcal{E}_1}.
\end{align}
We further decompose the first term $d_{i}\mathcal{F}_\omega(i_\delta)$. Note that  for an embedding $i$, using \eqref{invariant_torus3}, we have
\begin{align*}
 \mathcal{F}_\omega(G_\delta(i))&= \D_\omega (G_\delta(i)) - X_{H_{\zeta}}(G_{\delta}(i)) = DG_{\delta}(i)[\D_\omega i] - X_{H_{\zeta}}(G_\delta(i)) \\
 & = DG_{\delta}(i)[\D_\omega i - DG_\delta(i)^{-1}[X_{H_{\zeta}}(G_\delta(i))]].
\end{align*}
Since $K(i) = H_{\zeta}(G_\delta(i))$ (see\eqref{modified_Hamiltonian_2}), and $G_\delta$ is symplectic, we have
\begin{align}\label{rooosdjiwojsdsx1line}
X_{K}(i) = DG_{\delta}(i)^{-1}[X_{H_{\zeta}}(G_\delta(i))].
\end{align}
Therefore, we have
\begin{align}\label{eq1}
\mathcal{F}_\omega(G_\delta(i)) = DG_{\delta}(i)[\D_\omega i - X_{K}(i)].
\end{align}
Differentiating the above equation at $(\bar{i})$ in the direction $\ihat\in X_{i}$ (see \eqref{space_reversibles} for $X_i$), we obtain
\begin{align*}
d_{i}\mathcal{F}_\omega(G_\delta(\bar{i}))[DG_\delta(\bar{i})[\ihat]]  & = D^2G_\delta(\bar{i})[\D_\omega \bar{i}-X_{K}(\bar{i}),\ihat] + DG_{\delta}(\bar{i})\circ L[\ihat] \\
& = D^2G_{\delta}(\bar{i})[DG_\delta^{-1}[\mathcal{F}_\omega(i_\delta)],\ihat] + DG_{\delta}(\bar{i})\circ L[\ihat],
\end{align*}
where $L$ is as in \eqref{linearized_XK} and the last equality follows from \eqref{eq1} with $G_\delta(\bar{i})=i_\delta$ (see \eqref{def_G_delta}). Replacing $\ihat$ by $DG_\delta(\bar{i})^{-1}(\ihat)$, we obtain
\begin{align}
d_{i}\mathcal{F}_\omega(i_\delta)[\ihat] & = \underbrace{D^2G_\delta(\bar{i})[DG_\delta^{-1}[\mathcal{F}_\omega(i_\delta)],DG_\delta^{-1}[\ihat]]}_{=:\mathcal{E}_2[\ihat]} + DG_\delta(\bar{i})\circ L[DG_\delta(\bar{i})^{-1}[\ihat]]\label{e21s2x3}\\
& = \mathcal{E}_2[\ihat] + DG_\delta(\bar{i})\circ L \circ D{G}_\delta(\bar{i})^{-1}[\ihat]. \label{e2}
\end{align}
 Therefore, it follows from \eqref{e1} and \eqref{e2} that 
\begin{align}\label{e3ishere11}
d_{i}\mathcal{F}_\omega(i_0) =\mathcal{E}_1 +\mathcal{E}_2 + \underbrace{DG_\delta(\bar{i})\circ (L-\mathbb{D}) \circ D{G}_\delta(\bar{i})^{-1}}_{=:\mathcal{E}_3} + DG_\delta(\bar{i})\circ \mathbb{D}\circ D{G}_\delta(\bar{i})^{-1}.
\end{align}
Therefore, it follows from the definition of $T(i_0)$ in \eqref{tes3s} that 
\begin{align}\label{error_1}
d_{i}\mathcal{F}_\omega(i_0)\circ T(i_0) - I = (\mathcal{E}_1 +\mathcal{E}_2+\mathcal{E}_3)\circ T(i_0).
\end{align}
Now we estimate the size of $\mathcal{E}_i$ for $i=1,2,3$ . 
For $\mathcal{E}_1$, it follows from the definitions in $\mathcal{F}_\omega$, $H_\zeta$ and $i_\delta$ in \eqref{invariant_torus3}, \eqref{normal_formaa2}  and \eqref{iso_tropic} that
\begin{align*}
(d_{i}\mathcal{F}_\omega(i_0) - d_{i}\mathcal{F}_\omega(i_\delta))[\ihat] &= (d_iX_P(i_0) - d_iX_P(i_\delta)[\ihat] \\
&= \int_{0}^{1}\frac{d}{dt}(d_iX_P(t(i_0 -i_\delta)+i_\delta)[\ihat])dt\\
& = \int_0^{1}\partial_yd_iX_P(t(i_0-i_\delta)+i_\delta)[\ihat]\cdot (y-y_\delta)dt.
\end{align*}
Therefore, using \eqref{interpolation_2s}, we obtain
\begin{align}
& \rVert \mathcal{E}_1[\ihat]\rVert_{s}^{\Lip(\gamma,\Omega_\infty)} \nonumber\\
& \le_{\mathtt{pe},s}  
\sup_{t\in [0,1]}\rVert \partial_yd_iX_P(t(i_0-i_\delta)+i_\delta)[\ihat]\rVert_{s}^{\Lip(\gamma,\Omega_\infty)} \rVert y-y_\delta\rVert_{s_0}^{\Lip(\gamma,\Omega_\infty)} \nonumber\\
& \ +\sup_{t\in [0,1]}\rVert \partial_yd_iX_P(t(i_0-i_\delta)+i_\delta)[\ihat]\rVert_{s_0}^{\Lip(\gamma,\Omega_\infty)} \rVert y-y_\delta\rVert_{s}^{\Lip(\gamma,\Omega_\infty)}.\label{e1estimate}
\end{align}
Note that 
\begin{align}\label{rintertsd}
t(i_0-i_\delta)+i_\delta = \bar{i} + \mathfrak{I}_0 + (t-1)(i_0-i_\delta) =:\bar{i} + \mathfrak{I}_t.
\end{align}
 For $\mathfrak{I}_t$ as defined above, it follows from \eqref{y_estimate} that
\begin{align*}
\rVert \mathfrak{I}_t\rVert_{s}^{\Lip(\gamma,\Omega_\infty)} & \le_{\mathtt{pe},s}  \rVert \mathfrak{I}_0\rVert_{s}^{\Lip(\gamma,\Omega_\infty)} + \gamma^{-1}(\rVert Z \rVert_{s+\mu_0}^{\Lip(\gamma,\Omega_\infty)} + \rVert Z \rVert_{s_0+\mu_0}^{\Lip(\gamma,\Omega_\infty)}\rVert \mathfrak{I}_0\rVert_{s+\mu_0}^{\Lip(\gamma,\Omega_\infty)})\\
& \le_{\mathtt{pe},s}  \rVert \mathfrak{I}_0\rVert_{s+\mu_0}^{\Lip(\gamma,\Omega_\infty)}  + \gamma^{-1}\rVert Z\rVert_{s+\mu_0}^{\Lip(\gamma,\Omega_\infty)},
\end{align*}
where the last inequality follows from \eqref{size_assumption_2}.
Thus, it follows from \eqref{tame5} and \eqref{inverse_estimate2} that
\begin{align*}
&\rVert \partial_y d_i X_P(\bar{i}+\mathfrak{I}_t)[\ihat] \rVert_{s}^{\Lip(\gamma,\Omega_\infty)}\\
& \le_{\mathtt{pe},s}  \epsilon^{2b-1}(\rVert \ihat\rVert_{s+\mu_0}^{\Lip(\gamma,\Omega_\infty)} + \rVert \mathfrak{I}_t\rVert_{s+\mu_0}^{\Lip(\gamma,\Omega_\infty)}\rVert \ihat\rVert_{s_0+\mu_0}^{\Lip(\gamma,\Omega_\infty)})  \\
& \le_{\mathtt{pe},s}  \epsilon^{2b-1} (\rVert \ihat\rVert_{s+\mu_0}^{\Lip(\gamma,\Omega_\infty)} + \left(\rVert \mathfrak{I}_0\rVert_{s+\mu_0}^{\Lip(\gamma,\Omega_\infty)} + \gamma^{-1}\rVert Z\rVert_{s+\mu_0}^{\Lip(\gamma,\Omega_\infty)}\right)\rVert \ihat\rVert_{s_0+\mu_0}^{\Lip(\gamma,\Omega_\infty)}),
\end{align*}
and
\begin{align*}
\rVert \partial_y d_i X_P(\bar{i}+\mathfrak{I}_t)[\ihat] \rVert_{s_0}^{\Lip(\gamma,\Omega_\infty)} \le_{\mathtt{pe}}  \epsilon^{2b-1}\rVert  \ihat\rVert_{s_0 +\mu_0}^{\Lip(\gamma,\Omega_\infty)},
\end{align*}
where  we used \eqref{size_assumption_2} in the last inequality. Hence, plugging this and \eqref{y_estimate} with \eqref{rintertsd} into \eqref{e1estimate}, we obtain
\begin{align}\label{e1_estimate}
\rVert \mathcal{E}_1&[\ihat]\rVert_{s}^{\Lip(\gamma,\Omega_\infty)} \nonumber\\
& \le_{\mathtt{pe},s}  \epsilon^{2b-1}\gamma^{-1}\left(\rVert Z\rVert_{s_0+\mu_0}^{\Lip(\gamma,\Omega_\infty)}\rVert \ihat\rVert_{s+\mu_0}^{\Lip(\gamma,\Omega_\infty)}\right. \nonumber \\
& \left.\qquad + (\rVert Z\rVert_{s+\mu_0}^{\Lip(\gamma,\Omega_\infty)} + \rVert Z\rVert_{s_0+\mu_0}^{\Lip(\gamma,\Omega_\infty)}\rVert \mathfrak{I}_0\rVert_{s+\mu_0}^{\Lip(\gamma,\Omega_\infty)})\rVert \ihat\rVert_{s_0+\mu_0}^{\Lip(\gamma,\Omega_\infty)}\right).
\end{align}
For $\mathcal{E}_2$ in \eqref{e21s2x3}, it is straightforward from Lemma~\ref{G_delta_estimate} that \begin{align}
\rVert \mathcal{E}_2[\ihat]\rVert_s^{\Lip(\gamma,\Omega_\infty)} &\le_{\mathtt{pe},s}  \rVert Z\rVert_{s_0+\mu_0}^{\Lip(\gamma,\Omega_\infty)}\rVert \ihat\rVert_{s+\mu_0}^{\Lip(\gamma,\Omega_\infty)} \nonumber \\
& + (\rVert Z\rVert_{s+\mu_0}^{\Lip(\gamma,\Omega_\infty)} + \rVert Z\rVert_{s_0+\mu_0}^{\Lip(\gamma,\Omega_\infty)}\rVert \mathfrak{I}_0\rVert_{s+\mu_0}^{\Lip(\gamma,\Omega_\infty)})\rVert \ihat\rVert_{s_0+\mu_0}^{\Lip(\gamma,\Omega_\infty)}.\label{D2G_estimate}
\end{align}
Hence, it follows  that $\mathcal{E}_2$ satisfies \eqref{e1_estimate} (even without the coefficient $\epsilon^{2b-1}\gamma^{-1}\gg 1$). Similarly, $\mathcal{E}_3$ in \eqref{e3ishere11} can be estimated, using Lemma~\ref{G_delta_estimate} and \eqref{difference_2}, as
\begin{align*}
\rVert \mathcal{E}_3[\ihat]\rVert_s^{\Lip(\gamma,\Omega_\infty)}& \le_{\mathtt{pe},s}  \rVert Z\rVert_{s_0+\mu_0}^{\Lip(\gamma,\Omega_\infty)}\rVert \ihat\rVert_{s+\mu_0}^{\Lip(\gamma,\Omega_\infty)} \\
& + (\rVert Z\rVert_{s+\mu_0}^{\Lip(\gamma,\Omega_\infty)} + \epsilon^{-1}\rVert Z\rVert_{s_0+\mu_0}^{\Lip(\gamma,\Omega_\infty)}\rVert \mathfrak{I}_0\rVert_{s+\mu_0}^{\Lip(\gamma,\Omega_\infty)})\rVert \ihat\rVert_{s_0+\mu_0}^{\Lip(\gamma,\Omega_\infty)},
\end{align*}
which satisfies the same estimate as in \eqref{e1_estimate}. Hence, $\mathcal{E}:=\mathcal{E}_1+\mathcal{E}_2+\mathcal{E}_3$ satisfies the estimate in \eqref{e1_estimate}, that is,
\begin{align}
\rVert \mathcal{E}&[\ihat]\rVert_{s}^{\Lip(\gamma,\Omega_\infty)} \nonumber\\
&\le_{\mathtt{pe},s}  \epsilon^{2b-1}\gamma^{-1}\left(\rVert Z\rVert_{s_0+\mu_0}^{\Lip(\gamma,\Omega_\infty)}\rVert \ihat\rVert_{s+\mu_0}^{\Lip(\gamma,\Omega_\infty)}\right. \nonumber \\
& \left.\qquad \qquad \qquad+ (\rVert Z\rVert_{s+\mu_0}^{\Lip(\gamma,\Omega_\infty)} + \rVert Z\rVert_{s_0+\mu_0}^{\Lip(\gamma,\Omega_\infty)}\rVert \mathfrak{I}_0\rVert_{s+\mu_0}^{\Lip(\gamma,\Omega_\infty)})\rVert \ihat\rVert_{s_0+\mu_0}^{\Lip(\gamma,\Omega_\infty)}\right).\label{esti2m2ate_again_freee1ints1}
\end{align}
Hence, \eqref{inverse_estimate1} follows from \eqref{inverse_estimate2} and \eqref{error_1}.
\end{proofprop}

\begin{remark}\label{labesdof_sdxxcxc}
As mentioned, the proof of Proposition~\ref{normal_inversion} will be achieved throughout Chapter~\ref{Linear_op}-\ref{reduction}. The proof of Proposition~\ref{approx_inverse} in this section tells us that if Proposition~\ref{normal_inversion} holds for some $\mu_{\mathtt{p},0}\ge 0$, then Proposition~\ref{approx_inverse} holds for some $\mu_{\mathtt{p}}>0$, which is possibly larger than $\mu_{\mathtt{p},0}$. In other words, once we fix $\mu_{\mathtt{p},0}$, then $\mu_{\mathtt{p}}$ can be fixed, depending on $\mu_{\mathtt{p},0}$ and $\mathtt{p}$. Since both of $\mu_{\mathtt{p}},\mu_{\mathtt{p,0}}$ depend only on $\mathtt{p}$, we see that there exists a  constant $\mu_{\mathtt{p},2}=\mu_{\mathtt{p},2}(\mu_{\mathtt{p},0},\mathtt{p})>0$  such that 
\begin{align}\label{mu0andmu1}
\mu_{\mathtt{p},0} < \mu_{\mathtt{p}} <  \mu_{\mathtt{p},0} + \mu_{\mathtt{p},2}.
\end{align}
\end{remark}

\section{Linearized system at an invariant torus}
Before we close this chapter, we digress briefly to study the linearized Hamiltonian system at $i_0$, assuming that $i_0$ is a solution, that is $\mathcal{F}_\omega(i_0(\varphi))=0$.
If $i_0$ is an invariant torus\index{invariant torus}, that is, $\mathcal{F}_\omega(i_0)=0$, then the linearized Hamiltonian system at $i_0(\omega t)$ can be conjugated to the linear system with the vector field $d_i X_K(\bar{i})$. To see this more precisely, we see from \eqref{rooosdjiwojsdsx1line}  that 
\begin{align}\label{rjjsdsdww2}
DG_\delta(i) X_K(i) = X_{H_\zeta}(G_\delta(i)).
\end{align} Therefore, differentiating it in the direction $\ihat$ at the trivial embedding $\bar{i}$,  we have
\begin{align}\label{porisdsd2s}
D^2G_\delta(\bar{i})[X_K(\bar{i}),\ihat] + DG_\delta(\bar{i})\circ d_i X_K(\bar{i})[\ihat] = d_i X_{H_\zeta}(G_\delta(\bar{i}))\circ DG_\delta(\bar{i})[\ihat].
\end{align}
Furthermore, if $i_0$ is invariant, then using $i_\delta\overset{\eqref{y_estimate}}=i_0$ and $G_\delta(\bar{i})=i_\delta$, which follows from \eqref{def_G_delta}, we see that
\begin{align}
d_i X_{H_\zeta}(i_0)\circ DG_\delta(\bar{i})[\ihat] &= DG_\delta(\bar{i})\circ d_i X_K(\bar{i})[\ihat] + D^2G_\delta(\bar{i})[X_K(\bar{i}),\ihat].\label{wlfkfahsxh1}
\end{align}
Also, differentiating $G_\delta(\bar{i}(\varphi))=i_0(\varphi)$ in $\varphi$, we see that
\begin{align}\label{wlfjsd2sd}
DG_\delta(\bar{i}(\varphi))[(\omega,0,0)] = \D_\omega i_0(\varphi),\text{ for all $\varphi\in\mathbb{T}^\nu$.}
\end{align}
Now, consider the linearized system
\begin{align}\label{linearsdsds}
\dot{I}(t):=(\dot{\theta}(t),\dot{y}(t),\dot{z}(t)) = d_i X_{H_\zeta}(i_0(\omega t))[{\theta}(t),{y}(t),z(t)]=d_i X_{H_\zeta}(i_0(\omega t))[I(t)],
\end{align}
which is the linearized Hamiltonian system of $H_{\zeta}$ at the quasiperiodic solution $i_0(\omega t)$.
Then,  defining $A(t):=DG_\delta(\bar{i}(\omega t))^{-1}[I(t)]=(DG_\delta(\omega t,0,0))^{-1}[I(t)]$, we see that 
\begin{align*}
\frac{d}{dt}\left( DG_\delta(\bar{i}(\omega t))[{A}(t)] \right) &= D^2G_\delta(\bar{i}(\omega t))[(\omega,0,0),A(t)]  + DG_\delta(\bar{i}(\omega t))[\dot{A}(t)] \\
& = \dot{I}(t) = d_i X_{H_\zeta}(i_0(\omega t))[I(t)]\\
& = d_i X_{H_\zeta}(i_0(\omega t))\circ DG_\delta (\bar{i}(\omega t))[A(t)],
\end{align*}
which shows that
\begin{align}\nonumber
DG_\delta(\bar{i}(\omega t))[\dot{A}(t)] & =  d_i X_{H_\zeta}(i_0(\omega t))\circ DG_\delta (\bar{i}(\omega t))[A(t)] \\
& \qquad - D^2G_\delta(\bar{i}(\omega t))[(\omega,0,0),A(t)].\label{linear11sds}
\end{align}
Plugging \eqref{wlfjsd2sd} into \eqref{linear11sds}, we see that
\begin{align*}
DG_\delta(\bar{i}(\omega t))[\dot{A}(t)]  &=  DG_\delta(\bar{i}(\omega t))\circ d_i X_{K}(\bar{i}(\omega t))[A(t)] \\
& \ + D^2G_\delta(\bar{i}(\omega t))[X_K(\bar{i}(\omega t)),A(t)] \\
& \ -D^2G_\delta(\bar{i}(\omega t))[(\omega,0,0),A(t)],
\end{align*}
while we have
\begin{align*}
X_K(\bar{i}(\omega t)) - (\omega, 0 ,0) &= (DG_\delta(\bar{i}(\omega t)))^{-1}\circ (DG_\delta(\bar{i}(\omega t))\circ X_K(\bar{i}(\omega t)) \\
& -DG_\delta(\bar{i}(\omega t))[\omega,0,0] )\\
&\overset{\eqref{rjjsdsdww2},\eqref{wlfkfahsxh1}}= X_{H_\zeta}(G_\delta(\bar{i}(\omega t))) - D_\omega i_0(\omega t)\\
& \overset{G_\delta(\bar{i}) = i_0}= X_{H_\zeta}(i_0(\omega t)) - \D_\omega i_0(\omega t)\\
&=0.
\end{align*}
Hence, combining this with \eqref{linear11sds}, we see that $DG_\delta(\bar{i}(\omega t))[\dot{A}(t)] =  DG_\delta(\bar{i}(\omega t))\circ d_i X_{K}(\bar{i}(\omega t))[A(t)]$, that is, $
\dot{A}(t) = d_i X_{K}(\bar{i}(\omega t))[A(t)]$. Therefore, using the expression of $d_i X_K$ in \eqref{rjoosdwiwnds1s2s2} and using \eqref{difference_1}, we obtain the following:
\begin{lemma}\label{linearlsnsdosd}
Suppose $i_0$ satisfies $\mathcal{F}_\omega(i_0) = 0$ and consider the linear equation $I(t)$ given by $
\dot{I}(t) = d_i X_{H_\zeta}(i_0(\omega t))[I(t)].$
Then $A(t):=DG_\delta(\bar{i}(\omega t))^{-1}[I(t)]$ satisfies 
\begin{align}\label{rqqwsjspwsd2}
\dot{A}(t) = \colvec{ 0 & K_{20}(\omega t) & K_{11}(\omega t)^{T} \\ 0 & 0 & 0 \\ 0 & K_{11}(\omega t) & \partial_x\circ K_{02}(\omega t)} A(t).
\end{align}
\end{lemma}
\color{black}

\chapter{Linearized operator in the normal directions}\label{Linear_op}

In Chapter~\ref{Linear_op} and \ref{reduction}, we aim to prove Proposition~\ref{normal_inversion}. In this chapter, we will derive more explicit formula for the operator $\mathcal{L}_\omega=\D_\omega - \partial_xK_{02}(\psi)$\index{$\mathcal{L}_\omega$} in \eqref{third_equation}. More precisely, in view of \eqref{rksksd1kpsxx},  we need to find a linear operator $K_{02}(\psi):H_{S^\perp}\mapsto H_{S^\perp}$ such that
\begin{align}\label{K_{20}_1}
(K_{02}(\psi)[w],w)_{L^2(\mathbb{T})} = \frac{d^2}{dt^2}K(\psi,0,tw)|_{t=s=0}, \text{ for all }w\in H_{S^\perp}.
\end{align}
From now on, we use $\varphi$, instead of $\psi$, for the "time" variable.
From \eqref{modified_Hamiltonian_2}, \eqref{rescaled_hamiltonian} and Proposition~\ref{normal_form_prop11231}, we can easily see that for an embedding $i(\varphi)$, 
\begin{align}\nonumber
K(i(\varphi)) & =H_\zeta(G_\delta(i(\varphi)))=\epsilon^{-2b} H\circ U_\zeta\circ G_\delta(i(\varphi)) \\
& = \epsilon^{-2b}\mathcal{H}\circ \Phi^{WB}\circ U_\zeta\circ G_\delta(i(\varphi)).\label{whatklooklike}
\end{align}

We split this chapter into several sections. In Section~\ref{hohohohoho22}, we study some useful properties of the coefficients generated by $\Phi^{WB}\circ U_\zeta(i_\delta)$, at which we linearize the gSQG equation.  From \ref{item3_normal} in  Proposition~\ref{normal_form_prop11231}, we already know that  $\Phi^{WB}$ changes only a finite number of Fourier modes. This will generate smoothing operators in the expression of $\nabla_z^2 K$, which will be studied in Section~\ref{Finite dimensional operators}. In Section~\ref{linisdw}, we derive the explicit form of the Hessian $\nabla_z^2 K$.

\section{Homogeneous expansion}\label{hohohohoho22}
In this section, we investigate the structure of the coefficients and symbols arising in the linearized operator.

 To begin with, let us denote (recall $U_\zeta$ from \eqref{actionssx})\index{$\bar{v}$},
\begin{align}
u&:=u(\varphi):=\Phi^{WB}(U_\zeta(i_\delta(\varphi))),\label{u_def1}\\
\bar{v}&:=\bar{v}(\varphi,x):=\sum_{j_k\in S}\sqrt{j_{k}\zeta_{k}}e^{\ii \mathtt{l}(j_k)\cdot \varphi}e^{\ii j_k x}=:\sum_{j_k\in S}\bar{v}_{j_k}(\varphi,x),\label{norm_def_vbar}
\end{align}
where $\mathtt{l}:S\mapsto \mathbb{Z}^{\nu}$ such that \index{$\mathtt{l}$}
\begin{align}\label{def_lll}
\mathtt{l}(j_k) := \mathtt{e}_{k}, \text{ and }\mathtt{l}(-j_k):=-\mathtt{l}(j_k) \text { for $j_k\in S^{+}=\left\{ j_1,\ldots,j_{\nu}\right\}$},
\end{align}
and $\mathtt{e}_k:=\left(0,\ldots,1,\ldots,0\right)$ denotes\index{$\mathtt{e}_k$} the $k$-th vector of the canonical basis of $\mathbb{R}^\nu$.
We also denote the average operators:
\begin{align}\label{mean_def_52}
M_x(q):=\frac{1}{2\pi}\int_{\mathbb{T}}q(\omega,\varphi,x)dx, \quad M_{\varphi,x}(q):=\frac{1}{(2\pi)^{\nu}}\int_{\mathbb{T}^\nu}M_x(q)(\omega,\varphi)d\varphi.
\end{align}

\begin{lemma}\label{norm_phi_expansion1}
Let $u$ be as in \eqref{u_def1}. Then $u$ has the expansion:
\begin{align}\label{norm_phi_expansion2}
u = \epsilon \bar{v} +3\epsilon^2\partial_x\Pi_{S^\perp}K_2(\bar{v},\bar{v}) + q,
\end{align}
where $K_2$ is as in \eqref{K12def}. And we have that
\begin{align}
\rVert q \rVert^{\Lip(\gamma,\Omega_1)}_{s} & \le_{\mathtt{pe},s} \epsilon^{3} + \epsilon\rVert \mathfrak{I}_{\delta}\rVert^{\Lip(\gamma,\Omega_1)}_{s+\mu_0}, \nonumber\\
\rVert d_i q(i_0)[\ihat] \rVert_s & \le_{\mathtt{pe},s} \epsilon\left( \rVert \ihat \rVert_{s+\mu_0} + \rVert \mathfrak{I}_\delta \rVert_{s+\mu_0}\rVert \ihat \rVert_{s_0+\mu_0}\right). \label{jjsdsd2sdsz}
\end{align}
Furthermore, it holds that $M_x(q)=0$.
\end{lemma}
\begin{proof}
From Lemma~\ref{WB_expansion11},  \eqref{u_def1}, \eqref{actionssx} and \eqref{iso_tropic} , we have that
\begin{align*}
u & = \Phi^{WB}(\epsilon v_\epsilon(\theta_0,y_\delta) + \epsilon^bz_0) \\
& = \epsilon v_\epsilon + \epsilon^b z_0 + 6\epsilon^{1+b}\partial_x\Pi_SK_1( v_\epsilon,z_0) + 3\epsilon^2\partial_x\Pi_{S^\perp}K_2(v_\epsilon,v_\epsilon) + \Phi^{WB}_{\ge 3}(U_\zeta(i_\delta)).
\end{align*}
Let $q=q(\varphi,x)$ be
\begin{align}
q & :=\epsilon\left(v_\epsilon - \bar{v} \right) + 3\epsilon^2\partial_x\Pi_{S^\perp}\left( K_2(v_\epsilon,v_\epsilon) -K_2(\bar{v},\bar{v}) \right) \nonumber \\
& \qquad + \epsilon^bz_0 + 6\epsilon^{1+b}\partial_x\Pi_SK_1( v_\epsilon,z_0) + \Phi^{WB}_{\ge3}(U_\zeta(i_\delta)). \label{qes2rlaqkq}
\end{align}
Clearly, we have the expression in \eqref{norm_phi_expansion2}. Noticing that $|v_\epsilon(\theta_0,y_\delta)-\bar{v}| = O(\mathfrak{I}_\delta)$, more precisely (using \eqref{actionssx}),
\begin{equation}\label{vvvbbb}
\begin{aligned}
 \rVert v_\epsilon - \bar{v}\rVert^{\Lip(\gamma,\Omega_1)}_{s}&\le_{\mathtt{pe},s} \rVert \mathfrak{I}_\delta\rVert^{\Lip(\gamma,\Omega_1)}_s,\\
  \rVert d_i(v_\epsilon - \bar{v})(i_0)[\ihat]\rVert_{s}^{\Lip(\gamma,\Omega_1)}&\le_{\mathtt{pe},s} \rVert\ihat\rVert_{s}^{\Lip(\gamma,\Omega_1)}+\rVert\mathfrak{I}_\delta\rVert_{s}^{\Lip(\gamma,\Omega_1)}\rVert\ihat\rVert_{s_0}^{\Lip(\gamma,\Omega_1)},
\end{aligned}
\end{equation} and  using \eqref{size_assumption_3}, the estimate \eqref{jjsdsd2sdsz} follows straightforwardly (for example, in the estimate for $q$ in \eqref{jjsdsd2sdsz}, $\epsilon^3$ comes from $\Phi^{WB}_{\ge3}(U_\zeta(i_\delta))$ in \eqref{qes2rlaqkq} and $\epsilon\rVert \mathfrak{I}_{\delta}\rVert^{\Lip(\gamma,\Omega_1)}_{s+\mu_0}$ comes from $\epsilon\left(v_\epsilon - \bar{v} \right)$. All the other terms in \eqref{qes2rlaqkq} have smaller contribution in the estimates). Since $\Phi^{WB}(f)\in L^2_0$, for $f\in L^2_0$, we have that $M_x(u)=0$. Thus \eqref{norm_phi_expansion2} implies  $M_x(q)=0$, since $\overline{v}$ and $\partial_x\Pi_{S^\perp}K_2(\bar{v},\bar{v})$ have zero average.
\end{proof}

\subsection{Homogeneous expansion of the coefficients: Class $\mathfrak{C}_1$ and $\mathfrak{C}_2$}\label{hohohoss22s}
Now, we will define some classes of functions depending on the embedding $i_0$. The motivation of such classes is that Lemma~\ref{norm_phi_expansion1} tells us that $u(x),u(x)-u(y)$, where $u$ are in \eqref{actionssx}, is contained in those classes (see Lemma~\ref{properties_weights}).

\begin{definition}\label{homogeneous_expansion1}
We say that a function $a=a(\omega,\varphi,x)\in \mathfrak{C}_1(i_0)$, if  $a$ admits an expansion\index{$\mathfrak{C}_{1}$}:
\begin{align}\label{aaajsjdw}
a(\omega,\varphi,x)=\epsilon a_1(\omega,\varphi,x)+\epsilon^2a_2(\omega,\varphi,x) + \epsilon^3a_3(\omega,\varphi,x)+q(\omega,\varphi,x),
\end{align}
where\index{$M_x$}\index{$M_{\varphi,x}$}
\begin{enumerate}[label=(\arabic*)]
\item\label{asdsdelgisd} Each $a_i$ is of the form:
\begin{align}\label{homcoeeff}
a_i(\omega,\varphi,x)=\sum_{j_{k_1},\ldots,j_{k_i}\in S}C_{j_{k_1},\ldots,j_{k_i}}\bar{v}_{j_1}(\varphi,x)\cdots \bar{v}_{j_i}(\varphi,x),\text{ for $i=1,2,3$},
\end{align}
for some constants $C_{j_{k_1},\ldots,j_{k_i}}$, which do not depend on any of $i_0,\omega, \varphi,x,\epsilon$ (note that $\bar{v}$ depends on $\zeta$ (see \eqref{norm_def_vbar}) and thus $\omega$ by \eqref{xi_omega_dependent}).
\item \label{qqqsds22}$q$ in \eqref{aaajsjdw} satisfies that, assuming \eqref{size_assumption_3} for some $\mu_0=\mu_0(\mathtt{p})$, 
\begin{align}
\rVert q \rVert^{\Lip(\gamma,\Omega_1)}_{s} & \le_{\mathtt{pe},s} \epsilon^{3} + \epsilon\rVert \mathfrak{I}_{\delta}\rVert^{\Lip(\gamma,\Omega_1)}_{s+\mu_0}, \nonumber \\
\rVert d_i q(i)[\ihat] \rVert & \le_{\mathtt{pe},s} \epsilon\left( \rVert \ihat \rVert_{s+\mu_0} + \rVert \mathfrak{I}_\delta \rVert_{s+\mu_0}\rVert \ihat \rVert_{s_0+\mu_0}\right). \label{q_estimate_11}
\end{align}
Furthermore, $M_x(q)-M_{\varphi,x}(q)$ satisfies
\begin{equation}\label{q_estimate_112}
\begin{aligned}
&\rVert M_x(q)-M_{\varphi,x}(q)\rVert^{\Lip(\gamma,\Omega)}_{s}\le_{\mathtt{pe},s} \epsilon^4 + \epsilon^{2}\rVert\mathfrak{I}_\delta\rVert_{s+\mu_0}^{\Lip(\gamma,\Omega_1)},\\
& \rVert d_i(M_x(q)-M_{\varphi,x}(q)(i_0))[\ihat]\rVert_s\le_{\mathtt{pe},s} \epsilon^2\left(\rVert \ihat \rVert_{s+\mu_0} + \rVert \mathfrak{I}_\delta \rVert_{s+\mu_0}\rVert \ihat \rVert_{s_0+\mu_0} \right).
\end{aligned}
\end{equation}
\end{enumerate}

We also say that a function $a=a(\omega,\varphi,x,y)\in \mathfrak{C}_2(i_0)$, if $a$ admits an expansion\index{$\mathfrak{C}_2$}:
\begin{align}\label{aaajsjdw2}
a(\omega,\varphi,x,y)=\epsilon a_1(\omega,\varphi,x,y)+\epsilon^2a_2(\omega,\varphi,x,y) + \epsilon^3a_3(\omega,\varphi,x,y)+q(\omega,\varphi,x,y),
\end{align}
where
\begin{enumerate}[label=(\arabic*)]
\item \label{2tsksdysxxx}Each $a_i$ is of the form:
\begin{equation}\label{homcoeeff2}
\begin{aligned}
a_i(\omega,\varphi,x,y)=\sum_{j_{k_1},\ldots,j_{k_i}\in S}C_{j_{k_1},\ldots,j_{k_i}}(x-y)\bar{v}_{j_1}(\varphi,x)\cdots \bar{v}_{j_i}(\varphi,x),
\end{aligned}
\end{equation}
\text{ for $i=1,2,3$}, for some smooth functions $C_{j_{k_1},\ldots,j_{k_i}}(\cdot)$, which do not depend on any of $i_0,\omega, \varphi,\epsilon$, but $(x-y)$.

\item  \label{qsdwd}$q$ in \eqref{aaajsjdw2} satisfies the same estimates as in \eqref{q_estimate_11} (where the norm, $\rVert \cdot \rVert_{s}$, is replaced by $\rVert \cdot \rVert_{H^{s}(\mathbb{T}^{\nu+2})}$). Furthermore, denoting $Q:=q(\omega,\varphi,x,x-y)$, $M_{x}(Q)-M_{\varphi,x}(Q)$ satisfies
\begin{equation}\label{q_estimate461}
\begin{aligned}
&\rVert M_{x}(Q)-M_{\varphi,x}(Q)\rVert_{s}^{\Lip(\gamma,\Omega_1)}\le_{\mathtt{pe},s} \epsilon^4 + \epsilon^2\rVert\mathfrak{I}_\delta\rVert^{\Lip(\gamma,\Omega_1)}_{s+\mu_0},\\
&\rVert d_i(M_{x}(Q)-M_{\varphi,x}(Q))(i_0)[\ihat]\rVert_s\le_{\mathtt{pe},s} \epsilon^2\left(\rVert \ihat \rVert_{s+\mu_0} + \rVert \mathfrak{I}_\delta \rVert_{s+\mu_0}\rVert \ihat \rVert_{s_0+\mu_0} \right).
\end{aligned}
\end{equation} 
\end{enumerate}
\end{definition}
In the following lemma, we list some properties of the functions in $\mathfrak{C}_1(i_0),\mathfrak{C}_2(i_0)$. The proof is straightforward from its definition, therefore we omit it for simplicity.
\begin{lemma}\label{properties_weights}
We have:
 \begin{enumerate}[label=(\arabic*)]
 \item \label{slsd22sd2sdhw}If $a\in \mathfrak{C}_n(i_0)$, for $n=1,2$, it holds that
 \[
 \rVert a \rVert^{\Lip(\gamma,\Omega)}_{s}\le_{\mathtt{pe},s} \epsilon\left(1 + \rVert \mathfrak{I}_\delta\rVert^{\Lip(\gamma,\Omega)}_{s+\mu_0} \right).
 \]
 \item \label{item_homog_weight3}$u \in \mathfrak{C}_1(i_0)$ where $u$ is as in \eqref{u_def1}.
 \item \label{item_homog_weight4} $\mathfrak{C}_1(i_0),\mathfrak{C}_2(i_0)$ are closed under addition, multiplication, differentiation in $x,y,\varphi$.
  \item \label{item_homog_weight42}Denoting $u_1(x)=u_1(\omega,\varphi,x),\ u_2(y)=u_2(\omega,\varphi,y)\in \mathfrak{C}_1(i_0)$, $u_1(x)u_2(y)$, then $u_1(x)+u_2(y),\ J(u_1)(x,y)\in \mathfrak{C}_2(i_0)$, where $J(u_1)(x,y) = \frac{(u_1(x)-u_1(y))^2}{2-2\cos(x-y)}$.
 \item Let  $F:\mathbb{R}^3\mapsto \mathbb{R}$ be an analytic function such that $F(0,0,0)=0$.  We have that $F(u(x),u(y),J(u)(x,y))\in \mathfrak{C}_2(i_0)$, .\label{item_homog_weight5}
 \item If  $u\in \mathfrak{C}_1(i_0)$ and $a\in \mathfrak{C}_2(i_0)$, then
 $\Lambda^{\alpha-1}_{a}u,\Upsilon^{\alpha-3}_au\in \mathfrak{C}_1(i_0)$ (see \eqref{fractional_1} for the definitions of $\Lambda^{\alpha-1}_a$ and $\Upsilon^{\alpha-3}_a$.).\label{item_homog_weight6}
 \end{enumerate}
 \end{lemma}

 \begin{lemma}\label{meanisbetter}
 Let $a\in \mathfrak{C}_1(i_0)$. Then, $M_x(a)$ has a decomposition:
 \begin{align*}
 M_x(a)(\omega,\varphi)=\epsilon^2\mathtt{m}(\omega) + r(\omega) +\tilde{ q}(\omega,\varphi),\quad \int_{\mathbb{T}^\nu} \tilde{q}(\omega,\varphi)d\varphi = 0,
 \end{align*}
 such that $\mathtt{m}(\omega),r(\omega)$ do not depend on $\varphi,x$, and  $\mathtt{m}$ does not even depend on $i_0$ and 
 \begin{align}
 |\mathtt{m}|^{\Lip(\gamma,\Omega_1)} &\le_{\mathtt{pe}} 1,\label{const111}\\
 |r|^{\Lip(\gamma,\Omega_1)}&\le_{\mathtt{pe}} \epsilon^{7-4b},\quad |d_ir(i_0)[\ihat]|\le_{\mathtt{pe}} \epsilon\rVert \ihat|_{s_0+\mu_0},\label{const112}\\
  \rVert \tilde{q}\rVert_{s}^{\Lip(\gamma,\Omega_1)}
&\le_{\mathtt{pe},s} \epsilon^4 + \epsilon^{2}\rVert\mathfrak{I}_\delta\rVert_{s+\mu_0}^{\Lip(\gamma,\Omega_1)},\nonumber \\ \rVert d_i\tilde{q}(i_0)[\ihat]\rVert_{s}
& \le_{\mathtt{pe},s} \epsilon^2\left(\rVert \ihat \rVert_{s+\mu_0} + \rVert \mathfrak{I}_\delta \rVert_{s+\mu_0}\rVert \ihat \rVert_{s_0+\mu_0}\right).\label{const113}
 \end{align}
  \end{lemma}
  \begin{proof}
  By definition of $\mathfrak{C}_1(i_0)$, $a$ can be written as
  \[
  a=\epsilon a_1 + \epsilon^2 a_2 + \epsilon^3 a_3 + q,
  \]
  where $a_1,a_2,a_3$ have the form as in \eqref{homcoeeff} and $q$ satisfies the estimates in \ref{qqqsds22} in Definition~\ref{homogeneous_expansion1}. Clearly, we have
  \begin{align*}
  M_x(a) & = \epsilon M_x(a_1) + \epsilon^2 M_x(a_2) + \epsilon^3 M_x(a_3) + M_x(q) \\
  & = \epsilon^2 \underbrace{M_x(a_2)}_{=:\mathtt{m}} + \underbrace{\epsilon^3 M_x(a_3) + M_{\varphi,x}(q)}_{=:r} + \underbrace{M_x(q)-M_{\varphi,x}(q)}_{=:\tilde{q}},
  \end{align*}
  where the last equality follows from the fact that $0$ is not in the tangential sites therefore $\int_{\mathbb{T}}a_1dx = 0$ (see \eqref{tan_site}).
  
  By its definition, we have $
  \int_{\mathbb{T}^\nu}\tilde{q}d\varphi = 0.$
  Furthermore, \eqref{const111} follows from the structure of $a_2$ in \eqref{homcoeeff}: Recalling \eqref{norm_def_vbar} and \eqref{def_lll}, we have that $M_x(a_2) = M_{\varphi,x}(a_2)$.  \eqref{const112} follows from the  structure of $a_3$ in \eqref{homcoeeff} and the estimate for $q$ in \eqref{q_estimate_11} with \eqref{size_assumption_3}. The last estimate \eqref{const113} for $\tilde{q}$ is trivial from  \eqref{q_estimate_112}.
 \end{proof}

\subsection{Homogeneous expansion of the symbols: Class $\mathfrak{S}^m$}\label{symbol2sdsd}

\begin{definition}\label{symbolssssef}
Let $\mathfrak{a}=\mathfrak{a}(\omega,\varphi,x,\xi)$ be a symbol that depends on $\omega,\varphi$, that is, $a(\omega,\varphi,\cdot,\cdot)\in \mathcal{S}^m$\index{$\mathfrak{S}^m$}, for some $m\in \mathbb{R}$ (see \eqref{pdo_section} for the definition of $\mathcal{S}^m$). We say $\mathfrak{a}\in \mathfrak{S}^m(i_0)$ if $\mathfrak{a}$ admits an expansion:
\[
\mathfrak{a}=\epsilon\mathfrak{a}_1 + \epsilon^2\mathfrak{a}_2 + \epsilon^3\mathfrak{a}_3 + \mathfrak{q},
\]
with the following properties:
\begin{enumerate}[label=(\arabic*)]
\item \label{dksddefofssd}Each $\mathfrak{a}_i$ is of the form,
\begin{equation}\label{symbol_forms121}
\begin{aligned}
\mathfrak{a}_i(\omega,\varphi,x,\xi)&=\sum_{j_{k_1},\ldots,j_{k_i}\in S}C_{j_{k_1},\ldots,j_{k_i}}(\xi)\bar{v}_{j_1}(\varphi,x)\cdots \bar{v}_{j_i}(\varphi,x),
\end{aligned}
\end{equation}
\text{ for $i=1,2,3$}, for some smooth symbols $C_{j_{k_1},\ldots,j_{k_i}}(\xi)\in \mathcal{S}^m$, which do not depend on any of $i_0,\omega, \varphi,x,\epsilon$ (note that $\bar{v}$ depends on $\zeta$ (see \eqref{norm_def_vbar}) and thus $\omega$ by \eqref{xi_omega_dependent}).

\item  \label{sdsdsdsdsd2}$\mathfrak{q}=\mathfrak{q}(\omega,\varphi,\cdot,\cdot)\in \mathcal{S}^m$, satisfies  (under the assumption that \eqref{size_assumption_3} holds for some $\mu_0= \mu_0(\mathtt{p},\eta_0)$, 
\begin{align}
| \mathfrak{q}|^{\Lip(\gamma,\Omega_1)}_{m,s,\eta_0} & \le_{\mathtt{pe},s,\eta_0} \epsilon^{3} + \epsilon\rVert \mathfrak{I}_{\delta}\rVert^{\Lip(\gamma,\Omega_1)}_{s+\mu_0}, \nonumber \\
 | d_i \mathfrak{q}(i_0)[\ihat] |_{m,s,\eta_0} & \le_{\mathtt{pe},s,\eta_0} \epsilon\left( \rVert \ihat \rVert_{s+\mu_0} + \rVert \mathfrak{I}_\delta \rVert_{s+\mu_0}\rVert \ihat \rVert_{s_0+\mu_0}\right).\label{q_estimate_21}
\end{align}
Furthermore, $M_x(\mathfrak{q})-M_{\varphi,x}(\mathfrak{q})$ satisfies
\begin{equation}\label{q_estimate_212}
\begin{aligned}
&| M_x(\mathfrak{q})-M_{\varphi,x}(\mathfrak{q})|^{\Lip(\gamma,\Omega)}_{m,s,\eta_0}\le_{\mathtt{pe},s,\eta_0}\epsilon^4 + \epsilon^{2}\rVert\mathfrak{I}_\delta\rVert_{s+\mu_0}^{\Lip(\gamma,\Omega_1)},\\
&| d_i(M_x(\mathfrak{q})-M_{\varphi,x}(\mathfrak{q})(i_0))[\ihat]|_{m,s,\eta_0}\le_{\mathtt{pe},s,\eta_0}\epsilon^2\left(\rVert \ihat \rVert_{s+\mu_0} + \rVert \mathfrak{I}_\delta \rVert_{s+\mu_0}\rVert \ihat \rVert_{s_0+\mu_0} \right).
\end{aligned}
\end{equation}
\end{enumerate}
\end{definition}

\begin{remark}\label{eta0sdfix}
In what follows, we allow the loss of derivatives $\mu_0$ to depend on $\eta_0$, for example, \ref{sdsdsdsdsd2} of Definition~\ref{symbolssssef}. We note that $\eta_0$ will be bounded depending  on only $\mathtt{p}$, therefore we still have $\mu_0$ depend on only $\mathtt{p}$. For simplicity, we will not trace the explicit dependence of $\mu_0$ on $\eta_0$.
\end{remark}
The motivation of the above definition is that such symbols naturally arise from the coefficients in $\mathfrak{C}_2(i_0)$ in the linearized operator (see Lemma~\ref{sd22sd}). 
\begin{lemma}\label{sd22sd}
If $a=a(\omega,\varphi,x,y)\in \mathfrak{C}_2(i_0)$, then there exists $\mathfrak{a}\in \mathfrak{S}^{\alpha-3}(i_0)$ such that $\Upsilon^{\alpha-3}_a=Op^W(\mathfrak{a})$.
\end{lemma}
\begin{proof}
We first claim that
\begin{align}\label{121wesdsdsd}
\Upsilon^{\alpha-3}_a = Op(\mathfrak{b}),\text{ for some $\mathfrak{b}\in \mathfrak{S}^{\alpha-3}.$}
\end{align}
 By Definition~\ref{homogeneous_expansion1} for the class $\mathfrak{C}_2$, we have that
  \begin{align}\label{adexpenaios}
a=\sum_{i=1}^{3}\epsilon^ia_i + q,
\end{align}
which is  the decomposition of $a$ as in \eqref{homcoeeff2}. We aim to find a symbol $\mathfrak{b}\in \mathfrak{S}^{\alpha-3}$ such that 
\begin{align}\label{zmaosd2x}
Op(\mathfrak{b}) = \sum_{i=1}^3\epsilon^i \Upsilon^{\alpha-3}_{a_i} + \Upsilon^{\alpha-3}_{q}.
\end{align}
For each $i=1,2,3$, \ref{2tsksdysxxx} of Definition~\ref{homogeneous_expansion1} tells us that 
\begin{align}\label{ai1msdx}
a_{i}(\omega,\varphi,x,y) = \sum_{j_{k_1},\ldots,j_{k_i}\in S}C_{j_{k_1},\ldots,j_{k_i}}(x-y)\bar{v}_{j_1}(\varphi,x)\cdots \bar{v}_{j_i}(\varphi,x),
\end{align}
for some function $y\mapsto C_{j_{k_1},\ldots,j_{k_i}}(y)$ that does not depend on any of $i_0,\omega,\varphi,\epsilon$.
Recalling $S_{\Upsilon}$ from \eqref{def_ofsyma}, we define
\begin{align}\label{defbiqb}
\mathfrak{b}_{i}:=S_{\Upsilon}(a_i),\quad \mathfrak{q}_{\mathfrak{b}} = S_{\Upsilon}(q),\quad \mathfrak{b}:=\sum_{i=1}^3 \epsilon^i \mathfrak{b}_i + \mathfrak{q}_{\mathfrak{b}}.  
\end{align}
From the construction of $S_{\Upsilon}$ in \eqref{keypart}-\eqref{def_ofsyma}, we have that
\begin{align*}
\mathfrak{b}_i&=\sum_{j_{k_1},\ldots,j_{k_i}\in S} S_{\Upsilon}(C_{j_{k_1},\ldots,j_{k_i}})(\xi)\bar{v}_{j_1}(\varphi,x)\cdots \bar{v}_{j_i}(\varphi,x).
\end{align*}
Note that Since $ S_{\Upsilon}(C_{j_{k_1},\ldots,j_{k_i}})(\xi)$ does not depend on any of $i_0,\omega, \varphi,x,\epsilon$. From Lemma~\ref{lemmasd1sdxcxc}, \eqref{defbiqb} tells us that $S_{\Upsilon}(C_{j_{k_1},\ldots,j_{k_i}})(\xi)\in \mathcal{S}^{\alpha-3}$ and
\begin{align*}
Op(\mathfrak{b}) = \Upsilon^{\alpha-3}_{a}.
\end{align*}
To finish the proof of the claim that $\mathfrak{b}\in \mathfrak{S}^{\alpha-3}$, we need to show that $\mathfrak{q}_{\mathfrak{b}}$ satisfies the estimates in \eqref{q_estimate_21} and \eqref{q_estimate_212}. Thanks to the estimates for $q$  in \ref{qsdwd} of Definition~\ref{homogeneous_expansion1}, the estimates \eqref{q_estimate_21} and \eqref{q_estimate_212} for $\mathfrak{q}_{\mathfrak{b}}$ follows from \eqref{estsxzsdxc}.

Therefore, we have $\mathfrak{b}\in \mathcal{S}^{\alpha-3}$ satisfying \eqref{121wesdsdsd}. To find $\mathfrak{a}\in \mathfrak{S}^{\alpha-3}$ such that $Op^W(\mathfrak{a})=\Upsilon^{\alpha-3}_a = Op(\mathfrak{b})$, we use \ref{sjdsdjclasiss21} of Lemma~\ref{sdt2osd1}, then the result follows immediately.
\end{proof}

 As we will see later (especially in Lemma~\ref{fffs22}), the linearized operator $\mathcal{L}_\omega$ can be thought of as a linearized gSQG at $u$, up to a smoothing operator. Thanks to Lemma~\ref{norm_phi_expansion1},  we can see that  the terms of size $O(\epsilon), O(\epsilon^2)$ solely depend on $\bar{v}$, but not on $i_0$. The following definition will be useful to analyze the terms of size $O(\epsilon^2)$. 
 \begin{definition}\label{buildingbobo}
 We say\index{$\mathfrak{B}^m_1$}\index{$\mathfrak{B}^m_2$}\index{$\mathfrak{B}^m_3$} $\mathfrak{a}=\mathfrak{a}(\omega,\varphi,x,\xi)\in \mathfrak{B}^m_{1}$, $\mathfrak{b}=\mathfrak{b}(\omega,\varphi,x,\xi)\in \mathfrak{B}^m_{2}$, $\mathfrak{c}=\mathfrak{c}(\omega,\varphi,x,\xi)\in \mathfrak{B}^m_{3}$ if 
 \begin{equation}\label{f_k2expansion}
 \begin{aligned}
 \mathfrak{a}&=\sum_{j_k\in S}C_{j_{k}}(\xi)\bar{v}_{j_k}(\varphi,x),\\
  \mathfrak{b}&=\sum_{j_{k_1},j_{k_2}\in S}C_{j_{k_1},j_{k_2}}(\xi)\bar{v}_{j_{k_1}}(\varphi,x)\bar{v}_{j_{k_2}}(\varphi,x),\\
  \mathfrak{c}&=\sum_{j_{k_1},j_{k_2},j_{k_3}\in S}C_{j_{k_1},j_{k_2},j_{k_3}}(\xi)\bar{v}_{j_{k_1}}(\varphi,x)\bar{v}_{j_{k_2}}(\varphi,x)\bar{v}_{j_{k_3}}(\varphi,x),
\end{aligned}
\end{equation}
for some smooth symbols $C_{j_{k_1},\ldots,j_{k_i}}(\xi)\in \mathcal{S}^m$ for some $m\in \mathbb{R}$, which do not depend on any of $i_0,\varphi,x,\epsilon$ (note that $\zeta$ depend on $\omega$ by \eqref{xi_omega_dependent})
 \end{definition}
  \color{black}

 \color{black}

\section{Finite dimensional operators}\label{Finite dimensional operators}
We collect useful properties of the finite dimensional operators\index{finite dimensional operator} that are arising in the linearized operator.
\begin{definition}
We say that a $\varphi$-dependent linear operator $C^\infty_x\mapsto C^\infty_x$ is a finite dimensional operator if
\begin{align}\label{finite_3}
R(\varphi)[h]  = \sum_{|j| \le C}\int_0^1(h,g_j(\varphi,\tau))_{L^2_x}\chi_j(\varphi,\tau)d\tau,
\end{align}
for some $C>0$ and $g_j(\varphi,\tau),\chi_j(\varphi,\tau)\in C^\infty_{x}$.
\end{definition}
Clearly, $R$ is a smoothing operator\index{smoothing operator} (in the variable $x$) in the sense that $R(\varphi)[h]$ is as smooth as $g_j$ and $\chi_j$. Its  tame constants described in Subsection~\ref{s_decay_norm}  can be estimated in terms of $g$ and $\chi$:
\begin{lemma}\label{s_decay_finite} Let $R$ be a finite dimensional operator, where $g=(g_j)$ and $\chi=(\chi_j)$ depend on $\omega\in \Omega_1$. Then, the Lip-$0$-tame estimates of  $R$ (see Subsection~\ref{s_decay_norm}) satisfies
\begin{align*}
& \mathfrak{M}^\gamma_{R}(0,s)\\& \quad \le_{\mathtt{pe},s} \left(\sup_{\tau\in[0,1]}\rVert g(\tau) \rVert^{\Lip(\gamma,\Omega_1)}_s\rVert \chi(\tau)\rVert^{\Lip(\gamma,\Omega_1)}_{s_0}+ \rVert g(\tau) \rVert^{\Lip(\gamma,\Omega_1)}_{s_0}\rVert \chi(\tau)\rVert^{\Lip(\gamma,\Omega_1)}_{s}\right).
\end{align*}
Furthermore, if $g$ and $\chi$ depend on the embedding $i$,  then
\begin{align*}
&\mathfrak{M}_{d_iR(i)[\ihat]}(0,s)\\
&\quad  \le_{\mathtt{pe},s} \sup_{\tau\in[0,1]}\left(\rVert d_ig(i)[\ihat] \rVert_{s_0}\rVert \chi(i) \rVert_s+\rVert d_ig(i)[\ihat] \rVert_s\rVert \chi_j(i) \rVert_{s_0}\right.\\ & \left. \qquad \qquad \qquad +\rVert d_i\chi(i)[\ihat] \rVert_{s_0}\rVert g(i) \rVert_s+\rVert d_i\chi(i)[\ihat] \rVert_{s}\rVert g(i) \rVert_{s_0}\right)
\end{align*}
\end{lemma}
\begin{proof}
The proof is straightforward from the definition of the Lipschitz tame operators and using Lemma~\ref{interpolation_2s} and \eqref{bncal}.
\end{proof}

We define a class of finite dimensional operators, which will appear in our reduction procedure in Chapter~\ref{reduction}.
\begin{definition}\label{remiander_class12}
We say $R\in \mathfrak{R}(i_0)$,\index{$\mathfrak{R}$} if $R$ is a linear combination of the finite dimensional operators of the form $\epsilon^2R_1+R_2$, where $R_1$ is of the form,
\begin{align}\label{fisnsdw}
R_1[h]=\Pi_{S^\perp}Op^W(\mathfrak{a})\Pi_{S}(Op^W(\mathfrak{b})h),
\end{align}
for some $\mathfrak{a},\mathfrak{b}\in \mathfrak{B}^m_1$ for some $m\ge 0$ (see Definition~\ref{buildingbobo}) and $R_2$ is of the form,
\begin{align}\label{j2j2j2j2j2jssds}
R_2(\varphi)[h]  = \sum_{|j| \le C}\int_0^1(h,g_j(\varphi,\tau))_{L^2_x}\chi_j(\varphi,\tau)d\tau,
\end{align}
for some $g_j,\chi_j$ such that
\begin{equation}\label{sjdjsjsj}
\begin{aligned}
&\sup_{\tau\in[0,1]}\rVert g(\tau) \rVert^{\Lip(\gamma,\Omega_1)}_s\rVert \chi(\tau)\rVert^{\Lip(\gamma,\Omega_1)}_{s_0}+ \rVert g(\tau) \rVert^{\Lip(\gamma,\Omega_1)}_{s_0}\rVert \chi(\tau)\rVert^{\Lip(\gamma,\Omega_1)}_{s}\\
& \qquad \qquad\le_{\mathtt{pe},s}\epsilon^3+\epsilon^{2}\rVert \mathfrak{I}_\delta \rVert_{s+\mu_0}^{\Lip(\gamma,\Omega_1)},\\
&\sup_{\tau\in[0,1]}\left(\rVert d_ig(i_0)[\ihat] \rVert_{s_0}\rVert \chi(i_0) \rVert_s+\rVert d_ig(i_0)[\ihat] \rVert_s\rVert \chi_j(i_0) \rVert_{s_0}\right. \\ & \qquad \qquad \qquad \left.+\rVert d_i\chi(i_0)[\ihat] \rVert_{s_0}\rVert g(i_0) \rVert_s+\rVert d_i\chi(i_0)[\ihat] \rVert_{s}\rVert g(i_0) \rVert_{s_0}\right)\\
& \qquad \qquad  \le_{\mathtt{pe},s} \epsilon^2\rVert \ihat \rVert_{s+\mu_0} + \epsilon^{2b-1}\rVert \mathfrak{I}_\delta \rVert_{s+\mu_0}\rVert \ihat \rVert_{s+\mu_0}.
\end{aligned}
\end{equation}
\end{definition} 
Note that  the operator in \eqref{fisnsdw} is indeed a finite dimensional operator of the form in \eqref{finite_3}, because 
\begin{align*}
Op^W(\mathfrak{a})\Pi_S(Op^W(\mathfrak{b})h) & =\sum_{j\in S}(Op^W(\mathfrak{b})h,e^{\ii j x})_{L^2_x}Op^W(\mathfrak{a})[e^{\ii j x}] \\
& = \sum_{j\in S}(h, Op^W(\overline{\mathfrak{b}})[e^{\ii j x}])_{L^2_x}Op^W(\mathfrak{a})[e^{\ii j x}],
\end{align*}
where the last equality follows from \eqref{transfposdsdsd}.

\color{black}

\begin{lemma}\label{funitesd2sd}
Let $\Phi^{\pm}=\Phi,\Phi^{-1}$ be a symplectic transformations on $C^\infty_{\varphi,x}$ with the estimates:
\begin{equation}\label{transformation_estimate_clas2s}
\begin{aligned}
&\rVert (\Phi^{\pm}-I) h \rVert_{s}^{\Lip(\gamma,\Omega_1)} \le_{\mathtt{pe},s} \epsilon\left( \rVert h \rVert^{\Lip(\gamma,\Omega_1)}_{s+\mu_0} + \rVert \mathfrak{I}_\delta \rVert^{\Lip(\gamma,\Omega_1)}_{s+\mu_0}\rVert h \rVert^{\Lip(\gamma,\Omega_1)}_{s_0+\mu_0}\right)\\
&\rVert d_i\Phi^{\pm}(i_0)h[\ihat] \rVert_{s} \le_{\mathtt{pe},s} \left(\rVert h\rVert_{s+\mu_0} + \rVert \mathfrak{I}_\delta \rVert_{s+\mu_0}\rVert h \rVert_{s_0+\mu_0}\right) \rVert \ihat \rVert_{s_0+\mu_0} + \rVert h \rVert_{s_0+\mu_0}\rVert \ihat \rVert_{s+\mu_0}.
\end{aligned}
\end{equation}
If $\mathcal{R}\in \mathfrak{R}(i_0)$, then we have $\mathcal{R}\Phi,\Phi\mathcal{R}\in \mathfrak{R}(i_0)$.
\end{lemma}
\begin{proof}
We will show that if $\mathcal{R}=\epsilon^2 R_1 + R_2$ where $R_1$ and $R_2$ are of the form in \eqref{fisnsdw} and \eqref{j2j2j2j2j2jssds}, then  $\mathcal{R}\Phi\in \mathfrak{R}(i_0)$.  The same results for $\Phi\mathcal{R}$ and  a linear combination of such operators  follow in a similar way.

  We have that
\begin{align*}
 R_1h &= \Pi_{S^\perp}Op^W(\mathfrak{a})\Pi_{S}Op^W(\mathfrak{b})h, \text{ with $\mathfrak{a},\mathfrak{b}\in \mathfrak{B}^m_1$ for some $m\in \mathbb{R}$},\\
 R_2h& = \sum_{|j|\le C}\int_0^{1}(h,g_j(\tau))_{L^2}\chi_j(\tau)d\tau,
\end{align*}
as described in Definition~\ref{remiander_class12}. Therefore it follows that
\[
\mathcal{R}\Phi h = \epsilon^2R_1\Phi h + R_2 \Phi h = \epsilon^2 R_1 h + \underbrace{\left(\epsilon^2 R_1 (\Phi-I)h + R_2\Phi h \right)}_{=:\mathcal{R}_{2}}.
\]
It suffices to show that $\mathcal{R}_2$ can be written in the form \eqref{j2j2j2j2j2jssds} with the desired estimates in \eqref{sjdjsjsj}. Using $\Pi_S h = \sum_{j\in S}(h,e^{\ii j x})_{L^2}e^{\ii j x}$, we have
\begin{align*}
\epsilon^2R_1(\Phi - I)h &=\epsilon^2\sum_{j\in S} \Pi_{S^\perp}Op^W(\mathfrak{a})[(Op^W(\mathfrak{b})(\Phi-I)h,e^{\ii j x})_{L^2(\mathbb{T})}e^{\ii j x}]\\
& =\epsilon^2 \Pi_{S^\perp}\sum_{j\in S}(h,\left(Op^W(\mathfrak{b})(\Phi-I)\right)^{T}[e^{\ii j x}])_{L^2(\mathbb{T})}Op^W(\mathfrak{a})[e^{\ii j x}]\\
& = \epsilon^2\Pi_{S^\perp}\sum_{j\in S}(h,Op^W(\overline{\mathfrak{b}})(\Phi^T - I)[e^{\ii j x}])_{L^2}Op^W(\mathfrak{a})[e^{\ii j x}]\\
& = \Pi_{S^\perp}\sum_{j\in S}(h,\underbrace{\epsilon^2Op^W(\overline{\mathfrak{b}})(\Phi^T - I)[e^{\ii j x}]}_{=:\tilde{g}_j})_{L^2}\underbrace{Op^W(\mathfrak{a})[e^{\ii j x}]}_{=:\tilde{\chi}_j},
\end{align*}
where the third equality follows from \eqref{adjoint_weyl}. Note that the last expression is of the form  \eqref{j2j2j2j2j2jssds}. Also, using that $\Phi$ is symplectic and \eqref{symplectic73}, we have
\[
\tilde{g}_j =\epsilon^2(Op^W(\overline{\mathfrak{b}})(\partial_x^{-1}\Phi^{-1}\partial_x - I)) = \epsilon^2(Op^W(\overline{\mathfrak{b}})(\partial_x^{-1}(\Phi^{-1}-I)\partial_x ))
\]Using the definition of $\mathfrak{B}_1$ and the estimates \eqref{transformation_estimate_clas2s}, it follows straightforwardly that $\tilde{g}_j,\tilde{\chi}_j$ satisfy the estimates \eqref{sjdjsjsj}.
\end{proof}

\section{Linearized operator in the normal directions}\label{linisdw}

 In view of \eqref{K_{20}_1} and \eqref{whatklooklike}, we look at the contribution of each transformation to the Hessian\index{Hessian} of the Hamiltonian $K$.

\begin{lemma}\cite[Lemma 7.4]{Baldi-Berti-Montalto:KAM-quasilinear-kdv}\label{normal_direction1}
 Let  $H$ be as in \eqref{normal_form_f}. There exists a linear operator $R(\varphi)=R(i_0)(\varphi):H_{S^\perp}\mapsto H_{S^\perp}$ such that 
\[
K_{02}(\varphi)[w] = \Pi_{S^\perp}\nabla_z^2 H(U_\zeta(i_\delta(\varphi)))[w] + R(\varphi)[w],
\]
where $R(\varphi)$ is a finite dimensional operator in \eqref{finite_3} such that
$g_j(\varphi)=g_j(i_0)(\varphi)$ and $\chi_j(\varphi)=\chi(i_0)(\varphi)$ satisfy
\begin{align}
\rVert g_j \rVert^{\Lip(\gamma,\Omega_1)}_s\rVert \chi_j\rVert^{\Lip(\gamma,\Omega_1)}_{s_0}+ \rVert g_j \rVert^{\Lip(\gamma,\Omega_1)}_{s_0}\rVert \chi_j\rVert^{\Lip(\gamma,\Omega_1)}_{s} \le_s \epsilon^{b+1}\rVert \mathfrak{I}_\delta\rVert^{\Lip(\gamma,\Omega_1)}_{s+\mu_0}\label{finite_1}
\end{align}
and
\begin{align}
&\rVert d_ig_j(i_0)[\ihat] \rVert_{s_0}\rVert \chi_j(i_0) \rVert_s+\rVert d_ig_j(i_0)[\ihat] \rVert_s\rVert \chi_j(i_0) \rVert_{s_0} \nonumber \\
& \quad +\rVert d_i\chi_j(i_0)[\ihat] \rVert_{s_0}\rVert g_j(i_0) \rVert_s+\rVert d_i\chi_j(i_0)[\ihat] \rVert_{s}\rVert g_j(i_0) \rVert_{s_0}\nonumber\\
& \quad \le_s \epsilon^{b+1}\rVert \ihat \rVert_{s+\mu_0} + \epsilon^{2b-1}\rVert \mathfrak{I}_\delta \rVert_{s+\mu_0}\rVert \ihat \rVert_{s_0+\mu_0}.\label{finite_2}
\end{align}
\end{lemma}

\begin{remark}\label{finitesdsd}
Among the finite dimensional operators arising in $\nabla_z^2H$, the largest contribution comes form the Hessian of $\mathfrak{H}_1$ in \eqref{sjdjsjs93jx} (see Lemma~\ref{fffs22}). Note that $\nabla_z \mathfrak{H}_1$ is a finite dimensional operator, since $z\mapsto \mathfrak{H}_1(v+z)$ depends on only a finite number of Fourier modes of $z$, because of the projection $\Pi_S$ in the integral.
\end{remark}

Now, we consider the contribution of $\Phi^{WB}$. We first recall the following lemma:
\begin{lemma}\cite[Lemma 7.1]{Baldi-Berti-Montalto:KAM-quasilinear-kdv}\label{normal_direction2}
Let $\mathcal{H}$ be a Hamiltonian of class $C^2(H^1_0(\mathbb{T}),\mathbb{R})$ and consider a map $\Phi(f) := f+ \Psi(f)$ satisfying  $\Psi(f) = \Pi_E \Psi(\Pi_E f)$, for all $f\in H^1_0(\mathbb{T})$, for some finite dimensional subspace  $E$  as in \eqref{finite_dim}. Then
\begin{align*}
\nabla^2 (\mathcal{H}\circ \Phi)(f)[h] = \nabla^2\mathcal{H}(\Phi(f))[h] + R_{\mathcal{H}}(f)[h],
\end{align*}
where $R_{\mathcal{H}}(f)$ is a finite dimensional operator in \eqref{finite_3}. For the explicit form of  $R_{\mathcal{H}}(f)$, denoting  
\begin{align*}
h(x)& :=\sum_{j\in \mathbb{Z}} h_j e^{\ii j x},\quad \nabla \Psi(f)[\cdot] := \Pi_E \nabla \Psi(f) [\Pi_E\cdot],\\ \nabla^2 \Psi(f)[\cdot,\cdot]& = \Pi_E \nabla^2 \Psi(f)[\Pi_E \cdot , \Pi_E (\cdot)],
\end{align*} (which follow from $\Psi(f)=\Pi_E \Psi(\Pi_E f)$), we have that  
\begin{align}
R_{\mathcal{H},0}(f)[h] & = \sum_{|j| < C}h_j \chi_j(x),\quad \chi_j:=\nabla^2\mathcal{H}(\Phi(f))[\Pi_E\nabla \Psi(f)[e^{\ii j x}]], \label{R_0_form}\\
R_{\mathcal{H},1}(f)[h] & = \sum_{|j| < C}h_j \chi_j(x),\quad \chi_j:=\Pi_EA_j^{T}[\Pi_E\nabla \mathcal{H}(\Phi(f))], \nonumber \\
\text{ where } A_j[\cdot]& :=\nabla^2\Psi(f)[e^{\ii j x}, \Pi_E\cdot],\label{R_1_form}\\
R_{\mathcal{H},2}(f)[h] & = \sum_{|j| < C}h_j \chi_j,\quad \chi_j:=\Pi_E (\nabla \Psi(f))^T\Pi_E\nabla^2\mathcal{H}(\Phi(f))[e^{\ii j x}]. \label{R_2_form}
\end{align}
\end{lemma}
Recall from Proposition~\ref{normal_form_prop11231} that $H$ in \eqref{normal_form_f} is given by $H=\mathcal{H}\circ \Phi^{WB}$, where $\Phi^{WB}=f+\Psi$ and $\Psi$ satisfies $\Psi = \Pi_{E}\Psi\circ \Pi_E$ for a finite dimensional space $E$ of the form in \eqref{finite_dim}. 
\begin{lemma}\label{x2h2h2h2h}
Let $\mathcal{H}$ be as in Proposition~\ref{expansion_1}. There exist linear operators $R_1,R_2$ of the form \eqref{finite_3} such that 
\begin{align}\label{r1sded}
\nabla_z^2\left(\mathcal{H}\circ \Phi^{WB}\right)(U_\zeta(i_\delta)) = \nabla_z^2\mathcal{H}(\Phi^{WB}(U_\zeta(i_\delta))) + R_1(\varphi) + R_2(\varphi),
\end{align}
where
\begin{align*}
R_1(\varphi)[h] &=\epsilon^2\nabla_z^2\mathfrak{H}_1(v_\epsilon(\theta_0(\varphi),y_\delta(\varphi)))[h], \text{ where $\mathfrak{H}_1$ is as in \eqref{sjdjsjs93jx}},
\end{align*}
and $R_2$, which can be written as $R_2(\varphi)[h]=\sum_{|j|\le C}(h,e^{\ii j x})_{L^2}\chi_j(\varphi)$, satisfies
\begin{align}
\rVert \chi_j\rVert^{\Lip(\gamma,\Omega_1)}_{s} \le_{\mathtt{pe},s} \epsilon^{7-3b} + \epsilon^{1+b}\rVert \mathfrak{I}_\delta\rVert^{\Lip(\gamma,\Omega_1)}_{s+\mu_0}\label{finite_4}
\end{align}
and
\begin{align}\label{finite_5}
\rVert d_i\chi_j(i_0)[\ihat] \rVert_{s} \le_{\mathtt{pe},s} \epsilon^2\rVert \ihat \rVert_{s+\mu_0} + \epsilon^{2b-1}\rVert \mathfrak{I}_\delta \rVert_{s+\mu_0}\rVert \ihat \rVert_{s+\mu_0}.
\end{align}

\end{lemma}
\begin{proof}
From Lemma~\ref{normal_direction2}, it follows that there exists a finite dimensional operator $R$ such that
\[
\nabla_z^2(\mathcal{H}\circ\Phi^{WB})(U_\zeta(i_\delta)) = \nabla_z^2\mathcal{H}(\Phi^{WB}(U_\zeta(i_\delta))) + R(\varphi),
\]
where $R$ is of the form
\begin{align}\label{sd2222}
R(\varphi)[h] =\sum_{|j|\le C}(h,e^{\ii j x})_{L^2}\chi_j(\varphi).
\end{align}
We recall from Proposition~\ref{normal_form_prop11231}, we have that
\begin{align}\label{hhss22}
\mathcal{H}\circ\Phi^{WB}(f) =H(f)= H_2+H_3+H_4+H_{\ge 5},
\end{align}
where $H_{k}$ is homogeneous of degree $k$. First (denoting $f=v+z$ the tangential and normal component of $f$), we notice from \eqref{normal_formsf} that  $\nabla_z^2( H_2 + H_3)(U_\zeta(i_\delta))$ does not have a finite dimensional operator. In $H_4$, we can find from \ref{sdsd22sd} of Proposition~\ref{normal_form_prop11231} and \eqref{expansion_2} that $\mathcal{H}_{4,2},\mathfrak{H}_2$ do not induce a finite dimensional operator, while $\mathfrak{H}_1$ in \eqref{sjdjsjs93jx32x} gives a finite dimensional operator and we have 
\[
R_1(\varphi):=\nabla_z^2\mathfrak{H}_1(U_\zeta(i_\delta))[h] = \epsilon^2\nabla_z^2\mathfrak{H}_1(v_\epsilon(\theta_0(\varphi),y_\delta(\varphi)))[h].
\]
All the other homogeneous terms which can yield a nonzero Hessian in $z$ are of the form (recall the notation \eqref{multilinear52}):
\begin{align}\label{gform}
G(f)=G(v+z)=
\begin{cases}
R(vz^3), \text{ or }R(z^4) & \text{ from quartic terms,}\\
R(v^{n-k}z^{k}),& \text{ for $n=5$, and $k\ge 2$, or $n\ge 6$} \\& \text{ from higher order terms},
\end{cases}
\end{align}
which follows from \ref{sdsd22sd} of Proposition~\ref{normal_form_prop11231}.
We will show that a finite dimensional operator $\tilde{R}$ of the form in \eqref{sd2222}, given as $\tilde{R}[h]=\nabla_z^2G(U_\zeta(i_\delta))$, satisfies the estimates \eqref{finite_1} and \eqref{finite_2}. 

 We will prove the case where $G(f)=R(vz^3)$ case only, since the other terms can be treated in the same way. Thanks to \eqref{item3_normal_3} in Proposition~\ref{normal_form_prop11231}, there exists a $\mu_0$-regular trilinear map $A$ (see Definition~\ref{musmooth}) such that (recalling \eqref{actionssx})
 \begin{align*}
\tilde{ R}(\varphi)[h]& =\nabla_z^2G(U_\zeta(i_\delta))[h] = \epsilon^{1+b}A(v_\epsilon(\theta_0(\varphi),y_\delta(\varphi)),z_0(\varphi),h) \\
& = A(v_\epsilon(\theta_0(\varphi),y_\delta(\varphi)),z_0(\varphi),\Pi_Eh),
 \end{align*}
 for some finite dimensional space $E$, where the last equality follows from the assumption that $\tilde{R}$ in \eqref{sd2222} takes only a finite number of Fourier modes. This implies that $\chi_j$ is given by
 \[
 \chi_j (\varphi):=\epsilon^{1+b}A(v_\epsilon(\theta_0(\varphi),y_\delta(\varphi)),z_0(\varphi),e^{\ii j x}). 
 \]
 Since $A$ is $\mu_0$-regular, recalling its definition from Definition~\ref{musmooth} and using \eqref{size_assumption_3}, the estimates \eqref{finite_1} and \eqref{finite_2} follow straightforwardly.
\end{proof}

Combining Lemma~\ref{normal_direction1} and Lemma~\ref{x2h2h2h2h}, we have the following:
\begin{lemma}\label{fffs22}
We have that
\begin{align}\label{rksdsd1cxcxaa}
K_{02}(\varphi) = \Pi_{S^\perp}\nabla_z^2\mathcal{H}(\Phi^{WB}(U_\zeta(i_\delta))) + \epsilon^2\nabla_z^2\mathfrak{H}_1(\bar{v}(\varphi))[h] + R(\varphi),
\end{align}
where a finite dimensional operator $R$ in the form in \eqref{j2j2j2j2j2jssds} satisfies the estimates in \eqref{sjdjsjsj}.
\end{lemma}

\begin{proof}
Thanks to Lemma~\ref{normal_direction1}, Lemma~\ref{x2h2h2h2h} and \eqref{normal_form_f}, we have
\begin{align}
K_{02}(\varphi)[h] &= \Pi_{S^\perp}\nabla_z^2 H(U_\zeta(i_\delta))[h] +R_1(\varphi)\nonumber\\
& = \Pi_{S^\perp}\nabla_z^2\mathcal{H}(\Phi^{WB}(U_\zeta(i_\delta)))[h] + R_1(\varphi)[h] + R_2(\varphi)[h] \nonumber \\
& + \epsilon^2\nabla_z^2\mathfrak{H}_1(v_\epsilon(\theta_0(\varphi),y_\delta(\varphi)))[h],
\end{align}
for some  $R_1$ satisfying the estimates in \eqref{finite_1} and \eqref{finite_2}, and $R_2$ satisfying \eqref{finite_4} and \eqref{finite_5}. Note that the estimates \eqref{finite_1}, \eqref{finite_2}, \eqref{finite_4} and \eqref{finite_5} are stronger than \eqref{sjdjsjsj}, therefore, $R_1$ and $R_2$ satisfy \eqref{sjdjsjsj}.
Therefore, to finish the proof, it suffices to prove that
\[
R_*(\varphi):= \epsilon^2\left(\nabla_z^2\mathfrak{H}_1(v_\epsilon(\theta_0(\varphi),y_\delta(\varphi))) - \nabla_z^2\mathfrak{H}_1(\bar{v}(\varphi))\right),
\]
satisfies the estimate \eqref{sjdjsjsj}, which will give us \eqref{rksdsd1cxcxaa} for $R:=R_1+R_2+R_*$. Recalling $\mathfrak{H}_1$ from \eqref{sjdjsjs93jx}, we see that there exists a $\mu_0$-regular (for some $\mu_0\ge 0$. See Definition~\ref{musmooth}) trilinear operator $A=A(f_1,f_2,f_3)$ and a finite dimensional space $E$ of the form in \eqref{finite_dim} such that
\[
\nabla_z^2\mathfrak{H}_1(f)[h] = A(\Pi_S f,\Pi_{S}f,\Pi_E h).
\]
Hence,
\begin{align}
R_*[h] &= \epsilon^2\left( A(v_\epsilon(\theta_0,y_\delta),v_\epsilon(\theta_0,y_\delta),\Pi_Eh) - A(\bar{v},\bar{v},\Pi_E h)\right) \nonumber \\
& =: \sum_{|j|\le C}(h,e^{\ii j x})_{L^2}\chi_j,\label{sdsd1xcps2xcx}
\end{align}
where
\begin{align}\label{dhsk222}
\chi_j(\varphi,x) =\epsilon^2\left( A(v_\epsilon(\theta_0,y_\delta),v_\epsilon(\theta_0,y_\delta),e^{\ii j x}) - A(\bar{v},\bar{v},e^{\ii j x})\right).
\end{align}
Then, from \eqref{vvvbbb}, which says $|v_\epsilon - \bar{v}| = O(\epsilon^{6-4b} + |\mathfrak{I}_\delta|)$,  it follows that $g_j:=e^{\ii j x}$ and $\chi_j$ defined in \eqref{dhsk222} satisfy the estimates \eqref{sjdjsjsj} (the estimate for $d_i\chi(i_0)$ in \eqref{sjdjsjsj} follows straightforwardly, since $A$ is a trilinear map). Hence, $R_*$ in \eqref{sdsd1xcps2xcx}  satisfies  the estimates \eqref{sjdjsjsj}.
\end{proof}

 Now we state the main result of this section.
 \begin{proposition}\label{linearized_opspdosodw}
 The linear operator $\mathcal{L}_\omega$ in \eqref{third_equation} has the form\index{$\mathcal{L}_\omega$}:
 \begin{align}\label{linearisd2}
 \mathcal{L}_\omega h &= \D_\omega - \Pi_{S^\perp}\partial_xM[h] + \mathcal{R},\text{ with } M[h]:=\Lambda^{\alpha-1}_{-\frac{1}{2} + a_1} + \left(\frac{T_\alpha}{4}+a_2 \right)h + \Upsilon^{\alpha-3}_{a_3}
 \end{align}
 where
 \begin{enumerate}[label=(\arabic*)]
 \item \label{linsi1}$a_1,a_3\in \mathfrak{C}_2(i_0)$ and $a_2\in \mathfrak{C}_1(i_0)$ (see Definition~\ref{homogeneous_expansion1} for the sets $\mathfrak{C}_1(i_0),\mathfrak{C}_2(i_0)$). Furthermore, $a_1,a_2,a_3$ are real-valued and 
 \begin{align}
 a_1(\varphi,x,y) & = a_1(-\varphi,-x,-y),\quad a_3(\varphi,x,y) = a_3(-\varphi,-x,-y),\nonumber \\
 a_2(\varphi,x) & = a_2(-\varphi,-x).\label{reverssss}
 \end{align}
 \item \label{linsi2}$\mathcal{R}\in \mathfrak{R}(i_0)$. Specifically,
 \begin{align}\label{sd11cantake}
 \mathcal{R}(\varphi):=\epsilon^2\partial_x \nabla_z^2\mathfrak{H}_1(\bar{v}(\varphi))[h] + \partial_xR(\varphi),
 \end{align}
 where $R$ is a finite dimensional operator satisfying \eqref{sjdjsjsj}.
 \item \label{linsi3}$\mathcal{L}_\omega$ is real, reversible (see Definition~\ref{reversible_operators2s} for a reversible operator) and Hamiltonian (see Definition~\ref{hamiltonian_operatior230} for a Hamiltonian operator).
 
 \item \label{linsi4}  $\mathcal{L}$ is $\frac{2\pi}{\mathtt{M}}$-translation invariant and $a_1,a_2,a_3$ are also $\frac{2\pi}{\mathtt{M}}$-translation invariant, that is,
 \begin{align}\label{alsd2aem}
 \rho_{\mathtt{M}}(a_2)(\varphi,x) = a_2(\varphi,x),\quad \rho_{\mathtt{M}}(a_i)(\varphi,x,y) = a_i(\varphi,x,y),\text{ for $i=1,3$}, 
 \end{align}
where  $\rho_{\mathtt{M}}(a)(\varphi,x,y):=a(\varphi,x+\frac{2\pi}{\mathtt{M}},y+\frac{2\pi}{\mathtt{M}})$.
 
 \end{enumerate}
 \begin{proof}
 We first prove the expression of $\mathcal{L}_\omega$ in \eqref{linearisd2} and then give proofs for \ref{linsi1}-\ref{linsi4}. For $u(x)=u(\varphi,x)$, defined in \eqref{u_def1},
 we compute $\nabla_z^{2}\mathcal{H}(u)[h]$.  From Proposition~\ref{expansion_1}, it is clear that there exist  $F^1,F^2,F^3:\mathbb{R}^3\mapsto \mathbb{R}$, which are real analytic near $(0,0,0)$ such that $F^i(0,0,0)=0$ for $i=1,2,3$ and 
\begin{equation}\label{hessian1}
\begin{aligned}
& \nabla^2_z \mathcal{H}(u)[h] = -\frac{1}{2}\Lambda^{\alpha-1}h + \frac{T_\alpha}{4}h \\
& \  + \int_{\mathbb{T}}(2-2\cos(x-y))^{-\frac{\alpha}{2}}F^1(u(x),u(y),J(u)(x,y))(h(x)-h(y))dy\\
& \ + \int_{\mathbb{T}}(2-2\cos(x-y))^{-\frac{\alpha}{2}}F^2(u(x),u(y),J(u)(x,y))(u(x)-u(y))dyh(x)\\
& \ + \int_{\mathbb{T}}(2-2\cos(x-y))^{1-\frac{\alpha}{2}}F^3(u(x),u(y))h(y)dy,
\end{aligned}
\end{equation}
where $J(u)(x,y):=\frac{(u(x)-u(y))^2}{2-2\cos(x-y)}$.
 Let us denote
\begin{equation}\label{def_as}
\begin{aligned}
&a_1(x,y) := F^1(u(x),u(y),J(u)(x,y)),\\
&a_2(x) =  \int_{\mathbb{T}}(2-2\cos(x-y))^{-\frac{\alpha}{2}}F^2(u(x),u(y),J(u)(x,y))(u(x)-u(y))dy,\\
&a_3(x,y) := F^3(u(x),u(y)),
\end{aligned}
\end{equation}
so that we have
\begin{align}\label{hessian23}
\nabla_z^2\mathcal{H}(u)[h] = \Lambda^{\alpha-1}_{-\frac{1}{2} + a_1}h + \left(\frac{T_\alpha}{4}+a_2\right)h + \Upsilon^{\alpha-3}_{a_3}h = M[h].
\end{align}
Note that $a_1,a_2,a_3$ depend on $\omega,\varphi$ and $i_0$ as well, since so does $u$ (see \eqref{u_def1}).  Then, it follows from Lemma~\ref{fffs22} and \eqref{u_def1} that
\begin{align}\label{k2k2k2k2k2ksss}
K_{02} = M[h] + \epsilon^2\nabla_z^2\mathfrak{H}_1(\bar{v}(\varphi))[h] + R(\varphi),
\end{align}
where $R$ satisfies the estimates in \eqref{sjdjsjsj}.
From the definition of $\mathcal{L}_\omega$ in \eqref{third_equation}, we have
\[
\mathcal{L}_\omega = \D_\omega -\Pi_{S^\perp}\partial_xMh + \partial_x\mathcal{R},\quad \mathcal{R}:=\epsilon^2\nabla_z^2\mathfrak{H}_1(\bar{v}(\varphi))[h] + R(\varphi),
\]
which gives the expression of $\mathcal{L}_\omega$ in \eqref{linearisd2}.

We prove \ref{linsi2} first and then move on to \ref{linsi1}, \ref{linsi3} and \ref{linsi4}.

\vspace{0.5\baselineskip}
\noindent\textit{Proof of \ref{linsi2}.}
The decomposition $\mathcal{R}= \epsilon^2\nabla_z^2\mathfrak{H}_1(\bar{v}(\varphi))[h] + R(\varphi)$, where $R$ satisfies \eqref{sjdjsjsj}, has just been proved above. Noticing $\nabla_z^2\mathfrak{H}_1(\bar{v}(\varphi))[h]$ is of the form in \eqref{fisnsdw} (see \eqref{sjdjsjs93jx} for the definition of $\mathfrak{H}_1$), it implies that $\mathcal{R}\in \mathfrak{R}(i_0)$.

\vspace{0.5\baselineskip}
\noindent\textit{Proof of \ref{linsi1}.}
It follows immediately from Lemma~\ref{properties_weights} that $a_1,a_3\in \mathfrak{S}_2(i_0)$ and $a_2\in \mathfrak{S}_1(i_0)$. Furthermore, $a_1,a_2,a_3$ are real, since $F^{1},F^2,F^3$ in \eqref{def_as} are real-valued functions. To prove the symmetry in \eqref{reverssss},  note that $u\in X$, that is,
\[
u(-\varphi,-x)=u(\varphi,x),
\]
which follows from \eqref{u_def1}, Lemma~\ref{isotropicisrever}, \ref{item3_normal_2} of Proposition~\ref{normal_form_prop11231}, and \eqref{actionssx}. Hence, \eqref{reverssss} follows straightforwardly from \eqref{def_as}.

\vspace{0.5\baselineskip}
\noindent\textit{Proof of \ref{linsi3}.}
Clearly, $K$ is a real-valued Hamiltonian (see \eqref{modified_Hamiltonian_2}) since $H_\zeta$ is a real-valued Hamiltonian.
$K$ is a reversible Hamiltonian, since $H_\zeta$ is reversible (Lemma~\ref{rsksdtimerevsd}), $G_\delta$ is reversibility preserving (Lemma~\ref{Gsers}) and $K=H_\zeta\circ G_\delta$, which is given in \eqref{modified_Hamiltonian_2}.  From the definition of $\mathcal{L}_\omega$, and Lemma~\ref{reversible_sdesd}, we see that $\mathcal{L}_\omega$ is reversible. Furthermore, since $K_{02}$ is a symmetric operator (because it is a Hessian of a Hamiltonian, thanks to \eqref{K_{20}_1}), $\mathcal{L}_\omega$ (see \eqref{third_equation}) is a Hamiltonian operator.

\vspace{0.5\baselineskip}
\noindent\textit{Proof of \ref{linsi4}.} Lemma~\ref{mfisodswwdsd2s2} imples that $K_{02}$ in \eqref{K_{20}_1} is  $\frac{2\pi}{\mathtt{M}}$-translation invariant, and therefore so is $\mathcal{L}_\omega$ in \eqref{third_equation}, since $\D_\omega$ does not destroy this property. Since $i_\delta$ is $\frac{2\pi}{\mathtt{M}}$-translation invariant, it follows from \ref{item3_normal_2} of Proposition~\ref{normal_form_prop11231}, \eqref{uusd11} and \eqref{u_def1} that $u$ is $\frac{2\pi}{\mathtt{M}}$-translation invariant. Then, \eqref{alsd2aem} follows immediately from \eqref{def_as}.
 \end{proof}
 \end{proposition}
 
 \section{Structure of the operator of size $O(\epsilon^2)$}
 In Proposition~\ref{linearized_opspdosodw}, we decomposed the linear operator $\mathcal{L}_\omega$ into a sum (up to $\D_\omega$) of a finite dimensional operator $\mathcal{R}$ and a "non"-finite dimensional operator $M$. Before we close this chapter, we specify the terms in $\mathcal{L}_\omega$ of size $O(\epsilon^2)$ for the purpose of the analysis of the eigenvalues of $\partial_xM + \mathcal{R}$. The results in this section will be used in Section~\ref{linearbbss22}. 
 We say that for linear operators $L_1,L_2:C^\infty_{\varphi,x}\mapsto C^\infty_{\varphi,x}$, and $a>0$,
 \begin{align}\label{saytheyareequal}
 L_1=L_2, \text{ up to $O(\epsilon^a)$}, 
 \end{align}
 if there exists $C,\mu_0>0$ such that $\rVert(L_1-L_2)[h]\rVert_{{s_0}}\le C\epsilon^a\rVert h \rVert_{s_0+\mu_0}$, for all $h\in C^\infty_{\varphi,x}$.
 
 Arguing as in the proofs of Lemma~\ref{x2h2h2h2h} and Lemma~\ref{fffs22}, we see that the terms of size $O(\epsilon^2)$ depend only on $H_{2,2},H_{3,2}$ and $H_{4,2}$ in \ref{sdsd22sd} of Proposition~\ref{normal_form_prop11231}, that is, (for each fixed $\omega$)
 \[
 \nabla_z^2 H(U_\zeta(i_\delta)) =\nabla_z^2H_{2,2}(\overline{v})+ \epsilon\nabla_z^2 H_{3,2}(\bar{v}) + \epsilon^2\nabla_z^2\mathcal{H}_{4,2}(\bar{v}), 
 \]
 up to $O(\epsilon\rVert\mathfrak{I_\delta}\rVert_{s_0+\mu_0}) = O(\epsilon^{7-4b}$).
 Furthermore, \eqref{finite_1} in Lemma~\ref{normal_direction1} says that (note that $|\mathfrak{I}_\delta| = O(\epsilon^{6-4b})$ thanks to Remark~\ref{Insteadofimbd}),\[
 K_{02}=\Pi_{S^\perp}\nabla_z^2H(U_\zeta(i_\delta)) \text{  up to $O(\epsilon^{7-4b})$.}\]
 Therefore, we have that (for each fixed $\omega$),
 \begin{align}\label{epsilon2ssss}
 \mathcal{L}_\omega = \D_\omega - \partial_x K_{02} = \D_\omega - \left(\partial_x\nabla_z^2H_{2,2}(\bar{v}) +\epsilon  \partial_x\nabla_z^2H_{3,2}(\bar{v}) +\epsilon^2 \partial_x\nabla_z^2H_{4,2}(\bar{v})\right), 
 \end{align}
 up to $O(\epsilon^{7-4b})$,
 where $\overline{v}$ is defined in \eqref{norm_def_vbar}.
 
 \color{black}
 

 \chapter{Symplectic transformations}\label{proper_transform_9chapter}
In Chapter~\ref{reduction}, we will conjugate the operator $\mathcal{L}_\omega$ with reversible transformations $\Phi$ to reduce each coefficient of $M$ to a constant.
Before we start the conjugations, we study useful properties of the transformations that we will use throughout the chapter.
\begin{definition}\label{typesoftransformations}
\begin{enumerate}[label= (\arabic*)]
\item \label{type1} We say that $\Phi[h]$ is type (1)\index{type (1) transformation}\index{type (2) transformation}\index{type (3) transformation}\index{type (4) transformation} if it is  a time-$1$ flow map determined by a Hamiltonian PDE in $H_{S^\perp}$:
\begin{align}\label{type1flow}
\partial_t u = \Pi_{S^\perp}\partial_x(b(t)\Pi_{S^\perp}u), \text{ where $b(t)=b(t,\omega,\varphi,x):=\frac{\beta(\omega,\varphi,x)}{1+t\beta_x(\omega,\varphi,x)}$,}
\end{align}
for some $\beta\in \mathfrak{C}_1(i_0)$, where $\beta$ is real-valued and  $\mathfrak{C}_1(i_0)$ is defined in Definition~\ref{homogeneous_expansion1}.
\item \label{type2}We say that $\Phi[h]$ is type (2) if it is  a time-$1$ flow map determined by a Hamiltonian PDE in $H_{S^\perp}$:
\begin{align}\label{type1flow2s}
\partial_t u = \Pi_{S^\perp}\partial_x Op^W(\mathfrak{a})\Pi_{S^\perp}u, \text{ for some real-valued symbol $\mathfrak{a}$.} 
\end{align}

\item \label{typ21e2}We say that $\Phi[h]$ is type (3) if 
\begin{align}\label{sd3jwwj22}
\Phi[h] = h(\varphi+\omega p_1(\varphi),x),
\end{align}
for some real-valued function $p_1(\varphi) = p_1(\omega,\varphi)$, independent of $x$. 

\item \label{type3}We say that $\Phi[h]$ is type (4) if $\Phi[h] = \sum_{j\ne 0}h_j(\varphi)e^{\ii j (x + p_2(\varphi,j))}$ for some  real-valued symbol $p_2(\varphi,\xi) = p_2(\omega,\varphi,\xi)$ that depends on $i_0$ and $\omega$ but not on $x$.
\end{enumerate}
\end{definition}

Note that all of the transformations above are symplectic, more precisely, they satisfy \eqref{symplectic73}, since each of them is generated by a Hamiltonian PDE. For type (3) and (4), we have more explicit expressions for the transformations. In the rest of this section, we will study useful properties of the transformations of types (1) and (2), which will be used in Chapter~\ref{reduction}.

\section{Properties of the flow type (1)}\label{sub1sd}
We study the flow map $\Phi(t)$ determined by the PDE in \eqref{type1flow}. We first recall some properties of  $\Psi(t)$, which is the flow of the PDE:
\begin{align}\label{psei_type1}
u_t = \partial_x (b(t)u)=:\partial_x \mathcal{B}(t)u,
\end{align}
where $b(t)$ is as described in \eqref{type1flow}.  We recall that $\Psi(t)$ has an explicit form:
\begin{align}\label{psdisd2sd2sd}
\Psi(t)h(\varphi,x) = (1+t \beta_x)h(\varphi,x+t \beta(\varphi,x)).
\end{align}

\subsection{Invertibility}
Since $\beta\in \mathfrak{C}_1(i_0)$, it follows from \ref{slsd22sd2sdhw} in Lemma~\ref{properties_weights} that
\begin{align}\label{sdk2sjsjds}
\rVert \beta \rVert^{\Lip(\gamma,\Omega)}_{s}\le_{\mathtt{pe},s} \epsilon\left(1 + \rVert \mathfrak{I}_\delta\rVert^{\Lip(\gamma,\Omega)}_{s+\mu_0} \right).
\end{align}
Specifically, Remark~\ref{Insteadofimbd} implies that $\rVert \beta \rVert^{\Lip(\gamma,\Omega)}_{s_0}\le_{\mathtt{pe}}  \epsilon$. Therefore, Lemma~\ref{change_ofsd2} tells us that there exists $\tilde{\beta}(t)=\tilde{\beta}(t,\varphi,x)$ such that
\begin{align}\label{ftsdsd2}
B^{-1}(\varphi,x) = x+ \tilde{\beta}(t,\varphi,x),\text{ where }B(\varphi,x):=x+t\beta(\varphi,x),
\end{align}
with the estimate
\begin{align}\label{tildebest2}
\rVert \tilde{\beta(t)}\rVert^{\Lip(\gamma,\Omega_1)}_{s}\le_{\mathtt{pe},s}  \rVert \beta\rVert^{\Lip(\gamma,\Omega_1)}_{s+\mu_{0}}, 
\end{align} for some $\mu_0\ge0$, possibly larger than $\mu_0$ in above.
Using \eqref{sdk2sjsjds}, we have
\begin{align}\label{tildebest3}
\rVert \tilde{\beta}(t)\rVert^{\Lip(\gamma,\Omega_1)}_{s}\le_{\mathtt{pe},s} \epsilon\left(1 + \rVert \mathfrak{I}_\delta\rVert^{\Lip(\gamma,\Omega)}_{s+\mu_0} \right), 
\end{align}
for some $\mu_0\ge0$, possibly larger than $\mu_0$ in above.

The invertibility of $\Psi(t)$ follows immediately from the following result:
\begin{lemma}\cite[Section 8.1]{Baldi-Berti-Montalto:KAM-quasilinear-kdv}\label{invsdwsd}
If $\rVert \beta \rVert_{W^{1,\infty}(\mathbb{T}^{\nu+1})}\le \frac{1}{2}$, then the map $\Psi(t)$ in \eqref{psdisd2sd2sd} is invertible and
\[
\Psi(\tau)^{-1}h(\varphi,x) = (1+\tilde{\beta}_x(t,\varphi,x)) h(\varphi,x+\tilde{\beta}(t,{\varphi,x}),\text{ where $\tilde{\beta}$ is as in \eqref{ftsdsd2}.}
\]
\end{lemma}

\begin{lemma}\label{sed2sd}
We have that $\tilde{\beta}(t) \in \mathfrak{C}_1(i_0)$, for each $t\in [0,1]$.
\end{lemma}
\begin{proof}
Note that a fixed $t$ does not play a role in the proof, therefore in what follows, we omit it in the notations.

The lemma can be proven simply expanding $\tilde{\beta}$ by Taylor series. Indeed, \eqref{ftsdsd2} implies that 
\begin{align}\label{sh222}
\tilde{\beta}(\varphi,x)=-\beta(x+\tilde{\beta}(\varphi,x)),
\end{align} hence, using \eqref{elementary_taylor}, we get
\begin{align}
\tilde{\beta}(\varphi,x)&=-\underbrace{\beta(\varphi,x)}_{=:A_1} - \underbrace{\beta_x(\varphi,x)\tilde{\beta}({\varphi,x})}_{=:A_2} - \underbrace{\frac{1}{2}\beta_{xx}(\varphi,x)(\tilde{\beta}(\varphi,x))^2 }_{=:A_3} \nonumber\\
& \quad -  \underbrace{(\tilde{\beta}(\varphi,x))^3\int_0^1\beta_{xxx}(\varphi,x+s\tilde{\beta}(\varphi,x))\frac{(1-s)^2}{2!}ds}_{=:A_4}. \label{sj2ssjsjs22}
\end{align}
Using that $\beta\in \mathfrak{C}_1(i_0)$, let us write
\begin{align}\label{11wesxcxc}
\beta = \epsilon a_{1,\beta}+\epsilon^2 a_{2,\beta}+\epsilon^3 a_{3,\beta}+q_{\beta},
\end{align}
where $a_{i,\beta}$ are of the form in \eqref{homcoeeff} and $q_\beta$ satisfies the estimates \eqref{q_estimate_11} and \eqref{q_estimate_112}.  We look for an expansion of $\tilde{\beta}$:
\begin{align}\label{besd22sd2}
\tilde{\beta} = \epsilon a_{1,\tilde{\beta}}+\epsilon^2 a_{2,\tilde{\beta}}+\epsilon^3 a_{3,\tilde{\beta}}+ q_{\tilde{\beta}},
\end{align}
such that $a_{i,\tilde{\beta}}$ are of the form in \eqref{homcoeeff} and $q_{\tilde{\beta}}$ satisfies the estimates \eqref{q_estimate_11} and \eqref{q_estimate_112}. To simplify the notations, we write for a quantity $Q$, 
\begin{align}\label{notaesd11}
Q=O(\epsilon^k), \text{ for $k\ge 0$, if $\rVert Q \rVert^{\Lip(\gamma,\Omega_1)}_{s}\le_{\mathtt{pe},s} \epsilon^k(1+\rVert \mathfrak{I}_\delta \rVert^{\Lip(\gamma,\Omega_1)}_{s+\mu_0})$ for some $\mu_0>0$}.
\end{align}
 Note that thanks to Lemma~\ref{interpolation_2s} and \eqref{tildebest3}, the term of size $O(\epsilon)$ in \eqref{sj2ssjsjs22} is  from $A_1$, especially $a_{1,\beta}$ in \eqref{11wesxcxc} thus
 \begin{align}\label{a2sd22}
 \epsilon a_{1,\tilde{\beta}}:=-\epsilon a_{1,\beta},
 \end{align} so that 
 \[
 \tilde{\beta} - \epsilon a_{1,\tilde{\beta}} = O(\epsilon^2).
 \]
 Again, we see the terms of size $O(\epsilon^2)$ in \eqref{sj2ssjsjs22}. Defining
 \begin{align}\label{a2sd222}
 \epsilon^2 a_{2,\tilde{\beta}}:=-\epsilon^2\left(a_{2,\beta} + (a_{1,\beta})_xa_{1,\tilde{\beta}} \right),
 \end{align}
 we see that $\tilde{\beta}-\epsilon a_{1,\tilde{\beta}} - \epsilon^2 a_{2,\tilde{\beta}} = O(\epsilon^3)$.
 Similarly, for $a_{3,\tilde{\beta}}$ defined as
  \begin{align}\label{a2sd223}
 \epsilon^3 a_{3,\tilde{\beta}}:=-\epsilon^3\left(a_{3,\beta} + ((a_{2,\beta})_xa_{1,\tilde{\beta}}) +(a_{1,\beta})_x(a_{2,\tilde{\beta}}) +\frac{1}{2}(a_{1,\beta})_{xx}(a_{1,\tilde{\beta}})^2\right),
 \end{align}
we get (using the fact that $q_{\beta}$ satisfies the estimates \eqref{q_estimate_11} and \eqref{q_estimate_112}),
 \begin{align}\label{shdwsd}
\rVert \tilde{\beta}-\epsilon a_{1,\tilde{\beta}} - \epsilon^2 a_{2,\tilde{\beta}} - \epsilon^{3}a_{3,\tilde{\beta}} + q_{\beta}\rVert^{\Lip(\gamma,\Omega_1)}_{s} = O(\epsilon^4)
 \end{align}
 Therefore, writing
 \begin{align}\label{qtilds22}
 q_{\tilde{\beta}}:= - q_{\beta} +(\tilde{\beta}-\epsilon a_{1,\tilde{\beta}} - \epsilon^2 a_{2,\tilde{\beta}} - \epsilon^{3}a_{3,\tilde{\beta}} + q_{\beta}),
 \end{align}
 we have that \eqref{a2sd22}, \eqref{a2sd222} and \eqref{a2sd223} imply that each $a_{i,\tilde{\beta}}$ are of the form in \eqref{homcoeeff} and \eqref{shdwsd} implies that (recalling the notation in \eqref{notaesd11}) $q_{\tilde{\beta}}$ satisfies 
 \begin{align*}
 \rVert q_{\tilde{\beta}} \rVert^{\Lip(\gamma,\Omega_1)}_{s} & \le_{\mathtt{pe},s} \epsilon^{3} + \epsilon\rVert \mathfrak{I}_{\delta}\rVert^{\Lip(\gamma,\Omega_1)}_{s+\mu_0},\\
 \rVert M_x(q_{\tilde{\beta}})-M_{\varphi,x}(q_{\tilde{\beta}})\rVert^{\Lip(\gamma,\Omega)}_{s}& \le_{\mathtt{pe},s} \epsilon^4 + \epsilon^{2}\rVert\mathfrak{I}_\delta\rVert_{s+\mu_0}^{\Lip(\gamma,\Omega_1)}.
 \end{align*}
To finish the proof that $\tilde{\beta}\in \mathfrak{C}_1(i_0)$, it suffices to prove the estimates  in \eqref{q_estimate_11} and \eqref{q_estimate_112} for $d_iq_{\tilde{\beta}}$. Since $d_iq_{\beta}$ satisfies such estimates, and $a_{i,\tilde{\beta}}$ are independent of $i_0$, we see from \eqref{qtilds22} that  it suffices to show that
\begin{equation}\label{besjsdhsd}
\begin{aligned}
\rVert d_i \tilde{\beta}(i)[\ihat] \rVert_s &\le_{\mathtt{pe},s} \epsilon\left( \rVert \ihat \rVert_{s+\mu_0} + \rVert \mathfrak{I}_\delta \rVert_{s+\mu_0}\rVert \ihat \rVert_{s_0+\mu_0}\right),\\
  \rVert d_i(M_x(\tilde{\beta})-M_{\varphi,x}(\tilde{\beta})(i_0))[\ihat]\rVert_s&\le_{\mathtt{pe},s} \epsilon^2\left(\rVert \ihat \rVert_{s+\mu_0} + \rVert \mathfrak{I}_\delta \rVert_{s+\mu_0}\rVert \ihat \rVert_{s_0+\mu_0} \right).
\end{aligned}
\end{equation}
Note that the fact that $a_{i,\beta}$ is independent of $i_0$ implies that $\beta$ satisfies the same estimates in \eqref{besjsdhsd}, for $\tilde{\beta}$, replaced by $\beta$. 
 Towards the proof, we denote
 \begin{align}\label{notysdx}
 \mathcal{A}_{\tilde{\beta}}[h](\varphi,x):=h(\varphi,x+\tilde{\beta}(\varphi,x)),\text{ for $h\in C^\infty_{\varphi,x}$},
 \end{align}
 so that \eqref{ftsdsd2} (neglecting $t$, and denoting the dependence on $i_0$) implies
 \begin{align}\label{sdk2j2}
\tilde{\beta}(i_0)= - \mathcal{A}_{\tilde{\beta}}(i_0)[\beta(i_0)],\text{ and }\mathcal{A}_{\tilde{\beta}}(i_0)[\beta_x(i_0)] = -\frac{\tilde{\beta}(i_0)}{1+\tilde{\beta}(i_0)}.
 \end{align}
  Since $d_i(\mathcal{A}_\beta(i_0)[\ihat])[h](\varphi,x) = d_i\beta(i_0)[\ihat](\varphi,x)\mathcal{A}_{\beta}(i_0)[h_x](\varphi,x)$, which can be seen easily from \eqref{notysdx}, it follows from \eqref{sdk2j2} that 
 \begin{align*}
 d_i\tilde{\beta}(i_0)[\ihat] & = -d_i(\mathcal{A}_{\tilde{\beta}}(i_0)[\ihat])[\beta] -\mathcal{A}_{\tilde{\beta}}(i_0)[d_i\beta(i_0)[\ihat]] \\
 & = - d_{i}\tilde{\beta}(i_0)[\ihat]\mathcal{A}_{\tilde{\beta}}(i_0)[\beta_x(i_0)]- \mathcal{A}_{\tilde{\beta}}(i_0)[d_i\beta(i_0)[\ihat]].
 \end{align*}
 Therefore, we get
 \[
 d_i\tilde{\beta}(i_0)[\ihat] = -\frac{\mathcal{A}_{\tilde{\beta}}[d_i\beta(i_0)[\ihat]]}{1+\mathcal{A}_{\tilde{\beta}(i_0)}[\beta_x(i_0)]}.
 \]
 Since $\frac{1}{1+\mathcal{A}_{\tilde{\beta}}(\beta_x)} - 1 = O(\epsilon^2)$, and  $d_i\beta(i_0)[\ihat]$ satisfies \eqref{besjsdhsd}, we have that $\tilde{\beta}$ satisfies \eqref{besjsdhsd} as well.
 \end{proof}

\subsection{Adjoint operator}
Now, we study the adjoint operators \index{Adjoint operator}$\Psi(t)^{T}$ and $\left(\Psi(t)^{-1}\right)^T$. From the explicit expressions of $\Psi(t),\Psi(t)^{-1}$ in \eqref{psdisd2sd2sd} and Lemma~\ref{invsdwsd}, it follows straightforwardly from the change of variables that (see \cite[Section 8.1]{Baldi-Berti-Montalto:KAM-quasilinear-kdv}),
\begin{align}\label{adjoints}
\Psi(t)^{T}h(\varphi,x) = h(\varphi, x+\tilde{\beta}(t,\varphi,x)),\quad \left(\Psi(t)^{-1}\right)^{T}h(\varphi,x) = h(\varphi,x+t{\beta}(\varphi,x)).
\end{align}\subsection{Expansion of $\Psi(t)$}

\begin{lemma}\label{expasd111}
Let $\psi(t)\in \left\{\Psi(t),\Psi(t)^T,\Psi(t)^{-1},\left(\Psi(t)^T\right)^{-1}\right\}$. For each $\psi(t)$, there exist a symbol $\mathfrak{b}\in \mathfrak{B}^1_1$ (see definition~\ref{buildingbobo}) and an operator $R_{\psi}(t)$ such that
\begin{align}\label{exdsd21}
\psi(t) = I + \epsilon\tau Op^W(\mathfrak{b}) + R_{\psi}(\tau),
\end{align}
where $\mathcal{R}(t)\in \left\{R_{\Psi}(t),R_{\Psi^*}(t),R_{\Psi^{-1}}(t),R_{\left(\Psi^{-1}\right)^{T}}(t) \right\}$ satisfies
\begin{align}
\sup_{t\in [0,1]}\rVert \mathcal{R}(t)h\rVert_s^{\Lip(\gamma,\Omega_1)}  & \le_{\mathtt{pe},s} \left(\epsilon^2 \rVert h\rVert_{s+\mu_0}^{\Lip(\gamma,\Omega_1)} + \epsilon^2\rVert \mathfrak{I}_\delta \rVert_{s+\mu_0}^{\Lip(\gamma,\Omega_1)}\rVert h\rVert_{s_0+\mu_0}^{\Lip(\gamma,\Omega_1)}  \right),\label{rest1sd11}\\
\sup_{t\in [0,1]}\rVert d_i\mathcal{R}(t)(i_0)[\ihat]h\rVert_{s_0}&\le_{\mathtt{pe},s}   \epsilon \rVert \ihat \rVert_{s_0+\mu_0} \rVert h \rVert_{s_0+\mu_0}.\label{rest1sd112}
\end{align}
\end{lemma}
\begin{proof}
We will prove the case where $\psi(t)=\Psi(t)$ only. The proof is based on the explicit expression of $\Psi(t)$ in \eqref{psdisd2sd2sd} and the fact that $\beta\in \mathfrak{C}_1(i_0)$. The other operators $\Psi(t)^{-1},\Psi(t)^{T}$ and $\left(\Psi(t)^{T}\right)^{-1}$ can be proved in the same way using the expressions in \eqref{adjoints}, Lemma~\ref{invsdwsd} and Lemma~\ref{sed2sd}.

 To show \eqref{rest1sd11}, using the Taylor expansion in \eqref{elementary_taylor},  we write (see \eqref{psdisd2sd2sd})
 \begin{align}
 \Psi(t)h &= h(\varphi,x+t\beta(\varphi,x)) + t\beta_x(\varphi,x)h(\varphi,x + t\beta(\varphi,x))\nonumber\\
 & = h(\varphi,x) + t\beta h_x(\varphi,x) + (t\beta)^2\int_0^{1}h_{xx}(\varphi,x+st\beta(\varphi,x))(1-s)ds\nonumber\\
 & \ + t\beta_x\left(h(\varphi,x) + t\beta\int_0^{1}h_x(\varphi,x+st\beta(\varphi,x))ds \right)\nonumber\\
 & =: h(\varphi,x) + t\partial_x (\beta h ) +R_{\Psi,1}(t)[h], \label{sj2nsdsd2}
 \end{align}
 where 
 \begin{align}\nonumber
 R_{\Psi,1}(t)[h]& :=(t\beta)^2\int_0^{1}h_{xx}(\varphi,x+st\beta(\varphi,x))(1-s)ds \\
 & + t^2\beta_x\beta\int_0^{1}h_x(\varphi,x+st\beta(\varphi,x))ds.
 \label{dsd2sd}
 \end{align}
 Furthermore, using that $\beta\in \mathfrak{C}_1(i_0)$, we have an expansion
 \[
 \beta = \sum_{i=1}^{3}\epsilon^i a_i + q,
 \]
 where $a_i,q$ are as described in Definition~\ref{homogeneous_expansion1}. Hence, we can write from \eqref{sj2nsdsd2} that
 \begin{align}\label{qisequalt}
 \Psi(t) h = h(\varphi,x) + t\epsilon\partial_x(a_1 h) + \underbrace{t\partial_x\left(\left(\sum_{i=2}^3\epsilon^i a_i + q\right)h\right)}_{=:R_{\Psi,2}(t)[h]} + R_{\Psi,1}(t)[h].
 \end{align}
 Since $a_1$ is of the form in \eqref{homcoeeff} for $i=1$, we can find $\mathfrak{b}\in \mathfrak{B}^1_1$ such that $\partial_x(a_1 h)=: Op^W(\mathfrak{b}_1)h$, therefore we have
 \[
 \Psi(t) = I + \epsilon t Op^W(\mathfrak{b}) + R_{\Psi}(t),\quad R_{\Psi}(t):=R_{\Psi,1}(t)+R_{\Psi,2}(t).
 \]
 
 Now, it suffices to show that $R_{\Psi}$ satisfies the estimates in \eqref{rest1sd11} and \eqref{rest1sd112}.
 For $R_{\Psi,1}$ defined in \eqref{dsd2sd},  using \eqref{interpolation_2s} and \ref{sd2sd} of Lemma~\ref{change_ofsd2},  we have that for some $\mu_0\ge0$, 
 \begin{align*}
 \rVert R_{\Psi,1}(t)[h] \rVert^{\Lip(\gamma,\Omega_1)}_{s} & \le_{\mathtt{pe},s}  \rVert \beta\rVert^{\Lip(\gamma,\Omega_1)}_{s+\mu_0}\rVert \beta \rVert^{\Lip(\gamma,\Omega_1)}_{s_0+\mu_0} \rVert h \rVert^{\Lip(\gamma,\Omega_1)}_{s_0+\mu_0} \\
 & \quad + \left(\rVert \beta \rVert^{\Lip(\gamma,\Omega_1)}_{s_0+\mu_0}\right)^2\rVert h \rVert^{\Lip(\gamma,\Omega_1)}_{s+\mu_0}. 
 \end{align*}
 Using \eqref{sdk2sjsjds} and \eqref{size_assumption_2}, we obtain that $R_{\Psi,1}$ satisfies \eqref{rest1sd11}. To show \eqref{rest1sd112} for $R_{\Psi,1}$, we differentiate \eqref{dsd2sd} in $i_0$ to see that
\begin{align*}
\rVert d_i(R_{\Psi}(t)(i_0)[\ihat])h\rVert_{s_0}& \le_{\mathtt{pe},s}  \rVert d_i\beta(i_0)[\ihat]\rVert_{s_0+\mu_0}\rVert \beta \rVert_{s_0+\mu_0}\rVert h\rVert_{s_0+\mu_0}\\& \le \epsilon^2 \rVert \ihat \rVert_{s_0+\mu_0} \rVert h \rVert_{s_0+\mu_0},
\end{align*}
where the last inequality follows from $\beta\in \mathfrak{C}_1(i_0)$, 
which proves \eqref{rest1sd112} (with actually a higher power on $\epsilon$).
For $R_{\Psi,2}$ in \eqref{qisequalt}, we use that $q$ satisfies the estimates in \eqref{q_estimate_11} and $a_i$ is of the form in \eqref{homcoeeff}, then the result follows straightforwardly.

\end{proof}

\subsection{Finite dimensional operators}
\begin{lemma}\cite[Lemma A.4]{Feola-Giuliani:quasiperiodic-water-waves}\label{fhsjdwd2sd}
Let $\Phi(t),\Psi(t)$ be the flow maps determined by the PDEs in \eqref{type1flow} and \eqref{psei_type1}, respectively. Let $R(t):=\Phi(1)\Pi_{S^\perp} - \Pi_{S^\perp}\Psi(1)\Pi_{S^\perp}$. Then $R(t)$ is a finite dimensional operator, more precisely, $R(t)\in \mathfrak{R}(i_0)$ (see Definition~\ref{fisnsdw}).
\end{lemma}
\begin{proof}
As stated in \cite[Lemma A.4]{Feola-Giuliani:quasiperiodic-water-waves}, $R(t):=\Phi(t)\Pi_{S^\perp} - \Pi_{S^\perp}\Psi(t)\Pi_{S^\perp}$ is a finite dimensional operator. We aim to show that $R(t)\in \mathfrak{R}(i_0)$, following the strategy of \cite[Lemma C.1]{Feola-Giuliani-Procesi:kam-tori-degasperis-procesi}. We will show $R(1)\in \mathfrak{R}(i_0)$, since one can straightforwardly show that $R(1)\in \mathfrak{R}(i_0)$ implies  $R(t)\in \mathfrak{R}(i_0)$, by reparametrizing the variable $t$.

We define a operator $U(t)$ by
\begin{align}\label{usdwsd}
\frac{d}{d\tau} U(t) = -(\Psi(\tau))^{-1}\partial_x\mathcal{B}(t)\Pi_{S}\Psi(t)U(t),\, U(0)=I,\text{ where $\partial_x\mathcal{B}(t)$ is as in \eqref{psei_type1}}.
\end{align}
One can easily check that $
\Phi(t) = \Pi_{S^\perp}\Psi(t)U(t). $ Indeed, the left hand side, $\Pi_{S^\perp}\Psi(t)U(t)$, solves
\begin{align*}
\frac{d}{dt}\left(\Pi_{S^\perp}\Psi(t)U(t) \right) &= \Pi_{S^\perp}\left(\frac{d}{dt}\Psi(t)\right)U(t) + \Pi_{S^\perp}(\Psi(t))\left(\frac{d}{dt}U(t)\right)\\
& = \Pi_{S^\perp}\left(\partial_x \mathcal{B}(t)\Psi(t) - \Psi(t)\Psi(t)^{-1}\partial_x\mathcal{B}(t)\Pi_{S}\Psi(t)U(t) \right)\\
& = \Pi_{S^\perp}\partial_x \mathcal{B}(t)\Pi_{S^\perp}\Psi(t)U(t),
\end{align*}
which solves the same equation \eqref{type1flow} as  $\Phi(t)$.
Therefore, we have 
\begin{align}\label{s2sd2sj2sd2}
\Phi(t)\Pi_{S^\perp} = \Pi_{S^\perp}\Psi(t)\Pi_{S^\perp} +{ \Pi_{S^\perp}\Psi(t)(U(t)-I)\Pi_{S^\perp}}.
\end{align}
Therefore, it suffices to show that $\Pi_{S^\perp}\Psi(1)(U(1)-I)\Pi_{S^\perp}\in \mathfrak{R}(i_0)$. Towards the proof, let us denote $Z(t):=\partial_x \mathcal{B}(t)\Pi_{S}\Psi(t)$, so that the equation for $U(t)$ in \eqref{usdwsd} can be written as
\begin{align}\label{Useqxxzz}
\frac{d}{dt}U(t) = -(\Psi(t))^{-1}Z(t)U(t).
\end{align}
For the operator $Z(t)$, we have that
\begin{align}
Z(t)u & = \partial_x\mathcal{B}(t)\Pi_{S}\Psi(t)u = \partial_x\mathcal{B}(t)\left(\sum_{j\in S}(\Psi(t)u,e^{\ii j x })_{L^2(\mathbb{T})}e^{\ii j x} \right) \nonumber \\
& = \sum_{j\in S}(u,\underbrace{\Psi(t)^*[e^{\ii j x}]}_{=: g_j(t)})_{L^2(\mathbb{T})}\underbrace{\partial_x\mathcal{B}(t)[e^{\ii j x}]}_{=:\chi_j(t)}. \label{psisd2sdsd2}
\end{align}
Now, we write $U(t)$ as (for $u\in H_{S^\perp}$),
\begin{align*}
U(t)u-u &\overset{\eqref{Useqxxzz}}= \int_0^{t}\frac{d}{d\tau}U(\tau)ud\tau = \int_0^{t}-(\Psi(\tau))^{-1}Z(\tau)U(\tau)ud\tau\\
& \overset{\eqref{psisd2sdsd2}}= -\sum_{j\in S}\int_0^{t} (U(\tau)u,g_j(\tau))_{L^2(\mathbb{T})}(\Psi(\tau))^{-1}[\chi_j(\tau)] d\tau\\
& = -\sum_{j\in S}\int_{0}^{t}(u,{U(\tau)^{T}g_j(\tau)})_{L^2(\mathbb{T})}{(\Psi(\tau))^{-1}[\chi_j(\tau)]} d\tau.
\end{align*}
Therefore, we obtain  (using $g_j$ and $\chi_j$ in \eqref{psisd2sdsd2} and  $u\in H_{S^\perp}$),
\begin{equation}\label{sdjwsdsd}
\begin{aligned}
\Pi_{S^\perp}\Psi(1)(U(1) - I)u &= \Pi_{S^\perp}\sum_{j\in S}\int_0^{1} (u,\tilde{g}_j(t))\tilde{\chi}_j(t)dt,\\
\tilde{g}_j(t):=\Pi_{S^\perp}U(t)^{T}\Psi(t)^*[e^{\ii j x}],&\quad\tilde{\chi}_j(t):= -\Psi(1)\Psi(t)^{-1}\partial_x\mathcal{B}(t)[e^{\ii j x}], \, \text{ for $j\in S$}.
\end{aligned}
\end{equation}

To finish the proof, we need to estimate $\tilde{g}_j,\tilde{\chi}_j$. We need the following lemma.
\begin{lemma}\label{k22sd}
$\partial_x\mathcal{B}(t),U(t)^T$ have  expansions as in \eqref{exdsd21}. That is, there exist  symbols $\mathfrak{b}_{\mathcal{B}},\mathfrak{b}_U\in \mathfrak{B}^1_1$  and operators $R_{\mathcal{B}}(t),R_{U}(t)$ such that 
\begin{align*}
\partial_x\mathcal{B}(t) &= \epsilon Op^W(\mathfrak{b}_{\mathcal{B}}) + R_{\mathcal{B}}(t),\\
U(t)^T &= I + \epsilon t Op^W(\mathfrak{b}_U)\Pi_{S} + R_U(t),
\end{align*}
and $R_{\mathcal{B}}(t),R_{U}(t)$ satisfy the estimates in \eqref{rest1sd11} and \eqref{rest1sd112}. 
\end{lemma}
\begin{proof}
The results can be proved in the same way as in the proof of Lemma~\ref{expasd111}, using the Taylor expansions of the operators $\mathcal{B}(t)$ and $U(t)^T$.
\end{proof}

Thanks to Lemma~\ref{expasd111} and  Lemma~\ref{k22sd}, we have $\mathfrak{b}_1,\mathfrak{b}_2\in \mathfrak{B}^1_1$, operators $R_g(\tau),\ R_{\chi}(\tau)$ and functions $f_1(t),f_2(t)$ that depend on $t$ only such that  such that
\begin{align*}
U(t)^T\Psi(t)^*[e^{\ii j x}] &= e^{\ii j x} + \epsilon f_1(t)Op^W(\mathfrak{b}_1)[e^{\ii j x}] + R_g(\tau)[e^{\ii j x}],\\
-\Psi(1)\Psi(t)^{-1}\partial_x\mathcal{B}(t)[e^{\ii j x}] &= e^{\ii j x} + \epsilon f_2(t) Op^W(\mathfrak{b}_2)[e^{\ii j x}] + R_{\chi}(\tau)[e^{\ii j x}],
\end{align*}
where $ R_g,R_{\chi}$ satisfy the estimates \eqref{rest1sd11} and \eqref{rest1sd112}. Plugging this into \eqref{sdjwsdsd}, we obtain
\begin{align*}
& \Pi_{S^\perp}\Psi(1)(U(1) - I)u \\
&= \epsilon^2\underbrace{\Pi_{S^\perp}\sum_{j\in S} (u,Op^W(\mathfrak{b}_1)[e^{\ii j x}])_{L^2(\mathbb{T})}Op^W(\mathfrak{b}_2)[e^{\ii j x}]\left(\int_0^{1}f_1(t)f_{2}(t)dt\right)}_{=:\mathcal{R}_1[u]}\\
& \ + \underbrace{\Pi_{S^\perp}\sum_{j\in S}\int_0^{1}(u,\epsilon f_1(t)Op^W(\mathfrak{b}_1)[e^{\ii j x}])_{L^2(\mathbb{T})}R_{\chi}(t)[e^{\ii j x}]dt}_{\mathcal{R}_{2,1}[u]}\\
& \ + \underbrace{\Pi_{S^\perp}\sum_{j\in S}\int_0^1 (u,R_g(\tau)[e^{\ii j x}])_{L^2(\mathbb{T})}\left(\epsilon f_2(t)Op^{W}(\mathfrak{b}_2))+R_{\chi}(t)\right)[e^{\ii j x}]dt}_{=:\mathcal{R}_{2,2}[u]}.
\end{align*}
Using that  $R_g,R_{\chi}$ satisfy the estimates \eqref{rest1sd11} and \eqref{rest1sd112}, it follows straightforwardly that $\mathcal{R}_{2,1},\mathcal{R}_{2,2}$ are of the form  \eqref{j2j2j2j2j2jssds} with the estimates \eqref{sjdjsjsj}. We also have $\mathcal{R}_1$ is of the form  \eqref{fisnsdw}. Therefore, $\Pi_{S^\perp}\Psi(1)(U(1) - I)$ in \eqref{s2sd2sj2sd2} satisfies $\Pi_{S^\perp}\Psi(1)(U(1) - I)\in \mathfrak{R}(i_0)$. This finishes the proof.
\end{proof}

\begin{lemma}\label{reversibility}
 If $\beta(-\varphi,-x)=-\beta(\varphi,x)$, then $\Phi(t)$ is a reversibility preserving operator. 
\end{lemma}
\begin{proof}
From $b(t)$ in \eqref{type1flow}, it is clear that $b(t)(-\varphi,-x) =- b(t)(\varphi,x)$ for each $t$. In the equation \eqref{type1flow}, this implies that each space $X,Y$ is invariant under the vector field
\[
u\mapsto \Pi_{S^\perp}\partial_x(b(t)\Pi_S^{\perp}u).
\]
Therefore, the flow map is reversibility preserving. 
\end{proof}

\section{Properties of the flow type (2)}\label{study_hard_learn_more}
We study the flow map $\Phi(t)$ determined by the PDE in \eqref{type1flow2s}. We first recall some properties of  $\Psi(t)$, which is the flow of the PDE:
\begin{equation}\label{psei_type12}
\begin{aligned}
u_t &= \partial_x Op^W(\mathfrak{a})u=:\partial_x\mathcal{A}u,
\end{aligned}
\end{equation}
for some real-valued symbol $\mathfrak{a}\in \mathcal{S}^m$, for some $m<0$,
that is 
\begin{align}\label{sdkksdsdsd0293}
\frac{d}{dt}\Psi(\tau) &= \partial_x \mathcal{A}\Psi(t), \quad \Psi(0) = I.
\end{align}
Since the symbol $\mathfrak{a}$ is independent of $t$, the equation \eqref{sdkksdsdsd0293} is autonomous, therefore, we have
\begin{align}\label{commuetesd}
\Psi(t)^{-1} = \Psi(-t),\quad \Psi(t)\partial_x\mathcal{A} = \partial_x \mathcal{A}\Psi(t),\quad \Psi(t_1+t_2) = \Psi(t_1)\Psi(t_2).
\end{align}
Again, \eqref{adjoint_weyl} tells us that $\mathcal{A}$ is symmetric. In the following proposition, we collect tame estimates of the flow map $\Psi(t)$. The proof can be done following the same strategy in \cite[Appendix A]{Berti-Montalto:quasiperiodic-standing-gravity-capillary} word
by word. 

\begin{proposition}\cite[Appendix A]{Berti-Montalto:quasiperiodic-standing-gravity-capillary}\label{wjsd2sdsdesd}
Let  us fix the constants $\mathtt{b}_0,\mathtt{S},k$ so that
\begin{align}\label{bsd2sd}
\mathtt{b}_0\ge 0,\ \mathtt{S}\in(s_0,\infty),\ |k|\le \mathtt{b}_0 .
\end{align} 
We denote
\[
\Psi_k(\tau) := |D|^k \Psi(\tau)|D|^{-k}.
\]
There exist $\mu_0=\mu_0(\mathtt{p}, \mathtt{b}_0),\eta_0=\eta_0(\mathtt{b}_0)$ and $\delta=\delta(\mathtt{S})>0$, such that if,
\begin{align}\label{smallsdsimpsd}
|\mathfrak{a}|^{\Lip(\gamma,\Omega)}_{m,s_0+\mu_0,\eta_0}\le \delta,
\end{align}
the following hold:
\begin{enumerate}[label=(\arabic*)]
\item \label{swlsWkgoqhwk}For all $s,\mathtt{b}\in \mathbb{R}$ such that $s\in [0,\mathtt{S}]$ and $0\le \mathtt{b}\le s_0+\mathtt{b}_0$, it holds that
\begin{align}
\sup_{\tau\in [0,1]}\rVert \partial_{\varphi_j}^{\mathtt{b}}(\Psi_k(\tau))|D|^{-\mathtt{b}(1+m)}h\rVert_{s}&\le_{s,\mathtt{b}_0,m} \rVert h\rVert_{s} + |\mathfrak{a}|_{m,s+\mu_0,\eta_0}\rVert h \rVert_{s_0}, \label{sdsdsdkkk111111}\\
\sup_{\tau\in [0,1]}\rVert |D|^{-\mathtt{b}(1+m)}\partial_{\varphi_j}^{\mathtt{b}}(\Psi_k(\tau))h\rVert_{s}&\le_{s,\mathtt{b}_0,m} \rVert h\rVert_{s} + |\mathfrak{a}|_{m,s+\mu_0,\eta_0}\rVert h \rVert_{s_0} \label{sdsdsdkkk1111112}
\end{align}
\item \label{sddwsdsdsd}Given $\omega_1,\omega_2\in \Omega$,  we denote
   \[
   \Delta_{12}A(\omega):=\frac{{A(\omega_1)-A(\omega_2)}}{|\omega_1-\omega_2|},
   \]
   for an $\omega$-dependent object $A$ (symbols, operators,... etc). For all $s,\mathtt{b}\in \mathbb{R}$ such that $s\in [0,\mathtt{S}]$ and $0\le \mathtt{b}\le s_0+ \mathtt{b}_0$, it holds that
   \begin{align}
& \gamma\sup_{\omega_1,\omega_2\in \Omega, \tau\in [0,1]}\rVert \partial_{\varphi_j}^{\mathtt{b}}(\Delta_{12}\Psi_k(\tau))|D|^{-(1+\mathtt{b})(1+m)}h\rVert_{s}\nonumber \\&\le_{s,\mathtt{b}_0,m} \rVert h\rVert_{s} + |\mathfrak{a}|^{\Lip(\gamma,\Omega)}_{m,s+\mu_0,\eta_0}\rVert h \rVert_{s_0}, \label{sdsdsdkkk111113}\\
& \gamma\sup_{\omega_1,\omega_2\in \Omega, \tau\in [0,1]}\rVert |D|^{-(1+\mathtt{b})(1+m)}\partial_{\varphi_j}^{\mathtt{b}}(\Delta_{12}\Psi(\tau))h\rVert_{s} \nonumber \\&\le_{s,\mathtt{b}_0,m} \rVert h\rVert_{s} + |\mathfrak{a}|^{\Lip(\gamma,\Omega)}_{m,s+\mu_0,\eta_0}\rVert h \rVert_{s_0}, \label{sdsdsdkkk1121113}
\end{align}
\item \label{rlajsd2sdsd}In case $\mathfrak{a}$ depends on embeddings $i$, we have that for all $s,\mathtt{b}\in \mathbb{R}$ such that $s\in [0,\mathtt{S}]$ and $0\le \mathtt{b}\le_{s,\mathtt{b}_0,m} s_0+\mathtt{b}_0$, it holds that
\begin{align*}
& \sup_{\tau\in [0,1]}\rVert \partial_{\varphi_j}^{\mathtt{b}}\left(d_i\Psi_k(\tau)(i)[\ihat]\right)|D|^{-(1+b)(1+m)}\rVert_{s} \\
&  \le_{s,\mathtt{b}_0,m} |d_i\mathfrak{a}(i)[\ihat]|_{m,s_0+\mu_0,\eta_0}\rVert h \rVert_{s} +|d_i\mathfrak{a}(i)[\ihat]|_{m,s+\mu_0,\eta_0}\rVert h \rVert_{s_0} \\
& \sup_{\tau\in [0,1]}\rVert  |D|^{-(1+b)(1+m)}\partial_{\varphi_j}^{\mathtt{b}}\left(d_i\Psi_k(\tau)(i)[\ihat]\right)\rVert_{s} \\
&\le_{s,\mathtt{b}_0,m}|d_i\mathfrak{a}(i)[\ihat]|_{m,s_0+\mu_0,\eta_0}\rVert h \rVert_{s} +|d_i\mathfrak{a}(i)[\ihat]|_{m,s+\mu_0,\eta_0}\rVert h \rVert_{s_0}.
\end{align*}
\end{enumerate} 
\end{proposition}

\begin{lemma}\label{tnsjdsinver}
Under the assumption \eqref{smallsdsimpsd}, the operators $\Psi(\tau)^{-1}, \Psi(\tau)^{T},\left(\Psi(\tau)^{-1}\right)^{T}$ satisfy the same estimates as in Proposition~\ref{wjsd2sdsdesd}.
\end{lemma}
\begin{proof}
Since the evolution equation \eqref{psei_type12} is autonomous,  $\Psi(-\tau)=\Psi(\tau)^{-1}$ solves the same equation \eqref{psei_type12} with $\mathcal{A}\to-\mathcal{A}$. Therefore, Propsition~\ref{wjsd2sdsdesd} applies to $\Psi(\tau)^{-1}$. For $\Psi(\tau)^T$, it follows from \eqref{symplectic73} that $\Psi(\tau)^{T} = \partial_x^{-1}\Psi(\tau)^{-1}\partial_x$. Therefore, using that $\Psi(\tau)^{-1}$ satisfies the estimates in the proposition, it follows straightforwardly, noticing that $|D|^{-1}\partial_x,\ |D|\partial_x^{-1}$ are isomorphisms in between $H^s_x$.
\end{proof}

\subsection{Finite dimensional operator}
As in Lemma~\ref{fhsjdwd2sd}, we will study the smoothing effect of the operator $\Phi(1)\Pi_{S^\perp} - \Pi_{S^\perp}\Psi(1)\Pi_{S^\perp}$. 

\begin{lemma}\label{finitesd2sd}
$\Phi(1)\Pi_{S^\perp} - \Pi_{S^\perp}\Psi(1)\Pi_{S^\perp}$ is a finite dimensional operator. More precisely,  if the symbol $\mathfrak{a}\in \mathcal{S}^m$ in \eqref{psei_type12} satisfies
\begin{align}\label{sdsdsdaaaaasdkk}
\mathfrak{a} =\epsilon \mathfrak{b}_1 + \epsilon^2\mathfrak{b}_2 + \mathfrak{q},
\end{align}
for some $\mathfrak{b}_1\in \mathfrak{B}^m_1$, $\mathfrak{b}_2\in \mathfrak{B}^m_2$ (see Definition~\ref{buildingbobo}) and $\mathfrak{q}$ such that
\begin{equation}\label{qdsd2sd}
\begin{aligned}
|\mathfrak{q}|^{\Lip(\gamma,\Omega_1)}_{m,s,0}&\le_{\mathtt{pe},s}\gamma^{-1}\left( \epsilon^5+\epsilon^3\rVert \mathfrak{I}_\delta\rVert^{\Lip(\gamma,\Omega_1)}_{s+\mu_0}\right),\\
|d_i\mathfrak{q}(i_0)[\ihat]|_{m,s,0}&\le_{\mathtt{pe},s}\epsilon^3\gamma^{-1}\left( \rVert \ihat \rVert_{s+\mu_0} + \rVert \mathfrak{I}_\delta\rVert_{s+\mu_0}\rVert \ihat \rVert_{s_0+\mu_0}\right).
\end{aligned}
\end{equation}
Then,  $\Phi(1)\Pi_{S^\perp} - \Pi_{S^\perp}\Psi(1)\Pi_{S^\perp}\in \mathfrak{R}(i_0)$.
\end{lemma}
\begin{proof}
 Note that assuming the expansion of $\mathfrak{a}$ in \eqref{sdsdsdaaaaasdkk}, the transformation $\Psi(\tau)$ admits the same expansion as in Lemma~\ref{expasd111}. Then, the proof of the lemma can be proved by following the same argument word by word as in the proof of Lemma~\ref{fhsjdwd2sd}.
\end{proof}

\subsection{Approximate solutions of a Heisenberg equation}
Let $m:=1-\alpha <0$ (see \eqref{psei_type12} for $m$) and let us consider\index{Heisenberg equation} a real-valued symbol $\mathfrak{p}(x,\xi)=\mathfrak{p}(\omega,\varphi,x,\xi)\in \mathcal{S}^{m'}$ for some $m'\in \mathbb{R}$ and  $P(\tau$), defined as a solution to a Heisenberg equation:
\begin{align}\label{psd2kjde}
\frac{d}{d\tau}P(\tau) = [P(\tau),\mathcal{A}]_x,\quad P(0)=Op^W(\mathfrak{p}), \text{ where $
[A,B]_x := A\partial_x B - B\partial_x A$.}
\end{align} In the rest of this section, we will aim to find an approximate solution to the equation \eqref{psd2kjde}, using the argument introduced in \cite{Berti-Montalto:quasiperiodic-standing-gravity-capillary}. In the following proposition, we will use  the operator $\star$, defined in \eqref{syysdsdsd}.

\begin{proposition}\label{apposdsd}\cite[Lemma 6.15, 6.16]{Berti-Montalto:quasiperiodic-standing-gravity-capillary}
Let $P(\tau)$ be a solution to the Heisenberg equation \eqref{psd2kjde} and let $\mathtt{N}_\alpha\in \mathbb{N}$ be  fixed. Then we have an expansion $P(\tau) = Q(\tau) + W(\tau)$ such that
\begin{enumerate}[label=(\arabic*)]
\item \label{qsdisasymbol} $Q(\tau)$ is a pseudo differential operator with a symbol $\mathfrak{q}_{\le\mathtt{N}_\alpha}(\tau)$, that is $Q(\tau) = Op^W(\mathfrak{q}_{\le \mathtt{N}_\alpha}(\tau))$. The symbol $\mathfrak{q}_{\le \mathtt{N}_\alpha}$ is given by
\begin{align}\label{defsdq1sd}
\mathfrak{q}_{\le \mathtt{N}_\alpha}(\tau) &= \sum_{n=0}^{\mathtt{N}_\alpha} \frac{1}{n!}\tau^{n}\mathfrak{q}_n,\quad \mathfrak{q}_0:=\mathfrak{p},\quad \mathfrak{q}_{n}:=\mathfrak{q}_{n-1}\star \mathfrak{a},\text{ for $n=1,\ldots, \mathtt{N}_\alpha$.}
\end{align}
\item \label{sdsdsjjjsdd} For $W(\tau)$, we have
\begin{align}\label{jsd2sdxxq123s}
 W(\tau)  = \frac{1}{\mathtt{N}_\alpha!}\int_0^\tau \Psi(\tau-t)^{T}Op^W(\mathfrak{q}_{\mathtt{N}_\alpha}\star \mathfrak{a})\Psi(\tau-t)t^{\mathtt{N}_\alpha}dt,
\end{align}
\end{enumerate}
\end{proposition}
\begin{proof}
For a fixed constant $\mathtt{N}_\alpha\in \mathbb{N}$, we see that the sequence of symbols $\mathfrak{q}_n$ defined in \eqref{defsdq1sd} satisfies
 \begin{align}\label{qtsd2sdseq}
 \mathfrak{q}_n\in \mathcal{S}^{m'+n(1-\alpha)}.
 \end{align}
Indeed, \eqref{qtsd2sdseq} follows from \ref{symocsdwd} in Lemma~\ref{compandkskd2sd}. 
  Now we derive the expression for $W(\tau) = P(\tau) -Q(\tau)$.
 
We first observe that the symbol $\mathfrak{q}_{\le \mathtt{N}_\alpha}(\tau)$ satisfies (from  \eqref{defsdq1sd}),
\[
\frac{d}{d\tau}\mathfrak{q}_{\le \mathtt{N}_\alpha}(\tau) = \sum_{n=1}^{\mathtt{N}_\alpha}\mathfrak{q}_{n}\frac{\tau^{n-1}}{(n-1)!} = \sum_{n=0}^{\mathtt{N}_\alpha-1}\left(\mathfrak{q}_{n}\star\mathfrak{a}\right)\frac{\tau^{n}}{n!} = \mathfrak{q}_{\le \mathtt{N}_\alpha}(\tau)\star \mathfrak{a} - \frac{\tau^{\mathtt{N}_\alpha}}{\mathtt{N}_\alpha!}\mathfrak{q}_{\mathtt{N}_\alpha}\star\mathfrak{a}
\] 
Hence, $Q(\tau)$ solves 
\begin{align}\label{qeasd2sdsd}
\frac{d}{d\tau}Q(\tau) = [Q(\tau),Op^W(\mathfrak{a})]_x - \frac{\tau^{\mathtt{N}_\alpha}}{\mathtt{N}_\alpha!}Op^W(\mathfrak{q}_{\mathtt{N}_\alpha}\star\mathfrak{a}),\quad Q(0)= Op^W(\mathfrak{p}).
\end{align}
Therefore, using \eqref{psd2kjde}, we see that $W(\tau):=P(\tau)-Q(\tau)$ solves
\[
\frac{d}{d\tau}W(\tau) = [W(\tau),\mathcal{A}]_x +  \frac{\tau^{\mathtt{N}_\alpha}}{\mathtt{N}_\alpha!}Op^W(\mathfrak{q}_{\mathtt{N}_\alpha}\star\mathfrak{a}),\quad W(0)= 0.
\]
Now, let us denote $V(\tau):=\Psi(\tau)\partial_xW(\tau)\Psi(-\tau)$. Then, using \eqref{psei_type12} and \eqref{commuetesd}, we have
\begin{align*}
\frac{d}{d\tau}V(\tau) &=\frac{\tau^{\mathtt{N}_\alpha}}{\mathtt{N}_\alpha!}\Psi(\tau)\partial_xOp^W(\mathfrak{q}_{\mathtt{N}_\alpha}\star\mathfrak{a}) \Psi(-\tau),\quad V(0) = 0.
\end{align*}
Integrating in $\tau$, we obtain 
\[
\partial_x W(\tau) = \frac{1}{\mathtt{N}_\alpha!}\int_0^\tau \Psi(t -\tau)\partial_x Op^W( \mathfrak{q}_{\mathtt{N}_\alpha}\star\mathfrak{a})\Psi(\tau-t) t^{\mathtt{N}_\alpha}dt.
\]
Using that $\Psi(\tau)$ is symplectic and \eqref{symplectic73} holds, we obtain
\begin{align*}
W(\tau) & =\frac{1}{\mathtt{N}_\alpha!} \int_0^{\tau}\partial_x^{-1}\Psi(t-\tau)\partial_x Op^W(\mathfrak{q}_{\mathtt{N}_\alpha}\star\mathfrak{a})\Psi(\tau-t)t^{\mathtt{N}_\alpha}dt \\
& = \frac{1}{\mathtt{N}_\alpha!}\int_0^\tau \Psi(\tau-t)^{T}Op^W(\mathfrak{q}_{\mathtt{N}_\alpha}\star\mathfrak{a})\Psi(\tau-t)t^{\mathtt{N}_\alpha}dt.
\end{align*}
Therefore, we obtain the desired result \eqref{jsd2sdxxq123s}.
\end{proof}


\chapter{Reduction to a constant coefficients operator}\label{reduction}
In this chapter, we aim to reduce the linear operator $\mathcal{L}_\omega$ in Proposition~\ref{linearized_opspdosodw} to a constant coefficients operator to find its inverse, that is, we will prove Proposition~\ref{normal_inversion} in Section~\ref{loidneks1s2}.  Throughout Sections~\ref{change_of_the_space_variables}-\ref{taksd2tamesd}, we will conjugate $\mathcal{L}_\omega$ to a constant coefficients operator up to a Lip-$0$-modulo tame operator, using symplectic transformations that are reversibility preserving and $\frac{2\pi}{\mathtt{M}}$-translation invariance preserving (See Proposition~\ref{modulut2sosdtame}). In Section~\ref{rkppsdsdwdwdx1kaksd}, we will restrict the operator to $\frac{2\pi}{\mathtt{M}}$-translation invariant functions (the functions whose Fourier modes are supported only on the multiples of $\mathtt{M}$), and proceed to obtain the full reduction. 

\section{Change of the space variables}\label{change_of_the_space_variables}
The goal of this section is to make the coefficient of the highest order term $\Lambda^{\alpha-1}_{-\frac{1}{2} + a_1}$ in \eqref{linearisd2} independent of the variable $x$. The result of this section is summarized in the following proposition.
\begin{proposition}\label{toohard_2_3}
There exists a linear transformation $\Phi_1:H_{S^\perp}\mapsto H_{S^\perp}$ such that
\begin{equation}\label{linear_111111}
\begin{aligned}
\mathcal{L}^1[h]&:=(\Phi_1)^{-1}\mathcal{L}_\omega\Phi_1[h] = \D_\omega h -  \Pi_{S^\perp}\partial_xM_1[h] + R_1[h],\\
M_1[h] &:= {b}_1\Lambda^{\alpha-1}h + b_2 h + \Upsilon^{\alpha-3}_{{b}_3}h,
\end{aligned}
\end{equation}
satisfies the following:
\begin{enumerate}[label=(\arabic*)]
\item \label{sd2item1}${b}_1$ is independent of $x$ and there exist constants $\mathtt{m}_{\alpha,1}(\omega),\mathtt{m}_{\alpha,2}(\omega)$ and  $\tilde{q}_\alpha$, such that 
 \begin{align}\label{hsdjsddsd}
 b_1(\varphi)=-\frac{1}{2} + \epsilon^2\mathtt{m}_{\alpha,1} + \mathtt{m}_{\alpha,2} + \tilde{q}_\alpha,\quad \text{ with }\int_{\mathbb{T}^\nu}\tilde{q}_\alpha(\omega,\varphi)d\varphi = 0, 
 \end{align}
 and 
 \begin{equation}\label{qtildesd2dsd}
\begin{aligned}
 |\mathtt{m}_{\alpha,1}|^{\Lip(\gamma,\Omega_1)} &\le_{\mathtt{pe}} 1,\\
 |\mathtt{m}_{\alpha,2}|^{\Lip(\gamma,\Omega_1)}&\le_{\mathtt{pe}} \epsilon^{7-4b},\quad |d_i\mathtt{m}_{\alpha,2}(i_0)[\ihat]|\le_{\mathtt{pe}} \epsilon\rVert \ihat|_{s_0+\mu_0},\\
  \rVert \tilde{q}_\alpha\rVert_{s}^{\Lip(\gamma,\Omega_1)}
&\le_{\mathtt{pe},s} \epsilon^4 + \epsilon^{2}\rVert\mathfrak{I}_\delta\rVert_{s+\mu_0}^{\Lip(\gamma,\Omega_1)},\\
 \rVert d_i\tilde{q}_\alpha(i_0)[\ihat]\rVert_{s}
& \le_{\mathtt{pe},s} \epsilon^2\left(\rVert \ihat \rVert_{s+\mu_0} + \rVert \mathfrak{I}_\delta \rVert_{s+\mu_0}\rVert \ihat \rVert_{s_0+\mu_0} \right). 
\end{aligned}
\end{equation}
\item \label{sd2item2} ${b}_2-\frac{T_\alpha}{4}\in\mathfrak{C}_1(i_0)$ and ${b}_3\in \mathfrak{C}_2(i_0)$ and
\begin{align}\label{symmetred}
b_1(-\varphi) =   b_1(\varphi),\quad b_2(-\varphi,-x)= b_2(\varphi,x),\quad b_3(-\varphi,-x,-y)=b_3(\varphi,x,y).
\end{align}
\item \label{sd2item3}$R_1\in \mathfrak{R}(i_0)$.
\item \label{sd2item4} $\Phi_1$ is real, reversibility preserving and and symplectic. Therefore, $\mathcal{L}^1$ is a reversible and Hamiltonian.
\item \label{sd2item5}$\mathcal{L}^1$ and  $\Phi_1$ are $\frac{2\pi}{\mathtt{M}}$-translation invariance preserving and 
\begin{align}\label{b2b3arem}
\rho_{\mathtt{M}}(b_2) = b_2, \quad \rho_{\mathtt{M}}(b_3)=b_3.
\end{align}
\end{enumerate}
\end{proposition}
\begin{proof}
The proof will be completed using the lemmas presented in the rest of the section. The expression \eqref{toohard_2_3} and the item \ref{sd2item1} will follow from Lemma~\ref{space_reduction12} (note that $\partial_x(b_1)=0$ in \eqref{newposd22}, since $b_1$ does not depend on $x$) with the estimates for $\mathtt{m}_{\alpha,1}, \mathtt{m}_{\alpha,2},\tilde{q}_\alpha$ given in Lemma~\ref{beta_choice}. The items \ref{sd2item2}-\ref{sd2item4} will follow from Lemma~\ref{beta_choice}, Lemma~\ref{revs1},  Lemma~\ref{realshd2sds}, Lemma~\ref{revs2} and Lemma~\ref{sd2sdsdsddsd}.  Lastly, if $i_0$ is $\frac{2\pi}{\mathtt{M}}$-translation invariant, then \ref{linsi4} of Proposition~\ref{linearized_opspdosodw} tells us that our definition of $\beta$ in \eqref{const_coeff_1232} is also $\frac{2\pi}{\mathtt{M}}$-translation invariant. Hence, it follows from Lemma~\ref{minvariantce} that the flow map $\Phi_1$ defined by \eqref{tsd2222sd21} is $\frac{2\pi}{\mathtt{M}}$-translation invariance preserving, since the PDE in \eqref{tsd2222sd21} is a Hamiltonian PDE associated to $\frac{2\pi}{\mathtt{M}}$-translation invariant Hamiltonian: $u\mapsto \int_{\mathbb{T}}u(x)\Pi_{S^\perp}(b(\tau,\varphi,x)u(x))dx$.   Using \eqref{alsd2aem}, and \eqref{remainder_1111}, one can easily see that \eqref{b2b3arem} holds true.
\end{proof}
In order to find a symplectic transformation in Proposition~\ref{toohard_2_3}, we will use the transformation of type (1), described in Definition~\ref{type1flow}. Let $\Phi_1$ be the time-$1$ flow map of the PDE:
\begin{align}\label{tsd2222sd21}
\partial_\tau u &= \Pi_{S^\perp}\partial_x(b(\tau,\varphi,x)u), \\
 b(\tau,\varphi,x)&:=\frac{\beta(\varphi,x)}{1+\tau\beta_x(\varphi,x)}, \nonumber
\end{align}
for some $\beta\in \mathfrak{C}_1(i_0)$ to be chosen later (see Lemma~\ref{beta_choice}), and let $\Psi_1$ be the time-$1$ flow map of the PDE:
\begin{align}\label{tsd2222sd2}
\partial_\tau u = \partial_x(b(\tau,\varphi,x)u).
\end{align}
Our analysis in this section is based on the properties of $\Phi_\tau,\Psi_\tau$, studied in Section~\ref{sub1sd}.

First of all, it follows from \eqref{psdisd2sd2sd} that 
\begin{align}\label{A_def}
\Psi_1h(\varphi,x):=(1+\beta_x(\varphi,x))h(\varphi,x+\beta(\varphi,x)) =: \partial_xB(\varphi,x)h\circ B(\varphi,x),
\end{align} where $B(\varphi,x) =x+\beta(\varphi,x)$.  Since $\Psi_1$ is a symplectic transformation in $L^2_0$ (see \cite[Remark 3.3]{Baldi-Berti-Montalto:KAM-quasilinear-airy}), it holds that $\sigma(\Psi_1[f],\Psi_1[g])=\sigma(f,g)$ for all $f,g\in L^2_0$ (see \eqref{symplectic} for the definition of $\sigma$), therefore,
\begin{align}\label{symplectic_conju_withx}
\Psi_1\partial_x\Psi_1^T = \partial_x.
\end{align}
For $\Psi_1^{-1},\Psi_1^{T}$ , we use Lemma~\ref{invsdwsd}  and \eqref{adjoints} to see that
\begin{align}\label{transpose_inverse}
 \Psi_1^{-1}[h] = (1+\tilde{\beta}_x(\varphi,x))h(x+\tilde{\beta}(\varphi,x))=\partial_x\left( B^{-1}\right) h \circ B^{-1}, \quad \Psi^T h = h\circ B^{-1},
\end{align}
where $B^{-1}:x\mapsto x+\tilde{\beta}(\varphi,x)$ is the inverse of $B:x\mapsto x+\beta(\varphi,x)$.
To compute the conjugation of  $\mathcal{L}_\omega$ with $\Phi_1$, we denote, for a given $a=a(x,y)\in C^\infty(\mathbb{T}^2)$,
\begin{align}
B_2[a](x,y)&:=a(B(x),B(y)),\quad B_2^{-1}[a](x,y):=a(B^{-1}(x),B^{-1}(y)), \nonumber \\
\tilde{a}(x)& :=a(x,x).\label{reduc_1_notations}
\end{align}

Now we compute the conjugation of $\mathcal{L}_\omega$ with $\Phi_1$.
\subsection{Computing the conjugation}
\begin{lemma}\label{space_reduction12}
We have that
\begin{equation}\label{newposd22}
\begin{aligned}
&(\Phi_1)^{-1}\mathcal{L}_\omega\Phi_1[h] = \D_\omega h -  \Pi_{S^\perp}\partial_xM_1[h] + R_1[h],\\
&M_1[h] := {b}_1\Lambda^{\alpha-1}h + {b}_2 h + \frac{1}{2}\partial_x({b}_1)\partial_x\Upsilon^{\alpha-3}h + \Upsilon^{\alpha-3}_{b_3}h,
\end{aligned}
\end{equation}
where the new coefficients $b_1,b_2,b_3$ are given by
\begin{equation}\label{remainder_1111}
\begin{aligned}
{b}_1 & = \Psi_1^T\left[(\partial_xB)^{\alpha}\left(-\frac{1}{2} + \tilde{a_1}\right)\right],\\
{b}_2 & = \Psi_1^{T}\left[\Lambda^{\alpha-1}_{-\frac{1}{2}+a_1}\partial_xB +\left(\frac{T_\alpha}{4}+a_2 \right)\partial_xB - \D_\omega \beta \right] \\
& + \int_{\mathbb{T}}(2-2\cos(x-y))^{-\frac{\alpha}{2}}(a_{1,B}(x,y)-a_{1,B}(x,x))dy,\\
{b}_3 & =  \left(\frac{2-2\cos(B^{-1}(x)-B^{-1}(y))}{2-2\cos(x-y)} \right)^{1-\frac{\alpha}{2}} B_2[a_3](x,y) \\
& - \frac{a_{1,B}(x,y)-a_{1,B}(x,x)-\partial_ya_{1,B}(x,x)\sin(x-y)}{2-2\cos(x-y)},\\
a_{1,B}(x,y) & :=\left(\frac{2-2\cos(B^{-1}(x)-B^{-1}(y))}{2-2\cos(x-y)} \right)^{-\frac{\alpha}{2}}B_2^{-1}[a_1](x,y),
\end{aligned}
\end{equation}
with a new remainder $R_1$ given by
\begin{align}\label{newremaind223}
R_1[h] = \Phi_1^{-1}\left(\mathcal{R}_\omega - \mathcal{R}_M\right)h + \Phi_1^{-1}\mathcal{R}\Phi_1 h,
\end{align}
where $\mathcal{R}_\omega,\mathcal{R}_M$ are defined in \eqref{conju_with_dw1} and \eqref{rmdefs1}.
\end{lemma}
\begin{proof}
From Proposition~\ref{linearized_opspdosodw}, we have that for all $h(\varphi,x)\in C^\infty_{\varphi,x}$, such that $h(\varphi,\cdot)\in H_{S^\perp}$, 
\begin{align}\label{22jd2232}
(\Phi_1)^{-1}\mathcal{L}_\omega\Phi_1[h] = (\Phi_1)^{-1}\D_\omega\Phi_1[h] - (\Phi_1)^{-1}\Pi_{S^\perp}\partial_xM\Phi_1[h] + (\Phi_1)^{-1}\mathcal{R}\Phi_1[h]
\end{align}

\textbf{Conjugation of $\D_\omega$.}
Noticing that $\Phi_1$ depends on $\varphi$, we have that $
\D_\omega \Phi_1[h] = \Phi_{1}\D_\omega h  + \D_\omega(\Phi_1)[h],
$ hence,
\begin{align}
\D_\omega\Phi_1&=\Phi_1 \D_\omega +\D_\omega (\Phi_1)\nonumber\\
& = \Phi_1 \D_\omega + \Pi_{S^\perp}\D_\omega (\Psi_1) + \D_\omega (\Phi_1 - \Pi_{S^\perp}\Psi_1). \nonumber
\end{align}
From the definition of $\Psi_1$ in \eqref{A_def}, it follows that $\D_\omega(\Psi_1)[h] =\Psi_1\partial_x\left(\Psi_1^{T}[\D_\omega\beta(\varphi,\cdot)]h \right) $. Therefore, we have 
\begin{align}
\D_\omega\Phi_1 h&= \Phi_1\D_\omega h + \Pi_{S^\perp}\Psi_1\partial_x\left(\Psi_1^{T}[\D_\omega\beta(\varphi,\cdot)]h \right) + \D_\omega (\Phi_1 - \Pi_{S^\perp}\Psi_1)\nonumber\\
& = \Phi_1 \D_\omega h + \Pi_{S^\perp}\Psi_1\Pi_{S^\perp}\partial_x\left(\Psi_1^{T}[\D_\omega\beta(\varphi,\cdot)]h \right) \nonumber \\
& + \Pi_{S^\perp}\Psi_1\Pi_{S}\partial_x\left(\Psi_1^{T}[\D_\omega\beta(\varphi,\cdot)]h \right) + \D_\omega (\Phi_1 - \Pi_{S^\perp}\Psi_1)\nonumber\\
& =  \Phi_1 \D_\omega h  + \Phi_1\Pi_{S^\perp}\partial_x\left(\Psi_1^{T}[\D_\omega\beta(\varphi,\cdot)]h \right) + \mathcal{R}_\omega[h]\nonumber\\
\mathcal{R}_\omega[h]& := \left( \Pi_{S^\perp}\Psi_1\Pi_{S^\perp} - \Phi_1\Pi_{S^\perp}\right)\partial_x\left(\Psi_1^{T}[\D_\omega\beta(\varphi,\cdot)]h \right) \nonumber \\
& + \Pi_{S^\perp}\Psi_1\Pi_{S}\partial_x\left(\Psi_1^{T}[\D_\omega\beta(\varphi,\cdot)]h \right) + \D_\omega (\Phi_1 - \Pi_{S^\perp}\Psi_1). \label{conju_with_dw1}
\end{align}
Therefore, we obtain
\begin{align}\label{conju_with_dw}
\Phi_1^{-1}\D_\omega \Phi_1 = \D_\omega h + \Pi_{S^\perp}\partial_x\left(\Psi_1^{T}[\D_\omega\beta(\varphi,\cdot)]h \right) + \Phi_1^{-1}\mathcal{R}_\omega[h].
\end{align}
\textbf{Conjugation of $\Pi_{S^\perp}\partial_xM$.}
We compute
\begin{align}\label{red_phi_2}
\partial_xM\Phi_1\Pi_{S^\perp} &=\partial_x M  \Psi_1\Pi_{S^\perp} + \partial_x M \left( \Phi_1\Pi_{S^\perp} - \Psi_1\Pi_{S^\perp}\right)\nonumber\\
& = \partial_x M  \Psi_1\Pi_{S^\perp} + \partial_x M \left( \Phi_1\Pi_{S^\perp} -\Pi_{S^\perp} \Psi_1\Pi_{S^\perp}\right) - \partial_x M \Pi_{S}\Psi_1\Pi_{S^\perp}.
\end{align}
Let us focus on $\partial_x M \Psi_{1}\Pi_{S^\perp}$ first.  From $M$ in \eqref{linearisd2}, we have that for $h\in H_{S^\perp}$,
\begin{align*}
\Psi_1^T M\Psi_1h = \Psi_1^T \Lambda^{\alpha-1}_{-\frac{1}2+a_1} \Psi_1h + \Psi_1^T \left[\left(\frac{T_\alpha}{4} + a_2\right) \Psi_1h\right] + \Psi_1^T \Upsilon^{\alpha-3}_{a_3}\Psi_1h.
\end{align*}
Applying Lemma~\ref{conj_lem1}, we obtain
\begin{align}\label{tileMdef}
\Psi_1^TM\Psi_1[h] = \mathtt{a}_1\Lambda^{\alpha-1}h + \mathtt{a}_2h + \frac{1}{2}\partial_x(\mathtt{a}_1)\partial_x\Upsilon^{\alpha-3}h + \Upsilon^{\alpha-3}_{\mathtt{a}_3}h=:\tilde{M}[h],
\end{align}
where
\begin{align}
\mathtt{a}_1(\varphi,x)&= \Psi_1^T\left[(\partial_xB)^{\alpha}\left(-\frac{1}{2} + \tilde{a_1}\right)\right],\nonumber\\
\mathtt{a}_2(\varphi,x) &= \Psi_1^{T}\left[\Lambda^{\alpha-1}_{-\frac{1}{2}+a_1}\partial_xB +\left(\frac{T_\alpha}{4}+a_2 \right)\partial_xB\right] \nonumber \\
& + \int_{\mathbb{T}}(2-2\cos(x-y))^{-\frac{\alpha}{2}}(a_{1,B}(x,y)-a_{1,B}(x,x))dy,\label{b22sd2}\\
\mathtt{a}_3(\varphi,x,y) & = \left(\frac{2-2\cos(B^{-1}(x)-B^{-1}(y))}{2-2\cos(x-y)} \right)^{1-\frac{\alpha}{2}}B_2[a_3](x,y) \nonumber \\
& - \frac{a_{1,B}(x,y)-a_{1,B}(x,x)-(\partial_ya_{1,B})(x,x)\sin(x-y)}{2-2\cos(x-y)},\nonumber
\end{align}
where
\[a_{1,B}(x,y):=\left(\frac{2-2\cos(B^{-1}(x)-B^{-1}(y))}{2-2\cos(x-y)} \right)^{-\frac{\alpha}{2}}B_2^{-1}[a_1](x,y).\]
Recalling \eqref{symplectic_conju_withx}, it follows that $\partial_xM\Psi_1=\Psi_1\partial_x\Psi_1^TM\Psi_1 = \Psi_1 \partial_x \tilde{M}$, where $\tilde{M}$ is as in \eqref{tileMdef}. Therefore,  
\begin{align}
\Pi_{S^\perp}\partial_xM\Psi_1 & =\Pi_{S^\perp}\Psi_1\partial_x\tilde{M}\nonumber\\
& = \Pi_{S^\perp}\Psi_{1}\Pi_{S^\perp}\partial_x\tilde{M} + \Pi_{S^\perp}\Psi_1\Pi_S \partial_x \tilde{M}\nonumber\\
& = \Phi_1\Pi_{S^\perp} \partial_x \tilde{M} \nonumber\\
& \ +(\Pi_{S^\perp}\Psi_1 \Pi_{S^\perp} - \Phi_1\Pi_{S^\perp})\partial_x \tilde{M} + \Psi_1\Pi_S \partial_x \tilde{M},\label{2j2j2jsd2}
\end{align}  Hence, it follows from \eqref{red_phi_2} that
\begin{align}\label{rmdefs1}
\Pi_{S^\perp}\partial_x M \Phi_1\Pi_{S^\perp}& = \Pi_{S^\perp}\partial_x M  \Psi_1\Pi_{S^\perp} +\Pi_{S^\perp} \partial_x M \left( \Phi_1\Pi_{S^\perp} -\Pi_{S^\perp} \Psi_1\Pi_{S^\perp}\right) \nonumber \\ & - \Pi_{S^\perp}\partial_x M \Pi_{S}\Psi_1\Pi_{S^\perp}\nonumber\\
&= \Phi_1\Pi_{S^\perp} \partial_x \tilde{M} + (\Pi_{S^\perp}\Psi_1 \Pi_{S^\perp} - \Phi_1\Pi_{S^\perp})\partial_x \tilde{M} + \Psi_1\Pi_S \partial_x \tilde{M} \nonumber\\
& \  + \Pi_{S^\perp} \partial_x M \left( \Phi_1\Pi_{S^\perp} -\Pi_{S^\perp} \Psi_1\Pi_{S^\perp}\right) - \Pi_{S^\perp}\partial_x M \Pi_{S}\Psi_1\Pi_{S^\perp}\nonumber\\
& =:  \Phi_1\Pi_{S^\perp} \partial_x \tilde{M} + \mathcal{R}_{M}.
\end{align}
Then, we  obtain,
\begin{align}\label{conju_withM}
\Phi_1^{-1}\Pi_{S^\perp}\partial_x M \Phi_1 = \Pi_{S^\perp}\partial_x\tilde{M} + \Phi_1^{-1}\mathcal{R}_M.
\end{align}
Now, we plug the computations obtained in \eqref{conju_with_dw} and \eqref{conju_withM} into \eqref{22jd2232}, and obtain that for $h\in H_{S^\perp}$,
\begin{align*}
\Phi_1^{-1}\mathcal{L}_\omega \Phi_1 h  &= \D_\omega h - \Pi_{S^\perp}\partial_x \underbrace{\left( - \Psi_1^{T}[\D_\omega\beta(\varphi,\cdot)]h + \tilde{M}h \right)}_{=:M_1 [h]} \\ &  + \underbrace{\Phi_1^{-1}\left(\mathcal{R}_\omega - \mathcal{R}_M\right)h + \Phi_1^{-1}\mathcal{R}\Phi_1 h}_{=:R_1[h]} \\
& = \D_\omega h -\Pi_{S^\perp}\partial_x M_1[h] + R_1[h].
\end{align*}
Recalling the coefficients of $\tilde{M}$ from \eqref{tileMdef}, combining  the coefficients $\mathtt{a}_1,\mathtt{a}_2,\mathtt{a_3}$  with $-\Psi_1^T[\D_\omega\beta]$ in $M_1$, we obtain \eqref{newposd22}.
\end{proof}

\subsection{Analysis of the new coefficients $b_1,b_2,b_3$}

We need to choose $\beta$ to make $b_1$ in \eqref{newposd22} independent of the variable $x$ (recalling $B$ from \eqref{A_def}). 
\begin{lemma}\label{beta_choice}
Let $b_1$ be as in Lemma~\ref{space_reduction12}, that is, 
\begin{align}\label{b1ded2232}
 b_1 = \mathcal{A}^{T}\left[(\partial_xB)^\alpha\left(-\frac{1}{2} + \tilde{a}_1\right)\right],
 \end{align} 
where $B,\tilde{a}_1$ are as in \eqref{A_def} and \eqref{reduc_1_notations}. We have the following:
\begin{enumerate}[label=(\arabic*)]
\item \label{seccsdw2}We can choose $\beta$ so that
$b_1$ and $\beta$ satisfy
\begin{align}
\beta(\varphi,x) & :=(\partial_x)^{-1}\left(\left(\frac{1}{2}-\tilde{a}_1\right)^{-\frac{1}{\alpha}}(-b_1(\varphi))^{\frac{1}{\alpha}}-1 \right), \nonumber \\b_1(\varphi) & :=-\left( \frac{1}{2\pi}\int_\mathbb{T}\left(\frac{1}{2}-\tilde{a}_1(\varphi,x)\right)^{-\frac{1}{\alpha}}dx\right)^{-\alpha}.\label{const_coeff_1232}
\end{align}
 \item \label{sd2220ks} There exist $\mathtt{m}_{\alpha,1}(\omega),\mathtt{m}_{\alpha,2}(\omega)$, and   $\tilde{q}_\alpha$, such that 
 \begin{align}\label{decomsd2s}
 b_1(\varphi)=-\frac{1}{2} + \epsilon^2\mathtt{m}_{\alpha,1} + \mathtt{m}_{\alpha,2} + \tilde{q}_\alpha(\varphi),\quad \text{ with }\int_{\mathbb{T}^\nu}\tilde{q}_\alpha(\omega,\varphi)d\varphi = 0,
 \end{align}
 with the estimates,
\begin{align*}
 |\mathtt{m}_{\alpha,1}|^{\Lip(\gamma,\Omega_1)} &\le_{\mathtt{pe}} 1,\\
 |\mathtt{m}_{\alpha,2}|^{\Lip(\gamma,\Omega_1)}&\le_{\mathtt{pe}} \epsilon^{7-4b},\quad |d_i\mathtt{m}_{\alpha,2}(i_0)[\ihat]|\le_{\mathtt{pe}} \epsilon\rVert \ihat|_{s_0+\mu_0},\\
  \rVert \tilde{q}_\alpha\rVert_{s}^{\Lip(\gamma,\Omega_1)}
&\le_{\mathtt{pe},s} \epsilon^4 + \epsilon^{2}\rVert\mathfrak{I}_\delta\rVert_{s+\mu_0}^{\Lip(\gamma,\Omega_1)},
\\ \rVert d_i\tilde{q}_\alpha(i_0)[\ihat]\rVert_{s}
& \le_{\mathtt{pe},s} \epsilon^2\left(\rVert \ihat \rVert_{s+\mu_0} + \rVert \mathfrak{I}_\delta \rVert_{s+\mu_0}\rVert \ihat \rVert_{s_0+\mu_0}\right). 
\end{align*}
\item \label{betasd22} We have $\beta\in \mathfrak{C}_1(i_0)$ Furthermore, we have $b_2\in \mathfrak{C}_1(i_0)$ and $b_3\in \mathfrak{C}_2(i_0)$ where $b_2 - \frac{T_\alpha}4,b_3$ are as in Lemma~\ref{space_reduction12}.
\end{enumerate}
\end{lemma}
\begin{proof}
\vspace{0.5\baselineskip}
\noindent\textit{Proof of (1).}
 Note that \eqref{const_coeff_1232} can be proved in the same way as in \cite[Section 3]{Baldi-Berti-Montalto:KAM-quasilinear-airy}. Indeed,  taking $\left(\Psi_1^T\right)^{-1}$ in \eqref{b1ded2232} and recalling that $B'=1+\beta_x$ (see definition of $B$ in \eqref{A_def}), we are led to find $b_1(\varphi)$ and $\beta(\varphi,x)$ such that 
\begin{align*}
(1+\beta_x(\varphi,x))^\alpha\left(-\frac{1}{2} + \tilde{a_1}(\varphi,x)\right) = b_{1}(\varphi),
\end{align*}
which is equivalent to 
\[
(1+\beta_x)^{\alpha} = \left(  \frac{1}{2} - \tilde{a}_1\right)^{-1}\left(-b_{1}\right) \iff \beta_x= \left(\frac{1}{2} -\tilde{a}_1 \right)^{-\frac{1}{\alpha}} \left(-b_1\right)^{\frac{1}{\alpha}} -1.
\]
To find $\beta$, we have to take $\partial_x^{-1}$, which forces us to choose $b_1$ so that the very right-hand side has zero-average. This gives us a unique choice of $b_1$ and $\beta$ as given in \eqref{const_coeff_1232}.

  \vspace{0.5\baselineskip}
\noindent\textit{Proof of (2)}
  To see \ref{sd2220ks}, recall that $a_1\in \mathfrak{C}_2(i_0)$ thanks to \ref{linsi1} of Proposition~\ref{linearized_opspdosodw}. Recalling the definition of $\tilde{a}_1$ in  \eqref{reduc_1_notations}, it easily follows from the definition of $\mathfrak{C}_1(i_0),\mathfrak{C}_2(i_0)$ in Definition~\ref{homogeneous_expansion1} that  $ \tilde{a}_1\in \mathfrak{C}_1(i_0)$. Now, using the Taylor expansion of $x\mapsto \left(\frac{1}{2} - x\right)^{-\frac{1}{\alpha}}$ near $x=0$, it is easy to see from $b_1$ in \eqref{const_coeff_1232} that 
  \[
  (-b_1(\varphi))^{-\frac{1}{\alpha}} = \frac{1}{2\pi}\int_{\mathbb{T}} \left( \frac{1}{2} - \tilde{a}_1(\varphi,x)\right)^{-\frac{1}{\alpha}}dx = \left(\frac{1}{2}\right)^{-\frac{1}{\alpha}} + M_x(f),
  \] for some $f\in \mathfrak{C}_1(i_0)$.
  Again, using the Taylor expansion of $x\mapsto  \left(\left(\frac{1}2\right)^{-\frac{1}{\alpha}} + x\right)^{-\alpha}$ and Lemma~\ref{meanisbetter}, we obtain \eqref{decomsd2s} with the desired estimates described in the lemma.
  
\vspace{0.5\baselineskip}
\noindent\textit{Proof of \ref{betasd22}} This follows straightforwardly from the choice of $\beta$ in \eqref{const_coeff_1232}, item \ref{sd2220ks} and the fact that $\tilde{a}_1\in \mathfrak{C}_1(i_0)$.
\end{proof}

\subsection{Analysis of the transformation $\Phi_1$}

\begin{lemma}\label{align2sd2sd24s2}
$\Phi^{\pm1} = \Phi_1,\Phi_1^{-1}$ satisfies 
\begin{equation}\label{transformation_estimate_class}
\begin{aligned}
&\rVert (\Phi^{\pm}-I) h \rVert_{s}^{\Lip(\gamma,\Omega_1)} \le_{\mathtt{pe},s} \epsilon\left( \rVert h \rVert^{\Lip(\gamma,\Omega_1)}_{s+\mu_0} + \rVert \mathfrak{I}_\delta \rVert^{\Lip(\gamma,\Omega_1)}_{s+\mu_0}\rVert h \rVert^{\Lip(\gamma,\Omega_1)}_{s_0+\mu_0}\right)\\
&\rVert d_i\Phi^{\pm}(i_0)h[\ihat] \rVert_{s} \le_{\mathtt{pe},s} \left(\rVert h\rVert_{s+\mu_0} + \rVert \mathfrak{I}_\delta \rVert_{s+\mu_0}\rVert h \rVert_{s_0+\mu_0}\right) \rVert \ihat \rVert_{s_0+\mu_0} + \rVert h \rVert_{s_0+\mu_0}\rVert \ihat \rVert_{s+\mu_0}.
\end{aligned}
\end{equation}
\end{lemma}
\begin{proof}
Using that $\Phi_1=\Psi_1+ R$ on $H_{S^\perp}$ for a finite dimensional operator $R\in \mathfrak{R}(i_0)$, which follows from Lemma~\ref{fhsjdwd2sd}, we see that the  above estimates for $\Phi_1$ follow from the definition of $\mathfrak{R}(i_0)$ in Definition~\ref{remiander_class12} and Lemma~\ref{expasd111}. For $\Phi_1^{-1}$, we recall the following lemma:
\begin{lemma}\cite[Lemma 8.5]{Baldi-Berti-Montalto:KAM-quasilinear-kdv}\label{2jsd2e}
It holds that
\begin{align*}
\rVert \Phi_1^{\pm }h\rVert_s^{\Lip(\gamma,\Omega_1)}& \le_{\mathtt{pe},s}  \rVert h\rVert^{\Lip(\gamma,\Omega_1)}_{s+\mu_0} + \rVert \mathfrak{I}_\delta \rVert^{\Lip(\gamma,\Omega_1)}_{s+\mu_0}\rVert h \rVert_{s_0+\mu_0}^{\Lip(\gamma,\Omega_1)},\\
\rVert d_i\Phi_1^{\pm}(i_0)h[\ihat] \rVert_{s}  & \le_{\mathtt{pe},s} \left(\rVert h\rVert_{s+\mu_0} + \rVert \mathfrak{I}_\delta \rVert_{s+\mu_0}\rVert h \rVert_{s_0+\mu_0}\right) \rVert \ihat \rVert_{s_0+\mu_0} + \rVert h \rVert_{s_0+\mu_0}\rVert \ihat \rVert_{s+\mu_0}.
\end{align*}
\end{lemma}
With the above lemma, we differentiate (denoting $\Phi_\tau,\Psi_\tau$ be the time-$\tau$ flows of the PDEs  \eqref{tsd2222sd21} and \eqref{tsd2222sd2} respectively), $I=\Phi_\tau\circ \left(\Phi_{\tau}\right)^{-1}$ in $\tau$ to see that
\[
\frac{d}{d\tau}\Phi_\tau^{-1} =  -\Phi_\tau^{-1} \partial_x b(\tau)\Phi_\tau^{-1}, \quad \Phi_0^{-1} = I.
\]
Using the Taylor expansion near $\tau = 0$, we have
\[
\Phi_\tau^{-1} = I -\int_{0}^{\tau}\Phi_s^{-1}\partial_x b(\tau)\Phi_s^{-1}ds,
\]
therefore,
\begin{align}\label{2jsjsj2}
\rVert\left( \Phi_1^{-1} - I \right)h\rVert^{\Lip(\gamma,\Omega_1)}_{s} \le \sup_{\tau\in [0,1]}\rVert\Phi_{\tau}^{-1}\partial_x b(\tau)\Phi_\tau^{-1}h\rVert^{\Lip(\gamma,\Omega_1)_{s}}.
\end{align}
From $b(\tau)$ in \eqref{tsd2222sd21} with $\beta\in \mathfrak{C}_1(i_0)$ (see \ref{betasd22} in Lemma~\ref{beta_choice}),  it follows that (see \ref{slsd22sd2sdhw} of Lemma~\ref{properties_weights})
\[
\rVert b(\tau)\rVert^{\Lip(\gamma,\Omega_1)}_s\le_{\mathtt{pe},s} \epsilon\left(1 + \rVert \mathfrak{I}_\delta\rVert^{\Lip(\gamma,\Omega_1)}_{s+\mu_0} \right).
\] Hence, Lemma~\ref{2jsd2e} and \eqref{2jsjsj2} implies the first estimate in \eqref{transformation_estimate_class} for $\Phi_1^{-1}$ (note that Lemma~\ref{2jsd2e} is written for the time-$1$ flow, however, one can always reparametrize the time variable $\tau$ to apply the same estimates to $\Phi_\tau^{\pm}$ in \eqref{2jsjsj2} for $\tau\in[0,1]$). 
 The second estimate in \eqref{transformation_estimate_class} for $d_i(\Phi^{\pm})$ follows immediately from Lemma~\ref{2jsd2e}.
\end{proof}

\begin{lemma}\label{revs1}
$\Phi_1$ is  reversibility preserving and $b_1,b_2,b_3,\beta$ satisfy,
\begin{align*}
b_1(-\varphi) & =b_1(\varphi),\quad b_2(-\varphi,-x)=b_2(\varphi,x),\\
b_3(-\varphi,-x,-y) & =b_3(\varphi,x,y),\quad \beta(-\varphi,-x)=-\beta(\varphi,x).
\end{align*}
\end{lemma}
\begin{proof}
It is clear from $\beta,b_1$ in \eqref{const_coeff_1232} and the property of $a_1$ in \eqref{reverssss} that
\[
b_1(-\varphi) = b_1(\varphi),\quad \beta(-\varphi,-x) = \beta(\varphi,x).
\]
From Lemma~\ref{reversibility}, it follows that $\Phi_1$ is a reversibility preserving map. The results for $b_2,b_3$ follow straightforwardly from their expressions in Lemma~\ref{space_reduction12} with the definitions of $B,B_2$ and $\Psi_1^T$ in \eqref{transpose_inverse} and \eqref{reduc_1_notations} with the property of $a_i$ in \eqref{reverssss}.
\end{proof}

\begin{lemma}\label{realshd2sds}
$\Phi_1$ is a real operator and $b_1,b_2,b_3$ are real-valued.
\end{lemma}
\begin{proof}
Using that  $a_1,a_2,a_3$ in Proposition~\ref{linearized_opspdosodw} are real-valued, it follows from \eqref{const_coeff_1232} that $\beta$ is real. This clearly implies $\Phi_1$ is a real operator. Since $\mathcal{L}_\omega$ is real (see \eqref{linsi3} of Proposition~\ref{linearized_opspdosodw}), $\mathcal{L}^1$ is also real operator and the coefficients $b_2,b_3$ are real-valued as well.
\end{proof}

\begin{lemma}\label{revs2}
$\Phi_1$ is a symplectic transformation, therefore, $\mathcal{L}^1$ is  a Hamiltonian operator.
\end{lemma}
\begin{proof}
This follows from the fact that $\Phi_1$, defined in \eqref{tsd2222sd21}, is a time-$1$ flow of the Hamiltonian flow in $H_{S^\perp}$ and  Lemma~\ref{Hamiltonian_operators_under_symplectic}. 
\end{proof}

\subsection{Analysis of the remainder $R_1$}

\begin{lemma}\label{sd2sdsdsddsd}
$R_1\in \mathfrak{R}(i_0)$.
\end{lemma}
\begin{proof}
From \eqref{newremaind223} (also see \eqref{conju_with_dw1}, and  \eqref{rmdefs1}), it follows that each operator in $R_1$ involves a finite dimensional operator (see Lemma~\ref{fhsjdwd2sd}). Among all the terms, we only show that $\Phi_1^{-1}\mathcal{R}\Phi_1\in \mathfrak{R}(i_0)$ only, since the other terms can be done in the same way.  But this follows from $\mathcal{R}\in \mathfrak{R}(i_0)$, Lemma~\ref{align2sd2sd2} and Lemma~\ref{funitesd2sd}, that $\Phi_1^{-1}\mathcal{R}\Phi_1\in \mathfrak{R}(i_0)$. 
 \end{proof}


\section{Reparametrization of time}\label{reparmedumbple}
In this section, we consider the transformation of type (3):
\begin{align}\label{time_rep_transforms}
\Phi_2[h](\varphi,x):=h(\varphi + \omega p_1(\varphi),x), \quad \Phi_2^{-1}[h](\varphi,x) = h(\varphi + \tilde{p}_1(\omega,\varphi),x)\end{align}
for some $p_1:\mathbb{T}^\nu \mapsto \mathbb{R}$ and $\varphi\mapsto \varphi+\tilde{p}_1(\omega,\varphi)$, the inverse of $\varphi\mapsto \omega p_1(\varphi)$. We will find $p_1$ so that  the coefficient of $\Lambda^{\alpha-1}$ of the operator $\left(\Phi^2\right)^{-1}\mathcal{L}^1_\omega \Phi^2$ is constant. More precisely, we summarize the main result in the following proposition:
\begin{proposition}\label{prop_time_rep_92}
There exist a transformation $\Phi_2:H_{S^\perp}\mapsto H_{S^\perp}$ and $\rho$ such that
\begin{equation}\label{sdlinear2sd2}
\begin{aligned}
\mathcal{L}^2[h]&:=\frac{1}{\rho}(\Phi_2)^{-1}\mathcal{L}^1\Phi_2[h] = \D_\omega h -  \Pi_{S^\perp}\partial_xM_2[h] + R_2[h],\\
M_2[h] &:= \mathfrak{m}_\alpha(\omega)\Lambda^{\alpha-1}h + d_2 h + \Upsilon^{\alpha-3}_{d_3}h,
\end{aligned}
\end{equation}
satisfy the following:
\begin{enumerate}[label=(\arabic*)]
\item \label{itnsdwd} $\mathtt{m}_\alpha=-\frac{1}{2} + \epsilon^2 \mathtt{m}_{\alpha,1} + {\mathtt{m}}_{\alpha,2}$, where $\mathtt{m}_{\alpha,1},{\mathtt{m}}_{\alpha,2}$ are as in  \ref{sd2item1} in Proposition~\ref{toohard_2_3}
\item \label{sdsdwdsd} $d_2-\frac{T_\alpha}4 \in \mathfrak{C}_{1,p_1}(i_0)$ and $d_3\in \mathfrak{C}_{2,p_1}(i_0)$ (see Definition~\ref{class_TT} for the definitions of $\mathfrak{C}_{1,p_1},\mathfrak{C}_{2,p_1}$), and
\[
d_2(-\varphi,-x) = d_2(\varphi,x),\quad d_3(-\varphi,-x,-y)=d_3(\varphi,x,y).
\]
\item \label{2k232sd2} $R_2\in \mathfrak{R}(i_0)$.
\item \label{2k232sd2d}$\Phi_2$ is real and reversibility preserving, hence $\mathcal{L}^2$ is real and  reversible. Furthermore, it is Hamiltonian.
\item \label{2k232sd2d612} $\mathcal{L}^2$ and $\Phi_2$ are $\frac{2\pi}{\mathtt{M}}$-translation invariance preserving and 
\begin{align}\label{b2b322arem}
\rho_{\mathtt{M}}(d_2) = d_2, \quad \rho_{\mathtt{M}}(d_3)=d_3.
\end{align}
\end{enumerate}
\end{proposition}
\begin{proof}
The proposition is proved using the results obtained in the rest of this section. 
The expression of $\mathcal{L}^2$ is obtained in Lemma~\ref{phi_property_time} with $d_1=\mathtt{m}_\alpha$ in Lemma~\ref{2ksd2}. \ref{itnsdwd} is given in \eqref{sd2d2}. \ref{sdsdwdsd} is proved in \eqref{lsd2de2d3} and Lemma~\ref{2sk2sd2sd34sd}. \ref{2k232sd2} follows from Lemma~\ref{remisn2sd22}. \ref{2k232sd2d} follows from Lemma~\ref{2sk2sd2sd34sd}, Lemma~\ref{realsdsldsld2sd} and Lemma~\ref{hamiltosd2}. Lastly, \ref{2k232sd2d612} follows trivially since the transformation $\Phi_2$ in \eqref{time_rep_transforms} does not change the variable $x$ and therefore, \eqref{b2b322arem} follows from \eqref{remiander_time_rep} and \eqref{b2b3arem}.
\end{proof}

 Recall from \cite[Subsection 8.2]{Baldi-Berti-Montalto:KAM-quasilinear-kdv} that
\begin{align}\label{phi_property_time}
\Phi_2^{-1}\D_\omega \Phi_2 = \rho(\varphi)\D_\omega, \quad \rho(\varphi):=\Phi_2^{-1}(1+\D_\omega p_1(\varphi)).
\end{align}
As Proposition~\ref{toohard_2_3}, we start with computing the conjugation of the linear operator.
\subsection{Computing the conjugation}
\begin{lemma}\label{lemma_time_rep_1}
We have that
\begin{align}\label{2sd2sdsdkk0}
\mathcal{L}^2:=&\frac{1}{\rho(\varphi)}(\Phi_2)^{-1}\mathcal{L}^1\Phi_2[h] = \D_\omega - \Pi_{S^\perp}\partial_xM_2[h] + R_2[h],\\
&M_{2}[h] = d_1\Lambda^{\alpha-1} + d_2h + \Upsilon^{\alpha-3}_{d_3}h,\nonumber
\end{align}
where 
\begin{align}
d_1(\varphi) & = \frac{1}{\rho(\varphi)}b_1(\varphi+\tilde{p}_1(\omega,\varphi)), \quad d_2(\varphi,x) = \frac{1}{\rho(\varphi)}{b}_2(\varphi+\tilde{p}_1(\omega,\varphi),x), \nonumber \\
d_3(\varphi,x,y) & = \frac{1}{\rho(\varphi)}b_3(\varphi+\tilde{p}_1(\omega,\varphi),x,y),\label{coeff_time_rep}\\
R_2[h] & =\frac{1}{\rho(\varphi)}\Phi_2^{-1}R_1\Phi_2[h].\label{remiander_time_rep}
\end{align}
\end{lemma}

\begin{proof}
From Proposition~\ref{toohard_2_3}, we have that for all $h(\varphi,x)\in C^\infty_{\varphi,x}$, such that $h(\varphi,\cdot)\in H_{S^\perp}$, 
\begin{align}\label{22jd22322s}
(\Phi_2)^{-1}\mathcal{L}^1\Phi_2[h] = (\Phi_2)^{-1}\D_\omega\Phi_2[h] - (\Phi_2)^{-1}\Pi_{S^\perp}\partial_xM_1\Phi_2[h] + (\Phi_2)^{-1}R_1\Phi_2[h]
\end{align}

Let us  compute the conjugation of $\Pi_{S^\perp}\partial_x M_1$. Since $\Phi_2$ in \eqref{time_rep_transforms} commutes with $\Pi_{S^\perp}$ and $\partial_x$, it follows straightforwardly that (recalling $M_1$ from \eqref{linear_111111})
\begin{align}
\rho(\varphi)\Pi_{S^\perp}\partial_x M_2 & :=\Phi_2^{-1}\Pi_{S^\perp}\partial_x M_1\Phi_2 = \Pi_{S^\perp}\partial_x \Phi_2^{-1}M_1\Phi_2 \nonumber\\
& = \Pi_{S^\perp}\partial_x \left( \Phi_2^{-1}[b_1]\Lambda^{\alpha-1} + \Phi_2^{-1}[b_2] +\Upsilon^{\alpha-3}_{\Phi_2^{-1}[b_3]}\right).
\label{2sd2m222sd}
\end{align}
Plugging this into \eqref{22jd22322s} and using \eqref{phi_property_time}, the result follows immediately.
\end{proof}

\subsection{Analysis of the new coefficients $d_1,d_2,d_3$}

\begin{lemma}\label{2ksd2}
Letting 
\begin{align}\label{choice_22}
p_1(\varphi):=-\D_\omega^{-1}\left(1-\frac{{b}_1(\varphi)}{\mathtt{m}_\alpha} \right), \quad \mathtt{m}_\alpha:=\frac{1}{(2\pi)^{\nu}}\int_{\mathbb{T}^\nu}{b}_1(\varphi)d\varphi,
\end{align}  we have that $d_1=\mathtt{m}_\alpha$, which does not depend on $\varphi$.
\end{lemma}
\begin{proof}
From Lemma~\ref{lemma_time_rep_1}, we have that $d_1(\varphi) = \frac{1}{\rho(\varphi)}b_1(\varphi+\tilde{p}_1(\omega,\varphi))$. Hence, we need to choose $p_1$ and $\rho$ so that 
\begin{align}\label{choice_sd2sd}
b_1(\varphi+\tilde{p}_1(\omega,\varphi)) =\rho(\varphi)\mathtt{m}_\alpha,
\end{align}
 for some $\mathtt{m}_\alpha$ that does not depend on $\varphi$. From \eqref{phi_property_time}, we see that this can be achieved by  choosing $p_1$ and $\mathtt{m}_\alpha$ so that 
\begin{align}\label{m1choice}
b_1(\varphi) = (1+\D_\omega p_1(\varphi))\mathtt{m}_\alpha,
\end{align}
which gives us \eqref{choice_sd2sd} with 
\begin{align}\label{rhsd2sd}
\rho(\varphi):= 1 + \D_\omega p_1(\varphi + \tilde{p}_1(\omega,\varphi)).
\end{align}
Integrating \eqref{m1choice} in $\varphi$ and using that $\int_{\mathbb{T}^\nu}\D_\omega p_1(\varphi)d\varphi =0$,  it is clear that  $\mathtt{m}_\alpha$ and $p_1$ in \eqref{choice_22}  satisfy \eqref{m1choice}.
\end{proof}

For $p_1,\rho$ and $\mathtt{m}_\alpha$  chosen in Lemma~\ref{choice_22}, we have the following estimates:

\begin{lemma}\label{rep_time_m_estimates126}
We have that 
\begin{align}
&\rVert \D_\omega p_1 \rVert^{\Lip(\gamma,\Omega_1)}_s \le_{\mathtt{pe},s} \epsilon^4 + \epsilon^2\rVert \mathfrak{I}_\delta\rVert^{\Lip(\gamma,\Omega_1)}_{s+\mu_0}, \label{rep_time_est3}\\
& \rVert d_i \left(\D_\omega p_1\right)(i_0)[\ihat]\rVert_{s} \le_{\mathtt{pe},s} \epsilon^2\left( \rVert \ihat \rVert_{s+\mu_0} + \rVert \mathfrak{I}_\delta\rVert_{s+\mu_0}\rVert \ihat \rVert_{s_0+\mu_0}\right),
\label{rep_time_est4}\\
&\rVert p_1 \rVert^{\Lip(\gamma,\Omega_1)}_s \le_{\mathtt{pe},s} \gamma^{-1}\left(\epsilon^4 + \epsilon^2\rVert \mathfrak{I}_\delta\rVert^{\Lip(\gamma,\Omega_1)}_{s+\mu_0}\right), \label{rep_time_est1}\\
& \rVert d_i p_1(i)[\ihat]\rVert_{s} \le_{\mathtt{pe},s} \epsilon^2\gamma^{-1}\left( \rVert \ihat \rVert_{s+\mu_0} + \rVert \mathfrak{I}_\delta\rVert_{s+\mu_0}\rVert \ihat \rVert_{s_0+\mu_0}\right),\label{rep_time_est2}
\end{align}
\end{lemma}
\begin{proof}
\textbf{Estimates \eqref{rep_time_est3} and \eqref{rep_time_est4}.}
It follows from \eqref{m1choice} that 
\begin{align}\label{sd2d2}
\D_\omega p_1 = \frac{1}{\mathtt{m}_\alpha}\left(\mathtt{m}_\alpha-{b}_1 \right),
\end{align} 
while  \ref{sd2item1} in Proposition~\ref{toohard_2_3} tells us that
$\mathtt{m}_\alpha$ and $\mathtt{m}_\alpha-{b}_1$ are
\begin{equation}\label{m2sd2sjjd}
\begin{aligned}
\mathtt{m}_\alpha &=  - \frac{1}{2} + \epsilon^2\mathtt{m}_{\alpha,1} + \mathtt{m}_{\alpha,2},\quad \mathtt{m}_\alpha-b_1 = \tilde{q}_\alpha.
\end{aligned} 
\end{equation}
Then the estimates \eqref{rep_time_est3} and \eqref{rep_time_est4} follow from the estimates in \eqref{qtildesd2dsd} and \eqref{sd2d2}.

\textbf{ Estimates \eqref{rep_time_est1} and \eqref{rep_time_est2} .}
The estimates \eqref{rep_time_est1} and \eqref{rep_time_est2} follow from \eqref{domega_inverse1} and \eqref{domega_inverse}.
 \end{proof}

With the estimates for $p_1$ in Lemma~\ref{rep_time_m_estimates126}, the estimates for $\tilde{p}_1$ determined by \eqref{time_rep_transforms} are given in the next lemma:

\begin{lemma}\label{p_tilde_estimate11}
It holds that
\begin{align}
&\rVert \tilde{p}_1 \rVert^{\Lip(\gamma,\Omega_1)}_s \le_{\mathtt{pe},s} \gamma^{-1}\left(\epsilon^4 + \epsilon^2\rVert \mathfrak{I}_\delta\rVert^{\Lip(\gamma,\Omega_1)}_{s+\mu_0}\right), \label{rep_time_est7}\\
& \rVert d_i \tilde{p}_1(i_0)[\ihat]\rVert_{s} \le_{\mathtt{pe},s} \epsilon^2\gamma^{-1}\left( \rVert \ihat \rVert_{s+\mu_0} + \rVert \mathfrak{I}_\delta\rVert_{s+\mu_0}\rVert \ihat \rVert_{s_0+\mu_0}\right).\label{rep_time_est8}
\end{align}
\end{lemma}
\begin{proof}
The invertibility of the map $\varphi\mapsto \varphi + \omega p_1(\varphi)$ follows simply from the inverse function theorem. With the estimates for $p_1$ in Lemma~\ref{rep_time_m_estimates126}, the estimates \eqref{rep_time_est7} and \eqref{rep_time_est8} follow from the same proof as in Lemma~\ref{sed2sd}, using the Taylor expansion of $
\tilde{p}_1(\omega, \varphi) = \omega p_1(\varphi +\tilde{p}_1(\omega,\varphi)).$
\end{proof}

\begin{lemma}\label{rho_estimate_time_rep}
It holds that
\begin{align*}
&\rVert \rho-1 \rVert^{\Lip(\gamma,\Omega_1)}_{s}\le_{\mathtt{pe},s} \epsilon^4 + \epsilon^2\rVert \mathfrak{I}_\delta\rVert^{\Lip(\gamma,\Omega_1)}_{s+\mu_0},\\
&\rVert d_i \rho(i_0)[\ihat]\rVert_{s} \le_{\mathtt{pe},s} \epsilon^2\left( \rVert \ihat \rVert_{s+\mu_0} + \rVert \mathfrak{I}_\delta\rVert_{s+\mu_0}\rVert \ihat \rVert_{s_0+\mu_0}\right).
\end{align*}
\end{lemma}
\begin{proof}
It follows from the definition of $\rho$ in \eqref{rhsd2sd} and the estimates for $p_1$ and $\tilde{p}_1$ given in Lemma~\ref{rep_time_m_estimates126} (especially \eqref{rep_time_est3} and \eqref{rep_time_est4}) and Lemma~\ref{p_tilde_estimate11}.
\end{proof}
 

Recalling that $\gamma=\epsilon^{2b}$ from \eqref{frequency_set2}, $d_2$ and $d_3$ in \eqref{coeff_time_rep} are not in $\mathfrak{C}_1(i_0)$ and $\mathfrak{C}_2(i_0)$. For this reason, we denote by  $\mathfrak{C}_{1,p_1}(i_0), \mathfrak{C}_{2,p_1}(i_0)$, the classes of the coefficients to include $d_2,d_3$.

\begin{definition}\label{class_TT}
We say\index{$\mathfrak{C}_{1}$} a function $a=a(\omega,\varphi,x)\in \mathfrak{C}_{1,p_1}(i_0)$, if $a(\omega,\varphi,x)=\mathtt{a}(\omega,\varphi+\tilde{p_1}(\omega,\varphi),x)$ for some $\mathtt{a}\in \mathfrak{C}_1(i_0)$ where $\tilde{p}_1$ is in Lemma~\ref{p_tilde_estimate11}. Similarly, we say $a=a(\omega,\varphi,x,y)\in \mathfrak{C}_{2,p_1}(i_0)$, if $a(\omega,\varphi,x,y) = \mathtt{a}(\omega,\varphi+\tilde{p}_1(\omega,\varphi),x,y)$ for some $\mathtt{a}\in \mathfrak{C}_2(i_0)$.\index{$\mathfrak{C}_2$}\end{definition}

\begin{lemma}\label{lsd2de2d3}
We have that $d_2-\frac{T_\alpha}4\in \mathfrak{C}_{1,p_1}(i_0)$ and $d_3\in \mathfrak{C}_{2,p_1}(i_0)$.
\end{lemma}
\begin{proof}
We prove $d_2\in \mathfrak{C}_{1,p_1}(i_0)$ only, since $d_3\in \mathfrak{C}_{2,p_1}(i_0)$ can be proved in the same way. 

 Note that since $\varphi\mapsto \varphi + \tilde{p}_1(\omega,\varphi)$ is the inverse map of $\varphi\mapsto \varphi + \omega p_1(\varphi)$, we can write $d_2$ in \eqref{coeff_time_rep} as
 \begin{align}\label{2jsjd23232}
 d_2(\omega,\varphi,x) ={a}(\omega,\varphi+\tilde{p}_1(\omega,\varphi),x),\text{ where }{a}(\omega,\varphi,x) = \frac{1}{\rho(\varphi+\omega p_1)}b_2(\omega,\varphi,x).
 \end{align}
 In view of  item \ref{sd2item2} of Proposition~\ref{toohard_2_3}, which tells us  $b_2-\frac{T_\alpha}4 \in \mathfrak{C}_1(i_0)$, we can rewrite  $a(\omega,\varphi,x)$ as
 \[
 a(\omega,\varphi,x) = \frac{1}{\rho(\varphi +\omega p_1)}\left(b_2(\omega,\varphi,x)-\frac{T_{\alpha}}{4}\right) + \left(\frac{1}{\rho(\varphi+\omega p_1(\varphi))} - 1\right)\frac{T_\alpha}{4} + \frac{T_\alpha}4.
 \]
 Plugging this into \eqref{2jsjd23232}, we have
 \begin{equation}\label{2j2jsd2232sd2}
 \begin{aligned} d_2(\omega,\varphi,x) & - \frac{T_\alpha}4   = \mathtt{a}(\omega,\varphi +\tilde{p}_1(\omega,\varphi),x), \text{ where }\\
 \mathtt{a}(\omega,\varphi,x) & = \frac{1}{\rho(\varphi + \omega p_1(\varphi))}\left(b_2(\omega,\varphi,x) -\frac{T_\alpha}4\right) + \frac{T_\alpha}4\left(\frac{1}{\rho(\varphi + \omega p_1(\varphi))} -1 \right).
 \end{aligned}
 \end{equation}
 Recalling the definition of $\mathfrak{C}_{1,p_1}(i_0)$ from Definition~\ref{class_TT}, it suffices to show that $\mathtt{a}\in \mathfrak{C}_1(i_0)$. However, this follows straightforwardly from the fact that $b_2-\frac{T_\alpha}{4}\in \mathfrak{C}_1(i_0)$ (see \ref{sd2item2} in Proposition~\ref{toohard_2_3}) and the estimates for $\rho$ and $p_1$ in Lemma~\ref{rep_time_m_estimates126} and Lemma~\ref{rho_estimate_time_rep}.
\end{proof}

\subsection{Analysis of the transformation $\Phi_2$}

\begin{lemma}\label{align2sd2sd22}
$\Phi^{\pm1} = \Phi_2,\Phi_2^{-1}$ satisfies 
\begin{equation}\label{transformation_estimate_class2}
\begin{aligned}
\rVert (\Phi^{\pm}-I) h \rVert_{s}^{\Lip(\gamma,\Omega_1)} &\le_{\mathtt{pe},s} \epsilon\left( \rVert h \rVert^{\Lip(\gamma,\Omega_1)}_{s+\mu_0} + \rVert \mathfrak{I}_\delta \rVert^{\Lip(\gamma,\Omega_1)}_{s+\mu_0}\rVert h \rVert^{\Lip(\gamma,\Omega_1)}_{s_0+\mu_0}\right)\\
\rVert d_i\Phi^{\pm}(i_0)h[\ihat] \rVert_{s} &\le_{\mathtt{pe},s} \left(\rVert h\rVert_{s+\mu_0} + \rVert \mathfrak{I}_\delta \rVert_{s+\mu_0}\rVert h \rVert_{s_0+\mu_0}\right) \rVert \ihat \rVert_{s_0+\mu_0} + \rVert h \rVert_{s_0+\mu_0}\rVert \ihat \rVert_{s+\mu_0}\end{aligned}
\end{equation}
\end{lemma}
\begin{proof}
The estimates for $\Phi^{\pm1}$ follow from their definitions in \eqref{time_rep_transforms} and the estimates for $p_1,\tilde{p}_1$ in Lemma~\ref{rep_time_m_estimates126} and Lemma~\ref{p_tilde_estimate11}.
\end{proof}

\begin{lemma}\label{2sk2sd2sd34sd}
$\Phi_2$ is reversibility preserving. Furthermore,
\begin{align}\nonumber
\tilde{p}_1(\omega,-\varphi) & = -\tilde{p}_1(\omega,\varphi),\quad \rho(-\varphi) = \rho(\varphi),\\ d_2(-\varphi,-x) & =d_2(\varphi,x),\quad d_3(-\varphi,-x,-y)=d_3(\varphi,x,y). \label{sd2sd2txkk}
\end{align}
\end{lemma}
\begin{proof}
First, we see from \eqref{symmetred} and our choice of $p_1$ in \eqref{choice_22} that $\D_\omega p_1(-\varphi) = \D_\omega p_1(\varphi)$, therefore,
\begin{align}\label{2sd2sd2sd2}
p_1(-\varphi) = -p_1(\varphi).
\end{align}
Recalling the definition of the reversibility preserving operators in Definition~\ref{reversible_operators2s}, we have (recall the spaces $X,Y$ from \eqref{space_reversibles})
\begin{align*}
\Phi_2h(-\varphi,-x) &= h(-\varphi + \omega p_1(-\varphi),-x) \\
& = h(\varphi +\omega p_1(\varphi),x) = \Phi_2h(\varphi,x),\text{ for $h\in X$,}\\
\Phi_2h(-\varphi,-x) &= h(-\varphi + \omega p_1(-\varphi),-x) \\
& = -h(\varphi +\omega p_1(\varphi),x) = -\Phi_2h(\varphi,x),\text{ for $h\in Y$}.
\end{align*}
Therefore, $\Phi_2$ is reversibility preserving. Using the definition of $\tilde{p}_1$ in \eqref{time_rep_transforms}, \eqref{2sd2sd2sd2} implies that $\tilde{p}_1(\omega,-\varphi) = -\tilde{p}_1(\omega,\varphi).$ Therefore, it follows from \eqref{rhsd2sd} that $\rho(-\varphi)=\rho(\varphi)$. With these symmetry of $\tilde{p}_1$ and $\rho$, \eqref{sd2sd2txkk} follows from the definitions of $d_2,d_3$ in  \eqref{coeff_time_rep} and \eqref{symmetred}.
\end{proof}

\begin{lemma}\label{realsdsldsld2sd}
$\Phi_2$ is real and  $\rho,p_1,\tilde{p}_1$ are  real-valued functions.
\end{lemma}
\begin{proof}
Since $b_1$ is real-valued (see Lemma~\ref{realshd2sds}), it follows from \eqref{hsdjsddsd} and \eqref{m2sd2sjjd} that $\mathtt{m}_\alpha$ is real-valued. Therefore, it follows from \eqref{m1choice} and \eqref{rhsd2sd} that $p_1$ and $\rho$ are real-valued.  Using the definition of $\tilde{p}_1$, it is clear that $\tilde{p}_1$ is real-valued as well.
\end{proof}
\begin{lemma}\label{hamiltosd2}
$\mathcal{L}^2$ is a Hamiltonian operator.
\end{lemma}
\begin{proof}
Recalling the definition of Hamiltonian operators from Definition~\ref{hamiltonian_operatior230}, we only need to show that (from $\mathcal{L}^2$ in \eqref{2sd2sdsdkk0}) $M_2$ and $\partial_x^{-1}{R}_2$ are symmetric operators in $H_{S^\perp}$ for each fixed  $\varphi$. We show this for $M_2$ only since $R_2$ can be managed in the same way. From the definition of $M_2$ in \eqref{2sd2m222sd}, it suffices to show that $
M_2(\varphi)=\frac{1}{\rho(\varphi)}\Phi_2^{-1}M_1\Phi_2$ 
is a symmetric operator. Clearly, $\Phi_2$ reparametrizes the variable  $\varphi$ only therefore, $\Phi_2,\Phi_2^{-1}$ are symmetric operators on $H_{S^\perp}$ for each fixed $\varphi$. Since $M_1$ is also symmetric (see \ref{sd2item4} in Proposition~\ref{toohard_2_3}) and $\rho$ is real-valued (Lemma~\ref{realsdsldsld2sd}), $M_2$ is symmetric.
\end{proof}

\subsection{Analysis of the remainder $R_2$}

\begin{lemma}\label{remisn2sd22}
$R_2\in \mathfrak{R}(i_0)$.
\end{lemma}
\begin{proof}
Recalling that $R_1\in \mathfrak{R}_1(i_0)$ from \eqref{sd2item3} of Proposition~\ref{toohard_2_3}, it follows from Lemma~\ref{funitesd2sd} and Lemma~\ref{align2sd2sd22} that  $\rho R_2 \in \mathfrak{R}(i_0)$  (see \eqref{remiander_time_rep} for the definition of $R_2$). Thanks to the estimates of $\rho$ in Lemma~\ref{rho_estimate_time_rep}, the result follows straightforwardly from the definition of $\mathfrak{R}(i_0)$.
\end{proof}

\section{Egorov method}\label{sjdwjjjwdsegosdeoeo}
As\index{Egorov method} a result of Proposition~\ref{time_rep_transforms}, we obtained a conjugated linear operator $\mathcal{L}^2$ whose top order term has a constant coefficient, $\mathtt{m}_\alpha$, which does not depend on $\varphi,x$. In this section, we aim to make the coefficients of the lower order terms constant, by means of Egorov method. We organize this section into three parts. In the first part,  we express the linear operator $\mathcal{L}^2$ in terms of symbols. In the second and the third parts, we will remove the dependence on the variables $x$ and $\varphi$ respectively.
\color{black}

\subsection{Symbolic expression for the linear operator}\label{symbolic_1231_deulores}
We recall that the non-constant coefficients $d_2,d_3$ of $\mathcal{L}^2$ in \eqref{sdlinear2sd2} are in $\mathfrak{C}_{1,p_1}(i_0),\mathfrak{C}_{2,p_1}(i_0)$ respectively. We first define a class of symbols that arise from the coefficients in $\mathfrak{C}_{1,p_1}(i_0),\mathfrak{C}_{2,p_1}(i_0)$.
\begin{definition}\label{sdjjsdsdsddsgsd}
We say that a symbol $\mathfrak{a}=\mathfrak{a}(\omega,\varphi,x,\xi)\in \mathfrak{S}^m_{p_1}(i_0)$, if $\mathfrak{a}(\omega,\varphi,x,\xi) = \mathfrak{b}(\omega,\varphi + \tilde{p}_1(\omega,\varphi),x,\xi)$ for some $\mathfrak{b}\in \mathfrak{S}^m(i_0)$ (see Definition~\ref{symbolssssef}\index{$\mathfrak{S}^m$} for $\mathfrak{S}^m(i_0)$), for $\tilde{p}_1$ in \eqref{p_tilde_estimate11}.
\end{definition}

We collect useful properties of the symbols in $\mathfrak{S}^m_{p_1}(i_0)$.
\begin{lemma}\label{smysd2}
Let $a=a(\omega,\varphi,x,y)\in \mathfrak{C}_{2,p_1}(i_0)$. Then there exists $\mathfrak{a}\in \mathfrak{S}^{\alpha-3}_{p_1}(i_0)$ such that $Op^W(\mathfrak{a}) = \Upsilon^{\alpha-3}_{a}$.
\end{lemma}
\begin{proof}
By definition of $\mathfrak{C}_{2,p_1}(i_0)$, there exists a function $\mathtt{a}\in \mathfrak{C}_2(i_0)$ such that
\[
{a}(\omega,\varphi,x,y) = \mathtt{a}(\omega,\varphi+\tilde{p}_1(\omega,\varphi),x,y).
\]
Thanks to Lemma~\ref{sd22sd}, there exists a symbol $\mathfrak{b}\in \mathfrak{S}^{\alpha-3}(i_0)$ such that 
\begin{align}\label{2skkksd22}
Op^W(\mathfrak{b}) = \Upsilon_{\mathtt{a}}^{\alpha-3}.
\end{align}
We set
\begin{align}\label{sd2sdkjkk2sd}
\mathfrak{a}(\omega,\varphi,x,\xi):=\mathfrak{b}(\omega,\varphi+\tilde{p}_1(\omega,\varphi),x,\xi).
\end{align}
From the definition of $\mathfrak{S}_{p_1}^m(i_0)$, it is clear that $\mathfrak{a}\in\mathfrak{S}_{p_1}^{\alpha-3}(i_0)$.
Then $Op^W(\mathfrak{a}) =\Upsilon^{\alpha-3}_{a}$ follows trivially by reparametrizing $\varphi$ in \eqref{2skkksd22}. 
\end{proof}

As in Lemma~\ref{meanisbetter}, we can estimate the averaged symbol:
\begin{lemma}\label{jjsdsdwsdj0orsd}
Let $\mathfrak{a}\in \mathfrak{S}^m_{p_1}(i_0)$. Then, ${M}_x(\mathfrak{a})$ has a decomposition:
\begin{align}\label{symbsddefs}
{M}_x(\mathfrak{a}) = \epsilon^2 \mathfrak{m}(\omega,\xi) +\mathfrak{r}(\omega,\xi) + \tilde{\mathfrak{q}} (\omega,\varphi,\xi), \quad \int_{\mathbb{T}^\nu}\tilde{\mathfrak{q}}(\omega,\varphi,\xi)d\varphi = 0,
\end{align}
such that the symbols $\mathfrak{m},\mathfrak{r}$ are independent of $\varphi,x$ and $\mathfrak{m}$ does not depend on $i_0$. Furthermore, we have the estimates:
\begin{align}
|\mathfrak{m}|^{\Lip(\gamma,\Omega_1)}_{m,0,\eta_0}&\le_{\mathtt{pe},\eta_0} 1, \label{kksdsd22jjsd}\\
|\mathfrak{r}|^{\Lip(\gamma,\Omega_1)}_{m,0,\eta_0}&\le_{\mathtt{pe},\eta_0} \epsilon^{7-4b},\quad |d_i\mathfrak{r}(i_0)[\ihat]|_{m,0,\eta_0}\le_{\mathtt{pe},\eta_0} \epsilon\rVert \ihat \rVert_{s_0+\mu_0},\label{ksd2spp},\\
|\tilde{\mathfrak{q}}|^{\Lip(\gamma,\Omega_1)}_{m,s,\eta_0}&\le_{\mathtt{pe},s,\eta_0}  \epsilon^4 + \epsilon^2 \rVert \mathfrak{I}_\delta \rVert_{s+\mu_0}^{\Lip(\gamma,\Omega_1)},\nonumber \\ |d_i(\tilde{\mathfrak{q}})(i_0)[\ihat]|_{m,s,\eta_0} & \le_{\mathtt{pe},s,\eta_0} \epsilon^2 \left( \rVert \ihat \rVert_{s+\mu_0} + \rVert \mathfrak{I}_\delta\rVert_{s+\mu_0}\rVert \ihat \rVert_{s_0+\mu_0}\right).\label{vktmxkss}
\end{align}
\end{lemma}
\begin{proof}
By Definition~\ref{sdjjsdsdsddsgsd} and Definition~\ref{symbolssssef}, one can find $\mathfrak{a}_1,\mathfrak{a}_2,\mathfrak{a}_3,\mathfrak{q}_1$ such that
\begin{align}\label{krhcnwkdsd}
\mathfrak{a}(\varphi,x,\xi) = \epsilon\mathfrak{a}_1(\varphi,x,\xi)  + \epsilon^2\mathfrak{a}_2(\varphi,x,\xi)  + \epsilon^3\mathfrak{a}_3(\varphi,x,\xi)  + \mathfrak{q}_{1}(\varphi ,x,\xi) , 
\end{align}
such that
\begin{equation}\label{sj2tksdsdmmzzz}
\begin{aligned}
\mathfrak{a}_i(\varphi,x,\xi) & = \sum_{j_{k_1},\ldots,j_{k_i}\in S}C_{j_{k_1},\ldots,j_{k_i}}(\xi)\sqrt{j_{k_1}\zeta_{k_1}}\cdots \sqrt{j_{k_i}\zeta_{k_i}} \\
& \qquad \qquad \qquad  \times e^{\ii ((\mathtt{l}(j_{k_1})+\cdots+\mathtt{l}(j_{k_i}))\cdot \left(\varphi +\omega \tilde{p}_1(\varphi)\right) +\ii j_k x)},\\
\mathfrak{q}_1(\varphi,x,\xi) &= \mathfrak{q}(\varphi + \tilde{p}_1(\varphi),x,\xi),
\end{aligned}
\end{equation}
where $\mathfrak{q}(\varphi,x,\xi)$ satisfies the estimates in \eqref{q_estimate_21} and \eqref{q_estimate_212}. As in Lemma~\ref{meanisbetter}, we decompose $M_x(\mathfrak{a})$ as (using \eqref{krhcnwkdsd}) 
\begin{align*}
M_x(\mathfrak{a}) &=  \epsilon M_x(\mathfrak{a}_1) + \epsilon^2 M_x(\mathfrak{a}_2) + \epsilon^3M_x(\mathfrak{a}_3) + M_x(\mathfrak{q}_1) \\
& = \epsilon^2 M_{x}(\mathfrak{a}_2) + \epsilon^3M_x(\mathfrak{a}_3) + M_{x}(\mathfrak{q}_1)\\
& = \epsilon^2 \underbrace{M_{x}(\mathfrak{a}_2) }_{=:\mathfrak{m}(\omega,\xi)}+\underbrace{\epsilon^3M_x(\mathfrak{a}_3)  + M_{\varphi,x}(\mathfrak{q}_{1})}_{=:\mathfrak{r}(\omega,\xi)} + \underbrace{ M_{x}(\mathfrak{q}_1) - M_{\varphi,x}(\mathfrak{q}_{1})}_{=:\tilde{\mathfrak{q}}}.
\end{align*}
Indeed, recalling the definition of $\mathtt{l}(j)$ from \eqref{def_lll}, we see that the zero-th Fourier mode of $\mathfrak{a}_i$ for $i=1,2,3$ in \eqref{sj2tksdsdmmzzz} does not depend on $\varphi$, hence
 $\mathfrak{m}$ and $\mathfrak{r}$ are independent of $\varphi,x$. Since $C_{j_1,j_2}(\xi)$ is independent of  $i_0$, $\mathfrak{m}$ does not depend on $i_0$ either. Recalling the dependence between $\zeta$ and $\omega$ from \eqref{normal_formaa3}, we have \eqref{kksdsd22jjsd}.

To see \eqref{ksd2spp}, it is clear from $\mathfrak{a}_3$ in \eqref{sj2tksdsdmmzzz} that $\epsilon^3M_x(\mathfrak{a}_3)$ does not depend on $\varphi,x$ and satisfies \eqref{ksd2spp} (indeed, it satisfies  better estimates: $|\epsilon^3M_x(\mathfrak{a}_3)|_{m,0,\eta_0}^{\Lip(\gamma,\Omega_1)}=O(\epsilon^3)$ and $d_i\mathfrak{a}_3(i_0)=0$). Furthermore, it follows from the estimates of  $\tilde{p}_1$ in Lemma~\ref{p_tilde_estimate11}  and $\mathfrak{q}$, satisfying \eqref{q_estimate_21}, \eqref{q_estimate_212}  that $M_{\varphi,x}(\mathfrak{q}_1)$ satisfies \eqref{ksd2spp} and $\tilde{\mathfrak{q}}$ satisfies \eqref{vktmxkss}.
\end{proof}

\begin{lemma}\label{sdjjsdocmosds}
Let $\mathfrak{a}\in \mathfrak{S}^{m}_{p_1}(i_0)$, and $\mathfrak{b}\in \mathfrak{S}^{m'}_{p_1}(i_0)$. Then, $\mathfrak{a}\mathfrak{b}, \mathfrak{a}\star \mathfrak{b} \in \mathfrak{S}^{m+m'}_{p_1}(i_0)$.
\end{lemma}
\begin{proof}
It suffices to show that if $\mathfrak{a}\in \mathfrak{S}^m(i_0)$ and $\mathfrak{b}\in \mathfrak{S}^{m'}(i_0)$, then $\mathfrak{a}\mathfrak{b},\mathfrak{a}\star\mathfrak{b}\in \mathfrak{S}^{m+m'}(i_0)$. This follows straightforwardly  from the definition of $\mathfrak{S}^{m}(i_0)$ in Definition~\ref{symbolssssef}, Lemma~\ref{product_symbols} and Lemma~\ref{compandkskd2sd}.
\end{proof}

\begin{lemma}\label{sizeofthesymbolsdsd}
Let $\mathfrak{a}\in \mathfrak{S}^m_{p_1}(i_0)$. Then, it holds that
\begin{align}
|\mathfrak{a}|^{\Lip(\gamma,\Omega_1)}_{m,s,\eta_0}&\le_{\mathtt{pe},s,\eta_0} \epsilon(1+ \rVert \mathfrak{I}_\delta\rVert^{\Lip(\gamma,\Omega_1)}_{s+\mu_0}), \label{rmsidsdw}\\
\rVert d_i\mathfrak{a}(i_0)[\ihat]\rVert_{m,s,\eta_0}&\le_{\mathtt{pe},s,\eta_0} \epsilon^3\gamma^{-1}\left( \rVert \ihat \rVert_{s+\mu_0} + \rVert \mathfrak{I}_\delta\rVert_{s+\mu_0}\rVert \ihat \rVert_{s_0+\mu_0}\right).\label{rmsid2s}
\end{align}
Furthermore, $\mathfrak{a}$ admits an expansion 
\[
\mathfrak{a}=\epsilon\mathfrak{b}_1 + \epsilon^2 \mathfrak{b}_2 +\mathfrak{q}_{\mathfrak{a}},
\]
for some $\mathfrak{b}_1\in \mathfrak{B}^m_1,\mathfrak{b}_2\in \mathfrak{B}^m_2$ and $\mathfrak{q}_{\mathfrak{a}}$ satisfies
\begin{equation}\label{qdsd2sd2k}
\begin{aligned}
|\mathfrak{q}_{\mathfrak{a}}|^{\Lip(\gamma,\Omega_1)}_{m,s,\eta_0}&\le_{\mathtt{pe},s,\eta_0}\gamma^{-1}\left( \epsilon^5+\epsilon^3\rVert \mathfrak{I}_\delta\rVert^{\Lip(\gamma,\Omega_1)}_{s+\mu_0}\right),\\
|d_i\mathfrak{q}_{\mathfrak{a}}(i_0)[\ihat]|_{m,s,\eta_0}&\le_{\mathtt{pe},s,\eta_0}\epsilon^3\gamma^{-1}\left( \rVert \ihat \rVert_{s+\mu_0} + \rVert \mathfrak{I}_\delta\rVert_{s+\mu_0}\rVert \ihat \rVert_{s_0+\mu_0}\right).
\end{aligned}
\end{equation}
\end{lemma}
\begin{proof}
The estimates follow straightforwardly from Definition~\ref{sdjjsdsdsddsgsd} and the estimates for $\tilde{p}_1$ in Lemma~\ref{p_tilde_estimate11}.
\end{proof}
With the symbols in $\mathfrak{S}_{p_1}^m(i_0)$, we can rewrite $\mathcal{L}^2$ as follows:
\begin{lemma}\label{rewrite}
For $\mathcal{L}^2$ in Proposition~\ref{prop_time_rep_92},  there exists a symbol $\mathfrak{d}_{0}\in\mathfrak{S}^{0}_{p_1}(i_0)$ such that
\begin{equation}\label{newl2sds}
\begin{aligned}
\mathcal{L}^2[h] &= \D_\omega h - \Pi_{S^\perp}\partial_x M_2[h] + R_2[h],\\
M_2 &= Op^W(\mathtt{m}_\alpha m_{1,\alpha}(\xi) + \frac{T_\alpha}{4} + \mathfrak{d}_0(\varphi,x,\xi))=:Op^W(\mathfrak{p}_{M_2}(\varphi,x,\xi)),
\end{aligned}
\end{equation}
such that 
\begin{enumerate}[label=(\arabic*)]
\item \label{resd2sd} $\mathfrak{d}_0$ is a real-valued symbol.
\item \label{sdjjjjsdsd} $\mathfrak{d}_0$ is a reversible symbol.
\item \label{sdjsjdjsds2sdsd} $Op^W(\mathfrak{d}_0)$ is a real operator.
\item \label{sdsd1122askkssdsd}  $\mathfrak{d}_0$ is $\frac{2\pi}{\mathtt{M}}$-translation invariant, that is,
\begin{align}\label{dd1sdxx}
\rho_{\mathtt{M}}(\mathfrak{d}_0)(\varphi,x,\xi) := \mathfrak{d}_0(\varphi,x+\frac{2\pi}{\mathtt{M}},\xi) = \mathfrak{d}_0(\varphi,x,\xi).
\end{align}
\end{enumerate}
\end{lemma}
\begin{proof}
Recalling $M_2$ from \eqref{prop_time_rep_92}, we need choose $\mathfrak{d}_0$ so that
\begin{align}\label{sd21xxxxsd}
M_2 = \mathtt{m}_\alpha \Lambda^{\alpha-1} + \frac{T_\alpha}{4} + \left(d_2-\frac{T_\alpha}{4}\right) + \Upsilon^{\alpha-3}_{d_3} = Op^W\left(\mathtt{m}_\alpha m_{1,\alpha}(\xi) + \frac{T_\alpha}4 + \mathfrak{d}_0\right).
\end{align}
Thanks to the definition of $m_{1,\alpha}$ in Proposition~\ref{nichts2}, it suffices to choose $\mathfrak{d}_0$ such that
\begin{align}\label{sd22x0o2ssd}
\left( d_2 - \frac{T_\alpha}{4} \right) + \Upsilon^{\alpha-3}_{d_3} = Op^W(\mathfrak{d}_0).
\end{align}
It follows from Lemma~\ref{smysd2} and \ref{sdsdwdsd} of Proposition~\ref{prop_time_rep_92} that such a symbol $\mathfrak{d}_0$ exists and $\mathfrak{d}_0\in \mathfrak{S}^0_{p_1}(i_0)$. Since $M_2$ is  a symmetric operator, \eqref{adjoint_weyl} tells us that $\mathfrak{d}_0$ is real-valued. To prove item \ref{sdjjjjsdsd}, we see from the  symmetry of $d_2,d_3$ in \ref{sdsdwdsd} of Proposition~\ref{prop_time_rep_92} that  $\partial_x \left( d_2 - \frac{T_\alpha}{4}+ \Upsilon^{\alpha-3}_{d_3} \right)  $ is a reversible operator. Hence, $\partial_x Op^W(\mathfrak{d}_0)$ is a reversible operator (see Definition~\ref{sjdwdsd}).
 To see item \ref{sdjsjdjsds2sdsd}, it follows from \ref{2k232sd2d} of Proposition~\ref{prop_time_rep_92} that $Op^W(\mathfrak{d}_0)$ is a real operator.  Lastly, \ref{sdsd1122askkssdsd} follows straightforwardly from the properties of  $d_2$ and $d_3$ in \eqref{b2b322arem}.
\end{proof}

In the rest of this section, we will introduce two natural\index{$\mathtt{b}_0$}\index{$\mathtt{N}_\alpha$} numbers $\mathtt{N}_\alpha$ and $\mathtt{b}_0$, defined to be the smallest integers such that
\begin{align}\label{nsjdjwdsdnsde}
\mathtt{b}_0>  6(\tau+1),\quad \mathtt{N}_\alpha > \max\left\{\frac{3-\alpha}{\alpha-1}, \frac{(s_0+\mathtt{b}_0+2)(2-\alpha) + 2}{\alpha-1} ,3\right\}.
\end{align}
\begin{remark}\label{lsd1sdsdsdsdwds}
In the rest of the paper, the loss of derivatives $\mu_0$ will depend on $\mathtt{N}_\alpha$ and $\mathtt{b}_0$. Since $\mathtt{b}_0$ is completely determined by $\tau$, which is fixed in \eqref{frequency_set2}, we still have $\mu_0$ depend on only $\mathtt{p}$, see Remark~\ref{derv_loss}.
\end{remark}

\subsection{Elimination of the dependence on $x$}\label{subsuhamgqiwkinsl1}
The goal of this section is to make the symbol $\mathfrak{d}_0$ in \eqref{newl2sds} independent of the variable $x$, conjugating $\mathcal{L}^2$ with a flow map of type (2) defined in Definition~\ref{typesoftransformations}.  The result of this section is summarized in the following proposition.

\begin{proposition}\label{induction_egorov}
There exists a symplectic transformation $\Phi_3:H_{S^\perp}\mapsto H_{S^\perp}$ such that
\begin{equation}\label{sdlinesdar2sd11xx2}
\begin{aligned}
\mathcal{L}^3[h]&:=(\Phi_3)^{-1}\mathcal{L}^2\Phi_3[h] = \D_\omega h -  \Pi_{S^\perp}\partial_xM_3[h] +\partial_x\Pi_{S^\perp} W_0 + R_3[h],\\
M_3[h] &:= Op^W\left(\mathtt{m}_\alpha m_{1,\alpha}(\xi) +\frac{T_\alpha}4  + \sum_{k=0}^{\mathtt{N}_\alpha}M_x(\mathfrak{d}_k) + \mathfrak{r}_{-2}\right),
\end{aligned}
\end{equation}
satisfies the following:
\begin{enumerate}[label=(\arabic*)]
\item \label{rhshdosd1}$\mathfrak{d}_k\in \mathfrak{S}^{k(1-\alpha)}_{p_1}$ for $k=0,\ldots,\mathtt{N}_\alpha$ and $\mathfrak{r}_{-2}\in \mathfrak{S}^{-2}_{p_1}$.
\item \label{rhshdosd2} $\mathfrak{d}_k$ and $\mathfrak{r}_{-2}$ are real-valued, reversible symbols (see Definition~\ref{sjdwdsd}) 
\item \label{rhshdosd13} $W_0$ is defined in \eqref{sdkksdsds} and  $R_3\in \mathfrak{R}(i_0)$.
\item \label{rhshdosd14}$\Phi_3$ is real, reversibility preserving and  symplectic. Therefore, $\mathcal{L}^3$ is a real, reversible and Hamiltonian operator.

\item\label{egorsd1ppsd} $\mathcal{L}^3$ and $\Phi_3$ are  $\frac{2\pi}{\mathtt{M}}$-translation invariance preserving and 
\begin{align}\label{b2b3ss1sxare22m}
\rho_{\mathtt{M}}(\mathfrak{r}_{-2}) = \mathfrak{r}_{-2}.
\end{align}

\end{enumerate}
\end{proposition}
\begin{proof}
The proof will be presented throughout the subsection. The expression of $\mathcal{L}^3$ in \eqref{sdlinesdar2sd11xx2} is given in \eqref{new_sldwsdsdsds}, where the symbol of $M_3$ is given in \eqref{sdmsd2ksybsdsd} and \eqref{qformsd2s}. Items \ref{rhshdosd1} and \ref{rhshdosd2} are the results of Proposition~\ref{symobssd}. Item \ref{rhshdosd13} is in Lemma~\ref{ftdsdsd} and item \ref{rhshdosd14} is in Lemma~\ref{finsdsds} and \ref{weartnlsdsd} of Proposition~\ref{symobssd}. Lastly, for item \ref{egorsd1ppsd}, it follows from Lemma~\ref{mfodlsd211sdegosd} that $\mathfrak{a}$, determined by \eqref{asdansitsd} is $\frac{2\pi}{\mathtt{M}}$-translation invariant symbol. Hence the map $\Phi_3$, determined by a Hamiltonian PDE in \eqref{2jsd232},  gives us a $\frac{2\pi}{\mathtt{M}}$-translation invariance preserving transformation. Therefore $\mathcal{L}^3$ is $\frac{2\pi}{\mathtt{M}}$-translation invariance preserving.  \eqref{b2b3ss1sxare22m} is already explicitly stated in Lemma~\ref{mfodlsd211sdegosd}.
\end{proof}

In order to find a symplectic transformation in Proposition~\ref{induction_egorov}, we  use a symplectic transformation of type (2) in Definition~\ref{typesoftransformations}. We denote by $\Phi_3(\tau)$, the flow map of the PDE:
\begin{align}\label{2jsd232}
\partial_\tau u = \Pi_{S^\perp}\partial_x Op^W(\mathfrak{a})u,\text{ for some  real-valued symbol $\mathfrak{a}\in \mathfrak{S}^{1-\alpha}_{p_1}(i_0)$.}
\end{align}
As in Section~\ref{toohard_2_3}, we denote by $\Psi_3(\tau)$, the flow map of the PDE:
\begin{align}\label{k2k92dsd}
\partial_\tau u = \partial_x Op^{W}(\mathfrak{a})u.
\end{align}
When $\tau=1$, we simply denote
\[
\Phi_3 := \Phi_3(1),\quad \Psi_3:=\Psi_3(1).
\]
From \eqref{k2k92dsd}, we have that
\begin{align}
\frac{d}{d\tau}\Psi_3(\tau) & = \partial_x Op^W(\mathfrak{a})\Psi_3(\tau),\quad \frac{d}{d\tau}\Psi_3(\tau)^{-1} = -\Psi_3(\tau)^{-1}\partial_x Op^W(\mathfrak{a}),\nonumber \\
\Psi_3(0) & =\Psi_3(0)^{-1} = I.\label{rlaclsdsd}
\end{align}
Since the equation \eqref{k2k92dsd} is autonomous, we have that
\[
\Psi_3(\tau)\partial_xOp^W(\mathfrak{a}) = \partial_xOp^W(\mathfrak{a})\Psi_3(\tau).
\]
Our analysis in this section is based on the properties of $\Phi_3(\tau),\Psi_3(\tau)$, studied in Subsection~\ref{study_hard_learn_more}.

\subsubsection{Computing the conjugation}

\begin{lemma}\label{conjo_evg}
We have that
\begin{equation}\label{consj2sdsd}
\begin{aligned}
\mathcal{L}^3&:=\Phi_3^{-1}\mathcal{L}^2\Phi_3 = \D_\omega + \Pi_{S^\perp} \Psi_3^{-1}\D_\omega (\Psi_3) - \Pi_{S^\perp}\Psi_3^{-1}\partial_x M_2 \Psi_3 + R_3,\\
R_3&:=R_{3,1}  + R_{3,2}
\end{aligned}
\end{equation}
where $R_{3,1},R_{3,2}$ are defined in \eqref{resj2j2j} and \eqref{resj2j2j2}.
\end{lemma}
\begin{proof}
We compute
\begin{align}
\Phi_3^{-1}\mathcal{L}^2\Phi_3 &= \Phi_3^{-1}\mathcal{L}^2\Pi_{S^\perp}\Psi_3 + \Phi_3^{-1}\mathcal{L}^2 \left(\Phi_3\Pi_{S^\perp} - \Pi_{S^\perp}\Psi_3\Pi_{S^\perp}\right) \nonumber\\
& = \Phi_3^{-1}\mathcal{L}^2 \Psi_3 - \Phi_3^{-1}\mathcal{L}^2\Pi_{S}\Psi_3 + \Phi_3^{-1}\mathcal{L}^2 \left(\Phi_3\Pi_{S^\perp} - \Pi_{S^\perp}\Psi_3\Pi_{S^\perp}\right) \nonumber\\
&  = \Phi_3^{-1}\Pi_{S^{\perp}}\mathcal{L}^2 \Psi_3 + \Phi_s^{-1}\Pi_{S}\mathcal{L}^2\Psi_3 - \Phi_3^{-1}\mathcal{L}^2\Pi_{S}\Psi_3 \nonumber \\
& + \Phi_3^{-1}\mathcal{L}^2 \left(\Phi_3\Pi_{S^\perp} - \Pi_{S^\perp}\Psi_3\Pi_{S^\perp}\right) \nonumber\\
& = \Pi_{S^\perp}\Psi_3^{-1}\Pi_{S^\perp}\mathcal{L}^2 \Psi_{3} + (\Phi_3^{-1}\Pi_{S^\perp} - \Pi_{S^\perp}\Psi_3^{-1}\Pi_{S^\perp})\mathcal{L}^2\Psi_3+ \Phi_s^{-1}\Pi_{S}\mathcal{L}^2\Psi_3 \nonumber\\
& \  - \Phi_3^{-1}\mathcal{L}^2\Pi_{S}\Psi_3 + \Phi_3^{-1}\mathcal{L}^2 \left(\Phi_3\Pi_{S^\perp} - \Pi_{S^\perp}\Psi_3\Pi_{S^\perp}\right) \nonumber\\
& = \Pi_{S^\perp}\Psi_3^{-1}\mathcal{L}^2 \Psi_3 - \Pi_{S^\perp}\Psi_3^{-1}\Pi_S \mathcal{L}^2\Psi_3 + (\Phi_3^{-1}\Pi_{S^\perp} - \Pi_{S^\perp}\Psi_3^{-1}\Pi_{S^\perp})\mathcal{L}^2\Psi_3\nonumber\\
& \  + \Phi_s^{-1}\Pi_{S}\mathcal{L}^2\Psi_3 - \Phi_3^{-1}\mathcal{L}^2\Pi_{S}\Psi_3 + \Phi_3^{-1}\mathcal{L}^2 \left(\Phi_3\Pi_{S^\perp} - \Pi_{S^\perp}\Psi_3\Pi_{S^\perp}\right) \nonumber\\
& =: \Pi_{S^\perp}\Psi_3^{-1}\mathcal{L}^2 \Psi_3 + R_{3,1},\label{rsd2sd2sd2}
\end{align}
where
\begin{equation}\label{resj2j2j}
\begin{aligned}
R_{3,1} &:= - \Pi_{S^\perp}\Psi_3^{-1}\Pi_S \mathcal{L}^2\Psi_3 + (\Phi_3^{-1}\Pi_{S^\perp} - \Pi_{S^\perp}\Psi_3^{-1}\Pi_{S^\perp})\mathcal{L}^2\Psi_3+ \Phi_s^{-1}\Pi_{S}\mathcal{L}^2\Psi_3\\
& \  - \Phi_3^{-1}\mathcal{L}^2\Pi_{S}\Psi_3 + \Phi_3^{-1}\mathcal{L}^2 \left(\Phi_3\Pi_{S^\perp} - \Pi_{S^\perp}\Psi_3\Pi_{S^\perp}\right).
\end{aligned}
\end{equation}
In \eqref{rsd2sd2sd2}, we further decompose $\Psi_{3}^{-1}\mathcal{L}^2 \Psi_3$ as (using \eqref{sdlinear2sd2}),
\begin{align}
& \Pi_{S^\perp}\Psi_{3}^{-1}\mathcal{L}^2\Psi_3\nonumber \\
&= \Pi_{S^\perp}\Psi_3^{-1}\D_\omega \Psi_3 - \Pi_{S^\perp}\Psi_3^{-1}\Pi_{S^\perp}\partial_x M_2\Psi_3 + \Pi_{S^\perp}\Psi_3^{-1}R_2 \Psi_3\nonumber\\
& = \Pi_{S^\perp} \Psi_3^{-1}\D_\omega \Psi_3 - \Pi_{S^\perp}\Psi_3^{-1}\partial_x M_2 \Psi_3 + \Pi_{S^\perp}\Psi_3^{-1}\Pi_S\partial_x M_2 \Psi_3 +\Pi_{S^\perp}\Psi_3^{-1}R_2\Psi_3\nonumber\\
& =: \Pi_{S^\perp}\Psi_3^{-1}\D_\omega \Psi_3 - \Pi_{S^\perp}\Psi_3^{-1}\partial_x M_2 \Psi_3 + R_{3,2},\label{rlacsd2sd}
\end{align}
where
\begin{align}\label{resj2j2j2}
R_{3,2}:=  \Pi_{S^\perp}\Psi_3^{-1}\Pi_S\partial_x M_2 \Psi_3 +\Pi_{S^\perp}\Psi_3^{-1}R_2\Psi_3.
\end{align}
Therefore, writing $\D_\omega \Psi_3 = \D_\omega(\Psi_3) + \Psi_3 \D_\omega$ in \eqref{rlacsd2sd} and plugging it into \eqref{rsd2sd2sd2}, the result follows.
\end{proof}

Now, we analyze  the operators $\Psi_3^{-1}\partial_x M_2 \Psi_3 $ and $ \Psi_3^{-1}\D_\omega (\Psi_3)$. We first observe that each operator can be seen as a solution to a Heisenberg equation.
\subsubsection{Analysis of $\Psi_3^{-1}\partial_x M_2 \Psi_3$.}
We denote
\begin{align}\label{sdkpdefsdxk}
M_3(\tau):=\partial_x^{-1}\Psi_3(\tau)^{-1}\partial_x M_2 \Psi_3(\tau),
\end{align}
so that 
\[
\partial_x M_3(\tau) = \Psi_3(\tau)^{-1}\partial_x M_2 \Psi_3(\tau).
\]
Using \eqref{rlaclsdsd}, we have that
\begin{align}\label{sdsynsdesss}
\frac{d}{d\tau}M_3(\tau) = M_3(\tau)\partial_x Op^W(\mathfrak{a}) - Op^W(\mathfrak{a})\partial_x M_3(\tau), \quad M_3(0)=M_2.
\end{align}
Therefore, we have (recalling the definitions of $\mathfrak{p}_{M_2}$ from \eqref{newl2sds} and the notation $[\cdot,\cdot]_x$ from \eqref{commutator}) that $M_3(\tau)$ solves the Heisenberg equation\index{Heisenberg equation}:
\begin{align}\label{m2solvesheisenberg}
\begin{cases}
\frac{d}{d\tau}M_3(\tau) = [M_3(\tau),Op^W(\mathfrak{a})]_x, \\
M_3(0) = Op^W(\mathfrak{p}_{M_2}).
\end{cases}
\end{align}
Hence, we see from Proposition~\ref{apposdsd} that for the fixed $\mathtt{N}_\alpha$ in \eqref{nsjdjwdsdnsde}, 
\begin{equation}\label{eqssdjk2sd}
\begin{aligned}
M_3(1) &= Q_{M_2} + W_{M_2},\\
Q_{M_2} &= Op^W(\mathfrak{q}_{M_2,\le \mathtt{N}_\alpha}),\quad \mathfrak{q}_{M_2,\le \mathtt{N}_\alpha} := \sum_{n=0}^{\mathtt{N}_\alpha} \frac{1}{n!}\mathfrak{q}_{M_2,n},\\
\mathfrak{q}_{M_2,0}&:=\mathfrak{p}_{M_2},\quad \mathfrak{q}_{M_2,n}:=\mathfrak{q}_{M_2,n-1}\star \mathfrak{a},\text{ for $n=1,\ldots, \mathtt{N}_\alpha$},\\
W_{M_2} &=\frac{1}{\mathtt{N}_\alpha!} \int_0^1 \Psi_3(1-t)^{T}Op^W(\mathfrak{q}_{M_2,\mathtt{N}_\alpha}\star\mathfrak{a})\Psi_3(1-t)t^{\mathtt{N}_\alpha}dt.
\end{aligned}
\end{equation}

\subsubsection{Analysis of $\Psi_3^{-1}\D_\omega (\Psi_3)$.}
\begin{lemma}\label{conju_with_momega_1}(Conjugation with $\D_\omega$) We have that
\begin{align}\label{Q_equation}
\partial_x^{-1}\Psi_3^{-1}\D_\omega (\Psi_3) = \int_0^1S(\tau)d\tau, 
\end{align}
where $S(\tau)$ is a solution to the Heisenberg equation:
\begin{align}\label{sdsdshd2sd}
\begin{cases}
 \frac{d}{d\tau}S(\tau) = [S(\tau),Op^W(\mathfrak{a})]_x,\\
 S(0)=Op^W(\D_\omega \mathfrak{a}).
 \end{cases}
\end{align}
\end{lemma}
\begin{proof}
 Following the computations in \cite[Lemma 6.17]{Berti-Montalto:quasiperiodic-standing-gravity-capillary}, we have that
\begin{align*}
&  \frac{d}{d\tau}\left(\Psi_3(\tau)^{-1}\D_\omega (\Psi_3(\tau))\right) \\
&=-\partial_x\mathcal{A}\Psi_3(\tau)^{-1}\D_\omega(\Psi_3(\tau)) + \Psi_3(\tau)^{-1}\D_\omega(\partial_x\mathcal{A}\Psi_3(\tau))\\
& = -\partial_x\mathcal{A}\Psi_3(\tau)^{-1}\D_\omega(\Psi_3(\tau))  + \Psi_3(\tau)^{-1}\partial_x\D_\omega(\mathcal{A})\Psi_3(\tau) \\
& + \Psi_3(\tau)^{-1}\partial_x\mathcal{A}\D_\omega(\Psi_3(\tau))\\
& = -\Psi_3(\tau)^{-1}\partial_x\mathcal{A}\D_\omega\Psi_3(\tau) +\Psi_3(\tau)^{-1}\partial_xOp^W(\D_\omega\mathfrak{a})\Psi_3(\tau)\\
& + \Psi_3(\tau)^{-1}\partial_x\mathcal{A}\D_\omega(\Psi_3(\tau))\\
& =\Psi_3(\tau)^{-1}\partial_xOp^W(\D_\omega\mathfrak{a})\Psi_3(\tau)\\
& =: \partial_xS(\tau).
\end{align*}
where the third equality follows from \eqref{commuetesd}. $S(\tau)$ is defined to be
\begin{align}\label{stasd2sdx}
S(\tau):=\partial_x^{-1}\Psi_3(\tau)^{-1}\partial_x Op^W(\D_\omega \mathfrak{a})\Psi_3(\tau).
\end{align}
 Then, we have that
\[
\partial_x^{-1}\Psi_3^{-1}\D_\omega (\Psi_3) = \int_0^{1}S(\tau)d\tau,
\]
where we used  $\D_\omega(I) = 0$ and $\Psi_3(0) = I$. To show that $S(\tau)$ satisfies \eqref{sdsdshd2sd}, one can follow the same computations as in $M_2$, in \eqref{sdkpdefsdxk}, since the only difference is that $M_2=Op^W(\mathfrak{p}_{M_2})$ in \eqref{sdkpdefsdxk} is replaced by $Op^W(\D_\omega \mathfrak{a})$ in \eqref{stasd2sdx}.
\end{proof}
As in \eqref{eqssdjk2sd}, we use Proposition~\ref{apposdsd} and obtain the following:
\begin{equation}\label{ksddsdsds}
\begin{aligned}
 & Q_\omega + W_\omega  = \partial_x^{-1}\Psi_3(\tau)^{-1}\D_\omega (\Psi_3(\tau)),\\
Q_\omega& = Op^W(\mathfrak{q}_{\omega,\le \mathtt{N}_\alpha}),\quad \mathfrak{q}_{\omega, \le N_\alpha} := \sum_{n=0}^{\mathtt{N}_\alpha-2}\frac{1}{(n+1)!}\mathfrak{q}_{\omega,n},\\
\mathfrak{q}_{\omega,0} &:= \D_\omega\mathfrak{a},\quad \mathfrak{q}_{\omega,n}:=\mathfrak{q}_{\omega,n-1}\star \mathfrak{a}, \text{ for $n=1,\ldots,\mathtt{N}_\alpha-2$},\\
W_\omega &=\frac{1}{(\mathtt{N}_\alpha-2)!}\int_0^1\int_0^\tau \Psi_3(\tau-t)^{T}Op^W(\mathfrak{q}_{\omega,\mathtt{N}_\alpha-2}\star \mathfrak{a})\Psi_3(\tau-t)t^{\mathtt{N}_\alpha-2}dtd\tau.
\end{aligned}
\end{equation}

Plugging \eqref{ksddsdsds}, \eqref{eqssdjk2sd} into \eqref{consj2sdsd}, we obtain
\begin{align}\label{new_sldwsdsdsds}
\mathcal{L}^3 = \D_\omega -\partial_x \Pi_{S^\perp}\underbrace{\left(Q_{M_2} -Q_{\omega} \right)}_{=:M_3} + \partial_x\Pi_{S^\perp}\underbrace{\left(W_\omega -  W_{M_2}\right)}_{=:W_0} + R_3.
\end{align}

Next, we will choose $\mathfrak{a}$ to remove the dependence on the variable $x$ in the pseudo differential operator $Q_{M_2} -Q_{\omega}$. The term $W_\omega -  W_{M_2}$, as well as $R_3$, will be treated as a remainder. We postpone the estimates for $W_\omega -  W_{M_2}$ to Section~\ref{taksd2tamesd}.

\subsubsection{Choice of $\mathfrak{a}$ to eliminate the dependence on $x$.}
 \begin{lemma}\label{rlacqksd}
If a symbol $\mathfrak{d}$ satisfies $\mathfrak{d}\in \mathfrak{S}^m_{p_1}(i_0)$, then $\frac{\mathfrak{d}}{\mathtt{m}_\alpha\partial_\xi \lambda_{\alpha}}\in \mathfrak{S}^{m+(1-\alpha)}_{p_1}(i_0)$.
\end{lemma}
\begin{proof}
Thanks to Definition~\ref{sdjjsdsdsddsgsd}, it suffices to show that if $\mathfrak{d}\in \mathfrak{S}^{m}(i_0)$, then $\frac{\mathfrak{d}}{\mathtt{m}_\alpha\partial_\xi\lambda_{\alpha}}\in \mathfrak{S}^{m+(1-\alpha)}(i_0)$. This follows from Definition~\ref{symbolssssef}, \ref{extension_ofoperators3} of  Proposition~\ref{nichts2}, Lemma~\ref{product_symbols} and the property of $\mathtt{m}_\alpha$ in \ref{itnsdwd} of Proposition~\ref{sdlinear2sd2}.
\end{proof}

Recalling $\mathfrak{p}_{M_2}$ from \eqref{newl2sds}, we have that
\begin{equation}\label{sdmsd2ksybsdsd}
\begin{aligned}
M_3 & \overset{\eqref{new_sldwsdsdsds}}=Q_{M_2}- Q_{\omega} = Op^W(\mathfrak{q}),\\
\mathfrak{q}&\overset{\eqref{eqssdjk2sd},\eqref{ksddsdsds}}{:=}\mathfrak{q}_{M_2,\le \mathtt{N}_\alpha} - \mathfrak{q}_{\omega,\le \mathtt{N}_\alpha} \\
&= \mathfrak{p}_{M_2} + \mathfrak{p}_{M_2}\star \mathfrak{a} + \sum_{n=2}^{\mathtt{N}_\alpha}\left(\frac{1}{n!} \mathfrak{q}_{M_2,n} - \frac{1}{(n-1)!}\mathfrak{q}_{\omega,n-2}\right) \\
& \overset{\eqref{newl2sds}}= \mathtt{m}_\alpha m_{1,\alpha}(\xi) +\frac{T_\alpha}4 + \mathfrak{d}_0 +  \mathfrak{p}_{M_2}\star \mathfrak{a} + \sum_{n=2}^{\mathtt{N}_\alpha}\left( \frac{1}{n!} \mathfrak{q}_{M_2,n} - \frac{1}{(n-1)!} \mathfrak{q}_{\omega,n-2}\right).
\end{aligned}
\end{equation}
We use the notations of $\left\{ \cdot,\cdot\right\}_x$ and  $r_{\mathfrak{a},\mathfrak{b}}$ defined in \eqref{commutator_pseudo} to expand $\mathfrak{p}_{M_2}\star \mathfrak{a}$. Let us write
\begin{align}\label{qsdessds}
\mathfrak{q} & =\mathtt{m}_\alpha m_{1,\alpha}(\xi) +\frac{T_\alpha}4  + \mathfrak{d}_0 + \left\{ \mathfrak{p}_{M_2},\mathfrak{a}\right\}_x + r_{\mathfrak{p}_{M_2},\mathfrak{a}} + \sum^{\mathtt{N}_\alpha}_{n=2}\left( \frac{1}{n!} \mathfrak{q}_{M_2,n} - \frac{1}{(n-1)!} \mathfrak{q}_{\omega,n-2}\right)\nonumber\\
&\overset{\eqref{newl2sds}}= \mathtt{m}_\alpha m_{1,\alpha}(\xi) +\frac{T_\alpha}4  + \mathfrak{d}_0 + \left\{ \mathtt{m}_\alpha m_{1,\alpha}(\xi) +\frac{T_\alpha}4 ,\mathfrak{a}\right\}_x + \left\{ \mathfrak{d}_0,\mathfrak{a}\right\}_x+r_{\mathfrak{p}_{M_2},\mathfrak{a}} \nonumber \\
& + \sum^{\mathtt{N}_\alpha}_{n=2}\left( \frac{1}{n!} \mathfrak{q}_{M_2,n} - \frac{1}{(n-1)!} \mathfrak{q}_{\omega,n-2}\right)\nonumber\\
& \overset{\eqref{commutator_pseudo}}= \mathtt{m}_\alpha m_{1,\alpha}(\xi) +\frac{T_\alpha}4 + \mathtt{m}_\alpha\partial_\xi \lambda_{\alpha}(\xi)\partial_x\mathfrak{a} + \underbrace{\mathfrak{d}_0}_{\in \mathcal{S}^0} + \underbrace{\frac{T_\alpha}{4}\partial_x\mathfrak{a} + \left\{ \mathfrak{d}_0,\mathfrak{a}\right\}_x}_{\in \mathcal{S}^{1-\alpha}}  \nonumber \\
& + \sum^{\mathtt{N}_\alpha}_{n=2}\underbrace{\left( \frac{1}{n!} \mathfrak{q}_{M_2,n} - \frac{1}{(n-1)!} \mathfrak{q}_{\omega,n-2}\right)}_{\in \mathcal{S}^{(n-1)(1-\alpha)}} + \underbrace{r_{\mathfrak{p}_{M_2},\mathfrak{a}}}_{\in \mathcal{S}^{-2}}.
\end{align}
In order to choose $\mathfrak{a}$, we make an ansatz:
\begin{align}\label{asdansitsd}
\mathfrak{a}(\omega,\varphi,x,\xi) = \sum_{k=0}^{\mathtt{N}_\alpha}\mathfrak{a}_k(\omega,\varphi,x,\xi),\text{ for some $\mathfrak{a}_k\in \mathfrak{S}^{(k+1)(1-\alpha)}_{p_1}$}.
\end{align}
We use the notation: For $k= 0\ldots \mathtt{N}_\alpha$,
\begin{align}\label{lowershbosdsd}
\mathfrak{a}_{\le k}:= \sum_{n=0}^{k}\mathfrak{a}_n,\quad \mathfrak{a}_{> k}:= \sum_{n=k+1}^{\mathtt{N}_\alpha}\mathfrak{a}_n,
\end{align}
and we define
\begin{align}\label{qsdmsdsd}
\begin{cases}
\mathfrak{q}_{M_2,0,\le -1}:= 0,\\
\mathfrak{q}_{M_2,0,\le k} := \mathfrak{q}_{M_{2},0} = \mathfrak{p}_{M_2},\\
\mathfrak{q}_{M_2,n,\le k} := \mathfrak{q}_{M_2,n-1,\le k}\star \mathfrak{a}_{\le k} & \text{ for $n\ge 1$,}
\end{cases} \nonumber \\
 \begin{cases}
\mathfrak{q}_{\omega,0,\le -1}:=0,\\
\mathfrak{q}_{\omega,0,\le k}:= \mathfrak{q}_{\omega,0} = \D_\omega \mathfrak{a}_{\le k},\\
\mathfrak{q}_{\omega,n,\le k}:= \mathfrak{q}_{\omega,n-1,\le k}\star \mathfrak{a}_{\le k} & \text{ for $n\ge 1$.}
\end{cases}
\end{align}
From \eqref{qsdmsdsd}, it is clear that  each of $\mathfrak{q}_{M_{2},n,\le k}$ and $\mathfrak{q}_{\omega,n,\le k}$ collects the terms in $\mathfrak{q}_{M_2}$ and $\mathfrak{q}_{\omega,n}$ in \eqref{eqssdjk2sd} and \eqref{ksddsdsds} that depend on $\mathfrak{a}_{\le k}$ only, out of $\mathfrak{a}_1,\ldots,\mathfrak{a}_{\mathtt{N}_\alpha}$. Therefore, we have 
\[
\mathfrak{q}_{M_2,n}=\mathfrak{q}_{M_2,n,\le \mathtt{N}_\alpha},\quad \mathfrak{q}_{\omega,n}=\mathfrak{q}_{\omega,n,\le \mathtt{N}_{\alpha}},  \text{ for all $n=0,\ldots \mathtt{N}_\alpha$.}
\]
Then, we can write $\mathfrak{q}_{M_2,n}$ and $\mathfrak{q}_{\omega,n}$ as
\begin{align}
\mathfrak{q}_{M_2,n} & =\sum_{k=1}^{\mathtt{N}_{\alpha}+1}(\mathfrak{q}_{M_2,n,\le (k-1)} - \mathfrak{q}_{M_2,n,\le (k-2)}),\nonumber \\ \mathfrak{q}_{\omega,n} & =\sum_{k=1}^{\mathtt{N}_\alpha + 1}\left(\mathfrak{q}_{\omega,n,\le (k-1)}-\mathfrak{q}_{\omega,n,\le (k-2)} \right).\label{qm2decpm}
\end{align}
Using the decompositions of $\mathfrak{a},\mathfrak{q}_{M_2,n},\mathfrak{q}_{\omega,n}$ in \eqref{asdansitsd} and \eqref{qm2decpm}, we decompose the terms in \eqref{qsdessds} as
\begin{align*}
\mathtt{m}_\alpha\partial_{\xi}\lambda_\alpha(\xi)\partial_x\mathfrak{a}& = \mathtt{m}_\alpha\partial_{\xi}\lambda_\alpha(\xi)\partial_x\mathfrak{a}_0 + \sum_{k=1}^{\mathtt{N}_\alpha}\mathtt{m}_\alpha\partial_{\xi}\lambda_\alpha(\xi)\partial_x\mathfrak{a}_k\\
\frac{T_\alpha}4\partial_x\mathfrak{a} +\left\{ \mathfrak{d}_0, \mathfrak{a}\right\}_x &= \sum_{k=1}^{\mathtt{N}_\alpha}\left(\frac{T_\alpha}4\partial_x\mathfrak{a}_{k-1} + \left\{ \mathfrak{d}_0,\mathfrak{a}_{k-1}\right\}_x\right) + \left(\frac{T_\alpha}4\partial_x\mathfrak{a}_{\mathtt{N}_\alpha} + \left\{ \mathfrak{d}_0,\mathfrak{a}_{\mathtt{N}_\alpha}\right\}_x\right),
\end{align*}
\begin{align*}
& \frac{1}{n!}\mathfrak{q}_{M_2,n} - \frac{1}{(n-1)!}\mathfrak{q}_{\omega,n-2} \\
&=  \sum_{k=1}^{\mathtt{N}_\alpha} \left(\frac{1}{n!}\left(\mathfrak{q}_{M_2,n,\le (k-1)}-\mathfrak{q}_{M_2,n,\le (k-2)}\right) -\frac{1}{(n-1)!}(\mathfrak{q}_{\omega,n-2,\le (k-1)}-\mathfrak{q}_{\omega,n-2,\le (k-2)})  \right)\\
& \ +\frac{1}{n!}\left(\mathfrak{q}_{M_2,n,\le \mathtt{N}_\alpha}-\mathfrak{q}_{M_2,n,\le \mathtt{N}_\alpha-1}\right) -\frac{1}{(n-1)!}(\mathfrak{q}_{\omega,n-2,\le \mathtt{N}_\alpha}-\mathfrak{q}_{\omega,n-2,\le \mathtt{N}_\alpha-1}) .
\end{align*}
Hence, we have
\begin{align}\label{fsrdsdsd}
\frac{T_\alpha}4\partial_x\mathfrak{a} +\left\{ \mathfrak{d}_0, \mathfrak{a}\right\}_x + \sum_{n=2}^{\mathtt{N}_\alpha}\left(\frac{1}{n!}\mathfrak{q}_{M_2,n} - \frac{1}{(n-1)!}\mathfrak{q}_{\omega,n-2}\right)  = \sum_{k=1}^{\mathtt{N}_\alpha} \mathfrak{d}_k +\mathfrak{r}_{\mathtt{N}_\alpha},
\end{align}
where
\begin{align}
\mathfrak{d}_k&:=\frac{T_\alpha}4\partial_x\mathfrak{a}_{k-1} + \left\{ \mathfrak{d}_0,\mathfrak{a}_{k-1}\right\}_x \nonumber \\
& \ + \sum_{n=2}^{\mathtt{N}_\alpha}\left(\frac{1}{n!}\left(\mathfrak{q}_{M_2,n,\le (k-1)}-\mathfrak{q}_{M_2,n,\le (k-2)}\right)\right. \nonumber \\
& \qquad \qquad \qquad \qquad \qquad \left.-\frac{1}{(n-1)!}(\mathfrak{q}_{\omega,n-2,\le (k-1)}-\mathfrak{q}_{\omega,n-2,\le (k-2)})  \right),\label{qqtidlsddeds} \\
\mathfrak{r}_{\mathtt{N}_\alpha}&:=\frac{T_\alpha}4\partial_x\mathfrak{a}_{\mathtt{N}_\alpha} + \left\{ \mathfrak{d}_0,\mathfrak{a}_{\mathtt{N}_\alpha}\right\}_x \nonumber \\
& + \sum_{n=2}^{\mathtt{N}_\alpha}\left(\frac{1}{n!}\left(\mathfrak{q}_{M_2,n,\le \mathtt{N}_\alpha}-\mathfrak{q}_{M_2,n,\le \mathtt{N}_\alpha-1}\right) -\frac{1}{(n-1)!}(\mathfrak{q}_{\omega,n,\le \mathtt{N}_\alpha}-\mathfrak{q}_{\omega,n,\le \mathtt{N}_\alpha-1}) \right).\label{rnalpha2}
\end{align}
Plugging \eqref{fsrdsdsd} and \eqref{asdansitsd} into \eqref{qsdessds}, we obtain
\begin{equation}\label{descsdsd}
\begin{aligned}\mathfrak{q} &= \mathtt{m}_\alpha m_{1,\alpha}(\xi) +\frac{T_\alpha}4 + \left( \mathtt{m}_\alpha\partial_{\xi}\lambda_\alpha(\xi)\partial_x\mathfrak{a}_0 +\mathfrak{d}_0\right) \\
& \ + \sum_{k=1}^{\mathtt{N}_\alpha}\left(\mathtt{m}_\alpha\partial_{\xi}\lambda_\alpha(\xi)\partial_x\mathfrak{a}_k + \mathfrak{d}_k \right) + r_{\mathfrak{p}_{M_2},\mathfrak{a}} + \mathfrak{r}_{\mathtt{N}_\alpha}.
\end{aligned}
\end{equation}

\begin{lemma}\label{rlaksd00232}
For $k\in \left\{ 1 , \ldots , \mathtt{N}_\alpha\right\}$,   the following hold:
\begin{enumerate}[label=(\arabic*)]
\item \label{lsd2sdsd1} If $\mathfrak{a}_i \in \mathfrak{S}^{(i+1)(1-\alpha)}_{p_1}(i_0)$, for $i=0\ldots, k-1$, then, $\mathfrak{q}_{M_2,n,\le (k-1)}-\mathfrak{q}_{M_2,n,\le (k-2)},\ \mathfrak{q}_{\omega,n-2,\le (k-1)}-\mathfrak{q}_{\omega,n-2,\le (k-2)}\in\mathfrak{S}^{k(1-\alpha)}_{p_1}(i_0)$ for $n=\left\{ 2,\ldots,\mathtt{N}_\alpha\right\}$. 
\item  \label{lsd2sdsd2} If $\mathfrak{a}_i$ are real-valued for $i=1,\ldots k-1$, then $\mathfrak{q}_{M_2,n,\le (k-1)},\ \mathfrak{q}_{\omega,n-2,\le k-1}$ are real-valued  for $n=\left\{ 2,\ldots,\mathtt{N}_\alpha\right\}$.
\item  \label{lsd2sdsd3} If $\mathfrak{a}_i$ are  reversibility preserving symbols  for $i=1,\ldots k-1$, then $\mathfrak{q}_{M_2,n,\le (k-1)}$ and $ \mathfrak{q}_{\omega,n-2,\le k-1}$ are reversible symbols for $n=\left\{ 2,\ldots,\mathtt{N}_\alpha\right\}$.

Similarly, if $Op^W(\mathfrak{a}_i)$ are real operators for $i=1,\ldots,k-1$, then $Op^W(\mathfrak{q}_{M_2,n,\le (k-1)})$ and $ Op^W(\mathfrak{q}_{\omega,n-2,\le k-1})$ are real operators.
\item \label{sdsd2kkksdjrhsdnwd} If $\mathfrak{a}_i \in \mathfrak{S}^{(i+1)(1-\alpha)}_{p_1}(i_0)$ are real-valued, reversibility preserving and $Op^W(\mathfrak{a}_i)$ are real operators  for $i=1\ldots, k-1$, then    $\mathfrak{d}_k\in \mathfrak{S}^{k(1-\alpha)}_{p_1}(i_0)$ is real-valued, reversible and $Op^W(\mathfrak{d}_k)$ is a real operator. 
\end{enumerate}
\end{lemma}
\begin{proof}
We first observe that
\begin{align}\label{rkawkxnlrla}
(\mathtt{m}_\alpha m_{1,\alpha})\star \mathfrak{b} \in \mathfrak{S}^{m+(\alpha-1)}_{p_1}(i_0), \text{ if $\mathfrak{b}\in \mathfrak{S}^m_{p_1}$}.
\end{align}
Indeed, thanks to Proposition~\ref{nichts2}, we have that $m_{1,\alpha}\in \mathcal{S}^{\alpha-1}(i_0)$. Using the definition of $\mathfrak{S}^{m}_{p_1}$,  and the expansions of $\mathtt{m}_\alpha$ in \ref{itnsdwd} of Proposition~\ref{prop_time_rep_92},  \eqref{rkawkxnlrla} follows straightfowardly. Now we prove each item in the proposition separate.

\vspace{0.5\baselineskip}
\noindent\textit{Proof of \ref{lsd2sdsd1}.}
We prove only
\begin{align}\label{clsdisdsd}
\mathfrak{q}_{M_2,n,\le (k-1)}-\mathfrak{q}_{M_2,n,\le (k-2)}\in\mathfrak{S}^{k(1-\alpha)}_{p_1}(i_0),
\end{align} since $\mathfrak{q}_{\omega,n-2,\le (k-1)}-\mathfrak{q}_{\omega,n-2,\le (k-2)}\in\mathfrak{S}^{k(1-\alpha)}_{p_1}$ can be proved in the same way.

We prove this by induction.
If $n=2$, then it follows from definition of  $\mathfrak{q}_{M_2,n,\le k}$ in \eqref{qsdmsdsd} that
\begin{align}
& \mathfrak{q}_{M_2,2, \le (k-1)} - \mathfrak{q}_{M_2,2,\le (k-2)} = \nonumber \\
& \begin{cases}
\left( \mathfrak{p}_{M_2}\star \mathfrak{a}_0\right)\star \mathfrak{a}_0, & \text{ if $k=1$,}\\
 (\mathfrak{p}_{M_2}\star \mathfrak{a}_{k-1}) \star \mathfrak{a}_{\le (k-1)} +(\mathfrak{p}_{M_2}\star \mathfrak{a}_{\le (k-1)})\star \mathfrak{a}_{k-1} & \text{ if $k\ge 2$}.
\end{cases}\label{sjj22sdsd}
\end{align}
 Recalling $\mathfrak{p}_{M_2}$ from \eqref{newl2sds}, we see from \eqref{rkawkxnlrla} and Lemma~\ref{sdjjsdocmosds} that 
\begin{align}\label{tnseoqhrma}
\mathfrak{p}_{M_2}\star\mathfrak{a}_{k-1}&\in \mathfrak{S}^{(\alpha-1) + k(1-\alpha)}_{p_1}=\mathfrak{S}^{(k-1)(1-\alpha)}_{p_1}(i_0),
\end{align}
which implies (again by Lemma~\ref{sdjjsdocmosds}) 
\begin{align}\label{sjdsspp12d2sdsdd}
(\mathfrak{p}_{M_2}\star \mathfrak{a}_{k-1}) \star \underbrace{\mathfrak{a}_{\le (k-1)}}_{\in \mathfrak{S}^{1-\alpha}_{p_1}(i_0)}\in \mathfrak{S}^{k(1-\alpha)}_{p_1}(i_0).
\end{align}
 For $k=1$, this already gives \eqref{clsdisdsd} for $n=2$, thanks to \eqref{sjj22sdsd}. For $k\ge 2$, one can easily see from \eqref{tnseoqhrma} and \eqref{lowershbosdsd} that 
 \begin{align}\label{invasion_war}
 \mathfrak{p}_{M_2}\star \mathfrak{a}_{\le (k-1)}\in \mathfrak{S}^{(\alpha-1) + (1-\alpha)}_{p_1}(i_0)=\mathfrak{S}^{0}_{p_1}(i_0),
 \end{align} therefore, $(\mathfrak{p}_{M_2}\star \mathfrak{a}_{\le (k-1)})\star \mathfrak{a}_{k-1}\in \mathfrak{S}^{0 + k(1-\alpha)}_{p_1}(i_0)$. Together with \eqref{sjdsspp12d2sdsdd} and \eqref{sjj22sdsd}, we obtain \eqref{clsdisdsd} for $n=2$.

 Now, we assume \eqref{clsdisdsd} holds for  some $n\ge 2$ and we aim to prove it for $n+1$. Again, the definition of $\mathfrak{q}_{M_{2},n,\le k}$ in \eqref{qsdmsdsd} gives us that 
 \begin{align}
 \mathfrak{q}_{M_2,n+1,\le (k-1)}-\mathfrak{q}_{M_2,n+1,\le (k-2)} & =\mathfrak{q}_{M_2,n,\le (k-1)}\star \mathfrak{a}_{k-1} \nonumber \\
 & + (\mathfrak{q}_{M_2,n,\le (k-1)} - \mathfrak{q}_{M_2,n, \le (k-2)})\star \mathfrak{a}_{\le (k-2)}.
 \label{gkqltlqkfsha}
 \end{align}
 For the first term, \eqref{invasion_war} and  \eqref{qsdmsdsd} give us  that $\mathfrak{q}_{M_2,n,\le(k-1)} \in \mathfrak{S}^{ (n-1)(1-\alpha)}_{p_1}(i_0)$. Hence again it follows from Lemma~\ref{sdjjsdocmosds} that $\mathfrak{q}_{M_2,n,\le (k-1)}\star \mathfrak{a}_{k-1} \in \mathfrak{S}_{p_1}^{(n-1)(1-\alpha) + k(1-\alpha)}(i_0)\subset \mathfrak{S}^{k(1-\alpha)}_{p_1}(i_0)$. For the second term in  \eqref{gkqltlqkfsha}, we use the induction hypothesis \eqref{clsdisdsd} for $n$, and Lemma~\ref{sdjjsdocmosds} to obtain  $(\mathfrak{q}_{M_2,n,\le (k-1)} - \mathfrak{q}_{M_2,n, \le (k-2)})\star \mathfrak{a}_{\le (k-2)}\in \mathfrak{S}^{k(1-\alpha) + (1-\alpha)}_{p_1}(i_0)\subset \mathfrak{S}^{k(1-\alpha)}_{p_1}(i_0)$. Hence we obtain \eqref{clsdisdsd} for $n+1$ and this finishes the proof.

\vspace{0.5\baselineskip}
\noindent\textit{Proof of \ref{lsd2sdsd2}.}
We prove $\mathfrak{q}_{M_2,n,\le (k-1)}$ only as above. From \ref{2k232sd2d} of Proposition~\ref{prop_time_rep_92}, $\mathcal{L}^2$ is Hamiltonian operator, hence $M_2$ is a symmetric operator. Hence, it follows from \eqref{adjoint_weyl} and \eqref{newl2sds} that  $\mathfrak{p}_{M_2}$ is a real-valued symbol. Using \eqref{adjoint_weyl} and  the definition of the operator $\star$ in \eqref{syysdsdsd}, we have that $
\mathfrak{p}_{M_2}\star \mathfrak{a}_{\le (k-1)}$ is a real-valued symbol. Hence, recalling the definition of $\mathfrak{q}_{M_2,n,\le (k-1)}$ in \eqref{qsdmsdsd}, it follows immediately that $\mathfrak{q}_{M_2,n,\le (k-1)}$ is a real-valued symbol for any $n\ge 1$.

\vspace{0.5\baselineskip}
\noindent\textit{Proof of \ref{lsd2sdsd3}.}
We prove $\mathfrak{q}_{M_2,n,\le (k-1)}$ only, since the same results for $\mathfrak{q}_{M_2,n,\le (k-1)}$ follow in the same way.  Let us prove the reversibility first. Since $\xi\mapsto m_{1,\alpha}(\xi)$ is even (\ref{extension_ofoperators2} of Proposition~\ref{nichts2}), Lemma~\ref{rewrite} tells us that $\mathfrak{p}_{M_2}$ is a reversible symbol. $\mathfrak{a}_i$ are reversibility preserving symbols for $i=1,\ldots,k-1$, then $\mathfrak{a}_{\le k-1}$ is also a reversibility preserving symbol, hence, it follows from Lemma~\ref{reversible_reversbsd1} that $\mathfrak{p}_{M_2}\star \mathfrak{a}_{\le k-1}$ is a reversible symbol.
Thus, recalling the definition of $\mathfrak{q}_{M_2,n,\le(k-1)}$ from \eqref{qsdmsdsd}, we see that $\mathfrak{q}_{M_2,n,\le(k-1)}$ is obtain by the operator $\star$ of a reversible symbol and a reversibility preserving symbol. Hence, it follows from Lemma~\ref{reversible_reversbsd1} that $\mathfrak{q}_{M_2,n,\le(k-1)}$ is a reversible symbol.

 Similarly, \ref{2k232sd2d} of Proposition~\ref{prop_time_rep_92} tells us that $\mathcal{L}^2$ is a real operator, hence $Op^W(\mathfrak{p}_{M_2})$ is a real operator. Therefore, $Op^W(\mathfrak{p}_{M_2}\star \mathfrak{a}_{\le k-1})$ is a composition of real operators. Again, recalling the definition of $\mathfrak{q}_{M_2,n,\le(k-1)}$ from \eqref{qsdmsdsd}, we see that $\mathfrak{q}_{M_2,n,\le(k-1)}$ is a symbol of a composition of real operators. This implies $Op^W(\mathfrak{q}_{M_2,n,\le(k-1)})$ is a real operator.

\vspace{0.5\baselineskip}
\noindent\textit{Proof of \ref{sdsd2kkksdjrhsdnwd}.} The proof follows straightforwardly from Lemma~\ref{rewrite}, the definition of $\mathfrak{d}_k$ in \eqref{qqtidlsddeds} and the items \ref{lsd2sdsd1}-\ref{lsd2sdsd3} that we just proved.
\end{proof}

\begin{proposition}\label{symobssd} There exist symbols $\mathfrak{a}$, $\mathfrak{d}_k$ for $k=0,\ldots \mathtt{N}_\alpha$, $\mathfrak{r}_{-2}$ such that  $\mathfrak{a}\in \mathfrak{S}^{1-\alpha}_{p_1}(i_0)$, $\mathfrak{d}_{k}\in \mathfrak{S}^{k(1-\alpha)}_{p_1}(i_0)$  and  $\mathfrak{r}_{-2}\in \mathfrak{S}^{-2}_{p_1}(i_0)$ such that
\begin{enumerate}[label=(\arabic*)]
\item \label{itsdsd111}The symbol $\mathfrak{q}$ has the form:
\begin{align}\label{qformsd2s}
\mathfrak{q} =  \mathtt{m}_\alpha m_{1,\alpha}(\xi) +\frac{T_\alpha}4  + \sum_{k=0}^{\mathtt{N}_\alpha}M_x(\mathfrak{d}_k) + \mathfrak{r}_{-2}.
\end{align}
\item \label{itsdsd112}$\mathfrak{a},\mathfrak{d}_k$, for $k=0,\ldots,\mathtt{N}_\alpha$ and $\mathfrak{r}_{-2}$ are real-valued symbols.
\item \label{itsdsd113} $\mathfrak{a}$ is a reversibility preserving symbol, and $\mathfrak{d}_k$, for $k=1,\ldots \mathtt{N}_\alpha$, and $\mathfrak{r}_{-2}$ are reversible symbols.
\item \label{weartnlsdsd} Each of $Op^W(\mathfrak{a}),Op^W(\mathfrak{d}_k)$ and $\Phi_3$ are real operators.
\end{enumerate}
\end{proposition}
\begin{proof}
\vspace{0.5\baselineskip}
\noindent\textit{Proof of \ref{itsdsd111}.} In view of $\mathfrak{q}$ in \eqref{descsdsd}, we choose $\mathfrak{a}$ to be (see \eqref{asdansitsd})
\begin{align}\label{sdsdkkachosd}
\mathfrak{a}_0:= - \frac{\partial_x^{-1}(\mathfrak{d}_0- M_x(\mathfrak{d}_0))}{\mathtt{m}_\alpha \partial_\xi \lambda_\alpha(\xi)},\quad \mathfrak{a}_k:=-\frac{\partial_x^{-1}((\mathfrak{d}_{k}- M_x(\mathfrak{d}_k))}{\mathtt{m}_\alpha\partial_\xi\lambda_\alpha(\xi)},\quad \mathfrak{r}_{-2}:=r_{\mathfrak{p}_{M_2},\mathfrak{a}} + \mathfrak{r}_{\mathtt{N}_\alpha},
\end{align}
which gives us \eqref{qformsd2s}. Note that $\mathfrak{a}_k$ is well-defined, since $\mathfrak{d}_k$ in \eqref{qqtidlsddeds} depends only on $\mathfrak{a}_0,\ldots,\mathfrak{a}_{k-1}$ (by the construction in \eqref{qsdmsdsd}). Using Lemma~\ref{rlacqksd}, Lemma~\ref{rewrite} and \ref{sdsd2kkksdjrhsdnwd} in Lemma~\ref{rlaksd00232}, we see that $\mathfrak{a}_{k}\in \mathfrak{S}^{(k+1)(1-\alpha)}_{p_1}(i_0)$ and $\mathfrak{d}_{k}\in \mathfrak{S}^{k(1-\alpha)}_{p_1}(i_0)$ for $k\ge 0$. In order to see that $\mathfrak{r}_{-2}\in \mathfrak{S}^{-2}_{p_1}(i_0)$, we use the definition of $\mathfrak{r}_{\mathtt{N}_\alpha}$ from \eqref{rnalpha2}, and \ref{lsd2sdsd2} of  Lemma~\ref{rlaksd00232}, which yield $\mathfrak{r}_{\mathtt{N}_\alpha}\in \mathfrak{S}^{(\mathtt{N}_\alpha +1)(1-\alpha)}_{p_1}(i_0)$. Thanks to the choice of $\mathtt{N}_\alpha$ in \eqref{nsjdjwdsdnsde}, we have $(\mathtt{N}_\alpha +1)(1-\alpha) \le -2$. This tells us that
\begin{align}\label{rpsdjj2sd22s}
\mathfrak{r}_{\mathtt{N}_\alpha}\in \mathfrak{S}^{-2}_{p_1}(i_0).
\end{align} Also,   $r_{\mathfrak{p}_{M_2},\mathfrak{a}}\in \mathfrak{S}^{-2}_{p_1}(i_0)$ follows straightforwardly from its definition in \eqref{commutator_pseudo}, since $\mathfrak{p}_{M_2}\in \mathfrak{S}^{\alpha-1}_{p_1}(i_0)$ and $\mathfrak{a}\in \mathfrak{S}^{1-\alpha}_{p_1}(i_0)$. Together with  \eqref{rpsdjj2sd22s}, we see that $\mathfrak{r}_{-2}\in \mathfrak{S}^{-2}_{p_1}(i_0)$.

\vspace{0.5\baselineskip}
\noindent\textit{Proof of \ref{itsdsd112} and \ref{itsdsd113}.}
It follows from \ref{sdsd2kkksdjrhsdnwd} of Lemma~\ref{rlaksd00232}, \ref{resd2sd},\ref{sdjjjjsdsd} of Lemma~\ref{rewrite} and \eqref{sdsdkkachosd} that $\mathfrak{a},\mathfrak{d}_k$ satisfy the desired properties. To see $\mathfrak{r}_{-2}$ is a real-valued reversible symbol, we see that $\mathfrak{a}$ being real-valued and reversibility preserving implies that the transformation $\Phi_3$, determined by the PDE \eqref{2jsd232} is symplectic and reversibility preserving. Therefore, it follows from \ref{2k232sd2d} of Proposition~\ref{prop_time_rep_92} and \eqref{consj2sdsd} implies that $\mathcal{L}^3$ is a reversible Hamiltonian operator. Therefore, $M_3$ in \eqref{new_sldwsdsdsds} is a symmetric operator and $\partial_xM_3$ is a reversible operator. Recalling $\mathfrak{q}$ from \eqref{sdmsd2ksybsdsd}, we see that $\mathfrak{q}$ is real-valued and reversible. In \eqref{qformsd2s}, since we already have all the other symbols in $\mathfrak{q}$ are real-valued reversible symbols, we have that $\mathfrak{r}_{-2}$ is also a real-valued reversible symbol.

\vspace{0.5\baselineskip}
\noindent\textit{Proof of \ref{weartnlsdsd}.}
The result for $Op^W(\mathfrak{a})$ and $Op^W(\mathfrak{d}_k)$ follow from \ref{sdsd2kkksdjrhsdnwd} of Lemma~\ref{rlaksd00232}, \ref{sdjsjdjsds2sdsd} of Lemma~\ref{rewrite} and \eqref{sdsdkkachosd}. Since $\Phi_3$ is time-$1$ flow map determined by the PDE in \eqref{k2k92dsd}, $Op^W(\mathfrak{a})$ being a real operator implies $\Phi_3$ is a real operator as well.
\end{proof}

\subsubsection{Analysis of the transformation $\Phi_3$.}
For  the symbol $\mathfrak{a}$ chosen  in Proposition~\ref{symobssd}, we study some properties of the flow map $\Phi_3(\tau),\Psi(\tau)$, determined by \eqref{2jsd232} and \eqref{k2k92dsd} respectively.

\begin{lemma}\label{align2sd2sd2}
$\Phi^{\pm} = \Phi_3,\Phi_3^{-1},\Psi_3,\Psi_3^{-1}$ satisfies 
\begin{align}
&\rVert (\Phi^{\pm}-I) h \rVert_{s}^{\Lip(\gamma,\Omega_1)} \le_{\mathtt{pe},s} \epsilon\left( \rVert h \rVert^{\Lip(\gamma,\Omega_1)}_{s+\mu_0} + \rVert \mathfrak{I}_\delta \rVert^{\Lip(\gamma,\Omega_1)}_{s+\mu_0}\rVert h \rVert^{\Lip(\gamma,\Omega_1)}_{s_0+\mu_0}\right)\label{transformation_estimate_clae2sss}\\
&\rVert d_i\Phi^{\pm}(i_0)h[\ihat] \rVert_{s} \le_{\mathtt{pe},s} \left(\rVert h\rVert_{s+\mu_0} + \rVert \mathfrak{I}_\delta \rVert_{s+\mu_0}\rVert h \rVert_{s_0+\mu_0}\right) \rVert \ihat \rVert_{s_0+\mu_0} + \rVert h \rVert_{s_0+\mu_0}\rVert \ihat \rVert_{s+\mu_0}.\label{transformation_estimate_clae2sss2}
\end{align}
\end{lemma}

\begin{proof}
We prove $\Phi=\Psi_3$ only since the other operators can be proved in the same way (see Lemma~\ref{tnsjdsinver}).

Let us consider the estimates for $\Psi_3 -I$ first, and then move on to the estimates with the dependence on the embedding $i_0$.  

It follows from Lemma~\ref{sizeofthesymbolsdsd} that 
\begin{align}
|\mathfrak{a}|^{\Lip(\gamma,\Omega_1)}_{\alpha-1,s,\eta_0}&\le_{\mathtt{pe},s,\eta_0}\epsilon(1+ \rVert \mathfrak{I}_\delta\rVert^{\Lip(\gamma,\Omega_1)}_{s+\mu_0}),\nonumber \\
\quad \rVert d_i\mathfrak{a}(i_0)[\ihat]\rVert_{\alpha-1,s,\eta_0}& \le_{\mathtt{pe},s,\eta_0} \epsilon^3\gamma^{-1}\left(1 + \rVert \mathfrak{I}_\delta \rVert_{s+\mu_0}\right).\label{asdwsdsds}
\end{align}
for some $\mu_0 >0$. Hence, by assuming $\mu_0$ in \eqref{size_assumption_3} is large enough and $\epsilon$ is small enough,  we see that \eqref{smallsdsimpsd} is satisfied (recall our choice of $\mathtt{b}_0$  in \eqref{nsjdjwdsdnsde}, hence $\mu_0$ in \eqref{size_assumption_3} can be chosen depending on $\mathtt{p}$ and $\epsilon$ can be chosen depending on $s$, as well as on $\mathtt{pe}$).
It follows from Proposition~\ref{wjsd2sdsdesd}, Lemma~\ref{lipdeopsonh} and \eqref{asdwsdsds} that $\Psi_3$ (choosing $\mathtt{b}=0$ in \eqref{sdsdsdkkk111111} and \eqref{sdsdsdkkk111113}),  
\begin{align}\label{rlajsd2sd}
\rVert \Psi_3(\tau)h\rVert_{s}^{\Lip(\gamma,\Omega_1)}\le_{\mathtt{pe},s} \rVert h \rVert_{s+\mu_0}^{\Lip(\gamma,\Omega_1)} + \rVert \mathfrak{I}_\delta \rVert_{s+\mu_0}^{\Lip(\gamma,\Omega_1)}\rVert h\rVert_{s_0+\mu_0}.
\end{align}
Then, we can expand $\Psi_3(\tau)$ in \eqref{rlaclsdsd} using its Taylor expansion in $\tau$ to obtain
\[
\Psi_3(\tau) = I + \int_0^{\tau}\partial_x Op^W(\mathfrak{a})\Psi_3(t)dt.
\]
Hence \eqref{transformation_estimate_clae2sss} follows from \eqref{asdwsdsds}, Lemma~\ref{sdt2osd1} and \eqref{rlajsd2sd}.
For \eqref{transformation_estimate_clae2sss2}, the result follows from \ref{rlajsd2sdsd} in Proposition~\ref{wjsd2sdsdesd}, with $\mathtt{b}=0$.
\end{proof}

\begin{lemma}\label{finsdsds}
$\Phi_3$ is a symplectic transformation. Furthermore, $\Phi_3$ is reversibility preserving.
\end{lemma}
\begin{proof}
From Proposition~\ref{symobssd}, our choice of $\mathfrak{a}$  in \eqref{2jsd232} is real-valued. Hence, \eqref{2jsd232} is a Hamiltonian equation, associated to the Hamiltonian $u\mapsto \frac{1}{2}\int_{\mathbb{T}}uOp^W(\mathfrak{a})[u] dx$, with the symplectic form $\sigma$ in \eqref{symplectic}, restricted to $H_{S^\perp}$. Hence its flow map $\Phi_3$ is a symplectic transformation.
\end{proof}

\begin{lemma}\label{mfodlsd211sdegosd}
  $\mathfrak{d}_k, \mathfrak{a}_k$ and $\mathfrak{r}_{-2}$ are  $\frac{2\pi}{\mathtt{M}}$-translation invariant. 
\end{lemma}
\begin{proof}
It follows from \eqref{dd1sdxx} that $\mathfrak{d}_0$ is also $\frac{2\pi}{\mathtt{M}}$-translation invariant. Recalling $\mathfrak{p}_{M_2}$ in \eqref{newl2sds}, it is clear that $\mathfrak{p}_{M_2}$ is also  $\frac{2\pi}{\mathtt{M}}$-translation invariant. Therefore, using the construction of $\mathfrak{q}_{M_2,n,\le k}$ and $\mathfrak{q}_{\omega,n,\le k}$ in \eqref{qsdmsdsd} and the definition of $\mathfrak{d}_k$ in \eqref{descsdsd} that if $\mathfrak{a}_{i}$ is $\frac{2\pi}{\mathtt{M}}$-translation invariant for each $i-1,\ldots, k-1$ then $\mathfrak{a}_{k}$ defined in \eqref{sdsdkkachosd} is also $\frac{2\pi}{\mathtt{M}}$-translation invariant.  In \eqref{sdsdkkachosd}, $\mathfrak{a}_0$ is $\frac{2\pi}{\mathtt{M}}$-translation invariant since so is $\mathfrak{d}_0$, hence this property easily propagates to all $k$. This tells us that $\mathfrak{d}_k$ and $\mathfrak{a}_k$ are $\frac{2\pi}{\mathtt{M}}$-translation invariant for all $k$. For $\mathfrak{r}_{-2}$, we recall its definition in \eqref{sdsdkkachosd}, that is, $\mathfrak{r}_{-2} = r_{\mathfrak{p}_{M_2},\mathfrak{a}}+\mathfrak{r}_{\mathtt{N}_\alpha}$, where $r_{\mathfrak{p}_{M_2},\mathfrak{a}}$ and  $\mathfrak{r}_{\mathtt{N}_\alpha}$ are  from \eqref{rnalpha2} and \eqref{qsdessds}. Since we already have that $\mathfrak{a}_k,\mathfrak{q}_{M_2,n,\le k}$ and $\mathfrak{q}_{\omega,n,\le k}$ are all $\frac{2\pi}{\mathtt{M}}$-translation invariant for all $k$, so is $\mathfrak{r}_{\mathtt{N}_\alpha}$. Recalling $r_{\mathfrak{a},\mathfrak{b}}$ from \eqref{commutator_pseudo}, we see that $r_{\mathfrak{p}_{M_2},\mathfrak{a}}$ is a lower order symbol arising from a composition of two $\frac{2\pi}{\mathtt{M}}$-translation invariance preserving operators. Therefore, it follows straightforwardly that $r_{\mathfrak{p}_{M_2},\mathfrak{a}}$ is also $\frac{2\pi}{\mathtt{M}}$-translation invariant.
\end{proof}
\subsubsection{Analysis of the remainders  $W_0$  and  $R_3$.}
We postpone the detailed estimates for $W_0$ to Section~\ref{taksd2tamesd}. Recalling $W_0$ from \eqref{new_sldwsdsdsds}, \eqref{eqssdjk2sd} and \eqref{ksddsdsds}, we have that 
\begin{equation}\label{sdkksdsds}
\begin{aligned}
W_0 &= \frac{1}{(\mathtt{N}_\alpha-2)!}\int_0^1\int_0^\tau \Psi_3(\tau-t)^{T}Op^W(\mathfrak{q}_{\omega,\mathtt{N}_\alpha-2}\star \mathfrak{a})\Psi_3(\tau-t)t^{\mathtt{N}_\alpha-2}dtd\tau \\
& \ -  \frac{1}{\mathtt{N}_\alpha!}\int_0^1 \Psi_3(1-t)^{T}Op^W(\mathfrak{q}_{M_2,\mathtt{N}_\alpha}\star\mathfrak{a})\Psi_3(1-t)t^{\mathtt{N}_\alpha}dt.
\end{aligned}
\end{equation}
In the next lemma, we estimate the symbols $\mathfrak{q}_{\omega,\mathtt{N}_\alpha-2}\star \mathfrak{a}$ and $\mathfrak{q}_{M_2,\mathtt{N}_\alpha}\star\mathfrak{a}$.
\begin{lemma}\label{sjjdsdsymbolsd}
Let $\mathfrak{q}_{\mathtt{N}_\alpha}$ be $\mathfrak{q}_{\omega,\mathtt{N}_\alpha-2}\star \mathfrak{a}$  or $  \mathfrak{q}_{M_2,\mathtt{N}_\alpha}\star\mathfrak{a}$. Then $\mathfrak{q}_{\mathtt{N}_\alpha}\in \mathcal{S}^{\mathtt{N}_\alpha(1-\alpha)}$ and 
\begin{equation}\label{rlaksd1sdsd2}
\begin{aligned}
|\mathfrak{q}_{\mathtt{N}_\alpha}\star \mathfrak{a} |^{\Lip(\gamma,\Omega_1)}_{\mathtt{N}_\alpha(1-\alpha),s,\eta_0} & \le_{\mathtt{pe},s,\eta_0} \epsilon^{\mathtt{N}_\alpha}(1 + \rVert\mathfrak{I}_\delta\rVert_{s + \mu_0}^{\Lip(\gamma,\Omega_1)}),\\
 |d_i(\mathfrak{q}_{\mathtt{N}_\alpha}\star \mathfrak{a})(i_0)[\ihat]|_{\mathtt{N}_\alpha(1-\alpha),s,\eta_0} & \le_{\mathtt{pe},s,\eta_0} \epsilon^{\mathtt{N}_\alpha + 2}\gamma^{-1}\left( \rVert \ihat \rVert_{s+\mu_0} + \rVert \mathfrak{I}_\delta\rVert_{s+\mu_0}\rVert \ihat \rVert_{s_0+\mu_0}\right). 
\end{aligned}
\end{equation}
\end{lemma}
\begin{proof}
We prove the case where $\mathfrak{q}_{\mathtt{N}_\alpha}= \mathfrak{q}_{M_2,\mathtt{N}_\alpha}\star \mathfrak{a}$ only, since the same proof for  $\mathfrak{q}_{\omega,\mathtt{N}_\alpha-2}\star \mathfrak{a}$ can be easily adapted. 

Since $\mathfrak{a}\in \mathfrak{S}^{1-\alpha}_{p_1}(i_0)$ (see Proposition~\ref{symobssd}), it follows from \eqref{rkawkxnlrla}, the definition of $\mathfrak{p}_{M_2}$ in \eqref{newl2sds} and the definition of $\mathfrak{q}_{M_2,1}$ in \eqref{eqssdjk2sd} that $\mathfrak{q}_{M_2,1}\in \mathfrak{S}^{0}_{p_1}$. From Lemma~\ref{sizeofthesymbolsdsd}, we have that
\begin{align*}
|\mathfrak{q}_{M_2,1}|^{\Lip(\gamma,\Omega_1)}_{0,s,\eta_0}&\le_{\mathtt{pe},s,\eta_0} \epsilon(1+ \rVert \mathfrak{I}_\delta\rVert^{\Lip(\gamma,\Omega_1)}_{s+\mu_0}), \\
\rVert d_i\mathfrak{q}_{M_2,1}(i_0)[\ihat]\rVert_{0,s,\eta_0}&\le_{\mathtt{pe},s,\eta_0} \epsilon^3\gamma^{-1}\left( \rVert \ihat \rVert_{s+\mu_0} + \rVert \mathfrak{I}_\delta\rVert_{s+\mu_0}\rVert \ihat \rVert_{s_0+\mu_0}\right).
\end{align*}
Thanks to the above estimates, the desired estimates  \eqref{rlaksd1sdsd2} follow straightforwardly from \eqref{size_assumption_3}, 
 \ref{symocsdwd} in  Lemma~\ref{compandkskd2sd} and the definition of $\mathfrak{q}_{M_2,n}$ in \eqref{eqssdjk2sd}.
\end{proof}

\begin{lemma}\label{ftdsdsd}
$R_3\in \mathfrak{R}(i_0)$.
\end{lemma}
\begin{proof}
Thanks to Lemma~\ref{sizeofthesymbolsdsd} and Lemma~\ref{finitesd2sd}, we have that $\Phi_3\Pi_{S^\perp} - \Pi_{S^\perp}\Psi_3\Pi_{S^\perp}$ and  $\Phi_3^{-1}\Pi_{S^\perp} - \Pi_{S^\perp}\Psi_3^{-1}\Pi_{S^\perp}$ are finite dimensional operators in $\mathfrak{R}(i_0)$.
Recalling $R_3$  from \eqref{consj2sdsd}, \eqref{resj2j2j} and \eqref{resj2j2j2}, we see that   each operator in $R_3$ involves a finite dimensional operator, $\Pi_S$, $\Phi_3\Pi_{S^\perp} - \Pi_{S^\perp}\Psi_3\Pi_{S^\perp}$  or $\Phi_3^{-1}\Pi_{S^\perp} - \Pi_{S^\perp}\Psi_3^{-1}\Pi_{S^\perp}$. Then, the proof of the lemma can be accomplished as in Lemma~\ref{sd2sdsdsddsd}.
 \end{proof}

\subsection{Elimination of the dependence on $\varphi$}\label{subsuhamgqiwkinsl2}
The goal of this section is to make the symbol $\mathfrak{d}_k$ in Proposition~\ref{induction_egorov} independent of the variable $\varphi$ by conjugating $\mathcal{L}^3$ with a flow map of type (4). That is, a transformation $\Phi_4$ will be defined as
\begin{align}\label{phi4defsd2wsd}
\Phi_4h = \sum_{j\in \mathbb{Z}}e^{\ii j \mathfrak{p}_2(\varphi,j)}h_j(\varphi)e^{\ii j x},\quad \Phi_4^{-1}h = \sum_{j\in \mathbb{Z}}e^{-\ii j \mathfrak{p}_2(\varphi,j)}h_j(\varphi)e^{\ii j x},
\end{align}
for some real-valued symbol $\mathfrak{p}_2(\varphi,\xi)\in \mathcal{S}^{1-\alpha}$.
The choice of a symbol $\mathfrak{p}_2$ completely determines the transformation $\Phi_4$. Since $\mathfrak{p}_2$ is real-valued, we have $|e^{\ii j \mathfrak{p}_2(\varphi,j)}| = 1$, hence the transformation $\Phi_4$ is well-defined. 

 Note that $\Phi_4$ can be thought of as a time-$1$ flow map of a Hamiltonian PDE:
 \begin{align}\label{hamisdsd}
 \partial_\tau u =\Pi_{S^\perp} \partial_x Op^W(\mathfrak{p}_2)u,
 \end{align}
therefore, the well-definedness and its tame estimates follow from Proposition~\ref{wjsd2sdsdesd}. However, since the symbol $\mathfrak{p}_2$ is independent of the variable $x$, the analysis is much simpler.

 Before we start computing the conjugation, let us study how a pseudo differential operator can be conjugated by the transformation $\Phi_4$.  For a given symbol $\mathfrak{d} = \mathfrak{d}({\varphi,x,\xi})$, we define an operator:
\begin{align}
T_{\mathfrak{p}_2}[\mathfrak{d}](x,\xi) & :=\sum_{k\in \mathbb{Z}}\widehat{\mathfrak{d}}(k,\xi)e^{\ii((\xi-\frac{k}{2})\mathfrak{p}_2(\varphi,\xi - \frac{k}{2}) - (\xi + \frac{k}2)\mathfrak{p}_2(\varphi,\xi+\frac{k}2))}e^{\ii k x}, \nonumber\\
\text{ where }\widehat{\mathfrak{d}}(\varphi,k,\xi) & :=\frac{1}{2\pi}\int_\mathbb{T}\mathfrak{d}(\varphi,x,\xi)e^{\ii k x}dx. \label{def_Tppsd2}
 \end{align}
 Note that if $\mathfrak{d}\in \mathcal{S}^{m}$, then $T_{\mathfrak{p}_2}[\mathfrak{d}] \in \mathcal{S}^m$ as well \textcolor{black}{ (see Lemma~\ref{jjj2}).} With this notation $T_{\mathfrak{p}_2}$ and \eqref{phi4defsd2wsd}, the conjugation of a pseudo differential operator $Op^W(\mathfrak{d})$ can be computed as (using \eqref{weyl_quant}),
 \begin{align}\label{sdsdohssd}
 \Phi_4^{-1}Op^W(\mathfrak{d})\Phi_4 = Op^W(T_{\mathfrak{p}_2}[\mathfrak{d}]).
 \end{align}

 The result of this subsection is summarized in the next proposition:
 \begin{proposition}\label{rlaqkqdpsanfsdfmf}
 There exists a symplectic transformation $\Phi_4:H_{S^\perp}\mapsto H_{S^\perp}$ such that
\begin{equation}\label{sdlinesdar2sd2sj}
\begin{aligned}
\mathcal{L}^4[h]&:=(\Phi_4)^{-1}\mathcal{L}^3\Phi_4[h] = \D_\omega h -  \Pi_{S^\perp}\partial_xM_4[h] +\partial_x\Pi_{S^\perp} W_1 + R_4[h],\\
M_4[h] &:= Op^W\left(\mathtt{m}_\alpha m_{1,\alpha}(\xi) +\frac{T_\alpha}4  + \mathfrak{m}_{\le0}(\omega,\xi) + \mathfrak{r}_{-2,1}\right),
\end{aligned}
\end{equation}
satisfies the following:
\begin{enumerate}[label=(\arabic*)]
\item \label{rhshdosd122} The symbol $\mathfrak{m}_{\le 0}\in \mathcal{S}^0$ is real, reversible, independent of $\varphi,x$ and admit expansions:
\begin{align}\label{docsd2ksdsdsd}
\mathfrak{m}_{\le 0}(\omega,\xi) = \epsilon^2 \mathfrak{m}_{\le 0,1} + \mathfrak{m}_{\le 0,2},
\end{align}
satisfying
\begin{align}\label{m109symbolmsx}
|\mathfrak{m}_{\le 0,1}|^{\Lip(\gamma,\Omega_1)}_{0,0,\eta_0}&\le_{\eta_0} 1,\\
|\mathfrak{m}_{\le 0,2}|^{\Lip(\gamma,\Omega_1)}_{0,0,\eta_0}&\le_{\eta_0} \epsilon^{7-4b},\quad |d_i\mathfrak{m}_{\le 0,2}(i_0)[\ihat]|_{0,0,\eta_0}\le_{\eta_0} \epsilon\rVert \ihat \rVert_{s_0+\mu_0}.
\end{align}
\item \label{rhshdosd222} $\mathfrak{r}_{-2,1}$ is a real-valued reversible symbol.
\item \label{rhshdosd1322} $W_1$ is defined in \eqref{kk2dsdsd} and  $R_4\in \mathfrak{R}(i_0)$.
\item \label{rhshdosd1422}$\Phi_4$ is real, reversibility preserving and  symplectic. Therefore, $\mathcal{L}^4$ is a real, reversible and Hamiltonian.
\item\label{egor1sd122ppsd}$\mathcal{L}^4$ and $\Phi_4$ are $\frac{2\pi}{\mathtt{M}}$-translation invariance preserving and $\mathfrak{r}_{-2,1}$ is also  $\frac{2\pi}{\mathtt{M}}$-translation invariant.
\end{enumerate}
 \end{proposition}
 
 \begin{proof}
 For the expression for $\mathcal{L}^4$, we use  \eqref{kk2dsdsd} and Lemma~\ref{p2sdsdkcssd} so that we obtain \eqref{sdlinesdar2sd2sj} for
 \[
 \mathfrak{m}_{\le 0}:= M_{\varphi,x}(\sum_{k=0}^{\mathtt{N}_\alpha}\mathfrak{d}_k). 
 \]
 The property of $\mathfrak{m}_{\le 0}$ in item \ref{rhshdosd122} follows from Lemma~\ref{jjsdsdwsdj0orsd}. The property of $R_4$ in item \ref{rhshdosd222} will be proved in Lemma~\ref{sjrhcnskwdsd}. \ref{rhshdosd1322} follows from Lemma~\ref{remisn22sd2ss2}.  Also Lemma~\ref{sjhjsdsda2sdijwhat} tells us that $\Phi_4$ is a real, reversibility preserving symplectic operator. For \ref{egor1sd122ppsd}, it is clear that $\Phi_4$ in \eqref{phi4defsd2wsd} is  $\frac{2\pi}{\mathtt{M}}$-translation invariance preserving, since it does not change the variable $x$. Recalling \eqref{def_Tppsd2}, we see that $\mathfrak{r}_{-2,1}$ in \eqref{kk2dsdsd} is also  $\frac{2\pi}{\mathtt{M}}$-translation invariant, since so is $\mathfrak{r}_{-2}$ (see \eqref{b2b3ss1sxare22m}).
 \end{proof}

\subsubsection{Computing the conjugation}
\begin{lemma}\label{consydasd}
We have that
\begin{equation}\label{kk2dsdsd}
\begin{aligned}
\mathcal{L}^4&:=\Phi_4^{-1}\mathcal{L}^3\Phi_4 = \D_\omega - \Pi_{S^\perp}\partial_x M_4 + \partial_x\Pi_{S^\perp}W_1 + R_4,\\
M_4&:= Op^W\left( \mathtt{m}_\alpha m_{1,\alpha}(\xi) + \frac{T_\alpha}4  + \sum_{k=0}^{\mathtt{N}_\alpha}M_{\varphi,x}(\mathfrak{d}_k) + \mathfrak{r}_{-2,1}\right) \\
& + Op^W\left(\sum_{k=0}^{\mathtt{N}_\alpha}\left(M_x(\mathfrak{d}_k) - M_{\varphi,x}(\mathfrak{d}_k) \right)- \D_\omega \mathfrak{p}_2\right),\\
\mathfrak{r}_{-2,1}& := T_{\mathfrak{p}_2}[\mathfrak{r}_{-2}],\quad W_1 := \Phi_4^{-1}W_0\Phi_4,\quad  R_4 := \Phi_4^{-1}R_3 \Phi_4.
\end{aligned}
\end{equation}
\end{lemma}
\begin{proof}
We compute, recalling $\mathcal{L}^3$ from \eqref{sdlinesdar2sd11xx2} that
\begin{align}\label{kts00dsd}
\Phi_4^{-1}\mathcal{L}^3\Phi_4 = \Phi_4^{-1}\D_\omega \Phi_4 h - \Phi_4^{-1}\Pi_{S^\perp}\partial_x M_3\Phi_4 + \Phi_4^{-1}\partial_x \Pi_{S^\perp}W_0\Phi_4 + \Phi_4^{-1}R_3\Phi_{4},
\end{align}
and compute each term separately:

\textbf{Conjugation of $\D_\omega$.} From \eqref{phi4defsd2wsd}, we have
\[
\D_\omega \Phi_4h = \Phi_4Op^W(\ii \xi \D_\omega \mathfrak{p}_2)h + \Phi_4 \D_\omega h,
\]
therefore,
\begin{align}\label{dswodsdcons}
\Phi_4^{-1}\D_\omega \Phi_4 = Op^W(\ii \xi \D_\omega \mathfrak{p}_2) + \D_\omega = \partial_x Op^W(\D_\omega \mathfrak{p}_2) + \D_\omega.
\end{align}

\textbf{Conjugation of $\Pi_{S^\perp}\partial_x M_3$.} 
From \eqref{phi4defsd2wsd}, it is clear that $\Phi_4$ commutes with $\Pi_{S^\perp}$ and $\partial_x$. Recalling $M_3$ from \eqref{sdlinesdar2sd11xx2}, we have
\begin{align}
\Phi_4^{-1}\Pi_{S^\perp}\partial_x M_3\Phi_4
& = \Pi_{S^\perp}\partial_x \Phi_4^{-1}M_3 \Phi_4\nonumber \\
& = \Pi_{S^\perp}\partial_x \Phi_{4}^{-1}Op^W\left(\mathtt{m}_\alpha m_{1,\alpha}(\xi) +\frac{T_\alpha}4  + \sum_{k=0}^{\mathtt{N}_\alpha}M_x(\mathfrak{d}_k) + \mathfrak{r}_{-2}\right)\Phi_4\nonumber \\
& = \Pi_{S^\perp}\partial_x \left(Op^W\left(\mathtt{m}_\alpha m_{1,\alpha}(\xi) + \frac{T_\alpha}4 \right) \right. \nonumber \\ & \qquad \qquad \left.+ \sum_{k=0}^{\mathtt{N}_\alpha}\Phi_4^{-1}Op^W(M_x(\mathfrak{d}_k))\Phi_4 + \Phi_4^{-1}Op^W(\mathfrak{r}_{-2})\Phi_4 \right).
\label{rksdsdsddkk22}
\end{align}
To compute the conjugation of symbols with $\Phi_4$, we use \eqref{sdsdohssd} to see that
\begin{align*}
\Phi_4^{-1}Op^W(M_x(\mathfrak{d}_k))\Phi_4 & = Op^W(T_{\mathfrak{p}_2}[M_x(\mathfrak{d}_k)]) = Op^W(M_x(\mathfrak{d}_k)),\\
\Phi_4^{-1}Op^W(\mathfrak{r}_{-2})\Phi_4  & = Op^W(T_{\mathfrak{p}_2}[\mathfrak{r}_{-2}]).
\end{align*}
Plugging this into  \eqref{rksdsdsddkk22},  we get
\begin{align}\label{conju1lastnot}
\Phi_4^{-1}\Pi_{S^\perp}\partial_x M_3\Phi_4 = \Pi_{S^\perp}\partial_x Op^W\left( \mathtt{m}_\alpha m_{1,\alpha}(\xi) + \frac{T_\alpha}4  + \sum_{k=0}^{\mathtt{N}_\alpha}M_x(\mathfrak{d}_k) + T_{\mathfrak{p}_2}[\mathfrak{r}_{-2}]\right).
\end{align}

Plugging \eqref{dswodsdcons} and \eqref{conju1lastnot} into \eqref{dswodsdcons} and decomposing $M_x(\mathfrak{d}_k) = M_{\varphi,x}(\mathfrak{d}_k) + (M_{x}(\mathfrak{d}_k) - M_{\varphi,x}(\mathfrak{d}_k))$, we get  the expression for $\mathcal{L}^4$ in \eqref{kk2dsdsd}.
\end{proof}

\subsubsection{Choice of $\mathfrak{p}_2$}
In view of \eqref{kk2dsdsd}, let us denote
\begin{align}\label{dle0isdefinedhere}
\mathfrak{d}_{\ge 0}:=\sum_{k=0}^{\mathtt{N}_\alpha}\mathfrak{d}_k.
\end{align}
Since each $\mathfrak{d}_k\in \mathfrak{S}^{k(1-\alpha)}_{p_1}$, which follows from \ref{rhshdosd1} of  Proposition~\ref{induction_egorov}, we have 
\begin{align}\label{sdsdsd0gesd}
\mathfrak{d}_{\ge 0} \in \mathfrak{S}^0_{p_1}.
\end{align}
We choose $\mathfrak{p}_2$ to be
\begin{align}\label{defosdjsdj2sdsd}
\D_\omega \mathfrak{p}_2 = M_x(\mathfrak{d}_{\ge 0}) -M_{\varphi,x}(\mathfrak{d}_{\ge 0}),\quad \mathfrak{p}_2:=\D^{-1}_\omega( M_x(\mathfrak{d}_{\ge 0}) -M_{\varphi,x}(\mathfrak{d}_{\ge 0}) ).
\end{align}

\begin{lemma}\label{p2sdsdkcssd}
There exists  a real-valued symbol $\mathfrak{p}_2\in \mathcal{S}^0$ such that
\begin{align}\label{sdvashusds}
\sum_{k=0}^{\mathtt{N}_\alpha}\left(M_x(\mathfrak{d}_k) - M_{\varphi,x}(\mathfrak{d}_k) \right) = \D_\omega \mathfrak{p}_2.
\end{align}
Furthermore, $\mathfrak{p}_2$ satisfies 
\begin{equation}\label{p2esd1cxcwrworkdsd}
\begin{aligned}
|\D_\omega \mathfrak{p}_2|^{\Lip(\gamma,\Omega_1)}_{0,s,\eta_0}&\le_{\mathtt{pe},s,\eta_0} \epsilon^4 + \epsilon^2 \rVert \mathfrak{I}_\delta \rVert_{s+\mu_0}^{\Lip(\gamma,\Omega_1)},\\
|d_i(\D_\omega\mathfrak{p}_2)(i_0)[\ihat]|_{0,s,\eta_0}&\le_{\mathtt{pe},s,\eta_0}\epsilon^2\left( \rVert \ihat \rVert_{s+\mu_0} + \rVert \mathfrak{I}_\delta\rVert_{s+\mu_0}\rVert \ihat \rVert_{s_0+\mu_0}\right),\\
|\mathfrak{p}_2|^{\Lip(\gamma,\Omega_1)}_{0,s,\eta_0}&\le_{\mathtt{pe},s,\eta_0}\gamma^{-1}\left(\epsilon^4 + \epsilon^2 \rVert \mathfrak{I}_\delta \rVert_{s+\mu_0}^{\Lip(\gamma,\Omega_1)}\right),\\
 |d_i\mathfrak{p}_2(i_0)[\ihat]|_{0,s,\eta_0}&\le_{\mathtt{pe},s,\eta_0}\epsilon^2\gamma^{-1}\left( \rVert \ihat \rVert_{s+\mu_0} + \rVert \mathfrak{I}_\delta\rVert_{s+\mu_0}\rVert \ihat \rVert_{s_0+\mu_0}\right).
\end{aligned}
\end{equation} 
\end{lemma}
\begin{proof}
From our choice of $\mathfrak{p}_2$ in \eqref{defosdjsdj2sdsd}, \eqref{sdvashusds} follows trivially. Since $\mathfrak{d}_{\ge 0}$ is a real-valued symbol (\ref{rhshdosd1} of  Proposition~\ref{induction_egorov}), $\mathfrak{p}_2$ is also real-valued. Using \eqref{sdsdsd0gesd}, Lemma~\ref{jjsdsdwsdj0orsd}, \eqref{domega_inverse1} and \eqref{domega_inverse}, the desired estimates in \eqref{p2esd1cxcwrworkdsd} for $\mathfrak{p}_2$ follow straightforwardly. 
\end{proof}

\subsubsection{Analysis of the symbol $\mathfrak{r}_{-2,1}$}

\begin{lemma}\label{jjj2}
 Let $\mathfrak{a}(x,\xi)=\mathfrak{a}(\omega,\varphi,x,\xi)\in \mathcal{S}^m$ be a symbol.  With the choice of $\mathfrak{p}_2$ described in Lemma~\ref{p2sdsdkcssd},  $T_{\mathfrak{p}_2}[\mathfrak{a}]\in \mathcal{S}^m$ and there exists $\mu_0=\mu_0(\mathtt{p},\eta)$ such that
\begin{align}
|T_{\mathfrak{p}_2}[\mathfrak{a}]|^{\Lip(\gamma,\Omega_1)}_{m,s,\eta_0}&\le_{\mathtt{pe},s,\eta_0} |\mathfrak{a}|^{\Lip(\gamma,\Omega_1)}_{m,s+\mu_0,\eta_0} +\epsilon^2\gamma^{-1}|\mathfrak{a}|^{\Lip(\gamma,\Omega_1)}_{m,s_0+\mu_0,\eta_0}\rVert \mathfrak{I}_\delta\rVert_{s+\mu_0}^{\Lip(\gamma,\Omega_1)} ,\label{estimate12s2}\\
|d_i T_{\mathfrak{p}_2}[\mathfrak{a}](i_0)[\ihat]|_{m,s,\eta_0}& \le_{\mathtt{pe},s,\eta_0} |\mathfrak{a}|_{m,s+\mu_0,\eta_0}\rVert \ihat \rVert_{s_0+\mu_0} \nonumber \\
& +\epsilon^2\gamma^{-1} |\mathfrak{a}|_{m,s_0+\mu_0,\eta_0}(\rVert \ihat \rVert_{s+\mu_0} + \rVert \ihat \rVert_{s_0+\mu_0}\rVert \mathfrak{I}_\delta \rVert_{s+\mu_0}),
\label{estsdsd22}
\end{align}
\end{lemma}
\begin{proof}
Let us denote $\overline{\mathfrak{p}}(\varphi,\xi):= \xi \mathfrak{p}_2(\varphi,\xi)$. Since $\mathfrak{p}_2\in \mathcal{S}^0$, we have $\overline{\mathfrak{p}}\in \mathcal{S}^1$, and \eqref{p2esd1cxcwrworkdsd} gives us that 
\begin{align}
|\overline{\mathfrak{p}}|^{\Lip(\gamma,\Omega_1)}_{1,s,\eta_0}&\le_{\mathtt{pe},s,\eta_0}\gamma^{-1}\left(\epsilon^4 + \epsilon^2 \rVert \mathfrak{I}_\delta \rVert_{s+\mu_0}^{\Lip(\gamma,\Omega_1)}\right),\label{poverlineestsd1}\\
 |d_i\overline{\mathfrak{p}}(i_0)[\ihat]|_{1,s,\eta_0}&\le_{\mathtt{pe},s,\eta_0}\epsilon^2\gamma^{-1}\left( \rVert \ihat \rVert_{s+\mu_0} + \rVert \mathfrak{I}_\delta\rVert_{s+\mu_0}\rVert \ihat \rVert_{s_0+\mu_0}\right).\label{poverlineestsd2}
\end{align}
We will prove \eqref{estimate12s2} only, since \eqref{estsdsd22} can be done following the same argument, using \eqref{poverlineestsd2}.

We denote (see \eqref{def_Tppsd2} to recall the definition of $T_{\mathfrak{p}_2}$),
\begin{align}\label{sdjwjwj}
A(\varphi,k,\xi):=\int_{-\frac{k}2}^{\frac{k}2}\partial_{\xi}\overline{\mathfrak{p}}(\varphi,\xi+t)dt,\text{ so  that }T_{\mathfrak{p}_2}[\mathfrak{a}](\varphi,k,\xi)=\widehat{\mathfrak{a}}(\varphi,k,\xi)e^{\ii A(\varphi,k,\xi)}.
\end{align}
Now, let us fix $k,\xi$. It is clear that  for each $\eta\in [0,\eta_0]$, 
\begin{align*}
\rVert \partial_{\xi}^{\eta}A(\cdot, k,\xi)\rVert_{H^s_{\varphi}}^{\Lip(\gamma,\Omega_1)}& \le_{\eta} \int_{-\frac{k}2}^{\frac{k}2}\rVert \partial_\xi^{\eta+1}\overline{\mathfrak{p}}(\cdot,\xi+t)\rVert^{\Lip(\gamma,\Omega_1)}_{H^s_{\varphi}}dt \\
& \le_s \int_{-\frac{k}2}^{\frac{k}2}|\overline{\mathfrak{p}}|^{\Lip(\gamma,\Omega_1)}_{1,s,\eta+1}\langle \xi + t\rangle^{-\eta}dt.
\end{align*}
Therefore, using that
\[
\int_{-\frac{k}2}^{\frac{k}2} \langle \xi + t \rangle^{-\eta}dt \le_\eta \begin{cases}
\langle \xi \rangle^{-\eta}, & \text{ if $|\xi| > k$,}\\
1, & \text{ if  $|\xi| < k$,}
\end{cases} \le_\eta \langle k \rangle^{\eta}\langle \xi \rangle^{-\eta},
\]
we see that 
\begin{align}\label{sjdjwsadboud}
\rVert \partial_{\xi}^{\eta}A(\cdot,k,\xi)\rVert_{H^s_\varphi}^{\Lip(\gamma,\Omega_1)}\le_{s,\eta} |\overline{\mathfrak{p}}|^{\Lip(\gamma,\Omega_1)}_{1,s,\eta+1}\langle k\rangle^{\eta}\langle \xi \rangle^{-\eta}.
\end{align}
We claim that 
\begin{align}\label{jjsdksdlws}
 \rVert \partial_{\xi}^{\eta}\left(e^{\ii A(\cdot,k,\xi)}\right)\rVert^{\Lip(\gamma,\Omega_1)}_{H^s_{\varphi}}\le_{s,\eta_0}  \langle k \rangle^{\eta}\langle \xi \rangle^{-\eta} (1+|\overline{\mathfrak{p}}|^{\Lip(\gamma,\Omega_1)}_{1,s+\mu_0,\eta+1}), \text{ for all $0\le \eta\le \eta_0$}
 \end{align}
 for some $\mu_0=\mu_0(\mathtt{p},\eta_0)$. Assuming the claim \eqref{jjsdksdlws} for a moment, let us prove \eqref{estimate12s2}. From $T_{\mathfrak{p}_2}[\mathfrak{a}]$ in \eqref{sdjwjwj}, we have
 \begin{align}
 \rVert  \partial_{\xi}^{\eta}T_{\mathfrak{p}_2}[\mathfrak{a}]\rVert_{s}^2 &\le_s \rVert \partial_\xi^{\eta}T_{\mathfrak{p}_2}[\mathfrak{a}]\rVert_{H^s_xL^2_\varphi}^2 + \rVert\partial_{\xi}^\eta T_{\mathfrak{p}_2}[\mathfrak{a}]\rVert_{L^2_xH^s_\varphi}^2\nonumber\\
 & \le \sum_{\eta_1+\eta_2=\eta}C_{\eta_1,\eta_2}\sum_{k\in\mathbb{Z}}\left(\rVert \partial_{\xi}^{
 \eta_1}\widehat{\mathfrak{a}}(\cdot,k,\xi)\partial_{\xi}^{\eta_2}e^{\ii A(\cdot,k,\xi)}\rVert_{L^2_\varphi}^2\langle k \rangle^{2s} \right. \nonumber \\
 & \left.\qquad \qquad \qquad \qquad \qquad + \rVert \partial_{\xi}^{\eta_1}\widehat{\mathfrak{a}}(\cdot,k,\xi)\partial_{\xi}^{\eta_2}e^{\ii A(\cdot,k,\xi)}\rVert_{H^s_{\varphi}}^2 \right).\label{rkawkrnr}
 \end{align}
 Using the Sobolev embedding theorem, \eqref{jjsdksdlws}, \eqref{poverlineestsd1} and \eqref{size_assumption_3}, we have
 \begin{align}\nonumber
 \rVert \partial_{\xi}^{\eta_1}\widehat{\mathfrak{a}}(\cdot,k,\xi)\partial_{\xi}^{\eta_2}e^{\ii A(\cdot,k,\xi)}\rVert_{L^2_\varphi}^2 & \le \rVert  \partial_{\xi}^{\eta_1}\widehat{\mathfrak{a}}(\cdot,k,\xi)\rVert^2_{L^2_{\varphi}}\rVert \partial_{\xi}^{\eta_2}e^{\ii A(\cdot,k,\xi)}\rVert_{H^{s_0}_{\varphi}}^2 \\& \le\rVert  \partial_{\xi}^{\eta_1}\widehat{\mathfrak{a}}(\cdot,k,\xi)\rVert^2_{L^2_{\varphi}} \langle k\rangle^{2\eta_2}\langle \xi \rangle^{-2\eta_2}. \label{rlaclrnr}
 \end{align}
 Also, similarly, we have (using \eqref{interpolation_2s}),
 \begin{align*}
 & \rVert \partial_{\xi}^{\eta_1}\widehat{\mathfrak{a}}(\cdot,k,\xi)\partial_{\xi}^{\eta_2}e^{\ii A(\cdot,k,\xi)}\rVert_{H^s_{\varphi}}^2 \\
 &\le_s \rVert  \partial_{\xi}^{\eta_1}\widehat{\mathfrak{a}}(\cdot,k,\xi)\rVert^2_{H^s_{\varphi}}\rVert\partial_{\xi}^{\eta_2}e^{\ii A(\cdot,k,\xi)}\rVert_{H^{s_0}_{\varphi}}^2 + \rVert  \partial_{\xi}^{\eta_1}\widehat{\mathfrak{a}}(\cdot,k,\xi)\rVert^2_{H^{s_0}_{\varphi}}\rVert\partial_{\xi}^{\eta_2}e^{\ii A(\cdot,k,\xi)}\rVert_{H^{s}_{\varphi}}^2 \\
 & \le\rVert  \partial_{\xi}^{\eta_1}\widehat{\mathfrak{a}}(\cdot,k,\xi)\rVert^2_{H^s_{\varphi}}\langle k \rangle^{2\eta_2}\langle \xi \rangle^{-2\eta_2} \\
 & \qquad \qquad \qquad + \rVert  \partial_{\xi}^{\eta_1}\widehat{\mathfrak{a}}(\cdot,k,\xi)\rVert^2_{H^{s_0}_{\varphi}}\langle k \rangle^{2\eta_2}\langle \xi \rangle^{-2\eta_2}(1+|\overline{\mathfrak{p}}|^2_{1,s+\mu_0,\eta_2+1})\\
 &\le \langle k \rangle^{2\eta_2}\langle \xi \rangle^{-2\eta_2}\left(\rVert \partial_{\xi}^{\eta_1}\widehat{\mathfrak{a}}(\cdot,k,\xi)\rVert^2_{H^s_{\varphi}}  + \rVert \partial_{\xi}^{\eta_1}\widehat{\mathfrak{a}}(\cdot,k,\xi)\rVert^2_{H^{s_0}_{\varphi}}|\overline{\mathfrak{p}}|_{1,s+\mu_0,\eta_2+1}^2\right)
 \end{align*}
 Plugging this and \eqref{rlaclrnr} into \eqref{rkawkrnr}, we get
 \begin{align*}
  \rVert & \partial_{\xi}^{\eta}T_{\mathfrak{p}_2}[\mathfrak{a}]\rVert_{s}^2\\
  & \le \sum_{\eta_1+\eta_2=\eta}C_{\eta_1,\eta_2}\langle \xi \rangle^{-2\eta_2} \sum_{k\in \mathbb{Z}}\langle k \rangle^{2\eta_2}\\
  & \times\left( \rVert  \partial_{\xi}^{\eta_1}\widehat{\mathfrak{a}}(\cdot,k,\xi)\rVert^2_{L^2_{\varphi}}\langle k \rangle^{2s}  +\rVert \partial_{\xi}^{\eta_1}\widehat{\mathfrak{a}}(\cdot,k,\xi)\rVert^2_{H^s_{\varphi}}  + \rVert \partial_{\xi}^{\eta_1}\widehat{\mathfrak{a}}(\cdot,k,\xi)\rVert^2_{H^{s_0}_{\varphi}}|\overline{\mathfrak{p}}|_{1,s+\mu_0,\eta_2+1}^2 \right)\\
  & \le_{\eta}\sum_{\eta_1+\eta_2=\eta}C_{\eta_1,\eta_2} \langle \xi\rangle^{-2\eta_2} \\
  & \times \left( \rVert \partial_{\xi}^{\eta_1}\mathfrak{a}\rVert^2_{H^{s+\eta_2}_xL^2_{\varphi}}  +\rVert \partial_{\xi}^{\eta_1}\mathfrak{a}\rVert^2_{H^{\eta_2}_xH^{s}_{\varphi}} + \rVert \mathfrak{a}\rVert^2_{H^{\eta_2}_xH^{s_0}_\varphi}|\overline{\mathfrak{p}}|_{1,s+\mu_0,\eta_2+1}^2\right)\\
  & \le \sum_{\eta_1+\eta_2=\eta}C_{\eta_1,\eta_2}\langle \xi \rangle^{-2\eta_2}\left( \rVert \partial_{\xi}^{\eta_1}\mathfrak{a}\rVert^2_{s+\mu_0} + \rVert \partial_{\xi}^{\eta_1}\mathfrak{a}\rVert^2_{s_0+\mu_0}|\overline{\mathfrak{p}}|_{1,s+\mu_0,\eta_2+1}^2\right)\\
  &\le \sum_{\eta_1+\eta_2=\eta}C_{\eta_1,\eta_2}\langle \xi \rangle^{-2\eta_2}\langle \xi \rangle^{2m-2\eta_1}\left( |\mathfrak{a}|^2_{m,s+\mu_0,\eta_1} + |\mathfrak{a}|^2_{m,s_0+\mu_0,\eta_1}|\overline{\mathfrak{p}}|_{1,s+\mu_0,\eta_2+1}^2\right)\\
  & \le_\eta\langle \xi \rangle^{2m - 2\eta}\left( |\mathfrak{a}|^2_{m,s+\mu_0,\eta} + |\mathfrak{a}|^2_{m,s_0+\mu_0,\eta}|\overline{\mathfrak{p}}|_{1,s+\mu_0,\eta+1}^2\right).
 \end{align*}
 Therefore, we obtain
 \begin{align*}
 |T_{\mathfrak{p}_2}[\mathfrak{a}]|_{m,s,\eta}&\le_{\mathtt{pe},s,\eta} |\mathfrak{a}|_{m,s+\mu_0,\eta} + |\mathfrak{a}|_{m,s_0+\mu_0,\eta}|\overline{\mathfrak{p}}|_{1,s+\mu_0,\eta+1}\\
 &\le_{\mathtt{pe},s,\eta}  |\mathfrak{a}|_{m,s+\mu_0,\eta} +\gamma^{-1}\epsilon^2 |\mathfrak{a}|_{m,s_0+\mu_0,\eta}\rVert \mathfrak{I}_\delta\rVert_{s+\mu_0},
 \end{align*}
 where the last inequality follows from \eqref{poverlineestsd1}.
 
 For the Lipschitz dependence of $T_{\mathfrak{p}_2}[\mathfrak{a}]$, denoting $\Delta_{12}X:=X(\omega_1)-X(\omega_2)$, for a symbol $X$, we see from $T_{\mathfrak{p}_2}[\mathfrak{a}]$ in \eqref{sdjwjwj} that
 \[
 \Delta_{12}T_{\mathfrak{p}_2}[\mathfrak{a}](\varphi,k,\xi) = \Delta_{12}\widehat{\mathfrak{a}}(\varphi,k,\xi)e^{\ii A(\varphi,k,\xi)} + \widehat{\mathfrak{a}}(\omega_2,\varphi,k,\xi)\Delta_{12}e^{\ii A(\varphi,k,\xi)}.
 \]
 Then, using the same computations as above word by word, it follows straightforwardly that
 \begin{align*}
 \rVert \Delta_{12}T_{\mathfrak{p}_2}[\mathfrak{a}]\rVert_{m,s,\eta} &\le |\Delta_{12}\mathfrak{a}|_{m,s+\mu_0,\eta}|\overline{\mathfrak{p}}|_{1,s_0+\mu_0,\eta+1} + |\Delta_{12}\mathfrak{a}|_{m,s_0+\mu_0,\eta}|\overline{\mathfrak{p}}|_{1,s+\mu_0,\eta+1}\\
 &+ |\mathfrak{a}|_{m,s+\mu_0,\eta}|\Delta_{12}\overline{\mathfrak{p}}|_{1,s_0+\mu_0,\eta+1} + |\mathfrak{a}|_{m,s_0+\mu_0,\eta}|\Delta_{12}\overline{\mathfrak{p}}|_{1,s+\mu_0,\eta+1}.
 \end{align*}
 Plugging \eqref{poverlineestsd1} and taking the supremum over $\omega_1,\omega_2$, we obtain the desired estimates for the Lipschitz norm in \eqref{estimate12s2}.
 
\vspace{0.5\baselineskip}
\noindent\textit{Proof of the claim \eqref{jjsdksdlws}.}

For each $s\ge 0$, we claim that
\begin{align}\label{jsjdsmxjdsds}
\rVert e^{\ii A(\cdot,k,\xi)}\rVert_{H^s_{\varphi}}\le_{\mathtt{pe},s} \left(1+ \rVert A(\cdot,k,\xi)\rVert_{H^{s+\mu_0}_{\varphi}}\right),
\end{align}
for some $\mu_0$ depending on $s_0$.
The proof can be achieved using the usual induction argument in $s\ge0$. For $s=0$, the above claim holds trivially since $A$ is real-valued. Assuming the claim holds for $s\ge0$, we have
\begin{align*}
& \rVert \partial_{\varphi}^{s+1}e^{\ii A(\cdot,k,\xi)}\rVert_{L^2_{\varphi}} \\
&= \rVert \partial_{\varphi}^{s}(\ii \partial_{\varphi}(A) e^{\ii A})\rVert_{L^2_{\varphi}}\\
& \le_s \sum_{s_1+s_2=s}C_{s_1,s_2}\rVert \partial_{\varphi}^{s_1+1}(A) (\partial_{\varphi}^{s_2}e^{\ii A})\rVert_{L^2_{\varphi}} \\
& \le \sum_{s_1 +s_2=s}C_{s_1,s_2}\rVert\partial_{\varphi}^{s_1+1}A\rVert_{L^{\infty}_{\varphi}}\rVert \partial_{\varphi}^{s_2}e^{\ii A}\rVert_{L^2_{\varphi}}\\
&\le  \sum_{s_1 +s_2=s}C_{s_1,s_2} \rVert A\rVert_{H^{s_1+1+s_0}_{\varphi}}\rVert e^{\ii A}\rVert_{H^{s_2}_{\varphi}}\\
&\le  \sum_{s_1 +s_2=s}C_{s_1,s_2}\left(C(\delta) \rVert A \rVert_{H^{s_1+s_2+1+s_0}_{\varphi}}\rVert e^{\ii A}\rVert_{L^2} + \delta \rVert A\rVert_{H^{s_0}_{\varphi}}\rVert e^{\ii A}\rVert_{H^{s_1+s_2+1}_{\varphi}}\right)\\
& \le_s C(\delta)\rVert {A}\rVert_{H^{s+1+s_0}_{\varphi}} + C_{s}\delta\rVert A\rVert_{H^{s_0}_{\varphi}}\rVert e^{\ii A}\rVert_{H^{s+1}_{\varphi}},
\end{align*}
for all $\delta>0$, where the third inequality follows from the usual Sobolev embedding, and the fourth inequality follows from Lemma~\ref{intesdsdsd}. Since $\rVert e^{\ii A}\rVert_{L^2_{\varphi}}\le_{\nu} 1$, the above inequality shows that
\begin{align}\label{sdjkksdsd}
(1-C_{s}\delta \rVert A\rVert_{H^{s_0}_{\varphi}})\rVert e^{\ii A}\rVert_{H^{s+1}_{\varphi}}  \le C(\delta)\rVert A\rVert_{H^{s+1+s_0}_{\varphi}}.
\end{align}
Note that from \eqref{sjdjwsadboud} with $\eta:=0$, we have
\[
\rVert A \rVert_{H^{s_0}_{\varphi}} \le_{s_0,\eta} |\overline{\mathfrak{p}}|_{1,s_0,1}\overset{\eqref{poverlineestsd1},\eqref{size_assumption_3}}\le_{\mathtt{pe},\eta} \gamma^{-1}(\epsilon^4 + \epsilon^2 \epsilon^{6-4b})\le1, 
\]
for sufficiently small $\epsilon>0$,
where the last inequality follows from \eqref{frequency_set2}, which says $\gamma=\epsilon^{2b}$. Hence, we can choose $\delta$ depending on $s$ so that \eqref{sdjkksdsd} gives us  \eqref{jsjdsmxjdsds} for $s+1$. 

 Now, we consider the dependence on $\omega$. We claim that
 \begin{align}\label{clajsd2sd}
 \rVert  e^{\ii A(\cdot,k,\xi)}\rVert_{H^s_{\varphi}}^{\Lip(\gamma,\Omega_1)}\le C_s   \left( 1 + |\overline{\mathfrak{p}}|_{1,s+\mu_0,1}^{\Lip(\gamma,\Omega_1)}\right).
 \end{align}
 Denoting $\Delta_{12}A(\varphi,x,\xi):=A(\omega_1,\varphi,x,\xi) -A(\omega_2,\varphi,x,\xi)$ for $\omega_1,\omega_2\in \Omega_1$, we have
 \[
 \Delta_{12}e^{\ii A(\varphi,k,\xi)} = e^{\ii A(\omega_2)}\left(e^{\ii \Delta_{12}A} - 1\right) = \int_0^{1}\Delta_{12}A e^{\ii A(\omega_2)}e^{\ii t \Delta_{12} A}dt.
 \]
 Using \eqref{jsjdsmxjdsds} and \eqref{bncal} that
 \begin{align*}
 \rVert &\Delta_{12}e^{\ii A}\rVert_{H^s_{\varphi}}\\
 & \ \le_s\sup_{t\in[0,1]}\left( \rVert \Delta_{12}A\rVert_{H^s_{\varphi}}\rVert e^{\ii A(\omega_2)}\rVert_{H^{s_0}_{\varphi}}\rVert e^{\ii t \Delta_{12}A}\rVert_{H^{s_0}_\varphi}\right. \\
 & \qquad \qquad \qquad \left.+ \rVert \Delta_{12}A\rVert_{H^{s_0}_{\varphi}}\left( \rVert e^{\ii A}\rVert_{H^{s}_{\varphi}}\rVert e^{\ii t\Delta_{12}A }\rVert_{H^{s_0}_{\varphi}} +  \rVert e^{\ii A}\rVert_{H^{s_0}_{\varphi}}\rVert e^{\ii t\Delta_{12}A }\rVert_{H^{s}_{\varphi}}\right)\right)\\
 & \le  \rVert \Delta_{12}A\rVert_{H^s_{\varphi}} + \rVert \Delta_{12}A\rVert_{H^{s_0}_{\varphi}}\left(1+ \sup_{\omega\in \Omega_1}\rVert A(\omega)\rVert_{H^{s+\mu_0}_{\varphi}}\right),
 \end{align*}
 for some $\mu_0$.
 Hence, using \eqref{sjdjwsadboud} with $\eta=0$ and \eqref{jsjdsmxjdsds}, we obtain \eqref{clajsd2sd}.
 
 Now, we aim to prove the claim~\eqref{jjsdksdlws}
 for some $\mu_0=\mu_0(\mathtt{p},\eta_0)$.
We argue by induction as before. For $\eta=0$, we already have \eqref{clajsd2sd}. Assuming \eqref{jjsdksdlws} holds true for $\eta\ge 0$, let us show it for $\eta+1$. Clearly, we have
\begin{align*}
\partial_{\xi}^{\eta+1}\left(e^{\ii A(\varphi,k,\xi)} \right)& = \partial_{\xi}^{\eta}(\partial_{\xi}A(\varphi,k,\xi)e^{\ii A(\varphi,k,\xi)}) \\
& = \sum_{\eta_1+\eta_2=\eta}\partial_{\xi}^{\eta_1+1}A(\varphi,k,\xi)\partial_{\xi}^{\eta_2}e^{\ii A(\varphi,k,\xi)}.
\end{align*}
Therefore, using \eqref{interpolation_2s} and \eqref{sjdjwsadboud}, we have
\begin{align*}
&\rVert \partial_{\xi}^{\eta_1+1}A(\varphi,k,\xi)\partial_{\xi}^{\eta_2}e^{\ii A(\varphi,k,\xi)}\rVert^{\Lip(\gamma,\Omega_1)}_{H^s_{\varphi}} \\
& \le_{\eta_1,\eta_2}\rVert \partial_{\xi}^{\eta_1+1}A(\cdot,k,\xi)\rVert^{\Lip(\gamma,\Omega_1)}_{H^s_{\varphi}}\rVert \partial_{\xi}^{\eta_2}e^{\ii A(\cdot,k,\xi)}\rVert^{\Lip(\gamma,\Omega_1)}_{H^{s_0}_{\varphi}} \\
& \qquad \qquad + \rVert \partial_{\xi}^{\eta_1+1}A(\cdot,k,\xi)\rVert^{\Lip(\gamma,\Omega_1)}_{H^{s_0}_{\varphi}}\rVert \partial_{\xi}^{\eta_2}e^{\ii A(\cdot,k,\xi)}\rVert^{\Lip(\gamma,\Omega_1)}_{H^{s}_{\varphi}}\\
& \le \langle k \rangle^{\eta_1+1}\langle \xi \rangle^{-(\eta_1+1)}|\overline{\mathfrak{p}}|_{1,s,\eta_1+1}\rVert \partial_{\xi}^{\eta_2}e^{\ii A(\cdot,k,\xi)}\rVert^{\Lip(\gamma,\Omega_1)}_{H^{s_0}_{\varphi}}  \\
& \qquad \qquad + \langle k \rangle^{\eta_1+1}\langle \xi \rangle^{-(\eta_1+1)}|\overline{\mathfrak{p}}|_{1,s_0,\eta_1+1}\rVert \partial_{\xi}^{\eta_2}e^{\ii A(\cdot,k,\xi)}\rVert^{\Lip(\gamma,\Omega_1)}_{H^{s}_{\varphi}}\\
& \le \langle k \rangle^{\eta + 1}\langle \xi \rangle^{-(\eta+1)}\left(|\overline{\mathfrak{p}}|_{1,s,\eta_1+1}\left(1+|\overline{\mathfrak{p}}|_{1,s_0+\mu_0,\eta_2+1}\right)\right.  \\
& \qquad \qquad \qquad \qquad \qquad \left.+ |\overline{\mathfrak{p}}|_{1,s_0,\eta_1+1}\left( 1+ |\overline{\mathfrak{p}}|_{1,s+\mu_0,\eta_2+1}\right) \right)\\
&\le \langle k \rangle^{\eta+1}\langle \xi \rangle^{-(\eta+1)} \left(1 + |\overline{\mathfrak{p}}|_{1,s+\mu_0,\eta_0+N} \right),
\end{align*}
where  the third follows from our induction hypothesis of \eqref{jjsdksdlws}, and the last inequality follows from \eqref{poverlineestsd1} and \eqref{size_assumption_3}. This gives \eqref{jjsdksdlws} for $\eta +1$.
\end{proof}

\begin{lemma}\label{r2sd2sjrevsd}
$\mathfrak{r}_{-2,1}$ admits an expansion:
\begin{align}\label{rlsd1b1b2}
\mathfrak{r}_{-2,1} = \epsilon\mathfrak{b}_1 + \epsilon^2 \mathfrak{b}_2  +\mathfrak{r}_{-2,\le 3},
\end{align}
for some $\mathfrak{b}_1\in \mathfrak{B}_1^{-2}$, $\mathfrak{b}_2\in \mathfrak{B}_2^{-2}$. The symbol $\mathfrak{r}_{-2,\le 3}$ satisfies the estimates:
\begin{equation}\label{lastsxhsdsamurai}
\begin{aligned}
&|\mathfrak{r}_{-2,\le3}|^{\Lip(\gamma,\Omega_1)}_{-2,s,\eta_0} \le_{\mathtt{pe},s,\eta_0}  \gamma^{-1}\left(\epsilon^5 + \epsilon^3\rVert \mathfrak{I}_\delta\rVert^{\Lip(\gamma,\Omega_1)}_{s+\mu_0}\right),\\
&|d_i\mathfrak{r}_{-2,\le3}(i_0)[\ihat]|_{-2,s,\eta_0} \le_{\mathtt{pe},s,\eta_0} \epsilon^3\gamma^{-1}\left( \rVert \ihat \rVert_{s+\mu_0} + \rVert \mathfrak{I}_\delta\rVert_{s+\mu_0}\rVert \ihat \rVert_{s_0+\mu_0}\right).
\end{aligned}
\end{equation}
\end{lemma}

\begin{proof}
We first consider $\mathfrak{r}_{-2}$. Since $\mathfrak{r}_{-2}\in \mathfrak{S}^{-2}_{p_1}(i_0)$, which follows from \ref{rhshdosd1} in Proposition~\ref{induction_egorov}, Definition~\ref{sdjjsdsdsddsgsd} tells us that
\begin{align*}
\mathfrak{r}_{-2}(\varphi,x,\xi) &= \epsilon \mathfrak{b}_1(\varphi+\tilde{p}_1(\varphi),x,\xi) + \epsilon^2 \mathfrak{b}_2(\varphi+\tilde{p}_1(\varphi),x,\xi) \\
& + \epsilon^3\mathfrak{b}_3(\varphi+\tilde{p}_1(\varphi),x,\xi) + \mathfrak{q}(\varphi+\tilde{p}_1(\varphi,x,\xi)),
\end{align*}
where $\mathfrak{b}_i\in \mathfrak{B}^{-2}_{i}$ for $i=1,2,3$ and $\mathfrak{q}$ satisfies the estimates \eqref{q_estimate_21} and \eqref{q_estimate_212}.s
Let us denote
\[
T_{p_1}[\mathfrak{a}](\varphi,x,\xi):=\mathfrak{a}(\varphi+\tilde{p}_1(\varphi),x,\xi), \text{ for a symbol $\mathfrak{a}$},
\]
so that the above expansion for $\mathfrak{r}_{-2}$ can be written as
\[
\mathfrak{r}_{-2}  =\epsilon T_{p_1}[\mathfrak{b}_1] + \epsilon^2 T_{p_1}[\mathfrak{b}_2] + \epsilon^3 T_{p_1}[\mathfrak{b}_3] +  T_{p_1}[\mathfrak{q}]
\]
Therefore, recalling  $\mathfrak{r}_{-2,1}$ from \eqref{kk2dsdsd}, we have
\begin{align*}
& \mathfrak{r}_{-2,1}\\
& =T_{\mathfrak{p}_2}[\mathfrak{r}_{-2}] =  \epsilon T_{\mathfrak{p}_2}[T_{p_1}[\mathfrak{b}_1]]+\epsilon^2 T_{\mathfrak{p}_2}[T_{p_1}[\mathfrak{b}_2]]+\epsilon^3 T_{\mathfrak{p}_2}[T_{p_1}[\mathfrak{b}_3]]+T_{\mathfrak{p}_2}[T_{p_1}[\mathfrak{q}]]\\
& = \epsilon \mathfrak{b}_1 + \epsilon^2 \mathfrak{b}_2 \\
& \ + \underbrace{\epsilon\left( T_{\mathfrak{p}_2}[T_{p_1}[\mathfrak{b}_1]] - \mathfrak{b}_1\right)+\epsilon^2 \left( T_{\mathfrak{p}_2}[T_{p_1}[\mathfrak{b}_2]]-\mathfrak{b}_2\right)+ \epsilon^3 T_{\mathfrak{p}_2}[T_{p_1}[\mathfrak{b}_3]]+T_{\mathfrak{p}_2}[T_{p_1}[\mathfrak{q}]]}_{=: \mathfrak{r}_{-2,\le 3}}.
\end{align*}
Hence, it suffices to show that $\mathfrak{r}_{-2,\le 3}$ satisfies the estimates \eqref{lastsxhsdsamurai}. We will prove that  $\epsilon \left(T_{\mathfrak{p}_2}[T_{p_1}[\mathfrak{b}_1]]-\mathfrak{b}_1\right)$ satisfies \eqref{lastsxhsdsamurai} only, because the other terms are even smaller.

Recalling the definition of $\mathfrak{B}^{-2}_{1}$ from Definition~\ref{buildingbobo} and recalling $\overline{v}$ from \eqref{norm_def_vbar}, we have
\begin{align*}
& \epsilon  \left(T_{\mathfrak{p}_2}[T_{p_1}[\mathfrak{b}_1]]-\mathfrak{b}_1\right)\\
& = \epsilon\sum_{j_k\in S}C_{j_k}(\xi)\sqrt{j_k \zeta_k}\left(e^{\ii (\mathtt{l}(j_k)\cdot( \varphi+\tilde{p}_1(\varphi)) +j_k (x + \mathfrak{p}_2(\varphi,j_k)))} -e^{\ii (\mathtt{l}(j_k)\cdot \varphi + j_k x)}\right)\\
& = \epsilon\sum_{j_k\in S}C_{j_k}(\xi)\sqrt{j_k \zeta_k}e^{\ii (\mathtt{l}(j_k)\cdot \varphi + j_k x)}\left(e^{\ii (\mathtt{l}(j_k)\cdot(\tilde{p}_1(\varphi)) +j_k \mathfrak{p}_2(\varphi,j_k))} -1\right)\\
&=\epsilon\sum_{j_k\in S}C_{j_k}(\xi)\sqrt{j_k \zeta_k}e^{\ii (\mathtt{l}(j_k)\cdot \varphi + j_k x)}\\
& \times \int_0^{1}\left(\ii (\mathtt{l}(j_k)\cdot(\tilde{p}_1(\varphi)) +j_k \mathfrak{p}_2(\varphi,j_k))\right)e^{t\left(\ii (\mathtt{l}(j_k)\cdot(\tilde{p}_1(\varphi)) +j_k \mathfrak{p}_2(\varphi,j_k))\right)}dt
\end{align*}
Then, using the estimates for $\mathfrak{p}_2$ and $\tilde{p}_1$ from \eqref{p2esd1cxcwrworkdsd} and  Lemma~\ref{p_tilde_estimate11} and using Lemma~\ref{jjj2} that
\begin{align*}
& |\epsilon\left(T_{\mathfrak{p}_2}[T_{p_1}[\mathfrak{b}_1]]-\mathfrak{b}_1\right)|^{\Lip(\gamma,\Omega_1)}_{-2,s,\eta_0} \\
&\le_{\mathtt{pe},s,\eta_0} \epsilon \sup_{j_k\in S}|C_{j_k}|_{-2,0,\eta_0}\left( \rVert \tilde{p}_1\rVert^{\Lip(\gamma,\Omega_1)}_{s+\mu_0} + |\mathfrak{p}_2|_{0,s+\mu_0,\eta_0}^{\Lip(\gamma,\Omega_1)} \right)\\
&\le_{\mathtt{pe},s,\eta_0} \epsilon \gamma^{-1}\left(\epsilon^4 + \epsilon^2 \rVert \mathfrak{I}_\delta\rVert^{\Lip(\gamma,\Omega_1)}_{s+\mu_0} \right),
\end{align*}
and
\begin{align*}
|\epsilon d_i&\left(\left(T_{\mathfrak{p}_2}[T_{p_1}[\mathfrak{b}_1]]-\mathfrak{b}_1\right)\right)(i_0)[\ihat]|_{-2,s,\eta_0}\\
& \le_{\mathtt{pe},s,\eta_0}  \epsilon\left(\epsilon^2\gamma^{-1}\left( \rVert \ihat \rVert_{s+\mu_0} + \rVert \mathfrak{I}_\delta\rVert_{s+\mu_0}\rVert \ihat \rVert_{s_0+\mu_0}\right)\right).
\end{align*}
Therefore, $\epsilon \left(T_{\mathfrak{p}_2}[T_{p_1}[\mathfrak{b}_1]]-\mathfrak{b}_1\right)$ satisfies \eqref{lastsxhsdsamurai}.
\end{proof}

\begin{lemma}\label{veryexahsysdwds}
$\mathfrak{p}_2$ is a reversibility preserving symbol.
\end{lemma}

\begin{proof}
It follows from \ref{rhshdosd2} of Proposition~\ref{induction_egorov} and \eqref{dle0isdefinedhere} that $\mathfrak{d}_{\ge 0}$ is a reversible symbol, therefore, $M_x(\mathfrak{d}_{\ge 0}) - M_{\varphi,x}(\mathfrak{d}_{\ge 0})$ is also a reversible symbol. Hence recalling our choice of $\mathfrak{p}_2$ from \eqref{defosdjsdj2sdsd}, we have that $\mathfrak{p}_2$ is a reversibility preserving symbol.
\end{proof}

\begin{lemma}\label{sjrhcnskwdsd}
$\mathfrak{r}_{-2,1}$ is real-valued.  It is also a  reversible symbol.
\end{lemma}
\begin{proof}
The lemma follows straightforwardly  from Lemma~\ref{p2sdsdkcssd}, which says $\mathfrak{p}_2$ is a real-valued symbol, \ref{rhshdosd2} of Proposition~\ref{induction_egorov}, which says $\mathfrak{r}_{-2}$ is a real-valued symbol, and the definition of $\mathfrak{r}_{-2,1}$ in \eqref{kk2dsdsd}.
\end{proof}

\begin{lemma}\label{sjhjsdsda2sdijwhat}
$\Phi_4$ is a symplectic, reversibility preserving and  real operator.
\end{lemma}
\begin{proof}
From Lemma~\ref{p2sdsdkcssd}, we have that $\mathfrak{p}_2$ is a real-valued symbol, hence, the PDE in \eqref{hamisdsd} is a Hamiltonian PDE. Therefore its flow map is a symplectic transformation. Thanks to Lemma~\ref{veryexahsysdwds}, we have that the flow map is also reversibility preserving.
 From  \ref{weartnlsdsd} of Proposition~\ref{symobssd} and \ref{dle0isdefinedhere}, we have that $Op^W(\mathfrak{d}_{\ge 0})$ is a real operator. Recalling our choice of $\mathfrak{p}_2$ from \eqref{defosdjsdj2sdsd}, it is straightforward that $Op^W(\mathfrak{p}_2)$ is also a real operator. Therefore,  $\Phi_4$, the flow map determined by the PDE \eqref{hamisdsd}, is a real operator.
\end{proof}

\subsubsection{Analysis of the transformation $\Phi_4$}

\begin{lemma}\label{align2sd2sd2sx2}
$\Phi^{\pm} = \Phi_4,\Phi_4^{-1}$ satisfies 
\begin{equation}\label{transformation_estimate_cl22ass2}
\begin{aligned}
&\rVert (\Phi^{\pm}-I) h \rVert_{s}^{\Lip(\gamma,\Omega_1)} \le_{\mathtt{pe},s,\eta_0}  \epsilon\left( \rVert h \rVert^{\Lip(\gamma,\Omega_1)}_{s+\mu_0} + \rVert \mathfrak{I}_\delta \rVert^{\Lip(\gamma,\Omega_1)}_{s+\mu_0}\rVert h \rVert^{\Lip(\gamma,\Omega_1)}_{s_0+\mu_0}\right)\\
&\rVert d_i\Phi^{\pm}(i_0)h[\ihat] \rVert_{s} \le_{\mathtt{pe},s,\eta_0}  \left(\rVert h\rVert_{s+\mu_0} + \rVert \mathfrak{I}_\delta \rVert_{s+\mu_0}\rVert h \rVert_{s_0+\mu_0}\right) \rVert \ihat \rVert_{s_0+\mu_0} + \rVert h \rVert_{s_0+\mu_0}\rVert \ihat \rVert_{s+\mu_0}.
\end{aligned}
\end{equation}
\end{lemma}
\begin{proof}
The estimates for $\Phi^{\pm1}$ follow from their definitions in \eqref{phi4defsd2wsd} and the estimates for $\mathfrak{p}_2$ in Lemma~\ref{p2sdsdkcssd}.
\end{proof}

\subsubsection{Analysis of the remainder  $R_4$}
\begin{lemma}\label{remisn22sd2ss2}
$R_4\in \mathfrak{R}(i_0)$. 
\end{lemma}
\begin{proof}
Recalling $R_4$ from \eqref{kk2dsdsd}, we see that the result follows from Lemma~\ref{align2sd2sd2sx2}, Lemma~\ref{funitesd2sd} and \ref{rhshdosd13} in Proposition~\ref{induction_egorov}.
\end{proof}

\section{Linear Birkhoff normal form}\label{linearbbss22}

In this section, we aim to remove the terms of size $O(\epsilon)$, $O(\epsilon^2)$ in the linear operator $\mathcal{L}^4$ obtained in \eqref{sdlinesdar2sd2sj}. \index{Linear Birkhoff normal form}

\subsection{Linear Birkhoff normal form: Step 1}\label{firststeplbnf}
Using Lemma~\ref{r2sd2sjrevsd} and \eqref{sdlinesdar2sd2sj}, we have
\begin{align}
\mathcal{L}^4 & = \D_\omega -\Pi_{S^\perp}\partial_xOp^W\left( \mathtt{m}_\alpha m_{1,\alpha}(\xi) + \left(\frac{T_\alpha}{4}+\mathfrak{m}_{\le 0}(\xi)\right) + \epsilon\mathfrak{b}_1 + \epsilon^2\mathfrak{b}_2
+\mathfrak{r}_{-2,\le3} \right) \nonumber \\& + \partial_x\Pi_{S^\perp}W_1 + R_4,\label{sdsfeesd0923}
\end{align}
for some $\mathfrak{b}_1\in \mathfrak{B}^{-2}_1$ and $ \mathfrak{b}_2\in \mathfrak{B}^{-2}_2$. 

In order to remove $\mathfrak{b}_1$, we use a transformation $\Phi_5$ defined by
\begin{align}\label{phi5deff}
\Phi_{5}:=e^{\epsilon A_1},\quad \text{ where $A_1:=\partial_x\Pi_{S^\perp}Op^W({\rho_1})$},
\end{align}
for some symbol $\rho_1=\rho_1(\omega,\varphi,\cdot,\cdot)\in \mathcal{S}^{-1-\alpha}$, which will be chosen later (see Proposition~\ref{rhochosedefine}).  Note that $\Phi_5$ is a time-$1$ flow map of a pseudo PDE:
\begin{align}\label{skdsdpdesx}
\frac{d}{d\tau}u = \partial_x \Pi_{S^\perp}\epsilon Op^W(\rho_1)u.
\end{align}
Since $\partial_x Op^W(\rho_1) \in \mathcal{S}^{-\alpha}$,  $e^{\epsilon A_1}:=\sum_{n=0}^{\infty}\frac{1}{n!}(\epsilon A_1)^n$ is well-defined {(see Lemma~\ref{sdsdsdsd222kdkehsdjs})}.

 As in Section~\ref{study_hard_learn_more}, we will decompose $\Phi_5$ into a finite dimensional operator and $\Psi_5$, which is a time-$1$ flow map of the pseudo-PDE:
 \begin{align}\label{psi5defsseconddate}
 \partial_\tau u = \epsilon\partial_x Op^W(\rho_1)u,\text{ that is, }\frac{d}{d\tau}\Psi_5(\tau) = \epsilon\partial_xOp^W(\rho_1)\Psi_5(\tau),\quad \Psi_5(0)=I.
 \end{align}

 The main result is below:
\begin{proposition}\label{linearstep17}
There exists a symplectic transformation $\Phi_5:H_{S^\perp}\mapsto H_{S^\perp}$ such that
\begin{equation}\label{sdlinesdar2sd2s2j}
\begin{aligned}
\mathcal{L}^5[h]&:=(\Phi_5)^{-1}\mathcal{L}^4\Phi_5[h] = \D_\omega h -  \Pi_{S^\perp}\partial_xM_5[h] +\partial_x\Pi_{S^\perp} W_2 + R_5[h],\\
M_5[h] &:= Op^W\left(\mathtt{m}_\alpha m_{1,\alpha}(\xi) +\left(\frac{T_\alpha}4  + \mathfrak{m}_{\le0}(\omega,\xi)\right) +\epsilon^2\tilde{\mathfrak{b}}_2 + {\mathfrak{r}}_{-2,\le 3,*}\right),
\end{aligned}
\end{equation}
satisfies the following:

\begin{enumerate}[label=(\arabic*)]
\item \label{item22sx692} $\tilde{\mathfrak{b}}_2\in \mathfrak{B}^{-2}_2$ and ${\mathfrak{r}}_{-2,\le3,*}$ satisfies
\begin{equation}\label{sdsdtmddndkQk}
\begin{aligned}
&|\mathfrak{r}_{-2,\le3,*}|^{\Lip(\gamma,\Omega_1)}_{-2,s,\eta_0} \le_{\mathtt{pe},s,\eta_0}  \gamma^{-1}\left(\epsilon^5 + \epsilon^3\rVert \mathfrak{I}_\delta\rVert^{\Lip(\gamma,\Omega_1)}_{s+\mu_0}\right),\\
&|d_i\mathfrak{r}_{-2,\le3,*}(i_0)[\ihat]|_{-2,s,\eta_0} \le_{\mathtt{pe},s,\eta_0}   \epsilon^3\gamma^{-1}\left( \rVert \ihat \rVert_{s+\mu_0} + \rVert \mathfrak{I}_\delta\rVert_{s+\mu_0}\rVert \ihat \rVert_{s_0+\mu_0}\right).
\end{aligned}
\end{equation}
\item \label{item262221} $\tilde{\mathfrak{b}}_2$ is a real-valued and reversible symbol.
\item \label{item262222} $W_{2}:=\Psi_{5}^{T}W_1\Psi_5$, and  $R_5\in \mathfrak{R}(i_0)$.
\item \label{rlsd2ssss}$\Phi_5$ is real, reversibility preserving and symplectic. Therefore, $\mathcal{L}^5$ is  real, reversible Hamiltonian.

\item\label{egor1sd1211x2ppsd} $\mathcal{L}^5$ and $\Phi_5$ are $\frac{2\pi}{\mathtt{M}}$-translation invariance preserving and $\tilde{\mathfrak{b}}_{-2}$ is also  $\frac{2\pi}{\mathtt{M}}$-translation invariant.
\end{enumerate}
\end{proposition}
\begin{proof}
It follows from Lemma~\ref{linearbirss2} and Proposition~\ref{rhochosedefine} that
\[
\Phi_{5}^{-1}\mathcal{L}^4\Phi_5 = \D_\omega -\mathcal{B}_0 + \mathcal{Q}_1 + \mathcal{W}.
\]
Using \eqref{def_bssss} for $\mathcal{B}_0$,  Lemma~\ref{fisd2sd002sd} for $\mathcal{W}$ and Lemma~\ref{ibelieveshedoes} for $\mathcal{Q}_1$, we obtain the expression of $\mathcal{L}^5$ given in \eqref{sdlinesdar2sd2s2j} ($R_5$ is defined to be $R_5:=R_{5,1}+R_{5,2}$ where $R_{5,1},\ R_{5,2}$ are as in Lemma~\ref{fisd2sd002sd} and Lemma~\ref{ibelieveshedoes} respectively). Item \ref{item22sx692} follows from Lemma~\ref{ibelieveshedoes}. From Lemma~\ref{symbsrerer}, with \ref{rhshdosd1422} of Proposition~\ref{rlaqkqdpsanfsdfmf}, we obtain \ref{rlsd2ssss}. Since $\mathcal{L}^5$ is Hamiltonian, $\tilde{\mathfrak{b}}_2$ is real-valued (since $M_5$ is a symmetric operator (see \eqref{adjoint_weyl})). Also since $\mathcal{L}^5$ is reversible, therefore $\tilde{\mathfrak{b}}_2$ is a reversible symbol. This gives a proof for item \ref{item262221}. Item \ref{item262222} is a result from  Lemma~\ref{fisd2sd002sd} and Lemma~\ref{ibelieveshedoes}, which gives us $R_5:=R_{5,1}+R_{5,2}\in \mathfrak{R}(i_0)$. For item \ref{egor1sd1211x2ppsd}, we see from Proposition~\ref{rhochosedefine} that $\rho_1\in \mathfrak{B}_1^{-1-\alpha}$, that is, $x\mapsto \rho_1(\varphi,x,\xi)$ is supported on $\mathtt{M}$-th Fourier modes, since $S$ is a set of $\mathtt{M}$-multiples of integers (see \eqref{tan_site} and Definition~\ref{buildingbobo}). Therefore, the flow map generated by the PDE in \eqref{skdsdpdesx} is $\frac{2\pi}{\mathtt{M}}$-translation invariant.
\end{proof}

From now on, we aim to prove the lemmas that are used to prove the above proposition. In view of \eqref{sdsfeesd0923}, let us denote
\begin{equation}\label{def_bssss}
\begin{aligned}
\mathcal{B}_0&:=\Pi_{S^\perp}\partial_xOp^W\left(\mathtt{m}_\alpha m_{1,\alpha}(\xi) + \left(\frac{T_\alpha}4 + \mathfrak{m}_{\le 0}(\xi)\right)\right),\quad 
\mathcal{B}_1:= \Pi_{S^\perp}\partial_xOp^W(\mathfrak{b}_1),\\
 \mathcal{B}_2&=\Pi_{S^\perp}\partial_xOp^W(\mathfrak{b}_2),\quad \mathcal{Q}:=-\Pi_{S^\perp}\partial_x Op^W(\mathfrak{r}_{-2,\le 3}) + \Pi_{S^\perp}\partial_x W_1+ R_4.
\end{aligned}
\end{equation}
so that $\mathcal{L}^4$ in \eqref{sdsfeesd0923} can be written as
\begin{align}\label{l4_linear_B}
\mathcal{L}^4 = \D_\omega - (\mathcal{B}_0+ \epsilon \mathcal{B}_1 + \epsilon^2 \mathcal{B}_2) + \mathcal{Q}.
\end{align}

In view of \ref{itnsdwd} of Proposition~\ref{prop_time_rep_92} and \eqref{docsd2ksdsdsd}, we further split $\mathcal{B}_0$ in \eqref{def_bssss} as
\begin{align}\label{b01andb02sd}
\mathcal{B}_0& = \Pi_{S^\perp}\partial_xOp^W\left(\mathtt{m}_\alpha m_{1,\alpha}(\xi) + \left(\frac{T_\alpha}4 + \mathfrak{m}_{\le 0}(\xi)\right)\right)\nonumber\\
&= \Pi_{S^\perp}\partial_xOp^W\left(-\frac{1}{2}m_{1,\alpha}(\xi) + \frac{T_\alpha}4\right) \nonumber \\
& + \Pi_{S^\perp}\partial_xOp^W\left(\left(\mathtt{m}_{\alpha}+\frac{1}{2}\right)m_{1,\alpha}(\xi) + \mathfrak{m}_{\le 0}(\xi)\right)\nonumber\\
&=:\mathcal{B}_{0,1}+\mathcal{B}_{0,2}.
\end{align}

\subsubsection{Computing the conjugation}
 We expand $\Phi_5$ as
 \begin{equation}\label{expansionofphi_5}
\begin{aligned}
\Phi_5 &=:I + \epsilon\mathcal{A}_1 =: I + \epsilon A_1 +\epsilon^2\mathcal{A}_2 =: I + \epsilon A_1 + \frac{\epsilon^2 A_1^2}{2} +\epsilon^3 \mathcal{A}_3, \\
\Phi_5^{-1} &=: I + \epsilon\tilde{\mathcal{A}}_1 =:I - \epsilon A_1 + \epsilon^2\tilde{\mathcal{A}}_2 =: I - \epsilon A_1 + \frac{\epsilon^2A_1^2}{2} +\epsilon^3 \tilde{\mathcal{A}}_3.
\end{aligned}
\end{equation}

We first specify the operators of size $O(\epsilon)$:
\begin{lemma}\label{linearbirss2}
We have that
\begin{align}\label{sdsdxckksd1spqs}
\Phi_{5}^{-1}\mathcal{L}^4\Phi_5 = \D_\omega - \mathcal{B}_0 + \epsilon (\D_\omega(A_1)-[\mathcal{B}_{0,1},A_1]-\mathcal{B}_1)+  \mathcal{Q}_1 + \mathcal{W}.
\end{align}
where 
\begin{equation}\label{q1def2252323}
\begin{aligned}
\mathcal{Q}_1 & :=  \epsilon (\D_\omega - \D_{\overline{\omega}})(A_1)  - \epsilon[\mathcal{B}_{0,2},A_1]\\
& \ + \left( \epsilon^2 \D_\omega(\mathcal{A}_2) + (\Phi_5^{-1}-I)\left(\epsilon \D_\omega(A_1) + \epsilon^2\D_\omega(\mathcal{A}_2) \right)\right) \\
& \ + \epsilon^2\left( A_1\mathcal{B}_0A_1  + \tilde{\mathcal{A}_2}(\mathcal{B}_0 + \epsilon \mathcal{B}_0A_1) -\Phi_5^{-1}\mathcal{B}_0\mathcal{A}_2\right) \\
& \ -\left( \epsilon^2\mathcal{B}_1\mathcal{A}_1 + \epsilon^2\tilde{\mathcal{A}}_1\mathcal{B}_1\Phi_5\right) - \epsilon^2\Phi_5^{-1}\mathcal{B}_2 \Phi_5  -\Phi_5^{-1}\Pi_{S^\perp}\partial_x Op^W(\mathfrak{r}_{-2,\le 3})\Phi_5  \\& \ +  \Phi_5^{-1}R_4\Phi_5 \\
\mathcal{W} &:=\Phi_5^{-1}\Pi_{S^\perp}\partial_x W_1 \Phi_5.
\end{aligned}
\end{equation}

\end{lemma}

\begin{proof}
Using \eqref{l4_linear_B}, we have that 
\begin{align*}
\Phi_5^{-1} \mathcal{L}^4\Phi_5 = \Phi_5^{-1}\D_\omega \Phi_5 - \left( \Phi_5^{-1}\mathcal{B}_0\Phi_5 +    \epsilon\Phi_5^{-1}\mathcal{B}_1 \Phi_5 +  \epsilon^2\Phi_5^{-1}\mathcal{B}_2 \Phi_5 \right)+  \Phi_5^{-1}\mathcal{Q} \Phi_5.
\end{align*}

\textbf{Conjugation with $\D_\omega.$}
Using \eqref{expansionofphi_5}, we have 
\begin{align}
& \Phi_{5}^{-1}\D_\omega \Phi_5 \nonumber \\
&= \Phi_{5}^{-1}\left( \D_\omega (\Phi_5) + \Phi_5 \D_\omega\right) \nonumber\\
& = \D_\omega + \Phi_{5}^{-1}\left(\epsilon \D_\omega(A_1) + \epsilon^2\D_\omega(\mathcal{A}_2) \right)\nonumber\\
& = \D_\omega + \epsilon \D_\omega(A_1) + \left( \epsilon^2 \D_\omega(\mathcal{A}_2) + (\Phi_5^{-1}-I)\left(\epsilon \D_\omega(A_1) + \epsilon^2\D_\omega(\mathcal{A}_2) \right)\right)\nonumber\\
&=\D_\omega + \epsilon \D_{\overline{\omega}}(A_1)  + \epsilon (\D_\omega - \D_{\overline{\omega}})(A_1) \nonumber \\
& \qquad + \left( \epsilon^2 \D_\omega(\mathcal{A}_2) + (\Phi_5^{-1}-I)\left(\epsilon \D_\omega(A_1) + \epsilon^2\D_\omega(\mathcal{A}_2) \right)\right).\label{esdjsd2sd}
\end{align}

\textbf{Conjugation with $\mathcal{B}_0.$}
Using \eqref{expansionofphi_5} again, we have that
\begin{align*}
\Phi_5^{-1}\mathcal{B}_0\Phi_5 &= \Phi_{5}^{-1}\mathcal{B}_0(I + \epsilon A_1 + \epsilon^2 \mathcal{A}_2)\\
& = \Phi_5^{-1}\mathcal{B}_0(I + \epsilon A_1) + \epsilon^2\Phi_5^{-1}\mathcal{B}_0\mathcal{A}_2\\
& = \mathcal{B}_0 + \epsilon[\mathcal{B}_0,A_1] - \epsilon^2\left( A_1\mathcal{B}_0A_1  - \tilde{\mathcal{A}_2}(\mathcal{B}_0 + \epsilon \mathcal{B}_0A_1) -\Phi_5^{-1}\mathcal{B}_0\mathcal{A}_2\right).
\end{align*}

\textbf{Conjugation with $\mathcal{B}_1.$}
Similarly, \begin{align*}
\epsilon\Phi_5^{-1} \mathcal{B}_1\Phi_5 = \epsilon \mathcal{B}_1 + \epsilon^2\mathcal{B}_1\mathcal{A}_1 + \epsilon^2\tilde{\mathcal{A}}_1\mathcal{B}_1\Phi_5.
\end{align*}

\textbf{Conjugation with $\mathcal{B}_2.$}
We do not rewrite it since it is already $O(\epsilon^2)$.

\textbf{Conjugation with $\mathcal{Q}$.}
Recalling $\mathcal{Q}$ from \eqref{def_bssss}, we have
\begin{align*}
\Phi_5^{-1}\mathcal{Q}\Phi_5 & = -\Phi_5^{-1}\Pi_{S^\perp}\partial_x Op^W(\mathfrak{r}_{-2,\le 3})\Phi_5 + \Phi_5^{-1}\Pi_{S^\perp}\partial_x W_1 \Phi_5 + \Phi_5^{-1}R_4\Phi_5.
\end{align*}

Summing everything up, we obtain
\begin{align}
\Phi_5^{-1} \mathcal{L}^4\Phi_5  & = \D_\omega -\mathcal{B}_0 + \epsilon\left( \D_{\overline{\omega}}(A_1) - [\mathcal{B}_0, A_1] - \mathcal{B}_1 \right) \nonumber \\ & \qquad +\epsilon (\D_\omega - \D_{\overline{\omega}})(A_1) +  Q_1 + \mathcal{W}, \label{linear_bb_conju}
\end{align}
where $Q_1$ and $\mathcal{W}$ are 
\begin{equation}\label{q1dgef225232x3}
\begin{aligned}
{Q}_1 & := \left( \epsilon^2 \D_\omega(\mathcal{A}_2) + (\Phi_5^{-1}-I)\left(\epsilon \D_\omega(A_1) + \epsilon^2\D_\omega(\mathcal{A}_2) \right)\right) \\
& \ + \epsilon^2\left( A_1\mathcal{B}_0A_1  - \tilde{\mathcal{A}_2}(\mathcal{B}_0 + \epsilon \mathcal{B}_0A_1) -\Phi_5^{-1}\mathcal{B}_0\mathcal{A}_2\right) \\
& \ -\left( \epsilon^2\mathcal{B}_1\mathcal{A}_1 + \epsilon^2\tilde{\mathcal{A}}_1\mathcal{B}_1\Phi_5\right) \\
& \ - \epsilon^2\Phi_5^{-1}\mathcal{B}_2 \Phi_5  -\Phi_5^{-1}\Pi_{S^\perp}\partial_x Op^W(\mathfrak{r}_{-2,\le 3})\Phi_5  +  \Phi_5^{-1}R_4\Phi_5 \\
\mathcal{W} &:=\Phi_5^{-1}\Pi_{S^\perp}\partial_x W_1 \Phi_5.
\end{aligned}
\end{equation}
Recalling $\mathcal{B}_{0,1}$ and $\mathcal{B}_{0,2}$ from \eqref{b01andb02sd}, we have
\[
[\mathcal{B}_0,A_1] = [\mathcal{B}_{0,1},A_1] + [\mathcal{B}_{0,2},A_1].
\]
Plugging this into \eqref{linear_bb_conju} and defining $\mathcal{Q}_1 := Q_1 +  \epsilon (\D_\omega - \D_{\overline{\omega}})(A_1)  - \epsilon[\mathcal{B}_{0,2},A_1]$, we obtain \eqref{sdsdxckksd1spqs} with \eqref{q1def2252323}.
\end{proof}

\subsubsection{Choice of $\rho_1$}
We aim to choose $A_1$ by choosing the symbol $\rho_1$ (see \eqref{phi5deff}), to eliminate the terms of size $O(\epsilon)$ in \eqref{linear_bb_conju}, that is,
\begin{align}\label{homological_eq1}
\D_{\overline{\omega}}(A_1) - [\mathcal{B}_{0,1},A_1] - \mathcal{B}_1 = 0. 
\end{align}

\begin{proposition}\label{rhochosedefine}
There exists $\rho_1=\rho_1(\omega,\varphi,x,\xi)\in \mathfrak{B}^{-1-\alpha}_1$ such that
\begin{enumerate}[label=(\arabic*)]
\item \label{LBNF_item1} With $A_1$ defined as in \eqref{phi5deff}, we have that
\begin{align}\label{sjrlaclWlskdsd}
\D_{\overline{\omega}}(A_1) - [\mathcal{B}_{0,1},A_1] - \mathcal{B}_1 = 0. 
\end{align}
\item \label{LBNF_item3}$\rho_1$ is a real-valued and reversibility preserving symbol. Furthermore, $\partial_xOp^W(\rho_1)$ is a real operator.
\end{enumerate}
\end{proposition}
\begin{proof}

\vspace{0.5\baselineskip}\noindent\textit{Proof of \ref{LBNF_item1}.}  
Using \eqref{phi5deff} and \eqref{b01andb02sd}, it is easy to see that \eqref{homological_eq1} is equivalent to (see \eqref{commutator} for the notation $[\cdot,\cdot]_x$),
\begin{align}\label{homological_eq21s2}
\Pi_{S^\perp}\left( Op^W(\D_{\overline{\omega}} \rho_1) - \left[Op^W\left(-\frac{1}{2}m_{1,\alpha}(\xi) + \frac{T_\alpha}4\right),Op^W(\rho_1)\right]_x- Op^W(\mathfrak{b}_1)\right) = 0.
\end{align}
Recalling  $\mathfrak{b}_1\in \mathfrak{B}^{-2}_{1}$, which follows from Lemma~\ref{r2sd2sjrevsd}, and the definition of $\mathfrak{B}^{-2}_1$ from Definition~\ref{buildingbobo}, we can write $\mathfrak{b}_1$ as
\begin{align}\label{rlaxcxcc1}
\mathfrak{b}_1(\omega,\varphi,x,\xi)  = \sum_{j_k\in S} C_{\mathfrak{b},j_k}(\xi)\sqrt{j_k\zeta_k}e^{\ii (\mathtt{l}(j_k)\cdot \varphi + j_k x)},\text{ for some $C_{\mathfrak{b},j_k}\in \mathcal{S}^{-2}$.}
\end{align}
Using \eqref{usual_weyl}, let us denote by  $\mathfrak{b}_{1,s}$, a symbol such that $Op^W(\mathfrak{b}_1) = Op(\mathfrak{b}_{1,s})$, that is,
\begin{align}\label{bstandsxcsd}
\widehat{\mathfrak{b}}_{1,s}(\omega,\varphi,x,\xi):= \sum_{j_k\in S}\underbrace{C_{\mathfrak{b},j_k}(\xi+\frac{j_k}2)}_{=:C_{\mathfrak{b}_{1,s},j_k}(\xi)}\sqrt{j_k\zeta_k}e^{\ii (\mathtt{l}(j_k)\cdot \varphi + j_k x)}.
\end{align}
Since $S$ is a finite set (see \eqref{tan_site}), it is clear that $C_{\mathfrak{b}_{1,s},j_k}\in \mathcal{S}^{-2}$, therefore, $\mathfrak{b}_{1,s}\in \mathfrak{B}^{-2}_1$. 
We look for a symbol $\rho_{1,s}\in \mathfrak{B}^{-1-\alpha}_1$ such that
\begin{align}\label{stasdxcgsd}
\Pi_{S^\perp}\left(Op(\D_{\overline{\omega}}\rho_{1,s}) -\left[Op\left(-\frac{1}{2}m_{1,\alpha}(\xi) + \frac{T_\alpha}4\right),Op(\rho_{1,s})\right]_x-Op(\mathfrak{b}_{1,s})\right) = 0.
\end{align}
Once we find $\rho_{1,s}$, then we can define $\widehat{\rho}_{1}(\omega,\varphi,k,\xi):=\widehat{\rho}_{1,s}(\omega,\varphi,k,\xi - k)$ for $k\in S$, so that \eqref{usual_weyl} tells us  $Op^W(\rho_1)=Op(\rho_{1,s})$, therefore, $\rho_1$ is a solution to \eqref{homological_eq21s2}. 

 To find $\rho_{1,s}$ solving \eqref{stasdxcgsd}, we use the notations in \eqref{fourier_ssd2s} and \eqref{lambdadeffsd} to see that \eqref{stasdxcgsd} is equivalent to
 \begin{align}\nonumber
 \ii \left(\overline{\omega}\cdot l - \left(\left(-\frac{1}{2}\lambda_\alpha(k+j)-  (k+j) \frac{T_\alpha}4\right) - \left( -\frac{1}{2}\lambda_{\alpha}(j)-  j\frac{T_\alpha}4\right)\right)  \right)&\widehat{\rho}_{1,s}^{\varphi,x}(\omega,l,k,j)\\
 & =\widehat{\mathfrak{b}}^{\varphi,x}_{1,s}(\omega,l,k,j),\label{skjxc3ssdxcx}
 \end{align}
 for all $j\in S^\perp$. We set
\begin{align}\label{jjsd1xcxcxswihscsd}
{\rho}_{1,s}(\omega,\varphi,x,\xi):=\sum_{j_k\in S}2\ii C_{\mathfrak{b}_{1,s},j_k}(\xi)\kappa_{-j_k,1-\alpha}(\xi)\sqrt{j_k\zeta_k}e^{\ii (\mathtt{l}(j_k)\cdot \varphi + j_k x)},
\end{align}
where $\xi\mapsto \kappa_{-j_k,1-\alpha}(\xi)$ is the multiplier obtained in Lemma~\ref{kaoxixcc2xc}. Since $C_{\mathfrak{b}_{1,s,j_k}}\in \mathcal{S}^{-2}$ and $\kappa_{-j_k,1-\alpha}\in \mathcal{S}^{1-\alpha}$, we have that $\rho_{1,s}\in \mathfrak{B}^{-1-\alpha}_1$. Now, it suffices to show that $\rho_{1,s}$ solves \eqref{skjxc3ssdxcx}. Since $\rho_{1,s}$ and $\mathfrak{b}_{1,s}$ are supported on a finite number of Fourier modes as seen in \eqref{jjsd1xcxcxswihscsd}, and \eqref{rlaxcxcc1}, we have that  both sides in \eqref{skjxc3ssdxcx} are equal to $0$, if $(l,k)\ne (\mathtt{l}(j_k),j_k)$ for some $j_k\in S$. On the other hand, if $(l,k)= (\mathtt{l}(j_k),j_k)$ for some $j_k\in S$, then $\overline{\omega}$ in \eqref{linear_frequency_aa} tells us that (also see $\lambda_\alpha$ in \eqref{lambdadeffsd})
\begin{align}\label{heworkstoomuch}
\overline{\omega}\cdot l& - \left(\left(-\frac{1}{2}\lambda_\alpha(k+j)-  (k+j) \frac{T_\alpha}4\right) - \left( -\frac{1}{2}\lambda_{\alpha}(j)-  j\frac{T_\alpha}4\right)\right)\nonumber\\
& =\frac{1}{2}\left(\lambda_{\alpha}(j+j_k) - \lambda_\alpha(j) - \lambda_{\alpha}(j_k) \right),
\end{align}
hence, the left-hand side of \eqref{skjxc3ssdxcx} is given by 
\begin{align*}
LHS & = \frac{\ii}{2}\left(\lambda_{\alpha}(j+j_k) - \lambda_\alpha(j) - \lambda_{\alpha}(j_k) \right)\times 2\ii C_{\mathfrak{b}_{1,s}j_k}(j)\kappa_{-j_{k},1-\alpha}(j) \\
& = C_{\mathfrak{b}_{1,s},j_k}(j)=\widehat{\mathfrak{b}}_{1,s}^{\varphi,x}(\omega,\mathtt{l}(j_k),j_k,j),
\end{align*}
for all $j\in S^\perp$, where the second equality follows from \eqref{rjjsdsdjxcxcd}. Therefore, \eqref{jjsd1xcxcxswihscsd} solves \eqref{skjxc3ssdxcx} and thus, \eqref{stasdxcgsd}. 

\vspace{0.5\baselineskip}\noindent\textit{Proof of \ref{LBNF_item3}.} 
We will prove that $\rho_1$ can be chosen to be a real-valued symbol only. This property follows from the fact that $\rho_1$ is a solution to \eqref{homological_eq21s2}, while, $Op^W\left(-\frac{1}{2}m_{1,\alpha}(\xi) + \frac{T_\alpha}4\right)$ and $Op^W(\mathfrak{b}_1)$ are symmetric operators, which follows from \eqref{transfposdsdsd} with the fact that the symbols $-\frac{1}{2}m_{1,\alpha}(\xi) + \frac{T_\alpha}4$ and $\mathfrak{b}_1$ are real-valued, which follows from Proposition~\ref{nichts2} for $m_{1,\alpha}(\xi)$ and Lemma~\ref{sjrhcnskwdsd} with \eqref{rlsd1b1b2} for $\mathfrak{b}_1$). Indeed, Proposition~\ref{nichts2} and  Lemma~\ref{sjrhcnskwdsd} with \eqref{rlsd1b1b2}  tell us that $Op^W\left(-\frac{1}{2}m_{1,\alpha}(\xi) + \frac{T_\alpha}4\right)$ and $Op^W(\mathfrak{b}_1)$ are real and reversibility preserving operators, the other desired properties of $\rho_1$ follow in the same way.

 To show that $\rho_1$ can be chosen to be a real-valued symbol, we see that $\rho_1$ is a solution to \eqref{homological_eq21s2}, therefore, \eqref{transfposdsdsd} tells us that (using that $Op^W\left(-\frac{1}{2}m_{1,\alpha}(\xi) + \frac{T_\alpha}4\right)$ and $Op^W(\mathfrak{b}_1)$ are symmetric operators), $Op^W(\overline{\rho_1}$) solve \eqref{homological_eq21s2} as well. Since the homological equation \eqref{homological_eq21s2} is linear, we see that $\frac{1}{2}(\rho_1 + \overline{\rho_1})$ solves \eqref{homological_eq21s2}. Hence, if necesary, we can replace $\rho_1$ by $\frac{1}{2}(\rho_1 + \overline{\rho_1})$, which gives a real-valued symbol. Clearly, if $\rho_1\in \mathfrak{B}^{-1-\alpha}_1$, then $\frac{1}{2}(\rho_1 + \overline{\rho_1})\in \mathfrak{B}^{-1-\alpha}_1$ as well, which follows immediately from its definition in Definition~\ref{buildingbobo}.
 
 In view of \eqref{llsd1sdxcxcxcmd19sd2}, we can replace $\rho_1$ by $\frac{1}{2}(\rho_1(\varphi,x,\xi)+\overline{\rho_1}(\varphi,x,-\xi))$ and after than, again replace it by $\frac{1}{2}(\rho_{1}(\varphi,x,\xi) - \rho_1(-\varphi,-x,-\xi))$, which eventually gives us a symbol $\rho_1$ with all the desired properties.
 \end{proof}

 \subsubsection{Analysis of $\Phi_5$}
 
  \begin{lemma}\label{sdsdsdsd222kdkehsdjs}
For $k=1,2,3$, there exist $\rho_{k,*}(\tau)\in \mathcal{S}^{-k\alpha}$ such that  
\[
\Psi_5(\tau) - I = \epsilon Op(\rho_{1,*}(\tau)) + \frac{1}{2}\epsilon^2 Op(\rho_{2,*}(\tau)) + \epsilon^3Op(\rho_{3,*}(\tau)),\] and $\rho_{k,*}$ satisfy
\begin{align}
\sup_{\tau\in [0,1]}|\rho_{k,*}(\tau)|^{\Lip(\gamma,\Omega_1)}_{-k\alpha,s,\eta_0}&\le_{s,\mathtt{pe},\eta_0} 1,\label{rmsidgkse}\\
\sup_{\tau\in [0,1]}|d_i\rho_{k,*}(\tau)(i_0)[\ihat]|_{-k\alpha,s,\eta_0}& = 0. \label{rmsidgkse2}
\end{align}
 \end{lemma}
 \begin{proof}
 Let us pick $\rho_{1,1}\in \mathcal{S}^{-\alpha}$ such that
 \begin{align}\label{rh2s11}
 Op(\rho_{1,1}) = \partial_x Op^W(\rho_1),
 \end{align}
 so that \eqref{psi5defsseconddate} gives us that
 \begin{align}
 \Psi_5(\tau)  - I  & = \epsilon\tau Op(\rho_{1,1}) + \frac{(\epsilon\tau)^2}{2}Op(\rho_{1,1})\circ Op(\rho_{1,1}) \nonumber \\
 & +  (\epsilon\tau)^3\sum_{n=3}^{\infty}(\epsilon\tau)^{n-3} \frac{1}{n!}\left(Op(\rho_{1,1})\right)^n.
 \label{sj209393}
 \end{align}
From \eqref{usual_weyl} and Proposition~\ref{rhochosedefine}, which tells us that $\rho_{1}\in \mathfrak{B}^{-1-\alpha}_1$, hence $\rho_{1,1}\in \mathfrak{B}^{-\alpha}_1$. This gives us that 
 \begin{align}\label{rksh2sdsd}
 |\rho_{1,1}|^{\Lip(\gamma,\Omega_1)}_{-\alpha,s,\eta_0}\le C_1(\mathtt{pe},s,\eta_0),\quad |d_i\rho_{1,1}(i_0)[\ihat]|_{-\alpha,s,\eta_0}= 0,
 \end{align}
 for some constant $C_1(\mathtt{pe},s,\eta_0)>0$, where the estimate for $d_i\rho_{1,1}$ follows from the fact that $\rho_{1,1}$ does not depend on the embedding $i_0$, since $\rho_{1,1}\in \mathfrak{B}^{-\alpha}_1$ (see Definition~\ref{buildingbobo}). Also, this property tells us that
  \[
 \widehat{\rho}_{1,1}^{\varphi,x}(l,k,\xi) = 0, \text{ if $|l| + |k| \ge N_{S}$,}
 \]
 for some ${N}_{S}>0$, depending on the choice of the tangential sites in \eqref{tan_site}.
 Note that choosing $\rho_{1,*}$ and $\rho_{2,*}$ so that
 \begin{align}\label{rlajc2xx}
 \rho_{1,*}(\tau):=\tau\rho_{1,1},\quad Op(\rho_{2,*}(\tau)) = \tau^2Op(\rho_{1,1})\circ Op(\rho_{1,1}),
 \end{align}
 \eqref{rksh2sdsd} gives us \eqref{rmsidgkse} and \eqref{rmsidgkse2} for $k=1,2$, using \eqref{compandkskd2sd} and Lemma~\ref{sdtosd}.
 
 We focus on $\rho_{3,*}(\tau)$ from now on. For $k=3$,  we will prove \eqref{rmsidgkse} only, since \eqref{rmsidgkse2} is trivial ($\Phi_5$ is independent of the embedding, since so is $\rho_{1}$.)
 
 Let us define a sequence of symbols $\rho_{1,n+1}(\tau)$ by
 \[
 Op(\rho_{1,n+1}(\tau)):= \tau Op(\rho_{1,n}(\tau))\circ Op(\rho_{1,1}),\text{ for $n\ge 1$},
 \]
 so that we have
 \begin{align}\label{ksdpslieesd}
  (\epsilon\tau)^3\sum_{n=3}^{\infty}(\epsilon\tau)^{n-3} \frac{1}{n!}\left(Op(\rho_{1,1})\right)^n = \epsilon^3 \sum_{n=3}^{\infty} \epsilon^{n-3}\frac{1}{n!}Op(\rho_{1,n}(\tau)).
 \end{align}

  Note that using Lemma~\ref{finitemodes_foudsd}, there exists $C_2(\mathtt{pe},s,\eta_0)>1$ such that
 \begin{align}
  \nonumber &|\rho_{1,n+1}(\tau)|^{\Lip(\gamma,\Omega_1)}_{-(n+1)\alpha,s,\eta_0}\\
  &\le C_2(\mathtt{pe},s,\eta_0)\left( |{\rho_{1,n}}(\tau)|^{\Lip(\gamma,\Omega_1)}_{-n\alpha,s_0,\eta_0}|\rho_{1,1}|^{\Lip(\gamma,\Omega_1)}_{-\alpha,s,\eta_0} + |{\rho_{1,n}}(\tau)|^{\Lip(\gamma,\Omega_1)}_{-n\alpha,s,\eta_0}|\rho_{1,1}|^{\Lip(\gamma,\Omega_1)}_{-\alpha,s_0,\eta_0}\right)\label{fintsdsd2sdx}
 \end{align}
 Without loss of generality, we can assume that 
 \[
 s\mapsto C_1(\mathtt{pe},s,\eta_0), C_2(\mathtt{pe},s,\eta_0)\]
 {are increasing.}
 
 We claim that 
 \begin{align}\label{claim2sd2wlruqs}
 \begin{cases}
 |\rho_{1,n}(\tau)|^{\Lip(\gamma,\Omega_1)}_{-n\alpha,s_0,\eta_0}\le \left(2C_1(\mathtt{pe},s_0,\eta_0)C_2(\mathtt{pe},s_0,\eta_0)\right)^{n},\\
 |\rho_{1,n}(\tau)|^{\Lip(\gamma,\Omega_1)}_{-n,\alpha,s,\eta_0}\le \left(4C_1(\mathtt{pe},s,\eta_0)C_2(\mathtt{pe},s,\eta_0)\right)^{n} ,\text{ for $s\ge s_0$.}
 \end{cases}
 \end{align}
 We prove $s=s_0$ first. For $n=1$, the desired estimate follows immediately from \eqref{rksh2sdsd}. Assuming the estimate holds for $n$, let us prove it for $n+1$. Using \eqref{fintsdsd2sdx}, we have
 \begin{align*}
 |\rho_{1,n+1}(\tau)|^{\Lip(\gamma,\Omega_1)}_{-(n+1)\alpha,s_0,\eta_0}&\le 2C_2(\mathtt{pe},s_0,\eta_0)|{\rho_{1,n}}(\tau)|^{\Lip(\gamma,\Omega_1)}_{-n\alpha,s_0,\eta_0}|\rho_{1,1}|_{-\alpha,s_0,\eta_0}\\
 &  \le 2C_1(\mathtt{pe},s_0,\eta_0)C_2(\mathtt{pe},s_0,\eta_0)|\rho_{1,n}(\tau)|^{\Lip(\gamma,\Omega_1)}_{-n\alpha,s_0,\eta_0}\\
 &\le \left(2C_1(\mathtt{pe},s_0,\eta_0)C_2(\mathtt{pe},s_0,\eta_0)\right)^{n+1},
 \end{align*}
 where the second inequality follows from \eqref{rksh2sdsd} and the last inequality follows from the induction hypothesis. Therefore, we prove the estimate \eqref{claim2sd2wlruqs} for $s=s_0$.
 
 Now, we prove \eqref{claim2sd2wlruqs} for $s\ge s_0$. Again the estimate for $n=1$ holds, thanks to \eqref{rksh2sdsd}. For $n+1$, we  again use \eqref{fintsdsd2sdx} and \eqref{claim2sd2wlruqs} for $s=s_0$ and obtain
 \begin{align*}
 |\rho_{1,n+1}(\tau)|^{\Lip(\gamma,\Omega_1)}_{-(n+1)\alpha,s,\eta_0}&\le \left(2C_1(\mathtt{pe},s,\eta_0)C_2(\mathtt{pe},s,\eta_0)\right)^{n+1} \\
 & \qquad + C_1(\mathtt{pe},s,\eta_0)C_2(\mathtt{pe},s,\eta_0)|{\rho_{1,n}}(\tau)|^{\Lip(\gamma,\Omega_1)}_{-n\alpha,s,\eta_0}\\
 & \le (2^{n+1}+4^{n})\left(C_1(\mathtt{pe},s,\eta_0)C_2(\mathtt{pe},s,\eta_0)\right)^{n+1}\\
 & \le \left(4C_1(\mathtt{pe},s,\eta_0)C_2(\mathtt{pe},s,\eta_0)\right)^{n+1},
 \end{align*}
 where the second inequality follows from the induction hypothesis. Therefore, we have proved \eqref{claim2sd2wlruqs} holds for $n+1$.

Hence,  using Stirling's formula it follows that there exist a constant $C_3(\mathtt{pe},s,\eta_0)$ and a symbol $\rho_{3,*}(\tau)\in \mathcal{S}^{-3\alpha}$ such that 
\[
\sum_{n=1}^{\infty}\epsilon^{n-3}\frac{1}{n!}Op(\rho_{1,n}(\tau)) = Op(\rho_{3,*}(\tau)),\quad |\rho_{3,*}(\tau)|^{\Lip(\gamma,\Omega_1)}_{-\alpha,s,\eta_0}\le_{\mathtt{pe},s,\eta_0} 1,
\]
which is the desired result \eqref{rmsidgkse} for $k=3$, in view of \eqref{ksdpslieesd} and \eqref{sj209393}.
  \end{proof}

\begin{lemma}\label{fjsdinfsd}
$\Phi_5 \Pi_{S^\perp} - \Pi_{S^\perp}\Psi_5 \Pi_{S^\perp}\in \mathfrak{R}(i_0)$. 
\end{lemma}
\begin{proof}
The result follows from  Lemma~\ref{finitesd2sd}, the construction of $\Phi_5,\Psi_5$ in \eqref{skdsdpdesx}  and \eqref{psi5defsseconddate} and  $\rho_{1}\in \mathfrak{B}^{-1-\alpha}_1$, which follows from  Proposition~\ref{rhochosedefine}.
\end{proof}

\begin{lemma}\label{align2sd2sd2sssrx2}
$\Phi^{\pm} = \Phi_5,\Phi_5^{-1}$ satisfies 
\begin{equation}\label{transformation_es2stimatess_class2}
\begin{aligned}
&\rVert (\Phi^{\pm}-I) h \rVert_{s}^{\Lip(\gamma,\Omega_1)} \le_{\mathtt{pe},s} \epsilon \rVert h \rVert^{\Lip(\gamma,\Omega_1)}_{s}\\
&\rVert d_i\Phi^{\pm}(i_0)h[\ihat] \rVert_{s} =0.\end{aligned}
\end{equation}
\end{lemma}
\begin{proof}
It follows straightforwardly from Lemma~\ref{sdsdsdsd222kdkehsdjs} that $\Psi_5,\Psi_5^{-1}$ satisfy the estimates in \eqref{transformation_es2stimatess_class2}. Thanks to Lemma~\ref{fjsdinfsd}, $\Phi_5^{\pm}$ satisfies \eqref{transformation_es2stimatess_class2} for all $h$ such that $h(\varphi,\cdot)\in H_{S^\perp}$. Also, using that $\rho_1$ is independent of the embedding $i_0$, we  have that $d_i\Phi^{\pm}(i_0)=0$.
\end{proof}

 \begin{lemma}\label{symbsrerer}
 $\Phi_5$ is a  symplectic, reversibility preserving and real operator. 
 \end{lemma}
 \begin{proof}
 Using that $\rho_1$ is real-valued, reversibility preserving and $Op^W(\rho_1)$ is a real operator, which follows from Proposition~\ref{rhochosedefine}, the same proof as in Lemma~\ref{sjhjsdsda2sdijwhat} gives the desired result.
 \end{proof}
 
 \subsubsection{Analysis of  $\mathcal{W}$} 
 \begin{lemma}\label{fisd2sd002sd}
 $\mathcal{W}$ can be written as $\mathcal{W} =\Pi_{S^\perp}\partial_x\Psi_{5}^{T} W_1 \Psi_5 + R_{5,1}$ for some finite dimensional operator $R_{5,1}$ of the form in \eqref{j2j2j2j2j2jssds}.
 \end{lemma}
 \begin{proof}
 We write $\mathcal{W}$ in \eqref{q1def2252323} as
 \begin{align}\label{owsdpaosjd2}
 \mathcal{W}&= \Phi_5^{-1} \Pi_{S^\perp}\partial_x W_1\Phi_5 = \Phi_5^{-1} \Pi_{S^\perp}\partial_x W_1\Pi_{S^\perp}\Psi_5 + \Phi_5^{-1} \Pi_{S^\perp}\partial_x W_1 \left(\Phi_5 - \Pi_{S^\perp}\Psi_5\right)\nonumber\\
 & = \Phi_5^{-1} \Pi_{S^\perp}\partial_x W_1\Psi_5 - \Phi_5^{-1} \Pi_{S^\perp}\partial_x W_1\Pi_S\Psi_5 + \Phi_5^{-1} \Pi_{S^\perp}\partial_x W_1 \left(\Phi_5 - \Pi_{S^\perp}\Psi_5\right)\nonumber\\
 & = \Pi_{S^\perp}\Psi_{5}^{-1}\Pi_{S^\perp}\partial_x W_1 \Psi_5 + (\Phi_5^{-1}\Pi_{S^\perp} - \Pi_{S^\perp}\Psi_5^{-1} \Pi_{S^\perp})\partial_x W_1 \Psi_5\nonumber\\
 & \ - \Phi_5 \Pi_{S^\perp}\partial_x W_1\Pi_S\Psi_5 + \Phi_5 \Pi_{S^\perp}\partial_x W_1 \left(\Phi_5 - \Pi_{S^\perp}\Psi_5\right)\nonumber\\
 & = \Pi_{S^\perp}\Psi_5^{-1} \partial_x W_1 \Psi_5  - \Pi_{S^\perp}\Psi_5^{-1} \Pi_{S}\partial_x W_1 \Psi_5+ (\Phi_5\Pi_{S^\perp} - \Pi_{S^\perp}\Psi_5^{-1} \Pi_{S^\perp})\partial_x W_1 \Psi_5\nonumber\\
 & \ - \Phi_5^{-1} \Pi_{S^\perp}\partial_x W_1\Pi_S\Psi_5 + \Phi_5^{-1} \Pi_{S^\perp}\partial_x W_1 \left(\Phi_5 - \Pi_{S^\perp}\Psi_5\right)\nonumber\\
 & = \Pi_{S^\perp}\Psi_{5}^{-1}\partial_x W_1 \Psi_5 + R_{5,1},
 \end{align}
where
\begin{equation}\label{r51dsdsds1sd}
 \begin{aligned}
 R_{5,1} &:= -\Pi_{S}\Psi_5^{-1}\partial_xW\Psi_5- \Pi_{S^\perp}\Psi_5^{-1} \Pi_{S}\partial_x W_1 \Psi_5\\
 & \ + (\Phi_5^{-1}\Pi_{S^\perp} - \Pi_{S^\perp}\Psi_5^{-1} \Pi_{S^\perp})\partial_x W_1 \Psi_5\\
 & \ - \Phi_5^{-1} \Pi_{S^\perp}\partial_x W_1\Pi_S\Psi_5 + \Phi_5^{-1} \Pi_{S^\perp}\partial_x W_1 \left(\Phi_5 - \Pi_{S^\perp}\Psi_5\right).
 \end{aligned}
 \end{equation}
 It follows from \eqref{symplectic73} and \eqref{owsdpaosjd2} that $\mathcal{W} = \Pi_{S^\perp}\partial_x\Psi_{5}^{T} W_1 \Psi_5 + R_{5,1}$.
 
 Now, we show that $R_{5,1}$ is a finite dimensional operator of the form \eqref{j2j2j2j2j2jssds}.
  Thanks to Lemma~\ref{fjsdinfsd}, we can see that each term in $R_{5,1}$ has a finite dimensional operator of the form $\Phi_5 \Pi_{S^\perp} - \Pi_{S^\perp}\Psi_5 \Pi_{S^\perp}$ or $\Pi_S$. Among all the terms in $R_{5,1}$, we will prove the lemma only for the third term in \eqref{r51dsdsds1sd}, that is,  $(\Phi_5^{-1}\Pi_{S^\perp} - \Pi_{S^\perp}\Psi_5^{-1} \Pi_{S^\perp})\partial_x W_1 \Psi_5$, since the other terms can be treated in a similar way.
  
    Note that  $\Phi_5\Pi_{S^\perp} - \Pi_{S^\perp}\Psi_5 \Pi_{S^\perp} \in \mathfrak{R}(i_0)$, which follows from Lemma~\ref{fjsdinfsd}. Since the equations \eqref{skdsdpdesx} and \eqref{psi5defsseconddate} that determine $\Phi_5$ and $\Psi_5$ are autonomous, we see that $\Phi_5^{-1}$ and $\Psi_5^{-1}$ are the time-$1$ flow maps of their reversed PDEs, therefore $\Phi_5^{-1}\Pi_{S^\perp} - \Pi_{S^\perp}\Psi_5^{-1} \Pi_{S^\perp} \in \mathfrak{R}(i_0)$.  Recalling the definition of $\mathfrak{R}(i_0)$ from  Definition~\ref{remiander_class12}, we can decompose it into
    \[
    \Phi_5^{-1}\Pi_{S^\perp} - \Pi_{S^\perp}\Psi_5^{-1} \Pi_{S^\perp} = \epsilon^2 R_{5,1,1} + R_{5,1,2},
    \]
    where  $R_{5,1,2}$ is already of the form \eqref{j2j2j2j2j2jssds} and there exist $\mathfrak{a}_{5,1},\mathfrak{b}_{5,1}\in \mathfrak{B}^m_{1}$ such that
    \begin{align*}
    \epsilon^2R_{5,1,1}[h] &= \Pi_{S^\perp} \epsilon^2Op^W(\mathfrak{a}_{5,1})\Pi_S\left( Op^W(\mathfrak{b}_{5,1}h\right)\\
    & = \sum_{j\in S}\Pi_{S^\perp} \epsilon^2Op^W(\mathfrak{a}_{5,1})\left(\frac{1}{2\pi}(h,\left(Op^W(\mathfrak{b}_{5,1})\right)^{T}[e^{\ii j x}])_{L^2}e^{\ii j x}\right)\\
    & = \sum_{j\in S}(h,g_j)_{L^2}\chi_j,
    \end{align*}
    where
    \[
    g_j:=Op^W(\mathfrak{b}_{5,1})^{T}[e^{\ii j x}],\quad \chi_j :=  \epsilon^2Op^W(\mathfrak{a}_{5,1})[e^{\ii j x}].
    \]
    Clearly, using $\mathfrak{a}_{5,1},\mathfrak{b}_{5,1}\in \mathfrak{B}^m_1$ for some $m\ge 0$, we have that 
    \begin{align}\nonumber
   \rVert g_j \rVert^{\Lip(\gamma,\Omega_1)}_s\rVert \chi_j\rVert^{\Lip(\gamma,\Omega_1)}_{s_0}+ \rVert g \rVert^{\Lip(\gamma,\Omega_1)}_{s_0}\rVert \chi\rVert^{\Lip(\gamma,\Omega_1)}_{s} & \le\epsilon^2, \\
      \rVert d_ig_j(i_0)[\ihat] \rVert_{s} = \rVert d_i\chi_j(i_0)[\ihat] \rVert_{s} &= 0.\label{gchiest2sdsdzxc}
    \end{align}
   Hence, we get
   \[
   (\Phi_5^{-1}\Pi_{S^\perp} - \Pi_{S^\perp}\Psi_5^{-1} \Pi_{S^\perp})\partial_x W_1 \Psi_5 = \sum_{j\in S}(h,\underbrace{(\partial_x W_1\Psi_5)^T[g_j])}_{=:\tilde{g}_j}\chi_j.
   \] 
  Then, using Lemma~\ref{sdsdsdsd222kdkehsdjs}, recalling  $W_1$ from \eqref{kk2dsdsd} and \eqref{sdkksdsds}, it follows from the estimates for $\Psi_3,\Phi_4$ in Lemma~\ref{align2sd2sd2}, Lemma~\ref{align2sd2sd2sx2}, and the estimates for $\mathfrak{q}_{\omega,\mathtt{N}_\alpha\star \mathfrak{a}},\mathfrak{q}_{M_2,\mathtt{N}_\alpha}\star \mathfrak{a}$ in Lemma~\ref{sjjdsdsymbolsd} that $\tilde{g}_j$ and $\chi_j$ satisfy the estimates in \eqref{sjdjsjsj}. This finishes the proof.
 \end{proof}

 \subsubsection{Analysis of $\mathcal{Q}_1$} Now, we analyze the structure of the operator $\mathcal{Q}_1$ in \eqref{q1def2252323}.
   \begin{lemma}\label{ibelieveshedoes}
 $\mathcal{Q}_1$ is of the form:
 \begin{align}\label{formsdsjsjsd}
 \mathcal{Q}_1 = \Pi_{S^\perp}Op(\tilde{\mathfrak{b}}_2) + \Pi_{S^\perp}Op(\mathfrak{r}_{-2,\le 3,*}) + {R}_{5,2}
 \end{align}
 for some $\tilde{\mathfrak{b}}_2,\mathfrak{r}_{-2,\le 3,*}\in \mathcal{S}^{-2}$ and a finite dimensional operator $R_{5,2}$ such that 
 \begin{enumerate}[label=(\arabic*)]
\item \label{vkdhffhlql} $\tilde{\mathfrak{b}}_2\in \mathfrak{B}^{-2}_{2}$.
\item \label{vkdhffhlql2} $\mathfrak{r}_{-2,\le 3,*}$ satisfies
 \begin{align}
 &|\mathfrak{r}_{-2,\le3,*}|^{\Lip(\gamma,\Omega_1)}_{\alpha-3,s,\eta_0} \le_{\mathtt{pe},s,\eta_0}  \gamma^{-1}\left(\epsilon^5 + \epsilon^3\rVert \mathfrak{I}_\delta\rVert^{\Lip(\gamma,\Omega_1)}_{s+\mu_0}\right),\label{estiatefors20sd}\\
&|d_i\mathfrak{r}_{-2,\le3,*}(i_0)[\ihat]|_{\alpha-3,s,\eta_0} \le_{\mathtt{pe},s,\eta_0} \epsilon^3\gamma^{-1}\left( \rVert \ihat \rVert_{s+\mu_0} + \rVert \mathfrak{I}_\delta\rVert_{s+\mu_0}\rVert \ihat \rVert_{s_0+\mu_0}\right).\label{estiatefors20sd1}
 \end{align}
 \item \label{sd2sgfsd}$R_{5,2}\in \mathfrak{R}(i_0)$.
\end{enumerate}
\end{lemma}

\begin{proof}
Among all the terms of $\mathcal{Q}_1$ in \eqref{q1def2252323}, we prove the lemma for $\epsilon^2\D_\omega(\mathcal{A}_2)$ only, since all the other terms can be treated in the same way.

Clearly, it suffices to show that  $\epsilon^2\mathcal{A}_2$  has the decomposition in \eqref{formsdsjsjsd}, since all the properties in \ref{vkdhffhlql}, \ref{vkdhffhlql2} and \ref{sd2sgfsd} still hold after taking $\D_\omega$.
Recalling $\mathcal{A}_2$ from \eqref{expansionofphi_5} and $\Phi_5$ from \eqref{phi5deff}, we have that
\begin{align*}
\epsilon^2\mathcal{A}_2\Pi_{S^\perp} &= \int_0^1 (\epsilon A_1)^2\Phi_5(\tau)(1-\tau)d\tau \Pi_{S^\perp} \\
& = (\epsilon A_1)^2 \int_0^1 \Psi_5(\tau)(1-\tau)d\tau\Pi_{S^\perp} + (\epsilon A_1)^2\int_0^{1}\Pi_{S}\Psi_5(\tau)(1-\tau)d\tau\Pi_{S^\perp} \\
& \ + (\epsilon^2 A_1)^2 \int_0^1 (\Phi_5(\tau)\Pi_{S^\perp} - \Pi_{S^\perp}\Psi_5(\tau)\Pi_{S^\perp})(1-\tau)d\tau.
\end{align*}
From Lemma~\ref{sdsdsdsd222kdkehsdjs} and Lemma~\ref{fjsdinfsd}, we see that the last two operators are finite dimensional operators in class $\mathfrak{R}(i_0)$ (Note that Lemma~\ref{fjsdinfsd} for the time-$1$ flow, but we can always reparametrize the time variable $\tau$ so that $\Phi_5(\tau)\Pi_{S^\perp} - \Pi_{S^\perp}\Psi_5(\tau)\Pi_{S^\perp}\in \mathfrak{R}(i_0)$). Therefore,
 \begin{align}\label{jjsd2sds}
 \epsilon^2\mathcal{A}_2\Pi_{S^\perp} = (\epsilon A_1)^2 \int_0^1 \Psi_5(\tau)(1-\tau)d\tau\Pi_{S^\perp},
 \end{align}
 up to a finite dimensional operator in $\mathfrak{R}(i_0)$. Similarly, using  Lemma~\ref{sdsdsdsd222kdkehsdjs}, we can find $\rho_*(\tau)\in \mathcal{S}^{-\alpha}$ such that
 \begin{align*}
 \Psi_5(\tau)  =  I + \epsilon Op(\rho_*(\tau)), & \quad \sup_{\tau\in [0,1]}|\rho_{*}(\tau)|^{\Lip(\gamma,\Omega_1)}_{-\alpha,s,\eta_0}\le_{\mathtt{pe},s,\eta_0} 1,\\
\sup_{\tau\in [0,1]}|d_i\rho_{*}(\tau)(i_0)[\ihat]|_{-\alpha,s,\eta_0}& =0. 
 \end{align*}
Therefore, denoting all the finite dimensional operators in the right-hand side of \eqref{jjsd2sds} by $R_{5,2}$, we obtain
\begin{align} \nonumber
\epsilon^2\mathcal{A}_2 & = \epsilon^2 \partial_xOp(\rho_1)\partial_x Op(\rho_1) \\
& \ +  \epsilon^3 \partial_xOp(\rho_1)\partial_x Op(\rho_1)\int_0^{1}(1-\tau)Op(\rho_*(\tau))d\tau + {R}_{5,2},\label{rkfxpp11sd}
\end{align}
for some $R_{5,2}\in \mathfrak{R}(i_0)$.
Therefore, using $\rho_1\in \mathfrak{B}^{-1-\alpha}_1$  and the composition formulae of symbols in \eqref{compose_289}, it follows straightforwardly that
\[
\epsilon^2\mathcal{A}_2= Op(\tilde{\mathfrak{b}}_2) + Op(\mathfrak{r}_{-2,\le 3,*}) + R_{5,2},
\]
for some $\tilde{\mathfrak{b}}_2\in \mathfrak{B}^{-2\alpha}_2\subset \mathfrak{B}^{-2}_{2}$ and $\mathfrak{r}_{-2,\le 3,*}$ satisfying the estimates \eqref{estiatefors20sd} and \eqref{estiatefors20sd1} (indeed with a better estimates since $\rho_*$ does not depend on the embedding $i_0$ and the factor $\epsilon^3$ in \eqref{rkfxpp11sd}).
\end{proof}

\subsection{Structure of the operator of size $O(\epsilon^2)$}
Before eliminating the operators of size $O(\epsilon^2)$ in $\mathcal{L}^5$, we investigate the structure of such terms in \eqref{sdlinesdar2sd2s2j}.

Recalling the definition of $\mathfrak{R}(i_0)$ from Definition~\ref{remiander_class12}, \ref{item262222} in Proposition~\ref{linearstep17} tells us that we can pick $\mathfrak{c}_1,\mathfrak{c}_2\in \mathfrak{B}^{m}_{1}$ for some $m\ge 0$ such that 
\[
R_5 = \Pi_{S^{\perp}}\partial_xOp^W(\mathfrak{c}_1)\Pi_{S}Op^W(\mathfrak{c}_2)  + \mathcal{R}_5,
\] 
where $\mathcal{R}_5$ is of the form in \eqref{j2j2j2j2j2jssds}.
Therefore, we write $\mathcal{L}^5$ as
\begin{align}\label{hlaemsksd}
\mathcal{L}^5 = \D_\omega - \Pi_{S^\perp}\partial_x Op^W\left(\mathtt{m}_\alpha m_{1,\alpha}(\xi) +\left(\frac{T_\alpha}4  + \mathfrak{m}_{\le0}(\omega,\xi)\right) \right) -\epsilon^2 \partial_x\mathcal{C}_2 + \mathcal{R}_{\le 3},
\end{align}
where
\begin{equation}\label{cdsddef}
\begin{aligned}
\mathcal{C}_2 &= \Pi_{S^{\perp}}\left( Op^W(\tilde{\mathfrak{b}}_2) - Op^W(\mathfrak{c}_1)\Pi_{S}Op^W(\mathfrak{c}_2) \right),\\
& \qquad \qquad \qquad \qquad \qquad \qquad \qquad \quad \text{ for some $\mathfrak{c}_1,\mathfrak{c}_2\in \mathfrak{B}^m_1$ for some $m\ge0$,} \\
\mathcal{R}_{\le 3} & = \Pi_{S^\perp}\partial_x Op^W(\mathfrak{r}_{-2,\le 3,*}) + \partial_x \Pi_{S^\perp}W_2 + \mathcal{R}_{5},\text{ for some $\mathcal{R}_5$ of the form in \eqref{j2j2j2j2j2jssds}.}
\end{aligned}
\end{equation}

The main result of this section is summarized in the following proposition:
\begin{proposition}\label{sd2sdsd}
$\mathcal{L}^5$ has the form:
\begin{equation}\label{l5idenitsd}
\begin{aligned}\mathcal{L}^5 & = \D_\omega - \Pi_{S^\perp}\partial_x Op^W\left(\mathtt{m}_\alpha m_{1,\alpha}(\xi) + \left(\frac{T_\alpha}4 + \mathfrak{m}_{\le 0}(\omega,\xi) + \epsilon^2\mathfrak{m}_{\mathfrak{b}}(\omega,\xi) \right) \right)\\
& \   - \epsilon^2\Pi_{S^\perp}\partial_xOp^W(\mathfrak{b}_0) + \mathcal{R}_{\le 3},
\end{aligned}
\end{equation}
and satisfies the following:
\begin{enumerate}[label=(\arabic*)]
\item \label{rjsdsd} $\mathfrak{m}_{\mathfrak{b}}(\omega,\xi):=M_x(\tilde{\mathfrak{b}})$ is an $\omega$-dependent Fourier multiplier and it is a reversible symbol. Furthermore, we have
\begin{align}\label{idernfisdsd2}
\mathfrak{m}_{\mathfrak{b}}(\xi)+ \mathtt{m}_{\alpha,1}(\omega)m_{1,\alpha}(\xi)+\mathfrak{m}_{\le 0,1}(\omega,\xi) = \mathfrak{F}_{\alpha-1}(\omega,\xi),
\end{align}
where $\mathfrak{F}_{\alpha-1}$ is as given in \eqref{identification_gasd2sx} ($\mathtt{m}_{\alpha,1}$ and $\mathfrak{m}_{\le 0,1}$ are given in Proposition~\ref{prop_time_rep_92} and Proposition~\ref{rlaqkqdpsanfsdfmf}). 
\item \label{sdsdsd2jsds} $\mathfrak{b}_0\in \mathfrak{B}^{-2}_{2}$ is a real-valued and reversible symbol such that $M_x(\mathfrak{b}_0) = 0$.
\end{enumerate}
\end{proposition}
\begin{proof}
Denoting 
\begin{align}\label{sdsd1losc}
\mathfrak{m}_{\mathfrak{b}}:=M_x(\tilde{\mathfrak{b}}_2),\quad \mathfrak{b}_0=\tilde{\mathfrak{b}}_2 - M_x(\tilde{\mathfrak{b}}_2),
\end{align}
the result follows from \eqref{hlaemsksd} and Lemma~\ref{identification_prop}. Note that $\mathfrak{m}_{\mathfrak{b}}$ is a reversible symbol since $\tilde{\mathfrak{b}}_2$ is  reversible symbol (see Proposition~\ref{linearstep17}). $\mathfrak{b}_0$ being real-valued and reversible in \ref{sdsdsd2jsds} also follows from \ref{item262221} of Proposition~\ref{linearstep17}. $M_x(\mathfrak{b}_0)=0$ by its definition in \eqref{sdsd1losc}.
\end{proof}

In the rest of this section, we aim to prove the necessary lemmas to obtain the above proposition.

\subsubsection{Structure of $\mathcal{L}^5$.}
In the next lemma, we recall the notation in \eqref{saytheyareequal} to compare two linear operators.
\begin{lemma}\label{rlaqkaadssjxxsd}
We have that
\begin{enumerate}[label=(\arabic*)]
\item \label{rprkawkxnlrla} $\mathcal{R}_{\le 3}= 0$, up to $O(\epsilon^{9-6b})$.
\item \label{rprkawkxnlrla2} $Op^W\left(\mathtt{m}_\alpha(\omega) m_{1,\alpha}(\xi) +\left(\frac{T_\alpha}4  + \mathfrak{m}_{\le0}(\omega,\xi)\right) \right)  = Op^W(-\frac{1}{2}m_{1,\alpha}(\xi) + \frac{T_\alpha}4) + \epsilon^2Op^W(\mathtt{m}_{\alpha,1}(\omega)m_{1,\alpha}(\xi) + \mathfrak{m}_{\le 0,1}(\omega,\xi))$, up to $O(\epsilon^{9-6b})$.
\end{enumerate}
\end{lemma}
\begin{proof}
 For \ref{rprkawkxnlrla}, \eqref{cdsddef} tells us that  we need to show that $Op^W(\mathfrak{r}_{-2,\le 3,*}) = 0$, $W_2=0$ and $\mathcal{R}_5=0$ up to $O(\epsilon^{9-6b})$. For $\mathfrak{r}_{-2,\le 3,*}$,   it follows from \eqref{sdsdtmddndkQk},  \eqref{size_assumption_3} and \eqref{frequency_set2} that $Op^W(\mathfrak{r}_{-2,\le 3,*}) = 0$, up to $O(\gamma^{-1}(\epsilon^5 + \epsilon^3\epsilon^{6-2b}\gamma^{-1})) = O({\epsilon^{9-6b}})$. For $\mathcal{R}_5$, we already know from \eqref{cdsddef} that $\mathcal{R}_5$ is of the form in \eqref{j2j2j2j2j2jssds}. Thanks to the estimates in \eqref{sjdjsjsj} and Lemma~\ref{s_decay_finite}, we have that $\mathcal{R}_5=0$ up to $O(\epsilon^3 + \epsilon^2|\mathfrak{I}_\delta|)\overset{\eqref{size_assumption_3}}=O(\epsilon^{8-4b})=O(\epsilon^{9-6b})$. where the last equality follows from $b\in (1,1+1/12)$, hence $9-6b\le 8-4b$. For $W_2$, we recall from \ref{item262222} of Proposition~\ref{linearstep17} and \eqref{kk2dsdsd} and \eqref{sdkksdsds} that
 \begin{align}\label{W2isdsd2}
 W_2 &=\frac{1}{(\mathtt{N}_\alpha-2)!}\int_0^1\int_0^\tau (\Psi_3(\tau-t)\Phi_4\Psi_5)^T Op^W(\mathfrak{q}_{\omega,\mathtt{N}-2}\star\mathfrak{a}) \nonumber \\
 & \qquad \qquad \qquad \qquad \qquad \qquad \qquad \qquad \times \Psi_3(\tau-t)\Phi_4\Psi_5t^{\mathtt{N}_\alpha-2}dtd\tau\nonumber\\
 & \ -\frac{1}{\mathtt{N}_\alpha!}\int_0^1(\Psi_3(1-t)\Phi_4\Psi_5)^T Op^W(\mathfrak{q}_{M_2,\mathtt{N}_\alpha}\star \mathfrak{a})\Psi_3(1-t)\Phi_4\Psi_5t^{\mathtt{N}_\alpha}dt.
 \end{align}
  From Lemma~\ref{align2sd2sd2}, Lemma~\ref{align2sd2sd2sx2}, and Lemma~\ref{align2sd2sd2sssrx2}, we have that $(\Psi_3(\tau-t)\Phi_4\Psi_5) = I$ up to $O(\epsilon)$, and so is $(\Psi_3(\tau-t)\Phi_4\Psi_5)^{T}$, thanks to the fact that all the transformations are symplectic and  \eqref{symplectic73}. While, the symbols $\mathfrak{q}_{\omega,\mathtt{N}-2}\star\mathfrak{a}$ and $\mathfrak{q}_{M_2,\mathtt{N}_\alpha}\star \mathfrak{a}$ satisfy \eqref{rlaksd1sdsd2}, which tells us that
 $W_2 = O(\epsilon^{\mathtt{N}_\alpha})$. Recalling $\mathtt{N}_\alpha$  from \eqref{nsjdjwdsdnsde}, we have $W_2=0$ up to $O(\epsilon^{9-6b})$.

\ref{rprkawkxnlrla2} follows from the estimates for $\mathtt{m}_\alpha$ and $\mathfrak{m}_{\le 0}$ in \eqref{qtildesd2dsd} and \eqref{docsd2ksdsdsd} and $\epsilon^{7-4b}\le \epsilon^{9-6b}$ for sufficiently small $\epsilon >0$.
\end{proof}
Lemma~\ref{rlaqkaadssjxxsd} tells us that denoting 
\begin{equation}\label{sdsdsdenotignd}
\begin{aligned}
D:=Op^W\left(-\frac{1}{2}m_{1,\alpha}(\xi)+\frac{T_\alpha}4\right),\quad \mathcal{T}:=-\mathcal{C}_2 - Op^W(\mathtt{m}_{\alpha,1}m_{1,\alpha}(\xi)+\mathfrak{m}_{\le 0,1}(\omega,\xi)),
\end{aligned}
\end{equation}
\eqref{hlaemsksd} can be written as 
\begin{align}\label{sdjsdsfahfsd}
\mathcal{L}^5 = \left(\D_\omega -\partial_x D\right) +\epsilon^2\partial_x\mathcal{T},\text{ up to $O(\epsilon^{9-6b})$.}
\end{align}

\subsubsection{Structure of $\mathcal{L}^5$ from the evolution of $\mathcal{L}_\omega$.} Now, we analyze the structure of $\mathcal{L}^5$ from the expression in \eqref{epsilon2ssss} and compare it with \eqref{sdjsdsfahfsd} to specify an explicit formula for $\mathcal{T}$.

We recall from \eqref{epsilon2ssss} that
\begin{align}\label{ppwosdsd2}
 \mathcal{L}_\omega =  \D_\omega - \partial_x\nabla_z^2H_{2,2}(\bar{v}) -\epsilon  \partial_x\underbrace{\nabla_z^2H_{3,2}(\bar{v})}_{=:S_1} - \epsilon^2 \partial_x\underbrace{\nabla_z^2H_{4,2}(\bar{v})}_{=:S_2}, \text{ up to $O(\epsilon^{7-4b})$.}
\end{align}
Recalling $H_{2,2}$ from \eqref{normal_formsf} and $m_{1,\alpha}$ from \eqref{multiplier_2}, we see that $\nabla_{z}^2H_{2,2} = D$, where $D$ is as in \eqref{sdsdsdenotignd}, therefore,
\begin{align}\label{klsdd2sdsd}
\mathcal{L}_\omega = \left( \D_\omega - \partial_x D\right) - \epsilon \partial_x S_1 - \epsilon^2 \partial_x S_2,\text{ up to $O(\epsilon^{9-6b})$}.
\end{align}
From Propositions~\ref{toohard_2_3}, \ref{prop_time_rep_92}, \ref{induction_egorov}, \ref{rlaqkqdpsanfsdfmf} and \ref{linearstep17}, we have that
\[
\mathcal{L}^{5} = \Phi_{5}^{-1}\Phi_4^{-1}\Phi_3^{-1}\frac{1}{\rho}\Phi_2^{-1}\Phi_{1}^{-1}\mathcal{L}_\omega\Phi_1\Phi_2\Phi_3\Phi_4\Phi_5.
\]
From Lemma~\ref{rho_estimate_time_rep} and \eqref{size_assumption_3} that $h\mapsto \frac{1}{\rho}h = I$, up to $O(\epsilon^{8-4b})=O(\epsilon^{9-6b})$, therefore,
\begin{align}\label{sdsd2dkjgsx}
\mathcal{L}^{5} = \left(\Phi_1\Phi_2\Phi_3\Phi_4\Phi_5\right)^{-1}\mathcal{L}_\omega\Phi_1\Phi_2\Phi_3\Phi_4\Phi_5,\text{ up to $O(\epsilon^{9-6b}$).}
\end{align}

\begin{lemma}\label{sdrsdus2}
$\mathcal{L}^5$ is of the form:
\begin{equation}\label{lconsdsdforms}
\begin{aligned}
\mathcal{L}^{5} &= \D_\omega - \partial_x D -\epsilon(\partial_xS_1 - [\D_\omega - \partial_xD,\mathfrak{A}_1]) \\
& \ - \epsilon^2\left(\partial_xS_2 -[\mathfrak{A}_1,\partial_xS_1] -[\D_\omega - \partial_xD , \mathfrak{A}_2] - \frac{1}{2}[[\D_\omega - \partial_x D, \mathfrak{A}_1], \mathfrak{A}_1]\right),
\end{aligned}
\end{equation}
up to $O(\epsilon^{9-6b})$, where $\mathfrak{A}_i=\Pi_{S^\perp}Op(\mathfrak{t}_i)$ for some $\mathfrak{t}_i\in \mathfrak{B}^m_{i}$ for some $m\in \mathbb{R}$. 

\end{lemma}

\begin{proof}
First, we consider an operator $L$ is of the form:
\[
L = \D_\omega - \partial_xD + \epsilon T_1 + \epsilon^2 T_2,
\]
for some bounded operators   $T_1,T_2:H^{s_0}_{\varphi,x}\mapsto H^{s_0+\mu_0}_{\varphi,x}$ for sufficiently large $\mu_0\ge 0$. 
we claim that the conjugations with $\Phi_2,\Phi_4$ do not change its structure, that is,
\[
\Phi_i^{-1}L\Phi_i = L ,\text{ up to $O(\epsilon^{9-6b})$.}
\]
For Lemma~\ref{rep_time_m_estimates126} and Lemma~\ref{p2sdsdkcssd}, we see that (using \eqref{size_assumption_3})
\[
\Phi_i = I, \text{ up to $O(\gamma^{-1}\epsilon^2 \rVert \mathfrak{I}_\delta\rVert_{s_0+\mu_0}) =O(\epsilon^{8-6b})$, for $i=2,4$.}
\]
Therefore, the conjugation of $\Phi_{i}^{-1}\epsilon T_1\Phi_i - \epsilon T_1=\Phi_{i}^{-1}\epsilon^2 T_2\Phi_i - \epsilon^2 T_2 = 0$ up to $O(\epsilon^{9-6b})$.  Furthermore, since $\Phi_{2},\Phi_4$  commute with $\partial_x D$ (see \eqref{time_rep_transforms} and \eqref{phi4defsd2wsd}), the conjugation of $\partial_xD$ does not  change. Furthermore, the estimates $\D_\omega p_1$ and $\D_\omega \mathfrak{p}_2$ in \eqref{rep_time_est3} and \eqref{p2esd1cxcwrworkdsd} tell us that
\[
\Phi_i^{-1}\D_\omega \Phi_i = \D_\omega , \text{ up to $O(\epsilon^{8-4b})=O(\epsilon^{9-6b})$ for $i=2,4$},
\]
therefore the conjugation of $\D_\omega -\partial_x D$ with $\Phi_2,\Phi_4$ does not change $L$.

Therefore, in view of \eqref{sdsd2dkjgsx}, we have
\begin{align}\label{dhfsd2sd}
\mathcal{L}^5 = \left(\Phi_1\Phi_3\Phi_5\right)^{-1}\mathcal{L}_\omega \Phi_1\Phi_3\Phi_5,\text{ up to $O(\epsilon^{9-6b})$.}
\end{align}

We observe that each of  $\Phi_i$ for $i=1,3,5$ satisfies (see \eqref{tsd2222sd21}, \eqref{2jsd232} and \eqref{skdsdpdesx})
\begin{align}\label{sdsdjdsds12s23}
\frac{d}{d\tau}\Phi_i(\tau) = \mathfrak{A}(\tau)\Phi_i(\tau),\quad \Phi_i(0) = I,
\end{align}
for some operator $\mathfrak{A}(\tau)$ of the form
\begin{align}\label{differentsx1}
\mathfrak{A}(\tau) = \epsilon \underbrace{\Pi_{S^\perp}Op(\mathfrak{t}_{1})}_{=:\mathfrak{A}_1} + \epsilon^2\underbrace{\Pi_{S^\perp}Op(\mathfrak{t}_2)}_{=:\mathfrak{A}_2}, \text{ up to $O(\epsilon^{9-6b})$,}\end{align} for some 
$\mathfrak{t}_{1},\in \mathfrak{B}^{1}_{1}$ and $\mathfrak{t}_{2}\in \mathfrak{B}^{1}_{2}$.
Indeed, for $\Phi_1$, we expand  $b(\tau)$ in \eqref{tsd2222sd21} using the Taylor expansion in the variable $\tau$ at $\tau=0$, and use the fact that $\beta\in \mathfrak{C}_1(i_0)$ to see that $\Phi_1$ satisfies the above property. For $\Phi_3,\Phi_5$, we can see it from \eqref{2jsd232} and \eqref{skdsdpdesx}, while the corresponding operator $\mathfrak{A}(\tau)$ does not depend on $\tau$, since the equations in \eqref{2jsd232} and \eqref{skdsdpdesx} are autonomous.
Since we already observed that \eqref{sdsdjdsds12s23} is well-posed (Section~\ref{study_hard_learn_more}), the usual Taylor expansion gives us that for possibly different symbols $\mathfrak{t}_1,\mathfrak{t}_2$ (thus, different $\mathfrak{A}_1,\mathfrak{A}_2$) from the ones in  \eqref{differentsx1} but still in $\mathfrak{B}^m_{1},\mathfrak{B}^m_{2}$ for some $m\in \mathbb{R}$, respectively,
\begin{align*}
\Phi_1\Phi_3\Phi_5 &= I+\epsilon \mathfrak{A}_1 + \epsilon^2\left(\frac{1}{2}(\mathfrak{A}_1)^2 + \mathfrak{A}_2\right),\\
(\Phi_1\Phi_3\Phi_5)^{-1} &= I - \epsilon\mathfrak{A}_{1} + \epsilon^2\left(\frac{1}{2}(\mathfrak{A}_1)^2 - \mathfrak{A}_2 \right),
\end{align*}
 up to $O(\epsilon^{9-6b})$.
 Therefore, using \eqref{klsdd2sdsd} and \eqref{dhfsd2sd}, we see that
 \begin{align*}
\mathcal{L}^5 & = (\D_\omega-\partial_x D) - \epsilon(\partial_x S_1  -[\D_\omega  - \partial_x D,\mathfrak{A}_1] ) \\
& \ - \epsilon^2 \left(\partial_x S_2 - [ \mathfrak{A}_1,\partial_xS_1] - [\D_\omega-\partial_x D,\mathfrak{A}_2] - \frac{1}{2}[[\D_\omega-\partial_x D,\mathfrak{A}_1],\mathfrak{A}_1]\right),
 \end{align*}
which is in the form  \eqref{lconsdsdforms}.
\end{proof}

\begin{lemma}\label{lsd2bingchandler}
We have that 
\begin{align}\label{dkwlrrkfRlfdlajfek}
\mathcal{L}^5 = \D_\omega - \partial_x D - \epsilon^2\left(\partial_xS_2 -\frac{1}{2}[\mathfrak{A}_1,\partial_xS_1] -[\D_{\overline{\omega}} - \partial_xD ,\mathfrak{A}_2] \right), \text{up to $O(\epsilon^{9-6b})$,}
\end{align}
where  $\mathfrak{A}_2 = \Pi_{S^\perp}Op(\mathfrak{t}_2)$, for some $\mathfrak{t}_2\in \mathfrak{B}^m_2$ and $\mathfrak{A}_1 = \Pi_{S^\perp}Op(\mathfrak{t}_1)$ with $\mathfrak{t}_1\in \mathfrak{B}^m_1$ for some $m\in \mathbb{R}$. $\mathfrak{t}_2$ is determined uniquely by
\begin{align}\label{uniqsdssods}
\partial_x S_1 = [\D_{\overline{\omega}}- \partial_x D, \mathfrak{A}_1]
\end{align}
\end{lemma}
\begin{proof}
Comparing $O(\epsilon)$ terms in  \eqref{sdjsdsfahfsd} and \eqref{lconsdsdforms}, and using that $\D_\omega - \D_{\overline{\omega}} = 0$ up to $O(\epsilon^{2})$, which follows from \eqref{xi_omega_dependent} and \eqref{frequency_set1}, we have \eqref{uniqsdssods}. Therefore,
\[
\frac{1}{2}[[\D_\omega - \partial_x D, \mathfrak{A}_1], \mathfrak{A}_1] = \frac{1}{2}[\partial_x S_1 , \mathfrak{A}_1].
\]
Plugging this into \eqref{lconsdsdforms}, we obtain \eqref{dkwlrrkfRlfdlajfek}.
\end{proof}

In what follows, we will find explicit solution to \eqref{uniqsdssods}.

\begin{lemma}\label{sjsjxcxcwsdsd}
Let $\mathfrak{A}_1$ be a unique solution to \eqref{uniqsdssods}. Then there exist a Fourier multiplier $\mathfrak{F}_{\alpha-1}(\xi)$ and a symbol $\mathfrak{D}_0\in \mathfrak{B}^{\alpha-1}_2(\varphi,x,\xi)$ such that 
\begin{align}\label{rksdstsdsd}
\partial_x S_2 - \frac{1}{2}[\mathfrak{A}_1,\partial_x S_1] - [\D_{\overline{\omega}} - \partial_xD ,\mathfrak{A}_2]  = \Pi_{S^\perp}\partial_x Op^W(\mathfrak{F}_{\alpha-1} + \mathfrak{D}_0),
\end{align}
such that $M_x(\mathfrak{D}_0) = 0$ and $\mathfrak{F}_{\alpha-1}(\xi) = \mathfrak{F}_{\alpha-1}(\omega,\xi)$ can be  explicitly computed
\begin{equation}\label{identification_gasd2sx}
\begin{aligned}
\mathfrak{F}_{\alpha-1}(\omega,\xi) &=  \sum_{k=1}^{\nu}\left( \frac{12}{\pi}H_{4,j_k,-j_k,\xi,-\xi} +\frac{18}{\pi^2}\left(\frac{(\xi - j_k)H_{3,-j_k,\xi,-(\xi - j_k)}H_{3,j_k,\xi-j_k,-\xi}}{\lambda_\alpha(j_k)+\lambda_\alpha(\xi-j_k)-\lambda_\alpha(\xi)} \right. \right.\\
&\left. \left. +\frac{(\xi + j_k)H_{3,j_k,\xi,-(\xi + j_k)}H_{3,-j_k,\xi + j_k,-\xi}}{-\lambda_\alpha(j_k)+\lambda_\alpha(\xi+j_k)-\lambda_\alpha(\xi)} \right)\right) j_k\zeta_k.
\end{aligned}
\end{equation}
\end{lemma}
\begin{proof}
We first find a symbol $\mathfrak{t}_1$, that solves \eqref{uniqsdssods}. 

\textbf{Solution to \eqref{uniqsdssods}.}
Let us find a symbol $\mathfrak{t}_1\in \mathfrak{B}^m_1$ such that $\mathfrak{A}_1 = \Pi_{S^\perp}Op(\mathfrak{t}_1)$ satisfies \eqref{uniqsdssods}. By Definition~\ref{buildingbobo} and recalling $\overline{v}$ from \eqref{norm_def_vbar}, we have
\begin{align}\label{t1des2sxz}
\mathfrak{t}_1 = \sum_{j_k\in S}C_{j_k}(\xi)\sqrt{j_k\zeta_{k}}e^{\ii \mathtt{l}(j_k)\cdot \varphi + \ii j_k x}.
\end{align}
Let us first compute $\partial_x S_1$ using symbolic notations. 
 In view of the definition of $S_1$ in \eqref{ppwosdsd2}, we see from \eqref{rlaksd2sdsd} that 
 \begin{align}\label{h32sd2sdsd}
 H_{3,2}(f) = \sum_{\substack{j_1+j_2+j_3=0,\\ j_1\in S,\ j_2,j_3\in S^\perp}}3H_{3,j_1,j_2,j_3}f_{j_1}f_{j_2}f_{j_3}.
 \end{align}
 Using our convention for a Fourier expansion in \eqref{sjd2sjjjsd} and the $L^2$-gradient in \eqref{sjdsdsdxcxc12s}, we determine $\nabla_z^2H(\overline{v})$ as a unique linear map such that
  \begin{align}\label{2rlaclqhrmsd2}
 (\nabla_z^2H_{3,2}(\overline{v})[h],g)_{L^2} = \frac{d^2}{dtds}H_{3,2}(\overline{v}+th+sg)\bigg|_{t=s=0},\text{ for all $h,g\in H_{S^\perp}$. }
 \end{align}
 The left-hand side of \eqref{2rlaclqhrmsd2} can be computed as 
 \begin{align}\label{dsds3sd}
 (\nabla_z^2H_{3,2}(\overline{v})[h],g)_{L^2} = \sum_{j\in \mathbb{Z}}2\pi\left(\nabla_z^2H_{3,2}(\overline{v})[h]\right)_jg_{-j},
 \end{align}
while the right-hand side can be computed from \eqref{h32sd2sdsd} as
\[
\frac{d^2}{dtds}H_{3,2}(\overline{v}+th+sg)\bigg|_{t=s=0} = \sum_{j\in \mathbb{Z}}\left(\sum_{j_1+j_2 = j}6H_{3,j_1,j_2,-j}\overline{v}_{j_1}h_{j_2}\right)g_{-j}.
\]
Combining this with \eqref{dsds3sd}, we get, using \eqref{sjd2sjjjsd},
\[
\nabla_z^2H_{3,2}(\overline{v})[h]= \frac{1}{2\pi}\sum_{k\in S^\perp}\sum_{j_1+j_2 = k}6H_{3,j_1,j_2,-k}\overline{v}_{j_1}h_{j_2}e^{\ii k x}.
\]
Therefore, we have 
 \begin{align}\label{rjsd2s0022sd}
 \partial_x S_1[h]& = \partial_x \nabla_z^2H_{3,2}(\overline{v})[h] =\frac{1}{2\pi} \sum_{k\in S^\perp}\sum_{\substack{j_1+j_2=k,\\ j_1\in S,j_2\in S^\perp}}6\ii k H_{3,j_1,j_2,-k}\overline{v}_{j_1}h_{j_2}e^{\ii k x}\nonumber\\
  &= \Pi_{S^\perp}\sum_{k\in \mathbb{Z}}\sum_{\substack{j_1+j_2=k,\\ j_1\in S,j_2\in S^\perp}}\frac{3\ii}{\pi} k H_{3,j_1,j_2,-k}\overline{v}_{j_1}h_{j_2}e^{\ii k x}\nonumber\\
  & =: \Pi_{S^\perp}Op(\mathfrak{s}_1)[h],
 \end{align}
 where (using $\overline{v}$ in \eqref{norm_def_vbar}),
 \begin{align}\label{rlaqkqdprlaclwjs123s}
 {\mathfrak{s}_1}(x,\xi):=\sum_{j_k\in S}\frac{3\ii}{\pi} (j_k + \xi)H_{3,j_k,\xi,-(\xi + j_k)}\sqrt{j_k\zeta_k}e^{\ii \mathtt{l}(j_k)\cdot \varphi + \ii j_k x}.
 \end{align}
 Now, we compute the right-hand side of  \eqref{uniqsdssods}.  Using \eqref{t1des2sxz}, we have
 \begin{align}\label{rzlmdptsd}
 [\D_{\overline{\omega}} - \partial_x D, Op( \mathfrak{t}_1)] &= \D_{\overline{\omega}}(Op(\mathfrak{t}_1) )- [\partial_x D,Op(\mathfrak{t}_1)]\nonumber\\
 & = Op(\D_{\overline{\omega}}\mathfrak{t}_1)-[\partial_x D, Op(\mathfrak{t}_1)]\nonumber\\
 & = Op(\mathfrak{t}),
 \end{align}
 where  (using $\overline{\omega}$, and $D$ in \eqref{linear_frequency_aa} and \eqref{sdsdsdenotignd}),
  \begin{align}
 \mathfrak{t}&=\sum_{j_k\in S}C_{j_k}(\xi)\sqrt{j_k\zeta_k}\left(\ii \overline{\omega}_k-\left( \ii(\xi + j_k)\left(-\frac{1}{2}m_{1,\alpha}(\xi + j_k)+\frac{T_\alpha}{4}\right) \right. \right. \nonumber \\
 & \qquad \qquad \qquad \qquad \qquad \qquad \qquad \qquad \left. \left. - \ii \xi\left( -\frac{1}{2}m_{1,\alpha}(\xi)+\frac{T_\alpha}4\right) \right)\right)e^{\ii \mathtt{l}(j_k)\cdot \varphi + \ii j_k x}\nonumber\\
 & = \sum_{j_k\in S}\frac{1}{2}\ii C_{j_k}(\xi)\sqrt{j_k\zeta_k}\left( \lambda_{\alpha}(\xi+j_k) - \lambda_\alpha(\xi)-\lambda_\alpha(j_k)\right)e^{\ii \mathtt{l}(j_k)\cdot \varphi + \ii j_k x},\label{tdsd2sys}
 \end{align}
 where the second equality can be obtained from the definition of $\lambda_\alpha$ in \eqref{lambdadeffsd}. Therefore, plugging \eqref{rjsd2s0022sd} and \eqref{rzlmdptsd} into \eqref{uniqsdssods}, we have
 \[
\Pi_{S^\perp}Op(\mathfrak{s}_1) =  \partial_x S_1 = \Pi_{S^\perp}[\D_{\overline{\omega}}-\partial_x D,Op(\mathfrak{t}_1)] =\Pi_{S^\perp}Op( \mathfrak{t}),
 \]
 which gives (using \eqref{rlaqkqdprlaclwjs123s} and \eqref{tdsd2sys}),
 \begin{align}\label{ctdefs11sd}
 C_{j_k}(\xi) = \frac{6}{\pi}\frac{(j_k+\xi)H_{3,j_k,\xi,-(\xi + j_k)}}{\lambda_{\alpha}(\xi + j_k)-\lambda_\alpha(\xi) - \lambda_\alpha(j_k)}.
 \end{align}
 Therefore, we obtain from \eqref{t1des2sxz} that
 \begin{align}\nonumber
 \mathfrak{A}_1&=\Pi_{S^\perp}Op(\mathfrak{t}_1),\\
 \mathfrak{t}_1(x,\xi) &:= \frac{6}{\pi}\sum_{j_k\in S}\frac{(j_k+\xi)H_{3,j_k,\xi,-(\xi + j_k)}}{\lambda_{\alpha}(\xi + j_k)-\lambda_\alpha(\xi) - \lambda_\alpha(j_k)}\sqrt{j_k\zeta_{k}}e^{\ii \mathtt{l}(j_k)\cdot \varphi + \ii j_k x}.\label{dkwlsWkgkqltlqkftoRl}
 \end{align}
 
 \textbf{Computing $\partial_x S_2 - \frac{1}{2}[\mathfrak{A}_1,\partial_x S_1]$.}
 First,  using \eqref{ppwosdsd2} and \eqref{sjdjsjs93jx}, we have that for all $h\in H_{S^\perp}$,
 \begin{align*}
 \partial_x S_2[h] &=\partial_x\nabla_{z}^2\mathcal{H}_{4,2}(\overline{v})[h] + \partial_x\nabla_{z}^2\mathfrak{H}_{1}(\overline{v})[h] + \partial_x\nabla_{z}^2\mathfrak{H}_{2}(\overline{v})[h].
 \end{align*}
 Therefore, we have that
 \begin{align}\label{rsdsdgkqltlqkftoRL}
 \partial_x S_2 -\frac{1}{2}[\mathfrak{A}_1,\partial_x S_1] =\left(\partial_x \nabla_z^2\mathfrak{H}_1(\overline{v}) - \frac{1}{2}[\mathfrak{A}_1,\partial_x S_1]\right) + \partial_x\nabla_{z}^2\mathcal{H}_{4,2}(\overline{v})  + \partial_x\nabla_{z}^2\mathfrak{H}_{2}(\overline{v})
 \end{align}
 
 Let us compute $\partial_x\nabla_{z}^2\mathcal{H}_{4,2}(\overline{v})  + \partial_x\nabla_{z}^2\mathfrak{H}_{2}(\overline{v})$ first. 
 From \eqref{rlaksd2sdsd1} and \eqref{rlaksd2sdsd4}, we have
 \begin{align}\label{notajokeson}
 \mathcal{H}_{4,2}(f) + \mathfrak{H}_2(f) &= \sum_{\substack{j_1+j_2+j_3+j_4=0,\\ j_1,j_2\in S,\ j_3,j_4\in S^\perp\\ j_1+j_2 =0, j_3+j_4 =0}}6H_{4,j_1,j_2,j_3,j_4}f_{j_1}f_{j_2}f_{j_3}f_{j_4}\nonumber\\
 & \ +  \sum_{\substack{j_1+j_2+j_3+j_4=0,\ j_1,j_2\in S,\\ j_3,j_4\in S^\perp,\ j_1+j_2\ne 0} }G_{j_1,j_2,j_3,j_4}f_{j_1}f_{j_2}f_{j_3}f_{j_4}\nonumber\\
 & =: A_1(f) + A_2(f).
 \end{align}
 where 
 \[
 G_{j_1,j_2,j_3,j_4}:=6H_{4,j_1,j_2,j_3,j_4}-\frac{9}2(j_3+j_4)\frac{H_{3,j_1,j_2,-(j_1+j_2)}H_{3,j_3,j_4,-(j_3+j_4)}}{j_1\kappa_{j_1}+j_2\kappa_{j_2}-(j_1+j_2)\kappa_{j_1+j_2}},
 \]
  where we used \ref{tangent_2}, which implies that $j_1+j_2\in S^\perp\iff j_1+j_2\ne 0$, for $j_1,j_2\in S$. 
 For the first term $A_1$, we have (using the property of our tangential sites $S$ in \eqref{tan_site}),
 \begin{align*}
A_1(f) = 12\sum_{j_1\in S^+, j_3 \in S^{\perp}}H_{4,j_1,-j_1,j_3,-j_3}f_{j_1}f_{-j_1}f_{j_3}f_{-j_3}.
 \end{align*}
 Therefore, for all $h,g\in H_{S^\perp}$,
 \begin{align}\nonumber
 (\nabla_z^2A_1(\overline{v})[h],g)_{L^2} & = \frac{d^2}{dtds}A_1(\overline{v}+th+sg)\bigg|_{t=s=0} \\
 & = 24\sum_{j_1\in S^+,\ j_3\in S^\perp}H_{4,j_1,-j_1,j_3,-j_3}|\overline{v}_{j_1}|^2h_{j_3}g_{-j_3}.\label{tnseoqRma}
 \end{align}
 Defining a Fourier multiplier $\mathfrak{D}_1(\xi)$ by   (recalling our convention \eqref{sjd2sjjjsd})
 \begin{align}\label{d21sdsy}
 \mathfrak{D}_1(\xi):=\frac{12}{\pi}\sum_{j_1\in S^\perp}H_{4,j_1,-j_1,\xi,-\xi}|\overline{v}_{j_1}|^2,
 \end{align}
 we have 
 \[
   (Op^W(\mathfrak{D}_1)h,g)_{L^2} = 2\pi\sum_{j_3\in S^\perp}(Op^W(\mathfrak{D}_1)h)_{j_3}g_{-j_3} = (\nabla_z^2A(\overline{v})[h],g)_{L^2},
   \]
    hence,
 \begin{align}\label{wksclrnrtn}
 \nabla_z^2 A_1(\overline{v})[h] =  \Pi_{S^\perp}Op^W(\mathfrak{D}_1)h .
 \end{align}
 For $A_2$ in \eqref{notajokeson}, a similar argument as above shows that we can find a symbol $\mathfrak{D}_2(\varphi,x,\xi)$ such that
 \[
 \nabla_z^2A_2(\overline{v})[h] = \Pi_{S^\perp}Op^W(\mathfrak{D}_2)h,\quad \int_{\mathbb{T}}\mathfrak{D}_2(\varphi,x,\xi)dx = 0.
 \]
 Hence, plugging this and \eqref{wksclrnrtn} into \eqref{notajokeson}, we obtain
 \begin{align}\label{djehdjfjqmf}
 \partial_x\nabla_{z}^2\mathcal{H}_{4,2}(\overline{v})  + \partial_x\nabla_{z}^2\mathfrak{H}_{2}(\overline{v}) = \Pi_{S^\perp}\partial_xOp^W(\mathfrak{D}_1 + \mathfrak{D}_2).
 \end{align}

 Now, we compute $\partial_x \nabla_z^2\mathfrak{H}_1(\overline{v}) - \frac{1}{2}[\mathfrak{A}_1,\partial_x S_1]$ in \eqref{rsdsdgkqltlqkftoRL}. Using \eqref{rlaksd2sdsd3}, we have that for all $h,g\in H_{S^\perp}$,
 \begin{align*}
( \nabla_z^2&\mathfrak{H}_1(\overline{v})[h],g)_{L^2} =\frac{d^2}{dtds}\mathfrak{H}_1(\overline{v}+th+sg)|_{t=s=0}\\
& =  -9\sum_{\substack{j_1+j_2+j_3+j_4=0,\\ j_1,j_2\in S,\ j_3,j_4\in S^\perp, \\ j_2+j_3\in S}}(j_2+j_3)\frac{H_{3,-(j_2+j_3),j_2,j_3}H_{3,j_1,j_4,-(j_1+j_4)}}{j_1\kappa_{j_1}+j_4\kappa_{j_4}-(j_1+j_4)\kappa_{j_1+j_4}}\overline{v}_{j_1}\overline{v}_{j_2}h_{j_3}g_{j_4}\\
&  \ -9\sum_{\substack{j_1+j_2+j_3 + j_4=0,\\ j_1,j_2\in S,\ j_3,j_4\in S^\perp, \\ j_2+j_3\in S}}(j_1+j_4)\frac{H_{3,-(j_1+j_4),j_1,j_4}H_{3,j_2,j_3,-(j_2+j_3)}}{j_2\kappa_{j_2}+j_3\kappa_{j_3}-(j_2+j_3)\kappa_{j_2+j_3}}\overline{v}_{j_1}\overline{v}_{j_2}h_{j_3}g_{j_4},
 \end{align*}
 where we switched the indices $j_3$ with $j_4$ and $j_1$ with $j_2$ so that the Fourier mode of $h$ is denoted by the subindex $j_3$.
 Hence,
 \begin{align}\label{h1sdsdsxxxx}
 \partial_x& \nabla_z^2\mathfrak{H}_1(\overline{v})[h] \nonumber\\
 &= -\frac{9}{2\pi}\sum_{\substack{j_1+j_2+j_3 = j, \\ j_1,j_2\in S,\ j_3\in S^\perp,\\ j_2+j_3\in S}}(j_2+j_3)(j_1+j_2+j_3)\ii \nonumber \\
 & \qquad \qquad \qquad \times \frac{H_{3,-(j_2+j_3),j_2,j_3}H_{3,j_1,-(j_1+j_2+j_3),j_2+j_3}}{j_1\kappa_{j_1}+(j_2+j_3)\kappa_{j_2+j_3}-(j_1+j_2+j_3)\kappa_{j_1+j_2+j_3}}\overline{v}_{j_1}\overline{v}_{j_2}h_{j_3}e^{\ii j x}\nonumber\\
 & \  + \frac{9}{2\pi}\sum_{\substack{j_1+j_2+j_3=j,\\ j_1,j_2\in S,\ j_3,j_4\in S^\perp, \\ j_2+j_3\in S}}(j_2+j_3)(j_1+j_2+j_3)\ii \nonumber \\
 & \qquad \qquad \qquad \qquad \qquad\times \frac{H_{3,j_2+j_3,j_1,-(j_1+j_2+j_3)}H_{3,j_2,j_3,-(j_2+j_3)}}{j_2\kappa_{j_2}+j_3\kappa_{j_3}-(j_2+j_3)\kappa_{j_2+j_3}}\overline{v}_{j_1}\overline{v}_{j_2}h_{j_3}e^{\ii j x},
 \end{align}
 where we used $j_1+j_2+j_3 = j$ in the second summation.

Now,  we  compute $[\mathfrak{A}_1,\partial_x S_1] = \Pi_{S^\perp}\left(Op(\mathfrak{t}_1)\Pi_{S^\perp}Op(\mathfrak{s}_1) - Op(\mathfrak{s}_1)\Pi_{S^\perp}Op(\mathfrak{t}_1)\right)$. 
 Using \eqref{rlaqkqdprlaclwjs123s} and \eqref{dkwlsWkgkqltlqkftoRl}, let us denote
 \begin{align}\nonumber
 \widehat{\mathfrak{s}}_1(j_k,\xi)&:=\frac{3\ii}{\pi}(\xi + j_k)H_{3,j_k,\xi,-(\xi+j_k)}\overline{v}_{j_k},\\
 \widehat{\mathfrak{t}}_1(j_k,\xi)&:=\frac{6}{\pi}\frac{(j_k+\xi)H_{3,j_k,\xi,-(\xi+j_k)}}{\lambda_{\alpha}(\xi + j_k)-\lambda_\alpha(\xi) - \lambda_\alpha(j_k)}\overline{v}_{j_k},\label{s1sdskts1}
 \end{align}
which are the Fourier coefficients of $\mathfrak{s}_1$ and $\mathfrak{t}_1$. Therefore, we have that for all $h\in H_{S^\perp}$,
\begin{equation}\label{tjsd2wsdsds}
\begin{aligned}
& \left(Op(\mathfrak{t}_1)\Pi_{S^\perp}Op(\mathfrak{s}_1) - Op(\mathfrak{s}_1)\Pi_{S^\perp}Op(\mathfrak{t}_1)\right)[h] \\
&\ =\sum_{\substack{j_1+j_2+j_3=j,\\ j_1,j_2\in S, j_3\in S^\perp,\\ j_1+j_3\in S^\perp}}\left( \widehat{\mathfrak{t}}_1(j_2,j_1+j_3)\widehat{\mathfrak{s}}_1(j_1,j_3) - \widehat{\mathfrak{s}}_1(j_2,j_1+j_3)\widehat{\mathfrak{t}}_1(j_1,j_3)\right)h_{j_3}e^{\ii j x}\\
& =\frac{18\ii}{\pi^2} \sum_{\substack{j_1+j_2+j_3=j,\\ j_1,j_2\in S, j_3\in S^\perp,\\ j_1+j_3\in S^\perp}}\frac{(j_1+j_3)(j_1+j_2+j_3)H_{3,j_2,j_1+j_3,-(j_1+j_2+j_3)}H_{3,j_1,j_3,-(j_1+j_3)}}{\lambda_\alpha(j_1+j_2+j_3)-\lambda_\alpha(j_1+j_3)-\lambda_\alpha(j_2)}\\
& \qquad \qquad \qquad \qquad \qquad \qquad \qquad \qquad \qquad \qquad \qquad\qquad \qquad \times \overline{v}_{j_1}\overline{v}_{j_2}h_{j_3}e^{\ii j x}. \\
& \ - \frac{18\ii}{\pi^2} \sum_{\substack{j_1+j_2+j_3=j,\\ j_1,j_2\in S, j_3\in S^\perp,\\ j_1+j_3\in S^\perp}}\frac{(j_1+j_3)(j_1+j_2+j_3)H_{3,j_2,j_1+j_3,-(j_1+j_2+j_3)}H_{3,j_1,j_3,-(j_1+j_3)}}{\lambda_\alpha(j_1+j_3)-\lambda_\alpha(j_1)-\lambda_\alpha(j_3)} \\
& \qquad \qquad \qquad \qquad \qquad \qquad \qquad \qquad \qquad \qquad \qquad\qquad \qquad \times \overline{v}_{j_1}\overline{v}_{j_2}h_{j_3}e^{\ii j x}.
\end{aligned}
\end{equation}
We will combine this with $\partial_x \nabla_z^2 \mathfrak{H}_1$ in \eqref{h1sdsdsxxxx}. 
Recall the definitions of $\kappa_j$ in \eqref{moreexplicits2} and $\lambda_\alpha$ in \eqref{lambdadeffsd} so that we have
\[
j\kappa_j  = - \frac{\pi}2 \lambda_\alpha(j) + \frac{\pi T_\alpha}4j.
\]
Plugging this into \eqref{h1sdsdsxxxx} and and switching $j_1$ and $j_2$ in \eqref{tjsd2wsdsds}, we see that the summands in each summation are exactly the same, and 
\begin{align*}
&\partial_x \nabla_z^2\mathfrak{H}_1(\overline{v})[h] -\frac{1}{2}{[\mathfrak{A}_1,\partial_x S_1]}\\
&=\frac{9\ii}{\pi^2}\sum_{\substack{j_1+j_2+j_3 = j, \\ j_1,j_2\in S,\ j_3\in S^\perp}}(j_2+j_3)(j_1+j_2+j_3) \frac{H_{3,j_2,j_3,-(j_2+j_3)}H_{3,j_1,j_2+j_3,-(j_1+j_2+j_3)}}{\lambda_\alpha(j_1) + \lambda_\alpha({j_2+j_3}) - \lambda_\alpha(j_1+j_2+j_3)}\\
& \qquad \qquad \qquad \qquad \qquad \qquad \qquad \qquad \qquad \qquad \qquad\qquad \qquad \times \overline{v}_{j_1}\overline{v}_{j_2}h_{j_3}e^{\ii j x}. \\
& \ -\frac{9\ii}{\pi^2}\sum_{\substack{j_1+j_2+j_3 = j, \\ j_1,j_2\in S,\ j_3\in S^\perp}}(j_2+j_3)(j_1+j_2+j_3)\frac{H_{3,j_2,j_3,-(j_2+j_3)}H_{3,j_1,j_2+j_3,-(j_1+j_2+j_3)}}{\lambda_\alpha(j_2) + \lambda_\alpha(j_3) - \lambda_\alpha(j_2+j_3)}   \\
& \qquad \qquad \qquad \qquad \qquad \qquad \qquad \qquad \qquad \qquad \qquad\qquad \qquad \times \overline{v}_{j_1}\overline{v}_{j_2}h_{j_3}e^{\ii j x}.
\end{align*}
It is straightforward to see that denoting a symbol $\mathfrak{D}$ as
\begin{equation}\label{dsd2sdjgsdsdss}
\begin{aligned}
\mathfrak{D}(\varphi,x,\xi)&:=\sum_{j_1,j_2\in S}D_{j_1,j_2}(\xi)\overline{v}_{j_1}(\varphi,x)\overline{v}_{j_2}(\varphi,x),\\
D_{j_1,j_2}(\xi)&:=\frac{9}{\pi^2}\left(\xi + \frac{j_2-j_1}{2}\right)H_{3,j_2,\xi-\frac{j_1+j_2}2,-(\xi+\frac{j_2-j_1}2)}H_{3,j_1,\xi+\frac{j_2-j_1}2, -(\xi + \frac{j_1+j_2}2)}\\
& \times \left( \frac{1}{\lambda_\alpha(j_1) + \lambda_\alpha(\xi + \frac{j_2-j_1}2 ) - \lambda_\alpha(\xi + \frac{j_1+j_2}2))}\right. \\
& \qquad \qquad \qquad \qquad \left. -\frac{1}{\lambda_\alpha(j_2) + \lambda_\alpha(\xi - \frac{j_1+j_2}2) - \lambda_\alpha(\xi + \frac{j_2-j_1}2)} \right),
\end{aligned}
\end{equation}
we have
\[
\partial_x \nabla_z^2\mathfrak{H}_1(\overline{v}) -\frac{1}{2}{[\mathfrak{A}_1,\partial_x S_1]} = \Pi_{S^\perp}\partial_x Op^W(\mathfrak{D}).
\]
As in \eqref{djehdjfjqmf}, we can decompose $\mathfrak{D}$ as $M_x(\mathfrak{D}) + (\mathfrak{D} - M_x(\mathfrak{D}))$, more precisely,
\begin{align*}
\mathfrak{D}(\varphi,x,\xi) &= \sum_{j_1\in S^+}\left(D_{j_1,-j_1} + D_{-j_1,j_1}\right)(\xi)|\overline{v}_{j_1}(\varphi,x)|^2 \\
& \ + \sum_{j_1+j_2 \ne 0}D_{j_1,j_2}(\xi)\overline{v}_{j_1}(\varphi,x)\overline{v}_{j_2}(\varphi,x)\\
& = \mathfrak{D}_3(\xi) + 
\mathfrak{D}_4(\varphi,x,\xi) \\
\mathfrak{D}_3(\xi) & =\sum_{j_1\in S^+}\frac{18}{\pi^2}\left(\frac{(\xi - j_1)H_{3,-j_1,\xi,-(\xi - j_1)}H_{3,j_1,\xi-j_1,-\xi}}{\lambda_\alpha(j_1)+\lambda_\alpha(\xi-j_1)-\lambda_\alpha(\xi)} \right.\\
& \qquad \qquad \qquad \left.+\frac{(\xi + j_1)H_{3,j_1,\xi,-(\xi + j_1)}H_{3,-j_1,\xi + j_1,-\xi}}{-\lambda_\alpha(j_1)+\lambda_\alpha(\xi+j_1)-\lambda_\alpha(\xi)} \right)|\overline{v}_{j_1}|^2 \\
\mathfrak{D}_4(\varphi,x,\xi) & = \sum_{j_1+j_2\ne 0, j_1,j_2\in S}D_{j_1,j_2}(\xi)\overline{v}_{j_1}(\varphi,x)\overline{v}_{j_2}(\varphi,x),
\end{align*}
so that we have 
\begin{align}\label{rksdsecodnsdw}
\partial_x \nabla_z^2\mathfrak{H}_1(\overline{v}) -\frac{1}{2}{[\mathfrak{A}_1,\partial_x S_1]} = \partial_x Op^W(\mathfrak{D}_3 + \mathfrak{D}_4),
\end{align}
 where $\int_{\mathbb{T}}\mathfrak{D}_4(\varphi,x,\xi)dx = 0$. Plugging this and \eqref{djehdjfjqmf} into  \eqref{rsdsdgkqltlqkftoRL} and using \eqref{norm_def_vbar}, we get 
\begin{equation}\label{kjjjd2sdsd}
\begin{aligned}
 \partial_x S_2 -\frac{1}{2}[\mathfrak{A}_1,\partial_x S_1] &= \Pi_{S^\perp}\partial_x Op^W(\mathfrak{F}_{\alpha-1}(\xi) +\mathfrak{D}_5),
 \end{aligned}
\end{equation}
for some symbol $\mathfrak{D}_5\in \mathfrak{B}^{\alpha-1}_2$ such that $M_x(\mathfrak{D}_5)=0$ and $\mathfrak{F}_{\alpha-1}:=\mathfrak{D}_1 + \mathfrak{D}_3$, which is \eqref{identification_gasd2sx}.

\textbf{Contribution from $[\D_{\overline{\omega}} - \partial_xD ,\mathfrak{A}_2]$.} From Lemma~\ref{sdrsdus2}, we see that there exists a symbol $\mathfrak{t}_2\in \mathfrak{B}^1_{2}$ that defines $\mathfrak{A}_2$ in such a way that
\[
\mathfrak{A}_2 = \Pi_{S^\perp}Op(\mathfrak{t}_2).
\]

Since  $\D_{\overline{\omega}} - \partial_x D$ commutes with $\Pi_{S^\perp}$ and $Op(M_x(\mathfrak{t}_2))$,  it is straightforward that there is a symbol $\mathfrak{D}_6\in \mathfrak{B}^m_{2}$ for some $m\in \mathbb{R}$ such that 
\[
[\D_{\overline{\omega}} - \partial_xD ,\mathfrak{A}_2] = \Pi_{S^\perp}\partial_xOp^W(\mathfrak{D}_6),\quad M_x(\mathfrak{D}_6) = 0.
\]
Combining this with \eqref{kjjjd2sdsd}, we can denote $\mathfrak{D}_0:=\mathfrak{D}_5+\mathfrak{D}_6$ so that we obtain \eqref{rksdstsdsd}, with $M_x(\mathfrak{D}_0)=0$.
\end{proof}

Now, we specify the structure of $\mathcal{C}_2$ in \eqref{hlaemsksd}.
\begin{lemma}\label{identification_prop}
$\mathcal{C}_2$ in \eqref{cdsddef} satisfies
\begin{align}\label{jjsd2sdsdsd}
\mathcal{C}_2 =\Pi_{S^\perp} Op^W(\tilde{\mathfrak{b}}).
\end{align}
Furthermore,  we have that
\begin{align}\label{jstindefsdx}
M_x(\tilde{\mathfrak{b}}) + \mathtt{m}_{\alpha,1}m_{1,\alpha}(\xi)+\mathfrak{m}_{\le 0,1}(\omega,\xi) = \mathfrak{F}_{\alpha-1}(\omega,\xi),
\end{align}
where $\mathfrak{F}_{\alpha-1}$ is as given in \eqref{identification_gasd2sx}.
\end{lemma}
\begin{proof}
From \eqref{sdjsdsfahfsd}, \eqref{sdsdsdenotignd} and \eqref{cdsddef}, we have that
\begin{align*}
\mathcal{L}^5 & = \D_\omega - \partial_x D \\
& - \epsilon^2 \Pi_{S^\perp}\partial_x \left( Op^W(\tilde{\mathfrak{b}}_2+\mathtt{m}_{\alpha,1}m_{1,\alpha}(\xi)+\mathfrak{m}_{\le 0,1}(\omega,\xi)) - Op^W(\mathfrak{c}_1)\Pi_{S}Op^W(\mathfrak{c}_2)  \right), 
\end{align*}
up to $O(\epsilon^{9-6b})$.
On the other hand, Lemma~\ref{lsd2bingchandler} and Lemma~\ref{sjsjxcxcwsdsd} tells us that
\[
\mathcal{L}^5 = \D_\omega - \epsilon^2\partial_x D-\Pi_{S^\perp}\partial_xOp^W(\mathfrak{F}_{\alpha-1}(\omega,\xi) + \mathfrak{D}_0), \text{ up to $O(\epsilon^{9-6b})$.}
\]
Equating the two expressions of $\mathcal{L}^5$, we have
\begin{align*}
Op^W(\tilde{\mathfrak{b}}_2 + \mathtt{m}_{\alpha,1}m_{1,\alpha}(\xi)+\mathfrak{m}_{\le 0,1}(\omega,\xi)) & + Op^W(\mathfrak{c}_1)\Pi_{S}Op^W(\mathfrak{c}_2 ) \\
& = Op^W(\mathfrak{F}_{\alpha-1}(\omega,\xi) + \mathfrak{D}_0),
\end{align*} up to $O(\epsilon^{7-6b})$.
Note that the operators on the both sides are independent of $\epsilon>0$, while, $7-6b >0$, thanks to the range of $b$ in \eqref{parametersets1231}. Therefore, for sufficiently small $\epsilon>0$, we have exact equality and therefore, $\mathfrak{c}_1 = \mathfrak{c}_2=0$, which gives \eqref{jjsd2sdsdsd}, and 
\[
\tilde{\mathfrak{b}}_2 + \mathtt{m}_{\alpha,1}m_{1,\alpha}(\xi)+\mathfrak{m}_{\le 0,1}(\omega,\xi) = \mathfrak{F}_{\alpha-1}(\omega,\xi) + \mathfrak{D}_0.
\]
Using Lemma~\ref{sjsjxcxcwsdsd}, which implies $M_x(\mathfrak{D}_0)=0$ and taking the average in $x$ on both sides in the above equality, we get \eqref{jstindefsdx}.

\end{proof}

\subsection{Linear Birkhoff normal form: Step 2}\label{firststeplbnf2}
In this subsection, we eliminate the term $\mathfrak{b}_0$ in $\mathcal{L}^5$ in \eqref{sd2sdsd}.  
In order to remove $\mathfrak{b}_0$, we use a transformation $\Phi_6$ defined by
\begin{align}\label{phi6deff}
\Phi_{6}:=e^{\epsilon^2 A_{1,*}},\quad \text{ where $A_{1,*}:=\partial_x\Pi_{S^\perp}Op^W({\rho_2})$},
\end{align}
for some symbol $\rho_2=\rho_2(\omega,\varphi,\cdot,\cdot)\in \mathcal{S}^{-1-\alpha}$, which will be chosen later (see {Proposition~\ref{choiceofrho2common}}).  Note that $\Phi_6$ is a time-$1$ flow map of a pseudo PDE:
\begin{align}\label{skdsdssdpdesx}
\frac{d}{d\tau}u = \partial_x \Pi_{S^\perp}\epsilon^2Op^W(\rho_2)u.
\end{align}
Since $\partial_x Op^W(\rho_2) \in \mathcal{S}^{-\alpha}$,  $e^{\epsilon^2 A_{1,*}}:=\sum_{n=0}^{\infty}\frac{1}{n!}(\epsilon^2 A_{1,*})^n$ is well-defined {(see {Lemma~\ref{sdsdsdsd222kdkehsdjs2})}}.

 As in Section~\ref{study_hard_learn_more} (see Lemma~\ref{finitesd2sd}), we will decompose $\Phi_6$ into a finite dimensional operator and $\Psi_6$, which is the time-$1$ flow map of the pseudo-PDE:
 \begin{align}\label{psi6defsseconddat2e}
 \partial_\tau u = \epsilon^2\partial_x Op^W(\rho_2)u,\text{ that is, }\frac{d}{d\tau}\Psi_6(\tau) = \epsilon^2\partial_xOp^W(\rho_2)\Psi_6(\tau),\quad \Psi_6(0)=I.
 \end{align}

The main result is as follows:
\begin{proposition}\label{linearstep172sd}
There exists a symplectic transformation $\Phi_6:H_{S^\perp}\mapsto H_{S^\perp}$ such that
\begin{equation}\label{sdlinesdar2sd2ssdxcxcs2j}
\begin{aligned}
\mathcal{L}^6[h]&:=(\Phi_6)^{-1}\mathcal{L}^5\Phi_6[h] \\
& = \D_\omega h -  \Pi_{S^\perp}\partial_xM_6[h] +\Pi_{S^\perp}\partial_x Op^W(\mathfrak{r}_{-2,\le3,\sharp}) + \partial_x\Pi_{S^\perp} W_3 + R_6[h],\\
M_6[h] &:=  Op^W\left(\mathtt{m}_\alpha m_{1,\alpha}(\xi) + \left(\frac{T_\alpha}4 + \mathfrak{m}_{\le 0}(\omega,\xi) + \epsilon^2\mathfrak{m}_{\mathfrak{b}}(\omega,\xi) \right) \right),
\end{aligned}
\end{equation}
satisfies the following:

\begin{enumerate}[label=(\arabic*)]
\item \label{rjsdj2s2221sdsddsd} The symbol $\mathtt{m}_\alpha m_{1,\alpha}(\xi) + \left(\frac{T_\alpha}4 + \mathfrak{m}_{\le 0}(\omega,\xi) + \epsilon^2\mathfrak{m}_{\mathfrak{b}}(\omega,\xi) \right) $ is reversible
\item \label{item22sx62s92} $\mathfrak{r}_{-2,\le3,\sharp}$ satisfies
\begin{equation}\label{sdsdtm2sdsdddndkQk}
\begin{aligned}
&|\mathfrak{r}_{-2,\le3,\sharp}|^{\Lip(\gamma,\Omega_1)}_{-2,s,\eta_0} \le_{\mathtt{pe},s,\eta_0}    \gamma^{-1}\left(\epsilon^5 + \epsilon^3\rVert \mathfrak{I}_\delta\rVert^{\Lip(\gamma,\Omega_1)}_{s+\mu_0}\right),\\
&|d_i\mathfrak{r}_{-2,\le3,\sharp}(i_0)[\ihat]|_{-2,s,\eta_0} \le_{\mathtt{pe},s,\eta_0}   \epsilon^3\gamma^{-1}\left( \rVert \ihat \rVert_{s+\mu_0} + \rVert \mathfrak{I}_\delta\rVert_{s+\mu_0}\rVert \ihat \rVert_{s_0+\mu_0}\right).
\end{aligned}
\end{equation}
\item \label{item262gf222} $W_{3}:=\Psi_{6}^{T}W_2\Psi_6$, and  $R_6$ is of the form in \eqref{j2j2j2j2j2jssds}.
\item \label{rlsd2ssss2ss}$\Phi_6$ is real, reversibility preserving and symplectic. Therefore, $\mathcal{L}^6$ is real, reversible and Hamiltonian.

\item\label{egppsd}  $\mathcal{L}^6$  and $\Phi_6$ are $\frac{2\pi}{\mathtt{M}}$-translation invariance preserving.
\end{enumerate}
\end{proposition}

\begin{proof}
The expression for $\mathcal{L}^6$ given in \eqref{sdlinesdar2sd2ssdxcxcs2j} follows from Lemma~\ref{skdsdwdsd}, Proposition~\ref{choiceofrho2common}, Lemma~\ref{fisd2sd002sd1} and Lemma~\ref{wobsdand1q}, denoting $R_6:=R_{6,1}+R_{6,2}$.  To see \ref{rjsdj2s2221sdsddsd}, note that $\xi\mapsto m_{1,\alpha}(\xi)$ is even (see \eqref{explicit_multiplier}) Also, $\mathfrak{m}_{\le_0}$ is a reversible symbol, which follows from Proposition~\ref{rlaqkqdpsanfsdfmf}. From Proposition~\ref{sd2sdsd}, it follows that $\mathfrak{m}_\mathfrak{b}$ is also reversible, therefore \ref{rjsdj2s2221sdsddsd} follows immediately.   \ref{item22sx62s92} is proved in Lemma~\ref{lastbs2sd}. \ref{item262gf222} also follows from Lemma~\ref{fisd2sd002sd1} and Lemma~\ref{lastbs2sd}. \ref{rlsd2ssss2ss} is a result of Lemma~\ref{symbsrerer2}. Lastly, for \ref{egppsd}, we note from \ref{egor1sd1211x2ppsd} of Proposition~\ref{linearstep17} that $\mathcal{L}^5$ is $\frac{2\pi}{\mathtt{M}}$-translation invariance preserving. Therefore, the symbol $\mathfrak{b}_0$ in \eqref{l5idenitsd} is also $\frac{2\pi}{\mathtt{M}}$-translation invariant, which implies the Fourier coefficients of  $x\mapsto \mathfrak{b}_0(\varphi,x,\xi)$ are supported only on the modes which are the multiplies of  $\mathtt{M}$.  This property holds for $\rho_2$, since $\rho_2\in \mathfrak{B}^{-1-\alpha}_2$ (Proposition~\ref{choiceofrho2common}) therefore it is supported on the Fourier modes that are multiples of $\mathtt{M}$.  Hence, $\Phi_6$ in \eqref{phi6deff} is $\frac{2\pi}{\mathtt{M}}$-translation invariance preserving, therefore so is $\mathcal{L}^6$.
\end{proof}

From now on, we aim to prove the lemmas used to prove the above proposition. We see from the expansions of $\mathtt{m}_{\alpha}$ and $\mathfrak{m}_{\le 0}$ from Proposition~\ref{prop_time_rep_92} and Proposition~\ref{rlaqkqdpsanfsdfmf}, and \eqref{idernfisdsd2} that 
\begin{align*}
& \mathtt{m}_\alpha m_{1,\alpha}(\xi) + \frac{T_\alpha}4 + \mathfrak{m}_{\le 0}(\omega,\xi) + \epsilon^2\mathfrak{m}_{\mathfrak{b}}(\omega,\xi)  \\
& := \left(-\frac{1}{2}m_{1,\alpha}(\xi) + \frac{T_\alpha}{4}\right) + \epsilon^2\mathfrak{F}_{\alpha-1}(\omega,\xi) + (\mathtt{m}_{\alpha,2}m_{1,\alpha}(\xi) + \mathfrak{m}_{\le 0 , 2}).
\end{align*}
We denote
\begin{equation}\label{bdefsd2zpsdsayno}
\begin{aligned}
\mathcal{B}_{0,*}&:=\partial_x \Pi_{S^\perp}Op^W\left(-\frac{1}{2}m_{1,\alpha}(\xi) + \frac{T_\alpha}{4}\right),\quad \mathcal{B}_{1,*}:=\partial_x \Pi_{S^\perp}Op^W\left(\mathfrak{F}_{\alpha-1}(\omega,\xi)\right)\\
\mathcal{B}_{2,*}&:=\partial_x \Pi_{S^\perp}Op^W\left(\mathtt{m}_{\alpha,2}m_{1,\alpha}(\xi) + \mathfrak{m}_{\le 0 ,2}\right),\quad \mathcal{B}_{3,*}:=\partial_x \Pi_{S^\perp}Op^W\left(\mathfrak{b}_0\right),
\end{aligned}
\end{equation}
so that 
\begin{align}\label{j2sdwouldyouracherlsds}
\partial_x\Pi_{S^\perp}Op^W\left(\mathtt{m}_\alpha m_{1,\alpha}(\xi) + \frac{T_\alpha}4 + \mathfrak{m}_{\le 0}(\omega,\xi) + \epsilon^2\mathfrak{m}_{\mathfrak{b}}(\omega,\xi)\right) = \mathcal{B}_{0,*} + \epsilon^2 \mathcal{B}_{1,*} + \mathcal{B}_{2,*}.
\end{align}
Therefore, \eqref{l5idenitsd} can be written as 
\begin{align}\label{l6sd2syoulikedtodoingit}
\mathcal{L}^5 = \D_\omega -\mathcal{B}_{0,*} - \epsilon^2\mathcal{B}_{3,*} - \epsilon^2 \mathcal{B}_{1,*} -\mathcal{B}_{2,*} + \mathcal{R}_{\le 3}.
\end{align}

\subsubsection{Computing the conjugation}
As in subsection~\ref{firststeplbnf}, we expand $\Phi_6,\Phi_6^{-1}$ as
\begin{equation}\label{expansionofphi_6}
\begin{aligned}
\Phi_6 &=:I + \epsilon^2\mathcal{A}_{1,*} =: I + \epsilon^2 A_{1,*} +\epsilon^4\mathcal{A}_{2,*},\\
\Phi_6^{-1} &=: I + \epsilon^2\tilde{\mathcal{A}}_{1,*} =:I - \epsilon^2 A_{1,*} + \epsilon^4\tilde{\mathcal{A}}_{2,*}.
\end{aligned}
\end{equation}

\begin{lemma}\label{skdsdwdsd}
We have that
\begin{equation}\label{conjsd2ss20sdd}
\begin{aligned}
\Phi_6^{-1}\mathcal{L}^5\Phi_6 &= \D_\omega - \Pi_{S^\perp}\partial_x Op^W\left(\mathtt{m}_\alpha m_{1,\alpha}(\xi) + \left(\frac{T_\alpha}4 + \mathfrak{m}_{\le 0}(\omega,\xi) + \epsilon^2\mathfrak{m}_{\mathfrak{b}}(\omega,\xi) \right) \right) \\
& \ +\epsilon^2\left( \D_{\overline{\omega}}(A_{1,*}) - [\mathcal{B}_{0,*},A_{1,*}] - \mathcal{B}_{3,*}\right) +\mathcal{W}_*  + \mathcal{Q}_{1,*},
\end{aligned}
\end{equation}
where
\begin{equation}\label{wobsdand1q}
\begin{aligned}
\mathcal{W}_{*} &:=\Phi_6^{-1} \partial_x \Pi_{S^\perp}W_2\Phi_6,\\
\mathcal{Q}_{1,*} &:=\epsilon^2\left( \D_{{\omega}}(A_{1,*}) - \D_{\overline{\omega}}(A_{1,*})\right) + \epsilon^4 (\D_\omega(\mathcal{A}_{2,*}) + \tilde{\mathcal{A}}_{1,*}(\D_\omega(A_{1,*}) \\ & \ + \epsilon^2 \D_\omega(\mathcal{A}_{2,*})))
  - \epsilon^4 (\tilde{\mathcal{A}}_{2,*}\mathcal{B}_{0,*} + \tilde{\mathcal{A}}_{1,*}\mathcal{B}_{0,*}A_{1,*} + \Phi_6^{-1}\mathcal{B}_{0,*}\mathcal{A}_{2,*})\\
& \ - \epsilon^4\left( \tilde{\mathcal{A}}_{1,*}\mathcal{B}_{3,*}  + \Phi_6^{-1}\mathcal{B}_{3,*}\mathcal{A}_{1,*}\right) - \epsilon^4\left(\tilde{\mathcal{A}}_{1,*}\mathcal{B}_{1,*} + \Phi_6^{-1}\mathcal{B}_{1,*}\mathcal{A}_{1,*} \right)\\
& \ - \epsilon^2\left(\tilde{\mathcal{A}}_{1,*}\mathcal{B}_{2,*} + \Phi_6^{-1}\mathcal{B}_{2,*}\mathcal{A}_{1,*} \right)+  \Phi_{6}^{-1}\Pi_{S^\perp}\partial_x Op^W(\mathfrak{r}_{-2,\le 3,*})\Phi_6 + \Phi_6^{-1} \mathcal{R}_{5}\Phi_6^{-1} 
\end{aligned}
\end{equation}
\end{lemma}
\begin{proof}
Using \eqref{l6sd2syoulikedtodoingit}, we see that 
\begin{align}
\Phi_6^{-1}\mathcal{L}^{5}\Phi_6 & = \Phi_6^{-1}\D_\omega\Phi_6 - \Phi_6^{-1}\mathcal{B}_{0,*}\Phi_6 - \epsilon^2\Phi_6^{-1}\mathcal{B}_{3,*}\Phi_6 \nonumber \\
& -  \epsilon^2\Phi_6^{-1}\mathcal{B}_{1,*}\Phi_6 -\Phi_6^{-1}\mathcal{B}_{2,*}\Phi_6  + \Phi_6^{-1}\mathcal{R}_{\le 3}\Phi_6.\label{rlarhrsdw2}
\end{align}

\textbf{Conjugation with $\D_\omega.$}
Using \eqref{expansionofphi_6}, we have that
\begin{align}
\Phi_6^{-1}\D_\omega \Phi_6 &= \Phi_6^{-1}\left(\D_\omega(\Phi_6) + \Phi_6\D_\omega \right) \nonumber \\
& = \D_\omega + \Phi_6^{-1}(\epsilon^2 \D_\omega(A_{1,*}) + \epsilon^4 \D_\omega(\mathcal{A}_{2,*}))\nonumber\\
& = \D_\omega + \epsilon^2 \D_\omega(A_{1,*}) \nonumber \\
& + (\epsilon^4 \D_\omega(\mathcal{A}_{2,*}) + \epsilon^2\tilde{\mathcal{A}}_{1,*}(\epsilon^2 \D_\omega(A_{1,*}) + \epsilon^4 \D_\omega(\mathcal{A}_{2,*})))\nonumber\\
& = \D_\omega + \epsilon^2 \D_{\overline{\omega}}(A_{1,*}) + \epsilon^2\left( \D_{{\omega}}(A_{1,*}) - \D_{\overline{\omega}}(A_{1,*})\right) \nonumber \\
& + \epsilon^4 (\D_\omega(\mathcal{A}_{2,*}) + \tilde{\mathcal{A}}_{1,*}(\D_\omega(A_{1,*}) + \epsilon^2 \D_\omega(\mathcal{A}_{2,*}))).\label{youjustdidalittledance}
\end{align}

\textbf{Conjugation with $\mathcal{B}_{0,*}.$}
Again, \eqref{expansionofphi_6} gives us that
\begin{align}
\Phi_6^{-1} \mathcal{B}_{0,*}\Phi_6 &= \Phi_6^{-1}(\mathcal{B}_{0,*} + \epsilon^2 \mathcal{B}_{0,*}A_{1,*} + \epsilon^4\mathcal{A}_{2,*})\nonumber\\
& = \mathcal{B}_{0,*}  - \epsilon^2A_1\mathcal{B}_{0,*} + \epsilon^4\tilde{\mathcal{A}}_{2,*}\mathcal{B}_{0,*} + \epsilon^2\mathcal{B}_{0,*}A_1 \nonumber \\
& + \epsilon^4\tilde{\mathcal{A}}_{1,*}\mathcal{B}_{0,*}A_{1,*}+ \epsilon^4\Phi_6^{-1}\mathcal{A}_{2,*}\nonumber\\
& = \mathcal{B}_{0,*} + \epsilon^2[\mathcal{B}_{0,*},A_{1,*}] + \epsilon^4 (\tilde{\mathcal{A}}_{2,*}\mathcal{B}_{0,*} + \tilde{\mathcal{A}}_{1,*}\mathcal{B}_{0,*}A_{1,*} + \Phi_6^{-1}\mathcal{B}_{0,*}\mathcal{A}_{2,*}).\label{rlaqkqdpqkqakfdkajrwk}
\end{align}

\textbf{Conjugation with $\mathcal{B}_{3,*}.$}
Similarly, we have
\begin{align}
\epsilon^2\Phi_6^{-1}\mathcal{B}_{3,*}\Phi_6 & = \epsilon^2\Phi_6^{-1}(\mathcal{B}_{3,*} + \epsilon^2\mathcal{B}_{3,*}\mathcal{A}_{1,*})\nonumber\\
& = \epsilon^2 \mathcal{B}_{3,*} + \epsilon^4\left( \tilde{\mathcal{A}}_{1,*}\mathcal{B}_{3,*}  + \Phi_6^{-1}\mathcal{B}_{3,*}\mathcal{A}_{1,*}\right).
\end{align}

\textbf{Conjugation with $\mathcal{B}_{1,*}.$}
Similarly, we have
\begin{align}
\epsilon^2 \Phi_6^{-1}\mathcal{B}_{1,*}\Phi_6 = \epsilon^2 \mathcal{B}_{1,*} + \epsilon^4\left(\tilde{\mathcal{A}}_{1,*}\mathcal{B}_{1,*} + \Phi_6^{-1}\mathcal{B}_{1,*}\mathcal{A}_{1,*} \right).
\end{align}

\textbf{Conjugation with $\mathcal{B}_{2,*}.$}
Similarly, we have
\begin{align}\label{sdsrksdsdwsdd}
\Phi_6^{-1}\mathcal{B}_{2,*}\Phi_6^{-1} =  \mathcal{B}_{2,*} + \epsilon^2\left(\tilde{\mathcal{A}}_{1,*}\mathcal{B}_{2,*} + \Phi_6^{-1}\mathcal{B}_{2,*}\mathcal{A}_{1,*} \right)
\end{align}

\textbf{Conjugation with $\mathcal{R}_{\le 3}.$}
Recalling $\mathcal{R}_{\le 3}$ from \eqref{cdsddef}, we have
\begin{align}\label{jjsd2k4dfdf23}
\Phi_6^{-1}\mathcal{R}_{\le 3}\Phi_6 &= \Phi_6^{-1} \partial_x \Pi_{S^\perp}W_2\Phi_6 +  \Phi_{6}^{-1}\Pi_{S^\perp}\partial_x Op^W(\mathfrak{r}_{-2,\le 3,*})\Phi_6 + \Phi_6^{-1} \mathcal{R}_{5}\Phi_6^{-1}.
\end{align}

Plugging each conjugation into  \eqref{rlarhrsdw2}, we obtain that
\begin{align*}
\mathcal{L}^{5} & = \D_\omega - (\mathcal{B}_{0,*} + \epsilon^2 \mathcal{B}_{1,*} + \mathcal{B}_{2,*}) \\
& + \epsilon^2 \left( \D_{\overline{\omega}}(A_{1,*}) - [\mathcal{B}_{0,*},A_{1,*}] - \mathcal{B}_{4,*}\right) + \mathcal{Q}_{1,*} + \mathcal{W}_{*},
\end{align*}
where $\mathcal{Q}_{1,*}$ and $\mathcal{W}_{*}$ are as in \eqref{wobsdand1q}.  Recalling \eqref{j2sdwouldyouracherlsds}, we obtain \eqref{conjsd2ss20sdd}.
\end{proof}

\subsubsection{Choice of $\rho_2$} In view of \eqref{conjsd2ss20sdd}, we will choose $\rho_2$ in \eqref{phi6deff} so that 
\begin{align}\label{whihsd2xxx2}
\D_{\overline{\omega}}(A_{1,*}) - [\mathcal{B}_{0,*},A_{1,*}] - \mathcal{B}_{3,*} = 0,
\end{align}
\begin{proposition}\label{choiceofrho2common}
There exists $\rho_2=\mathfrak{B}_2^{-1-\alpha}$ such that
\begin{enumerate}[label=(\arabic*)]
\item \label{LBNF_ite2m1} With $A_{1,*}$ defined as in \eqref{phi6deff}, we have that
\begin{align}\label{sjrlaclWlsk12dsd}
\D_{\overline{\omega}}(A_{1,*}) - [\mathcal{B}_{0,*},A_{1,*}] - \mathcal{B}_{3,*} = 0. 
\end{align}
\item \label{LBNF_ite2m3}$\rho_2$ is a real-valued and reversibility preserving symbol. Furthermore, $\partial_xOp^W(\rho_2)$ is a real operator.
\end{enumerate}
\end{proposition}
\begin{proof}
Using \eqref{phi6deff} and \eqref{bdefsd2zpsdsayno}, we have that the above equation is equivalent to 
\begin{align}\label{whihsd2xxx22}
\Pi_{S^\perp}\left(Op^W(\D_{\overline{\omega}}(\rho_2)) - [Op^W\left(-\frac{1}{2}m_{1,\alpha}(\xi) + \frac{T_\alpha}4\right),Op^W(\rho_2)]_x - Op^W(\mathfrak{b}_0)\right) = 0.
\end{align}
We aim to choose $\rho_2$ satisfying  \eqref{whihsd2xxx22}. First, it follows from \ref{sdsdsd2jsds} of  Proposition~\ref{sd2sdsd} that $\mathfrak{b}_0\in \mathfrak{B}^{-2}_2$, hence, Definition~\ref{buildingbobo} tells us that there exists a smooth symbol $C_{j_{k_1}j_{k_2}}(\xi)\in \mathcal{S}^{-2}$ such that
\begin{align}\label{bo2sdsd}
\mathfrak{b}_0(\varphi,x,\xi) = \sum_{\substack{j_{k_1},j_{k_2}\in S, \\ j_{k_1}+j_{k_2} \ne 0}}C_{j_{k_1},j_{k_2}}(\xi)\overline{v}_{j_{k_1}}\overline{v}_{j_{k_2}}e^{\ii \left(\mathtt{l}(j_{k_1})+\mathtt{l}(j_{k_2})\right)\cdot \varphi + \ii (j_{k_1}+j_{k_2})x}.
\end{align}
As in the proof of Proposition~\ref{rhochosedefine}, let us choose, using \eqref{usual_weyl}, $\mathfrak{b}_{0,s}\in \mathfrak{B}^{-2}_2$ so that $Op^W(\mathfrak{b}_0)=Op(\mathfrak{b}_{0,s})$. That is,
\begin{align}\label{tqtxxx}
\mathfrak{b}_{0,s}(\varphi,x,\xi) = \sum_{\substack{j_{k_1},j_{k_2}\in S, \\ j_{k_1}+j_{k_2} \ne 0}}\underbrace{C_{j_{k_1},j_{k_2}}\left(\xi+\frac{j_{k_1}+j_{k_2}}{2}\right)}_{=:C_{j_{k_1},j_{k_2},s}(\xi)}\overline{v}_{j_{k_1}}\overline{v}_{j_{k_2}}e^{\ii \left(\mathtt{l}(j_{k_1})+\mathtt{l}(j_{k_2})\right)\cdot \varphi + \ii (j_{k_1}+j_{k_2})x}
\end{align}
Again, by Definition~\ref{buildingbobo}, $\xi\mapsto C_{j_{k_1},j_{k_2},s}(\xi)$ is a Fourier multiplier in $\mathcal{S}^{-2}$, that does not depend on $i_0,\omega,\varphi,x,\epsilon$. Therefore, in view of \eqref{whihsd2xxx22}, it suffices to find $\rho_{2,s}\in \mathfrak{B}^{-1-\alpha}_2$ such that 
\begin{align}\label{equivsd1xx}
\Pi_{S^\perp}\left(Op(\D_{\overline{\omega}}\rho_{2,s}) - [Op\left(-\frac{1}{2}m_{1,\alpha}(\xi) + \frac{T_\alpha}4\right),Op(\rho_{2,s}))]_x - Op(\mathfrak{b}_{0,s}) \right)= 0.
\end{align}
Using the notation in \eqref{fourier_ssd2s}, \eqref{tqtxxx} tells us that it suffices to find $\rho_{2,s}$ such that (see the proof of Proposition~\ref{rhochosedefine})  for all $j\in S^\perp$,
\begin{align}\label{njjxcxx2xcx}
\widehat{\rho}_{2,s}^{\varphi,x}(l,k,j) =
\begin{cases}
\frac{\widehat{\mathfrak{b}}_{0,s}^{\varphi,x}(l,k,j)}{\ii\delta_{ljk}}, & \text{if $(l,k)\in \mathcal{J}_{S}$},\\
 0, & \text{otherwise,}
\end{cases}
\end{align}
where
\begin{align}
\mathcal{J}_S&:=\left\{(\mathtt{l}(j_{k_1}+j_{k_2}),j_{k_1}+j_{k_2}):  j_{k_1},j_{k_2}\in S,\ j_{k_1}+j_{k_2}\ne 0, \right. \nonumber \\
&\left. \qquad \qquad \qquad \qquad \qquad \qquad \ \ j_{k_1}+j_{k_2} + j \in S^\perp \text{ for all $j\in S^\perp$} \right\},\label{delsdxcmasd}\\
\delta_{ljk}&:=\left(\overline{\omega}\cdot l - \left(\left(-\frac{1}{2}\lambda_\alpha(k+j)-  (k+j) \frac{T_\alpha}4\right) - \left( -\frac{1}{2}\lambda_{\alpha}(j)-  j\frac{T_\alpha}4\right)\right)\right).\label{delsdxcmasd2}
\end{align}
 Plugging  $(l,k)=(\mathtt{l}(j_{k_1}+j_{k_2}),j_{k_1}+j_{k_2})$ in \eqref{njjxcxx2xcx}, we see that the denominator in the left-hand side can be written as (see \eqref{linear_frequency_aa} and \eqref{def_lll} for $\overline{\omega}$ and $\mathtt{l}(j_{k_1}+j_{k_2}$)),
 \begin{align*}
\frac{1}{\delta_{ljk}}=\frac{2}{\left(\lambda_{\alpha}(j+j_{k_1}+j_{k_2})-\lambda_\alpha(j) - \lambda_\alpha(j_{k_1}+j_{k_2}) \right)}=-2\kappa_{-(j_{k_1}+j_{k_2}),1-\alpha}(j),
 \end{align*}
 where the last equality follows from \eqref{rjjsdsdjxcxcd}.
 Therefore, in view of \eqref{tqtxxx} and Lemma~\ref{kaoxixcc2xc}, we can choose
 \begin{align}\label{rho2sdsd}
 \rho_{2,s}(\varphi,x,\xi)& :=\sum_{\substack{j_{k_1},j_{k_2}\in S, \\ j_{k_1}+j_{k_2} \ne 0}}2\ii C_{j_{k_1},j_{k_2},s}(\xi)\kappa_{-(j_{k_1}+j_{k_2}),1-\alpha}(\xi) \nonumber \\
 & \qquad \qquad \qquad \times \overline{v}_{j_{k_2}}e^{\ii \left(\mathtt{l}(j_{k_1})+\mathtt{l}(j_{k_2})\right)\cdot \varphi + \ii (j_{k_1}+j_{k_2})x}.
 \end{align}
Clearly, $\rho_{2,s}$ satisfies \eqref{njjxcxx2xcx}.  Since $C_{j_{k_1},j_{k_2},s}(\xi)\in \mathcal{S}^{-2}$ and $\kappa_{-(j_{k_1}+j_{k_2}),1-\alpha}\in \mathcal{S}^{1-\alpha}$, we have $\rho_{2,s}\in \mathfrak{B}^{-1-\alpha}_2$.
 For the proof of \ref{LBNF_ite2m3}, the proof follows in the exactly same way as in the proof of \ref{LBNF_item3} of Proposition~\ref{rhochosedefine}. Indeed, comparing \eqref{whihsd2xxx22} and \eqref{homological_eq21s2}, we see that $\mathfrak{b}_0$ is also a real-valued, reversible symbol, as proved in \ref{sdsdsd2jsds} of Proposition~\ref{sd2sdsd}. Also $Op^W(\mathfrak{b}_0)$ is a real operator, since so is $\mathcal{L}^5$ (see \eqref{l5idenitsd} and \ref{rlsd2ssss} of Proposition~\ref{linearstep17}).
\end{proof}

 \subsubsection{Analysis of $\Phi_6$}
 
  \begin{lemma}\label{sdsdsdsd222kdkehsdjs2}
For $k=1,2$, there exist $\tilde{\rho}_{k,*}(\tau)\in \mathcal{S}^{-k\alpha}$ such that  
\[
\Psi_6(\tau) - I = \epsilon^2 Op(\tilde{\rho}_{1,*}(\tau)) + \frac{1}{2}\epsilon^4 Op(\tilde{\rho}_{2,*}(\tau)) ,\] and $\tilde{\rho}_{k,*}$ satisfy
\begin{align}
\sup_{\tau\in [0,1]}|\tilde{\rho}_{k,*}(\tau)|^{\Lip(\gamma,\Omega_1)}_{-k\alpha,s,\eta_0}&\le_{\mathtt{pe},s,\eta_0} 1,\label{rmsidgkse22}\\
\sup_{\tau\in [0,1]}|d_i\tilde{\rho}_{k,*}(\tau)(i_0)[\ihat]|_{-k\alpha,s,\eta_0}&=0\label{rmsidgkse222}
\end{align}
 \end{lemma}
\begin{proof}
The proof is identical to Lemma~\ref{sdsdsdsd222kdkehsdjs}.
\end{proof} 

\begin{lemma}\label{fjsdinf6sd}
$\Phi_6 \Pi_{S^\perp} - \Pi_{S^\perp}\Psi_6 \Pi_{S^\perp}\in \mathfrak{R}(i_0)$. 
\end{lemma}
\begin{proof}
The proof is identical to Lemma~\ref{fjsdinfsd}.
\end{proof}

\begin{lemma}\label{align2sd2sd2sssrx26}
$\Phi^{\pm1} = \Phi_6,\Phi_6^{-1}$ satisfies 
\begin{equation}\label{transformation_es2stimatess_cla2sss2}
\begin{aligned}
&\rVert (\Phi^{\pm}-I) h \rVert_{s}^{\Lip(\gamma,\Omega_1)} \le_{\mathtt{pe},s} \epsilon^2 \rVert h \rVert^{\Lip(\gamma,\Omega_1)}_{s}\\
&\rVert d_i\Phi^{\pm}(i_0)h[\ihat] \rVert_{s} =0
\end{aligned}
\end{equation}
\end{lemma}
\begin{proof}
The proof can be done following the proof of Lemma~\ref{align2sd2sd2sssrx2}, word by word, using Lemma~\ref{sdsdsdsd222kdkehsdjs2}.
\end{proof}

 \begin{lemma}\label{symbsrerer2}
 $\Phi_6$ is a symplectic, reversibility preserving, real operator. 
 \end{lemma}
 \begin{proof}
 The proof is identical to Lemma~\ref{symbsrerer}, using \ref{LBNF_ite2m3} of Proposition~\ref{choiceofrho2common}.
 \end{proof}

 \subsubsection{Analysis of  $\mathcal{W}_*$} 
 
\begin{lemma}\label{fisd2sd002sd1}
  $\mathcal{W}_*$ in \eqref{wobsdand1q} can be written as $\mathcal{W}_* =\Pi_{S^\perp}\partial_x\Psi_{6}^{T} W_2 \Psi_6 + R_{6,1}$ for some finite dimensional operator $R_{6,1}$ of the form in \eqref{j2j2j2j2j2jssds}.
 \end{lemma}
 \begin{proof}
 The proof is identical to Lemma~\ref{fisd2sd002sd}.
 \end{proof}

 \subsubsection{Analysis of $\mathcal{Q}_{1,*}$}
 \begin{lemma}\label{lastbs2sd}
 $\mathcal{Q}_{1,*}$ in \eqref{wobsdand1q} is of the form:
 \begin{align*}
 Q_{1,*} = \Pi_{S^\perp}Op(\mathfrak{r}_{-2,\le 3, \sharp}) + R_{6,2},
 \end{align*}
 where the symbol  $\mathfrak{r}_{-2,\le 3, \sharp}$ satisfies
 \begin{align}
 &|\mathfrak{r}_{-2,\le3,\sharp}|^{\Lip(\gamma,\Omega_1)}_{\alpha-3,s,\eta_0} \le_{\mathtt{pe},s,\eta_0}  \gamma^{-1}\left(\epsilon^5 + \epsilon^3\rVert \mathfrak{I}_\delta\rVert^{\Lip(\gamma,\Omega_1)}_{s+\mu_0}\right),\label{estiatefors20sd2}\\
&|d_i\mathfrak{r}_{-2,\le3,\sharp}(i_0)[\ihat]|_{\alpha-3,s,\eta_0} \le_{\mathtt{pe},s,\eta_0} \epsilon^3\gamma^{-1}\left( \rVert \ihat \rVert_{s+\mu_0} + \rVert \mathfrak{I}_\delta\rVert_{s+\mu_0}\rVert \ihat \rVert_{s_0+\mu_0}\right),\label{estiatefors20sd12}
 \end{align}
 and $R_{6,2}$ is a finite dimensional operator of the form in \eqref{j2j2j2j2j2jssds}.
\end{lemma}
\begin{proof}
The proof can be completed following the proof of Lemma~\ref{ibelieveshedoes}, using Lemma~\ref{sdsdsdsd222kdkehsdjs2} and Lemma~\ref{align2sd2sd2sssrx26}.
\end{proof}

\section{Lipschitz tame estimates for the remainders}\label{taksd2tamesd}
 Before we obtain the full reduction of the linear operator, we first derive tame estimates of the remainder in $\mathcal{L}^6$ obtained in Proposition~\ref{linearstep172sd}. We denote
 \begin{align}\label{r0definsd2sd}
 \mathcal{R}_0:=\partial_x Op^W(\mathfrak{r}_{-2,\le3,\sharp}) + \partial_xW_3 + R_6[h].
 \end{align}
 In the following lemma, we recall $\mathtt{b}_0\in \mathbb{N}$ is fixed as in \eqref{nsjdjwdsdnsde}.
 \begin{lemma}\label{tamesd2sd}
Given $\mathtt{S}\gg s_0$, we have the following:
\begin{enumerate}[label=(\arabic*)]
\item \label{tameestsd1sd} For $s\in [s_0,\mathtt{S}]$, and $\vec{\mathtt{b}}\in \mathbb{Z}^\nu$ such that $|\mathtt{b}|\in [0,s_0+\mathtt{b}_0]$, the operators  $\partial_{\varphi}^{\vec{\mathtt{b}}}(\mathcal{R}_0),[\partial_{\varphi}^{\vec{\mathtt{b}}}(\mathcal{R}_0),\partial_x]$ are Lip-$0$-tame operators and 
\begin{align}\label{dsd2sdsd}
\mathfrak{M}^\gamma_{\partial_{\varphi}^{\vec{\mathtt{b}}}(\mathcal{R}_0)}(0,s),\mathfrak{M}^\gamma_{[\partial_{\varphi}^{\vec{\mathtt{b}}}(\mathcal{R}_0),\partial_x]}(0,s)&\le_{\mathtt{pe},\mathtt{S}}\gamma^{-1}\left(\epsilon^5 + \epsilon^3\rVert \mathfrak{I}_\delta\rVert^{\Lip(\gamma,\Omega_1)}_{s+\mu_0}\right).
\end{align}
\item \label{tameestsd1sd2}  For $s\in [s_0,\mathtt{S}]$, and $\vec{\mathtt{b}}\in \mathbb{Z}^\nu$ such that $|\mathtt{b}|\in [0,s_0+\mathtt{b}_0]$, the operators  $\partial_{\varphi}^{\vec{\mathtt{b}}}(d_i(\mathcal{R}_0)(i_0)[\ihat]),[\partial_{\varphi}^{\vec{\mathtt{b}}}(d_i(\mathcal{R}_0)(i_0)[\ihat]),\partial_x]$ are $0$-tame operators and
\begin{align} \nonumber
\sup_{ \vec{\mathtt{b}}\in\mathbb{N}_0^\nu,\ |\vec{\mathtt{b}}|  \le s_0 + \mathtt{b}_0}\mathfrak{M}_{\mathcal{Q}_{\vec{\mathtt{b}}}}(0,s_0) & \le_{\mathtt{pe},\mathtt{S}}\epsilon^3\gamma^{-1}\rVert\ihat\rVert_{s_0+\mu_0}, \\
\mathcal{Q}_{\vec{\mathtt{b}}}& \in \left\{ \partial_{\varphi}^{\vec{\mathtt{b}}}(d_i(\mathcal{R}_0)(i_0)[\ihat]),[\partial_{\varphi}^{\vec{\mathtt{b}}}(d_i(\mathcal{R}_0)(i_0)[\ihat]),\partial_x]\right\}.\label{desiresithappensd}
\end{align}
\end{enumerate}
 \end{lemma}
 \begin{proof}
 \vspace{0.5\baselineskip}\noindent\textit{Proof of \ref{tameestsd1sd}.} We claim that 
 \begin{align}
 \mathfrak{M}^\gamma_{\partial_x\partial_{\varphi}^{\vec{\mathtt{b}}}(\mathcal{R}_0)}(0,s)&\le_{\mathtt{pe},\mathtt{S}}\gamma^{-1}\left(\epsilon^5 + \epsilon^3\rVert \mathfrak{I}_\delta\rVert^{\Lip(\gamma,\Omega_1)}_{s+\mu_0}\right) , \label{tamealmostover1} \\
 \mathfrak{M}^\gamma_{\partial_{\varphi}^{\vec{\mathtt{b}}}(\mathcal{R}_0)\partial_x}(0,s)&\le_{\mathtt{pe},\mathtt{S}}\gamma^{-1}\left(\epsilon^5 + \epsilon^3\rVert \mathfrak{I}_\delta\rVert^{\Lip(\gamma,\Omega_1)}_{s+\mu_0}\right).\label{tamealmostover2}
 \end{align}
 Let us assume for a moment that \eqref{tamealmostover1} and \eqref{tamealmostover2} hold. Since $\partial_x^{-1}$ is clearly a Lip-$0$-tame operator, the tame estimates for $\partial_{\varphi}^{\vec{b}}(\mathcal{R}_0) = \partial_x^{-1}(\partial_x\partial_{\varphi}^{\vec{b}}(\mathcal{R}_0))$ follow from \eqref{tamealmostover1} and Lemma~\ref{consdsdpcosdldsx}. Furthermore, it follows immediately from Definition~\ref{lip_tame} that, for two operators $A,B$,  $\mathfrak{M}^\gamma_{A+B}(0,s)\le \mathfrak{M}^\gamma_{A}(0,s) + \mathfrak{M}^\gamma_{B}(0,s)$. Hence \eqref{tamealmostover1} and \eqref{tamealmostover2} give us the tame estimates \eqref{dsd2sdsd} for $[\partial_{\varphi}^{\vec{\mathtt{b}}}(\mathcal{R}_0),\partial_x]$. For the proof of \eqref{tamealmostover1} and \eqref{tamealmostover2}, we prove \eqref{tamealmostover1} only, since the other estimate can be proved in the same way.
 
 In order to show \eqref{tamealmostover1}, we need to show that, in view of \eqref{r0definsd2sd},
 \begin{align}\nonumber
\mathfrak{M}^\gamma_{\partial_x^2Op^W(\partial_{\varphi}^{\vec{\mathtt{b}}}(\mathfrak{r}_{-2,\le3,\sharp}))}(0,s), & \ \mathfrak{M}^\gamma_{\partial_x\partial_{\varphi}^{\vec{\mathtt{b}}}(R_6)}(0,s), \ \mathfrak{M}^\gamma_{\partial_x^2\partial_{\varphi}^{\vec{\mathtt{b}}}(W_3) }(0,s) \\
 & \le_{\mathtt{pe},s} \gamma^{-1}\left(\epsilon^5 + \epsilon^3\rVert \mathfrak{I}_\delta\rVert^{\Lip(\gamma,\Omega_1)}_{s+\mu_0}\right).
 \label{rjsd2zxcxczc}
 \end{align}
 For $\partial_x^2Op^W(\partial_{\varphi}^{\vec{\mathtt{b}}}(\mathfrak{r}_{-2,\le3,\sharp}))$, it follows from \eqref{sdsdtm2sdsdddndkQk} and \ref{compositsd} of Lemma~\ref{compandkskd2sd} that we can find a symbol $r_0\in \mathcal{S}^0$ such that $\partial_x^2 Op^W(\partial_{\varphi}^{\vec{\mathtt{b}}}(\mathfrak{r}_{-2,\le3,\sharp})) = Op(r_0)$ and $r_0$ satisfies (Recalling Remark~\ref{eta0sdfix}),
\begin{align}
&|r_{0}|^{\Lip(\gamma,\Omega_1)}_{0,s,0} \le_{s,\mathtt{pe}}  \gamma^{-1}\left(\epsilon^5 + \epsilon^3\rVert \mathfrak{I}_\delta\rVert^{\Lip(\gamma,\Omega_1)}_{s+\mu_0}\right),\label{talkedtoyou}\\
&|d_ir_0(i_0)[\ihat]|_{0,s,0} \le_{s,\mathtt{pe}} \epsilon^3\gamma^{-1}\left( \rVert \ihat \rVert^{\Lip(\gamma,\Omega_1)}_{s+\mu_0} + \rVert \mathfrak{I}_\delta\rVert_{s+\mu_0}\rVert \ihat \rVert_{s_0+\mu_0}\right).\label{talkedtoyou2}
\end{align}
 Hence the desired tame estimates follow from \eqref{talkedtoyou} and Lemma~\ref{symbosldtame}.
 
  For the operator $\partial_x\partial_{\varphi}^{\vec{\mathtt{b}}}(R_6)$, it follows from \ref{item262gf222} of Proposition~\ref{linearstep172sd} that we can pick $g_j(\tau),\chi_j(\tau)$  such that 
  \[
   R_6[h] = \sum_{|j|\le C} \int_0^1(h,g_j(\tau))_{L^2}\chi_j(\tau)d\tau,
   \]
    and satisfy
  \begin{align}
 &\sup_{\tau\in[0,1]}\rVert g(\tau) \rVert^{\Lip(\gamma,\Omega_1)}_s\rVert \chi(\tau)\rVert^{\Lip(\gamma,\Omega_1)}_{s_0}+ \rVert g(\tau) \rVert^{\Lip(\gamma,\Omega_1)}_{s_0}\rVert \chi(\tau)\rVert^{\Lip(\gamma,\Omega_1)}_{s}\nonumber \\
 & \le_{\mathtt{pe},s}\epsilon^3+\epsilon^{2}\rVert \mathfrak{I}_\delta \rVert_{s+\mu_0}^{\Lip(\gamma,\Omega_1)},\label{gsdwsdgchies}\\
 &\sup_{\tau\in[0,1]}\left(\rVert d_ig(i_0)[\ihat] \rVert_{s_0}\rVert \chi(i_0) \rVert_s+\rVert d_ig(i_0)[\ihat] \rVert_s\rVert \chi_j(i_0) \rVert_{s_0}\right. \\
 & \qquad \left.+\rVert d_i\chi(i_0)[\ihat] \rVert_{s_0}\rVert g(i_0) \rVert_s+\rVert d_i\chi(i_0)[\ihat] \rVert_{s}\rVert g(i_0) \rVert_{s_0}\right)\nonumber\\
& \quad \le_{\mathtt{pe},s} \epsilon^2\rVert \ihat \rVert_{s+\mu_0} + \epsilon^{2b-1}\rVert \mathfrak{I}_\delta \rVert_{s+\mu_0}\rVert \ihat \rVert_{s+\mu_0}.\label{gsdwsdgchies2}
  \end{align}
  It is clear that
  \[
  \partial_x\partial_\varphi^{\vec{\mathtt{b}}}(R_6)h = \int_0^1 \sum_{\vec{b}_1+\vec{b}_2=\vec{\mathtt{b}}}\sum_{|j|\le C}(h,\underbrace{\partial_\varphi^{\vec{b}_1}g_j(\tau)}_{=:g_{b_1,b_2,j}(\tau)})_{L^2_x}\underbrace{C_{\vec{b}_1,\vec{b}_2} \partial_x\partial_{\varphi}^{\vec{b}_2}(\chi_j(\tau))}_{=:\chi_{b_1,b_2,j}(\tau)}d\tau
  \]
  Clearly, $g_{b_1,b_2,j}(\tau),\chi_{b_1,b_2,j}(\tau)$ satisfy the estimates in \eqref{gsdwsdgchies} with possibly larger $\mu_0$. Hence, the tame estimate \eqref{rjsd2zxcxczc} for $\partial_x\partial_{\varphi}^{\vec{\mathtt{b}}}(R_6)$ follows from Lemma~\ref{s_decay_finite}.
  
  Lastly, we deal with $W_3$. It follows from \ref{item262gf222} of Proposition~\ref{linearstep172sd} and \eqref{W2isdsd2} that
  \begin{equation}\label{jjdsdwosdsd}
  \begin{aligned}
  & \partial_x^2W_3 =\frac{1}{(\mathtt{N}_\alpha+1)!} \int_0^1\int_0^\tau \partial_x^2\left(\Psi_3(\tau-t)\Phi_4\Psi_5\Psi_6\right)^{T}Op^W(\mathfrak{q}_{\omega,\mathtt{N}_\alpha-2}\star \mathfrak{a}) \\
  & \qquad \qquad \qquad  \qquad \qquad \qquad \qquad \qquad \qquad \times \Psi_3(\tau-t)\Phi_4\Psi_5\Psi_6 t^{\mathtt{N}_\alpha -2}dtd\tau\\
  & \  -  \frac{1}{\mathtt{N}_\alpha!}\int_0^1 \partial_x^2\left(\Psi_3(1-t)\Phi_4\Psi_5\Psi_6\right)^{T}Op^W(\mathfrak{q}_{M_2,\mathtt{N}_\alpha}\star\mathfrak{a})\Psi_3(1-t)\Phi_4\Psi_5\Psi_6 t^{\mathtt{N}_\alpha}dt.
  \end{aligned}
  \end{equation}
 
 Note that it follows from \eqref{rlaclsdsd}, \eqref{hamisdsd}, \eqref{psi5defsseconddate} and \eqref{psi6defsseconddat2e} that the transformations $\Psi_3,\Phi_4,\Psi_5,\Psi_6$ are generated by the symbols $\mathfrak{a}\in \mathcal{S}^{1-\alpha},\ \mathfrak{p}_2\in \mathcal{S}^{1-\alpha},\ \epsilon\rho_1, \epsilon^2\rho_2\in\mathcal{S}^{-1-\alpha}\subset \mathcal{S}^{1-\alpha}$. Furthermore, since $\mathfrak{a}\in \mathfrak{S}^{1-\alpha}_{p_1}$ (see \eqref{2jsd232}), using Lemma~\ref{sizeofthesymbolsdsd} and also using Lemma~\ref{p2sdsdkcssd}, Proposition~\ref{rhochosedefine} and  Proposition~\ref{choiceofrho2common}  that for $
 \mathfrak{u}\in \left\{ \mathfrak{a},\mathfrak{p}_2,\epsilon \rho_1, \epsilon^2 \rho_2\right\},$
 we have the estimates
 \begin{align}
|\mathfrak{u}|^{\Lip(\gamma,\Omega_1)}_{1-\alpha,s,\eta_0}&\le_{\mathtt{pe},s,\eta_0} \epsilon(1+ \rVert \mathfrak{I}_\delta\rVert^{\Lip(\gamma,\Omega_1)}_{s+\mu_0}), \label{rmsidsdw2}\\
\rVert d_i\mathfrak{u}(i_0)[\ihat]\rVert_{1-\alpha,s,\eta_0}&\le_{\mathtt{pe},s,\eta_0} \epsilon^3\gamma^{-1}\left( \rVert \ihat \rVert_{s+\mu_0} + \rVert \mathfrak{I}_\delta\rVert_{s+\mu_0}\rVert \ihat \rVert_{s_0+\mu_0}\right).\label{rmsidsd22s}
\end{align}
especially, \eqref{size_assumption_3} tells us that  
\begin{align}\label{kgkscenlqsd}
|\mathfrak{u}|^{\Lip(\gamma,\Omega_1)}_{1-\alpha,s_0+\mu_0,\eta_0}\le_{\mathtt{pe},s,\eta_0} \epsilon.
\end{align}
  Therefore, we can apply Proposition~\ref{wjsd2sdsdesd} for sufficiently small $\epsilon>0$, depending on $\mathtt{S}$ (Recalling Remark~\ref{eta0sdfix}, we choose $\eta_0=\eta_0(\mathtt{b}_0)$ in \eqref{kgkscenlqsd}). Then, \ref{swlsWkgoqhwk} and \ref{sddwsdsdsd} of Proposition~\ref{wjsd2sdsdesd} tell us that
  \begin{align}\nonumber
  \sup_{|k|\le \mathtt{b}_0, |\vec{b}|\le |\vec{\mathtt{b}}|} &\mathfrak{M}^\gamma_{|D|^k \partial_{\varphi}^{\vec{{b}}}(\psi) |D|^{-k-|\vec{b}|(2-\alpha)}}(0,s), \\
  & \ \mathfrak{M}^\gamma_{|D|^{k-|\vec{b}|(2-\alpha)} \partial_{\varphi}^{\vec{{b}}}(\psi) |D|^{-k}}(0,s)\le_{\mathtt{pe},s} (1+ \rVert \mathfrak{I}_\delta\rVert^{\Lip(\gamma,\Omega_1)}_{s+\mu_0}),\label{finaltameksd2}
  \end{align}
  for $\psi\in \left\{ \Psi_3(\tau),\Phi_4,\Psi_5,\Psi_6\right\}$ for $\tau\in [0,1]$. Furthermore, it follows from Lemma~\ref{sjjdsdsymbolsd} and \ref{compositsd} of Lemma~\ref{compandkskd2sd} that denoting a symbol $\mathfrak{q}_{*}$
  \begin{align} \nonumber
Op^W(\mathfrak{q}_{*} ) \in &\left\{|D|^{(|\vec{b}_1|+1)(2-\alpha)+2}Op^W(\mathfrak{q}_\mathfrak{a})|D|^{(|\vec{b}_2|+1)(2-\alpha)}\right. \\
& \left.: \mathfrak{q}_\mathfrak{a}=\mathfrak{q}_{\omega,\mathtt{N}_\alpha-2}\star \mathfrak{a}\text{ or }\mathfrak{q}_{M_2,\mathtt{N}_\alpha}\star \mathfrak{a},\ \vec{b}_1+\vec{b}_2 = \vec{\mathtt{b}}\right\},\label{ericcartman}
  \end{align}
  we have that $\mathfrak{q}_*\in \mathcal{S}^{\mathtt{N}_\alpha(1-\alpha)+\left((|\vec{\mathtt{b}}| +2)(2-\alpha) + 2\right)}\subset \mathcal{S}^{0}$ (see \eqref{nsjdjwdsdnsde} and recall that $|\mathtt{b}|\le s_0+\mathtt{b}_0$), and 
\begin{align}
|\mathfrak{q}_{*} |^{\Lip(\gamma,\Omega_1)}_{0,s,0} & \le_{s,\mathtt{pe}} \epsilon^{\mathtt{N}_\alpha}(1 + \rVert\mathfrak{I}_\delta\rVert_{s + \mu_0}^{\Lip(\gamma,\Omega_1)}),\label{ssssss2020iwassds}\\
 |d_i(\mathfrak{q}_{*}(i_0)[\ihat]|_{0,s,0} & \le_{s,\mathtt{pe}} \epsilon^{\mathtt{N}_\alpha + 2}\gamma^{-1}\left( \rVert \ihat \rVert_{s+\mu_0} + \rVert \mathfrak{I}_\delta\rVert_{s+\mu_0}\rVert \ihat \rVert_{s_0+\mu_0}\right). \label{ssssss2020iwassds2}
\end{align}
   Therefore, taking $\partial_{\varphi}^{\vec{\mathtt{b}}}$ in \eqref{jjdsdwosdsd}, we see, using \eqref{finaltameksd2}, \eqref{ssssss2020iwassds} and Lemma~\ref{consdsdpcosdldsx}, that
  \[
 \mathfrak{M}^\gamma_{\partial_x^2\partial_\varphi^{\vec{\mathtt{b}}}(W_3)}\le_{\mathtt{pe},s} \epsilon^{\mathtt{N}_\alpha+1}(1 + \rVert\mathfrak{I}_\delta\rVert_{s + \mu_0}^{\Lip(\gamma,\Omega_1)}),
  \]  
  which gives us \eqref{rjsd2zxcxczc}, since $\mathtt{N}_\alpha \ge 2$, thanks to \eqref{nsjdjwdsdnsde}.
  
  \vspace{0.5\baselineskip}\noindent\textit{Proof of \ref{tameestsd1sd2}.} The proof of \ref{tameestsd1sd2} goes similarly as in \ref{tameestsd1sd}, hence we briefly explain the outline of the proof. As above, we claim that
   \begin{align}
 \mathfrak{M}_{\partial_x\partial_{\varphi}^{\vec{\mathtt{b}}}(d_i(\mathcal{R}_0)(i_0)[\ihat])}(0,s_0)&\le_{\mathtt{pe}}\epsilon^3\gamma^{-1}\rVert\ihat\rVert_{s_0+\mu_0}, \label{tamealmostover21} \\
 \mathfrak{M}_{\partial_{\varphi}^{\vec{\mathtt{b}}}(d_i(\mathcal{R}_0)(i_0)[\ihat])\partial_x}(0,s_0)&\le_{\mathtt{pe}}\epsilon^3\gamma^{-1}\rVert\ihat\rVert_{s_0+\mu_0}.\label{tamealmostover22}
 \end{align}
 Once we have the above estimates, then the similar argument described above, using that $\partial_x^{-1}$ is $0$-tame operator, gives us the desired result \eqref{desiresithappensd}. As mentioned in \ref{tameestsd1sd}, the proofs for \eqref{tamealmostover21} and \eqref{tamealmostover22} are similar. For \eqref{tamealmostover21}, in view of \eqref{r0definsd2sd}, we are led to prove
 \begin{align}
  \mathfrak{M}_{\partial_x^2Op^W(\partial_{\varphi}^{\vec{\mathtt{b}}}(d_i\mathfrak{r}_{-2,\le3,\sharp}(i_0)[\ihat]))}(0,s_0), & \ \mathfrak{M}_{\partial_x\partial_{\varphi}^{\vec{\mathtt{b}}}(d_iR_6(i_0)[\ihat])}(0,s_0), \nonumber \\
  & \ \mathfrak{M}_{\partial_x^2\partial_{\varphi}^{\vec{\mathtt{b}}}(d_iW_3(i_0)[\ihat]) }(0,s_0) \le_{\mathtt{pe}} \epsilon^3\gamma^{-1}\rVert\ihat\rVert_{s_0+\mu_0}.\label{rjsdsdsd}
 \end{align}
For each operator in above, one can follow the same argument as in the proof of \ref{tameestsd1sd},  using the estimates \eqref{talkedtoyou2}, \eqref{gsdwsdgchies2}, \eqref{rmsidsd22s} and \eqref{ssssss2020iwassds2}, with \eqref{size_assumption_3}.
 \end{proof}

 Recalling Remark~\ref{derv_loss} and Remark~\ref{Insteadofimbd}, we see that the loss of derivatives, $\mu_0$, arising in the reduction steps in Chapter~\ref{reduction} is finite number depending on only $\mathtt{p}$, fixed as in \eqref{parametersets1231}, assuming that $\epsilon$ is small enough depending on $s,\mathtt{pe}$. We summarize the results that we obtain through Sections~\ref{change_of_the_space_variables}-~\ref{linearbbss22} and Lemma~\ref{tamesd2sd}.
 
 \begin{proposition}\label{modulut2sosdtame}
There exists a constant 
 \begin{align}\label{constandusd2}
 \mu_{\mathtt{p},1} \ge 0, \text{ depending only on $\mathtt{p}$ in \eqref{parametersets1231}},
 \end{align}
 such that  if \eqref{size_assumption_2} is satisfied with $\mu =\mu_{\mathtt{p},1}$ and for some $\mathtt{C}>0$, then the followings hold:
$\mathcal{L}^6$ in \eqref{sdlinesdar2sd2ssdxcxcs2j} can be written as
 \begin{align}\label{semifinallinearoperator}
L_0:= \mathcal{L}^6 = \D_\omega -D_0 + \mathcal{R}_0,
 \end{align}
 where the operators $D_0$ and $\mathcal{R}_0$ satisfy
 \begin{enumerate}[label=(\arabic*)]
 \item \label{diagonalpart_1} $D_0$ is a diagonal operator between $H_{S^\perp}$ and $H_{S^\perp}$ such that 
 \begin{align}
 D_0 &= Op^W(d_0(\omega,\xi)),\quad d_0({\omega,\xi}):=\ii\left( \mathtt{m}_\alpha(\omega)\lambda_\alpha(\xi) +  \xi \frac{T_\alpha}4 + \mathfrak{m}_{1}(\omega,\xi) \right),\label{d0defsd2}\\
 \mathfrak{m}_{1} &:= \xi\left(\mathfrak{m}_{\le 0}(\omega,\xi) + \epsilon^2\mathfrak{m}_{\mathfrak{b}}(\xi)\right)\in \mathcal{S}^1, \label{d0defsd22}
 \end{align}
 and  $\mathtt{m}_{\alpha}=-\frac{1}{2}+\epsilon^2\mathtt{m}_{\alpha,1}+\mathtt{m}_{\alpha,2}$ satisfies \begin{align}
|\mathtt{m}_{\alpha,1}|^{\Lip(\gamma,\Omega_1)} & \le_{\mathtt{pe}} 1,\quad |\mathtt{m}_{\alpha,2}|^{\Lip(\gamma,\Omega_1)}\le_{\mathtt{pe}} \epsilon^{7-4b}, \nonumber \\
 |d_i\mathtt{m}_{\alpha,2}(i_0)[\ihat]| & \le_{\mathtt{pe}} \epsilon\rVert \ihat\rVert_{s_0+\mu_{\mathtt{p},1}},\label{nusd2sdsd}
 \end{align} and  $\mathfrak{m}_1(\omega,\xi)$ satisfies
 \begin{align}
 |\mathfrak{m}_1|^{\Lip(\gamma,\Omega_1)}_{1,0,2}&\le_{\mathtt{pe}}  \epsilon^2,\quad |d_i\mathfrak{m}_1(i_0)[\ihat]|_{1,0,2}\le_{\mathtt{pe}} \epsilon \rVert \ihat \rVert_{s_0+\mu_{\mathtt{p},1}}.\label{rejsmsdim1}
 \end{align}
 
  \item \label{remainderpart_1} Given $\mathtt{S}\gg s_0$, there exists $\epsilon_0=\epsilon_0(\mathtt{pe},\mathtt{S})$ such that if $\epsilon \in (0,\epsilon_0)$, then $\mathcal{R}_0$ is a Lip-$0$-modulo tame operator up $\mathtt{S}$ and $d_i\mathcal{R}_0(i_0)[\ihat]$ is a $0$-tame operator  and satisfies the following estimates:
  \begin{align}
  \sup_{\vec{b}\in \mathbb{N}_0^{\nu},\ |\vec{b}|\le \mathtt{b}_0 }\mathfrak{M}^{\sharp,\gamma}_{\partial_{\varphi}^{\vec{b}}(\mathcal{R}_0)}(0,s)&\le_{\mathtt{pe},s} \gamma^{-1}\left(\epsilon^5 + \epsilon^3\rVert \mathfrak{I}_0\rVert^{\Lip(\gamma,\Omega_1)}_{s+\mu_{\mathtt{p},1}}\right),\label{fromyourbag1}\\
  \sup_{ \vec{{b}}\in\mathbb{N}_0^\nu,\ |\vec{{b}}|\le \mathtt{b}_0}\mathfrak{M}^\sharp_{\partial_\varphi^{\vec{b}}(d_i\mathcal{R}_0(i_0)[\ihat])}(0,s_0)&\le_{\mathtt{pe}}\epsilon^3\gamma^{-1}\rVert\ihat\rVert_{s_0+\mu_{\mathtt{p},1}}, \label{fromyourbag2}\end{align}
  for all $s\in [s_0,\mathtt{S}]$.
  \item \label{reversiblitsd} $ D_0$ is a real, reversible operator. Also, $\mathcal{R}_0$ is a real and  reversible operator.

  \item \label{reversibility_presdsx1d} There exist reversibility preserving real operators $\Phi_{1-6,1},\Phi_{1-6,2}:H_{S^\perp}\mapsto H_{S^\perp}$ such that
 \begin{align}\label{lomegatol2}
 L_0=\Phi_{1-6,1}\mathcal{L}_\omega \Phi_{1-6,2},
 \end{align}
 and $\Phi:=\Phi_{1-6,1},\Phi_{1-6,2}$ satisfies (assuming $\epsilon \in (0, \epsilon_0(\mathtt{pe},\mathtt{S}$)))
 \begin{align}\label{phiestfinal1}
 \rVert (\Phi-I) h \rVert_{s}^{\Lip(\gamma,\Omega_1)} \le_{\mathtt{pe},s} \epsilon\left( \rVert h \rVert^{\Lip(\gamma,\Omega_1)}_{s+\mu_{\mathtt{p},1}} + \rVert \mathfrak{I}_0 \rVert^{\Lip(\gamma,\Omega_1)}_{s+\mu_{\mathtt{p},1}}\rVert h \rVert^{\Lip(\gamma,\Omega_1)}_{s_0+\mu_{\mathtt{p},1}}\right),
 \end{align}
 for all $s\in[s_0,\mathtt{S}]$.
 
\item\label{egpp2sd}  $\mathcal{R}_0$ is $\frac{2\pi}{\mathtt{M}}$-translation invariant and $\Phi_{1-6,1},\Phi_{1-6,2}$ are $\frac{2\pi}{\mathtt{M}}$-translation invariance preserving.
 \end{enumerate}
 \end{proposition}
 \begin{proof}
 Fix $\mu_{\mathtt{p},1}$ be the largest $\mu_0$, arising in Chapter~\ref{reduction}, which certainly depends on only $\mathtt{p}$. We simply rewrite $\mathcal{L}^6$ in \eqref{sdlinesdar2sd2ssdxcxcs2j} using $d_0$ and $\mathcal{R}_0$,  defined in \eqref{d0defsd2} and \eqref{r0definsd2sd} to obtain the expression \eqref{semifinallinearoperator}.
 
\vspace{0.5\baselineskip}\noindent\textit{Proof of \ref{diagonalpart_1}.} \eqref{nusd2sdsd} can be found in \eqref{qtildesd2dsd}. To see \eqref{rejsmsdim1},  notice that \ref{rjsdsd} of Proposition~\ref{sd2sdsd} tells us that  the symbol $\mathfrak{m}_{\mathfrak{b}}(\xi)$ is a Fourier multiplier that does not depend on $\omega$ or $i_0$, and it satisfies $\mathfrak{m}_{\mathfrak{b}}\in \mathcal{S}^{-2}$.
 Then, \eqref{rejsmsdim1} immediately follows from the definition of $\mathfrak{m}_1$ in \eqref{d0defsd22} with the estimates for $\mathfrak{m}_{\le 0}$ in  \ref{rhshdosd122} of Proposition~\ref{rlaqkqdpsanfsdfmf}.

  \vspace{0.5\baselineskip}\noindent\textit{Proof of \ref{remainderpart_1}.}
  \ref{remainderpart_1} follows from Lemma~\ref{tamesd2sd} and Lemma~\ref{modulsdo2ds}, replacing $\mathfrak{I}_\delta$ by $\mathfrak{I}_0$ using \eqref{isjkxcoksd1sd}. 
  
 \vspace{0.5\baselineskip}\noindent\textit{Proof of \ref{reversiblitsd}.} It follows from \ref{rlsd2ssss2ss} of Proposition~\ref{linearstep172sd} that   $\mathcal{L}^6$ is a  real and reversible operator, therefore so are $D_0$ and $\mathcal{R}_0$.

  \vspace{0.5\baselineskip}\noindent\textit{Proof of \ref{reversibility_presdsx1d}.}  We denote
    \[
    \Phi_{1-6,1}:=\Phi_6\Phi_5\Phi_4\Phi_3\frac{1}{\rho}\Phi_2\Phi_1,\quad \Phi_{1-6,2}:=\Phi_1\Phi_2\Phi_3\Phi_4\Phi_5\Phi_6.
    \]
    Then the expression \eqref{lomegatol2} follows from \eqref{linear_111111}, \eqref{sdlinear2sd2}, \eqref{sdlinesdar2sd11xx2}, \eqref{sdlinesdar2sd2sj}, \eqref{sdlinesdar2sd2s2j} and \eqref{sdlinesdar2sd2ssdxcxcs2j}, with  \eqref{semifinallinearoperator}. The estimates \eqref{phiestfinal1} follows from Lemma~\ref{align2sd2sd24s2}, \ref{align2sd2sd22}, \ref{align2sd2sd2}, \ref{align2sd2sd2sx2}, \ref{align2sd2sd2sssrx2}, \ref{align2sd2sd2sssrx26}, where we can replace $\mathfrak{I}_\delta$ by $\mathfrak{I}_0$, using \eqref{isjkxcoksd1sd}.

  \vspace{0.5\baselineskip}\noindent\textit{Proof of \ref{egpp2sd}.}   We see from \ref{egppsd} of Proposition~\ref{linearstep172sd} that $\mathcal{L}^6$ is $\frac{2\pi}{\mathtt{M}}$-translation invariance preserving, therefore so is $\mathcal{R}_0$. $\frac{2\pi}{\mathtt{M}}$-translation invariance preserving property of $\Phi_{1-6,1},\Phi_{1-6,2}$ follows from (5) of Proposition~\ref{toohard_2_3}, \ref{prop_time_rep_92}, \ref{induction_egorov}, \ref{rlaqkqdpsanfsdfmf}, \ref{linearstep17}, \ref{linearstep172sd}.
 \end{proof}
  In view of Proposition~\ref{normal_inversion}, we  need to find the inverse image of $Y_{\perp}\cap X_{\mathtt{M}}$ in $X_{\perp}\cap X_{\mathtt{M}}$. Especially, we can restrict the operator $L_0$ in \eqref{semifinallinearoperator} to the space of functions whose Fourier modes are supported on only $\mathtt{M}$-multiples. More precisely, we consider $L_0$ as a linear operator between $H_{S^\perp}\cap X_{\mathtt{M}}$. Therefore, without loss of generality, we have that (using \eqref{mnorlam}),
 \begin{align}\label{rkjspspdwjskds}
(\mathcal{R}_0)^{j_1}_{j_2}(l) = 0,\text{ if $j_1,j_2\notin \left\{ \mathtt{M}j\in \mathbb{Z}\cap S^\perp: j\in \mathbb{Z}\right\}=S_{\mathtt{M}}^\perp$.}
 \end{align}

 \section{KAM reducibility and inversion of $\mathcal{L}_\omega$}\label{rkppsdsdwdwdx1kaksd}
 \subsection{The homological equation}
\index{KAM reducibility}
 \begin{lemma}\label{homols22d2sd}
 Fix $N\ge  1$, $\mathtt{S}\gg s_0$ and let $\mu_{\mathtt{p},1},\mathtt{b}_0$ be as in \eqref{constandusd2} and \eqref{nsjdjwdsdnsde}. Let $\Omega$ be a subset of $\Omega_1$. Let an $\omega$-dependent operator $\mathcal{R}$ and an $\omega$-dependent Fourier multiplier $r(\omega,j)$ be well-defined for  $\omega\in \Omega$.  We make the following assumptions:
 \begin{enumerate}[label=(A\arabic*)]
 \item \label{rosdus2sd1} $\mathcal{R}$ is a  Lip-$0$-modulo-tame operator up  to $\mathtt{S}$ and satisfies
  \begin{align}\label{rauspsd2sd}
 (\mathcal{R})^{j}_{k}(l) = 0, \text{ of $j,k\notin S_{\mathtt{M}}^\perp$.}
 \end{align} 
 \item\label{rosdus2sd2} $r(\omega,j)$ satisfies 
  \begin{align}\label{kjjsdnksdsdxc}
 |r|^{\Lip(\gamma,\Omega)}_{0,0,0}\le_{\mathtt{pe}} \epsilon^{9}\gamma^{-3}, \quad |d_ir(i_0)|^{\sup}_{0,0,0}\le_{\mathtt{pe}}\epsilon^3\gamma^{-1}\rVert \ihat \rVert_{s_0+\mu_{\mathtt{p},1}},
 \end{align}
where
\[
 |r|^{\Lip(\gamma,\Omega)}_{0,0,0}:=\sup_{\omega,\omega_1,\omega_2\in \Omega,\ j\in \mathbb{Z}\backslash \left\{ 0 \right\}} |r(\omega,j)| +  \gamma \frac{|r(\omega_1,j)-r(\omega_2,j)|}{|\omega_1-\omega_2|}.
\] 
 \end{enumerate}
 We denote 
\begin{align}
d_r(\omega,j)&:= d_0(\omega,j) + r(\omega,j),\text{ where $d_0$ is as in \eqref{d0defsd2}},\quad D_r:= Op^W(d_r),\label{diagonasd2}\\
\tilde{\Omega}&:=\left\{ \omega\in \Omega: |\ii\omega\cdot l - (d_r(\omega,j) - d_r(\omega,k))|\ge \frac{\gamma}{2} |l|^{-\tau}|\lambda_\alpha(j)-\lambda_\alpha(k)|,\right. \nonumber\\
& \qquad \qquad \left.\ \forall |l|\le N,\ \forall j,k\in S_{\mathtt{M}}^\perp\cup\left\{ 0 \right\}\right\},\label{setfofsdfreq}
\end{align}
(we recall $\tau = \frac{1}{\alpha-1}+\nu + 2$ from \eqref{xi_omega_dependent}).
Then, for each $\omega\in \tilde{\Omega}$, the homological equation for $\Psi$:
\begin{align}
\D_\omega(\Psi) - [D_r, \Psi] + \Pi_{\le N}\mathcal{R} & =[\mathcal{R}], \nonumber \\
\text{ where }([\mathcal{R}])^{j}_k(l) & := \begin{cases} \mathcal{R}^j_j(0), & \text{ if $(l,j,k)=(0,j,j)$, $j\in S_{\mathtt{M}}^\perp$}\\ 0, & \text{ otherwise},
\end{cases}\label{homological_sayit}
\end{align}
has a  unique solution $\Psi$ such that 
\begin{enumerate}[label=(\arabic*)]
 \item \label{item1sdkams1}  For $\vec{b}\in \mathbb{N}_0^{\nu}$ such that $|\vec{b}|\le \mathtt{b}_0$, it holds that for $\omega\in \tilde{\Omega}$ and $s\in [s_0,\mathtt{S}]$,  
\begin{align}
\mathfrak{M}^{\sharp,\gamma}_{\partial_\varphi^{\vec{b}}(\Psi)}(0,s) & \le_{\mathtt{pe}}\epsilon^2\gamma^{-2}N^{2\tau+1}\mathfrak{M}^{\sharp,\gamma}_{\partial_\varphi^{\vec{b}}(\mathcal{R})}(0,s),\label{homoglogisd1}\\
  \mathfrak{M}^{\sharp}_{\partial_\varphi^{\vec{b}}(d_i\Psi(i_0)[\ihat])}(0,s_0) & \le_{\mathtt{pe}}\left(\epsilon^3\gamma^{-3}N^{2\tau+1}\mathfrak{M}^{\sharp}_{\partial_\varphi^{\vec{b}}(\mathcal{R})}(0,s_0)\rVert \ihat \rVert_{s_0 +\mu_{\mathtt{p},1}} \right. \nonumber \\
  & \left. \qquad \qquad + \gamma^{-1}N^{\tau}\mathfrak{M}^{\sharp}_{\partial_{\varphi}^{\vec{b}}(d_i\mathcal{R}(i_0)[i\hat])}(0,s_0)\right),\label{homoglogisd2}
\end{align}
\item \label{item2sdksdsd} If $D_r$ and $\mathcal{R}$ are reversible, then $\Psi$ is reversibility preserving.
\item \label{itemforrealsdsd} If $D_r$ and $\mathcal{R}$ are  real operators, then $\Psi$ is a real opeartor.
\item \label{itemforrealsdsd1}  $\Psi$ satisfies
\begin{align}\label{rhhsdsdppwsdqxsd}
(\Psi)^{j}_k(l) = 0, \text{ if  $j,k\notin S^\perp_{\mathtt{M}}$}
\end{align} and $\Psi$ is   $\frac{2\pi}{\mathtt{M}}$-translation invariance preserving.
\end{enumerate}
 \end{lemma}
 \begin{proof}
We first solve the equation \eqref{homological_sayit} and prove the properties of the solution in \ref{item1sdkams1}-\ref{itemforrealsdsd1}.   

We write \eqref{homological_sayit} as
 \begin{align}\label{defsd22scheesed}
 \underbrace{\left(\ii \omega\cdot  l - (d_r(\omega,j) - d_r(\omega,k))\right)}_{=:\delta_{ljk}(\omega)} \Psi^{k}_j(l) + \mathcal{R}^k_j(l) = ([\mathcal{R}])^k_j(l).
 \end{align}
 Therefore, we define\index{homological equation}
 \begin{align}\label{solutiontothehom}
 \Psi^k_j(l):=\begin{cases}
 -\frac{\mathcal{R}^k_j(l)}{ \delta_{ljk}(\omega)}, & \text{ if $(l,j,k)\ne (0,j,j)$ for $j\in S_{\mathtt{M}}^\perp$},\\
 0, & \text{ otherwise.}
 \end{cases}
 \end{align} 
 Clearly, $\Psi$ defined \eqref{solutiontothehom} is a unique solution to \eqref{defsd22scheesed}.
 
 \textbf{Estimates for $\delta_{ljk}(\omega)$: Supremum in $\omega$.}
 First, we see from the definition of $\tilde{\Omega}$ in \eqref{setfofsdfreq} that
 \begin{align}\label{deltaforsup}
\sup_{\omega \in \Omega} \left|\frac{1}{\delta_{ljk}(\omega)}\right| \le 2\gamma^{-1}N^{\tau}|\lambda_\alpha(j) - \lambda_\alpha(k)|^{-1}.
 \end{align}
 
  \textbf{Estimates for $\delta_{ljk}(\omega)$: Lipschitz dependence on $\omega$.}
 For $\omega_1,\omega_2\in \tilde{\Omega}$, we will estimate the dependence on $\omega$ of $\delta_{ljk}$, that is,
 \begin{align}\label{kksd2sdsdpropose}
 \frac{1}{\delta_{ljk}(\omega_1)}- \frac{1}{\delta_{ljk}(\omega_2)} = \frac{\delta_{ljk}(\omega_2)-\delta_{ljk}(\omega_1)}{\delta_{ljk}(\omega_1)\delta_{ljk}(\omega_2)}.
 \end{align}
 Clearly, it follows from the definition of $\delta_{ljk}$ in \eqref{defsd22scheesed} that
 \begin{align} \nonumber 
 \delta_{ljk}(\omega_2)-\delta_{ljk}(\omega_1)& \le |\omega_1-\omega_2|N \\
 & \ + |(d_r(\omega_1,j)-d_r(\omega_2,j)) - (d_r(\omega_1,k)-d_r(\omega_2,k))|
 \label{ksd2sdsds}
 \end{align}
 From \eqref{diagonasd2} and \eqref{d0defsd2}, we have that 
 \begin{equation}\label{rsdsdsd2sdsd}
 \begin{aligned}
 |(d_r(\omega_1,j)-d_r(\omega_2,j)) &- (d_r(\omega_1,k)-d_r(\omega_2,k))|\\
 &  \le |(\mathtt{m}_\alpha(\omega_1)-\mathtt{m}_\alpha(\omega_2))(\lambda_\alpha(j)-\lambda_\alpha(k))| \\
 & \ + \left| \int_{j}^{k} \partial_\xi\mathfrak{m}_{1}(\omega_1)(\xi)-\partial_\xi\mathfrak{m}_{1}(\omega_2)(\xi)d\xi \right|\\
 & \ + |r(\omega_1,j)-r(\omega_2,j)| +  |r(\omega_1,k)-r(\omega_2,k)|. 
 \end{aligned}
 \end{equation}
 Using the estimates for $\mathtt{m}_\alpha$ in\eqref{nusd2sdsd} and using \eqref{rejsmsdim1} and \eqref{kjjsdnksdsdxc}, we have
 \begin{equation}\label{lsdjsd2sdsd12}
 \begin{aligned}
 |(\mathtt{m}_\alpha(\omega_1)-\mathtt{m}_\alpha(\omega_2))(\lambda_\alpha(j)-\lambda_\alpha(k))|  &\le_{\mathtt{pe}}\epsilon^2 \gamma^{-1} |\omega_1-\omega_2||\lambda_\alpha(j)-\lambda_\alpha(k)|,\\
 \left| \int_{j}^{k} \partial_\xi\mathfrak{m}_{1}(\omega_1)(\xi)-\partial_\xi\mathfrak{m}_{1}(\omega_2)(\xi)d\xi \right| & \le_{\mathtt{pe}}\epsilon^2\gamma^{-1}|\omega_1-\omega_2| |j-k|,\\
 |r(\omega_1,j)-r(\omega_2,j)| +  |r(\omega_1,k)-r(\omega_2,k)|& \le_{\mathtt{pe}}\epsilon^2\gamma^{-1} |\omega_1-\omega_2|
 \end{aligned}
 \end{equation}
 Hence, using \eqref{pheobe1}, we have that for $j\ne k$, 
 \begin{align*}
 \left| (d_r(\omega_1,j)-d_r(\omega_2,j)) \right. & \left.- (d_r(\omega_1,k)-d_r(\omega_2,k))\right| \\
 & \le_{\mathtt{pe}}\gamma^{-1}\epsilon^2|\omega_1-\omega_2| |\lambda_\alpha(j)-\lambda_\alpha(k)|.
 \end{align*}
 Plugging this into  \eqref{ksd2sdsds} and \eqref{kksd2sdsdpropose}, and using \eqref{deltaforsup}, we obtain
 \begin{align}\label{sdsdsdimsdwsx}
 \gamma \left| \frac{1}{\delta_{ljk}(\omega_1)}- \frac{1}{\delta_{ljk}(\omega_2)} \right|\le_{\mathtt{pe}} |\omega_1-\omega_2| \epsilon^2\gamma^{-2} N^{2\tau + 1}|\lambda_\alpha(j)-\lambda_\alpha(k)|^{-1}.
 \end{align}
 
  \textbf{Estimates for $\delta_{ljk}(\omega)$: Dependence on $i_0$.} For a fixed $\omega\in \tilde{\Omega}$ (we will omit its notation), we have that
    \begin{align}\label{imgonnatakeoff}
  \left| d_i\left( \frac{1}{\delta_{ljk}}\right)(i_0)[\ihat] \right| = \left|\frac{d_i\delta_{ljk}(i_0)[\ihat]}{(\delta_{ljk})^2}\right|.
  \end{align}
  Following the same computations as in \eqref{rsdsdsd2sdsd}, (one can simply replace the difference on $\omega$ by $d_i$), it is straightforward to see that (using the estimates in \eqref{nusd2sdsd}, \eqref{rejsmsdim1}, \eqref{kjjsdnksdsdxc} and Lemma~\ref{melnikove_23s}),
  \begin{align*}
  \left|d_i\delta_{ljk}(i_0)[\ihat] \right|& \le   \left| d_i(\mathfrak{m}_\alpha(i_0)[\ihat])(\lambda_\alpha(j)-\lambda_\alpha(k))\right| + \left| \int_{j}^{k} \partial_\xi(d_i\mathfrak{m}_1(i_0)[\ihat])(\xi)d\xi \right| \\
  & + |d_ir(i_0)[\ihat](j)| + |d_ir(i_0)[\ihat](k)|\\
  & \le_{\mathtt{pe}}\epsilon^3\gamma^{-1} \rVert \ihat \rVert_{s_0+\mu_{\mathtt{p},1}} | \lambda_\alpha(j)-\lambda_\alpha(k)|.
    \end{align*}
    Therefore, \eqref{imgonnatakeoff} and \eqref{deltaforsup} give us that
    \begin{align}\label{pibsdonimbddsd}
    \left| d_i\left( \frac{1}{\delta_{ljk}}\right)(i_0)[\ihat] \right| \le C\epsilon^3\gamma^{-3}N^{2\tau}| \lambda_\alpha(j)-\lambda_\alpha(k)|^{-1}\rVert \ihat \rVert_{s_0+\mu_{\mathtt{p},1}}.
    \end{align}
    
    \vspace{0.5\baselineskip}\noindent\textit{Proof of \ref{item1sdkams1}.}
    Plugging \eqref{deltaforsup}, \eqref{sdsdsdimsdwsx} and \eqref{pibsdonimbddsd} into \eqref{solutiontothehom}, it is straightforward that (using Lemma~\ref{melnikove_23s}, which tells us that $|\lambda_\alpha(j)-\lambda_\alpha(k)|^{-1}\le_{\mathtt{pe}}1$, unless $j=k$), 
    \begin{align}
    |\Psi^k_j(l)| &\le_{\mathtt{pe}} \gamma^{-1}N^{\tau}|\mathcal{R}^k_j(l)|,\label{taem1hanson2}\\
    \gamma\left| \Psi^k_j(l)(\omega_1)-\Psi^k_j(l)(\omega_2)\right| & \le_{\mathtt{pe}} |\omega_1-\omega_2|\epsilon^2\gamma^{-2}N^{2\tau+1} |R^k_j(l)| \nonumber \\
    & + N^\tau |\mathcal{R}^k_j(l)(\omega_1)-\mathcal{R}^k_j(l)(\omega_2)| , \label{taem1hanson3}\\ 
      |(d_i\Psi(i_0)[\ihat])^k_j(l)& \le_{\mathtt{pe}} \epsilon^3\gamma^{-3}N^{2\tau}|\mathcal{R}^k_j(l)|\rVert \ihat \rVert_{s_0+\mu_0} \nonumber \\
      & +  \gamma^{-1}N^\tau |\left(d_i\mathcal{R}(i_0)[\ihat]\right)^k_j(l)| \label{taem1hanson23}
    \end{align}
    
Towards the estimates \eqref{homoglogisd1},  we pick $h\in C^\infty_{\varphi,x}$, and $\vec{b}\in \mathbb{N}_0^{\nu}$ such that $|\vec{b}|\le \mathtt{b}_0$. We have that (recall the notations \eqref{matrix_rep_2} and \eqref{maojisd2sd} for a linear operator),
    \begin{align*}
    & \rVert \underline{ \partial_\varphi^{\vec{b}}(\Psi)} h \rVert_{s}^2 \\
    &\le \sum_{(l,j)\in\mathbb{Z}^{\nu+1}}\langle l,j\rangle^{2s}\left(\sum_{(l',k)\in \mathbb{Z}^{\nu+1}}\left|(\partial_\varphi^{\vec{b}}(\Psi))^k_j(l-l')\right| |h_k(l')|\right)^2\\
    & = \sum_{(l,j)\in\mathbb{Z}^{\nu+1}}\langle l,j\rangle^{2s}\left(\sum_{(l',k)\in \mathbb{Z}^{\nu+1}}|(l-l')^{\vec{b}}|\left|\Psi^k_j(l-l')\right| |h_k(l')|\right)^2\\
    &\overset{\eqref{taem1hanson2}}{\le}_{\mathtt{pe}} \gamma^{-2}N^{2\tau}\sum_{(l,j)\in\mathbb{Z}^{\nu+1}}\langle l,j\rangle^{2s}\left(\sum_{(l',k)\in \mathbb{Z}^{\nu+1}}|(l-l')^{\vec{b}}|R^{k}_j(l-l') |h_k(l')|\right)^2\\
    & \overset{\eqref{maojisd2sd1},\eqref{maojisd2sd}}{=} \gamma^{-2}N^{2\tau}\rVert \underline{\partial_\varphi^{\vec{b}}(\mathcal{R})}\underline{h}\rVert_{s}^2\\
    & \overset{\text{Definition~\ref{def_tame1}}}\le\gamma^{-2}N^{2\tau}\left(\mathfrak{M}^\sharp_{\partial_{\varphi}^{\vec{b}}(\mathcal{R})}(0,s)\right)^2\rVert \underline{h}\rVert^2_s\\
    & \overset{\eqref{endedupspending}}{=} \gamma^{-2}N^{2\tau}\left(\mathfrak{M}^\sharp_{\partial_{\varphi}^{\vec{b}}(\mathcal{R})}(0,s)\right)^2\rVert {h}\rVert^2_s.
    \end{align*}   
   Therefore, we have
   \begin{align}\label{modulotsd11sd2}
   \mathfrak{M}^\sharp_{\partial_{\varphi}^{\vec{b}}(\Psi)}(0,s)\le_{\mathtt{pe}} \gamma^{-1}N^{\tau}\mathfrak{M}^\sharp_{\partial_\varphi^{\vec{b}}(\mathcal{R})}(0,s).
   \end{align}
   Similarly,  following the same computations but using \eqref{taem1hanson3} and \eqref{taem1hanson23}, one can easily see the desired estimates \eqref{homoglogisd1} and \eqref{homoglogisd2}.

   \vspace{0.5\baselineskip}\noindent\textit{Proof of \ref{item2sdksdsd} and \ref{itemforrealsdsd}}
    To see items \ref{item2sdksdsd} and \ref{itemforrealsdsd}, in view of \eqref{lsdrjdmlRmxdlqhdlsekx}, let us denote involutions $P_{rev},P_{real}$ acting on linear operators by
    \[
    (P_{rev}\mathcal{A})^{j_1}_{j_2}(l):=\overline{ \mathcal{A}^{-j_1}_{-j_2}(-l)},\quad \left(P_{real}\mathcal{A}\right)^{j_1}_{j_2}(l) := \overline{(\mathcal{A})^{j_2}_{j_1}(l)}.
    \]
    Clearly, $P_{rev},P_{real}$ are involutions, since $P_{rev}\circ P_{rev} = I$ and $P_{real}\circ P_{real} = I$. Recalling that $D_r$ and $\mathcal{R}$ are reversible, it follows immediately, taking the involutions in the homological equation \eqref{homological_sayit} that $-P_{rev}\Psi$ and $P_{real}\Psi$ are also solutions to \eqref{homological_sayit}. Since the solution is unique, we have that $\Psi  = -P_{rev}\Psi = P_{real}\Psi$. Therefore, $\Psi$ is reversibility preserving and real.
     \end{proof}
     
     \vspace{0.5\baselineskip}\noindent\textit{Proof of \ref{itemforrealsdsd1}.} 
It is clear that $\Psi$ is $\frac{2\pi}{\mathtt{M}}$-translation invariant if and only if
     \begin{align}\label{sd2cxcxc3mmma}
    (\Psi)^{j_1}_{j_2}(l) = 0, \text{ if $j_1 - j_2 \notin \left\{ \mathtt{M}j: j\in \mathbb{Z}\right\}$},
     \end{align}
    Therefore, the result follows immediately from \eqref{rauspsd2sd} and \eqref{solutiontothehom}.
 
 \subsection{Choice of the increasing sequence $N_n$}
 We fix
 \begin{equation}\label{parameterssd1sd}
 \begin{aligned}
 N_n&:=N_0^{\chi^n},\text{ for $n\ge0$},\quad \chi:=\frac{3}{2},\quad N_{-1}:=1,\quad N_0:=\left(\epsilon\gamma^{-1}\right)^{\rho},\\
 \sigma_1 &:= 3(2\tau +1),\quad \mu_{\mathtt{p},0}:=2\mu_{\mathtt{p},1} + 2\tau +1\\
   \rho^{-1}&:=\max\left\{ 22(2b-1)(\sigma_1+2\tau+1),\frac{2\tau(2b-1)}{9-8b},\right. \\
   & \qquad \qquad \qquad \left.  \frac{(2b-1)(2\tau+1)}{2b-2},\frac{40(\mu_{\mathtt{p},0} + \mu_{\mathtt{p},2})(2b-1)}{5-4b}\right\}.
 \end{aligned}
 \end{equation}
 \begin{remark}\label{jsdjksdksd}
 Note that the above sequence is well-defined in the sense that $\tau$ is fixed in \eqref{frequency_set2}, $b$ is fixed in \eqref{parametersets1231}, which gives $\sigma_1$ as above. $\mu_{\mathtt{p},1}$ is fixed in \eqref{constandusd2}, which gives us $\mu_{\mathtt{p},0}$ and this gives us $\mu_{\mathtt{p},2}$ as described in Remark~\ref{labesdof_sdxxcxc}, therefore $\rho$ is well-defined.
 \end{remark}
 
One can easily check that for all sufficiently small $\epsilon>0$, depending on $\mathtt{p}$, and for all $n\ge0$,
 \begin{align}
 \mathtt{b}_0&\overset{\eqref{nsjdjwdsdnsde}}\ge 2\sigma_1,\label{kjsoxcxcj}\\
\epsilon^{\frac{1}{22}} N_n^{2\tau +1}N_{n-1}^{-2\sigma_1} &\le N_n^{-\sigma_1},\label{para1chs12s}\\
 N_{n}^{-\mathtt{b}_0}N_{n-1}&\le N_n^{1-\sigma_1}N_n^{-\sigma_1} \le N_0^{1-\sigma_1}N_n^{-\sigma_1}\label{para1chs1}\\
\epsilon^{11-10b} N_{n}^{2\tau+1}N_{n-1}^{-2\sigma_1}&=\epsilon^{1-\frac{10}{11}} \underbrace{\epsilon^{10(\frac{1}{11} - (b-1))}}_{\le 1,\ \because\eqref{parametersets1231}}N_{n}^{2\tau+1}N_{n-1}^{-2\sigma_1} \nonumber \\
& \le \epsilon^{\frac{1}{22}}\epsilon^{\frac{1}{22}}N_n^{2\tau+1}N_{n-1}^{-2\sigma_1}\overset{\eqref{para1chs12s}}\le \epsilon^{\frac{1}{22}}N_n^{-\sigma_1}\label{para1chs12}.
 \end{align}

 \subsection{KAM reducibility}
 We take $L_0,D_0,\mathcal{R}_0,d_0$ given\index{KAM reducibility} in Proposition~\ref{modulut2sosdtame}. For $n\ge 0$, we define inductively: 
 \begin{align}
 r_0&:= 0,\quad \Omega^\gamma_0(i_0):=\Omega_1,\label{sddpooo1}\\
 \phi_{n}&:=I+\psi_n,\  \phi_n^{-1}=:I-\tilde{\psi}_n, \nonumber
 \end{align}
 where $\psi_n$ solves \begin{align}\D_\omega(\psi_n) - [D_n,\psi_n]+\Pi_{\le N_n}\mathcal{R}_n = [\mathcal{R}_n] \in H_{S^\perp}\cap X_{\mathtt{M}},
 \label{sddpooo2}
 \end{align}
 \begin{align}
r_{n+1}&:=r_{n+1}(\omega,j):=\begin{cases}
-[\mathcal{R}_n(\omega)]^j_j, & \text{ if $j\in S_{\mathtt{M}}^\perp$},\\
0, & \text{ if $j\ne S_{\mathtt{M}}^\perp$},
\end{cases}  \text{ for  $\omega\in \Omega^\gamma_n(i_0)$},\label{sddpooo4}\\
 D_{n+1}&:=Op^W(d_{n+1}),\quad d_{n+1}:=d_n + r_{n+1} = d_0 + \sum_{k=0}^{n+1}r_k =: d_0 + r_{\le {n+1}},\label{sddpooo5}\\
\mathcal{R}_{n+1}&:=\tilde{\psi}_n\circ Op^W(r_n) + \phi_n^{-1}  \circ \left( \Pi_{\ge N_n}\mathcal{R}_n + \mathcal{R}_n\psi_{n}\right).\label{sddpooo6}
 \end{align}
 Furthermore, we consider a sequence of sets:\index{$\Omega^\gamma_n$}
 \begin{align}  &\Omega^\gamma_{n+1}(i_0)\nonumber\\
 & \quad :=\left\{ \omega\in \Omega_n^\gamma(i_0): |\ii \omega\cdot l - (d_n(\omega,j) - d_n(\omega,k))|\ge \frac{1}{2}\gamma |l|^{-\tau} |\lambda_\alpha(j)-\lambda_\alpha(k)|,\right. \nonumber\\
 & \qquad \qquad \qquad \qquad \qquad \qquad \qquad \qquad \qquad \qquad \left.\ \forall |l|\le N_{n},\ \forall j,k\in S_{\mathtt{M}}^\perp\cup\left\{0\right\}\right\}.\label{sddpooo3}
\end{align}
 Trivially, we have that 
 \begin{align}\label{2gammareason}
 \Omega^{2\gamma}_{n+1}(i_0)\subset \Omega^{\gamma}_{n+1}(i_0).
 \end{align}
 The motivation is that we aim to conjugate $L_n=\D_\omega - D_n +\mathcal{R}_n$ with a transformation $\phi_n=I + \psi_n, \phi^{-1}_n=I+\tilde{\psi}_n$ to reduce the size of the remainder. One can easily show that for each $n\ge 0$,
 \begin{align}\label{consd2sdsdsd}
 &\phi_n^{-1}L_n\phi_n \nonumber\\
 &= \D_\omega + \underbrace{\left(- D_n + \left( \D_\omega(\psi_n) -[D_n,\psi_n] + \Pi_{\le N_n}\mathcal{R}_n\right)\right)}_{= - D_{n+1}}\nonumber\\
 &\ + \underbrace{\tilde{\psi}_n\circ\left( \D_\omega(\psi_n) -[D_n,\psi_n] + \Pi_{\le N_n}\mathcal{R}_n \right) + \Pi_{\ge N_n}\mathcal{R}_n + \mathcal{R}_n\psi_n}_{=\mathcal{R}_{n+1}}.
 \end{align}
 Therefore, we simply choose $\psi_n$ to be a solution to the corresponding homological equation in \eqref{sddpooo2}, while the diagonal part is updated. Since our initial remainder in this scheme, $\mathcal{R}_0$, is already a $0$-tame operator, we can repeat this procedure infinitely many times without any loss of derivatives (since a composition of $0$-tame operators is still a $0$-tame operator). The parameters in \eqref{parameterssd1sd} are carefully chosen so that we can close the infinite number of conjugations.  
  Whether they are well-defined or not will be checked in Proposition~\ref{kaksdsddsd}. Note that for all $n\ge0$ such that
\begin{align}\label{kjjsdnksdsd23xc}
|r_{\le n}|^{\Lip(\gamma,\Omega^\gamma_{n}(i_0))}_{0,0,0}\le_{\mathtt{pe}}  \epsilon^{9}\gamma^{-3}, \quad |d_ir_{\le n}(i_0)|^{\sup}_{0,0,0}\le_{\mathtt{pe}}   \epsilon^3\gamma^{-1}\rVert \ihat \rVert_{s_0+\mu_{\mathtt{p},1}},
 \end{align}
 Lemma~\ref{homols22d2sd} tells us that $\psi_n=\psi_n(\omega)$ exists in \eqref{sddpooo2} for $\omega\in \Omega_n^\gamma(i_0)$. Once we have $\psi_n$, then assuming $\mathfrak{M}^{\sharp,\gamma}_{\psi_n}(0,s_0)$ is small enough, Lemma~\ref{kinverson_modulsd} gives us the existence of $\tilde{\psi}_n$ in \eqref{sddpooo2}, which allows us to have $\mathcal{R}_{n+1}$ and $D_{n+1}$ in \eqref{sddpooo6} and \eqref{sddpooo5} well-defined for each $\omega\in \Omega^\gamma_{n}(i_0)$. Indeed, we have the following:
 \begin{proposition}\label{kaksdsddsd}
 Given $\mathtt{S}\gg s_0$, let us denote 
 \begin{equation}\label{defofmm2ksdsd}
 \begin{aligned}
\mathbb{M}_0(s)&:=\sup_{\vec{b}\in \mathbb{N}_0^\nu, |\vec{b}|\le \mathtt{b}_0}\mathfrak{M}^{\sharp,\gamma}_{\partial_{\varphi}^{\vec{b}}\mathcal{R}_0}(0,s),\text{ for $s\in [s_0,\mathtt{S}]$.}
 \end{aligned}
 \end{equation}
 Assuming \eqref{size_assumption_2} for $\mu=\mu_{\mathtt{p},1}$, there exists $\epsilon_0(\mathtt{p},\mathtt{S})>0$ such that for all $\epsilon\in(0,\epsilon_0)$, we have that, for every $n\ge 0$:
 \begin{enumerate}[label=(\arabic*)]

 \item \label{kamrdu3}(Estimates for $\mathcal{R}_n$) $\mathcal{R}_n$ is well defined for $\omega\in \Omega^\gamma_{n}(i_0)$ as a real, reversible operator such that
 \begin{align}\label{sjowk1sdsds1}
 (\mathcal{R}_n)^{j}_k(l) = 0 , & \text{ if $j,k\notin S^\perp_{\mathtt{M}}$.}
 \end{align}
 Furthermore, it satisfies the following estimates: For all $s\in [s,\mathtt{S}]$, 
 \begin{align}
 \mathfrak{M}^{\sharp,\gamma}_{\mathcal{R}_n}(0,s)&\le \mathbb{M}_0(s) N_{n-1}^{-\sigma_1},\label{sdjsds11ram1}\\
   \sup_{\vec{b}\in \mathbb{N}^\nu_0,\ |\vec{b}| = \mathtt{b}_0}\mathfrak{M}^{\sharp,\gamma}_{\partial_\varphi^{\vec{b}}(\mathcal{R}_n)}(0,s)&\le \mathbb{M}_0(s) N_{n-1}, \label{sdjsds11ram2}\\
   \mathfrak{M}^{\sharp}_{d_i(\mathcal{R}_n)(i_0)[\ihat]}(0,s_0)&\le_{\mathtt{pe}} \epsilon^3\gamma^{-1}N_{n-1}^{-\sigma_1}\rVert \ihat \rVert_{s_0+\mu_{\mathtt{p},1}},\label{sdjsds11ram3}\\
    \sup_{\vec{b}\in \mathbb{N}^\nu_0,\ |\vec{b}| = \mathtt{b}_0}\mathfrak{M}^{\sharp}_{\partial_\varphi^{\vec{b}}(d_i\mathcal{R}_n(i_0)[\ihat])}(0,s_0)&\le_{\mathtt{pe}}  \epsilon^3\gamma^{-1}N_{n-1}\rVert \ihat \rVert_{s_0+\mu_{\mathtt{p},1}},\label{sdjsds11ram4}.
 \end{align}

 \item \label{kamrdu4}(Estimates for $D_{n+1}$) $D_{n+1}$ is well defined (that is, $r_{n+1}$ is well-defined) for $\omega\in \Omega^\gamma_{n}(i_0)$ as a diagonal, real and reversible operator. Furthermore, it satisfies the following estimates: 
 \begin{align}
|r_{n+1}|^{\Lip(\gamma,\Omega^\gamma_n(i_0))}_{0,0,0} & \le_{\mathtt{pe}} \mathbb{M}_0(s_0) N_{n-1}^{-{\sigma_1}},\label{rsdj1sd11}\\
|r_{\le n+1}|^{\Lip(\gamma,\Omega^\gamma_n(i_0))}_{0,0,0} & \le_{\mathtt{pe}}\mathbb{M}_0(s_0),\label{jjsd2sdsdsdsdsdsdsdsdsdsdsds}\\
|d_ir_{n+1}(i_0)[\ihat]|^{\sup}_{0,0,0} & \le_{\mathtt{pe}}\epsilon^3\gamma^{-1}N_{n-1}^{-{\sigma_1}}\rVert \ihat \rVert_{s_0+\mu_\mathtt{p}},\label{rsdj1sd12}\\
|d_ir_{\le {n+1}}(i_0)[\ihat]|^{\sup}_{0,0,0} & \le_{\mathtt{pe}}\epsilon^3\gamma^{-1}\rVert \ihat \rVert_{s_0+\mu_\mathtt{p}}. \label{sdsdsd2sd1sxf}
 \end{align} 
     \item \label{kamrdu2}(Estimates for $\psi_{n}$) $\psi_n$ is well defined for $\omega\in \Omega^\gamma_{n+1}(i_0)$ as a real, reversibility preserving operator such that
     \begin{align}\label{sjowk1sdsds2}
     (\psi_n)^{j}_k(l) = 0 , & \text{ if $j,k\notin S^\perp_{\mathtt{M}}$.}
     \end{align}
     Furthermore, it satisfies the following estimate: For all $s\in [s,\mathtt{S}]$,
 \begin{align}
  \mathfrak{M}^{\sharp,\gamma}_{\psi_{n}}(0,s)&\le_{\mathtt{pe}}\epsilon^2\gamma^{-2} N_{n-1}^{-\sigma_1}N_n^{2\tau +1}\mathbb{M}_0(s)\label{psdjsds11ram1}
   \end{align}
  \item \label{kamrdu5} Let $i_1,i_2$ be two embeddings satisfying \eqref{size_assumption_2}, \eqref{sjwojsdsdiwosdsd} for $\mu_{\mathtt{p},1}$ and some $\mathtt{C}>0$. There exists  $C=C(\mathtt{pe})>0$ such that  if   $C\epsilon^3\gamma^{-1}N_{n-1}^{\tau}\rVert i_1-i_2\rVert_{s_0+\mu_{\mathtt{p},1}}\le \delta$, for some $\delta\in(0,\gamma/2)$, then we have that
 \begin{align}\label{tlwoclaosd}
 \implies \Omega_n^{2\gamma}(i_1)\subset \Omega^{2(\gamma-\delta)}_{n}(i_2),
 \end{align}
 for sufficiently small $\epsilon$ depending on $\mathtt{pe}$.

 \end{enumerate}
Lastly for all $n\ge 0$ and $\omega\in \cap_{n\ge 0}\Omega_n^{\gamma}(i_0)$, we have that  $L_{n}:=\D_\omega - D_n + \mathcal{R}_n$ satisfies
 \begin{align}\label{consjdsdsd2sd}
 L_{n+1} = \phi_{n}^{-1}L_n\phi_n.
 \end{align}
    Furthermore, $L_n$ is real and reversible and $\frac{2\pi}{\mathtt{M}}$-translation invariant.
 \end{proposition}
\begin{proof}
We will prove items \ref{kamrdu3}-\ref{kamrdu5} by induction. Once \ref{kamrdu3}-\ref{kamrdu5} are verified, then we will prove \eqref{consjdsdsd2sd}, and the reversibility, reality and $\frac{2\pi}{\mathtt{M}}$-translation invariance preserving properties of $L_n$.

\vspace{0.5\baselineskip}\noindent\textit{Proof of \ref{kamrdu3}-\ref{kamrdu5} for $n=0$.}
   For item \ref{kamrdu3}, $\mathcal{R}_0$ is already well-defined in Proposition~\ref{modulut2sosdtame} for $\omega\in \Omega_1=\Omega_0^\gamma(i_0)$. Furthermore, \ref{reversiblitsd} of Proposition~\ref{modulut2sosdtame} tells us that $\mathcal{R}_0$ is real and reversible. \eqref{sjowk1sdsds1}  follows from \eqref{rkjspspdwjskds}. 
   
   For the estimates, recalling from \eqref{parameterssd1sd} that $N_{-1}=1$, the estimates \eqref{sdjsds11ram1} and \eqref{sdjsds11ram2} hold trivially, thanks to the definition of $\mathbb{M}_0(s)$ in \eqref{defofmm2ksdsd}. The estimates \eqref{sdjsds11ram3} and \eqref{sdjsds11ram4} follow from \eqref{fromyourbag2}.    

 For item \ref{kamrdu4}, recalling $r_1$ from \eqref{sddpooo4}, we see that $r_1$ is well-defined, since so is $\mathcal{R}_0$, as proved above. Since $\mathcal{R}_0, D_0$ is real and reversible (see \ref{reversiblitsd} of Proposition~\ref{modulut2sosdtame}), it follows straightforwardly from the definition of $r_1$ in \eqref{sddpooo4} and \eqref{lsdrjdmlRmxdlqhdlsekx} that $Op^W(r_{1})$ is also real and reversible. For the estimates, using $N_{-1}=1$ and using  the estimates in item \ref{kamrdu3} for $n=0$ and Lemma~\ref{stsdjwdsymbodlsd}, we see that  the estimates \eqref{rsdj1sd11}-\eqref{sdsdsd2sd1sxf} follow immediately.  
  
 For item \ref{kamrdu2}, we apply Lemma~\ref{homols22d2sd} to obtain $\psi_0$ from \eqref{sddpooo2} (since  \ref{rosdus2sd1} follows from \eqref{sjowk1sdsds1} for $\mathcal{R}_0$ and $r_0=0$, hence \ref{rosdus2sd2} holds trivially). Then for $\omega\in \Omega^\gamma_0(i_0)$, we see that $\psi_0(\omega)$ is well-defined. From \ref{item2sdksdsd}, \ref{itemforrealsdsd} and \ref{itemforrealsdsd1} of Lemma~\ref{homols22d2sd}, we see that $\Psi$ is real, reversibility preserving and  satisfies \eqref{rhhsdsdppwsdqxsd}. Also, the estimate \eqref{psdjsds11ram1} follows from \eqref{homoglogisd1} with $\vec{b} = 0,\mathcal{R}=\mathcal{R}_0$ and $N=N_{0}$.
  
 For \ref{kamrdu5}, \eqref{tlwoclaosd} holds trivially since $\Omega_0^{2\gamma}(i_0)$ actually does not depend on the embedding $i_0$, see \eqref{sddpooo1}.

Now, we assume that \ref{kamrdu3}-\ref{kamrdu2} hold true for some $n\ge 0$ and aim to prove \ref{kamrdu3} for $n+1$.

\textbf{Inductive step for \ref{kamrdu3}.} 
First, we check the well-definedness and reality, reversibility and \eqref{sjowk1sdsds1} for $\mathcal{R}_{n+1}$. We notice that our induction hypotheses for  \eqref{jjsd2sdsdsdsdsdsdsdsdsdsdsds} and \eqref{sdsdsd2sd1sxf} imply that the condition \eqref{kjjsdnksdsdxc} is satisfied for the homological equation for $\psi_{n}$ in \eqref{sddpooo2}, which guarantees that $\psi_{n}$ is well-defined for $\omega\in \Omega^\gamma_{n+1}(i_0)$. We use the induction hypotheses that $\mathcal{R}_n$ and $D_{n}$ are real, reversible and $\frac{2\pi}{\mathtt{M}}$-translation invariant to see that \ref{item2sdksdsd} of Lemma~\ref{homols22d2sd}  implies that  $\Psi_{n}$ is a   real and  reversibility and $\frac{2\pi}{\mathtt{M}}$-translation invariance preserving operator, and so are $\phi_n,\phi_n^{-1}$.  In view of \eqref{consd2sdsdsd}, we have that $\D_\omega -  D_{n+1}  + \mathcal{R}_{n+1}$ is real, reversible. Thanks to the induction hypotheses for \ref{kamrdu4}, $D_{n+1}$ is real and reversible. Therefore, $\mathcal{R}_{n+1}$ is real and reversible. Since $\mathcal{R}_n$ and $\Psi_n$ satisfy \eqref{sjowk1sdsds1} and \eqref{sjowk1sdsds2}, it follows from the definition of $\mathcal{R}_{n+1}$ in \eqref{sddpooo6} that $\mathcal{R}_{n+1}$ satisfies \eqref{sjowk1sdsds1} as well.

 Towards the estimates in \ref{kamrdu3},  we first derive the estimates for $\psi_n$, $\tilde{\psi}_n$, $r_n$, $\Pi_{\ge N_n}\mathcal{R}_n$ and $\mathcal{R}_n\psi_n$.
 
 \textbf{Estimates for $\psi_n$.}
  Applying  \eqref{homoglogisd1} and \eqref{homoglogisd2} to \eqref{sddpooo2}, we see that
 \begin{align*}
 \mathfrak{M}^{\sharp,\gamma}_{\psi_n}(0,s) &\le_{\mathtt{pe}}\epsilon^2\gamma^{-2}N_n^{2\tau+1}\mathfrak{M}^{\sharp,\gamma}_{\mathcal{R}_n}(0,s),\\
\sup_{\vec{b}\in \mathbb{N}^\nu_0,\ |\vec{b}|=\mathtt{b}_0} \mathfrak{M}^{\sharp,\gamma}_{\partial_\varphi^{\vec{b}}(\psi_n)}(0,s) &\le_{\mathtt{pe}}\epsilon^2\gamma^{-2}N_n^{2\tau+1}\sup_{\vec{b}\in \mathbb{N}^\nu_0,\ |\vec{b}|=\mathtt{b}_0}\mathfrak{M}^{\sharp,\gamma}_{\partial_\varphi^{\vec{b}}(\mathcal{R}_n)}(0,s),\\
 \mathfrak{M}^{\sharp}_{d_i\Psi_n(i_0)[\ihat]}(0,s_0) & \le_{\mathtt{pe}}\left(\epsilon^3\gamma^{-3}N_n^{2\tau+1}\mathfrak{M}^{\sharp}_{\mathcal{R}_n}(0,s_0)\rVert \ihat \rVert_{s_0 +\mu_\mathtt{p}} \right.\\
 & \qquad \qquad \qquad \left.+ \gamma^{-1}N_n^{\tau}\mathfrak{M}^{\sharp}_{d_i\mathcal{R}_n(i_0)[i\hat]}(0,s_0)\right),
 \end{align*}
 and 
 \begin{align*}
&\sup_{\vec{b}\in \mathbb{N}^\nu_0,\ |\vec{b}|=\mathtt{b}_0}\mathfrak{M}^{\sharp}_{\partial_\varphi^{\vec{b}}(d_i\Psi_n(i_0)[\ihat])}(0,s_0)\nonumber\\
& \quad  \le_{\mathtt{pe}}\left(\epsilon^3\gamma^{-3}N_n^{2\tau+1}\sup_{\vec{b}\in \mathbb{N}^\nu_0,\ |\vec{b}|=\mathtt{b}_0}\mathfrak{M}^{\sharp}_{\partial_\varphi^{\vec{b}}(\mathcal{R}_n)}(0,s_0)\rVert \ihat \rVert_{s_0 +\mu_\mathtt{p}}\right. \\
& \qquad \qquad \qquad \left.+ \gamma^{-1}N_n^{\tau}\sup_{\vec{b}\in \mathbb{N}^\nu_0,\ |\vec{b}|=\mathtt{b}_0}\mathfrak{M}^{\sharp}_{\partial_{\varphi}^{\vec{b}}(d_i\mathcal{R}_n(i_0)[i\hat])}(0,s_0)\right).
 \end{align*}
 Using the induction hypotheses for $n$, we can use the estimates in \ref{kamrdu3} and we get
 \begin{align}
 \mathfrak{M}^{\sharp,\gamma}_{\psi_n}(0,s) &\le_{\mathtt{pe}}\epsilon^2\gamma^{-2}N_n^{2\tau +1}N_{n-1}^{-\sigma_1}\mathbb{M}_0(s),\label{jsdj2sdsd1sd1}\\
 \sup_{\vec{b}\in \mathbb{N}^\nu_0,\ |\vec{b}|=\mathtt{b}_0} \mathfrak{M}^{\sharp,\gamma}_{\partial_\varphi^{\vec{b}}(\psi_n)}(0,s) &\le_{\mathtt{pe}}\epsilon^2\gamma^{-2}N_{n}^{2\tau +1}N_{n-1}\mathbb{M}_0(s),\label{jsdj2sdsd1sd2}
 \end{align}
  \begin{align}
\nonumber & \mathfrak{M}^{\sharp}_{d_i\psi_n(i_0)[\ihat]}(0,s_0) \\
& \qquad \qquad \qquad \le_{\mathtt{pe}}\left( \epsilon^3\gamma^{-3}N_{n}^{2\tau+1}N_{n-1}^{-\sigma_1}\mathbb{M}_0(s_0) + \epsilon^3\gamma^{-2}N_n^\tau N_{n-1}^{-\sigma_1}\right)\rVert \ihat \rVert_{s_0+\mu_{\mathtt{p},1}},\label{jsdj2sdsd1sd3}
  \end{align}
  \begin{align}
  \nonumber & \sup_{\vec{b}\in \mathbb{N}^\nu_0,\ |\vec{b}|=\mathtt{b}_0}\mathfrak{M}^{\sharp}_{\partial_\varphi^{\vec{b}}(d_i\psi_n(i_0)[\ihat])}(0,s_0) \\
  & \qquad \qquad \le_{\mathtt{pe}}\left( \epsilon^3\gamma^{-3}N_{n}^{2\tau+1}N_{n-1}\mathbb{M}_0(s_0) + \epsilon^3\gamma^{-2}N_n^\tau N_{n-1}\right)\rVert \ihat \rVert_{s_0+\mu_{\mathtt{p},1}}.\label{jsdj2sdsd1sd4}
 \end{align}
 We see from \eqref{jsdj2sdsd1sd1} that
\begin{align*}
\mathfrak{M}^{\sharp,\gamma}_{\psi_n}(0,s_0) &\le_{\mathtt{pe}}\epsilon^2\gamma^{-2}N_n^{2\tau +1}N_{n-1}^{-\sigma_1}\mathbb{M}_0(s_0),
\end{align*}
while \eqref{fromyourbag1} and \eqref{defofmm2ksdsd} give us that
\begin{align}\label{m01sd2est}
\mathbb{M}_0(s_0)\le_{\mathtt{pe}}\gamma^{-1}\left( \epsilon^5 + \epsilon^3\rVert \mathfrak{I}_\delta\rVert^{\Lip(\gamma,\Omega_1)}_{s_0+\mu_\mathtt{p}}\right) \overset{\eqref{size_assumption_2}}\le_{\mathtt{pe}} \epsilon^{9-6b}.
\end{align}
Thus, we have that
\begin{align}\nonumber
\mathfrak{M}^{\sharp,\gamma}_{\psi_n}(0,s_0)&\overset{\eqref{jsdj2sdsd1sd1}}\le_{\mathtt{pe}}\epsilon^2\gamma^{-2}N_n^{2\tau +1}N_{n-1}^{-\sigma_1}\mathbb{M}_0(s_0) \\
& \overset{\eqref{m01sd2est}} \le_{\mathtt{pe}} \epsilon^{11-10b}N_{n}^{2\tau+1}N_{n-1}^{-\sigma_1} \overset{\eqref{para1chs12}}\le_{\mathtt{pe}}\epsilon^{\frac{1}{22}}
 \label{sjdjsddrhcnwkdWlp}
\end{align}
Furthermore,   it follows from \eqref{jsdj2sdsd1sd3},\eqref{jsdj2sdsd1sd4} with \eqref{m01sd2est} that (using $\epsilon^{3}\gamma^{-3}\mathbb{M}_0(s_0)\overset{\eqref{frequency_set2},\eqref{parametersets1231}}\le \epsilon^{3-4b}=\epsilon^3\gamma^{-2}$)
\begin{equation}\label{psiisds2sd}
\begin{aligned}
 \mathfrak{M}^{\sharp}_{d_i{\psi}_n(i_0)[\ihat]}(0,s_0)&\le_{\mathtt{pe}}\epsilon^{3-4b}N_n^{2\tau+1}N_{n-1}^{-\sigma_1}\rVert \ihat \rVert_{s_0 + \mu_{\mathtt{p},1}},\\
 \sup_{\vec{b}\in \mathbb{N}^\nu_0,\ |\vec{b}|=\mathtt{b}_0}\mathfrak{M}^{\sharp}_{\partial_\varphi^{\vec{b}}(d_i\psi_n(i_0)[\ihat])}(0,s_0) &\le_{\mathtt{pe}}\epsilon^{3-4b}N_n^{2\tau+1}N_{n-1}\rVert \ihat \rVert_{s_0 + \mu_{\mathtt{p},1}}.
\end{aligned}
\end{equation}

\textbf{Estimates for $\tilde{\psi}_n$.}
 For $\tilde{\psi}_n$, 

Thanks to \eqref{sjdjsddrhcnwkdWlp}, we see from Lemma~\ref{kinverson_modulsd}  that  we can choose $\epsilon$ sufficiently small, depending on $\mathtt{pe}$ so that $\tilde{\psi}_n$ in \eqref{sddpooo2} is well-defined and satisfies
\begin{align}
\mathfrak{M}^{\sharp,\gamma}_{\tilde{\psi}_n}(0,s)
&\le_{\mathtt{pe}}\mathfrak{M}^{\sharp,\gamma}_{{{\psi}_n}}(0,s)\overset{\eqref{jsdj2sdsd1sd1}}{\le}_{\mathtt{pe}}\epsilon^2\gamma^{-2}N_n^{2\tau +1}N_{n-1}^{-\sigma_1}\mathbb{M}_0(s),\label{ptildesd1}\end{align}
\begin{align}
\sup_{\vec{b}\in \mathbb{N}^\nu_0, |\vec{b}|=\mathtt{b}_0}&\mathfrak{M}^{\sharp,\gamma}_{\partial_\varphi^{\vec{b}}(\tilde{\psi}_n)}(0,s) \nonumber \\
&\le_{\mathtt{pe}}\sup_{\vec{b}\in \mathbb{N}^\nu_0,\ |\vec{b}|=\mathtt{b}_0}\left(\mathfrak{M}^{\sharp,\gamma}_{\partial_\varphi^{\vec{b}}(\psi_n)}(0,s) +\mathfrak{M}^{\sharp,\gamma}_{\partial_\varphi^{\vec{b}}(\psi_n)}(0,s_0)\mathfrak{M}^{\sharp,\gamma}_{\psi_n}(0,s) \right)\nonumber\\
&\overset{\eqref{jsdj2sdsd1sd2}}\le_{\mathtt{pe}}\epsilon^2\gamma^{-2}N_n^{2\tau+1}N_{n-1}\mathbb{M}_0(s)\left( 1 + \epsilon^2\gamma^{-2}N_{n}^{2\tau+1}N_{n-1}^{-\sigma_1}\mathbb{M}_0(s_0)\right)\nonumber\\
&\overset{\eqref{m01sd2est}}\le_{\mathtt{pe}}\epsilon^2\gamma^{-2}N_n^{2\tau+1}N_{n-1}\mathbb{M}_0(s)(1+ \underbrace{\epsilon^{11-10b}N_{n}^{2\tau+1}N_{n-1}^{-\sigma_1}}_{\le 1,\ \because \eqref{para1chs12}})\nonumber\\
& \le_{\mathtt{pe}}\epsilon^2\gamma^{-2}N_n^{2\tau+1}N_{n-1}\mathbb{M}_0(s). \label{ptildesd2}
\end{align}
Therefore, for $\phi_n^{-1}$ in \eqref{sddpooo2}, we have that 
\begin{equation}\label{psd22323dsd}
\begin{aligned}
\mathfrak{M}^{\sharp,\gamma}_{\phi_n^{-1}}(0,s)&\le 1 + \mathfrak{M}^{\sharp,\gamma}_{\tilde{\psi}_n}(0,s) \overset{\eqref{ptildesd1}}\le 1+c_{\mathtt{pe}}\epsilon^2\gamma^{-2}N_n^{2\tau +1}N_{n-1}^{-\sigma_1}\mathbb{M}_0(s) ,\\
\sup_{\vec{b}\in \mathbb{N}^\nu_0,\ |\vec{b}|=\mathtt{b}_0} \mathfrak{M}^{\sharp,\gamma}_{\phi_n^{-1}}(0,s)&\le_{\mathtt{pe}}  \sup_{\vec{b}\in \mathbb{N}^\nu_0,\ |\vec{b}|=\mathtt{b}_0}\mathfrak{M}^{\sharp,\gamma}_{\partial_\varphi^{\vec{b}}(\tilde{\psi}_n)}(0,s) \\&\overset{\eqref{ptildesd2}} \le_{\mathtt{pe}}\epsilon^2\gamma^{-2}N_n^{2\tau+1}N_{n-1}\mathbb{M}_0(s).
\end{aligned}
\end{equation}
Furthermore, using that $\phi_n^{-1}\circ \phi_n = I$, which implies $d_i(\phi_n^{-1})(i_0)[\ihat] = -\phi_n^{-1}\circ d_i(\psi_n)(i_0)[\ihat]\circ \phi_n^{-1}$ and Lemma~\ref{consdsdpcosdldsx},  we see from \eqref{psd22323dsd} and \eqref{psiisds2sd} that
\begin{equation}\label{qkrwoalssd1}
\begin{aligned}
\mathfrak{M}^{\sharp}_{d_i(\phi_n^{-1})(i_0)[\ihat]}(0,s_0)&=\mathfrak{M}^{\sharp}_{d_i(\tilde{\psi}_n)(i_0)[\ihat]}(0,s_0) \\&\le_{\mathtt{pe}}\epsilon^{3-4b}N_n^{2\tau+1}N_{n-1}^{-\sigma_1}\rVert \ihat \rVert_{s_0 + \mu_{\mathtt{p},1}},\\
\sup_{\vec{b}\in \mathbb{N}^\nu_0,\ |\vec{b}|=\mathtt{b}_0}\mathfrak{M}^{\sharp}_{\partial_{\varphi}^{\vec{b}}(d_i(\phi_n^{-1})(i_0)[\ihat])}(0,s_0) &= \sup_{\vec{b}\in \mathbb{N}^\nu_0,\ |\vec{b}|=\mathtt{b}_0}\mathfrak{M}^{\sharp}_{\partial_{\varphi}^{\vec{b}}(d_i(\tilde{\psi}_n)(i_0)[\ihat])}(0,s_0) \\&\le_{\mathtt{pe}}\epsilon^{3-4b}N_n^{2\tau+1}N_{n-1}\rVert \ihat \rVert_{s_0 + \mu_{\mathtt{p},1}}.
\end{aligned}
\end{equation}

 \textbf{Estimates for $Op^W(r_n)$.}
 Using Lemma~\ref{rjsdj2jsdsd}, we have that (in the applications of the estimates \eqref{rsdj1sd11} and \eqref{rsdj1sd12} below,  we use the fact that $\Omega_{n}^\gamma(i_0)\subset \Omega^{\gamma}_{n-1}(i_0)$, which follows from the definition of the set of frequencies in \eqref{sddpooo3}), 
 \begin{equation}\label{cjdrnrwkd}
\begin{aligned}
\mathfrak{M}_{Op^W(r_n)}^{\sharp,\gamma}(0,s)&\le_{\mathtt{pe},s}|r_{n}|^{\Lip(\gamma,\Omega^{\gamma}_{n}(i_0))}_{0,0,0} \overset{\eqref{rsdj1sd11}}{\le}_{\mathtt{pe},s}\mathbb{M}_0(s_0) N_{n-1}^{-\sigma_1},\\
\mathfrak{M}^{\sharp}_{d_i(Op^W(r_n))[\ihat]}(0,s_0)&\le_{\mathtt{pe},s}|d_ir_n(i_0)[\ihat]|^{\sup}_{0,0,0}\overset{\eqref{rsdj1sd12}}{\le}_{\mathtt{pe},s}\epsilon^3\gamma^{-1}N_{n-1}^{-\sigma_1}\rVert \ihat \rVert_{s_0+\mu_{\mathtt{p},1}}.
\end{aligned}
\end{equation}

\textbf{Estimates for $\Pi_{\ge N_n}\mathcal{R}_n$.}
For $ \Pi_{\ge N_n}\mathcal{R}_n$ in \eqref{sddpooo6}, it follows from Lemma~\ref{smoothing1123} that
\begin{equation}\label{ehlswkdWlro}
\begin{aligned}
 \mathfrak{M}^{\sharp,\gamma}_{\Pi_{\ge N_{n}}\mathcal{R}_n}(0,s)&\le  N_{n}^{-\mathtt{b}_0}   \sup_{\vec{b}\in \mathbb{N}^\nu_0,\ |\vec{b}| = \mathtt{b}_0}\mathfrak{M}^{\sharp,\gamma}_{\partial_\varphi^{\vec{b}}(\mathcal{R}_n)}(0,s) \\ & \overset{\eqref{sdjsds11ram2}}\le \mathbb{M}_0(s) N_n^{-\mathtt{b}_0}N_{n-1},\\
   \sup_{\vec{b}\in \mathbb{N}^\nu_0,\ |\vec{b}| = \mathtt{b}_0}\mathfrak{M}^{\sharp,\gamma}_{\partial_\varphi^{\vec{b}}(\Pi_{\ge N_{n}}\mathcal{R}_n)}(0,s)&\le  \sup_{\vec{b}\in \mathbb{N}^\nu_0,\ |\vec{b}| = \mathtt{b}_0}\mathfrak{M}^{\sharp,\gamma}_{\partial_\varphi^{\vec{b}}(\mathcal{R}_n)}(0,s) \\
   & \overset{\eqref{sdjsds11ram2}}\le \mathbb{M}_0(s) N_{n-1},\\
   \mathfrak{M}^{\sharp}_{d_i(\Pi_{\ge N_{n}}\mathcal{R}_n)(i_0)[\ihat]}(0,s_0)&\le N_{n}^{-\mathtt{b}_0} \sup_{\vec{b}\in \mathbb{N}^\nu_0,\ |\vec{b}| = \mathtt{b}_0}\mathfrak{M}^{\sharp}_{\partial_\varphi^{\vec{b}}(d_i(\mathcal{R}_n)(i_0)[\ihat])}(0,s_0)\\
   & \overset{\eqref{sdjsds11ram4}}\le_{\mathtt{pe}} \epsilon^3\gamma^{-1}N_{n-1}N_n^{-\mathtt{b}_0}\rVert \ihat \rVert_{s_0+\mu_\mathtt{p}},\\
     \sup_{\substack{\vec{b}\in \mathbb{N}^\nu_0 \\ |\vec{b}| = \mathtt{b}_0}}\mathfrak{M}^{\sharp}_{\partial_\varphi^{\vec{b}}(d_i(\Pi_{\ge N_{n}}\mathcal{R}_n)(i_0)[\ihat])}(0,s_0)
    &\le_{\mathtt{pe}}   \sup_{\vec{b}\in \mathbb{N}^\nu_0,\ |\vec{b}| = \mathtt{b}_0}\mathfrak{M}^{\sharp}_{\partial_\varphi^{\vec{b}}(d_i(\mathcal{R}_n)(i_0)[\ihat])}(0,s_0) \\
    & \overset{\eqref{sdjsds11ram4}}\le_{\mathtt{pe}} \epsilon^3\gamma^{-1}N_{n-1}\rVert \ihat \rVert_{s_0+\mu_{\mathtt{p},1}}.
\end{aligned}
\end{equation}
    
\textbf{Estimates for $\mathcal{R}_n\psi_n$.}
    For $\mathcal{R}_n\psi_n$ in \eqref{sddpooo6}, we use \eqref{chanrul2} to obtain
    \begin{align}
     \mathfrak{M}^{\sharp,\gamma}_{\mathcal{R}_n\psi_n}(0,s)&\overset{\eqref{jsdj2sdsd1sd1},\eqref{sdjsds11ram1}}\le_{\mathtt{pe}}\epsilon^2\gamma^{-2}N_{n}^{2\tau+1}N_{n-1}^{-2\sigma_1}\mathbb{M}_0(s_0)\mathbb{M}_0(s)\nonumber\\
     & \overset{\eqref{m01sd2est}}\le_{\mathtt{pe}}\epsilon^{11-10b}N_{n}^{2\tau+1}N_{n-1}^{-2\sigma_1}\mathbb{M}_0(s),\label{sdjsds11ram1321}
     \end{align}
     and
     \begin{align}
      &\sup_{\vec{b}\in \mathbb{N}^\nu_0,\ |\vec{b}| = \mathtt{b}_0}  \mathfrak{M}^{\sharp,\gamma}_{\partial_\varphi^{\vec{b}}(\mathcal{R}_n\psi_n)}(0,s) \nonumber\\
      & \quad \le_{\mathtt{pe}}\left(    \sup_{\vec{b}\in \mathbb{N}^\nu_0,\ |\vec{b}| = \mathtt{b}_0}\mathfrak{M}^{\sharp,\gamma}_{\partial_{\varphi}^{\vec{b}}(\mathcal{R}_n)}(0,s)\mathfrak{M}^{\sharp,\gamma}_{\psi_n}(0,s_0)\right. \nonumber \\ & \qquad \qquad\qquad\left.+   \sup_{\vec{b}\in \mathbb{N}^\nu_0,\ |\vec{b}| = \mathtt{b}_0}\mathfrak{M}^{\sharp,\gamma}_{\partial_{\varphi}^{\vec{b}}(\mathcal{R}_n)}(0,s_0)\mathfrak{M}^{\sharp,\gamma}_{\psi_n}(0,s) \right.\nonumber \\ & \qquad \qquad\qquad\qquad\left.+   \mathfrak{M}^{\sharp,\gamma}_{\mathcal{R}_n}(0,s)   \sup_{\vec{b}\in \mathbb{N}^\nu_0,\ |\vec{b}| = \mathtt{b}_0}\mathfrak{M}^{\sharp,\gamma}_{\partial_{\varphi}^{\vec{b}}(\psi_n)}(0,s_0)\right. \nonumber \\ & \qquad \qquad\qquad\qquad\qquad \left.+  \mathfrak{M}^{\sharp,\gamma}_{\mathcal{R}_n}(0,s_0)   \sup_{\vec{b}\in \mathbb{N}^\nu_0,\ |\vec{b}| = \mathtt{b}_0}\mathfrak{M}^{\sharp,\gamma}_{\partial_{\varphi}^{\vec{b}}(\psi_n)}(0,s)\right)\nonumber\\
& \quad  \overset{\eqref{sdjsds11ram1},\eqref{sdjsds11ram2},\eqref{jsdj2sdsd1sd1},\eqref{jsdj2sdsd1sd2}}\le_{\mathtt{pe}}\epsilon^2\gamma^{-2}N_{n}^{2\tau+1}N_{n-1}^{1-\sigma_1}\mathbb{M}_0(s_0)\mathbb{M}_0(s)\nonumber\\
& \quad  \overset{\eqref{m01sd2est}}\le_{\mathtt{pe}}\epsilon^{11-10b}N_{n}^{2\tau+1}N_{n-1}^{1-\sigma_1}\mathbb{M}_0(s).\label{rjsj2dsd2211sd}
\end{align}
Similarly, we have
\begin{align}
&\mathfrak{M}^{\sharp}_{d_i(\mathcal{R}_n\psi_{n})(i_0)[\ihat]}(0,s_0)\nonumber\\
 &\le_{\mathtt{pe}}\left( \mathfrak{M}^{\sharp}_{d_i(\mathcal{R}_n)(i_0)[\ihat]}(0,s_0)\mathfrak{M}^{\sharp}_{\psi_{n}}(0,s_0)   +  \mathfrak{M}^{\sharp}_{d_i(\psi_n)(i_0)[\ihat]}(0,s_0)\mathfrak{M}^{\sharp}_{\mathcal{R}_{n}}(0,s_0)\right) \nonumber \\
& \overset{\substack{\eqref{sdjsds11ram3},\eqref{sdjsds11ram1} \\ \eqref{sjdjsddrhcnwkdWlp},\eqref{psiisds2sd}}}{\le_{\mathtt{pe}}} \left( \epsilon^{3}\gamma^{-1}\epsilon^{11-10b} + \epsilon^{3-4b}\mathbb{M}_0(s_0)\right)\epsilon^{11-10b}N_{n}^{2\tau+1}N_{n-1}^{-2\sigma_1}\rVert \ihat \rVert_{s_0+\mu_{\mathtt{p},1}}\nonumber \\
& \overset{\eqref{m01sd2est}}{\le_{\mathtt{pe}}}\epsilon^{3}\gamma^{-1} \epsilon^{11-10b}N_{n}^{2\tau+1}N_{n-1}^{-2\sigma_1}\rVert \ihat \rVert_{s_0+\mu_{\mathtt{p},1}},
 \end{align}
 and using \eqref{sdjsds11ram1}, \eqref{sdjsds11ram2}, \eqref{sdjsds11ram3}, \eqref{sdjsds11ram4}, \eqref{jsdj2sdsd1sd1}, \eqref{jsdj2sdsd1sd2}, \eqref{psiisds2sd} and \eqref{m01sd2est},
 \begin{align}\label{rjskds2sda11}
  & \sup_{\vec{b}\in \mathbb{N}^\nu_0,\ |\vec{b}| = \mathtt{b}_0}\mathfrak{M}^{\sharp}_{\partial_\varphi^{\vec{b}}(d_i(\mathcal{R}_n\psi_n)(i_0)[\ihat])}(0,s_0)\nonumber\\
   & \quad \overset{\eqref{chanrul2}}{\le_{\mathtt{pe}}}\left(  \sup_{\vec{b}\in \mathbb{N}^\nu_0,\ |\vec{b}| = \mathtt{b}_0}\mathfrak{M}^{\sharp}_{\partial_\varphi^{\vec{b}}(d_i(\mathcal{R}_n)(i_0)[\ihat])}(0,s_0) \mathfrak{M}^{\sharp}_{\psi_n}(0,s_0) \right.\nonumber \\
   & \qquad \qquad  + \sup_{\vec{b}\in \mathbb{N}^\nu_0,\ |\vec{b}| = \mathtt{b}_0}\mathfrak{M}^{\sharp}_{\partial_\varphi^{\vec{b}}(\mathcal{R}_n)}(0,s_0) \mathfrak{M}^{\sharp}_{d_i\psi_n(i_0)[\ihat]}(0,s_0) \nonumber\\
  &\qquad \qquad \qquad \qquad   +  \sup_{\vec{b}\in \mathbb{N}^\nu_0,\ |\vec{b}| = \mathtt{b}_0}\mathfrak{M}^{\sharp}_{\partial_\varphi^{\vec{b}}(d_i(\psi_n)(i_0)[\ihat])}(0,s_0) \mathfrak{M}^{\sharp}_{\mathcal{R}_n}(0,s_0) \nonumber  \\
  & \left.\qquad \qquad \qquad \qquad \qquad \qquad + \sup_{\vec{b}\in \mathbb{N}^\nu_0,\ |\vec{b}| = \mathtt{b}_0}\mathfrak{M}^{\sharp}_{\partial_\varphi^{\vec{b}}(\psi_n)}(0,s_0) \mathfrak{M}^{\sharp}_{d_i\mathcal{R}_n(i_0)[\ihat]}(0,s_0)  
   \right)\nonumber\\
   &\quad  \ \le_{\mathtt{pe}} \epsilon^{3}\gamma^{-1}\epsilon^{11-10b}N_{n}^{2\tau+1}N_{n-1}^{1-\sigma_1}\rVert \ihat \rVert_{s_0+\mu_{\mathtt{p},1}}.
 \end{align}
  
 Lastly, using Lemma~\ref{consdsdpcosdldsx}, \eqref{m01sd2est} and the definition of $\mathcal{R}_{n+1}$ in \eqref{sddpooo6}, we see from the above estimates \eqref{psiisds2sd}, \eqref{ptildesd1}, \eqref{ptildesd2}, \eqref{qkrwoalssd1}, \eqref{cjdrnrwkd}, \eqref{ehlswkdWlro}, \eqref{sdjsds11ram1321}, \eqref{rjsj2dsd2211sd} and \eqref{rjskds2sda11}  that
 \begin{equation}\label{rkawkqufsd2}
 \begin{aligned}
  \mathfrak{M}^{\sharp,\gamma}_{\mathcal{R}_{n+1}}(0,s)&\le_{\mathtt{pe},s}\left(\epsilon^{11-10b}N_{n}^{2\tau+1}N_{n-1}^{-2\sigma_1} + N_{n}^{-\mathtt{b}_0}N_{n-1}\right)\mathbb{M}_0(s),\\
   \sup_{\substack{\vec{b}\in \mathbb{N}^\nu_0 \\  |\vec{b}| = \mathtt{b}_0}}\mathfrak{M}^{\sharp,\gamma}_{\partial_\varphi^{\vec{b}}(\mathcal{R}_{n+1})}(0,s)&\le_{\mathtt{pe},s}\left(\epsilon^{11-10b}N_{n}^{2\tau+1}N_{n-1}^{1-\sigma_1} + N_{n-1} \right)\mathbb{M}_0(s), \\
   \mathfrak{M}^{\sharp}_{d_i(\mathcal{R}_{n+1})(i_0)[\ihat]}(0,s_0)\le_{\mathtt{pe}}&(\epsilon^{11-10b}N_{n}^{2\tau+1}N_{n-1}^{-2\sigma_1} + N_n^{-\mathtt{b}_0}N_{n-1}) \epsilon^3\gamma^{-1}\rVert \ihat \rVert_{s_0+\mu_{\mathtt{p},1}},\\
    \sup_{\substack{\vec{b}\in \mathbb{N}^\nu_0 \\  |\vec{b}| = \mathtt{b}_0}}\mathfrak{M}^{\sharp}_{\partial_\varphi^{\vec{b}}(d_i\mathcal{R}_{n+1}(i_0)[\ihat])}(0,s_0)&\le_{\mathtt{pe}}(\epsilon^{11-10b}N_{n}^{2\tau +1}N_{n-1}^{1-\sigma_1} + N_{n-1})\epsilon^3\gamma^{-1}\rVert \ihat \rVert_{s_0+\mu_{\mathtt{p},1}}.
 \end{aligned}
 \end{equation}
  Using \eqref{parameterssd1sd}, \eqref{para1chs1} and \eqref{para1chs12}, it is clear that 
  \begin{align*}
  \epsilon^{11-10b}N_{n}^{2\tau+1}N_{n-1}^{-2\sigma_1} + N_{n}^{-\mathtt{b}_0}N_{n-1} &\le_{\mathtt{pe}} (\epsilon^{\frac{1}{22}} + N_0^{1-\sigma_1})N_{n}^{-\sigma_1},\\
  \epsilon^{11-10b}N_{n}^{2\tau+1}N_{n-1}^{1-\sigma_1} + N_{n-1}  &=  \left(\epsilon^{11-10b}N_{n}^{2\tau+1}N_{n-1}^{-\sigma_1} + 1\right)\frac{N_{n-1}}{N_{n}}N_{n}\\
  & \overset{\eqref{para1chs12}}\le \left(\epsilon^{\frac{1}{22}} + 1\right)\frac{N_{-1}}{N_0}N_n\\
  & \le 2N_0^{-1}N_n.
  \end{align*}
 Plugging these into  \eqref{rkawkqufsd2}, we have that
 \begin{align*}
 \mathfrak{M}^{\sharp,\gamma}_{\mathcal{R}_{n+1}}(0,s)&\le_{\mathtt{pe},s}(\epsilon^{\frac{1}{22}} + N_0^{1-\sigma_1})N_{n}^{-\sigma_1}\mathbb{M}_0(s),\\
   \sup_{\vec{b}\in \mathbb{N}^\nu_0,\ |\vec{b}| = \mathtt{b}_0}\mathfrak{M}^{\sharp,\gamma}_{\partial_\varphi^{\vec{b}}(\mathcal{R}_{n+1})}(0,s)&\le_{\mathtt{pe},s} N_0^{-1}N_n\mathbb{M}_0(s), \\
   \mathfrak{M}^{\sharp}_{d_i(\mathcal{R}_{n+1})(i_0)[\ihat]}(0,s_0)&\le_{\mathtt{pe}}(\epsilon^{\frac{1}{22}} + N_0^{1-\sigma_1})N_{n}^{-\sigma_1} \epsilon^3\gamma^{-1}\rVert \ihat \rVert_{s_0+\mu_{\mathtt{p},1}},\\
    \sup_{\vec{b}\in \mathbb{N}^\nu_0,\ |\vec{b}| = \mathtt{b}_0}\mathfrak{M}^{\sharp}_{\partial_\varphi^{\vec{b}}(d_i\mathcal{R}_{n+1}(i_0)[\ihat])}(0,s_0)&\le_{\mathtt{pe}}N_0^{-1}N_n\epsilon^3\gamma^{-1}\rVert \ihat \rVert_{s_0+\mu_{\mathtt{p},1}}.
 \end{align*}
    Using \eqref{parameterssd1sd},  we can choose $\epsilon$ sufficiently small depending on $\mathtt{pe}$ and $\mathtt{S}$, so that $N_0^{-1}$ and $N_0^{1-\sigma_1}$ are as small as we need. This gives the desired estimates \eqref{sdjsds11ram1}-\eqref{sdjsds11ram4}, for $n+1$.

\textbf{Inductive step for \ref{kamrdu4}.} 
Assuming \ref{kamrdu3}, \ref{kamrdu4}, and \ref{kamrdu2} are true for $n$, we aim to prove \ref{kamrdu4} for $n\to n+1$.  
Using \eqref{sddpooo4}, we see that  $Op^W(r_{n+2})$ is also real and reversible, since   $\mathcal{R}_{n+1}$ is real and reversible and $Op^W(r_{n+2})$ is nothing but a restriction to the diagonal part of $\mathcal{R}_{n+1}$. For the estimates, we use  Lemma~\ref{stsdjwdsymbodlsd} and \eqref{sddpooo4} to see that  
\begin{align*}
|r_{n+2}|^{\Lip(\gamma,\Omega^{\gamma}_{n+1}(i_0))}_{0,0,0}&\le_{\mathtt{pe}} \mathfrak{M}^{\sharp,\gamma}_{\mathcal{R}_{n+1}}(0,s_0),\\
|d_ir_{n+2}(i_0)[\ihat]|^{\Lip(\gamma,\Omega^{\gamma}_{n+1}(i_0))}_{0,0,0}&\le_{\mathtt{pe}} \mathfrak{M}^{\sharp,\gamma}_{d_i\mathcal{R}_{n+1}(i_0)[\ihat]}(0,s_0),
\end{align*}
where the second estimate follows from  the relation between $r_{n+1}$ and $\mathcal{R}_n$ being linear. Hence, using \eqref{sdjsds11ram1}, \eqref{sdjsds11ram3} for $n\to n+1$, which have been already proved above, and \eqref{m01sd2est}, we have
\begin{align*}
|r_{n+2}|^{\Lip(\gamma,\Omega^{\gamma}_{n+1}(i_0))}_{0,0,0}&\le_{\mathtt{pe}}\epsilon^{9-6b}N_{n}^{-\sigma_1},\\
|d_ir_{n+2}(i_0)[\ihat]|^{\Lip(\gamma,\Omega^{\gamma}_{n+1}(i_0))}_{0,0,0}&\le_{\mathtt{pe}}\epsilon^{3}\gamma^{-1}N_{n}^{-\sigma_1}\rVert \ihat \rVert_{s_0+\mu_{\mathtt{p},1}},
\end{align*}
and this gives \eqref{rsdj1sd11} and \eqref{rsdj1sd12}. Recalling the definitions of $r_{\le n+2}$ from \eqref{sddpooo5} and noting that $\sum_{n\ge 0}N_{n-1}^{-\sigma_1} < \infty$, the estimates for $r_{\le n+2}$ in \eqref{jjsd2sdsdsdsdsdsdsdsdsdsdsds} and \eqref{sdsdsd2sd1sxf}  follow immediately.

\textbf{Inductive step for \ref{kamrdu2}.} 
By the induction hypothesis, we already know that $D_{n+1}$ is real and reversible, and the estimates \eqref{jjsd2sdsdsdsdsdsdsdsdsdsdsds} and \eqref{sdsdsd2sd1sxf} hold true for $r_{\le n+1}$. Furthermore, we have already proved that $\mathcal{R}_{n+1}$ is also real, reversible and satisfying \eqref{sjowk1sdsds1}. Therefore, we apply Lemma~\ref{homols22d2sd}  to the homological equation in \eqref{sddpooo2}, which gives us that $\psi_{n+1}$ is well-defined for $\omega\in \Omega_{n+2}^{\gamma}(i_0)$ and it is real, reversibility preserving and satisfies \eqref{sjowk1sdsds2}. Furthermore, \eqref{homoglogisd1}with the estimates \eqref{sdjsds11ram1} for $\mathcal{R}_{n+1}$, which we have already proved, gives us the estimate \eqref{psdjsds11ram1} for $\psi_{n+1}$. 

\textbf{Inductive step for \ref{kamrdu5}.} 
 Using the induction hypothesis, we have that
\[
\Omega_{n+1}^{2\gamma}(i_1)\overset{\eqref{sddpooo3}}\subset \Omega^{2\gamma}_{n}(i_1)\overset{\eqref{tlwoclaosd}}\subset \Omega^{2(\gamma-\delta)}_n(i_2) \overset{\because\ \delta<\gamma/2}\subset \Omega_n^{\gamma}(i_2)\overset{\eqref{sddpooo3}}\subset \Omega_{n-1}^{\gamma}(i_2).
\]
Therefore denoting the dependence of $d_n$ on the embedding, we have that  $r_{n}(i_2)(\omega,j)$ and $r_{n+1}(i_1)(\omega,j)$ are well-defined for $\omega\in \Omega_{n+1}^{2\gamma}(i_1)$ (indeed, $r_{n+1}(i_1)$ is well-defined on $\Omega^{\gamma}_{n}(i_1)\overset{\eqref{sddpooo3}}\supset \Omega^{2\gamma}_{n}(i_1)$ and $r_{n}(i_2)$ is well-defined on $\Omega_{n-1}^{\gamma}(i_2))$. For such $\omega$, we have that
\begin{align}
|&\left(d_{n}(i_2)(\omega,j)-d_n(i_2)(\omega,k) \right)-\left(d_{n}(i_1)(\omega,j)-d_n(i_1)(\omega,k) \right)|\nonumber\\
& \le |\left(d_{0}(i_2)(\omega,j)-d_0(i_2)(\omega,k) \right)-\left(d_{0}(i_1)(\omega,j)-d_0(i_1)(\omega,k) \right)| \nonumber \\
& +  2\sup_{j\in S^\perp}|r_{\le n}(i_1)(\omega,j)-r_{\le n}(i_2)(\omega,j)|.\label{morningtrain}
\end{align}
From \eqref{d0defsd22}, we have that
\begin{align*}
|&\left(d_{0}(i_2)(\omega,j)-d_0(i_2)(\omega,k) \right)-\left(d_{0}(i_1)(\omega,j)-d_0(i_1)(\omega,k) \right)| \\
& \le|(\mathtt{m}_\alpha(i_2(\omega))(\omega) - \mathtt{m}_\alpha(i_1(\omega))(\omega))(\lambda_\alpha(j)-\lambda_\alpha(k))|\\
& \ + |\left( \mathfrak{m}_1(i_2)(\omega,j)-\mathfrak{m}_1(i_2)(\omega,k)\right) - \left( \mathfrak{m}_1(i_1)(\omega,j)-\mathfrak{m}_1(i_1)(\omega,k)\right)|\\
& \overset{\eqref{nusd2sdsd}}\le_{\mathtt{pe}}\epsilon^2 \rVert i_1-i_2\rVert_{s_0+\mu_{\mathtt{p},1}}|\lambda_\alpha(j)-\lambda_\alpha(k)| \\
& \qquad + \int_j^k \left | \partial_\xi\mathfrak{m}_1(i_1)(\omega,\xi) - \partial_\xi\mathfrak{m}_1(i_2)(\omega,\xi) \right|d\xi\\
&\overset{ \eqref{rejsmsdim1} }\le_{\mathtt{pe}}\left( \epsilon^2 \rVert i_1-i_2\rVert_{s_0+\mu_{\mathtt{p},1}}|\lambda_\alpha(j)-\lambda_\alpha(k)| + \epsilon^2|k-j| \rVert i_1-i_2\rVert_{s_0+\mu_{\mathtt{p},1}}  \right)\\
& \overset{\text{ Lemma~\ref{melnikove_23s}}}\le_{\mathtt{pe}} \epsilon^2 \rVert i_1-i_2\rVert_{s_0+\mu_{\mathtt{p},1}}|\lambda_\alpha(j)-\lambda_\alpha(k)|.
\end{align*}
Plugging this into \eqref{morningtrain} and using \eqref{sdsdsd2sd1sxf} for $r_{\le n}$, we obtain
\begin{align*}
|&\left(d_{n}(i_2)(\omega,j)-d_n(i_2)(\omega,k) \right)-\left(d_{n}(i_1)(\omega,j)-d_n(i_1)(\omega,k) \right)|\nonumber\\
&\le_{\mathtt{pe}} \epsilon^2 \rVert i_1-i_2\rVert_{s_0+\mu_{\mathtt{p},1}}|\lambda_\alpha(j)-\lambda_\alpha(k)| + \epsilon^3\gamma^{-1}\rVert i_1 - i_2\rVert_{s_0 + \mu_{\mathtt{p},1}}\\
&\le_{\mathtt{pe}}\epsilon^3\gamma^{-1}\rVert i_1 - i_2\rVert_{s_0 + \mu_{\mathtt{p},1}}|\lambda_\alpha(j)-\lambda_\alpha(k)|\\
&\le C({\mathtt{pe}})\epsilon^3\gamma^{-1}\rVert i_1-i_2\rVert_{s_0+\mu_{\mathtt{p},1}}|\lambda_\alpha(j)-\lambda_\alpha(k)|,
\end{align*}
for some $C(\mathtt{pe})$.
Therefore, if $$C(\mathtt{pe})\epsilon^3\gamma^{-1}N_{n}^{\tau}\rVert i_1-i_2\rVert_{s_0+\mu_{\mathtt{p},1}}\le \delta,$$ then we have that
\begin{align}\nonumber
|\left(d_{n}(i_2)(\omega,j)-d_n(i_2)(\omega,k) \right)&-\left(d_{n}(i_1)(\omega,j)-d_n(i_1)(\omega,k) \right)|\\
& \le \delta N_n^{-\tau}|\lambda_\alpha(j)-\lambda_\alpha(k)|.\label{sdskxdnwnfmfcxc}
\end{align}
Thus, for $\omega\in \Omega_{n+1}^{2\gamma}(i_1)$ and $|l|\le N_n$,
\begin{align*}
|&\ii \omega - (d_{n}(i_2)(\omega,j) - d_n(i_2)(\omega,k))|\\
& \ge |\ii \omega - (d_{n}(i_1)(\omega,j) - d_n(i_1)(\omega,k))| \\
& \qquad \qquad - |\left(d_{n}(i_2)(\omega,j)-d_n(i_2)(\omega,k) \right)-\left(d_{n}(i_1)(\omega,j)-d_n(i_1)(\omega,k) \right)|\\
& \overset{\because\omega\in \Omega_{n+1}^{2\gamma}(i_1) }\ge  \gamma|l|^{-\tau}|\lambda_\alpha(j)-\lambda_\alpha(k)| \\
& \qquad \qquad -|\left(d_{n}(i_2)(\omega,j)-d_n(i_2)(\omega,k) \right)-\left(d_{n}(i_1)(\omega,j)-d_n(i_1)(\omega,k) \right)|\\
& \overset{\eqref{sdskxdnwnfmfcxc}}\ge (\gamma-\delta)|l|^{-\tau}|\lambda_\alpha(j)-\lambda_\alpha(k)|.
\end{align*}
which implies $\omega \in \Omega_{n+1}^{2(\gamma-\delta)}(i_2)$, thus, $\Omega_{n+1}^{2\gamma}(i_1)\subset \Omega^{2(\gamma-\delta)}_{n+1}(i_2)$.

\vspace{0.5\baselineskip}\noindent\textit{Proof of \eqref{consjdsdsd2sd}.}  From the above induction steps, it follows that $D_{n+1},\mathcal{R}_n,\psi_n, \tilde{\psi}_n$ are all well-defined for $\omega\in \cap_{n\ge 0}\Omega_n^{\gamma}(i_0)$. Then,  \eqref{consjdsdsd2sd} follows from \eqref{consd2sdsdsd}. 

Lastly, since $\phi_n$ is real, reversibility and $\frac{2\pi}{\mathtt{M}}$-translation invariance preserving and $L_0$ is real and reversible and $\frac{2\pi}{\mathtt{M}}$-translation invariant , $L_n$ is real, reversible and $\frac{2\pi}{\mathtt{M}}$-translation invariant as well.  
\end{proof}

 \begin{lemma}\label{kskd22sdsds}
For $\omega \in \cap_{n\ge 0}\Omega^{\gamma}_{n}(i_0)$, let us define (note that each $\phi_n$ depends on $\omega$)
\begin{align}\label{pusjdwsd1ss}
\Phi_{7,0}:=\phi_0,\quad \Phi_{7,n+1}:=\Phi_{7,n}\circ \phi_{n+1}.
\end{align}
Then given $\mathtt{S}\gg s_0$, there exists $\epsilon_0(\mathtt{pe},\mathtt{S})>0$ such that if $\epsilon\in(0,\epsilon_0)$, then there exists  an operator $\Phi_{7,\infty}$ such that 
\begin{align}\label{strongconvers}
\lim_{n\to \infty} \rVert (\Phi_{7,\infty} - \Phi_{7,n})h\rVert_s = 0, \text{ for all $h\in C^\infty_{\varphi,x}$ and $s\in [s_0,\mathtt{S}]$}.
\end{align} Furthermore, $\Phi_{7,\infty}-I$ is Lip-$0$-modulo-tame operator up to $\mathtt{S}$ with the estimates
\begin{align}\label{finallsd1sdsd}
\mathfrak{M}^{\sharp,\gamma}_{\Phi_{7,\infty}-I}(0,s)\le_{\mathtt{pe},\mathtt{S}}\epsilon^2\gamma^{-2}N_0^{2\tau+1}\mathbb{M}_0(s).
\end{align}
Lastly, the same results hold as above for $\phi_{n}{-1}$, that is, $\Phi_{7,\infty}^{-1}$, the inverse of $\Phi_{7,\infty}$, satisfies
\[
\lim_{n\to\infty}\rVert (\Phi_{7,\infty}^{-1} - \Phi_{7,n}^{-1})h\rVert_s =0,\text{ and }\mathfrak{M}^{\sharp,\gamma}_{\Phi_{7,\infty}{-1}-I}(0,s)\le_{\mathtt{S},\mathtt{pe}}\epsilon^2\gamma^{-2}N_0^{2\tau+1}\mathbb{M}_0(s),
\]
for $s\in [s_0,\mathtt{S}]$.
 \end{lemma}
 \begin{proof}
 In the proof, we will denote by $c_{\mathtt{pe}}$ a constant that varies from line to line, depending on $\mathtt{pe}$ only. We will prove the lemma for $\Phi_{7,\infty}$ only, since the same argument applies to $\Phi_{7,\infty}^{-1}$. 
 
   First, we derive a tame constant of  $\Phi_{7,n}$
  We claim that 
  \begin{align}\label{claimfors_0023212}
  \mathfrak{M}^{\sharp,\gamma}_{\Phi_{7,n}}(0,s_0)\le c_{\mathtt{pe}}^{n+1} \text{ for some $c_{\mathtt{pe}}\ge1$ for all $n\ge 0$}.
  \end{align}
   When $n=0$, it follows from \eqref{psdjsds11ram1} that
   \begin{align*}
    \mathfrak{M}^{\sharp,\gamma}_{\Phi_{7,0}}(0,s_0) &\le 1 +  \mathfrak{M}^{\sharp,\gamma}_{\psi_0}(0,s_0)\\
    & \le 1 + c_{\mathtt{pe}}\epsilon^2\gamma^{-2}N_{-1}^{-\sigma_1}N_0^{2\tau+1}\mathbb{M}_0(s_0)\\
    & \overset{\eqref{m01sd2est},\ N_{-1}=1}\le 1 + c_\mathtt{pe}\epsilon^{11-10b}N_{-1}^{-2\sigma_1}N_0^{2\tau+1}\\
    &\overset{\eqref{para1chs12}}\le c_\mathtt{pe}.
   \end{align*}
   Assuming \eqref{claimfors_0023212} for $n\ge 0$, Lemma~\ref{consdsdpcosdldsx} tells us that
   \begin{align*}
   \mathfrak{M}^{\sharp,\gamma}_{\Phi_{7,n+1}}(0,s_0) &=2\mathfrak{M}^{\sharp,\gamma}_{\Phi_{7,n}}(0,s_0)\mathfrak{M}^{\sharp,\gamma}_{\phi_{n+1}}(0,s_0),\\
   &\overset{\eqref{claimfors_0023212}}\le 2c_{\mathtt{pe}}^{n+1}( 1 + \mathfrak{M}^{\sharp,\gamma}_{\psi_{n+1}}(0,s_0))\\
   &\overset{\eqref{psdjsds11ram1}}\le 2c_{\mathtt{pe}}^{n+1}\left(1 +c_{\mathtt{pe}}\epsilon^2\gamma^{-2}N_{n}^{-\sigma_1}N_{n+1}^{2\tau+1}\mathbb{M}_0(s_0) \right)\\
   &\overset{\eqref{m01sd2est}}\le 2c_{\mathtt{pe}}^{n+1}(1+ c_{\mathtt{pe}}\epsilon^{11-10b}N_{n}^{-\sigma_1}N_{n+1}^{2\tau+1})\\
   &\overset{\eqref{para1chs12}}\le2c_{\mathtt{pe}}^{n+1}\left( 1+\epsilon^{\frac{1}{22}}c_{\mathtt{pe}}\right)
   \end{align*}
    Therefore, assuming that $c_{\mathtt{pe}}$ is large enough depending on $\mathtt{pe}$, we obtain \eqref{claimfors_0023212} for $n+1$.
   
   Now, we claim that 
   \begin{align}\label{sd2sdsdasd}
   \mathfrak{M}^{\sharp,\gamma}_{\Phi_{7,n}}(0,s)\le c_{\mathtt{pe}}^{n+1}(1+\epsilon^2\gamma^{-2}N^{-\sigma_1}_{n-1}N_n^{2\tau+1}\mathbb{M}_0(s)),
   \end{align}
   for some constant $c_{\mathtt{pe}}>0$.
   When $n=0$, the proof is identical to the case when $s=s_0$ as above. For $n\ge 0$, Lemma~\ref{consdsdpcosdldsx} yields that
   \begin{align*}
   \mathfrak{M}^{\sharp,\gamma}_{\Phi_{7,{n+1}}}(0,s) &\le \mathfrak{M}^{\sharp,\gamma}_{\Phi_{7,n}}(0,s)\mathfrak{M}^{\sharp,\gamma}_{\phi_{n+1}}(0,s_0)+\mathfrak{M}^{\sharp,\gamma}_{\Phi_{7,n}}(0,s_0)\mathfrak{M}^{\sharp,\gamma}_{\phi_{n+1}}(0,s)\\
   & \overset{\eqref{claimfors_0023212}}\le  \mathfrak{M}^{\sharp,\gamma}_{\Phi_{7,{n}}}(0,s)(1+ \mathfrak{M}^{\sharp,\gamma}_{\psi_{{n+1}}}(0,s_0)) +c_{\mathtt{pe}}^{n+1}(1+\mathfrak{M}^{\sharp,\gamma}_{\psi_{n+1}}(0,s))\\
   &\overset{\eqref{psdjsds11ram1},\eqref{sd2sdsdasd}}\le 4c_{\mathtt{pe}}^{n+1}(1+\epsilon^2\gamma^{-2}N_{n}^{-\sigma_1}N_{n+1}^{2\tau+1}\mathbb{M}_0(s)), \\
      \end{align*}
      where we used $\mathfrak{M}^{\sharp,\gamma}_{\psi_{{n+1}}}(0,s_0) \le 1$ for sufficiently small $\epsilon>0$, which follows from \eqref{sjdjsddrhcnwkdWlp}. Assuming that $c_{\mathtt{pe}}\ge 4$, we obtain \eqref{sd2sdsdasd} for $n+1$.
   Therefore, we have that  for some $c_{\mathtt{pe}}>0$, 
  \begin{align}\nonumber
 &\mathfrak{M}^{\sharp,\gamma}_{\Phi_{7,n}\circ \psi_{n+1}} (0,s)\\
 &\overset{\text{Lemma~\ref{consdsdpcosdldsx}}}\le  \mathfrak{M}^{\sharp,\gamma}_{\Phi_{7,n}} (0,s) \mathfrak{M}^{\sharp,\gamma}_{ \psi_{n+1}} (0,s_0) +   \mathfrak{M}^{\sharp,\gamma}_{\Phi_{7,n}} (0,s_0) \mathfrak{M}^{\sharp,\gamma}_{ \psi_{n+1}} (0,s)\nonumber\\
  & \overset{\eqref{sd2sdsdasd}, \eqref{psdjsds11ram1}}\le 2c_\mathtt{pe}^{n+1}\epsilon^2\gamma^{-2}N_n^{-\sigma_1}N_{n+1}^{2\tau+1} \nonumber \\
  & \qquad \qquad \qquad \times \left( \mathbb{M}_0(s_0) + (1 + \underbrace{\epsilon^2\gamma^{-2}N^{-\sigma_1}_{n-1}N_n^{2\tau+1}\mathbb{M}_0(s_0)}_{\le 1, \because \eqref{m01sd2est}, \eqref{para1chs12}})\mathbb{M}_0(s)\right),\nonumber\\
  &\overset{\mathbb{M}_0(s_0)\le \mathbb{M}_0(s)}\le 6c_{\mathtt{pe}}^{n+1}\epsilon^2\gamma^{-2}N_n^{-\sigma_1}N_{n+1}^{2\tau+1} \mathbb{M}_0(s)\nonumber \\
  &\qquad \le 6c_{\mathtt{pe}}^{n+1}\epsilon^2 \gamma^{-2}N_n^{-\frac{\sigma_1}2}\mathbb{M}_0(s),\label{ericcartma22n}
  \end{align} 
  where the last inequality follows from
  \[
  N_n^{-\sigma_1}N_{n+1}^{2\tau+1}\overset{\eqref{parameterssd1sd}} =  N_0^{\chi^n(3/2(2\tau+1)-\sigma_1)} \le N_0^{-\chi^n\frac{\sigma_1}{2}} = N_n^{-\frac{\sigma_1}{2}}.
  \]
Therefore, for all $n\ge0$ and $m\ge1$, noting that
\begin{align}\label{rksdopwogodnseos1}
\Phi_{7,n+m}-\Phi_{7,n} = \sum_{k=n}^{n+m-1}(\Phi_{7,k+1}-\Phi_{7,k}) \overset{\eqref{pusjdwsd1ss},\eqref{sddpooo2}}=  \sum_{k=n}^{n+m-1}\Phi_{7,k}\circ\psi_{k+1},
\end{align}
we see that
\begin{align}\label{sd2sdzxcxc}
\mathfrak{M}^{\sharp,\gamma}_{\Phi_{7,n+m}-\Phi_{7,n}}(0,s)\overset{\eqref{rksdopwogodnseos1}}{\le}\sum_{k=n}^{\infty}\mathfrak{M}^{\sharp,\gamma}_{\Phi_{7,k}\circ \psi_{k+1}}(0,s)\overset{\eqref{ericcartma22n}}\le 6\epsilon^2 \gamma^{-2}\mathbb{M}_0(s)\sum_{k=n}^{\infty}c_{\mathtt{pe}}^{k+1}N_k^{-\frac{\sigma_1}2}.
\end{align}
  For each $\omega \in \omega \in \cap_{n\ge 0}\Omega^{\gamma}_{n}(i_0)$, for which the $\phi_n$'s are well-defined for all $n\ge 0$ (see \ref{kamrdu2} of Proposition~\ref{kaksdsddsd}), we define $\Phi_{7,\infty}=\Phi_{7,\infty}(\omega)$ to be
  \begin{align}\label{limitsd1sd}
  \Phi_{7,\infty}h :=\lim_{n\to \infty} \Phi_{7,n}h, \text{ for $h\in C^\infty_{\varphi,x}$,}
  \end{align}
   in the strong limit in $H^s_{\varphi,x}$. Indeed, we need to show that for each $h\in C^\infty_{\varphi,x}$, the sequence $\Phi_{7,\infty}h$ converges strongly in $H^s_{\varphi,x}$. Using the notation in \eqref{reminiscings1}, it follows from \eqref{sd2sdzxcxc} that  for $\omega,\omega_1,\omega_2\in  \cap_{n\ge 0}\Omega^{\gamma}_{n}(i_0)$,
      \begin{align}
   &\rVert (\Phi_{7,n+m}(\omega) - \Phi_{7,n}(\omega))h \rVert_s  \nonumber \\
   & \le \sum_{k=0}^1\rVert (\underline{\Delta^k_{1,2}\Phi_{7,n+m}} - \underline{\Delta^k_{12}\Phi_{7,n}})h\rVert_s\nonumber\\
   & \le \mathfrak{M}^{\sharp,\gamma}_{\Phi_{7,n+m}-\Phi_{7,n}}(0,s)\rVert h \rVert_{s_0} + \mathfrak{M}^{\sharp,\gamma}_{\Phi_{7,n+m}-\Phi_{7,n}}(0,s_0)\rVert h \rVert_{s}\nonumber\\
   & \le_{\mathtt{pe}}  \epsilon^2\gamma^{-2}(\mathbb{M}_0(s_0)\rVert h \rVert_s + \mathbb{M}_0(s)\rVert h \rVert_{s_0}) \sum_{k=n}^{\infty}c_{\mathtt{pe}}^{k+1}N_k^{-\frac{\sigma_1}2}.  \label{cosd1dssds}
   \end{align}
   Using $\sigma_1$ and $N_n$ in \eqref{parameterssd1sd}, one can easily check that $\sum_{k=0}^{\infty} c_{\mathtt{pe}}^{k+1}N_k^{-\frac{\sigma_1}2} <\infty$ for sufficiently small $\epsilon>0$, since $N_n$ increases to $\infty$ double-exponentially in $n$. Therefore, $n\mapsto \Phi_{7,n}(\omega)h$ is a Cauchy sequence in $H^s_{\varphi,x}$, hence, \eqref{limitsd1sd} is well-defined. Furthermore, the above estimates also show that $n\mapsto \underline{\Delta^k_{12}\Phi_{7,n}}h$ is also a Cauchy sequence for $k=0,1$. Therefore, using the definition of Lip-$0$-modulo tame operators, one can easily check that $\lim_{n\to \infty}\underline{\Delta^k_{12}\Phi_n} h = \underline{\Delta^k_{12}\Phi_{\infty}}h$ in $H^s_{\varphi,x}$ strongly, therefore, plugging $n=0$ and letting $m\to\infty$ in \eqref{cosd1dssds} to obtain that
   \begin{align*}
   \rVert \left(\underline{\Delta^k_{12}\Phi_{7,\infty}} - \underline{\Delta^k_{12}\Phi_{7,0}}\right)h\rVert_s \le c_\mathtt{pe}\epsilon^2\gamma^{-2}(\mathbb{M}_0(s_0)\rVert h \rVert_s + \mathbb{M}_0(s)\rVert h \rVert_{s_0}).
   \end{align*}
 Therefore, taking the supremum in $\omega,\omega_1,\omega_2$ over $\cap_{n\ge 0}\Omega^{\gamma}_{n}(i_0)$, we obtain
 \begin{align}\label{iirpsdsd1sd}
 \mathfrak{M}^{\sharp,\gamma}_{\Phi_{7,\infty}-\Phi_{7,0}}(0,s)\le_{\mathtt{pe}}\epsilon^2\gamma^{-2}\mathbb{M}_0(s).
 \end{align}
 Thus, we get
  \begin{align*}
 \mathfrak{M}^{\sharp,\gamma}_{\Phi_{7,\infty} - I} (0,s)&\le  \mathfrak{M}^{\sharp,\gamma}_{\Phi_{7,\infty} - \Phi_{7,0}}(0,s) +  \mathfrak{M}^{\sharp,\gamma}_{\Phi_{7,0} - I}(0,s) \\
 & \overset{\eqref{sddpooo2}}\le  \mathfrak{M}^{\sharp,\gamma}_{\Phi_{7,\infty} - \Phi_{7,0}}(0,s) +  \mathfrak{M}^{\sharp,\gamma}_{\psi_0} (0,s)\\
 &\overset{\eqref{psdjsds11ram1},\eqref{iirpsdsd1sd}}\le c_{\mathtt{pe}}\epsilon^2\gamma^{-2}\mathbb{M}_0(s) + c_{\mathtt{pe}}\epsilon^2\gamma^{-2}N_0^{2\tau+1}\mathbb{M}_0(s)\\
&\le c_{\mathtt{pe}}\epsilon^2\gamma^{-2}N_0^{2\tau+1}\mathbb{M}_0(s),
 \end{align*}
 which proves \eqref{finallsd1sdsd}.
       \end{proof}

   Now, we denote\index{$d_0$}\index{$d_\infty$}\index{$r_\infty$}\index{$\Omega^{4\gamma}_\infty$}
   \begin{equation}\label{sjdjsd2sdsddenoteasd}    
   \begin{aligned}   r_{\infty}(0,j)&:= \sum_{n=0}^{\infty}r_n,\quad d_{\infty}(i_0)(\omega,j):=d_0(\omega,j) + r_{\infty}(\omega,j),\\ 
   \Omega^{4\gamma}_{\infty}(i_0)&:=\left\{ \omega\in \Omega_1: |\ii \omega - (d_\infty(i_0)(\omega,j) - d_{\infty}(i_0)(\omega,k))| \right.\\
   & \left.\qquad \qquad \qquad \quad \ge 2\gamma |l|^{-\tau}|\lambda_\alpha(j)-\lambda_\alpha(k)|,\  \forall l\ne0, \forall j,k\in S_{\mathtt{M}}^\perp\cup\left\{0\right\}\right\}.
\end{aligned}
\end{equation}

\begin{lemma}\label{jsdsd2reamsd2sd}We have that

\[
|r_{\infty}|^{\Lip(\gamma, \cap_{n\ge 0}\Omega^{\gamma}_{n}(i_0))}_{0,0,0} \le_{\mathtt{pe}}\mathbb{M}_0(s_0),\quad |d_ir_{\infty}(i_0)[\ihat]|^{\sup}_{0,0,0} \le_{\mathtt{pe}}\epsilon^3\gamma^{-1}\rVert \ihat \rVert_{s_0+\mu_{\mathtt{p},1}}.
\]

We also have that 
\begin{align}\label{javisd1sds}
\Omega^{4\gamma}_\infty(i_0)\subset \cap_{n\ge 0}\Omega^\gamma_n(i_0),
\end{align}
therefore $r_{\infty}$ is well-defined for all $\omega \in  \Omega^{4\gamma}_{\infty}(i_0)$.
\end{lemma}
 \begin{proof}
 The estimates for $r_\infty$ follow immediately from \ref{kamrdu4} of Proposition~\ref{kaksdsddsd}.

For each  $\omega\in \Omega^{4\gamma}_{\infty}(i_0)$, we have $\omega\in \Omega_1\subset \Omega^{\gamma}_0(i_0)$ by definition. For $n\ge 0$, we have that  that (recalling the definition of $\Omega^\gamma_{n+1}(i_0)$ from \eqref{sddpooo3}) for $|l|\le N_n$ and $j\ne k$,
 \begin{align}
& |\ii \omega\cdot l - (d_n(\omega,j) - d_n(\omega,k))| \nonumber \\
&\ge |\ii \omega\cdot l - (d_\infty(\omega,j) - d_\infty(\omega,k))| - 2\sum_{k\ge n+1}|r_{k}|^{\sup}_{0,0,0} \nonumber\\
& \overset{\eqref{sjdjsd2sdsddenoteasd},\eqref{rsdj1sd11},\eqref{m01sd2est}}\ge 2\gamma |l|^{-\tau}|\lambda_\alpha(j)-\lambda_\alpha(k)| - c_{\mathtt{pe}}\epsilon^{9-6b}\sum_{k\ge n+1}N_{k-2}^{-\sigma_1}\nonumber\\
&\overset{\text{Lemma~\ref{melnikove_23s}}}\ge \gamma |l|^{-\tau}|\lambda_\alpha(j)-\lambda_\alpha(k)|\left(2 - c_{\mathtt{pe}}\epsilon^{9-8b}|l|^{\tau}\sum_{k\ge n}N_{k-1}^{-\sigma_1} \right).\label{rjdmlRxdlsd}
 \end{align}
 Let us claim 
 \begin{align}\label{clisdqhsds}
 c_{\mathtt{pe}}\epsilon^{9-8b}|l|^{\tau}\sum_{k\ge n}N_{k-1}^{-\sigma_1} \le 1, \text{ for $|l|\le N_n, \ n\ge 0$.}
 \end{align}
 When $n=0$, we have (using $|l|\le N_n=N_0$),
 \begin{align*}
 c_{\mathtt{pe}}\epsilon^{9-8b}|l|^{\tau}\sum_{k\ge 0}N_{k-1}^{-\sigma_1} \le c_{\mathtt{pe}}\epsilon^{9-8b}N_0^{\tau} \overset{\eqref{parameterssd1sd}}\le c_{\mathtt{pe}}\epsilon^{9-8b}(\epsilon^{1-2b})^{\rho\tau} \le c_{\mathtt{pe}}\epsilon^{(9-8b)/2},
 \end{align*}
 where the last inequality follows from \eqref{parameterssd1sd}, which says $\rho^{-1}\ge \frac{2\tau(2b-1)}{9-8b}$. This proves \eqref{clisdqhsds} for $n=0$.
  For $n\ge 1$, one can also use \eqref{parameterssd1sd} to see that
  \begin{align*}
  \sum_{k\ge n}N_{k-1}^{-\sigma_1} & = N_{n-1}^{-\sigma_1}\left( 1 + \sum_{k\ge n+1}\frac{N_{k-1}^{-\sigma_1}}{N_{n-1}^{-\sigma_1}}\right) \\
  & = N_{n-1}^{-\sigma_1}\left(1 + \sum_{k\ge n+1}N_0^{-\sigma_1\chi^{k-1}(1-\chi^{n-k})}\right) \le 2N_{n-1}^{-\sigma_1},
  \end{align*}
  for sufficiently small $\epsilon>0$. Therefore,
 \begin{align*}
c_{\mathtt{pe}}\epsilon^{9-8b}|l|^{\tau}\sum_{k\ge n}N_{k-1}^{-\sigma_1}& \le 2c_{\mathtt{pe}}\epsilon^{9-8b}N_n^{\tau}N_{n-1}^{-\sigma_1} \\
& \le 2c_{\mathtt{pe}}\epsilon^{9-8b}N_0^{\chi^{n-1}(\frac{3}2\tau-\sigma_1)} \le c_{\mathtt{pe}}\epsilon^{9-8b},
 \end{align*}
 for small $\epsilon$, which gives \eqref{clisdqhsds} for $n\ge 1$. Plugging \eqref{clisdqhsds} into \eqref{rjdmlRxdlsd}, we get
  \[
 |\ii \omega\cdot l - (d_n(\omega,j) - d_n(\omega,k))|\ge \gamma |l|^{-\tau}|\lambda_\alpha(j)-\lambda_\alpha(k)|, \text{ for all $|l|\le N_n$},
 \]
 which implies $\omega \in \Omega_{n+1}^\gamma(i_0)$ for all $n\ge 0$. This proves \eqref{javisd1sds}.
 \end{proof}

We summarize the main results of this section:
 \begin{proposition}\label{kamresdsd}
 Given $\mathtt{S}\gg s_0$, there exists $\epsilon_0(\mathtt{pe},\mathtt{S})>0$ such that if $\epsilon\in (0,\epsilon_0)$ and \eqref{size_assumption_2} holds for $\mu=\mu_{\mathtt{p},1}$ and for some $\mathtt{C}>0$, there exist a real, reversibility preserving operator $\Phi_7=\Phi_7(\omega)$  and a diagonal operator $D_\infty(\omega):=\text{diag}_{j\in S_\mathtt{M}^\perp}d_\infty(\omega,j)$ that depend on $\omega\in \Omega_1$ satisfying
 \begin{enumerate}[label=(\arabic*)]
 \item\label{kamredu11} (Eigenvalues\index{Eigenvalues}) $d_{\infty}(\omega,j)= d_\infty(i_0(\omega))(\omega,j) = d_{0}(\omega,j) + r_{\infty}(\omega,j)$, where $r_\infty(\omega,j)$ satisfies
 \begin{align}\label{shsjd2sressdsdsd}
 |r_{\infty}|^{\Lip(\gamma, \Omega_1)}_{0,0,0} \le_{\mathtt{pe}} \epsilon^{9-6b} ,\quad |d_ir_{\infty}(i_0)[\ihat]|^{\sup}_{0,0,0} \le_{\mathtt{pe}}\epsilon^3\gamma^{-1}\rVert \ihat \rVert_{s_0+\mu_{\mathtt{p},1}}. 
 \end{align}
 \item\label{kamredu12} (Conjugation\index{Conjugation})  $\Phi=\Phi_{7},\Phi_7^{-1}$ satisfies that (assuming $\epsilon\in (0,\epsilon_0(\mathtt{pe},\mathtt{S}))$)
 \begin{align}\label{oosdsdboudncsd}
\rVert \Phi h\rVert_s^{\Lip(\gamma,\Omega_1)} \le_{\mathtt{pe},\mathtt{S}} \rVert h\rVert_{s}^{\Lip(\gamma,\Omega_1)} +\epsilon^{7}\gamma^{-4} \rVert \mathfrak{I}_0\rVert_{s+\mu_{\mathtt{p},1}}^{\Lip(\gamma,\Omega_1)}\rVert h\rVert_{s_0},\text{ for $s\in [s_0,\mathtt{S}]$}.
 \end{align}
 \item \label{kamredu13}(Diagonalization\index{Diagonalization}) For $\omega\in \Omega^{4\gamma}_{\infty}(i_0)$, where
 \begin{align*}
 \Omega^{4\gamma}_{\infty}(i_0):=\left\{ \omega\in \Omega_1: \right. & |\ii \omega - (d_\infty(\omega,j) - d_{\infty}(\omega,k))|\\
 & \left.\ge 2\gamma |l|^{-\tau}|\lambda_\alpha(j)-\lambda_\alpha(k)|, \forall l\ne0, \forall j,k\in S_{\mathtt{M}}^\perp\cup\left\{ 0 \right\}\right\},
 \end{align*}
 We have that
 \begin{align}\label{liabssdw2sss2sd2}
 \mathcal{L}_{\infty}:= \Phi_{7}(\omega)^{-1}L_0(\omega)\Phi_7(\omega) = \D_\omega - D_\infty(\omega).
 \end{align}
 Also, $\mathcal{L}_{\infty}$ is real and reversible.
 \item \label{kamredu213}($\frac{2\pi}{\mathtt{M}}$-translation invariance) $\Phi_7$ is $\frac{2\pi}{\mathtt{M}}$-translation invariance preserving.
 \end{enumerate}
 \end{proposition}
\begin{proof}
We already have $r_{\infty},\Phi_{7,\infty},\Phi_{7,\infty}^{-1}$, which are well-defined for $\omega\in \Omega^{4\gamma}_{\infty}(i_0)$ (Lemma~\ref{javisd1sds}). Using Kirszbraun's Theorem \cite{Kirszbraun:zusammenziehende-lipschitzsche-transformationen}, we can extend them to $\omega\in \Omega_1$ with the same Lipschitz constant.  denoting $\Phi_7:=\Phi_{7,\infty}$, where $\Phi_{7,\infty}$ is as in Lemma~\ref{kskd22sdsds}, we see that $\Phi_7^{\pm}$ are real and reversibility preserving, since they are compositions of real, reversibility preserving operators.

Item  \ref{kamredu11} follows from Lemma~\ref{jsdsd2reamsd2sd}.

For item \ref{kamredu12},  we obtain from \eqref{finallsd1sdsd}, \eqref{defofmm2ksdsd} and  \eqref{fromyourbag1} that\index{$\mathcal{L}_\infty$}
\[
\mathfrak{M}^{\sharp,\gamma}_{\Phi_7 - I}(0,s)\le_{\mathtt{pe},\mathtt{S}} \epsilon^2\gamma^{-2}N_0^{2\tau+1}\gamma^{-1}\left(\epsilon^5 + \epsilon^3 \rVert \mathfrak{I}_0\rVert_{s+\mu_{\mathtt{p},1}}^{\Lip(\gamma,\Omega_1)}\right).
\]
Therefore, using \eqref{size_assumption_2},
\begin{align}\label{rkkxcxcppsds11}
\mathfrak{M}^{\sharp,\gamma}_{\Phi_7 - I}(0,s_0)\le_{\mathtt{pe}}\epsilon^{11-10b}N_0^{2\tau+1} \overset{\eqref{para1chs12},n=0}\le 1.
\end{align}
Also, we have 
\begin{align}\label{jsdk2seves}
\mathfrak{M}^{\sharp,\gamma}_{\Phi_7 - I}(0,s) &\le_{\mathtt{pe},\mathtt{S}} \epsilon^{7-4b}N_0^{2\tau+1} + \epsilon^{5}\gamma^{-3}N_0^{2\tau+1}\rVert \mathfrak{I}_0\rVert^{\Lip(\gamma,\Omega_1)}_{s+\mu_{\mathtt{p},1}}\nonumber\\
&\overset{\eqref{rkkxcxcppsds11}}\le 1+ \epsilon^{5}\gamma^{-3}N_0^{2\tau+1}\rVert \mathfrak{I}_0\rVert^{\Lip(\gamma,\Omega_1)}_{s+\mu_{\mathtt{p},1}}\nonumber\\
&\le 1 +  \epsilon^{7-4b}\rVert \mathfrak{I}_0\rVert_{s+\mu_{\mathtt{p},1}}^{\Lip(\gamma,\Omega_1)},
\end{align}
where the last inequality follows from \eqref{parameterssd1sd}, which says $\rho^{-1}\ge \frac{(2b-1)(2\tau+1)}{2b-2}$, therefore,
\begin{align*}
\epsilon^{5}\gamma^{-3}N_0^{2\tau+1}& = \epsilon^{5}\gamma^{-3}\epsilon^{(1-2b)\rho(2\tau+1)}\le \epsilon^{5}\gamma^{-3}\epsilon^{2-2b} \overset{\eqref{frequency_set2}}= \epsilon^7\gamma^{-4}.
\end{align*}
Hence, \eqref{oosdsdboudncsd} for $\Phi=\Phi_7$ follows from \eqref{rkkxcxcppsds11} and \eqref{jsdk2seves}. The proof for  $\Phi=\Phi_7^{-1}$ is identical.

 Item \ref{kamredu13} follows from \eqref{consjdsdsd2sd} with the strong convergence of $\Phi_{7,n}$ to $\Phi_{7}$, which have been proved in Lemma~\ref{kskd22sdsds}. $\mathcal{L}_\infty$ is real and reversible, since so is $L_0$ (see \eqref{semifinallinearoperator} and \ref{rlsd2ssss2ss} of Proposition~\ref{linearstep172sd}) and $\Phi_7$ is real and reversibility preserving. 
 
 Lastly, since $\Phi_7$ is composition of $\frac{2\pi}{\mathtt{M}}$-translation invariance preserving, so is $\Phi_7$.
\end{proof}

\section{Invertibility of $\mathcal{L}_\omega$: Proof of Proposition 7.14}\label{loidneks1s2}
\begin{proofprop}{normal_inversion}
Let us fix $\mathtt{S}\gg s_0$ and 
\begin{align}\nonumber
\mu_{\mathtt{p},0}& :=2\mu_{\mathtt{p},1} + 2\tau +1 \ (\text{recall $\mu_{\mathtt{p,1}}$ from \eqref{parameterssd1sd}}),\\
   \Omega_{\infty}(i_0) & :=\Omega^{4\gamma}_\infty(i_0), \text{ see \eqref{sjdjsd2sdsddenoteasd}}. \label{finally_muisfie}
\end{align}
It follows from \ref{reversibility_presdsx1d} of\index{$\Omega_\infty$} Proposition~\ref{modulut2sosdtame} and \ref{kamredu13} of Proposition~\ref{kamresdsd} that 
denoting 
\[
\Phi_{1,\infty}:=\Phi_7^{-1}\circ\Phi_{1-6,1},\quad \Phi_{2,\infty}:=\Phi_{1-6,2}\circ \Phi_7.\]
 we have
\begin{align}\label{dixcxc2sdcs}
\Phi_{1,\infty}\mathcal{L}_\omega \Phi_{2,\infty} = \mathcal{L}_\infty,
\end{align}
Applying  \eqref{phiestfinal1} and \eqref{oosdsdboudncsd} (we replace $\mathtt{S}$ in the application of  \eqref{phiestfinal1} and \eqref{oosdsdboudncsd},  by $\mathtt{S}+2\mu_{\mathtt{p},1}+2\tau+1$,  where the later $\mathtt{S}$ was fixed in the beginning of this proof, hence we can choose $\epsilon$ sufficiently small depending on only $\mathtt{pe}$ and $\mathtt{S}$ to obtain the following estimates for $\Phi_{i,\infty}$) and $i = 1,2$:
\begin{align}\label{rrsdxcxc}
\rVert (\Phi_{i,\infty})^{\pm}h\rVert^{\Lip(\gamma,\Omega_\infty)}_s \le_{\mathtt{pe},s} \rVert h\rVert^{\Lip(\gamma,\Omega_\infty)}_{s+\mu_{\mathtt{p},1}} + \epsilon^{7}\gamma^{-4}\rVert \mathfrak{I}_0\rVert^{\Lip(\gamma,\Omega_\infty)}_{s+\mu_{\mathtt{p},1}}\rVert h\rVert_{s_0+\mu_{\mathtt{p},1}}^{\Lip(\gamma,\Omega_\infty)},
\end{align}
for all $s\in[s_0,\mathtt{S}+2\mu_{\mathtt{p},1}+2\tau+1]$.

 Now, it suffices to prove the estimate \eqref{rjjsxcxcxcxc2sdpp}, assuming \eqref{size_assumption_2} holds, for some $\mathtt{C}>0$ and  $\mu=\mu_{\mathtt{p},0}$, which is fixed in \eqref{finally_muisfie}. For all $\omega\in \Omega_\infty(i_0)$, it follows from \ref{kamredu13} of Proposition~\ref{kamresdsd} that 
\begin{align}\label{firstmelnikov}
|\ii \omega -d_\infty(\omega,j)| \ge 2\gamma |l|^{-\tau}|{\lambda}_\alpha(j)|, \text{ for all $j\in S_{\mathtt{M}}^\perp$.}
\end{align}
For such $\omega$, we have that
\begin{align}\label{rlllsrjwsdxx}
\mathcal{L}_\infty(\omega)^{-1}g = \sum_{l\in \mathbb{Z}^\nu,\ j\in \mathbb{Z}\backslash\left\{ 0 \right\}}\frac{1}{\ii \omega \cdot l - d_\infty(\omega,j)}\widehat{g}_j(l)e^{\ii (l\cdot \varphi + j x)}, \text{ for $g(\varphi,\cdot)\in H_{S^\perp}\cap X_{\mathtt{M}}$}
\end{align}
hence, \eqref{firstmelnikov} implies that
\begin{align}\label{rlasd1sdxxxxcxc}
\rVert \mathcal{L}_\infty(\omega)^{-1}g\rVert_s \le 2 \gamma^{-1} \rVert  g \rVert_{s+{\tau}}.
\end{align}
For $\omega_1,\omega_2\in \Omega_{\infty}(i_0)$, it follows from \eqref{shsjd2sressdsdsd}, \eqref{d0defsd2}, \eqref{nusd2sdsd} and \eqref{rejsmsdim1} that, denoting $\delta(\omega,l,j):=\ii \omega - d_{\infty}(\omega,j)$,
\begin{align*}
& \left|\frac{1}{\delta(\omega_1,l,j)} - \frac{1}{\delta(\omega_2,l,j)}\right| \\
&\le \frac{|\omega_1-\omega_2||l| + |d_\infty(\omega_1,j)-d_\infty(\omega_2,j)|}{|\delta(\omega_1,l,j)||\delta(\omega_2,l,j)|} \\
&\overset{\eqref{firstmelnikov}}\le 4\gamma^{-2}|\omega_1-\omega_2||l|^{2\tau}|\lambda_\alpha(j)|^{-2}\left( |l| + \frac{|d_\infty(\omega_1,j)-d_\infty(\omega_2,j)|} {|\omega_1-\omega_2|}\right)\\
&\le_{\mathtt{pe}}\gamma^{-2}|\omega_1-\omega_2||\lambda_\alpha(j)|^{-2}|l|^{2\tau}\left( |l| + \epsilon^2\lambda_\alpha(j) \right)\\
&\le_{\mathtt{pe}}\gamma^{-2}|\omega_1-\omega_2|\left(|l|^{2\tau +1} + 1 \right).
\end{align*}
where we used \eqref{pheobe1}, which implies $|\lambda_\alpha(j)|^{-2}\le |j|^{2\alpha}$ for the third inequality (since $g(\varphi,\cdot)\in H_{S^\perp}$, we exclude the case where $j=0$). Hence, using \eqref{rlllsrjwsdxx}, we see that 
\begin{align*}
\gamma \frac{1}{|\omega_1-\omega_2|}\rVert \left(\mathcal{L}_\infty(\omega_1)^{-1}-\mathcal{L}_\infty(\omega_2)^{-1}\right)g\rVert_s\le_{\mathtt{pe},s}\gamma^{-1}\rVert g\rVert_{s+2\tau+1}.
\end{align*}
Combining this with \eqref{rlasd1sdxxxxcxc}, we have
\begin{align}\label{linverse1123}
\rVert \mathcal{L}_\infty^{-1}g \rVert_{s}^{\Lip(\gamma,\Omega_\infty(i_0))} \le_{\mathtt{pe},s}\gamma^{-1}\rVert g \rVert_{s+2\tau +1}^{\Lip(\gamma,\Omega_\infty(i_0))},\text{ for all $s\ge 0$.}
\end{align}
Clearly, $\mathcal{L}_\omega [f] =g$ is equivalent to (from \eqref{dixcxc2sdcs}) 
\begin{align}\label{rijspwiiiwnsdax}
f = \Phi_{2,\infty}^{-1}\mathcal{L}_\infty^{-1}\Phi_{1,\infty}^{-1}g.
\end{align}
Since $\Phi_{1,\infty},\Phi_{2,\infty}$ are real, reversibility preserving and $\frac{2\pi}{\mathtt{M}}$-translation invariance preserving (see Proposition~\ref{kamresdsd} and \ref{reversibility_presdsx1d}, \ref{egpp2sd} of Proposition~\ref{modulut2sosdtame}), and $\mathcal{L}_\infty$ is real, reversible and diagonal (see \ref{kamredu13} of Proposition~\ref{kamresdsd}), we have that  $g\in Y_{\perp}\cap X_{\mathtt{M}}$ implies $f\in X_{\perp}\cap X_{\mathtt{M}}$. For the estimates, it follows from \eqref{rrsdxcxc} and \eqref{linverse1123} that  for all $s\in [s_0,\mathtt{S}]$,
\begin{align*}
\rVert f \rVert_{s}^{\Lip(\gamma,\Omega_\infty)} &\le_{\mathtt{pe},\mathtt{S}} \gamma^{-1}\left(\rVert g \rVert^{\Lip(\gamma,\Omega_\infty)}_{s+2\mu_{\mathtt{p},1} + 2\tau +1} +\epsilon^{7}\gamma^{-4}\rVert \mathfrak{I}_0\rVert^{\Lip(\gamma,\Omega_\infty)}_{s+2\mu_{\mathtt{p},1} + 2\tau +1}\rVert g \rVert_{s_0+2\mu_{\mathtt{p},1}+2\tau+1}\right)\\
&\overset{\eqref{finally_muisfie}}=\gamma^{-1}\left(\rVert g \rVert^{\Lip(\gamma,\Omega_\infty)}_{s+\mu_{\mathtt{p},0}} +\epsilon^{7}\gamma^{-4}\rVert \mathfrak{I}_0 \rVert^{\Lip(\gamma,\Omega_\infty)}_{s+\mu_{\mathtt{p},0}}\rVert g \rVert_{s_0+\mu_{\mathtt{p},0}}\right) ,
\end{align*}
which gives us \eqref{rjjsxcxcxcxc2sdpp}. 
\end{proofprop}

\chapter{Nash-Moser iteration}\label{NSmoes}
In this chapter, we aim to prove\index{Nash-Moser iteration} Theorem~\ref{main1}. We denote
\[
E_n:=\left\{\mathfrak{I}=(\Theta,y,z)\in C^\infty_\varphi\times C^\infty_\varphi\times C^\infty_{\varphi,x},\  \Theta=\Pi_n\Theta,\ y =\Pi_ny,\ \Pi_nz=z\right\},
\]
where
\begin{align*}
\Pi_n\Theta(\varphi) & =\sum_{|l|\le N_n}\Theta_le^{\ii l\cdot \varphi},\quad \Pi_ny(\varphi):=\sum_{|l|\le N_n}y_le^{\ii l\cdot \varphi},\\\Pi_nz(\varphi,x) & :=\sum_{|l|,|j|\le N_n}z_{l,j}e^{\ii l\cdot \varphi + \ii j x},
\end{align*}
where $N_n$ is as defined in \eqref{parameterssd1sd}. Given $\mathfrak{I}(\varphi):=( \Theta(\varphi),y(\varphi),z(\varphi))$, we denote $\Pi_n\mathfrak{I}:=(\Pi_n \Theta,\Pi_n y, \Pi_n z)$, so that  for $\Omega\subset \mathbb{R}^\nu$,
\begin{align}\nonumber
\rVert \Pi_n\mathfrak{I} \rVert_{s+\mu}^{\Lip(\gamma,\Omega)} & \le_{s,\mu} N_n^{\mu} \rVert \mathfrak{I}\rVert^{\Lip(\gamma,\Omega)}_{s},\\
\rVert (I-\Pi_n)\mathfrak{I}\rVert^{\Lip(\gamma,\Omega)}_{s} & \le_{s,\mu} N_n^{-\mu}\rVert\mathfrak{I}\rVert^{\Lip(\gamma,\Omega)}_{s+\mu},\text{ for all $s\ge0, \mu\ge0$.}\label{truncation_3212}
\end{align}

We define inductively\index{$\mathfrak{I}_n$} for $n\ge 0$,
\begin{align}
\mathfrak{I}_0 & :=(0,0,0),\quad \mathfrak{I}_{n+1}:=\mathfrak{I}_n-\Pi_{n}T(i_n)\Pi_{n}\mathcal{F}_\omega(i_n), \nonumber \\ i_n & :=i_{triv} + \mathfrak{I}_n, \quad 
\widehat{\mathfrak{I}}_{n+1}  :=\mathfrak{I}_{n+1}-\mathfrak{I}_{n},\quad \widehat{\mathfrak{I}}_0 := 0,
\label{iteration_step2}
\end{align}
where $T(i_n)$ is the approximate inverse constructed in Proposition~\ref{approx_inverse}. We also define\index{$\mathcal{G}_n$}
\begin{equation}\label{setoffreqences}
\begin{aligned}
&\mathcal{G}_0:=\left\{ \omega\in \Omega_\epsilon: |\omega\cdot l|\ge 2\gamma|l|^{-\tau},\quad \forall l\in \mathbb{Z}^\nu\backslash\left\{ 0 \right\} \right\}=\Omega_0,\\
&\mathcal{G}_{n+1}:=\left\{ \omega\in \mathcal{G}_n: |\ii \omega\cdot l - (d_\infty(i_n(\omega))(\omega,j)-d_\infty(i_n(\omega)(\omega,k)))|\right. \\
& \qquad \qquad \qquad \qquad \left.\ge \frac{2\gamma_n|{\lambda}_\alpha(j)-{\lambda}_\alpha(k)|}{\langle l \rangle^{\tau}},\quad \forall j,k\in S_{\mathtt{M}}^\perp\cup\left\{0\right\},\ l\in \mathbb{Z}^\nu\right\},\\
&\text{where $\gamma_n:=\gamma(1+2^{-n})$.}
\end{aligned}
\end{equation}
Let us denote
\begin{align}\label{rhoreplacesds}
A:=\frac{5-4b}{(2b-1)\rho}.
\end{align}
 Recalling $\rho$ from \eqref{parameterssd1sd}, which says $\rho^{-1}\ge \frac{40(\mu_{\mathtt{p,0}}+\mu_{\mathtt{p},2})(2b-1)}{5-4b}$, we see that
\begin{align}\label{sdsdxc2sdAAA}
A \ge 40(\mu_{\mathtt{p},0}+\mu_{\mathtt{p},2})\overset{\eqref{mu0andmu1}}\ge 40\mu_{\mathtt{p}}.
\end{align}
Since $\rho$ in \eqref{parameterssd1sd} depends on  $\mathtt{p},b$ only (see \eqref{parametersets1231} for $\mathtt{p}$), so does $A$. For such $A$, we fix
\begin{equation}\label{parambesd2sdjj22}
\begin{aligned}
 \sigma_3=\sigma_5:= \frac{9}{16}A ,\quad \sigma_4:= \frac{3}{4}A, \quad \mathtt{k}:=7A,\quad \mathtt{b}_1:=9A,\quad \mathtt{S}:=s_0 +2 \mu_{\mathtt{p}} + \mathtt{b}_1.
\end{aligned}
\end{equation}

\begin{proposition}\label{nashmoser2d}
There exists $\epsilon_0=\epsilon_0({\mathtt{p},b,\mathtt{M},S_0})>0$ such that for all $\epsilon\in (0,\epsilon_0)$, and for all $n\ge 0$, 
\begin{enumerate}[label=(\arabic*)]
\item \label{embeddingsdsd2sd2} $Z_{n}(\omega):=\mathcal{F}_\omega(i_{n}(\omega))$ is well-defined\index{$Z_n$} for $\omega\in \mathcal{G}_n$ and $Z_n(\omega)\in Y_i\cap X_{i,\mathtt{M}}$. It also satisfies
\begin{align}
\rVert Z_{n}\rVert_{s_0}^{\Lip(\gamma,\mathcal{G}_{n})}&\le_{\mathtt{p},b,\mathtt{M},S_0} \epsilon^{6-2b}N_{n-1}^{-\sigma_4},\label{jjsd2asdsa1}\\
\rVert Z_{n}\rVert_{s_0+\mathtt{b}_1-\mu_{\mathtt{p}}}^{\Lip(\gamma,\mathcal{G}_n)}&\le_{\mathtt{p},b,\mathtt{M},S_0} \epsilon^{6-2b}N_{n-1}^{\mathtt{k}},\label{jjsd2asdsa2}
\end{align}
Furthermore, we have 
\begin{align}
\rVert Z_{n}\rVert_{s_0 + \mu_{\mathtt{p}}}^{\Lip(\gamma,\mathcal{G}_{n})}&\le_{\mathtt{p},b,\mathtt{M},S_0} \epsilon^{6-2b}N_{n-1}^{-\sigma_5},\label{jjsd2asdsa3}
\end{align}

\item \label{embeddingsdsd2sd1} $\mathfrak{I}_{n+1}=\mathfrak{I}_{n+1}(\omega)$ is well-defined for $\omega\in \mathcal{G}_{n+1}$ and  $\mathfrak{I}_{n+1}(\omega) \in X_i\cap X_{i,\mathtt{M}}$.  It also satisfies
\begin{align}
\rVert \widehat{\mathfrak{I}}_{n+1}\rVert^{\Lip(\gamma,\mathcal{G}_{n+1})}_{s_0+2\mu_\mathtt{p}}&\le_{\mathtt{p},b,\mathtt{M},S_0} \epsilon^{6-4b}N_{n-1}^{-\sigma_3},\label{lowernorm_11232}\\
 \rVert \widehat{\mathfrak{I}}_{n+1}\rVert^{\Lip(\gamma,\mathcal{G}_{n+1})}_{s_0+2\mu_\mathtt{p}+\mathtt{b}_1}&\le_{\mathtt{p},b,\mathtt{M},S_0} \epsilon^{6-4b}N_{n}^{\mathtt{k}} \label{hihsdsd2sdsd}
\end{align}

\end{enumerate}
\end{proposition}

\begin{proof}
The proof will be completed by the usual induction argument. For simplicity, we denote
\begin{align}\label{nzzotzzations1zzdddd}
|A_n|_{s}:=\rVert A_n \rVert_{s}^{\Lip(\gamma,\mathcal{G}_n)},\quad A\in \left\{ \mathfrak{I}, \widehat{\mathfrak{I}}, Z\right\},\quad \mathtt{p_e}:=\left\{ \mathtt{p},b,\mathtt{M},S_0\right\}.
\end{align}

\vspace{0.5\baselineskip}\noindent\textit{Proof of \ref{embeddingsdsd2sd2},\ref{embeddingsdsd2sd1} for $n=0$.} For item \ref{embeddingsdsd2sd2}, we note that $i_0=i_{triv}\in X_i\cap X_{i,\mathtt{M}}$, which is well-defined for all $\omega\in \Omega_1=\mathcal{G}_0$, hence so is $Z_0$. And Lemma~\ref{function_spacess} and Lemma~\ref{sjdsd1sdsd} tell us that $Z_0\in Y_i\cap X_{i,\mathtt{M}}$. The estimates \eqref{jjsd2asdsa1}, \eqref{jjsd2asdsa2} and \eqref{jjsd2asdsa3} follow from \eqref{tame92} with $N_{-1}=1$ (see \eqref{parameterssd1sd}). Indeed, \eqref{tame92} tells us that
\begin{align}\label{sdjj22sdxcxcosd}
|Z_0|_s\le_{\mathtt{p_e}} \epsilon^{6-2b},\text{ for all $s\in [s_0,\mathtt{S}]$}.
\end{align}

For \ref{embeddingsdsd2sd1}, we prove \ref{embeddingsdsd2sd1} for $\mathfrak{I}_1$. Trivially, we have  $\mathfrak{I}_{0}(\omega)=0\in X_i\cap X_{i,\mathtt{M}}$ well-defined.  Furthermore, by \eqref{jjsd2asdsa3} for $Z_0$, which was just proved, we have that 
 \begin{align}\label{zaisdzjstsd}
|\mathfrak{I}_0|_{s_0+\mu_{\mathtt{p}}}=0,\quad  |Z_0|_{s_0+\mu_{\mathtt{p}}}\le_{\mathtt{p_e}} \epsilon^{6-2b}.
 \end{align}
 Therefore, \eqref{size_assumption_2} holds for $\mu=\mu_{\mathtt{p}}$ and for some $\mathtt{C}=\mathtt{C}_{\mathtt{p_e}}$. Then we can apply Proposition~\ref{approx_inverse} with $\mathtt{S}$ in \eqref{parambesd2sdjj22}, and $\mathtt{C}=\mathtt{C}_{\mathtt{p_e}}$, for sufficiently $\epsilon$ small depending only on  $\mathtt{p_e}$, since $\mathtt{S}$ in \eqref{parambesd2sdjj22} and $\mathtt{C}_{\mathtt{p_e}}$ are fixed depending on $\mathtt{p_e}$. Note that Proposition~\ref{approx_inverse} tells us that we can find an approximate inverse $T(i_0)$ for each $\omega \in \Omega_{\infty}(i_0)$. From the definition of $\mathcal{G}_{1}$ in \eqref{setoffreqences} and $\Omega_{\infty}(i_0)$ in \eqref{finally_muisfie}, it is clear that
 \begin{align}\label{seisxcgisanska}
 \mathcal{G}_{n+1}\subset \Omega_{\infty}(i_n).
 \end{align} 
 Therefore, for all $\omega \in \mathcal{G}_{1}$, $T(i_0(\omega)):Y_i\cap X_{i,\mathtt{M}}\mapsto X_i\cap X_{i,\mathtt{M}}$ is well-defined. Since $Z_0(\omega) \in Y_i\cap X_{i,\mathtt{M}}$, we see that $\mathfrak{I}_{1}$ given in \eqref{iteration_step2} is well-defined for $\omega\in \mathcal{G}_{1}$ and $\mathfrak{I}_{1}\in X_i\cap X_{i,\mathtt{M}}$. 
 
  Now, we prove \eqref{lowernorm_11232} and \eqref{hihsdsd2sdsd} for $\mathfrak{I}_1$. Towards the estimates \eqref{lowernorm_11232}, we use \eqref{inverse_estimate2} to obtain
 \begin{align}
|\widehat{\mathfrak{I}}_{1}|_{s_0+2\mu_{\mathtt{p}}} &= |\Pi_{0}T(i_{0})\Pi_{n}Z_{0}|_{s_0+2\mu_{\mathtt{p}}}\nonumber\\
&\le|T(i_{0})\Pi_{0}Z_{0}|_{s_0+2\mu_{\mathtt{p}}}\nonumber\\
&\overset{\eqref{inverse_estimate2}}{\le_{\mathtt{p_e}}}\gamma^{-1}(1+\epsilon^7\gamma^{-4}\underbrace{|\mathfrak{I}_{0}|_{s_0+2\mu_{\mathtt{p}}}}_{=0})|\Pi_{0}Z_{0}|_{s_0+2\mu_\mathtt{p}}\nonumber\\
&\overset{\eqref{sdjj22sdxcxcosd}}{\le_{\mathtt{p_e}}}\epsilon^{6-2b}\gamma^{-1},\label{rksd1ccssds0001}
 \end{align}
 which gives us \eqref{lowernorm_11232} for $\widehat{\mathfrak{I}}_{1}$. For \eqref{hihsdsd2sdsd},  we compute similarly,
 \begin{align}
 |\widehat{\mathfrak{I}}_{1}|_{s_0+2\mu_{\mathtt{p}}+\mathtt{b}_1} &= |\Pi_{0}T(i_{0})\Pi_{0}Z_{0}|_{s_0+2\mu_{\mathtt{p}}+\mathtt{b}_1}\nonumber\\
&\le|T(i_{0})\Pi_{0}Z_{0}|_{s_0+2\mu_{\mathtt{p}}+\mathtt{b}_1}\nonumber\\
&\overset{\eqref{inverse_estimate2}, \mathfrak{I}_0=0}{\le_{\mathtt{p_e}}} \gamma^{-1}|\Pi_{0}Z_{0}|_{s_0+2\mu_{\mathtt{p}}+\mathtt{b}_1} \nonumber\\
&\overset{\eqref{sdjj22sdxcxcosd}}{\le_{\mathtt{p_e}}}\epsilon^{6-2b}\gamma^{-1},\label{rksxxd1sds00012}
 \end{align}
which gives us \eqref{hihsdsd2sdsd} for $\widehat{\mathfrak{I}}_{1}$.

\vspace{0.5\baselineskip}\noindent\textit{Proof of \ref{embeddingsdsd2sd2} for $n+1$.} 
We aim to prove \ref{embeddingsdsd2sd2} for $Z_{n+1}$, assuming \ref{embeddingsdsd2sd2} and \ref{embeddingsdsd2sd1} hold for $Z_k$ and $\mathfrak{I}_{k+1}$ for $k=1,\ldots,n$ respectively. Since we $i_{n+1}(\omega)\in X_{i}\cap X_{i,\mathtt{M}}$, which is well-defined for $\omega\in \mathcal{G}_{n+1}$, $Z_{n+1}$ is well-defined for $\omega\in \mathcal{G}_{n+1}$ and Lemma~\ref{function_spacess} and Lemma~\ref{sjdsd1sdsd} tell us that $Z_{n+1}\in Y_i\cap X_{i,\mathtt{M}}$. Therefore, it suffices to prove \eqref{jjsd2asdsa1}, \eqref{jjsd2asdsa2} and \eqref{jjsd2asdsa3} for $Z_{n+1}$.
We first note that the induction hypotheses on \eqref{lowernorm_11232} and \eqref{hihsdsd2sdsd} tells us that
\begin{align}
 |\mathfrak{I}_{k+1}|_{s_0+2\mu_{\mathtt{p}}}&\overset{\eqref{iteration_step2}}\le \sum_{i=0}^{k}|\widehat{\mathfrak{I}}_{i+1}|_{s_0+2\mu_{\mathtt{p}}} \le_{\mathtt{p_e}} \epsilon^{6-4b}\sum_{i=0}^{\infty}N_{i-1}^{-\sigma_3}\le_{\mathtt{p_e}}\epsilon^{6-4b} ,\label{fortheloveof2xgod}\\ 
 |\mathfrak{I}_{k+1}|_{s_0+2\mu_{\mathtt{p}}+\mathtt{b}_1}& \le \sum_{i=0}^{k}|\widehat{\mathfrak{I}}_{i+1}|_{s_0+2\mu_{\mathtt{p}}+\mathtt{b}_1}\le_{\mathtt{p_e}} \epsilon^{6-4b}\sum_{i=0}^kN_{i}^{\mathtt{k}} \overset{\eqref{choices1xs}}{\le_{\mathtt{p_e}}}\epsilon^{6-4b}N_{k}^{\mathtt{k}}, \label{forthe2x1loveofgod2} 
\end{align}
for all $0\le k\le n$.

We compute that for all $s\in  [s_0,\mathtt{S}]$ (see $\mathtt{S}$ in \eqref{parambesd2sdjj22}),
\begin{align*}
|Z_{n+1}|_{s} &\le |Z_{n+1} - Z_n - d_i\mathcal{F}_\omega(i_{n})[\widehat{\mathfrak{I}}_{n+1}]|_{s} + |Z_n +  d_i\mathcal{F}_\omega(i_{n})[\widehat{\mathfrak{I}}_{n+1}]|_{s}\\
&\overset{\eqref{iteration_step2}}\le |Z_{n+1} - Z_n - d_i\mathcal{F}_\omega(i_{n})[\widehat{\mathfrak{I}}_{n+1}]|_{s}  +  | Z_n - d_i(\mathcal{F}_\omega)(i_n)[\Pi_n T(i_n)\Pi_n Z_n]|_{s}.
\end{align*}
We further decompose $ | Z_n - d_i(\mathcal{F}_\omega)(i_n)[\Pi_n T(i_n)\Pi_n Z_n]|_{s}$ into \begin{align}
 & | Z_n - d_i(\mathcal{F}_\omega)(i_n)[\Pi_n T(i_n)\Pi_n Z_n]|_{s} \\
 &\le  |Z_n - d_i(\mathcal{F}_\omega)(i_n)[T(i_n)\Pi_n Z_n] |_s + |d_i(\mathcal{F}_\omega)(i_n) (I-\Pi_n)T(i_n)\Pi_n Z_n|_s\nonumber\\
 &\le |(I-\Pi_n)Z_n|_s + |(I-d_i(\mathcal{F}_\omega)(i_n)T(i_n))\Pi_n Z_n|_s\nonumber\\
 & \ +  |d_i(\mathcal{F}_\omega)(i_n) (I-\Pi_n)T(i_n)\Pi_n Z_n|_s.\nonumber\\
 &  = |(I-\Pi_n)Z_n|_s  + |\Pi_n \left(I-d_i(\mathcal{F}_\omega)(i_n)T(i_n))\Pi_n Z_n\right)|_s\nonumber\\
 & \ +  |(I-\Pi_n)d_i(\mathcal{F}_\omega)(i_n)T(i_n))\Pi_n Z_n|_s +  |d_i(\mathcal{F}_\omega)(i_n) (I-\Pi_n)T(i_n)\Pi_n Z_n|_s.\nonumber
\end{align}
Therefore, we have
\begin{equation}\label{intd1sdz}
\begin{aligned}
|Z_{n+1}|_{s} & \le |Z_{n+1} - Z_n - d_i\mathcal{F}_\omega(i_{n})[\widehat{\mathfrak{I}}_{n+1}]|_{s}\\
& \ +  |(I-\Pi_n)Z_n|_s\\
& \   + |\Pi_n \left(I-d_i(\mathcal{F}_\omega)(i_n)T(i_n))\Pi_n Z_n\right)|_s\\
 & \ +  |(I-\Pi_n)d_i(\mathcal{F}_\omega)(i_n)T(i_n))\Pi_n Z_n|_s \\
 & +  |d_i(\mathcal{F}_\omega)(i_n) (I-\Pi_n)T(i_n)\Pi_n Z_n|_s.
\end{aligned}
\end{equation}
We estimate the low/high norms of each term separately. We recall from Remark~\ref{derv_loss} that
\begin{align}\label{muoisduxcs}
\mu_{\mathtt{p}}\ge \mu_1,\mu_2,\mu_3.
\end{align}

For $|Z_{n+1} - Z_n - d_i\mathcal{F}_\omega(i_{n})[\widehat{\mathfrak{I}}_{n+1}]|_{s_0}$, we compute
\begin{align}
& |Z_{n+1} - Z_n - d_i\mathcal{F}_\omega(i_{n})[\widehat{\mathfrak{I}}_{n+1}]|_{s_0} \nonumber \\
&\overset{\eqref{fffs2x},\eqref{muoisduxcs}}{\le_{\mathtt{p_e}}} \epsilon \left( 1+ (|\mathfrak{I}_{n+1}|_{s_0+\mu_{\mathtt{p}}}+|\mathfrak{I}_n|_{s_0+\mu_{\mathtt{p}}}) \right)|\widehat{\mathfrak{I}}_{n+1}|_{s_0+\mu_{\mathtt{p}}}^2\nonumber\\
&\overset{\eqref{fortheloveof2xgod}}\le 2\epsilon |\widehat{\mathfrak{I}}_{n+1}|_{s_0+\mu_{\mathtt{p}}}^2 \nonumber\\
& \overset{\eqref{lowernorm_11232}}{\le_{\mathtt{p_e}}} \epsilon^{13-8b}N_{n-1}^{-2\sigma_3}\nonumber\\
&\overset{\eqref{choicess6}}\le \epsilon^{6-2b}N_n^{-\sigma_4}.\label{lovwe1}
\end{align}
Similarly, we have
\begin{align}
& |Z_{n+1} - Z_n - d_i\mathcal{F}_\omega(i_{n})[\widehat{\mathfrak{I}}_{n+1}]|_{s_0+\mathtt{b}_1-\mu_{\mathtt{p}}} \nonumber \\
&\overset{\eqref{fffs2x},\eqref{muoisduxcs}}{\le_{\mathtt{p_e}}} \epsilon\left( |\widehat{\mathfrak{I}}_{n+1}|_{s_0+\mathtt{b}_1}|\widehat{\mathfrak{I}}_{n+1}|_{s_0+\mu_\mathtt{p}}+ (|\mathfrak{I}_{n+1}|_{s_0+\mathtt{b}_1}+|\mathfrak{I}_n|_{s_0+\mathtt{b}_1})|\widehat{\mathfrak{I}}_{n+1}|_{s_0+\mu_\mathtt{p}}^2 \right)\nonumber\\
&\overset{\eqref{lowernorm_11232},\eqref{hihsdsd2sdsd},\eqref{forthe2x1loveofgod2}}{\le_{\mathtt{p_e}}} \epsilon^{13-8b}\left( N_{n}^{\mathtt{k}}N_{n-1}^{-\sigma_3} + \epsilon^{6-4b}N_{n-1}^{-2\sigma_3}N_{n}^{\mathtt{k}}\right)\nonumber\\
& \le_{\mathtt{p_e}} \epsilon^{6-2b}N_n^{\mathtt{k}}. \label{dhodlfjgrp}
\end{align}

For $  |(I-\Pi_n)Z_n|_{s_0}$, we have
\begin{align}
  |(I-\Pi_n)Z_n|_{s_0} & \overset{\eqref{truncation_3212}}\le N_n^{-(\mathtt{b}_1 - \mu_{\mathtt{p}})}|Z_n|_{s_0 + \mathtt{b}_1-\mu_{\mathtt{p}}} \nonumber \\& \overset{\eqref{jjsd2asdsa2}}{\le_{\mathtt{p_e}}} \epsilon^{6-2b}N_{n-1}^{\mathtt{k}}N_n^{-(\mathtt{b}_1-\mu_{\mathtt{p}})} \overset{\eqref{choicess2}}\le \epsilon^{6-2b} N_n^{-\sigma_4}. \label{dlsdusdlfksrp}
\end{align}
For   $|(I-\Pi_n)Z_n|_{s_0+\mathtt{b}_1-\mu_{\mathtt{p}}}$, we simply estimate it as
\begin{align}
 |(I-\Pi_n)Z_n|_{s_0+\mathtt{b}_1-\mu_{\mathtt{p}}}\le  |Z_n|_{s_0+\mathtt{b}_1-\mu_{\mathtt{p}}} \overset{\eqref{jjsd2asdsa2}}{\le_{\mathtt{p_e}}} \epsilon^{6-2b}N_n^{\mathtt{k}}.\label{dlsdusdlfksrp2}
\end{align}

For $ |\Pi_n \left(I-d_i(\mathcal{F}_\omega)(i_n)T(i_n))\Pi_n Z_n\right)|_{s_0}$, we have that
\begin{align}
 &|\Pi_n \left(I-d_i(\mathcal{F}_\omega)(i_n)T(i_n))\Pi_n Z_n\right)|_{s_0} \nonumber \\
 &\le  |\left(I-d_i(\mathcal{F}_\omega)(i_n)T(i_n))\Pi_n Z_n\right)|_{s_0}\nonumber\\
 & \overset{\eqref{inverse_estimate1}}{\le_{\mathtt{p_e}}} \epsilon^{2b-1}\gamma^{-2}\left(1 +\epsilon^7\gamma^{-4}|\mathfrak{I}_n|_{s_0+\mu_{\mathtt{p}}} \right)|Z_n|_{s_0+\mu_\mathtt{p}}|\Pi_n Z_n|_{s_0+\mu_{\mathtt{p}}}\nonumber\\
 &\overset{\eqref{fortheloveof2xgod},\gamma=\epsilon^{2b}}{\le_{\mathtt{p_e}}}\epsilon^{2b-1}\gamma^{-2}(1+\epsilon^{13-12b})|Z_n|_{s_0+\mu_{\mathtt{p}}}^2 \nonumber \\
 &\overset{\eqref{jjsd2asdsa3},\eqref{parametersets1231}}{\le_{\mathtt{p_e}}} \epsilon^{11-6b}N_{n-1}^{-2\sigma_5}\nonumber\\
 &\overset{\eqref{choicess3}}\le \epsilon^{6-2b}N_{n}^{-\sigma_4}.\label{tofindajsdj}
\end{align}
For $ |\Pi_n \left(I-d_i(\mathcal{F}_\omega)(i_n)T(i_n))\Pi_n Z_n\right)|_{s_0+\mathtt{b}_1-\mu_{\mathtt{p}_1}}$, we have
\begin{align}\label{rlaqkqdprhcnwkdakfdl2}
 |\Pi_n &\left(I-d_i(\mathcal{F}_\omega)(i_n)T(i_n))\Pi_n Z_n\right)|_{s_0+\mathtt{b}_1-\mu_{\mathtt{p}_1}}\nonumber\\
 &\overset{\eqref{truncation_3212}}\le N_n^{\mu_{\mathtt{p}}}|\left(I-d_i(\mathcal{F}_\omega)(i_n)T(i_n))\Pi_n Z_n\right)|_{s_0+\mathtt{b}_1-2\mu_{\mathtt{p}}}\nonumber\\
 &\overset{\eqref{inverse_estimate1}}{\le_{\mathtt{p_e}}} N_n^{\mu_{\mathtt{p}}}\epsilon^{2b-1}\gamma^{-2}\left(|Z_n|_{s_0+\mu_{\mathtt{p}}}|\Pi_nZ_n|_{s_0+\mathtt{b}_1-\mu_\mathtt{p}} + (|Z_n|_{s_0+\mathtt{b}_1-\mu_{\mathtt{p}}} \right.\nonumber \\
 & \qquad \qquad \qquad \qquad  \left.+ \epsilon^7\gamma^{-4}|Z_n|_{s_0+\mu_{\mathtt{p}}}|\mathfrak{I}_n|_{s_0+\mathtt{b}_1-\mu_{\mathtt{p}}})|\Pi_nZ_n|_{s_0+\mu_{\mathtt{p}}} \right)\nonumber\\
 &\overset{\gamma=\epsilon^{2b}}{\le_{\mathtt{p_e}}} N_n^{\mu_{\mathtt{p}}}\epsilon^{-2b-1}|Z_n|_{s_0+\mu_{\mathtt{p}}}\left(|Z_n|_{s_0+\mathtt{b}_1-\mu_\mathtt{p}}  + \epsilon^7\gamma^{-4}|Z_n|_{s_0+\mu_{\mathtt{p}}}|\mathfrak{I}_n|_{s_0+\mathtt{b}_1-\mu_{\mathtt{p}}}\right)\nonumber\\
 &\overset{\eqref{jjsd2asdsa3},\eqref{forthe2x1loveofgod2},\eqref{jjsd2asdsa2}}{\le_\mathtt{p_e}} \epsilon^{5-4b}N_n^{\mu_{\mathtt{p}}}N_{n-1}^{-\sigma_5}\left(\epsilon^{6-2b}N_{n-1}^{\mathtt{k}} + \epsilon^{6-2b}\underbrace{\epsilon^{13-12b}}_{\le 1, \because \eqref{parametersets1231}}N_{n-1}^{-\sigma_5}N_{n}^{\mathtt{k}}\right)\nonumber\\
 &\overset{\eqref{choicess4}}{\le_{\mathtt{p_e}}} \epsilon^{6-2b}N_n^{\mathtt{k}}.
\end{align}

For $|(I-\Pi_n)d_i(\mathcal{F}_\omega)(i_n)T(i_n))\Pi_n Z_n|_{s_0}$,  We see from \eqref{tame200zwanzig}, \eqref{tame6} and \eqref{muoisduxcs} that
\begin{align}\label{kkksdsxcxx222xc11}
|d_i\mathcal{F}_{\omega}(i_n)[\ihat]|_{s}\le_{\mathtt{p_e}} |\ihat|_{s+\mu_{\mathtt{p}}} + |\mathfrak{I}_n|_{s+\mu_{\mathtt{p}}}|\ihat|_{s_0+\mu_\mathtt{p}}.
\end{align}
Then, we compute
\begin{align}
|(I&-\Pi_n)d_i(\mathcal{F}_\omega)(i_n)T(i_n))\Pi_n Z_n|_{s_0}\nonumber\\
&\overset{\eqref{truncation_3212}}\le N_{n}^{-(\mathtt{b}_1-\mu_\mathtt{p})}|d_i(\mathcal{F}_\omega)(i_n)T(i_n))\Pi_n Z_n|_{s_0+\mathtt{b}_1-\mu_{\mathtt{p}}}\nonumber\\ 
&\overset{\eqref{kkksdsxcxx222xc11}}{\le_{\mathtt{p_e}}} N_n^{-(\mathtt{b}_1-\mu_{\mathtt{p}})}\left( | T(i_n)\Pi_n Z_n|_{s_0+\mathtt{b}_1} + |\mathfrak{I}_n|_{s_0+\mathtt{b}_1}|T(i_n)\Pi_nZ_n|_{s_0+\mu_{\mathtt{p}}}\right)\label{notnotis}.
\end{align}
From \eqref{inverse_estimate2}, it follows that
\begin{align}
|T(i_n)\Pi_{n}Z_n|_{s_0+\mathtt{b}_1}&\le_{\mathtt{p_e}} \gamma^{-1}\left(|\Pi_nZ_n|_{s_0+\mathtt{b}_1+\mu_{\mathtt{p}}} + \epsilon^7\gamma^{-4}|\mathfrak{I}_{n}|_{s_0+\mathtt{b}_1+\mu_{\mathtt{p}}}|\Pi_nZ_n|_{s_0+\mu_{\mathtt{p}}}\right)\nonumber\\
&\le\gamma^{-1}\left(N_n^{2\mu_{\mathtt{p}}}|Z_n|_{s_0+\mathtt{b}_1-\mu_{\mathtt{p}}} + \epsilon^7\gamma^{-4}|\mathfrak{I}_{n}|_{s_0+\mathtt{b}_1+\mu_{\mathtt{p}}}|\Pi_nZ_n|_{s_0+\mu_{\mathtt{p}}}\right)\nonumber\\
&\overset{\eqref{jjsd2asdsa2},\eqref{jjsd2asdsa3},\eqref{forthe2x1loveofgod2}}{\le_{\mathtt{p_e}}} \epsilon^{6-2b}\gamma^{-1}\left(N_n^{2\mu_{\mathtt{p}}}N_{n-1}^{\mathtt{k}} + \epsilon^{13-12b}N_{n-1}^{\mathtt{k}}N_{n-1}^{-\sigma_5} \right)\nonumber\\
&\overset{\eqref{parametersets1231}}{\le_{\mathtt{p_e}}} \epsilon^{6-2b}\gamma^{-1}(N_n^{2\mu_{\mathtt{p}}}N_{n-1}^{\mathtt{k}} + N_{n-1}^{\mathtt{k}})\nonumber\\
&\overset{\eqref{choicess5}}\le 2\epsilon^{6-2b}N_{n}^{\mathtt{k}}.\label{thrhrltrnr}
\end{align}
Similarly,
 \begin{align*}
&|\mathfrak{I}_n|_{s_0+\mathtt{b}_1}|T(i_n)\Pi_nZ_n|_{s_0+\mu_{\mathtt{p}}} \\&\overset{\eqref{inverse_estimate2}}{\le_{\mathtt{p_e}}} |\mathfrak{I}_n|_{s_0+\mathtt{b}_1+\mu_{\mathtt{p}}}\gamma^{-1}\left( |\Pi_nZ_n|_{s_0+2\mu_{\mathtt{p}}} + \epsilon^{7}\gamma^{-4}|\mathfrak{I}_n|_{s_0+2\mu_{\mathtt{p}}}|\Pi_nZ_n|_{s_0+\mu_{\mathtt{p}}}\right)\\
&\overset{\eqref{truncation_3212},\eqref{forthe2x1loveofgod2}}{\le_{\mathtt{p_e}}}  |\mathfrak{I}_n|_{s_0+\mathtt{b}_1+\mu_{\mathtt{p}}}\gamma^{-1}\left(1+\epsilon^{13-12b} \right)N_n^{2\mu_{\mathtt{p}}}|Z_n|_{s_0}\\
&\overset{\eqref{jjsd2asdsa1},\eqref{forthe2x1loveofgod2}}{\le_{\mathtt{p_e}}} \epsilon^{6-4b}N_{n-1}^{\mathtt{k}}\gamma^{-1}\epsilon^{6-2b}N_{n}^{2\mu_{\mathtt{p}}}N_{n-1}^{-\sigma_4}\\
&\overset{\eqref{choicess5}}\le \epsilon^{6-2b}N_{n}^{\mathtt{k}}.
\end{align*}
Therefore, together with \eqref{thrhrltrnr}, we obtain
\begin{align}\label{jjsd11009sd2}
\left( | T(i_n)\Pi_n Z_n|_{s_0+\mathtt{b}_1} + |\mathfrak{I}_n|_{s_0+\mathtt{b}_1}|T(i_n)\Pi_nZ_n|_{s_0+\mu_{\mathtt{p}}}\right) \le_{\mathtt{p_e}}  \epsilon^{6-2b}N_{n}^{\mathtt{k}}.
\end{align}
Plugging this  into \eqref{notnotis}, we have
\begin{align}\label{sk111sdsdssd}
|(I&-\Pi_n)d_i(\mathcal{F}_\omega)(i_n)T(i_n))\Pi_n Z_n|_{s_0}\le_{\mathtt{p_e}} N_{n}^{-\mathtt{b}_1+\mu_{\mathtt{p}}}\epsilon^{6-2b}N_{n}^{\mathtt{k}}\overset{\eqref{choicess2}}\le\epsilon^{6-2b}N_n^{-\sigma_4}.
\end{align}

For $|(I-\Pi_n)d_i(\mathcal{F}_\omega)(i_n)T(i_n))\Pi_n Z_n|_{s_0+\mathtt{b}_1-\mu_\mathtt{p}}$, we have
\begin{align}
&|(I-\Pi_n)d_i(\mathcal{F}_\omega)(i_n)T(i_n))\Pi_n Z_n|_{s_0+\mathtt{b}_1-\mu_\mathtt{p}}\nonumber \\
&\le |d_i(\mathcal{F}_\omega)(i_n)T(i_n))\Pi_n Z_n|_{s_0+\mathtt{b}_1-\mu_\mathtt{p}}\nonumber\\
&\overset{\eqref{kkksdsxcxx222xc11}}{\le_{\mathtt{p_e}}} \left( | T(i_n)\Pi_n Z_n|_{s_0+\mathtt{b}_1} + |\mathfrak{I}_n|_{s_0+\mathtt{b}_1}|T(i_n)\Pi_nZ_n|_{s_0+\mu_{\mathtt{p}}}\right)\nonumber\\
&\overset{\eqref{jjsd11009sd2}}{\le_{\mathtt{p_e}}} \epsilon^{6-2b}N_{n}^{\mathtt{k}}.\label{rPrlakwsd}
\end{align}

For $|d_i(\mathcal{F}_\omega)(i_n) (I-\Pi_n)T(i_n)\Pi_n Z_n|_{s_0}$, we have
\begin{align}
 |&d_i(\mathcal{F}_\omega)(i_n) (I-\Pi_n)T(i_n)\Pi_n Z_n|_{s_0}\nonumber\\
 & \overset{\eqref{kkksdsxcxx222xc11}}{\le_{\mathtt{p_e}}} \left( |(I-\Pi_n)T(i_n)\Pi_n Z_n|_{s_0+\mu_\mathtt{p}} + |\mathfrak{I}_n|_{s_0+\mu_\mathtt{p}}|(I-\Pi_n)T(i_n)\Pi_n Z_n |_{s_0+\mu_\mathtt{p}}\right)\nonumber\\
 &\overset{\eqref{truncation_3212}}{\le_{\mathtt{p_e}}} N_{n}^{-(\mathtt{b}_1-\mu_{\mathtt{p}})}\left(|T(i_n)\Pi_n Z_n|_{s_0+\mathtt{b}_1}  + |\mathfrak{I}_n|_{s_0+\mu_{\mathtt{p}}}|T(i_n)\Pi_nZ_n|_{s_0+\mathtt{b}_1}\right)\nonumber\\
 &\overset{\eqref{inverse_estimate2}}{\le_{\mathtt{p_e}}} N_n^{-\mathtt{b}_1+\mu_{\mathtt{p}}}\Big( |T(i_n)\Pi_n Z_n|_{s_0+\mathtt{b}_1}  \nonumber \\
 & \qquad \qquad \left. + |\mathfrak{I}_n|_{s_0+\mu_{\mathtt{p}}}\gamma^{-1}(|\Pi_nZ_n|_{s_0+\mathtt{b}_1+\mu_{\mathtt{p}}} + \epsilon^{7}\gamma^{-4}|\mathfrak{I}_n|_{s_0+\mathtt{b}_1+\mu_{\mathtt{p}}}|\Pi_nZ_n|_{s_0+\mu_\mathtt{p}})\right)\label{sd1rlaqkqakfdl1}
\end{align}
We see each term in $|\mathfrak{I}_n|_{s_0+\mu_{\mathtt{p}}}\gamma^{-1}(|\Pi_nZ_n|_{s_0+\mathtt{b}_1+\mu_{\mathtt{p}}} + \epsilon^{7}\gamma^{-4}|\mathfrak{I}_n|_{s_0+\mathtt{b}_1+\mu_{\mathtt{p}}}|\Pi_nZ_n|_{s_0+\mu_\mathtt{p}})$ separately. For the first term, we have that
\begin{align}\label{rkawkxkdrnranf}
|\mathfrak{I}_n|_{s_0+\mu_{\mathtt{p}}}\gamma^{-1}|\Pi_n Z_n|_{s_0+\mathtt{b}_1+\mu_{\mathtt{p}}}&\overset{\eqref{truncation_3212},\eqref{fortheloveof2xgod}}{\le_{\mathtt{p_e}}} \gamma^{-1}\epsilon^{6-4b}N_n^{2\mu_{\mathtt{p}}}|Z_n|_{s_0+\mathtt{b}_1-\mu_{\mathtt{p}}}\nonumber\\
&\overset{\eqref{jjsd2asdsa2}}{\le_{\mathtt{p_e}}}  \epsilon^{6-2b}\epsilon^{6-4b}\gamma^{-1}N_n^{2\mu_{\mathtt{p}}}N_{n-1}^{\mathtt{k}}\nonumber\\
&\overset{\eqref{choicess5}}\le \epsilon^{6-2b}N_n^{\mathtt{k}}.
\end{align}
For the second term,  we have
\begin{align}
& \gamma^{-1}|\mathfrak{I}_n|_{s_0+\mu_{\mathtt{p}}} \epsilon^{7}\gamma^{-4}|\mathfrak{I}_n|_{s_0+\mathtt{b}_1+\mu_{\mathtt{p}}}|\Pi_nZ_n|_{s_0+\mu_\mathtt{p}}\nonumber \\&\overset{\eqref{fortheloveof2xgod},\eqref{forthe2x1loveofgod2},\gamma=\epsilon^{2b}}{\le_{\mathtt{p_e}}}\gamma^{-1}\underbrace{\epsilon^{19-16b}}_{\le 1,\ \because\eqref{parametersets1231}}N_{n-1}^{\mathtt{k}}|\Pi_nZ_n|_{s_0+\mu_\mathtt{p}}\nonumber\\
&\overset{\eqref{jjsd2asdsa3}}{\le_{\mathtt{p_e}}} \epsilon^{6-2b}\gamma^{-1}N_{n-1}^{\mathtt{k}}\nonumber\\
&\overset{\eqref{choicess5}}\le\epsilon^{6-2b}N_n^{\mathtt{k}}. \label{rkawkxkdWlro}
\end{align}
Plugging \eqref{rkawkxkdWlro}, \eqref{rkawkxkdrnranf} and \eqref{thrhrltrnr} into \eqref{sd1rlaqkqakfdl1}, we obtain
\begin{align}\label{rjsdrkawkskfk}
|&d_i(\mathcal{F}_\omega)(i_n) (I-\Pi_n)T(i_n)\Pi_n Z_n|_{s_0}\le_{\mathtt{p_e}} \epsilon^{6-2b}N_n^{\mathtt{k}+\mu_{\mathtt{p}}-\mathtt{b}_1}\overset{\eqref{choicess2}}\le \epsilon^{6-2b}N_n^{-\sigma_4}.
\end{align}

For $|d_i(\mathcal{F}_\omega)(i_n) (I-\Pi_n)T(i_n)\Pi_n Z_n|_{s_0+\mathtt{b}_1-\mu_{\mathtt{p}}}$, we have
\begin{align}
|&d_i(\mathcal{F}_\omega)(i_n) (I-\Pi_n)T(i_n)\Pi_n Z_n|_{s_0+\mathtt{b}_1-\mu_{\mathtt{p}}}\nonumber\\
& \overset{\eqref{kkksdsxcxx222xc11}}{\le_{\mathtt{p_e}}} \left( |(I-\Pi_n)T(i_n)\Pi_n Z_n|_{s_0+\mathtt{b}_1} + |\mathfrak{I}_n|_{s_0+\mathtt{b}_1}|(I-\Pi_n)T(i_n)\Pi_n Z_n |_{s_0+\mu_\mathtt{p}}\right)\nonumber\\
&\le|T(i_n)\Pi_n Z_n|_{s_0+\mathtt{b}_1} +  |\mathfrak{I}_n|_{s_0+\mathtt{b}_1}| T(i_n)\Pi_n Z_n |_{s_0+\mu_\mathtt{p}}\nonumber\\
&\overset{\eqref{jjsd11009sd2}}\le \epsilon^{6-2b}N_{n}^{\mathtt{k}}.\label{gkqlrotlqkftoRl}
\end{align}

Then, plugging \eqref{lovwe1}, \eqref{dlsdusdlfksrp}, \eqref{tofindajsdj}, \eqref{sk111sdsdssd} and \eqref{rjsdrkawkskfk} into \eqref{intd1sdz}, we obtain
\begin{align}\label{znplus1l}
|Z_{n+1}|_{s_0} \le_{\mathtt{p_e}} \epsilon^{6-2b}N_n^{-\sigma_4}.
\end{align}
Similarly, plugging \eqref{dhodlfjgrp}, \eqref{dlsdusdlfksrp2}, \eqref{rlaqkqdprhcnwkdakfdl2}, \eqref{rPrlakwsd}, and \eqref{gkqlrotlqkftoRl} into \eqref{intd1sdz}, we obtain
\begin{align}\label{znplus1l2}
|Z_{n+1}|_{s_0 + \mathtt{b}_1-\mu_{\mathtt{p}}}\le_{\mathtt{p_e}} \epsilon^{6-2b}N_n^{\mathtt{k}}.
\end{align}
Therefore, we prove \eqref{jjsd2asdsa1} and \eqref{jjsd2asdsa2} for $Z_{n+1}$. 

Now we aim to prove \eqref{jjsd2asdsa3} for $Z_{n+1}$. We simply apply Lemma~\ref{GNint} and obtain
\begin{align*}
|Z_{n+1}|_{s_0+\mu_{\mathtt{p}}}& \le |Z_{n+1}|_{s_0}^{1-\frac{\mu_{\mathtt{p}}}{\mathtt{b}_1-\mu_{\mathtt{p}}}}|Z_{n+1}|_{s_0+\mathtt{b}_1-\mu_{\mathtt{p}}}^{\frac{\mu_{\mathtt{p}}}{\mathtt{b}_1-\mu_{\mathtt{p}}}} \\& \overset{\eqref{znplus1l},\eqref{znplus1l2}}{\le_{\mathtt{p_e}}} \epsilon^{6-2b}N_n^{-\sigma_4\left(1-\frac{\mu_{\mathtt{p}}}{\mathtt{b}_1-\mu_{\mathtt{p}}}\right)+\frac{\mathtt{k}\mu_{\mathtt{p}}}{\mathtt{b}_1-\mu_{\mathtt{p}}}}\overset{\eqref{choicess71}}\le \epsilon^{6-2b}N_n^{-\sigma_5},
\end{align*}
which proves \eqref{jjsd2asdsa3} for $Z_{n+1}$.

\vspace{0.5\baselineskip}\noindent\textit{Proof of \ref{embeddingsdsd2sd1} for $n+1$.} 
We assume \ref{embeddingsdsd2sd2}, \ref{embeddingsdsd2sd1} hold true for $Z_{n+1}$, which was just proved above, and $\widehat{\mathfrak{I}}_k$ for $k\le n+1$ and will aim to prove \ref{embeddingsdsd2sd1} for $\widehat{\mathfrak{I}}_{n+2}$.

 Since $\mathfrak{I}_{n+1}(\omega)\in X_i\cap X_{i,\mathtt{M}}$ is already well-defined by the induction hypotheses for $\omega \in \mathcal{G}_{n+1}$, using \eqref{fortheloveof2xgod}, \eqref{forthe2x1loveofgod2}, and \eqref{jjsd2asdsa3} for $Z_{n+1}$, we have
 \begin{align}\label{ezveryhgsdxxxcgoodd}
|\mathfrak{I}_{n+1}|_{s_0+2\mu_{\mathtt{p}}}\le_{\mathtt{p_e}}\epsilon^{6-4b},\quad |\mathfrak{I}_{n+1}|_{s_0+2\mu_{\mathtt{p}}+\mathtt{b}_1}\le \epsilon^{6-4b}N_{n-1}^{\mathtt{k}},\quad |Z_{n+1}|_{s_0+\mu_{\mathtt{p}}}\le_{\mathtt{p_e}}\epsilon^{6-2b}.
 \end{align}
therefore, \eqref{size_assumption_2} holds for $\mu=\mu_{\mathtt{p}}$ and for some $\mathtt{C}=\mathtt{C}_{\mathtt{p_e}}$. Then we can apply Proposition~\ref{approx_inverse} with $\mathtt{S}$ in \eqref{parambesd2sdjj22}, and $\mathtt{C}=\mathtt{C}_{\mathtt{p_e}}$, for sufficiently small depending on $\mathtt{p_e}$, since $\mathtt{S}$ in \eqref{parambesd2sdjj22} and $\mathtt{C}_{\mathtt{p_e}}$ are fixed depending on $\mathtt{p_e}$. Note that Proposition~\ref{approx_inverse} tells us that we can find an approximate inverse $T(i_{n+1})$ for all $\omega \in \Omega_{\infty}(i_{n+1})$. From definition of $\mathcal{G}_{n+2}$ in \eqref{setoffreqences} and $\Omega_{\infty}(i_{n+1})$ in \eqref{finally_muisfie}, it is clear that
 \begin{align}\label{seisxcgisansk22a}
 \mathcal{G}_{n+2}\subset \Omega_{\infty}(i_{n+1}).
 \end{align} 
 Therefore, for all $\omega \in \mathcal{G}_{n+2}$, $T(i_{n+1}(\omega)):Y_i \cap X_{i,\mathtt{M}} \mapsto X_i \cap X_{i,\mathtt{M}}$ is well-defined. Since $Z_{n+1}(\omega) \in Y_i\cap X_{i,\mathtt{M}}$ is well-defined for $\omega\in \mathcal{G}_{n+1}\subset\mathcal{G}_{n+2}$, we see that $\mathfrak{I}_{n+2}$ in \eqref{iteration_step2} is well-defined for $\omega\in \mathcal{G}_{n+2}$ and $\mathfrak{I}_{n+2}\in X_i\cap X_{i,\mathtt{M}}$.
 
  Now, it suffices to show the estimates \eqref{lowernorm_11232} and \eqref{hihsdsd2sdsd} for $\widehat{\mathfrak{I}}_{n+2}$. Towards the estimates \eqref{lowernorm_11232}, we use \eqref{inverse_estimate2} to obtain
 \begin{align}
 |\widehat{\mathfrak{I}}_{n+2}|_{s_0+2\mu_{\mathtt{p}}} 
&= |\Pi_{n+1}T(i_{n+1})\Pi_{n+1}Z_{n+1}|_{s_0+2\mu_{\mathtt{p}}}\nonumber\\
&\overset{\eqref{truncation_3212}}\le N_{n+1}^{\mu_{\mathtt{p}}}|T(i_{n+1})\Pi_{n+1}Z_{n+1}|_{s_0+\mu_{\mathtt{p}}}\nonumber\\
&\overset{\eqref{inverse_estimate2}}\le_{\mathtt{p_e}} N_{n+1}^{\mu_{\mathtt{p}}}\gamma^{-1}(1+\epsilon^7\gamma^{-4}|\mathfrak{I}_{n+1}|_{s_0+2\mu_{\mathtt{p}}})|\Pi_{n+1}Z_{n+1}|_{s_0+2\mu_\mathtt{p}}\nonumber\\
&\overset{\eqref{fortheloveof2xgod},\eqref{parametersets1231},\eqref{truncation_3212}}{\le_{\mathtt{p_e}} }N_{n+1}^{3\mu_{\mathtt{p}}}\gamma^{-1}(1+\underbrace{\epsilon^{7}\gamma^{-4}\epsilon^{6-4b}}_{\overset{\gamma=\epsilon^{2b}}=\epsilon^{13-12b}\overset{\eqref{parametersets1231}}\le1})|Z_{n+1}|_{s_0}\nonumber\\
&\overset{\eqref{znplus1l}}\le_{\mathtt{p_e}}\epsilon^{6-2b}\gamma^{-1}N_{n+1}^{3\mu_{\mathtt{p}}}N_{n}^{-\sigma_4}\nonumber\\
&\overset{\eqref{choicess31}}\le  \epsilon^{6-2b}\gamma^{-1} N_{n}^{-\sigma_3},\label{rksd1sds0001}
 \end{align}
 which gives us \eqref{lowernorm_11232} for $\widehat{\mathfrak{I}}_{n+2}$. For \eqref{hihsdsd2sdsd},  we compute similarly,
 \begin{align}
 & |\widehat{\mathfrak{I}}_{n+2}|_{s_0+2\mu_{\mathtt{p}}+\mathtt{b}_1} \nonumber \\
 &= |\Pi_{n+1}T(i_{n+1})\Pi_{n+1}Z_{n+1}|_{s_0+2\mu_{\mathtt{p}}+\mathtt{b}_1}\nonumber\\
&\overset{\eqref{truncation_3212}}\le N_{n+1}^{4\mu_{\mathtt{p}}}|T(i_{n+1})\Pi_{n+1}Z_{n+1}|_{s_0-2\mu_{\mathtt{p}}+\mathtt{b}_1}\nonumber\\
&\overset{\eqref{inverse_estimate2}}{\le_{\mathtt{p_e}}} N_{n+1}^{4\mu_{\mathtt{p}}}\gamma^{-1}(|\Pi_{n+1}Z_{n+1}|_{s_0-\mu_{\mathtt{p}}+\mathtt{b}_1} + \epsilon^7\gamma^{-4}|\mathfrak{I}_{n+1}|_{s_0-\mu_{\mathtt{p}}+\mathtt{b}_1}|\Pi_{n+1}Z_{n+1}|_{s_0+\mu_{\mathtt{p}}})\nonumber\\
&\overset{\eqref{truncation_3212}}\le  N_{n+1}^{5\mu_{\mathtt{p}}}\gamma^{-1}(|Z_{n+1}|_{s_0-\mu_{\mathtt{p}}+\mathtt{b}_1} + \epsilon^7\gamma^{-4}|\mathfrak{I}_{n+1}|_{s_0-\mu_{\mathtt{p}}+\mathtt{b}_1}|Z_{n+1}|_{s_0})\nonumber\\
&\overset{\eqref{forthe2x1loveofgod2},\eqref{znplus1l}, \eqref{znplus1l2}}{\le_{\mathtt{b}}}  \gamma^{-1}\epsilon^{6-2b}N_{n+1}^{5\mu_{\mathtt{p}}}\left(N_{n}^{\mathtt{k}} +\underbrace{\epsilon^{13-12b}}_{\le 1,\because \eqref{parametersets1231}}N_{n}^{\mathtt{k}}N_{n}^{-\sigma_4}\right)\nonumber\\
&\le \gamma^{-1}\epsilon^{6-2b}N_{n+1}^{5\mu_{\mathtt{p}}}N_n^{\mathtt{k}}\nonumber\\
&\overset{\eqref{choicess09}}\le \gamma^{-1}\epsilon^{6-2b}N_{n+1}^{\mathtt{k}},\label{rksd1sds00012}
 \end{align}
which gives us \eqref{hihsdsd2sdsd} for $\widehat{\mathfrak{I}}_{n+2}$. 
\end{proof}

\begin{lemma}\label{usesequensd}
 With the choice of $\mathtt{b}_1,\sigma_3,\sigma_4,\sigma_5,\mathtt{k}$ in \eqref{parambesd2sdjj22} and the sequence $N_n$ in \eqref{parameterssd1sd}, we have the following: For all sufficiently small $\epsilon$ depending on $\mathtt{p},b$, 
\begin{align}
\gamma^{-1}N_n^{2\mu_{\mathtt{p}}}N_{n-1}^{\mathtt{k}}&\le N_n^{\mathtt{k}},\label{choicess5}\\
\epsilon^{7-6b}N_{n-1}^{-2\sigma_3}&\le N_n^{-\sigma_4},\label{choicess6}\\
N_{n}^{\mathtt{k}}N_{n}^{-(\mathtt{b}_1-\mu_{\mathtt{p}})}&\le N_n^{-\sigma_4} \label{choicess2},\\
\epsilon^{5-4b}N_{n-1}^{-2\sigma_5} &\le N_n^{-\sigma_4}, \label{choicess3}\\
\epsilon^{5-4b}N_n^{\mu_\mathtt{p}}N_{n-1}^{-\sigma_5}&\le 1, \label{choicess4}\\
N_{n+1}^{3\mu_{\mathtt{p}}}N_{n}^{-\sigma_4}&\le N_{n}^{-\sigma_3}, \label{choicess31}\\
N_{n+1}^{5\mu_{\mathtt{p}}}N_{n}^{\mathtt{k}}&\le N_{n+1}^{\mathtt{k}},\label{choicess09}\\
-\sigma_4\left(1-\frac{\mu_{\mathtt{p}}}{\mathtt{b}_1-\mu_{\mathtt{p}}}\right)+\frac{\mathtt{k}\mu_{\mathtt{p}}}{\mathtt{b}_1-\mu_{\mathtt{p}}} &\le -\sigma_5, \label{choicess71}\\
\sum_{k=0}^nN_{k}^{\mathtt{k}} &\le 2N_{n}^{\mathtt{k}} ,\label{choices1xs}
\end{align}
for all $n\ge 0$.\end{lemma}
\begin{proof}
Recalling $\gamma$, $N_{-1}$, $N_0$ and $N_n$ from \eqref{frequency_set2} and \eqref{parameterssd1sd}, we have
\begin{align}\label{N_0andjr}
\gamma^{-1} = N_0^{\frac{2b}{(2b-1)\rho}},\quad N_{-1} =1,\quad N_0=\epsilon^{(1-2b)\rho}, \quad \chi=3/2,\quad N_n=N_0^{(\chi)^{n}},\text{ for $n\ge 0$.}
\end{align}

\vspace{0.5\baselineskip}\noindent\textit{Proof of \eqref{choicess5}.}
When $n=0$, we see that
\begin{align}\nonumber
\gamma^{-1} N_n^{2\mu_{\mathtt{p}}}N_{n-1}^{\mathtt{k}}\le N_n^{\mathtt{k}}  & \overset{n=0}\iff \gamma^{-1}N_0^{2\mu_{\mathtt{p}}} \le N_0^{\mathtt{k}}\\
&\overset{\eqref{N_0andjr}} \iff N_0^{\frac{2b}{(2b-1)\rho}+2\mu_{\mathtt{p}}}\le N_0^{\mathtt{k}}  \iff \mathtt{k}\ge \frac{2b}{(2b-1)\rho}+2\mu_{\mathtt{p}}.\label{tlqkfwhwrkxek1}
\end{align}
When $n\ge 1$, we see that
\begin{align*}
\gamma^{-1} N_n^{2\mu_{\mathtt{p}}}N_{n-1}^{\mathtt{k}}\le N_n^{\mathtt{k}} & \iff \frac{2b}{(2b-1)\rho}+2\mu_{\mathtt{p}}\chi^{n} + \mathtt{k}\chi^{n-1}\le \mathtt{k}\chi^{n}\\ & \overset{\chi =3/2} \Longleftarrow \mathtt{k}\ge \frac{4b}{(2b-1)\rho} +6\mu_{\mathtt{p}}.
\end{align*}
Therefore, for any $n\ge 0$, it suffices to show that
\begin{align}\label{sjjsdxcxcxc1}
\mathtt{k}\ge \frac{4b}{(2b-1)\rho} +6\mu_{\mathtt{p}}.
\end{align}
Since $
\frac{4b}{(2b-1)\rho} \overset{\eqref{rhoreplacesds}}= A\frac{4b}{5-4b}\overset{b\in(1,1+1/12)}\le \frac{13}{2}A,$ it suffices to show that
\begin{align}\label{jrkdpwissd2s2}
\mathtt{k}\ge \frac{13}{2}A + 6\mu_\mathtt{p}.
\end{align}
From \eqref{parambesd2sdjj22}, we have $\mathtt{k}=7A$, therefore 
\[
\mathtt{k}-\left(\frac{13}2 A +6\mu_{\mathtt{p}}\right) = \frac{1}{2}A - 6\mu_{\mathtt{p}}\overset{\eqref{sdsdxc2sdAAA}}\ge 0.
\]
Therefore, \eqref{jrkdpwissd2s2} holds true.

\vspace{0.5\baselineskip}\noindent\textit{Proof of \eqref{choicess6}.}
When $n=0$, we see that
\begin{align*}
\epsilon^{7-6b}N_{-1}^{-2\sigma_3}\le N_0^{-\sigma_4} & \iff \epsilon^{7-6b}\le N_0^{-\sigma_4}\overset{N_0=(\epsilon\gamma^{-1})^\rho=\epsilon^{\rho(1-2b)}}{\iff }\epsilon^{7-6b}\le \epsilon^{(2b-1)\rho\sigma_4} \\& \iff \sigma_4\le \frac{7-6b}{(2b-1)\rho} \iff \sigma_4 \le\frac{5-4b}{(2b-1)\rho}\frac{7-6b}{5-4b}\\
&\overset{\eqref{rhoreplacesds}}\iff \sigma_4 \le A\frac{7-6b}{5-4b}\overset{b\in(1,1+1/12)}{\Longleftarrow}\sigma_4 \le \frac{3}{4}A
\end{align*}
which holds true thanks to $\sigma_4$ in \eqref{parambesd2sdjj22}.
When $n\ge 1$, we see that
\begin{align*}
\epsilon^{7-6b}N_{n-1}^{-2\sigma_3}\le N_n^{-\sigma_4}& \Longleftarrow N_{n-1}^{-2\sigma_3}\le N_n^{-\sigma_4} \\
& \iff N_0^{-2\sigma_3\chi^{n-1}}\le N_0^{-\sigma_4\chi^{n}}\iff \frac{3}{2}\sigma_4\le 2\sigma_3,
\end{align*}
which holds true thanks to $\sigma_4,\sigma_3$ in \eqref{parambesd2sdjj22}.

\vspace{0.5\baselineskip}\noindent\textit{Proof of \eqref{choicess2}.}
We see that
\[
N_{n}^{\mathtt{k}}N_{n}^{-(\mathtt{b}_1-\mu_{\mathtt{p}})}\le N_n^{-\sigma_4}  \iff \sigma_4\le \mathtt{b}_1-\mathtt{k}-\mu_\mathtt{p}.
\]
From \eqref{parambesd2sdjj22}, we see that
\[
\mathtt{b}_1-\mathtt{k}-\mu_{\mathtt{p}}-\sigma_4 = \frac{5}{4}A-\mu_{\mathtt{p}} \overset{\eqref{sdsdxc2sdAAA}}\ge 0.
\]
Therefore, \eqref{choicess2} holds true.

\vspace{0.5\baselineskip}\noindent\textit{Proof of \eqref{choicess3}.}
When $n=0$, we see that
\[
\epsilon^{5-4b}\le N_{0}^{-\sigma_4}\overset{\eqref{N_0andjr}}\iff \epsilon^{5-4b} \le  \epsilon^{(2b-1)\rho\sigma_4} \iff \sigma_4 \le \frac{5-4b}{(2b-1)\rho}\overset{\eqref{rhoreplacesds}}\iff \sigma_4 \le A,
\]
which holds true thanks to $\sigma_4$ in \eqref{parambesd2sdjj22}.
When $n\ge 1$, we see that
\begin{align*}
\epsilon^{5-4b}N_{n-1}^{-2\sigma_5}\le N_n^{-\sigma_4}& \Longleftarrow N_{n-1}^{-2\sigma_5}\le N_n^{-\sigma_4} \\& \iff N_0^{\sigma_4\chi^{n}}\le N_0^{2\sigma_5 \chi^{n-1}} \iff \frac{3}{2}\sigma_4 \le 2\sigma_5,
\end{align*}
which holds true, thanks to $\sigma_4,\sigma_5$ in \eqref{parambesd2sdjj22}.

\vspace{0.5\baselineskip}\noindent\textit{Proof of \eqref{choicess4}.}
When $n=0$, we see that
\begin{align*}
\epsilon^{5-4b}N_0^{\mu_{\mathtt{p}}}\le 1 &\overset{\eqref{N_0andjr}}\iff \epsilon^{5-4b + \mu_{\mathtt{p}}\rho(1-2b)} \le 1\iff \rho\mu_{\mathtt{p}}(2b-1)\le 5-4b \\&\iff \frac{1}{\rho}\ge \frac{(2b-1)\mu_{\mathtt{p}}}{5-4b}
\overset{\eqref{rhoreplacesds}}\iff A\ge \mu_{\mathtt{p}},
\end{align*}
which holds true thanks to \eqref{sdsdxc2sdAAA}.
When $n\ge 1$, we see that
\[
\epsilon^{5-4b}N_{n}^{\mu_{\mathtt{p}}}N_{n-1}^{-\sigma_5}\le 1 \Longleftarrow N_{n}^{\mu_{\mathtt{p}}}N_{n-1}^{-\sigma_5}\le 1 \iff N_0^{\mu_{\mathtt{p}}\chi^{n} -\sigma_5\chi^{n-1}}\le 1 \iff \sigma_5 \ge \frac{3}{2}\mu_{\mathtt{p}},
\]
which holds also true, thanks to the choice of $\sigma_5$ in \eqref{parambesd2sdjj22} and the estimate for $A$ in \eqref{sdsdxc2sdAAA}.

\vspace{0.5\baselineskip}\noindent\textit{Proof of \eqref{choicess31}.}
We see that
\begin{align*}
N_{n+1}^{3\mu_{\mathtt{p}}}N_{n}^{-\sigma_4}\le N_n^{-\sigma_3}&\iff N_0^{3\mu_{\mathtt{p}}\chi^{n+1} - \sigma_4 \chi^{n}}\le N_0^{-\sigma_3\chi^n} \\&\overset{\chi =3/2}\iff \frac{9}2\mu_{\mathtt{p}} - \sigma_4 \le - \sigma_3 \iff \sigma_4 - \sigma_3 \ge \frac{9}2\mu_{\mathtt{p}}.
\end{align*}
From $\sigma_4,\sigma_3$ in \eqref{parambesd2sdjj22}, we have
\[
\sigma_4-\sigma_3 - \frac{9}2\mu_{\mathtt{p}} = \frac{3}{16}A-\frac{9}2{\mu_{\mathtt{p}}}\overset{\eqref{sdsdxc2sdAAA}}\ge \frac{15}{2}\mu_{\mathtt{p}}-\frac{9}{2}\mu_{\mathtt{p}}\ge 0.
\]
Therefore, \eqref{choicess31} holds true.

\vspace{0.5\baselineskip}\noindent\textit{Proof of \eqref{choicess09}.}
We see that
\begin{align*}
N_{n+1}^{5\mu_{\mathtt{p}}}N_n^{\mathtt{k}}\le N_{n+1}^{\mathtt{k}}&\iff N_0^{5\mu_{\mathtt{p}}\chi^{n+1} + \mathtt{k}\chi^{n}}\le N_0^{\mathtt{k}\chi^{n+1}}\\
&\iff 5\mu_{\mathtt{p}}\chi^{n+1} + \mathtt{k}\chi^{n} \le \mathtt{k}\chi^{n+1}\overset{\chi=3/2}\Longleftarrow 15\mu_{\mathtt{p}}\le \mathtt{k}.
\end{align*}
From \eqref{parambesd2sdjj22}, we have
\[
\mathtt{k} - 15\mu_{\mathtt{p}} = 7A -15\mu_{\mathtt{p}}\overset{\eqref{sdsdxc2sdAAA}}\ge 0.
\]
Therefore, \eqref{choicess09} holds true.

\vspace{0.5\baselineskip}\noindent\textit{Proof of \eqref{choicess71}.}
Plugging $\sigma_4,\sigma_5$ and $\mathtt{k}$ from \eqref{parambesd2sdjj22} into \eqref{choicess71}, it suffices to show that
\[
-(1-\frac{\mu_{\mathtt{p}}}{\mathtt{b}_1-\mu_{\mathtt{p}}}) + \frac{7\mu_{\mathtt{p}}}{\mathtt{b}_1-\mu_{\mathtt{p}}}\le -\frac{3}{4},
\]
which is equivalent to
\[
\frac{8\mu_{\mathtt{p}}}{\mathtt{b}_1-\mu_{\mathtt{p}}} \le \frac{1}{4}.
\]
This is equivalent to $33\mu_{\mathtt{p}}\le \mathtt{b}_1$, while we have $
\mathtt{b}_1 \overset{\eqref{parambesd2sdjj22}}= 9A \overset{\eqref{sdsdxc2sdAAA}}\ge 360\mu_{\mathtt{p}}.$
Therefore,  \eqref{choicess71} holds true.

\vspace{0.5\baselineskip}\noindent\textit{Proof of \eqref{choices1xs}.}
Recalling $N_n$ from \eqref{N_0andjr}, we see that it suffices to show that
\[
\sum_{k=0}^{n-1}N_k^{\mathtt{k}}\le N_n^{\mathtt{k}},\text{ for all $n\ge 1$.}
\]
This can be proved by the usual induction argument. If $n=1$, the result follows immediately since $n\mapsto N_n$ is an increasing sequence. For $n+1$, we have that
\[
\sum_{k=0}^{n}N_k^{\mathtt{k}} = \sum_{k=0}^{n-1}N_k^{\mathtt{k}}  + N_{n}^{\mathtt{k}} \le 2N_n^{\mathtt{k}},
\]
where the last inequality follows from the induction hypothesis. Therefore, it suffices to show that $2N_{n}^{\mathtt{k}}\le N_{n+1}^{\mathtt{k}}$, in other words (recalling $N_n$ from \eqref{N_0andjr})
\[
2N_0^{\mathtt{k}\chi^{n}}\le N_0^{\mathtt{k}\chi^{n+1}}, \text{ for all $n\ge 0$.} 
\]
Since $\mathtt{k}>0$, $\chi=3/2$ and $n\ge 0$, this certainly holds for sufficiently large $N_0$, that is, for sufficiently small $\epsilon$ (see $N_0$ in \eqref{N_0andjr}).
\end{proof}

\section{Measure of the frequency set}\label{measure_estimate2ssd2}
In this section, we estimate the measure of the set $\mathcal{G}_n$ in \eqref{setoffreqences}. 
Recall from \eqref{iteration_step2} and \ref{embeddingsdsd2sd1} of Proposition~\ref{nashmoser2d} that $i_n(\omega)$ is a well-defined reversible, $\frac{2\pi}{\mathtt{M}}$-translation invariant embedding for each $\omega\in \mathcal{G}_n$ for all $n \ge 0$. It follows from \eqref{lowernorm_11232} and \eqref{jjsd2asdsa3} that (using the sequence $\sum_{n=-1}^{\infty} N_n^{-\sigma}$ is summable, and arguing as in \eqref{fortheloveof2xgod}), each $i_n$ for all $n\ge0$ satisfies 
\begin{align}\label{rosisidsd22sd}
\rVert \mathfrak{I}_n\rVert^{\Lip(\gamma,\mathcal{G}_n)}_{s_0+2\mu_{\mathtt{p}}}\le_{\mathtt{p_e}} \epsilon^{6-4b},\quad \rVert Z_n\rVert^{\Lip(\gamma,\mathcal{G}_n)}_{s_0+\mu_{\mathtt{p}}}\le_{\mathtt{p_e}} \epsilon^{6-2b},
\end{align}
 where $\mathtt{p_e}=\left\{\mathtt{p},b,\mathtt{M},S_0\right\}$ (recall from \eqref{nzzotzzations1zzdddd}).
Specifically, $i_n$ satisfies \eqref{size_assumption_2} with $\Omega_1=\mathcal{G}_n$, $\mu=\mu_{\mathtt{p}}$ and for some $\mathtt{C}=\mathtt{C}(\mathtt{p_e})$. Since $\mu_{\mathtt{p}}>\mu_{\mathtt{p},1}$ (see \eqref{finally_muisfie} and \eqref{mu0andmu1}), we apply Proposition~\ref{kamresdsd} and \ref{diagonalpart_1} of Proposition~\ref{modulut2sosdtame} with $\mathtt{S}$ (fixed in \eqref{parambesd2sdjj22}) to see that the eigenvalues
\begin{align}
d_\infty(i_n(\omega))(\omega,j) &= d_0(i_n(\omega))(\omega,j) + r_{\infty}(i_n(\omega))(\omega,j),\label{d11secxoxpxossxx}\\
 d_0(i_n(\omega))(\omega,j)&=\ii\left( \mathtt{m}_\alpha(i_n(\omega))\lambda_\alpha(j) +  j \frac{T_\alpha}4 + \mathfrak{m}_{1}(i_n(\omega))(\omega,j) \right),\nonumber
\end{align}
are well-defined and satisfy the decompositions:
\begin{equation}\label{decoprjdmlRmxdlqhdlsekxxx}
\begin{aligned}
\mathtt{m}_{\alpha}(i_n(\omega))&=-\frac{1}{2}+\epsilon^2\mathtt{m}_{\alpha,1}(\omega)+\mathtt{m}_{\alpha,2}(i_n(\omega)),\\
\mathfrak{m}_{1}(\omega,\xi) &:= \xi\left(\mathfrak{m}_{\le 0}(i_n(\omega))(\omega,\xi) + \epsilon^2\mathfrak{m}_{\mathfrak{b}}(\omega,\xi)\right)\\
& \overset{\eqref{docsd2ksdsdsd}} = \xi\left(\epsilon^2(\mathfrak{m}_{\le 0 ,1}(\omega,\xi)+\mathfrak{m}_{\mathfrak{b}}(\omega,\xi)) + \mathfrak{m}_{\le 0 ,2}(i_n(\omega))(\omega,\xi) \right),
\end{aligned}
\end{equation}
with the estimates (see \eqref{qtildesd2dsd}, \eqref{rejsmsdim1}, \eqref{m109symbolmsx} and \eqref{shsjd2sressdsdsd} and recall $\mathtt{p_e}$ from \eqref{nzzotzzations1zzdddd}),
\begin{align}
|\mathtt{m}_{\alpha,1}|^{\Lip(\gamma,\mathcal{G}_n)} &\le_{\mathtt{p_e}} 1,\quad |\mathtt{m}_{\alpha,2}|^{\Lip(\gamma,\mathcal{G}_n)}\le_{\mathtt{p_e}} \epsilon^{7-4b},\nonumber \\
& |\mathtt{m}_{\alpha,2}(i_{n}) -\mathtt{m}_{\alpha,2}(i_{m}) | \le_{\mathtt{p_e}} \epsilon\rVert i_n-i_m\rVert_{s_0+\mu_{\mathtt{p}}},\label{measure_est_11s1}\\
|\mathfrak{m}_1|^{\Lip(\gamma,\mathcal{G}_n)}_{1,0,2}&\le_{\mathtt{p_e}}  \epsilon^2,\quad |\mathfrak{m}_1(i_n) - \mathfrak{m}_1(i_m)|_{1,0,2}\le_{\mathtt{p_e}} \epsilon \rVert i_n-i_m \rVert_{s_0+\mu_{\mathtt{p}}},\label{measure_est_11s2}\\
|\mathfrak{m}_{\le 0,1}|^{\Lip(\gamma,\Omega_1)}_{0,0,2}&\le_{d} 1,\quad |\mathfrak{m}_{\le 0,2}|^{\Lip(\gamma,\Omega_1)}_{0,0,2}\le_{\mathtt{p_e}} \epsilon^{7-4b},\nonumber \\ & |\mathfrak{m}_{\le 0,2}(i_n) -\mathfrak{m}_{\le 0,2}(i_m)  |_{0,0,2}\le_{\mathtt{p_e}} \epsilon\rVert i_n - i_m \rVert_{s_0+\mu_\mathtt{p}},\label{measure_est_11s3}\\
 |r_{\infty}|^{\Lip(\gamma, \mathcal{G}_n)}_{0,0,0} &\le_{\mathtt{p_e}} \epsilon^{9-6b} ,\quad |r_{\infty}(i_n) - r_{\infty}(i_m)|_{0,0,0} \le_{\mathtt{p_e}}\epsilon^3\gamma^{-1}\rVert i_n - i_m \rVert_{s_0+\mu_{\mathtt{p}}},\label{measure_est_11s4}
\end{align}
for all $n,m\in\mathbb{N}_0$, where we used the fact that the inequalities in the estimates depend on $\mathtt{pe}$ (see \eqref{psd1sdsdsd}), while $\mathtt{C}$ depends on $\mathtt{p_e}$, which is defined in \eqref{nzzotzzations1zzdddd}, hence the inequalities depend on only $\mathtt{p_e}$. Furthermore, \eqref{idernfisdsd2} tells us that 
\begin{align}\label{seamsxcondresauesc}
\mathfrak{m}_{\mathfrak{b}}(i_n(\omega))(\omega,\xi)+ \mathtt{m}_{\alpha,1}(\omega)m_{1,\alpha}(\xi)+\mathfrak{m}_{\le 0,1}(i_n(\omega))(\omega,\xi) = \mathfrak{F}_{\alpha-1}(\omega,\xi).
\end{align}
Hence, plugging \eqref{seamsxcondresauesc} and \eqref{decoprjdmlRmxdlqhdlsekxxx} into \eqref{d11secxoxpxossxx}, we have (rearranging terms according to the power of $\epsilon$)
\begin{equation}\label{seamxxsxcondresauesc1}
\begin{aligned}
d_\infty(i_n(\omega))(\omega,j) &= \ii \left( -\frac{1}{2}\lambda_{\alpha}(j) +j\frac{T_\alpha}{4}  + \epsilon^2 j\mathfrak{F}_{\alpha-1}(\omega,j) \right.\\
 &  \ + \left(\mathtt{m}_{\alpha,2}(i_n(\omega)) \lambda_\alpha(j) + j\mathfrak{m}_{\le 0,2}(i_n(\omega))(\omega, j) \right)\Big)\\
 & \ +r_{\infty}(i_n(\omega))(\omega,j),
\end{aligned}
\end{equation}
if $j\in S^\perp$,
while,  $r_{\infty}(i_n(\omega))(\omega,0)\overset{\eqref{sddpooo4},\eqref{sjdjsd2sdsddenoteasd}}=0$ and $d_0(i_n(\omega))(\omega,0)\overset{\eqref{d0defsd2},\eqref{d0defsd22}}=0$, therefore, \eqref{d11secxoxpxossxx} yields that
\begin{align}\label{rkkxcxcaaj1s}
d_\infty(i_n(\omega))(\omega,0) = 0.
\end{align}

In view of the definition of $\mathcal{G}_n$, let us define
\begin{align}\nonumber
R_{ljk}(i_n):=\left\{\omega \in \mathcal{G}_n : |\ii \omega \cdot l \right.& - (d_\infty(i_n(\omega))(\omega,j) - d_\infty(i_n(\omega))(\omega,k) )| \\
& \left. < 2 \gamma_n \langle l \rangle^{-\tau} |\lambda_\alpha(j)-\lambda_\alpha(k)| \right\}.
\label{define_rlkj}
\end{align}
From \eqref{setoffreqences}, it is clear that
\begin{align}\label{rkaksxcxskfk}
\mathcal{G}_n\backslash\mathcal{G}_{n+1} = \cup_{l\in \mathbb{Z}^\nu,\ j,k\in S_{\mathtt{M}}^\perp\cup\left\{ 0 \right\}} R_{ljk}(i_n),\text{ for all $n\ge 0$.}
\end{align}
In view of \eqref{rkaksxcxskfk}, we will assume in the rest of the proof that 
\begin{align}\label{jkasxcmfoldandnonjexc}
j,k\subset S_{\mathtt{M}}^\perp \cup\left\{ 0 \right\},\text{ and }j\ne k.
\end{align}
\begin{lemma}\label{rkawkskfk182s1}
For small enough $\epsilon>0$ depending on $\mathtt{p_e}$, it holds that $R_{ljk}(i_{n+1})\subset R_{ljk}(i_n)$, if $|l|<N_{n}$ for $n\ge 0$. 
\end{lemma}
\begin{proof}
Let $\omega \subset R_{ljk}(i_{n+1})$, that is,
\begin{align}\label{omega11xxxsd}
|\ii \omega \cdot l - ( d_{\infty}(i_{n+1}(\omega))(\omega,j) - d_{\infty}(i_{n+1}(\omega))(\omega,k))| < 2\gamma_{n+1}\langle l \rangle^{-\tau}|\lambda_\alpha(j) - \lambda_\alpha(k)|.
\end{align}
Hence, for $|l|< N_n$, we have
\begin{align} \nonumber
& |\ii \omega \cdot l - ( d_{\infty}(i_{n}(\omega))(\omega,j) - d_{\infty}(i_{n}(\omega))(\omega,k)))|  \\ & \le |\ii \omega \cdot l - ( d_{\infty}(i_{n+1}(\omega))(\omega,j) - d_{\infty}(i_{n+1}(\omega))(\omega,k))| \nonumber \\
& \  + |( d_{\infty}(i_{n+1}(\omega))(\omega,j) - d_{\infty}(i_{n+1}(\omega))(\omega,k)) \nonumber \\
& \qquad \qquad \qquad \qquad \qquad -  ( d_{\infty}(i_{n}(\omega))(\omega,j) - d_{\infty}(i_{n}(\omega))(\omega,k))|\nonumber\\
&\overset{\eqref{omega11xxxsd}}< 2\gamma_{n+1}\langle l\rangle^{-\tau}|\lambda_\alpha(j) - \lambda_\alpha(k)| \nonumber\\
& \ \qquad + |( d_{\infty}(i_{n+1}(\omega))(\omega,j) - d_{\infty}(i_{n+1}(\omega))(\omega,k)) \nonumber \\ & \qquad \qquad \qquad \qquad \qquad-  ( d_{\infty}(i_{n}(\omega))(\omega,j) - d_{\infty}(i_{n}(\omega))(\omega,k))|.\label{gamman_1sx}
\end{align}
We claim that 
\begin{align}
& |( d_{\infty}(i_{n+1}(\omega))(\omega,j) - d_{\infty}(i_{n+1}(\omega))(\omega,k)) \nonumber \\
& -  ( d_{\infty}(i_{n}(\omega))(\omega,j) - d_{\infty}(i_{n}(\omega))(\omega,k))| \le 2^{-n}\gamma\langle l \rangle^{-\tau}|\lambda_\alpha(j) - \lambda_\alpha(k)|,\label{cliamx1meaysr21}
\end{align}
for $|l| < N_n$.
Assuming \eqref{cliamx1meaysr21} holds, plugging this into \eqref{gamman_1sx}, we see that
\begin{align}\nonumber
|\ii \omega \cdot l &- ( d_{\infty}(i_{n}(\omega))(\omega,j) - d_{\infty}(i_{n}(\omega))(\omega,k)))|  \\
& <2 \langle l \rangle^{-\tau}|\lambda_\alpha(j) - \lambda_\alpha(k)| (\gamma_{n+1} + 2^{-n-1}\gamma)\nonumber\\
&  = 2  \langle l \rangle^{-\tau}|\lambda_\alpha(j) - \lambda_\alpha(k)|  \gamma_n,\label{axcxcxc11cxpos}
\end{align}
where the last inequality follows from $\gamma_n = \gamma(1+2^{-n})$ (see \eqref{setoffreqences}). Clearly, \eqref{axcxcxc11cxpos} and \eqref{define_rlkj} imply that $\omega \in R_{ljk}(i_n)$.

To finish the proof, let us prove \eqref{cliamx1meaysr21}. From \eqref{d11secxoxpxossxx}, we have
\begin{align*}
|( & d_{\infty}(i_{n+1}(\omega))(\omega,j) 
- d_{\infty}(i_{n+1}(\omega))(\omega,k)) \\
& \qquad \qquad \qquad \qquad \qquad \qquad \qquad-  ( d_{\infty}(i_{n}(\omega))(\omega,j) - d_{\infty}(i_{n}(\omega))(\omega,k))|  \\
& \le |(\mathtt{m}_{\alpha}(i_{n+1}(\omega)) - \mathtt{m}_{\alpha}(i_{n}(\omega)))(\lambda_\alpha(j) - \lambda_\alpha(k))| \\
& \ + |(\mathfrak{m}_1(i_{n+1}(\omega))(\omega,j) -\mathfrak{m}_1(i_{n+1}(\omega))(\omega,k) ) \\
& \qquad \qquad \qquad \qquad \qquad \qquad \qquad - (\mathfrak{m}_1(i_{n}(\omega))(\omega,j) -\mathfrak{m}_1(i_{n}(\omega))(\omega,k) )|\\
& \ + 2\sup_{ j\in S^\perp\left\{ 0 \right\}}|r_{\infty}(i_n(\omega))(\omega,j) -r_{\infty}(i_{n+1}(\omega))(\omega,j) |\\
& = \left|(\mathtt{m}_{\alpha}(i_{n+1}(\omega)) - \mathtt{m}_{\alpha}(i_{n}(\omega)))(\lambda_\alpha(j) - \lambda_\alpha(k))\right|  \\
& + \left|\int_{k}^j \partial_{\xi}\mathfrak{m}_{1}(i_{n+1}(\omega))(\omega,\xi) -  \partial_{\xi}\mathfrak{m}_{1}(i_{n}(\omega))(\omega,\xi) d\xi \right|\\
& \ +2\sup_{ j\in S^\perp\left\{ 0 \right\}}|r_{\infty}(i_n(\omega))(\omega,j) -r_{\infty}(i_{n+1}(\omega))(\omega,j) |\\
& \overset{\eqref{decoprjdmlRmxdlqhdlsekxxx}}= | (\mathtt{m}_{\alpha,2}(i_{n+1}(\omega)) -\mathtt{m}_{\alpha,2}(i_{n}(\omega)))(\lambda_\alpha(j) - \lambda_\alpha(k))| \\
& \ + \left|\int_k^j \partial_\xi\left(\xi \mathfrak{m}_{\le 0 ,2}(i_{n+1}(\omega))(\omega,\xi) - \xi \mathfrak{m}_{\le 0 ,2}(i_{n}(\omega))(\omega,\xi) \right) d\xi \right| \\
& + 2\sup_{ j\in S^\perp\left\{ 0 \right\}}|r_{\infty}(i_n(\omega))(\omega,j) -r_{\infty}(i_{n+1}(\omega))(\omega,j) |\\
& \overset{\eqref{measure_est_11s1},\eqref{measure_est_11s3},\eqref{measure_est_11s4}}{\le_{\mathtt{p_e}}}  \epsilon^{3}\gamma^{-1} \rVert i_{n+1}(\omega)-i_{n}(\omega)\rVert_{s_0+\mu_{\mathtt{p}}}\left( |\lambda_\alpha(j)-\lambda_\alpha(k)| + |j-k| + 1\right) \\
& \overset{\eqref{pheobe1},\eqref{jkasxcmfoldandnonjexc}}{\le_{\mathtt{p_e}}}  \epsilon^{3}\gamma^{-1}\rVert i_{n+1}(\omega)-i_{n}(\omega)\rVert_{s_0+\mu_{\mathtt{p}}} |\lambda_\alpha(j)-\lambda_\alpha(k)| \\
& \overset{\eqref{lowernorm_11232}}\le C_{\mathtt{p_e}} \epsilon^{9-6b}N_{n-1}^{-\sigma_3} |\lambda_\alpha(j)-\lambda_\alpha(k)|,
\end{align*}
for some constant $C_\mathtt{p_e}$, where the the third equality follows from the fact that $\mathtt{m}_{\alpha,1},\mathfrak{m}_{\le 0,1}$ and $\mathfrak{m}_{\mathfrak{b}}$ are independent of  embeddings. Comparing this to the claim \eqref{cliamx1meaysr21}, it suffices to show that
\[
C_\mathtt{p_e}\epsilon^{9-6b}N_{n-1}^{-\sigma_3} \le 2^{-n}\gamma \langle l \rangle^{-\tau},\text{ for $|l|\le N_{n}$ for all $n\ge0$.}
\]
From $\gamma=\epsilon^{2b}$ (see \eqref{frequency_set2}) and $\langle l \rangle = 1+|l|$ (see \eqref{gbxx11ssx}), it suffices to show that 
\begin{align}\label{syixcxc}
C_\mathtt{p_e}\epsilon^{9-8b}N_{n-1}^{-\sigma_3}\le 2^{-n}(1+N_n)^{-\tau}, \text{ for all $n\ge0$.}
\end{align}
We consider two cases: $n=0$ and $n\ge 1$. If $n=0$, then it follows from the definition of $N_n$ in \eqref{parameterssd1sd} that we need to show
\[
C_\mathtt{p_e}\epsilon^{9-8b}\le \epsilon^{(2b-1)\rho\tau}, \text{ for some $C_\mathtt{p_e}$, possibly smaller than the one in \eqref{parameterssd1sd}}.
\]
Since \eqref{parameterssd1sd} tells us that $\rho^{-1}\ge \frac{2\tau (2b-1)}{9-8b}$, which imples $2(b-1)\rho\tau\le \frac{9-8b}{2}$. Therefore,  it is enough to show $C_\mathtt{p_e}\epsilon^{9-8b} \le  \epsilon^{\frac{9-8b}{2}}$, which is equivalent to $C_{\mathtt{p_e}}\le \epsilon^{-\frac{9-8b}{2}}$. Since $0<b-1<\frac{1}{2}$ (see \eqref{parametersets1231}), this holds true for sufficiently small $\epsilon$ depending on $\mathtt{p_e}$.

 Now we turn to $n\ge 1$. Again from \eqref{parameterssd1sd} for $N_n$, we see that \eqref{syixcxc} can be implied by showing 
 \[
 C_{\mathtt{p_e}}\epsilon^{9-8b}N_0^{-\sigma_3\chi^{n-1}} \le 2^{-n}N_0^{-\tau\chi^{n}},
 \]
 which is equivalent to (using $\chi=\frac{3}{2}$ in \eqref{parameterssd1sd}) $C_\mathtt{p_e} \epsilon^{9-8b}2^n\le N_0^{\chi^{n-1}(\sigma_3-\frac{3}2\tau)}$.  Note that
 \[
 \sigma_3 \overset{\eqref{parambesd2sdjj22}} = \frac{3}{4}A \overset{\eqref{sdsdxc2sdAAA}}\ge15 \mu_{\mathtt{p},0}\overset{\eqref{parameterssd1sd}}\ge 15(2\tau +1),
 \]
 therefore, $\sigma_3-\frac{3}{2}\tau \ge 15$, and we are led to show $C_\mathtt{p_e} \epsilon^{9-8b}2^{n}\le N_0^{15\left( 3/2\right)^{n-1}} $ for all $n\ge 1$. This inequality obviously true for sufficiently small $\epsilon$ since the right-hand side increases double-exponentially. Therefore, \eqref{syixcxc} is proved.
\end{proof}

\begin{lemma}\label{measurelemma1}
If $|l| <  N_n$, then $R_{ljk}(i_{n+1})=\emptyset$, for $n\ge 0$.
\end{lemma}
\begin{proof}
By definition in \eqref{define_rlkj}, we have that $R_{ljk}(i_{n+1})\subset \mathcal{G}_{n+1}$, while \eqref{rkaksxcxskfk} tells us that $\mathcal{G}_{n+1}\cap R_{ljk}(i_n) = \emptyset$. Therefore, for $|l| < N_n$, Lemma~\ref{rkawkskfk182s1} implies that  $
R_{ljk}(i_{n+1}) \subset \mathcal{G}_{n+1}\cap R_{ljk}(i_n) = \emptyset$, which proves the lemma.
\end{proof}

\begin{lemma}\label{rkawkskfk182s2}
 There exists $C_\mathtt{p_e}>0$ such that if $\max\left\{|j|,|k|\right\}\ge C_\mathtt{p_e}(|l|^{\frac{1}{\alpha-1}} + 1)$, then $R_{ljk}(i_n) = \emptyset$ for $n\ge 0$.
\end{lemma}
\begin{proof}
Let us assume that $R_{ljk}(i_{n})\ne \emptyset$ for some 
\begin{align}\label{rlkxcxcjqqsd2xx1}
\max\left\{ |j|,|k|\right\} \ge N (|l|^{\frac{1}{\alpha-1}} + 1),
\end{align}
for some $N$. Assuming $N$ is sufficiently large, depending on $\mathtt{p_e}$, we will derive a contradiction, which yields that $R_{ljk}(i_n)$ has to be empty.

For $\omega \in R_{ljk}(i_n)$,  \eqref{define_rlkj} tells us that
\[
|\ii \omega \cdot l - (d_\infty(i_n(\omega))(\omega,j) - d_\infty(i_n(\omega))(\omega,k) )| < 2 \gamma_n \langle l \rangle^{-\tau} |\lambda_\alpha(j)-\lambda_\alpha(k)|,
\]
that is,
\begin{align}\label{jsijustsent}
 |d_\infty(i_n(\omega))(\omega,j) - d_\infty(i_n(\omega))(\omega,k)| -2 \gamma_n \langle l \rangle^{-\tau} |\lambda_\alpha(j)-\lambda_\alpha(k)| < |\omega\cdot l |.
\end{align}
From \eqref{d11secxoxpxossxx} and \eqref{decoprjdmlRmxdlqhdlsekxxx}, we have
\begin{align}
& |d_\infty(i_n(\omega))(\omega,j) - d_\infty(i_n(\omega))(\omega,k)| \nonumber \\
& \ge  \left| -\frac{1}{2}+\epsilon^2\mathtt{m}_{\alpha,1}(\omega)+\mathtt{m}_{\alpha,2}(i_n(\omega))\right| |\lambda_\alpha(j)-\lambda_\alpha(k)| -\frac{T_\alpha}4|j-k| \nonumber\\
& \ - |\mathfrak{m}_{1}(i_n(\omega))(\omega,j)-\mathfrak{m}_{1}(i_n(\omega))(\omega,k)|  - |r_{\infty}(i_n(\omega))(\omega,j)-r_{\infty}(i_n(\omega))(\omega,k)|\nonumber\\
&\overset{\eqref{measure_est_11s1},\eqref{measure_est_11s2},\eqref{measure_est_11s3},\eqref{measure_est_11s4}}{\ge_{\mathtt{p_e}}}  \left(\frac{1}{2} - \epsilon^2\right)|\lambda_\alpha(j)-\lambda_\alpha(k)| - c_{\mathtt{p_e}}|j-k|\nonumber\\
& \ge_{\mathtt{p_e}} |\lambda_\alpha(j)-\lambda_\alpha(k)|\left( \frac{1}{2}- c_\mathtt{p_e}\frac{|j-k|}{|\lambda_\alpha(j)-\lambda_\alpha(k)|}\right)\nonumber\\
&\overset{\eqref{pheobe1}}{\ge_{\mathtt{p_e}}}|\lambda_\alpha(j)-\lambda_\alpha(k)|\left( \frac{1}{2} - \frac{c_{\mathtt{p_e}}}{|j|^{\alpha-1}+|k|^{\alpha-1}}\right).\label{lsixcxc11ca}
\end{align}
Therefore, for sufficiently large $N$ in \eqref{rlkxcxcjqqsd2xx1} depending on $\mathtt{p_e}$, we have
\[
|d_\infty(i_n(\omega))(\omega,j) - d_\infty(i_n(\omega))(\omega,k)|\ge C_\mathtt{p_e}|\lambda_\alpha(j) - \lambda_\alpha(k)|.
\]
Plugging this into \eqref{jsijustsent}, we have
\[
|\omega\cdot l | \ge |\lambda_\alpha(j)-\lambda_\alpha(k)| \left( C_\mathtt{p_e} - 2\gamma_n\langle l \rangle^{-\tau} \right),
\]
 Recalling $\gamma_n$ from \eqref{setoffreqences} and using \eqref{pheobe1}, we can choose $\epsilon$ small enough depending on $\mathtt{p_e}$ so that 
\begin{align}\label{coxx1s}
|\omega\cdot l | \ge_{\mathtt{p_e}} |j-k|\left( |j|^{\alpha-1} + |k|^{\alpha-1}\right) \overset{\eqref{rlkxcxcjqqsd2xx1},\eqref{jkasxcmfoldandnonjexc}}\ge C_{\mathtt{p_e}}N^{\alpha-1}|l| .
\end{align}
On the other hand, \eqref{rpishsnsdiwsxcsd} and \eqref{frequency_set1} tells us 
\begin{align}\label{rkxcxcxcxxx11p111ccxc}
|\omega|\le_{\mathtt{p_e}} 1,\text{ for all $\omega\in \Omega_\epsilon$.}
\end{align}  Noting that $R_{ljk}(i_n)\subset \Omega_\epsilon$, we have $|\omega\cdot l| \le_{\mathtt{p_e}} |l|$. Combining this with \eqref{coxx1s}, we get
\[
C_{\mathtt{p_e}}N^{\alpha-1}|l| \le |l|,
\]
which is a contradiction if $N$ is sufficiently large depending on $\mathtt{p_e}$.
\end{proof}

By Lemma~\ref{measurelemma1} and Lemma~\ref{rkawkskfk182s2}, we are led to estimate the measure of $R_{ljk}(i_n)$, only for $(l,j,k)$ such that $l \ge N_n$ and $|j|,|k|\le_{\mathtt{p_e}}\langle l\rangle^{\frac{1}{\alpha-1}}$, since otherwise, $R_{ljk}(i_n) = \emptyset$. For such $(l,j,k)$, we now aim to estimate the measure of $R_{ljk}(i_n)$.

\begin{lemma}\label{rkhxcxc11} Let the tangential sites $S$ in \eqref{tan_site} be chosen so that it satisfies the hypotheses \ref{hypos1s1}, \ref{hypothsdj2} and \ref{hypothsdj22} in Section~\ref{rpoisdsd1sd}. Then, 
\begin{align}\label{rkxcxccc1111kk}
|R_{ljk}(i_n)| \le_{\mathtt{p_e}} \epsilon^{2(\nu -1)}\gamma \langle l \rangle^{-\tau}, \text{ for all $n\ge 0$ and $j,k\in S_{\mathtt{M}}^\perp\cup\left\{ 0 \right\}$ such that $j\ne k$.}
\end{align}
\end{lemma}

We will postpone the proof of Lemma~\ref{rkhxcxc11} to the end of this section.

\begin{proposition}\label{set_measure_last}
Let the tangential sites $S$ in \eqref{tan_site} be chosen so that they satisfy the hypotheses \ref{hypos1s1}, \ref{hypothsdj2} and \ref{hypothsdj22} in Section~\ref{rpoisdsd1sd}. Then,  for sufficiently small $\epsilon$, depending on $\mathtt{p_e}$, it holds that 
\begin{align}\label{porisdsdws2}
|\cap_{n\ge 0}\mathcal{G}_n | \ge |\Omega_\epsilon| - C_{\mathtt{p_e}}\epsilon^{2(\nu - 1)}\gamma.
\end{align}
\end{proposition}
\begin{proof}
We have
\begin{align}\label{dhswodsd0}
|\mathcal{G}_0\backslash\mathcal{G}_1|&\overset{\eqref{rkaksxcxskfk}}\le  \sum_{l\in \mathbb{Z}^\nu,\ j,k\in S^\perp_{\mathtt{M}}\cup \left\{ 0 \right\}, j\ne k}|R_{ljk}(i_0)|\nonumber\\
& \overset{Lemma~\ref{rkawkskfk182s2}}\le \sum_{l\in \mathbb{Z}^\nu, |j|,|k| \le C_\mathtt{p_e} \left(|l|^{\frac{1}{\alpha-1}}+1 \right)}|R_{ljk}(i_0)|\nonumber\\
&\overset{Lemma~\ref{rkhxcxc11}}{\le_{\mathtt{p_e}}}\sum_{l\in \mathbb{Z}^\nu}\epsilon^{2(\nu-1)}\gamma \langle l \rangle^{-\tau}\langle l \rangle^{\frac{2}{\alpha-1}}\nonumber\\
&\overset{\eqref{frequency_set2}}{\le_\mathtt{p_e}}\epsilon^{2(\nu-1)}\gamma\sum_{l\in\mathbb{Z}^\nu}\langle l \rangle^{-\nu-2}\nonumber\\
&\le_{\mathtt{p_e}} \epsilon^{2(\nu-1)}\gamma.
\end{align}
For all $n\ge 0$,
\begin{align*}
|\mathcal{G}_{n+1}\backslash \mathcal{G}_{n+2}&|\overset{\eqref{rkaksxcxskfk}}\le \sum_{l\in \mathbb{Z}^\nu, j,k\in S_{\mathtt{M}}^\perp\cup\left\{0\right\},j\ne k}|R_{ljk}(i_{n+1})|\\
&\overset{Lemma~\ref{measurelemma1}}\le \sum_{|l| \ge N_n,\ j,k\in S_{\mathtt{M}}^\perp\cup\left\{0\right\},j\ne k}|R_{ljk}(i_{n+1})|\\
&\overset{Lemma~\ref{rkawkskfk182s2}}\le \sum_{|l|\ge N_n, \ |j|,|k| \le C_\mathtt{p_e}\left(|l|^{\frac{1}{\alpha-1}}+1\right)}|R_{ljk}(i_{n+1})|\\
& \overset{Lemma~\ref{rkhxcxc11}}{\le_{\mathtt{p_e}}} \sum_{|l|\ge N_n} \epsilon^{2(\nu - 1)}\gamma  \langle l\rangle^{-\tau+\frac{2}{\alpha-1}}\\
&\overset{\eqref{frequency_set2}}\le_{\mathtt{p_e}} \epsilon^{2(\nu -1)}\gamma \sum_{l\in \mathbb{Z}^\nu,\ |l| \ge N_n}\langle l \rangle^{-\nu - 2}\\
& \le_{\mathtt{p_e}}\epsilon^{2(\nu -1)}\gamma N_n^{-1}.
\end{align*}
Therefore, recalling $\mathcal{G}_n$ from \eqref{setoffreqences}, we have that  
\begin{align*}
|\cap_{n\ge 0}\mathcal{G}_n| & = |\mathcal{G}_0|-|\mathcal{G}_0\backslash \mathcal{G}_{1}| - \sum_{n\ge 1}|\mathcal{G}_n\backslash \mathcal{G}_{n+1}|  \\
& \ge |\Omega_0|-C_\mathtt{p_e}\epsilon^{2(\nu-1)}\gamma \left( 1 + \sum_{n\ge0}N_n^{-1}\right)\ge |\Omega_0| - C_\mathtt{p_e}\epsilon^{2(\nu-1)}\gamma,
\end{align*}
where $N_n$ is given in \eqref{parameterssd1sd} and it is a double-exponentially increasing sequence depending on $\mathtt{p_e}$. Recalling $\Omega_\epsilon, \Omega_0$ from \eqref{frequency_set1} and \eqref{frequency_set2}, and using \eqref{rpishsnsdiwsxcsd} it is classical that
\[
|\Omega_\epsilon\backslash \Omega_0|\le_{\mathtt{p_e}} \epsilon^{2(\nu-1)}\gamma,
\]
therefore,
\begin{align*}
|\cap_{n\ge 0}\mathcal{G}_n|  & \ge |\Omega_0| - C_\mathtt{p_e}\epsilon^{2(\nu-1)}\gamma \\
& \ge |\Omega_\epsilon| - |\Omega_\epsilon\backslash \Omega_0| - C_\mathtt{p_e}\epsilon^{2(\nu-1)}\gamma \ge |\Omega_\epsilon| - C_\mathtt{p_e}\epsilon^{2(\nu-1)}\gamma,
\end{align*}
which proves \eqref{porisdsdws2}.
\end{proof}

\subsection{Proof of Lemma~\ref{rkhxcxc11}.}

 We  notice from \eqref{identification_gasd2sx}  that $\zeta\mapsto \mathfrak{F}_{\alpha-1}(\omega,\xi)$ is linear. Indeed, it follows from \eqref{identification_gasd2sx} and \eqref{eigenvalusdsd1} that 
 \begin{align}\label{mathfrakFxe}
  \mathfrak{F}_{\alpha-1}&(\omega,\xi) = \vec{D}(\xi)\cdot \vec{\zeta} \overset{\eqref{xi_omega_dependent}}=\frac{\pi}{6\epsilon^2}\vec{D}(\xi)\cdot\mathbb{A}^{-1}(\omega-\overline{\omega}),\text{ for $\xi\in S^\perp$.}
 \end{align}
Then we see from \eqref{mathfrakFxe}, \eqref{eigenvalusdsd2}, \eqref{seamxxsxcondresauesc1} and \eqref{rkkxcxcaaj1s} that 
 \begin{align}\label{rlajjxcxcx2yysssd}
 d_{\infty}(i_n(\omega))(\omega,j)& = \ii \left( -\frac{1}{2}\lambda_\alpha(j) + j\frac{T_\alpha}{4} +\frac{j\pi}{6}\vec{D}(j)\cdot\mathbb{A}^{-1}(\omega-\overline{\omega}) \right)\nonumber \\
& \ + \ii  \left(\mathtt{m}_{\alpha,2}(i_n(\omega)) \lambda_\alpha(j) + j\mathfrak{m}_{\le 0,2}(i_n(\omega))(\omega, j) \right) \nonumber \\
& \ +r_{\infty}(i_n(\omega))(\omega,j).
  \end{align}
 Therefore, we can decompose  $ \omega\cdot l -\frac{1}{\ii}(d_\infty(i_n(\omega))(\omega,j) - d_\infty(i_n(\omega))(\omega,k))$ as 
 \begin{align}
 \phi(\omega)&:=\omega\cdot l -\frac{1}{\ii}(d_\infty(i_n(\omega))(\omega,j) - d_\infty(i_n(\omega))(\omega,k)) \nonumber \\
 & = a_{jk} + \vec{b}_{ljk}\cdot \omega + q_{jk}(\omega)\label{ab11xxcx1}\\
 a_{jk}&:= -\left(-\frac{1}{2}(\lambda_\alpha(j)-\lambda_\alpha(k)) + \frac{T_\alpha}4(j-k))\right) + \frac{\pi}6\left((j\vec{D}(j)-k\vec{D}(k))\cdot\mathbb{A}^{-1}\overline{\omega} \right)\label{ab11xxcx2}\\
 \vec{b}_{ljk}&:=l - \frac{\pi}6\left(j\left(\mathbb{A}^{-1}\right)^{T}\vec{D}(j)-k\left(\mathbb{A}^{-1}\right)^{T}\vec{D}(k)\right),\label{ab11xxcx3}\\
 q_{jk}(\omega)&:=-\mathtt{m}_{\alpha,2}(i_n(\omega))(\lambda_\alpha(j) - \lambda_\alpha(k)) \nonumber \\
 & - (j\mathfrak{m}_{\le 0 ,2}(i_n(\omega))(\omega,j)-k\mathfrak{m}_{\le 0 ,2}(i_n(\omega))(\omega,k)) \nonumber \\
 & \ - \frac{1}{\ii}\left( r_{\infty}(i_n(\omega))(\omega,j) - r_{\infty}(i_n(\omega))(\omega,k)\right).\label{ab11xxc4}
 \end{align}
 From \eqref{ab11xxcx2} and our hypothesis on the tangential sites in \ref{hypothsdj2} (see \eqref{poisosuduwoqs}), we have that
 \begin{enumerate}[label=$(\mathtt{H}\arabic*)$]
   \setcounter{enumi}{1}
 \item \label{hypo11} There exists a constant $C_{\mathtt{H}2}(\mathtt{p_e})>0$ such that  if $j,k\in S_\mathtt{M}^\perp\left\{ 0 \right\}$, $j\ne k$ and $\max\left\{|j|,|k|\right\}\ge C_{\mathtt{H}2}$, then
 \begin{align}
|j-k||(|j|^{\alpha-1} + |k|^{\alpha-1})&\le_{\mathtt{p_e}}  |a_{jk}|.\label{ajsj1estimax1}
 \end{align}
\end{enumerate}

Lemma~\ref{rkhxcxc11} will be proved by using several  auxiliary lemmas that will be proved below. If $\max\left\{|j|,|k|\right\}\ge C_{\mathtt{H}2}$, then \eqref{rkxcxccc1111kk} follows from Lemma~\ref{lemma12232xsetmas}. Assuming $\max\left\{|j|,|k|\right\} <  C_{\mathtt{H}2}$, then Lemma~\ref{smalljklargel} tells us that \eqref{rkxcxccc1111kk} holds for $|l|\ge L$ for some $L>0$. If $\max\left\{|j|,|k|\right\} <  C_{\mathtt{H}2}$ and $|l|\le L$, then the estimate \eqref{rkxcxccc1111kk} follows from Lemma~\ref{smalljksmalllbutsamedec2}.

\begin{lemma}\label{lemma12232xsetmas}
For $j,k\in S_\mathtt{M}^\perp\left\{ 0 \right\}$, $j\ne k$ and $\max\left\{|j|,|k|\right\}\ge C_{\mathtt{H}2}$, it holds that
\begin{align}\label{estimate11clast1}
|R_{ljk}(i_n)|\le_{\mathtt{p_e}} \epsilon^{2(\nu-1)}\gamma \langle l\rangle^{-\tau},
\end{align}
\end{lemma}
\begin{proof}
Let $l,j,k$ be fixed so that $j\ne k$ and $\max\left\{ |j|,|k|\right\} \ge C_{\mathtt{H}2}$. For $\omega \in R_{ljk}(i_n)$, we see from \eqref{ab11xxcx1} and \eqref{define_rlkj} that $|\phi(\omega)|< 2\gamma_n\langle l \rangle^{-\tau}|\lambda_\alpha(j) - \lambda_\alpha(k)|$, therefore,
\begin{align*}
|a_{j,k}|  - |\vec{b}_{ljk}\cdot \omega | - |q_{jk}(\omega)| & < 2\gamma_n\langle l \rangle^{-\tau}|\lambda_\alpha(j) - \lambda_\alpha(k)| \\
& < 2\gamma_n|\lambda_\alpha(j) - \lambda_\alpha(k)| \overset{\eqref{setoffreqences}}\le 4\gamma |\lambda_\alpha(j) - \lambda_\alpha(k)|,
\end{align*}
which implies
\begin{align}\label{bestimat1sxx}
 |\vec{b}_{ljk}\cdot \omega |  \ge |a_{j,k}| - |q_{jk}(\omega)|  -  4\gamma|\lambda_\alpha(j) - \lambda_\alpha(k)|.
\end{align}
For $|q_{jk}(\omega)|$, it follows from its definition in \eqref{ab11xxc4} and the estimates in \eqref{measure_est_11s1}, \eqref{measure_est_11s3} and \eqref{measure_est_11s4} that
\begin{align}\label{qdsxxcxc2ekdtlsgksxp}
|q_{jk}|^{\Lip(\gamma,\mathcal{G}_n)}\le_{\mathtt{p_e}} \epsilon^{9-6b}|\lambda_\alpha(j)-\lambda_\alpha(k)|,
\end{align}
therefore , plugging this and \eqref{ajsj1estimax1} into \eqref{bestimat1sxx} and using \eqref{pheobe1}, we obtain
\[
|\vec{b}_{ljk}\cdot \omega| \ge (C_\mathtt{p_e} - (4\gamma +\epsilon^{9-6b}))|\lambda_\alpha(j)-\lambda_\alpha(k)| \overset{\gamma=\epsilon^{2b}}\ge \frac{C_\mathtt{p_e}}2|\lambda_\alpha(j)-\lambda_\alpha(k)|,
\]
for sufficiently small $\epsilon>0$. Meanwhile, we see from \eqref{rkxcxcxcxxx11p111ccxc} that $
|\vec{b}_{ljk}\cdot \omega| \le_{\mathtt{p_e}} |\vec{b}_{ljk}|.$
Hence, we have that
\begin{align}\label{vecbislargess}
|\vec{b}_{ljk}| \ge C_{\mathtt{p_e}}|\lambda_\alpha(j)-\lambda_\alpha(k)|.
\end{align}
Towards the estimate of $|R_{ljk}(i_n)|$, we define
\begin{equation}\label{decompsdxcxcx}
\begin{aligned}
 b^{\perp}&:=\left\{ \vec{w}\in \mathbb{R}^\nu : \vec{w}\cdot \vec{b}_{ljk} =0, \text{ and $s \frac{\vec{b}_{ljk}}{|\vec{b}_{ljk}|} + \vec{w}\in R_{ljk}(i_n)$ for some $s\in \mathbb{R}$ }\right\},\\
 I_b(\vec{w})&:=\left\{ s\in \mathbb{R}: s \frac{\vec{b}_{ljk}}{|\vec{b}_{ljk}|} + \vec{w}\in R_{ljk}(i_n)\right\},
\end{aligned}
\end{equation}
so that each $\omega\in R_{ljk}(i_n)$ can be uniquely decomposed as 
\begin{align}\label{decopomegs}
\omega = s\frac{\vec{b}_{ljk}}{|\vec{b}_{ljk}|} + \vec{w}=: \omega(s,\vec{w}),\text{ for some $\vec{w}\in b^{\perp}$ and $s\in I_b(\vec{w})$}.
\end{align}
Denoting $\Psi(s):=\phi(s \frac{\vec{b}_{ljk}}{|\vec{b}_{ljk}|} + \vec{w})$ for each fixed $\vec{w}\in b^\perp$, we have that for each $\vec{w}\in b^\perp$ and $s_1,s_2\in I_b(\vec{w})$,
\begin{align*}
 &|\Psi(s_1)-\Psi(s_2)|\\
 & = |\phi(\omega(s_1,\vec{w}))-\phi(\omega(s_2,\vec{w}))| \\
 & \overset{\eqref{ab11xxcx1}}= |\vec{b}_{ljk}\cdot(\omega(s_1,\vec{w}) - \omega(s_2,\vec{w})) + q_{jk}(\omega(s_1,\vec{w})) - q_{jk}(\omega(s_2,\vec{w}))|\\
 &\overset{\eqref{decopomegs}}= \left|(s_1 - s_2)|\vec{b}_{ljk}| + \gamma\frac{q_{jk}(\omega(s_1,\vec{w})) - q_{jk}(\omega(s_2,\vec{w}))}{|\omega(s_1,\vec{w}) - \omega(s_2,\vec{w})|}\cdot \frac{|\omega(s_1,\vec{w}) - \omega(s_2,\vec{w})|}{\gamma}\right|\\
 &\overset{\eqref{qdsxxcxc2ekdtlsgksxp},\eqref{constnas1}}{\ge_{\mathtt{p_e}}} |s_1-s_2||\vec{b}_{ljk}|  - \epsilon^{9-6b}|\lambda_\alpha(j)-\lambda_\alpha(k)|\gamma^{-1}|s_1-s_2|\\
 &\overset{\eqref{vecbislargess},\gamma=\epsilon^{2b}}{\ge_{\mathtt{p_e}}}|s_1-s_2||\lambda_\alpha(j)-\lambda_\alpha(k)|- \epsilon^{9-8b}|\lambda_\alpha(j)-\lambda_\alpha(k)||s_1-s_2|\\
 & \ge_{\mathtt{p_e}}|s_1-s_2||\lambda_\alpha(j)-\lambda_\alpha(k)|,
\end{align*}
for all sufficiently small $\epsilon$, depending on $\mathtt{p_e}$. On the other hand,  we also have (using $\omega(s_1,\vec{w}),\omega(s_2,\vec{w})\in R_{ljk}(i_n)$, by \eqref{decompsdxcxcx})
\begin{align*}
|\Psi(s_1)-\Psi(s_2)|& = |\phi(\omega(s_1,\vec{w}))-\phi(\omega(s_2,\vec{w}))| \\
& \overset{\eqref{ab11xxcx1},\eqref{define_rlkj}}< 4\gamma_n \langle l \rangle^{-\tau}|\lambda_\alpha(j)-\lambda_\alpha(k)|\\
&\overset{\eqref{setoffreqences}}<8\gamma \langle l \rangle^{-\tau}|\lambda_\alpha(j)-\lambda_\alpha(k)|.
\end{align*}
Therefore, we have that
\[
|s_1-s_2|\le_{\mathtt{p_e}} \gamma \langle l \rangle^{-\tau}, \text{ for all $\vec{w}\in b^\perp$, and  $s_1,s_2\in I_b(\vec{w})$,}
\]
in other words, 
\begin{align}\label{setofIestimate}
|I_{b}(\vec{w})| \le_{\mathtt{p_e}} \gamma\langle l \rangle^{-\tau},\text{ for all $\vec{w}\in b^\perp$.}
\end{align}
For the set $b^\perp$ in \eqref{decompsdxcxcx}, it follows from $R_{ljk}(i_n)\subset \mathcal{G}_n\subset \Omega_\epsilon$ (see \eqref{rkaksxcxskfk}, \eqref{setoffreqences})  and the definition of $b^\perp$ in \eqref{decompsdxcxcx} that 
\[
b^\perp \subset \left\{ \vec{w}\in\mathbb{R}^\nu : \vec{w}\cdot \vec{b}_{ljk} =0,\text{ and $s \frac{\vec{b}_{ljk}}{|\vec{b}_{ljk}|} + \vec{w}\in \Omega_\epsilon$ for some $s\in \mathbb{R}$ }\right\}.
\]
Since $\Omega_\epsilon$ in \eqref{frequency_set1} is a $\nu$-dimensional cube with side-length $O(\epsilon^2)$ (see \eqref{rpishsnsdiwsxcsd}), we have that $|b^{\perp}| \le_{\mathtt{p_e}} \epsilon^{2(\nu -1)}$. Lastly, it follows from \eqref{decompsdxcxcx} that
\[
|R_{ljk}(i_n)| =\int_{b^\perp} \int_{I_b(\vec{w})} 1 ds d\vec{w} \overset{\eqref{setofIestimate}}{\le_{\mathtt{p_e}}}\gamma\langle l \rangle^{-\tau} |b^\perp| \le_{\mathtt{p_e}}\epsilon^{2(\nu - 1)}\gamma\langle l \rangle^{-\tau},
\]
which proves  \eqref{estimate11clast1}.
\end{proof}

\begin{lemma}\label{smalljklargel}
There exists $L=L(\mathtt{p_e})>0$ such that if $|l| \ge L$ and  $|j|,|k|< C_{\mathtt{H}2}$, then
\begin{align}\label{rkxcxcskshisd1}
|R_{ljk}(i_n)|\le_{\mathtt{p_e}} \epsilon^{2(\nu - 1)}\gamma \langle l \rangle^{-\tau}.
\end{align} 
\end{lemma}
\begin{proof}
First we see from \eqref{mathfrakFxe} and \eqref{seamsxcondresauesc} that $\xi\mapsto (\vec{D}(\xi))_k$ for each $k=1,\ldots,\nu$ is a Fourier multiplier in $\mathcal{S}^{\alpha-1}$, therefore,
\begin{align}\label{jjjsss11sdssdsossdow1}
|j\vec{D}(j) -k\vec{D}(k)|\le_{\mathtt{p_e}} |j-k|(|j|^{\alpha-1} + |k|^{\alpha-1}),\text{ for all $j,k\in S_{\mathtt{M}}^\perp\cup\left\{ 0 \right\}$}.
\end{align}
Then, it follows from \eqref{rpishsnsdiwsxcsd} that
\begin{align}
|(|j\vec{D}(j) - k\vec{D}(k))\cdot \mathbb{A}^{-1}\vec{w}|&\le_{\mathtt{p_e}}|\vec{w}||j-k|(|j|^{\alpha-1} +|k|^{\alpha-1}),\label{ajsj1estimax2}
\end{align}
for all $\vec{w}\in \mathbb{R}^\nu$, $j,k\in S^\perp\left\{ 0 \right\}$.

Given $\omega_1,\omega_2\in R_{ljk}(i_n)$, we have
\begin{align}\label{rkshxxxcxc1}
&\left|\left(d_\infty(i_n(\omega_1))(\omega_1,j)-d_\infty(i_n(\omega_1))(\omega_1,j) \right)\right.\nonumber \\& \qquad \qquad \qquad \qquad \qquad \qquad \left.-\left(d_\infty(i_n(\omega_1))(\omega_1,j)-d_\infty(i_n(\omega_1))(\omega_1,j) \right)\right|\nonumber\\
&  \overset{\eqref{rlajjxcxcx2yysssd}}\le\left| \frac{\pi}{6}\left(j\vec{D}(j)-k\vec{D}(k) \right)\cdot\mathbb{A}^{-1}(\omega_1-\omega_2) \right| \nonumber \\
& \ +\left| (\mathtt{m}_{\alpha,2}(i_n(\omega_1))-\mathtt{m}_{\alpha,2}(i_n(\omega_2)))(\lambda_\alpha(j) - \lambda_\alpha(k))\right|\nonumber\\
& \ + \left| \int_k^j \partial_\xi \left(\xi\mathfrak{m}_{\le 0,2}(i_n(\omega_1))(\omega_1,\xi)-\xi\mathfrak{m}_{\le 0,2}(i_n(\omega_2))(\omega_2,\xi) \right)d\xi \right| \nonumber \\
& \ + 2|r_\infty(i_n(\omega_1))-r_\infty(i_n(\omega_2))|_{0,0,0}\nonumber\\
& \overset{\eqref{measure_est_11s1},\eqref{measure_est_11s3},\eqref{measure_est_11s4}}{\le_\mathtt{p_e}} \left| \frac{\pi}{6}\left(j\vec{D}(j)-k\vec{D}(k) \right)\cdot\mathbb{A}^{-1}(\omega_1-\omega_2) \right|  + \epsilon^{9-6b}\gamma^{-1}|\omega_1-\omega_2|\nonumber\\
& \overset{\eqref{ajsj1estimax2}}{\le_{\mathtt{p_e}}}|j-k|(|j|^{\alpha-1} + |k|^{\alpha-1})|\omega_1-\omega_2| + \epsilon^{9-8b}|\omega_1-\omega_2| \nonumber\\
& \le_{\mathtt{p_e}}|j-k|(|j|^{\alpha-1} + |k|^{\alpha-1})|\omega_1-\omega_2|,
\end{align}
for sufficiently small $\epsilon$, depending on $\mathtt{p_e}$.

For a fixed $0\ne l \in \mathbb{R}^\nu$, we define
\begin{align*}
l^\perp &:= \left\{ \vec{w}\in \mathbb{R}^\nu : \vec{w}\cdot l =0, s\widehat{l} + \vec{w}\in R_{ljk}(i_n),\text{ for some $s\in\mathbb{R}$, where $\widehat{l}:=\frac{l}{|l|}$}\right\},\\
I_l(\vec{w})&:=\left\{ s\in\mathbb{R}: s\widehat{l} + \vec{w}\in R_{ljk}(i_n)\right\},
\end{align*}
so that each $\omega\in R_{ljk}(i_n)$ can be uniquely written as
\begin{align*}
\omega = \omega(s,\vec{w}):=  s \widehat{l} +  \vec{w},\text{ for some $s\in I_{l}(\vec{w})$ and $\vec{w}\in l^\perp$.}
\end{align*}
For a fixed $\vec{w}\in l^\perp$, let us denote
\[
\Psi(s):=\phi(s\widehat{l} + \vec{w}), \text{ for $s \in I_l(\vec{w})$}.
\]
Denoting $\omega_1=\omega(s_1,\vec{w}),\omega_2=\omega(s_2,\vec{w})$, we have that
\begin{align*}|\Psi(s_1)-\Psi(s_2)| & = |\phi(\omega_1) - \phi(\omega_2)|\nonumber\\
& \overset{\eqref{ab11xxcx1}}\ge |(\omega_1-\omega_2)\cdot l|\nonumber \\
& \ -\left|\left(d_\infty(i_n(\omega_1))(\omega_1,j)-d_\infty(i_n(\omega_1))(\omega_1,j) \right)\right. \nonumber\\
& \left.\qquad \qquad \qquad -\left(d_\infty(i_n(\omega_1))(\omega_1,j)-d_\infty(i_n(\omega_1))(\omega_1,j) \right)\right|\nonumber\\
& \overset{\eqref{rkshxxxcxc1}}\ge |l| |s_1-s_2| -c_\mathtt{p_e}|j-k|(|j|^{\alpha-1} + |k|^{\alpha-1})|s_1-s_2|\nonumber\\
& = ||l| - c_{\mathtt{p_e}}|j-k|(|j|^{\alpha-1} + |k|^{\alpha-1})| |s_1-s_2|\nonumber\\
& \ge (L-\mathtt{c}_d)|s_1-s_2|,
\end{align*}
where the last inequality follows from $|l| \ge L$ and $|j|,|k| <  C_{\mathtt{H}2}=C_{\mathtt{H}2}(\mathtt{p_e})$. Therefore we can choose $L$ large enough depending on $\mathtt{p_e}$ so that 
\begin{align}\label{shutfupakxcx}
|\Psi(s_1)-\Psi(s_2)|  \ge \frac{L}{2}|s_1-s_2|,\text{ for all $s_1,s_2\in I_{l}(\vec{w})$, for all $\vec{w}\in l^\perp$.}
\end{align}
On the other hand, it follows from $\omega_1,\omega_2\in R_{ljk}(i_n)$ and \eqref{define_rlkj} that
\[
|\Psi(s_1)-\Psi(s_2)| = |\phi(\omega_1) - \phi(\omega_2)| \overset{\gamma_n=(1+2^{-n}\gamma)}< 2\gamma  \langle l \rangle^{-\tau}|\lambda_\alpha(j) - \lambda_\alpha(k)|\le_{\mathtt{p_e}} \gamma \langle l \rangle^{-\tau}.
\]
Combining this with \eqref{shutfupakxcx}, we have 
\[
|I_l(\vec{w})|\le_{\mathtt{p_e}} \gamma\langle l \rangle^{-\tau},\text{ for all $\vec{w}\in l^\perp$.}
\]
As in the proof of Lemma~\ref{lemma12232xsetmas}, we have $R_{ljk}(i_n)\subset \Omega_\epsilon$, hence,
$
|l^\perp|\le_{\mathtt{p_e}} \epsilon^{2(\nu-1)},$ therefore
\[
|R_{ljk}(i_n)| = \int_{l^\perp}\int_{I_l(\vec{w})} 1 ds d\vec{w} \le_{\mathtt{p_e}}\gamma \langle l \rangle^{-\tau}|l^\perp|\le_{\mathtt{p_e}}\epsilon^{2(\nu -1)}\gamma \langle l \rangle^{-\tau},
\]
which proves \eqref{rkxcxcskshisd1}.
\end{proof}

\begin{lemma}\label{smalljksmalllbutsamedec}
 Let  $L$ be fixed as in Lemma~\ref{smalljklargel} and we recall the definition of $W$ from \eqref{defoflsksdjwwww}:
 \[
 W(j)=-\frac{1}{2}\lambda_\alpha(j) + \frac{T_\alpha}{4}j.
 \] 
 If $|l|\le L$,  $|j|,|k|< C_{\mathtt{H}2}$, and
\[
\delta:=|\overline{\omega}\cdot l - (W(j) - W(k))|\ne 0,
\]
then $R_{ljk}(i_n) = \emptyset$ for all sufficiently small $\epsilon$ depending on $\mathtt{p_e}$ and $\delta$. 
\end{lemma}
\begin{proof}
For each $\omega\in \Omega_\epsilon$, it follows from \eqref{frequency_set1} and \eqref{rpishsnsdiwsxcsd} that 
\begin{align}\label{omegadivffs}
|\omega - \overline{\omega}|\le_{\mathtt{p_e}} \epsilon^2.
\end{align}
Furthermore, using \eqref{d11secxoxpxossxx}, \eqref{decoprjdmlRmxdlqhdlsekxxx}, we see that
\begin{align*}
& |d_\infty(i_n(\omega))(\omega,j)-\ii W(j)|\\ &= |(\epsilon^2 \mathtt{m}_{\alpha,1}(\omega)+\mathtt{m}_{\alpha,2}(i_n(\omega)))\lambda_\alpha(j) + \mathfrak{m}_{1}(i_n(\omega))(\omega,j)+r_{\infty}(i_n(\omega))(\omega,j)|\\
&\overset{\eqref{measure_est_11s1},\eqref{measure_est_11s2},\eqref{measure_est_11s4}}{\le_{\mathtt{p_e}}}\epsilon^2|\lambda_\alpha(j)| \le_{\mathtt{p_e}} \epsilon^2,
\end{align*}
where the last inequality follows from $|j|< C_{\mathtt{H}2}$. Therefore, if $\omega\in R_{ljk}(i_n)$, then
\begin{align*}
2\gamma_n\langle l \rangle^{-\tau}|\lambda_\alpha(j) - \lambda_\alpha(k)| &\overset{\eqref{define_rlkj}}> |\ii\omega\cdot l - (d_{\infty}(i_n(\omega))(\omega,j)-d_{\infty}(i_n(\omega))(\omega,k))|\\
&  \ge |\overline{\omega}\cdot l - (W(j) - W(k))|  - |(\omega - \overline{\omega})\cdot l |\\
& \ -|d_\infty(i_n(\omega))(\omega,j)- \ii W(j)|\\
& \ -|d_\infty(i_n(\omega))(\omega,k)- \ii W(k)|\\
& \ge \delta - c_{\mathtt{p_e}}\epsilon^2 (|l| + 1)\\
&\ge \delta - c_{\mathtt{p_e}}\epsilon^2,
\end{align*}
where the last inequality follows from that $|l|\le L$ and $L$ in Lemma~\ref{smalljklargel} depends on only $\mathtt{p_e}$. Hence, for small enough $\epsilon$ depending on $\delta$ and $\mathtt{p_e}$, we have
\[
\frac{\delta}2  < 2\gamma_n\langle l \rangle^{-\tau}|\lambda_\alpha(j) - \lambda_\alpha(k)| \overset{\gamma_n= (1+2^{-n})\gamma,\ \eqref{pheobe1}}{\le_{\mathtt{p_e}}}\gamma (|j|^\alpha + |k|^{\alpha})\le_{\mathtt{p_e}}\gamma\le \epsilon^{2b},
\]
which is a contradiction for sufficiently small $\epsilon$, therefore $R_{ljk}(i_n)$ is empty.
\end{proof}

\begin{lemma}\label{smalljksmalllbutsamedec2}Let  $L$ be fixed as in Lemma~\ref{smalljklargel}. If $|l|\le L$,  $j,k< C_{\mathtt{H}2}$ then,
 $|R_{ljk}(i_n)|\le_{\mathtt{p_e}}\epsilon^{2(\nu -1)}\gamma \langle l \rangle^{-\tau}$ for all sufficiently small $\epsilon$ depending on $\mathtt{p_e}$.
\end{lemma}
\begin{proof}
Let us suppose $l,j,k$ are chosen according to the assumptions of the lemma. If $|\overline{\omega}\cdot l - (W(j) - W(k))|\ne 0$, then $R_{ljk}(i_n)=\emptyset$ which follows from Lemma~\ref{smalljksmalllbutsamedec}.  Therefore, let us assume that
\begin{align}\label{degenerxxx11}
|\overline{\omega}\cdot l - (W(j) - W(k))|= 0.
\end{align}
We first write $d_\infty(i_n)$ in \eqref{rlajjxcxcx2yysssd} using $W(j)$ in \eqref{defoflsksdjwwww} as
\begin{align*}
d_\infty(i_n(\omega)) &(\omega,j) \\
& = \ii \left( W(j) +\frac{j\pi}{6}\vec{D}(j)\cdot\mathbb{A}^{-1}(\omega-\overline{\omega}) \right)\\
 &  \ + \ii \left(\mathtt{m}_{\alpha,2}(i_n(\omega)) \lambda_\alpha(j) + j\mathfrak{m}_{\le 0,2}(i_n(\omega))(\omega, j) \right)+r_{\infty}(i_n(\omega))(\omega,j)\\
 &=:\ii (W(j) + \frac{j\pi}{6}\vec{D}(j)\cdot\mathbb{A}^{-1}(\omega-\overline{\omega}) )  + \tilde{r}_{\infty}(\omega,j).
\end{align*}
Using \eqref{measure_est_11s1}, \eqref{measure_est_11s3}, \eqref{measure_est_11s4}, we have that for each $\omega_1,\omega_2\in R_{ljk}(i_n)$,
\begin{align}\label{rxcxcxxcqmcxc111}
|\tilde{r}_\infty(\omega_1,j) - \tilde{r}_\infty(\omega_2,j)|\le_{\mathtt{p_e}}  \epsilon^{9-6b}\gamma^{-1}|\omega_1-\omega_2| = \epsilon^{9-8b}|\omega_1-\omega_2|,
\end{align}
for all $j\in S_{\mathtt{M}}^\perp$ such that $j< C_{\mathtt{H}2}\max_{j\in S^+}|j|$.  Therefore, the expression of $\phi$ in \eqref{ab11xxcx1} gives us that for $\omega\in R_{ljk}(i_n)$,
\begin{align*}
\phi(\omega) &= \omega \cdot l - \frac{1}{\ii} (d_\infty(i_n(\omega))(\omega,j) - d_\infty(i_n(\omega))(\omega,k))\\
& = \omega\cdot l - (W(j) - W(k) + \frac{\pi}{6}(j\vec{D}(j)-k\vec{D}(k))\cdot\mathbb{A}^{-1}(\omega-\overline{\omega})  ) \\
& - \frac{1}{\ii}(\tilde{r}_\infty(\omega,j)-\tilde{r}_\infty(\omega,k))\\
& =  \omega\cdot l \\
& - (W(j) - W(k) +\frac{\pi}{6}\frac{(j\vec{D}(j)-k\vec{D}(k))}{W(j) -W(k)}\cdot\mathbb{A}^{-1}(\omega-\overline{\omega})(W(j) - W(k)))\\
&\ - \frac{1}{\ii}(\tilde{r}_\infty(\omega,j)-\tilde{r}_\infty(\omega,k))\\
&\overset{\eqref{degenerxxx11}}= \left(\omega-\overline{\omega}\right)\cdot l - \frac{\pi}{6}\frac{(j\vec{D}(j)-k\vec{D}(k))}{W(j) -W(k)}\cdot\mathbb{A}^{-1}(\omega-\overline{\omega})\overline{\omega}\cdot l \\
& - \frac{1}{\ii}(\tilde{r}_\infty(\omega,j)-\tilde{r}_\infty(\omega,k))\\
&\overset{\eqref{mathbsk2sdx}}= \left((\mathbb{A}^{-1})^{T}\mathbb{A}^T(\omega-\overline{\omega})\right)\cdot l  - \left(\left(\mathbb{A}^{-1}\right)^{T} \frac{\pi}6\mathbb{B}_{j,k}(\omega-\overline{\omega})\right)\cdot l \\
& - \frac{1}{\ii}(\tilde{r}_\infty(\omega,j)-\tilde{r}_\infty(\omega,k))\\
& = \left(\mathbb{A}^{-1}\right)^{T}(\mathbb{A}^T - \frac{\pi}6\mathbb{B}_{j,k})(\omega- \overline{\omega}) \cdot l- \frac{1}{\ii}(\tilde{r}_\infty(\omega,j)-\tilde{r}_\infty(\omega,k))\\
&\overset{\eqref{mathbsk2sdx}}=\left( \mathbb{A}^{-1}\right)^T \mathbb{C}_{j,k}(\omega-\overline{\omega})\cdot l - \frac{1}{\ii}(\tilde{r}_\infty(\omega,j)-\tilde{r}_\infty(\omega,k))\\
&=\mathbb{C}_{j,k}^{T}\mathbb{A}^{-1}l \cdot (\omega - \overline{\omega})- \frac{1}{\ii}(\tilde{r}_\infty(\omega,j)-\tilde{r}_\infty(\omega,k)).
\end{align*}
Therefore, for $\omega_1,\omega_2\in R_{ljk}(i_n)$, denoting
\begin{align}\label{lac2c}
l_{ac}:=\frac{\mathbb{C}_{j,k}^{T}\mathbb{A}^{-1}l}{|\mathbb{C}_{j,k}^{T}\mathbb{A}^{-1}l|},
\end{align}
we have
\begin{align}
\phi(\omega_1)-\phi(\omega_2) & = l_{ac}\cdot(\omega_1-\omega_2) \nonumber \\
& - \frac{1}{\ii}\left((\tilde{r}_\infty(\omega_1,j)-\tilde{r}_\infty(\omega_1,k)) - (\tilde{r}_\infty(\omega_2,j)-\tilde{r}_\infty(\omega_2,k))\right)
\label{rkxcxcxcs22xxx}
\end{align}
As in the proof of Lemma~\ref{lemma12232xsetmas}, we denote
\begin{align*}
l_{ac}^\perp &:= \left\{ \vec{w}\in \mathbb{R}^\nu : \vec{w}\cdot l_{ac} =0, \text{ and } sl_{ac} + \vec{w}\in R_{ljk}(i_n),\text{ for some $s\in\mathbb{R}$}\right\},\\
I_{l_{ac}}(\vec{w})&:=\left\{ s\in\mathbb{R}: sl_{ac} + \vec{w}\in R_{ljk}(i_n)\right\}
\end{align*}
For each $\vec{w}\in l_{ac}^\perp$ and $s_1,s_2\in I_{l_{ac}}(\vec{w})$, it follows from \eqref{rkxcxcxcs22xxx} that denoting $\Psi(s):=\phi(s l_{ac} + \vec{w})$,
\begin{align*}
|\Psi(s_1)-\Psi({s_2})| &\ge |\mathbb{C}_{j,k}^{T}\mathbb{A}^{-1}l||s_1-s_2| \\
& - |\left((\tilde{r}_\infty(\omega_1,j)-\tilde{r}_\infty(\omega_1,k)) - (\tilde{r}_\infty(\omega_2,j)-\tilde{r}_\infty(\omega_2,k))\right)|\\
&\overset{\eqref{rxcxcxxcqmcxc111}, \eqref{rpishsnsdiwsxcsd}, \eqref{oosdhwwdwd1ssd}}{\ge_{\mathtt{p_e}}}|s_1-s_2|, 
\end{align*}
for sufficiently small $\epsilon$, depending on $\mathtt{p_e}$. On the other hand, that $\omega_1,\omega_2\in R_{ljk}(i_n)$ implies 
\[
|\Psi(s_1)-\Psi({s_2})|  = |\phi(\omega_1)-\phi(\omega_2)|\overset{\eqref{define_rlkj}}{\le_{\mathtt{p_e}}}\gamma \langle l \rangle^{-\tau}|\lambda_\alpha(j) - \lambda_\alpha(k)| \le_{\mathtt{p_e}}\gamma \langle l \rangle^{-\tau},
\]
where the last inequality follows from that  $j,k< C_{\mathtt{H}2}$.
This implies $|I_{l_{ac}}(\vec{w})|\le_{\mathtt{p_e}} \gamma \langle l \rangle^{-\tau}$ for all $\vec{w}\in  l^\perp_{ac}$. As in the proof of Lemma~\ref{smalljklargel}, we have $|l_{ac}^\perp|\le_{\mathtt{p_e}} \epsilon^{2(\nu -1)}$, therefore,
\[
|R_{ljk}(i_n)| = \int_{l^\perp_{ac}}\int_{I_{l_{ac}}(\vec{w})} 1 ds d\vec{w} \le_{\mathtt{p_e}}\epsilon^{2(\nu -1)}\gamma \langle l \rangle^{-\tau},
\]
which finishes the proof.
\end{proof}


\chapter{Proof of Theorem 6.5}\label{Final_chapter_3}
\begin{proofthm}{main1} Let us fix $\mathtt{p}$ as in \eqref{parametersets1231} and fix $S_0^+,\mathtt{M}$ and $S$ so that \ref{tangent_2},\ref{tangent_1} and the hypotheses \ref{hypos1s1}-\ref{hypothsdj22} hold. Let $\mathcal{C}_\epsilon:=\cap_{n\ge 0}\mathcal{G}_n$. Then, Proposition~\ref{nashmoser2d} tells us that the sequence of embeddings $i_n$, constructed in \eqref{iteration_step2}, is well-defined for $\omega\in \mathcal{C}_\epsilon$ for sufficiently small $\epsilon$, depending on $\mathtt{p},b,\mathtt{M},S_0$ and satisfies that for each $\omega\in \mathcal{C}_\epsilon$,
\begin{enumerate}[label=(\arabic*)]
\item\text{$\mathcal{F}_\omega(i_n(\omega)) \to 0$ as $n\to \infty$ in $H^{s_0}_{\varphi,x}$, (see \eqref{jjsd2asdsa1})}.\label{ressults1}
\item \text{$i_{n}(\omega)=i_{triv} + \mathfrak{I}_n(\omega)$ is Cauchy sequence in $H^{s_0 + 2\mu_{\mathtt{p}}}_{\varphi,x}$ (see \eqref{lowernorm_11232} and \eqref{iteration_step2})}\label{ressults2}.
\item  $\rVert \mathfrak{I}_n(\omega)\rVert_{s_0+2\mu_{\mathtt{p}}}\le_{\mathtt{p},\mathtt{M},S_0,b}\epsilon^{6-4b}$ (see \eqref{rosisidsd22sd}).\label{ressults4}
\item $i_n(\omega)\in X_i\cap X_{i,\mathtt{M}}$, that is, $i_n(\omega)$ is a reversible and $\frac{2\pi}{\mathtt{M}}$-translation invariant, (see \ref{embeddingsdsd2sd1} of Proposition~\ref{nashmoser2d}).\label{ressults3}
\end{enumerate}
Therefore, for $\omega\in \mathcal{C}_\epsilon$, we have that $i_\infty(\omega):=\lim_{n\to\infty}i_n(\omega)$ is  reversible and $\frac{2\pi}{\mathtt{M}}$-translation invariant, and satisfies that $\mathcal{F}_\omega(i_\infty(\omega))=0$ with the  estimate
\begin{align}\label{lastsdsdw2s2s}
\rVert \mathfrak{I}_\infty(\omega)\rVert_{s_0+2\mu_{\mathtt{p}}}\le_{\mathtt{p},\mathtt{M},S_0,b} \epsilon^{6-4b},\text{ where $\mathfrak{I}_\infty(\omega):=i_\infty(\omega)-i_{triv}$.}
\end{align}
From Proposition~\ref{set_measure_last} and \eqref{rosisidsd22sd}, we see that
\[
|\mathcal{C}_\epsilon| \ge |\Omega_\epsilon| - C_{\mathtt{p},\mathtt{M},S_0,b}\epsilon^{2(\nu-1)}\gamma,
\]
Furthermore \eqref{frequency_set1} and \eqref{rpishsnsdiwsxcsd} imply that $|\Omega_\epsilon|\ge c_{\mathtt{p},\mathtt{M},S_0,b}\epsilon^{2\nu}$. Therefore,
\[
\lim_{\epsilon\to 0}\frac{|\mathcal{C}_\epsilon|}{|\Omega_\epsilon|} \ge 1 -C_{\mathtt{p},\mathtt{M},S_0,b}\lim_{\epsilon\to 0}\gamma\epsilon^{-2} \overset{\eqref{frequency_set2},\gamma=o(\epsilon^{2})}=1.
\]
Since $\mathcal{C}_\epsilon\overset{\eqref{setoffreqences}}\subset \Omega_\epsilon$, we obtain \eqref{rpsdsetfsd}.

To finish the proof, it suffices to show that $i_\infty(\omega)$ is linearly stable.  Let $I(t):=(\dot{\theta}(t),\dot{y}(t),\dot{z}(t))$ be $\frac{2\pi}{\mathtt{M}}$-translation invariant and solve the linearized Hamiltonian system at $i_\infty(\omega)(\omega t)$, that is, 
\begin{align}\label{risdsdsdsd2222sd}
\dot{I}(t) = d_i X_{H_\zeta}(i_\infty(\omega)(\omega t))[I(t)].
\end{align}
Since $H_\zeta$ is $\frac{2\pi}{\mathtt{M}}$-translation invariant Hamiltonian, one can easily see that $z(t)$ is also $\frac{2\pi}{\mathtt{M}}$-translation invariant, that is $z(t)\in X_{\mathtt{M}}$. 
Since $i_\infty(\omega)$ solves $\mathcal{F}_\omega(i_\infty(\omega))=0$, Lemma~\ref{linearlsnsdosd} tells us that
\begin{align}\label{alsdjsdw}
A(t):=(\psi(t),\eta(t),w(t)):=DG_\delta(\bar{i}(\omega t))^{-1}[I(t)]
\end{align}
solves 
\begin{align}\label{rjsdeuqsdx}
\colvec{\dot{\psi}(t) \\ \dot{\eta}(t) \\ \dot{w}(t)} = \colvec{K_{20}(\omega t)[\eta(t)] + K_{11}(\omega t)^{T}[w(t)] \\ 0 \\ \partial_x K_{02}(\omega t)[w(t)] + K_{11}(\omega t)[\eta(t)]}.
\end{align}
Clearly, $|\psi(t)|$ is bounded since $\psi(t)\in\mathbb{T}^\nu$, and $\eta(t)$ does not evolve in time, hence $\eta(t)=\eta(0)$.  Therefore, if we prove that the Sobolev norm of $w$ is bounded uniformly in time, then the estimate for $DG_\delta$ in Lemma~\ref{G_delta_estimate} and \eqref{alsdjsdw} implies the linear system\index{linear system} \eqref{risdsdsdsd2222sd} is stable (see Remark~\ref{jrjsdpppwwdds}).
Denoting $f(\omega t):=K_{11}(\omega t)\eta(0)$, we see that the evolution of $w$ is described by 
\begin{align}\label{riiisdsdwdw}
\dot{w}(t) -\partial_x K_{02}(\omega t)[w(t)] =f(\omega t).
\end{align}
  Lemma~\ref{Gsers} and \eqref{alsdjsdw} imply that $w(t)\in X_{\mathtt{M}}$ and Lemma~\ref{mfisodswwdsd2s2} implies $f(\omega t)\in X_{\mathtt{M}}$. Indeed, in Chapter~\ref{reduction}, (specially in Section~\ref{loidneks1s2}), we proved that there exist $\varphi$-dependent transformations $\Phi_{1,\infty}(\varphi),\Phi_{2,\infty}(\varphi)$ between the spaces of $\frac{2\pi}{\mathtt{M}}$-translation invariant functions in $H_{S^\perp}$ such that
 \begin{align}
 &\Phi_{1,\infty}(\varphi)\circ (\D_\omega - \partial_{x}K_{02}(\varphi))\circ \Phi_{2,\infty}(\varphi) =\mathcal{L}_\infty= \D_\omega - D_\infty, \text{ cf. \eqref{dixcxc2sdcs},\eqref{liabssdw2sss2sd2}},\label{rkksddw1q1}\\
 &\text{ $D_\infty=\text{diag}_{j\in S^\perp_{\mathtt{M}}}d_\infty(\omega,j)$ is real and reversible (see \ref{kamredu13} of Proposition~\ref{kamresdsd})}.\label{rkssdwdwdw1}
 \end{align}
 Using \eqref{rkksddw1q1} and \eqref{riiisdsdwdw}, and denoting 
  \begin{align}\label{rj2221jsdsd2}
  v(t):=\left(\Phi_{2,\infty}(\omega t)\right)^{-1}w(t),\quad g(\varphi):=\Phi_{1,\infty}(\varphi)[f(\varphi)],
  \end{align}
  we see that $v$ satisfies
  \begin{align}\label{rjjsodwksdsdj}
  \dot{v}(t) - D_\infty v(t) = g(\omega t), \quad v(0)=v_0 = \Phi_{2,\infty}(0)^{-1}w(0).
  \end{align}
  whose solution is given by
  \begin{align*}
  v(t)&  = v(t,x) = \sum_{j\in S_{\mathtt{M}}^\perp}\widehat{v}_j(t)e^{\ii j x},\\ \widehat{v}_j(t) & :=(\widehat{v_0})_je^{d_\infty(\omega,j)t} + \sum_{l\in \mathbb{Z}^\nu}\frac{\widehat{g}_j(l)}{\ii \omega \cdot l - d_\infty(\omega,j)}e^{\ii \omega\cdot l t}.
  \end{align*}
 From \eqref{rkssdwdwdw1}, we have that $d_\infty(\omega,j)$ is purely imaginary. Hence, the Sobolev norm\index{Sobolev norm} of $v$ is uniformly bounded, therefore  the Sobolev norm of $w$ does not grow in time.
\end{proofthm}

\appendix


\chapter{Conjugation with the operators $\Lambda^{\alpha-1}_a,\Upsilon^{\alpha-3}_a$}\label{Appeo100021}
In this chapter, we compute the conjugation of the operators $\Lambda^{\alpha}_a,\Upsilon^{\alpha-3}$ with the transformation $\Psi_1$ in \eqref{A_def}. We recall that $B(x):=x+\beta(x)$ and $\Psi_1^{T}[h](x) = h(B^{-1}(x))$. We simply denote $\Psi:=\Psi_1$ and  use the notation from \eqref{reduc_1_notations}.
\begin{lemma}\label{conj_lem1}
Given $a=a(x,y)\in C^{\infty}(\mathbb{T}^2)$, let us denote $\tilde{a}(x):=a(x,x).$ We have that 
\begin{align*}
&\Psi^T\Lambda^{\alpha-1}_{a}\Psi h =  a_1\Lambda^{\alpha-1}h + a_2h + \frac{1}{2}\partial_x(a_1)\partial_x\Upsilon^{\alpha-3}h -\Upsilon^{\alpha-3}_{a_3}h,\\
&\Psi^{T}[a\mathcal{A}[h]]=\Psi^{T}[aB']h,\\
&\Psi^T\Upsilon^{\alpha-3}_a\Psi[h] = \Upsilon^{\alpha-3}_{a_4}h,
\end{align*}
where 
\begin{align*}
&a_1=\Psi^{T}[(B')^{\alpha}\tilde{a}],\\
&a_2=\Psi^{T}[\Lambda^{\alpha-1}_aB']  + \int (2-2\cos(x-y))^{-\frac{\alpha}{2}}(\mathfrak{a}_1(x,y)-\mathfrak{a}_1(x,x))dy \\
&a_3=\frac{\mathfrak{a}_1(x,y)-\mathfrak{a}_1(x,x) -\partial_y\mathfrak{a}_1(x,x)\sin(x-y) }{2-2\cos(x-y)},\\
&\mathfrak{a}_1(x,y) = \left( \frac{2-2\cos(B^{-1}(x)-B^{-1}(y))}{2-2\cos(x-y)}\right)^{-\frac{\alpha}{2}}B_2^{-1}[a](x,y),\\
&a_4=\left(\frac{2-2\cos(B^{-1}(x)-B^{-1}(y))}{2-2\cos(x-y)} \right)^{1-\frac{\alpha}{2}}B_2^{-1}[a](x,y).
\end{align*}
\end{lemma}

\begin{proof}
The conjugation $\Psi^T[a\Psi[h]]$ and $\Psi^T\Upsilon^{\alpha-3}_a\Psi[h]$ follows in the same way as below. Therefore, we compute $\Psi^T\Lambda^{\alpha-1}_a\Psi[h]$ only.

Since $\Psi^T\Lambda^{\alpha-1}_a\Psi[h](x) = \Lambda^{\alpha-1}_a(\Psi h)(B^{-1}(x))$, direct computations show that
\begin{align}\label{conju_112121}
\Lambda^{\alpha-1}_{a}&(\Psi h)(B^{-1}(x))\nonumber \\
&= \int(2-2\cos(B^{-1}(x)-B^{-1}(y)))^{-\frac{\alpha}{2}} B_2^{-1}[a](x,y) \nonumber \\
& \qquad \qquad \times (B'(B^{-1}(x))h(x)-B'(B^{-1}(y))h(y))\frac{1}{B'(B^{-1}(y))}dy\nonumber\\
& = \int(2-2\cos(B^{-1}(x)-B^{-1}(y)))^{-\frac{\alpha}{2}}B_2^{-1}[a](x,y) \\
& \qquad \qquad \qquad \qquad \qquad \times \left(\frac{B'(B^{-1}(x))}{B'(B^{-1}(y))} -1 \right)dyh(x)\nonumber\\
& \ + \int(2-2\cos(B^{-1}(x)-B^{-1}(y)))^{-\frac{\alpha}{2}}B_2^{-1}[a](x,y)(h(x)-h(y))dy\nonumber\\
& =: A_1 + A_2.
\end{align}
For $A_1$, it follows from the change of variables $y\mapsto  B(y)$ that
\begin{align}\label{conj_A1}
A_1 = \Psi^{T}[\Lambda^{\alpha-1}_aB'] h.
\end{align}
For $A_2$,  we denote
\begin{align}\label{a1_appendix_def}
\mathfrak{a}_1(x,y):=\left( \frac{2-2\cos(B^{-1}(x)-B^{-1}(y))}{2-2\cos(x-y)}\right)^{-\frac{\alpha}{2}}B_2^{-1}[a](x,y),
\end{align}
so that $A_2=\Lambda^{\alpha-1}_{\mathfrak{a}_1}h$.
Writing 
\[
\mathfrak{a}_1(x,y)=\mathfrak{a}_1(x,x) + \left(\mathfrak{a}_1(x,y)-\mathfrak{a}_1(x,x)\right),
\]
we have that
\begin{equation}\label{conj_62}
\begin{aligned}
A_2 &= \mathfrak{a}_1(x,x)\Lambda^{\alpha-1}h + \int (2-2\cos(x-y))^{-\frac{\alpha}{2}}(\mathfrak{a}_1(x,y)-\mathfrak{a}_1(x,x))dyh(x) \\
& \ - \int (2-2\cos(x-y))^{-\frac{\alpha}{2}}(\mathfrak{a}_1(x,y)-\mathfrak{a}_1(x,x))h(y)dy.
\end{aligned}
\end{equation}
For the last integral, we write 
\begin{align*}
\mathfrak{a}_1(x,y)-\mathfrak{a}_1(x,x)&=\partial_y\mathfrak{a}_1(x,x)\sin(x-y) \\
& + \left( \mathfrak{a}_1(x,y)-\mathfrak{a}_1(x,x) -\partial_y\mathfrak{a}_1(x,x)\sin(x-y) \right),
\end{align*}
and $\mathfrak{a}_2(x,y):=\frac{\mathfrak{a}_1(x,y)-\mathfrak{a}_1(x,x) -\partial_y\mathfrak{a}_1(x,x)\sin(x-y) }{2-2\cos(x-y)}$, so that 
\begin{align*}
& \int (2-2\cos(x-y))^{-\frac{\alpha}{2}}(\mathfrak{a}_1(x,y)-\mathfrak{a}_1(x,x))h(y)dy \\
& = \partial_y\mathfrak{a}_1(x,x)\frac{1}{2-\alpha}\partial_x\Upsilon^{\alpha-3}h +\Upsilon^{\alpha-3}_{\mathfrak{a}_2}h. 
\end{align*}
Therefore, from \eqref{conju_112121}, \eqref{conj_A1}, \eqref{conj_62}, we have that
\begin{equation}\label{conju_247}
\begin{aligned}
& \Psi^T\Lambda^{\alpha-1}_a\Psi h \\
&  = \mathfrak{a}_1(x,x)\Lambda^{\alpha-1}h + \Big( \Psi^{T}[\Lambda^{\alpha-1}_aB']  \\
& \qquad \qquad \qquad \left.+ \int (2-2\cos(x-y))^{-\frac{\alpha}{2}}(\mathfrak{a}_1(x,y)-\mathfrak{a}_1(x,x))dy\right) h\\
& \ -  \partial_y\mathfrak{a}_1(x,x)\frac{1}{2-\alpha}\partial_x\Upsilon^{\alpha-3}h  - \Upsilon^{\alpha-3}_{\mathfrak{a}_2}h.
\end{aligned}
\end{equation}
Direct computations from \eqref{a1_appendix_def} shows that
\begin{align*}
&\mathfrak{a}_1(x,x) = \Psi^{T}[(B')^{\alpha}\tilde{a}], \quad \text{ where $\tilde{a}(x):=a(x,x)$},\\
&\partial_y\mathfrak{a}_1(x,x) = \frac{1}{2}\partial_x \tilde{\mathfrak{a}}_1(x), \quad \text{ where $\tilde{\mathfrak{a}}_1(x):=\mathfrak{a}_1(x,x)$}.
\end{align*}
Plugging this into  \eqref{conju_247}, the result follows.
\end{proof}

\chapter{Analysis of the non-resonance condition}\label{applanxxx}
In this appendix, we aim to prove Proposition~\ref{freq_amp}. Recall from \eqref{moreexplicits2} and \eqref{mutsd2xcxcxc} that 
\begin{equation}\label{mu1ltip1liersss}
\begin{aligned}
m^\circ_{1,\alpha}(j)&=\int_{\mathbb{T}}(2-2\cos(\theta))^{-\frac{\alpha}2}(1-e^{-\ii j \theta})d\theta,\\ T_{\alpha}&= \int_{\mathbb{T}}(2-2\cos(\theta))^{1-\frac{\alpha}2},\\
\lambda_\alpha(j)&=jm_{1,\alpha}(j)\\
m^\circ_{3,\alpha}(j,k)&=m^\circ_{4,\alpha}(j,-j,k,-k) \\
& = \int_{\mathbb{T}}(2-2\cos x)^{-1-\frac{\alpha}2}(2-2\cos(jx))(2-2\cos(kx))dx,\\
\kappa_j &= -\frac{\pi}2\left( m^\circ_{1,\alpha}(j) - \frac{T_\alpha}2\right).
\end{aligned}
\end{equation}
We first derive expressions of $\mathbb{A}$ and $\vec{D}(\xi)$ in terms of $\lambda_\alpha,m_{1,\alpha}$ and $m_{3,\alpha}$ and analyze the dominating terms in their expressions. Lastly, we will prove Proposition~\ref{freq_amp}.

\section{Expressions of $\mathbb{A}$ and $\vec{D}(\xi)$}
\subsection{Expression of $\mathbb{A}$.}

Plugging \eqref{mu1ltip1liersss} into \eqref{moreexplicits2}, we have
\begin{equation}\label{skd2sdsd0123}
\begin{aligned}
H_{4,j,-j,k,-k} &= -\frac{\alpha(\alpha+4)\pi}{192}(m^\circ_{1,\alpha}(j) + m^\circ_{1,\alpha}(k)) \\
& \ - \frac{\alpha^2\pi}{384}(m^\circ_{1,\alpha}(0) +m^\circ_{1,\alpha}(j+k)+m^\circ_{1,\alpha}(j-k)) \\
& \ +\frac{\alpha\pi}{192}m^\circ_{3,\alpha}(j,k) + \frac{\alpha(\alpha+2)\pi T_\alpha}{192}\\
& = -\frac{\alpha \pi}{192}\left(({\alpha+4})(m^\circ_{1,\alpha}(j) + m^\circ_{1,\alpha}(k)) \right. \\
& \left. \ \qquad \qquad + \frac{\alpha}{2}(m^\circ_{1,\alpha}(0) + m^\circ_{1,\alpha}(j+k) + m^\circ_{1,\alpha}(j-k))\right. \\
& \qquad \qquad \qquad \left.- m^\circ_{3,\alpha}(j,k)-(\alpha+2)T_\alpha \right),\\
H_{3,j,k,-(j+k)} &= \frac{\alpha \pi}{24}(m^\circ_{1,\alpha}(j)+m^\circ_{1,\alpha}(k) + m^\circ_{1,\alpha}(j+k) - T_\alpha),\\
H_{3,j,-k,-(j-k)} &= \frac{\alpha \pi}{24}(m^\circ_{1,\alpha}(j)+m^\circ_{1,\alpha}(k) + m^\circ_{1,\alpha}(j-k) - T_\alpha).
\end{aligned}
\end{equation}
Therefore, plugging this into \eqref{h4h52s}, we have that if $j=k$,
\begin{align*}
&H^{(3)}_{4,j,-j,j,-j} \\
&= H_{4,j,-j,j,-j} - \frac{3}{4}\frac{j (H_{3,j,j,-2j})^2}{j\kappa_{2j}-j\kappa_j} \\
& = -\frac{\alpha \pi}{192}\left(2(\alpha+4)m^\circ_{1,\alpha}(j) + \frac{\alpha}2(m^\circ_{1,\alpha}(2j) + 2m^\circ_{1,\alpha}(0)) - m^\circ_{3,\alpha}(j,j) -(\alpha+2)T_\alpha \right) \\
& +\frac{\alpha^2\pi}{2\cdot 192}\frac{(2m^\circ_{1,\alpha}(j) + m^\circ_{1,\alpha}(2j) - T_\alpha)^2}{m^\circ_{1,\alpha}(2j) - m^\circ_{1,\alpha}(j)}\\
& = \frac{\alpha \pi}{192}\left(\frac{\alpha}{2}\frac{(2m^\circ_{1,\alpha}(j) + m^\circ_{1,\alpha}(2j) - T_\alpha)^2}{m^\circ_{1,\alpha}(2j) - m^\circ_{1,\alpha}(j)}  \right.
\\& \qquad  - \left(2(\alpha+4)m^\circ_{1,\alpha}(j) + \frac{\alpha}2(m^\circ_{1,\alpha}(2j) + 2m^\circ_{1,\alpha}(0)) - m^\circ_{3,\alpha}(j,j) -(\alpha+2)T_\alpha \right) \Bigg).
\end{align*}
If $j\ne k$, we have
\begin{align*}
& H^{(3)}_{4,j,-j,k,-k} \\
& \hspace{-0.005cm}=H_{4,j,-j,k,-k} - \frac{3}{4}\left( \frac{(j+k) (H_{3,j,k,-(j + k)})^2}{(j + k)\kappa_{j+k} -j\kappa_j -k\kappa_{k}} + \frac{(j - k) (H_{3,j,-k,-(j - k)})^2}{(j-k)\kappa_{j-k} -j\kappa_j +k\kappa_{k}} \right)\\
& \hspace{-0.005cm} = -\frac{\alpha \pi}{192}\left(({\alpha+4})(m^\circ_{1,\alpha}(j) + m^\circ_{1,\alpha}(k)) + \frac{\alpha}{2}(m^\circ_{1,\alpha}(0) + m^\circ_{1,\alpha}(j+k) + m^\circ_{1,\alpha}(j-k)) \right. \\
& \qquad \qquad \qquad \qquad \qquad \qquad \qquad \qquad \qquad \qquad \qquad \left.  - m^\circ_{3,\alpha}(j,k)-(\alpha+2)T_\alpha \right)\\
& \ +\frac{\alpha^2 \pi}{2\cdot 192}\left(\frac{(j+k)(m^\circ_{1,\alpha}(j)+m^\circ_{1,\alpha}(k) + m^\circ_{1,\alpha}(j+k) - T_\alpha)^2}{(j+k)m_{1,\alpha}(j+k) - jm^\circ_{1,\alpha}(j) - km^\circ_{1,\alpha}(j)} \right. \\& \left.\qquad \qquad \qquad \qquad \qquad + \frac{(j-k)(m^\circ_{1,\alpha}(j)+m^\circ_{1,\alpha}(k) + m^\circ_{1,\alpha}(j-k) - T_\alpha)^2}{(j-k)m^\circ_{1,\alpha}(j-k) - jm^\circ_{1,\alpha}(j) + km^\circ_{1,\alpha}(j)} \right)\\
& \hspace{-0.005cm}= \frac{\alpha \pi}{192}\left( \frac{\alpha}{2}\left(\frac{(j+k)(m^\circ_{1,\alpha}(j)+m^\circ_{1,\alpha}(k) + m^\circ_{1,\alpha}(j+k) - T_\alpha)^2}{(j+k)m^\circ_{1,\alpha}(j+k) - jm^\circ_{1,\alpha}(j) - km^\circ_{1,\alpha}(j)} \right. \right. \\& \left.\qquad \qquad \qquad \qquad \qquad + \frac{(j-k)(m^\circ_{1,\alpha}(j)+m^\circ_{1,\alpha}(k) + m^\circ_{1,\alpha}(j-k) - T_\alpha)^2}{(j-k)m^\circ_{1,\alpha}(j-k) - jm^\circ_{1,\alpha}(j) + km^\circ_{1,\alpha}(j)} \right)  \\
&  - \left(({\alpha+4})(m^\circ_{1,\alpha}(j) + m^\circ_{1,\alpha}(k)) + \frac{\alpha}{2}(m^\circ_{1,\alpha}(0) + m^\circ_{1,\alpha}(j+k) + m^\circ_{1,\alpha}(j-k)) \right. \\ & \quad \qquad \qquad \qquad \qquad \qquad \qquad \qquad \qquad \qquad \qquad \left.- m^\circ_{3,\alpha}(j,k)-(\alpha+2)T_\alpha \right) \Bigg).
\end{align*}
Therefore, it follows from \eqref{amplitude_modulation1} that for $i,k\in \left\{ 1,\ldots,\nu\right\}$ and the tangential sites in \eqref{tan_site}, $S^+=\left\{  0 < j_1 < \cdots < j_\nu\right\}$,

\begin{align}\label{lxcxc11sdsd}
\mathbb{A}^{i}_k &= \frac{\alpha \pi}{192} A(j_i,j_k),\text{ for $i,k=1,\ldots,\nu$}
\end{align}
where
\begin{align*}
A&(j,j)  =j^2\left( \frac{\alpha}{2}\frac{(2m^\circ_{1,\alpha}(j) + m^\circ_{1,\alpha}(2j) - T_\alpha)^2}{m^\circ_{1,\alpha}(2j)- m^\circ_{1,\alpha}(j)} \right. \\
&  - \left(2(\alpha+4)m^\circ_{1,\alpha}(j) + \frac{\alpha}2(m^\circ_{1,\alpha}(2j) + 2m^\circ_{1,\alpha}(0)) - m^\circ_{3,\alpha}(j,j) -(\alpha+2)T_\alpha \right)\Bigg),
\end{align*}
and for $j\ne k$, 
\begin{align*}
A& (j,k)  = 2jk\left(\frac{\alpha}{2}\frac{(j+k)(m^\circ_{1,\alpha}(j)+m^\circ_{1,\alpha}(k) + m^\circ_{1,\alpha}(j+k) - T_\alpha)^2}{(j+k)m^\circ_{1,\alpha}(j+k) - jm^\circ_{1,\alpha}(j) - km^\circ_{1,\alpha}(k)} \right.\\
& \qquad \qquad + \frac{\alpha}{2}\frac{(j-k)(m^\circ_{1,\alpha}(j)+m^\circ_{1,\alpha}(k) + m^\circ_{1,\alpha}(j-k) - T_\alpha)^2}{(j-k)m^\circ_{1,\alpha}(j-k) - jm^\circ_{1,\alpha}(j) + km^\circ_{1,\alpha}(k)}   \\
&  \quad - \left(({\alpha+4})(m^\circ_{1,\alpha}(j) + m^\circ_{1,\alpha}(k)) + \frac{\alpha}{2}(m^\circ_{1,\alpha}(0) + m^\circ_{1,\alpha}(j+k) + m^\circ_{1,\alpha}(j-k))\right. \\
& \qquad \qquad \left.- m^\circ_{3,\alpha}(j,k)-(\alpha+2)T_\alpha \right)\Big). 
\end{align*}
To simplify the notation, let us denote
\begin{align}
\mathtt{B}_1(j,k)&:=(\alpha+4)(m^\circ_{1,\alpha}(j)+m^\circ_{1,\alpha}(k)) \nonumber \\
& + \frac{\alpha}2(m^\circ_{1,\alpha}(0) + m^\circ_{1,\alpha}(j+k) + m^\circ_{1,\alpha}(j-k))-(\alpha+2)T_\alpha,\label{b1def}\\
\mathtt{B}_2(j,k)&:=\frac{(j+k)(m^\circ_{1,\alpha}(j)+m^\circ_{1,\alpha}(k) + m^\circ_{1,\alpha}(j+k) - T_\alpha)^2}{(j+k)m^\circ_{1,\alpha}(j+k) - jm^\circ_{1,\alpha}(j) - km^\circ_{1,\alpha}(k)} \nonumber \\
& - \frac{(j-k)(m^\circ_{1,\alpha}(j)+m^\circ_{1,\alpha}(k) + m^\circ_{1,\alpha}(j-k) - T_\alpha)^2}{jm^\circ_{1,\alpha}(j)  - (j-k)m^\circ_{1,\alpha}(j-k) - km^\circ_{1,\alpha}(k)},\label{b2def}\\
\mathtt{B}_3(j)&:= \frac{\alpha }{2}\frac{(2m^\circ_{1,\alpha}(j) + m^\circ_{1,\alpha}(2j) - T_\alpha)^2}{m^\circ_{1,\alpha}(2j) - m^\circ_{1,\alpha}(j)} ,\label{b3def}
\end{align}
 so that we can write $A$ as 
 \begin{align*}
 A(j,k) = \begin{cases}
j^2 \mathtt{B}_3(j)- j^2\mathtt{B}_{1}(j,j) + j^2m^\circ_{3,\alpha}(j,j), & \text{ if $j=k$},\\
2jk\left(\frac{\alpha}2 \mathtt{B}_2(j,k) - \mathtt{B}_1(j,k) \right) + 2 jkm_{3,\alpha}^\circ(j,k),& \text{ if $j\ne k$.}
 \end{cases}
 \end{align*}
Therefore, defining
\begin{equation}\label{ejjxcxcxcsdsd11}
\begin{aligned}
A_1(j,k) & :=\begin{cases}
j^2m^\circ_{3,\alpha}(j,j), & \text{ if $j=k$,}\\
2jk m^\circ_{3,\alpha}(j,k) & \text{ if $j\ne k$},
\end{cases} \\
 A_2(j,k) & :=
\begin{cases}
  j^2\mathtt{B}_3(j)- j^2\mathtt{B}_{1}(j,j)  & \text{ if $j=k$},\\
2jk\left(\frac{\alpha}2 \mathtt{B}_2(j,k) - \mathtt{B}_1(j,k) \right) ,& \text{ if $j\ne k$,}
 \end{cases}
\end{aligned}
\end{equation}
we have $A(j,k) =A_1(j,k) + A_{2}(j,k)$. Therefore \eqref{lxcxc11sdsd} gives that
\begin{align}\label{Amatrixexpression}
\mathbb{A}^i_k &= \frac{\alpha \pi}{192}A_1(j_i,j_k) + \frac{\alpha \pi}{192}A_2(j_i,j_k)\\
&=:\frac{\alpha\pi}{192}\left(\mathbb{A}_1\right)^i_k + \frac{\alpha\pi}{192}\left(\mathbb{A}_2\right)^i_k,\text{ for $i,k\in\left\{ 1,\ldots,\nu\right\}$ for $S^+=\left\{ 0<j_1<\cdots<j_\nu\right\}$.}\nonumber
\end{align}
Note that using the evenness of $j\mapsto m^\circ_{1,\alpha}(j)$ and $m_{3,\alpha}(j,k)=m_{3,\alpha}(k,j)$, which easily follow from \eqref{mu1ltip1liersss}, we see that 
\begin{align}\label{Aeushhsdsdx}
\mathbb{A}^T = \mathbb{A},\quad \mathbb{A}_1^T = \mathbb{A}_1, \quad \mathbb{A}_2^T = \mathbb{A}_2.
\end{align}
\subsection{Expression of $\vec{D}(\xi)$.}

Using \eqref{skd2sdsd0123}, we have that for $\xi\in S^\perp$ and $j_k\in S^+=\left\{ 0 < j_1 < \cdots < j_\nu\right\}$ for $k=1,\ldots,\nu$, we have 
\begin{align*} 
&\frac{12}\pi H_{4,j_k,-j_k,\xi,-\xi} \\
&= -\frac{\alpha}{16}\left(({\alpha+4})(m^\circ_{1,\alpha}(\xi) + m^\circ_{1,\alpha}(j_k))\right. \\
& \left.\qquad + \frac{\alpha}{2}(m^\circ_{1,\alpha}(0) + m^\circ_{1,\alpha}(\xi+j_k) + m^\circ_{1,\alpha}(\xi-j_k)) - m^\circ_{3,\alpha}(\xi,j_k)-(\alpha+2)T_\alpha \right),
\end{align*}
and using $\lambda_\alpha(\xi)=\xi m^\circ_{1,\alpha}(\xi)$ for all $\xi\in \mathbb{Z}\backslash\left\{ 0 \right\}$ (see \eqref{lambdadeffsd} and \eqref{def_of23sm}),
\begin{align*}
\frac{18}{\pi^2}&\frac{(\xi + j_k)(H_{3,j_k,\xi, -(\xi +j_k)})^2}{\lambda_\alpha(\xi +j_k) -\lambda_\alpha(\xi) - \lambda_\alpha(j_k)} \\
& = \frac{\alpha^2}{32}\frac{(\xi +j_k)(m^\circ_{1,\alpha}(\xi)+m^\circ_{1,\alpha}(j_k) + m^\circ_{1,\alpha}(\xi+j_k) - T_\alpha)^2}{(\xi + j_k)m_{1,\alpha}^\circ -\xi m_{1,\alpha}^\circ(\xi)- j_km^\circ_{1,\alpha}(j_k)},\\
\frac{18}{\pi^2}&\frac{(\xi - j_k)(H_{3,j_k,\xi, -(\xi -j_k)})^2}{\lambda_\alpha(\xi) -\lambda_\alpha(\xi-j_k) - \lambda_\alpha(j_k)} \\
& = \frac{\alpha^2}{32}\frac{(\xi -j_k)(m^\circ_{1,\alpha}(\xi)+m^\circ_{1,\alpha}(j_k) + m^\circ_{1,\alpha}(\xi-j_k) - T_\alpha)^2}{\xi m^\circ_{1,\alpha}(\xi) -(\xi-j_k)m^\circ_{1,\alpha}(\xi-j_k) - j_k m^\circ_{1,\alpha}(j_k)}.
\end{align*}
Plugging this into \eqref{eigenvalusdsd1}, we have that for $k=1,\ldots,\nu$, $S^+=\left\{ 0 < j_1 < \cdots < j_\nu \right\}$ and $\xi\in S^\perp$, 
\begin{align*}
(\vec{D}(\xi))_k& = \frac{\alpha}{16}m^\circ_{3,\alpha}(\xi,j_k)j_k,\\
& -j_k\frac{\alpha}{16}\Big(({\alpha+4})(m^\circ_{1,\alpha}(\xi) + m^\circ_{1,\alpha}(j_k))  \\
& \left.\qquad \qquad + \frac{\alpha}{2}(m^\circ_{1,\alpha}(0) + m^\circ_{1,\alpha}(\xi+j_k) + m^\circ_{1,\alpha}(\xi-j_k)) -(\alpha+2)T_\alpha \right)\\
& + \frac{\alpha^2}{32} j_k \left(\frac{(\xi +j_k)(m^\circ_{1,\alpha}(\xi)+m_{1,\alpha}(j_k) + m^\circ_{1,\alpha}(\xi+j_k) - T_\alpha)^2}{(\xi + j_k)m_{1,\alpha}^\circ -\xi m_{1,\alpha}^\circ(\xi)- j_km^\circ_{1,\alpha}(j_k)} \right. \\ & \qquad \qquad \qquad \left.-  \frac{(\xi -j_k)(m^\circ_{1,\alpha}(\xi)+m^\circ_{1,\alpha}(j_k) + m^\circ_{1,\alpha}(\xi-j_k) - T_\alpha)^2}{\xi m^\circ_{1,\alpha}(\xi) -(\xi-j_k)m^\circ_{1,\alpha}(\xi-j_k) - j_k m^\circ_{1,\alpha}(j_k)}\right).
\end{align*}
Hence, using the notation in \eqref{b1def} and \eqref{b2def}, we write $\vec{D}(\xi)_k$ as
\begin{align}\label{desstimaxcxcsd}
(\vec{D}(\xi))_k = \begin{cases}
\frac{\alpha}{16}m^\circ_{3,\alpha}(\xi,j_k)j_k - \frac{\alpha}{16}j_k\mathtt{B}_1(\xi,j_k) + \frac{\alpha^2}{32}j_k\mathtt{B}_2(\xi,j_k), & \text{ if $\xi\in S^\perp$,}\\
0, & \text{ if $\xi = 0$},
\end{cases}
\end{align}
for $k=1,\ldots,\nu$, $S^+=\left\{ 0 < j_1 < \cdots < j_\nu \right\}$.

Note that using the evenness of $j\mapsto m^\circ_{1,\alpha}(j)$ and $m_{3,\alpha}(j,k)=m_{3,\alpha}(k,j)$, which easily follow from \eqref{mu1ltip1liersss}, we see that 
\begin{align}\label{evenessofvectsx}
\vec{D}(-\xi) = \vec{D}(\xi).
\end{align}

\section{Estimates for $\mathtt{B}_1,\mathtt{B}_2, \mathtt{B}_3$ and $m^\circ_{3,\alpha}$}
From \eqref{Amatrixexpression} and \eqref{desstimaxcxcsd} we see that $\mathbb{A}$ and $\vec{D}(\xi)$ are completely determined by the choice of the tangential sites $S^+$ in \eqref{tan_site}, that is, by the choice of $S_0^+$ and $\mathtt{M}$. In order to reduce the complexity of the  computation, we will take  $\mathtt{M}$ large so that the nature of $\mathbb{A}$ and $\vec{D}(\xi)$ is determined by the asymptotic behavior of $m^\circ_{1,\alpha}$ and $m_{3,\alpha}^\circ$. To do so, we first prove the necessary estimates for $\mathtt{B}_1,\mathtt{B}_1$ and $m^\circ_{3,\alpha}$ in this section.

\begin{lemma}\label{lambdaamdmwes}
For $\alpha\in (1,2)$ and $j,j_1,j_2,k\in \mathbb{Z}\backslash\left\{ 0 \right\}$ such that $j_1\ne j_2,\ |j|,|j_1|,|j_2|\ne |k|$, we have that
\begin{align}
|m_{1,\alpha}(j+k)-m_{1,\alpha}(j)|\le_\alpha |k||j|^{\alpha-2},\label{p1qosiixx1}\\
|2\lambda_{\alpha}(j) - \lambda_{\alpha}(j+k)-\lambda_{\alpha}(j-k)|&\le_\alpha |k|^2| j |^{\alpha-2},\label{m1sess3}
\end{align}
and
\begin{align}
 |m_{1,\alpha}(j_1+k) - m_{1,\alpha}(j_2 +k)|&\le_\alpha |j_1-j_2|(|j_1| +|j_2|)^{\alpha-2},\label{m1sess234}\\
|(m_{1,\alpha}(j_1+k)-m_{1,\alpha}(j_1))&-(m_{1,\alpha}(j_2+k)-m_{1,\alpha}(j_2))|\nonumber\\
&\le_{\alpha} |j_1-j_2|(|j_1| + | j_2|)^{\alpha-3}|k|,\label{p1qosiixx2}\\
|(2\lambda_{\alpha}(j_1) - \lambda_{\alpha}(j_1+k)-\lambda_{\alpha}(j_1-k)) &- (2\lambda_{\alpha}(j_2) - \lambda_{\alpha}(j_2+k)-\lambda_{\alpha}(j_2-k))|\nonumber\\
&\le_\alpha |j_1-j_2|(| j_1| + | j_2 |)^{\alpha-3}|k|^2.\label{m2sess3}
\end{align}
\end{lemma}
\begin{proof}
Since the proof for each estimate is similar, we will prove \eqref{m2sess3} only.  Let us denote
\[
T:=|(2\lambda_{\alpha}(j_1) - \lambda_{\alpha}(j_1+k)-\lambda_{\alpha}(j_1-k)) - (2\lambda_{\alpha}(j_2) - \lambda_{\alpha}(j_2+k)-\lambda_{\alpha}(j_2-k))|.
\]

We claim that
\begin{align}\label{climasx22}
T \le_\alpha |j_1-j_2|(| j_1|^{\alpha-3} + | j_2 |^{\alpha-3})|k|^2.
\end{align}
Let us assume the claim for a moment and derive the estimate \eqref{m2sess3}. If $|j_1| \ge 2|j_2|$, then it follows from \eqref{m1sess3} that
\[
T \le_\alpha |k|^2 |j_1|^{\alpha-2} ,
\]
while, 
\[
|j_1-j_2|(| j_1| + | j_2 |)^{\alpha-3}|k|^2\ge_\alpha |j_1|^{\alpha-2}|k|^2,
\]
therefore, we obtain \eqref{m2sess3}. Similarly, we have \eqref{m2sess3} in case $|j_1| \le \frac{1}{2}|j_2|$. If $\frac{1}{2}|j_2| \le |j_1|\le 2|j_2|$, then clearly the estimate \eqref{climasx22} is equivalent to \eqref{m2sess3}. Therefore, it suffices to prove \eqref{climasx22} to obtain \eqref{m2sess3}.

Towards \eqref{climasx22}, we see from Proposition~\ref{nichts2}  that $\xi\mapsto \lambda_\alpha(\xi)$ is odd. Therefore without loss of generality, we can assume that
\begin{align}\label{amspsdjxx1}
j_1>0, \ |j_1|>|j_2|, \ k>0,
\end{align}
since the estimate for this case can cover the all possible cases of $j_1,j_2,k$ in the estimate \eqref{m1sess3}. We further decompose into more cases: $j_1>0>j_2$ and $j_1>j_2>0$:

\textbf{Case $j_1>0>j_2$:} In this case, \eqref{m1sess3} and the triangular inequality yield that 
\[
T \le |k|^2(|j_1|^{\alpha-2}+|j_2|^{\alpha-2}) \overset{\eqref{amspsdjxx1}}\le 2 |j_1|^{\alpha-2}|k|^2 \overset{j_1>0>j_2}\le_\alpha |j_1-j_2| |j_1|^{\alpha-3}|k|^2,
\]
which gives \eqref{climasx22}.

For the case $j_1>j_2>0$, we divide it into three cases: $0<j_2<j_1\le 2k,\ 0< j_2\le 2k < j_1 ,\ 0< 2k < j_2 < j_1$.

\textbf{Case $0<j_2<j_1\le 2k$:}
In this case,  we have
\[
T =\left|\int_{j_2}^{j_1} 2\partial_\xi \lambda_\alpha(\xi) - \partial_\xi\lambda_\alpha(\xi + k) - \partial_\xi \lambda_\alpha(\xi - k) d\xi\right|,
\]
while each term of the integrand is bounded by $C_\alpha (|\xi|^{\alpha-1} + |k|^{\alpha-1})\le_\alpha |k|^{\alpha-1}$. Therefore, we have $
T\le |j_1-j_2||k|^{\alpha-1}$. Using that $|k|^{\alpha-1} = |k|^{\alpha-3}|k| \le_\alpha |j_1|^{\alpha-3}|k|$, the estimate \eqref{climasx22} follows.

\textbf{Case $0< j_2\le 2k < j_1$:}
We compute
\begin{align*}
T &= \left|\int_{j_2}^{j_1} 2\partial_\xi \lambda_\alpha(\xi) - \partial_\xi\lambda_\alpha(\xi + k) - \partial_\xi \lambda_\alpha(\xi - k) d\xi\right|\\
&\le \left| \int_{j_2}^{2k}2\partial_\xi \lambda_\alpha(\xi) - \partial_\xi\lambda_\alpha(\xi + k) - \partial_\xi \lambda_\alpha(\xi - k) d\xi\right| \\
& + \left| \int_{2k}^{j_1}2\partial_\xi \lambda_\alpha(\xi) - \partial_\xi\lambda_\alpha(\xi + k) - \partial_\xi \lambda_\alpha(\xi - k) d\xi\right|\\
&\le_\alpha |2k - j_2||k|^{\alpha-1} +\left|\int_{2k}^{j_1}\int_0^k\int_{-t}^t \partial_{\xi\xi\xi}\lambda_\alpha(\xi +u)dudtd\xi\right|\\
&\le_\alpha  |2k - j_2||k|^{\alpha-1} +\left| \int_{2k}^{j_1}\int_0^k\int_{-t}^{t}|\xi|^{\alpha-3}du dt d\xi\right|\\
&\le_\alpha |2k - j_2||k|^{\alpha-1}  + |j_1-2k| |k|^{\alpha-1}\\
& \le_\alpha |j_1-j_2| |k|^{\alpha-1}.
\end{align*}
Using that $ 0 < j_2 < 2k$, we have $|j_1-j_2| |k|^{\alpha-1}\le_\alpha |j_1-j_2| |j_2|^{\alpha-3}|k|^2$, therefore \eqref{climasx22} follows.

\textbf{Case $0< 2k < j_2 < j_1$:}
In this case, we have
\begin{align*}
T & = \left| \int_{j_2}^{j_1}\int_0^{k}\int_{-t}^t \partial_{\xi\xi\xi}\lambda\alpha(\xi +u)dudtd\xi\right| \\\
& \le_\alpha \int_{j_2}^{j_1}\int_0^{k}\int_{-t}^t \langle \xi +u \rangle^{\alpha-3} dudtd\xi \le_\alpha |j_1-j_2||j_2|^{\alpha-3}|k|^2,
\end{align*}
therefore we obtain \eqref{climasx22}.
\end{proof}

 \begin{lemma}
For $\alpha\in (1,2)$ and $j,j_1,j_2,k\in \mathbb{Z}\backslash\left\{ 0 \right\}$ such that $j_1\ne j_2,\ |j|,|j_1|,|j_2|\ne |k|$, we have that 
\begin{align}
|\mathtt{B}_1(j,k)|& \le_\alpha (|j|^{\alpha-1} + |k|^{\alpha-1}), \text{ for all $j,k\in\mathbb{Z},k\ne 0$}\label{b1estimate11sx}\\
|\mathtt{B}_1(j_1,k) - \mathtt{B}_{1}(j_2,k)|& \le_\alpha  |j_1-j_2|(|j_1| + | j_2|)^{\alpha-2}.\label{bestimate11sx1}
\end{align}
\end{lemma}
\begin{proof}In the expression of $\mathtt{B}_1$ in \eqref{b1def}, we can replace $m^\circ_{1,\alpha}(j,k)$ by $m_{1,\alpha}$ since $j,k\ne 0$. Then the estimates \eqref{b1estimate11sx} and \eqref{bestimate11sx1} follow immediately from the fact that $\xi\mapsto m_{1,\alpha}(\xi) \in \mathcal{S}^{\alpha-1}$ and the estimate \eqref{m1sess234}.
\end{proof}

Now, we estimate $\mathtt{B}_2$. We further decompose $\mathtt{B}_2$ as
\begin{align}
\mathtt{B}_2(j,k) &= \mathtt{B}_{21}(j,k) + \mathtt{B}_{22}(j,k),\label{b1deco1}\\
B_{21}
(j,k) & := \frac{\mathtt{B}_{21}^N(j,k)}{\mathtt{B}^D_{21}(j,k)}\nonumber\\
\mathtt{B}_{21}^N(j,k)&= (j+k)(m_{1,\alpha}(j)+m_{1,\alpha}(k) + m_{1,\alpha}(j+k) - T_\alpha)^2 \nonumber \\
& \ -(j-k)(m_{1,\alpha}(j)+m_{1,\alpha}(k) + m_{1,\alpha}(j-k) - T_\alpha)^2 ,\nonumber\\
\mathtt{B}_{21}^D(j,k)&= 
{\lambda^\circ_\alpha(j+k) - \lambda^\circ_\alpha(j) - \lambda^\circ_\alpha(k)} \label{b1deco2} \\
\mathtt{B}_{22}(j,k) & :=\frac{\mathtt{B}^N_{22}(j,k)}{\mathtt{B}^D_{22}(j,k)}.\nonumber \\
\mathtt{B}_{22}^N(j,k) & = (j-k)(m_{1,\alpha}(j)+m_{1,\alpha}(k) + m_{1,\alpha}(j-k) - T_\alpha)^2 \nonumber \\
& \ \times (2\lambda_\alpha^\circ(j) - \lambda_\alpha^\circ(j-k) - \lambda_\alpha^\circ(j+k)) \nonumber \\
\mathtt{B}_{22}^D(j,k)  & = (\lambda^\circ_\alpha(j+k) - \lambda^\circ_\alpha(j) - \lambda^\circ_\alpha(k))(\lambda^\circ_\alpha(j) -\lambda^\circ_\alpha(j-k)- \lambda^\circ_\alpha(k)) \label{b1deco3}
\end{align}

\begin{lemma}\label{rotlqkftocx}
For $\alpha\in (1,2)$ and $j,j_1,j_2,k\in \mathbb{Z}\backslash \left\{ 0 \right\}$ such that $|j|,|j_1|,|j_2|\ne |k|$, it holds that 
\begin{align}
|\mathtt{B}^N_{21}(j,k)|&\le_\alpha  |k|(|j|^{2\alpha-2} + |k|^{2\alpha-2}),\label{resumoxcqq1}\\
|\mathtt{B}^N_{21}(j_1,k)-\mathtt{B}^N_{21}(j_2,k)|&\le_\alpha|k||j_1-j_2| (|j_1|+|j_2|)^{\alpha-2}(|j_1|^{\alpha-1}+|j_2|^{\alpha-1}+|k|^{\alpha-1}) \label{resumoxcqq2}\\
|\mathtt{B}^N_{22}(j,k)|&\le_\alpha |k|^2|j|^{\alpha-2} (|j|^{2\alpha-1} + |k|^{2\alpha-1}),\label{resumoxcqq3}\\
|\mathtt{B}^N_{22}(j_1,k)-\mathtt{B}^N_{22}(j_1,k)|&\le_\alpha|k|^2 |j_1-j_2| (|j_1|+|j_2|+|k|)^{2\alpha-1}(|j_1|+|j_2|)^{\alpha-3},\label{resumoxcqq4}\\
|\mathtt{B}^D_{21}(j,k)|&\le_\alpha |k|(|j|^{\alpha-1} + |k|^{\alpha-1}),\label{resumoxcqq5}\\
|\mathtt{B}^D_{21}(j_1,k) - \mathtt{B}^D_{21}(j_2,k)|&\le_\alpha|j_1-j_2||k| (|j_1|+|j_2|)^{\alpha-2},\label{resumoxcqq6}\\
|\mathtt{B}^D_{22}(j,k)|&\le_\alpha |k|^2(|j| + |k|)^{2\alpha-2},\label{resumoxcqq7}\\
|\mathtt{B}^D_{22}(j_1,k)-\mathtt{B}^D_{22}(j_2,k)|&\le_\alpha|j_1-j_2||k|^2(|j_1|+|j_2|)^{\alpha-2}(|j_1|+|j_2|+|k|)^{\alpha-1}.\label{resumoxcqq8}
\end{align}
\end{lemma}
\begin{proof}
\vspace{0.5\baselineskip}\noindent\textit{Proof of \eqref{resumoxcqq1}, \eqref{resumoxcqq2}.}
Applying \eqref{p1qosiixx1} and \eqref{p1qosiixx2} to 
\begin{align*}
m_{1,\alpha}(j+k)-m_{1,\alpha}(j-k) = (m_{1,\alpha}(j+k) - m_{1,\alpha}(j)) + (m_{1,\alpha}(j) - m_{1,\alpha}(j-k))
\end{align*}
we see that
\begin{align}
 |m_{1,\alpha}(j+k)-m_{1,\alpha}(j-k)|&\le_\alpha |k||j|^{\alpha-2},\label{m1sess1}
 \end{align}
 \begin{align}
 |(m_{1,\alpha}(j_1+k)-m_{1,\alpha}(j_1-k))&-(m_{1,\alpha}(j_2+k)-m_{1,\alpha}(j_2-k))| \nonumber \\
 &\le_{\alpha} |j_1-j_2|(|j_1| + | j_2|)^{\alpha-3}|k|.\label{m1sess2}
\end{align}

Let
\begin{equation}\label{l1ssd11whatsapp}
\begin{aligned}
L_1(j,k)&:=(j+k)(m_{1,\alpha}(j)+m_{1,\alpha}(k) + m_{1,\alpha}(j+k) - T_\alpha)^2\\
L_2(j,k)&:=m_{1,\alpha}(j)+m_{1,\alpha}(k)  - T_\alpha.
\end{aligned}
\end{equation}
Since $\xi\mapsto m_{1,\alpha}(\xi)\in \mathcal{S}^{\alpha-1}$ is even, we have
\begin{align}
|L_1(j,k)|&\le_\alpha |j|^{2\alpha-1} + |k|^{2\alpha-1},\label{l1eess111s1}\\
|L_2(j,k)|&\le_\alpha |j|^{\alpha-1} + |k|^{\alpha-1},\label{l1eess111s2}\\
|L_2(j,k)-L_2(j,-k)|& = 0,\label{l1eess111s3}\\
|L_2(j_1,k) - L_2(j_2,k)|&\overset{\eqref{m1sess234}}{\le_\alpha} |j_1-j_2|(|j_1| + |j_2|)^{\alpha-2}. \label{l1eess111s4}
\end{align}
Since $\mathtt{B}^N_{21}(j,k) = L_1(j,k) - L_1(j,-k)$, we expand
\begin{align}\label{rksdsdssxxxb21}
\mathtt{B}^N_{21}(j,k)&= (j+k)(L_2(j,k)+m_{1,\alpha}(j + k))^2 -  (j - k)(L_2(j,k)+m_{1,\alpha}(j - k))^2\nonumber\\
& = j((L_2(j,k)+m_{1,\alpha}(j+k))^2 - (L_2(j,k)+m_{1,\alpha}(j-k))^2) \nonumber\\
& \ +k ((L_2(j,k)+m_{1,\alpha}(j+k))^2 + (L_2(j,k)+m_{1,\alpha}(j-k))^2)\nonumber\\
& = j\underbrace{(2 L_2(j,k) + m_{1,\alpha}(j+k)+m_{1,\alpha}(j-k))}_{=:L_3(j,k)} \nonumber \\
& \qquad \qquad \qquad \qquad \qquad \times \underbrace{(m_{1,\alpha}(j+k) - m_{1,\alpha}(j-k))}_{=:L_4(j,k)}\nonumber\\
& \ +2kL_2(j,k)^2 + \underbrace{2kL_2(j,k)(m_{1,\alpha}(j+k) +m_{1,\alpha}(j-k))}_{=:L_5(j,k)} \nonumber \\
& + \underbrace{k(m_{1,\alpha}(j+k)^2+m_{1,\alpha}(j-k)^2)}_{=:L_6(j,k)}.
\end{align}
We have that
\begin{align}
|L_3(j,k)|&\overset{\eqref{l1eess111s2}}{\le_\alpha} |j|^{\alpha-1} + |k|^{\alpha-1},\label{qkfrjfmadl1}\\
|L_3(j_1,k) - L_3(j_2,k)|&\overset{\eqref{l1eess111s4},\eqref{m1sess234}}{\le_\alpha} |j_1-j_2|(|j_1| + |j_2|)^{\alpha-2} ,\label{qkfrjfmadl2}\\
|L_4(j,k)|&\overset{\eqref{m1sess1}}{\le_\alpha}|k| |j|^{\alpha-2},\label{qkfrjfmadl3}\\
|L_4(j_1,k) - L_4(j_2,k)|&\overset{\eqref{m1sess2}}{\le_\alpha}|j_1-j_2| (|j_1| + |j_2|)^{\alpha-3}|k|\label{qkfrjfmadl4}\\
|L_5(j,k)|&\overset{\eqref{l1eess111s2}}{\le_\alpha}|k|( |j|^{2\alpha-2} + |k|^{2\alpha-2}),\label{qkfrjfmadl5}\\
|L_5(j_1,k) - L_5(j_2,k)|&\overset{\eqref{l1eess111s4}\eqref{m1sess234}}{\le_\alpha}|k|  |j_1-j_2|(|j_1| + |j_2|)^{\alpha-2}(|j_1|^{\alpha-1} + |k|^{\alpha-1}),\label{qkfrjfmadl6}\\
|L_6(j,k)|&\le_\alpha |k| (|j|^{2\alpha-2} + |k|^{2\alpha-2}),\label{qkfrjfmadl7}\\
|L_6(j_1,k) - L_6(j_2,k)|&\overset{\eqref{m1sess234}}{\le_\alpha}|k|(|j_1|^{\alpha-1}+|j_2|^{\alpha-1}+|k|^{\alpha-1}) |j_1-j_2|(|j_1| + |j_2|)^{\alpha-2}\label{qkfrjfmadl8}
\end{align}
Therefore,
\begin{align*}
|\mathtt{B}^N_{21}(j,k)|&\overset{\eqref{rksdsdssxxxb21}}=|j||L_3(j,k)||L_4(j,k)|+2k|L_2(j,k)|^2 + |L_5(j,k)|+|L_6(j,k)|\\
&\overset{\eqref{qkfrjfmadl1},\eqref{qkfrjfmadl4},\eqref{qkfrjfmadl5},\eqref{qkfrjfmadl7},\eqref{l1eess111s2}}{\le_\alpha}|k|(|j|^{2\alpha-2} + |k|^{\alpha-1}|j|^{\alpha-1} + |k|^{2\alpha-2})\\
&\le_\alpha |k|(|j|^{2\alpha-2} + |k|^{2\alpha-2}),
\end{align*}
where the last inequality follows from Young's inequality. Therefore, the estimate \eqref{resumoxcqq1} follows. Towards the estimate \eqref{resumoxcqq2},  we have that 

\begin{align*}
&|\mathtt{B}^N_{21}(j_1,k) - \mathtt{B}^N_{21}(j_2,k)|\\
&\overset{\eqref{rksdsdssxxxb21}}\le |j_1-j_2| |L_3(j_1,k)||L_4(j_1,k)| \\
& \ + |j_2||L_3(j_1,k)-L_3(j_2,k)||L_4(j_1,k)| \\
& \ + |j_2||L_3(j_2,k)||L_4(j_1,k)-L_4(j_2,k)|\\
& \ + |L_2(j_1,k)+L_2(j_2,k)||L_2(j_1,k) - L_2(j_2,k)|\\
&\ + |L_5(j_1,k)-L_5(j_2,k)|+|L_6(j_1,k)-L_6(j_2,k)|\\
&\overset{\eqref{l1eess111s2},\eqref{qkfrjfmadl1}-\eqref{qkfrjfmadl8}}{\le_\alpha}|k| |j_1-j_2|(|j_1| + |j_2|)^{\alpha-2}(|j_1|^{\alpha-1}+|j_2|^{\alpha-1}+|k|^{\alpha-1}),
\end{align*}
which proves  \eqref{resumoxcqq2}.
\color{black}

\vspace{0.5\baselineskip}\noindent\textit{Proof of \eqref{resumoxcqq3}, \eqref{resumoxcqq4}.} Let
\begin{align}\label{l3dexxxx}
L_7(j,k) := 2\lambda_\alpha^\circ(j) - \lambda_\alpha^\circ(j-k) - \lambda_\alpha^\circ(j+k).
\end{align}
 From \eqref{b1deco3} and \eqref{l1ssd11whatsapp}, we have
\begin{align}\label{1px}
\mathtt{B}^N_{22}(j,k) = L_1(j,-k)L_7(j,k),
\end{align}
while it holds that
\begin{align}
& |L_1(j_1,k) - L_1(j_2,k)|\nonumber \\ &\overset{\eqref{l1ssd11whatsapp}}\le |j_1-j_2||L_2(j_1,k) + m_{1,\alpha}(j_1+k)|^2\nonumber\\
& \ +|j_2+k|\left(L_2(j_1,k) +L_2(j_2,k) + m_{1,\alpha}(j_1+k)+m_{1,\alpha}(j_2+k) \right)\nonumber\\
& \times \left(L_2(j_1,k) -L_2(j_2,k) + m_{1,\alpha}(j_1+k)-m_{1,\alpha}(j_2+k) \right)\nonumber\\
& \overset{\eqref{l1eess111s2},\eqref{l1eess111s4},\eqref{m1sess234}}{\le_\alpha} |j_1-j_2|(|j_1| + |j_2| + |k|)^{2\alpha-2}.\label{rllxxxcxcxcxc1}
\end{align}
From \eqref{1px}, we have
\begin{align}\label{B22xcxcxc1}
|\mathtt{B}^N_{22}(j,k)|{\le_\alpha}|L_1(j,-k)||L_7(j,k)|\overset{\eqref{l1eess111s1},\eqref{m1sess3}}{\le_\alpha}|k|^2|j|^{\alpha-2} (|j|^{2\alpha-1} + |k|^{2\alpha-1}),
\end{align}
which gives \eqref{resumoxcqq3}.
Towards the estimate \eqref{resumoxcqq4}, let us assume, without loss of generality, that $|j_1|\ge |j_2|$. Then we have that
\begin{align}\label{fndksh2xxx}
& |\mathtt{B}^N_{22}(j_1,k)-\mathtt{B}^N_{22}(j_2,k)| \nonumber \\
&\le |L_1(j_1,-k) - L_1(j_2,-k)||L_7(j_1,k)| + |L_1(j_2,-k)||L_7(j_1,k) - L_7(j_2,k)|\nonumber\\
&\overset{\eqref{rllxxxcxcxcxc1},\eqref{l1eess111s1},\eqref{m1sess3},\eqref{m2sess3}}{\le_\alpha}|k|^2|j_1-j_2|\nonumber\\
& \times (|j_1|^{\alpha-2}(|j_1|+|j_2|+|k|)^{2\alpha-2} +(|j_1|+|j_2|+|k|)^{2\alpha-1}(|j_1|+|j_2|)^{\alpha-3})\nonumber\\
&\le_\alpha |k|^2 |j_1-j_2| (|j_1|+|j_2|+|k|)^{2\alpha-1}((|j_1|+|j_2|)^{\alpha-3} + |j_1|^{\alpha-3}).
\end{align}
Since $j_1,j_2\ne 0$, and $|j_1|\ge |j_2|$, we have $|j_1|\ge_\alpha |j_1|+|j_2|$, therefore, \eqref{fndksh2xxx} implies \eqref{resumoxcqq4}.

\vspace{0.5\baselineskip}\noindent\textit{Proof of \eqref{resumoxcqq5}, \eqref{resumoxcqq6}.}
From \eqref{b1deco2}, \eqref{lambdadeffsd1}, and \eqref{def_of23sm}, we see that
\begin{align*}
\mathtt{B}^N_{22}(j,k) = \underbrace{j(m_{1,\alpha}(j+k) -m_{1,\alpha}(j))}_{=:L_8(j,k)}  + \underbrace{k(m_{1,\alpha}(j+k) - m_{1,\alpha}(k))}_{=:L_9(j,k)}.
\end{align*}
From \eqref{p1qosiixx1}, we have
\begin{align}
|L_8(j,k)|\le_\alpha |k||j|^{\alpha-1},\label{rlappqos1}.
\end{align}
We claim that 
\begin{align}
|L_8(j_1,k) - L_8(j_2,k)|\le_\alpha  |j_1-j_2||k| (|j_1|+|j_2|)^{\alpha-2}.\label{js2d2x2cx2cxc1}
\end{align}
Indeed, assuming $|j_1|\ge |j_2|$, without loss of generality, we have
\begin{align*}
|L_8(j_1,k) - L_8(j_2,k)|&\le_\alpha |j_1-j_2||m_{1,\alpha}(j_1+k)-m_{1,\alpha}(j_1)|\\
&\  + |j_2||m_{1,\alpha}(j_1+k)-m_{1,\alpha}(j_2+k)+m_{1,\alpha}(j_1)-m_{1,\alpha}(j_2)|\\
&\overset{\eqref{p1qosiixx1},\eqref{p1qosiixx2}}{\le_\alpha} |j_1-j_2||k|(|j_1|^{\alpha-2} + (|j_1|+|j_2|)^{\alpha-2})\\
&\le_\alpha |j_1-j_2||k| (|j_1|+|j_2|)^{\alpha-2},
\end{align*}
which gives \eqref{js2d2x2cx2cxc1}.
For $L_9$, we compute 
\begin{align}
|L_9(j,k)|&\le_\alpha |k|(|j|^{\alpha-1} + |k|^{\alpha-1}),\label{xmflqldjf1}\\
|L_9(j_1,k) - L_9(j_2,k)|&\overset{\eqref{m1sess234}}\le_\alpha |k||j_1-j_2|(|j_1|+|j_2|)^{\alpha-2}.\label{xmflqldjf2}
\end{align}
Therefore, we have
\begin{align*}
|\mathtt{B}^N_{22}(j,k)|&\le_\alpha|k|(|j|^{\alpha-1} + |k|^{\alpha-1}),\\
|\mathtt{B}^N_{22}(j_1,k)-\mathtt{B}^N_{22}(j_2,k)|&\le_\alpha  |j_1-j_2||k| (|j_1|+|j_2|)^{\alpha-2},
\end{align*}
which gives \eqref{resumoxcqq5}, and \eqref{resumoxcqq6}.

\vspace{0.5\baselineskip}\noindent\textit{Proof of \eqref{resumoxcqq7}, \eqref{resumoxcqq8}.}
From \eqref{b1deco2} and \eqref{b1deco3}, we see that
\begin{align}\label{tlqkfgk1x}
\mathtt{B}^D_{22}(j,k)=- \mathtt{B}^D_{21}(j,k)\mathtt{B}^D_{21}(j,-k).
\end{align}
Therefore, we have
\begin{align*}
|\mathtt{B}^D_{22}(j,k)|&\overset{\eqref{resumoxcqq5}}{\le_\alpha}|k|^2(|j| + |k|)^{2\alpha-2},\\
|\mathtt{B}^D_{22}(j_1,k) -\mathtt{B}^D_{22}(j_2,k)|&\overset{\substack{\eqref{resumoxcqq5} \\ \eqref{resumoxcqq6}}}{\le_\alpha}|j_1-j_2||k|^2(|j_1|+|j_2|)^{\alpha-2}(|j_1|+|j_2|+|k|)^{\alpha-1},
\end{align*}
which gives \eqref{resumoxcqq7} and \eqref{resumoxcqq8}.
\end{proof}

\begin{lemma}
For $\alpha\in (1,2)$ and $j,j_1,j_2,k\in \mathbb{Z}\backslash \left\{ 0 \right\}$ such that $|j|,|j_1|,|j_2|\ne |k|$, it holds that 
\begin{align}
|\mathtt{B}_{21}(j,k)|&\le_\alpha \frac{|k|}{\min\left\{ |j|,|k|\right\}} (|j|+|k|)^{\alpha-1},\label{rlaqkqskfkcs1}\\
|\mathtt{B}_{21}(j_1,k) -\mathtt{B}_{21}(j_2,k)|&\le_\alpha |j_1-j_2|\left(\frac{|k|}{\min\left\{ |j_1|,|j_2|,|k|\right\}}\right)^2 \nonumber \\
& \qquad \qquad \times\left(\frac{(|j_1|+|j_2|)^{\alpha-2}(|j_1|+|j_2|+|k|)^{2\alpha-2}}{(|j_1|+|k|)^{\alpha-1}(|j_2|+|k|)^{\alpha-1}} \right),\label{rlaqkqskfkcs2}\\
|\mathtt{B}_{22}(j,k)|&\le_\alpha \left(\frac{|k|}{\min\left\{ |j|,|k|\right\}} \right)^2|j|^{\alpha-2}(|j| +|k|),\label{rlaqkqskfkcs3}\\
|\mathtt{B}_{22}(j_1,k) -\mathtt{B}_{22}(j_2,k)|&\le_\alpha |j_1-j_2|\left( \frac{|k|}{\min\left\{ |j_1|,|j_2|,|k|\right\}}\right)^4 \nonumber \\
& \qquad \qquad \times \left( \frac{(|j_1|+|j_2|+|k|)^{4\alpha-3}(|j_1|+|j_2|)^{\alpha-3}}{(|j_1|+|k|)^{2\alpha-2}(|j_2|+|k|)^{2\alpha-2}}\right).\label{rlaqkqskfkcs4}
\end{align}
\end{lemma}

\begin{proof}
From Lemma~\ref{integralsd2sd}, we have
\begin{align*}
|\mathtt{B}^D_{21}(j,k)|& \ge_\alpha (|j|+|k|)^{\alpha-1}\min\left\{ |j|,|k|\right\},\\
|\mathtt{B}^D_{22}(j,k)|& \overset{\eqref{tlqkfgk1x}}{\ge_{\alpha}}(|j|+|k|)^{2\alpha-2}\min\left\{ |j|,|k|\right\}^2.
\end{align*}
Therefore, the results  \eqref{rlaqkqskfkcs1}-\eqref{rlaqkqskfkcs4} follow straightforwardly from Lemma~\ref{rotlqkftocx}.
\end{proof}

\begin{lemma}
For $\alpha\in (1,2)$ and $j\in \mathbb{N}$, we have that
\begin{align}\label{lemxcxcb3}
|\mathtt{B}_3(j)|\le_\alpha |j|^{\alpha-1}.
\end{align}
\end{lemma}
\begin{proof}
We first recall the definition of $\mathtt{B}_3$ from \eqref{b3def} and estimate the denominator/numerator separately. For the denominator,
Recalling $\lambda_{\alpha-1}$ from \eqref{def_inmu}, it follows from \eqref{explicit_multiplier} that
\[
m_{1,\alpha}^\circ(2j) - m_{1,\alpha}^\circ(j) = \mathcal{C}_\alpha (\lambda_{\alpha-1}(2j) - \lambda_{\alpha-1}(j)),
\]
where $ \mathcal{C}_\alpha$ is a positive constant given in \eqref{rkaksd2sdsdmental}. For $j\in \mathbb{N}$, we see that
\begin{align}\label{ndxcxcxcusd1}
|m_{1,\alpha}^\circ(2j) - m_{1,\alpha}^\circ(j)| =  \mathcal{C}_\alpha \int_{j}^{2j}\partial_\xi \lambda_{\alpha-1}(\xi)d\xi\overset{\eqref{sdsds11taxdsd}}{\ge_\alpha} \int_{j}^{2j}\xi^{\alpha-2}d\xi \ge_\alpha j^{\alpha-1}.
\end{align}
For the numerator, we see that  $\xi\mapsto m_{1,\alpha}(\xi)$ is a Fourier multiplier in $\mathcal{S}^{\alpha-1}$ and $m_{1,\alpha}(j)=m_{1,\alpha}^\circ(j)$ (see \eqref{def_of23sm}), therefore we have
\begin{align}\label{ndxcxcxcusd2}
|2m_{1,\alpha}^\circ(j) +m_{1,\alpha}^\circ(2j) - T_{\alpha}|\le_\alpha j^{\alpha-1}.
\end{align}
Therefore, plugging \eqref{ndxcxcxcusd1} and \eqref{ndxcxcxcusd2} into \eqref{b3def}, we obtain \eqref{lemxcxcb3}.

\end{proof}

Now, we start simplifying $m_{3,\alpha}^\circ(j,k)$ in \eqref{mu1ltip1liersss}. Since  $j\mapsto m_{3,\alpha}^\circ(j,k), k \mapsto m_{3,\alpha}^\circ(j,k) $ are even, we will focus on $j,k>0$.
\begin{lemma}
For $\alpha\in (1,2)$ and $j,k\in \mathbb{Z}\backslash \left\{ 0 \right\}$, we have that
\begin{align}
m_{3,\alpha}^\circ(j,k) & =\frac{\mathcal{C}_{\alpha}}{\alpha(\alpha+1)}\left( |j+k|^{\alpha+1} + |j-k|^{\alpha+1} - 2(|j|^{\alpha+1} +|k|^{\alpha+1})\right) \nonumber \\
& +m_{3,\alpha,r}(j,k),\label{m3deccsdx}
\end{align}
where $\mathcal{C}_{\alpha}$ is as in \eqref{Caohsdsdc}, and $m_{3,\alpha,r}(j,k)$ satisfies
\begin{align}\label{estimate1m3alax3}
|m_{3,\alpha,r}(j,k)|\le_\alpha |k|^2\langle j \rangle^{\alpha-2} + |k|^\alpha.
\end{align} 
Also, if $j_1,j_2\in \mathbb{Z}$ and $j_1,j_2\ne k$, then
\begin{align}\label{estimate1m3alax4}
|m_{3,\alpha,r}(j_1,k) - m_{3,\alpha,r}(j_2,k)|\le_\alpha (|j_1| +| j_2|)^{\alpha-2} |k|^2.
\end{align}
\end{lemma}
\begin{proof}
Since $j,k\mapsto m_{3,\alpha}^\circ(j,k)$ is even and symmetric, that is $m_{3,\alpha}^\circ(j,k)=m^\circ_{3,\alpha}(k,j)$, we can assume, without loss of generality, that $j,k\in \mathbb{N}$ and $j>k$.

For a natural number $k$, we  decompose $(2-2\cos k x)$ as 
\begin{align*}
2-2\cos kx &= (1-e^{\ii k x})(1-e^{-\ii k x}) = (1-e^{\ii x})(1-e^{-\ii x})\sum_{n_1,n_2=0}^{k-1} e^{\ii (n_1-n_2) x} \\
& = (2-2\cos x)\sum_{n_1,n_2=0}^{k-1} e^{\ii (n_1-n_2) x} \\
& =(2-2\cos x)\sum_{n_1,n_2=0}^{k-1}\cos((n_1-n_2)x).
\end{align*}
Recalling $m^\circ_{3,\alpha}$ from \eqref{mu1ltip1liersss}, we have
\begin{align*}
& m_{3,\alpha}^\circ(j,k) \\
&= \sum_{n_1,n_2=0}^{k-1}\int_{\mathbb{T}}(2-2\cos x)^{-\frac{\alpha}2}(2-2\cos j x)\cos((n_1-n_2)x)dx\\
& = \sum_{n_1,n_2=0}^{k-1}\int_{\mathbb{T}}(2-2\cos x)^{-\frac{\alpha}2}((1-\cos((j+n_1-n_2)x)) \\& 
 \qquad \qquad \qquad \qquad \qquad  + (1-\cos((j-n_1+n_2)x))   - (2-2\cos((n_1-n_2)x)))dx\\
& = \frac{1}{2}  \sum_{n_1,n_2=0}^{k-1}(m_{1,\alpha}(j+(n_1-n_2)) + m_{1,\alpha}(j-(n_1-n_2)) - 2m_{1,\alpha}(n_1-n_2))\\
& = \sum_{n_1,n_2=0}^{k-1}(m_{1,\alpha}(j+n_1-n_2) - m_{1,\alpha}(n_1-n_2)),
\end{align*}
where the second last equality follows from the definition of $m_{1,\alpha}^\circ$ in \eqref{mu1ltip1liersss} and \eqref{def_of23sm}, and the last equality follows from evenness of $m_{1,\alpha}$ and exchange $n_1$ and $n_2$ in the summation.
Now, we recall $\lambda_{\alpha-1}$ from \eqref{def_inmu} so that using explicit formula for $m_{1,\alpha}^\circ$ in \eqref{explicit_multiplier}, we have
\begin{align*}
m_{3,\alpha}^\circ(j,k) = \mathcal{C}_{\alpha}  \sum_{n_1,n_2=1}^{k-1} (\lambda_{\alpha-1} (j+ (n_1-n_2)) - \lambda_{\alpha-1} (n_1-n_2)), 
\end{align*}
In view of \eqref{jsdguasd1sd}, we write
 \begin{align}\label{m1disjxc11}
 m_{3,\alpha}^\circ(j,k) =\mathcal{C}_{\alpha} \left( f_1(j,k) + f_2(j,k)-f_2(0,k) + f_3(j,k)-f_3(0,k)\right),
 \end{align}
 where
 \begin{align}
f_1(j,k)&:= \int_0^k\int_0^k (j+s-t)^{\alpha-1}  - |s-t|^{\alpha-1}dsdt,\label{f1defx1}\\
f_2(j,k)&:=\sum_{n_1,n_2=0}^{k-1}(\lambda_{\alpha-1}(j+ n_1-n_2)  - |j+n_1-n_2|^{\alpha-1}, \label{f1defx2}\\
f_3(j,k)&:= \sum_{n_1,n_2=0}^{k-1}|j+n_1-n_2|^{\alpha-1}  -  \int_0^k\int_0^k|j+s-t|^{\alpha-1}  dsdt. \label{f1defx3}
\end{align}
$f_1$ can be explicitly computed as
\begin{align}\label{f2exxx1s111}
f_1(j,k) = \frac{1}{\alpha(\alpha+1)}\left( (j+k)^{\alpha+1} + (j-k)^{\alpha+1} - 2(j^{\alpha+1} +k^{\alpha+1})\right), \text{ for $j\ge k$}.
\end{align}
For $f_2$, it follows from \eqref{jsdguasd1sd} that (denoting $n:=n_1-n_2$),
\[
|\lambda_{\alpha-1}(j + n)- |j+n|^{\alpha-1}|\le_\alpha \langle j+n\rangle^{\alpha-2}, \text{ for $j\ne n\in \mathbb{N}$}.
\]
Therefore, we obtain
\begin{align}\label{f2exxx1s112}
|f_2(j,k)|\le_\alpha \begin{cases}
\sum_{n_1,n_2=0}^{k-1} \langle j+n_1-n_2\rangle^{\alpha-2}\le_\alpha |k|^2\langle j \rangle^{\alpha-2}, &\text{ for $j\ge k$, $j,k\in \mathbb{N}$,}\\
\sum_{n_1,n_2=0}^{k-1} \langle n_1-n_2\rangle^{\alpha-2}\le_\alpha |k|^\alpha,&\text{ for $j=0$, $k\in\mathbb{N}$.}
\end{cases}
\end{align}
For $f_3$, we have
\begin{align*}
f_3(j,k)&= \sum_{n_1,n_2=0}^{k-1} \int_0^1\int_0^1|j+n_1-n_2|^{\alpha-1} - |j+n_1-n_2+s-t|^{\alpha-1}dsdt.
\end{align*}
Using
\[
|j+n_1-n_2|^{\alpha-1} - |j+n_1-n_2+s-t|^{\alpha-1} \le_\alpha \langle j + n_1-n_2\rangle^{\alpha-2},
\]
for $j,n_1,n_2\in \mathbb{N}$, $s,t\in [0,1]$, we see that
\begin{align}\label{f2exxx1s113}
|f_3(j,k)|\le_\alpha \begin{cases}
|k|^2\langle j \rangle^{\alpha-2}, &\text{ for $j\ge k$, $j,k\in \mathbb{N}$,}\\
|k|^\alpha,&\text{ for $j=0$, $k\in\mathbb{N}$.}
\end{cases}
\end{align}
Therefore, denoting $m_{3,\alpha,r}(j,k):=\mathcal{C}_{\alpha} (f_2(j,k)-f_2(0,k) + f_3(j,k)-f_3(0,k))$, the decomposition of $m_{3,\alpha}$ as in \eqref{m3deccsdx}  follows from \eqref{f2exxx1s111} and \eqref{m1disjxc11}. The estimates \eqref{estimate1m3alax3} and \eqref{estimate1m3alax4}  follow immediately from \eqref{f2exxx1s112} and \eqref{f2exxx1s113}.
\end{proof}
\begin{lemma}\label{prosdjsdowsds1}
For $\alpha\in (1,2)$ and $j\in \mathbb{Z}$, we have that
\begin{align}\label{wksjdsdsxcsd}
W(j)=-\frac{\mathcal{C}_{\alpha}}2j |j|^{\alpha-1} + W_r(j),
\end{align}
where $W_r(j)$ satisfies
\begin{align}
|W_r(j)&|\le_\alpha |j|, \label{wxxrxxcxcxw2}\\ 
|W_r(j) - W_r(k)|&\le_\alpha |j-k|. \label{wxxrxxcxcxw3}
\end{align}
\end{lemma}
\begin{proof}
It follows from \eqref{def_of23sm},\eqref{lambdadeffsd},\eqref{lambdadeffsd1},\eqref{alpha_1aslongas} and \eqref{def_inmu} that $\lambda_\alpha(j) = \mathcal{C}_{\alpha} j \lambda_{\alpha-1}(j) - \mathcal{C}_{\alpha}\lambda_{\alpha-1}(0)j$. Therefore, we have
\begin{align}
W(j) &\overset{\eqref{defoflsksdjwwww}}= -\frac{1}{2}\lambda_\alpha(j) + \frac{T_\alpha}4j\nonumber \\
&=-\frac{\mathcal{C}_{\alpha}}{2} j\lambda_{\alpha-1}(j) + \left(\frac{\mathcal{C}_{\alpha}}{2}\lambda_{\alpha-1}(0) + \frac{T_\alpha}{4}\right)j\nonumber\\
& =-\frac{\mathcal{C}_{\alpha}}2 j|j|^{\alpha-1}+ W_r(j),\label{wrdefixchsdsd}
\end{align}
where
\begin{align}\label{wrdeisx}
W_r(j):=\frac{\mathcal{C}_{\alpha}}2j(|j|^{\alpha-1}-\lambda_{\alpha-1}(j)) + \left(\frac{\mathcal{C}_{\alpha}}{2}\lambda_{\alpha-1}(0) + \frac{T_\alpha}{4}\right)j.
\end{align}
Therefore, it suffices to show that $W_r$ in \eqref{wrdeisx} satisfies the estimates in \eqref{wxxrxxcxcxw2} and \eqref{wxxrxxcxcxw3}. These estimates follow straightforwardly from \eqref{jsdguasd1sd}.
\end{proof}

\section{Asymptotics of $\mathbb{A}$ and $\vec{D}(\xi)$ for large $\mathtt{M}$}
\subsection{Notation}
 We recall from \eqref{s00sd1sd}, \eqref{tan_site} and \eqref{mnorlam} that
 \begin{align}\label{stanxcxcx11}
 S_0^+=\left\{0< \mathtt{j}_1< \cdots < \mathtt{j}_\nu : \mathtt{j}_1,\ldots, \mathtt{j}_\nu\in \mathbb{N}\right\},&\quad  S^+=\left\{0< j_1 < \cdots < j_\nu: j_i = \mathtt{M}\mathtt{j}_i\right\},\\
  S=\left\{ \pm j : j\in S^+\right\},&\quad S_{\mathtt{M}}^\perp :=\left\{ \mathtt{M}j\in \mathbb{Z}\backslash\left\{ 0 \right\}: \mathtt{M}j\notin S\right\}.
 \end{align}
We also denote by $D_{S_0^+}$ the diagonal matrix whose $i$-th diagonal entry is $\mathtt{j}_i$, that is,
\begin{align}\label{dkxcxc111sds}
\left(D_{S_0^+}\right)^{i}_i:= \mathtt{j}_i,\quad \left(D_{S_0^+}\right)^{i}_j = 0, \text{ for $i,j=1,\ldots,\nu$ and $i\ne j$.}
\end{align}
Given a matrix $B\in \mathbb{R}^{\nu\times \nu}$, we denote
\begin{align}\label{sposx}
|B|:=\max_{i,j=1,\ldots,\nu}|B^i_j|.
\end{align}
We also denote
\begin{align}
(\vec{u})_k := 1,\text{ for $k=1,\ldots,\nu$}.\label{vecudefsd}
\end{align}
Given $S_0^+$, we denote\index{$\mathbb{A}_{S_0^+}$} by $\mathbb{A}_{S_0^+}$ be a $\nu\times \nu$ symmetric matrix defined as
\begin{align}\label{ansxs2hd1sdcs}
\left(\mathbb{A}_{S_0^+}\right)^i_k:=
\begin{cases}
 (2^{\alpha+1}-4)\mathtt{j}_i^{\alpha+3}, & \text{ if $i=k$,}\\
  2\mathtt{j}_i\mathtt{j}_k\left( (\mathtt{j}_i+\mathtt{j}_k)^{\alpha+1} + (\mathtt{j}_i-\mathtt{j}_k)^{\alpha+1} - 2(\mathtt{j}_i^{\alpha+1} +\mathtt{j}_k^{\alpha+1})\right), & \text{ if $i>k$,}\\
  (\mathbb{A}_{S_0^+})^k_i, & \text{ if $i<k$.}
 \end{cases}
\end{align}

In this section, we will study the asymptotics of the hypotheses \ref{hypos1s1}-\ref{hypothsdj22} in Section~\ref{rpoisdsd1sd} for large $\mathtt{M}$.

\subsection{Asymptotic analysis for the hypothesis \ref{hypos1s1}}
\begin{lemma}\label{Infsohsjdxcxc}
 $\mathbb{A}$ in \eqref{amplitude_modulation1} admits a decomposition
\begin{align}\label{dkhhxcxc11sdsd}
\mathbb{A} = \frac{\mathcal{C}_\alpha \pi\mathtt{M}^{\alpha+3}}{192(\alpha+1)}\left(\mathbb{A}_{S^+_0} + \mathbb{B}_1\right),
\end{align}
where $\mathbb{B}_1\in \mathbb{R}^{\nu\times \nu}$ satisfies
\begin{align}\label{bksjdexx1}
|\mathbb{B}_1|\le_{\alpha,S_0^+}\mathtt{M}^{-1}.
\end{align}
\end{lemma}
\begin{proof}
We recall the decomposition of $\mathbb{A}$ from \eqref{Amatrixexpression}:
\begin{align}\label{amsjxcxcxc2}
\mathbb{A} = \frac{\alpha\pi}{192}\mathbb{A}_1 + \frac{\alpha\pi}{192}\mathbb{A}_2.
\end{align}
We claim that 
\begin{align}\label{clamixa122sc}
|\mathbb{A}_2|\le_{\alpha,S^+_0} \mathtt{M}^{\alpha+1}.
\end{align}
In view of the definition of $\mathbb{A}_2$ in \eqref{Amatrixexpression}, we estimate $A_2(j_i,j_k)$: For $j_1,\ldots,j_\nu\in S^+$, it holds that for $i,k=1,\ldots,\nu$,
\begin{align}
|\mathtt{B}_1(j_i,j_k)|&\overset{\eqref{b1estimate11sx}}{\le_{\alpha,S_0^+}} \mathtt{M}^{\alpha-1},\quad |\mathtt{B}_{2}(j_i,j_k)|\overset{\eqref{b1deco1},\eqref{rlaqkqskfkcs1},\eqref{rlaqkqskfkcs3}}{\le_{\alpha,S_0^+}}\mathtt{M}^{\alpha-1}, \nonumber \\
|\mathtt{B}_3(j_i)| & \overset{\eqref{lemxcxcb3}}{\le_{\alpha,S_0^+}}\mathtt{M}^{\alpha-1}.\label{22xxaestima3}
\end{align}
Plugging these estimates into $A_2$ in \eqref{ejjxcxcxcsdsd11}, we have $|A_2(j_i,j_k)|\le_{\alpha,S_0^+}\mathtt{M}^{\alpha+1}$ for all $j_i,j_k\in S^+$. Therefore, plugging this into  the definition of $\mathbb{A}_2$ in \eqref{Amatrixexpression}, we obtain \eqref{clamixa122sc}.

 Now, we estimate $\mathbb{A}_1$. In view of the definition of $\mathbb{A}_1$ in \eqref{Amatrixexpression}, \eqref{ejjxcxcxcsdsd11} and the decomposition of $m_{3,\alpha}^\circ$ in \eqref{m3deccsdx}, we have
 \begin{align}
 \mathbb{A}_1&= \frac{\mathcal{C}_\alpha}{\alpha(\alpha+1)}\mathbb{A}_{11}+\mathbb{A}_{12},\text{ where $\mathcal{C}_{\alpha}=-\frac{2\pi\Gamma(1-\alpha)}{\Gamma\left(\frac{\alpha}{2}\right)\Gamma\left(1-\frac{\alpha}{2}\right)}$,}\label{a1furtherdocms}
 \end{align}
 where $\mathbb{A}_{11},\mathbb{A}_{12}$ are symmetric matrices defined by
 \begin{align}
 \left(\mathbb{A}_{11}\right)^i_k& := \begin{cases}
 (2^{\alpha+1}-4)j_i^{\alpha+3}, & \text{ if $i=k$,}\\
  2j_ij_k\left( (j_i+k_i)^{\alpha+1} + (j_i-k_i)^{\alpha+1} - 2(j_i^{\alpha+1} +k_i^{\alpha+1})\right), & \text{ if $i>k$,}\\
 \end{cases}\label{a1furtherdocms2}\\
  \left(\mathbb{A}_{12}\right)^i_k &:=
  \begin{cases}
  j_i^2m_{3,\alpha,r}(j_i,j_i), & \text{ if $i=k$},\\
  2j_ij_km_{3,\alpha,r}(j_i,j_k) & \text{ if $i >k $,}
  \end{cases}\label{a1furtherdocms3}
 \end{align}
 for $j_1,\ldots,j_\nu\in S^+$. 
From \eqref{estimate1m3alax3} and \eqref{stanxcxcx11}, we have
\begin{align}\label{a21sisdmxx1}
|\mathbb{A}_{12}|\le_{\alpha,S^+_0}\mathtt{M}^{\alpha+2}.
\end{align}
For $\mathbb{A}_{11}$, it follows from \eqref{ansxs2hd1sdcs} and \eqref{stanxcxcx11} that
\begin{align}\label{Alsjjxca1sd}
\mathbb{A}_{11}=\mathtt{M}^{\alpha+3}\mathbb{A}_{S^+_0}.
\end{align}
Thus, we have
\begin{align}
\mathbb{A} &\overset{\eqref{amsjxcxcxc2}}=\frac{\alpha\pi}{192}(\mathbb{A}_1 + \mathbb{A}_2) \overset{\eqref{a1furtherdocms}}=\frac{\mathcal{C}_\alpha \pi}{192(\alpha+1)} \left( \mathbb{A}_{11} + \frac{\alpha(\alpha+1)}{\mathcal{C}_\alpha}(\mathbb{A}_{12} +\mathbb{A}_2)\right) \nonumber \\
& \overset{\eqref{Alsjjxca1sd}}{=:}\frac{\mathcal{C}_\alpha \pi\mathtt{M}^{\alpha+3}}{192(\alpha+1)}\left(\mathbb{A}_{S^+_0} + \mathbb{B}_1\right),\label{decosgsdhsoghx1A}
\end{align}
where 
\begin{align}
\mathbb{B}_1:=\frac{\alpha(\alpha+1)}{\mathcal{C}_\alpha\mathtt{M}^{3+\alpha}}(\mathbb{A}_{12} +\mathbb{A}_2).
\end{align}
Then, the estimate \eqref{bksjdexx1} follows from \eqref{a21sisdmxx1} and \eqref{clamixa122sc}.
\end{proof}

\begin{lemma}\label{inversitjbhsksdwsz}
If $\mathbb{A}_{S_0^+}$ is invertible, then there exists $\mathtt{M}_0=\mathtt{M}(\alpha,\nu,S_0^+)$ such that if $\mathtt{M}\ge \mathtt{M}_0$, then   $\mathbb{A}$ is invertible and satisfies
\begin{align}\label{B2dsdx11sx2}
\mathbb{A}^{-1} = 
\frac{192(\alpha+1)}{\mathcal{C}_\alpha \pi\mathtt{M}^{\alpha+3}}\left(\mathbb{A}_{S^+_0}^{-1} + \mathbb{B}_2\right)
\end{align}
where $\mathbb{B}_2\in \mathbb{R}^{\nu\times \nu}$ satisfies
\begin{align}\label{B2dsdx11sx}
|\mathbb{B}_2|\le_{\alpha,\nu,S_0^+}\mathtt{M}^{-1}.
\end{align}
for some constants $c_1=c_1(\alpha,\nu),c_2=c_2(\nu,\alpha,S^+_0)>0$.
\end{lemma}
\begin{proof}
In view of \eqref{dkhhxcxc11sdsd}, we need to show that $\mathbb{A}_{S^+_0}+\mathbb{B}_1$ is invertible. Using the invertibility of $\mathbb{A}_{S_0^+}$, we can write
\[
\mathbb{A}_{S^+_0}+\mathbb{B}_1 = \mathbb{A}_{S_0^+}\left( I + \left(\mathbb{A}_{S_0^+}\right)^{-1}\mathbb{B}_1 \right),
\]
and \eqref{bksjdexx1} gives us that 
\begin{align}\label{rsdsdsdspwiggxcinters}
\left|\left(\mathbb{A}_{S_0^+}\right)^{-1}\mathbb{B}_1 \right|\le_{\alpha,\nu,S_0^+}\mathtt{M}^{-1}.
\end{align}
Therefore, for sufficiently large $\mathtt{M}$, Gershgorin circle theorem tells us that the matrix $I + \left(\mathbb{A}_{S_0^+}\right)^{-1}\mathbb{B}_1$ is invertible and we can find a matrix $\mathbb{B}_2$ such that
\[
(\mathbb{A}_{S^+_0}+\mathbb{B}_1)^{-1}=\mathbb{A}_{S^+_0}^{-1}+\mathbb{B}_2,
\]
that is,  \eqref{B2dsdx11sx2} is satisfied. The estimate in \eqref{B2dsdx11sx} follows from the estimate \eqref{rsdsdsdspwiggxcinters} with the Neumann series of $\left(I+\left(\mathbb{A}_{S_0^+}\right)^{-1}\mathbb{B}_1\right)^{-1}$.
\end{proof}

\subsection{Asymptotic analysis for the hypothesis \ref{hypothsdj2}}
We will first decompose $\vec{D}(\xi)$ depending on the dependence on $\mathtt{M}$. 
\begin{lemma}\label{sdjjjxcxcsdsdqd}
For each $\xi\in S_{\mathtt{M}}^\perp$, $\vec{D}(\xi)$ admits a decomposition,
\begin{align}\label{decompoasd}
\vec{D}(\xi) =\frac{\alpha}{16}\vec{m}_{3,1}(\xi) + \vec{m}_{4}(\xi), 
\end{align}
where 
\begin{align}
(\vec{m}_{3,1}(\xi))_k&:=m_{3,\alpha}^\circ(\xi,j_k)j_k,\text{ for $j_k\in S^+$ and  for $k=1,\ldots,\nu$,}\label{m3dscsxcsw}
\end{align}
and 
$\vec{m}_{4}(\xi)$ satisfies
\begin{align}
|\vec{m}_{4}(\xi)|&\le_{\alpha,S^+_0,\nu} \mathtt{M}|\xi|^{\alpha-1},\label{mrmsidx1}\\
|\vec{m}_{4}(\xi_1) - \vec{m}_4(\xi_2)|&\le_{\alpha,S^+_0,\nu}\mathtt{M}|\xi_1-\xi_2|(|\xi_1|+|\xi_2|)^{\alpha-2},\text{ for $\xi_1,\xi_2\in S_{\mathtt{M}}^\perp$.} \label{mrmsidx2}
\end{align}
\end{lemma}
\begin{proof}
 Given $\xi\in S_{\mathtt{M}}^\perp$, let us denote by $i_\xi$ the integer such that 
 \begin{align}\label{xidsd2sdsdw1sd}
 \mathtt{M}i_\xi =\xi.
 \end{align}
 In view of \eqref{desstimaxcxcsd}, we define
 \begin{align}
 (\vec{m}_{4}(\xi))_k&:= - \frac{\alpha}{16}j_k\mathtt{B}_1(\xi,j_k) + \frac{\alpha^2}{32}j_k\mathtt{B}_2(\xi,j_k),\label{m3dscsxcsw2}
 \end{align}
 for $j_1,\ldots,j_k\in S^+$, so that we have the decomposition as in \eqref{decompoasd}.  Now, it suffices to prove the estimates \eqref{mrmsidx1} and \eqref{mrmsidx2}.
 
   We first prove \eqref{mrmsidx1}. Indeed, we have
   \begin{align}
   |\mathtt{B}_1(\xi,j_k)|&\overset{\eqref{b1estimate11sx}}{\le_\alpha} (|\xi|^{\alpha-1} + |j_k|^{\alpha-1})\overset{\eqref{xidsd2sdsdw1sd},\eqref{stanxcxcx11}}=\mathtt{M}^{\alpha-1}(|i_\xi|^{\alpha-1}+|\mathtt{j}_k|^{\alpha-1})\nonumber \\
   & \le_{\alpha,S_0^+}\mathtt{M}^{\alpha-1}|i_\xi|^{\alpha-1}\overset{\eqref{xidsd2sdsdw1sd}}=|\xi|^{\alpha-1}.\label{rxkkxsxdsxxcxc1}
   \end{align}
   Similarly, we have
   \begin{align*}
   |\mathtt{B}_{21}(\xi,j_k)|&\overset{\eqref{rlaqkqskfkcs1}}{\le_\alpha}|\xi|^{\alpha-1}+|j_k|^{\alpha-1}\le_{\alpha,S^+_0}|\xi|^{\alpha-1},\\
   |\mathtt{B}_{22}(\xi,j_k)|&\overset{\eqref{rlaqkqskfkcs3}}{\le_\alpha}|\xi|^{\alpha-2}(|\xi| +|j_k|)\le_{\alpha,S^+_0}|\xi|^{\alpha-1},
   \end{align*}
   which implies
   \begin{align}\label{rxkkxsxdsxxcxc2}
   |\mathtt{B}_{2}(\xi,j_k)|\overset{\eqref{b1deco1}}{\le_{\alpha,S^+_0}}|\xi|^{\alpha-1}.
   \end{align}
  Plugging \eqref{rxkkxsxdsxxcxc1} and \eqref{rxkkxsxdsxxcxc2} into \eqref{m3dscsxcsw2}, we obtain \eqref{mrmsidx1}.
  
  Towards \eqref{mrmsidx2}, we denote by $i_{\xi,1},i_{\xi_2}$ the integers such that
  \begin{align}\label{xjasd1xs1}
  \mathtt{M}i_{\xi,1}=\xi_1,\quad \mathtt{M}i_{\xi,2}=\xi_2.
  \end{align}
We have that
  \begin{align}
  |\mathtt{B}_{1}(\xi_1,j_k) - \mathtt{B}_{1}(\xi_2,j_k)|&\overset{\eqref{bestimate11sx1}}{\le_\alpha}|\xi_1-\xi_2|(|\xi_1|+|\xi_2|)^{\alpha-2}.\label{bhsjd111sdxcxc}
  \end{align}
  For $\mathtt{B}_{21}$, we consider to cases: 1) $|i_{\xi,1}|\ge 2|i_{\xi,2}|$ or $|i_{\xi,2}|\ge 2|i_{\xi,1}|$ and $2) \frac{1}{2}|i_{\xi,2}|\le |i_{\xi,1}|\le 2|i_{\xi,2}|$. In the first case, assuming $|i_{\xi,1}|\ge 2|i_{\xi,2}|$, without loss of generality, we have
  \begin{align*}
  |\mathtt{B}_{21}(\xi_1,j_k) - \mathtt{B}_{21}(\xi_2,j_k)|&\le |\mathtt{B}_{21}(\xi_1,j_k)| + |\mathtt{B}_{21}(\xi_2,j_k)|\\
  &\overset{\eqref{rlaqkqskfkcs1}}{\le_{\alpha,S_0^+}} |\xi_1|^{\alpha-1} + |\xi_2|^{\alpha-1}\le_{\alpha,S^+_0} |\xi_1-\xi_2|(|\xi_1| + |\xi_2|)^{\alpha-2},
  \end{align*}
  where the last inequality follows from that $|i_{\xi,1}| \ge 2|i_{\xi,2}|$ implies $|\xi_1|\ge 2|\xi_2|$. For the second case, assuming $\frac{1}{2}|i_{\xi,2}|\le |i_{\xi,1}|\le 2|i_{\xi,2}|$, we have
  
  \begin{align*}
  & |\mathtt{B}_{21}(\xi_1,j_k) - \mathtt{B}_{21}(\xi_2,j_k)|\\
  &\overset{\eqref{rlaqkqskfkcs2},\eqref{xjasd1xs1},\eqref{stanxcxcx11}}{\le_{\alpha}} \mathtt{M}^{\alpha-2}|\xi_1-\xi_2|\left(\frac{(|i_{\xi,1}|+|i_{\xi,2}|)^{\alpha-2}(|i_{\xi,1}| + |i_{\xi,2}| + |\mathtt{j}_k|)^{2\alpha-2}}{(|i_{\xi,1}| + |\mathtt{j}_k|)^{\alpha-1}(|i_{\xi,2}| + |\mathtt{j}_k|)^{\alpha-1}} \right)\\
  &\qquad \le_{\alpha,S^+_0}\mathtt{M}^{\alpha-2}|\xi_1-\xi_2|(|i_{\xi,1}| + |i_{\xi,2}|)^{\alpha-2}\\
  &\qquad \overset{\eqref{xjasd1xs1}}= |\xi_1-\xi_2|(|\xi_1| + |\xi_2|)^{\alpha-2}.
  \end{align*}
  Therefore, in either case, we have
  \begin{align}\label{b221jjsdsxcxc}
    |\mathtt{B}_{21}(\xi_1,j_k) - \mathtt{B}_{21}(\xi_2,j_k)|\le_{\alpha,S^+_0}|\xi_1-\xi_2|(|\xi_1| + |\xi_2|)^{\alpha-2}.
  \end{align}
  Similarly, if $|i_{\xi,1}\ge 2|i_{\xi,2}|$ or $|i_{\xi_2}|\ge 2|i_{\xi,1}|$, we have
  \begin{align*}
    |\mathtt{B}_{22}(\xi_1,j_k) - \mathtt{B}_{22}(\xi_2,j_k)| & \le   |\mathtt{B}_{22}(\xi_1,j_k)| + |\mathtt{B}_{22}(\xi_2,j_k)| \\
    & \overset{\eqref{rlaqkqskfkcs3}}\le_{\alpha}|\xi_1|^{\alpha-1}+|\xi_2|^{\alpha-1}\le_{\alpha}|\xi_1-\xi_2|(|\xi_1| + |\xi_2|)^{\alpha-2},
  \end{align*}
  and if $\frac{1}{2}|i_{\xi,2}|\le |i_{\xi,1}|\le 2|i_{\xi,2}|$, we have
  \begin{align*}
    |\mathtt{B}_{22}(\xi_1,j_k) - \mathtt{B}_{22}(\xi_2,j_k)|& \overset{\eqref{rlaqkqskfkcs4}}{\le_{\alpha,S^+_0}}\mathtt{M}^{\alpha-2}|\xi_1-\xi_2|(|i_{\xi,1}|+|i_{\xi,2}|)^{\alpha-2} \\
    & \le_\alpha |\xi_1-\xi_2|(|\xi_1| + |\xi_2|)^{\alpha-2}.
  \end{align*}
  Therefore, in either case, we have $ |\mathtt{B}_{22}(\xi_1,j_k) - \mathtt{B}_{22}(\xi_2,j_k)|\le_{\alpha,S^+_0}|\xi_1-\xi_2|(|\xi_1| + |\xi_2|)^{\alpha-2}$. Plugging this and \eqref{b221jjsdsxcxc} into \eqref{b1deco1}, we have
  \begin{align}\label{kksdsxcx11s2sd2}
  |\mathtt{B}_{2}(\xi_1,j_k) - \mathtt{B}_{2}(\xi_2,j_k)|\le_{\alpha,S^+_0}|\xi_1-\xi_2|(|\xi_1| + |\xi_2|)^{\alpha-2}.
  \end{align}
  Then, the estimate \eqref{mrmsidx2} follows straightforwardly from \eqref{m3dscsxcsw2}, \eqref{bhsjd111sdxcxc} and \eqref{kksdsxcx11s2sd2}.
\end{proof}

In view of \eqref{m3deccsdx}, we can further decompose $\vec{m}_{3,1}$ for sufficiently large $\xi$:
\begin{lemma}\label{xisdesxcsmm3}
For $\xi\in S_{\mathtt{M}}^\perp$ such that $|\xi|\ge \sup\left\{j_k: j_k\in S^+ \right\}$, $\vec{m}_{3,1}(\xi)$ admits a decomposition,
\begin{align}\label{m3sdjsdjxcxcsds}
\vec{m}_{3,1}(\xi) = \frac{\mathcal{C}_\alpha}{\alpha(\alpha+1)}\vec{m}_{3,2}(\xi) + \vec{m}_{3,3}(\xi),
\end{align}
where
\begin{align}\label{m032sdsdedosjvsd}
(\vec{m}_{3,2}(\xi))_k = (|\xi + j_k|^{\alpha+1} + |\xi - j_k|^{\alpha+1} - 2(|\xi|^{\alpha+1} + |j_k|^{\alpha+1}))j_k,
\end{align}
and $\vec{m}_{3,3}(\xi)$ satisfies
\begin{align}
|(\vec{m}_{3,3}(\xi))_k|&\le_{\alpha,S^+_0} \mathtt{M}^2|\xi|^{\alpha-1}, \label{sdjjwkkksdxcxsd}\\
|(\vec{m}_{3,3}(\xi_1))_k-(\vec{m}_{3,3}(\xi_2))_k|&\le_{\alpha,S^+_0}\mathtt{M}^2|\xi_1-\xi_2|(|\xi_1|+|\xi_2|)^{\alpha-2} . \label{sdjjwkkksdxcxsd1}
\end{align}
\end{lemma}
\begin{proof}
For $\xi\in S^\perp$ such that $|\xi|\ge \sup\left\{j_k: j_k\in S^+ \right\}$, we define (using that $m_{3,\alpha}^\circ(\xi,k)=m_{3,\alpha}^\circ(-\xi,k)$, which follows from its definition in \eqref{mu1ltip1liersss}),
\begin{align}\label{sdjjxcxc11}
(\vec{m}_{3,3}(\xi))_k:=m_{3,\alpha,r}(\xi,j_k)j_k,
\end{align}
so that \eqref{m3dscsxcsw}, \eqref{m032sdsdedosjvsd} and \eqref{m3deccsdx} yield the decomposition \eqref{m3sdjsdjxcxcsds}. Therefore, it suffices to show that  $\vec{m}_{3,3}$ in \eqref{sdjjxcxc11}  satisfies the estimates in \eqref{sdjjwkkksdxcxsd} and \eqref{sdjjwkkksdxcxsd1}. 
Since $\xi\mapsto m_{3,\alpha}^\circ(\xi,k),\vec{m}_{3,2}(\xi)$ are even, we only need to show the estimate for $\xi>0$. For positive $\xi,\xi_1,\xi_2 \in S^\perp$ such that $\xi,\xi_1,\xi_2\ge  \sup\left\{|j|: j\in S \right\}$, we have 
\begin{align*}
|m_{3,\alpha,r}(\xi,j_k)|&\overset{\eqref{estimate1m3alax3}}{\le_\alpha}|j_k|^2|\xi|^{\alpha-2}  + |j_k|^\alpha  \\
 & \overset{\eqref{xidsd2sdsdw1sd},\eqref{stanxcxcx11}}{\le_{\alpha,S^+_0}}\mathtt{M}^\alpha |i_{\xi}|^{\alpha-2} =\mathtt{M}^2 |\xi|^{\alpha-2}\le \mathtt{M}|\xi|^{\alpha-1},\\
|m_{3,\alpha,r}(\xi_1,j_k)-m_{3,\alpha,r}(\xi_2,j_k)|&\overset{\eqref{estimate1m3alax4}}{\le_{\alpha,S^+_0}}\mathtt{M}^2|\xi_1+\xi_2|^{\alpha-2} \\
 & = \mathtt{M}\frac{\mathtt{M}}{|\xi_1-\xi_2|}|\xi_1-\xi_2|(\xi_1+\xi_2)^{\alpha-2}\\
  & \le \mathtt{M}|\xi_1-\xi_2|(\xi_1+\xi_2)^{\alpha-2},
\end{align*}
where the two  last equalities in the estimates follow from that $\xi,\xi_1,\xi_2$ are $\mathtt{M}$-multiples of natural numbers. Therefore, plugging these estimates into \eqref{sdjjxcxc11}, we obtain the estimates \eqref{sdjjwkkksdxcxsd} and \eqref{sdjjwkkksdxcxsd1}.
\end{proof}

\begin{lemma}\label{jsdsdxkajw11}
For $\xi\in S_{\mathtt{M}}^\perp$ such that $|\xi|\ge 2\sup\left\{j_k: j_k\in S^+ \right\}$, $\vec{m}_{3,2}(\xi)$ admits a decomposition,
\begin{align}\label{tlqkfrowhwrkxcxza}
(\vec{m}_{3,2}(\xi))_k =  \alpha(\alpha+1)j_k^3|\xi|^{\alpha-1}+ (\vec{m}_{3,r}(\xi))_k,
\end{align}
where $(\vec{m}_{3,r}(\xi))_k$ satisfies
\begin{align}
|(\vec{m}_{3,r}(\xi))_k|&\le_\alpha j_k^4|\xi|^{\alpha-2}, \label{rlqsjsdsxc1}\\
| (\vec{m}_{3,r}(\xi_1))_k - (\vec{m}_{3,r}(\xi_2))_k|&\le_\alpha j_k^4|\xi_1-\xi_2|(|\xi_1|+|\xi_2|)^{\alpha-3}, \label{rlqsjsdsxc2}
\end{align}  for $|\xi_1|,|\xi_2| \ge 2\sup\left\{|j_k|: j_k\in S \right\}$.
\end{lemma}
\begin{proof}
In view of \eqref{m032sdsdedosjvsd}, we rewrite $j_k(|\xi + j_k|^{\alpha+1} + |\xi - j_k|^{\alpha+1} - 2(|\xi|^{\alpha+1} + |j_k|^{\alpha+1}))$ as
\begin{align}\label{def_fsjsd1mmees}
& j_k(|\xi + j_k|^{\alpha+1} + |\xi - j_k|^{\alpha+1} - 2(|\xi|^{\alpha+1} + |j_k|^{\alpha+1}))\nonumber\\
& = \alpha(\alpha+1)j_k^3|\xi|^{\alpha-1} + \alpha(\alpha+1)(\alpha-1)j_k|\xi|^{\alpha+1}\int_0^{\frac{j_k}{|\xi|}}\int_{-t}^{t}\int_0^s (1+u)^{\alpha-2}dudsdt\nonumber\\
&=:\alpha(\alpha+1)j_k^3|\xi|^{\alpha-1} + (\vec{m}_{3,r}(\xi))_k,
\end{align}
which gives us the decomposition in \eqref{tlqkfrowhwrkxcxza}. It suffices to check the estimates \eqref{rlqsjsdsxc1} and \eqref{rlqsjsdsxc2}.
For \eqref{rlqsjsdsxc1}, we have that (since $|\xi|\ge 2\sup\left\{|j_k|: j_k\in S \right\}$, the integrand is strictly bounded),
\begin{align}\label{skksdxcxc11}
\left|\int_0^{\frac{j_k}{|\xi|}}\int_{-t}^{t}\int_0^s (1+u)^{\alpha-2}dudsdt \right|\le_\alpha \int_0^{\frac{j_k}{|\xi|}}\int_{-t}^t\int_0^s 1 dudsdt =\frac{1}{3}j_k^3|\xi|^{-3}.
\end{align}
Plugging this into the definition of $(\vec{m}_{3,r}(\xi))_k$ in \eqref{def_fsjsd1mmees}, we obtain \eqref{rlqsjsdsxc1}.
To show \eqref{rlqsjsdsxc2}, we first  assume that $\xi_1$ and $\xi_2$ have opposite signs. Then, we have
\begin{align}\label{jjsdxcx111xc}
| (\vec{m}_{3,r}(\xi_1))_k - (\vec{m}_{3,r}(\xi_2))_k|\overset{\eqref{rlqsjsdsxc1}}{\le_\alpha} j_k^4(|\xi_1|^{\alpha-2} + |\xi_2|^{\alpha-2}),
\end{align}
while it holds that
\begin{align}\label{rjjjsdszcxc1}
|\xi_1-\xi_2|(|\xi_1|+|\xi_2|)^{\alpha-3} \ge_{\alpha} (|\xi_1| + |\xi_2|)(|\xi_1| + |\xi_2|)^{\alpha-3} \ge_\alpha(|\xi_1| + |\xi_2|)^{\alpha-2},
\end{align}
 therefore we obtain \eqref{rlqsjsdsxc2}.
Now, we assume that $\xi_1,\xi_2>0$. Without loss of generality, we can assume that $\xi_1>\xi_2>0$. If $\xi_1 \ge 2\xi_2$, we still have  \eqref{jjsdxcx111xc}, and the inequalities in \eqref{rjjjsdszcxc1} hold true as well, therefore we obtain \eqref{rlqsjsdsxc2}. Lastly, we assume 
\begin{align}\label{auusdsxcxc1}
0<\xi_2<\xi_1<2\xi_2.
\end{align} Then we have
  \begin{align*}
&|(\vec{m}_{3,r}(\xi_1))_k - (\vec{m}_{3,r}(\xi_2))_k|\\
&\qquad \overset{\eqref{def_fsjsd1mmees}}{\le_\alpha}j_k(\xi_1^{\alpha+1}-\xi_2^{\alpha+1})\int_0^{\frac{j_k}{|\xi_2|}}\int_{-t}^{t}\int_0^s (1+u)^{\alpha-2}dudsdt\\
& \qquad \ + j_k\xi_2^{\alpha+1}\int_{\frac{j_k}{|\xi_1|}}^{\frac{j_k}{|\xi_2|}}\int_{-t}^{t}\int_0^s (1+u)^{\alpha-2}dudsdt\\
& \qquad \overset{\eqref{skksdxcxc11}}{\le_\alpha} j_k |\xi_1-\xi_2| (|\xi_1|^{\alpha} + |\xi_2|^{\alpha})j_k^3|\xi_2|^{-3} +j_k\xi_2^{\alpha+1}\int_{\frac{j_k}{|\xi_1|}}^{\frac{j_k}{|\xi_2|}} t^2 dt\\
&\qquad \overset{\eqref{auusdsxcxc1}}{\le_\alpha} j_k^4|\xi_1-\xi_2||\xi_1+\xi_2|^{\alpha-3} + j_k^4\xi_2^{\alpha+1}\left(\frac{1}{\xi_2^3} - \frac{1}{\xi_1^3}\right)\\
&\qquad \le_\alpha j_k^4|\xi_1-\xi_2||\xi_1+\xi_2|^{\alpha-3},
\end{align*}
which gives \eqref{rlqsjsdsxc2}.
\end{proof}

In view of  \eqref{poisosuduwoqs}, let us denote 
\begin{align}\label{Eajhxcxc1}
\mathtt{A}(\xi):= W(\xi) - \frac{\pi}6\xi\vec{D}(\xi)\cdot\mathbb{A}^{-1}\overline{\omega},\text{ for $\xi\in S_{\mathtt{M}}^\perp\cup\left\{ 0\right\}$}.
\end{align}
We  simplify the expression of $\mathtt{A}$:
\begin{lemma}\label{kkkrpsosdsd12p2p2p}
Suppose $\mathbb{A}_{S_0^+}$ in \eqref{ansxs2hd1sdcs} is invertible. For $\xi\in S_{\mathtt{M}}^\perp$ such that $|\xi|\ge 2 \max\left\{ |j_k|: j_k\in S\right\}$, and for sufficiently large $\mathtt{M}$, depending on $\alpha,\nu,S_0^+$, we have
\begin{align}\label{mathasexxsd}
\mathtt{A}(\xi) = -\frac{\mathcal{C}_\alpha}{2}\left( 1-2\alpha(\alpha+1)I_{S_0^+}\right)\xi|\xi|^{\alpha-1} + m_{\mathtt{A},r}(\xi),
\end{align}
where
\begin{align}I_{S_0^+} &:=\vec{u}\cdot\left(\left(D_{S_0^+}\right)^{-\alpha} \cdot \mathbb{A}_{S^+_0} \left(D_{S^+_0}\right)^{-3}\right)^{-1}\vec{u}, \quad (\vec{u})_k = 1, \text{ for $k=1,\ldots, \nu$},\label{Idjsdxc1}
\end{align}
and $m_{\mathtt{A},r}(\xi)$ satisfies
\begin{align}
|m_{\mathtt{A},r}(\xi)|&{\le_{\alpha,\nu,S_0^+}} \mathtt{M}^{1-\alpha}|\xi|^{\alpha} + \mathtt{M}|\xi|^{\alpha-1},\label{skkksdhsd21}\\
|m_{\mathtt{A},r}(\xi_1) - m_{\mathtt{A},r}(\xi_2)|&{\le_{\alpha,\nu,S_0^+}}|\xi_1-\xi_2|(\mathtt{M}^{1-\alpha}(|\xi_1|+|\xi_2|)^{\alpha-1} + \mathtt{M}(|\xi_1|+|\xi_2|)^{\alpha-2}). \label{skkksdhsd22}
\end{align}
\end{lemma}
\begin{proof}
Recalling $W$ from \eqref{defoflsksdjwwww}, we have
\begin{align}\label{rjskscjxcksdsq}
-\frac{1}{2}\lambda_\alpha(\xi) + \frac{T_\alpha}4 \xi = W(\xi) \overset{\eqref{wksjdsdsxcsd}}= -\frac{\mathcal{C}_\alpha}2\xi|\xi|^{\alpha-1} +W_r(\xi).
\end{align}

  We choose a sufficiently  large $\mathtt{M}$ so that Lemma~\ref{inversitjbhsksdwsz} gives us that $\mathbb{A}$ is invertible. In this case, it follows from \eqref{B2dsdx11sx2}  that
\begin{align}\label{Absdsdxc}
\mathbb{A}^{-1} = \frac{192(\alpha+1)}{\mathcal{C}_\alpha \pi\mathtt{M}^{\alpha+3}}\left(\mathbb{A}_{S^+_0}^{-1} + \mathbb{B}_2\right),
\end{align}
where 
\begin{align}\label{bk2jsdsdad}
|\mathbb{B}_2| \overset{\eqref{B2dsdx11sx}}\le_{\alpha,\nu,S_0^+}\mathtt{M}^{-1}.
\end{align}
Also, we recall $\overline{\omega}$ from \eqref{linear_frequency_aa} so that we can write $
(\overline{\omega})_k\overset{\eqref{defoflsksdjwwww}}=W(j_k),$ for $j_1,\ldots,j_\nu\in S^+$. In view of \eqref{wksjdsdsxcsd}, we write $\overline{\omega}$ as (using the notations in \eqref{stanxcxcx11} and \eqref{dkxcxc111sds})
\begin{align}
\overline{\omega}&=-\frac{\mathcal{C}_\alpha}{2}\mathtt{M}^\alpha \left(\left(D_{S^+_0}\right)^\alpha\vec{u} + \vec{w}\right),\label{delfjscxcxcs}
\end{align}
where
\begin{align}
(\vec{u})_k=1, \text{ for $k=1,\ldots,\nu$},\quad (\vec{w})_k=-\frac{2 \mathtt{M}^{-\alpha}}{C_\alpha} W_r(j_k),\text{ for $j_k\in S^+$.} \label{delfjscxcxcs2}
\end{align}
For the size of $\vec{w}$, we have
\begin{align}\label{shhxcxca11}
|\vec{w}|\overset{\eqref{wxxrxxcxcxw2}}\le_{\alpha} \mathtt{M}^{-\alpha}\sum_{k=1}^{\nu}|j_k|\overset{\eqref{stanxcxcx11}}\le_{\alpha,S_0^+} \mathtt{M}^{1-\alpha}.
\end{align}
Then, we compute
\begin{align}\label{rksjxck11p0s2}
\mathbb{A}^{-1}\overline{\omega}& \overset{\eqref{Absdsdxc},\eqref{delfjscxcxcs}}= -\frac{96(\alpha+1)}{\pi\mathtt{M}^{3}}\left(\mathbb{A}_{S^+_0}^{-1} + \mathbb{B}_2\right)\left(\left(D_{S^+_0}\right)^\alpha\vec{u} + \vec{w}\right)\nonumber\\
&= -\frac{96(\alpha+1)}{\pi\mathtt{M}^{3}}\left(\mathbb{A}_{S^+_0}^{-1} \left(D_{S^+_0}\right)^\alpha\vec{u} + \vec{v} \right),
\end{align}
where 
\begin{align}
\vec{v}:=\left(\mathbb{A}_{S_0^+}\right)^{-1}\vec{w} + \mathbb{B}_2\left( D_{S^+_0}\right)^{\alpha}\vec{u} + \mathbb{B}_2\vec{w}.
\end{align}
For $\vec{v}$, we estimate its size as,
\begin{align}
|\vec{v}| &\overset{\eqref{bk2jsdsdad},\eqref{shhxcxca11}}\le_{\alpha,\nu,S_0^+} \mathtt{M}^{1-\alpha}.\label{vsiz2s}
\end{align}

For $\vec{D}(\xi)$, we simplify it only for $|\xi|\ge 2 \max\left\{ j_k: j_k\in S^+\right\}$, since in the other case, its expression is not necessary. Assuming $|\xi|\ge 2 \max\left\{ j_k: j_k\in S^+\right\}$, we have
\begin{align}\label{tlqkfhshxc22}
\vec{D}(\xi)&\overset{\eqref{decompoasd}}=\frac{\alpha}{16}(\vec{m}_{3,1}(\xi))_k +  (\vec{m}_{4}(\xi))_k\nonumber\\
&\overset{\eqref{m3sdjsdjxcxcsds}}= \frac{\mathcal{C}_\alpha}{16(\alpha+1)}(\vec{m}_{3,2}(\xi))_k + \frac{\alpha}{16}(\vec{m}_{3,3}(\xi))_k +(\vec{m}_{4}(\xi))_k\nonumber\\
&\overset{\eqref{tlqkfrowhwrkxcxza}}= \frac{\alpha  \mathcal{C}_\alpha}{16}j_k^3|\xi|^{\alpha-1} +\frac{\mathcal{C}_\alpha}{16(\alpha+1)}(\vec{m}_{3,r}(\xi))_k+ \frac{\alpha}{16}(\vec{m}_{3,3}(\xi))_k +(\vec{m}_{4}(\xi))_k,\nonumber\\
&=: \frac{\alpha  \mathcal{C}_\alpha}{16}\Bigg(j_k^3|\xi|^{\alpha-1} \nonumber  \\
 & \left.\qquad \qquad +\left(\frac{\alpha  \mathcal{C}_\alpha}{16}\right)^{-1}\left(\frac{\mathcal{C}_\alpha}{16(\alpha+1)}(\vec{m}_{3,r}(\xi))_k+ \frac{\alpha}{16}(\vec{m}_{3,3}(\xi))_k +(\vec{m}_{4}(\xi))_k\right)\right).
\end{align}
Therefore, using the notation for $D_{S_0^+}$ and $\vec{u}$ in \eqref{dkxcxc111sds} and \eqref{vecudefsd}, we can write
\begin{align}\label{dkksd1sdsd}
\vec{D}(\xi)=\frac{\alpha \mathcal{C}_\alpha\mathtt{M}^3}{16}\left(|\xi|^{\alpha-1}(D_{S_0^+})^{3}\vec{u} + \vec{m}_5(\xi)\right),\text{ for $|\xi|\ge 2 \max\left\{ j_k: j_k\in S^+\right\}$},
\end{align}
where
\[
(\vec{m}_5(\xi))_k:=\frac{1}{\mathtt{M}^3}\left(\frac{\alpha  \mathcal{C}_\alpha}{16}\right)^{-1}\left(\frac{\mathcal{C}_\alpha}{16(\alpha+1)}(\vec{m}_{3,r}(\xi))_k+ \frac{\alpha}{16}(\vec{m}_{3,3}(\xi))_k +(\vec{m}_{4}(\xi))_k\right)
\]
For the estimates of $\vec{m}_5$, we have that
\begin{align}\label{rjjjsdmu1}
|\vec{m}_5(\xi)|&\le_\alpha \mathtt{M}^{-3}(|\vec{m}_{3,r}(\xi)| + |\vec{m}_{3,3}(\xi)| + |\vec{m}_4(\xi)|)\nonumber\\
& \overset{\eqref{rlqsjsdsxc1}}{\le_{\alpha,\nu,S_0^+}} \mathtt{M}|\xi|^{\alpha-2}+\mathtt{M}^{-3}|\vec{m}_{3,3}(\xi)| + \mathtt{M}^{-3}|\vec{m}_4(\xi)|\nonumber\\
&\overset{\eqref{sdjjwkkksdxcxsd}}{\le_{\alpha,S_0^+,\nu}}\mathtt{M}|\xi|^{\alpha-2} + \mathtt{M}^{-1}|\xi|^{\alpha-1} +\mathtt{M}^{-3} |\vec{m}_4(\xi)|\nonumber\\
&\overset{\eqref{mrmsidx1}}{\le_{\alpha,S_0^+,\nu}}\mathtt{M}|\xi|^{\alpha-2} + \mathtt{M}^{-1}|\xi|^{\alpha-1} .
\end{align}
For $\xi_1,\xi_2\in S_{\mathtt{M}}^\perp$ such that $|\xi_1,\xi_2|\ge2 \max\left\{ j_k: j_k\in S^+\right\}$, we also have
\begin{align}
& |\vec{m}_5(\xi_1)-\vec{m}_5(\xi_2)| \nonumber\\
&\le_\alpha \mathtt{M}^{-3}\left(|\vec{m}_{3,r}(\xi_1)-\vec{m}_{3,r}(\xi_2)| + |\vec{m}_{3,3}(\xi_1)-\vec{m}_{3,3}(\xi_2)| + |\vec{m}_4(\xi_1)-\vec{m}_4(\xi_2)|\right)\nonumber\\
&\overset{\eqref{rlqsjsdsxc2},\eqref{sdjjwkkksdxcxsd1},\eqref{mrmsidx2}}{\le_{\alpha,\nu,S_0^+}}\mathtt{M}|\xi_1-\xi_2|(|\xi_1|+|\xi_2|)^{\alpha-3} +\mathtt{M}^{-1}|\xi_1-\xi_2|(|\xi_1|+|\xi_2|)^{\alpha-2}.\label{rjjjsdmu2}
\end{align}

We combine \eqref{rksjxck11p0s2} and \eqref{dkksd1sdsd} and obtain (for $|\xi|\ge 2 \max\left\{ j_k: j_k\in S^+\right\}$),
\begin{align}\label{jjjsdad1hsdsd}
\frac{\pi}6\xi\vec{D}(\xi)\cdot\mathbb{A}^{-1}\overline{\omega} &= - {\alpha(\alpha+1) \mathcal{C}_\alpha}\xi \left(|\xi|^{\alpha-1}(D_{S_0^+})^{3}\vec{u} + \vec{m}_5(\xi)\right)\cdot \left(\mathbb{A}_{S^+_0}^{-1} \left(D_{S^+_0}\right)^\alpha\vec{u} + \vec{v} \right)\nonumber\\
&=:-\alpha(\alpha+1) \mathcal{C}_\alpha  I_{S_0^+} \xi|\xi|^{\alpha-1} + m_{6}(\xi),
\end{align}
where
\begin{align}
 I_{S_0^+} &:=(D_{S_0^+})^{3}\vec{u} \cdot \mathbb{A}_{S^+_0}^{-1} \left(D_{S^+_0}\right)^\alpha\vec{u},\label{Isjxcs11}\\
 m_6(\xi)&:=- {\alpha(\alpha+1) \mathcal{C}_\alpha}\left( \xi |\xi|^{\alpha-1}(D_{S_0^+})^{3}\vec{u} \cdot \vec{v} +  \xi \vec{m}_5(\xi)\cdot\left(\mathbb{A}_{S^+_0}^{-1} \left(D_{S^+_0}\right)^\alpha\vec{u} + \vec{v} \right)\right).\label{kksjsxcxjcs1}
\end{align}
For the size of $m_6$, we have
\begin{align}
|m_6(\xi)|&\overset{\eqref{rjjjsdmu1},\eqref{vsiz2s}}{\le_{\alpha,\nu,S_0^+}} \mathtt{M}^{1-\alpha}|\xi|^{\alpha},\label{skkksdhsd1}\\
|m_6(\xi_1) - m_6(\xi_2)|&\overset{\eqref{vsiz2s},\eqref{rjjjsdmu2}}{\le_{\alpha,\nu,S_0^+}}|\xi_1-\xi_2|(\mathtt{M}^{1-\alpha}(|\xi_1|+|\xi_2|)^{\alpha-1} + \mathtt{M}(|\xi_1|+|\xi_2|)^{\alpha-2}). \label{skkksdhsd2}
\end{align}
Combining \eqref{jjjsdad1hsdsd} with \eqref{rjskscjxcksdsq} and \eqref{Eajhxcxc1}, we obtain \eqref{mathasexxsd}, with
\[
m_{\mathtt{A},r}(\xi):= W_r(\xi) - m_6(\xi).
\]
Then the estimates \eqref{skkksdhsd21} and \eqref{skkksdhsd22} follow from \eqref{skkksdhsd1}, \eqref{skkksdhsd2}, \eqref{wxxrxxcxcxw2} and \eqref{wxxrxxcxcxw3}.
\end{proof}

\subsection{Asymptotic analysis for the hypothesis \ref{hypothsdj22}}
We  derive a more explicit expression for the matrices $\mathbb{B}_{\xi_1,\xi_2}$ and $\mathbb{C}_{\xi_1,\xi_2}$ in \eqref{mathbsk2sdx}.  We denote:
\begin{align}
&\text{For  $\xi,\xi_1,\xi_2\in S_{\mathtt{M}}^\perp$, we denote by $i_\xi,i_{\xi,1},i_{\xi,2}$, the integers} \nonumber \\
& \text{ such that $\mathtt{M}i_{\xi}=\xi,\mathtt{M}i_{\xi,1}=\xi_1,\mathtt{M}i_{\xi,2}=\xi_2$, }\label{ixusdkksdw}\\
&f(x,y):=|x+y|^{1+\alpha}+|x-y|^{1+\alpha}-2(|x|^{1+\alpha} + |y|^{1+\alpha}).\label{fesjsdxc}
\end{align}

\begin{lemma}
For $\xi_1,\xi_2\in S_{\mathtt{M}}^\perp\cup\left\{ 0 \right\}$ such that $\xi_1\ne \xi_2$, we have that
\begin{align}
\mathbb{C}_{\xi_1,\xi_2} &= \frac{ \mathcal{C}_\alpha \pi \mathtt{M}^{\alpha+3}}{192(\alpha+1)}(\mathbb{A}_{S_0^+} -2 (D_{S_0^+})^{3}\vec{F}(\xi_1,\xi_2)\vec{u}^{T}(D_{S_0^+})^\alpha) + W_{r,*}(\xi_1,\xi_2),\label{rksjjsjdsdsd1}
\end{align}
where $\vec{F}$ is defined as
\begin{align}\label{vecosdF}
(\vec{F}(\xi_1,\xi_2))_k:= \begin{cases}
\frac{ f(i_{\xi,1},\mathtt{j}_k) }{|i_{\xi,1}|^{\alpha-1}\mathtt{j}_k^2}, & \text{ if $\xi_1\ne 0, \xi_2=0$},\\
(\vec{F}(\xi_2,\xi_1))_k, &\text{ if $\xi_1=0,\xi_2\ne0$},\\
\frac{\left(i_{\xi,1} f(i_{\xi,1},\mathtt{j}_k)-i_{\xi,2} f(i_{\xi,2},\mathtt{j}_k)\right)}{(i_{\xi,1}|i_{\xi,1}|^{\alpha-1}-i_{\xi,2}|i_{\xi,2}|^{\alpha-1})\mathtt{j}_k^2},&\text{ if $\xi_1,\xi_2\ne 0$}.
\end{cases}
\end{align}
Furthermore, $W_{r,*}(\xi_1,\xi_2)$ satisfies
\begin{align}\label{remshshsdh1s}
|W_{r,*}(\xi_1,\xi_2)|\le_{\alpha,S_0^+,\nu}\mathtt{M}^4.
\end{align}
\end{lemma}
\begin{proof}
For $\xi\in S_{\mathtt{M}}^\perp$, we derive expressions for $\vec{D}(\xi)$, $W(\xi)$ and $\overline{\omega}$ in \eqref{mathbsk2sdx} without assuming $|\xi|$ is dominating the size of the tangential site $S$. 

 We will derive expressions for (see \eqref{mathbsk2sdx} and \eqref{Aeushhsdsdx}),
 \begin{align}\label{sjjsjdwiwd1}
 \mathbb{C}_{\xi_1,\xi_2}=\mathbb{A} -\frac{\pi}6(\mathbb{B})_{\xi_1,\xi_2}.
 \end{align}
 For $\xi\vec{D}(\xi)$ and $k=1,\ldots,\nu$ and $\xi\in S_{\mathtt{M}}^\perp$, we have 
 \begin{align}\label{djhhsdiixcxcsw1}
\xi(\vec{D}(\xi))_k&\overset{\eqref{decompoasd}}=\frac{\alpha}{16}(\xi\vec{m}_{3,1}(\xi))_k + (\xi\vec{m}_4(\xi))_k\nonumber\\
& \overset{\eqref{m3dscsxcsw}}=\frac{\alpha}{16}\xi m^\circ_{3,\alpha}(\xi,j_k)j_k +  (\xi\vec{m}_4(\xi))_k\nonumber\\
&\overset{\eqref{m3deccsdx},\eqref{fesjsdxc}}=\frac{\mathcal{C}_\alpha}{16(\alpha+1)}\xi f(\xi,j_k)j_k +\xi m_{3,\alpha,r}(\xi,j_k) +  (\xi\vec{m}_4(\xi))_k\nonumber\\
&\overset{\eqref{ixusdkksdw},\eqref{stanxcxcx11}}=\frac{\mathcal{C}_\alpha\mathtt{M}^{\alpha+3}}{16(\alpha+1)} i_\xi f(i_\xi,\mathtt{j}_k)\mathtt{j}_k + \underbrace{\xi m_{3,\alpha,r}(\xi,j_k) + (\xi\vec{m}_4(\xi))_k}_{=:(\vec{m}_{7}(\xi))_k}.
 \end{align}
Therefore, for $\xi_1,\xi_2\in S_{\mathtt{M}}^\perp$, we have
\begin{align}
& \xi_1(\vec{D}(\xi_1))_k -\xi_2(\vec{D}(\xi_2))_k  \\
& = \frac{\mathcal{C}_\alpha\mathtt{M}^{\alpha+3}}{16(\alpha+1)}\left( i_{\xi,1} f(i_{\xi,1},\mathtt{j}_k)\mathtt{j}_k-i_{\xi,2} f(i_{\xi,2},\mathtt{j}_k)\mathtt{j}_k\right)  +(\vec{m}_7(\xi_1)-\vec{m}_{7}(\xi_2))_k.\label{djhhsdiixcxcsw2}
\end{align}
 For $\vec{m}_7$ in \eqref{djhhsdiixcxcsw1}, we estimate it as
 \begin{align}
 |\vec{m}_7(\xi)|&\overset{\eqref{estimate1m3alax3},\eqref{mrmsidx1}}{\le_{\alpha,S_0^+,\nu}} \mathtt{M}^{1+\alpha}|i_\xi|^{\alpha-1},\label{rmsdwsxc}\\
 |\vec{m}_7(\xi_1) -\vec{m}_7(\xi_2)|&\overset{\eqref{estimate1m3alax4},\eqref{mrmsidx2}}{\le_{\alpha,S_0^+,\nu}}\mathtt{M}^{1+\alpha}|i_{\xi,1}-i_{\xi,2}|(|i_{\xi,1}|+|i_{\xi,2}|)^{\alpha-1}.\label{xcxcxc11sloww}
 \end{align} 
 Now, we consider $W(\xi)$. From \eqref{wksjdsdsxcsd} and \eqref{wxxrxxcxcxw2}, we have
 \begin{align}
\frac{1}{W(\xi)} & = \frac{1}{-\frac{\mathcal{C}_\alpha}{2}\xi|\xi|^{\alpha-1} + W_r(\xi)}=:-\frac{2}{\mathcal{C}_\alpha}\frac{1}{\xi|\xi|^{\alpha-1}} + W_{r,1}(\xi)\nonumber \\
& =-\frac{2}{\mathcal{C}_\alpha}\mathtt{M}^{-\alpha}\frac{1}{i_{\xi}|i_{\xi}|^{\alpha-1}}+W_{r,1}(\xi),\label{sdsdsddpqodlqkfrj1}
 \end{align}
 where $W_{r,1}(\xi)$ satisfies
 \begin{align}\label{skksdxcxc1}
 |W_{r,1}(\xi)|\le_{\alpha,S_0^+,\nu} |\xi|^{1-2\alpha} =\mathtt{M}^{1-2\alpha}|i_{\xi}|^{1-2\alpha}.
 \end{align}
 Similarly, we have
 \begin{align}\label{qkfrjfskdcxc1}
 \frac{1}{W(\xi_1)-W(\xi_2)}&=-\frac{2}{\mathcal{C}_\alpha}\frac{1}{\xi_1|\xi_1|^{\alpha-1}-\xi_2|\xi_2|^{\alpha-1}}+W_{r,2}(\xi_1,\xi_2) \nonumber\\
 &= -\frac{2}{\mathcal{C}_\alpha}\mathtt{M}^{-\alpha}\frac{1}{i_{\xi,1}|i_{\xi,1}|^{\alpha-1} - i_{\xi,2}|i_{\xi,2}|^{\alpha-1}} + W_{r,2}(\xi_1,\xi_2),
 \end{align}
 where (using \eqref{wxxrxxcxcxw3})
 \begin{align}\label{skcxskksd1xsx2xc}
 |W_{r,2}(\xi_1,\xi_2)|\le_{\alpha,S_0^+,\nu}\mathtt{M}^{1-2\alpha}\frac{1}{|i_{\xi,1}-i_{\xi,2}| (i_{\xi,1}|+|i_{\xi,2}|)^{2\alpha-2}}.
 \end{align}
Therefore, we have
\begin{align}
\frac{\xi(\vec{D}(\xi))_k}{W(\xi)}&\overset{\eqref{djhhsdiixcxcsw1},\eqref{sdsdsddpqodlqkfrj1}}=-\frac{\mathtt{M}^3}{8(\alpha+1)}\frac{ i_\xi f(i_\xi,\mathtt{j}_k)\mathtt{j}_k }{i_{\xi}|i_{\xi}|^{\alpha-1}} + (\vec{W}_{r,3}(\xi))_k,\label{rkksdspqksdx1}\\
\frac{\xi_1(\vec{D}(\xi_1))_k -\xi_2(\vec{D}(\xi_2))_k}{W(\xi_1)-W(\xi_2)}&\overset{\substack{\eqref{djhhsdiixcxcsw2} \\ \eqref{qkfrjfskdcxc1}}}=-\frac{ \mathtt{M}^{3}}{8(\alpha+1)}\frac{\left(i_{\xi,1} f(i_{\xi,1},\mathtt{j}_k)-i_{\xi,2} f(i_{\xi,2},\mathtt{j}_k)\right)\mathtt{j}_k}{\xi_1|\xi_1|^{\alpha-1}-\xi_2|\xi_2|^{\alpha-1}} \nonumber \\
& + (\vec{W}_{r,4}(\xi_1,\xi_2))_k,\label{rkksdspqksdx2}
\end{align}
 where $(\vec{W}_{r,3}(\xi))_k$ and $(\vec{W}_{r,4}(\xi_1,\xi_2))_k$ satisfy
 \begin{align}
 |\vec{W}_{r,3}(\xi)|&\le_{\alpha,\nu,S_0^+} \mathtt{M}^{4-\alpha}|i_\xi|^{1-\alpha},\label{ksdsd1sdsd1}\\
  |\vec{W}_{r,4}(\xi_1,\xi_2)|&\le_{\alpha,\nu,S_0^+}\mathtt{M}^{4-\alpha},\label{ksdsd1sdsd2}
 \end{align}
 which come from the contributions of \eqref{xcxcxc11sloww} and \eqref{skcxskksd1xsx2xc}.
\end{proof}

Therefore, combining these with \eqref{delfjscxcxcs}, we have
\begin{align}
\frac{\xi(\vec{D}(\xi))_i}{W(\xi)}(\overline{\omega})_k&=\frac{ \mathcal{C}_\alpha\mathtt{M}^{\alpha+3}}{16(\alpha+1)} \frac{ i_\xi f(i_\xi,\mathtt{j}_i)\mathtt{j}_i }{i_{\xi}|i_{\xi}|^{\alpha-1}} \mathtt{j}_k^{\alpha} + (W_{r,5}(\xi))^{i}_k,\label{sjdsjdsdad1}\\
\frac{\xi_1(\vec{D}(\xi_1))_i -\xi_2(\vec{D}(\xi_2))_i}{W(\xi_1)-W(\xi_2)}(\overline{\omega})_k&=\frac{ \mathcal{C}_\alpha\mathtt{M}^{\alpha+3}}{16(\alpha+1)}\frac{\left(i_{\xi,1} f(i_{\xi,1},\mathtt{j}_i)-i_{\xi,2} f(i_{\xi,2},\mathtt{j}_i)\right)\mathtt{j}_i}{\xi_1|\xi_1|^{\alpha-1}-\xi_2|\xi_2|^{\alpha-1}} \mathtt{j}_k^{\alpha} \nonumber\\
& \ + (W_{r,6}(\xi))^i_k,
\end{align}
where $ (W_{r,5}(\xi))^{i}_k,  (W_{r,6}(\xi))^{i}_k$ satisfy (from \eqref{ksdsd1sdsd1},\eqref{ksdsd1sdsd2} and \eqref{shhxcxca11})
\begin{align}
|W_{r,5}(\xi))^{i}_k|&\le_{\alpha,S_0^+,\nu}\mathtt{M}^4|i_{\xi}|^{1-\alpha}\le \mathtt{M}^4,\label{hhsddsadsdwq1}\\
|W_{r,6}(\xi))^{i}_k|&\le_{\alpha,S_0^+,\nu}\mathtt{M}^4.\label{hhsddsadsdwq2}
\end{align}
Therefore, using \eqref{dkhhxcxc11sdsd}, \eqref{bksjdexx1}, \eqref{sjjsjdwiwd1} and \eqref{mathbsk2sdx} we obtain 

\begin{align}
\left(\mathbb{C}_{\xi,0}\right)^i_k &=\frac{ \mathcal{C}_\alpha \pi \mathtt{M}^{\alpha+3}}{192(\alpha+1)} \left(\left(\mathbb{A}_{S_0^+}\right)^i_k -2\frac{ i_\xi f(i_\xi,\mathtt{j}_i)\mathtt{j}_i }{i_{\xi}|i_{\xi}|^{\alpha-1}} \mathtt{j}_k^{\alpha}\right) + \left(W_{r,7}(\xi)\right)^i_k,\label{kksdppwfinslds1}\\
\left(\mathbb{C}_{\xi_1,\xi_2}\right)^i_k  & = \frac{ \mathcal{C}_\alpha \pi \mathtt{M}^{\alpha+3}}{192(\alpha+1)} \left(\left(\mathbb{A}_{S_0^+}\right)^i_k - 2\frac{\left(i_{\xi,1} f(i_{\xi,1},\mathtt{j}_i)-i_{\xi,2} f(i_{\xi,2},\mathtt{j}_i)\right)\mathtt{j}_i}{i_{\xi,1}|i_{\xi,1}|^{\alpha-1}-i_{\xi,2}|i_{\xi,2}|^{\alpha-1}} \mathtt{j}_k^{\alpha} \right) \nonumber \\ & + \left(W_{r,8}\right)^i_k,\label{kksdppwfinslds2}
\end{align}
where
\begin{align}
| \left(W_{r,7}(\xi)\right)^i_k|, | \left(W_{r,8}(\xi)\right)^i_k|\le_{\alpha,S_0^+,\nu} \mathtt{M}^4.\label{slllsds1w111}
\end{align}
Therefore, using the notation $\vec{F}$ in \eqref{vecosdF}, we obtain \eqref{vecosdF}. The estimate \eqref{remshshsdh1s} follows from \eqref{kksdppwfinslds1} and \eqref{kksdppwfinslds2}.

\section{Choice of tangential sites}
The goal of this section is to prove Proposition~\ref{vefiisdkswjsjdsdwsd}, which shows that there are infinitely many choices of $S_0^+$ and $\mathtt{M}$ so that the hypotheses described in Section~\ref{rpoisdsd1sd} are satisfied.

\begin{lemma}\label{ksjjsdsd1B12s}
For each $2\le\nu\in \mathbb{N}$ and $\alpha\in (1,2)$, there exist infinite number of choices of $S^+_0=\left\{0<\mathtt{j}_1<\cdots<\mathtt{j}_\nu\right\}$ such that 
\begin{enumerate}[label=(\arabic*)]
\item \label{measure_ef0} The tangential sites $S$ constructed by \eqref{tan_site} satisfy \ref{tangent_2} and \ref{tangent_1}.
\item \label{measure_ef1} $\mathbb{A}_{S^+_0}$ is invertible.
\item \label{measure_ef2}$1-2\alpha(\alpha+1)I_{S_0^+} \ne 0$. 
\item \label{measure_ef3}For all $\xi_1,\xi_2\in S_\mathtt{M}^\perp\cup\left\{ 0 \right\}$ such that $\xi_1\ne \xi_2$, it holds that
\begin{align}\label{jjsususdosod1}
\min_{\vec{U}\in \mathbb{R}^{\nu},\ |\vec{U}|=1}\left|\vec{U}\cdot \left( \mathbb{A}_{S_0^+} -2 (D_{S_0^+})^{3}\vec{F}(\xi_1,\xi_2)\vec{u}^{T}(D_{S_0^+})^\alpha\right)\vec{U}\right|\ge_{\alpha,\nu,S_0^+}1.
\end{align}
\end{enumerate}
\end{lemma}
\begin{proof}
We look for $S_0^+$ such that $\mathtt{j}_1=1$, that is, 
\begin{align}\label{forms1}
S_0^+= \left\{ 1 < \mathtt{j}_2 < \cdots < \mathtt{j}_\nu\right\}.
\end{align}
We denote 
\begin{align}\label{rationslsdj1s}
R_{S_0^+}:=\max_{i < k}\left\{ \frac{\mathtt{j}_i}{\mathtt{j}_k}\right\} < 1.
\end{align}
We will show that if $S_0^+$ of the form \eqref{forms1} is chosen so that $R_{S_0^+}$ is sufficiently small, then \ref{measure_ef1}-\ref{measure_ef3} are satisfied. From the conditions \ref{tangent_2} and \ref{tangent_1}, and the asymptotic properties of $\kappa_j$ that follow from \ref{defoflsksdjwwww} and Lemma~\ref{prosdjsdowsds1}, it is straightforward to see that such a choice of $S_0^+$ with small $R_{S_0^+}$ can be made so that \ref{measure_ef0} holds true as well. Therefore, we will focus on the proof of \ref{measure_ef1}-\ref{measure_ef3}.
 
 \vspace{0.5\baselineskip}\noindent\textit{Proof of \ref{measure_ef1}.}
From \eqref{ansxs2hd1sdcs} and \eqref{fesjsdxc}, we rewrite $\mathbb{A}_{S_0^+}$ as
\begin{align}\label{rkksdsdsdsaldksgotj}
\mathbb{A}_{S_0^+}=( 2^{\alpha + 1}- 4 )(D_{S_0^+})^{\frac{\alpha+3}{2}}\left(I + \frac{2}{2^{\alpha+1}-4} L\right)(D_{S_0^+})^{\frac{\alpha+3}{2}},
\end{align}
where
\begin{align}\label{lssjsksjddw2}
(L)^i_{k} = \begin{cases}
0, & \text{ if $i=k$,}\\
f\left(\frac{\mathtt{j}_i}{\mathtt{j}_k},\frac{\mathtt{j}_k}{\mathtt{j}_i}\right), & \text{ if $i\ne k$.}
\end{cases}
\end{align}
Note that for $0\ne x\in\mathbb{R}$, recalling $f$ from \eqref{fesjsdxc}, we see that
\[
f(x,x^{-1})= |x+x^{-1}|^{1+\alpha} + |x-x^{-1}|^{1+\alpha} - 2(x^{1+\alpha}+x^{-(1+\alpha)})
\]
satisfies
\begin{align}
|f(x,x^{-1})|\le_\alpha
\begin{cases}
|x|^{\alpha-3}, & \text{ if $|x|>1$},\\
|x|^{3-\alpha}, & \text{ if $|x|<1$.}
\end{cases}
\end{align}
Therefore, \eqref{lssjsksjddw2} tells us that 
\begin{align}\label{liksohsds1}
|L^i_k|\le_\alpha 
\begin{cases}
\left(\frac{\mathtt{j}_i}{\mathtt{j}_{k}}\right)^{3-\alpha}\le R_{S_0^+}^{3-\alpha}, & \text{ if $i > k$},\\
\left(\frac{\mathtt{j}_k}{\mathtt{j}_{i}}\right)^{3-\alpha}\le R_{S_0^+}^{3-\alpha}, & \text{ if $i < k$}.
\end{cases}
\end{align}
Since $2^{\alpha+1}-4\ne 0$ for $\alpha>1$, we have that $I + \frac{2}{2^{\alpha+1}-4} L$ is invertible, if $R_{S_0^+}$ is sufficiently small depending on $\nu$. Therefore from \eqref{rkksdsdsdsaldksgotj}, we see that  $\mathbb{A}_{S_0^+}$ is invertible.

 \vspace{0.5\baselineskip}\noindent\textit{Proof of \ref{measure_ef2}.}
 We first derive more careful estimates for the inverse of $I + \frac{2}{2^{\alpha+1}-4} L$. Using a Neumann series, we write
 \begin{align}\label{ltidlejssd2}
 \left(I + \frac{2}{2^{\alpha+1}-4} L\right)^{-1} =I + \sum_{n=1}^{\infty}\left(-\frac{2}{2^{\alpha+1}-4}\right)^n L^{n}=: I - \tilde{L}
  \end{align}
  We claim that  for all $n\ge 1$,
  \begin{align}\label{inmpiisdsdsd1}
  \begin{cases}
  |({L}^n)^{i}_k|\le \left(3\nu R_{S_0^+}^{\frac{3-\alpha}{2}}\right)^{n-1}\left(\frac{\mathtt{j}_i}{\mathtt{j}_k}\right)^{\frac{3-\alpha}{2}}, &\text{ if $k\ge i$},\\
  |({L}^n)^{i}_k|\le \left(3\nu R_{S_0^+}^{\frac{3-\alpha}{2}}\right)^{n-1} \left(\frac{\mathtt{j}_k}{\mathtt{j}_i}\right)^{\frac{3-\alpha}{2}}, &\text{ if $k\le i$}.
  \end{cases}
  \end{align}
  Since $L$ is a symmetric matrix, we can assume that $k\ge i$ to prove the claim. If $n=1$, then the estimate follows immediately from \eqref{liksohsds1} and \eqref{lssjsksjddw2}. Assuming \eqref{ltidlejssd2} holds for some $n\ge 1$, we estimate $(L^{n+1})^i_k$ for $k\ge i$ as 
   \begin{align}\label{kfkdgosd1}
   |(L^{n+1})^i_k| &\le \sum_{j=1}^\nu|(L^n)^j_k||L^i_j|\nonumber\\
   &=\sum_{1\le j\le i}|(L^n)^j_k||L^i_j|+\sum_{ i<j\le k}|(L^n)^j_k||L^i_j|+\sum_{k <j\le \nu}^\nu|(L^n)^j_k||L^i_j|\nonumber\\
   &\overset{\eqref{inmpiisdsdsd1},\eqref{liksohsds1}}\le \left(3\nu R_{S_0^+}^{\frac{3-\alpha}{2}}\right)^{n-1}\left(\sum_{1\le j\le i} \left(\frac{\mathtt{j}_j}{\mathtt{j}_k} \right)^{\frac{3-\alpha}{2}}\left(\frac{\mathtt{j}_j}{\mathtt{j}_i}\right)^{3-\alpha}\right.\nonumber\\
   & \left.+ \sum_{i < j\le k} \left(\frac{\mathtt{j}_j}{\mathtt{j}_k} \right)^{\frac{3-\alpha}{2}}\left(\frac{\mathtt{j}_i}{\mathtt{j}_j}\right)^{3-\alpha} + \sum_{k < j\le \nu} \left(\frac{\mathtt{j}_k}{\mathtt{j}_j} \right)^{\frac{3-\alpha}{2}}\left(\frac{\mathtt{j}_i}{\mathtt{j}_j}\right)^{3-\alpha}\right)
     \end{align}
Using that $0 <\mathtt{j}_1 < \cdots < \mathtt{j}_\nu$, we have 
\begin{align*}
\sum_{1\le j\le i} \left(\frac{\mathtt{j}_j}{\mathtt{j}_k} \right)^{\frac{3-\alpha}{2}}\left(\frac{\mathtt{j}_j}{\mathtt{j}_i}\right)^{3-\alpha}&\le \nu R_{S_0^+}^{3-\alpha}\left( \frac{\mathtt{j}_i}{\mathtt{j}_k}\right)^{\frac{3-\alpha}2},\\
\sum_{i < j\le k} \left(\frac{\mathtt{j}_j}{\mathtt{j}_k} \right)^{\frac{3-\alpha}{2}}\left(\frac{\mathtt{j}_i}{\mathtt{j}_j}\right)^{3-\alpha} &=\sum_{i < j\le k} \left(\frac{\mathtt{j}_i}{\mathtt{j}_k} \right)^{\frac{3-\alpha}{2}}\left(\frac{\mathtt{j}_j}{\mathtt{j}_i}\right)^{\frac{3-\alpha}2}\left(\frac{\mathtt{j}_i}{\mathtt{j}_j}\right)^{3-\alpha}\le \nu R_{S_0^+}^{\frac{3-\alpha}2}\left(\frac{\mathtt{j}_i}{\mathtt{j}_k}\right)^{\frac{3-\alpha}2},\\
\sum_{k < j\le \nu} \left(\frac{\mathtt{j}_k}{\mathtt{j}_j} \right)^{\frac{3-\alpha}{2}}\left(\frac{\mathtt{j}_i}{\mathtt{j}_j}\right)^{3-\alpha} &=\sum_{k < j\le \nu} \left(\frac{\mathtt{j}_i}{\mathtt{j}_k}\right)^{\frac{3-\alpha}2}\left(\frac{\mathtt{j}_i\mathtt{j}_k^2}{\mathtt{j}_j^{3}} \right)^{\frac{3-\alpha}2}\le \nu R_{S_0^+}^{\frac{3-\alpha}2}\left(\frac{\mathtt{j}_i}{\mathtt{j}_k}\right)^{\frac{3-\alpha}2}.
\end{align*}
 Plugging these estimates into \eqref{kfkdgosd1}, we obtain \eqref{inmpiisdsdsd1} for $n+1$, therefore, the claim is proved.
 
  From the definition of $\tilde{L}$ in \eqref{ltidlejssd2}, we see that \eqref{inmpiisdsdsd1} implies
  \begin{align}\label{shhsdasdjwdsd}
  |\tilde{L}^i_k|\le_{\alpha,\nu}R_{S_0^+}^{\frac{3-\alpha}2}\min\left\{ \left( \frac{\mathtt{j}_i}{\mathtt{j}_k}\right)^{\frac{3-\alpha}2}, \left( \frac{\mathtt{j}_k}{\mathtt{j}_i}\right)^{\frac{3-\alpha}2} \right\}.
  \end{align}
  
 Recalling from \eqref{Idjsdxc1},  we have that
 \begin{align*}
 I_{S_0^+} &= \vec{u}\cdot \left((D_{S_0^+})^{-\alpha}\mathbb{A}_{S_0^+}(D_{S_0^+})^{-3}\right)^{-1}\vec{u}\\
 & = \vec{u}\cdot \left((D_{S_0^+})^{\frac{3-\alpha}2}(D_{S_0^+})^{-\frac{\alpha+3}{2}}\mathbb{A}_{S_0^+}(D_{S_0^+})^{-\frac{\alpha+3}2}(D_{S_0^+})^{\frac{\alpha-3}2}\right)^{-1}\vec{u}\\
 & = \vec{u}\cdot (D_{S_0^+})^{\frac{3-\alpha}2} \left( (D_{S_0^+})^{-\frac{\alpha+3}{2}}\mathbb{A}_{S_0^+}(D_{S_0^+})^{-\frac{\alpha+3}2}\right)^{-1} (D_{S_0^+})^{\frac{\alpha-3}2}\vec{u}\\
 & \overset{\eqref{rkksdsdsdsaldksgotj}}= (D_{S_0^+})^{\frac{3-\alpha}2}\vec{u}\cdot (( 2^{\alpha + 1}- 4 )I + 2 L)^{-1} (D_{S_0^+})^{\frac{\alpha-3}2}\vec{u}\\
 &\overset{\eqref{ltidlejssd2}}= (D_{S_0^+})^{\frac{3-\alpha}2}\vec{u} \cdot \left(( 2^{\alpha + 1}- 4 )^{-1}(I-\tilde{L})\right) (D_{S_0^+})^{\frac{\alpha-3}2}\vec{u}\\
 & =( 2^{\alpha + 1}- 4 )^{-1}\left( |\vec{u}|^2 +  (D_{S_0^+})^{\frac{3-\alpha}2}\vec{u} \cdot \tilde{L} (D_{S_0^+})^{\frac{\alpha-3}2}\vec{u}\right)\\
 &\overset{\eqref{Idjsdxc1}}=( 2^{\alpha + 1}- 4 )^{-1}\left(\nu +\sum_{i,k}\mathtt{j}_k^{\frac{3-\alpha}2}\tilde{L}^i_k\mathtt{j}_i^{\frac{\alpha-3}2}\right)\\
 &\overset{ \eqref{shhsdasdjwdsd}}\ge ( 2^{\alpha + 1}- 4 )^{-1}\left(\nu - c_{\alpha,\nu,S_0^+} R_{S_0^+}^{\frac{3-\alpha}2}\right),
 \end{align*}
 for some constant $c_{\alpha,\nu,S_0^+}>0$. Therefore, choosing $S_0^+$ so that $R_{S_0^+}^{\frac{3-\alpha}2}$ is sufficiently small, we have $$I_{S_0^+}>( 2^{\alpha + 1}- 4 )^{-1}\frac{\nu}2 \overset{\nu \ge 2}\ge( 2^{\alpha + 1}- 4 )^{-1}.$$ Since $\alpha\in(1,2)$, we see that
 \[
 2\alpha(\alpha+1)I_{S_0^+} > 2\alpha(\alpha+1)(2^{\alpha+1}-4)^{-1}\ge 3,
 \]
  which implies that $1-2\alpha(\alpha+1)I_{S_0^+} \ne 0$.

 \vspace{0.5\baselineskip}\noindent\textit{Proof of \ref{measure_ef3}.}
From the definition of $f$ in \eqref{fesjsdxc}, the following properties can be verified straightforwardly: For each $j\in \mathbb{N}$,
 \begin{align}
& f(\xi,j) > 0, \text{ for all $\xi,j\ne 0$},\label{fprossd1}\\
& \xi\mapsto \frac{f(\xi,j)}{|\xi|^{\alpha-1}} \text{ is even and monotone increasing in $\xi>0$,}\label{fprossd2}\\
&  \frac{f(\xi,j)}{|\xi|^{\alpha-1}j^2} \le_{\alpha} 1\text{ and }\frac{\xi_1 f(\xi_1,j) - \xi_2 f(\xi_2,j)}{(\xi_1|\xi_1|^{\alpha-1} -\xi_2|\xi_2|^{\alpha-1})j^2}\le_\alpha 1, \text{ for all $\xi,\xi_1,\xi_2\ne 0$} \label{fprossd3}
 \end{align}
With the above properties of $f$, we will estimate $(\vec{F}(\xi_1,\xi_2))_k$ in \eqref{vecosdF}. For the case where one of $\xi_1,\xi_2$ equals to $0$, we can assume without loss of generality that $\xi_1\ne 0$ and $\xi_2=0$, since $\vec{F}(\xi_1,\xi_2)=\vec{F}(\xi_2,\xi_1)$. In such case, we have
\begin{align}\label{xssd22111}
(\vec{F}(\xi_1,\xi_2))_k \overset{\eqref{vecosdF}}= \frac{f(i_{\xi,1},\mathtt{j}_k)}{{|i_{\xi,1}|^{\alpha-1}\mathtt{j}_k^2}}\overset{\eqref{fprossd2}}\ge \frac{f(2,\mathtt{j}_k)}{\mathtt{j}_k^2},
\end{align}
where the last equality follows from $|{i}_{\xi,1}|> \mathtt{j}_1=1$ for $\xi_1\in S^\perp$. For an upper bound, we also have
\begin{align}\label{xssd22112}
|(\vec{F}(\xi_1,\xi_2))_k| =\left| \frac{f(i_{\xi,1},\mathtt{j}_k)}{{|i_{\xi,1}|^{\alpha-1}\mathtt{j}_k^2}}\right|\overset{\eqref{fprossd3}}{\le_\alpha} 1.
\end{align}
Now, let us assume $\xi_1,\xi_2\ne 0$. In this case, it follows from the explicit formula for $\vec{F}$ in \eqref{vecosdF} and \eqref{fprossd2} that we can without loss of generality that $\xi_1>|\xi_2|>0$. For such $\xi_1,\xi_2$, it follows from \eqref{fprossd1} and  \eqref{fprossd2} that
\begin{align}
(\vec{F}(\xi_1,\xi_2))_k &\overset{\eqref{vecosdF}}=\frac{\left(i_{\xi,1} f(i_{\xi,1},\mathtt{j}_k)-i_{\xi,2} f(i_{\xi,2},\mathtt{j}_k)\right)}{(i_{\xi,1}|i_{\xi,1}|^{\alpha-1}-i_{\xi,2}|i_{\xi,2}|^{\alpha-1})\mathtt{j}_k^2}  \nonumber \\ & \ge \frac{i_{\xi,1} f(i_{\xi,1},\mathtt{j}_k)}{i_{\xi,1}|i_{\xi,1}|^{\alpha-1}\mathtt{j}_k^2}\overset{\eqref{fprossd2}}\ge\frac{f(2,\mathtt{j}_k)}{\mathtt{j}_k^2},\label{jjrjsdssdsdsd1}
\end{align}
where the first inequality follows from the elementary fact that $\frac{b-d}{a-c}\ge \frac{b}{a}$ for $a,b,c,d>0$ such that $a>c$ and  $\frac{b}{a}>\frac{d}{c}$.
For an upper bound, we also have
\begin{align}\label{jjrjsdssdsdsd2}
|(\vec{F}(\xi_1,\xi_2))_k| \overset{\eqref{vecosdF}}=\left|\frac{\left(i_{\xi,1} f(i_{\xi,1},\mathtt{j}_k)-i_{\xi,2} f(i_{\xi,2},\mathtt{j}_k)\right)}{(i_{\xi,1}|i_{\xi,1}|^{\alpha-1}-i_{\xi,2}|i_{\xi,2}|^{\alpha-1})\mathtt{j}_k^2}\right|\overset{\eqref{fprossd3}}{\le_\alpha}1.
\end{align}
Therefore, combining the above bounds for all $\xi_1,\xi_2\in S^\perp\left\{ 0 \right\}$ such that $\xi_1\ne \xi_2$, we have
\begin{align}
(\vec{F}(\xi_1,\xi_2))_k &\overset{\eqref{xssd22111},\eqref{jjrjsdssdsdsd1}}\ge  \frac{f(2,\mathtt{j}_k)}{\mathtt{j}_k^2} >0, \label{sjsjdjsdad1sxfF}\\
|(\vec{F}(\xi_1,\xi_2))_k|&\overset{\eqref{xssd22112},\eqref{jjrjsdssdsdsd2}}{\le_\alpha}1.\label{sjsjdjsdad1sxfF2}
\end{align}
Towards \eqref{jjsususdosod1}, we write
\begin{align}
&\mathbb{A}_{S_0^+} -2 (D_{S_0^+})^{3}\vec{F}(\xi_1,\xi_2)\vec{u}^{T}(D_{S_0^+})^\alpha\nonumber\\
&= (D_{S_0^+})^{\frac{\alpha+3}2} \nonumber \\
& \qquad \times \left((D_{S_0^+})^{-\frac{\alpha+3}2}\mathbb{A}_{S_0^+}(D_{S_0^+})^{-\frac{\alpha+3}2} -2 (D_{S_0^+})^{\frac{3-\alpha}{2}}\vec{F}(\xi_1,\xi_2)\vec{u}^{T}(D_{S_0^+})^{\frac{\alpha-3}2} \right)(D_{S_0^+})^{\frac{\alpha+3}2}\nonumber\\
&\overset{\eqref{rkksdsdsdsaldksgotj}}=(2^{\alpha+1}-4)(D_{S_0^+})^{\frac{\alpha+3}2} \nonumber \\
& \qquad \qquad \times \underbrace{\left(I +\frac{2}{2^{\alpha+1}-4}L  -\frac{2}{2^{\alpha+1}-4} (D_{S_0^+})^{\frac{3-\alpha}{2}}\vec{F}(\xi_1,\xi_2)\vec{u}^{T}(D_{S_0^+})^{\frac{\alpha-3}2} \right)}_{=:\mathbb{B}_{S_0^+}(\xi_1,\xi_2)} (D_{S_0^+})^{\frac{\alpha+3}2}\label{Bsisdjsddef}
\end{align}
Using the norm in \eqref{sposx},  it follows from   \eqref{liksohsds1} and \eqref{sjsjdjsdad1sxfF2}  that
\begin{align}\label{lsdjsdsd1sd}
|\mathbb{B}_{S_0^+}(\xi_1,\xi_2)|\le_{\alpha,\nu,S_0^+} 1
\end{align}
where the upper bound is  independent of $\xi_1,\xi_2$. This implies that the largest eigenvalue of $\mathbb{B}_{S_0^+}$ is bounded. Furthermore, we recall the matrix determinant Lemma:
\begin{align}
\text{det}(A + \vec{U}\vec{V}^T) = (1+ \vec{V}\cdot A^{-1}\vec{U})\text{det}(A),\text{ for $\vec{U},\vec{V}\in \mathbb{R}^\nu,\ A\in\mathbb{R}^{\nu\times \nu}$}.
\end{align}
Applying this to $\mathbb{B}_{S_0^+}$, we get
\begin{align}
&\text{det}(\mathbb{B}_{S_0^+}(\xi_1,\xi_2))\nonumber\\
& = \left(1 - \frac{2}{2^{\alpha+1}-4}(D_{S_0^+})^{\frac{\alpha-3}2}\vec{u}\cdot \left( I+\frac{2}{2^{\alpha+1}-4}L\right)^{-1} (D_{S_0^+})^{\frac{3-\alpha}{2}}\vec{F}(\xi_1,\xi_2)\right) \nonumber \\
& \qquad \qquad \qquad \qquad \qquad \qquad \qquad \qquad \qquad \qquad \qquad \qquad \times \text{det}\left(I+\frac{2}{2^{\alpha+1}-4}L\right)\nonumber\\
&\overset{\eqref{ltidlejssd2}}= \left(1 - \frac{2}{2^{\alpha+1}-4}(D_{S_0^+})^{\frac{\alpha-3}2}\vec{u}\cdot (I-\tilde{L})(D_{S_0^+})^{\frac{3-\alpha}{2}}\vec{F}(\xi_1,\xi_2)\right)\text{det}\left(I+\frac{2}{2^{\alpha+1}-4}L\right)\label{rkksdjvs11s2ss}
\end{align}
We can further compute
\begin{align}
\frac{2}{2^{\alpha+1}-4}(D_{S_0^+})^{\frac{\alpha-3}2}\vec{u}\cdot (D_{S_0^+})^{\frac{3-\alpha}{2}}\vec{F}(\xi_1,\xi_2) &\overset{\eqref{vecudefsd}}= \frac{2}{2^{\alpha+1}-4}\sum_{k=1}^{\nu}(\vec{F}(\xi_1,\xi_2))_k\nonumber\\
& \overset{\eqref{sjsjdjsdad1sxfF}}\ge \frac{2}{2^{\alpha+1}-4}\sum_{k=1}^{\nu}\frac{f(2,\mathtt{j}_k)}{\mathtt{j}_k^2}\nonumber\\
&\overset{\eqref{fprossd1},\eqref{forms1}}\ge\frac{2}{2^{\alpha+1}-4}{f(2,1)}\nonumber\\
&\overset{\eqref{fesjsdxc}}=\frac{2(3^{1+\alpha}-2^{2+\alpha}-1)}{2^{\alpha+1}-4}\nonumber\\
&\ge 3,\label{rkjsdsjdshdww11}
\end{align}
where the last inequality can be straightforwardly proved using that $\alpha\in (1,2)$. Also, we have
\begin{align}
\left|\frac{2}{2^{\alpha+1}-4}(D_{S_0^+})^{\frac{\alpha-3}2}\vec{u}\cdot \tilde{L}(D_{S_0^+})^{\frac{3-\alpha}{2}}\vec{F}(\xi_1,\xi_2)\right|& \le_{\alpha}\sum_{k,i=1}^{\nu}\mathtt{j}_i^{\frac{\alpha-3}2}|\tilde{L}^{i}_k|\mathtt{j}_k^{\frac{3-\alpha}{2}}|(\vec{F}(\xi_1,\xi_2))_k|\nonumber\\
&\overset{\eqref{shhsdasdjwdsd},\eqref{sjsjdjsdad1sxfF2}}{\le_{\alpha,\nu}}R_{S_0^+}^{\frac{3-\alpha}2}\nonumber\\
& \le1,\label{smallsdjsjds}
\end{align}
for sufficiently small $R_{S_0^+}$.
Plugging this and \eqref{rkjsdsjdshdww11} into \eqref{rkksdjvs11s2ss} and using that $I+\frac{2}{2^{\alpha+1}-4}L$ is invertible, we see that 
\begin{align}\label{sdjsjdnotns}
|\text{det}(\mathbb{B}_{S_0^+}(\xi_1,\xi_2))| \ge_{\alpha,\nu} 1,
\end{align}
independently of $\xi_1,\xi_2$. Combining this with \eqref{lsdjsdsd1sd}, we see that 
\begin{align*}
\min_{\vec{U}\in \mathbb{R}^{\nu},\ |\vec{U}|=1}\left|\vec{U}\cdot \mathbb{B}_{S_0^+}(\xi_1,\xi_2)\vec{U}\right|\ge_{\alpha,\nu,S_0^+}1.
\end{align*}
In view of \eqref{Bsisdjsddef}, this implies \eqref{jjsususdosod1}. 
\end{proof}

\begin{proposition}\label{vefiisdkswjsjdsdwsd}
For each $2\le\nu\in \mathbb{N}$ and $\alpha\in (1,2)$, there exist  infinite number of choices of $S_0^+$ such that for sufficiently large $\mathtt{M}$ depending on $\alpha,\nu,S_0^+$, the following hold:
\begin{enumerate}[label=(\arabic*)]
\item \label{rkjjsodsdwewsd0} The tangential sites $S$ constructed by \eqref{tan_site} satisfies \ref{tangent_2} and \ref{tangent_1}.
\item \label{rkjjsodsdwewsd1} \ref{hypos1s1} holds true. That is, $\mathbb{A}$ in \eqref{amplitude_modulation1} is invertible.
\item \label{rkjjsodsdwewsd2}  \ref{hypothsdj2} holds true. That is, there exists a constant $C_{\mathtt{H}2}({S}_0^+,\nu,\alpha,\mathtt{M})>0$ such that  if $\xi_1,\xi_2\in S_{\mathtt{M}}^\perp\cup \left\{ 0 \right\}$ and  $\max\left\{|\xi_1|,|\xi_2|\right\}\ge C_{\mathtt{H}2}$, then
 \begin{align}
 &|\xi_1 - \xi_2| (|\xi_1|^{\alpha-1} + |\xi_2|^{\alpha-1}) \nonumber \\
 & \le_{\alpha,\nu,S} \left|\left(W(\xi_1) - \frac{\pi}6\xi_1\vec{D}(\xi_1)\cdot\mathbb{A}^{-1}\overline{\omega}\right)-\left(W(\xi_2) - \frac{\pi}6\xi_2\vec{D}(\xi_2)\cdot\mathbb{A}^{-1}\overline{\omega}\right)\right|.\label{ajsj21esti2ma22x12}
 \end{align}
\item \label{rkjjsodsdwewsd3} \ref{hypothsdj22} holds true. That is, $\text{det}(\mathbb{C}_{\xi_1,\xi_2})\ne 0$ for all $\xi_1,\xi_2\in S_{\mathtt{M}}^\perp\cup\left\{ 0 \right\}$.
\end{enumerate}
\end{proposition}
\begin{proof}
Let us choose $S_0^+$ such that \ref{measure_ef0}-\ref{measure_ef3} in Lemma~\ref{ksjjsdsd1B12s} are satisfied. For such a choice of $S_0^+$, item \ref{rkjjsodsdwewsd0} follows trivially from \ref{measure_ef0} of Lemma~\ref{ksjjsdsd1B12s}.

\vspace{0.5\baselineskip}\noindent\textit{Proof of \ref{rkjjsodsdwewsd1}.}
 \ref{rkjjsodsdwewsd1} follows from  \ref{measure_ef1} of Lemma~\ref{ksjjsdsd1B12s} and \eqref{inversitjbhsksdwsz}.

\vspace{0.5\baselineskip}\noindent\textit{Proof of \ref{rkjjsodsdwewsd2}.}
It follows from \eqref{Eajhxcxc1} that \eqref{ajsj21esti2ma22x12} is equivalent to
\begin{align}\label{ajsj1estima2x12}
|\xi_1-\xi_2|(|\xi_1|^{\alpha-1} + |\xi_2|^{\alpha-2})\le_{\alpha,\nu,S}\left| \mathtt{A}(\xi_1)-\mathtt{A}(\xi_2)\right|.
\end{align}

From \ref{measure_ef2} in Lemma~\ref{ksjjsdsd1B12s} and Lemma~\ref{kkkrpsosdsd12p2p2p}, there exist two constants  $C_1(\alpha,\nu,S_0^+)\ne 0\ ,C_2(\alpha,\nu,S_0^+)> 0$ such that for $|\xi|\ge 2 \text{max}_{S^+}:=2\max\left\{ j_k: j_k\in S^+\right\}$,
\begin{align}
\mathtt{A}(\xi) &= C_1\xi|\xi|^{\alpha-1} + m_{\mathtt{A},r}(\xi),\label{rjjsdsdsdsd1}\\
|m_{\mathtt{A},r}(\xi)|&{\le}C_2\left( \mathtt{M}^{1-\alpha}|\xi|^{\alpha} + \mathtt{M}|\xi|^{\alpha-1}\right),\label{skkksdhsd2123}\\
|m_{\mathtt{A},r}(\xi_1) - m_{\mathtt{A},r}(\xi_2)|&{\le}C_2|\xi_1-\xi_2|(\mathtt{M}^{1-\alpha}(|\xi_1|+|\xi_2|)^{\alpha-1} + \mathtt{M}(|\xi_1|+|\xi_2|)^{\alpha-2}). \label{skkksdhsd2223}
\end{align}
We also see from \eqref{Eajhxcxc1} and \eqref{jjjsss11sdssdsossdow1} that for all $\xi,\xi_1,\xi_2\in S_\mathtt{M}^\perp\cup\left\{ 0 \right\}$,
\begin{align}
|\mathtt{A}(\xi)|&\le C_3|\xi|^{\alpha},\label{rlkksdowow1}
\end{align}
for some $C_3=C_3(\alpha,\nu,S_0^+,\mathtt{M})>0$. In addition, it is straightforward to see that 
\begin{align}\label{alpsdoosd1}
|\xi_1|\xi_1|^{\alpha-1} - \xi_2|\xi_2|^{\alpha-1}|\ge c_\alpha |\xi_1-\xi_2|(|\xi_1|^{\alpha-1} + |\xi_2|^{\alpha-1}),
\end{align}
for some $c_\alpha>0$, that depends on only $\alpha\in (1,2)$. Without loss of generality, we can assume that \begin{align}
c_\alpha<1.\label{csdjsjd1}
\end{align}

To prove \eqref{ajsj1estima2x12}, let us choose  $\mathtt{M}$ and $C_{\mathtt{H}2}$ so that 
\begin{align}\label{cjjjsdiidwda}
 c_\alpha|C_1| - C_2\mathtt{M}^{1-\alpha}\ge \frac{c_\alpha|C_1|}{2}, \quad C_{\mathtt{H}2}\ge \max\left\{\frac{4C_2\mathtt{M}}{|C_1|c_\alpha}, 2 \text{max}_{S^+}, \frac{16C_3\text{max}_{S+}}{|C_1|c_\alpha}\right\},
\end{align}
which is possible since $C_1$ and $C_2$ depend on $\alpha,\nu,S_0^+$ only.
Note that without loss of generality, we can assume that $|\xi_1| \ge |\xi_2|$, so that
\begin{align}\label{xisd1alsdwd}
|\xi_1|\ge C_{\mathtt{H}2}.
\end{align} Then we consider two cases: 1) $|\xi_2|\le 2 \text{max}_{S^+}$ and 2) $|\xi_2|\ge 2 \text{max}_{S^+}$.

 In the first case, we have that
 \begin{align}\label{skjsdsdxi2222s}
|\mathtt{A}(\xi_2)|\overset{\eqref{rlkksdowow1}}\le  C_3|\xi_2|^{\alpha}\le 2\text{max}_{S^+}C_3.
\end{align}
while $\mathtt{A}(\xi_1)$ can be estimated as
 \begin{align}
 |\mathtt{A}(\xi_1)| &\overset{\eqref{rjjsdsdsdsd1},\eqref{xisd1alsdwd}}{\ge} |C_1||\xi_1|^{\alpha} - |m_{\mathtt{A},r}(\xi_1)| \nonumber\\
 &\overset{\eqref{skkksdhsd2123},\eqref{csdjsjd1}}\ge \left(c_\alpha|C_1|-C_2\mathtt{M}^{1-\alpha}\right)|\xi_1|^{\alpha} - C_2\mathtt{M}|\xi_1|^{\alpha-1}\nonumber\\
 & \overset{\eqref{cjjjsdiidwda},\eqref{xisd1alsdwd}}{\ge} \frac{c_\alpha|C_1|}{2}|\xi_1|^{\alpha}- \frac{c_\alpha|C_1|}{4}|\xi_1|^{\alpha}\nonumber\\
 &\ge \frac{c_\alpha|C_1|}{4}|\xi_1|^\alpha.\label{largemsdjxck1}
 \end{align}
 Therefore, we have
 \begin{align} \nonumber
 |\mathtt{A}(\xi_1)-\mathtt{A}(\xi_2)|& \ge |\mathtt{A}(\xi_1)|-|\mathtt{A}(\xi_2)| \\& \overset{\eqref{skjsdsdxi2222s},\eqref{largemsdjxck1}}\ge \frac{c_\alpha|C_1|}{4}|\xi_1| - 2\text{max}_{S^+}C_3 \overset{\eqref{cjjjsdiidwda},\eqref{xisd1alsdwd}}\ge \frac{c_\alpha|C_1|}{8}|\xi_1|,\label{tlksdmsgosd}
 \end{align}
 while, it holds that
 \begin{align}
 |\xi_1-\xi_2|(|\xi_1|^{\alpha-1} +|\xi_2|^{\alpha-1})\overset{{|\xi_1|\ge|\xi_2|}}\le 4|\xi_1|^{\alpha}.
 \end{align}
 Combining this with \eqref{tlksdmsgosd}, we obtain
 \begin{align}\label{firstscxcs}
  |\mathtt{A}(\xi_1)-\mathtt{A}(\xi_2)|\ge \frac{c_\alpha|C_1|}{16} |\xi_1-\xi_2|(|\xi_1|^{\alpha-1} +|\xi_2|^{\alpha-1}),\text{ if $|\xi_2|\le 2\text{max}_{S^+}$.}
 \end{align}
 
 In the second case, assuming $|\xi_2|\ge 2\text{max}_{S+}$, we have
 \begin{align}
 &|\mathtt{A}(\xi_1)-\mathtt{A}(\xi_2)| \nonumber \\
 &\overset{\eqref{rjjsdsdsdsd1}}\ge |C_1||\xi_1|\xi_1|^{\alpha-1} - \xi_2|\xi_2|^{\alpha-1}| - |m_{\mathtt{A},r}(\xi_1)-m_{\mathtt{A},r}(\xi_2)|\nonumber\\
 &\overset{\eqref{skkksdhsd2223}}{\ge} |C_1||\xi_1|\xi_1|^{\alpha-1} - \xi_2|\xi_2|^{\alpha-1}| \nonumber \\
 & \qquad \qquad - C_2|\xi_1-\xi_2|(\mathtt{M}^{1-\alpha}(|\xi_1|+|\xi_2|)^{\alpha-1} + \mathtt{M}(|\xi_1|+|\xi_2|)^{\alpha-2})\nonumber\\
 &\overset{\eqref{alpsdoosd1},\eqref{cjjjsdiidwda}}\ge   \frac{|C_1|c_\alpha }{2}|\xi_1-\xi_2|(|\xi_1|^{\alpha-1} + |\xi_2|^{\alpha-1}) -C_2 \mathtt{M}(|\xi_1|+|\xi_2|)^{\alpha-2}\nonumber\\
 &\ge \frac{|C_1|c_\alpha }{2}|\xi_1-\xi_2|(|\xi_1|^{\alpha-1} + |\xi_2|^{\alpha-1}) - \frac{C_2\mathtt{M}}{|\xi_1|}(|\xi_1|^{\alpha-1} + |\xi_2|^{\alpha-1})\nonumber\\
 &\overset{\eqref{cjjjsdiidwda},\eqref{xisd1alsdwd}}\ge \frac{|C_1|c_\alpha}{4}|\xi_1-\xi_2|(|\xi_1|^{\alpha-1} + |\xi_2|^{\alpha-1}).\label{qkfrjsdjsdiq}
  \end{align}
Combining this with \eqref{firstscxcs}, we obtain
\begin{align}\label{sksdsdsdw11koisi}
|\mathtt{A}(\xi_1) - \mathtt{A}(\xi_2)|\ge \frac{|C_1|c_\alpha}{16}|\xi_1-\xi_2|(|\xi_1|^{\alpha-1} + |\xi_2|^{\alpha-1}),
\end{align}
if $\max\left\{|\xi_1|,|\xi_2|\right\} \ge C_{\mathtt{H}2}$ and $\xi_1,\xi_2\in S_\mathtt{M}^\perp\left\{ 0 \right\}$.
This proves \eqref{ajsj1estima2x12}.

\vspace{0.5\baselineskip}\noindent\textit{Proof of \ref{rkjjsodsdwewsd3}.}
From \eqref{jjsususdosod1},and  \eqref{rksjjsjdsdsd1}, we see that for some $C_1(\alpha,\nu,S_0^+),C_2(\alpha,\nu,S_0^+)>0$, 
\begin{align*}
\min_{\vec{U}\in \mathbb{R}^\nu,\ |\vec{U}|=1} \left|\vec{U} \cdot \mathbb{C}_{\xi_1,\xi_2} \vec{U}\right| \ge C_{1}\mathtt{M}^{\alpha+3} - |\vec{U}\cdot W_{r,*}\vec{U}|\\
\ge  C_{1}\mathtt{M}^{\alpha+3} - C_2\mathtt{M}^4.
\end{align*}
Since $\alpha>1$, we can choose $\mathtt{M}$ large enough, depending on $\alpha,\nu,S_0^+$ so that 
\[
\min_{\vec{U}\in \mathbb{R}^\nu,\ |\vec{U}|=1} \left|\vec{U} \cdot \mathbb{C}_{\xi_1,\xi_2} \vec{U}\right| \ge \frac{C_1}{2}{\mathtt{M}}^{\alpha+3}.
\]
Therefore, $\text{det}(\mathbb{C}_{\xi_1,\xi_2} )\ne 0$. 
\end{proof}

\newpage
\thispagestyle{empty}

\backmatter

\addtocontents{toc}{{\string\vskip1\baselineskip}}%
\addtocontents{toc}{\string\enlargethispage{-24pt}}%

%
%

\bibliographystyle{abbrv}
\bibliography{references}

\begin{tabular}{l}
\textbf{Javier G\'omez-Serrano}\\
{Department of Mathematics} \\
{Brown University} \\
{314 Kassar House, 151 Thayer St.} \\
{Providence, RI 02912, USA} \\
{Email: javier\_gomez\_serrano@brown.edu} \\ \\
\textbf{Alexandru D. Ionescu}\\
{Department of Mathematics}\\
{Princeton University} \\
{705 Fine Hall, Washington Rd,}\\
{Princeton, NJ 08544, USA}\\
{e-mail: aionescu@math.princeton.edu}\\ \\
\textbf{Jaemin Park} \\
{Department of Mathematics} \\
{Yonsei University} \\
{50 Yonsei-Ro, Seodaemun-Gu} \\
{Seoul, 03722, South Korea} \\
{Email: jpark776@yonsei.ac.kr} \\ \\
\end{tabular}
\end{document}